\documentclass[12pt,onecolumn]{book}	
\usepackage[top=1in, bottom=1in, left=1in, right=1in]{geometry}
\usepackage{graphicx}
\usepackage{amssymb}
\usepackage{amsmath}  				
\usepackage{epstopdf}  				
\usepackage{graphicx,float,wrapfig} 	
\usepackage[font=small,labelfont=bf,textfont=sf,
labelsep=quad]{caption}				

\makeatletter
\newcommand{\getcitenumber}[1]{
  \@ifundefined{b@#1}
    {\hbox{\reset@font\bfseries ?}}
    {\csname b@#1\endcsname}}
\makeatother

\usepackage[superscript,space]{cite}	

\usepackage{setspace}   

\usepackage [english]{babel}
\usepackage [english = american]{csquotes}
\MakeOuterQuote{"}

 \usepackage[protrusion=false,expansion=false]{microtype}

\newcommand{\mb}[0]{\mbox}
\newcommand{\ov}[1]{\overline{#1}}
\newcommand{\un}[1]{\underline{#1}}
\newcommand{\notequiv}[0]{\equiv\hti{-2}/\hti{1}}

\newcommand{\varss}[0]{\begin{picture}(0,0)(0,0)\put(0,8){\vector(1,0){20}}
\end{picture}}
\newcommand{\ssec}[1]{\subsection{#1}}
\newcommand{\sssec}[1]{\subsubsection{#1}}
\newcommand{\hti}[1]{\hspace*{#1ex}}
\newcommand{\htn}{\hspace{8ex}}
\newcommand{\hth}{\hspace*{8ex}}
\newcommand{\htp}[1]{\hspace*{#1pt}}
\newcommand{\enumb}{\begin{enumerate}\setcounter{enumi}{-1}}
\newcommand{\enume}{\end{enumerate}}
\newcommand{\enumdb}{\begin{enumerate}\setcounter{enumii}{-1}}
\newcommand{\enumde}{\end{enumerate}}
\newcommand{\itemb}{\begin{itemize}}
\newcommand{\iteme}{\end{itemize}}
\newcommand{\ve}[2]{\left(\begin{array}{c}#1\\#2\end{array}\right)}
\newcommand{\vet}[3]{\left(\begin{array}{c}#1\\#2\\#3\end{array}\right)}

\newcommand{\ma}[4]{\left(\begin{array}{cc}#1&#2\\#3&#4\end{array}\right)}

\newcommand{\mat}[9]{\left(\begin{array}{ccc}#1&#2&#3\\#4&#5&#6\\
		#7&#8&#9\end{array}\right)}
\newcommand{\matb}[2]{\left#1\begin{array}{#2}}
\newcommand{\mate}[1]{\end{array}\right#1}
\newcommand{\umod}[1]{\: \underline{mod} \:#1}
\newcommand{\js}[2]{\left(\begin{array}{c}#1\\ \hline #2\end{array}\right)}

\newcommand{\bigx}[0]{$ {\large$ \times\hspace{-10pt}\times$ }$}
\newcommand{\para}[0]{\:/\!/\:}
\newcommand{\paras}[1]{\:/\!/_{#1}\:}
\newcommand{\incid}[0]{\:\iota\:}
\newcommand{\notinc}[0]{\;\iota\hspace{-8pt}-}

\newcommand{\smv}[0]{\hspace{-2pt}\vline\hspace{-2pt}}

\begin{document}

\title{Finite Euclidean and Non-Euclidean Geometries}
\author
{R. De Vogelaere$^{1}$\\
\\
\normalsize{$^{1}$Department of Mathematics}\\
\normalsize{University of California, Berkeley, CA}\\
}
\date{}  

\maketitle

\noindent {\Huge \bf Foreword}
\vspace{10mm}

The author of this monograph was my father, Professor Ren\'e De Vogelaere. He received his PhD in Mathematics in 1948 from the University Louvain, Belgium. Shortly after graduation, he immigrated to Canada and taught at l'Universit\'e Laval in Quebec, followed by Notre Dame in South Bend, Indiana and then the University of California, Berkeley, where he spent most of his career. He studied and taught a wide range of subjects, including differential equations, numerical analysis, number theory, group theory, and Euclidean geometry, to mention a few.  

Georges Lema\^itre, the founder of the "Big Bang" theory, was my father's thesis advisor and lifelong mentor. He was often a guest in our home, and at these meetings he encouraged my father to study astronomy and planetary motion.  After earning his doctorate degree, Professor Lema\^itre spent a year working with Arthur Eddington. Professor Eddington postulated that there were a finite number of protons in the universe. This is known as Eddington's number.  

Ren\'e spent much of his career modeling the continuous world with discrete, finite numbers.  In the late 70's he asked himself: what if the world was discrete rather than continuous? Would the proofs found in different mathematical branches still work?  That is when his research in finite geometry began, culminating in this monograph, to which he dedicated the last 10 years of his life.  In my family, I was the only one who had studied math at the graduate level, and so I was uniquely qualified to share in the excitement of his discoveries and the number of theorems he was able to prove. He was like an archeologist having found a new field of dinosaur bones---discovering something new, then examining and documenting it.  He taught classes on his findings, and wrote many papers (see the bibliography).  He did not publish his book; there were too many exciting theorems to prove, which were much more interesting to him than working with a publisher.

Upon his passing in 1991, I inherited his unfinished book.  I worked with a good friend and past classmate, Michael Thwaites, to try to compile the book written in LaTeX.  But life was busy with family and work. It wasn't easy stepping into my father's shoes to complete this very involved task.  Throughout the following 25 years, I looked for a way to preserve the book and disseminate its knowledge.  Eventually technology and the right person came together.  One late evening I was discussing my father's Finite Geometry book with William Gilpin.  He has just completed his PhD in Physics from Stanford University.  He knew LaTeX very well, and was able to assemble all the files.  He also knew of the Cornell's arXiv and recommended posting it there.  It is the perfect place to store Professor De Vogelaere's magnum opus.

\vspace{8mm}

I would like to thank:
\begin{itemize}
\item My mother, Elisabeth De Vogelaere, who made it possible for my father to dedicate his career to mathematics, which he loved
\item Arthur Eddington, for inspiring my father
\item Georges Lema\^itre, for inspiring and teaching my father
\item The University of California, Berkeley, for providing the facilities and allowing him time to do his research, as well as for archiving the work he did over the 43 years of his career.
\item Michael Thwaites, for helping me get started on the book, and for encouraging me to continue the work
\item My wife Cynthia Haines, for carefully keeping and storing the computer disks, files and papers all these years
\item My daughter Beth, for finding William Gilpin
\item William Gilpin, for his extreme generosity of time to rapidly assemble the book and for facilitating having it stored at Cornell's arXiv.
\item My siblings, Helene, Andrew, and Gabrielle for their patience and faith that this would happen!
\item And Cornell's arXiv, for being there to disseminate knowledge.
\end{itemize}

\vspace{8mm}
\null\hfill Charles De Vogelaere \\
\null\hfill Mountain View, CA \\
\null\hfill August 2019 
\clearpage

\setcounter{chapter}{-1}
\setcounter{tocdepth}{2} 
\tableofcontents

\chapter{Preface}



\def\Bbb#1{#1}
\def\ov#1{\overline{#1}}



\subsection*{Purpose}
The purpose of this book and of others that are in progress
is to give an exposition of Geometry from a point of view which in
some sense complements Klein's Erlangen program. The emphasis is on
extending the classical Euclidean geometry to the finite case, but it goes
way beyond that.

\subsection*{Plan}
In this preface, after a brief introduction, which gives the main theme, and
was presented in some details at the first Berkeley Logic Colloquium of Fall 1989, I present the
main results, according to a synthetic view of the subject, rather that
chronologically. First, some variation on the axiomatic treatment of
projective geometry, then new results on quaternionian geometry, then
results in geometry over the reals which are
generalized over arbitrary fields, then those which depend on properties
of finite fields, then results in finite mechanics. The role of the
computer, which was essential for these inquires is briefly surveyed.
The methodology to obtain illustrations by drawings is described.
The interaction between Teaching and Research is then given. I end with
a table which enumerates enclosed additional material which constitutes
a small but representative part of what I have written.

\subsection*{Introduction}
My inquiry started with rethinking Geometry, by
examining first, what could be preserved among the properties of Euclidean 
geometry when the field of reals is replaced by a finite field. This led
me to a separation of the notions concerned with the distance between 2 points 
and the angle between an ordered pair of lines, into two sets,
those concerned with equality and those concerned with measure.
Properties relating to equality are valid for a Pappian geometry, whatever the
underlying field,  those pertaining to measure require specifying the field.\\
I have also come to the conclusion that the more fruitful approach to the
axiomatic of Euclidean geometry is to reduce it to that of Projective geometry
followed by a preference of certain elements, namely the isotropic points
on the ideal line. This preference can be presented alternately by choosing 2 
points relatively
to a triangle of coordinates, namely the barycenter and the orthocenter.
The barycenter is used to define the ideal line, the orthocenter is then
used to define the fundamental involution of this line, for which the
isotropic points are the (imaginary) fixed points. This program extends to
all non-Euclidean geometries.\\
The preference method, which I call the "Berkeley Program", can be considered 
as the synthetic equivalent of the group theoretical relations between
geometries, as advocated in Felix Klein's Erlangen program.\\
When I refer to Euclidean geometry, I always mean that the set of points
and lines of the geometry of Euclid have been completed by the ideal line
and the ideal points on that line.\\[20pt]

\subsection*{Axiomatic of projective geometry} 
{\bf Projective Geometry}\\
{\em Axiomatic.}\\
The approach, used by Artzy, has the advantage of giving the equivalence between
the synthetic axioms and the algebraic axioms, at each stage
of the axiomatic development: for perspective planes, Veblen-Wedderburn planes,
Moufang planes, Desarguesian planes, Pappian planes, ordered planes, and
finally, projective planes. I have revised it,
to give a uniform treatment (particularly lacking at the intermediate step of
the Veblen-Wedderburn plane, in which, for instance, vectors are introduced by
Artzy and others, to prove commutativity of addition) and by giving, for all 
proofs, explicit, rather than implicit constructions, together with drawings.\\
{\em Notation.}\\
The Theorems of Desargues, Pappus and Pascal play an important role
in synthetic proofs in Projective geometry. A notation has been introduced
for the repeated use of these theorems and their converse, in an efficient
and unambiguous way. A notation for configurations has been introduced, which 
further helps in distinguishing non isomorphic configurations.\\[20pt]

\subsection*{Desarguesian geometry}
{\bf Quaternionian Geometry.}\\
{\em With Relative Preference of 2 Points.}\\
A quaternionian plane is a well known, particularly important, example of a 
Desarguesian plane. I have introduced in it,
the relative preference of 2 points, the barycenter and the
cobarycenter and have obtained several Theorems, which in the
sub-projective planes of the geometry correspond to Theorems in involutive
geometry which are associated with the circumcircle and with the point of 
Lemoine. But these Theorems cannot be considered as simple generalizations.
For instance, in the involution on the ideal line, defined by the circumcircular
polarity, which corresponds to a circumcircle, the direction of a side
and that of the comedian, which generalizes an altitude, are not corresponding 
elements, although these
correspond to each other, in the sub-projective planes. Moreover, what I call
the Lemoine polarity degenerates in the sub-projective planes into all the
lines through the point of Lemoine.  The proofs given are all algebraic.
These investigations are just the beginning of what should become a very
rich field of inquiries.

{\bf Finite Quaternionian Geometry.}\\
The Theorems in quaternionian geometry were conjectured using a geometry
whose points and lines are represented by 3 homogeneous coordinates in the ring
of finite quaternions over ${\Bbb Z}_p$. In the corresponding plane, the 
axioms of allignment are not allways satisfied. If they are, Theorems and
proofs for the quaternionian plane extend to the finite case.\\[20pt]

\subsection*{Pappian geometry over arbitrary fields}
{\bf Pappian Geometry.}\\
This can be considered as a projective geometry over an arbitrary field.\\
{\em On Steiner's Theorem.}\\
Pappus' Theorem is one of the fundamental axioms
of Projective geometry. If the 3 points on one of the lines are permuted,
we obtain 6 Pappian lines which pass 3 by 3 through 2 points, this is the
Theorem of Jakob Steiner. By duality, we can obtain from these, 6 points on 2
lines. That these 2 lines are the same as the original ones is a new Theorem.
Detailed computer analysis of the mapping in special cases leads to conjectures
in which twin primes appear to play a role.\\
{\em Generalization of Wu's Theorem.}\\
I obtained some 80 new Theorems in Pappian geometry, generalizing a Theorem,
in projective geometry, of Wen-Tsen Wu, related to conics through 6 Pascal
points of 6 points on a conic, I have obtained a computer proof for all of 
these Theorems by means of a single program, which includes convincing checks,
and then succeeded in obtaining a synthetic proof for each of these
Theorems, using several different patterns and approaches including duality
and symmetry.  These proofs have benefited from the projective geometry
notation. Drawings have been made for a large number of these Theorems
which have suggested 2 new Theorems and a (solid) Conjecture.  Many of the
Theorems can be considered as Theorems in Euclidean geometry, (only one of which
was known, the Theorem of Brianchon-Poncelet), others can be considered as
Theorems in Affine or in Galilean geometry.\\
{\em Generalization of Euclidean Theorems.}\\
The Theorems, given for involutive geometry, can be considered, alternately,
as Theorems in Pappian geometry, because they involve only the
preference of 2 elements of the projective plane and not additional axioms.

{\bf Involutive Geometry.}\\
I call involutive plane, a Pappian plane in which I prefer 2 points
relative to a triangle, $M$, the barycenter and ${\ov M}$, the orthocenter.
$M$ allows for the definition of
the ideal line, ${\ov M}$ allows, subsequently, for the definition of the
fundamental involution on that line.\\
{\em Generalization of Theorems in Euclidean and Minkowskian Geometry over
Arbitrary Fields.}\\
In this, which constitues the more extensive part of my research, I have
generalized, when the involution is elliptic, a very large number of Theorems
in Euclidean geometry, namely those which are characterized by not using the
measure of distance and of angles and not involving elements whose construction
leads to more than one solution. When the fundamental involution is hyperbolic,
each of the Theorems gives a corresponding Theorem in the geometry
of Hermann Minkowski.\\
{\em Symmetry and Duality.}\\
The barycenter and orthocenter have a symmetric role for many Theorems of 
Euclidean geometry, the
line of Euler and the circle of Brianchon-Poncelet being the simpler examples.
This has been systematically exploited, to almost double the number of
Theorems known in that part of Euclidean geometry which involves congruence
and not measure. Duality can also be extended to Euclidean geometry by
associating to $M$ and $\ov{M}$, the ideal line and the orthic line and 
vice-versa.
This also has been systematically exploited to help me, in obtaining
constructions of new elements, and should be helpful in future constructions.\\
{\em Notation.}\\
A set of notations was introduced, to allow for a compact description of
some 1006 definitions, 1073 conclusions and for the corresponding proofs. The
counts correspond to one form of counting, other forms give higher numbers.
All these Theorems are valid for any Pappian plane
and give directly both statement and new proofs in both Euclidean and 
Minkowskian geometry.\\
{\em The Geometry of the Triangle.}\\
During the period 1870 to 1900, there was an explosion of results in what has
been called the geometry of the triangle, prepared by Theorems due to
Leonhard Euler, Jean Poncelet, Charles Brianchon, Emile Lemoine and others.
 The synthesis of the subject
was never successfully accomplished, not only because of the wealth of
Theorems, but because of the difficulty of insuring that elements defined
differently were in fact, in general, distinct. The proofs, used in involutive 
geometry, not only throw a new light on the reason for the
explosive number of results for the geometry of the triangle but also
gives a exhaustive synthetic view of the subject.\\
{\em Diophantine Equations.}\\
Because an algebraic expression of the homogeneous coordinates of points and 
lines and the coefficients for conics is given in terms of polynomials in 3 
variables, a large number of particular results on diophantine equations in 3
variables are  implicitly obtained in these investigations.\\
{\em Construction with the Ruler only.}\\
In all of the classical investigations, the most extensive one being that of 
Henri Lebesgue, the impression is given that the compass is indispensable
for most constructions in geometry.  More than half of the Theorems for
which a count is given above, can be characterized as using the ruler only.
Implicit, in this part of my Research, is, that many constructions, which
usually or by necessity were assumed to require the compass, in fact
need the ruler only, the simplest one is that for constructing the perpendicular
to a line. The more remarkable one is that the circles of Apollonius
can be constructed with the ruler only. These are defined as the circles which 
have as diameter the intersections of the bisectrices of an angle of a triangle 
with the opposite sides. It is this reduction to construction with the
ruler alone, which allows for the straigthforward proofs which constitutes
a major success of these investigations.\\
{\em Construction with the Ruler and Compass.}\\
The construction with compass can be envisioned as follows. Given $M$ and 
$\ov{M}$,  by finding the intersection of 2 circles centered at 2 of
the vertices of a triangle with the adjacent sides and by constructions
with the ruler, we can construct the bissectrices of these angles, the
incenter (center of the inscribed circle) and the point of Joseph Gergonne
(the common intersection of the lines through a vertex of the triangle and the
point of tangency of the inscribed circle with the opposite side). From these,
a very large number of other points, lines and circles can be constructed with
the ruler only, for instance, the point of Karl Feuerbach, the excribed circles
and the circles of Spieker. One can therefore, in the framework of involutive 
geomatry, prefer instead of $M$ and $\ov{M}$, the incenter $I$ and the
point of Gergonne $J$. Starting from $I$ and $J$, we can construct $M$ and 
$\ov{M}$, using the ruler alone. This allows to extend the proof methodology
considerably, allowing the generalization to arbitrary fields of Theorems
involving elements whose construction, in the
classical case, would requires the compass.\\
{\em Cubics.}\\
Very little has been written on the construction of cubics by the ruler.
Starting with the work of Herman Grassmann of R. Tucker and of Ian Barbilian,
 I have obtaining a few results in this direction, one of which, incidentally,
gives a illustration of the procedure of construction with the compass as I am 
envisioning it, which is much simpler than those involving bissectrices.

{\bf Galilean geometry.}\\
When the fundamental involution is parabolic and when the field is the field of
reals, the geometry is called Galilean, because its group is the group of 
Galilean transformations of classical mechanics. Extending to the Pappian case 
and starting from the definitions and conclusions of involutive geometry, I have
made  appropriate modifications to obtain
Theorems which are valid in Galilean geometry, but I have not yet completed the
careful check that is required to insure the essential accuracy.
Again a very large number of Theorems have been obtained, which are new, even 
in the case of the field of reals.

{\bf Polar Geometry.}\\
The extension, to $n$ dimensions, can be obtained using an appropriate adaptation
of the algebra of Herman Grassmann.  A first set of Theorems has been 
obtained in the
case of 3 dimensions, again for a Pappian space over arbitrary fields, in which
preference is given to one plane, the ideal plane and one quadric. These
Theorems generalize Theorems on the tetrahedron due to E. Prouhet, Carmelo
Intrigila and Joseph Neuberg.  The special case of the orthogonal tetrahedron
has also been studied in a way which puts in evidence the reasons behind
many of the Theorems obtained in this case.

{\bf Non-Euclidean Geometry.}\\
The beginning of the preference approach to obtain new results in non-Euclidean
geometry was started in January 1982. The confluence, in the case of a finite 
field, of the geometries of Janos Bolyai and of Nikolai Lobachevsky was then 
explored. A new point, called the center of a triangle was discovered and
its properties were proven.\\[20pt]

\subsection*{Pappian geometry over finite fields}
{\bf The Case of Finite Fields.}\\
All the results given for involutive geometry and in the following sections
are true, irrespective of fields. In what follows, we describe results
for finite fields.

{\bf Projective Geometry.}\\
{\em Representation on Pythagorean and Archimedean solids.}\\
Fernand Lemay has shown how to represent the projective planes corresponding to
the Galois fields, 2, 3 and 5 respectively on the
tetrahedron, the cube (or octahedron) and the dodecahedron (or icosahedron).
I have shown, that if we choose instead of the Pythagorean solids, the
Archimedean ones, the results extend to $2^2$ and the 5-gonal antiprism
and to $3^2$ and the truncated dodecahedron. I have studied also the
corresponding representations of the conics on the dodecahedron. This is
useful for the representation on it of the finite non-Euclidean Geometry
associated with $GF(5)$.

{\bf Involutive Geometry.}\\
{\em Partial Ordering.}\\
In the case of finite fields, ordering and therefore the notions of limits
and continuity are not present.  
By using Farey sets or, alternately, by using a symmetry property of the 
continued fraction algorithm, I have introduced partial ordering in 
${\Bbb Z}_p$. If only, the properties of order have to be preserved which are 
related to the additive inverse and multiplicative inverse, then
a Theorem of Mertens allows me to estimate the cardinality of the
ordered subset of ${\Bbb Z}_p$ by .61 $p$, when $p$ is large. The
cardinality is decreased logarithmicaly, by a factor 2, for each additional
operation of addition and multiplication, for which order needs to be
preserved.\\
{\em Orthogonal polynomials.}\\
Orthogonal polynomials can be defined in a straightforward way in ${\Bbb Z}_p$.
For those I have studied, it turns out, that the classical scaling used
in defining the classical orthogonal polynomials, there is a symmetry which
is exibited in each case, with the exception of those of Charles Hermite. 
In this case, by using an alternate scaling, with different expressions
for the polynomials of even and odd degree, symmetry can also be obtained.\\
{\em Finite Trigonometry.}\\
Ones the measure of angles between an ordered pair of non ideal lines
and the measure of the square of the distance between two ordinary points
has been defined, it is straightforward to obtain the trigonometric functions
in ${\Bbb Z}_p$. There are in fact, for each prime $p$, two sets of 
trigonometric 
functions, one corresponding to the circular ones, one to the hyperbolic ones.
The proofs required, depend on the existence of primitive roots, in the case 
corresponding to Minkowskian geometry, and on a generalization
to the Galois field $GF(p^2)$ in the case corresponding to Euclidean geometry.\\
{\em Finite Riccati Functions.}\\
The functions of Vincenzo Riccati, which are generalization of the
trigonometric functions have been defined and studied in the finite case.
They enable the definition of a Riccati geometry. An invariant defines
distances, the addition formulas, which correspond to multiplication of
associated Toeplitz matrices, define addition of angles.
 This again should be a fruitful field of inquiry.\\
{\em Finite Elliptic Functions.}\\
After I conjectured that the Theorem of Poncelet on polygons
inscribed to a conic and circumscribed to an other conic extended to the
finite case, I knew that Finite Elliptic functions could be defined in the
finite case, because I had learned from Georges Lema\^{\i}tre the
relation between Theorems on elliptic functions and the Theorem of Poncelet.
The functions I defined, correspond to the functions $sn$, $cn$ and $dn$
of Karl Jacobi. After I found that John Tate had defined the Weierstrass type
of finite elliptic functions I established the relation between the 2.\\
{\em Construction with the compass.}\\
In the case of finite fields, the points $I$ and $J$ will only exist if 2
and therefore all angles of the triangle are even.  Prefering $I$ and $J$
instead of $M$ and $\ov{M}$, insures that the triangle is even.

{\bf Isotropic Geometry.}\\
Many of the Theorems in involutive and polar geometry do not apply to the case
of fields of characteristic 2, because the diagonal points of a complete 
quadrilateral are collinear, because every conics has all its tangents incident 
to a single point and because in the algebraic formulations, 2, which occurs in 
many of 
the algebraic expressions involved in corresponding proofs of involutive
geometry is to be replaced by 0. I call isotropic plane, a Pappian plane,
with field of characteristic 2 and with the relative preference of 2 points,
$M$, the barycenter and, $O$, the center. The orthocenter does not exist when
the characteristic is 2 because each line can be considered as perpendicular to 
itself. The difference sets of J. Singer, called selectors by Fernand Lemay,
were an essential tool in these investigations. In an honor Thesis, Mark 
Spector, now a Graduate Student in Physics at 
M.I.T. wrote a program to check the consistency of the notation in the 
statements of the Theorems and the accuracy of the proofs. He obtained
new results. My results on cubics are not retained in his honors Thesis.\\
Some of the results in isotropic geometry were anticipated by the work of
J. W. Archbold, Lawrence Graves, T. G. Ostrom and D. W. Crowe.\\[20pt]

\subsection*{Finite mechanics and simplectic integration}
I was asked to participate in a discussion, Spring 1988, at Los Alamos,
on the field of simplectic integration which I originated in 1955.
Simplectic integration methods are methods of numerical integration
which preserve the properties of canonical or simplectic transformations.
It then occured to me, that these methods were precisely what was needed
to extend to the finite case the solution of problems in Mechanics.
I had searched for a solution to this problem since I obtained, as first
example, the solution, using finite elliptic functions, for the motion in
${\Bbb Z}_p$ of the pendulum with large amplitude, as well as the polygonal 
harmonic motion, whose study was suggested by a Theorem of John Casey, and
led to an equation similar to Kepler's equation.\\
More specifically, whenever the classical Hamiltonian describing a motion
has no singularities, a set of difference equations can be produced
whose solutions at successive steps have the properties associated with 
simplectic transformations. To confirm the solidity of this approach,
I studied, in detail, the bifurcation properties for one particular 
Hamiltonian.
The study can be made in a more complete fashion than in the classical
case and requires a much simpler analysis using the $p$-adic analysis of
Kurt Hensel.\\[20pt]

\subsection*{The role of the computer for conjectures and verification}
The computer was an essential tool in the conjecture part of the Research 
described above, in the verification of the order of the statements and to
insure the consistency of the notation used in the statements of the Theorems
as well as in the verification of the proofs. 
In particular, the Theorem refered to in the Steiner section was conjectured
from examples from finite geometry. All of the Theorems generalizing Wu's
Theorem were conjectured by examining, in detail, one appropriately
chosen example, for a single finite field. Many Theorems in involutive geometry
and all the Theorems in quaternionian geometry were so conjectured and the
methodology used was such that almost all conjectures could be proven.
The remaining ones could easily be disposed of, by a counterexample or
algebraically. The only exception are the conjectures, indicated in the
section on Steiner's Theorem, which refer to twin primes.\\[20pt]

\subsection*{Illustrations by drawings}
Responding to natural requests for figures which illustrate the many
Theorems obtained, I have also prepared a large number of drawings.
These have been done for the case of the field of reals and therefore in the
framework of classical Euclidean geometry. These are created by means of a 
VMS-BASIC program, which constructs a POSTSCRIPT file, for any set a data, including
points, lines, conics and cubics. The position of the labels of points and
lines can be adjusted by adding the appropriate information to the data file
in order to position the labels properly. One such illustration was chosen by George
Bergman, for this years poster on "Graduate opportunities in Mathematics for 
minority and women students".\\[20pt]

\subsection*{Interaction between research and teaching}
These 2 obligations are for me very closely intertwined, my specific
contributions to teaching are given in a separate document. The
conjecture aspect of my research was exclusively dependent on VMS-BASIC
programs which were a natural extension of programs which I wrote for my
classes. Many of the proofs are dependent on material contained in
notes I prepared for students while teaching courses not related to my
original specialty of Numerical Analysis and of Ordinary Differential
Equations.\\
Many results have been presented in courses, a few, in Computation
Mathematics, (Math. 100), Abstract Algebra (Math. 113) and Number Theory
(Math. 115), a large number,  in a seminar on Geometry, 2 years
ago, and in Foundations of Geometry (Math. 255), Fall 1989.\\[20pt]

\subsection*{Notes and publication}
The scope of the results and their constant interaction during the years
made it impractical to publish incrementally without slowing down considerably
the pace of the inquiry. I have only given a brief overview in 1983 and in
 1986.
\begin{quotation} {\em Finite Euclidean and non-Euclidean Geometry with
	application to the finite Pendulum and the polygonal harmonic
	Motion.  A first step to finite Cosmology.}  The Big Bang
	and Georges Lema\^{\i}tre, Proc. Symp. in honor of 50 years after his
	initiation of Big-Bang Cosmology, Louvain-la-Neuve, Belgium,
	October 1983., D. Reidel Publ. Co, Leyden, the Netherlands. 341-355.

	{\em G\'{e}om\'{e}trie Euclidienne finie.  Le cas p premier
	impair.} La Gazette des Sciences Math\'{e}matiques du Qu\'{e}bec, Vol.
	10, Mai 1986.

	{\em Basic Discoveries in Mathematics using a Computer.}
	Symposium on Mathematics and Computers, Stanford, August 1986.
\end{quotation}

\subsection*{A short guide to the reader.}

The reader may want to start directly with Chapter II and to read sections
of the introductory Chapter as needed. He may perhaps wish to read the
section on a model of finite Euclidean geometry with the framework of
classical geometry, if he wishes to be more confortable about the generalization
of the Euclidean notions to the finite case. If at some stage the readers
wants a more tourough axiomatic treatment it will want to read the section on
axiomatic of the first Chapter.

Chapter II is written in terms of finite projective geometry associated to the
prime $p$, but, except in obvious places, all definitions and Theorem
apply to Pappian planes over arbitrary fields. Among the new results,
included in this Chapter, are, a Theorem related to the Steiner-Pappus Theorem,
considerations on a "general conic", a description of the Theorems of
Steiner, Kirkmanm Cayley and Salmon in terms of permutation maps.
After describing the representation of the finite projective planes for
p = 2, 3 and 5 on Pythagorian solids, the generalization to the
projective plane of order $p^2$ on the truncated dodecahedron is given
as well as that of the plane of order on the antiprism.  Difference sets
involving non primitive polynomials are studied which allow a definition of the
notion of distance for affine as well as other planes.

Attention is also drawn to B\'{e}zier curves, which have not yet entered
the classical repertoire of Projective Geometry. These are used extensively
in the computer drawing of curves and surfaces.

	One of the reason for the historical delay of extended the Euclidean 
notions associated with distance between points and angle between lines is
the lack of early distinction between equality and measure.  Equality is
a simpler notion which can be dealt with over arbitrary fields, while
measure requires greater care.  This is examplified by the comment on finite
projective geometries by O'Hara and Ward, p. 289.
\begin{quote}
Their analytic treatment involves the theory of numbers, and, in particular the 
theory of numerical congruences; it may be assumed that the synthetic treatment
of them is correspondingly complicated.
\end{quote}

It is my fondest hope that some of the material on
finite geometry will be assimilated to form the basis of renewal of the teaching
of geometry at the high school level, combined with a well-thought related
use of computers at that level.


\chapter{MAIN HISTORICAL DEVELOPMENTS}



\def\Bbb#1{#1}
\def\ov#1{\overline{#1}}



\setcounter{section}{-1}
	\section{Introduction.}
In this chapter, I give the main historical developments in Mathematics which
have a bearing on the generalization of Euclidean Geometry to the finite case
and to non Euclidean Geometries.\\
What could be consider as the first contribution to Mathematics which covers
number theory, geometry and trigonometry is a tablet in the Plimpton
collection, this is briefly described and discussed in a note at the end of the 
Chapter. The key to the treatment of geometry and its use of continuity dates
from the discovery of the irrationals by the school of Pythagoras. This is
commented upon to suggest an alternative which is consistent with finite
Euclidean
geometry. I thought it would be handy for many readers to have at hand the
definitions and postulates of Euclid, as well as a brief description of his
13 books, if only to see how we have travelled in getting a more precise
description of concepts and theorems in geometry. Distances play an essential,
if independent role, in the development of geometry, until recently,
after some comments on the subject, I give some post Euclidean theorems
involving distnaces on the sides of a triangle due to Menelaus and Ceva.
The geometry of the triangle, which has played an important historical role,
is illustrated by theorems due to Euler, Brianchon and Poncelet, Feuerbach,
Lemoine and Schr\"{o}ter.\\
I then review quickly some of the major developments in projective geometry due 
to Menaechmus, Apollonius, Desargues, Pascal, MacLaurin, Carnot, Poncelet,
Gergonne and Chasles. In the next section, I start the process of going back
from projective, to affine, to involutive, to Euclidean geometry.\\
I then review the algebraization of geometry starting with Descartes and 
Poncelet and ending with James Singer, who spured by a paper of Veblen and
MacLagan-Wedderburn, introduced the notion of difference sets which allows
the representation of every point and line in a finite Pappian plane by an
integer, allowing an easy determination of incidence, without
coordinatization.\\
This is followed by a section on trigonometry which gives the Lambert formulas
valid in the case of finite fields.\\
The section on algebra is for the reader which has been away from the subject 
for some time. It includes algorithms to solve linear diophantine equations
and to obtain the representation of numbers as sum of 2 squares, the
definition of primitive roots and the application to the extraction of
square roots in a finite field, contrasting with the solution of the
school of Pythagoras.\\
The section on Farey sets includes original material on
partial ordering of distances, which at least suggest that the essential
notion of ordering in the classical case can be extended to the finite case.\\
Definition of complex and quaternion integers, loops, groups,
Veblen-Wedderburn systems and ternary rings are given as a preparation for
the section on axiomatic. The important relevant contributions of Klein,
Gauss, Weierstrass, Riemann, Hermite and Lindenbaum are then recalled.\\
The subject of elliptic functions and the application of geometry to mechanics
has lost, at the present time,
the great interest it had during last century. Because this too generalizes
to the finite case and because this is not now part of the Mathematics
curriculum, I have a long section introducing one of its components,
the motion of the pendulum to introduce elliptic integrals, the
elliptic functions of Jacobi as well as his theta functions, ending
with the connection given first by Lagrange between spherical trigonometry
and elliptic functions.\\
To add credibility to the existence of non Euclidean geometries, models
were divised to give models within the framework of Euclidean geometry.
The next section gives a model of finite Euclidean geometry also within
this framework. It can be used as an introduction to the subject.\\ 
The axiomatic of geometry in the next section is done using a uniform treatment,
and explicit constructions. It includes a plane which is, like the Moufang
plane, intermediate
between the Veblen-Wedderburn plane and the Desaguesian plane. The
geometry of Lenz-Barlotti of type I.1 discovered by Veblen and
MacLagan-Wedderburn and studied by Hughes is an example of this intermediate
plane.

	\section{Before Euclid.}

	\ssec{The Babylonians and Plimpton 322.}\label{sec-Sbab}

	\setcounter{subsubsection}{-1}
	\sssec{Introduction.}
	Besides estimating areas and volumes, the Babylonians
	had a definite interest in so called Pythagorian triples, integers $a,$
	$b$ and $c$ such that $a^2 = b^2 + c^2.$

	In tablet 322 of the Plimpton library collection from Columbia
	University, dated 1900 to 1600 B.C., a table gives, with 4 errors, and
	in hexadesimal notation, 15 values of\\
	\hth$  a,$ $b,$ and $(\frac{a}{c})^2 = sec^2(B),$\\
	corresponding to angles varying fairly regularly from near $45^{\bullet}$
	to near $32^{\bullet}$. (See Note \ref{sec-plimpton}).

	It is still debated if their interest was purely arithmetical or was
	connected with geometry (See Note \ref{sec-babyl}).

	\ssec{The Pythagorean school.}\label{sec-Spyth}

	That the ratio of the length of the sides of a triangle
	is equal to the ratio of 2 integers was first contradicted by the
	counterexample of an isosceles right triangle $A_0,$ $A_1,$ $A_2,$ with
	right angle at $A_0$ and with sides $a_1$ and hypothenuse $a_0.$
	The theorem of Pythagoras states that\\
	\hth$a_0^2 = a_1^2 + a_1^2 = 2 a_1^2,$ $a_0 > a_1 > 0.$\hfill(1)\\
	If $a_0$ and $a_1$ are positive integers, it follows from the fact
	that the square of an odd integer is odd and that of an even integer
	is even, and from (1), that $a_0^2$ and therefore $a_0$ is even,
	therefore $a_0 = 2 a_2$ and\\
	\hth    $a_1^2 = 2 a_2^2,$ $a_1 > a_2 > 0.$\hfill(2)\\
	The argument can be repeated indefinitely and an infinite sequence
	of decreasing positive integers is obtained,\\
	\hth$     a_0 > a_1 >  \ldots  > a_n >  \ldots  > 0.$\hfill(3)\\
	But this contradicts the fact that only a finite number of positive
	integers exist which are less than $a_0.$\\
	Geometrically, the proof follows from the following figure:

\begin{picture}(200,60)(-160,0)
\put(-40,0){\line(1,1){40}}
\put(-20,0){\line(1,1){20}}
\put(-10,0){\line(1,1){10}}
\put(10,0){\line(-1,-1){10}}
\put(20,0){\line(-1,-1){20}}
\put(40,0){\line(-1,-1){40}}

\put(-40,0){\line(1,-1){40}}
\put(-20,0){\line(1,-1){20}}
\put(-10,0){\line(1,-1){10}}
\put(10,0){\line(-1,1){10}}
\put(20,0){\line(-1,1){20}}
\put(40,0){\line(-1,1){40}}

\put(40,40){\line(0,-1){80}}
\put(20,20){\line(0,-1){40}}
\put(10,10){\line(0,-1){20}}
\put(-40,40){\line(0,-1){80}}
\put(-20,20){\line(0,-1){40}}
\put(-10,10){\line(0,-1){20}}

\put(40,40){\line(-1,0){80}}
\put(20,20){\line(-1,0){40}}
\put(10,10){\line(-1,0){20}}
\put(-40,-40){\line(1,0){80}}
\put(-20,-20){\line(1,0){40}}
\put(-10,-10){\line(1,0){20}}

\put(-6,44){$a_0$}
\put(-6,24){$a_2$}
\put(20,24){$a_1$}
\put(8,14){$a_3$}
\end{picture}\\[40pt]

	This argument has been refined through the ages, by a careful
	construction of the integers, see for instance the Appendix by Professor
	A. Morse in Professor J. Kelley's book on Topology, by an analysis of
	their divisibility properties (see the Theorem of Aryabatha) and by
	their ordering properties (the well ordering axiom of the integers).
	What is implicit in the geometry considered by the Greeks, after
	Pythagoras, is that the circle with center $A_1$ and radius $a_0$
	meets the line through $A_1$ and
	$A_0$ at a point, but this assumption is not made explicitely.  From it
	follows the existence of points on the line corresponding to the
	irrational $\sqrt{2}$ and also the existence of the unrelated
	irrationals, $\sqrt{3},$ \ldots, $\sqrt{17},$  \ldots, more generally,
	$\sqrt{p},$ for $p$ prime, eventually this lead Euclid to consider
	that the set of points on each line forms a continuous set.\\
	Moreover the theorem of Pythagoras assumes the axiom on parallels of
	Euclid.\\
	In finite affine geometry, I will keep the axiom of parallels but
	assume that the number of points on each line is finite.  In finite
	Euclidean Geometry most of the
	notions of ordinary Euclidean geometry are preserved, the measure of
	angles presents no difficulties and the measure of distances requires
	the introduction of one irrational.  On the other hand circles meet
	half of the lines through their center in 2 points and the other half
	in no point and $\sqrt{2}$ need not be irrational.
	See \ref{sec-sqr2rat}.



	\section{Euclidean Geometry.}

	\ssec{Euclid.(3-th Century B.C.)}\label{sec-Seuclid}
	The greek geometer Euclid (300 B.C) constructed a careful theory of
	geometry based on the primary notions of points, lines and planes and on
	a set of axioms, the last one being the axiom on parallels.\\
	His first 3 books are devoted to a study of the triangle, of the circle
	and of similitude.\\
	I will list here the definitions, postulates and common notions as
	translated by Heath, p. 153 to 155:

	\sssec{Definitions.}
	\enumb
	\item  A {\em point} is that which has no parts.
	\item  A {\em line} is breadthless length.
	\item  The extremities of a line are points.
	\item  A {\em straight line} is a line which lies evenly with the
		points on itself.
	\item  A {\em surface} is that which has length and breath only.
	\item  The extremities of a surface are lines.
	\item  A {\em plane surface} is a surface which lies evenly with the
		straight lines on itself.
	\item  A {\em plane angle} is the inclination to one another of two
		lines in a plane which meet one another and do not lie in a
		straight line.
	\item  And when the lines containing the angle are straight, the
		{\em angle is} called {\em rectilinear}.
	\item When a straight line set up on a straight line makes the adjacent
		angles equal to one another, each of the equal {\em angles} is
		{\em right}, and the straight line standing on the other is
		called a {\em perpendicular} to that on which it stands.
	\item An {\em obtuse angle} is greater than the right angle.
	\item An {\em acute angle} is an angle less than a right angle.
	\item A {\em boundary} is that which is an extremity of anything.
	\item A {\em figure} is that which is contained by any boundary or
		boundaries.
	\item A {\em circle} is a plane figure contained by one line such that
		all the straight lines falling upon it from one point among
		those lying within the figure are equal to one another.
	\item And the point is called the {\em center of the circle}.
	\item A {\em diameter} of the circle is any straight line through the
		center 	and terminated in both directions by the circumference
		of the circle, and such a straight line also {\em bisects the
		circle}.
	\item A {\em semicircle} is the figure contained by the diameter and the
		circumference cut off by it.  And the center of the semicircle
		is the same as that of the circle.
	\item {\em Rectilineal figures} are those which are contained by
		straight lines, trilateral figures being those contained by
		three, quadrilateral those contained by four, and multilateral
		those contained by more than four straight lines.
	\item Of trilateral figures, an {\em equilateral triangle} is that
		which has its three sides equal, an {\em isosceles triangle}
		that which has two of its sides alone equal, and a {\em scalene
		triangle} that which has its three sides unequal.
	\item Further, of trilateral figures, a {\em right-angled triangle} is
		that which has a right angle, an {\em obtuse-angled triangle}
		that which has an obtuse angle, and an {\em acute-angled
		triangle} that which has three angles acute.
	\item Of quadrilateral figures, a {\em square} is that which is both
		equilateral and right-angled; an {\em oblong} that which is
		right-angled but not equilateral; a {\em rhombus} that which is
		equilateral but not right-angled; and a {\em rhomboid} that
		which has opposites sides and angles equal to one another but
		is neither equilateral or right-angled. And let quadrilaterals
		other than these be called {\em trapezia}.
	\item {\em Parallel straight lines} are straight lines which, being in
		the same plane and being produced indefinitely in both
		directions, do not meet one another in either direction.
	\enume

	\sssec{Postulates.}
	Let the following be postulated.
	\enumb
	\item  To draw a straight line from one point to any point.
	\item  To produce a finite straight line continuously in a straight
		line.
	\item  To desribe a circle with any center and distance.
	\item  That all right angles are equal to one another.
	\item  That, if a straight line falling on two straight lines make the
		interior angles on the same side less than two right angles,
		the two straight lines, if produced indefinitely, meet on that
		side on which are the angles less than the two right angles.
	\enume

	\sssec{Common notions.}
	\enumb
	\item {\em Things which are equal to the same thing are also equal to
		one another.}
	\item {\em If equals be added to equals, the wholes are equal.}
	\item {\em If equals be subtracted from equals, the remainders are
		equal.}
	\item {\em Things which coincide with one another are equal to one
		another.}
	\item {\em The whole is greater than the part.}
	\enume

	\sssec{Short description of the Books of Euclid.}
	The work of Euclid consists of 13 books which contain
	propositions which are either theorems proving properties of
	geometrical figures or theorems concerned with proving that certain
	figures can be constructed.  It also consists of a study of
	integers, rationals and reals.\\
	\begin{itemize}
	\item Book 1 is devoted mainly to congruent figures, area of triangles
	and culminates with the Theorem of Pythagoras (Proposition 47).
	\item Book 2 is concerned with construction of which the following is
	typical, determine $P$ on $A B$ such that $AP^2 = AB . BP$.
	\item Book 3 studies in detail circles, tangent to circles, tangent
	circles.
	\item Book 4 constructs polygons inscribed and outscribed to circles.
	\item Book 5 gives the theory of proportions.
	\item Book 6 applies the theory of proportions to geometrical figures.
	\item Book 7 studies integers, their greatest common divisor
	(Proposition 2) and their least common multiple (Proposition 34).
	\item Book 8 studies proportional numbers.
	\item Book 9 studies geometrical progression, in Proposition 20,
	the proof that the number of primes is infinite is given.
	\item Book 10 studies the commensurables and incommensurables.
	\item Book 11 is on 3 dimensional or solid geometry.
	\item Book 12 studies similar figures in solid geometry.
	\item Book 13 studies properties of pentagons and decagons as well as
	the regular solids.
	\end{itemize}

	\sssec{Comment.}
	These definitions, postulates and axioms have been discussed
	since the time of Euclid.  The reader is urged to study some of these
	discussion, for instance those in the book of Heath.  Already Proclus
	(see Paul van Eecke) criticizes Postulate 5, and claim that it should be
	proven.  Let me only observe here that except for the notion of being on
	the same side and the notion of continuity, which are absent from finite
	Euclidean geometry, in some sense all of the definitions and postulates
	given above are valid in finite Euclidean geometry.  It should be
	stressed that the expression ``produced indefinitely" (eis apeiron)
	cannot be translated by ``to infinity" (see Heath, p. 190).\\
	Heath observes also (p. 234) that Euclid implies that ``straight lines
	and circles determine by their intersections other points in addition to
	those given" and that ``the existence of such points of intersection
	must be postulated".  He concludes that ``the deficiency can only be
	made good by the Principle of Continuity" and proceed by giving the
	axioms of Killing.\\
	We will see that the alternate route of finite Euclidean geometry
	disposes of the problem quite differently and that some figures cannot
	always be constructed.\\
	It will also be seen that the great emphasis given to distance between
	points and angles of two straight lines and their equality are
	notions which we will derive from more basic notions and that following
	the point of view adopted since the 19-th century no attempt will be
	made to define points and lines, as in Euclid, but we will give instead
	properties that they possess.  In this connection the critique of
	Laurent, H., 1906, p.69 is of interest:

	\begin{quotation}Euclid and Legendre have imagined that the word
	`distance' has a meaning
	and they believed that the proofs using superposition have a `logical'
	value.  Moreover few of the present day geometers have observed that
	Legendre and Euclid have erred.  And that is, I believe, one
	of the more curious psychological phenomenons that for more than two
	thousand years one does geometry without realizing that its fundamental
	propositions have no sense from a `logical' point of view
	\footnote
	{Euclide et Legendre se sont figur\'{e} que le mot `distance' avait un
	sens et ils ont cru que les d\'{e}monstrations par superposition avaient
	une valeur `logique'.  D'ailleurs peu de g\'{e}om\`{e}tres aujourd'hui
	s'apercoivent que Legendre et Euclide ont divagu\'{e}.  Et c'est l\`{a},
	\`{a} mon avis, un des ph\'{e}nom\`{e}nes psychologiques les plus
	curieux que, depuis plus de deux milles ans, on fait de le g\'{e}ometrie
	sans s'apercevoir que ses propositions fondamentales n'ont aucun sens
	au point de vue `logique'.}.
	\end{quotation}

	I will now state a few theorems which play an important role in Part II.
	a few of which are not in Euclid or Legendre.


	\sssec{Definition.}
	The {\em altitude} through $A_0$ is the line through $A_0$ which is
	perpendicular to $A_1 A_2.$  The {\em foot of the altitude} through
	$A_0$ is the point $H_0$ on the altitude and on $A_1 A_2.$

	\sssec{Theorem.}
	{\em The altitudes through $A_0,$ $A_1$ and $A_2$ are concurrent in
	$H.$}

	\sssec{Definition.}
	The point H is called the {\em orthocenter.}

	\sssec{Theorem.}
	{\em Let $M_0$ be the mid-point of $A_1 A_2,$ let $M_1$ be the
	mid-point of
	$A_2 A_0$ and let $M_2$ be the mid-point of $A_0 A_1,$ then $A_0 M_0,$
	$A_1 M_1$ and $A_2 M_2$ are concurrent in $M.$}

	\sssec{Definition.}
	The point $M$ is called the {\em barycenter} or {\em center of mass}.

	\sssec{Definition.}
	$m_0$ is the {\em mediatrix} of $A_1 A_2$ if $m_0$ passes through
	$M_0$ and is perpendicular to $A_1 A_2.$
	
	\sssec{Theorem. [Euclid, Book 4, Proposition 5.]}
	{\em The mediatrices $m_0,$ $m_1$ and $m_2$ are concurrent in $O.$}

	\sssec{Definition.}
	The point $O$ is called the {\em center of the circumcircle of the
	triangle} $A_0 A_1 A_2.$

	\sssec{Theorem. [Euler]}\label{sec-teul}
	{\em The points $H,$ $M$ and $O$ are on the same line $e.$}

	\sssec{Definition.}
	The line $e$ is called the {\em line of Euler}.

	\sssec{Comment.}
	The usual proof of \ref{sec-teul} is geometric.  The proof given by
	Euler is entirely algebraic.  It is based on an expression for the
	distances of $HG,$ $HO$ and $OG$ in terms of the sides of the triangle.
	Let $a,$ $b$ and $c$ be the sides of the triangle.  Let\\
	\hth$p = a + b + c, q = bc + ca + ab, r = a b c,$\\
	(the symmetric functions of $a,$ $b$ and $c$).\\
	The area $A$ is given by  $AA = \frac{1}{16}(-p^4 + 4ppq - 8pr)$
	\footnote{I use here the notation of Euler and of mathematicians before
	the middle of the 19th century, namely $AA$ for $A.A.$}.\\
	Euler obtains\\
	\hth$  HM\:HM = \frac{1}{4}\frac{rr}{AA} - \frac{4}{9} (pp - 2q),$\\
	\hth$  HO\:HO = \frac{9}{16} \frac{rr}{AA} - (pp - 2q),$\\
	\hth$  MO\:MO = \frac{1}{16} \frac{rr}{AA} - \frac{1}{9} (pp -2q).$\\
	Therefore $MO = \frac{1}{2}HM = \frac{3}{2}HO$ and $HO = HM + MO$
	therefore the points $H,$ $M$ and $O$ are collinear.\\
	If $I$ is the center of the inscribed circle, Euler determines also
	$HI,$ $GI$ and $IO.$

	\sssec{Theorem. [Euclid, Book 3, Propositions 35 and 36.]}
	\label{sec-powercircle}
	{\em If 2 lines through $M,$ not on a circle meet that circle,
	the first one in $A$ and $B,$ the second one in $C$ and $D,$ then\\
	\hth$| MA|  | MB|  = | MC|  | MD| .$}

	\sssec{Theorem.}
	{\em Let $A_0 A_1 A_2$ be a triangle and $H_0 H_1 H_2$ be the feet of
	the perpendiculars from the vertices to the opposite sides.
	then $A_0 H_0$ bisects the angle $H_1 H_0 H_2.$}

	Proof: If $H$ is the orthocenter, the quadrangle $H H_1 A_2 H_0$
	can be inscribed in a circle and therefore the angles
	$A_0 H_0 H_1$ and $H_2 A_2 A_0$ are equal.  Similarly the angles
	$H_2 H_0 A_0$
	and $A_0 A_1 H_1$ are equal, but the angles $H_2 A_2 A_0$ and
	$A_0 A_1 H_1$ are
	equal because there sides are perpendicular, therefore $A_0 H_0 H_1$
	and $A_0 A_1 H_1$ are equal.

	\sssec{Definition.}
	The triangle $H_0 H_1 H_2$ is called the {\em orthic triangle.}

%
%
%
	\ssec{Menelaus (about 100 A.D) and Ceva (1647-1734?).}
	\label{sec-Smence}

	\setcounter{subsubsection}{-1}
	\sssec{Introduction.}
	The following theorems give a metric characterization of
	three points on the sides of a triangle which are collinear or
	which are
	such that the line joining these points to the opposite vertex are
	concurrent.  For these theorems, an orientation is provided on each of
	the sides and therefore the distances have a sign.  The theorems are as
	follows:

	\sssec{Theorem. [Menelaus]}
	{\em If $X_0$ is on $a_0,$ $X_1$ is on $a_1$ and $X_2$
	is on $a_2,$
	then the points $X_0,$ $X_1$ and $X_2$ are collinear iff\\
	\hth$| A_1 X_0| | A_2 X_1| | A_0 X_2|
		= | A_2 X_0| | A_0 X_1| | A_1 X_2| .$}

	\sssec{Theorem. [Ceva]}
	{\em If $X_0$ is on $a_0,$ $X_1$ is on $a_1$ and $X_2$ is on $a_2,$ then
	the lines $A_0 X_0,$ $A_1 X_1$ and $A_2 X_2$ are concurrent iff\\
	\hth$ | A_1 X_0| | A_2 X_1| | A_0 X_2|
	= - | A_2 X_0| | A_0 X_1| | A_1 X_2| .$}

	The following theorem
	is a direct consequence of the theorem of Ceva.
	Theorem \ldots see Coxeter, I believe I saw it later ???

	\sssec{Theorem.}
	{\em Let $X$ be a point not on the sides of a triangle $A_0 A_1 A_2,$
	let $X_0,$ $X_1,$ $X_2$ be the intersection of $X A_0$ with $A_1 A_2,$
	of $X A_1$ with $A_2 A_0,$ and of $X A_2$ with $A_0 A_1,$,
	let $Y_0,$ $Y_1,$ $Y_2$ be the other intersection of the circle through
	$X_0,$ $X_1$ and $X_2$ with the sides of the triangle,
	then $A_0 Y_0,$ $A_1 Y_1$ and $A_2 Y_2$ are concurrent.}

	Proof:  If we eliminate $| A_i X_j|$  from the relation of Ceva and from
	the relations\\
	\hth$| A_0 X2| | A_0 Y_2| = | A_0 X1| | A_0 Y_1| $\\
	\hth$	| A_1 X0| | A_1 Y_0| = | A_1 X2| | A_1 Y_2| $\\
	\hth$	| A_2 X1| | A_2 Y_1| = | A_2 X0| | A_2 Y_0| $\\
	obtained from Theorem \ref{sec-powercircle}, we obtain\\
	\hth$| A_1 Y_0| | A_2 Y_1| | A_0 Y_2|
		= - | A_2 Y_0| | A_0 Y_1| | A_1 Y_2| .$\\
	Therefore by the Theorem of Ceva, the lines 
	$A_0 Y_0,$ $A_1 Y_1$ and $A_2 Y_2$ are concurrent.

	\ssec{Euler (1707-1783) and Feuerbach (1800-1834).}
	\label{sec-Seulfe}

	\sssec{Inroduction.}
	The geometry of the triangle has its origin in the
	following theorems.

	\sssec{Theorem.}\label{sec-tbaryc}
	{\em The 3 medians of a triangle meet at a point called the
	barycenter or, in mechanics, the center of mass.}

	\sssec{Theorem.}\label{sec-torthoc}
	{\em The 3 altitudes of a triangle meet at a point called the
	orthocenter.}

	\sssec{Theorem.}\label{sec-tcenterc}
	{\em The 3 mediatrices of a triangle meet at a point which is
	the center of the circumcircle.}

	\sssec{Theorem.}
	{\em The 3 bisectrices of a triangle meet at a point which is
	the center of the inscribed circle.}

	\sssec{Theorem. [Euler]}\label{sec-teulerhgo}
	{\em The points $H,$ $G$ and $O$ are on a line, called the
	line of Euler, moreover\\
	\hth$	|HG = 2 |GO|$ and $|HO| = 3 |GO|.$}

	The proof of Euler is algebraic.  He determines the distance $HG,$ $GO$
	and $HO$ in terms of the length of the sides of the triangle.  Other
	distances are also determined in the same paper.

	\sssec{Theorem. [Brianchon and Poncelet]}\label{sec-tbrponc}
	{\em The mid-points of the sides of
	a triangle, the feet of the altitudes and the mid-points of the
	segments joining the orthocenter to the vertices of the triangle are on
	a circle, called the circle of Brianchon-Poncelet.  It is also called
	the 9 point circle or the circle of Feuerbach, who discovered it
	independently, and improperly the circle of Euler.}

	\sssec{Theorem. [Feuerbach]}\label{sec-tfeuer}
	{\em The circle of Brianchon-Poncelet is tangent to
	the inscribed circle and to the three excribed circles, the point of
	tangency for the inscribed circle is called the point of Feuerbach.}

	The proof given by Feuerbach is algebraic and trigonometric in
	character.  It expresses distances in terms of the length of the sides
	and of the trigonometric functions of the angles of the triangle.

	\ssec{The Geometry of the Triangle.  Lemoine (1840-1912).}
	\label{sec-Striangle}

	\setcounter{subsubsection}{-1}
	\sssec{Introduction.}
	An interesting development of Euclidean geometry occured
	during the 19-th century, known under the name of the geometry of the
	triangle.  The activity in this area was most intense during the
	period 1870-1900.  A large number of elementary results were obtained
	especially in Belgium and France, but also in England, Germany and
	elsewhere.  Strictly speaking, the Theorem of Euler of
	\ref{sec-teulerhgo} can be considered as the first important new result
	in this connection since Euclid.  Others which prepared the way were
	the theorems of Brianchon-Poncelet of \ref{sec-tbrponc} and the Theorem
	of Feuerbach of \ref{sec-tfeuer}.  A few theorems
	will be extracted from the long list.

	\sssec{Theorem. [Schr\"{o}ter]}
	{\em If $a\times b$ denotes the point on $a$ and $b$ and
	$A\times B$ denotes the line through $A$ and $B$,\\
	Let}\\
	\hth$F_0 := (M_1 \times H_2) \times (M_2 \times H_1),$\\
	\hth$F_1 := (M_2 \times H_0) \times (M_0 \times H_2),$\\
	\hth$F_2 := (M_0 \times H_1) \times (M_1 \times H_0).$\\
	\hth$G_0 := (M_1 \times M_2) \times (H_1 \times H_2),$\\
	\hth$G_1 := (M_2 \times M_0) \times (H_2 \times H_0),$\\
	\hth$G_2 := (M_0 \times M_1) \times (H_0 \times H_1).$
	\enumb
	\item {\em $F_0,$ $F_1$ and $F_2$ are on the line $e$ of Euler.}
	\item {\em $A_0 \times G_0,$ $A_1 \times G_1$ and $A_2 \times G_2$
	are parallel and are perpendicular to $e.$}
	\item {\em $A_0,$ $F_0,$ $G_1$ and $G_2$ are collinear, and so are
	$A_1,$ $F_1,$ $G_2$ and $G_0$ as well as
	$A_2,$ $F_2,$ $G_0$ and $G_1$.}
	\item {\em $G_0,$ $G_1,$ $G_2$ are the vertices of a triangle conjugate
	to the circle of Brianchon-Poncelet.}
	\item {\em $M_0 \times G_0,$ $M_1 \times G_1$ and $M_2 \times G_2$
	pass through the same point $S.$}
	\item {\em $H_0 \times G_0,$ $H_1 \times G_1$ and $H_2 \times G_2$
	pass through the same point $S'.$}
	\item {\em $S$ and $S'$ are on the circle of Brianchon-Poncelet.}
	\item {\em $S$ and $S'$ are on the polar of $H$ with respect to the
	triangle $A_0,$	$A_1,$ $A_2.$ }
	\enume

	$S$ and $S'$ are called the {\em points of Schr\"{o}ter.}

	The proof of this Theorem published by Schr\"{o}ter in ``Les Nouvelles
	Annales de Math\'{e}matiques" in 1864, was obtained by several people.
	The published proof is that of a student of Sainte-Barbe, L. Lacachie.
	The Theorem is generalized to Projective Geometry in III.D8.1,D8.2,C8.0.
	It is stated in finite involutive geometry in III.\ref{sec-tinvschr}.


\setcounter{section}{1}
	\section{Projective Geometry.}

	\ssec{The preparation. Menaechmus (about 340 B.C.), Apollonius
	(260? B.C - 200? B.C.), Pappus (300 - ?).}\label{sec-Smenap}

	The projective geometry has its source in the discovery of the conic
	sections, the ellipse, the parabola and the hyperbola, which is
	ascribed by Proclus to the Greek mathematician Menaechmus, a pupil of
	Plato and Eudoxus.
	The conic sections were studied by Aristaeus the Elder, Euclid,
	Archimedes, Pappus of Alexandria and finally by Apollonius of Perga.
	The conics are defined as the intersection of a (circular) cone by a
	plane not passing through its vertex.
	If we make a cut of the cone with a
	plane through the vertex we obtain
	two lines $c_1$ and $c_2.$  Line $a$ is the cut
	of a plane giving an hyperbola, line $b$
	is the cut of a plane giving $a$
	parabola, line $c$ gives an ellipse and
	line $d$ gives the special case of a circle.\\
\begin{picture}(250,140)(-50,-20)
\put(50,64){$V$}
\thicklines
\put(5,45){\line(3,1){190}}
\put(175,93){$c_1$}
\put(5,75){\line(3,-1){190}}
\put(175,23){$c_2$}
\thinlines
\put(5,50){\line(1,0){190}}
\put(190,44){$a$}
\put(95,0){\line(0,1){100}}
\put(98,10){$d$}
\put(90,100){\line(1,-2){50}}
\put(110,65){$b$}
\put(5,65){\line(3,-1){190}} \put(175,1){$c$}
\end{picture}\\

	Among the many contributions of Pappus I will cite the discovery
	that the anharmonic ratio of 4 points is unchanged after projection,
	where the anharmonic ratio of $A$, $B$, $C$ and $D$
	is $\frac{dist(C,A)\:dist(D,B}{dist(C,B)\:dis(D,A)}$. This is a
	fundamental property in geometry.

	The important notion of point at infinity can be traced to Kepler,
	in 1604, and Desargues, in 1639 (see Heath, I, p. 193).
	This leads to the notion of the extended Euclidian plane which
	contains besides the ordinary points, the directions, each one is what
	is what is common to the set of parallel line, and the set of
	all directions, or line at infinity.

	\ssec{G\'{e}rard Desargues (1593-1661) and Blaise Pascal
	(1623-1662).}\label{sec-Sdespa}

	\setcounter{subsubsection}{-1}
	\sssec{Introduction.}
	The extensive study of the conics by Apollonius was
	eventually taken up again by Pascal.  One of his many new results is
	Theorem \ref{sec-tpascal} which allows the construction
	\ref{sec-cpascal} and \ref{sec-cmaclaurin} of a
	conic using the ruler only.  The second construction is attributed to
	MacLaurin.  But the two constructions are closely related to each
	other as will be seen.  The Theorem of Pascal was generalized to $n$
	dimension by Arthur Buchheim in 1984.

	\sssec{Notation.}
	I will introduce in III.\ref{sec-defnot} detailed notations which allow
	a compact description of constructions.  For instance,\\
	\hth$	a_0 := A_1 \times A_2$\\
	means that the line $a_0$ is defined as the line through the 2 points
	$A_1$ and $A_2.$

	\sssec{Theorem. [Pascal]}\label{sec-tpascal}
	{\em Given the points $A_0,$ $A_1,$ $A_2,$ $A_3,$ $A_4$ and $A_5.$\\
	Let $P_0$ be the point common to $A_0 \times A_1$ and $A_3 \times A_4,$
	let $P_1$ be the point common to $A_1 \times A_2$ and $A_4 \times A_5,$
	let $P_2$ be the point common to $A_2 \times A_3$ and $A_5 \times A_0.$
	A necessary and sufficient condition for $A_0,$ $A_1,$ $A_2,$ $A_3,$
	$A_4$ and $A_5$ to be on the same conic is that $P_0,$ $P_1$ and $P_2$
	be collinear.}

	This theorem leads to 2 construction of conics.

	\sssec{Construction.[Pascal]}\label{sec-cpascal}
	Given 5 points $A_0,$ $A_1,$ $A_2,$ $A_3$ and $A_4.$
	To each line through $A_4$ corresponds a point $A_5$ on the conic.\\
	\hth$a_0 := A_0 \times A_1,$ $a_1 := A_1 \times A_2,$
	$a_2 := A_2 \times A_3,$ $a_3 := A_3 \times A_4,$\\
	\hth$P_0 := a_0 \times a_3,$ $a_4$ is an arbitrary line through $A_4,$\\
	\hth$P_1 := a_1 \times a_4,$ $e := P_0 \times P_1, P_2 := e \times a_2,
	$\\
	\hth$a_5 := A_0 \times P_2,$ $A_5 := a_4 \times a_5.$

	\sssec{Construction. [MacLaurin]}\label{sec-cmaclaurin}
	\vspace{-19pt}\hspace{172pt}\footnote{as stated by Braikenridge}\\[8pt]
	If the sides of a triangle pass through three fixed points, and two
	vertices trace straight lines, the third vertex will trace a conic
	through two of the given points.

	The proof follows from Pascal's Theorem.  The construction
	can be given in the following explicit form:\\
	\hth$A_0,$ $A_1,$ $A_2,$ $A_3,$ $A_4$ are 5 given points.\\
	To each line $l$ through $P_0$
	will correspond a point $A_5$ on the conic.\\
	\hth$a_0 := A_0 \times A_1,$ $a_1 := A_1 \times A_2,$
	$a_2 := A_2 \times A_3,$ $a_3 := A_3 \times A_4,$\\
	\hth$P_0 := a_0 \times a_3,$ 
	$P_1 := l \times a_1,$ $P_2 := l \times a_2,$\\
	$a_4 := A_4 \times P_1,$
	$a_5 := A_0 \times P_2,$ $A_5 := a_4 \times a_5.$\\
	The triangle is $\{P_1,P_2,A_5\},$ $P_1$ is on $a_1,$ $P_2$ is on $a_2,$
	$P_1 \times P_2$ passes through $P_0,$ $P_1 \times A_5$ passes through
	$A_4,$ $P_2 \times A_5$ passes through $A_0.$

	\sssec{Comment.}
	Pascal would not have easily accepted a finite geometry.
	Indeed in his "Pens\'{e}es", he says (p. 567),
	\begin{quote}that there are no
	geometers which do not believe that space is infinitely divisible.
	\end{quote}
	Also discussing both the infinitely large and the infinitely small,
	he writes (p. 564)
	\begin{quote}In one word, whatever the motion, whatever the
	number, whatever the space, whatever the time, there is always one
	which is larger and one which is smaller, in such a way they they
	sustain each other between nothing and infinity, being always
	infinitely removed from those extremes.  All these truths cannot be
	proven, and still they are the foundations and the principles of
	geometry.
	\end{quote}

	\ssec{Lazare Carnot (1783-1823).}

	A contemporary of Poncelet, Carnot obtained many results of which
	the following is in the line of Manelaus and Ceva applied to conics.

	\sssec{Theorem. [Carnot]}
	{\em If a conic cuts the side $A \times B$ of a triangle $\{A,B,C\}$
	at $C_1$ and $C_2,$ and similarly the side $B\times C$ cut the conic
	at $A_1$ and $A_2$ and the side $C\times A$ at $B_1$ and $B_2$, then
	the oriented distances satisfy\\
	\hth$AC_1.AC_2.BA_1.BA_2.CB_1.CB_2 = AB_1.AB_2.BC_1.BC_2.CA_1.CB_2$}

	This is generalized to curves of degree $n$.\footnote{Eves p.358}
 
	\sssec{Theorem.}
	{\em Let $A_0 B_0 C_0$ be a triangle and $X$ be a point not on its
	sides,\\
	Let $A_0 \times X$ meet $A_1 \times A_2$ at $X0,$ $A_1 \times X$ meet
	$A_2 \times A0$ at $X_1$ and $A_2 \times X$ meet
	$A_0 \times A_1$ at $X_2.$  Let $Y_0$ be a point on $A_1 \times A_2,$
	$Y1$ be a point on $A_2 \times A_0$
	and $Y2$ be a point on $A_0 \times A_1,$ then a necessary and sufficient
	condition
	for $X_0,$ $X_1,$ $X_2,$ $Y_0,$ $Y_1,$ $Y_2$ to be on the same conic
	is that the lines $A_0 \times Y_0,$ $A_1 \times Y_1,$
	$A_2 \times Y_2$ be concurrent.}

	This is a consequence of the Theorem of Carnot.

	\ssec{Jean Poncelet (1788-1867).}\label{sec-Sponc}

	The work of Poncelet done while a prisoner of Russia at the end of
	Napoleon's campaign, was fundamental in isolating those properties of
	Euclidean geometry which are independent of the notions of distances
	and measure of angles and dependent only on incidence properties and
	appropriate axioms which involve only incidence.
	One of is celebrated Theorems is the following.

	\sssec{Theorem.}
	{\em If a $n$ sided polygon is inscribed in a conic and outscribed
	to an other conic, then if with start from any point on the first conic
	and draw a tangent to the second, then obtain the other intersection
	with the first conic and repeat the construction, the new polygon
	closes after $n$ steps.}

	There are many proofs of this Theorem.  The proofs which are done using
	the theory of elliptic functions, suggested to me that the Jacobi
	elliptic functions could be generalized to the finite case.

%
%
	\ssec{Joseph Gergonne (1771-1858).}\label{sec-Sgerg}
	Gergonne was the first to recognize the property of duality
	which plays a fundamental role in projective geometry.\footnote{Coxeter,
	p.13}

%
	\ssec{Michel Chasles (1793-1880).}\label{sec-Schas}
	Chasles greatest contribution to projective geometry, according to
	Coolidge \footnote{p.96.} is the study of the cross ratio also called
	anharmonic ratio.\footnote{See also Coxeter, p.165.}


{\tiny\footnotetext[1]{G13.TEX [MPAP], \today}}

\setcounter{section}{2}
	\section{Relation between Projective and Euclidean Geometry.}

	\setcounter{subsection}{-1}
	\ssec{Introduction.}

	Projective geometry is concerned only with those properties in
	geometry which are preserved under projection.  Euclidean, as well as
	non Euclidean geometry can be derived from projective geometry.
	The connection through transformation groups will be described in
	section \ref{sec-Sklein}.

	The first connection goes back to the work of Poncelet, but it is has
	been deemphasized in the teaching of the subject, except for the first
	step (affine geometry).  I will presently summarize this approach.
	Terms which are unknown to the reader, will be defined in the later
	Chapters.

	In projective geometry, no line is distinguished from any other, no
	point is similarly distinguished.  The main notions are those of
	incidence, perspectivity, projectivity, involution and polarity, the
	last notion leading naturally to conics.  Euclidean geometry can be
	considered as derived from projective geometry by choosing some elements
	in it and distinguishing them from all others. 	I will proceed in 3
	steps.

	\ssec{Affine Geometry.}

	\setcounter{subsubsection}{-1}
	\sssec{Introduction.}
	In this first step one line is distinguished.  This
	line is called the ideal line, or line at infinity.  When we do so, we
	obtain the so called affine geometry.  Points fall now into two
	categories, the ordinary points, which are not on the ideal line and
	the ideal points which are.  Lines fall in two categories, the ideal
	line and the others which we can call ordinary.  From the basic notion
	of parallelism  follow the derived notions of parallelogram, equality
	of vectors on the same line or on parallel lines, trapeze or rhombus,
	mid-point, barycenter, center of a conic, area of triangles.

	\sssec{Definition.}
	Two distinct ordinary lines are {\em parallel} iff their common
	point is an ideal point.

	\sssec{Definition.}
	A {\em vector} $\varss B,C$ is an ordered pair of points.

	\sssec{Definition.}
	If the lines $B \times C$ and $D \times E$ are parallel, the vectors
	$\varss B,C$ and $\varss D,E$ {\em are equal} iff the lines
	$B \times D$ and $C \times E$ are also parallel.

	\sssec{Definition.}
	If $B,$ $C,$ $D$ and $E$ are on the same line, the vectors $\varss
	B,C$ and $\varss D,E$ {\em are equal} iff there exists 2 points $F$
	and $G$ on a parallel line, such that $\varss B,C = \varss F,G$ and
	$\varss F,G = \varss D,E.$  This definition has,
	of course, to be justified.  It can be replaced by:
	$\varss B,C$ and $\varss D,E$ {\em are equal} iff there exists a
	parabolic projectivity,
	with the fixed point being the ideal point on the line, which
	associates $C$ to $B$ and $E$ to $D.$

	\sssec{Definition.}
	The {\em center} of a conic is the pole of the ideal line, in
	the polarity whose fixed points are the conic.

	\sssec{Definition.}
	Two points are {\em conjugate} iff one is on the polar of the
	other.

	\sssec{Theorem.}
	{\em Conjugate points on a given line determine an involution.}

	\sssec{Definition.}
	A {\em parabola} is a conic tangent to the ideal line.
	The point of tangency is called the {\em direction of the parabola}.

	\sssec{Example.}
	The parabola $y^2 = 4c x,$ in homogeneous coordinates is\\
	\hth$Y^2 = 4 c X Z.$\hfill(1)\\
	Its intersection with $Z = 0$ is $Y = 0.$  The parabola is tangent to
	$Z$ = 0 at (1,0,0).

	\sssec{Definition.}
	The {\em focus of a parabola} is the intersection of the
	ordinary tangents to the parabola from the isotropic points.  The
	{\em directrix} of the parabola is the polar of the focus.  The
	{\em axis} of the parabola is the line through the focus and the
	direction of the parabola.  The {\em vertex} of the parabola is the
	point of the parabola on its axis.

	\sssec{Example.}
	The tangent to the parabola at $(X_0, Y_0,1)$ is\\
	\hth$2c \:X - Y_0 \:Y + 2c\:X_0 \:Z = 0.$\hfill(2)\\
	It passes through the isotropic point $(1,i,0)$ if $2c = Y_0\:i,$ hence
	because of (1), $X_0 = -c.$  The tangent is therefore
	$X + Y\:i - c = 0.$  The tangent from the other isotropic point is
	$X - Y\:i - c = 0.$  They both intersect at $(c,0,1).$\\
	The polar is obtained by substituting in (2) this point for
	$(X_0,Y_0,1),$ this gives $X = -c\:Z.$  The axis is $Y = 0,$ the
	vertex is (0,0,1).

	\sssec{Comment.}
	The terminology can be changed by accepting as points and
	lines only those which are ordinary.  An ideal point is renamed a
	direction.  We obtain in this way, something which is closer to the
	terminology used by Euclid.

	\ssec{Involutive geometry.}

	\setcounter{subsubsection}{-1}
	\sssec{Introduction.}
	The second step consists in considering the involutions
	on the ideal line.  Among all the involutions we can distinguish one of
	them and call it the fundamental involution.  Three cases are possible,
	the involution may have 2 fixed ideal points, in which case it is
	called hyperbolic, one fixed point in which case it is called
	parabolic and no fixed point, in which case it is called elliptic.
	If we extend the projective geometry to the complex case, these ideal
	points then exist, but are not real.\\
	The elliptic case, which leads to Euclidean Geometry and the
	hyperbolic case which leads to the Geometry of Minkowski can be
	studied together.  The parabolic case, which leads to the Galilean
	Geometry is studied separately.\\
	Using the fundamental involution, either elliptic or hyperbolic,
	we can introduce the basic notion of perpendicularity and from it
	follow the derived notion of right triangle, rectangle, altitude,
	orthocenter, circle, equal segment, isosceles and equilateral
	triangles, center of circumcircle, Euler line, circle of
	Brianchon-Poncelet.

	In the alternate second step, one involution with 2 real fixed points
	is distinguished.  It is only if we stay with real projective geometry
	as opposed to complex projective geometry that the hyperbolic
	involutive geometry is distinct from the elliptic involutive geometry.
	Staying with real projective geometry, the notions which are introduced
	can be given the same name as in the elliptic involutive geometry, the
	definitions may differ slightly, but properties are quite analogous.

	\sssec{Definition.}
	When the fundamental involution has no real fixed points,
	I will call the geometry {\em elliptic involutive geometry}.\\
	When the fundamental involution has no real fixed points,
	I will call the geometry {\em hyperbolic involutive geometry}.

	\sssec{Definition.}
	The fixed points of the fundamental involution are called
	{\em isotropic points}.  Any ordinary line through an isotropic point
	is called an {\em isotropic line}.  Strictly speaking, the {\em ideal
	points} are those on the ideal line which are not isotropic, and the
	{\em ordinary lines} are those which are not isotropic.

	\sssec{Definition.}
	Two {\em lines are perpendicular} iff their ideal points are
	pairs of the fundamental involution.

	\sssec{Definition.}
	A conic is a {\em circle} iff the involution that the conic
	determines on the ideal line is the fundamental involution.

	\sssec{Theorem.}
	{\em A conic which passes through the 2 isotropic points is a circle .}

	\sssec{Definition.}
	A {\em segment} $[A B]$ is an unordered pair of points.

	\sssec{Definition.}
	The segment $[A B]$ and the segment $[C D]$ {\em are equal} iff
	the point $E$ constructed in such a way that $ACDE$ is a parallelogram,
	is such that $E$ and $B$ are on the same circle centered at $A.$

	\sssec{Definition.}
	The {\em center of a circle} is the intersection of the tangents
	to the circle at the isotropic points.

	\sssec{Comment.}
	A geometry could also be constructed in which the
	correspondence on the ideal line associates every point to one of them.
	This corresponds, using algebra, to the transformation\\
	\hth$	T(x) = \frac{a x + b}{c x - a},$ $aa + bc = 0.$\\
	This is the parabolic involutive geometry.

	Before leaving the subject of involutive geometry, I would like to make
	the following observation, which will be useful to understand
	terminology in non-Euclidean geometry.
	The step to construct non-Euclidean geometry from projective geometry,
	which correspond to involutive geometry, is to choose a particular
	conic as ideal, or set of ideal points.  In view of the fact that a
	line conic can
	degenerate in the set of lines passing through either one or the other
	of 2 points, we can observe that the ideal in the involutive geometry
	is such a degenerate conic.  This analogy will be pursued to define,
	using the ideal conic, notions in non-Euclidean geometry which are
	related to notions of Euclidean geometry and will help in an economy of
	terminology, but nothing more.

%
%
%

{\tiny\footnotetext[1]{G14.TEX [MPAP], \today}}

\setcounter{section}{3}
	\section{Analytic Geometry.}

	\ssec{Ren\'{e} Descartes (1596-1650)[La G\'{e}om\'{e}trie].}

	The prime motivation of Descartes when he wrote, "La G\'{e}om\'{e}trie"
	appears to have been a long standing problem, the determination of
	the locus of Pappus.\footnote{p. 8}\\
	In present day notation, given lines $l_i$ and angles
	$\alpha_i$, the problem is to determine the locus of
	a point $C$ and its $\alpha_i$ projections $U_i$ on $l_i$, such that
	the angle of $C\times U_i$ with $l_i$ is $\alpha_i$, and for instance,
	with i = 0,1,2,	3, such that\\
	\hth$|CU_0|\:|CU_2| = k\:|CU_1|\:|CU_3|$.\hfill(1)

	Descartes chooses as axis $l_0$ and $u_0 := C\times U_0$, he chooses
	also some
	orientation which allows him to associate to the points on these axis,
	some real number.  If $x := |U_0,A_1|,$ $y := |C,U_0|,$
	$a_i := |A_1,A_i|$, $i = 1,2,3,$ if $X_i$ are the intersection of $l_i$
	with $a$, then the prescribed angles imply by
	similarity\\
	\hth$\frac{|U_0X_i|}{|U_0A_i|} = \frac{b_i}{e},
		\frac{|CU_i|}{|CX_i|} = \frac{c_i}{e},$\\
	for some $b_i$, $c_i$ and unit of distance $e$.\\
	The distances $|CU_i|$ are linear functions of $x$ and $y$ and therefore
	replacing in (1) gives the equation of a conic through $A_1$.
	By symmetry, the conic passes through $A_3$, $B_1$ and $B_3$.\\
	Indeed,\\
	$|CU_0| = y,$,\\
	$U_0A_1 = x,$ $U_0X_1 = \frac{b_1}{e}x,$ $|CX_1| = |CU_0| + |U_0X_1| =
	y + \frac{b_1}{e},$ $CU_1 =  (y + \frac{b_1}{e})\frac{c_1}{e},$\\
	$|U_0A_2| = x + a_2,$ $U_0X_2 = (x + a_2)\frac{b_2}{e},$
	$CX_2 = y + (x + a_2)\frac{b_2}{e},$\\
	$CU_2 = (y + (x + a_2)\frac{b_2}{e})\frac{c_2}{e},$
	$CU_3 = (y + (x + a_3)\frac{b_3}{e})\frac{c_3}{e}.$\\

	Nowhere, in his work are the axis or arrows on them indicated 
	specifically or are the axis chosen at a right angle, except if
	convenient to solve the problem at hand.

	\ssec{After Descartes.}
	Using modern terminology, the problem posed by Descartes, was to
	construct an algebraic structure which is isomorphic to Euclidean
	geometry.  More precisely the problem is to obtain algebraic elements
	$P'$ which are in one to one correspondence with points $P,$
	algebraic elements $l'$ which are in one to one correspondence with
	lines $l,$ an algebraic relation $P' \cdot l' = 0$ associated to the
	incidence
	relation in geometry, $P$ is on $l$ or $l$ is through $P,$ written
	$P \cdot l = 0,$ such that if $l'$ corresponds to $l$ and $P'$ to $P,$
	$P' \cdot l' = 0$ if and only if $P \cdot l = 0.$\\
	Similar correspondences have to be given for perpendicularity, equality
	of angles and segments, measure of angle and segments, etc.
	Descartes' solution is to choose 2 lines $xx$ and $yy$
	in the Euclidean plane and to associate, if these are perpendicular,
	to a point $P$ the 2 real numbers $x$ and
	$y$ which are the distances from $P$ to $yy$ and to $xx.$\\
	\hth$	P' = \delta(P) = (x,y).$\\
	This correspondence is not one to one.  If $x, y \neq 0,$ there are
	four 	points which will give the same pair $(x, y).$
	To solve this problem a sign must be associated to the distances,
	corresponding to an orientation on the lines $xx$ and $yy.$  Usually,
	with $xx$ horizontal and $yy$ vertical, $x$ is positive to the right
	of $yy,$ $y$ is positive above $xx.$\\
	The distance between $(x,y)$ and $(x',y')$ of the points
	$P$ and $Q$ is given by\\
	\hth$	\sqrt{ (x'-x)^2 + (y'-y)^2 }.$\\
	To represent the lines, several choices are possible, one such choice,
	is the pair $[m,b],$ where $b$ is the (oriented) slope and $b$ it the
	distance from the intersection of the line with $yy,$ the so called
	$y$ intercept. 	In this case, if $(x, y)$ corresponds to the point
	$P$ and $[m, b]$ to the line $l,$ $(x, y)$ is on $[m, b]$ if and only
	if\\
	\hth$	y = mx + b.$\\
	Perpendicularity of $[m, b]$ and $[n, c]$ is defined by $m\: n = -1.$\\
	The difficulty of this representation is that lines perpendicular to
	$xx$ do not have a (finite) slope.  Reversing the role of $xx$ and
	$yy$ does not help.\\
	An other representation of lines that can be chosen, is to take the
	pair $\{l_0, l_1\}$ of the distances $l_1$ and $l_0$ from the origin
	to the intersection of the line $l$ with $xx$ and $yy,$\\
	\hth$	l' = \delta (l) = \{l_0,l_1\},$\\
	with the incidence property represented by the relation\\
	\hth$	l_0 x + l_1 y - l_0 l_1 = 0.$\\
	In particular, the points $(l_1,0)$ and $(0,l_0)$ are on $l'.$
	The perpendicularity property of $l'$ and $m'$ is represented by the
	bilinear relation\\
	\hth$	l_0 m_1 + l_1 m_0 = 0.$\\
	This again is not suitable because this representation fails for lines
	through the origin.\\
	The correspondence finally chosen by Descartes is a triple of real
	numbers $[a, b, c]$ which are obtained from $l_0,$ $l_1$ and
	$l_0$ $l_1$ by
	multiplication by some arbitrary non zero real $k.$\\
	\hth$	[a, b, c] = k [l_0, l_1, l_0 l_1].$\\
	For a line through the origin $c = 0,$ $\frac{b}{a}$ is the slope,
	A line parallel to $yy$ is represented by $(1,0,c)$ where $-c$ is the
	$x$ intercept.\\
	The incidence property is the familiar linear relation\\
	\hth$	a x + b y + c = 0.$\\
	But it is important to realize that the correspondence is not one to
	one.  The line is represented by the set of all triples corresponding to
	all the possible value of $k,$ a so called equivalence class, the
	numbers
	$a,$ $b,$ $c$ are called the homogeneous coordinates of the line.
	Perpendicularity of $[a, b, c]$ and $[a', b', c']$ is represented by\\
	\hth$	a a' + b b' + c c' = 0.$\\
	By analogy, one could represent points by a triple $(x, y, 1)$ or by
	any equivalent set $(X, Y, Z) = k(x, y, 1),$ $k\neq  0.$
	This implies that $Z = k\neq  0$ and $X = kx,$ $Y = ky$ or
	$x = \frac{X}{Z},$ $y = \frac{Y}{Z}.$\\
	The incidence property is then\\
	\hth$ a X + b Y + c Z = 0.\hfill(2)$\\
	$(X, Y, Z)$ are the so called homogeneous coordinates of an algebraic
	point.

	\ssec{Jean Poncelet (1788-1867).}

	Poncelet was one of the first to take full advantage of the fact that
	parallel lines define a
	direction, which can be called the point at infinity and that all the
	points at infinity can be considered to be on a line, the line at
	infinity.  This constitutes the decisive step towards the development
	of projective geometry.\\
	The algebraic points $(X, Y, 0),$ with $X$ and $Y$ not both 0,
	correspond to
	these new geometric points, They are all on the line [0,0,1] which is
	the line at infinity. The distance between the algebraic points
	$(X, Y, Z)$ and $(X', Y', Z'),$
	with $Z$ and $Z'\neq 0,$ is given by\\
	\hth$	\sqrt{ (\frac{X'}{Z'} - \frac{X}{Z})^2
	+ (\frac{Y'}{Z'} - \frac{Y}{Z})^2 }.$

	\sssec{Comment.}
	The extension of the Euclidian plane by adding the points at
	infinity and a line at infinity is distinct from the extension
	of the complex plane, in which to all the points $x + iy,$ $x$ and $y$
	real (and $i^2$ = -1), we add 1 point at infinity.  In a complex plane,
	all lines pass through the point at infinity.

	It is not the place here to review all the other basic formulas of
	analytic geometry.  However, there is an important consequence of the
	isomorphism between synthetic geometry and analytic geometry, which
	is implicit in the work of Poncelet and is associated to the properties
	of circles, which was the basis of Poncelet's method to obtain
	properties for conics in general.

	The equation of a circle of center $(a,b)$ and radius $R$ is\\
	\hth$	(x-1)^2 + (y-b)^2 = R^2,$\\
	or in homogeneous coordinates,\\
	\hth$	(X - aZ)^2 + (Y - bZ)^2 = R^2 Z^2.$\\
	The points on the circle and on the line at infinity $Z = 0$ satisfy\\
	\hth$	X^2 + Y^2 = 0,$\\
	which has no real solution.
	The introduction of complex numbers, whose use had become standard
	by the time of Poncelet, suggested the definition of a complex
	analytic geometry, with elements\\
	\hth$	(X,Y,Z) = k(X,Y,Z),$ $k,$ $X,$ $Y,$ $Z$ complex, $k\neq  0$
	and not all $X,$ $Y,$ $Z$ equal to zero, and with elements\\
	\hth$	(a,b,c) = k'(a,b,c),$ $k',$ $a,$ $b,$ $c$ complex,
	$ k'\neq  0$ and not all $a,$ $b,$ $c$ equal to zero,
	The incidence property being again (1).\\
	The complex elements which are not real correspond to new points and
	lines in synthetic geometry, the complex points and the complex lines.
	In this structure, $(1,i,0)$ and $(1,-i,0)$ are 2 points on the line at
	infinity which are also on every circle.  They are called isotropic
	points and play an essential role in both Euclidean geometry, extended
	to the complex and in what I call involutive geometry.


	\ssec{James Singer on Difference sets and finite projective
	Geometry.}

	\setcounter{subsubsection}{-1}
	\sssec{Introduction.}
	Inspired by the paper of Veblen and MacLagan-Wedderburn of 1907,
	Singer introduced in October 1934 (Singer, 1938, Baumert, 1971) the
	important concept of cyclic difference
	sets which allows for an arithmetization of projective geometry which
	is as close to the synthetic point of view as is possible. With this
	notion, it becomes possible to label points and hyperplanes in $N$
	dimensional projective geometry of order $p^k$.  With it, in the plane,
	it is not only
	trivial to determine all the points on a line, and lines incident to a
	point but also the lines through 2 points and points on 2 lines.\\
	Completely independently, one of my first students at the
	"Universit\'{e} Laval", Quebec City, made the important discovery
	that the regular polyhedra can be used as models for finite geometries
	associated with
	2, 3 and 5.  Then, he introduced the nomenclature of selector
	(s\'{e}lecteur) for the notion of cyclic difference sets, to construct
	an appropriate numbering of the points and lines on the polyhedra.
	The definition
	of selector function and selector correlation is implicit in his work.\\
	The notion of cyclic difference sets makes duality explicit through the
	correlation, which is the polarity when $p \geq  5,$ introduced by
	Fernand Lemay.

	After defining selector and selector function, I associate with them
	points and lines in the projective plane, represented by integers and
	give Singer's results which prove the existence of selectors using
	the notion of primitive polynomials, \ref{sec-dprimpol}.

	\ref{sec-tsing3} is a special case of what is needed to determine when
	an irreducible polynomial is a primitive polynomial.
	\footnote{Baumert, p. 101}
	\ref{sec-tsing4} gives a form of the primitive polynomial and the
	generator, so chosen that the polynomials whose coefficient
	define the homogeneous coordinates of points and lines satisfy the same
	4 term recurrence relation.

	\sssec{Definition.}
	Given a power $q=p^k$ of a prime $p$, a {\em selector} or
	{\em difference set} is a subset of $q+1$ distinct integers, such that
	their $q(q+1)$ differences modulo $n := q^2+q+1$ are all of the
	integers from 1 to $q^2+q$.

	\sssec{Example. [Singer.]}
	The following are selectors  with $q = p^k:$ \\
	For $p = 2:$  0, 1, 3, modulo 7.\\
	For $p = 3:$  0, 1, 3, 9, modulo 13.\\
	For $q = 2^2:$  0, 1, 4, 14, 16, modulo 21.\\
	For $p = 5:$  0, 1, 3, 8, 12, 18, modulo 31.\\
	For $p = 7:$  0, 1, 3, 13, 32, 36, 43, 52, modulo 57.\\
	For $q = 2^3:$  0, 1, 3, 7, 15, 31, 36, 54, 63, modulo 73.\\
	For $q = 3^2:$  0, 1, 3, 9, 27, 49, 56, 61, 77, 81, modulo 91.\\
	For $q = 11:$ 0, 0, 1, 3, 12, 20, 34, 38, 81, 88, 94,104, 109
	modulo 133.

	\sssec{Theorem.}
	{\em If $s_i,$ $i = 0$ to $p$ is a selector then, for any $j$,}
	\enumb
	\item{\em $s'_i	= a + k s_{i+j},$ is also a selector.}
	\enume

	The indices are computed modulo $q+1$ and the selector numbers,
	modulo $n.$

	\sssec{Definition.}
	If $a = 1$ and $k = -1,$ the selector $s'_i := 1 - s_i$ is called
	the {\em complementary selector} or {\em co-selector} of $s_i.$ 
	The selectors obtained using $k = 2,$ $\frac{1}{2},$ are called
	respectively {\em bi-selector}, {\em semi-selector}.

	\sssec{Example.}
	\enumb
	\item For $q = 4,$ other selectors are 10, 12, 17, 18, 21 and
		0, 1, 6, 8, 18.
	\item For $p = 7,$ if\\
	the selector is		 0, 1, 7,24,36,38,49,54,\\
	then\\
	the co-selector is	 0, 1, 4, 9,20,22,34,51,\\
	the bi-selector is	 0, 1, 5,27,34,37,43,45,\\
	the semi-selector is	 0, 1, 9,11,14,35,39,51.
	\enume


	\sssec{Definition.}
	The {\em selector function} $f$ is the function from ${\Bbb Z}_n$ to
	${\Bbb Z}_n$ \\
	\hth	$f(s_j-s_i) = s_i,$ $i\neq j.$

	\sssec{Theorem.}
	$f(j-i) -i = f(i-j) - j.$

	\sssec{Example.}
	For $p = 3,$ and $n = 13$, the selector function associated with
	the selector {0,1,3,9} is\\
	$\begin{array}{crrrrrrrrrrrr}
	i  & 1& 2& 3& 4& 5& 6& 7& 8& 9&10&11&12\\
	f(i)&0& 1& 0&-4&-4& 3&-4& 1& 0& 3& 3& 1
	\end{array}$

	\sssec{Definition.}
	{\em Given a selector, {\em points} in the projective plane associated
	with $q = p^k$, with $n = q^2+q+1$ elements are integers in
	${\Bbb Z}_n,$ and {\em lines} are integers in ${\Bbb Z}_n$
	followed by $ ^*$, with the incidence defined by\\ 	
	$a$ is on $b^*$ iff $f(a+b) = 0$ or $a+b = 0.$}

	\sssec{Theorem.}
	\enumb
	\item$ 	a \times b = (f(b-a)-a)^*.$ 
	\item 	$a^* \times b^* = f(b-a) - a.$
	\item{\em $a$ on $b^* \Rightarrow b$ on $a^*.$}
	\enume

	The Statements immediately reflect the duality in projective geometry.

	\sssec{Example.}
	For $p = 3,$ and the selector {0,1,3,9}, the lines and the points on 
	them are\\
	$\begin{array}{ccrrrrrrrrrrrrr}
	lines  &\vline&0^*&1^*&2^*&3^*&4^*&5^*&6^*&7^*&8^*&9^*&
		10^*&11^*&12^*\\
	\hline
	points&\vline& 0&12&11&10& 9& 8& 7& 6& 5& 4& 3& 2& 1\\
	on    &\vline& 1& 0&12&11&10& 9& 8& 7& 6& 5& 4& 3& 2\\
	^*    &\vline& 3& 2& 1& 0&12&11&10& 9& 8& 7& 6& 5& 4\\
	      &\vline& 9& 8& 7& 6& 5& 4& 3& 2& 1& 0&12&11&10\\
	\end{array}$

	\sssec{Theorem.}
	{\em In the projective geometry associated with	$q = p^k$
	and the selector $\{s_0,$ $s_1,$  \ldots  , $s_q\}$\\
	\hth$	i + s_0,$ $i + s_1,$  \ldots, $i + s_p$\\
	are the $q+1$ points on the line $-i^*,$ the addition being done modulo
	$n = q^2+q+1.$}

	\sssec{Definition. [Singer]}\label{sec-dprimpol}
	$P$ is a {\em primitive polynomial} in the Galois Field $GF(p^k)$ iff
	$P$ is of degree $k$ and $I^{p^{k-1}}$ is the smallest power of $I,$
	modulo $P$, which is identical to 1.

	\sssec{Example.}
	\enumb
	\item 	$I^3$	+ I + 1 = 0 is primitive in $GF(2^3).$ 
	\enume
	With 2 $\equiv$  0 modulo 2, we have, modulo P, $I^3$ = I + 1,
	$I^4$ = $I^2$ + I,\\
	$I^5$	= $I^2$ + I + 1, $I^6$ = $I^2$ + 1, $I^7$ = 1.

	It is well known that

	\sssec{Theorem.}
	{\em A primitive polynomial always exists.}

	\sssec{Theorem.}
	$P$ is a {\em primitive polynomial of degree} $m$ over the
	{\em Galois field} $GF(q),$  iff  $P$ is an irreducible polynomial of
	degree $m$ over $GF(q)$ and for a given primitive root $\rho$ of
	$GF(q^m),$ $P(\rho) = 0.$

	\sssec{Theorem. [Singer]}\label{sec-tsing3}
	{\em For each value of $q = p^k$, a selector can be obtained by choosing
	a primitive polynomial of degree 3 over $GF(q).$
	It is, with 0, the set of exponents of $I$ such that the coefficient
	of $I^2$ is 0.}

	\sssec{Example.}
	For $p = 3,$ $P = I^3 - I + 1,$
	$I^3 = I - 1,$ $I^4 = I^2 - I,$ $I^5 = I^2 - I + 1,$
	$I^6 = I^2 + I + 1,$ $I^7 = I^2 - I - 1,$ $I^8 = I^2 + 1,$
	$I^9 = I + 1,$ $I^{10} = I^2 + I,$
	$I^{11} = I^2 + I - 1,$ $I^{12} = I^2 - 1,$ $I^{13} = 1.$\\
	Therefore the selector is 0, 1, 3, 9.

	\sssec{Theorem.}\label{sec-tsing4}
	{\em Let the primitive polynomial be $P_3 := I^3 + b I - c$ and the
	generator be $G := I + g,$\\
	let $g' := 3g^2 + b,$ $h' = g^3 + bg + c,$ $h = \frac{1}{h'}$,\\
	let $J^{(n)} := h^{-n} G^{-n+1} * G^{-n},$ then}
	\enumb
	\item 0. $G^2 = I^2 + 2g I + g^2,$\\
	 1.	$G^3 = 3g I^2 + (3g^2-b) I + (g^3 + c),$\\
	 2.	$G^{-1} = h I^2 - gh I + (g^2+b)h,$ \\
	 3.	$G^{-2} = g'h^2 I^2 + (1-g'gh)h I + (-2gh + g'(g^2+b)h^2).$ 
	\item 0. $G^{n+3} = 3g G^{n+2} - g' G^{n+1} + h' G^n,$ \\
	 1.	$G^n = h(g' G^{n+1} - 3g G^{n+2} + G^{n+3}).$
	\item 0. $J^{(0)} = I^2,$ \\
	 1.	$J^{(1)} = g I^2 + I,$\\
	 2.	$J^{(2)} = (g^2-b) I^2 + 2g I + 1.$
	\item 0. $J^{(n+3)} = 3g J^{(n+2)} - g' J^{(n+1)} + h'J^{(n)},$\\
	 1.	$J^{(n)} = h(g' J^{(n+1)} - 3g J^{(n+2)} + J^{(n+3)}).$
	\enume

	In other words the 4 term recurrence relation is the same for the
	points associated to $G^n$ (1.) as for the lines associated to
	$J^{(n)}$ (3.).

	Proof:  0.0. is immediate.\\
	$G^3 = (I+g)^3,$ or because of $P_3$ we get 0.1.  Eliminating 1,
	$I$ and $I^2$ from
	$G^1,$	$G^2$ and $G^3$ gives $G^3 = 3g G^2 - g' G + h'.$
	Multiplying by $G^n$ gives 1.0 hence 1.1.\\
	From this recurrence relation is it easy to get 0.2. and 0.3.
	$J^{(n)} := h^{-n} G^{1-n} * G^{-n},$ this gives easily 2.0., 2.1, 2.2.  We
	should be careful not to scale.\\
	The definition of $J^{(n)}$ implies\\
	$J^{(n+3)} = h^{-n-3} G^{-n-2} * G^{-n-3}\\
	\hth	= h^{-n-2} G^{-n-2} * (g' G^{-n-2} - 3g G^{-n-1} + G^{-n})\\
	\hth	= 3g J^{(n+2)} - h^{-n-2} G^{-n} * G^{-n-2}\\
	\hth	= 3g J^{(n+2)} - h^{-n-1} G^{-n} \times (g' G^{-n-1} - 
		3g G^{-n} + G^{-n+1})\\
	\hth	= 3g J^{(n+2)} - g' J^{(n+1)} + h' J^{(n)} .$

	\sssec{Example.}
	$p = 5,$ $g = -2,$ $b = -1,$ $c = 2,$ $g' = 1,$ $h' = 1,$ $h = 1,$\\
	$P_3 = I^3 - I - 2,$ $G^3 = - G^2 - G + 1.$
	$J^{(3)} = - J^{(2)} - J^{(1)} + 1.$\\
	$\begin{array}{rcccrrcccrrccc}
	i &\vline&G^i	&J^{(i)}&\vline\vline&i &\vline&G^i	&J^{(i)}
		&\vline\vline&i &\vline&G^i	&J^{(i)}\\
	\cline{1-3}
	-2&\vline&(1,3,2)&[3,3,1]&\vline\vline& 9&\vline&(2,4,3)&[0,4,2]
		&\vline\vline&20&\vline&(4,2,4)&[3,2,4]\\
	-1&\vline&(1,2,3)&[4,2,1]&\vline\vline&10&\vline&(0,2,3)&[3,2,0]
		&\vline\vline&21&\vline&(4,4,0)&[4,4,4]\\
	 0&\vline&(0,0,1)&[1,0,0]&\vline\vline&11&\vline&(2,4,4)&[1,4,2]
		&\vline\vline&22&\vline&(1,1,3)&[4,1,1]\\
	 1&\vline&(0,1,3)&[3,1,0]&\vline\vline&12&\vline&(0,3,1)&[1,3,0]
		&\vline\vline&23&\vline&(4,2,1)&[0,2,4]\\
	 2&\vline&(1,1,4)&[0,1,1]&\vline\vline&13&\vline&(3,0,3)&[1,0,3]
		&\vline\vline&24&\vline&(4,1,1)&[0,1,4]\\
	 3&\vline&(4,3,4)&[3,3,4]&\vline\vline&14&\vline&(4,1,0)&[4,1,4]
		&\vline\vline&25&\vline&(3,3,1)&[4,3,3]\\
	 4&\vline&(0,2,0)&[0,2,0]&\vline\vline&15&\vline&(3,2,3)&[1,2,3]
		&\vline\vline&26&\vline&(2,3,4)&[1,3,2]\\
	 5&\vline&(2,1,0)&[2,1,2]&\vline\vline&16&\vline&(1,2,0)&[1,2,1]
		&\vline\vline&27&\vline&(4,0,1)&[0,0,4]\\
	 6&\vline&(2,0,4)&[1,0,2]&\vline\vline&17&\vline&(0,2,2)&[2,2,0]
		&\vline\vline&28&\vline&(2,0,1)&[3,0,2]\\
	 7&\vline&(1,1,1)&[2,1,1]&\vline\vline&18&\vline&(2,3,1)&[3,3,2]
		&\vline\vline&29&\vline&(1,3,2)&[3,3,1]\\
	 8&\vline&(4,0,0)&[4,0,4]&\vline\vline&19&\vline&(4,2,2)&[1,2,4]
		&\vline\vline&30&\vline&(1,2,3)&[4,2,1]\\
	\end{array}$
	The selector is {0,1,4,10,12,17}.  Line 1$^*$ is incident to
	points -1=30, 0,3,9,11 and 16.

{\tiny\footnotetext[1]{G15.TEX [MPAP], \today}}

\setcounter{section}{4}
	\section{Trigonometry and Spherical Trigonometry.}

	\ssec{Aryabatha I (476-?).}

        The first known table of trigonometric functions corresponds\\
	\hth$crd(\alpha ) = 2 sin(\frac{1}{2}\alpha )$\\
        and to $\alpha  = 0$ to $90^o$ step $15'$, using two sexadesimal places.
        for instance, $crd(36^o) = 2 sin(18^o) = ;37, 4,55.$  (See \ldots).\\
        The trigonometric functions were first defined as ratios of the sides
        of a triangle by Rh\"{a}ticus, who constructed 10 place tables for
        $sin,$ $cos,$ $tan,$ $cot,$ $sec$ and $cosec,$ in increments of $10'',$
        and 15 place tables for $sin,$ with first second and third difference.
        They were edited by Piticus.

	\ssec{Jean Henri Lambert (1728-1777).}
	Lambert gives, in 1770 (I, 190-191), the values of the
	trigonometric function sine for arguments in units $\frac{\pi}{60}$.\\
	These require $s3 = \sqrt{3},$ $s2 = \frac{\sqrt{2}}{2},$
	$s5 = \sqrt{5},$ $s5p = \sqrt{5+s5},$ $s5m = \sqrt{5-s5}.$\\
	His table can then be rewritten as follows:
	$sin(1)= \frac{-s3\:s5p + s5p + s2\:s3\:s5 + s2\:s5 - s2\:s3 - s2}{8},\\
	sin(2)  = \frac{2 s2\:s3\:s5m - s5 - 1}{8},\\
	sin(3)  = \frac{s5\: s2 + s2 - s5m}{4},\\
	sin(4)  = \frac{2 s2\: s5p - s3\: s5 + s3}{8},\\
	sin(5)  = \frac{s3\: s2 - s2}{2},\\
	sin(6)  = \frac{s5 - 1}{4},\\
	sin(7)  = \frac{s3\:s5m + s5m - s2\:s3\:s5 + s2\:s5 - s2\:s3 + s2}{8},\\
	sin(8)  = \frac{- 2 s2\: s5m + s3\: s5 + s3}{8},\\
	sin(9)  = \frac{- s5\: s2 + s2 + s5p}{4},\\
	sin(10) = \frac{1}{2}
	sin(11) = \frac{s3\:s5p - s5p + s2\:s3\:s5 + s2\:s5 - s2\:s3 - s2}{8},\\
	sin(12) = \frac{1}{2}s2\: s5m,\\
	sin(13)= \frac{-s3\:s5m + s5m + s2\:s3\:s5 + s2\:s5 + s2\:s3 + s2}{8},\\
	sin(14) = \frac{2 s2\: s3\: s5p - s5 + 1}{8},\\
	sin(15) = s2,\\
	sin(16) = \frac{2 s2\: s5p + s3\: s5 - s3}{8},\\
	sin(17) = \frac{s3\:s5m + s5m + s2\:s3\:s5 - s2\:s5 + s2\:s3 - s2}{8},\\
	sin(18) = \frac{s5 + 1}{4},\\
	sin(19) = \frac{s3\:s5p + s5p - s2\:s3\:s5 + s2\:s5 + s2\:s3 - s2}{8},\\
	sin(20) = \frac{3s}{2},\\
	sin(21) = \frac{s5\:s2 - s2 + s5p}{4},\\
	sin(22) = \frac{2 s2\:s3\:s5m + s5 + 1}{8},\\
	sin(23) = \frac{s3\:s5m - s5m + s2\:s3\:s5 + s2\:s5 + s2\:s3 + s2}{8},\\
	sin(24) = \frac{1}{2}s2\: s5p,\\
	sin(25) = \frac{s3\: s2 + s2}{2},\\
	sin(26) = \frac{2 s2\: s3\: s5p + s5 - 1}{8},\\
	sin(27) = \frac{s5\: s2 + s2 + s5m}{4},\\
	sin(28) = \frac{2 s2\: s5m + s3\: s5 + s3}{8},\\
	sin(29) = \frac{s3\:s5p + s5p + s2\:s3\:s5 - s2\:s5 - s2\:s3 + s2}{8},\\
	sin(30) = 1.$

	These tables are given here, because they can be used in the case of
	finite fields for appropriate values of p.

	\ssec{Menelaus  of Alexandria (about 100 A. D.)}

        The first appearance of a spherical triangle is in book I of
        Menelaus' treatise Sphaerica, known through its translation into Arabic.
        In it appears the first time a study of spherical triangles and of
        the formula for a spherical triangle $ABC$ with points $L,M,N$ on
	the sides corresponding to IV.\ldots ?\\
	\hth$sin(AN) sin(BL) sin(CM) = - sin(NB) sin(LC) sin(MA).$

	\ssec{al-Battani, or Albategnius (850?-929?).}
        The law of cosine for a spherical triangle was given by al-Battani,
        it will be generalized to finite non-Euclidean geometry in
        IV.\ldots 2.0.\\
        The formula, for a spherical right triangle, called Geber's Theorem,
        will be generalized in IV \ldots 1.1.


\setcounter{section}{6}
\setcounter{subsection}{1}
	\setcounter{subsubsection}{-1}
	\subsubsection{Introduction.}
	This section uses extensively, material learned
	from Professor George Lema\^{\i}tre, in his class on Analytical
	Mechanics, given to first year students in Engineering and in
	Mathematics and Physics, University of Louvain, Belgium, 1942.
	We first determine the differential equation for the pendulum 6.1.3.
	using the Theorem of Toricelli 6.1.1. , we then define the elliptic
	integral of the first kind and the elliptic functions of Jacobi
	6.1.5., we then derive the Landen transformation which relates
	elliptic functions with different parameters 6.1.10., use it to obtain
	the Theorem of Gauss which determines the complete elliptic integrals
	of the first kind from the arithmetico-geometric mean of its 2
	parameters 6.1.14. and obtain the addition formulas for the these
	functions 6.1.16. using the Theorem of Jacobi on pendular motions which
	differ by their initial condition 6.1.7.  We also derive the Theorem
	of Poncelet on the existence of infinitely many polynomials inscribed
	in one conic and circumscribed to another 6.1.9.  We state, without
	proof, the results on the imaginary period of the elliptic functions of
	Jacobi 6.1.19. and 6.1.20.  A Theorem of Lagrange is then given which
	relates identities for spherical trigonometry and those for elliptic
	function 6.1.23.  Finally we state the definitions and some results
	on the theta functions.
	Using this approach, the algebra is considerably simplified by using
	geometrical and mechanical considerations.

	\subsubsection{Theorem. [Toricelli]}
	{\em  If a mass moves in a uniform gravitational field
	its velocity $v$ is related to its height $h$ by
	\begin{enumerate}
	\setcounter{enumi}{-1}
	\item 	$v = \sqrt{2g(h_0-h)},$
	\end{enumerate}
	where $g$ is the gravitational constant and $h_0$ is a constant,
	corresponding to the height at which the velocity would be 0.}

	Proof:  The laws of Newtonian mechanics laws imply the conservation of
	energy.
	In this case the total energy is the sum of the kinetic energy
	$\frac{1}{2}$ $mv^2$
	and the potential energy $mgh,$ therefore\\
	\hth$	\frac{1}{2} mv^2 + mgh = mgh_0,$ for some $h_0.$

	\subsubsection{Definition.}
	A {\em circulatory pendular motion} is the motion of a mass $m$
	restricted to stay on a vertical frictionless circular track, whose
	total energy allows the mass to reach with positive velocity the
	highest point on the circle.  An {\em oscilatory pendular motion} is one
	for which the total energy is such that the highest point on the
	circle is not reached.  The mass in this case oscillates back and forth.
	The following Theorem gives the equation satisfied by a pendular
	motion.

	\subsubsection{Theorem.}
	{\em If a mass $m$ moves on a vertical circle of radius $R,$ with
	lowest point $A,$ highest point $B$ and center $O,$ its position $M$
	at time $t,$
	can be defined by $2\phi(t) = \angle(AOM)$ which satisfies
	\begin{enumerate}
	\setcounter{enumi}{-1}
	\item 	$D\phi = \sqrt{a^2-c^2sin^2\phi},$
	where
	\item 	$a^2 := 2gh_0\frac{1}{4R^2},$ $c^2 = \frac{g}{R},$
	for some $h_0.$
	\end{enumerate}}

	Proof:  If the height is measured from $A,$\\
	\hth$	h(t) = R - R cos(2\phi (t)) = 2R sin^2\phi (t),$\\
	the Theorem of Toricelli gives\\
	\hth$	R D(2\phi )(t) = v(t) = \sqrt{2gh_0 - 4gR sin^2\phi (t)},$\\
	hence 0.  The motion is circulatory if $h_0 > 2R$ or $a > c,$ it is
	oscilatory if $0 < h_0 , 2R$ or $c > a.$

	\subsubsection{Notation.}
	\begin{enumerate}
	\setcounter{enumi}{-1}
	\item 	$k := \frac{c}{a},$ $b^2 := a^2- c^2,$ $k' := \frac{b}{a},$
	\end{enumerate}

	\subsubsection{Definition.}
	If $a = 1,$ and we express $t$ in terms of $\phi(t),$
	\begin{enumerate}
	\setcounter{enumi}{-1}
	\item 	$t = \int_0^{\phi(t)}\frac{1}{\sqrt{1 - k^2sin^2}}.$
	The integral 0. is called the {\em incomplete elliptic integral of the
	first kind}. Its inverse function $\phi$  is usually noted
	\item 	$am(t),$ the {\em amplitude function},\\
	The functions
	\item 	$sn := sin \circ am,$ $cn := cos \circ am,$
	$dn := \sqrt{1 - k^2sn^2},$ \\
	are called the {\em elliptic functions of Jacobi}.
	\item 	$K := \int_0^{\frac{1}{2}\pi} \frac{1}{\sqrt{1 - k^2sin^2}}.$
	is called the {\em complete integral of the first kind},  it gives
	half the period, $\frac{K}{a},$ for the circular pendulum.
	The functions which generalize $tan,$ $cosec,$  \ldots are
	\item 	$ns := \frac{1}{sn}, nc := \frac{1}{cn}, nd := \frac{1}{dn},$
	\item	$sc := \frac{sn}{cn}, cd := \frac{cn}{dn}, ds := \frac{dn}{sn},$
	\item	$cs := \frac{cn}{dn}, dc := \frac{dn}{cn}, sd := \frac{sn}{dn}.$
	\end{enumerate}

	\subsubsection{Theorem.}
	{\em If
	\begin{enumerate}
	\setcounter{enumi}{-1}
	\item   $s_1 := sn(t_1),$ $c_1 = cn(t_1),$ $d_1 = dn(t_1)$ and
	\item   $s_2 := sn(t_2),$ $c_2 = cn(t_2),$ $d_2 = dn(t_2),$ \\
	we have
	\item 	$sn^2 + cn^2 = 1,$ $dn^2 + k^2 sn^2 = 1,$
	$dn^2 - k^2cn^2 = k'^2.$
	\item 	$1 - k^2s^2_1 s^2_2 = c^2_1 + s^2_1 d^2_2
	= c^2_2 + s^2_2 d^2_1.$
	\end{enumerate}}

	\subsubsection{Theorem. [Jacobi]}
	{\em Let $M(t)$ describes a pendular motion,
	Given the circle $\gamma$  which has the line $r$ at height $h_0$ as
	radical axis
	and is tangent to $AM(t_0),$ if $N(t)M(t)$ remains tangent to that
	circle,
	then $N(t)$ also describes a pendular motion, with $N(t_0) = A.$}

	Proof:  With the abbreviation $M = M(t),$ $N = N(t),$ let $N M$ meets
	$r$ at $D,$
	let $M',$ $N'$ be the projections of $M$ and $N$ on $r,$ let $T$ be
	the point of tangency of $MN$ with $\gamma$ ,
	\begin{enumerate}
	\setcounter{enumi}{-1}
	\item 	$DM \: DN = DT^2,$\\
	therefore
	\item $\frac{DT}{ND} = \frac{DM}{DT} = \frac{DT-DM}{ND-DT}
		= \frac{MT}{NT} = \sqrt{\frac{DT}{ND}\frac{DM}{DT}}
		= \sqrt{\frac{DM}{ND}} = \sqrt{\frac{M'M}{N'N}}$\\
	When $t$ is replaced by $t+\epsilon$,\\
	\item $	\frac{v_M}{v_N}
		= lim \frac{M(t+\epsilon ) - M(t)}{N(t+\epsilon) - N(t)}
		= lim \frac{M(t) T}{N(t+\epsilon) T} = \frac{MT}{NT},$
	
	because the triangles $T,M,M(t+\epsilon)$ and $T,N,N(t+\epsilon)$ are
	similar, because $\angle(T,N,N(t+\epsilon) = \angle(T,M(t+\epsilon),M)$
	as well as $\angle(M(t+\epsilon),T,M) = \angle(N(t+\epsilon),T,N).$\\
	Therefore
	\item $	\frac{v_M}{v_N} = \sqrt{\frac{M'M}{N'N}}.$
	\end{enumerate}
	The Theorem of Toricelli asserts that $v_M = \sqrt{2g M'M},$ this
	implies, as we have just seen, $v_N = \sqrt{2g N'N},$ therefore $N$
	describes the same pendular motion with a difference in the origin of
	the independent variable.

	\subsubsection{Corollary.}
	{\em If $M = B$ and $N = A,$ the line $M(t)\times N(t)$ passes through a
	fixed point $L$ on the vertical through $O$ called {\em point of
	Landen}.\\
	Moreover, if $b := BL$ and $a := LA,$ we have\\
	\hth$\frac{v_M}{v_N} = \frac{b}{a}$ and $h_0 = \frac{a^2}{a - b}.$}

	This follows at once from from 6.1.7.2. and 6.1.7.1.

	\subsubsection{Theorem. [Poncelet}
	{\em Given 2 conics $\theta$  and $\gamma$ , if a polygon $P_i,$ $i = 0$
	to $n,$ $P_n = P_0,$ is such that $P_i$ is on $\theta$ and
	$P_i\times P_{i+1}$ is tangent to $\gamma$ , then there exists
	infinitely many such polygons.\\
	Any such polygon is obtained by choosing $Q_0$ on $\theta$ drawing a
	tangent $Q_0Q_1$ to $\gamma$, with $Q_1$ on $\theta$ and successively
	$Q_i,$ such that $Q_i$ is on $\theta$ and $Q_{i-1}\times Q_i$ is
	tangent to $\gamma$, the Theorem asserts that $Q_n = Q_0.$}

	The proof follows at once from 6.1.7. after using projections which
	transform the circle $\theta$  and the circle $\gamma$  into the given
	conics.
	The Theorem is satisfied if the circle have 2 points in common or not.

	\subsubsection{Theorem.}
	{\em If $M(t)$ describes a circular pendular motion, then the
	mid-point $M_1(t)$ of $M(t)$ and $M(t+K)$ describes also a circular
	pendular motion.  More precisely, $M_1(t)$ is on a cicle with diameter
	$LO,$ with $LA = a,$ $LB = b,$ and if $\phi_1(t) = \angle(O,L,M_1(t),$
	\begin{enumerate}
	\setcounter{enumi}{-1}
	\item 	$t = \int_0^{\phi(t)} \frac{D\phi}{\Delta}
		= \frac{1}{2} \int_0^{\phi_1(t)}\frac{D\phi_1}{\Delta_1}.$\\
	where
	\item 	$\Delta^2 := a^2 cos^2\phi  + b^2 sin^2\phi$  and
		$\Delta^2_1 := a^2_1 cos^2\phi _1 + b^2_1 sin^2\phi _1, $\\
	where the relation between $\phi$ and $\phi_1$ is given by
	\item 	$tan(\phi_1-\phi ) = k' tan\phi$, or
	\item 	$sin(2\phi - \phi_1) = k_1 sin\phi_1,$\\
	with
	\item 	$a_1 := \frac{1}{2}(a+b),$ $b_1 := \sqrt{ab},$
		$c_1 := \frac{1}{2}(a-b),$ therefore
	\item 	$a = a_1 + c_1,$ $b = a_1-c_1,$ $c = 2\sqrt{a_1c_1}. $
	\end{enumerate}}

	Proof:  First, it follows from the Theorem of Toricelli that the
	velocity $v_A$ at $A$ and $v_B$ at $B$ satisfy\\
	\hth$v_A = \sqrt{2gh_0} = 2R a,$
		$v_B = \sqrt{2gh_0-2R} = \sqrt{4R^2a^2 - 4 c^2R^2} = 2R b,$\\
	therefore $\frac{BL}{LA} = \frac{b}{a}.$\\
	If $P$ is the projection of $L$ on $BM$ and $Q$ the projection of
	$L$ on $AM,$\\
	\hth$LM^2 = LP^2 + LQ^2 = a^2 cos^2\phi + b^2 sin^2\phi = \Delta ^2.$\\
	\hth$LQ = LM cos(\phi_1-\phi) = a cos \phi.$

	We can proceed algebraically. Differentiating 2. gives\\
	$a(1 + tan^2(\phi_1-\phi)) (D\phi_1 - D\phi) = b(1+tan^2\phi) D\phi,$ or
	$a(1 + tan^2(\phi_1-\phi)) D\phi_1
		= (a(1+tan^2(\phi_1-\phi) + b(1+tan^2\phi)) D\phi\\
	\hth    = (a+b + \frac{b^2}{a} tan^2\phi + b tan^2\phi) D\phi\\
	\hth    = (a+b)(1 + \frac{b}{a} tan^2\phi) D\phi\\
	\hth    = (a+b)(1 + tan\phi tan(\phi_1-\phi)) D\phi,$\\
	or\\
	\hth$\frac{a}{cos^2(\phi_1-\phi)} D\phi_1
	= 2a_1 \frac{cos(2\phi -\phi_1)}{cos\phi cos(\phi_1-\phi)} D\phi,$ or\\
	\hth$	\frac{D\phi}{\frac{a cos\phi}{cos(\phi_1-\phi)}}
	= \frac{D\phi_1}{2a_1 cos(2\phi -\phi_1)},$\\
	or because $LM = \Delta$ \\
	\hth$	\frac{D\phi}{\Delta} = \frac{D\phi_1}{2 \Delta_1}.$

	We can also proceed using kinematics.\\
	The velocity at $M$ is\\
	\hth$	v_M = 2R D\phi = 2R \Delta ,$\\
	If we project the velocity vector on a perpendicualr to $LM,$\\
	\hth$LM D\phi_1 = v_M cos(2\phi_1-\phi)
		= 2R cos(2\phi_1-\phi) \Delta \phi.$\\
	Therefore\\
	\hth$\frac{D\phi}{\Delta} = \frac{D\phi_1}{ 2R cos(2\phi_1-\phi )}
	= \frac{a_1}{2R} \frac{D\phi_1}{\Delta_1} = \frac{D\phi_1}{2\Delta_1}.$

	\subsubsection{Definition.}
	The transformation from $\phi$ to $\phi_1$ is called the {\em forward
	Landen transformation}.  The transformation from $\phi_1$ to
	$\phi$  is called the {\em backward Landen transformation}.

	\subsubsection{Comment.}
	The formulas 3. and 1. are the formulas which are used to
	compute $t$ from $\phi(t)$.  The formulas 4. and 2. are used to compute
	$\phi(t)$ from $t$.

	\subsubsection{Theorem. [Gauss]}
	{\em Given $a_0 > b_0 > 0,$ let
	\begin{enumerate}
	\setcounter{enumi}{-1}
	\item 	$a_{i+1} := \frac{1}{2} (a_i+b_i),$
	\item 	$b_{i+1} := \sqrt{a_ib_i},$
	\end{enumerate}
	then the sequence $a_i$ and $b_i$ have a common limit $a_{\infty }.$
	The sequence $a_i$ is monotonically decreasing and the sequence $b_i$ is
	monotonically increasing.}

	Proof:  Because\\
	\hth$a_i > a_{1+1},$ $b_{i+1} > b_i,$\\
	it follows that the sequence $a_i$ is bounded below by $b_0,$ the
	sequence $b_i$ is bounded above by $a_0,$ therefore both have a limit
	$a_{\infty }$ and $b_{\infty }.$ Taking the limit of 0. gives at once
	$a_{\infty } = b_{\infty }.$

	\subsubsection{Theorem.}
	{\em For the complete integrals we have}
	\begin{enumerate}
	\setcounter{enumi}{-1}
	\item 	$\frac{K}{a} = \int_0^{\frac{1}{2}\pi}
		 \frac{1}{\sqrt{a^2cos^2 + b^2sin^2}}
		= \frac{\frac{\pi}{2}}{a_{\infty}}.$

	Proof:  If $\phi(K) = \frac{\pi}{2}$, then $\phi_1(K) = \pi$,
	therefore
	\item $K = \int_0^{\frac{\pi}{2}}\frac{D\phi}{\Delta}
		= \int_0^{\pi} \frac{D\phi_1}{2 \Delta_1}
		= \frac{1}{2} \int_0^{\frac{\pi}{2}} \frac{D\phi_1}{\Delta_1}
	+ \frac{1}{2} \int_{\frac{\pi}{2}}^{\pi} \frac{D\phi_1}{\Delta_1}
		= \int_0^{\pi} \frac{D\phi_1}{\Delta_1}$
	\hth$	= \int_0^{\frac{\pi}{2}} \frac{D\phi_n}{\Delta_n}
		= \int_0^{\frac{\pi}{2}} \frac{1}{a_{\infty }}
		= \frac{\frac{\pi}{2}}{a_{\infty }}.$
	\end{enumerate}

	\subsubsection{Lemma.}
	\begin{enumerate}
	\setcounter{enumi}{-1}
	\item 	$c_2 = c_1 cn(t_1+t_2) + d_2s_1 sn(t_1+t_2),$
	\item 	$d_2 = d_1 dn(t_1+t_2) + k^2 s_1c_1 sn(t_1+t_2).$

	Proof:  We use the Theorem 6.1.7. of Jacobi.  Let $R$ be the radius of
	$\theta$ and $O$ its center, let $r$ be the radius of $\gamma$ and $O'$
	its center, let $s := OO'.$  Let $A,$ $N,$ $M',$ $M$ be the position of
	the mass at time 0, $t_1,$ $t_2,$ $t_1+t_2.$\\
	The lines $A \times M'$ and $N \times M$ are tangent to the same circle
	$\gamma$ at $T'$ and $T.$\\
	Let $X$ be the intersection of $O \times M$ and $O' \times T,$
	$2\phi := \angle(A,O,N),$
	\item $2\phi' := \angle(A,O,M),$\\
	we have $\angle(N,O,M) = 2(\phi'-\phi ),$ $\angle(M,X,T) = \phi'-\phi,$
	$\angle(T,O',O) = \phi'+\phi$.\\
	If we project $MOO'$ on $O'T,$\\
	\hth$	r = R cos(\phi'-\phi ) s cos(\phi'+\phi),$ or
	\item $r = (R+s) cos\phi  cos\phi' + (R-s) sin\phi sin\phi'.$\\
	\hth	$\phi = am t_1,$ $\phi' = am(t_1+t_2),$ \\
	\hth	$sin\phi' = sn(t_1+t_2),$ $cos\phi' = cn(t_1+t_2),$ \\
	\hth	$sin\phi = sn\: t_1 = s_1,$\\
	\hth	$cos\phi = cn\: t_1 = c_1,$\\
	\end{enumerate}
	when $t_1$ = 0,\\
	\hth$cos(\angle(A,B,M') = cn\:t_2 = c_2 = \frac{BM'}{AB}
		= \frac{O'T'}{AO'} = \frac{r}{R + s}$,\\
	the ratio of the velocities is\\
	\hth$\frac{v_{M'}}{v_A} = \frac{dn \:t_2}{dn\: 0} = d_2
		= \frac{TM'}{AT} = \frac{O'B}{AO'} = \frac{R - s}{R + s},$
	substituting in 2. gives 0.\\
	The proof of 1. is left as an exercise.

	\subsubsection{Theorem. [Jacobi]}
	\begin{enumerate}
	\setcounter{enumi}{-1}
	\item $\frac{sn\: u_1 cn\:  u_2 dn\: u_2 + sn\: u_2 cn\: u_1 dn\:u_1}
		{sn(u_1+u_2)}= 1 - k^2 sn^2 u_1 sn^2 u_2.$
	\item $\frac{cn\:u_1 cn\:u_2 - sn\:u_1 dn\:u_1 sn\:u_2 dn\:u_2}
		{cn(u_1+u_2) } = 1 - k^2 sn^2 u_1 sn^2 u_2. $
	\item $\frac{dn\:u_1 dn\:u_2 - k^2 sn\:u_1 sn\:u_2 cn\:u_1 cn\:u_2}
		{dn(u_1+u_2) } = 1- k^2 sn^2 u_1 sn^2 u_2. $

	Proof: Let $w = \frac{1}{1 - k^2 s^2_1 s^2_2}.$\\
	Let $s_1,$ $s_2,$  $\ldots$  denote $sn\:u_1,$ $sn\:u_2,$  $\ldots$,
	define $S$ and $C$ such that\\
	\hth$sn(u_1+u_2) = S w,$ $cn(u_1+u_2) = C w.$\\
	The 6.1.15.0. gives\\
	\hth$c_2 = c_1 C w + d_2 s_1 S w$ or
	\item $c_1 C w = - d_2 s_1 S w + c_2,$ \\
	6.1.6.2. gives\\
	\hth$S^2 w^2 + C^2 w^2 = 1,$\\
	eliminating $C$ gives the second degree equation in $S w$:\\
	\hth$(c^2_1 + d^2_2 s^2_1 (S w)^2 - 2 s_1 c_2 d_2 (S w) + c^2_2 - c^2_1
		 = 0,$\\
	one quarter of the discriminant is\\
	\hth$s^2_1 c^2_2 d^2_2 - (c^2_2 - c^2_1)(c^2_1 + d^2_2 s^2_1)\\
	\hth	= s^2_1 c^2_2 d^2_2 - c^2_1 c^2_2 + c^4_1 -s^2_1 c^2_2 d^2_2
		+ s^2_1 c^2_1 d^2_2\\
	\hth	= c^2_1(c^2_1 - c^2_2 + s^2_1 d^2_2) = c^2_1 s^2_2 d^2_1,$\\
	therefore\\
	\hth$	S w = (s_1 c_2 d_2 \pm  c_1 d_1 s_2)w.$\\
	One sign correspond to one tangent from $M$ to $\gamma$ , the other to
	the other
	tangent, therefore one corresponds to the addition, the other to the
	subtration formula.  From the special case $k = 0,$ follows that, by
	continuity, the $+$ sign should be used.  This gives 0., 1. follows from
	3, 2. is left as an exercise.
	\end{enumerate}

	\subsubsection{Corollary.}
	{\em
	\begin{enumerate}
	\setcounter{enumi}{-1}
	\item $	sn(u+K)  = cd(u),$ $cn(u+K) = -k' sd(u),$ $dn(u+K)  = k' nd(u).$
	\item $	sn(u+2K) = -sn(u),$ $cn(u+2K) = - cn(u),$ $dn(u+2K) = dn(u).$
	\item $	sn(u+4K) = sn(u),$ $cn(u+4K) = cn(u),$ $dn(u+4K) = dn(u).$
	\end{enumerate}}

	\subsubsection{Definition.}
	\hth$K'(k^2) = K(k'^2).$

	\subsubsection{Theorem.}{\em
	\begin{enumerate}
	\setcounter{enumi}{-1}
	\item 0. $k sn \circ I+iK' = sn,$\\
	 1.	 $i k cn \circ I+iK' = ds,$\\
	 2.	 $i dn \circ I+iK' = cs,$
	\item 0. $sn \circ I+2iK' = sn,$\\
	 1.	 $cn \circ I+2iK' = - cn,$\\
	 2.	 $dn \circ I+2iK' = - dn,$
	\end{enumerate}}

	\subsubsection{Theorem.}{\em
	\begin{enumerate}
	\setcounter{enumi}{-1}
	\item 	$sn$ has periods $4K$ and $2iK'$ and pole $\pm iK',$
	\item 	$cn$ has periods $4K$ and $4iK'$ and pole $\pm iK',$
	\item 	$dn$ has periods $2K$ and $4iK'$ and pole $\pm iK'.$
	\end{enumerate}}

	\subsubsection{Theorem.}{\em
	\begin{enumerate}
	\setcounter{enumi}{-1}
	\item $	k = 0 \Rightarrow  sn = sin,$ $cn = cos,$ $dn = \underline{1},$
	\item $	k = 1 \Rightarrow  sn = tanh,$ $cn = sech,$ $dn = sech.$
	\end{enumerate}}

	\subsubsection{Theorem. [Lagrange]}
	{\em From the addition formulas of elliptic functions
	we can derive those for a spherical triangle as follows.  Let}
	\begin{enumerate}
	\setcounter{enumi}{-1}
	\item   $u_1 + u_2 + u_3 = 2K,$\\
	{\em define}
	\item   $sin a := - sn u_1,$ $cos a := - cn u_1,$\\
	    $sin b := - sn u_2,$   $cos b := - cn u_2,$\\
	    $sin c := - sn u_3,$   $cos c := - cn u_3,$\\
	    $sin A := - k sn u_1,$ $cos A := - dn u_1,$\\
	    $sin B := - k sn u_2,$ $cos B := - dn u_2,$\\
	    $sin C := - k sn u_3,$ $cos C := - dn u_3,$\\
	{\em then to any formula for elliptic functions of $u_1,$ $u_2,$ $u_3,$
	corresponds
	a formula for a spherical triangle with angles $A,$ $B,$ $C$ and sides
	$a,$ $b,$ $c.$  For instance,}
	\item   $\frac{sin A}{sin a} = \frac{sin B}{sin b}
		= \frac{sin C}{sin c} = k.$
	\item  $  cos a = cos b cos c + sin b sin c cos A,$
	\item  $  cos A = - cos B cos C + sin B sin C cos a,$
	\item  $  sin B cot A = cos c cos B + sin c cot a.$

	Proof.  2. follows from the definition.  3. follows from
	\hth$c_2 = c_1 cn(t_1+t_2) + d_2 s_1 sn(t_1+t_2)$
	after interchanging $t_1$ and $t_2$ and using
	\item 0. $sn(t_1+t_2) = sn(2K-t_1-t_2) = sn\: t_3 = s_3,$ \\
	 1.	$cn(t_1+t_2) = - cn(2K-t_1-t_2) = - cn\: t_3 = - c_3,$\\
	 2.	$dn(t_1+t_2) = dn(2K-t_1-t_2) = dn \:t_3 = d_3,$\\
	similarly, 4. follows from\\
	\hth$c_2	= c_1 cn(t_1+t_2) + d_2 s_1 sn(t_1+t_2)$ \\
	after interchanging $t_1$ and $t_2$ and using 6 and 5. from\\
	\hth$sn\:t_2 dn\:t_1 = cn\:t_1 sn(t_1+t_2)
		- sn\:t_1 dn\:t_2 cn(t_1+t_2)$\\
	after division by $sn\:t_1.$
	\end{enumerate}

	\subsubsection{Definition.}
	Given the parameter q, called the nome,
	\begin{enumerate}
	\setcounter{enumi}{-1}
	\item $	q := e^{-\pi \frac{K'}{K}},$\\
	the functions
	\item $\theta_1 := 2 q^{\frac{1}{4}} \sum_{n=0}^{\infty}
		(-1)^n q^{n(n+1)} sin(2n+1)I$
	\item $\theta_2 := 2 q^{\frac{1}{4}} \sum_{n=0}^{\infty}
		q^{n(n+1)} cos(2n+1)I$
	\item $\theta_3 := 1 + 2 \sum_{n=1}^{\infty} q^{n^2 } cos 2nI$
	\item $\theta_4 := 1 + 2 \sum_{n=1}^{\infty} (-1)^n q^{n^2 } cos 2nI$
	are the {\em theta functions of Jacobi}.
	\end{enumerate}

	\subsubsection{Definition.}
	The functions, with $v = \pi  \frac{I}{2K)}$
	\begin{enumerate}
	\setcounter{enumi}{-1}
	\item 	$\theta_s := \frac{2K \theta_1 \circ v}{D \theta_1(0)},$
	$\theta_c := \frac{\theta_2 \circ v}{\theta_2(0)},$
	$\theta_d := \frac{\theta_3 \circ v}{\theta_3(0)},$
	$\theta_n := \frac{\theta_4 \circ v}{\theta_4(0)},$\\
	are called the {\em theta functions of Neville}.
	\end{enumerate}

	\subsubsection{Theorem.}
	{\em If $p,q$ denote any of $s,$ $c,$ $d,$ $n,$\\
	\hth$pq = \frac{\theta_p}{\theta_q}.$
	For instance\\
	\hth$sn = \frac{\theta_s}{\theta_n }
		= \frac{2K \theta_1 \circ v}{D\theta_1(0)}
		. \frac{\theta_4(0)}{\theta_4 \circ v}.$}

	\subsubsection{Theorem.}
	{\em The Landen transformation replaces the parameter $q,$ by $q^2.$}

	\subsubsection{References.}
	Jacobi, Fundamenta Nava Theoriae Funktionum Ellipticarum, 1829.\\
	Legendre, Trait\'{e} des fonctions elliptiques et des int\'{e}grales
	elliptiques. III, 1828.\\
	Gauss, Ostwald Klasiker?\\
	Landen John, Phil. Trans. 1771, 308.\\
	Abel, Oeuvres, 591.\\
	Bartky, Numerical Calculation of generalized complete integrals,
	    Rev. of Modern Physics, 1938, Vol. 10, 264-\\
	Lemaitre G. Calcul des integrales elliptiques, Bull. Ac. Roy. Belge,
	    Classe des Sciences, Vol. 33, 1947, 200-211.\\
	Fettis, Math. of Comp., 1965.\\
	Appell, Cours de Mechanique,

	\subsubsection{Notes}

	\hth$(Dy)^2	= C_0 (y^2 - A_0) (y^2 + B_0),$\\
	\hth$(Dz)^2	= C_1 (z^2 - A_1) (z^2 + B_1),$\\
	\hth$z = d (y + \frac{1}{y}),$ $l\neq 0,$ $d > 0.$

	The equations are compatible iff (l in the beginning of next expres.?)\\
	$d^2(1-	\frac{2}{y^2}) C_0(y^2 + A_0) (y^2 + B_0)
		= C_1(d^2(\frac{y^2+l}{y})^2 + A_1)(d^2(\frac{y^2+l}{y})^2
			 + A_1)$\\
	this requires $\sqrt{l}$ to be a root of one of the factore of the
	second member, let it be the second factor, this implies\\
	\hth$d^2 4 l  + B_1 = 0,$\\
	then, the second factor becomes,\\
	\hth$d^2(\frac{y^2+l}{y})^2 + B_1
		= d^2((\frac{y^2+l}{y})^2  - 4 K)
		= d^2(\frac{y^2-l}{y})^2$\\
	therefore $\sqrt{l}$ is a double root of the second memeber and\\
	\hth$C_0(y^4 + (A_0+B_0) y^2 + A_0B_0)
		= d^2 C_1(y^4 + (2 l + A_1)y^2 + l^2),$
	therefore\\
	\hth$C_0 = d^2C_1,$ $A_0B_0 = \frac{ B_1^2}{16 d^4},
		 A_0+B_0 =  \frac{A_1 - \frac{1}{2} B_1}{d^2},$\\
	For real transformations, $A_0B_0 > 0,$ if $j_0 = sign(B_0)$ and
	$j_1 = sgn(B_1),$\\
	\hth$B_1 = 4j_0d^2\sqrt{A_0B_0},$
		$A_1 = d^2(A_0+B_0+2j_1\sqrt{A_0B_0}$\\
	\hti{11}$= j_0(\sqrt{|A_0|} + j_0j_1\sqrt{|B_0|})^2.$\\
	If we want $A_1$ $B_1$ >0 then $j_0$ = $j_1.$

{\tiny\footnotetext[1]{G16.TEX [MPAP], \today}}

\setcounter{section}{5}
	\section{Algebra, Modular Arithmetic.}

	\setcounter{subsection}{-1}
	\ssec{Introduction.}

	Geometry can be handled synthetically, with little or no reference
	to algebra.  But it was discovered little by little that an underlying
	algebraic structure lurks behind geometry.  If we deal with a geometry
	with a finite number of points on each line, we have to deal with an
	underlying algebraic structure which involves a finite number of
	integers.  Such structure presented itself in connection with
	application of mathematics to astronomy (and astrology), in studying
	the relative motion of sun and moon and the relative motion of the
	planets, mainly Jupiter.
	If the smallest unit of time used is $t,$ the period of the sun around
	the earth, is $s.t,$ the position of the sun is the same after 2
	revolutions hence $2s$ is equivalent to $s$ and $2s+1$ is equivalent to
	$s+1$
	as well as 1.  This led to the notion of working modulo $s.$

	\ssec{The integers.}

	\sssec{Definition.}
	$p$ is a {\em prime}  iff  $p$ is an integer larger than 1, which
	is only divisible by 1 and $p.$

	\sssec{Restriction.}
	In the sequel, it is always assumed that \un{$p$ is odd.}

	\ssec{The integers modulo p.}

	\setcounter{subsubsection}{-1}
	\sssec{Introduction.}
	Although much of what I will do can be generalized, to the case of
	powers of primes, I will, for simplicity, restrict myself to the case of
	a prime $p.$

	\sssec{Definition.}
	The {\em integers modulo} $p$ are the integers $x$ satisfying\\
	\hth$	0 \leq  x < p.$

	The set of these integers is denoted ${\Bbb Z}_p.$
	The operations modulo $p$ are defined in terms of the operations on the
	integers as follows:

	\sssec{Definition.}
	\enumb
	\item If $x$ and $y$ are integers modulo $p,$ {\em addition modulo} $p,$
		denoted $+_p$ is defined as the least non negative remainder of
		the division of the integer $x+y$ by $p.$
	\item {\em Multiplication modulo} $p,$ denoted $._p,$ is defined as the
		least non negative remainder of the division of $x.y$ by $p.$
	\item {\em Subtraction modulo} $p,$ denoted $-_p,$ is defined as the
		inverse operation of addition, $c +_p b = a \implies
		a -_p b = c.$
	\item {\em Division modulo} $p,$ denoted $/_p,$ is defined as the
		inverse operation of addition, $c ._p b = a \implies
		a /_p b = c,$ provided $b\neq 0$.
	\enume

	\sssec{Convention.}
	As I will not use simultaneously 2 different primes,
	and as it will usually be clear from the context that the addition,
	multiplication, \ldots, are done modulo $p,$ I will replace
	$+_p$ by $+,$ \ldots.
	An alternate notation, useful when several different moduli are used,
	is to use\\
	\hth$	a + b \equiv  c \pmod{p}.$

	\sssec{Example.}
	We have, $0 +_5 3 = 3,$ $5 +_7 4 = 2,$ $5 +_{11} 6 = 0.$\\
	Modulo 7: $5 + 4 = 2,$ $5 - 4 = 1,$ $5 . 4 = 6,$ $5 / 4 = 3,$
		$5 + 0 = 5.$\\
	Modulo 7: the inverses of 1 through 6 are respectively $1,4,5,2,3,6.$
	Modulo 11: $9 + 5 = 3,$ $9 - 5 = 4,$ $9 . 5 = 1,$ $9 / 5 = 4,$
		$9 . 0 = 0.$

	\sssec{Comment.}
	Addition, subtraction and multiplication are easy to
	perform, moreover hand calculators and languages for microprocessors
	have functions which allow easy computations.  Division requires either
	a table of inverses or the inverses can be obtained, for large primes
	$p,$ using the Euclid-Aryabatha algorithm.
	\footnote
	{To appreciate this contribution of the Hindus, Aryabatha lived at 
	    the end of the fifth Century, while an equivalent algorithm was
	    only developed in the Western World by Bachet de Meziriac in 1624.}
	\footnote{Pulverizing $a$ is meant to convey what we would now express
	by finding the inverse of $a$ modulo $n.$}

	\sssec{Algorithm. [Euclid]}\label{sec-aeuc}
	Let $a \geq b > 0.$  We determine in succession\\
	\hth $a_0 := a,$ $a_1 := b,$ $q_1,$ $a_2,$ $q_2,$  $\ldots  , 
		a_n = 0 \ni \\
	0.\hti{6}a_{j-1} := a_j q_j + a_{j+1},$ $0 \leq  a_{j+1} < a_j.$

	\sssec{Algorithm. (Pulverizer of Aryabatha)}\label{sec-aary}
	Given $q_1,$ $q_2,$  $\ldots$  , $q_{n-1},$ determine\\
	\hth $b_{n-1} := 0, \: b_{n-2} := 1,\\
	0.\hti{6}b_{j-1} := b_j q_j + b_{j+1},$ for $j = n-2, \ldots , 1.$

	\sssec{Algorithm. (Continued fraction algorithm)}
	\label{sec-acfr}
	Given $q_1,$ $q_2,$  $\ldots$  , $q_{n-1},$ determine\\
	\hth $c_0 := 0,$ $c_1 := 1,$ $d_0 := 1,$ $d_1 := 0,$\\
	0.\hti{6}$c_{j+1} := c_j q_j + c_{j-1},$ for $j = 1,\ldots, n-1.$\\
	1.\hti{6}$d_{j+1} := d_j q_j + d_{j-1},$ for $j = 1,\ldots, n-1.$

	\sssec{Algorithm.}\label{sec-a1d11}
	\hth$u_j := c_j^2 + d_j^2.$\\
	\hth$v_j := c_j c_{j+1} + d_j d_{j+1}.$

	\sssec{Example.}\label{sec-e1}
{\small	Let $a = 10672$ and $b = 4147,$ \ref{sec-aeuc}, \ref{sec-aary},
	\ref{sec-acfr} and \ref{sec-a1d11} give\\
	$\begin{array}{rrrrrrrrr}
	 &&& \vline&\uparrow b_0 =& {\mathit 175}\: \vline&\downarrow \:c_0 =
		& {\bf 0} \:d_0 =&{\bf 1}\\
	a_0 =& {\bf 10672} =& {\bf 4147}.2 +& 2378\:     \vline&   
	 \: b_1 =&  {\mathit 68}  \:\vline&\:   c_1 =&   {\bf 1} \:d_1 =& {\bf 0}\\
	a_1 =&  4147 =& 2378.1 +& 1769\:  \vline&  \: b_2 =&  39\: \vline&\:
		 c_2 =& 2  \: d_2 =&   1\\
	a_2 =&  2378 =& 1769.1 +&  609\:  \vline&  \: b_3 =&  29\: \vline&\:
		 c_3 =&   3  \: d_3 =&   1\\
	a_3 =&  1769 =&  609.2 +&  551\:  \vline&  \: b_4 =&  10\: \vline&\:
		 c_4 =&   5  \: d_4 =&   2\\
	a_4 =&   609 =&  551.1 +&   58\:  \vline&  \: b_5 =&   9\: \vline&\:
		 c_5 =&  13  \: d_5 =&   5\\
	a_5 =&   551 =&   58.9 +&   29\:  \vline&  \: b_6 =&   {\bf 1}\: 
		\vline&\: c_6 =&  18  \: d_6 =&   7\\
	a_6 =&    58 =&   29.2 +&    0\:  \vline&\uparrow  \: b_7 =&  {\bf 0}\: 
		\vline&\: c_7 =& {\mathit 175}  \: d_7 =&  {\mathit 68}\\
	a_7 =&    {\mathit 29}\:\:\:\:&&&                                 
	&   \vline&\downarrow c_8 =& {\mathit 368}  \: d_8 =& {\mathit 143}
	\end{array}$\\
	\hth$n = {\mathit 8},$ $4147\cdot 175 - 10673\cdot 68 = 29$.\\
	\hth$u_7 = 35249, \: u_8 = 155873, \: v_7 = 74124.$\\
	The bold-faced number are initial values, the italicized numbers
	are final values.  Notice that all $a$'s have to be computed before the
	$b$'s are computed, but this is not so for the $c$'s and the $d$'s.\\
	For instance, for the line starting with $a_3,$ 1769 and 609 come from
	the preceding line, 2 is the quotient of the division of 1769 by 609
	and 551 is the remainder,\\
	$b_4 = q_5.b_5+b_6 = 1.9 + 1 = 10.\\
	c_4 = c_3.q_3+c_2 = 3.1+2 = 5,$ $d_4 = d_3.q_3+d_2 = 1.1+1 = 2.$\\
	Observe that 175 and 68 are obtained in 2 different ways.
}

	\sssec{Definition.}
	The {\em greatest common divisor} of $a$ and $b$ is the largest positive
	integer which divides $a$ and $b$, it is denoted $(a,b).$

	\sssec{Theorem.}\label{sec-t1d13}
	\enumb
	\item{\em The algorithm} \ref{sec-aeuc} {\em terminates in a finite
		number of steps.}
	\item$(a,b) = (a_0,a_1) = (a_1,a_2) = \ldots = (a_{n-1},a_n) = a_{n-1}.$
	\item$b_j a_{j-1} - b_{j-1} a_j = (-1)^{n-j} (a,b),$
		{\em in particular,} $b_0 b - b_1 a = (-1)^n (a,b).$
	\item$b_j < a_j.$ {\em in particular,} $b_1 < b,$ $b_0 < a.$
	\enume

	\sssec{Theorem.}\label{sec-t2d13}
	\enumb
	\item$a_0/a_1 = q_1 + 1/(q_2 + 1/(q_3 +  \ldots  + 1/q_{n-1})).$
	\item$b<\frac{a}{2}\implies c_i<c_{i+1}$, $i = 0,\ldots,n-1.$
	\item$a_ic_{i+1}+a_{i+1}c_i = a,$ $a_id_{i+1}+a_{i+1}d_i = b.$
	\item$b.c_i\equiv(-1)^{i+1}a_i\pmod{a},$
		$a.d_i\equiv(-1)^ia_i\pmod{b}.$
	\enume

	If $\frac{a}{2}\leq b < a$ then $q_1 = 1$ and $c_2 = c_1$.

	\sssec{Definition.}
	The second member in \ref{sec-t2d13}.0. is called a
	{\em terminating continued fraction}.

	\sssec{Theorem. [Symmetry property]}\label{sec-tcontfrsymm}
	{\em If $(a,b) = 1,$ $b<\frac{a}{2},$ and we repeat the algorithm
	with $a' := a$ and\\ $b' := \pm b^{-1}\pmod{a},$ $b' \leq \frac{a}{2},$
	then this algorithm terminates in the same number $n$ of steps and\\
	\hth$a'_j = c_{n-j},$ $c'_j = a_{n-j},$ $q'_j = q_{n-j}.$\\
	In particular, if $b^2\equiv -1\pmod{a}$ and $b<a$, then $n = 2n'+1$ is
	odd and}\\
	\hth$c_j = a_{n-j},$ $q_j = q_{n-j}$ and $a^2_{n'}+a^2_{n'+1} = a.$

	\sssec{Example.}
	$\begin{array}{rrrrrcrrrrrl}
	i&a_i&q_i&c_i&&\vline&i&a_i&q_i&c_i&\\
	\cline{1-5}\cline{7-11}
	0&378&&0&8&\vline&	0&65&&0&7\\
	1&143&2&1&7&\vline&	1&18&3&1&6\\
	2&82&1&2&6&\vline&	2&11&1&3&5\\
	3&61&1&3&5&\vline&	3&7&1&4&4\\
	4&21&2&5&4&\vline&	4&4&1&7&3\\
	5&19&1&13&3&\vline&	5&3&1&11&2\\
	6&2&9&18&2&\vline&	6&1&3&18&1\\
	7&1&2&175&1&\vline&	7&0&&65&0\\
	\cline{7-11}
	8&0&&368&0&\vline&	&c'_j&q'_j&a'_j&j&18^2+1\equiv 0 \pmod{65}.\\
	\cline{1-5}
	&c'_j&q'_j&a'_j&j&\vline&&&&&&7^2 + 4^2 = 65.
	\end{array}$
 
	\sssec{Theorem. [Euler]}
	{\em Every integer whose prime factors to an odd power are congruent to
	1 modulo 4, can be written as a sum of 2 squares and vice-versa.}

	\sssec{Example.}
	\hth$13 = 2^2 + 3^2,$ $52 = 4^2 + 6^2.$\\
	\hth$585 = 9^2 + 24^2 = 12^2 + 21^2.$

	\ssec{Quadratic Residues and Primitive Roots.}

	\sssec{Definition.}
	$n$ is a {\em quadratic residue} of $p$  iff  there exist a integer
	$x$ such that $x^2$ is congruent to $n$ modulo $p.$  We write, with
	Gauss, $n$ R $p.$  $n$ is a non residue of $p,$ if there are no integer
	whose square is congruent to $n$ modulo $p,$ and we write $n$ N $p.$

	\sssec{Theorem.}
	{\em The product of 2 quadratic residues or of 2 non residues is
	a quadratic residue.  The product of a quadratic residue by a non
	residue is a non residue.}

%
	\sssec{Theorem. [Fermat]}
	{\em If $a$ is not divisible by $p$ then $a^{p-1} \equiv 1 \pmod{p}.$}

	\sssec{Definition.}
	$g$ is a {\em primitive root} of $p$ iff 
	$a^i \equiv 1 \pmod{p}$ and $0<i<p \implies i = p-1.$
	In other words, $p-1$ is the smallest positive power of $g$ which
	is congruent to 1 modulo $p$.

	\sssec{Notation. [Euler]}
	$\phi(n)$ denotes the number of integers betweem 1 and $p,$ relatively 
	prime to $p$.

	\sssec{Theorem. [Gauss]}
	\enumb
	\item {\em There are $\phi(p-1)$ primitive roots of $p$.}
	\item {\em If $g$ is a primitive root of $p$, all primitive roots are
		$g^i$ with $(i,p-1) = 1.$}
	\enume

	\sssec{Example.}
	For $p = 13,$ 2 is a primitive root,\\
	$\begin{array}{ccccccccccccccc}
	i&0&1&2&3&4&5&6&7&8&9& 10& 11& 12& \pmod{12}\\
	\hline
	g^i&	 1& 2&4&-5&3& 6& -1&-2 &-4 &5& -3&-6&1& \pmod{13}
	\end{array}$

	The other primitive roots are $2^5 = 6$, $2^7 = -2$ and $2^{11} = -6.$\\
	The easiest method to obtain all inverses moudo $p$ is to first
	obtain a primitive root and then to use $g^i\:.\:(g^{p-1-i})^{-1} = 1$.

	\sssec{Theorem.}
	{\em If $\delta$  is a primitive root of $p,$ the square root of an
	integer can be unambiguously defined if we chose a particular primitive
	root.}

	It is sufficient to choose $a$ or $a\delta$, with
	$0\leq a<\frac{p-1}{2}.$

	\sssec{Examples.}
	Modulo 5, $\delta^2$ can be chosen equal to 2 or 3, with
	$\delta^2 = 3,$ we have\\
	$\begin{array}{cccccc}
	i	&0& 1& 2& 3& 4\\
	\hline
	\sqrt{i}&0& 1&2\delta& 1\delta&2.
	\end{array}$

	Modulo 7, $\delta^2 = 3$ can be chosen equal to 3 or 5, with
	$\delta^2 = 5,$ we have\\
	$\begin{array}{cccccccc}
	i	&0& 1& 2& 3& 4& 5& 6\\
	\hline
	\sqrt{i}&0& 1& 3&1\delta&2&2\delta& 3\delta
	\end{array}$

	\sssec{Theorem.}
	\enumb
	\item $p \equiv 1 \pmod{4} \Leftrightarrow -1 $R$p$ and $p$ odd.
	\item $p \equiv 1,-1 \pmod{8} \Leftrightarrow 2 $R$p$ and $p$ odd.
	\enume

	\sssec{Theorem.}\label{sec-sqr2rat}
	{\em $\sqrt{2}$ is rational in the field $Z_{17}.$}

	This follows at once from the following figure and the fact
	that the mid-point of the segment joining $(1,0)$ to $(0,1)$
	is $(-8,-8)$ when $p = 17.$ This figure originates with the geometric
	construction corresponding to the proof by the school of Pythagoras
	that there is no rational number whose square is 2.  In fact,
	$\sqrt{2} = \pm 6,$ when $p = 17.$ In the case of real numbers, a 
	corresponding figure corresponds to the geometric interpretation of
	the classical proof of the irrationality of
	$\sqrt{2}$, the squares becoming smaller and smaller. I suggest that
	the reader reflects on this, from a geometric point of view, together
	with the atomic structure of our Universe.

\begin{picture}(300,100)(-60,0)
\put(-40,0){\line(1,1){40}}
\put(-20,0){\line(1,1){20}}
\put(-10,0){\line(1,1){10}}
\put(10,0){\line(-1,-1){10}}
\put(20,0){\line(-1,-1){20}}
\put(40,0){\line(-1,-1){40}}

\put(-40,0){\line(1,-1){40}}
\put(-20,0){\line(1,-1){20}}
\put(-10,0){\line(1,-1){10}}
\put(10,0){\line(-1,1){10}}
\put(20,0){\line(-1,1){20}}
\put(40,0){\line(-1,1){40}}

\put(40,40){\line(0,-1){80}}
\put(20,20){\line(0,-1){40}}
\put(10,10){\line(0,-1){20}}
\put( 5, 5){\line(0,-1){10}}
\put(-40,40){\line(0,-1){80}}
\put(-20,20){\line(0,-1){40}}
\put(-10,10){\line(0,-1){20}}
\put( -5, 5){\line(0,-1){10}}

\put(40,40){\line(-1,0){80}}
\put(20,20){\line(-1,0){40}}
\put(10,10){\line(-1,0){20}}
\put( 5, 5){\line(-1,0){10}}
\put(-40,-40){\line(1,0){80}}
\put(-20,-20){\line(1,0){40}}
\put(-10,-10){\line(1,0){20}}
\put( -5, -5){\line(1,0){10}}

\put(-40,0){\line(1,0){80}}
\put(0,-40){\line(0,1){80}}

\put(40,40){\circle*{2}}
\put(-40,-40){\circle*{2}}
\put(40,-40){\circle*{2}}
\put(-56,-48){\makebox(0,0){(-8,-8)}}
\put(52,48){\makebox(0,0){(8,8)}}
\end{picture}\\[60pt]

\ssec{Non Linear Diophantine Equations and Geometry.}
\setcounter{subsubsection}{-1}
	\sssec{Introduction.}
	There has been, historically, a constant interplay between geometry and 
	diophantine equations, the former suggesting problems of the latter kind
	which also indicate the interest of having problems in geometry solved
	using integers only. As evidence I will give just one such problem
	considered by Euler\footnote{Opera Minora Collecta, II, (1778) 294-301,
	(1779) 362-365, (1782) 488-491}.

	\sssec{Definition.}
	The {\em median problem} consists in constructing a triangle
	with integer sides and medians.

	\sssec{Theorem.}
	{\em If $a_i$ are the length of the sides and $g_i$ are twice the
	length of the medians, then}
	\enumb
	\item $2a_{i+1}^2 + 2a_{i-1}^2 = a_i^2 + g_i^2.$
	\enume

	Proof:\\
	\hth$a_0^2+a_2^2 - 2a_0a_2cos(A_1) = a_1^2,$\\
	\hth$\frac{1}{4}a_0^2+a_2^2 - a_0a_2cos(A_1) = \frac{1}{4}g_0^2,$\\
	eliminating the terms involving the angle gives\\
	\hth$2a_1^2 + 2a_2^2 = a_0^2 + g_0^2.$

	\sssec{Theorem. [Euler]}
	{\em The solution of the preceding problem can be expressed in terms
	of 2 parameters $a$ and $b$, using}\\
	\hth$C = (4ab)^2,$ $D = (9a^2+b^2)(a^2+b^2)$,
		$F = 2(3a^2+b^2)(3a^2-b^2)$,\\
	\hth$a_0 = 2a(D-F),$ $a_1+a_2 = 2a(C+D),$ $a_1-a_2 = 2b(C-D),$\\
	\hth$g_0 = 2b(D+F),$ $\:g_1+g_2 = 6a(C-D),$ $\:g_1-g_2 = 2b(C+D),$

	\sssec{Theorem. [Euler]}\label{sec-teulmed}
	{\em An other solution can be obtained corresponding to $a' = b$ and
	$b' = 3a$ for which\\
	\hth$a'_i = g_i$ and $g'_i = 3a_i$.}

	An example is provided with the pair (1,2) giving the pair (2,3) in
	the Example. In fact we have the following

	\sssec{Theorem.}
	{\em If both $a$ and $b$ are not divisible by 3, then $3|g_i$ and
	3 does not divide $a_0$ therefore the preceding Theorems gives a
	solution and for this solution $b'$ is divisible by 3.}

	Indeed, $a^2 \equiv b^2 \equiv 1,$ $C \equiv F\equiv 1,$
	$D\equiv -1,$ $a_0\equiv -a$, $-a_1\equiv a_2\equiv b,$ $g_i\equiv 3.$
	It is therefore sufficient, if we use Theorem \ref{sec-teulmed}
	to consider all pairs in which one of the integers in the pair is
	divisible by 3. Similarly we have the following Theorem.

	\sssec{Theorem.}
	{\em If $a$ is not divisible by 3 and $b$ is divisible by 9, then
	$9|a_i$, $27|\hti{-1}/a_i$, $27|g_i$ and the solution is the same as
	that obtained from $a' = b/3$ and $b' = a$.}

	Proof:\\
	If $9|b$, then modulo 27, $C\equiv 0,$ $D\equiv 9$ and $F\equiv 18,$
	therefore $9|a(0)$ but $27|\hti{-1}/a(0)$ while $27|g(i)$.

	\sssec{Example.}
	The solutions for the pairs $(a,b) = $(1,3), (1,2), (2,3), (1,6), (3,1),
	(3,5) are given by Euler. Except for the pair (1,2), They are ordered
	by increasing maximum values of $a_i$.\\
	$\begin{array}{ccrrrrrrrrrrr}
	a&\vline&1&1&2&1&5&3& ^*4&7&3&5& ^*3\\
	b&\vline&3&2&3&6&3&1&3&3&5&6&2\\
	\hline
	a_0&\vline&3&158& 68&314&145&477& 184&1099&2547& 2690& 1926\\
	a_1&\vline&1&127& 85&159&207&277& 739& 810&2699& 5277& 3985\\
	a_2&\vline&2&131& 87&325&328&446&1077&1339&2704& 5953& 6101\\
	g_0&\vline&1&204&158&404&529&569&1838&1921&4765&10924&10124\\
	g_1&\vline&5&261&131&619&463&881&1357&2312&4507& 7583& 8123\\
	g_2&\vline&4&255&127&377&142&640&   5&1391&4498& 5893& 1399\\
	\end{array}$\\[10pt]
	$\begin{array}{ccrrrrrrrr}
	a&\vline&3&11&3&8& ^*7& ^*3&10& ^*3\\
	b&\vline&7& 3&4&3&6&11& 3&13\\
	\hline
	a_0&\vline& 8163&12287& 8874&18288&   42&40563& 59820& 75123\\
	a_1&\vline& 5050& 6416&13703&11663&15091& 4232& 32621&  6953\\
	a_2&\vline& 5897& 9897&14671&19105&20567&28531& 51439& 58580\\
	g_0&\vline& 7343&11281&26968&25838&36076& 4301& 61982& 36283\\
	g_1&\vline&13316&21370&20005&35537&24865&70006&106699&134543\\
	g_2&\vline&12227&16921&17827&23999& 5699&50125& 81481& 89174\\
	\end{array}$

	The pair (1,2) corresponds to a degenerate triangle (1+2=3).\\
	The pairs marked with $^*$ are solutions only in a geometry with
	complex coordinates because $a_{i+1}+a_{i-1}<a_i$ for some $i$.
	The other degenerate solutions are obtained by observing that, in
	Euler's proof, other solutions are obtained when $b^2 = a^2$ or $9a^2.$
	
	\ssec{Farey sets and Partial Ordering.}
\setcounter{subsubsection}{-1}
	\sssec{Introduction.}
	The basic idea is the following, a subset $T_1(n)$, $n>0$, of the set
	${\Bbb Z}_p$ can be placed into one to one correspondance with the set
	$H_n$ of irreducible rationals whose numerator, in modulus, and
	denominator are not larger than $n$, provided $2n^2-2n+1 < p$.
	The ordering $\leq$ in $H_n \in {\Bbb Q}$ induces an
	ordering in $T_1(n)$ such that $a\leq b$ and $b\leq c \implies a\leq c$
	and $-b\leq -a$. If, morever, $0\leq a$ then $b^{-1}\leq a^{-1}.$
	If order is to be preserved, when we do one addition or one 
	multiplication, we have to use $T_2(n) := T_1(n')$ instead of $T_1(n)$,
	with $n = 2n'^2.$ This insures that the sum or product of 2
	elements in $T_2(n)$ is in $T_1(n)$. $T_1$, $T_2$ are defined as the
	sets $T_1(n)$, $T_2(n)$ corresponding to the largest $n$.
	This can be repeated for a finite number of additions and 
	multiplications provided $p$ is large enough.\\
	$H_n$ is related to the Farey set $F_n$ which is its subset in [0,1].
	Farey sets have been used, for instance, by my colleague and friend
	Professor R. Sherman Lehman to factor medium sized numbers.\\
	The cardinality of the partially ordered set is estimated in
	\ref{sec-tcardt1}.\\
	The complement (${\Bbb Z}_p-H_n-\{\pm\sqrt{-1}\})$ can be partitioned
	into 4 sets $\epsilon$, $-\epsilon$, $\lambda$ and $-\lambda$
	which might play the role of the sets of smallest elements and the sets
	of largest elements as given in \ref{sec-dcontfreps} to
	\ref{sec-tcontfreps}.
	Given an integer $k$, we can determine the corresponding irreducible
	rationals, or in which of the small or large set $k$ belongs, using
	algorithm \ref{sec-achaint1}, which depends on the symmetry Theorem
	\ref{sec-tcontfrsymm}.
	We end by contrasting with the notion of continuity in the set of real
	numbers.

	\sssec{Definition.}
	A {\em Farey set} $F_n$ is the set of irreducible rationals
	$\frac{a_i}{b_i},$ in ascending order, between 0 and 1, whose
	numerator and denominator do not exceed $n.$\\
	A {\em Haros set} $H_n$ is the set of irreducible rationals
	$\frac{a_i}{b_i},$ in ascending order, between $-n$ and $n$, whose
	numerator, in modulus, and denominator do not exceed $n.$

	\sssec{Theorem. [Haros]}
	{\em If $\frac{a_i}{b_i}$ and $\frac{a_{i+1}}{b_{i+1}}$ are
	any 2 successive rationals of a Farey set $F_n,$ then}
	\enumb
	\item 	$a_{i+1} b_i - a_i b_{i+1} = 1.$
	\item 	The numerators and denominators of 2 successive rationals are
		relatively prime.
	\item 	$\frac{a_i}{b_i} = \frac{a_{i-1} + a_{i+1}}{b_{i-1} + b_{i+1}}.$
	\item   {\em The set $F_n$ can be constructed starting from
	$\frac{0}{1}$ and $\frac{1}{1}$
	    by inserting rationals using formula $2$ while the resulting
	    numerators and denominators of the second member are not larger
	    than $n.$}
	\enume

	For a proof see Hardy and Wright, p. 23 to 26.\\
	The set $H_n$ can be deduced from $F_n$, by multiplicative symmetry
	with respect to 1 and then by additive symmetry with respect to 0.
	It can also be obtained from $\frac{-n}{1}$ and $\frac{n}{1}$
	using formula 2, but reduction is required
	and the termination condition is not as simple as for the set $F_n$.

	\sssec{Definition.}
	A set $S$ is {\em partially ordered} by $\leq$ iff, with $a,b,c \in S$,
	\enumb
	\item $a \leq a$ for all $a$ in $S,$
	\item $a \leq b$ and $b \leq a \implies a = b.$
	\item $a \leq b$ and $b \leq  c \implies a \leq c.$
	\enume

	But, for any 2 distinct elements $a$ and $b$ in $S,$ we need not have
	$a \leq  b$ or\\
	\hth$b \leq a.$

	\sssec{Notation.}
	\hth$a < b$ if $a \leq b$ and $a\neq b.$

	\sssec{Definition.}\label{sec-dchaint0t1}
	We define the {\em set} $T_1(n)$ by:
	\enumb
	\item The set $T_1(n)$ := \{$ \frac{a_i}{b_i},$ $a_i$ and $b_i$
		relatively prime, $| a_i| \leq n,$ $0 < b_i \leq n$ \}.
	\enume

	\sssec{Theorem.}\label{sec-tpartialord}
	{\em If $0<n$ and $2n(n-1) + 1 < p$ or equivalently if}
	$0 < n < \frac{\sqrt{2p-1}+1}{2}$,
	\enumb
	\item {\em there is a bijection between the irreducible rationals in
		$H_n$ and the elements in the subset $T_1(n)\in {\Bbb Z}_p$.}
	\item {\em If the order in $T_1(n)$ is that induced by the order in
		$H_n\in{\Bbb Q}$ and if\\
		$x,y \in T_1(n),$ then}\\
	0. {\em the set $H_n$ is partialy ordered,}\\
	1. $x < y  \Rightarrow   -y < -x,$\\
	2. $0 < x < y  \Rightarrow   0 < 1/y < 1/x.$
	\enume

	Proof:  It is sufficient to prove, that under the given hypothesis,
	if $\frac{a_i}{b_i}$ and $\frac{a_j}{b_j}$ are any 2 distinct elements
	in ${\Bbb Q},$ they correspond to distinct elements of $T_1(n)$ in
	${\Bbb Z}_p.$  Indeed, if $\frac{r}{s} \equiv \frac{t}{u}$ then
	$ru-ts \equiv  0$ modulo $p,$
	but $| ru-ts|  \leq n^2 + (n-1)^2 = 2n(n-1) +1<p,$
	hence, by hypothesis, $\frac{r}{s} = \frac{t}{u}.$
	The bound cannot be improved for $T_1,$ because,
	$\frac{n}{n-1} \equiv -\frac{n-1}{n}$ if $n^2 + (n-1)^2 = p,$
	whose positive root is
	$\frac{\sqrt{2p-1}+1}{2},$ and the sequence of primes of the form
	$\frac{m^2+1}{2}$ is infinite.

	\sssec{Definition.}\label{sec-dchaint2}
	For a given $p,$ let $n_p$ be the largest positive integer such
	that
	\enumb
	\item 	$2n_p(n_p-1)+ 1 < p,$\\
	then\\
	 0. $T_1 := T_1(n_p),$\\
	 1. $T_2 := T_1([\sqrt{\frac{n_p}{2}}]).$
	\enume

	\sssec{Theorem.}\label{sec-tchaint2}
	{\em If $x, y, x', y' \in T_2,$}
	\enumb
	\item $x.y, x+y \in T_1$,
	\item $0 < x',$ $x < y  \Rightarrow x.x' < y.x'.$
	\item $x \leq  y,$ $x' \leq y' \Rightarrow x + x' \leq y + y'.$
	\enume

	Indeed, if $|a|,$ $|b|,$ $|c|,$ $|d| \leq m := [\sqrt{\frac{n_p}{2}}]$
	then $|ad + bc| \leq 2m^2$, $|ac|\leq m^2$, and $|bd| \leq m^2,$
	therefore if
	$x, y, x', y' \in T_2,$ $x + x'$ and $y + y' \in T_1(2m^2)
	= T_1(n_p) = T_1.$\\
	Of course, for multiplication only, we could replace $2m^2$ by $m^2$.

	\sssec{Example.}
	In this, and in other examples, I have chosen as representative of an
	element in ${\Bbb Z}_p$, that which is in modulus less than
	$\frac{p}{2}$.
	\enumb
	\item   For $p = 31,$ $n_{31} = 4,$ $T_1 = T_1(4)$ is\\
		$-4 < -3 < -2 < 14 < 9 < -1 < 7 < -11 < 15 < 10 < -8 < 0\\
		   < 8 < -10 < -15 < 11 < -7 < 1 < -9 < -14 < 2 < 3 < 4.$\\
	Indeed, the Farey set $F_4$ is\\
	    $\frac{0}{1} < \frac{1}{4} < \frac{1}{3} < \frac{1}{2}
		< \frac{2}{3} < \frac{3}{4} < \frac{1}{1},$\\
	the values in ${\Bbb Z}_{31}$ are\\
	\hth$    0   <  8  < -10 < -15 <  11 <  -7 <  1,$\\
	their inverses are\\
	\hth$	   4  >  3  >  2  > -14 >  -9 >  1.$\\
	$T_2$ = $T_1(1)$ is\\
	\hth$	-1  <  0  < 1 .$\\
	For one multiplication we could use $T'_2$ = $T_1(2)$ $(2^2 = 4)$
	which is\\
	\hth$	-2  < -1  < 15  < 0 < -15 < 1 < 2.$

	\item
{\small For $p = 617,$ $n_{617} = 18,$ the positive elements of $T_1$
	are\\
	$\begin{array}{rrrrrrrrrrrr}
	240<&-254<&270<&288<&-44<&95<&-257<&-56<&-185<&-137\\
	<\:\:\:109<&-77<&-41<&88<&190<&103<&-145<&-112<&193<&247\\
	<-132<&-274<&285<&218<&-154<&-82<&-168<&-34<&-176<&-36\\
	<\:\:\:\:\:62<&-237<&116<&206<&-290<&-220<&-224<&-231<&-142<&-171\\
	<-123<&73<&-51<&-264<&39<&69<&-280<&-47<&165<&-181\\
	<-308<&182<&-164<&48<&281<&-68<&-38<&265<&52<&-72\\
	<\:\:\:124<&172<&143<&232<&225<&221<&291<&-205<&-115<&238\\
	<\:\:-61<&37<&177<&35<&169<&83<&155<&-217<&-284<&275\\
	<\:\:\:133<&-246<&-192<&113<&146<&-102<&-189<&89<&42<&78\\
	<-108<&138<&186<&57<&258<&-94<&45<&-287<&-269<&255\\
	<-239<&1\:\:\:\:&\\<-253\:\:\:\:\\
	<\:\:\:271<&289<&-43<&96<&-256<&-55<&-184<&-136<&-76<&-40\\
	<\:\:-87<&191<&104<&-111<&248<&-131<&-273<&286<&-153<&-167\\
	<-175<&63<&-236<&207<&-223<&-230<&-141<&-122<&-50<&-263\\
	<\:\:\:\:\:70<&-279<&-307<&282<&-67<&266<&125<&233<&226<&-204\\
	<\:\:-60<&178<&156<&276<&-245<&-101<&90<&79<&139<&2\\
	<\:\:-75<&-86<&105<&249<&-152<&-174<&208<&-121<&-262<&-306\\
	<\:\:\:267<&126<&-203<&157<&-244<&-100<&3<&250<&-151<&209\\
	<-120<&-305<&127<&-202<&158<&4<&-150<&210<&-304<&-201\\
	<\:\:\:\:\:\:\:\:5<&211<&-303<&-200<&6<&-302<&7<&-301<&8<&-300\\
	<\:\:\:\:\:\:\:\:9<&10<&11<&12<&13<&14<&15<&16<&17<&18.
\end{array}$
}
	\enume

	The positive elements in $T_2 = T_1(3)$ are\\
	\hth$ 206<-308<-205<   1<-307<   2<   3.$

	\sssec{Theorem (Mertens).}
	\hth$\sum_{b=1}^n (\phi (b)) = \frac{3n^2}{\pi^2} + O(n log(n)),$

	{\em where the last notation implies that the error divided by
	$n log(n)$ is bounded as $n$ tends to infinity.}

	\sssec{Theorem.}\label{sec-tcardt1}
	{\em The number of terms in $T_1$ is of the order of\\
	\hth$	\frac{6}{\pi^2} p + O(p^{\frac{1}{2} } log(p)),$\\
	or approximately $0.6079 p.$}

	This follows at once from the fact that the number of irreducible
	rationals with denominator $b$ is $\phi(b),$ from 
	$T_1 = 4 \sum_{b=2}^n(\phi(b))+3$ from $p = 2 n^2 + O(n)$,
	from $\phi(1) = 1$ and from the Theorem of Mertens.\\
	For $p = 31,$ $23 = .74 p,$ for $p = 617,$ $405 = .656 p.$

	The following Theorem gives a method to determine if a given integer in
	${\Bbb Z}_p$ is in $T_1.$

	\sssec{Algorithm. [Modified continued fraction] }
	\label{sec-achaint1}
	Given $a_0 := p,$ let $n := n_p,$ $0 < a_1 := a < \frac{p}{2},$
	$c_0 := 0,$ $d_0 = 1,$ $c_1 := 1,$ $d_1 := 0,$ $i := 1$\\
	l:\hti{7}$q_i := a_{i-1} / a_i,$ $a_{i+1} = a_{i-1} - a_i q_i,$\\
	\hti{16}$c_{i+1} = c_{i-1} + c_i q_i,$ $d_{i+1} = d_{i-1} + d_i q_i,$\\
	\hti{16}$\un{if}\: a_{i+1} \geq  c_{i+1} \:\un{then}\:
		\un{begin}\:i := i+1; \un{goto}\:l\:\un{end},$\\
	\hth$\un{if}\: a_i<n \:\un{then}\: a \equiv
		\frac{(-1)^i a_{i+1}}{c_{i+1}} \pmod{p} \in T_1,$ \\
	\hth$\un{if}\: c_{i+1}<n \:\un{then}\:a \equiv
		\frac{(-1)^ic_i}{a_i} \pmod{p} \in T_1,$ \\
	\hth$\un{if}\: a_ia_{i+1}>c_ic_{i+1} \:\un{then}\: a \in
		(-1)^{i+1} \lambda,$\\
	\hth$\un{if}\: a_ia_{i+1}<c_ic_{i+1} \:\un{then}\: a \in
		(-1)^i \epsilon,$\\
	\hth$\un{if}\: a_ia_{i+1}=c_ic_{i+1} \:\un{then}\: a
		= - 1/a$.

	$i$ is therefore the largest index for which $a_i \geq c_i$.\\
	We observe that if we start with $a' := a$ and $b' := \pm b^{-1}
	\pmod{a}$, the $+$ sign is to be chosen when $n$ is even, and that by
	the symmetry property, when the algorithm stops, $c'_j\geq a'_j$,
	$c'_{j+1}<a'_{j+1}$, therefore $j = i+1$ and we have consistent 
	conditions.

	\sssec{Example.}
	For $p = 31,$ $n_{31} = 4,$
	\enumb
	\item   if $a = 14,$ the continued fraction algorithm gives\\
	$\begin{array}{rrrrrcrrrrrcrrrrr}
	i&a_i&q_i&c_i&&\vline&i&a_i&q_i&c_i&&\vline
	&i&a_i&q_i&c_i&\\
	\cline{1-17}
	0&31&&0&5&\vline& 0&31&&0&6&\vline&0&31&&0&3\\
	1&14&2&1&4&\vline& 1&12&2&1&5&\vline&\un{1}&\un{6}&5&\un{1}&2\\
	\un{2}&\un{3}&4&\un{2}&\un{3}&\vline&2&7&1&2&4&\vline&2&\un{1}&6&\un{5}
		&\un{1}\\
	3&2&1&9&2&\vline& \un{3}&\un{5}&1&\un{3}&\un{3} &\vline&3&0&&31&0\\
	4&1&2&11&1&\vline& 4&\un{2}&2&\un{5}&2&\vline\\
	5&0&&31&0&\vline& 5&1&2&13&1&\vline\\
	&&&&&\vline& 6&0&&31&0&\vline\\
	\cline{1-17}
	&c'_j&q'_j&a'_j&j&\vline&&c'_j&q'_j&a'_j&j&\vline&&c'_j&q'_j&a'_j&j
	\end{array}$
	\item For $a = 14,$ $i = 2,$ $14.2 - 31.1 = -3,$ $| -3|  \leq  4,$
	$14 \equiv  -\frac{3}{2} \pmod{31},$ which is in $T_1.$
	\item For $a = 12,$ $i = 3,$ $12.3 - 31.1 = 5 > 4,$
	$14 \equiv \frac{5}{3} \pmod{31},$ which is not in $T_1.$ But
	$12\in -\epsilon$ and $13\in -\lambda$.
	\item For $a = 6,$ $i = 1,$ $6.1 - 31.1 = -25,$ $| -25| >  4,$
	$6 \equiv -\frac{6}{1} \pmod{31},$ which is in $T_1.$ But
	$6\in \lambda$ and $5\in -\epsilon$.
	\item In conclusion, $\lambda  = \{ 6, -13\},$
		$\epsilon  = \{ -5,-12\}.$
	\enume

	\sssec{Example.}
	For $p = 617,$ the elements in $\lambda$  and, below them, their
	inverse in $\epsilon$, are given below.  Those in $-\lambda$  and
	$-\epsilon$ are obtained by replacing $x$ by $-x.$\\
{\small
	$\begin{array}{rrrrrrrrrrrrrr}
	\lambda:& 19& 20& 21& 23& 24& 25& 26& 28& 29& 30& 32& 53& 58\\
	\epsilon:&65&216&$-$235&161&180&$-$74&$-$261&$-$22&$-$234&144&135&163
		&$-$117\\
	\lambda:&64&80&85&91&92&$-$98&$-$99&106&107&$-$118&$-$119&128&129\\
	\epsilon:&$-$241&54&$-$196&278&$-$114&$-$170&$-$268&$-$227&173&183
		&$-$140&188&$-$110\\
	\lambda:&159&160&162&187&$-$195&$-$197&$-$198&$-$199&213&214&215&
		$-$228&$-$229\\
	\epsilon:&260&27&$-$179&33&$-$212&$-$166&$-$134&31&$-$84&$-$222&$-$66
		&46&$-$97\\
	\lambda:&$-$242&$-$243&251&252&259&$-$272&277&$-$293&$-$294&$-$296
		&$-$297&$-$298&$-$299\\
	\epsilon:&283&$-$292&59&$-$71&81&93&$-$49&219&149&$-$148&$-$295&147
		&130\\
	\end{array}$
}

	\sssec{Definition.}\label{sec-dcontfreps}
	Let $x \in \hti{-2}|\:\:\:\:T_1,$ let $a_i$ and $c_i$ be defined as in
	\ref{sec-achaint1}, with $b$ replaced by $x$ and let
	$a_{i+1}$ and $c_{i+1}$ be the next pair, let $a'_i$ and $c'_i,$
	$a'_{i+1}$ and $c'_{i+1},$ be the corresponding quadruple for
	$b' := \pm x^{-1},$ the sign so chosen that $b' < a/2$, if
	\enumb
	\item $c'_{i+1} = a_i,$ $a'_{i+1} = c_i,$ $c'_i = a_{i+1},$
		$a'_i = c_{i+1},$

	then
	\item 0. $ a_i a_{i+1} < c_i c_{i+1}$ and $i$ even $\Rightarrow
		x \in \epsilon$,\\
	 1. $ a_i a_{i+1} < c_i c_{i+1}$ and $i$ odd $\Rightarrow
		x \in -\epsilon$,
	\item 0. $ a_i a_{i+1} > c_i c_{i+1}$ and $i+1$ even $\Rightarrow
		x \in \lambda$,\\
	 1. $ a_i a_{i+1} > c_i c_{i+1}$ and $i+1$ odd $\Rightarrow
		x \in -\lambda$.

	and we have the partial ordering of these sets by $<<$,\\
	\hth$	-\lambda  << -\epsilon  << 0 << \epsilon  << \lambda .$
	\enume

	\sssec{Theorem.}\label{sec-tcontfreps}
	\enumb
	\item $x \in \epsilon \Rightarrow -x \in -\epsilon,$
		$1/x \in \lambda,$ $-1/x \in -\lambda.$
	\item {\em For a given $p$, all integers in the set $[0,p-1]$ are
		either in the set $T_1(n_p)$ or in one of the sets $\epsilon$,
		$-\epsilon,$ $\lambda,$ or $-\lambda,$ with the exception of
		$\pm \sqrt{-1}$ when $p\equiv -1 \pmod{4}.$}
	\enume

	We leave the proof as an exercise.

	\sssec{Theorem.}
	{\em If for all $\epsilon \in  (0,\epsilon_1)\:\exists\:\delta(\epsilon)
	> 0 \ni\\
	x_0-\delta(\epsilon) < x < x_0+\delta(\epsilon) \Rightarrow
	| f(x) - f(x_0) | < \epsilon$, then $f$ is continuous at $x_0.$}

	Indeed for the continuity criterium, we can choose $\delta(\epsilon)
	= \delta(\epsilon_1)$ for $\epsilon \geq \epsilon_1.$

	\sssec{Comment.}
	The preceding Theorem is implicit in most text.
	In the older texts, it is alluded to by adding in the definition of
	continuity the phrase "however small is $\epsilon$".  If we choose
	$\epsilon_1 = 10^{-100}$, say, and assume that for a given $f$ and
	$x_0$, the hypothesis of the preceding Theorem is satisfied, it follows
	that the continuity at $x_0$ depends only on the value of the function
	in the interval $(x_0-10^{-100},x_0+10^{-100}).$  If we now try to
	give an example from the world we live in, no meaning can be given to
	physical objects which have distances from each other less than
	$\epsilon_1.$
	The definition of continuity gives therefore problems of interpretation
	in Atomic Physics.  The same is true is Cosmology when the distances
	are of the order of the dimension of the Universe. Continuity requires
	the notion of ordered set. We need to apply the more general concept 
	of partialy ordered set, to allow for a criterium which test values
	which are small, but not too small, or large but not too large.  This
	is what is achieved using Farey sets.

	\ssec{Complex and quaternion integers.}

\setcounter{subsubsection}{-1}
	\sssec{Introduction.}
	Hamilton introduced the notion of quaternions, to try to
	generalize the notion of complex number for application to
	3 dimensional geometry.\\
	The elements are of the form $a + b {\bf i} + c{\bf j} + d{\bf k},$ 
	with $a,$ $b,$ $c,$ $d$ not all zero, with \\
	\hth${\bf j}.{\bf k} = -{\bf k}.{\bf j} = {\bf i},
	{\bf k}.{\bf i} = - {\bf i}.{\bf k} = {\bf j},
	{\bf i}.{\bf j} = - {\bf j}.{\bf i} = {\bf k},
	{\bf i}.{\bf i} = {\bf j}.{\bf j} = {\bf k}.{\bf k} = 0,$\\
	and real numbers commute with {\bf i}, {\bf j}, {\bf k}, addition
	of quaternions is commutative and is distributive over
	multiplication.

	\sssec{Definition.}
	Given a prime $p$ and a non quadratic residue $d,$ the set
	of {\em complex integers} ${\Bbb C}_p$ is the set\\
	\hth$a + b \delta,$ $a,b \in {\Bbb Z}_p,$ $\delta^2 = d.$ \\
	The operations are those of {\em addition},\\
	\hth$(a_0 + b_0 \delta ) + (a_1 + b_1 \delta ) = a_2 + b_2 \delta ,$\\
	where $a_2 := a_0 + a_1 \pmod{p},$ $b_2 := b_0 + b_1 \pmod{p}.$\\
	and of {\em multiplication},\\
	\hth$(a_0 + b_0 \delta ) . (a_1 + b_1 \delta ) = a_3 + b_3 \delta ,$\\
	where $a_3 := a_0.a_1 + b_0.b_1.d \pmod{p},$
	$b_3 := a_0.b_1 + a_1.b_0 \pmod{p}.$\\
	This is entirely similar to the introduction of complex numbers,\\
	\hth$	\delta^2 = d$ replacing $i^2 = -1.$

	\sssec{Example.}
	For $p = 5$ and $d = 2,$\\
	\hth$(1 + \delta) + (1 + 3 \delta) = 2 + 4 \delta,$\\
	\hth$(1 + \delta) . (1 + 3 \delta) = 2 + 4 \delta.$

	\sssec{Definition.}
	A {\em quaternion integer} is a quaternion with coefficients in
	${\Bbb Z}_p$.

	\sssec{Theorem.}
	{\em If $p\equiv 1,3\pmod{8}$, the quaternion integers are isomorphic
	to 2 by 2 matrices over ${\Bbb Z}_p$.}

	The isomorphisms is deduced from the correspondance\\
	\hth$ {\bf 1} \sim \ma{1}{0}{0}{1},$ ${\bf i} \sim \ma{1}{b}{b}{-1},$
	${\bf j} \sim \ma{-b}{1}{1}{b},$ ${\bf k} \sim \ma{0}{1}{-1}{0},$\\
	with $b^2 = -2$.  For instance, for $p = 11,$ $b = 3$, for $p = 17,$
	$b = 7$.

	\sssec{Theorem.}
	\enumb
	\item {\em The quaternions form a skew field (or division ring).}
	\item {\em The quaternion integers form a 
	non commutative ring with unity for which if a right inverse exists
	then it is also a left inverse.}
	\enume

	\ssec{Loops.}
	\sssec{Definition.}
	A {\em loop} $(L,+)$ is a non empty set of elements $L$ together with
	a binary operation " + " such that, if $l_1,$ $l_2$, $l_3$,
	are elements in $L$,
	\enumb
	\item   $l_1+l_2$ is a {\em well defined} element of $L.$
	\item   There exists a {\em neutral element} $e \in L,$ such that
		$e \:+\: l_1 = l_1 \:+\: e = l_1.$
	\item$l_1 \:+\: x = l_2$ has a unique solution $x \in L$, denoted
	$x = l_1\: \vdash \:l_2$ (or $x = l_1\: \setminus \: l_2,$ for $(L,.)$),
	\item$y + l_1 = l_2$ has a unique solution $y \in L$ denoted
	$y = l_2\: \dashv \:l_1,$ (or $y = l_2\: / \:l_1,$  for $(L,.)$)
	\enume

	\ssec{Groups.}

	\sssec{Definition.}
	A {\em group} $(G,.)$ is a non empty set of elements $G$
	together with an operation . such that
	\enumb
	\item   If $g_1$ and $g_2$ are any elements of $G,$ $g_1.g_2$ is a
	 {\em well defined} element of $G.$
	\item   The operation is {\em associative}, or for any elements
	$g_1,$ $g_2,$ $g_3$ of $G,$\\
	\hth$    (g_1.g_2).g_3 = g_1.(g_2.g_3).$
	\item   There exists a {\em neutral element} $e$ in $G,$ such that for
	 all elements\\
	\hth$    g \in G,$ $e.g = g.e = g.$
	\item   Every element $g$ of $G$ has an {\em inverse,} written
	$g^{-1},$ such that\\
	\hth$    g.g^{-1} = g^{-1}.g = e.$
	\enume

	\sssec{Notation.}
	If the operation is noted $+$ instead of $.$, the
	neutral element is called a zero and is noted 0.

	\sssec{Comment.}
	$(G,+)$ or $(G,.)$ is often abbreviated as $G,$ if the operation
	is clear from the context.

	\sssec{Theorem.}
	{\em In a group, the neutral element is unique and in element
	has only one inverse.}

	\sssec{Definition.}
	A group $(G,+)$ is {\em abelian} or {\em commutative}  iff
	for every element $g_1$ and $g_2$ of $G,$\\
	\hth$g_1.g_2 = g_2.g_1.$

	\sssec{Notation.}
	In a group $(G,+),$ we define
	\enumb
	\item $g = e,$ $1.g = g,$ $(n+1).g = n.g + g$ and $(-n).g = -(n.g)$
	where $n$ is any positive integer.
	\enume
	In a group $(G,.),$ we use instead of $0.g,$ $1.g$ and $n.g,$
	$ g^0$ $g^1$ and $g^n,$ where $n$ is any positive or negative integer.

	\sssec{Definition.}
	A {\em cyclic group} $(G,.)$ is a group for which there exist an
	element $g,$ called a generator of the group such that every
	element if $G$ is of the form $g^n.$ ($n.g$ if the operation is $+$).

	\sssec{Examples.}
	\enumb
	\item $({\Bbb Z},+)$ is a cyclic group, $1$ and $-1$ are generators.
	\item $({\Bbb Z}_p,+),$ $p$ prime, is a cyclic group, every element
	different from 0 is a generator.
	\item $({\Bbb Z}_n,+),$ $n$ composite, is an abelian group which is not
	cyclic.
	\item $({\Bbb Z}_p-\{0\},.),$ $p$ prime, is a cyclic group, any
	primitive root is a generator.
	\enume

	\ssec{Veblen-Wederburn system.}
	\sssec{Definition.}
	A {\em Veblem-Wederburn system} ($\Sigma,+,\cdot$), is a set $\Sigma$,
	containing at least the elements 0 and 1 which is such that for
	$a,$ $b$, $c \in \Sigma$,
	\enumb
	\item$\Sigma$ is closed under the binary operations " + " and
	$\:"\cdot"\:$,
	\item($\Sigma$,+) {\em is an abelian group,}
	\item($\Sigma-\{0\},\cdot)$ {\em is a loop,}
	\item $(a + b) \cdot c = a \cdot c \:+\:b\cdot c$,
	\item $a \cdot 0 = 0,$
	\item {\em is right distributive,}\\
	\hth$(a+b)\cdot c = a \cdot c + b\cdot c$.
	\item$a \neq b \Rightarrow x \cdot a = x \cdot b + c$ {\em has a unique
	solution.}
	\enume

	\sssec{Definition.}
	A {\em division ring} is a Veblen-Wederburn system which is left 
	distributive.

	\sssec{Definition.}
	A {\em alternative division ring} is a division ring for which for all
	elements $a \neq 0$, the right inverse $a^R$ and left inverse $a^L$ are
	equal, so that we can write it as $a^{-1}$ and such that for all $b$
	in the set\\
	\hth$(a\cdot b)\cdot b^{-1} = b^{-1}\cdot (b\cdot a) = a.$

	\sssec{Theorem.}
	{\em In an alternative division ring},\\
	\hth$(b\cdot a)\cdot a = b\cdot a^2,$ $a\cdot(a\cdot b) = a^2\cdot b.$

	\sssec{Definition.}
	The {\em Cayley numbers} or {\em octaves} consist of
	(${\bf p + \ov{q}e},+,\cdot$), with (see also Stevenson p. 379)
	\enumb
	\item {\bf p} and {\bf q} are quaternions over the reals,
	\item $({\bf p + \ov{q}e}) + ({\bf p' + \ov{q'}e}) =
		{\bf (p+p') + \ov{(q+q')}e},$
	\item $({\bf p + \ov{q}e}) \cdot ({\bf p' + \ov{q'}e})
		= {\bf (pp' - \ov{q'}q + \ov{q'p + q\ov{p'}}e},$
	\enume

	\sssec{Comment.}
	With ${\bf l}$ and ${\bf l'}$ denoting {\bf i} or {\bf j} or {\bf k},\\
	${\bf e} \cdot {\bf l} = - {\bf l} \cdot {\bf e} = -{\bf le},$\\
	${\bf l}^2 =  ({\bf l} \cdot {\bf e})^2 = -1,$\\
	${\bf e} \cdot ({\bf l e}) = -({\bf l e}) \cdot {\bf e} = {\bf l},$\\
	$({\bf l e}) \cdot {\bf l'} = - ({\bf l l' e}),$
	${\bf l} \neq {\bf l'} \Rightarrow ({\bf l e}) \cdot ({\bf l' e})
		=- {\bf l l'}.$

	\sssec{Definition.}
	The {\em conjugate} of an octave ${\bf o} = {\bf p} + {\bf \ov{q}e}$ is
	defined by\\
	\hth $\ov{\bf o} = \ov{\bf p}- {\bf \ov{q}e},$\\
	the {\em norm} of an octave is defined by\\
	\hth$N({\bf o}) = {\bf o} \cdot {\bf \ov{o}}$.

	\sssec{Theorem.}
	If ${\bf o} = {\bf p} + {\bf \ov{q}e}$, then
	\enumb
	\item $\ov{\ov{{\bf o}}} = {\bf o},$
	\item $\ov{{\bf o}\cdot{\bf o'}} = \ov{{\bf o'}}\cdot\ov{{\bf o}},$
	\item$N({\bf o}) = N({\bf \ov{o}}) = N({\bf p})+N({\bf q}),$
	\item$N({\bf o}) = 0$ iff ${\bf o} = 0,$
	\item$N({\bf o} \cdot {\bf o'}) = N({\bf o})\: N({\bf o'}).$
	\enume

	\sssec{Theorem.}
	\enumb
	\item{\em The octaves is an alternative division ring which is non
	associative.}
	\enume

	For instance, $({\bf i} \cdot {\bf j})\cdot {\bf e} = {\bf k}{\bf e},$
	and ${\bf i} \cdot ({\bf j} \cdot {\bf e}) = -{\bf k}{\bf e}.$

%
%
	\ssec{Ternary Rings.}
	\sssec{Definition.}
	A {\em ternary ring} $(\Sigma,*)$ is a set of elements $\Sigma$
	with at least 2 distinct elements 0 and 1, together with an ternary
	operation " $*$ " such that if $a_1$, $a_2$, $a_3$, $a_4$ are elements in
	$\Sigma$, then
	\enumb
	\item$a_1*a_2*a_3$ is a well defined element of $\Sigma$,
	\item$a_1 * 0 * a_2 = a_2,$
	\item$0 * a_1 * a_2 = a_2,$
	\item$1 * a_1 * 0 = a_1,$
	\item$a_1 * 1 * 0 = a_1,$
	\item$a_1 \neq a_2 \Rightarrow x * a_1 * a_3 = x * a_2 * a_4,$ has
	a unique solution $x \in \Sigma,$
	\item$a_1 * a_2 * y = a_3$ has a unique solution $y \in \Sigma,$
	\item$a_1 \neq a_2 \Rightarrow a_1 * x * y = a_3$ and $a_2 * x * y =
	a_4$ have a unique solution $(x,y)$, $x \in \Sigma,$ $y \in \Sigma.$
	\enume

	\sssec{Theorem.}
	\enumb
	\item $a \neq 0 \Rightarrow \exists a^R \ni a \cdot a^R = 1,$
	$a^R$ is called the right inverse of $a$.
	\item$ a \neq 0 \Rightarrow \exists a^L \ni a^L \cdot a = 1,$
	$a^L$ is called the left inverse of $a$.
	\enume

	\sssec{Definition.}
	The {\em addition} in a ternary ring is defined by\\
	\hth$a \:+\: b := a * 1 * b,$\\
	the {\em multiplication} in a ternary ring is defined by\\
	\hth$ a \:\cdot\: b := a * b * 0.$

	\sssec{Theorem.}
	{\em In a ternary ring} $(\Sigma,*)$,
	\enumb
	\item($\Sigma$,+) is a loop with neutral element 0.
	\item($\Sigma-\{0\},\cdot)$ is a loop with neutral element 1 and
	$a \cdot 0 = 0 \cdot a = 0.$.
	\enume

	\ssec{Felix Klein (1849-1925).  Transformation groups.}
	\label{sec-Sklein}

	The approach which has dominated the non axiomatic study of geometry
	during the last one hundred years has been influenced, almost
	exclusively\footnote{Diedonn\'{e} characterizes it as a "ligne de partage des eaux" in the reedition of the French translation}, by the
	celebrated Inaugural address given by Felix Klein,
	when he became Professor of the Faculty of Philosophy of University of
	Erlangen and a member of its senate in 1872. In it\footnote{Abhandlungen, p.460-497}, Klein states that Geometries are
	characterized by a
	subgroup of the projective group, with, for instance, the group of
	congruences characterizing the Euclidean Geometry. The success of
	this approach to the study of Geometry has been such that in may very
	well have led to the decline of the synthetic Research and Teaching.
	It is hopped that this work, with its underlying program, which
	I call the Berkeley program, will revitalize the subject from the
	high school level on.

	\ssec{Functions.}

	\sssec{Definition.}
	A {\em function} $f$ from a set $D$ to a set $R$ is a set of $ordered$
	pairs $(d_0,r_0),$ $d_0$ in $D,$ $r_0$ in $R,$ such that if to pairs
	have the same first elements, they have the same second element.
	We write $r_0 = f(d_0).$

	\sssec{Definition.}
	The {\em domain} of a function is the set $D'$ which is the union
	of all the first elements of the pairs, the {\em range} of a function
	is the set $R'$ which is the union of all the second pairs.

	\sssec{Definition.}
	A {\em function is one to one} or {\em bijective}  iff  for every pair
	$(d_0,r_0)$ $(d_1,r_1)$ such that $r_0 = r_1$ then $d_0 = d_1.$

	\sssec{Theorem.}
	{\em If a function is one to one, the set of pairs $(r_0,d_0)$ is
	a function $f^{-1}$ from $R$ to $D.$}

	\sssec{Definition.}
	The function $f^{-1}$ is called the {\em inverse} of $f.$

%
	\sssec{Definition.}
	Given 2 functions $f$ and $g$ such that the domain of $g$ is a subset
	of the range of $f,$ the {\em composition} $g \circ f$ is the function
	$(d_0,g(d_0)).$

	\sssec{Theorem.}
	{\em The composition is associative.  In other words,\\
	\hth$	(h \circ g) \circ f = h \circ (g \circ f).$}

	\ssec{Cyclotomic polynomials.  Constructibility with ruler and
	compass.}

	One of the most extensive type of problems in Euclidean Geometry
	is the constructibility of geometric figures using the ruler and the
	compass. The construction of regular polygons lead Gauss,
	in his celebrated Disquisitiones Arithmeticae of 1801, to the study
	of roots of cyclotomic polynomials and his discovery that the
	regular polygon with 17 sides is so constructible. More generally,
	this is the case whenever the number of sides has the form
	$2^n\prod F_j^{i_j}$, where the $F_j$'s are Fermat primes (of the form
	$2^k+1$)
	\footnote{Gauss gives, in n. 366, the polygons with number of sides
	less than 300, constructible with rule and compass, namely,
	$2,4,8,16,32,64,128,256, 3,6,12,24,48,96,192, 5,10,20,40,80,160,\\
	15,30,60,120,240, 17,34,68,136,272, 51,102,204, 85,170, 255, 257.$}.
	In so doing Gauss introduced, for the special case
	of cyclotomic equations, the method, which could be described as
	baby Galois Theory, which was generalized by Galois to the case of
	general polynomial equations. But in his case Gauss gives
	explicitely the various subgroups required to analyze completely the
	solution to the problem.\\
	The general problem of constructibility has been extensively studied,
	I will mention only here the work Emile Lemoine (1902), of Henri
	Lebesgue (1950) and of A. S. Smogorzhevskii (1961).
	In finite geometry, it would appear at first that the ruler is
	sufficient for all constructions because any point in the plane can be
	obtained from 4 points forming a complete quadrangle. But this
	interpretation should be rejected in favor of that which implies that
	the construction of geometric figures should be given completely
	independently of the prime, or power of prime, which characterizes the
	finite Euclidean geometry.  The impression is given that with the ruler
	very little can be constructed. One of the consequences of the results 
	of Chapter 3, is to demonstrate, that both in the finite and classical 
	case many more points, lines, circles, \ldots can be constructed
	with the ruler than heretofore assumed, and that it is a useful
	exercise to reduce the problem  of construction with the compass
	to that of a few points obtained with it and then the ruler alone.
	This is pursued extensively starting with the construction, first of
	the center of the inscribed circle. Of note is that the circle
	of Apollonius can be constructed with the ruler alone.

{\tiny\footnotetext[1]{G17.TEX [MPAP], \today}}

\setcounter{section}{6}
	\section{The real numbers.}

	\ssec{The arithmetization of analysis. [Karl Weierstrass
	(1815-1897) and Riemann (1826-1866)]}
	In his Introduction to the History of Mathematics, Eves
	\footnote{p.426} ascribed the beginning the arithmetization of
	analysis by Weierstrass and his followers to the problem
	presented by the existence (Riemann, 1874) of a continuous curve
	having no tangents at any of its points and that (Riemann) of a
	function which is continuous for all irrational values and discontinuous
	for all rational values in its domain of definition.

	\ssec{Algebraic and transcendental numbers.
	[Hermite (1822-1901) and Lindemann (1852-1939)]}

	\setcounter{subsubsection}{-1}
	\sssec{Introduction.}
	We have seen that the Pythagoreans discovered that if we
	want any circle centered at the origin and passing to a point with
	rational coordinates to intersect always the $x$ axis, irrationals have
	to be introduced.  If that was all that was desired, it would be
	sufficient to construct first the extension field\\
	\hth$	Q(\sqrt{2}) = \{u + v \sqrt{2}\},$\\
	where $u$ and $v$ are in $Q$ then\\
	\hth$	Q(\sqrt{2}(\sqrt{5}) = \{u1 + v1 \sqrt{5}\},$\\
	where $u1$ and $v1$ are in $Q(\sqrt{2}),$  \ldots.  The successive
	integers
	5, 13, 17 are all the primes congruent to 1 modulo 4, because
	of the result of Euler \ldots and because all we would need was to
	obtain the square root of integers which can be written as a sum of 2
	squares.\\
	But, in fact we would like that circles centered at the origin, through
	a point with coordinates in one of these extension fields also intersect
	the $x$ axis at a number in our system.  This requires the introduction
	of algebraic numbers:

	\sssec{Definition.}
	An {\em algebraic number} is one which can be obtained as the
	real solution of a polynomial with integer coefficients.\\
	A {\em transcendental number} is a real number which is not algebraic.

	\sssec{Example.}
	$\sqrt{2}$ is algebraic, being a root of $x^2 - 2 = 0,$\\
	An outstanding problems of the last part of the 19-th century was
	the following, is $\pi$, which is the limit of the ratio of the
	length of a regular polygon with $n$ sides to the diameter, algebraic
	or not.\\
	The proof that it was not algebraic was first give by Lindemann in
	1882, using an earlier result of Hermite of 1973, that $e$ is not
	algebraic.

\setcounter{section}{7}
\section{The pendulum and the elliptic functions.}
	\setcounter{subsection}{-1}
	\ssec{Introduction.}
	This section uses extensively, material learned
	from George Lema\^{\i}tre, in his class on Analytical
	Mechanics, given to first year students in Engineering and in
	Mathematics and Physics, University of Louvain, Belgium, 1942
	and from de la Vall\'{e}e Poussin in his class on elliptic
	functions in 1946.\\
	We first determine the differential equation for the pendulum
	\ref{sec-tpend} using the Theorem of Toricelli \ref{sec-ttor}, we then
	define the elliptic integral of the first kind and the elliptic
	functions of Jacobi \ref{sec-djac1} and \ref{sec-djac2}, we then
	derive the Landen
	transformation which relates elliptic functions with different
	parameters \ref{sec-tlandentr}, use it to obtain
	the Theorem of Gauss which determines the complete elliptic integrals
	of the first kind from the arithmetico-geometric mean of its 2
	parameters \ref{sec-tgaussell}. and obtain the addition formulas for
	the these functions \ref{sec-tjacell} using the Theorem of Jacobi on
	pendular motions which differ by their initial condition
	\ref{sec-tjacpend}.  We also derive the Theorem
	of Poncelet on the existence of infinitely many polynomials inscribed
	in one conic and circumscribed to another \ref{sec-tponc}.  We state,
	without
	proof, the results on the imaginary period of the elliptic functions of
	Jacobi \ref{sec-tjacell1} and \ref{sec-tjacell2}.  A Theorem of
	Lagrange is then given which
	relates identities for spherical trigonometry and those for elliptic
	function \ref{sec-tlagrell}.  Finally we state the definitions and
	some results on the theta functions.
	Using this approach, the algebra is considerably simplified by using
	geometrical and mechanical considerations.\\
	For references, see, Landen (1771), Legendre (1828), Jacobi (1829),
	Eisenstein (1847), Lagrange (Oeuvres), Gauss (Ostwald Klassiker),
	Abel (Oeuvres), Weierstrass (Werke), Cayley (1884), Emch (1901),
	Appell (1924), Bartky (1938), Lema\^{\i}tre (1947), Fettis (1965).

	\ssec{The pendulum.}
	\sssec{Theorem. [Toricelli]}\label{sec-ttor}
	{\em If a mass moves in a uniform gravitational field,
	its velocity $v$ is related to its height $h$ by}
	\enumb
	\item 	$v = \sqrt{2g(h_0-h)},$
	\enume
	{\em where $g$ is the gravitational constant and $h_0$ is a constant,
	corresponding to the height at which the velocity would be 0.}

	Proof:  The laws of Newtonian mechanics laws imply the conservation of
	energy.
	In this case the total energy is the sum of the kinetic energy
	$\frac{1}{2}$ $mv^2$
	and the potential energy $mgh,$ therefore\\
	\hth$	\frac{1}{2} mv^2 + mgh = mgh_0,$ for some $h_0.$

	\sssec{Definition.}
	A {\em circulatory pendular motion} is the motion of a mass $m$
	restricted to stay on a vertical frictionless circular track, whose
	total energy allows the mass to reach with positive velocity the
	highest point on the circle.  An {\em oscilatory pendular motion} is one
	for which the total energy is such that the highest point on the
	circle is not reached.  The mass in this case oscillates back and forth.
	The following Theorem gives the equation satisfied by a pendular
	motion.

	\sssec{Theorem.}\label{sec-tpend}
	{\em If a mass $m$ moves on a vertical circle of radius $R,$ with
	lowest point $A,$ highest point $B$ and center $O,$ its position $M$
	at time $t,$ can be defined by\\
	$2\phi(t) = \angle(AOM)$ which satisfies}
	\enumb
	\item 	$D\phi = \sqrt{a^2-c^2sin^2\phi},$
	{\em where}
	\item 	$a^2 := \frac{gh_0}{2R^2},$ $c^2 = \frac{g}{R},$
	{\em for some $h_0.$}
	\item 	$D^2\phi = -\frac{g}{2R}\:sin \circ (2\phi).$
	\enume

	Proof:  If the height of the mass is measured from $A,$\\
	\hth$	h(t) = R - R cos(2\phi (t)) = 2R sin^2\phi (t),$\\
	the Theorem of Toricelli gives\\
	\hth$	R D(2\phi )(t) = v(t) = \sqrt{2gh_0 - 4gR sin^2\phi (t)},$\\
	hence 0.\\
	The motion is circulatory if $h_0 > 2R$ or $a > c,$ it is
	oscilatory if $0 < h_0 < 2R$ or $c > a.$\\
	2, follows by squaring 0 and taking the	derivative.

	\sssec{Notation.}
	\enumb
	\item 	$k := \frac{c}{a},$ $b^2 := a^2- c^2,$ $k' := \frac{b}{a},$
		$m := k^2$.
	\enume

	\sssec{Theorem. [Jacobi]}\label{sec-tjacpend}
	{\em Let $M(t)$ describes a pendular motion. Given the circle $\gamma$
	which has the line $r$ at height $h_0$ as radical axis and is tangent
	to $AM(t_0),$ if $N(t)M(t)$ remains tangent to that circle, then $N(t)$
	describes the same pendular motion, with $N(t_0) = A.$}

	Proof:  With the abbreviation $M = M(t),$ $N = N(t),$ let $N M$ meets
	$r$ at $D,$
	let $M',$ $N'$ be the projections of $M$ and $N$ on $r,$ let $T$ be
	the point of tangency of $MN$ with $\gamma$,
	\enumb
	\item 	$DM \: DN = DT^2,$\\
	therefore
	\item $\frac{DT}{ND} = \frac{DM}{DT} = \frac{DT-DM}{ND-DT}
		= \frac{MT}{NT} = \sqrt{\frac{DT}{ND}\frac{DM}{DT}}
		= \sqrt{\frac{DM}{ND}} = \sqrt{\frac{M'M}{N'N}}$\\
	When $t$ is replaced by $t+\epsilon$,
	\item $	\frac{v_M}{v_N}
		= lim \frac{M(t+\epsilon ) - M(t)}{N(t+\epsilon) - N(t)}
		= lim \frac{M(t) T}{N(t+\epsilon) T} = \frac{MT}{NT},$
	
	because the triangles $T,M,M(t+\epsilon)$ and $T,N,N(t+\epsilon)$ are
	similar, because $\angle(T,N,N(t+\epsilon) = \angle(T,M(t+\epsilon),M)$
	as well as $\angle(M(t+\epsilon),T,M) = \angle(N(t+\epsilon),T,N).$\\
	Therefore
	\item $	\frac{v_M}{v_N} = \sqrt{\frac{M'M}{N'N}}.$
	\enume
	The Theorem of Toricelli asserts that $v_M = \sqrt{2g M'M},$ this
	implies, as we have just seen, $v_N = \sqrt{2g N'N},$ therefore $N$
	describes the same pendular motion with a difference in the origin of
	the independent variable.

	\sssec{Corollary.}
	{\em If $M = B$ and $N = A,$ the line $M(t)\times N(t)$ passes through a
	fixed point $L$ on the vertical through $O$.\\
	Moreover, if $b := BL$ and $a := LA,$ we have\\
	\hth$\frac{v_M}{v_N} = \frac{b}{a}$ and $h_0 = \frac{a^2}{a - b}.$}

	This follows at once from from \ref{sec-tjacpend}.2, and 1.

	\sssec{Definition.}
	The point $L$ of the preceding Corollary is called {\em point of
	Landen}.

	\sssec{Theorem. [Poncelet]}\label{sec-tponc}
	{\em Given 2 conics $\theta$  and $\gamma$ , if a polygon $P_i,$ $i = 0$
	to $n,$ $P_n = P_0,$ is such that $P_i$ is on $\theta$ and
	$P_i\times P_{i+1}$ is tangent to $\gamma$ , then there exists
	infinitely many such polygons.\\
	Any such polygon is obtained by choosing $Q_0$ on $\theta$ drawing a
	tangent $Q_0Q_1$ to $\gamma$, with $Q_1$ on $\theta$ and successively
	$Q_i,$ such that $Q_i$ is on $\theta$ and $Q_{i-1}\times Q_i$ is
	tangent to $\gamma$, the Theorem asserts that $Q_n = Q_0.$}

	The proof follows at once from \ref{sec-tjacpend}, after using
	projections which transform the circle $\theta$  and the circle
	$\gamma$  into the given conics.\\
	The Theorem is satisfied if the circle have 2 points in common or not.

	\sssec{Theorem.}\label{sec-tlandentr}
	{\em If $M(t)$ describes a circular pendular motion, then the
	mid-point $M_1(t)$ of $M(t)$ and $M(t+K)$ describes also a circular
	pendular motion.  More precisely, $M_1(t)$ is on a circle with diameter
	$LO,$ with $LA = a,$ $LB = b,$ and if\\
	$\phi_1(t) = \angle(O,L,M_1(t),$}
	\enumb
	\item 	$t = \int_0^{\phi(t)} \frac{D\phi}{\Delta}
		= \frac{1}{2} \int_0^{\phi_1(t)}\frac{D\phi_1}{\Delta_1}.$\\
	{\em where}
	\item 	$\Delta^2 := a^2 cos^2\phi  + b^2 sin^2\phi$  and
		$\Delta^2_1 := a^2_1 cos^2\phi _1 + b^2_1 sin^2\phi _1,$\\
	{\em where the relation between $\phi$ and $\phi_1$ is given by}
	\item 	$tan(\phi_1-\phi ) = k' tan\phi$, {\em or}
	\item 	$sin(2\phi - \phi_1) = k_1 sin\phi_1,$\\
	{\em with}
	\item 	$k' := \frac{b}{a}$, $k_1 := \frac{c_1}{a_1}$,
	\item	$a_1 := \frac{1}{2}(a+b),$ $b_1 := \sqrt{ab},$
		$c_1 := \frac{1}{2}(a-b),$ {\em therefore}
	\item 	$a = a_1 + c_1,$ $b = a_1-c_1,$ $c = 2\sqrt{a_1c_1}. $
	\enume

	Proof:  First, it follows from the Theorem of Toricelli that the
	velocity $v_A$ at $A$ and $v_B$ at $B$ satisfy\\
	\hth$v_A = \sqrt{2gh_0} = 2R a,$
		$v_B = \sqrt{2gh_0-2R} = \sqrt{4R^2a^2 - 4 c^2R^2} = 2R b,$\\
	therefore $\frac{BL}{LA} = \frac{b}{a}.$\\
	If $P$ is the projection of $L$ on $BM$ and $Q$ the projection of
	$L$ on $AM,$\\
	\hth$LM^2 = LP^2 + LQ^2 = a^2 cos^2\phi + b^2 sin^2\phi = \Delta ^2.$\\
	\hth$LQ = LM cos(\phi_1-\phi) = a cos \phi.$

	We can proceed algebraically. Differentiating 2. gives\\
	$a(1 + tan^2(\phi_1-\phi)) (D\phi_1 - D\phi) = b(1+tan^2\phi) D\phi,$\\
	or\\
	$a(1 + tan^2(\phi_1-\phi)) D\phi_1
	\hth	= (a(1+tan^2(\phi_1-\phi) + b(1+tan^2\phi)) D\phi\\
	\hth    = (a+b + \frac{b^2}{a} tan^2\phi + b tan^2\phi) D\phi\\
	\hth    = (a+b)(1 + \frac{b}{a} tan^2\phi) D\phi\\
	\hth    = (a+b)(1 + tan\phi tan(\phi_1-\phi)) D\phi,$\\
	or\\
	\hth$\frac{a}{cos^2(\phi_1-\phi)} D\phi_1
	= 2a_1 \frac{cos(2\phi -\phi_1)}{cos\phi cos(\phi_1-\phi)} D\phi,$ or\\
	\hth$	\frac{D\phi}{\frac{a cos\phi}{cos(\phi_1-\phi)}}
	= \frac{D\phi_1}{2a_1 cos(2\phi -\phi_1)},$\\
	or because $LM = \Delta$ \\
	\hth$	\frac{D\phi}{\Delta} = \frac{D\phi_1}{2 \Delta_1}.$

	We can also proceed using kinematics.\\
	The velocity at $M$ is\\
	\hth$	v_M = 2R D\phi = 2R \Delta ,$\\
	If we project the velocity vector on a perpendicualr to $LM,$\\
	\hth$LM D\phi_1 = v_M cos(2\phi_1-\phi)
		= 2R cos(2\phi_1-\phi) \Delta \phi.$\\
	Therefore\\
	\hth$\frac{D\phi}{\Delta} = \frac{D\phi_1}{ 2R cos(2\phi_1-\phi )}
	= \frac{a_1}{2R} \frac{D\phi_1}{\Delta_1} = \frac{D\phi_1}{2\Delta_1}.$

	\sssec{Definition.}
	The transformation from $\phi$ to $\phi_1$ is called the {\em forward
	Landen transformation}.\\
	The transformation from $\phi_1$ to
	$\phi$  is called the {\em backward Landen transformation}.

	These transformations have also been applied to the integrals of the
	second kind and of the third kind.

	\sssec{Comment.}
	The formulas 3. and 1. are the formulas which are used to
	compute $t$ from $\phi(t)$.  The formulas 4. and 2. are used to compute
	$\phi(t)$ from $t$.

	\sssec{Comment}
	Given the first order differential equations,\\
	\hth$(Dy)^2	= C_0 (y^2 + A_0) (y^2 + B_0),$\\
	\hth$(Dz)^2	= C_1 (z^2 + A_1) (z^2 + B_1),$\\
	with\\
	\hth$z = d (y + \frac{l}{y}),$ $l\neq 0,$ $d > 0.$\\
	These equations are compatible iff\\
	$d^2(1-	\frac{l}{y^2})^2 C_0(y^2 + A_0) (y^2 + B_0)
		= C_1(d^2(\frac{y^2+l}{y})^2 + A_1)(d^2(\frac{y^2+l}{y})^2
			 + B_1)$\\
	this requires $\sqrt{l}$ to be a root of one of the factors of the
	second member, let it be the second factor, this implies\\
	\hth$d^2 4 l  + B_1 = 0,$\\
	then, the second factor becomes,\\
	\hth$d^2(\frac{y^2+l}{y})^2 + B_1
		= d^2((\frac{y^2+l}{y})^2  - 4 l)
		= d^2(\frac{y^2-l}{y})^2$,\\
	therefore $\sqrt{l}$ is a double root of the second member and\\
	\hth$C_0(y^4 + (A_0+B_0) y^2 + A_0B_0)
		= d^2 C_1(y^4 + (2 l + \frac{A_1}{d^2})y^2 + l^2),$
	therefore\\
	\hth$C_0 = d^2C_1,$ $A_0B_0 = \frac{ B_1^2}{16 d^4},
		 A_0+B_0 =  \frac{A_1 - \frac{1}{2} B_1}{d^2},$\\
	For real transformations, $A_0B_0 > 0,$ if $j_0 = sign(B_0)$ and
	$j_1 = sign(B_1),$\\
	\hth$B_1 = 4j_1d^2\sqrt{A_0B_0},$
		$A_1 = d^2(A_0+B_0+2j_1\sqrt{A_0B_0}$\\
	\hti{11}$= j_0d^2(\sqrt{|A_0|} + j_0j_1\sqrt{|B_0|})^2.$\\
	If we want $A_1B_1 > 0$ then $j_0 = j_1.$

	\ssec{The elliptic integral and the arithmetico-geometric mean.}
	\setcounter{subsubsection}{-1}
	\sssec{Introduction.}
	Gauss began his investigations after he showed that the length of the 
	lemniscate could be computed from the arithmetico geometric mean of
	$\sqrt{2}$ and 1. More precisely, the lemniscate  is the curve
	$r^2=cos(2\theta)$, in polar coordinates. A quarter of its length
	is given by the integral
	\[\int_0^1 \frac{dr}{\sqrt{1-r^4}},\]
	which is easily deduced from the general formula for the square
	of the arc length in polar coordinates, $ds^2 = dr^2+r^2 (d\theta)^2.$\\
	Gauss observed that to 9 decimal places the integral was
	1.311028777 and so is $\frac{\pi/2}{agm(\sqrt{2},1}$,
	where $agm(a,b)$ denotes the arithmetico geometric mean of 2 numbers,
	defined below.

	\sssec{Theorem. [Gauss]}\label{sec-tgaussagm}
	{\em Given $a_0 > b_0 > 0,$ let}
	\enumb
	\item 	$a_{i+1} := \frac{1}{2} (a_i+b_i),$
	\item 	$b_{i+1} := \sqrt{a_ib_i},$
	\item {\em The sequences $a_i$ and $b_i$ have a common limit}
		$a_{\infty }.$
	\item {\em The sequence $a_i$ is monotonically decreasing and the
		sequence $b_i$ is monotonically increasing.}
	\enume

	Proof:  Because\\
	\hth$a_i > a_{1+1},$ $b_{i+1} > b_i,$\\
	it follows that the sequence $a_i$ is bounded below by $b_0,$ the
	sequence $b_i$ is bounded above by $a_0,$ therefore both have a limit
	$a_{\infty }$ and $b_{\infty }.$ Taking the limit of 0. gives at once
	$a_{\infty } = b_{\infty }.$

	\sssec{Definition.}\label{sec-agm}
	$a_\infty$ is called the {\em arithmetico-geometric mean} of $a_0$ and
	$b_0$.

	\sssec{Example.}
	With $a_0 = \sqrt{2}$ and $b_0 = 1$,\\
	$a_1 = 1.207106781,$ $b_1 = 1.189207115$,\\
	$a_2 = 1.198156948,$ $b_2 = 1.198123521$,\\
	$a_3 = 1.198140235,$ $b_2 = 1.198140235$.

	\sssec{Definition.}\label{sec-djac1}
	If $a = 1,$ and we express $t$ in terms of $\phi(t),$
	\enumb
	\item 	$t = \int_0^{\phi(t)}\frac{1}{\sqrt{1 - k^2sin^2}}.$
	This integral is called the {\em incomplete elliptic integral of the
	first kind}. Its inverse function $\phi$  is usually noted
	\item 	$am := \phi$, the {\em amplitude function},
	\item 	$K := \int_0^{\frac{\pi}{2}} \frac{1}{\sqrt{1 - k^2sin^2}}$
	is called the {\em complete integral of the first kind},  it gives
	half the period, $\frac{K}{a},$ for the circular pendulum.
	\enume

	\sssec{Theorem.}
	\enumb
	\item For the circulatory pendulum, the angle $2\phi$ between the
	lowest position of the mass and that at time $t$ is given by
	$\phi = am(at)$. The coordinates are $R\: sin(2\phi),$
	$R-R\: cos(2\phi)$.
	\item For the oscillatory pendulum, if the highest point is
	$2R sin^2(\alpha)$ above the lowest point, the angle $2\theta$ between
	the lowest position of the mass and that at time $t$ is given by
	$sin\theta = sin\phi\:sin\alpha$ where $\phi$ is given by
	$\phi = am(at,sin^2\alpha)$.
	\enume

	\sssec{Theorem.}\label{sec-tgaussell}
	{\em For the complete integrals we have}
	\enumb
	\item 	$\frac{K}{a} = \int_0^{\frac{\pi}{2}}
		 \frac{1}{\sqrt{a^2cos^2 + b^2sin^2}}
		= \frac{\frac{\pi}{2}}{a_{\infty}}.$
	\enume

	Proof:  If $\phi(K) = \frac{\pi}{2}$, then $\phi_1(K) = \pi$,
	therefore
	\begin{enumerate}
	\item $K = \int_0^{\frac{\pi}{2}}\frac{D\phi}{\Delta}
		= \int_0^{\pi} \frac{D\phi_1}{2 \Delta_1}
		= \frac{1}{2} \int_0^{\frac{\pi}{2}} \frac{D\phi_1}{\Delta_1}
	+ \frac{1}{2} \int_{\frac{\pi}{2}}^{\pi} \frac{D\phi_1}{\Delta_1}
		= \int_0^{\pi} \frac{D\phi_1}{\Delta_1}
		= \int_0^{\frac{\pi}{2}} \frac{D\phi_n}{\Delta_n}$\\
	\hti{3}$= \int_0^{\frac{\pi}{2}} \frac{1}{a_{\infty }}
		= \frac{\frac{\pi}{2}}{a_{\infty }}.$
	\enume

	\ssec{The elliptic functions of Jacobi.}
	\sssec{Definition.}\label{sec-djac2}
	The functions
	\enumb
	\item 	$sn := sin \circ am,$ $cn := cos \circ am,$
	$dn := \sqrt{1 - k^2sn^2},$ \\
	are called the {\em elliptic functions of Jacobi}.

	The functions which generalize $tan,$ $cosec,$  \ldots are
	\item 	$ns := \frac{1}{sn}, nc := \frac{1}{cn}, nd := \frac{1}{dn},$
	\item	$sc := \frac{sn}{cn}, cd := \frac{cn}{dn}, ds := \frac{dn}{sn},$
	\item	$cs := \frac{cn}{dn}, dc := \frac{dn}{cn}, sd := \frac{sn}{dn}.$
	\enume

	\sssec{Theorem.}\label{sec-tjacell0}
	{\em If}
	\enumb
	\item   $s_1 := sn(t_1),$ $c_1 = cn(t_1),$ $d_1 = dn(t_1)$ {\em and}
	\item   $s_2 := sn(t_2),$ $c_2 = cn(t_2),$ $d_2 = dn(t_2),$ \\
	{\em we have}
	\item 	$sn^2 + cn^2 = 1,$ $dn^2 + k^2 sn^2 = 1,$
		$dn^2 - k^2cn^2 = k'^2.$
	\item 	$1 - k^2s^2_1 s^2_2 = c^2_1 + s^2_1 d^2_2
	= c^2_2 + s^2_2 d^2_1.$
	\enume

	\sssec{Lemma.}\label{sec-ljac}
	\enumb
	\item 	$c_2 = c_1 cn(t_1+t_2) + d_2s_1 sn(t_1+t_2),$
	\item 	$d_2 = d_1 dn(t_1+t_2) + k^2 s_1c_1 sn(t_1+t_2).$
	\enume

	Proof:  We use the Theorem \ref{sec-tjacpend} of Jacobi.  Let $R$ be
	the radius of
	$\theta$ and $O$ its center, let $r$ be the radius of $\gamma$ and $O'$
	its center, let $s := OO'.$  Let $A,$ $N,$ $M',$ $M$ be the position of
	the mass at time 0, $t_1,$ $t_2,$ $t_1+t_2.$\\
	The lines $A \times M'$ and $N \times M$ are tangent to the same circle
	$\gamma$ at $T'$ and $T.$\\
	Let $X$ be the intersection of $O \times M$ and $O' \times T,$
	$2\phi := \angle(A,O,N),$
	\begin{enumerate}
	\setcounter{enumi}{1}
	\item $2\phi' := \angle(A,O,M),$\\
	we have $\angle(N,O,M) = 2(\phi'-\phi ),$ $\angle(M,X,T) = \phi'-\phi,$
	$\angle(T,O',O) = \phi'+\phi$.\\
	If we project $MOO'$ on $O'T,$\\
	\hth$	r = R cos(\phi'-\phi ) s cos(\phi'+\phi),$ or
	\item $r = (R+s) cos\phi  cos\phi' + (R-s) sin\phi sin\phi'.$\\
	\hth	$\phi = am t_1,$ $\phi' = am(t_1+t_2),$ \\
	\hth	$sin\phi' = sn(t_1+t_2),$ $cos\phi' = cn(t_1+t_2),$ \\
	\hth	$sin\phi = sn\: t_1 = s_1,$\\
	\hth	$cos\phi = cn\: t_1 = c_1,$
	\enume
	when $t_1$ = 0,\\
	\hth$cos(\angle(A,B,M') = cn\:t_2 = c_2 = \frac{BM'}{AB}
		= \frac{O'T'}{AO'} = \frac{r}{R + s}$,\\
	the ratio of the velocities is\\
	\hth$\frac{v_{M'}}{v_A} = \frac{dn \:t_2}{dn\: 0} = d_2
		= \frac{TM'}{AT} = \frac{O'B}{AO'} = \frac{R - s}{R + s},$
	substituting in 2. gives 0.\\
	The proof of 1. is left as an exercise.

	\sssec{Theorem. [Jacobi]}\label{sec-tjacell}
	\enumb
	\item $\frac{sn\: u_1 cn\:  u_2 dn\: u_2 + sn\: u_2 cn\: u_1 dn\:u_1}
		{sn(u_1+u_2)}= 1 - k^2 sn^2 u_1 sn^2 u_2.$
	\item $\frac{cn\:u_1 cn\:u_2 - sn\:u_1 dn\:u_1 sn\:u_2 dn\:u_2}
		{cn(u_1+u_2) } = 1 - k^2 sn^2 u_1 sn^2 u_2. $
	\item $\frac{dn\:u_1 dn\:u_2 - k^2 sn\:u_1 sn\:u_2 cn\:u_1 cn\:u_2}
		{dn(u_1+u_2) } = 1- k^2 sn^2 u_1 sn^2 u_2. $
	\enume

	Proof: Let $w = \frac{1}{1 - k^2 s^2_1 s^2_2}.$\\
	Let $s_1,$ $s_2,$  $\ldots$  denote $sn\:u_1,$ $sn\:u_2,$  $\ldots$,
	define $S$ and $C$ such that\\
	\hth$sn(u_1+u_2) = S w,$ $cn(u_1+u_2) = C w.$\\
	The \ref{sec-ljac}.0. gives\\
	\hth$c_2 = c_1 C w + d_2 s_1 S w$ or
	\begin{enumerate}
	\setcounter{enumi}{2}
	\item $c_1 C w = - d_2 s_1 S w + c_2,$
	\enume
	\ref{sec-tjacell0}.2. gives\\
	\hth$S^2 w^2 + C^2 w^2 = 1,$\\
	eliminating $C$ gives the second degree equation in $S w$:\\
	\hth$(c^2_1 + d^2_2 s^2_1 (S w)^2 - 2 s_1 c_2 d_2 (S w) + c^2_2 - c^2_1
		 = 0,$\\
	one quarter of the discriminant is\\
	\hth$s^2_1 c^2_2 d^2_2 - (c^2_2 - c^2_1)(c^2_1 + d^2_2 s^2_1)\\
	\hth	= s^2_1 c^2_2 d^2_2 - c^2_1 c^2_2 + c^4_1 -s^2_1 c^2_2 d^2_2
		+ s^2_1 c^2_1 d^2_2\\
	\hth	= c^2_1(c^2_1 - c^2_2 + s^2_1 d^2_2) = c^2_1 s^2_2 d^2_1,$\\
	therefore\\
	\hth$	S w = (s_1 c_2 d_2 \pm  c_1 d_1 s_2)w.$\\
	One sign correspond to one tangent from $M$ to $\gamma$ , the other to
	the other
	tangent, therefore one corresponds to the addition, the other to the
	subtration formula.  From the special case $k = 0,$ follows that, by
	continuity, the $+$ sign should be used.  This gives 0., 1. follows from
	3, 2. is left as an exercise.

	\sssec{Corollary.}
	\enumb
	\item $	sn(u+K)  = cd(u),$ $cn(u+K) = -k' sd(u),$ $dn(u+K)  = k' nd(u).$
	\item $	sn(u+2K) = -sn(u),$ $cn(u+2K) = - cn(u),$ $dn(u+2K) = dn(u).$
	\item $	sn(u+4K) = sn(u),$ $cn(u+4K) = cn(u),$ $dn(u+4K) = dn(u).$
	\enume

	\sssec{Definition.}
	\hth$K'(k^2) = K(k'^2).$

	\sssec{Theorem.}\label{sec-tjacell1}
	\enumb
	\item 0. $k sn \circ I+iK' = sn,$\\
	 1.	 $i k cn \circ I+iK' = ds,$\\
	 2.	 $i dn \circ I+iK' = cs,$
	\item 0. $sn \circ I+2iK' = sn,$\\
	 1.	 $cn \circ I+2iK' = - cn,$\\
	 2.	 $dn \circ I+2iK' = - dn,$
	\enume

	\sssec{Theorem.}\label{sec-tjacell2}
	\enumb
	\item {\em $sn$ has periods $4K$ and $2iK'$ and pole $\pm iK',$}
	\item {\em $cn$ has periods $4K$ and $4iK'$ and pole $\pm iK',$}
	\item {\em $dn$ has periods $2K$ and $4iK'$ and pole $\pm iK'.$}
	\enume

	\sssec{Theorem.}
	\enumb
	\item $	k = 0 \Rightarrow  sn = sin,$ $cn = cos,$ $dn = \underline{1},$
	\item $	k = 1 \Rightarrow  sn = tanh,$ $cn = sech,$ $dn = sech.$
	\enume

	\ssec{The theta functions of Jacobi.}

	\sssec{Definition.}
	Given the parameter $q,$ called the nome,
	\enumb
	\item $	q := e^{-\pi \frac{K'}{K}},$\\
	the functions
	\item $\theta_1 := 2 q^{\frac{1}{4}} \sum_{n=0}^{\infty}
		(-1)^n q^{n(n+1)} sin(2n+1)I$
	\item $\theta_2 := 2 q^{\frac{1}{4}} \sum_{n=0}^{\infty}
		q^{n(n+1)} cos(2n+1)I$
	\item $\theta_3 := 1 + 2 \sum_{n=1}^{\infty} q^{n^2 } cos 2nI$
	\item $\theta_4 := 1 + 2 \sum_{n=1}^{\infty} (-1)^n q^{n^2 } cos 2nI$
	are the {\em theta functions of Jacobi}.
	\enume

	\sssec{Definition.}
	The functions, with $v = \frac{\pi\:I}{2K}$
	\enumb
	\item 	$\theta_s := \frac{2K\: \theta_1 \circ v}{D \theta_1(0)},$
	$\theta_c := \frac{\theta_2 \circ v}{\theta_2(0)},$
	$\theta_d := \frac{\theta_3 \circ v}{\theta_3(0)},$
	$\theta_n := \frac{\theta_4 \circ v}{\theta_4(0)},$\\
	are called the {\em theta functions of Neville}.
	\enume

	\sssec{Theorem.}
	{\em If $p,q$ denote any of $s,$ $c,$ $d,$ $n,$\\
	\hth$pq = \frac{\theta_p}{\theta_q}.$\\
	For instance\\
	\hth$sn = \frac{\theta_s}{\theta_n }
		= \frac{2K \theta_1 \circ v}{D\theta_1(0)}
		. \frac{\theta_4(0)}{\theta_4 \circ v}.$}

	\sssec{Theorem.}
	{\em The Landen transformation replaces the parameter $q,$ by $q^2.$}

	\ssec{Spherical trigonometry and elliptic functions.}
	\sssec{Theorem. [Lagrange]}\label{sec-tlagrell}
	{\em From the addition formulas of elliptic functions,
	we can derive those for a spherical triangle as follows.  Let}
	\enumb
	\item   $u_1 + u_2 + u_3 = 2K,$\\
	{\em define}
	\item   $sin a := - sn u_1,$ $cos a := - cn u_1,$\\
	    $sin b := - sn u_2,$   $cos b := - cn u_2,$\\
	    $sin c := - sn u_3,$   $cos c := - cn u_3,$\\
	    $sin A := - k\: sn u_1,$ $cos A := - dn u_1,$\\
	    $sin B := - k\: sn u_2,$ $cos B := - dn u_2,$\\
	    $sin C := - k\: sn u_3,$ $cos C := - dn u_3,$\\
	{\em then to any formula for elliptic functions of $u_1,$ $u_2,$ $u_3,$
	corresponds
	a formula for a spherical triangle with angles $A,$ $B,$ $C$ and sides
	$a,$ $b,$ $c.$  For instance,}
	\item   $\frac{sin A}{sin a} = \frac{sin B}{sin b}
		= \frac{sin C}{sin c} = k.$
	\item  $  cos a = cos b\: cos c + sin b\: sin c\: cos A,$
	\item  $  cos A = - cos B\: cos C + sin B\: sin C\: cos a,$
	\item  $  sin B\: cot A = cos c\: cos B + sin c\: cot a.$
	\enume

	Proof.  2. follows from the definition.  3. follows from\\
	\hth$c_2 = c_1 cn(u_1+u_2) + d_2 s_1 sn(u_1+u_2)$
	after interchanging $u_1$ and $u_2$ and using
	\begin{enumerate}
	\setcounter{enumi}{5}
	\item 0. $sn(u_1+u_2) = sn(2K-u_1-u_2) = sn\: u_3 = s_3,$ \\
	 1.	$cn(u_1+u_2) = - cn(2K-u_1-u_2) = - cn\: u_3 = - c_3,$\\
	 2.	$dn(u_1+u_2) = dn(2K-u_1-u_2) = dn \:u_3 = d_3,$\\
	similarly, 4. follows from\\
	\hth$c_2	= c_1 cn(u_1+u_2) + d_2 s_1 sn(u_1+u_2)$ \\
	after interchanging $u_1$ and $u_2$ and using 6, and 5. from\\
	\hth$sn\:u_2 dn\:u_1 = cn\:u_1 sn(u_1+u_2)
		- sn\:u_1 dn\:u_2 cn(u_1+u_2)$\\
	after division by $sn\:u_1.$
	\enume

	\ssec{The p function of Weierstrass.}
	\setcounter{subsubsection}{-1}
	\sssec{Introduction.}
	Because it is not germane in this context, I will only mention briefly
	the important contribution of Weierstrass, which proved that all
	doubly periodic meromorphic functions can be expressed in terms of one 
	of them, the $p$ function. The addition formulas for this function
	anf for the Jacobi functions and many other properties generalize to
	the finite case (De Vogelaere, 1983)..

	\ssec{References.}
	\enumb
\item Abel, Niels Henrik, {\em Oeuvres Compl\`{e}tes}, Nouv. ed., publi\'{e}es
	aux frais de l'Etat norv\'{e}gien par L. Sylow et S. Lie, Christiania,
	Grondahl \& Son, 1881, Vol. 1,2. 
\item Appell, Paul Emile \& Dautheville, S., {\em Pr\'{e}cis de M\'{e}canique
	Rationnelle}, Paris, Gauthier-Villars, 1924, 721 pp.
\item Bartky, {\em Numerical Calculation of Generalized Complete Integrals},
	 Rev. of Modern Physics, 1938, Vol. 10, 264-
\item Cayley, Arthur, {\em On the Addition of Elliptic Functions}, Messenger of
	Mathematics, Vol. 14, 1884, 56-61.
\item Cayley, Arthur, {\em Note sur l'Addition des Fonctions Elliptiques},
	Crelle J., Vol. 41, 57-65.
\item De Vogelaere, Ren\'{e}, {\em Finite Euclidean and non-Euclidean Geometry
	with application to the Finite Pendulum and the Polygonal Harmonic
	Motion. A First Step to Finite Cosmology}. The Big Bang
	and Georges Lema\^{i}tre, Proc. Symp. in honor of 50 years after his
	initiation of Big-Bang Cosmology, Louvain-la-Neuve, Belgium,
	October 1983., D. Reidel Publ. Co, Leyden, the Netherlands. 341-355.
\item Eisenstein, Ferdinand Gotthold Max, {\em Mathematische Abhandlungen}, 
	besonders aus dem Gebiete der hoheren Arithmetik und den elliptischen
	Functionen, Mit einer Vorrede von C.F. Gauss. (Reprografischer 
	Nachdruk der Ausg., Berlin 1847.) Hildesheim, G. Olms, 1967.
\item Emch, {\em An Application of Elliptic Functions}, Annals of Mathematics,
	Ser. 2, Vol. 2, 1901. III.4.4.0.
\item Fettis, Henri E., Math. of Comp., 1965.
\item Gauss, Carl Friedrich, Ostwald Klassiker der Exakten Wissenschaften,
	Nr 3.
\item Jacobi, Karl Gustav Jakob, {\em Fundamenta Nova Theoriae Funktionum
	Ellipticarum}, 1829.
\item Landen, John, Phil. Trans. 1771, 308.
\item Lagrange, Joseph Louis, comte, {\em Oeuvres}, publi\'{e}es par les soins
	de m. J.-A. Serret, sous les auspices de Son Excellence le ministre de
	l'instruction publique, Paris, Gauthier-Villars, 1867-92.
\item Legendre, Adrien Marie, {\em Trait\'{e} des Fonctions Elliptiques et des
	Int\'{e}grales Euleriennes}, avec des tables pour en faciliter le calcul
	num\'{e}rique, Paris, Huzard-Courcier, Vol. 1-3, 1825-1828.
\item Lema\^{\i}tre, Georges, {\em Calcul des Int\'{e}grales Elliptiques}, Bull.
	Ac. Roy. Belge, Classe des Sciences, Vol. 33, 1947, 200-211.
\item Weierstrass, Karl, {\em Mathematische Werke}, Hildesheim, G. Olms, New
	York, Johnson Reprint 1967.
	\enume

	\ssec{Texts on and tables of elliptic Functions.}
	\enumb
\item Abramowitz, Milton \& Stegun, Irene A. Edit., {\em Handbook of
	Mathematical Functions}, U.S. Dept of Commerce, Nat. Bur. of Stand.,
	Appl. Math, Ser., number 55, 1964, 1046 pp.
\item Adams, Edwin Plimpton, Ed., {\em Smithsonian Mathematical Formulae and
	Tables of Elliptic Functions}, under the direction of Sir George
	Greenhill, 3d reprint, City of Washington, 1957, Its Smithsonian
	miscellaneous collections, v.74, no.1, Smithsonian Institution,
	Publication 2672.
\item Halphen, Georges Henri, {\em Trait\'{e} des Fonctions Elliptiques et
	de leurs Applications}, Paris, Gauthier-Villars, 1886-91.
\item Hancock, Harris, {\em Lectures on the Theory of Elliptic Functions}, v.
	1. 1st ed., 1st thousand, New York, J. Wiley, 1910. Dover Publ., 1958.
\item Jahnke, Eugen \& Emde, Fritz, {\em Tables of Functions with Formulae and
	Curves}, 4th ed., New York, Dover Publications, 1945.
\item Jahnke, Eugen \& Emde, Fritz \& Losch, Friedrich, {\em Tables of Higher
	Functions}, 6th ed. rev. by Losch, New York, McGraw-Hill, 1960.
\item Jordan, Camille, {\em Fonctions Elliptiques}, New York, Springer-Verlag,
	1981.
\item King, Louis Vessot, {\em On the Direct Numerical Calculations of
	Elliptic Functions and Integrals}, Cambridge, Eng., Univ. Press, 1924.
\item Lang, Serge, {\em Elliptic functions}, 2nd ed., New York,
	Springer-Verlag, 1987, Graduate texts in mathematics, 112.
\item Mittag-Leffler, Magnus Gustaf, {\em An Introduction to the Theory of
	Elliptic Functions}, Lancaster, Pa., 1923, Hamburg, Germany, Lutcke \&
	Wulff.
\item Neville, Eric Harold, {\em Elliptic Functions, a primer}, Prepared for
	publication by W. J. Langford, 1st d. Oxford, New York, Pergamon
	Press, 1971.
\item Neville, Eric Harold, {\em Jacobian Elliptic Functions}, Oxford, The
	Clarendon Pr., 1944.
\item Oberhettinger, Fritz Wilhelm \& Magnus, Wilhelm, {\em Anwendung der
	elliptischen Funktionen in Physik und Technik}, Berlin, Springer, 1949,
	Die Grundlehren der mathematischen Wissenschaften in
	Einzeldarstellungen, Bd. 55. 
\item Riemann, Bernhard, {\em Elliptische Functionen}, Mit zusatzen
	herausgegeben von Hermann Stahl \ldots Leipzig, B. G. Teubner, 1899.
\item Schuler, Max \& Gebelein, H., {\em Eight and Nine Place Tables of
	Elliptical Functions based on Jacobi Parameter $q$}, with an English
	text by Lauritz S. Larsen, Berlin, Springer-Verlag, 1955, XXIV+296 pp.
\item Spenceley, G.W. \& Spenceley, R.M., {\em Smithsonian Elliptic Functions
	Tables}, Washington, Smithsonian Institution 1947, Smithsonian
	miscellaneous collections v. 109.
\item Sturm, Charles Fran\c{c}ois, {\em Cours d'Analyse de l'Ecole
	Polytechnique}, revu et corrig\'{e} par E. Prouhet et augment\'{e} de
	la Th\'{e}orie \'{e}l\'{e}mentaire des fonctions elliptiques, par H.
	Laurent, 14. ed., rev. et mise au courant du nouveau programme de la
	licence, par A. de Saint-Germain, Paris, Gauthier-Villars, 1909.
\item Tannery, Jules \& Molk, Jules, {\em El\'{e}ments de la th\'{e}orie des
	fonctions elliptiques}, Paris, Gauthier-Villars, Vol. 1-4, 1893-1902.
\item T\"{o}lke, Friedrich, {\em Praktische Funktionenlehre}, Berlin, Springer,
	Vol. I to VI ab. Vol. 3 and 4 deal with the elliptic functions of
	Jacobi.
\item Tricomi, Francesco Giacomo, {\em Elliptische Funktionen}, \"{u}bers. und
	bearb. von Maximilian Krafft. Leipzig, Akademische
	Verlagsgesellschaft Geest \& Portig, 1948, Mathematik und ihre
	Anwendungen in Physik und Technik, Bd. 20.
\item Whittaker, Edmund Taylor \& Watson, G. N., {\em A Course of Modern
	Analysis}, an introduction to the general theory of infinite processes
	and of analytic functions, Cambridge, Eng., Univ. Pr., 1963, 606 pp.,
	(1927).
	\enume


\setcounter{section}{8}
	\section{Model of Finite Euclidean Geometry in Classical Euclidean
		Geometry.}

	\setcounter{subsection}{-1}
	\ssec{Introduction.}

	The purpose of this section is to give an informal introduction to
	finite Euclidean geometry for those familiar with classical Euclidean
	geometry and analytic geometry.\\
	The definitions of points and lines will be given in terms of
	equivalence classes. The Theorems will be derived from these
	definitions or can be derived from the classical Theorems.
	I will restrict myself to the 2 dimensional case and will not attempt
	to give the most general results.  In particular, I will assume that
	distances are defined in only one way.\\
	In this restricted framework there is one finite geometry for each
	prime integer $p.$  $p$ is assumed to be larger than 2, non degenerate
	circles require $p$ larger than $3$.  The examples correspond
	to small $p.$  The reader is encouraged to think of the implications
	when $p$ is very large, for instance of the order of $10^{32}$ say,
	and is looking at points with coordinates of the order of $10^8$ to
	$10^{20}$.

	\ssec{Points and lines in finite Euclidean geometry.}

	\sssec{Notation.}
	A point  $P'$  in Euclidean geometry will be denoted by its
	cartesian coordinates $((x,y))$ given between double parenthesis.  
	A line $l'$  will be denoted by the coefficients $[[a,b,c]]$ of its
	equation\\
	\hth$a x + b y + c = 0,$\\
	given between double brackets.
	These coefficients are not unique.  They can be replaced by\\
	\hth$[[k a, k b, k c]]$\\
	where $k$ is any real number different from 0.\\
	For the points  $P$  and lines  $l$  in finite geometry, I will use the
	same notation with single parenthesis and single brackets.

	\sssec{Definition.}
	Given a prime $p,$ if $x$ is an integer, $x \umod{p}$ denotes
	the smallest positive remainder of the division of $x$ by $p.$
	\begin{verse}
	    For instance, $28 \umod{13} = 2,$ $-5 \umod{11} = 6.$
	\end{verse}
	We observe that if  $x$ and $y$  are non negative integers less than
	$p,$ for any integers $l$ and $m$,\\
	\hth$x + l p \umod{p} = x,$ $y + m p \umod{p} = y.$

	\sssec{Definition.}
	Let $x$ and $y$ be integers.  For any integers $l$ and $m,$
	the points\\
	\hth$((x + l p, y + m p))$\\
	are called {\em equivalent points}.  A set of equivalent points is
	called a {\em point} $(x,y)$ in finite geometry.

	Let $a,$ $b,$ $c$ be integers, $a$ and $b$ not both zero.  For any
	integers $l,$ $m$  and  $n,$ the lines\\
	\hth$[[ k(a + l p), k(b + m p),$ $k(c + np) ]],$ $k \neq 0,$\\
	are called {\em equivalent lines}.  A set of equivalent lines is called
	a {\em line} in finite geometry.\\
	If  $P = (x + m p, y + n p)$ is on $l = [ a + k'p, b + m'p, c + n'p],$
	then\\
	\hth$(a + k'p)(x + m p) + (b + m'p)(y + n p) + (c + n'p) = 0$\\
	and therefore\\
	\hth$(a x + b y + c) \umod{p} = 0,$\\
	this is \ref{sec-incidence} below.  This method of reducing modulo $p$
	allows us to
	extend many of the properties of Euclidean geometry to the finite case.

	\sssec{Example.}
	Let $p = 7.$  The line $a = [[1,-1,-5]]$ is equivalent to the lines
	$a_0 = [[1,-1,-12]],$ $a_1 = [[1,-1,2]],$ $a_2 = [[1,-1,9]].$\\
	The line $b = [[1,2,-17]]$ is equivalent to the lines $b_0 =
	[[1,2,-3]],$
	$b_1 = [[1,2,-10]],$ $b_2 = [[1,2,-24]],$ $b_3 = [[1,2,-31]].$\\
	The intersection $P = ((9,4))$ of  $a$  and  $b$ is equivalent to the
	points all labelled $Q,$ $((2,4)),$ $((16,4)),$ $((2,11)),$ $((9,11)),$
	$((16,11))$.  Only one of the points equivalent to $P$ is in the
	domain\\
	\hth$0 \leq x,y < p,$ namely $((2,4)).$
	\pagebreak
	\begin{verbatim}
	^
	y
	.   .  Qb2  .   .   .   .   .   .  Qb3  .   .   .   .   .   .  Qb4  .

	.   a2  .   .   b2  .   .   .   a1  .   .   b3  .   .   .   a   .   .

	a2  .   .   .   .   .   b2  a1  .   .   .   .   .   b3  a   .   .   .

	.   b   .   .   .   .   a1  .   b2  .   .   .   .   a   .   b3  .   .

	.   .   .   b   .   a1  .   .   .   .   b2  .   a   .   .   .   .   b3
	__________________________
	.   .   .   .   a1  b   . |  .   .   .   .   a   b2  .   .   .   .   .
	                          | 
	b1  .   .   a1  .   .   . |  b   .   .   a   .   .   .   b2  .   .   a0
	                          | 
	.   .  Qb1  .   .   .   . |  .   .  Pb   .   .   .   .   .   .  Qb2  .
	                          | 
	.   a1  .   .   b1  .   . |  .   a   .   .   b   .   .   .   a0  .   .
	                          | 
	a1  .   .   .   .   .   b1|  a   .   .   .   .   .   b   a0  .   .   .
	                          | 
	.   b0  .   .   .   .   a |  .   b1  .   .   .   .   a0  .   b   .   .
	                          | 
	.   .   .   b0  .   a   . |  .   .   .   b1  .   a0  .   .   .   .   b
 x >
	\end{verbatim}
	\begin{center}
			{\em Equivalence of points and lines}.\\
	                        {\em Fig.0a, p = 7}.
	\end{center}

	In Fig.0a, I have not given those lines which are equivalent to $a$
	but have a different slope, if $R$ is any point on such a line which
	is in the lower right square it is either $Q$ or a point labelled
	$a$  or  $a_1.$\\
	In finite geometry, we do not distinguish points labelled $ a_1$  from
	those labelled  $a$ or the points labelled  $b_0$  and  $b_1$ from those
	labelled  $b.$  We have therefore Fig.0b below.
	\pagebreak
	\begin{verbatim}

	                    ^
	                    y
	                    __________________________
	                    .   .   .   .   a   b   . | 
	                                              | 
	                    b   .   .   a   .   .   . | 
	                                              | 
	                    .   .  Qb   .   .   .   . | 
	                                              | 
	                    .   a   .   .   b   .   . | 
	                                              | 
	                    a   .   .   .   .   .   b | 
	                                              | 
	                    .   b   .   .   .   .   a | 
	                                              | 
	                    .   .   .   b   .   a   . |     x >

	\end{verbatim}
	\begin{center}
		{\em Points and lines in finite Euclidean geometry}.\\
				{\em Fig.0b, p = 7}.
	\end{center}

	The point $Q = (2,4)$ is on the lines  $a = [1,6,2] = [4,3,1]$ and
	$b = [1,2,4] = [2,4,1]$ in the finite Euclidean geometry associated with
	$p = 7.$\\
	Observe that from one point on  $a,$ the others are obtained by moving
	one to the right and one up, for  $b,$ we move 2 to the right and one
	down.  Observe also what happens at the boundary using, if needed
	Fig. 0a.

	If we attempt to use the equivalence method when $p$ is not a prime,
	the situation for  6  points is typical.
	If  $a = [[1,1,-5]],$ $b = [[1,3,1]]$ and $c = [[1,3,4]],$ the points
	$P = ((2,3))$ and $Q = ((5,0))$ are both on the lines  $a$  and  $b,$
	while the lines  $a$  and  $c$  or their equivalent have no point in
	common with coordinates reduced modulo 6.

	\pagebreak
	\begin{verbatim}
	                        ^
	                        y
	                        a   .   b   .   .   c

	                        .   a   c   .   .   b

	                        .   .   bP  .   .   c

	                        .   .   c   a   .   b

	                        .   .   b   .   a   c

	                        .   .   c   .   .   bQ   x >

	\end{verbatim}
	\begin{center}
			  {\em 6  points per line}.\\
				{\em Fig.0c}.
	\end{center}

	\sssec{Comment.}
	When giving numerical examples, it is convenient to assume that $x,$
	$y,$ $a,$ $b,$ $c$ are non negative integers less than $p,$ and that the
	right most of the triplet $a,$ $b,$ $c$ which is non 0 is chosen to be
	1. This is always possible because, if $c$ is for instance different
	from $0$, then we can choose $k$ in such a way that,
	 $k c  \umod{p} = 1.$  This property requires $p$ to be a prime and was
	known, together with an algorithm to
	obtain  $k,$ by the Indian Astromomer-Mathematician Aryabatha, 5-th
	Century A. D. as well as by the Chinese, the date of the invention of
	their algorithm, called the  chiu-i, or search for 1 is not known.
	\begin{verse}
	    For instance, I will use, when $p = 11,$ $[4,3,1]$ instead of
	    $[5,1,4],$ $[4,1,0]$ instead of $[2,6,0]$ and $[1,0,0]$ instead of
		$[7,0,0].$
	\end{verse}
	The main advantage of this convention is that it insures a unique
	representation of the lines in finite geometry.

	Because {\em all computations for finite Euclidean geometry have to be
	done modulo $p$}, it is useful to have ready a table of multiples of
	$p,$ of inverses modulo $p$ and of squares modulo $p.$ Two such tables
	are given, the others should be completed.

	$\begin{array}{llcrrrrrrrrrrrrrrrrrrr}
	p =  7,\\
	 i    &\vline&  \:0&  1&  2&  3&  4&  5&  6\\
		\cline{1-9}
		p\:i  &\vline&   0&  7& 14& 21& 28& 35& 42\\
		\frac{1}{i}  &\vline&   - & 1&  4&  5&  2&  3&  6\\
		i^2  &\vline&   0&  1&  4&  2&  2&  4&  1\\

	p = 11,\\
	 i&\vline& 0&1&2&3&4&5&6&7&8&9& 10\\
		\cline{1-13}
		p\:i&\vline& 0& 11& 22& 33& 44& 55& 66& 77& 88& 99& 110\\
		\frac{1}{i}&\vline& -&1&6&4&3&9&2&8&7&5 &10\\
		i^2&\vline& 0&1&4&9&5&3&3&5&9&4&1\\

	p = 13,\\
	 i&\vline& 0&1&2&3&4&5&6&7&8&9&10& 11& 12\\
		\cline{1-15}
		p\:i&\vline& \\
		\frac{1}{i}&\vline& \\
		i^2&\vline&\\

	p = 17,\\
	 i&\vline& 0&1&2&3&4&5&6&7&8&9&10&11&12\\
		\cline{1-15}
		p\:i&\vline& \\
		\frac{1}{i}&\vline& \\
		i^2&\vline& \\
	\\[-5pt]
	 i&\vline& 13&14&15&16\\
		\cline{1-6}
		p\:i&\vline& \\
		\frac{1}{i}&\vline& \\
		i^2&\vline& \\

	p = 19,\\
	 i&\vline& 0&1&2&3&4&5&6&7&8&9&10&11&12\\
		\cline{1-15}
		p\:i&\vline \\
		\frac{1}{i}&\vline \\
		i^2&\vline \\
	\\[-5pt]
	 i&\vline& 13&14&15&16&17&18\\
		\cline{1-8}
		p\:i&\vline& \\
		\frac{1}{i}&\vline& \\
		i^2&\vline& \\
	\end{array}$

	\sssec{Theorem.}\label{sec-incidence}
	{\em A point $(x,y)$ is on a line $[a,b,c]$ if and only if\\
	\hth$(a x + b y + c) \umod{p} = 0.$}

	\begin{verse}
	    For instance, with $p = 11,$ from $((12,14))$ on the line
	$[[16,-10,4]],$ it follows that $(1,2)$ is on line $[5,1,4]$ or
	$[4,3,1].$
	\end{verse}

	\sssec{Theorem.}
	{\em There are $p^2$ points and $p^2 + p$ lines.}

	\sssec{Theorem.}
	{\em 2 distinct points determine a unique line.}

	Moreover, if $A = (A_0,A_1)$ and $B = (B_0,B_1),$ the line  $a$  through
	$A$  and  $B$ is $a = [A_1 - B_1, B_0 - A_0, A_0 B_1 - B_0 A_1].$
	\begin{verse}
	    For instance, for $p = 11,$ if $A = (9,8)$ and $B = (8,6),$
	$a = [2,10,1].$
	\end{verse}

	\sssec{Theorem.}\label{sec-lines}
	{\em 2 distinct lines have at most one point in common.\\
	Moreover, if $l = [l_0,l_1,l_2]$ and $m = [m_0,m_1,m_2],$
	let $d := l_0 m_1 - l_1 m_0,$ if $d$ is different from 0, then the
	point $P$ common to $l$ and $m$ is}\\
	\hth$P = ( \frac{l_1 m_2 - l_2 m_1}{d},$
		$\frac{l_2 m_0 - l_0 m_2}{d} ).$

	\begin{verse}
	    For instance, with $p = 11,$ if $l = [2,10,1]$ and
	$m = [9,9,1],$ $d = 5,$\\
	\hth$\frac{1}{5} \umod{11} = 9,$ $P = (1\:.\:9, 7\:.\:9) = (9,8).$
	\end{verse}

	\sssec{Definition.}\label{sec-parallel}
	If 2 lines have no points in common, they are called {\em parallel}.\\
	This will occur if  $d = 0,$ because of \ref{sec-lines}.
	\begin{verse}
	For instance, with $p = 11,$ $a = [2,10,1]$ is parallel to
	$b = [5,3,1].$
	\end{verse}

	The following figure gives also a representation of points in finite
	geometry.  The representative which is chosen is that with integer
	coordinates, non negative and less than $p$ $( 0 \leq x,y < p).$\\
	The reader is asked to ignore for now the information at the left
	of the figure.  The possible points are indicated with ".", a named
	point has its name just to the right of it.  All points on a line  $a$
	are indicated by replacing "." by $"a"$.  If 2 lines have a point in
	common, one of the 2 lines is chosen.  The other could be indicated
	by the reader, if he so desires.

	\pagebreak
	\sssec{Example.}
	\begin{verbatim}
	.         ^
	          y
	bDb       .   .   .   .   .   c   .   b   .   .   a

	.         .   .   .   c   .   .   .   .   bA  .   .

	.         .   c   .   .   .   .   aB  .   .   b   .

	.         .   .   .   .   a   .   .   .   .   .   cC

	.         b   .   a   .   .   .   .   .   cD  .   .

	.         a   b   .   .   .   .   c   .   .   .   .

	.         .   .   b   .   c   .   .   .   .   a   .

	.         .   .   c   b   .   .   .   a   .   .   .

	cDa       c   .   .   .   b   a   .   .   .   .   .

	.         .   .   .   a   .   b   .   .   .   c   .

	.         .   a   .   .   .   .   b   c   .   .   .	   x  >

	\end{verbatim}
	\begin{center}
			   {\em Points, lines and parallels}.\\
				{\em Fig. 1, p =  11}.
	\end{center}

	The points are $A = (8,9),$ $B = (6,8),$ $C = (10,7),$ the line
	$a = [10,2,1]$ is the line through  $A$  and  $B,$ it passes through the
	points (0,5), (2,6), (4,7), (6,8), (8,9), (10,10), (1,0), (3,1), (5,2),
	(7,3) and (9,4).  The line $b = [9,9,1]$ is the line through  $A$  and
	$C.$ The line $c = [3,5,1]$ passes through  $C,$ has no points in
	common with  $a$ and therefore is parallel to  $a.$  The line  $d,$
	which is not indicated on the picture, is parallel to  $b$  and passes
	 through $B.$  The point $D$  is on $c$  and $d$.\\
	As an exercise determine the coordinates of  $D$  and the other points
	of  $d.$

	\sssec{Notation.}
	To ease description of constructions in geometry, I have
	introduced the notation  $A \times B,$ for the line  $l$  through the
	distinct points $A$ and $B,$ and  $l \times m$, for the point $C$
	common to the distinct lines $l$ and $m.$

	\sssec{Comment.}
	If $p$ is very large, and the unit used for the representation is very
	 small, the Angstr\"{o}m = $10^{-8}$ cm, say, the points on a line will
	appear as we imagine them in the classical case.  But it is clear that
	they are not connected.  Connectedness is a property in classical
	Euclidean geometry, which has no counterpart in the finite case.
	Moreover, the finite case should, when fully understood, give a better
	model for a world which is atomic, whatever the smallest particle is and
	which is finite, whatever the size of the universe is.

	\ssec{Parallels, parallelograms, distance.}

	\setcounter{subsubsection}{-1}
	\sssec{Introduction.}
	Parallels have been defined in \ref{sec-parallel}.  In this section,
	I will give properties of parallel lines and define parallelograms.
	It is appropriate at this stage to define distances between points.
	In the finite case, the square of a distance is the appropriate basic
	concept, if we do not want to introduce "imaginaries".  Properties of
	the parallelogram allow us then to derive a construction for the
	mid-point of a segment.  The barycenter will be define in section
	\ref{sec-bary}.

	\sssec{Theorem.}
	{\em Given a line  $l$  and a point  $P$  not on  $l,$ there exists
	a unique line  $m$  through  $P$  parallel to $l.$\\
	Moreover, if $P = (P_0,P_1)$ and $l = [l_0,l_1,l_2]$ then\\
	\hth$m = [l_0,l_1,-(P_0 l_0 + P_1 l_1)],$}

	\sssec{Definition.}
	Given 3 points $A,$ $B,$ $C$ not on the same line, let $c$ be
	the line through $C$ parallel to the line $a$ through $A$ and $B,$ let
	$d$ be the line through $A$ parallel to the line $b$ trough $B$ and $C,$
	\{$A,$ $B,$ $C,$ $D$\} is called a {\em parallelogram}.  The lines
	$A \times C$ and $B \times D$ are called {\em diagonals}, their
	intersection is called the {\em center} of the parallelogram.

	\sssec{Comment.}
	In Euclidean geometry, opposites sides of a parallelogram are
	equal.  To generalize, I observe first, that distances in Euclidean
	geometry are always considered positive.  This is consistent with the
	distance of $A B$ equal to the distance of $B A.$  But when working
	modulo
	$p,$ we cannot introduce positive numbers, keeping the requirement that
	the product and sum of positive integers modulo $p$ is positive.  Also
	not every integer modulo $p$ has a square root hence we use the square
	of the distance instead.  To use a terminology reminiscent of that used
	in Euclid's time, I will say the square on $A B,$ for the square of the
	distance between  $A$  and  $B.$

	\sssec{Definition.}
	Given 2 points $A = (A_0,A_1)$ and $B = (B_0,B_1),$ the {\em square on}
	$A B$, denoted $(A B)^2$ is\\
	\hth$(B_0-A_0)^2 + (B_1-A_1)^2  \umod{p}.$

	\begin{verse}
	    For instance, for $p = 19,$ if $A = (8,11),$ $B = (12,9),$
	    the square on $A B$ is $(A B)^2 = (16 + 4) \umod{19} = 1.$
	    But, for $p = 13,$ the square on $A B,$ where $A = (1,2)$ and
	$B = (2,7)$
	    is $(1 + 25) \umod{13} = 0,$ therefore the square can be zero for
	    distinct points $A$ and $B.$  See \ref{sec-dist0}.
	\end{verse}

	\sssec{Definition.}
	2 segments $\{A C\}$ and $\{B D\}$ {\em are equal} if the square on $A C$
	equals the square on $B D.$\\
	I write $A C = B D$.

	\begin{verse}
	    For instance, with $p = 19,$ if $A = ( 7,11),$ $B = (11,9),$
	$C = (9,7)$
	    and $D = (9,13)$ then $(A C)^2 = 1,$ $(B D)^2 = 1,$ $(A B)^2 = 1,$
	    $(C D)^2 = 17.$ Therefore $A C = B D$.
	\end{verse}

	\sssec{Definition.}
	A point  $M$  on  the line through $A$ and $C$ such that the square on
	$A M$ is equal to the square on $C M$ is called the {\em mid-point}
	of $A C.$

	\sssec{Theorem.}
	{\em In a parallelogram \{$A,$ $B,$ $C,$ $D,$\} with $A \times B$
	parallel to $C \times D$
	and $A \times D$ parallel to $B \times C,$ the square on $A B$ is equal
	to the square on $C D,$ the square on $A D$ is equal to the square on
	$B C.$  The center  $M$
	is the midpoint of the diagonals $A \times C$ and $B \times D.$}

	Moreover, if $A = (A_0,A_1),$ $B = (A_0+B_0,A_1+B_1),$
	$D = (A_0+C_0,A_1+C_1),$\\
	then $C = (A_0+B_0+C_0,A_1+B_1+C_1),$ $M = (A_0+\frac{B_0+C_0}{2},$
	$A_1+\frac{B_1+C_1}{2}),$\\
	$(A B)^2 = (C D)^2 = B_0^2 + B_1^2,$
	$(A D)^2 = (B C)^2 = C_0^2 + C_1^2,$\\
	$(A M)^2 = (M  C)^2 = \frac{1}{4}((B_0+C_0)^2 +(B_1+C_1)^2).$

	\sssec{Example.}
	The given points are $A = (7,11),$ $B = (9,13),$ $C = (11,9).$
	$D = (9,7),$
	$M = (9,10).$  $(A B)^2 = (C D)^2 = 1,$ $(A D)^2 = (B C)^2 = 8.$
	$(A M)^2 = (M C)^2 = 5,$ $(B M)^2 = (M D)^2 = 9.$
	\pagebreak
	\begin{verbatim}

	.         ^
	          y
	.         . c . . . . . . . . . . . d a . b . .
	.Df       c . . . d . . b . . . . . a . . . . .
	.         . . . . . . . . . . . . a . d . . b c
	.         . . . . . d . . b . . a . . . . . c .
	.         . . . . . . . . . . a . . . . d c . b
	.         . . . . . . d . . bB. . . . . c . . .
	.         b . . . . . . . a . . . . . c . d . .
	.         . . . . . . . dA. . b . . c . . . . .
	.Dm       . b . . . . a . . .M. . c . . . . d .
	dDb       . . . . . a . . d . . cC. . . . . . .
	.         . . b . a . . . . . c . . . . . . . d
	.         . . . a . . . . . dD. . b . . . . . .
	.         d . a b . . . . c . . . . . . . . . .
	.         . a . . . . . c . . d . . b . . . . .
	.         a d . . b . c . . . . . . . . . . . .
	.         . . . . . c . . . . . d . . b . . . a
	.         . . d . c b . . . . . . . . . . . a .
	cDa       . . . c . . . . . . . . d . . b a . .
	.         . . c d . . b . . . . . . . . a . . .     x  >

	\end{verbatim}
	\begin{center}
		       {\em Parallelogram and mid-point of A B}.\\
				{\em Fig. 2, p =  19}.
	\end{center}

	\ssec{Perpendicularity.}

	\setcounter{subsubsection}{-1}
	\sssec{Introduction.}
	The perpendicularity of lines is defined.  Theorem \ref{sec-ortho}
	follows from the corresponding theorem in classical geometry and from
	the analytical property of perpendicular lines adapted modulo $p$.\\
	An application giving the orthocenter of a triangle is also given.

	\sssec{Definition.}
	Two lines $l = [l_0,l_1,l_2]$ and $m = [m_0,m_1,m_2]$ are
	{\em perpendicular} iff\\
	\hth  $l_0 m_0 + l_1 m_1 \umod{p} = 0.$

	\begin{verse}
	    For instance, with $p = 11,$ $l = [7,6,1]$ and $m = [3,2,1]$ are
	    perpendicular.
	\end{verse}

	\sssec{Theorem.}
	{\em If 2 lines $l$ and $m$ are perpendicular to the same line $a,$
	then they are parallel.}

	\begin{verse}
	    For instance, with $p = 11,$ $l = [9,1,1]$ and $m = [1,5,1]$ are
	    perpendicular to $a = [6,1,1]$ and are parallel.
	\end{verse}

	\sssec{Theorem.}\label{sec-ortho}
	{\em Given a triangle $A,$ $B,$ $C$ with sides $a,$ $b$ and $c,$ if
	$p$ is the perpendicular from $A$ to $a,$ $q$ is the perpendicular from
	$B$ to $b$ and $r$ the perpendicular from $C$ to $c,$ then the three
	lines $p,$ $q$ and $r$ have a point $H$ in common.}

	\sssec{Definition.}
	The lines $p,$ $q$ and $r$ of Theorem \ref{sec-ortho} are called
	altitudes, the point $H$ is called the {\em orthocenter} of the
	triangle \{$A,$ $B,$ $C$\}.

	\sssec{Example.}
	The given points are $A = (8,4),$ $B = (4,8),$ $C = (3,2),$ the sides
	are
	$a = [1,9,1],$ $b = [6,7,1],$ $c = [10,10,1],$ the altitudes are
	$p = [1,6,1],$ $q = [2,3,1],$ $r = [10,1,1]$ and the orthocenter is
	$H = (7,6).$\\
	As an exercise, indicate on the figure one of the sides and compare
	how points are derived from each other with those of the perpendicular
	line,  $c$  and  $r$ are the easiest, $b$  and  $q$ the more difficult.
	\pagebreak
	\begin{verbatim}
	.         ^
	          y
	.Dc       r   q   .   .   .   p   .   .   .   .   .

	.         p   .   .   .   .   .   .   .   q   .   r

	.Db       .   .   .   .   qB  .   p   .   .   r   .

	.         q   p   .   .   .   .   .   .   r   .   .

	.         .   .   .   .   .   .   .   qH  .   .   .

	pDp       .   .   p   q   .   .   r   .   .   .   .

	qDq       .   .   .   .   .   r   .   .   pA  .   q

	.         .   .   .   p   r   .   q   .   .   .   .

	.Da       .   .   q   rC  .   .   .   .   .   p   .

	rDr       .   .   r   .   p   .   .   .   .   q   .

	.         .   r   .   .   .   q   .   .   .   .   p	   x  >

	\end{verbatim}
	\begin{center}
			  {\em The orthocenter H of \{A, B, C\}}.\\
				{\em Fig. 3, p =  11}.
	\end{center}



\setcounter{section}{9}
\setcounter{subsection}{3}
	\ssec{Circles, tangents and diameters.}

	\setcounter{subsubsection}{-1}
	\sssec{Introduction.}
	Having the notion of distance, we can define a circle. Having a
	diameter $\{A,B\}$ we can define the tangent at $A$ as the
	perpendicular to $A \times B.$  The medians and barycenter are defined
	and the
	relation between the center of a circumcircle and the mediatrices of
	the sides is given.  The proofs depend on the following Theorem:

	Given a prime $p,$ there exists a circle ${\cal C}'$ of radius $r'$ in
	Euclidean geometry which contains representatives $P'$ of each point
	$P$ of a circle ${\cal C}$ of radius $r := r' \umod{p}$ in finite
	Euclidean geometry.

	\begin{verse}
	    For instance, when $p = 19,$ the circle ${\cal C}$ of radius $1$
	    contains the points $(0,1),$ $(3,7),$ $(2,4),$ and therefore the
	    points $(7,3),$ $(3,-7),$ $(7,-3),$ $(-3,7),$ $(-7,3),$ $(-3,-7),$
	    $(-7,-3),$ $(4,2),$ $(2,-4),$ $(4,-2),$ $(-2,4),$ $(-4,2),$
	    $(-2,-4),$ $(-4,-2),$ $(1,0),$ $(0,-1),$ $(-1,0)$ altogether $20$
	     points.\\
	    $0^2+75^2 = 60^2+45^2=21^2+72^2 = 75^2$, therefore, in Euclidean
	    Geometry, $((0,75)),$ $(60,45)),$ $((21,72))$ are on the circle 
	    ${\cal C}'$  of radius $r' = 75,$ moeover $r' \umod{19} = -1,$
	    $60 \umod{19} = 3,$ $45 \umod{19} = 7,$ $21 \umod{19} = 2,$
	    $72 \umod{19} = -4.$  Appropriate change of signs give
	    the other points, for instance ((-60,45)) corresponds to (-3,-7) is
	    also on ${\cal C}'$.
	\end{verse}

	We can also replace in the Theorem just quoted the radius  $r'$  by the
	radius square $r_2' = r'^2.$\\
	The following solutions are especially attractive, because the points
	on the circle in the Euclidian plane are also the representatives in
	the finite Euclidean plane.\\
	$\begin{array}{crrl}
	\hth	&p & r_2'&{\rm points}\\
		\cline{2-4}
		& 5&    1& (0,1)\\
		& 5&    2& (1,1)\\
		& 5&    3& (2,2)\\
		& 5&    4& (0,2)\\
		& 7&    5& (1,2)\\
		& 7&   13& (2,3)\\
		& 11&   25& (0,5), (3,4)\\
		& 13&   25& (0,5), (3,4)\\
		& 17&   65& (1,8), (4,7)
	\end{array}$

	\sssec{Definition.}
	Given a point $A$ and an integer $d,$ the points $P$ such that the
	square on $P A$ is equal to  $d$  are on a {\em circle of center}
	$A$ and {\em radius square} $d.$

	\sssec{Notation.}
	From here on, it is often more convenient to have the origin
	at the center of the figure.  We will then replace the condition\\
	\hth$       0 \leq x,y,a,b,c < p$\\
	by\\
	\hth$       -\frac{p}{2}  <  x,y,a,b,c  < \frac{p}{2}.$

	\sssec{Theorem.}
	{\em For a circle centered at the origin,\\
	if $(x,0)\:  ( x \neq 0 )$ is a point, so are $(-x,0),$ $(0,x),$
		$(0,-x),$\\
	if $(x,x) \: ( x \neq 0 )$ is a point, so are $(x,-x),$ $(-x,x),$
		$(-x,-x),$\\
	if $(x,y)$ is a point ($x \neq y$, both non zero), so are $(y,x),$
		$(-x,y),$ $(y,-x),$ $(x,-y),$
	\hth$(-y,x),$ $(-x,-y),$ $(-y,-x).$}

	\sssec{Example.}
	\begin{verbatim}
	.                               ^
	                                y
	.           .   .   .   .   .   .   .   .   .   .   .

	.           .   c   .   .   .   .   .   .   .   c   .

	.           .   .   .   .   c   .   c   .   .   .   .

	.           .   .   .   .   .   .   .   .   .   .   .

	.           .   .   c   .   .   .   .   .   c   .   .

	.           .   .   .   .   .   .A  .   .   .   .   .	x  >

	.           .   .   c   .   .   .   .   .   c   .   .

	.           .   .   .   .   .   .   .   .   .   .   .

	.           .   .   .   .   c   .   c   .   .   .   .

	.           .   c   .   .   .   .   .   .   .   c   .

	.           .   .   .   .   .   .   .   .   .   .   .

	\end{verbatim}
	\begin{center}
			      {\em Circle of center A}.
				{\em Fig. 4, p = 11}.
	\end{center}

	$A = (0,0),$ the points labelled $c$ are on a circle with center $A$ and
	with radius square 10.  The line [0,1,0] through  $A$  has no point in
	common with the circle, the line [1,-3,0] has 2 points in common with
	the circle, (3,1) and (3,-1), the line [1,-1,0] has also 2 points in
	common with the circle, (4,4) and (-4,-4).

	\sssec{Exercise.}
	Indicate on the Fig. 4 by  $r$  the points on a circle with radius
	square $3$.

	\sssec{Theorem.}
	{\em If $p+1$ is divisible by 4, there are $p+1$ points on the circle.
	Otherwise, there are $p-1$ points on the circle.}

	\sssec{Definition.}
	If a line $t$ through a point $P$ on a circle has no other
	points in common with the circle, it is called a {\em tangent} to the
	circle.

	\sssec{Theorem.}
	{\em If a line $l$ through $a$ point $P$ of a circle is not tangent to
	it, it intersects the circle at an other point $Q.$}

	\sssec{Definition.}
	A line through the center of a circle is called a {\em diameter}.

	\sssec{Theorem.}
	{\em Half of the diameters have 2 points in common with the
	circle, half of them of no points in common with the circle.}

	\sssec{Theorem.}
	{\em The tangent at a point  $A$  of a circle is perpendicular to
	the diameter passing through  $A.$}

	\sssec{Example.}
	Given the point $A = (6,6)$ and the radius square 5, the points
	labeled  $c$  are on the circle centered at  $A.$  The tangent  $t$  at
	$P = (4,7)$ is $[9,1,1].$  The point $C =  (2,3)$ is ont the tangent.
	The line $d = [3,6,1]$ is a diameter through $P.$
	\newpage
	\begin{verbatim}
	.         ^
	          y
	.         t   .   .   .   .   .   c   .   .   d   .

	dDd       d   .   .   .   .   t   .   .   .   .   .

	.         .   .   d   .   .   c   .   c   .   .   t

	.         .   .   .   .   tP  .   .   .   c   .   .

	tDt       .   .   c   .   .   .   dA  .   .   t   c

	.         .   .   .   t   c   .   .   .   d   .   .

	.         .   .   .   .   .   c   .   c   t   .   d

	.         .   d   tC  .   .   .   .   .   .   .   .

	.         .   .   .   d   .   .   c   t   .   .   .

	.         .   t   .   .   .   d   .   .   .   .   .

	.         .   .   .   .   .   .   t   d   .   .   .	   x  >

	\end{verbatim}
	\begin{center}
			  {\em Circle, tangent and diameter}.
				{\em Fig. 5, p =  11}.
	\end{center}

	\sssec{Exercise.}
	Determine the other point $Q$ on $d$ and the circle and the tangent at
	$Q.$

	\sssec{Theorem.}
	{\em If $B$ and $C$ are points on a diameter of a circle and on the
	circle and $A$ is an other point of the circle, $A \times B$ is
	perpendicular to $A \times C.$}

	\sssec{Definition.}
	The {\em medians} of a triangle are the lines joining the
	vertices to the mid-points of the opposite side.

	\sssec{Theorem.}\label{sec-tmedians}
	{\em The medians of a triangle have a point $G$  in common.}

	\sssec{Definition.}\label{sec-bary}
	The point $G$ of \ref{sec-tmedians} is called the {\em barycenter} of
	the triangle.

	\sssec{Definition.}
	The {\em anti-complementary triangle} $\{D,E,F\}$ has its side
	$E \times F$ through $A$ parallel to $B\times C$, and similarly for
	$E$ and $F.$

	\sssec{Theorem.}
	{\em The mid-points $M,$ $N,$ $O$ of the sides of the triangle
	$\{A,B,C\}$ are also the mid-points of $\{A D\},$ $\{B E\}$ and 
	$\{C F\}.$}

	\sssec{Example.}\label{sec-emedian}
	\begin{verbatim}
	.Da       ^
	          y
	.         . . . o m n . . . . . . . . . . . . .
	.         . . . n . . o . . . . . . . m . . . .
	o         . n . . . m . . . o . . . . . . . . .
	.         . . . . . . . . . . . . o . . m . . n
	.         . . . . . . mA. . . . . . . . oFnE. .
	.Dc       . . . . . . . . . . . . . . n . m . o
	.         . . o . . . . m . . . . n . . . . . .
	.         . . . . . oO. . . . nN. . . . . . m .
	.         . . . . . . . . oG. . . . . . . . . .
	m         . . . . . . n . . . . o . . . . . . m
	.         . . . . nB. . . . mM. . . . oC. . . .
	.         m . n . . . . . . . . . . . . . . o .
	.         n o . . . . . . . . m . . . . . . . .
	.Db       . m . . o . . . . . . . . . . . . n .
	.         . . . . . . . o . . . m . . . n . . .
	.         . . m . . . . . . . o . . n . . . . .
	n         . . . . . . . . . . . n mDo . . . . .
	.         . . . m . . . . . n . . . . . . o . .
	.         o . . . . . . n . . . . . m . . . . .     x  >

	\end{verbatim}
	\begin{center}
			     {\em Medians and barycenter}.
				 {\em Fig. 6a, p = 19}.
	\end{center}

	The given points are $A = (6,14),$ $B = (4,8),$ $C = (14,8).$\\
	The anti-complementary points are $D = (12,2),$ $E = (16,14),$
		$F = (15,14).$\\
	The mid-points are $M = (12,2),$ $N = (15,14),$ $O = (16,14).$\\
	The medians are $m = [16,8,1],$ $n = [8,3,1],$ $o = [13,1,0].$
	The barycenter is $G = (8,10).$

	\sssec{Exercise.}
	Indicate on Fig. 6a, the line  $F \times D$  through 2 mid-points and
	observe that $F \times D$ is parallel to  $C \times A.$

	\sssec{Definition.}
	The {\em mediatrix} of  $A B$  is the line through the mid-point
	of $A B$ perpendicular to  $A \times B.$

	\sssec{Theorem.}
	{\em The mediatrices of the sides of a triangle pass through the
	center of the circumcircle of the triangle.}

	\sssec{Example.}
	\begin{verbatim}
	.Da$      ^
	          y
	.         . q . r . . c . . p . . c . . . . . .
	.         . . . . . q r . . p . . . . . . . . .
	rDr       . . . . . . . . . rZ. . . . . . . . .
	qDq       . . . . . . . . . p . . r q . . . . .
	.         . . . . . . cA. . p . . c . . rF.Eq .
	.Dc       . . q . . . . c . p . c . . . . . . r
	.         . . r . . c q . . p . . . c . . . . .
	.         . c . . . rO. . . p qN. . . . . . c .
	.         . . . . . . . . r p . . . . q . . . .
	.         . . . . . . . . . p . r . . . . . . q
	.         . . . q cB. . . . pM. . . . cC. . . .
	.         . . . . . . . q . p . . . . . . . r .
	.         . r . . . . . . . p . q . . . . . . .
	.Db       . . . . c . . . . p . . . . c q . . .
	.         q . . . . . . r . p . . . . . . . . .
	.         . . . . q . . . . p r . . . . . . . .
	.         . c . . . . . . q p . . .Dr . . . c .
	.         . . . . . c . . . p . . q c . . r . .
	pDp       r . . . . . . c . p . c . . . . q . .     x  >

	\end{verbatim}
	\begin{center}
		     {\em Mediatrices and center of circumcircle}.
				 {\em Fig. 6b, p = 19}.
	\end{center}

	The given points are the same as in Example \ref{sec-emedian}.  The
	mediatrices
	$p = [2,0,1],$ $q = [13,14,1],$ $r = [13,1,0]$ pass through the center
	$Z = (9,16)$ of the circumcircle ${\cal C}$ of the triangle $\{A,B,C\}.$

	\sssec{Exercise.}
	Determine the radius square of the circle and check that
	$(A Z)^2 = (B Z)^2 = (C Z)^2.$\\
	Check that if  $Y$  is some point on  $q,$ $(A Y)^2 = (C Y)^2.$ 

	\sssec{Theorem.}
	{\em If $A,$ $B,$ $C$ and $D$ are points on a circle and  $A B$ is
	parallel to $C D$ then the square on $A C$ equals the square on $B D$
	and the square on $A D$ equals the square on $B C.$}


\setcounter{section}{9}
\setcounter{subsection}{4}
	\ssec{The ideal line, the isotropic points and the isotropic
	lines.}

	\setcounter{subsubsection}{-1}
	\sssec{Introduction.}
	It is now time to explain the points located at the
	left of each figure.  In classical geometry, the plane can be extended
	to contain elements which are not points but have similar properties.
	For instance, all lines which are parallel to a given line  $l$  have
	no points in common, but they all have the same direction.  A
	direction is also called a point at infinity or an ideal point.
	If we extend the Euclidean plane in this way we see that 2 points
	ideal or not determine a unique line, with the exception of 2
	ideal points.  To have no exceptions, we also introduce the line at
	infinity or ideal line, which contains all ideal points.
	This extended Euclidean plane, which is unfortunately not
	part of high school education, is a first step to the understanding
	of projective geometry.  Other notions which are known to those
	familiar with complex Euclidean geometry are the isotropic points, the
	isotropic lines and their properties.  These notions also extend to the
	finite case and, with the definition of distance used, give rise to
	real points when the prime is of the form $4k + 1$.  The distance
	between points which are not both ordinary is not defined.

	To represent points we will now use, as for lines, 3 coordinates,
	not all 0, and $(x,y,z)$ will not be considered distinct from
	$(k x, k y, k z),$ $k \neq 0$.\\
	The ordinary points $(x,y)$ will also be noted $(x,y,1)$ or
	$(k x, k y, k).$

	\sssec{Definition.}
	The {\em ideal line} is the line $[0,0,1],$ the {\em ideal points} or
	{\em directions} are the points $(P_0,P_1,0).$

	\sssec{Definition.}
	A point $P = (P_0,P_1,P_2)$ {\em is on a line} $l = [l_0,l_1,l_2]$ if\\
	\hth$       P_0 l_0 + P_1 l_1 + P_2 l_2  \umod{p} = 0.$

	\sssec{Theorem.}
	{\em All $p+1$ ideal points are on the ideal line.  There are
	$p^2 + p +1$ ordinary and ideal points and $p^2+p+1$ ordinary and ideal
	lines.}

	\sssec{Theorem.}
	{\em If 2 lines  $l$  and  $m$  are parallel, they have an ideal
	point in common or have the same direction.}

	Moreover, if $l = [l_0,l_1,l_2]$ this point is  $D_l = (l_1,-l_0,0)$
	and if $m = [m_0,m_1,m_2],$ then $d: = l_0 m_1 - l_1 m_0 = 0.$

	\sssec{Definition.}
	If 2 lines  $l$  and  $m$  are {\em perpendicular}, their ideal
	{\em points} or {\em directions are} said to be {\em perpendicular}.

	Moreover, if the direction of $l$ is  $D_l = (l_1,-l_0,0)$, that of
	$m$ is $D_m = (l_0,l_1,0).$

	\sssec{Comment.}
	The ideal points are represented to the left of the figures.
	(1,0,0) is at the top the other points are from the bottom up
	(0,1,0), (1,1,0), (2,1,0), (3,1,0), $\ldots$.

	\sssec{Example.}

	In Fig. 1, the point $D_a = (2,1,0)$ is the ideal point on
	$a = [10,2,1],$ and $c = [3,5,1],$ the point $D_b = (10,1,0)$ is the
	ideal point on $b = [9,9,1]$ and $d.$

	In Fig. 2, the points $D_a = (1,1,0),$ $D_b = (9,1,0),$
	$D_f = (17,1,0),$
	$D_m = (10,1,0)$ are respectively the ideal points on $a = [5,14,1],$
	$b = [3,11,1],$ $f = [17,15,1],$ $m = [14,12,1].$  $m$ is the mediatrix
	of $A C.$

	In Fig. 3, $D_a = (2,1,0)$ and $D_p = (5,1,0),$ $D_b = (8,1,0)$ and
	$D_q = (4,1,0),$ $D_c = (10,1,0)$ and $D_r = (1,1,0)$ are the direction
	of pairs of perpendicular lines.

	In Fig. 5, $D_d = (9,1,0)$ is the direction of the diameter
	$d = [3,6,1].$ $D_t = [6,1,0]$ is the direction perpendicular to $Dd$
	and of the tangent $t = [9,1,1].$

	\sssec{Definition.}
	The {\em isotropic points} are the ideal points $(i,0,1)$ and
	$(-i,0,1)$ where $i$ is a solution of  $i^2 + 1 = 0.$

	\sssec{Theorem.}
	{\em The isotropic points exist if $p$ is of the form $4k + 1$
	(or $p$ is congruent to 4 modulo 1), they do not, otherwise.}

	The proof of this result goes back to Euler.
	\begin{verse}
	    For instance, if $p = 5,$ $i = 2,$ if $p = 13,$ $i = 5,$
	if $p = 17,$ $i = 4.$
	\end{verse}

	\sssec{Definition.}\label{sec-drad2}
	In the extended Euclidean plane, $(X_0,X_1,X_2)$ is on the
	{\em circle} with {\em center} $(C_0,C_1,C_2)$ and {\em radius square}
	$R_2$ if
	\enumb
	\item 	$(X_0 - C_0 X_2)^2 + (X_1 - C_1 X_2)^2 = R_2 X_2^2.$

	If $X_2 = 1,$ we obtain the usual equation.
	\enume

	\sssec{Theorem.}
	{\em When the isotropic points exist, they are on each of the circles.}

	Indeed, if $X_0 = i,$ $X_1 = 1$ and $X_2 = 0,$ \ref{sec-drad2}0 becomes
	$i^2 + 1 = 0.$

	\sssec{Definition.}
	The {\em isotropic lines} are any ordinary line passing through
	an isotropic point.

	\sssec{Theorem.}
	{\em The isotropic lines are perpendicular to themselves.}

	\sssec{Theorem.}
	{\em The isotropic lines through the center of a circle are
	tangent to that circle at the isotropic point.}

	\sssec{Theorem.}\label{sec-dist0}
	{\em If $A$ and $B$ are ordinary points on the same isotropic line,
	the square on $A B$ is 0.}

	Indeed, if $A = (A_0,A_1,1)$ and $B = (B_0,B_1,1),$ the line
	$A \times B$ which is $[A_1-B_1,B_0-A_0,A_0B_1-A_1B_0]$ passes through
	$ (i,1,0)$ if\\
	\hth$       (A_1-B_1) i = A_0-B_0.$\\
	But the square on $A B$ is $(A_0-B_0)^2$ + $(A_1-B_1)^2
	       = (A_1-B_1)^2 (i^2 + 1) = 0.$

	\sssec{Comment.}
	Because of \ref{sec-dist0}, when $p$ is congruent to 1 modulo 4, it is
	possible
	for the square on $A B$ to be 0 for distinct points $A$ and $B.$

	\sssec{Example.}
	The circle ${\cal C}$ of center $A = (8,8)$ passes through $P = (4,10)$ 
	and through the isotropic points $J = (4,1,0)$ and $K = (13,1,0).$
	The isotropic lines through $A$ are $j = [5,14,1]$ and $k = [14,5,1].$
	\newpage
	\begin{verbatim}
	i         ^
	          y
	i         . . . . . . j . . . k . . . . . .
	i         . . j . . . . . . . . . . . k . .
	i         . k . . . . . c . c . . . . . j .
	cK        . . . . . k . . . . . j . . . . .
	i         . . . . . . c j . k c . . . . . .
	i         . . . j . . . . . . . . . k . . .
	i         k . . . cP. . . . . . . c . . . j
	i         . . c . k . . . . . . . j . c . .
	i         . . . . . . . . kA. . . . . . . .
	i         . . c . j . . . . . . . k . c . .
	i         j . . . c . . . . . . . c . . . k
	i         . . . k . . . . . . . . . j . . .
	cJ        . . . . . . c k . j c . . . . . .
	i         . . . . . j . . . . . k . . . . .
	i         . j . . . . . c . c . . . . . k .
	i         . . k . . . . . . . . . . . j . .
	i         . . . . . . k . . . j . . . . . .     x  >

	\end{verbatim}
	\begin{center}
		{\em Ideal line, isotropic points and isotropic lines}.
				{\em Fig. 7, p =  17}.
	\end{center}

	\ssec{Equality of angles and measure of angles.}

	\setcounter{subsubsection}{-1}
	\sssec{Introduction.}
	The definition of angles is the most difficult aspect of
	finite geometry.  To approach the subject, I will give a construction
	which obtains points on a circle which are equidistant.  If we obtain
	in this way all $2q = p+1$ or $p-1$ distinct points on the circle, then
	the smallest angle or unit angle can be defined.  The proof that this is
	always possible will be given.  The equidistance follows from that in
	Euclidean geometry, for the same construction, but to better illustrate,
	I will give an independent argument in \ref{sec-teqangles}
	Examples for $p = 11,$ 13
	and 17 have been chosen, because in these particular cases, the points
	on the circle in finite Euclidean geometry are also points on a circle
	in classical Euclidean geometry.

	\sssec{Construction.}\label{sec-ceqangles}
	Given a diameter of a circle with the points $A_0$ and
	$A_q$ on it and an other point $A_1$ on the circle, I will construct
	$A_2,$ $A_3,$  $\ldots$  , as follows.\\
	Let $C$ be the center of the circle,  $A_2$ is such that
	$A_q \times A_2$ is parallel to $C \times A_1,$ and such that
	$A_0 \times A_2$ is perpendicular to $C \times A_1$ 
	$(A_0 \times A_2$ is parallel to the tangent at $A_1$).
	Given some point $A_j$ different from $A_1,$ $A_j+1$ is such that
	$A_0 \times A_j+1$ is parallel to $A_1 \times A_j$ and
	$A_q \times A_j+1$ is perpendicular to $A_1 \times A_j.$  Using
	$j = 2,$ 3,  \ldots, we obtain $A_3,$ $A_4,$  \ldots.

	\sssec{Theorem.}\label{sec-teqangles}
	{\em The square on $A_j A_j+1$ is equal to the square on
	$A_0 A_1.$}

	The proof for $j = 1$ is as follows, let $C = (0,0),$ $A_0 = (r,0)$ and
	$A_q = (-r,0)$ and $A_j = (x_j,y_j),$ the square on $A_0 A_1$ is\\
	\hth$       (x_1 - r)^2 + y_1^2 = 2 r (r - x_1),$\\
	$A_0 \times A_2 = [ -y_2, x_2 - r, r y_2],$
	$A_q \times A_2 = [ -y_2, x_2 + r, -r y_2],$\\
	$C \times A_1 = [ y_1, - x_1, 0],$\\[10pt]
	parallelism requires
	\enumb
	\item $y_1 (x_2 + r) - x_1 y_2 = 0,$\\[10pt]
	perpendicularity requires
	\item $y_1 y_2 + x_1 (x_2 - r) = 0,$\\[10pt]
	therefore,\\
	\hth$ (A_1 A_2)^2 = (x_2 - x_1)^2 + (y_2 - y_1)^2 =$\\
	\hth$   2 (r^2 - x_1 x_2 - y_1 y_2) = 2 r (r - x_1)
		 = (a_0 A_1)^2,$\\
	because of 1.  Multiplying 0, by $x_1 (x_2-r)$ and 1, by $y_1 (x_2+r)$
	and subtracting gives\\
	\hth$       x_1 y_1 (x_2+r)(x_2-r) = x_1 y_2 y_1 y_2$\\
	or because $x_1 y_1$ is different from $0,$
	\item $x_2^2 + y_2^2 = r^2.$\\
	$A_2$ is therefore on the circle.
	\enume

	\sssec{Theorem.}
	{\em Given the construction \ref{sec-ceqangles},}
	\enumb
	\item $(A_j A_j+k)^2 = (A_0 A_k)^2,$
	\item $(j + m - k -l) \pmod{2q} = 0$ implies $A_j \times A_k$
		{\em is parallel to}\\
	\hth$       A_l \times A_m.$
	\enume

	\sssec{Definition.}
	Assume that the construction \ref{sec-ceqangles} gives all $2q$ points
	of the circle, let the direction of $A_q A_j$ be $I_j$ and that of
	the tangent at $A_q$ be $I_q,$ the set $I_0,$ $I_1,$  $\ldots$  , 
	$I_{2q}$ define a {\em scale} on the ideal line.

	\sssec{Definition.}
	Let the lines $l,$ $m,$ have directions $I_l,$ $I_m,$
	the {\em angle between} $l$ and $m$ is given by $(m - l) \pmod{2q}.$

	\sssec{Theorem.}
	{\em The sum of the angles of a triangle is  0  $\pmod{2q}$.}

	Indeed, if the directions of the sides  $a,$ $b,$ $c$ are $I_a,$
	$I_b,$ $I_c,$ the angles are  $(c - b) \pmod{2q},$
	$(a - c) \pmod{2q}$ and  $(b - a) \pmod{2q}.$

	\sssec{Theorem.}
	{\em If 2 angles of a triangle are even, the third angle is even.}

	\sssec{Definition.}
	A triangle is called {\em even} if 2 of its angles and therefore
	all its angles are even.

	\sssec{Example.}
	The points $A_0 = (10,5),$ $A_1 = (1,2),$ $A_2 = (2,1),$ $A_3 = (5,10),$
	$A_4 = (8,1),$ $A_5 = (9,2),$ $A_6 = (0,5),$ $A_7 = (9,8),$
	$A_8 = (8,9),$ $A_9 = (5,0),$ $A_a = (2,9),$ $A_b = (1,8)$ are on a
	circle centered at $C$
	with radius square 3.  These points have been obtained from
	$A_0,$ $A_1$ and $A_6$ by the construction \ref{sec-ceqangles}\\
	The angles can be determined using the scale defined by the ideal
	points $I_0 = (1,0,0),$ $I_1 = (7,1,0),$ $I_2 = (5,1,0),$
	$I_3 = (1,1,0),$
	$I_4 = (9,1,0),$ $I_5 = (8,1,0),$ $I_6 = (0,1,0),$ $I_7 = (3,1,0),$
	$I_8 = (2,1,0),$ $I_9 = (10,1,0),$ $I_a = (6,1,0),$ $I_b = (4,1,0).$
	If $i + l = j + k$ then $A(i) A(j)$ is parallel to $A(k) A(l),$
	for instance,
	$b$ or $A_8 A_7$ is parallel to $c$ or $A_9 A_6$.\\
	\newpage
	\begin{verbatim}
	.I0     ^
	        y
	cI9     .   .   .   .   .   .A3 c   b   .   .   .

	.I4     .   .   .Aa .   .   .   .   c   bA8 .   .

	.I5     .   .Ab .   .   .   .   .   .   c   bA7 .

	.I1     .   .   .   .   .   .   .   .   .   c   b

	.Ia     b   .   .   .   .   .   .   .   .   .   c

	.I2     cA6 b   .   .   .   .C  .   .   .   .   .A0

	.Ib     .   c   b   .   .   .   .   .   .   .   .

	.I7     .   .   c   b   .   .   .   .   .   .   .

	.I8     .   .A1 .   c   b   .   .   .   .   .A5 .

	.I3     .   .   .A2 .   c   b   .   .   .A4 .   .

	.I6     .   .   .   .   .   cA9 b   .   .   .   .	   x  >

	\end{verbatim}
	\begin{center}
		   {\em Angles and equidistant points on a circle}.
				{\em Fig. 8, p = 11}.
	\end{center}

	\sssec{Exercise.}
	Obtain using the construction \ref{sec-ceqangles}, the
	point $A_3$ from the point $A_2.$

	\sssec{Theorem.}\label{sec-bisect}
	{\em If the angle of 2 lines  $l$  and  $m$  is even there are
	two lines  $b_1$  and  $b_2$  which form an equal angle with  l  and  m.
	The lines  $b_1$  and  $b_2$ are perpendicular.}

	\sssec{Definition.}
	The lines  $b_1$  and  $b_2$  of \ref{sec-bisect} are called
	{\em bisectrices.}

	\sssec{Theorem.}
	{\em If a triangle is even, there exist 4 points $C_0,$ $C_1,$ $C_2,$
	$C_3,$ which are on 3 bisectrices, each passing by a different vertex
	of the triangle.}

	More precisely, if the 3 bisectrices $d,$ $e,$ $f,$ which pass
	respectively Through $A_0,$ $A_1,$ $A_2,$ are such that
	\enumb
	\item $(angle(d,a) + angle(e,b) + angle(f,c)) \pmod{2q} = q$\\
	then\\
	\hth	$d,$ $e$ and $f$ have a point in common.
	\enume

	\sssec{Definition.}
	The 4 points $C_0,$ $C_1,$ $C_2,$ $C_3$ are called {\em center of
	the tangent circles}.

	\sssec{Theorem.}
	{\em There exist a circle with center $C_i$ tangent to each of
	the sides of the triangle.}

	\newpage
	\sssec{Example.}
	\begin{verbatim}
	.                         ^
	                          y
	.         r . c . . . . .B0 . . . a . . . r
	c         . . . . cA1 . . . . . . . . . . .
	.         . . . . . . c . . . . . . a . . .
	r         . . . . . a . . c . . . . . . . .
	.         . .B1 . . . . . . . c . . . a . .
	.         . . . . . r a . . . . r c . . . .
	.         . . . . . r . . . . . r . . c a .
	.         . . . . . . . a . . . . . . . . cA0
	a         . c . . . . . . . . . . . . . . a     x  >
	.         . . . c . . . . a . . . . . . . .B2
	.         a . . . . c . . . . . . . . . . .
	.         r . . . . . . c . a . . . . . . r
	r         . aA2 . . . . r . c . . . . . . .
	.         . . . . . . r . . . a c . . . . .
	.         . . a . . . . . .I. . . . c . . .
	.         . . . . . . r . . . r a . . . c .
	.         c . . a . . . r . r . . . . . . .

	\end{verbatim}
	\begin{center}
		       {\em Bissectrices and inscribed circle}.
				 {\em Fig.9, p = 17}.
	\end{center}

	The given triangle is $A_0 = (8,1,1),$ $A_1 = (-4,7,1),$
	$A_2 = (-7,-4).$
	Its sides are $a = [2,1,1],$ $b = [-7,4,1],$ $c = [5,-7,1].$
	The bisectrices meet the circle at $B_0 = (-1,8),$ $B_1 = (-7,4),$
	$B_2 = (8,-1).$  They are $d = [8,3,1],$ $e = [-3,3,1],$ $f = [-4,3,1]$
	 and have
	the point $I = (0,-6)$ in common.  The tangent circle  $r$  has radius
	square  5.  Its points of contact with  $a,$ $b,$ $c$ are respectively
	(2,-5), (-3,3), (1,-4).  Only  $a$  and  $c$ are given on the Figure,
	not
	to clutter it.  The other centers of tangent circles are (-1,4), (3,-3)
	and (-2,5).

	\sssec{Exercise.}
	Determine that $b$ is tangent to the circle $r$ and that $(-1,4)$ is
	indeed a center of a tangent circle.

	\ssec{Finite trigonometry.}

	\sssec{Definition.}
	If $r = 1$ and the construction \ref{sec-ceqangles} gives all $2q$
	points $A_j = (x_j,y_j)$ of the circle with radius square 1.
	I define\\
	\hth$       {\mathit sin}(2j) := y_j, {\mathit cos}(2j) := x_j.$

	\sssec{Comment.}
	Because, in general, several points $A_1$ can be chosen, there
	are several distinct but related trigonometric functions  sine  and
	cosine.  Each corresponds to a different choice of the unit angle.
	This is similar to the real case in which many different units are used,
	those with angles in radians, degrees, grades, for instance.

	\sssec{Comment.}
	I will develop the properties of the trigonometric functions and
	obtain functions which can be considered as an analogue of the
	hyperbolic functions.  An efficient method to obtain them for
	large $p$ will also be given.

	\sssec{Example.}
	For $p = 11,$\\
	$\begin{array}{cccc}
	i& A_i&angle(i)-180 i&      A'_i\\
	\cline{1-4}
	1& (-4,-3)&  \:36^{\circ}87& ((-4,-3))\\
	2& (-3,-4)&  \:73^{\circ}74& ((\frac{7}{5},\frac{24}{5}))\\
	3& ( 0, 5)&  110^{\circ}61& ((\frac{44}{25},-\frac{117}{25}))\\
	4& ( 3,-4)&  147^{\circ}48& ((-\frac{527}{125},\frac{336}{125}))\\
	5& ( 4,-3)&  184^{\circ}35& ((\frac{3116}{625},\frac{237}{625}))\\
	6& (-5, 0)&  221^{\circ}22& ((-\frac{11753}{3125},-\frac{10296}{3125}))
	\end{array}$

	If $A_i = (-4,-3),$ then $4^2 + 3^2 = 5^2,$ $cos(i) = -\frac{4}{5}
	= -3$ and\\
	 $sin(i) = -\frac{3}{5} = -5.$

	For $p = 13,$\\
	$\begin{array}{cccc}
	i& A_i&angle(i)-180 i&      A'_i\\
	\cline{1-4}
	1& (-3,-4)&	\:53^{\circ}13&  ((-3,-4))\\
	2& (-4,-3)&	 106^{\circ}26&  ((-\frac{7}{5},\frac{24}{5}))\\
	3& ( 0, 5)&	 159^{\circ}39&  ((\frac{117}{25},-\frac{44}{25}))\\
	4& ( 4,-3)&	 212^{\circ}52&  ((-\frac{527}{125},-\frac{336}{125}))\\
	5& ( 3,-4)&	 265^{\circ}65&  ((\frac{237}{625},\frac{3116}{625}))\\
	6& (-5, 0)&	 318^{\circ}78&  ((\frac{11753}{3125},-\frac{10296}
	{3125}))
	\end{array}$

	For $p = 17,$\\
	$\begin{array}{crrc}
	i& A_i&angle(i)-180 i&      A'_i\\
	\cline{1-4}
	-1 &( 8, 1)& -7^{\circ}125&(( 8, 1))\\
	 1 &( 8,-1)&\:7^{\circ}250&(( 8,-1))\\
	 3 &(-4, 7)& 21^{\circ}375&((\frac{488}{64},-\frac{191}{65}))\\
	 5 &( 7,-4)& 35^{\circ}625&((\frac{27688}{4225},-\frac{19841}{4225}))\\
	 7 &(-1, 8)& 49^{\circ}875&((\frac{1426888}{274625},-\frac{1692991}
	{274625}))\\
	 9\\
	11\\
	13\\
	15
	\end{array}$

	\sssec{Exercise.}
	Continue the last table obtaining the missing values.

	\sssec{Exercise.}
	Obtain trigonometric functions for $p = 11$ and check the familiar
	identities
	\enumb
	\item$sin^2(x) + cos^2(y) = 1,$
	\item$sin(x+y) = sin(x) cos(y) + sin(y) cos(x),$
	\item$cos(x+y) = cos(x) cos(y) - sin(x) sin(y),$
	\enume

	\sssec{Notation.}
	In a finite field there is no ambiguity in defining $\pi := 2q.$

{\tiny\footnotetext[0]{G19.TEX [MPAP], \today}}
\setcounter{section}{9}
	\section{Axiomatic}

\setcounter{subsection}{-1}
	\ssec{Introduction to Axiomatic.}
	The axiomatic study of Geometry has a long history, starting with
	Euclid. Among the main earlier contributors are Giovanni Saccheri
	(1667-1733), Karl Gauss (1777-1855),
	Janos Bolyai (1802-1860), Nikolai Ivanovich Lobachevsky
	(1792-1856), de Tilly (1837-1906)\footnote
	{Blumenthal considers than in the paper of 1892, de Tilly makes a
	fundamental contribution by introducing n-point relations to
	characterize a space metrically.},
	Pieri, Carl Menger, Oswald Veblen (1880-1960), William Young
	(1863-1942),
	Julius Dedekind, Frederigo Enriques (1871-1946),
	I. Schur, David Hilbert (1862-1943), Marshall Hall and Alfred Tarski
	(1901-1983).\\
	To obtain a clear understanding of the relation between the synthetic
	and the algebraic point of view, an important step was the
	realization of the connection between the Axiom of Pappus and the
	commutativity of multiplication, first considered by Schur, in 1898,
	then by Hilbert in 1899 (p. 71), by Artin in 1957 and many others,
	see Artzy (1965), Hartshorne (1967).\\
	A detailed history of the developments concerning Finite Geometry can be
	obtained from the monumental work of Dembowsky, 1968 and Pickert
	(Chapter 12).\\
	For some authors, the word projective geometry as moved away from its
	original meaning, to become a synonym of incidence geometry. I will
	not follow that practice.

	What follows can be used to obtain a justification of the relation
	between the synthetic and algebraic axioms of Chapter II.
	With the exception of the proof of associativity and
	commutativity of addition, I have borrowed heavily from Artzy's book,
	which contains proofs not given here,
	increasing the formalism to prepare for eventual computarization.

	The axioms will progress from those of\\
	the perspective plane, with $(\Sigma,+,\cdot)$ a ternary ring,
	($A*B*C$) and $(\Sigma,+)$, $(\Sigma-\{0\},\cdot)$ are loops,\\ 
	to Veblen-Wedderburn plane, with $(\Sigma,+,\cdot)$ a quasifield,
	(linear, right distributivity) and $(\Sigma,+)$ an Abelian group,\\
	to Moufang plane, with $(\Sigma,+,\cdot)$ an alternative division ring
	(left distributivity, right and left inverse property),\\
	to Desarguesian plane, with $(\Sigma,+,\cdot)$ a skew field,
	(associativity of multiplication),\\
	to Pappian plane, with $(\Sigma,+,\cdot)$ a field, ( commutativity of
	multiplication),\\
	to Separable Pappian plane, with $(\Sigma,+,\cdot)$ an ordered field,\\
	to Continuous Pappian plane with $(\Sigma,+,\cdot)$ the field of 
	reals.\\
	The definitions of Desargues and Pappus configurations, given in Chapter
	II, will not be repeated here.

	\ssec{The Perspective Plane.}

\setcounter{subsubsection}{-1}
	\sssec{Introduction.}
	Marshall Hall and D.T. Perkins independently succeeded to construct
	an algebraic structure, called ternary ring, \ref{sec-dternary} to
	coordinatize \ref{sec-dprojter} the perspective plane\ref{sec-dpersp}.
	Theorem \ref{sec-ttern} shows that the first 4 conditions of the
	definition of a ternary ring are associated with the incidence
	property \ref{sec-tincidpersp}.3 and the others with Theorem
	\ref{sec-tincid1}. Theorem \ref{sec-tternloop} proves that the set of
	the ternary ring is a loop under addition and multiplication.

	\sssec{Axioms. [Of Allignment]}
	Given a set of elements called {\em points} and a set of elements called
	{\em lines} with the relation of {\em incidence}, such that
	\enumb
	\item 2 points are incident to one and only one line.
	\item 2 lines are incident to one and only one point.
	\item there exists at least 4 points, any 3 of which are not collinear,
	\enume
	we say that the {\em axioms of allignment} are satisfied.

	The terminology is that of Seidenberg, 1962, p. 56.

	\sssec{Definition.}\label{sec-dpersp}
	A {\em perspective} plane is a set of points and lines satisfying the
	axioms of alignment.  It is also called a {\em rudimentary projective
	plane}. (Artzy, p. 201.)

	\sssec{Theorem.}
	{\em Duality is satisfied in a perspective plane.}

	Menger gives a self dual set of equivalent axioms.

	\sssec{Definition.}
	Given a point $P$ and 2 lines $a$ and $b$ not incident to $P,$
	a {\em perspectivity} $\Pi(P,a,b)$ is the correspondance between
	$A_i\incid a$ and $B_i\incid b,$ with $B_i := (P\times A_i)\times b.$\\
	$\Pi^{-1}(P,a,b) := \Pi(P,b,a)$ is the inverse correspondance which
	associates to $B_i,$\\
	$A_i = (P\times B_i) \times a.$\\
	I will also use the notation $\Pi(P,A_i,B_i)$.\\
	A {\em projectivity} is a perspectivity or the composition of 2 or more
	perspectivities.

	\sssec{Theorem.}
	{\em $\Pi$ is a bijection.}

	\sssec{Definition.}
	Given a line $m$, we say that {\em l is m-parallel to $l'$} iff
	$l,$ $l'$ and $m$ are incident and we write\\
	\hth$l\paras{m}l'$ and $I^m_l := l\times m$.\\
	$I^m_l$ is called the {\em m-direction} of $l$.

	\sssec{Definition.}\label{sec-dtrpr}
	Given a line $m$, 2 points $A$ and $A'$, not on $m$ and a point $B$
	neither on $m$ nor on $a := A\times A'$, the
	{\em translation} ${\cal T}^{m,B}_{AA'}$ is the transformation which
	associates to $I,\:I$ if $I\incid m$ and to points $P$ neither on $a$,
	nor on $m$ the point
	$P' := (P\times I^m_{A\times A'})\times (A'\times I^m_{A\times P}),$
	and to $C\incid a$, $C' := (I^m_{B\times C}\times B')
	\times a,$ where $B' := {\cal T}^{m,B}_{AA'}(B).$

	\sssec{Definition. [Marshall Hall]}\label{sec-dternary}
	\vspace{-18pt}\hspace{182pt}\footnote{1943, also unpublished work by
	D. T. Perkins}\\[8pt]
	$(\Sigma,*)$ is a {\em ternary ring} iff $\Sigma$ is a set and $*$
	is an operation which associates to an ordered triple in the set
	an element in the set satisfying the following properties
	\enumb
	\item$A*0*C = C,$
	\item$0*B*C = C,$
	\item$1*B*0 = B,$
	\item$A*1*0 = A,$
	\item$A*B*X = D$ {\em has a unique solution X,}
	\item$B_1\neq B_2\implies X*B_1*C_1 = X*B_2*C_2$ {\em has a unique
	solution X,}
	\item$A_1\neq A_2\implies A_1\times X\times X' = D_1$ and $A_2\times
	X\times X' = D_2$
	{\em have a unique solution (X,X')}.
	\enume

	\sssec{Theorem.}
	{\em $X * 1 * 1 = 0$ has a unique solution $X$.}

	Proof: \ref{sec-dternary}.0 implies $X * 0 * 0 = 0,$ $1\neq 0$,
	the Theorem follows from \ref{sec-dternary}.5.

	\sssec{Definition.}\label{sec-dprojter}
	A perspective plane can be coordinatized as follows, (Fig. 20a')\\
	H0.0.\hti{3}$Q_0$, $Q_1$, $Q_2$, $U$, 4 points, no 3 of which are
		collinear,\\
	D0.0.\hti{3}$q_0 := Q_0 \times Q_1,$ $q_1 := Q_2 \times Q_0,$
		$m := q_2:= Q_1 \times Q_2$,\\
	D0.1.\hti{3}$v := Q_2 \times U$, $i := Q_0 \times U$,
		$V := i \times q_2$, $I := v \times q_0$,\\
	D0.2.\hti{3}$u := V \times I$.\\
	Let $\Sigma$ be the set of points on $q_2$, {\em distinct} from
	$Q_2$. Define $0 := Q_1$, $1 := V$.

	The point $Q_2$ is represented by $(\infty)$, $\infty$ being
	a new symbol.\\
	The points $Q$ on $q_2$, distinct from $Q_2$ are {\em represented} by
	the element $Q$ in $\Sigma$, placed between parenthesis, Q = (Q).\\
	A point $P$ not on $q_2$ is {\em represented} by a pair of elements
	$(P_0, P_1)$ in $\Sigma$ defined by (Fig. 20a'')\\
	\hth$P_0 := (((((P \times Q_2) \times i) \times Q_1) \times v) \times
	Q_0) \times q_2,$\\
	\hth$P_1 := (((P \times Q_1) \times v) \times Q_0) \times q_2.$

	In particular, if a point $A$ is on $q_0$, then its second coordinate
	$A_1 = 0$, we represent its first coordinate by $A$, if a point $C$ is
	on $q_1$, then its first coordinate $C_0 = 0$, we represent its second
	coordinate by $C$. Points on $v$ have first coordinate 1.  

	The line $q_2$ is represented by $[\infty]$,\\
	a line $l_0$ through $Q_2$ distinct from $q_2$ is {\em represented} by
	$[A]$, with $a \times q_0 = (A,0),$\\
	a line $m$ not through $Q_2$ is {\em represented} by the pair
	$M_0,M_1]$, where $(M_0)$ is the representation of the point $m \times
	q_2$ and where the point $m\times q_1$ on $q_1$ is $(0,M_1)$.

	Let $P := ((A,0)\times Q_2)\times ((B)\times (0,C)) = (A,Y).$
	$Y$ is a function of $A$, $B$ and $C$ which we denote by\\
	\hth$Y := A*B*C$.

	\sssec{Theorem.}
	{\em There is a bijection between the points $(A,0)$ on $q_0,$ $(0,A)$
	on $q_1 $ and $(A)$ on $q_2.$}

	Proof:
	We use the perspectivity $\Pi(Q_2,q_0,i),$ followed by $\Pi(Q_1,i,q_1),$
	or $\Pi(Q_1,q_1,u),$ followed by $\Pi(Q_0,u,q_2).$

	\sssec{Comment.}
	If $P = (P_0,P_1)$, all points on $P\times Q_2$ have the same first
	coordinate, $P_0$, in particular, $(P\times Q_2)\times q_0 = (P_0,0)$,
	all points on $P\times Q_1$ have the same second coordinate, $P_1$, in
	particular, $(P\times Q_1)\times q_1 = (0,P_1)$.

	In Euclidean Geometry, if $q_0$ is the $x$ axis, $q_1$ is the $y$ axis
	$q_2$ is the ideal line and $U = (1,1,1),$ $(A,B)$ corresponds to
	$(A,B,1),$ $(A)$ to $(1,A,0)$ which is the direction of lines with
	slope $A,$ $(\infty)$ to $(0,1,0)$ which is the direction of $y$ axis.
	The slope of the line joining the origin to $(A,B,1)$ is $\frac{B}{A}$.

	\sssec{Theorem.}\label{sec-tincidpersp}
	{\it The incidence, noted "$\iota$" satisfies},
	\enumb
	\item$(Q)\incid [\infty]$,
	\item$(P_0)\incid [P_0,P_1]$,
	\item$(P_0,P_1)\incid[P_0])$,
	\item$(P_0,P_1)\incid[M_0,M_1]$ iff $P_1 = P_0 * M_0 * M_1$.
	\enume

	\sssec{Theorem.}\label{sec-tincid1}
	\enumb
	\item$(R_0,R_1)\incid(Q_2\times(A,0))\implies R_0 = A$,
	\item$(S_0,S_1)\incid(Q_1\times(0,C))\implies S_1 = C$,
	\item$X := v\times (Q_0\times (B)) \implies X = (1,B)$,
	\item$(Y_0,Y_1)\incid(Q_0\times V)\implies Y_0 = Y_1$.
	\enume

	\sssec{Theorem.}\label{sec-ttern}
	{\em The pespective plane as coordinatized in} \ref{sec-dprojter}
	{\em satisfies the properties of a ternary ring.  In particular}
	\enumb
	\item {\em the unique solution $X$ of $A*B*X = D$ is the second
	coordinate of  the point} $((A,D)\times (B))\times q_1,$
	\item {\em with $B_1\neq B_2$, the unique solution $X$ of $X*B_1*C_1
	 = X*B_2*C_2$ is the first coordinate of the point}
	$((0,C_1)\times (B_1))\times ((0,C_2)\times (B_2)),$
	\item {\em with $B_1\neq B_2$, the unique solution $(X,X')$ of
	$A_1\times X\times X' = D_1$ and $A_2\times X\times X' = D_2$ is given
	by}\\
	\hth$X := ((A_1,D_1)\times (A_2,D_2))\times q_2,$
		$X' := ((A_1,D_1)\times (A_2,D_2))\times q_1,$
	\enume

	Proof:
	For 0. to 3. of \ref{sec-dprojter}, we consider the points\\
	$(R_0,R_1) := ((A,0)\times Q_2)\times ((0)\times (0,C)) = (A,A*0*C)
	= (A,C)$,\\
	$(S_0,S_1) := ((Q_0\times Q_2)\times ((B)\times (0,C)) = (0,0*B*C) =
	(0,C)$,\\
	$X := ((1,0)\times Q_2)\times ((B)\times (0,0)) = (1,1*B*0) = (1,B),$\\
	$(Y_0,Y_1) := ((A,0)\times Q_2)\times ((V)\times (0,0)) = (A,A*1*0)
	= (A,A)$.

	\sssec{Theorem.}
	{\em The pespective plane satisfies also the properties:}
	\enumb
	\item$X*B*C = D$ {\em has a unique solution, the first coordinate of 
	the point} $((0,D)\times q_1)\times ((B)\times(0,C)),$
	\item$A*X*C = D$ {\em has a unique solution, the coordinate of 
	the point} $((A,D)\times (0,C))\times q_2,$
	\enume

	Proof: For 0, $(X,Y)\incid[B,C],$ $(X,Y)\incid[0,D],$ therefore
	$Y = X*B*C = X*0*D = D$.  For 1, $(A,D)\incid[X,C]$ therefore
	$D = A*X*C$.

	\sssec{Example.}\label{sec-eprojter}
	$Q_0 = (0,0),$ $Q_1 = (0),$ $Q_2 = (\infty),$
	$U = (1,1)$, $V = (1),$ $I = (1,0)$, $u = [1,z]$ with $1*1*z = 0.$\\
	$q_0 = [0,0],$ $q_1 = [0],$ $q_2 = [\infty],$
	$v = [1],$ $i = [1,0],$\\
	Let (see Fig. 20a')\\
	D0.3,\hti{3}$J := u \times q_1$ $j := U \times Q_1,$
		$W := j \times q_1,$ $w := J \times Q_1$,\\
	D0.4.\hti{3}$T := v \times w$, $t := V \times W,$
		$R := t \times q_0,$ $r := T \times Q_0,$ $S := r \times q_2,$\\
		then, with S = (S),\\
	$J = (0,S)$, $j = [0,1],$ $W = (0,1),$ $w = [0,S],$
	$T = (1,S),$ $t = [1,1],$ $r = [S,0],$ $R = (y,0)$ with $y*1*1 = 0.$

	\sssec{Definition.}
	The dual coordinatization can also be chosen.  I will use the
	subscript $d$ to indicate the dual representation,\\
	The notation in the preceding example is chosen to allow the
	dual coordinatization using as elements of $\Sigma$ the lines through
	a given point $(\infty).$  We choose $(\infty)_d$ as $Q_2$.
	The line $q_2$ is represented by $[\infty]_d$,
	the line $l_0 := Q_2 \times (L_0,0)$ is {\em represented} by
	$[l_0]_d = [L_0],$ the line $n$ not through $Q_2$ is {\em represented}
	by $[n_0,n_1]_d,$  with\\
	\hth$n_0 := (((((l \times q_2) \times I) \times q_1 \times V) \times
	q_0) \times Q_2,$\\
	\hth$n_1 := (((l \times q_1 \times V) \times q_0) \times Q_2.$\\
	in this case $q_0 = [0,0]_d,$ $q_1 = [0]_d,$ $q_2 = [\infty]_d,$
	$u = [1,1]_d,$, $w = [0,1]_d,$ $i = [1,0]_d,$ but\\
	$j = [0,R]_d,$ with $t\times q_0 = (R,0)_d$ and $t = [1,R]_d.$
	The representation of points is done dually as in \ref{sec-dprojter},
	with $N$ {\em represented} by $(N_0,N_1)_d$, with
	$N\times Q_2 = [N_0]_d$ and $N\times Q_1 = [0,N_1]_d.$
	
	\sssec{Theorem.}\label{sec-tternloop}
	\enumb
	\item $(\Sigma,+)$ {\em is a loop, with 0 as neutral element,}
	\item$(\Sigma-\{0\},\cdot)$ {\em is a loop with 1 as neutral element.}
	\enume

	Proof:
	For the addition, the neutral element property follows from
	\ref{sec-dternary}.1 with $B = 1$ and from .3. The solution property
	follows from \ref{sec-dternary}.4 and .6 with $B = 1$.\\
	For the multiplication, the neutral element property follows from
	\ref{sec-dternary}.3 and 4. The solution property follows from
	\ref{sec-dternary}.4 and .5 with $C = 0$.

	\sssec{Theorem.}
	{\em If the number of elements in $\Sigma$ is a small number $n$,
	\enumb
	\item {\em If $n$ = 2,3,4,5, there is only one perspective plane,}
	\item {\em If $n$ = 6, there is no perspective plane,}
	\item {\em If $n$ = 14,21,22,30,33,38,42,46,54,57,62,66,69,70,77,
			78,86,93,94,\ldots,}
	\item {\em If $n$ = 10, there is no perspective plane,}
	\enume

	0, is easily settled, see II\\
	1, originates with the problem of the 36 officers, Euler (1782),
		was settled by Tarry (1900),\\
	2, depends on the next Theorem,\\
	3, has a long history, and was finaly proven, using computers, by
		Lam, Thiel and Swiercz (1989), see also Lam (1991).

	\sssec{Theorem. [Bruck and Ryser]}
	{\em If $n\equiv 1,2 \umod{4}$ and there are no integers $x,y$
	such that $x^2+y^2=n$ then there are no perpective plane of order $n$.}

	\sssec{Notation.}
	\hth$A + B := A*1*B$,\\
	\hth$A\cdot B := A*B*0$,\\
	\hth$A + (A\vdash B) = B,$ $(B\dashv A) + A = B.$\\
	\hth When $A\neq 0,$ $A \cdot (A\setminus B) = B,$ $(B/ A) \cdot A = B.$

	In the Euclidean case, the line joining the point $(A,A+B)$ to the
	point $(0,B)$ has slope 1 and the slope of the line joining $Q_0$ to
	$(A,A.B)$ is $C = A.B$.

	\sssec{Definition.} A ternary ring $(\Sigma,*)$ is {\em linear} iff
	for every $A$, $B$, $C$ in the set\\
	\hth$(A*B*0)*1*C = A*B*C.$

	\sssec{Theorem.}
	{\em If a ternary ring is linear then}\\
	\hth$A*B*C = A\cdot B + C.$

	\sssec{Axiom. [Fano]}
	The diagonal points of every quadrangle are not collinear.

	\sssec{Axiom. [N-Fano]}
	The diagonal points of every quadrangle are collinear.

	\sssec{Definition.}
	A {\em Fano plane} is a perspective plane which satisfies the N-Fano
	axiom.

	\sssec{Theorem.}
	In a Fano plane $A+A = 0$.

	Proof: For the quadrangle $Q_0 = (0,0),$ $X_A = (A,0),$ $Y_A =(0,A),$
	$A0 = (A,A),$ 2 of the diagonal points are on $q_2,$ therefore the
	third diagonal point is $V = (Q_0\times A0)\times q_2$, therefore
	$(A,A*1*A)$ coincides with $X_A$ and $A+A = 0$.
	
	\sssec{Exercise.}
	\enumb
	\item Prove that in a Fano plane $(A\:*\:B)\:*\:(A\cdot B) = 0.$
	\item Determine a subset of quadrangles with collinear diagonal points
		which justify the preceding property in a perspective plane.
	\item Same question for the property $A\:+\:A = 0.$
	\enume

	\sssec{Definition.}
	Two triangles $\{APQ\}$ and $\{A'P'Q'\}$ are {\em $m$-parallel} iff\\
	\hth$A\times P\paras{m}A'\times P'$,
	$A\times Q\paras{m}A'\times Q'$, $P\times Q\paras{m}P'\times Q'$.

	\sssec{Theorem.}
	{\em If $A\incid l$, $l' := A'\times I^m_l$ and $P\incid l$ then
	$P'\incid l'$.}

	In general, a line $n$ not through $A$ is not transformed into a line.
	For this to be so, if $P\incid n$ and $Q\incid n$, we want
	$P' := {\cal T}^m_{AB}(P)$ and $Q' := {\cal T}^m_{AA'}(Q)$ to be
	collinear with $I^m_{P\times Q}$.  This suggest the following 
	Definition.

	\sssec{Axiom. [Of Desargues]}
	In a perspective plane, given any 2 triangles $\{A_i,a_i\}$ and
	$\{B_i,b_i\}$,\\
	let $c_i := A_i\times B_i,$ and $C_i := a_i\times b_i,$
	incidence$(c_i,C) \implies$ incidence$(C_i,c)$.\\
	$C$ is called the {\em center}, $c$ is called the {\em axis} of the
	configuration.\\
	I write Desargues$(C,\{A_i\},\{B_i\};\langle C_i\rangle,c)$.

	\sssec{Axiom. [Elated Desargues]}
	The {\em Elated Desargues axiom} is the special case when we restrict
	Desargues' axiom to the case when the axis $c$ passes through the
	center $C$ of the configuration.  More specifically,
	$C\incid c$, and for the 2 triangles $\{A_i\}$ and $\{B_i\}$,\\
	let $C_i := (A_{i+1}\times A_{i-1})\times(B_{i+1}\times B_{i-1}),$\\
	$c_i := (A_i\times B_i),$ $c_i\incid C$, $i = 0,1,2,$
	incidence($A_0\times A_j,B_0\times B_j,c)$, $j = 1,2$,\\
	$\implies$ incidence($A_1\times A_2,B_1\times B_2,c)$. We write\\
	\hth Elated-Desargues$(C,\{A_i\},\{B_i\};\langle C_i\rangle,c).$

	The terminology comes from that in projective geometry, which calls
	elation, a collineation with an axis of fixed point and a center of
	fixed lines, with the center on the axis. This axiom is also called
	the minor Desargues axiom, see for instance Artzy, p. 210.

	\sssec{Theorem.}
	
	{\em Given 2 triangles $\{A_i\}$ and $\{B_i\}$, let
	$C_i := (A_{i+1}\times A_{i-1})\times(B_{i+1}\times B_{i-1}),$
	$C_i := A_i\times B_i,$ and $C := c_1\times c_2,$\\
	\hth$\langle C_i,c\rangle$ and $C\incid c\implies c_0\incid C.$
	We write}\\
	\hth Elated-Desargues$^{-1}(c,\{A_i\},\{B_i\};\langle \un{c_0},c_1,c_2
	\rangle,C)$

	Proof: Desargues$(C_0,\{A_1,B_1,C_2\},\{A_2,B_1,C_1\};\langle B_0,A_0,C
	\rangle,c)$.

	\ssec{Veblen-Wedderburn Planes.}

	\sssec{Definition.}
	A {\em Veblen-Wedderburn plane} is a perspective plane for which
	the elated Desargues axiom is satisfied on a specific line of the
	plane.

	\sssec{Comment.}
	In all the construction that follow, H0.0 and .1, D0.0 to .4,
	of \ref{sec-dprojter} and \ref{sec-eprojter} will be assumed, but not
	all these constructions are necessarily required.

	\sssec{Lemma. [For the linearity property.]}
	H1.0.\hti{3}$X_A = (A,0),$ $Y_C := (0,C)$, $(B),$ (See Fig. 21a)\\
	D1.0.\hti{3}$j_b := Q_0 \times B,$ $j'_b := Y_C \times B,$
		$j_1 := Q_0 \times V,$ $j'_1 := Y_C \times V,$\\
	D1.1.\hti{3}$x := X_A \times Q_2, K := x \times j_b,$
		$k_0 := K \times Q_1,$ $L := k_0 \times j_1,$\\
	D1.2.\hti{3}$K' := x \times j'_b,$ $c := L\times Q_2,$
		$L' := c \times j'_1,$ $k'_0 := L' \times K',$\\
	C1.0.\hti{3}$Q_1\incid k'_0.$

	{\em Moreover,}\\
	$K = (A,A\cdot B),$ $L = (A\cdot B,A\cdot B),$ $K' = (A,A*B*C),$
	$L' = (A\cdot B,A\cdot B+C),$\\
	C1.0 $\implies A*B*C = A\cdot B + C.$

	Proof:\\
	Elated-Desargues($Q_2,\{Q_0,K,L\},\{Y_C,K',L'\};\langle Q_1,V,B\rangle,
	q_2)\implies Q_1\incid k'_0.$

	\sssec{Lemma. [For the additive associativity law]}
	H1.0.\hti{3}$X_A,$ $X_B,$ $Y_C,$ (See Fig. 21b)\\
	D1.0.\hti{3}$a := X_A \times Q_2,$\\
	D1.1.\hti{3}$b := X_B \times Q_2,$ $B_1 := b \times i,$
		 $x_1 := B_1 \times Q_1,$ $Y_1 := x_1 \times q_1,$\\
	D1.2.\hti{3}$i_2 := Y_1 \times V,$ $A_1 := i_2 \times a,$
		 $x_3 := A_1 \times Q_1,$ $D_1 := x_3 \times i,$\\
	D1.3.\hti{3}$d := D_1 \times Q_2,$ $i_1 := Y_C \times V,$
		$D_2 := d \times i_1,$\\
	D1.4.\hti{3}$B_2 := i_1 \times b,$ $x_2 := B_2 \times Q_1,$
		$Y_2 := x_2 \times q_1,$\\
	D1.5.\hti{3}$i_3 := Y_2 \times V,$ $A_2 := i_3 \times a,$
		$x_4 := A_2 \times D_2,$\\
	D1.6.\hti{3}$e_1 := A_1 \times B_1,$ $e_2 := A_2 \times B_2,$
		$E := e_1\times e_2,$\\
	C1.0.\hti{3}$E\incid q_2,$\\
	C1.1.\hti{3}$Q_1\incid x_4,$

	{\em Moreover,}\\
	$B_1 = (B,B),$ $Y_B = (0,B),$ $A_1 = (A,A+B),$ $D_1 = (A+B,A+B),$\\
	$B_2 = (B,B+C),$ $Y_1 = (0,B+C),$ $A_2 = (A,A+(B+C)),$\\
		$D_2 = (A+B,((A+B)+C)$,\\
	C1.1. $\implies A+(B+C) = (A+B)+C.$

	Proof:\footnote{variant due to Michael Sullivan, October 24, 1989.}\\
	Elated-Desargues($Q_2,\{A_1,B_1,Y_1\};\{A_2,B_2,Y_2\};
	\langle Q_1,V,E\rangle,q_2),\\
	\implies$ Elated-Desargues$^{-1}(q_2,\{A_1,B_1,D_1\};\{A_2,B_2,D_2\};
		\langle V,Q_1,E\rangle,Q_2)\\
	\implies Q_1\incid x_4.$

	\sssec{Corollary.}
	{\em If 2 $m$-parallelograms $\{A_j\}$ and $\{B_j\},$ $j = 0,1,2,3,$ are
	such that\\
	$A_k\times B_k \paras{m} A_0\times B_0$, $k = 1,2,$ the same
	is true for $k = 3.$} (See Fig. 21e)\marginpar{21e?}

	The parallelograms for which the proof is given in the Lemma are\\
	$\{A_1,Y_B,Y_1,A_2\}$ and $\{D_1,B_1,B_2,D_2\}$.

	\sssec{Lemma. [For the right distributive law]}
	H1.0.\hti{3}$X_A = (A,0),$ $Y_1 = (0,B),$ $(C),$ (See Fig. 21c)\\
	D1.0.\hti{3}$x := X_A\times Q_2,$\\
	D1.1.\hti{3}$x_1 := Q_1\times Y_1,$ $B_1 := x_1\times i,$
		$i_1 := Y_1\times V,$ $A_1 := i_1\times x,$\\
	D1.2.\hti{3}$x_3 := A_1\times Q_1,$ $F_1 := x_3\times i,$
		$f := F_1\times Q_2,$ $c_1 := Q_0\times C,$\\
	D1.3.\hti{3}$b := B_1\times Q_2,$ $B_2 := b\times c_1,$
		$F_2 := f\times c_1,$\\
	D1.4.\hti{3}$x_2 := B_2\times Q_1,$ $Y_2 := x_2\times q_1,$
		$e_1 := Y_2\times A_1,$ $e_2 := B_2\times F_1,$\\
	D1.5.\hti{3}$E := e_1\times e_2,$\\
	D1.6.\hti{3}$c_2 := Y_2\times C,$ $A_2 := c_2\times x,$
		$x_4 := A_2\times F_2,$\\
	C1.0.\hti{3}$E\incid q_2$.\\
	C1.1.\hti{3}$Q_1\incid x_4$.

	{\em Moreover,}\\
	$B_1 = (B,B),$ $A_1 = (A,A+B),$ $F_1 = (A+B,A+B),$
		$B_2 = (B,B\cdot C),$\\
	$F_2 = (A+B,(A+B)\cdot C),$ $Y_2 = (0,B\cdot C),$
	$A_2 = (A,A*C*(B\cdot C)),$ and\\
	C1.1 $\implies (A+B)\cdot C = A\cdot C + B\cdot C.$

	Proof:\\
	Elated-Desargues($Q_1,\{Y_1,A_1,Y_2\},\{B_1,F_1,B_2\};
	\langle E,Q_2,V\rangle,q_2)\implies E\incid q_2.$\\
	Elated-Desargues$^{-1}(q_2,\{A_2,A_1,Y_2\},\{F_2,F_1,B_2\};
	\langle E,C,Q_2\rangle,Q_1)\implies Q_1\incid x_4.$\\
	Finally, from C1.0 follows $(A+B)\cdot C = A*C*(B\cdot C)$, but by
	linearity, the second member equals $A\cdot C + B\cdot C.$

	\sssec{Exercise.}
	Determine the identity corresponding to C1.0 or to the $m$-parallelism
	of $Y_2\times A_1$ and $B_2\times F_1$.

	\sssec{Lemma. [For the commutativity law]}
	H1.0.\hti{3}$A_1,$ $B_1,$ (See Fig. 21d)\\
	D1.0.\hti{3}$a_0 := A_1 \times Q_0,$ $A := a_0 \times q_2,$
		$a_2 := A \times B_1,$\\
	D1.1.\hti{3}$b_0 := B_1 \times Q_0,$ $B := b_0 \times q_2,$
		$b_2 := A_1 \times B,$ $D := b_2 \times a_2,$\\
	D1.2.\hti{3}$x_1 := A_1 \times Q_1,$ $Y_A := x_1 \times q_1,$
		$b_1 := Y_A \times B,$\\
	D1.3.\hti{3}$y_2 := B_1 \times Q_2,$ $B_2 := y_2 \times b_1,$
		$x_2 := B_2 \times D,$ $y_1 := A_1 \times Q_2,$\\
	D1.4.\hti{3}$x_3 := B_1 \times Q_1,$ $Y_B := x_3 \times q_1,$
		$a_1 := Y_B \times A,$ $A_2 := a_1 \times y_1,$\\
	C1.0.\hti{3}$Q_1\incid x_2,$\\
	C1.1.\hti{3}$A_2\incid x_2,$

	{\em Moreover,\\
	if $A_1 = (X_A,Y_A),$ and $B_1 = (X_B,Y_B),$ then}
	$A_2 = (X_A,Y_A+Y_B),$ $B_2 =(X_B,Y_B+Y_A),$\\
	C1.0 and .1 $\implies Y_A+Y_B = Y_B+Y_A$.

	Proof:\\
	Elated-Desargues($B,\{Q_0,Y_A,A_1\},\{B_1,B_2,D\};
	\langle Q_1,A,Q_2\rangle,q_2)\implies Q_1\incid x_2.$\\
	Elated-Desargues($A,\{Q_0,Y_B,B_1\},\{A_1,A_2,D\};
	\langle Q_1,B,Q_2\rangle,q_2)\implies A_2\incid x_2.$\\
	therefore $A_2$ and $B_2$ have the same second coordinate $Y$.\\
	Because $A_1\incid a_0\incid (A)$, $Y_A = X_A\cdot A,$ by construction
	and because of linearity,
	$A_2 = (X_A,X_A*A*Y_B) = (X_A,X_A\cdot A + Y_B)$
	similarly $Y_B = X_B\cdot B,$ and $B_2 = (X_B,X_B\cdot B + Y_A).$

	\sssec{Corollary.}
	If we make the same constructions as in the lemma with $A = B = J,$ then
\\
	$Q_1\incid(A_2\times B_2).$

	\sssec{Lemma. [Addition an Negation in Veblen-Wedderburn planes.]}
	H0.0.	$Y_A,$ $Y_B,(Fig. 21e)$\\
	D1.0.	$i_1 := Y_A\times V,$ $i_2 := Y_B\times V,$\\
	D1.1.	$x_1 := Y_A\times Q_1,$ $A_1 := x_1\times i,$
		$a := A_1\times Q_2,$ $A_2 := a\times i_2,$\\
	D1.2.	$x_3 := Y_B\times Q_1,$ $B_1 := x_3\times i,$
		$b := B_1\times Q_2,$ $B_2 := b\times i_1,$\\
	D1.3.	$x_2 := A_2\times B_2,$\\
	C1.0.	$Q_1\incid x_2,$\\
	D2.0.	$U_1 := x_1\times v,$ $c := U_1\times Q_0,$
		$A := c\times q_2,$\\
	D2.1.	$c- := Y_A\times I,$ $A- := c-\times q_2,$\\
	{\em Moreover,\\
	If $Y_A = (0,A)$ and $Y_B = (0,B),$ then} $A_1 = (A,A),$ $B_1 = (B,B),$
	$A_2 = (A,A+B),$ $B_2 = (B,B+A),$\\
	$U_1 = (1,A),$ $A = (A),$ $A- = (-A).$

	\sssec{Theorem.}\label{sec-tveblenw}
	{\em In a Veblen-Wedderburn plane, the ternary ring ($\Sigma$,*) is a
	quasifield in the terminology of Dembowski} (p. 129):
	\enumb
	\item$(\Sigma,*)$ is linear, $a * b * c = a \cdot b \:+\:c$,
	\item($\Sigma$,+) {\em is an abelian group,}
	\item($\Sigma-\{0\},\cdot)$ {\em is a loop,}
	\item ($\Sigma,*) = (\Sigma,+,\cdot)$ {\em is right distributive,}
	$(a+b)\cdot c = a \cdot c + b\cdot c$.
	\item$a \neq b \implies x \cdot a = x \cdot b + c$
	{\em has a unique solution.}
	\enume

	\sssec{Theorem.}
	{\em In a Veblen-Wedderburn plane with ideal line $m$,
	${\cal T}^{m,B}_{AA'}(C)$, $C\incid A\times A'$, is independent of $B$.
	We can therefore use ${\cal T}^m_{AA'}$ as notation for a translation.}
	
	\sssec{Definition.}
	$m$-equality is defined by\\
	\hth$[A,A'] =_m [P,P']$ iff $P' = {\cal T}^m_{AA'}(P).$

	\sssec{Theorem.}\label{sec-t3coord}
	{\em In a Veblen-Wedderburn plane we can use systematically 3
	coordinates as follows\\
	$(Q_2)$ is equivalent to $(0,1,0)$,\\
	$(P_0)$ is equivalent to $(1,P_0,0),$\\
	$(P_0,P_1)$ is equivalent to $(P_0,P_1,1),$\\
	$[q_2]$ is equivalent to $[0,0,1],$\\
	$[M_0]$ is equivalent to $[1,0,-M_0]$,\\
	$[M_0,M_1]$ is equivalent to $[M_0,-1,M_1].$\\
	A point $(P_0,P_1,P_2)$ is incident to the line $[l_0,l_1,l_2]$ iff\\
	\hth$P_0l_0+P_1l_1+P_2l_2 = 0.$}

	Proof:
	In the general case, because of linearity, a point $(P_0,P_1)$ is
	incident to the line $[M_0,M_1]$ if $P_1 = P_0\cdot M_0 + M_1,$ which we
	can rewrite\\
	\hth$P_0\cdot M_0 + P_1\cdot (-1) + 1 \cdot M_1.$\\
	The other correspondances can be verified using \ref{sec-tincidpersp}.

	\sssec{Theorem.}
	{\em In a Veblen-Wedderburn plane with ideal line $m$, $m$-equality is
	an equivalence relation.}

	\ssec{Moufang Planes.}

	\sssec{Definition.}
	A {\em Moufang plane} is a Veblen-Wedderburn plane in which the elated
	Desargues axiom is satisfied for every line in the plane. (See Fig. 3f).

	\sssec{Theorem.}
	{\em Duality is satisfied in a Moufang plane.}

	\sssec{Definition.}
	The {\em C-Desargues Configuration} is a Desargues Configuration, for
	which 2 corresponding sides intersect on the line joining the other
	vertices.  The point of intersection will be underlined.

	\sssec{Lemma.}
	{\em The Elated-Desargues Configuration for all lines in the planes
	implies the C-Desargues Configuration.}

	Proof: (See Fig 3f.)
	To prove C-Desargues$(C,\{A_i\},\{\un{B_0},B_1,B_2\};\langle C_i\rangle,
	c),$ we apply Elated-Desargues$^{-1}(c_0,\{A_1,B_1,C_2\},
	\{A_2,B_2,C_1\};\langle B_0,A_0,C\rangle,C_0)$.

	\sssec{Definition.}
	The {\em 1-Desargues Configuration} is a Desargues Configuration, for
	which the vertex of 1 triangle is on the side of the other, this
	vertex will be underlined.

	\sssec{Lemma.}
	{\em The Elated-Desargues Configuration for all lines in the planes
	implies the 1-Desargues Configuration.}

	Proof: (See Fig 3b.)
	To prove 1-Desargues$(C,\{A_i\},\{B_i\};\langle \un{C_0},C_1,C_2\rangle,
	c),$ we apply Elated-Desargues$(B_0,\{A_0,C_1,C_2\},\{C,B_2,B_1\};
	\langle C_0,A_1,A_2\rangle,a_0)$.

	\sssec{Theorem.}\label{sec-t1des}
	{\em In a Moufang plane}
	\enumb
	\item {\em the C-Desargues Theorem is true}.
	\item {\em the 1-Desargues Theorem is true}.
	\enume

	\sssec{Lemma. [For the left distributive law]}
	H1.0.\hti{3}$X_A = (A,0),$ $(B),$ $(C),$ (See Fig. 22a)\\
	D1.0.\hti{3}$a := X_A\times Q_2,$ $c_1 := Q_0\times C,$
		$U_1 := c_1\times u,$ $x_1 := U_1\times Q_1,$\\
	D1.1.\hti{3}$Y_1 := x_1\times q_1,$ $b_1 := Y_1\times B,$
		$U_2 := b_1\times u,$\\
	D1.2.\hti{3}$A_1 := a\times c_1,$ $x_2 := A_1\times Q_1,$
		$Y_2 := x_2\times q_1,$\\
	D1.3.\hti{3}$b_2 := Y_2\times B,$ $A_2 := b_2\times a,$
		$d := U_2\times Q_0,$\\
	C1.0.\hti{3}$A_2 \incid d.$

	{\em Moreover,}\\
	$U_1 = (1,C),$ $Y_1 = (0,C),$ $U_2 = (1,1*B*C),$ $A_1 = (A,A\cdot C),$
	$Y_2 = (0,A\cdot C),$ $A_2 = (A,A*B*(A\cdot C),$\\
	C1.0 $\implies A\cdot B)+(A\cdot C) = A*B*(A\cdot C) = A\cdot(1*B*C) =
		A\cdot((1\cdot B)+C) = A\cdot(B+C).$

	Proof:\\
	C-Desargues$(Q_0,\{Y_1,U_1,U_2\},\{Y_2,A_1,A_2\};
		\langle \un{Q_2},B,Q_1\rangle,q_2$)\\
	$\implies((UA_1\times A_2)\times (U_1\times U_2))\incid
	(Y_1\times Y_2).$

	\sssec{Lemma. [For the inverse property]}
	H1.0.\hti{3}$A,$ (See Fig. 22b)\\
	D1.0.\hti{3}$a := X_A \times Q_2,$ $ A_0 := a \times j,$
		$A_1 := a \times i,$\\
	D1.1.\hti{3}$a_1 := A_1 \times Q_1,$ $ A_2 := a_1 \times u,$
		$a_0 := A_0 \times Q_0,$ $ A_3 := a_0 \times u,$\\
	D1.2.\hti{3}$a_2 := A_2 \times Q_0,$ $ A_4 := a_2 \times j,$\\
	D1.4.\hti{3}$a_3 := A_3 \times Q_1,$ $ a_4 := A_4 \times Q_2,$
		$A_5 := a_3 \times a_4,$\\
	D1.5.\hti{3}$d_1 := A_0 \times A_2,$ $d_2 := A_3 \times A_4,$
		$E := d_1 \times d_2,$\\
	C1.0.\hti{3}$E\incid q_2,$\\
	C1.1.\hti{3}$A_5\incid i.$

	{\em Moreover,}\\
	$A_0 = (A,1),$ $A_1 = (A,A),$ $A_2 = (1,A),$ $A_4 = (A^L,1),$
	$A_3 = (1,A^R),$\\
	$A_5 = (A^L,A^R),$\\
	C1.1 $\implies A^L = A^R.$

	Proof:\\
	1-Desargues$(Q_0,\{\un{A_0},A_2,A_1\},\{A_3,A_4,U\};
	\langle Q_1,Q_2,E\rangle,q_2)\implies E\incid q_2.$\\
	1-Desargues$^{-1}(Q_0,\{\un{A_4},A_3,A_5\},\{A_2,A_0,U\};
	\langle Q_1,Q_2,E\rangle,q_2)\implies A_5\incid i.$

	\sssec{Notation.}
	If $B \neq 0,$ {\em we write $B^{-1} = B^R$.}

	\sssec{Lemma. [For the right inverse property]}
	H1.0.\hti{3}$ X_A,$ $ B,$ (See Fig. 22c)\\
	D1.0.\hti{3}$a := X_A \times Q_2,$ $ A_1 := a \times i,$
		$b := Q_0 \times B,$ $ C_2 := j \times b,$\\
	D1.1.\hti{3}$c := C_2 \times Q_2,$ $ C_1 := c \times i,$
		$x1 := C_1 \times Q_1,$ $ U_1 := x_1 \times u,$\\
	D1.2.\hti{3}$A_2 := a \times b,$ $x_3 := A_2 \times Q_1,$
		$AB_2 := x_3 \times i,$ $ab := AB_2 \times Q_2,$\\
	D1.3.\hti{3}$b' := U_1 \times Q_0,$ $AB_1 := ab \times b',$
		$x_2 := A_1 \times AB_1,$\\
	D1.4.\hti{3}$d := U_1 \times A_1,$ $ e := U \times A_2,$
		$S := d \times e,$\\
	C1.0.\hti{3}$S \incid q_0,$
	C1.1.\hti{3}$Q_1\incid x_2.$

	{\em Moreover,}\\
	$A_1 = (A,A),$ $C_2 = (B^{-1},1),$ $C_1 = B^{-1},B^{-1}),$
	$U_1 = (1,B^{-1}),$
	$A_2 = (A,A\cdot B),$ $AB_2 = (A\cdot B,A\cdot B),$
	$AB_1 = (A\cdot B,(A\cdot B)\cdot B^{-1}),$\\
	C1.1 $\implies (A\cdot B)\cdot B^{-1} = A.$

	Proof:\\
	1-Desargues($Q_2,\{U_1,A_1,C_1\},\{\un{U},A_2,C_2\};
	\langle Q_0,Q_1,S\rangle,q_0)\implies S\incid q_0.$\\
	1-Desargues$^{-1}(Q_2,\{U_1,\un{A_1},AB_1\},\{U,A_2,AB_2\};
	\langle Q_0,Q_1,S\rangle,q_0)\implies Q_1\incid x_2.$

	\sssec{Lemma. [For the left inverse property]}
	H1.0.\hti{3}$X_A,$ $B,$ (See Fig. 22d)\\
	D1.2.\hti{3}$b := Q_0 \times B,$ $U_3 := b \times u,$
		$x_3 := U_3 \times Q_1,$\\
	D1.3.\hti{3}$a := X_A \times Q_2,$ $A_1 := a \times j,$
		$b' := Q_0 \times A_1,$ $U_1 := b' \times u,$\\
	D1.4.\hti{3}$x_1 := U_1 \times Q_1,$ $C_1 := x_1 \times i,$
		$c := C_1 \times Q_2,$ $C_2 := c  \times x_3,$\\
	D1.5.\hti{3}$A_2 := a \times b,$ $x_2 := A_2 \times Q_1,$
		$U_2 := x_2 \times u,$\\
	D1.6.\hti{3}$ab := U_2 \times Q_0,$\\
	D1.7.\hti{3}$r_1 := U \times A_2,$ $r_2 := U_3 \times C_1,$
		$R := r_1 \times r_2,$\\
	C1.0.\hti{3}$R \incid q_2,$\\
	C1.1.\hti{3}$C_2 \incid ab,$

	{\em Moreover,}\\
	$A_1 = (A,1),$ $A_2 = (A,A\cdot B),$ $U_1 = (1,A^{-1}),$ $U_3 = (1,B),$
	$U_2 = (1,A\cdot B),$
	$C_1 = (A^{-1},A^{-1}),$ $C_2 = (A^{-1},A^{-1}\cdot (A\cdot B)).$\\
	C1.1 $\implies A^{-1}\cdot(A\cdot B) = B.$

	Proof:\\
	1-Desargues($Q_0,\{A_1,\un{U},A_2\},\{U_1,C_2,U_3\};
	\langle R,Q_2,Q_1\rangle,q_2)\\
	\implies$ 1-Desargues$^{-1}(q_2,\{U,U_2,A_2\},\{C_1,C_2,\un{U_3}\};
	\langle Q_1,R,Q_2\rangle,Q_0)\implies C_2\incid ab.$

	\sssec{Theorem.}
	{\em With the coordinatization of the plane as given in}
	\ref{sec-dprojter},
	\enumb
	\item {\em the ternary ring ($\Sigma,+,\cdot)$ is left distributive, or
}\\
	\hth$A \cdot (B+C) = A \cdot B + A \cdot C$.
	\item $B \neq 0 \implies  B^R = B^L = B^{-1},$
	\item $(A\cdot B)\cdot B^{-1} = B^{-1}\cdot(B\cdot A) = A$
	{\em for all $A$}.
	\enume
	{\em In other words, ($\Sigma,+,\cdot)$ is an alternative division
	ring.}

	\ssec{Desarguesian Planes.}

	\sssec{Definition.}
	A {\em Desarguesian plane} is a plane in which the Desargues Axiom
	is always satisfied.

	\sssec{Theorem.}
	{\em Duality is satisfied in a Desarguesian plane.}

	\sssec{Comment.}
	Instead of the Axiom of Desargues one can use the equivalent axiom
	of Reidemeister (See Theorem II.\ref{sec-treid} and Klingenberg, 1955).

	\sssec{Lemma. [For Associativity]}
	H1.0.\hti{3}$X_A,$ $B,$ $C,$ (See Fig. 23.)\\
	D1.0.\hti{3}$b := Q_0 \times B,$ $c := Q_0 \times C,$
		$U_1 := b \times u,$ $x_1 := U_1 \times Q_1,$\\
	D1.1.\hti{3}$D_1 := x_1 \times i,$ $d := D_1 \times Q_2,$
		$D_2 := d \times c,$ $x_2 := D_2 \times Q_1,$\\
	D1.2.\hti{3}$U_2 := x_2 \times u,$ $bc := U_2 \times Q_0,$\\
	D1.3.\hti{3}$a := X_A \times Q_2,$ $A_1 := a \times b,$
		$x_3 := A_1 \times Q_1,$ $AB_1 := x_3 \times i,$\\
	D1.4.\hti{3}$ab := AB_1 \times Q_2,$ $AB2 := ab \times c,$
		$x_4 := AB_2 \times Q_1,$ $A_2 := x_4 \times a,$\\
	D1.5.\hti{3}$r_1 := U_1 \times D_2,$ $r_2 := A_1 \times AB_2,$
		$R := r_1\times r_2,$\\
	C1.0.\hti{3}$A_2 \incid bc,$

	{\em Moreover,}\\
	$A_1 = (A,A\cdot B),$ $AB_1 = (A\cdot B,A\cdot B),$
	$AB_2 = (A\cdot B,(A\cdot B)\cdot C),$\\
	$A_2 = (A,A\cdot(B\cdot C)),$
	$U_1 = (1,B),$ $D_1 = (B,B),$ $D_2 = (B,B\cdot C),$
	$U_2 = (1,B\cdot C).$\\
	C1.0 $\implies A\cdot(B\cdot C) = (A\cdot B)\cdot C.$

	Proof:\\
	Desargues($Q_0,\{D_1,D_2,U_1\},\{AB_1,AB_2,A_1\};
	\langle R,Q_1,Q_2\rangle,q_2)\\
	\implies$ Desargues$^{-1}(q_2,\{U_2,U_1,D_2\},\{A_2,A_1,AB_2\};
	\langle R,Q_1,Q_2\rangle,Q_0)\implies A_2\incid bc.$

	\sssec{Theorem.}
	{\em With the coordinatization of the plane as given in
	\ref{sec-dprojter}},
	\enumb
	\item {\em ($\Sigma,\cdot)$ is associative,}\\
	\hth$A \cdot (B \cdot C) = (A \cdot B) \cdot C.$
	\enume
	In other words, ($\Sigma,+,\cdot)$ is a skew field.

	\sssec{Theorem.}
	{\em If a Desarguesian plane we use the coordinates  of
	\ref{sec-t3coord}, we can make them homogeneous by multiplying the
	coordinates of points to the left by the same element in the set
	$\Sigma,$ and those of lines to the right by the same element in the
	set $\Sigma.$}\\

	Associativity of multiplication is essential to allow for the left
	equivalence of points and the right equivalence of lines.

	\ssec{Pappian planes.}

	\sssec{Axiom. [Of Pappus]}
	In a perspective plane:
	If $A_i$ are 3 distinct points on a line $a$ and $B_i$ are 3 distinct
	points on a line $b$ and $C_i := (A_{i+1}\times B_{i-1}) \times
	(A_{i-1}\times B_{i+1}$ then incidence$(C_i)$.

	I write Pappus$(\{A_i\},\{B_i\};\{C_i\}).$

	\sssec{Definition.}
	A {\em Pappian plane} is a plane in which the Pappus Axiom
	is always satisfied.

	\sssec{Comment.}
	There are other axioms which are equivalent to that of Pappus.
	The Fundamental axiom and Axiom A (See Seidenberg, p. 25 and Chapter
	IV). The Fundamental axiom states that there is at most one
	projectivity which associates 3 given distinct collinear points into 
	3 given distinct collinear.  Axiom A states that if a projectivity which
	associates a line $l$ into a distinct line $l'$ leaves $l\times l'$
	invariant then it is a perspectivity.

	\sssec{Theorem.}
	{\em Duality is satisfied in a Pappian plane.}

	\sssec{Theorem.}
	{\em A Pappian plane is a Desarguesian plane.}

	\sssec{Lemma. [For Commutativity]}
	H1.0.\hti{3}$(A),$ $(B),$ (See Fig. 24)\\
	D1.0.\hti{3}$a := Q_0 \times A,$ $b := Q_0 \times B,$
		$U_1 := a \times u,$\\
	D1.1.\hti{3}$x_1 := U_1 \times Q_1,$ $C_1 := x_1 \times i,$
		$c := C_1 \times Q_2,$ $C_2 := c \times b,$\\
	D1.2.\hti{3}$U_2 := b \times u,$\\
	D1.3.\hti{3}$x_2 := U_2 \times Q_1,$ $D_1 := x_2 \times i,$
		$d := D_1 \times Q_2,$ $D_2 := d \times a,$\\
	D1.4.\hti{3}$x_3 := C_2 \times D_2,$\\
	C1.0.\hti{3}$D_2 \incid x_3,$

	{\em Moreover,}\\
	$U_1 = (1,A),$ $U_2 = (1,B),$ $C_1 = (A,A),$ $D_1 = (B,B),$
	$C_2 = (A,AB),$ $D_2 = (B,BA),$
	C1.0 $\implies A\cdot B = B\cdot A.$

	Proof:\\
	Pappus($\langle D_1,C_1,Q_0\rangle,\langle U_1,U_2,Q_2\rangle;
	\langle C_2,D_2,Q_2\rangle)\implies D_2\incid x_3.$

	\sssec{Theorem.}
	{\em With the coordinatization of the plane as given in
	\ref{sec-dprojter}},
	\enumb
	\item {\em $(\Sigma,\cdot)$ is commutative,}\\
	\hth$a \cdot b = b \cdot a.$
	\enume
	{\em In other words, $(\Sigma,+,\cdot)$ is a field.}

	\sssec{Theorem.}
	{\em The field of a Pappus-Fano plane has characteristic 2.
	Vice-versa if a field has characteristic 2, the corresponding Pappian
	plane satisfies the axiom N-Fano.}

	Proof:
	We have seen than in a Fano plane $A + A = 0$, for all $A \in \Sigma$,
	therefore the characteristic of the field is 2.
	To prove the converse,we choose as coordinates of the vertices of the
	quadrangle $A_0 = (1,0,0),$ $A_1 = (0,1,0),$ $A_2 = (0,0,1)$ and
	$M = (1,1,1),$ the diagonal elements are $M_0 = (0,1,1),$
	$M_1 = (1,0,1),$ $M_2 = (1,1,0),$ which are collinear iff $1 + 1 = 0.$

	\ssec{Separable Pappian Planes.}

	\sssec{Axiom. [Of separation]}\label{sec-asepar}
	In a perspective plane, if $A_i,$ $i = 0,1,2,3,4$ are distinct points on
	the same line:
	\enumb
	\item There are at least 4 points on a line.
	\item $\sigma(A_0,A_1\vert A_2,A_3) \implies
	\sigma(A_0,A_1\vert A_3,A_2)$ and
		$\sigma(A_3,A_2\vert A_0,A_1)$
	\item only one of the relations $\sigma(A_0,A_1\vert A_2,A_3)$,
	$\sigma(A_0,A_2\vert A_1,A_3),$ $\sigma(A_0,A_3\vert A_1,A_2)$ holds.
	\item $\sigma(A_0,A_1\vert A_2,A_3)$ and $\sigma(A_1,A_2\vert A_3,A_4)
	\implies \sigma(A_0,A_4\vert A_2,A_3).$
	\item $\Pi(P,A_j,A'_j),$ $j = 0,1,2,3,$ and $\sigma(A_0,A_1\vert
	A_2,A_3)\implies \sigma(A'_0,A'_1\vert A'_2,A'_3).$
	\enume

	\sssec{Definition.}
	A {\em separable Pappian plane} is a Pappian plane in which the
	separation axioms are satisfied.

	\sssec{Theorem.}\label{sec-tsepar}
	\enumb
	\item $\sigma(A_0,A_1\vert A_2,A_3) \implies
	\sigma(A_1,A_0\vert A_2,A_3)$, $\sigma(A_0,A_1\vert A_3,A_2)$,
	$\sigma(A_1,A_0\vert A_3,A_2)$,\\
	\hti{4}$\sigma(A_2,A_3\vert A_0,A_1)$, $\sigma(A_2,A_3\vert A_1,A_0)$
	$\sigma(A_3,A_2\vert A_0,A_1)$, $\sigma(A_3,A_2\vert A_1,A_0).$
	\item $\sigma(A_0,A_1\vert A_2,A_3) {\mathit and}
	\sigma(A_1,A_2\vert A_3,A_4)\implies \sigma(A_0,A_4\vert A_1,A_2)$,
	\enume

	\sssec{Notation.}
	When we use \ref{sec-asepar}.3 or \ref{sec-tsepar}.1, I will underline
	the element in each quadruple of point which is distinct, to ease
	the application of the axiom and write, for instance\\
	\hth$\sigma(\un{A_0},A_1\vert A_2,A_3)$ and $\sigma(A_1,A_2\vert
	A_3,\un{A_4}) \implies \sigma(A_0,A_4\vert A_2,A_3),$ or\\
	\hth$\sigma(A_3,A_2\vert A_1,\un{A_0})$ and $\sigma(A_2,A_1\vert
	A_3,\un{A_4}) \implies \sigma(A_0,A_4\vert A_2,A_1).$

	\sssec{Theorem.}
	{\em In a Pappus-Fano plane, given a harmonic quadrangle $A,B,C,D$,
	(See Fig. 2a''), $P,R\vert U,V$, where $P$, $R$ are diagonal points
	and $U$, $V$ are the intersection with $P\times R$ of the sides
	of the quadrangle which are not incident to $P$ or $R$.}

	Proof:\\
	$\Pi(C,\{P,U,R,V\},\{D,Q,B,V\}),$ $\Pi(A,\{D,Q,B,V\},\{R,U,P,V\}),$
	therefore\\
	$P,R\vert U,V \implies R,P\vert U,V$, while
	$P,U\vert R,V \implies R,U\vert P,V$,
	$P,V\vert U,R \implies R,V\vert U,P$,
	the last 2 conclusions are contradicted by \ref{sec-asepar}.2.

	\sssec{Corollary.}\label{sec-csepar}
	$(O,\infty\vert A,-A).$

	\sssec{Definition.}
	Given $A_i,$ $i = 0,1,2$ on a line $a$, a {\em segment} 
	$seg(A_0,A_1\:\setminus\:A_2)$ is the set of points $A \incid a$ such
	that $\sigma(A_0,A_1\vert A,A_2)$.

	\sssec{Lemma.}
	{\em If $A_i\in\Sigma$ and $\sigma(A_0,A_1\vert A_2,A_3)$,}
	\enumb
	\item $\sigma(P+A_0,P+A_1\vert P+A_2,P+A_3)$,
	\item $P \neq 0 \implies \sigma(P\cdot A_0,P\cdot A_1
		\vert P\cdot A_2,P\cdot A_3).$
	\item {\em More generally, if $\Pi$ is a projectivity which associates
		to\\
		$X$, $(A\cdot X + B)\cdot (C\cdot X + D),$ $A\cdot C\cdot
	D\neq 0$, $A\cdot D\neq B\cdot C$,\\
		then} $\sigma(\Pi(A_0),\Pi(A_1),\Pi(A_2),\Pi(A_3)).$
	\enume
	{\em The same properties hold if one of the $A_i$ is replaced by
	$\infty$ and we use $\infty + A = \infty$ and with $A\neq 0$,
	$\infty \cdot A = \infty.$}

	Proof:\\
	$\Pi(V,q_1,p) \circ \Pi(Q_1,p,q_1)$ transforms $(0,A_i)$ into
	$(P,P+A_i)$ into $(0,P+A_i).$\\
	$\Pi(Q_1,q_1,i) \circ \Pi(Q_2,i,Q_0\times(P))
	\circ \Pi(Q_2,i,Q_0\times(P))$ transforms $(0,A_i)$ into
	$(A_i,A_i)$ into $(A_i,P\cdot A_i)$ into $(0,P\cdot A_i).$
	The rest of the proof is left as an exercise.

	\sssec{Lemma.}
	{\em In a separable Pappian plane, the characteristic is not 2.}

	Proof:  If the characteristic was 2 and $A$ is different from 0, 1 and
	$\infty$,\\
	either $\sigma(0,1\vert A,\infty)$ or $\sigma(0,A\vert \infty,1)$ or
	$\sigma(0,\infty\vert 1,A)$.\\
	In the first case, adding 1 or $A$ gives $\sigma(\un{1},0\vert A+1,
	\infty)$ or $\sigma(\un{A},A+1\vert 0,\infty)$, combining gives
	$\sigma(1,A\vert 0,\infty)$ which contradicts $\sigma(0,1\vert A,
	\infty)$.  In the second case we add 1 or $A$ and in the third case
	we add 1 or $A+1$ and proceed similarly to show contradiction.

	\sssec{Definition.}
	$P$ is {\em positive}, or $P > 0$, iff $\sigma(0,\infty\vert -1,P)$.\\
	$P$ is {\em negative}, or $P < 0$, iff $-P > 0$ or iff $\sigma(0,
	\infty\vert -1,-P)$ or iff $\sigma(0,\infty\vert 1,P)$.

	\sssec{Theorem.}
	\enumb
	\item $1 > 0$.
	\item $A,B \in \Sigma$, $A >0$ {\em and} $B > 0 \implies A+B >0.$ 
	\item {\em $A \in \Sigma$, either $A = 0$ or $A >0$ or $-A >0$.}
	\item $A,B \in \Sigma$, $A >0$ {\em and} $B > 0 \implies A\cdot B >0.$ 
	\enume

	Proof:\\
	For 0, we use Corollary \ref{sec-csepar}.\\
	For 1, $A > 0 \implies \sigma(0,\infty\vert -1,A)$ by the projectivity
	which associates to $X$, $A-X-1$,\\
	\hti{3}4.~~$\sigma(A-1,\infty\vert A,\un{-1}),$\\
	$B > 0 \implies \sigma(0,\infty\vert -1,B) \implies$ (adding $A$)
	$\sigma(A,\infty\vert A-1,\un{A+B})\implies$ (combining with 4.)
	$\sigma(-1,\un{A+B}\vert A,\infty)$, with $\sigma(\un{0},\infty\vert -1,
	A) \implies \sigma(0,\un{A+B}\vert A,\infty),$ with $\sigma(0,\infty
	\vert\un{-1},A) \implies \sigma(-1,A+B\vert 0,\infty) \implies
	A+B > 0.$\\
	For 2, by the definition of $A > 0$ or $-A > 0$, it follows that
	$A$ is not 0.  $A > 0$ and $-A >0$ are also mutually
	exclusive, otherwize $A+(-A) = 0$ would be positive. If $A = -1,$ then
	$-A = 1 > 0$. It remains to examine for a given $A$ distinct from 0 and 
	-1, the 3 possibilities,\\
	$\sigma(0,-1\vert \un{A},\infty)$ and $\sigma(0,\infty\vert -1,\un{1})
	\implies \sigma(1,A\vert 0,\infty)\implies A < 0.$\\
	$\sigma(0,\un{A}\vert \infty,-1)$ and $\sigma(0,\infty\vert -1,\un{1})
	\implies \sigma(1,A\vert 0,\infty)\implies A < 0.$\\
	$\sigma(0,\infty\vert -1,A)\implies A > 0.$\\
	For 3, $\sigma(0,\infty\vert -1,B) \implies \sigma(0,\infty\vert -A,
	A\cdot B)$,\\
	$A > 0 \implies$ not $\sigma(0,\infty\vert -1,-A),$ therefore either
	$\sigma(0,-1\vert \infty,-A)$ or $\sigma(-1,\infty\vert 0,-A)$.
	In the first case, $\sigma(0,\infty\vert -A,\un{A\cdot B})$ and
	$\sigma(0,\un{-1}\vert \infty,-A)\implies \sigma(-1,A\cdot B\vert
	0,\infty).$\\
	In the second case, $\sigma(0,\infty\vert -A,\un{A\cdot B})$, and
	$\sigma(\un{-1},\infty\vert 0,-A)\implies \sigma(-1,A\cdot B\vert
	0,\infty)$.

	\sssec{Theorem.}
	{\em With the coordinatization of the separable Pappian plane as
	given in \ref{sec-dprojter}},
	\enumb
	\item {\em ($\Sigma,+,\cdot)$ is an ordered field.}
	\enume

	\ssec{Continuous Pappian or Classical Projective Planes.}

	\sssec{Axiom. [Of continuity]}
	Let ${\cal S} \subset(seg(A,C\:\setminus\:B)$, ${\cal S}$ non empty,
	$\exists L$ and $U \ni$ all $P \in {\cal S},$ $\sigma(AP\vert LU)
	\implies \exists G$ and $H$ $\ni$ $\sigma(LP\vert GU)$ and
	$\sigma(UP\vert HL)$.

	\sssec{Definition.}
	A {\em Continuous Pappian} or {\em Classical Projective Plane}
	is a separable plane for which the continuity axiom is satisfied.

	\sssec{Theorem.}
	{\em The field associated to a Continuous Pappian plane is the real
	field ${\Bbb R}$.} 

	\ssec{Isomorphisms of Synthetically and Algebraically defined Planes.}

	\sssec{Introduction.}
	We have seen that we can coordinatize the various perspective planes
	by ternary rings which have special properties.  The converse is also 
	true.  If a ternary rings has appropriate properties there exists
	a plane as defined above which is isomorphic to it.
	More specifically:

	\sssec{Theorem.}
	{\em There is an isomorphism between}
	\enumb
	\item {\em perspective planes and ternary rings $(\Sigma,*)$.}
	\item {\em Veblem-Wedderburn planes and ternary rings with the
	properties \ref{sec-tveblenw}.}
	\item {\em Moufang planes and alternative division rings.}
	\item {\em Desarguesian planes and skew fields.}
	\item {\em Papian planes and fields.}
	\enume

	\ssec{Examples of Perspective Planes.}

	\sssec{Definition.}
	A {\em Moulton plane} (1902) is the set of points in the Euclidean plane
	coordinatized with Cartesian coordinates and the lines,
	\enumb
	\item the ideal line, [0,0,1],
	\item the lines [$m,-1,n$], $m\leq 0,$
	\item the lines consisting of two parts, first, the subset of
	[$m,-1,n$], $m > 0,$ which is in the lower half plane or on the ideal
	line, second, the subset of [$m/2,-1,n$], $m > 0,$ which is in the
	upper half plane.
	\enume

	\sssec{Theorem.}
	\enumb
	\item{\em The Moulton plane is a perspective plane.}
	\item{\em The Moulton plane is not a Veblen-Wedderburn plane.}
	\enume

	Proof: See Artzy, p. 210.

	\sssec{Definition.}
	A {\em 2-Q plane} is defined like a quaternion plane with
	$ij = -ji = k$ replaced by $ij = -ji = 2k$.

	\sssec{Theorem.}
	\enumb
	\item{\em The 2-Q plane is a Veblen-Wedderburn plane.}
	\item{\em The 2-Q plane is not a Moufang plane.}
	\enume

	Proof: See Artzy, p. 226.

	\sssec{Definition.}
	A {\em Cayleyian plane} is defined like a quaternion plane, using
	Cayley numbers instead of quaternions.

	\sssec{Theorem.}
	\enumb
	\item{\em The Cayleyian plane is a Moufang plane.}
	\item{\em The Cayleyian plane is not a Desarguesian plane.}
	\enume

	Proof: See Artzy, p. 226.

	\sssec{Definition.}
	A {\em quaternion plane} is defined using quaternions as coordinates
	instead of real numbers.

	\sssec{Theorem.}
	\enumb
	\item{\em The quaternion plane is a Desarguesian plane.}
	\item{\em The quaternion plane is not a Pappian plane.}
	\enume

	Proof: See Artzy, p. 226.

	\sssec{Definition.}
	A {\em finite Pappian plane} is a Pappian plane for which the
	number of points on one line is finite. The field associated to it
	is therefore a finite field which is necessarily a Galois field
	$GF(p^k)$ with $p$ prime, the number of points being $p^k+1$.

	\sssec{Theorem.}
	\enumb
	\item{\em The finite Pappian plane is a Pappian plane.}
	\item{\em The finite Pappian plane is not a Separable Pappian plane.}
	\enume

	Proof: See Artzy, p. 210.

	\sssec{Definition.}
	If the field is the field of rationals, the Pappian plane is called
	the {\em Rational Pappian plane}.

	\sssec{Theorem.}
	\enumb
	\item{\em The Rational Pappian plane is a Separable Pappian plane.}
	\item{\em The Rational Pappian plane is not a Continous Pappian plane.}
	\enume

	Proof: See Artzy, p. 210.

	\sssec{Exercise.}
	Give a synthetic definition of a
	\enumb
	\item The rational Pappian plane.
	\item The quaternion plane.
	\item The Cayleyian plane.
	\item 2-Q plane.
	\enume
	It is clear how to proceed for the rational plane, immitating the
	definition of the rational numbers as equivalence classes of the
	integers.  It is not known to me how to solve the other exercises.

	\ssec{Collineations and Correlations in Perspective to Pappian Planes.}
\setcounter{subsubsection}{-1}
	\sssec{Introduction.}
	For collineations, correlations and polarities in finite planes, see
	Dembowski, section 3.3 and Chapter 4.

	\sssec{Definition.}
	Given \ref{sec-dtrpr},we say that the {\em vectors $\ov{AA'}$ and
	$\ov{PP'}$ are $m$-equal} and we write\\
	\hth$\ov{AA'} =_m \ov{PP'}.$

	\sssec{Definition.}
	In a Veblen-Wedderburn plane with ideal line $m$, the elements of
	the {\em set} ${\cal V}$ are the equivalence classes of $m-$equal
	vectors and the {\em addition of vectors} is defined by\\
	\hth$\ov{P_1P_2} + \ov{Q_2Q_3} := \ov{P_1P_3}$,\\
	where
	\enumb
	\item $Q_2 = P_2 \implies P_3 = Q_3,$
	\item $P_2$, $Q_2$, $Q_3$ non collinear $\implies P_3$ is the point
	defined by $\ov{P_2P_3} =_m \ov{Q_2Q_3}$,
	\item if $P_2$, $Q_2$, $Q_3$ collinear and $X$ is not on $Q_2\times Q_3
	\implies P_3$ is the point defined by $\ov{P_2P_3} =_m \ov{Q_2Q_3}$,
	$\ov{P_2P_3} =_m \ov{Q_2Q_3}$.
	\enume

	\sssec{Theorem.}
	{\em The addition of vectors is well defined and (${\cal V}$,+)
	is an abelian group.}

	This follows from the fact that any vector is equivalent to a vector
	$\ov{(0,0)(A,B)}$ for some $A$ and $B$ \ref{sec-tveblenw}.1.

	\sssec{Theorem.}
	{\em The translations in a Veblen-Wedderburn plane with ideal line $m$
	are collineations, in other words, the image of points $P$ on a
	fixed line $l$ are points $P'$ on a line $l'$. Each collineation is an
	elation with axis $m$ and center $m\times (A\times A')$.}

	\sssec{Theorem.}
	{\em For any line $n$, the $n$-translations in a Moufang plane are
	collineations, in other words, the image of points $P$ on a
	fixed line $l$ are points $P'$ on a line $l'$. Each collineation is an
	elation with axis $n$ and center $n\times (A\times A')$.}

	\sssec{Definition.}
	In a in Veblen-Wedderburn Plane the {\em pre correlation
	configuration} is defined as follows, (See Fig. 25)\\
	Hy0.\hti{4}$\{Q_i\}$, $u\incid Q_2$, $i,a,b,a' \incid Q_0$,\\
	De.\hti{4}$U_1 := u \times a,$ $D_2 := d \times a,$
		$x_1 := U_1\times Q_1,$ $C_1 := x1\times i,$\\
	De.\hti{4}$c := C_1\times Q_2,$ $C_2 := c\times b,$
		$U_2 := b\times u,$ $x_2 := U_2\times Q_1,$\\
	De.\hti{4}$D_1 := x2\times i,$ $d := D_1\times Q_2,$
		$D_2 := d\times a,$ $x_3 := D_2\times Q_1,$\\
	De.\hti{4}$D'_2 := a'\times x_3,$ $d' := D'_2\times Q_2,$
		$D'_1 := d'\times i,$ $x'_2 := D'_1\times Q_1,$\\
	De.\hti{4}$U'_2 := x'2\times u,$ $b' := U'_2\times Q_0,$
		$B' := b'\times q_2,$ $U'_1 := a'\times u,$\\
	De.\hti{4}$x'_1 := U'1\times Q_1,$ $C'_1 := x'1\times i,$
		$c' := C'_1\times Q_2,$ $C'_2 := c'\times b',$\\
	De.\hti{4}$x'_3 := C_2\times Q_1,$\\
	Hy1.\hti{4}$C'_2 \incid x'3,$

	Let $(A) = a\times q_2,$ $(B) = b\times q_2,$
	$(A') = a'\times q_2,$ $(B') = b'\times q_2,$\\
	then $D_2 = (B,B\cdot A,$ $D'_2 = (B',B'\cdot A'),$
	$C_2 = (A,A\cdot B,$ $C'_2 = (A',A'\cdot B'),$
	If $b'$ or $B'$ is chosen in such a way that $B\cdot A = B'\cdot A'$,
	the configuration requires $A\cdot B = A'\cdot B'$.
	This defines a correspondance $\gamma$ between $X = A\cdot B$ and
	$X' = B\cdot A$.

	\sssec{Exercise.}
	If we associate to $(Q),$ $[Q]$ and to $(P_0,P_1),$
	$[P_0\gamma,P_1\gamma],$ is the correspondance is a correlation?\\
	If not which of the axioms given below are required for the
	correspondance to be a correlation.

	\ssec{Three Nets in Perspective Geometry.}

	\sssec{Definition.}
	A {\em three net associated to the 3 points $A$, $B$, $C$} in a 
	perspective plane is the set of points $P$ in the plane and
	the set of lines $P\times A,$ $P\times B,$ $P\times C$.

	\sssec{Theorem.}
	{\em The coordinates of the lines of the three net associated to the
	points (0), (1), $(\infty)$ are $[0,P_0]$, $[1,P_1]$, $[P_2]$, where}\\
	\hth$P_0 := (((P \times (0))  \times v) \times Q_0) \times q_2$,\\
	\hth$P_1 := (((((P \times (1)) \times q_1))\times (0)) \times v) \times
	Q_0) \times q_2$,\\
	\hth$P_2 := (((((P\times (\infty))\times i)\times (0))\times v)\times
	Q_0)\times q_2$.

	\sssec{Lemma.}
	{\em Let $Y_A = (0,A),$ $Y_B = (0,B),$ then}\\
	$(((Q_0\times(1))\times(Y_A\times(0)))\times(\infty))\times
	(((Q_0\times(\infty))\times(Y_B\times(0)))\times(1)) = (A,A+B).$

	\sssec{Definition.}
	Given $Q'_0 = (F,F+G)$, the {\em F-G-sum of} $A$ and $B$, $A\oplus B$
	is defined by\\
	$(((Q'_0\times(1))\times(Y_A\times(0)))\times(\infty))\times
	(((Q'_0\times(\infty))\times(Y_B\times(0)))\times(1)) = (X,A\oplus B).$

	\sssec{Theorem.}
	\hth$A\oplus B = (A\dashv G)+(F\vdash B).$\\
	\hth$X = A\dashv G$.

	Proof:\\
	$X_A := ((Q'_0\times(1))\times(Y_A\times(0)) = (X,A)$ and $X+G = A$,
	therefore $X = A\dashv G$.\\
	$F_B := ((Q'_0\times(\infty))\times(Y_B\times(0)) = (F,B),$\\
	if $Y_Z := ((F,B)\times(1))\times q_1 = (0,Z),$ then $C = F+Z$ and
	$Z = F\vdash C,$\\
	finally $X_Y := (X_A\times(\infty))\times(F_B\times(1) =
	(X_A\times(\infty))\times(Y_Z\times(1) = (X,A\oplus B),$
	therefore $(A\oplus B) = X+Z,$ substituting for $X$ and $Z$
	gives the Theorem.

	\sssec{Theorem.}
	{\em $(\Sigma,\oplus)$ is a loop.\\
	The neutral element is $F+G$.\\
	The solutions of $A\oplus B = C$ are given by}\\
	\hth$A = (C\dashv (F\vdash B)) + G,$ $B = F + ((A\dashv G)\vdash C).$

	Proof:
	The solutions follow directly from the preceding Theorem, the neutral
	element property follows from\\
	\hth$(F+G)\oplus H = ((F+G)\dashv G)+(F\vdash H) = F+(F\vdash H) = H,$\\
	\hth$H\oplus(F+G) = (H\dashv G)+(F\vdash(F+G)) = (H\dashv G)+G = H$.

	\sssec{Theorem.}
	{\em The coordinates of the lines of the three net associated to the
	points $Q_0,$ $Q_1$, $Q_2$ are $[P_1,0]$, $[0,P_0]$, $[P_2]$, where}\\
	\hth$P_0 := (((P \times (0)) \times v) \times Q_0) \times q_2$,\\
	\hth$P_1 := (P \times Q_0) \times q_2$,\\
	\hth$P_2 := (((((P\times (\infty))\times i)\times (0))\times v)\times
	Q_0)\times q_2$.

	\sssec{Exercise.}
	Determine Theorems analogous to those associated with (0), (1) and
	$(\infty)$.  See Artzy, p. 206 and p.210, 15.

\newpage
	\ssec{Bibliography.}
	\enumb
\item Artin, Emil, {\em Geometric Algebra}, New York, Interscience Publishers,
	1957. Interscience tracts in pure and applied mathematics, no. 3.
\item Artzy, Rafael, {\em Linear Geometry}, Reading Mass., Addison-Wesley, 1965,
	273 pp.
\item Bolyai, Farkas, {\em Tentamen Juventutem Studiosam ein Elementa Mathiseos
	Parae, introducendi}, Maros-Vasarhely, 1829, see Smith D. E.
	p. 375.
\item Bolyai, Janos, {\em The Science Absolute of Space Independent of the Truth
	and Falsity of Euclid's Axiom XI}, transl. by Dr George Brus
	Halstead, Austin, Texas, The Neomon, Vol. 3, 71 pp, 1886.
\item Bolyai, Janos, {\em Appendix, the theory of space}, with introduction,
	comments, and addenda, edited by Ferenc Karteszi, supplement by Barna
	Szenassy, Amsterdam, New York, North-Holland, New York, Sole
	distributors for the U.S.A. and Canada, Elsevier Science Pub. Co.,
	1987, North-Holland mathematics studies, 138.
\item Bruck R. H. and Ryser H. J., {\em The non existence of certain finite 
	projective planes}, Can. J. Math., Vol. 1, 1949, 88-93.
\item Dedekind, Julius Wilhelm, {\em Stetigkeit und Irrationalen Zahlen}, 1872,
	see Smith D. E., p. 35. (I,9,p.1)
\item Dembowski, Peter, {\em Finite Geometries}, Ergebnisse der Mathematik und
	ihrer Grenzgebiete, Band 44, Springer, New-York, 1968, 375 pp.
\item Enriques, Federigo, {\em Lezioni di Geometria Proiettiva}, Bologna, 1904,
	French Transl., Paris 1930.
\item Enriques, Federigo, {\em Lessons in Projective Geometry}, transl. from
	the Italian by Harold R. Phalen. Annandale-on-Hudson, N.Y., printed by
	the translator, [1932?].
\item Euler, Leonhard, {\em Recherches sur une nouvelle esp\`{e}ce de
	quarr\'{e}s magiques}, Verh. Zeeuwsch. Genootsch. Wetensch. Vlissengen, 
	Vol. 9, 1782, 85-239.
\item Fano, Gino, {\em Sui Postulati Fondamentali della Geometria Proiettiva},
	Giorn. di mat., Vol. 30, 1892, 106-132. (PG($n,p$))
\item Hall, Marshall, Jr, {\em Projective Planes}, Trans. Amer. Math. Soc.,
	Vol. 54, 1943, 229-277.
\item Hartshorne, Robin C., {\em Foundation of Projective Geometry}, N. Y.
	Benjamin, 1967, 161 pp.
\item Hilbert, David, {\em Grundlagen der Geometrie}, 1899, tr. by E. J.
	Townsend, La Salle, Ill., Open Court Publ. Cp., 1962, 143 pp.
\item Hilbert, David, {\em The Foundations of Geometry}, authorized transl.
	by E.J. Townsend \ldots Chicago, The Open court publishing company;
	London, K. Paul, Trench, Trubner \& co., ltd., 1902.
\item Klingenberg, Wilhelm, {\em Beweis des Desargueschen Satzes aus der
	Reidemeisterfigur und Verwandte S\"{a}tze}. Abh. Math. Sem. Hamburg,
	Vol. 19, 1955, 158-175.
\item Klingenberg, Wilhelm, {\em Grundlagen der Geometrie}, Mannheim,
	Bibliographisches Institut, 1971, B. I.-Hochschul\-skripten 746-746a.
\item Lam, C. W. H., {\em The Search for a Finite Projective Planes of Order
	10}, Amer. Math. Monthly, Vo. 98, 1991, 305-318.
\item Lam, C. W. H., Thiel, L. H. \& Swiercz S., {\em The non existence of
	finite projective planes of order 10}, Can. J. of Mat., Vol. 41, 1989,
	1117-1123.
\item Lobachevskii, Nikolai Ivanovich, see Norden A., {\em Elementare
	Einfuhrung in die Lobachewskische Geometrie}, Berlin, VEB Deuscher
	Verlag der Wissenschaften, 1958, 259 pp.
\item Lobachevskii, Nikolai Ivanovich, {\em Geometrical Researches on the
	Theory of Parallels}, translated from the original by George Bruce
	Halsted, Austin, University of Texas, 1891.
\item Menger, Karl, {\em Untersunchungen \"{u}ber allgemeine Metrik}, Math. Ann.
	Vol. 100, 1928, 75-163.
\item Moufang, Ruth, {\em Alternatievk\"{o}rper und der Satz vom
	Vollst\"{a}ndigen Vierseit}, Abh. Math. Sem. Hamburg, Vol. 9, 1933,
	207-222.
\item Pickert. Gunter, {\em Projektive Ebenen}, Berlin, Springer, 1955, 343 pp.
\item Pieri, {\em Un Sistema di Postulati per la Geometria Proiettiva}, Rev.
	Math\'{e}m. Torino, Vol 6, 1896. See also Atti Torino, 1904, 1906.
\item Pieri, {\em I Principii della Geometria di Posizione, composti in Sistema
	Logico Deduttivo}, Mem. della Reale Acad. delle Scienze di Torino,
	serie 2, Vol.48, 1899, 1-62.
\item Reidemeister, Kurt, {\em Grundlagen der Geometrie}, Berlin, Springer,
	Grundl. der math. Wissens. in Einz., Vol. 32, 1968, (1930),
\item Saccheri, Giovanni Girolamo, {\em Euclides ab omni Naevo Vindicatus},
	Milan, 1732. tr. George Halstead, London Open Court Pr. 1920, 246 pp.
	See St\"{a}ckel.
\item Schur, Friedrich, {\em Grundlagen der Geometrie}, mit 63 figuren im
	text. Leipzig, Berlin, B. G. Teubner, 1909.
\item Schur, Issai, {\em Gesammelte Abhandlungen}, Hrsg. von Alfred Brauer u.
	Hans Rohrbach, Berlin, Heidelberg, New York: Springer, 1973.
\item Stackel, Paul Gustav, {\em Die Theorie der Parallellinien von Euklid bis
	auf Gauss, eine Urkundensammlung zur Vorgeschichte der
	nichteuklidischen Geometrie, in Gemeinschaft mit Friedrich Engel},
	hrsg. von Paul Stackel, New York, Johnson Reprint Corp., 1968,
	Bibliotheca mathematica Teubneriana, Bd. 41.
\item Tilly, Joseph Marie de, {\em Essai sur les Principes fondamentaux de
	G\'{e}om\'{e}trie et de M\'{e}canique}, Bruxelles, Mayolez, 1879, 192
	pp. Also, M\'{e}m. Soc. science phys. et natur. de Bordeaux, Vol III,
	Ser. 2, cahier 1.
\item Tilly, Joseph Marie de, {\em Essai de G\'{e}om\'{e}trie Analytique
	G\'{e}n\'{e}rale}, Bruxelles, 1892.\\
\item Tarry, G., {\em Le probl\`{e}me des 36 officiers}, C. R. Assoc. Franc.
	Av. Sci., Vol. 1 (1900), 122-123, Vol. 2 (1901), 170-203.

\item Veblen, Oswald \& Wedderburn, Joseph Henri MacLagan,
	{\em Non-Desarguesian and non-Pascalian Geometries}, Trans. Amer. Math.
	Soc., Vol. 8, 1907, 279-388.
\item Veblen, Oswald \& Young, John, {\em Projective Geometry}, Wesley,
	Boston, I, 1910, II, 1918. 
	\enume

\setcounter{section}{10}
	\section{Mechanics.}
	\vspace{-18pt}\hspace{102pt}\footnote{30.10.87}\\[1pt]

	\setcounter{subsection}{-1}
	\ssec{Introduction.}

	Geometry is to be the support of the description of phenomenon in the
	real world.  I will briefly review Newton's laws and 2 results to be
	generalized, the central force theorem of Hamilton and the motion of the
	pendulum.

	\ssec{Kepler (1571-1630).}

	\setcounter{subsubsection}{-1}
	\sssec{Introduction.}
	Among many of the contribution of Kepler those which perpetuate his
	name are his 3 laws of Mechanics and his equation discovered from 1605
	to 1621. The first and third law are in Astronomi Nova, the second law
	and his equation in section V of his Epitome.
	 $\ldots$  We also know that an ellipse can be generated by
	moving a segment of length $a+b$ with one point on an axis and the other
	point on a perpendicular axis.  \ref{sec-tkepl}.2 shows that the angle
	of the line is also the eccentric anomaly.?

	\sssec{Theorem.}\label{sec-tkepl}
	{\em If a point $P(x,y)$ is restricted to move on an ellipse
	with major axis $2a,$ minor axis $2b$ and eccentricity $e,$ and origin
	at a focus,}
	\enumb
	\item $x = a(cos \circ E - e),$ $y = b\: sin \circ E,$\\
	{\em where $E$ is the excentric anomaly.\\
	If$E(0) = 0,$ then $x(0) = ae,$ $y(0) = 0.$\\
	Let $I$ be the identity function, the motion preserves area  iff
	{\em Kepler's equation}}
	\item $I = E - e sin \circ E$\\
	{\em is satisfied.\\[5pt]
	If $v := \angle(A,F,P),$ called {\em true anomaly} then}
	\item 	$tan v = \frac{b sin E}{a(cos E - e)}$.
			
	{\em Finally, if the line through $P$ makes an angle $E$ with $AF,$
	and intersect the major axis at $L$ and the minor axis at $M,$}
	\item 	$PL = b,$ $PM = a.$\\[5pt]
	{\em Let $A$ be twice the area $(0,0),$ $(a,0),$ $(x,y)$ along the
	ellipse, divided by $ab.$\\
	Let $T$ be twice the area of the triangle $(0,0),$ $(x,y),$ $(x',y')$
	divided by $ab,$ then}
	\item 	$T = x y' - x'y\\
	\hti{2}	  = (cos \circ E - e) sin \circ E' - (cos \circ E' - e)
			 sin \circ E\\
	\hti{2}	  = sin \circ (E'-E) - e (sin \circ E' - sin \circ E),$\\[5pt]
	{\em If $E' = E + \Delta E,$ and $\Delta E$ is small, then}
	\item $\Delta A = (1 - e cos \circ E) \Delta E.$\\[5pt]
	Integrating gives
	\item $A = E - e \:sin \circ E.$
	\enume

	Therefore, if the area $A$ is a linear function, with a proper
	choice of the unit of time, $A = I$ and we have 1.  Vice-versa,
	if 1. is satisfied then comparing 6, and 1, gives $A = I$ and the area
	is proportional to the time.

	\ssec{Newton (1642-1727).}

	\ssec{Hamilton (1805-1865).}

	\sssec{Theorem. [Hamilton]}\label{sec-thamil}
	{\em Assuming Newton's law, if a mass is to move on an ellipse, under a
	force passing through a fixed point (central force),}
	\enumb
	\item[0.~~]{\em this force is proportional to the distance to the center
		and inversely proportional to the cube of the distance to
		the polar of the center of force.}
	\item[1.~~]{\em the relation between the eccentric anomaly $E$ and the
		time $t$ is given by}
	\hth	$a \:E(t) + c \:sin(E(t)) = C \:t.$

	{\em Consider a conic with major axis of length $2a$ on the $x$ axis,
	with minor axis of length $2b$ and with center at $(c,d),$
	the parametric representation is}
	\item[2.~~]$x = c + a \:cos \circ E,$ $y = d + b\: sin \circ E.$

	{\em The acceleration is}
	\item[3.0.]$D^2 x = -a\:cos \circ E (DE)^2 - a\:sin \circ E D^2 E,$
	\item  [1.] $D^2 y = -b\:sin \circ E (DE)^2 + b\:cos \circ E D^2 E,$

	If we accept Newton's law, the acceleration has to be in the
	direction of the force, if the force is $f \circ E\:g \circ E,$ where
	\item[4.~~]$(g \circ E)^2 = (c + a\:cos \circ E)^2 + (d + b\:sin
		\circ E)^2,$\\
	$g \circ E$ being the distance to the center of force,
	\item[5.0.]$D^2 x = f \circ E (c + a\:cos \circ E),$
	\item  [1.] $D^2 y = f \circ E (d + b\:sin \circ E),$
	\item[6 0.]$-a\:cos \circ E (DE)^2 - a\:sin \circ E D^2 E
		= f \circ E (c + a\:cos \circ E),$
	\item  [1.]$-b\: sin \circ E (DE)^2 + b\:cos \circ E D^2 E
		= f \circ E (d + b\:sin \circ E),$\\
	hence equating 3.0 and 5.0 as well as 3.1 and 5.1 we get 6.0 and
	6.1, the combinations\\
	\hth$(-b\:sin \circ E)\: 6.0. + (a\:cos \circ E)\:6.1.$ and\\
	\hth$(-b\:cos \circ E)\:6.0. - (a\:sin \circ E)\:6.1.$\\
	give
	\item[7.0.]$ab\:D^2 E = f \circ E (ad\:cos \circ E - bc\:sin \circ E),$
	\item  [1.]$ab\:(DE)^2 = - f \circ E (ad\:sin \circ E + bc\:cos \circ E
		+ ab).$

	Taking the derivative of this equation and subtracting $2 DE$ times 7.0.
	gives
	\item[8.~~]$D(f \circ E) (ad\:sin \circ E + bc\:cos \circ E + ab) +
		3 f \circ E (ad\:cos \circ E - bc\:sin \circ E) DE  = 0.$

	Integrating gives
	\item[9.~~]$f \circ E (ad\:sin \circ E + bc\:cos \circ E + ab)^3
		= - a C_1$\\
	for some constant $C_1,$
	but the polar of the origin is the line\\
	\hth$b^2c x + a^2d y -b^2c^2 - a^2d^2 + a^2b^2 = 0,$\\
	therefore the distance of $(x,y)$ to it is proportional to\\
	\hth$b^2c (c + a cos \circ E) + a^2d (d + b sin \circ E) -b^2c^2
		- a^2d^2 + a^2b^2$\\ 
	or to
	\item[10.~]$bc\:cos \circ E + ad\:sin \circ E + ab,$\\
	hence part 1 of the theorem.

	Replacing in 7.1. $f \circ E$ by its value gives\\
	\hth$(DE)^2 = \frac{C}{(ad\: sin \circ E + bc\: cos \circ E + ab)^2},$\\
	therefore $C_1$ must be positive.

	{\em Let $C_1 = C^2,$ then}
	\item[11.~] $(ad\:sin \circ E + bc\:cos \circ E + ab) DE = C,$\\
	{\em and we obtain a generalization of Kepler's equation}
	\item[12.~]$-ad\: cos \circ E + bc\:sin \circ E  + ab\:E  = C I,$

	{\em Let $e$ and $A$ be such that}
	\item[13.~]$bc = ab\:e\: cos(A), ad = ab\:e\: sin(A),$\\
	{\em then}
	\item[14.~]$e^2 = (\frac{c}{a})^2 + (\frac{d}{b})^2,$ 
	\item[15.~]$tan(A) = \frac{ad}{bc}.$

	{\em Let}
	\item[16.~]$F = E-A$ and $M = CI - ab\:A,$

	{\em then}
	\item[17.~]$e\:sin(F) + F = M.$
	\enume

	\sssec{Comment.}
	If the center of the conic is the center of force, $c = d = 0,$
	$f \circ E$ is a constant and the force is proportional to the
	distance.  If the center of force is on the conic,
	\ref{sec-thamil}.7 becomes\\
	\hth$f \circ E (ab)^3 (1 - cos \circ (E-E_0))^3 = - a C_1,$\\
	when the center of force is $c + a cos(E_0),$ $d + b sin(E_0).$\\
	When the conic is a circle,\\
	\hth$g \circ E^2 = 2 a^2(1 - cos \circ (E-E_0))^2,$ \\
	therefore, the force is inversely proportional to the 5-th power
	of the distance.?

	\sssec{Comment.}
	\ref{sec-thamil}.8 is proportional to $g \circ E$ if\\
	\hth$h^2((c	+ a \:cos \circ E)^2 + (d + b\: sin \circ E)^2) 
			= (a\:d\:sin \circ E + b\:c\:cos \circ E + a\:b)^2$ \\
	expanding will give terms in $cos^2,$ $sin\:cos,$ $sin,$ $cos$ and $1.
$\\
	The coefficient of $sin\:cos$ must be $0,$ hence $cd = 0.$\\
	Let $d = 0,$ the $sin$ term disappears and the coefficients of
	$1,$ $cos \circ E,$ $cos^2 \circ E$ give\\
	\hth$	h^2(c^2 + b^2) = a^2 b^2,$\\
	\hth$	h^2(2 a c) = 2 b^2 a c,$\\
	\hth$	h^2(a^2 - b^2) = b^2 c^2,$\\
	hence\\
	\hth$h = b$ and $a^2 = b^2 + c^2$\\
	or\\
	\hth$e := \frac{c}{a} = \sqrt{1-\frac{b}{a}^2}.$ 

	\sssec{Definition.}
	Given a curve $(x,y),$ the {\em hodograph of the curve} is the
	curve $(Dx,Dy).$

	\sssec{Comment.}
	The concept was first introduced by M\"{o}bius (Mechanik des
	Himmels, (1843)), the name was chosen by Hamilton when he gave,
	independently,	the definition in the Proc. Roy. Irish Acad., Vol. 3,
	(1845-1847) pp. 344-353.

	\sssec{Theorem.}
	{\em If the force is central, and the center is chosen as the
	origin, the hodograph of the hodograph is the original curve.}

	Indeed, the hodograph is $(D^2x,D^2y) = f \circ E (x,y).$

	\sssec{Theorem.}
	{\em If the central force obeys Newton's law, the hodograph of
	the ellipse, 0.0. is the circle}\\
	\hth$b^2((Dx)^2 + (Dy)^2) + 2a\:e\:C Dy - C^2 = 0,$

	The proof is straightforward, the verification using\\
	\hth$Dx = -a\: sin \circ E DE,$ $Dy = b\: cos \circ E DE$ and
	\ref{sec-thamil}.9, .16 is even simpler.

	If the equation of the circle is\\
	\hth$     (-R sin(G),$ $-k + R cos(G)),$\\
	equating to $(-a sin \circ E DE, b cos \circ E DE)$
	for $E = 0$ and $\pi$ and therefore $G = 0$ and $\pi,$ gives\\
	\hth$	b C = a b(-k + R)(1 + E),$ $-b C = a b(-k-R)(1-e),$
	therefore\\
	\hth$	-k + R = \frac{C}{a(1+e)}),$ $k+R = \frac{C}{a(1-e)},$\\
	hence\\
	\hth$	R = C (\frac{a}{b})^2,$ $k = R e.$\\
	moreover\\
	\hth$	cos(G) = e + \frac{b^2}{a^2} 
		\frac{cos \circ E}{1 + e\:cos \circ E}.$

	\ssec{Preliminary remarks extending mechanics to finite geometry.}

	\setcounter{subsubsection}{-1}
	\sssec{Introduction.}
	The generalization of classical mechanics to
	finite geometry turned out to be a thorny task.

	\sssec{Lemma.}
	{\em If $x_0,$ $y_0$ is a solution of\\
	\hth$	x (p+1) - y p = 1,$\\
	all solutions are given by\\
	\hth$	x = x_0 + k p, y = y_0 + k(p+1),$\\
	or\\
	\hth$	x \equiv x_0 \pmod{p}, y \equiv y_0 \pmod{p+1}.$}

	\sssec{Definition.}
	{\em Kepler's equation} associated to the prime $p$ is given by\\
	\hth$	(e\: sin \circ E) (p + 1) - (E - M) p = 1.$

	This definition can be justified as follows, first when $p$ is very
	large,
	we get the classical Kepler equation.  Moreover from Lemma 10.4.1.
	all solutions are such that $e\:sin \circ E$ are equal modulo $p$ and
	$E-M$ are equal modulo $p+1$ which are precisely the congruence
	relations for $e,$ $sin \circ E$ and for $E-M.$

	\sssec{Example.}
	For $p = 101,$ $\ldots$
	
	\sssec{Theorem.}
	{\em (Of the circular hodograph of Hamilton). $\ldots$}	

	\ssec{Eddington (18?-1944).  The cosmological constant.}

	Starting with the work of Edwin P. Hubble,(1934) there had been mounting
	observational astronomical evidence that the Universe is finite.
	This lead, Monseigneur Georges Lema\^{\i}tre to his hypothesis of
	the Primeval Atom and Sir Arthur Eddington to a possible a priori
	determination of the cosmical number $N = 3.68.2^{256} = 2.36 10^{79}.$
	In his article published in 1944, in the Proc. of the Camb. Phil. Soc.,
	he first describes the number ``picturesquely as the number of protons
	and electrons in the universe" and ``interprets it by the consideration
	of a distribution of hydrogen in equilibrium at zero temperature,
	because the presence of the matter produces a curvature in space, the
	curvature causes the space to close when the number of particles
	contained in it reaches the total $N$".

	If the work of Eddington would be reexamined today, protons and
	electron would
	probably be replaced by quarks, if it were to be reexamined at some
	time in the future some other particles might play the fundamental
	role.  In any case the lectures of Lema\^{\i}tre and the work of
	Eddington have been a primary motivation for my work on finite
	Euclidean and non-Euclidean geometry.  As will be examined in more
	details when application will be made to the finite pendulum, some
	elementary particle occupies a position and the possible positions
	are discrete, they do this at a certain time, but again the time is
	not a continuous function but a discrete monotonic function.  The
	fact that there are no infinitesimals in finite geometry may very
	well be related to the uncertainty principle of Heisenberg (1927).\\
	H. Pierre Noyes and ANPA

	\section{Description of Algorithms and Computers.}

%
	All the earlier proofs in Mathematics were constructive, these proofs
	not only showed the existence of objects, for instance the existence of
	the orthocenter of a triangle, where the 3 perpendiculars from the
	vertex to the opposite sides meet, but also how to construct that point,
	by giving an explicit construction for a perpendicular to a line from a
	point outside it.  Little by little mathematicians have used more and
	more proofs using non constructive arguments, which show the existence
	of the object in question, without giving a method of construction.
	Such proofs are essential when no finite construction is possible,
	and are considered by many as intellectually superior to a constructive
	proof when this one is possible.  In finite geometry, it is desirable to
	limit oneself to constructive proofs, although this is not always
	possible, at a given point in time.
	I will give 2 examples later, the proof of Aryabatha's theorem and the
	proof of the existence of primitive roots.
	Because in a finite geometry it is not easy to rely on tools such as
	the straightedge or the compass to experiment for the purpose of
	conjecturing theorems, it is useful if not necessary to rely on computer
	experiments.  Moreover, although the simpler algorithm were for
	centuries given in the vernacular language, see for instance the
	description of the so called Chinese remainder theorem by
	Ch'in Chiui-Shao, in Ulrich Libbrecht's translation, often the
	description avoids special cases or is ambiguous.  Careful description
	of algorithms started to appear with the advent of computers.
	\footnote{Already in 1957, Lema\^{\i}tre used precise descriptions
	to communicate by letter with a person doing his calcultations on
	a EUCLID mechanical calculator.}

	The first formula oriented language was FORTRAN which evolved to FORTRAN
	4 then FORTRAN 77. It was developed enpirically. ALGOL was developed in
	1958 and its syntax carefully defined in 1960 using the Backus 
	normal form to attempt to define a priori an algorithmic
	language with a carefully constructed block structure. Its immediate
	successors were ALGOL 68 and PASCAL. APL was developped by Iverson to
	describe carfully the logic of the hardware of computers. It was
	magistrally adapted for the programming of Mathematical problems.
	LISP and its family of languages were developed when a list structure is
	required. BASIC was created at Dartmouth, to allow all undergraduates
	to learn programming in a friendly environment. It is the language
	which has evolved the most since its early days especially by a small 
	group at the Digital Equipment Corporation. This is the language
	which I found most useful to discover mathematical conjectures because
	of the flexibility it offers in changing the program while in core
	and in examining easily, when needed, intermediate results without
	prior planning. MAXIMA and its family of languages, MABEL, MATHEMAICA
	and other recently developed languages are sure to play a more and
	more important role in discoveries.

	Elsewhere, I will describe some of the BASIC programs, that I have
	written to investigate new areas of Mathematics, as well as the style
	used in the program descriptions and in their documentation and use.
%

	\section{Notes.}

	\ssec{On Babylonian Mathematics.}\label{sec-babyl}

	Besides estimating areas and volumes, the Babylonians had a definite
	interest in so called Pythagorian triples, integers $a,$ $b$ and $c$
	such that $a^2 = b^2 + c^2.$  It is still debated if their interest was
	purely arithmetical or was connected with geometry.  On the one hand
	Neugebauer, states\\
	"It is easy to show that geometrical concepts play a very secondary
	part in Babylonean algebra, however extensively a geometrical
	terminology is used." (p. 41)\\
	However, more recent discoveries, let him state (p.46), that these
	"contributions lie in the direction of geometry".  One tablet computes
	the radius $r$ of a circle which circumscribes an isosceles triangle
	of sides $50,$ $50$ and $60$.  An other tablet gives the regular
	hexagon, and from this the approximation
	$\sqrt{3} = 1;45 (1 + \frac{45}{60})$
	can be deduced.  \ldots  $(\sqrt{2}= 1;25),$  \ldots  $\pi = 3;$
		$ 7,30 (3 \frac{1}{8}),  \ldots  ".$\\
	He also describes, with Sachs, the data contained
	in tablet 322 of the Plimpton library collection from Columbia
	University (see Neugebauer and Sachs, vii and 38-41) as clearly
	indicating a relationship with right triangles "with angles varying
	regularly between almost 45 degrees to almost 31 degrees", while
	Bruins interpretation of the same table is purely algebraic.
	In fact the variation although monotonic is not that regular and the
	last triangle corresponds to 31.84 degrees.\\
	Freiberg\\
	The tablet, dated 1900 to 1600 B.C., gives, with 4 errors, and in
	hexadesimal notation 15 values of\\
	\hth$	a, b,$ and $(\frac{a}{c})^2 = sec^2(B)$\\
	where B is the angle opposite b,\\
	from $2 49    1 59   [1 59] 15$         or $169  119  \frac{7155}{3600}
	$\\
	to   $  53      56   [1] 23 13 46 4_0$ or $\frac{17977600}{1296000}.$\\
	Where the values between brackets are reconstructed values and 56 should
	be corrected to 28.

	\ssec{On Plimpton 322, Pythagorean numbers in Babylonean
	Mathematics.}\label{sec-plimpton}

	The tablet gives in hexadesimal notation columns I, II, III and IV,
	except for the line labelled 11a in column IV.
	diff. is the difference between the numbers in column IV.
	The numbers in the second line give, in hexadesimal notation 
	$\frac{u}{v}$ and
	$\frac{v}{u},$ for instance $2;24 = 2+\frac{24}{60} = \frac{12}{5}.$\\
	$\begin{array}{rrrrrrrrl}
	IV&&&     III&&          II &              I\\
	&v& u& a& c& b&B& (\frac{a}{c})^2&diff.\\
	1& 5&12& 169& 120& 119& 44.7603& 1.9834\\
	&2;&24,&0;&25,\\
	2&27&64&4825&3456&3367& 44.2527& 1.9492& -.034244\\
	&2;&22,13,20,&0;&25,18,45,\\
	3&32&75&6649&4800&4601& 43.7873& 1.9188& -.030356\\
	&2;&20,37,30,&0;&25,36,\\
	4&54& 125& 18541& 13500& 12709& 43.2713& 1.8862& -.032554\\
	&2;&18,53,20,&0;&25,55,12,\\
	5& 4& 9&97&72&65& 42.0750& 1.8150& -.071240\\
	&2;&15,&0;&26,40,\\
	6& 9&20& 481& 360& 319& 41.5445& 1.7852& -.029815\\
	&2;&13,20,&0;&27,\\
	7&25&54&3541&2700&2291& 40.3152& 1.7200& -.065209\\
	&2;& 9,36,&0;&27,46,40,\\
	8&15&32&1249& 960& 799& 39.7703& 1.6927& -.027274\\
	&2;&8,&0;&28, 7,30,\\
	9&12&25& 769& 600& 481& 38.7180& 1.6427& -.050040\\
	&2;&5,&0;&28,48,\\
	 10&40&81&8161&6480&4961& 37.4372& 1.5861& -.056547\\
	&2;& 1,30,&0;&29,37,46,40,\\
	 11& 1& 2& 5& 4& 3& 36.8699& 1.5625& -.023623\\
	&2;&&0;&30,\\
	 11a& 64& 125& 19721& 16000& 11529& 35.7751& 1.5192& -.043290\\
	&1;&57,11,15,&0;&30,43,12,\\
	 12&25&48&2929&2400&1679& 34.9760& 1.4894& -.029793\\
	&1;&55,12,&0;&31,15,\\
	 13& 8&15& 289& 240& 161& 33.8550& 1.4500& -.039399\\
	&1;&52,30,&0;&32,\\
	 14&27&50&3229&2700&1771& 33.2619& 1.4302& -.019779\\
	&1;&51,6,40,&0;&32,24,\\
	 15& 5& 9& 106&90&56& 31.8908& 1.3872& -.043078 *\\
	&1;&48,&0;&33,20,
	\end{array}$

	There are 2 interpretations for the method of obtaining this table.
	The method of Neugebauer and Sachs, assumes the knowledge of
	the formulae\\
	\hth$	a = u^2 + v^2,$ $b = u^2 - v^2,$ $c = 2 u v.$\\
	It was proven later that all integer solutions of $a^2 = b^2 + c^2,$ 
	can be obtained from these formulae and that the values of $a,$ $b$ and
	$c$ are relatively prime if $u$ and $v$ are relatively prime and not
	both odd. They observe that $u$ and $v$ are always regular, it is,
	have only 2, 3 and
	5 as divisors, this implies that the reciprocals have a finite
	representation if we use hexadesimal notation.\\
	This point of view is confirmed if we observe that $u$ and $v$ are
	precisely all the regular numbers, which are relatively prime,
	satisfying
	\enumb
	\item $(\sqrt{2}-1) u < v < u <= 125,$
	\enume
	except for the added pair, 11a, $u = 125,$ $v = 64.$ The first condition
	corresponds to requiring that the triangle has an angle $B$ opposite $b$
	less than 45 degrees.  In this range, only one pair is such
	that $u$ and $v$ are both odd.  This is the pair $u = 9,$ v = 5, which
	gives $a = 106,$ $b = 56,$ $c = 90.$  The values $a = 53,$ $b = 45,$
	$c = 28,$ could have
	been obtained with $u = 7$ and $v = 2,$ but these numbers are not both
	regular.  It is interesting that one of the errors occurs for this pair,
	$a$ being divided by 2 but not $b.$\\
	The other point of view is presented by Bruins which claims that
	$a$ and $b$ are obtained from a subset of tables of reciprocals, which
	we could write $\frac{u}{v}$ and $\frac{v}{u},$ giving the values of $a$
	and $b,$ because of\\
	\hth$	(\frac{u}{v} + \frac{v}{u})^2
		= (\frac{u}{v} - \frac{v}{u})^2 + 2^2,$\\
	after removing the common factors, which are necessarily 2, 3 or 5.
	This would give the table for monotonically varying values of
	$\frac{a}{c}.$\\
	We have given the corresponding hexadesimal values of $\frac{u}{v}$ and
	$\frac{v}{u}$ on alternate lines.\\
	Condition 0. adds credibility to the point of view of Neugebauer
	and would strengthen the geometrical content of the table.
	A hope to get a deciding clue from one of the errors in the table is
	not easily fulfilled.  Indeed the second line gives for $a$ and $b,$
	11521 and 3367, instead of 4825 and 3367.\\
	One explanation, which I consider farfetched, is given by Gillings.
	He assumes that 11521 is obtained using $(64+27)^2 + 2*27*60.$
	This requires several errors, first to add before squaring, then to
	add $2*27*60,$ which is explained by Gillings by the use of\\
	\hth$	u^2 + v^2 = (u+v)^2 - 2uv$\\
	with $-$ replaced by $+$ and $v = 64$ replaced by v = $60.$
	An other explanation, only slightly less farfetched is to observe
	that, if we use Bruins approach, both numbers $2;22,13,20$ and
	$0;25,18,45$
	have to be divided 3 times by 5, (or multiplied by 12 in hexadesimal
	notation).  This gives for $a,$ $1,20,25$ in base $60.$  If we assume
	that the
	scribe wrote instead $1,20,,25,$ using a large space, rather than a
	small one, and multiplies by 12 twice more, we get 3,12,1.\\
	An other explanation could start by explaining why the scribe computed
	instead of\\
	\hth$ (64 (= 60+ 4) )^2 + (27 (= 24+3) )^2 = 4825,$\\
	\hth$(100 (= 60+40) )^2 + (39 (= 36+3) )^2 = 11521.$\\
	The argument could be decided if other tablets which continue this
	table are found.  The table Plimpt.tab, gives the values for angles
	less than 31.5 degrees, using criteria 0.\\
	There is an other minor controversy in the literature concerning
	the fact that the 1 in column IV is visible or not in the tablet.
	If the opinion is taken, which is contrary to Neugebauer, that 1
	is not there, column I is then $(\frac{b}{c})^2
	= tan^2(angle\:opposite\:b)
	= ( \frac{1}{2}(\frac{u}{v} - \frac{v}{u}))^2$ instead of
	$(\frac{a}{c})^2 = (\frac{1}{2}(\frac{u}{v} + \frac{v}{u}) )^2.$ 



	\begin{center}
	{\Huge CHAPTER I\\[40pt]
	FINITE PROJECTIVE\\[20pt] GEOMETRY}\\[40pt]
	\end{center}

\setcounter{section}{89}
	\section{Answers to problems and miscellaneous notes.}
	\ssec{Algebra and modular arithmetic.}
	\sssec{Example.}
	Modulo 7, the inverses of 1 through 6 are respectively\\
	\hth$	1,4,5,2,3,6.$

	\sssec{Answer to}
	\vspace{-18pt}\hspace{94pt}{\bf \ref{sec-}.}\\

Notes for section on axiomatic, Pieri (coxeter, p. 12), Menger (Coxeter, p. 14)
Dedekind (Coxeter, p 22), Enriques (Coxeter, p.22)

	The following does not work, leave for
	examination of other types, 1 where the triangles have sides through
	$Q_1$ and $Q_2$ may give something, see also Pickert p. 74,75,80

	\ssec{Linear Associative Planes.}
	\sssec{Axiom. [2-point Desargues]}
	The {\em 2-point Desargues axiom} is the special case when we restrict
	Desargues' axiom to the case when the center $C$ of the configuration
	is one of 2 given points $Q_1$ or $Q_2$ of the given axis $c$.
	More specifically,\marginpar{?}
	$C\incid c$, and for the 2 triangles $\{A_i\}$ and $\{B_i\}$,\\
	let $C_i := (A_{i+1}\times A_{i-1})\times(B_{i+1}\times B_{i-1}),$\\
	$c_i := (A_i\times B_i),$ $c_i\incid C$, $i = 0,1,2,$
	incidence($A_0\times A_j,B_0\times B_j,c)$, $j = 1,2$,\\
	$\implies$ incidence($A_1\times A_2,B_1\times B_2,c)$. We write\\
	\hth 2-point-Desargues$(C,\{A_i\},\{B_i\};\langle C_i\rangle,c).$

	\sssec{Theorem.}
	{\em Given 2 triangles $\{A_i\}$ and $\{B_i\}$, let
	$C_i := (A_{i+1}\times A_{i-1})\times(B_{i+1}\times B_{i-1}),$
	$C_i := A_i\times B_i,$ and $C := c_1\times c_2,$\\
	\hth$\langle C_i,c\rangle$ and $C\incid c\implies c_0\incid C.$
	We write}\\
	\hth 2-point-Desargues$^{-1}(c,\{A_i\},\{B_i\};\langle \un{c_0},c_1,c_2
	\rangle,C)$

	Proof: 2-point-Desargues$(C_0,\{A_1,B_1,C_2\},\{A_2,B_1,C_1\};\langle
	 B_0,A_0,C\rangle,c)$.

	\sssec{Definition.}
	A {\em linear associative plane} is a perspective plane for which
	the 2-point Desargues axiom is satisfied for 2 specific points on a
	specific line of the plane.

	If the line is $q_2$ and the points are $Q_1$ and $Q_2$, we have

	\sssec{Theorem.}\label{sec-tveblenw}
	{\em In a linear associative plane, the ternary ring ($\Sigma$,*) is a
	\ldots, more specifically:}
	\enumb
	\item$(\Sigma,*)$ is linear, $a * b * c = a \cdot b \:+\:c$,
	\item($\Sigma$,+) {\em is a group,}
	\item($\Sigma-\{0\},\cdot)$ {\em is a loop,}
	\item ($\Sigma,*) = (\Sigma,+,\cdot)$ {\em is right distributive,}
	$(a+b)\cdot c = a \cdot c + b\cdot c$.
	\item$a \neq b \implies x \cdot a = x \cdot b + c$
	{\em has a unique solution.}
	\enume

 before 
	\ssec{Veblen-Wedderburn Planes.}


\chapter{FINITE PROJECTIVE GEOMETRY}
{\tiny\footnotetext[1]{G20.TEX [MPAP], \today}}


\setcounter{section}{-1}
	\section{Introduction.}

	In Section 1, I give the axiomatic definition of synthetic projective
	geometry. In Section
	2, I give an algebraic model of projective geometry. Although I will
	use, whenever possible a synthetic proof, I will use extensively an
	algebraic proof to proceed more expeditiously, if not more elegantly.
	The reader is encouraged to replace these by the more satisfying
	synthetic proofs.  In Section 3, I discuss the geometric model of the
	projective plane of order 2, 3 and 5, discovered by Fernand Lemay and
	relate each model to classical configurations.

	\section{Synthetic Finite Projective Geometry.}

\setcounter{subsection}{-1}
	\ssec{Introduction.}
	Projective Geometry implies usually that when we write down the
	equivalent algebraic axioms, the underlying field is the field of reals.
	Most of the properties that I will discuss in this Chapter and in the
	next one are valid whatever the field chosen.  To deal with a set of
	Axioms which characterize the plane, in a simpler setting, I will assume
	instead that the field is finite.  See \ref{sec-ax2}.  Most properties
	generalize to any field.

	\ssec{Notation.}
	The objects or elements of plane projective geometry are points and
	lines.  The relation between points and lines is called incidence.
	A point and a line are incident if and only if the point is on the line
	or if the line passes through the point.\\
	Identifiers are sequences of letters and digits, starting with a
	letter.  If the first letter is a lower case letter, the identifier
	will denote a line.  If the first letter is an upper case letter, the
	identifier will denote a point.  If the line $ab$ is constructed as the
	line through the points $A$ and $B,$ we write\\
	\hth$		ab := A \times B.$\\
	If the point $A_0$ is constructed as the point on both $a_1$ and $a_2,$
	we write\\
	\hth$		A_0 := a_1 \times a_2.$\\
	The symbol "~:=~" pronounced "is defined as" indicates a definition
	of a new point or of a new line.  The symbol "$\times$" will be
	justified in \ref{sec-ntimes}.\\
	\hth$		A \cdot ab = 0,$ or $A\incid ab$,\\
	is an abbreviation for the statement "the point $A$ is on the line
	$ab$ ".\\
	\hth$		A \cdot ab\neq 0$ or $A\notinc ab$,\\
	is an abbreviation for the statement "the point $A$ is not on the
	line $ab$".\\
	\hth$		A = B,$ $x = y,$\\
	are abbreviations for "the points $A$ and $B$ or the lines $x$ and
	$y$", all previously defined, "are identical".\\
	\hth$\{A, B, C\}$ or $\{a, b, c\}$\\
	denotes a triangle with vertices
	$A,$ $B$ and $C$ or sides $a,$ $b$ and $c.$

	For Projective Geometry over fields we will use the following Axioms.

	\ssec{Axioms.}\label{sec-ax1}
	\hth{\em Of incidence and existence or of allignment}:
	\enumb
	\item Given 2 distinct points, there exists one and only one line
		incident to, or passing through, the 2 points.
	\item Given 2 distinct lines, there exists one and only one point
		incident to, or on, the 2 lines.
	\item There exists at least 4 points, any 3 of which are not collinear.

	\noindent\hth{\em Of Pappus:}
	\item Let $A_0,$ $A_1,$ $A_2$ be distinct points on $a,$\\
	 	let $B_0,$ $B_1,$ $B_2$ be distinct points on $b.$\\
	 Let $C_0$ be the intersection of $A_1 \times B_2$ and $A_2 \times 
	B_1$ or\\
	\hth$C_0 := (A_1 \times B_2) \times (A_2 \times B_1).$\\
	Similarly, let\\
	\hth$C_1 := (A_2 \times B_0) \times (A_0 \times B_2),$
	$C_2 := (A_0 \times B_1) \times (A_1 \times B_0),$\\
	then the points $C_0,$ $C_1,$ $C_2$ are collinear.  (Fig. 1a)
	\enume

	\sssec{Notation.}\label{sec-npappus}
	The subscript $i$ is usually restricted to the set \{0,1,2\} and 
	addition is then done modulo 3. I write\\
	\hth Pappus$(\langle A_i\rangle,\langle B_i\rangle;\langle C_i\rangle)$
	or more generally\\
	\hth Pappus$(\langle A_i\rangle[,a],\langle B_i\rangle[,b];
	\langle C_i\rangle[,c][,X])$.\\
	where "$\langle X_i\rangle$" indicate that the points $X_i$ are
	collinear,
	where the brackets indicate that what is between them need not be given,
	and where $X$, if written, is the intersection of $a$ and $b$.

	The axiom is trivially satisfied if $X$ is one of the points $A_i$ or 
	$B_i$.  If the axiom is used in proofs, it is always assumed that
	the points $A_i$ and $B_i$ are distinct from $X$.

	Any plane satisfying the allignment axioms and the axiom of Pappus
	is called a {\em Pappian plane}.
	For Projective Geometry over a specific field we will add one axiom
	or a set of associated axioms, for instance,
	for finite Projective Geometry over a ${\Bbb Z}_p^k$, we add

	\ssec{Axiom (the finite field).}\label{sec-ax2}
	\hth On the line $l$ there are exactly $p^k+1$ points, $p$ a prime.

	\sssec{Exercise.}
	Write down the appropriate existence axiom associated with
	the fields,
	\enumb
	\item${\Bbb R}$, classical Projective Geometry ,
	\item${\Bbb C}$, complex Projective Geometry,
	\item${\Bbb Q}$, rational Projective Geometry.
	\enume

	\ssec{Basic consequences.}

	\sssec{Theorem.}\label{sec-tinc}
	\enumb
	\item{\em Each line is incident to exactly $p^k+1$ points.}
	\item{\em Each point is incident to exactly $p^k+1$ lines.}
	\item{\em There are exactly $p^{2k} + p^k + 1$ points and lines.}
	\enume

	The proof is left as an exercise.

	\sssec{Corollary.}
	{\em There exists at least 4 lines, any 3 of which are not incident.}

	\sssec{Comment.}
	If, contrary to \ref{sec-ax1}.2, there is only one point $P$ not
	on the line $l,$ the geometry reduces to $l,$  to a \{ pencil \} of
	$p + 1$
	lines through $P,$ to $P$ and to a set of $p + 1$ points on $l.$
	The axiom of Pappus is satisfied vacuously because no 2 distinct
	lines contain 3 points each.

	\sssec{Definition.}
	The line through $C_0,$ $C_1$ and $C_2$, in the axiom of Pappus,
	is called the {\em Pappus line.}

	\sssec{Notation.}\label{sec-nsynthetic}
	I introduce in the next Chapter 
	a detailed notation for
	algebraic projective geometry.  An incomplete notation for the
	synthetic approach will now be introduced.  The purpose is to
	formalize the Theorems, without the details of the approach of
	Russell and Whitehead.\\
	\hth$\langle X_i\rangle$ or $(\langle X_i\rangle,x)$ indicates that the
	points $X_i$ are collinear and distinct, on $x$,\\
	\hth$\langle x_i\rangle$ or $(\langle x_i\rangle,X)$ indicates that the
	lines $x_i$ are incident and distinct, through $X$,\\
	\hth$\{X_i\}$ indicates that the points $X_i$ are distinct and
	not collinear, in other words form a triangle and similarly for the 
	sides, $\{x_i\}$.\\
	\hth$incidence(A,B,C[,l])$ or $incidence(A_j[,l])$, $j \in \{0,1,\ldots,
	k\},$ $k \leq 2$,\\
	is used to state that the points $A,$ $B,$ $C$ or the points $A_j$ are
	on the same line $l$. "$[,l]$" indicates that the name of the line need
	not be given explicitely.\\
	\hth$incidence(a,b,c[,L])$ or $incidence(a_j[,L])$\\
	is the corresponding statement for lines $a,b,c$ or $a_j$ incident to
	the point $L$.

	\noindent No.\hti{6}Pappus$(\langle A_i\rangle,\langle B_i\rangle;
	\langle C_i\rangle)$ and the corresponding axioms can be written,\\
	\hth in greater detail, as follows.\\
	Hy0.\hti{5}$\langle A_i\rangle.$\\
	Hy1.\hti{5}$\langle B_i\rangle.$\\
	De.\hti{6}$C_i := (A_{i+1} \times B_{i-1})
		\times (A_{i-1} \times B_{i+1}).$\\
	Co.\hti{6}$\langle C_i\rangle.$

	"No" is an abbreviation for "nomenclature" or "notation", "Hy",
	for "hypothesis", "De", for "Definition", "Co" for "conclusion".

	Notice that the order of the points is important.\\
	The reciprocal,\\
	\hth$Pappus^{-1}(\langle A_i\rangle,\langle C_i\rangle;\langle B_i
	\rangle)$\\
	exchanges Hy1. and Co. and follows from\\
	\hth Pappus$(\langle A_i\rangle,\langle C_i\rangle;\langle B_i\rangle)$.

	In a statement, different letters indicate different elements
	with no special relationship between them except as stated in the
	hypothesises "Hy".

	\sssec{Theorem.}
	Pappus$(\langle A_i\rangle,\langle B_i\rangle;\langle C_i\rangle)
	\implies$ Pappus$(\langle A_0,B_1,C_2\rangle,\langle B_0,C_1,A_2\rangle;
	\langle C_0,A_1,B_2\rangle).$

	\ssec{The Theorem of Desargues.}

	\sssec{Theorem. [Desargues]}\label{sec-tdes}
	Hy0.\hti{6}$\{A_i\},\:\{B_i\},$\\
	De0.\hti{6}$c_i := A_i\times B_i,$\\
	Hy1.\hti{6}$C\incid c_i,$\\
	De1.\hti{6}$a_i := A_{i+1}\times A_{i-1},$\\
	De2.\hti{6}$b_i := B_{i+1}\times B_{i-1},$\\
	De3.\hti{6}$C_i := a_i \times b_i,$\\
	Co.\hti{6}$(\langle C_i\rangle,c).$\\
	No.\hti{6}$Desargues(C,\{A_i\}[,\{a_i\}],\{B_i\}[,\{b_i\}];
	\langle C_i\rangle[,\langle c_i\rangle][,c])$.

	This is the notation for the following statements.\\
	\noindent{\em Given two triangles $\{A_0, A_1, A_2\}$ and $\{B_0, B_1,
	 B_2\},$ such that the lines $A_0 \times B_0,$ $A_1 \times B_1$ and
	$A_2 \times B_2$ have a point $C$ in common.  Let\\
	\hth$	C_0 := (A_1 \times A_2) \times (B_1 \times B_2),$
	$C_1 := (A_2 \times A_0) \times (B_2 \times B_0),$\\
	\hth$	C_2 := (A_0 \times A_1) \times (B_0 \times B_1).$\\
	Then $C_0,$ $C_1,$ $C_2$ are incident to the same line $c$}  (Fig. 3a).
	It is assumed that the triangles are distinct and that the lines $c_i$
	are distinct.

	This theorem can be proven using the incidence axioms in 3 dimensions.
	In 2 dimensions, it can be taken as an axiom or it can be derived
	from the axiom of Pappus, see \ref{sec-pdes}.  But the axiom of Pappus
	does not derive from the incidence axioms and the Theorem of Desargues
	taken as axiom.

	\sssec{Theorem.}
	{\em The axiom of incidence and the axiom of Pappus \ref{sec-ax1}.4.
	imply the Theorem of Desargues.}  See \ref{sec-pdes}.

	\ssec{Configurations.}

	\setcounter{subsubsection}{-1}
	\sssec{Introduction.}
	One of the characteristics of synthetic geometry is to
	start from a set of points and lines, to construct from them new points
	and lines and to extract known sets which have known properties.
	Hence, it is useful to describe some of the important sets, which
	are called {\em configurations}.  We have seen 2 such configurations.
	In that of Pappus, we have 9 points and 9 lines.
	In that of Desargues, we have 10 points and 10 lines.
	I will define here the complete quadrangle and the complete 
	quadrilateral configuration, the special Desargues configuration,
	as well as closely related  configurations.
	To characterize the configuration further, I will use the
	following notation:

	\sssec{Notation.}
	\hth$	10 * 3  \:\&\:  10 * 3, (11)$\\
	indicates that each of the 10 points are incident to 3 lines,
	that each of the 10 lines are incident to 3 points and that
	the construction requires 11 independent data elements
	(2 for a given point or line, 1 for a point on a given line or a
	line through a given point). Or\\
	\hth$	3 * 6 + 8 * 3  \:\&\:  12 * 3 + 3 * 2,$\\
	indicates that 3 points are incident to 6 lines, that 8 points are
	incident to 3 lines and that 12 lines are incident to 3 points and that
	3 lines are incident to 2 points.
	The order chosen is that of decreasing number of incident elements.

	The notation does not uniquely define the configuration but is a
	useful tool.

	\sssec{Definition.}
	A {\em confined configuration} is a configuration in the description
	of which "~$*$~2~" does not occur.\\
	Except for the triangle and the complete quadrangle or quadrilateral,
	I will restrict the word configuration to confined configuration and
	will use the adjective "non confined" otherwize.\\
	A {\em self dual type configuration} is one for which the information
	to the left of "~\&~" is the same as that to the right.

	It should not be confused with the notion of self dual configuration
	that will be introduced later.  A self dual configuration is a
	self dual type configuration but not vice-versa.

	\sssec{Theorem.}\label{sec-tpph}
	{\em The configuration of Pappus is of type $9 * 3  \:\&\:  9 * 3, (10).
	$  It can be viewed as a degenerate case of that of Pascal.  See
	\ref{sec-tpascal}. Hence the
	alternate name {\em Pappus-Pascal hexagon}:  If the alternate points of
	the hexagon $A_0,$ $A_1,$ $A_2,$ $A_3,$ $A_4,$ $A_5$ are on 2 lines, the
	three pairs of
	opposites sides of the hexagon meet in 3 collinear points $P_0,$ $P_1$
	and $P_2.$}

	The correspondence between this notation and that used in the
	Theorem of Pappus is:\\
	$\begin{array}{cccccccccc}
	\hth&A_0,&A_1,&A_2,&A_3,&A_4,&A_5,&P_0,&P_1,&P_2,\\
	    &B_2,&A_1,&B_0,&A_2,&B_1,&A_0,&C_0,&C_2,&C_1.
	\end{array}$

	\sssec{Theorem.}
	{\em The configuration of Desargues is of type $10 * 3  \:\&\:  10 * 
	3, (11).$\\
	It can also be viewed as consisting of 2 pentagons which are inscribed
	one into the other.  The points $P_0,$ $P_1,$ $P_2,$ $P_3,$ $P_4$ and
	the points $Q_0,$ $Q_1,$ $Q_2,$ $Q_3,$ $Q_4$ being such that $P_0$ is on
	$Q_0 \times Q_1,$
	$P_1$ is on $Q_1 \times Q_2,$
	$P_2$ is on $Q_2 \times Q_3,$ $P_3$ is on $Q_3 \times Q_4$ and $P_4$
	is on $Q_4 \times Q_0.$
	$Q_0$ is on $P_1 \times P_3,$ $Q_1$ is on $P_2 \times P_4,$ $Q_2$ is
	on $P_3 \times P_0,$
	$Q_3$ is on $P_4 \times P_1$ and $Q_4$ is on $P_0 \times P_2.$}

	The correspondence between this notation and that used in the Theorem
	of Desargues is:\\
	$\begin{array}{ccccccccccc}
	\hth&	P_0,&P_1,&P_2,&P_3,&P_4,&Q_0,&Q_1,&Q_2,&Q_3,&Q_4,\\
	&	B_1,&A_0,&B_2,&A_1,&C_1,&C_2,&B_0,&C,& A_2,&C_0.
	\end{array}$

	\sssec{Definition.}\label{sec-dcq}
	A {\em complete quadrangle} is a configuration consisting of
	4 points $A_0,$ $A_1,$ $A_2,$ $A_3,$ no 3 of which are on the same line
	and of the 6 lines through each pair of points:  
	$a_0 := A_0 \times A_1,$ $a_1 := A_0 \times A_2,$
	$a_2 := A_0 \times A_3,$
	$a_3 := A_2 \times A_3,$ $a_4 := A_3 \times A_1,$
	$a_5 := A_1 \times A_2.$ (Fig. 2a)\\
	It is of type\\
	\hth$4 * 3  \:\&\:  6 * 2, (8).$

	\sssec{Definition.}\label{sec-ddiagpoint}
	The 3 points $D_0 := a_0 \times a_3,$ $D_1 := a_1 \times a_4$ and
	$D_2 := a_2 \times a_5$ are called the {\em diagonal points} of the
	complete quadrangle.\\
	The lines $d_i$ joining the diagonal points are called {\em diagonal
	lines}.\\
	These form, together with the quadrangle configuration, the {\em
	completed non confined quadrangle configuration.} See Fig. 2a'.

	\sssec{Definition.}\label{sec-dconic2}
	Given a complete quadrangle, a {\em conic2 pseudo non confined
	configuration} is the
	sub configuration consisting of 3 of the points and the 3 lines
	joining these points to the 4-th one.  It is of type \\
	\hth$1 * 3 + 3 * 1 \:\&\: 1 * 3 + 3 * 1.$ (8)\\
	See \ref{sec-tconic2}.

	\sssec{Definition.}\label{sec-dcompleted}
	Given a complete quadrangle, a {\em completed quadrangle configuration}
	is the configuration consisting of the complete quadrangle,
	the diagonal points and the lines joining the diagonal points.

	\sssec{Theorem.}
	\enumb
	\item {\em If $p = 2$ the completed quadrangle configuration is of type}
\\
	\hti{16}$7 * 3  \:\&\: 7 * 3, (8).$\\
	\hth See \ref{sec-tproj2} and \ref{sec-tconic2}
	\item {\em If $p > 2,$ it is of type\\
	\hti{16}$	3 * 4 + 4 * 3  \:\&\:  6 * 3 + 3 * 2$ (8)\\
	and is not confined.}
	\enume

	\sssec{Definition.}
	A {\em complete $n$-angle} is a configuration consisting of
	$n$ points no 3 of which are on the same line and of the
	$\frac{n(n-1)}{2}$ lines through each pair of points.

	\sssec{Theorem.}
	{\em A complete 5-angle does not exist if $p < 5.$  Indeed, on the
	line through 2 of the points, we must have 3 other points which are the
	intersection with the 3 pairs of lines through the other 3 points.
	We must have therefore at least 5 points on each line.}

	\sssec{Exercise.}
	For which value of $p$ does a complete $n$-angle exist for
	$n > 5?$

	\sssec{Definition.}\label{sec-dcql}
	A {\em complete quadrilateral} is a configuration consisting of
	4 lines $a_0,$ $a_1,$ $a_2,$ $a_3,$ no 3 of which are incident to the
	same point and of the 6 points through each pair of lines:  
	$A_0 := a_0 \times a_1,$ $A_1 := a_0 \times a_2,$
	$A_2 := a_0 \times a_3,$
	$A_3 := a_2 \times a_3,$ $A_4 := a_3 \times a_1,$
	$A_5 := a_1 \times a_2.$ (Fig. 2b)\\
	It is of type\\
	\hth$6 * 2  \:\&\:  4 * 3, (8).$

	\sssec{Definition.}\label{sec-ddiagline}
	The 3 lines $A_0 \times A_3,$ $A_1 \times A_4$ and $A_2 \times A_5$
	are called the {\em diagonal lines} of the complete quadrilateral.\\
	The points joining the diagonal lines are called {\em diagonal points}.
	These together with the complete quadrilateral configuration form
	the {\em completed quadrilateral non confined configuration} (Fig. 2b').

	\sssec{Definition.}\label{sec-dsd}
	The {\em special Desargues configuration}, consists of 13 points
	and 13 lines obtained as follows.  $A_0,$ $A_1,$ $A_2,$ $C$ is a
	complete quadrilateral,\\
	\hth$	a_0 := A_1 \times A_2,$ $a_1 := A_2 \times A_0,$
	$a_2 := A_0 \times A_1,$\\
	\hth$	c_0 := C\: \times A_0,$ $c_1 := C\: \times A_1,$
	$c_2 := C\: \times A_2,$\\
	\hth$	B_0 := a_0 \times c_0,$ $B_1 := a_1 \times c_1,$
	$B_2 := a_2 \times c_2,$\\
	\hth$	b_0 := B_1 \times B_2,$ $b_1 := B_2 \times B_0,$
	$b_2 := B_0 \times B_1,$\\
	\hth$	C_0 := a_0 \times b_0,$ $C_1 := a_1 \times b_1,$
	$C_2 := a_2 \times b_2,$\\
	\hth$	r_0 := A_0 \times C_0,$ $r_1 := A_1 \times C_1,$
	$r_2 := A_2 \times C_2,$\\
	\hth$	R_0 := r_1 \times r_2,$ $R_1 := r_2 \times r_0,$
	$R_2 := r_0 \times r_1,$\\
	\hth$	c\: := C_1 \times C_2.$  (Fig. 3e')

	This configuration is also called the {\em quadrangle-quadrilateral
	configuration}.  The quadrangle is \{$R_i,$ $C$\} or \{$c_i,r_i$\},
	the quadrilateral is \{$b_i,$ $c$\} or \{$C_i,$ $B_i.$\}
	The diagonal points are $A_i$ and the diagonal lines, $a_i.$ 

	\sssec{Comment.}
	The dual construction can be obtained with the upper case letters
	exchanged for the lower case ones except for the exchange of $B_i$ and
	$r_i$ and $b_i$ and $R_i$.\\
	This configuration plays an essential role in Euclidean Geometry.  An
	example consist of a triangle $\{A_i\}$, $C$ the barycenter, $a_i$,
	the sides, $c_i,$ the medians, $B_i$, the mid-points, $b_i$, the sides
	of the complementary triangle, $C_i$, the directions of the sides,
	$r_i$, the sides of the anticomplementary triangle, $R_i$, its vertices,
	c, the ideal line.
 
	\sssec{Definition.}\label{sec-dconic3}
	Given a complete quadrangle-quadrilateral configuration, a {\em conic3
	pseudo non confined configuration} is the sub configuration consisting
	of the quadrangle\\
	\{$R_i,$ $C$\} and the quadrilateral \{$b_i,$ $c$\}. It is of type\\
	\hth$1 * 3 + 3 * 1 \:\&\: 1 * 3 + 3 * 1$ (8).\\
	See \ref{sec-tqqc3}.

	\sssec{Theorem. [Special Desargues]}\label{sec-tsd}
	\enumb
	\item {\em $C_0$ is on $c,$}
	\item {\em $R_0$ is on $c_0,$ $R_1$ is on $c_1,$ $R_2$ is on $c_2.$}
	\item {\em If $p = 3$ the special Desargues configuration is of type}\\
	\hti{16}$13 * 4  \:\&\: 13 * 4$ (8)\\
	\hth See \ref{sec-eqqc}\\
	{\em If $p > 3,$ it is of type}\\
	\hti{16}$	9 * 4 + 4 * 3  \:\&\:  9 * 4 + 4 * 3, (8).$
	\item {\em If we exclude $r_0,$ $r_1,$ $r_2,$ $R_0,$ $R_1,$ $R_2,$
	we obtain a special case of the Desargues configuration in which\\
	\hth$P_0$ is on $A_1 \times A_2,$ $P_1$ is on $A_2 \times A_0$ and
	$P_2$ is on $A_0 \times A_1.$}
	\enume

	The proof will be given in section \ref{sec-psd}.

	\sssec{Definition.}
	$c$ is called the {\em polar of $C$ with respect to the triangle}
	$\{A_0,$ $A_1,$ $A_2\}.$
	$C$ is called the {\em pole of $c$ with respect to the triangle}.

	\sssec{Notation.}\label{sec-nsd}
	Part of Definition \ref{sec-dsd} and Theorem \ref{sec-tsd} can be noted
	as follows.\\
	No.\hti{6}$Special\:Desargues(C,A_i;C_i,c).$\\
	De0.\hti{5}$a_i := A_{i+1} \times A_{i-1}.$\\
	De1.\hti{5}$B_i := a_i \times (C \times A_i).$\\
	De2.\hti{5}$C_i := a_i \times (B_{i+1} \times B_{i-1}).$\\
	Co.\hti{6}$(\langle C_i\rangle,c).$

	\sssec{Exercise.}\label{sec-eqqc}
	Construct the configuration starting from $R_i,$ $C,$ and prove
	the 4 incidence properties corresponding to \ref{sec-tsd} in this
	construction.

	\sssec{Exercise.}\label{sec-eqqc3}
	For $p = 3,$ prove that $B_i$ is on $r_i$ and $C$ is on $c.$
	See also \ref{sec-tqqc3}.

	For a connection between conics and the quadrangle-quadrilateral
	configuration, when $p = 3,$ see \ref{sec-tqqc3}.

	\ssec{Other Configurations.}

	\setcounter{subsubsection}{-1}
	\sssec{Introduction.}
	There exist 2 other configurations of type $9 * 3 \:\&\: 9 * 3,$
	these will be constructed and defined. Many special cases of Desargues 
	configurations will be defined, as well as the extended special 
	Desargues configuration and the dodecahedral configuration.  I end
	by making some comments on the complete triangle in the more
	general case of the perspective plane.

	\sssec{Definition.}
	H0.0.	A$_i$, M, d$_0$,\\
	H0.1.	A$_0 \incid $d$_0$,\\
	D0.0.	a$_i$ := A$_{i+1} \times $A$_{i-1}$,\\
	D1.0.	d$_1$ := M$ \times $A$_1$, d$_2$ := M$ \times $A$_2$,\\
	D1.1.	B$_i$ := d$_{i} \times $a$_i$,\\
	D1.2.	mm$_0$ := B$_1 \times $B$_2$, MA$_0$ := mm$_0 \times $a$_0$,\\
	D1.3.	eul := M$ \times $MA$_0$, C$_0$ := eul$ \times $d$_0$,\\
	D1.4.	c$_2$ := B$_2 \times $C$_0$, c$_1$ := B$_1 \times $C$_0$,\\
	D1.5.	C$_1$ := d$_1 \times $c$_2$, C$_2$ := d$_2 \times $c$_1$,
		c$_0$ := C$_1 \times $C$_2$,\\
	then\\
	C0.0.	B$_0 \incid $c$_0$, (Fig. 1c')

	This defines the {\em extended 2-Pappus Configuration}.

	\sssec{Definition.}
	The {\em 2-Pappus Pseudo Configuration} is the subset of the extended
	2-Pappus Configuration consisting of the point A$_i$, B$_i$, C$_i$ and
	of the lines a$_i$, b$_i$, c$_i$ (Fig. 1c).

	\sssec{Theorem.}
	{\em The extended 2-Pappus Configuration is of type}\\
	\hth$3 * 4 + 8 * 3 \:\&\: 3 * 4 + 8 * 3, (9).$\\
	{\em The 2-Pappus Pseudo Configuration is of type}\\
	\hth$9 * 3 \:\&\: 9 * 3, (9).$

	The proof is left as an exercise.

	\sssec{Definition.}
	H0.0.	A$_i$, M, c$_2$,\\
	H0.1.	A$_1 \incid $c$_2$,\\
	D0.0.	a$_i$ := A$_{i+1} \times $A$_{i-1}$,\\
	D1.0.	X$_0$ := c$_2 \times $a$_1$, \\
	D1.1.	ma$_1$ := M$ \times $A$_1$, ma$_2$ := M$ \times $A$_2$,
	      B$_1$ := ma$_{1} \times $a$_1$, B$_2$ := ma$_{2} \times $a$_2$,\\
	D1.2.	x$_0$ := X$_0 \times $B$_2$, X$_1$ := a$_0 \times $x$_0$,
		x$_1$ := X$_1 \times $M, C$_1$ := x$_1 \times $c$_2$,\\
	D1.3.	b$_0$ := B$_1 \times $B$_2$, C$_0$ := b$_0 \times $c$_2$,
		c$_0$ := A$_2 \times $C$_1$,\\
	D1.4.	c$_1$ := A$_0 \times $C$_0$, C$_2$ := c$_0 \times $c$_1$,\\
	D1.5.	b$_1$ := B$_2 \times $C$_1$, b$_2$ := B$_1 \times $C$_2$,
		B$_0$ := b$_1 \times $b$_2$,\\
	then\\
	C1.0.	B$_0 \incid $a$_00$, (Fig. 1d')

	This defines the {\em extended 1-Pappus Configuration}.

	\sssec{Definition.}
	The {\em 1-Pappus Pseudo Configuration} is the subset of the extended
	1-Pappus Configuration consisting of the point A$_i$, B$_i$, C$_i$ and
	of the lines a$_i$, b$_i$, c$_i$,  (Fig. 1d).

	\sssec{Theorem.}
	{\em The extended 1-Pappus Configuration is of type}\\
	\hth$1 * 5 + 4 * 4 + 7 * 3 \:\&\: 3 * 4 + 11 * 3, (9).$\\
	{\em The 1-Pappus Pseudo Configuration is of type}\\
	\hth$9 * 3 \:\&\: 9 * 3, (9).$

	\sssec{Definition.}
	There are many special cases of the Desargues configuration.
	\enumb
	\item 1-Desargues$(\{A_i\},\{\un{B_0},B_1,B_2\},\langle C_i\rangle),$
		in which $B_0\incid a_0$ (Fig. 3b).
	\item 2-Desargues$(\{A_0,\un{A_1},\un{A_2}\},\{B_i\},
		\langle C_i\rangle),$
		in which $A_1\incid b_1$ and $A_2\incid b_2$\\
		(Fig. 3c).
	\item 1-1-Desargues$(\{\un{A_0},A_1,A_2\},\{\un{B_0},B_1,B_2\},
		\langle C_i\rangle),$
		in which $A_0\incid b_0$ and $B_0\incid a_0$ (Fig. 3d).
	\item 3-Desargues$(\{A_i\},\{\un{B_i}\},\langle C_i\rangle),$
		in which $B_i\incid a_i$ (Fig. 3e).
	\item C-Desargues$(\{A_i\},\{B_i\},\langle \un{C_0},C_1,C_2\rangle),$
		in which $C_0\incid c_0$ (Fig. 3f).
	\item C-1-Desargues$(\{A_i\},\{B_i\},\langle \un{C_0},C_1,C_2\rangle),$
		in which $B_1\incid a_1$ and $C_2\incid c_2$\\
		(Fig. 3g).
	\item Elated-Desargues$(\un{C},\{A_i\},\{\un{B_0},B_1,B_2\},
		\langle C_i\rangle,c),$ in which $C\incid c$ (Fig. 3h).
	\enume
	In each case the additional incident point(s) is (are) underlined.

	\sssec{Definition.}\label{sec-desd}
	The {\em extended special Desargues} or {\em  extended
	quadrangle-quadrilateral configuration}, consists of 25 points and
	25 lines, those of \ref{sec-dsd} and\\
	\hth$	PQ_i := p_{i+1} \times q_{i-1},$
	$QP_i := q_{i+1} \times p_{i-1},$\\
	\hth$	QR_i := q_i \times r_i,	$
	$PR_i := p \times r_i,$\\
	\hth$	pq_i := P_{i+1} \times Q_{i-1},$
	$qp_i := Q_{i+1} \times P_{i-1},$\\
	\hth$	qr_i := Q_i \times R_i,	$ $pr_i := P \times R_i.$

	\sssec{Theorem.}\label{sec-tesd}
	{\em All 25 points are on the 6 lines $p_i,$ $r_i$ of the quadrangle
	$\{Q_i,P\}.$	 All 25 lines are on the 6 points $P_i,$ $R_i$ of the
	quadrilateral $\{q_i,p\}.$\\
	If $p = 5$ the extended special Desargues configuration is of type\\
	\hth$	10 * 6 + 15 * 4  \:\&\:  10 * 6 + 15 * 4.$\\
	If $p > 5$ it is of type\\
	\hth$	10 * 6 + 3 * 4 + 12 * 2  \:\&\:  10 * 6 + 3 * 4 + 12 * 2, (8).$}

	\sssec{Definition.}\label{sec-dconical}
	The {\em conical points} and {\em lines} of the extended
	quadrangle-quadrilateral configuration are the 6 points and 6 lines\\
	\hth$	AF_i := a_{i+1} \times pr_{i-1},$
	$FA_i := pr_{i+1} \times a_{i-1},$\\
	\hth$	af_i := A_{i+1} \times PR_{i-1},$
	$fa_i := PR_{i+1} \times A_{i-1},$

	\sssec{Theorem.}\label{sec-tesd1}
	\hth{\em $AF_i\cdot pq_{i+1} = FA_i\cdot qp_{i-1} = 0.$}

	Proof:  To show that $AF_0\cdot pq_1$ = 0, we can use the dual of
	Desargues' theorem applied to\\
	\hth$	\{p_0,p,p_1\} =    \{R_1,Q_2,R_0\}$\\
	and\\
	\hth$	\{a_1,p_2,r_1\} = \{Q_1,P_1,A_2\}$\\
	with axial points\\
	\hth$	A_0,R_2,A_1$ on the axis $a_2$\\
	and therefore central lines $QR_1,PQ_0,a_0$ on the center $AF_2.$ 

	\sssec{Comment.}\label{sec-cesd}
	We will see in \ref{sec-tconical} that the conical points
	are points on a conic.  The conic therefore appears in a
	natural way for $p = 5,$ in which case there are exactly 25 + 6 points
	and lines.  (The Pascal line of $N_1,M_2,N_0,M_1,N_2,M_0$ is
	$R_0,R_1,R_2.)$ 
	Although, in some sense, the conic exits already for $p = 2$ and
	$p = 3,$ see \ref{sec-dconic2}, \ref{sec-tqqc3}.

	\sssec{Definition.}
	In view of \ref{sec-tdodeca}, we define as the {\em dodecahedral
	configuration}, the configuration obtained by adding the 6 conical
	points to the extended special Desargues configuration.

	\sssec{Theorem.}
	{\em If $p = 5,$ the dodecahedral configuration is of type\\
	\hth$	25 * 6  \:\&\:  25 * 6.$\\
	If $p > 5,$ the dodecahedral configuration is of type\\
	\hth$13 * 6 + 3 * 4 + 12 * 3 + 3 * 2 \:\&\:
	13 * 6 + 3 * 4 + 12 * 3 + 3 * 2.$}

	Proof:  The first part follows from \ref{sec-tesd} and \ref{sec-tesd1}.
	For $p = 5,$ all the points and lines of the dodecahedral configuration
	are distinct and are all the points and lines of the corresponding
	finite projective geometry.  Any of the 6 conical points can be chosen
	to construct the extended special Desargues configuration.  Moreover,
	$pq_i$	contains also $PR_i,$ $QP_i$ and $FA_{i+1};$ $qp_i$ contains
	$PR_i,$ $PQ_i$ and $AF_{i-1};$ $qr_i$ contains $QP_i,$ $FA_i,$
	$QR_{i+1}$
	and $QR_{i-1};$ $pr_i$ contains $PQ_i,$ $QR_i,$ $FA_{i-1}$ and
	$AF_{i+1};$ $fa_i$ contains $QP_{i+1}$  $QR_i$ and $AF_i;$ $af_i$
	contains $PQ_{i-1},$ $QR_i$ and $FA_i.$	

	We leave, as an exercise, the proof of the following Theorem and the
	generalization of the definitions given therein.

	\sssec{Theorem.}
	{\em The dodecahedral configuration can be continued indefinitely.}

	Starting with $A_0$ = (1,0,0), $A_1$ = (0,1,0), $A_2$ = (0,0,1) and 
	$P = (1,1,1),$
	the coordinates of the points and lines obtained by replacing lower
	case letter by the corresponding upper case letter are the same, e.g.
	$p = [1,1,1].$\\
	These are\\
	\hth$	A_0 = (1,0,0),$ $R_0 = [0,1,-1],$ $P_0 = (0,1,1),$
	$Q_0 = (-1,1,1),$\\
	\hth$	PQ_0 = (-1,2,1),$       $QP_0 = (-1,1,2),$\\
	\hth$	QR_0 = (2,1,1),  $      $PR_0 = (-2,1,1),$\\
	\hth$	AF_0 = (2,0,-1),  $     $FA_0 = (2,-1,0),$\\
	More points are\\
	\hth$	PG_i := p_{i+1} \times g_{i-1},$
	$GP_i := g_{i+1} \times p_{i-1},$\\
	\hth$	AG_i := a_{i+1} \times g_{i-1},$
	$GA_i := g_{i+1} \times a_{i-1},$\\
	\hth$	QRQR_i := qr_{i+1} \times qr_{i-1},$
	$PQQP_i := PQ_i \times QP_i,$\\
	and the lines are defined similarly, e.g.\\
	\hth$	pg_i := P_{i+1} \times G_{i-1}.$\\
	We have\\
	\hth$	PG_0 = (-1,3,1),$       $GP_0 = (-1,1,3),$\\
	\hth$	AG_0 = (2,0,1),$        $GA_0 = (2,1,0),$\\
	\hth$	QRQR_0 = (-3,1,1),$    $PQQP_0 = (3,1,1).$\\
	and we have\\
	\hth$	GP_i\cdot qp_{i+1} = PG_i\cdot pq_{i-1} = 0,$\\
	\hth$	GA_i\cdot pq_i = AG_i\cdot qp_i = 0,$\\
	\hth$	QRQR_i\cdot r_i = 0.$\\
	Besides the conic\\
	\hth$	\{AF_i,FA_i\} = 2 (X_0^2 + X_1^2 + X_1^2)
		+ 5 (X_1 X_2 + X_2 X_0 + X_0 X_1) = 0,$\\
	there are many more, such as\\
	\hth$	\{PQ_i,QP_i\} = (X_0^2 + X_1^2 + X_1^2)
		+ 6 (X_1 X_2 + X_2 X_0 + X_0 X_1) = 0,$\\
	\hth$	\{PG_i,GP_i\} = (X_0^2 + X_1^2 + X_1^2)
		+ 11 (X_1 X_2 + X_2 X_0 + X_0 X_1) = 0,$\\
	\hth$	\{AG_i,GA_i\} = 2 (X_0^2 + X_1^2 + X_1^2)
		- 5 (X_1 X_2 + X_2 X_0 + X_0 X_1) = 0.$

	\sssec{Comment.}
	We started with the special Desargues configuration with 13 points
	(and lines) which are all of the points when $p = 3,$ the extended
	special Desargues configuration consists of adding 18 points and lines
	which are 31 distinct points and lines when p = $5.$
	It would appear that we could extend the construction in such a way
	that we get from the configuration with 31 points a configuration with
	57 points which would be all distinct when $p = 7,$ of 133 points which
	would be all distinct when $p = 11,$  $\ldots$  .
	But this is not possible.  For $p = 7,$ $(1,1,-1) \times (2,-1,0)$ gives
	(1,2,3) and by symmetry we get 5 other points but the points
	(0,1,3) give by symmetry $(3,0,1) = (1,0,-2)$ which has already been
	constructed.  Moreover, the point $(1,2,-3)$ gives by symmetry
	$(-3,1,2) = (1,2,-3)$ hence for $p = 7,$ the same point.  It is
	therefore
	not clear how to proceed in a systematic way.  This may be related to
	the fact that there are only 5 regular polyhedra which are associated
	to $p = 2,$ 3 and 5.  See Section 3.

	\sssec{Exercise.}\label{sec-epph}
	Rewrite the statement of Theorem \ref{sec-tpph}. in the form of a
	necessary and sufficient condition for $A_1,$ $A_3,$ $A_5$ to be
	collinear, given that $A_0,$ $A_2$ and $A_4$ are collinear.

	\sssec{Exercise.}\label{sec-einfcub}
	Let $\omega$ satisfy $\omega^2 + \omega + 1 = 0$.\\
	Let $P_0 = (0,1,-1),$ $Q_0 = (0,1,-\omega),$ $R_0 = (0,1,-\omega^2),$\\
	$P_1 = (-1,0,1),$ $Q_1 = (-\omega,0,1),$ $R_1 = (-\omega^2,0,1),$\\
	$P_2 = (1,-1,0),$ $Q_2 = (1,-\omega,0),$ $R_2 = (1,-\omega^2,0),$\\
	then, with\\
	$p = [1,1,1],$ $q = [1,\omega^2,\omega],$ $r = [1,\omega,\omega^2],$\\
	$p_0 = [1,0,0],$ $q[0] = [1,\omega,\omega],$
	$r_0 = [1,\omega^2,\omega^2],$
	\enumb
	\item $\omega^3 = 1,$
	\item incidence$(P_i,p)$, incidence$(Q_i,q)$, incidence$(R_i,r)$,
	\item incidence$(P_i,Q_i,R_i,p_i)$,
	\item incidence$(P_i,Q_{i+1},R_{i-1},q_i)$,incidence$(P_i,Q_{i-1},
	R_{i+1},r_i),$
	\item the configuration is therefore of type $9 * 4 \:\&\: 12 * 3.$
	\enume
	This configuration is that of the 9 inflection points of the cubic,
	$X_0^3+X_1^3+X_2^3 + k X_0X_1X_2 = 0$.

	\sssec{Comment.}
	Let {\bf\{}$A_0,A_1,A_2,A_3${\bf\}} be a complete quadrangle and
	$D_0,D_1,D_2$ be the diagonal points, several situation are possible
	in a perspective plane (See I).
	\enumb
	\item The diagonal points are always collinear, in this case, we have
	the N-Fano Configuration, N-Fano$({\bf\{}$A,B,C,D${\bf\}};\langle P,Q,R
	\rangle).$
	\item The diagonal points are never collinear,in this case, we have the
	Fano Configuration, Fano$({\bf\{}$A,B,C,D${\bf\}};\{P,Q,R\}).$
	\item The diagonal points are sometimes collinear,in this case, we have
	either the pseudo configuration, $({\bf\{}$A,B,C,D${\bf\}},\langle P,Q,R
	\rangle)$ or the pseudo configuration, $({\bf\{}$A,B,C,D${\bf\}},
	\{P,Q,R\}).$
	\enume
	 Notice the ";" in the first 2 cases.

	\ssec{Proof of the Theorem of Desargues.
	The hexagon of Pappus-Brianchon.
	The configuration of Reidemeister.}
	\sssec{Proof of the Theorem of Desargues.}\label{sec-pdes}

	Proof:  The proof that Theorem \ref{sec-tdes} follows from the axioms
	of incidence and of Pappus will now be given.\\
	Cronheim (1953)\footnote{Proc. Amer. Math. Soc., {\bf 4}, 219-221.}
	showed that the proof reduces to 2 cases. In the first
	one, a permutation of the indices 0, 1, 2 is chosen in such a way that
	$A_0\notinc b_0$ and $B_2\notinc a_2$. In the second one, except
	perhaps for an exchange of $A_i$ and $B_i$, $B_i\incid a_i.$

	In the first case (Hessenberg, 1905), let\\
	He1.0.\hti{3}$A_0\notinc b_0,$ $B_2\notinc a_2$,\\
	De1.0.\hti{3}$d := A_0\times B_2,$ $D := d\times  c_1,$
		$e := D\times C_2,$ $E := e\times b_1,$\\
	De1.1.\hti{3}$f := D \times C_0,$ $F := f\times a_1,$
		$G := a_2\times b_0,$ $g := F\times G.$\\
	De2.0.\hti{3}$X := d\times a_0,$ $Y := d\times b_2,$ $Z := d\times g$,\\
	The Pappus-Pascal hexagon $D,$ $A_1,$ $A_0,$ $A_2,$ $B_2,$ $C_0
	\implies G,$ $C$ and $F$ are collinear.\\
	The Pappus-Pascal hexagon $D,$ $B_1,$ $B_2,$ $B_0,$ $A_0,$ $C_2
	\implies G,$ $C$ and $E$ are collinear.\\
	Hence $A_0,$ $F,$ $E$, $D,$ $B_2,$ $G$ is a Pappus-Pascal
	hexagon and $C_0,$ $C_1$ and $C_2$ are collinear.
	It is easy to verify that, because of He1.0, $X$ is distinct from
	$D,A_0,B_2,A_2,C_0,A_1$, that $Y$ is distinct from $D,B_2,A_0,B_0,C_2,
	B_1$ and that $Z$ is distinct from $A_0,D,B_2,E,G,F.$

	This will be abbreviated as follows.\\
	Pr1.0.\hti{3}Pappus$(\langle D,A_0,B_2\rangle,d,
		\langle A_2,C_0,A_1\rangle,a_0;\langle G,C,F\rangle,X),$\\
	Pr1.1.\hti{3}Pappus$(\langle D,B_2,A_0\rangle,d,
		\langle B_0,C_2,B_1\rangle,b_2;\langle G,C,E\rangle,Y),$\\
	Pr1.2.\hti{3}$(\langle E,G,F\rangle,g),$\\
	Pr1.3.\hti{3}Pappus$(\langle A_0,D,B_2\rangle,d,\langle E,G,F\rangle,g;
		\langle C_i\rangle,Z),$\\
	Pr1.4.\hti{3}$\langle C_i\rangle.$

	In the second case (Cronheim, 1953), we have the 3-Desargues
	configuration (Fig. 3e), let\\
	De3.0.\hti{3}$r_2 := A_2 \times C_2,$
		$R_0 := c_0 \times r_2,$ $R_1 := c_1 \times r_2,$\\
	De4.0.\hti{3}$X := c_2\times b_2$,\\
	Pr3.0.\hti{3}Pappus($\langle C,A_2,B_2\rangle,c_2,
		\langle C_2,B_1,B_0\rangle,b_2;\langle C_0,A_0,R_1\rangle
		,r_0,X),$\\
	Pr3.1.\hti{3}Pappus($\langle C,A_2,B_2\rangle,c_2,
		\langle C_2,B_0,B_1\rangle,b_2;\langle C_1,A_1,R_0\rangle
		,r_1,X),$\\
	Pr3.2.\hti{3}Pappus($\langle R_0,B_0,A_0\rangle,c_0,
		\langle B_1,R_1,A_1\rangle,c_1;\langle C_0,C_1,C_2\rangle,c,C).$

	\sssec{Exercise.}
	Prove the Theorem of Cronheim, on the reduction to 2 cases, refered to
	in \ref{sec-pdes}.

	\sssec{Theorem. [Dual of Pappus]}\label{sec-tdupap}
	{\em If the alternate sides of the hexagon \{$a_0,$ $a_1,$ $a_2,$ $a_3,$
	$a_4,$ $a_5$\} pass through 2 points,
	three pairs of opposites points of the hexagon are on 3 lines $p_0,$
	$p_1$ and $p_2$ which pass through the same point.}
	See \ref{sec-Sdua}.  (Fig. 1b)

	We write dual-Pappus$(\langle a_2,a_0,a_4\rangle,
	\langle a_5,a_3,a_1\rangle,\langle p_0,p_1,p_2\rangle).$

	Proof:  Let $A_0 := a_0 \times a_1,$ $A_1 := a_1 \times a_2,$
	$A_2 := a_2 \times a_3,$ $A_3 := a_3 \times a_4,$
	$A_4 := a_4 \times a_5,$ $A_5 := a_5 \times a_0.$\\ 
	Let $B_0 := a_0 \times a_2,$ $B_1 := a_1 \times a_3,$
	$p_0 := A_0 \times A_3,$ $p_1 := A_1 \times A_4,$ 
	$p_2 := A_2 \times A_5,$ $B_2 := p_0 \times p_2.$\\
	By hypothesis, $B_0 \cdot a_4 = B_1 \cdot a_5 = 0.$
	$B_0,$ $A_2,$ $A_5,$ $B_1,$ $A_0,$ $A_3$
	is a Pappus-Pascal hexagon, therefore $A_1,$ $B_2$ and $A_4$ are
	collinear, in other words $p_1$ passes through $B_2.$

	\sssec{Definition.}
	The preceding configuration is a degenerate form of that
	of Brianchon.  I will call it the {\em Pappus-Brianchon hexagon}.
	The point common to $p_0,$ $p_1$ and $p_2$ is called the
	{\em Pappus point}.

	\sssec{Proof of the special Desargues Theorem.}
	\label{sec-psd}

	The proof of Theorem \ref{sec-tsd} is as follows:
	0. is a direct consequence of \ref{sec-pdes}.
	1. follows from the Axiom of Pappus \ref{sec-ax1}.4. applied to the
	points $P_1,$ $A_2,$ $R_1$ and $P_2,$
	$A_1,$ $R_2,$ proving that $Q_0,$ $P_0$ and $P$ are collinear.

	\sssec{Exercise.}
	The proof \ref{sec-psd} of Theorem \ref{sec-tsd} is only given in the
	general case.  Describe all the exceptional cases and give a proof for
	each case.

	\sssec{Definition.}
	The {\em Reidemeister configuration} consists of 11 points\\
	$A_0,$ $A_1,$ $A_2,$ $B_{00},$ $B_{11},$ $B_{22},$ $B_{33},$ $B_{01},$
	$B_{10},$ $B_{23},$ $B_{32},$\\
	and 15 lines,\\
	$a_0,$ $a_1,$ $a_2,$ $b_{00},$ $b_{01},$ $b_{02},$ $b_{03},$ $b_{10},$
	$b_{11},$ $b_{12},$ $b_{13},$ $b_{20},$ $b_{21},$ $b_{22},$ $b_{23}$:\\
	Let $A_0,$ $A_1,$ $A_2$ be a triangle, $a_0 := A_1 \times A_2,$
	$a_1 := A_2 \times A_0,$ $a_2 := A_0 \times A_1,$\\
	let $b_{00},$ $b_{01},$ $b_{02}$ be 3 lines through $A_0$ distinct from
	$a_1$ and $a_2,$\\
	let $B_{00},$ $B_{22}$ be points on $b_{01}$ not on $a_0,$
	$b_{10} := A_1 \times B_{00},$ $ b_{12} := A_1 \times B_{22},$
	$b_{20} := A_2 \times B_{00},$ $b_{22} := A_2 \times B_{22},$
	$B_{01} := b_{00} \times b_{10},$ $B_{23} := b_{00} \times b_{12},$
	$B_{10} := b_{02} \times b_{20},$ $B_{32} := b_{02} \times b_{22},$
	$b_{11} := A_1 \times B_{10},$ $b_{13} := A_1 \times B_{32},$
	$b_{21} := A_2 \times B_{01},$ $b_{23} := A_2 \times B_{23},$
	$B_{11} := b_{11} \times b_{21},$ $B_{33} := b_{13} \times b_{23},$
	$b_{03} := B_{11} \times B_{33}.$ (Fig. 11a)

	\sssec{Lemma.}
	Let $c_{00} := B_{01} \times B_{10},$ $c_{01} := B_{32} \times B_{23},$
	$C_0 := c_{00} \times c_{01},$ then\\
	incidence$(C_0,A_1,A_2)$.

	Proof:\\
	Desargues$(A_0,B_{00}\:B_{10}\:B_{01},B_{22}\:B_{32}\:B_{23};
	C_0\:A_1\:A_2,a_0)$.

	\sssec{Theorem. [Reidemeister]}\label{sec-treid}
	\enumb
	\item$A_0\cdot b_{03} = 0.$
	\item{\em The Reidemeister configuration is of type}\\
	\hth$	3 * 6 + 8 * 3  \:\&\:  12 * 3 + 3 * 2.$
	\enume

	Proof: After using the preceding Lemma, we use\\
	Desargues$(a_0,\{c_{00},b_{11},b_{21}\},\{c_{01},b_{13},b_{23}\};
	\langle b_{03},b_{00},b_{02}\rangle,A_0)$.

	\sssec{Theorem.}
	Let\\
	\hth$c_{02} := B_{11} \times B_{22},$\\
	\hth$c_{12} := B_{01} \times B_{32},$\\
	\hth$c_{22} := B_{23} \times B_{10},$\\
	then\\
	\hth incidence$(c_{i2},C)$.

	Proof:\\
	Desargues$^{-1}(a_2,\{b_{00},b_{12},c_{20}\},\{b_{02},b_{11},c_{21}\};
		\langle c_{02},c_{12},c_{22}\rangle,C)$,

	\sssec{Theorem.}
	Let\\
	\hth$c_{00} := B_{01} \times B_{10},$ $c_{01} := B_{32} \times B_{23},$
\\
	\hth$c_{10} := B_{10} \times B_{22},$ $c_{11} := B_{23} \times B_{11},$
\\
	\hth$c_{20} := B_{22} \times B_{01},$ $c_{21} := B_{11} \times B_{32},$
\\
	\hth$C_i := c_{i0} \times c_{i1},\\
	then\\
	\hth$incidence$(C_i,c)$.

	Proof: Using the preceding Theorem,\\
	Desargues$(C,\{B_{22},B_{01},B_{10}\},\{B_{11},B_{32},B_{23}\};
	\langle C_i\rangle,c).$

	\sssec{Exercise.}
	Let\\
	\hth$c'_{00} := B_{00}\times B_{11},$ $c'_{01} := B_{22}\times B_{33},$
\\
	\hth$c'_{10} := B_{00}\times B_{32},$ $c'_{11} := B_{01}\times B_{33},$
\\
	\hth$c'_{20} := B_{00}\times B_{23},$ $c'_{21} := B_{10}\times B_{33},$
\\
	then\\
	incidence$(C_i,C'_{i+1},C'_{i-1}).$

	\sssec{Definition.}
	The {\em extended Reidemeister configuration} consists of the points
	$A_0,$ $A_1,$ $A_2,$ $B_{jj},$ $j$ = 0,1,2,3, $B_{01},$ $B_{10},$
	$B_{23},$ $B_{32},$ $C$, $C_i$, $C'_i$, i = 0,1,2,
	and of the lines $a_0,$ $a_1,$ $a_2,$ $b_{ij},$ $i$ = 0,1,2,
	$j$ = 0,1,2,3, $c_{ik},$ $c'_{ik},$ $c'_i,$, $i,$ = 0,1,2, $k$ = 0,1,
	$c0,$ see Fig. 11f.

	\sssec{Exercise.}\label{sec-ereid}
	Prove
	\enumb
	\item[0.~~] that for given $A_0,$ $A_1,$ $A_2,$ $B_{01},$ $B_{10},$ the
		correspondance 	between $B_{22}$ and $B_{33}$ is a projectivity
		with center $AEul_0$ on $A_1 \times A_2$ (See
		\ref{sec-Sprojinv}).
	\item[1.~~] The lines $b_{01}$ and $b_{03}$ coincide if the point
		$B_{10}$ is on the conic through $B_{01}$ tangent at $A_1$ to
		$A_1 \times A_0$ and tangent at $A_2$ to $A_2 \times A_0,$
		represented by the matrix\\
		$\matb{[}{ccc}2&0&0\\0&0&-1\\0&-1&0\mate{]}$
	\item[2.~~] that if we permute cyclically $A_0,$ $A_1,$ $A_2,$, then
	\item[  0.] the lines $B_{00} \times B_{11}$ and the 2 other
		corresponding lines pass through the same point $K$.
	\item[  1.] the lines $A_0 \times B_{22}$ and the 2 other
		corresponding lines pass through the same point $P$. The same
		is true for the lines $A_0 \times B_{33}$ and the 2 other
		corresponding lines, giving $\ov{P}$.
	\enume

	This configuration, see Fig. 26b, which I will call the
	{\em K-Reidemeister configuration} is part of the Hexal configuration
	studied in Chapter III, with the correspondance\\
	$\begin{array}{cccccccccccc}
	\hti{4}&A_{i}&B_{01}&B_{10}&b_{00}&b_{10}&b_{21}&b_{02}&b_{11}&b_{20}
		&c_{00}&c'_{00}\\
	    &A_{i}&M&\ov{M}&ma_0&ma_1&ma_2&\ov{m}a_0&\ov{m}a_1
	&\ov{m}a_2&eul&mMa_0\\
	\cline{2-12}
		&B_{00}&B_{11}&b_{01}&b_{03}\\
	&Maa_{0}&\ov{M}a_0&cc_0&\ov{c}c_0
	\end{array}$

	\ssec{The extended Pappus configuration and a remarkable
	Theorem.}

	\setcounter{subsubsection}{-1}
	\sssec{Introduction.}
	If we permute in all possible way the 6 points of the
	Pappus configuration we obtain 6 Pappus lines.  I prove in Theorem
	\ref{sec-tpap1}.  that these pass 3 by 3 through 2 points.  We obtain
	therefore a dual configuration, which therefore determines 6 Pappus
	points, which are 3 by 3 on 2 lines.  I prove in Theorem \ref{sec-tpap2}
	that these lines are
	the 2 original ones of the Pappus configuration.  The points are not,
	in general the same.  The proof of the first Theorem is synthetic, I
	have no synthetic proof of the second Theorem.  The algebraic proof
	uses a notation introduced in Chapter III.  Special cases of this
	configuration have been studied, but because some of the results are
	still at the conjecture stage, these will not be discussed here, others
	are given as exercises.\\
	The term ``rotate the points
	$\ov{M}_0,$ $\ov{M}_1,$ $\ov{M}_2$" means that we take
	the even permutations of 
	$\ov{M}_0,$ $\ov{M}_1,$ $\ov{M}_2$, namely
	$\ov{M}_1,$ $\ov{M}_2,$ $\ov{M}_0$ and
	$\ov{M}_2,$ $\ov{M}_0,$ $\ov{M}_1$.\\
	The notation is explained, in details, in Chapter III.

	\sssec{Theorem. [Steiner (Pappus)]}\label{sec-tpap1}
	\vspace{-18pt}\hspace{180pt}\footnote{Steiner, Werke, I, p. 451}\\[8pt]
	{\em If we fix the points $M_0,$ $M_1,$ $M_2$ on $d$ and rotate the
	points $\ov{M}_0,$ $\ov{M}_1,$ $\ov{M}_2$ on
	$\ov{d},$ we obtain the 3 Pappus lines $m_0,$ $m_1,$ 
	$m_2.$	 These pass through the point $D$.  Similarly if we reverse the
	order of the points of $\ov{d}$ and rotate, we obtain 3 other
	Pappus lines, $\ov{m}_0,$	$\ov{m}_1,$ $\ov{m}_2.$
	These pass through the point $\ov{D}$.  In detail, let}\\
	H0.\hti{6}$M_i,$ $\ov{M}_i,$\\
	D0.\hti{6}$	a_i := M_i \times \ov{M}_i,$\\
	D1.\hti{6}$    b_i := M_{i+1} \times \ov{M}_{i-1},$
	$\ov{b}_i := \ov{M}_{i+1} \times M_{i-1},$\\
	D2.\hti{6}$	L_i := b_i \times \ov{b}_i,$\\
	D3.\hti{6}$    N_i := a_i \times b_i,$
	$\ov{N}_i := a_i \times \ov{b}_i,$\\
	D4.\hti{6}$	m_0 := L_1 \times L_2,$\\
	D5.\hti{6}$    m_1 := N_1 \times N_2,$
	$m_2 = \ov{N}_1 \times \ov{N}_2,$\\
	D6.\hti{6}$	D := m_1 \times m_2,$\\
	D7.\hti{6}$	Q_i := a_{i+1} \times a_{i-1},$\\
	D8.\hti{6}$    P_i := b_{i+1} \times b_{i-1},$
	$\ov{P}_i := \ov{b}_{i+1} \times \ov{b}_{i-1},$\\
	D9.\hti{6}$	\ov{m}_i := P_i \times \ov{P}_i,$\\
	D10.\hti{5}$	\ov{D} := \ov{m}_1 \times
	\ov{m}_2,$\\
	{\em then}\\
	C0.\hti{6}$m_0.L_0 = 0 (*).$\\
	C1.\hti{6}$    m_1.N_0 = m_2 \cdot \ov{N}_0 = 0 (*).$\\
	C2.\hti{6}$D.m_0 = 0.$\\
	C3.\hti{6}$\ov{m}_i \cdot Q_i = 0.$\\
	C4.\hti{6}$\ov{D}\cdot \ov{m}_0 = 0 (*).$  See Fig. 9,

	Proof:  A synthetic proof is as follows,
	C0, C1, C3, are direct consequences of Pappus' theorem applied to\\
	Pappus($\langle M_0,M_1,M_2\rangle,d,\langle\ov{M}_0,\ov{M}_1,\ov{M}_2,
	\rangle,\ov{d};\langle L_0,L_1,L_2\rangle,m_0)$\\
	Pappus($\langle M_2,M_0,M_1\rangle,d,\langle\ov{M}_1,\ov{M}_2,\ov{M}_0,
	\rangle,\ov{d};\langle N_0,N_1,N_2\rangle,m_0)$\\
	Pappus($\langle M_1,M_2,M_0\rangle,d,\langle\ov{M}_2,\ov{M}_0,\ov{M}_1,
	\rangle,\ov{d};\langle \ov{N}_0,\ov{N}_1,\ov{N}_2\rangle,m_0)$\\
	Pappus($\langle M_0,M_1,M_2\rangle,d,\langle\ov{M}_0,\ov{M}_1,\ov{M}_2,
	\rangle,\ov{d};\langle L_0,L_1,L_2\rangle,m_0)$\\
	$\ov{M}_2,$ $\ov{M}_0,$ $\ov{M}_1$ or 
	$\ov{M}_1,$ $\ov{M}_2,$ $\ov{M}_0$ 
	and $\ov{M}_2,$ $\ov{M}_1,$ $\ov{M}_0,$ or
	$\ov{M}_0,$ $\ov{M}_2,$ $\ov{M}_1$ or
	 $\ov{M}_1,$ $\ov{M}_0,$ $\ov{M}_2.$ 
	The triangles $L_i,$ $N_i,$ $\ov{N}_i$ have $\ov{m}_0$ as
	 axis of perspectivity
	for $i = 1$ and 2 therefore they have a center of perspectivity $D,$ by
	Desargues.  I also note that for $i = 2$ and 0 the axis is
	$\ov{m}_1$ and for
	$i = 0$ and 1 the axis is $\ov{m}_2.$  Hence C2.
	Symmetrically we get C4.

	For an algebraic proof, useful because of \ref{sec-tpap2}, let\\
	H0.\hti{6}$    M_0 = (0,1,-1), M_1 = (-1,0,1), M_2 = (1,-1,0)$\\
	H1.\hti{6}$\ov{M}_0 = (0,$m$_2,-$m$_1), \ov{M}_1 = (-$m$_2,0,$m$_0),$
	$\ov{M}_2 = ($m$_1,-$m$_0,0).$\\
	then\\
	P0.\hti{6}$     a_0 = [1,0,0].$\\
	P1.\hti{6}$    b_0 = [$m$_0,$m$_1,$m$_0], \ov{b}_0 = [$m$_0,$m$_0,$m$_2].$\\
	P2.\hti{6}$     L_0 = ($m$_0^2-$m$_1$m$_2,$m$_0($m$_2-$m$_0),-$m$_0($m$_0-$m$_1)).$\\
	P3.\hti{6}$    N_0 = (0,$m$_0,-$m$_1), \ov{N}_0 = (0,$m$_2,-$m$_0).$\\
	P4.\hti{6}$     m_0 = [$m$_0($m$_1+$m$_2),$m$_1($m$_2+$m$_0),$m$_2($m$_0+$m$_1)],$\\
	P5.\hti{6}$    m_1 = [$m$_0$m$_1,$m$_1$m$_2,$m$_2$m$_0],	m_2 = [$m$_2$m$_0,$m$_0$m$_1,$m$_1$m$_2],$\\
	P6.\hti{6}$D = ( $m$_1$m$_2($m$_0^2-$m$_1$m$_2),$m$_2$m$_0($m$_1^2-$m$_2$m$_0),$m$_0$m$_1($m$_2^2-$m$_0$m$_1)).$\\
	P7.\hti{6}$     Q_0 = (1,0,0).$\\
	P8.\hti{6}$    P_0 = ($m$_2($m$_1-$m$_2),$m$_2($m$_0-$m$_1),$m$_1($m$_2-$m$_0)),$\\
	\hti{9}$\ov{P}_0 = ($m$_1($m$_1-$m$_2),$m$_2($m$_0-$m$_1),$m$_1($m$_2-$m$_0)).$\\
	P9.\hti{6}$\ov{m}_0 = [0,$m$_1($m$_2-$m$_0),-$m$_2($m$_0-$m$_1)],$\\
	P10.\hti{5}$\ov{D} = ($m$_1$m$_2($m$_2-$m$_0)($m$_0-$m$_1),$m$_2$m$_0($m$_0-$m$_1)($m$_1-$m$_2),\\
	\hti{12}$m$_0$m$_1($m$_1-$m$_2)($m$_2-$m$_0)).$

	\sssec{Definition.}
	The configuration of Theorem \ref{sec-tpap1} which consists of 26
	points and 17 lines is called the {\em extended Pappus configuration}.

	It is of type $6 * 4 + 20 * 3  \:\&\:  9 * 6 + 6 * 4 + 2 * 3.$ (10)

	It can also be viewed, because of the synthetic proof as a multiple
	Desargues configuration, with 3 triangles perspective from $D$ and 3
	triangles perspective from $D'$ in which the axis of one are the
	concurrent lines of the other.

	\sssec{Definition. [Steiner]}
	The sub-configuration consisting of the points $L_i,$ $N_i,$
	$\ov{N}_i,$ $Q_i,$ $P_i,$ $\ov{P}_i,$ $D$,
	$\ov{D}$ and of the lines $a_i$, $b_i,$ $\ov{b}_i,$
	$m_i,$ $\ov{m}_i,$ is called the {\em Steiner configuration}.
	It is of type\\
	\hth$	20 * 3\: \&\: 15 * 4.$ (10)

	\sssec{Comment.}\label{sec-csteiner}
	Part of a dual of the extended configuration is described in
	sections 1, 3 and 4 of the involutive geometry of the triangle.
	The relation between the notations is as follows:\\
	$\begin{array}{llllllllll}
	d&\ov{d}&M_i&   \ov{M}_i&  a_i&   b_i&   \ov{b}_i& 
	 L_i&   N_i&   \ov{N}_i \\
	M&\ov{M}&ma_i&  \ov{m}a_i& A_i&   Maa_i& \ov{M}aa_i&
	mMa_i&cc_i&\ov{c}c_i\\
	\cline{1-10}
	m_0&	  m_1&   m_2&   D&      Q_i&   P_i&  
	\ov{P}_i&  \ov{m}_i&  \ov{D}\\
	K&    P& \ov{P}&  pp&a_i&pap_i&\ov{p}ap_i&
	Pap_i&pap
	\end{array}$

	In particular, $K$ is the point of Lemoine.
	On the other hand there is the following correspondence
	$\ov{m}_i$ and the dual of $abr1_i,$ $\ov{D}$
	and the dual of $Ste,$ which passes through $BRa$ and $Abr.$

	\sssec{Theorem.}\label{sec-tpap2}
	{\em If we make the dual construction starting with $m_0,$ $m_1,$ 
	$m_2$ on $D$ and $\ov{m}_i$ on $\ov{D},$
	the points $Ma_i$ dual of $m_i$ is on the original line $d$ and those
	$\ov{M}a_i$ dual of $\ov{m}_i$ is on the original line $\ov{d}$}:\\
	C5.\hti{6}$    Ma_i \cdot d = 0.$\\
	C6.\hti{6}$    \ov{M}a_i \cdot \ov{d} = 0.$ See Fig. 10,

	Proof:  An algebraic proof is as follows.
	\footnote{The reader will want to wait to check these algebraic
	manipulations until the notation has been explained.}\\
	P'0.\hti{4}$A_0 = (2$m$_1$m$_2($m$_1-$m$_2),$m$_2($m$_1+$m$_2)
		($m$_0-$m$_1),$m$_1($m$_1+$m$_2)($m$_2-$m$_0)).$\\
	\hth$	A_1 = ($m$_1$m$_2($m$_0-$m$_1),$m$_0($m$_1^2+$m$_2$m$_0
		-2$m$_0$m$_1),$m$_0$m$_1($m$_1-$m$_2)),$\\
	\hth$	A_2 = ($m$_1$m$_2($m$_2-$m$_0),$m$_2$m$_0($m$_1-$m$_2),
		-$m$_0($m$_2^2-2$m$_2$m$_0+$m$_0$m$_1)),$\\
	P'1.\hti{4}$   B_0 = ($m$_1$m$_2($m$_2-$m$_0),$m$_2$m$_0($m$_1-$m$_2),
		$m$_1($m$_2^2-2$m$_1$m$_2+$m$_0$m$_1)),$\\
	\hth$	B_1 = (-$m$_1($m$_0^2+$m$_1$m$_2-2$m$_0$m$_1),$m$_2$m$_0($m$_0
		-$m$_1),$m$_0$m$_1($m$_2-$m$_0)),$\\
	\hth$	B_2 = ($m$_2($m$_2+$m$_0)($m$_0-$m$_1),2$m$_2$m$_0($m$_2-$m$_0)
		,$m$_0($m$_1-$m$_2)($m$_2+$m$_0)),$\\
	\hth$\ov{B}_0 = ($m$_1$m$_2($m$_0-$m$_1),-$m$_2($m$_1^2
		-2$m$_1$m$_2+$m$_2$m$_0),$m$_0$m$_1($m$_1-$m$_2)),$\\
	\hth$\ov{B}_1 = ($m$_1($m$_2-$m$_0)($m$_0+$m$_1),$m$_0($m$_1
		-$m$_2)($m$_0+$m$_1),
	2$m$_0$m$_1($m$_0-$m$_1)),$\\
	\hth$\ov{B}_2 = ($m$_2($m$_0^2+$m$_1$m$_2-2$m$_2$m$_0),
		$m$_2$m$_0($m$_0-$m$_1),$m$_0$m$_1($m$_2-$m$_0)).$\\
	P'2.\hti{4}$l_0 = [2$m$_1^3$m$_2+2$m$_2^3$m$_1-$m$_2^3$m$_0
		-$m$_1^3$m$_0-5$m$_1^2$m$_2^2-$m$_2^2$m$_0^2$\\
	\hti{16}$-$m$_0^2$m$_1^2+$m$_0$m$_1$m$_2($m$_0+2$m$_1+2$m$_2),$\\
	\hti{12}m$_1($m$_2^2$m$_1-2$m$_1^2$m$_2+$m$_0^2$m$_2-2$m$_2^2$m$_0
		-2$m$_0^2$m$_1+$m$_1^2$m$_0+3$m$_0$m$_1$m$_2),$\\
	\hti{12}m$_2($m$_1^2$m$_2-2$m$_2^2$m$_1+$m$_0^2$m$_1-2$m$_1^2$m$_0
		-2$m$_0^2$m$_2+$m$_2^2$m$_0+3$m$_0$m$_1$m$_2),$\\
	\hth$	l_1 =    [ $m$_0($s$_{21}-6$m$_2$m$_1$m$_0),$
	\footnote{s$_{21}$ is the symmetric function in m$_i$, namely,
	m$_0^2($m$_1+$m$_2)+$m$_1^2($m$_2+$m$_0)+$m$_2^2($m$_0+$m$_1).$}\\
	\hti{12}m$_1(3$m$_0^3+4$m$_1^2$m$_0-5$m$_1$m$_0^2-2$m$_0^2$m$_2
		+$m$_0$m$_2^2+$m$_2^2$m$_1-2$m$_2$m$_1^2),$\\
	\hti{12}$-($m$_1+$m$_0)($m$_1$m$_0^2+$m$_0$m$_2^2+$m$_2$m$_1^2
		-2($m$_1^2$m$_0+$m$_0^2$m$_2+$m$_2^2$m$_1)\\
	\hti{16}+3$m$_2$m$_1$m$_0)],$\\
	\hth$	l_2 =    [ $m$_0($s$_{21}-6$m$_1$m$_2$m$_0),$\\
	\hti{12}$-($m$_2+$m$_0)($m$_2$m$_0^2+$m$_0$m$_1^2
		+$m$_1$m$_2^2-2($m$_2^2$m$_0+$m$_0^2$m$_1+$m$_1^2$m$_2)\\
	\hti{16}+3$m$_1$m$_2$m$_0),$\\
	\hti{12}m$_2(3$m$_0^3+4$m$_2^2$m$_0-5$m$_2$m$_0^2-2$m$_0^2$m$_1
		+$m$_0$m$_1^2+$m$_1^2$m$_2-2$m$_1$m$_2^2)],$\\
	P'3.\hti{4}$n_0 = [-($m$_1+$m$_2)($m$_1$m$_2^2+$m$_2$m$_0^2
	+$m$_0$m$_1^2-2($m$_1^2$m$_2+$m$_2^2$m$_0+$m$_0^2$m$_1)\\
	\hti{16}+3$m$_0$m$_1$m$_2),$\\
	\hti{12}m$_1(3$m$_2^3+4$m$_1^2$m$_2-5$m$_1$m$_2^2-2$m$_2^2$m$_0
		+$m$_2$m$_0^2+$m$_0^2$m$_1-2$m$_0$m$_1^2),$\\
	\hti{12}m$_2($s$_{21}-6$m$_0$m$_1$m$_2),$\\
	\hth$\ov{n}_0 =
		[($m$_1+$m$_2)(-($m$_1^2$m$_2+$m$_2^2$m$_0+$m$_0^2$m$_1)
		+2($m$_2^2$m$_1+$m$_0^2$m$_2+$m$_1^2$m$_0)\\
	\hti{16}-3$m$_0$m$_1$m$_2),
	$\\
	\hti{12} m$_1($s$_{21}-6$m$_0$m$_1$m$_2),$\\
	\hti{12} m$_2(3$m$_1^3-5$m$_1^2$m$_2+4$m$_2^2$m$_1-2$m$_2^2$m$_0
		+$m$_0^2$m$_2+$m$_0^2$m$_1-2$m$_0$m$_1^2)],$\\
	\hth$	l_2 = (l_{10},l_{12},l_{11}) ($m$_1,$m$_0,$m$_2).$\\
	\hth$	n_0 = (l_{12},l_{11},l_{10}) ($m$_2,$m$_1,$m$_0).$\\
	\hth$	n_1 = (l_{01},l_{02},l_{00}) ($m$_2,$m$_0,$m$_1).$\\
	\hth$	n_2 = (l_{11},l_{12},l_{10}) ($m$_2,$m$_0,$m$_1).$\\
	P'4.\hti{4}$   Ma_0 = (2$m$_0-$m$_1-$m$_2,2$m$_1-$m$_2-$m$_0,2$m$_2
		-$m$_0-$m$_1).$\\
	P'5.\hti{4}$   Ma_1 = (2$m$_1-$m$_2-$m$_0,2$m$_0-$m$_1-$m$_2,2$m$_2
		-$m$_0-$m$_1),$\\
	\hth$	Ma_2 = (2$m$_2-$m$_0-$m$_1,2$m$_1-$m$_2-$m$_0,2$m$_0-$m$_1
		-$m$_2).$\\
	P'7.\hti{4}$q_0 = [-$m$_0($m$_1^3($m$_2+$m$_0)+$m$_2^3($m$_0+$m$_1)
		-$m$_1^2$m$_2^2-2$m$_2^2$m$_0^2
	-2$m$_0^2$m$_1^2$\\
	\hti{16}$		+$m$_0$m$_1$m$_2(5$m$_0-2$m$_1-2$m$_2)),$\\
	\hti{12}$-$m$_1$m$_2(-($m$_1^2$m$_2+$m$_2^2$m$_0+$m$_0^2$m$_1)
		+2($m$_2^2$m$_1+$m$_0^2$m$_2+$m$_1^2$m$_0)\\
	\hti{16}-3$m$_0$m$_1$m$_2),$\\
	\hti{12} m$_1$m$_2(-2($m$_1^2$m$_2+$m$_2^2$m$_0+$m$_0^2$m$_1
		+($m$_2^2$m$_1+$m$_0^2$m$_2+$m$_1^2$m$_0)\\
	\hti{16}+3$m$_0$m$_1$m$_2)],$\\
	P'8.\hti{4}$   p_0 = [ -$m$_2$m$_0($s$_{21}-6$m$_0$m$_1$m$_2),$\\
	\hti{12}$-$m$_1($m$_2+$m$_0)(-($m$_1^2$m$_2+$m$_2^2$m$_0+$m$_0^2$m$_1)
	+2($m$_2^2$m$_1+$m$_0^2$m$_2+$m$_1^2$m$_0)\\
	\hti{16}-3$m$_0$m$_1$m$_2),$\\
	\hti{12}$-$m$_2($m$_0^3$m$_2-2$m$_0^3$m$_1+3$m$_1^2$m$_2^2
	+$m$_2^2$m$_0^2+4$m$_0^2$m$_1^2\\
	\hti{16}-$m$_0$m$_1$m$_2(2$m$_2+5$m$_1))],$\\
	\hth$\ov{p}_0 = [ -$m$_0$m$_1($s$_{21}-6$m$_0$m$_1$m$_2),$\\
	\hti{12}$-$m$_1($m$_0^3$m$_1-2$m$_0^3$m$_2+3$m$_1^2$m$_2^2
	+$m$_0^2$m$_1^2+4$m$_2^2$m$_0^2\\
	\hti{16}-$m$_0$m$_1$m$_2(2$m$_1+5$m$_2)),$\\
	\hti{12} m$_2($m$_0+$m$_1)(-2($m$_1^2$m$_2+$m$_2^2$m$_0+$m$_0^2$m$_1)
		+($m$_2^2$m$_1+$m$_0^2$m$_2+$m$_1^2$m$_0)\\
	\hti{16}+3$m$_0$m$_1$m$_2)
	].$\\
	P'9.\hti{4}$\ov{M}a_0 = ( $m$_1$m$_2(2$m$_1$m$_2-$m$_2$m$_0
		-$m$_0$m$_1),
	$m$_2$m$_0(2$m$_0$m$_1-$m$_1$m$_2-$m$_2$m$_0),$\\
	\hti{12}$	$m$_0$m$_1(2$m$_2$m$_0-$m$_0$m$_1-$m$_1$m$_2),$\\
	\hth$\ov{M}a_1 = ( $m$_1$m$_2(2$m$_2$m$_0-$m$_0$m$_1-$m$_1$m$_2),
		$m$_2$m$_0(2$m$_1$m$_2-$m$_2$m$_0-$m$_0$m$_1),$\\
	\hti{12}$	$m$_0$m$_1(2$m$_0$m$_1-$m$_1$m$_2-$m$_2$m$_0)),$\\
	\hth$\ov{M}a_2 = ( $m$_1$m$_2(2$m$_0$m$_1-$m$_1$m$_2-$m$_2$m$_0),
		$m$_2$m$_0(2$m$_2$m$_0-$m$_0$m$_1-$m$_1$m$_2),$\\
	\hti{12}$	$m$_0$m$_1(2$m$_1$m$_2-$m$_2$m$_0-$m$_0$m$_1)).$

	\sssec{Comment.}
	Continuing \ref{sec-csteiner} we have the following relation
	between the above notation and that in the involutive geometry of
	the triangle.\\
	$\begin{array}{llllllllll}
	d&\ov{d}&M_i&   \ov{M}_i\\
	pp&pap&K,P,P~&Pap_i\\
	\cline{1-4}
	a_i&   b_i&   \ov{b}_i\\ 
	kpa_0,\ov{t}pa_1,tpa_2&\ov{t}pa_2,tpa_0,kpa_1
	&tpa_1,kpa_2,\ov{t}pa_0\\
	\cline{1-4}
	L_i&   N_i&   \ov{N}_i \\
	Ttp_0,Tkp_1,Tk\ov{p}_2&Tk\ov{p}_1,Ttp_2,Tkp_0&Tkp_2,Tk\ov{p}_0,Ttp_1\\
	\cline{1-3}
	m_i&   D&  \ov{m}_i&  \ov{D}\\
	apa_i&M&\ov{a}pa_0,\ov{a}pa_2,\ov{a}pa_1&\ov{M}\\
	\cline{1-4}
	Q_i&P_i&\ov{P}_i\\
	\ov{T}tp_0,\ov{T}k\ov{p}_1,\ov{T}kp_2&
	\ov{T}k\ov{p}_2,\ov{T}kp_0,\ov{T}tp_1&\ov{T}kp_1,\ov{T}tp_2
	,\ov{T}k\ov{p}_0&
	\end{array}$
 
	\sssec{Definition.}
	The mapping which associates to the points $M_i$ and $\ov{M}_i,$
	the points $Ma_i$ and $\ov{M}a_i,$ is called the
	{\em Pappus-dual-Pappus mapping}.

	\sssec{Exercise.}
	If $a_0,$ $a_1$ and $a_2$ have a point in common, prove that the
	elements defined in \ref{sec-tpap1} and their dual defined in 
	\ref{sec-tpap2} determine a self-dual configuration and the points
	$Ma_i$ and $\ov{M}a_i$ 
	coincide, as a set, with the points $M_i$ and $\ov{M}_i.$  
	If $p > 5,$ there are 29 points and 29 lines.\\
	If $p = 5,$ there are 25 points and 25 lines, the type is\\
	\hth$10 * 6 + 4 * 5 + 11 * 4 \:\&\:  10 * 6 + 4 * 5 + 11 * 4.$\\
	If $p = 7,$ it is of type\\
	\hth$12 * 6 + 8 * 5 + 9 * 4  \:\&\:  12 * 6 + 8 * 5 + 9 * 4.$\\
	If $p > 7,$ it is of type\\
	\hth$12 * 6 + 4 * 5 + 1 * 4 + 12 * 3  \:\&\:  12 * 6 + 4 * 5 + 1 * 4 
	+ 12 * 3.$\\
	The configuration is therefore distinct from the extended special
	Desargues configuration of \ref{sec-desd}.\\
	Prove that the 6 points and lines left over are also on a conic, as in
	\ref{sec-cesd}.

	\ssec{Duality.}\label{sec-Sdua}
	\setcounter{subsubsection}{-1}
	\sssec{Introduction.}
	This important concept, prepared by the work of
	Maurolycus and Poncelet, was introduced by Joseph Diaz Gergonne.  We
	observe that if we join Theorem \ref{sec-tdupap} to the axioms
	\ref{sec-ax1} and to Theorems \ref{sec-tinc}, and then exchange the
	words line and point, we
	obtain the same statements in some other order.  Therefore in any
	result obtained, we can exchange the words line and point.

	\sssec{Definition.}
	The method of obtaining from a result an other result
	by exchange of the words line and point is called {\em duality}.
	In particular, the Theorem of Desargues \ref{sec-tdes}, becomes:

	\sssec{Theorem. [Dual of Desargues' Theorem]}
	\label{sec-tdudes}
	{\em Given two triangles
	\{$a_0,$ $a_1,$ $a_2$\} and \{$b_0,$ $b_1,$ $b_2,$\} such that the
	points $a_0 \times b_0,$ $a_1 \times b_1$ and $a_2 \times b_2$ are on
	the same line $c.$  Let $c_0: = (a_1 \times a_2) \times (b_1 \times
	b_2),$ $c_1 := (a_2 \times a_0) \times (b_2 \times b_0),$
	$c_2 := (a_0 \times a_1) \times (b_0 \times b_1).$
	Then $c_0,$ $c_1,$ $c_2$ are incident to the same point $C.$} Fig. 3a)\\
	This is the dual of Theorem \ref{sec-tdes}.

	\sssec{Comment.}
	Fig. 1a and 1b are dual of each other, so are Fig. 2a and 2b,
	Fig. 2a' and 2b', Fig. 9 and 10.\\
	Fig. 3a, Fig. 3e, Fig. 3h are self dual.

	\ssec{Complete quadrangles and homologic quadrangles.}

	\sssec{Theorem.}\label{sec-tcqu}
	{\em If 2 quadrangles \{$A_0,$ $A_1,$ $A_2,$ $A_3$\} and 
	\{$A'_0,$ $A'_1,$ $A'_2,$ $A'_3$\} are
	such that none of their points and none of their lines coincide and are
	such that 5 of their corresponding lines are on the same line $p,$ then
	the 6-th pair of lines intersect on $p.$}

	Proof:  Using the notation \ref{sec-dcq} and ``$'$" for the second
	quadrangle,
	let $B_k := a_k \times a'_k$, $k$ = 0 to 5
	and let $B_0,$ $B_1,$ $B_2,$ $B_3,$ $B_4,$
	be all on the line $p.$  Theorem \ref{sec-tdudes}, dual of Desargues
	can be applied to the
	triangles $\{A_0, A_2, A_3\}$ and $\{A'_0, A'_2, A'_3\}$
	 then to the triangles $\{A_0, A_3, A_1\}$
	and $\{A'_0, A'_3, A'_1\}$.  The consequence is that the
	lines $A_0 \times A'_0,$
	$A_2 \times A'_2,$ $A_3 \times A'_3$ have a point $P$ in common which
	is also on $A_1 \times A'_1.$
	Therefore the Theorem of Desargues can be applied to the triangles
	$\{A_0, A_1, A_2\}$ and $\{A'_0, A'_1, A'_2\}$ which implies that
	the lines $a_5$ and $a'_5$ intersect on the line $p.$\\
	Or using the synthetic notation, let $b_j := A_j\times A'_j$, $j$ = 0
	to 3\\
	Desargues$^-1(p,\{a_3,a_2,a_1\},\{A_0,A_2,A_3\},
		\{a'_3,a'_2,a'_1\};\{A'_0,A'_2,A'_3\};
		\langle b_0,b_2,b_3\rangle,P)$,\\
	Desargues$^-1(p,\{a_4,a_0,a_2\},\{A_0,A_3,A_1\},
		\{a'_4,a'_0,a'_2\};\{A'_0,A'_3,A'_1\};
		\langle b_0,b_3,b_1\rangle,,Q)$,\\
	$\implies P = (A_0\times A'_0)\times (A_3\times A'_3) = Q,$\\
	$\implies$Desargues($P,\{A_0,A_1,A_2\},\{a_5,a_0,a_1\},
		\{A'_0,A'_1,A'_2\},\{a'_5,a'_0,a'_1\};
		\langle B_5,B_0,B_1\rangle,p$),

	\sssec{Definitions.}
	The quadrangles of Theorem \ref{sec-tcqu} are said to be
	{\em homologic.}
	$p$ is called the {\em axis} and $P$ the {\em center of the homology}.

	\sssec{Corollary.}\label{sec-ccqu}
	{\em If two complete quadrangles with no points and lines in
	common are such that\\
	\hth$	K := a_0 \times a_3$ is on $a'_0$ and $a'_3,$
	$L := a_1 \times a_4$ is on $a'_1$ and $a'_4,$\\
	\hth$	M := a_2 \times a'_2$ is on $K \times L,$\\
	then\\
	\hth$	N := a_5 \times a'_5$ is also on $K \times L.$}

	\sssec{Construction.}
	Given three points $K,$ $L,$ $M$ on a line $p,$
	choose arbitrarily a point $A_0$ not on $p$ and a point $A_1$ on
	$A_0 \times K$ distinct from $A_0$ and $K.$  Define\\
	\hth$    A_3 := (A_1 \times L) \times (A_0 \times M),$
	$A_2 := (A_3 \times K) \times (A_0 \times L),$\\
	\hth$    N  := (A_1 \times A_2) \times (K \times L).$ See Fig. 2a''.

	It follows from \ref{sec-ccqu} that $N$ is independent of the choice of
	$A_0$ and $A_1.$

	\sssec{Definition.}
	$N$ is called the {\em harmonic conjugate of $M$ with respect to
	$K$ and $L.$}

	\sssec{Theorem.}
	If each line has $q+1$ points on it, let $l(n,q)$ denote the
	number of points on a complete $n$-angle,
	let $l^*(n,q)$ denote the number of points not on a complete n-angle,
	let $L(n,q)$ denote the number of complete $n$-angles,\\
	0.\hti{7}$     l(n,q) = n\frac{n-1}{2} q - n(n-3)\frac{n^2-3n+6}{8}.$\\
	1.\hti{7}$l^*(n,q) = q^2 - (n+1)\frac{n-2}{2} q 
	+ (n-2)\frac{n^3-4n^2+7n-4}{8}.$\\
	2.\hti{7}$     L(n+1,q) = \frac{1}{n+1} L(n,q) l^*(n,q).$\\
	3.\hti{7}$     l(n+1,q)-l(n,q) = n q - \frac{1}{2}(n-1)(n^2-2n+2).$

	Proof.  $l(n+1,q)$ is obtained from $l(n,q)$ by adding points on each of
	the $n$ lines through the new point $A_n$ and through one of the old points
	A say, plus the new point itself.  On each of the lines $A_n \times A,$
	we have $q+1$ points from which we have to subtract the points $A$ and
	$A_n$
	as well as the points on the $(n-1)\frac{n-2}{2}$ lines through each
	pair of the old points, $A$ excluded.  This gives\\
	\hth$	n(q+1 - \frac{1}{2}(n-1)(n-2) -2) + 1 
	= nq - \frac{1}{2}(n^3-3n^3+4n-2).$\\
	Using $l(1,q) = 1,$ 0. follows by induction, 1. follows from
	$l^*(n,q) = q^2+q+1 - l(n,q),$ 2. follows from the fact that to each
	complete $n$-angle and each point not on its sides is associated a
	complete $(n+1)$-angle each being counted $n+1$ times.

	\sssec{Exercise.}
	\enumb
	\item $l(I,q)$ is a polynomial of degree 4, its successive
		forward differences at 0 are 1, $q-1,$ 0 and $-3.$
	\item $l^*(\frac{5}{2}+x,q) = l^*(\frac{5}{2}-x,q)$
	\item $l^*(n,I)$ is a quadratic function.  Its discriminant
		is $-\frac{1}{4}(n-2)(n-3)(n^2-5n+2),$ its successive
		forward differences at 0 are 4, $-5,$ 6 and $-6.$\\
		The discriminant is negative if $n > 4.$
	\enume

	\sssec{Table.}
	$\begin{array}{lllrl}
	n &      l(n,q)&  l^*(n,q)&discr.&L(n,q)\\
	0 &      0&  q^2+q+1  &-3&1\\
	1 &      1&  (q+1)q &1&q^2+q+1\\
	2 &      q+1&  q^2      &0&\frac{1}{2!}q(q+1)(q^2+q+1)\\
	3 &      3q&  (q-1)^2  &0&\frac{1}{3!}q^3(q+1)(q^2+q+1)\\
	4 &      6q-5&  (q-2)(q-3)    &1&\frac{1}{4!}q^3(q^2-1)(q^3-1)\\
	5 &      10q-20&  q^2-9q+21&-3&\frac{1}{5!}q^3(q^2-1)(q^3-1)(q-2)(q-3)\\
	6 &      15q-54\\
	7 &      21q-119\\
	8 &      28q-230
	\end{array}$

	\sssec{Exercise.}
	Complete the last 3 lines of the preceding table.
	\ssec{Collineation and Correlation.}\label{sec-Scollin}

	\sssec{Definition.}\label{sec-dcollin}
	A {\em collineation} consists of a one to one function $\gamma$  from
	the set of points of the plane onto itself, such that all points on a
	line have their image also on a line and of the induced function
	$\gamma'$ from the set of lines of the plane onto itself.

	\sssec{Definition.}\label{sec-dcorrel}
	A {\em correlation} consists of a one to one function $\rho$  from
	the set of lines of the plane onto itself, such that all lines through
	a point have their image also through a point and of the induced
	function $\rho'$ from the set of points of the plane onto itself.

	\sssec{Theorem.}
	{\em If the geometry is of prime order, a collineation or a
	correlation is determined by the image of a complete quadrangle onto a
	complete quadrangle or quadrilateral.}  (See \ref{sec-tcoll})

	\ssec{Finite projective planes for small $p$.}

	\setcounter{subsubsection}{-1}
	\sssec{Introduction.}
	There is a well known, see for instance
	Stevenson, p. 72, or Dembowski, p. 144, 14.  that there is, up to
	isomorphism, only one plane satisfying the incidence axioms, the
	axiom of Pappus  and the finite field axiom \ref{sec-ax2}.
	In the general case, the proof will require a full knowledge of the
	material not only of section 1, but also of the existence of
	fundamental projectivities of order $p-1$ and $p+1.$
	The axiom of Pappus is not required for $p \leq 7$, as proven by
	MacInnes in 1907 for $p = 2,$ 3 and 5.
	For $p = 7,$ see Bose and Nair, 1941, Hall 1953, 1954b, Pierce, 1953,
	Pickert, 1955.

	\sssec{Theorem.}\label{sec-tproj2}
	{\em For $p = 2,$}
	\enumb
	\item {\em There exists, up to isomorphism, only one plane satisfying
	the incidence axioms.}
	\item {\em The diagonal points of a complete quadrangle configuration
	are collinear.}
	\enume

	Proof:  Assume that line [3] contains the points (0), (1) and (2).
	Let (3) be an other point.  Define line [0] as the line through (1)
	and (3), we abbreviate this as $[0] := (1) \times (3).$  Similarly,
	$[1] := (0) \times (3),$ $[2] := (2) \times (3).$  Let the third point
	on [0] be (5),
	on [1] be (4) and on [2] be (6).  Let $[4] := (4) \times (6),$
	$[5] := (5) \times (6),$ $[6] := (4) \times (5).$
	The incidence properties imply (0) is on [5], which we abbreviate
	$(0) \cdot [5] = 0,$ similarly $(1) \cdot [4] = 0$ and
	$(2) \cdot [6] = 0.$ This completes the incidence tables:\\
	line : Points on line\hti{6} Point : lines through Point\\
	$\begin{array}{rrrrcrrrr}
	  0 :&  1&  3&  5&\hti{12}&   0 :&  1&  3&  5\\
	  1 :&  0&  3&  4&	&   1 :&  0&  3&  4\\
	  2 :&  2&  3&  6&	&   2 :&  2&  3&  6\\
	  3 :&  0&  1&  2&	&   3 :&  0&  1&  2\\
	  4 :&  1&  4&  6&	&   4 :&  1&  4&  6\\
	  5 :&  0&  5&  6& 	&   5 :&  0&  5&  6\\
	  6 :&  2&  4&  5&	&   6 :&  2&  4&  5
	\end{array}$

	\sssec{Theorem.}\label{sec-tproj3}
	{\em For $p = 3,$}
	\enumb
	\item {\em there exists, up to isomorphism, only one plane satisfying
	the incidence axioms.}
	\enume

	Proof:  Assume that the line [4] contains the points (0), (1), (2)
	and (3).  Let (4) be a point not on [4].  Let $[0] := (1) \times (4),$
	$[1] := (0) \times (4),$ $[2] := (3) \times (4),$
	$[3] := (2) \times (4).$
	Let (7) and (10) be the other points on [0].  Let $[7] := (0) \times
	(10),$ $[10] := (0) \times (7),$ $[9] := (2) \times (10),$ $[11] := (2)
	\times (7),$ $[8] := (3) \times (10),$ $[12] := (3) \times (7).$ 
	Let $(9) := [2] \times [10],$ $(11) := [2] \times [7],$
	$(8) := [3] \times [10],$ $(12) := [3] \times [7],$
	$(5) := [1] \times [9],$ $(6) := [1] \times
	[8].$  Let $[5] := (1) \times (9),$ $[6] := (1) \times (8).$

	At this stage with have the following incidence table:\\
	\hti{1}line : Points on line\hti{4} Point : lines through Point\\
	$\begin{array}{rrrrrcrrrrr}
	  0 :&  1&  4&  7& 10&  &   0 :&  1&  4&  7& 10\\
	  1 :&  0&  4&  5&  6&	&   1 :&  0&  4&  5&  6\\
	  2 :&  3&  4&  9& 11&	&   2 :&  3&  4&  9& 11\\
	  3 :&  2&  4&  8& 12&	&   3 :&  2&  4&  8& 12\\
	  4 :&  0&  1&  2&  3&	&   4 :&  0&  1&  2&  3\\
	  5 :&  1&&     9&& 	&   5 :&  1&  &   9\\
	  6 :&  1&&     8&&	&   6 :&  1&  &   8\\
	  7 :&  0& 10& 11& 12&&     7 :&  0& 10& 11& 12\\
	  8 :&  3&  6&  &  10&&     8 :&  3&  6&  &  10\\
	  9 :&  2&  5&  &  10&&     9 :&  2&  5&  &  10\\
	 10 :&  0&  7&  8&  9&&    10 :&  0&  7&  8&  9\\
	 11 :&  2&&    7     &&&   11 :&  2&  &   7\\
	 12 :&  3&&    7     &&&   12 :&  3&  &   7
	\end{array}$

	It remains to complete the table using the incidence axioms:\\
	Line [8] contains (3), (6) and (10), but (3) is already on line [2]
	with (4) hence (4) cannot be on [8].  Similarly (3) excludes
	(9), (11), (0), (1), (2), (6), (7); (6) excludes (5) and (10) excludes
	(12).  The only point left is (8).\\
	Line [9] contains (2), (5) and (10), (2) excludes (4), (8), (12), (0),
	(1), (3), (7) and (10) excludes (11), (3), (6), only (9) remains.\\
	Line [5] contains (1) and (9), (1) excludes (4), (7), (10), (0), (2),
	(3), (8) and (10) excludes (11), (5), only (6) and (12) remain.\\
	Line [6] contains (1) and (8), (1) excludes (4), (7), (10), (0), (2),
	(3), (6), (9), (12), only (5) and (11) remain.\\
	Line [11] contains (2) and (7), (2) excludes (4), (8), (12), (0), (1),
	(3), (5), (9), (10), only (6) and (11) remain.\\
	Line [12] contains (3) and (7), (3) excludes (4), (9), (11), (0), (1),
	(2), (6), (8), (10), only (5) and (12) remain.\\
	This completes the incidence tables:\\
	line : Points on line\hti{6} Point : lines through Point\\
	$\begin{array}{rrrrrcrrrrr}
	  5 :&  1&  6&  9& 12&\hti{4}&     5 :&  1&  6&  9& 12\\
	  6 :&  1&  5&  8& 11&&     6 :&  1&  5&  8& 11\\
	  8 :&  3&  6&  8& 10&&     8 :&  3&  6&  8& 10\\
	  9 :&  2&  5&  9& 10&&     9 :&  2&  5&  9& 10\\
	 11 :&  2&  6&  7& 11&&    11 :&  2&  6&  7& 11\\
	 12 :&  3&  5&  7& 12&&    12 :&  3&  5&  7& 12
	\end{array}$


	\sssec{Exercise.}
	 $D$ and $\ov{D}$, are harmonic conjugates to d and $\ov{d}$.
	Comes from Steiner for conics.

	\sssec{Exercise.}
	Let\\
	\ldots complete this, change notation for $BM$ in D2.\\
	D0.\hti{6}$    AM_i := a_{i+1} \times m_{i+1},$\\
	D2.\hti{6}$    BM_i := m_{i+1} \times \ov{b}_{i-1},$\\
	D3.\hti{6}$    ab_i := BM_{i+1} \times AM_{i-1},$\\
	then\\
	C0.\hti{6}$    P_i.ab_i = 0.$

	Proof:  $ab_2$ is the axis of perspectivity of the triangles
	$N_0,$
	with center of perspectivity $M_0.$%

	\sssec{Comment.}
	A configuration associated to antipolarity in 3
	dimensions implies a configuration of 20 points and 22 lines
	in 2 dimensional geometry, see VI.\ref{sec-tconf2022}.

INTEGRATE THERE\\
	Dx.\hti{5}$    ce_i := M\ov{M}a_{i+1} \times
	M\ov{M}a_{i-1},$\\
	\hth$
	 \ov{M}\ov{M}a_{i-1},$\\
	Dy.\hti{5}$     PQ_i := ce_i \times \ov{c}e_i,$\\
	Px.\hti{5}$ce_0 = [$m$_1($m$_1-$m$_2),$m$_2($m$_0-$m$_1),
		$m$_1($m$_2-$m$_0)],$\\
	\hth$\ov{c}e_0 = [$m$_2($m$_1-$m$_2),$m$_2($m$_0-$m$_1),
		$m$_1($m$_2-$m$_0)],$\\
	Py.\hti{5}$PQ_0 = (0,$m$_1($m$_2-$m$_0),-$m$_2($m$_0-$m$_1)).$

	\sssec{Examples.}

	In the following examples we can replace $\gamma$  and $\gamma'$ by
	$\rho$ and $\rho'$.  $\rho$  composed with $\rho'$ gives a collineation
	$\sigma$. Properties and special cases of collineations and correlations
	will be discussed in \ref{sec-Scollin} and \ref{sec-Scorrel}.
	In these examples, the complete quadrangle in the domain is
	always (0), (1), (6), (12).  $t(i)$ denotes the smallest positive
	integer such that $(\gamma^t)(i) = i.$  $C_3 = C_1^3$ indicates that the
	function $\gamma$ 
	of the collineation $C_3$ corresponds to the function $\gamma$  of the
	collineation $C_1$ composed with itself 3 times.  The examples will be
	used in 1.8.12.

	For $p = 5,$ 4 points? 0,1,6,12\\
	$C_0$ 4 image points? 1,6,12,3$\\
	\begin{array}{crrrrrrrrrrrrrrrrrrrrrrrrrrrrrrr}
    i&      0& 1& 2& 3& 4& 5& 6& 7& 8& 9&10&11&12&13&14&15\\
\gamma(i)&  1& 6&21&11&26&16&12&22& 7&17&27& 2& 3& 0& 4& 5\\
\gamma'(i)& 5&10&18&26&14&22& 0& 2& 1& 4& 3& 9&17&30&13&21\\
  t(i)&    31&31&31&31&31&31&31&31&31&31&31&31&31&31&31&31\\
    i&     &16&17&18&19&20&21&22&23&24&25&26&27&28&29&30\\
\gamma(i)& &30&20&10&25&15&18&13& 8&28&23&24&29& 9&14&19\\
\gamma'(i)&& 7&11&20&24&28& 8&29&25&16&12& 6&23&15&27&19\\
  t(i)&    &31&31&31&31&31&31&31&31&31&31&31&31&31&31&31
	\end{array}$

	$C_1$ 4 image points? 1,6,12,5$\\
	\begin{array}{crrrrrrrrrrrrrrrrrrrrrrrrrrrrrrr}
    i&      0& 1& 2& 3& 4& 5& 6& 7& 8& 9&10&11&12&13&14&15\\
\gamma(i)&  1& 6&16&26&11&21&12&27&17& 7&22& 2& 5& 4& 0& 3\\
\gamma'(i)& 5&10&22&14&26&18& 0& 3& 4& 1& 2& 9&21&13&30&17\\
  t(i)&    24&24&24& 6&24&24&24& 6&24&24&24&24&24& 1& 6&24\\
    i&     &16&17&18&19&20&21&22&23&24&25&26&27&28&29&30\\
\gamma(i)& &30&15&25&10&20&18&23&28& 8&13&24&19&14& 9&29\\
\gamma'(i)&& 7&28&24&20&11& 8&12&16&25&29& 6&19&27&15&23\\
  t(i)&    & 6&24&24&24&24&24&24& 6&24&24&24&24&24&24& 6
	\end{array}$

	$C_2=C_1^2$ 4 image points? 6,12,5,21$\\
	\begin{array}{crrrrrrrrrrrrrrrrrrrrrrrrrrrrrrr}
    i&      0& 1& 2& 3& 4& 5& 6& 7& 8& 9&10&11&12&13&14&15\\
\gamma(i)&  6&12&30&24& 2&18& 5&19&15&27&23&16&21&11& 1&26\\
\gamma'(i)&18& 2&12&30& 6&24& 5&14&26&10&22& 1& 8&13&23&28\\
  t(i)&    12&12&12& 3&12&12&12& 3&12&12&12&12&12& 1& 3&12\\
    i&     &16&17&18&19&20&21&22&23&24&25&26&27&28&29&30\\
\gamma(i)& &29& 3&13&22&20&25&28&14&17& 4& 8&10& 0& 7& 9\\
\gamma'(i)&& 3&27&25&11& 9& 4&21& 7&29&15& 0&20&19&17&16\\
  t(i)&    & 3&12&12&12&12&12&12& 3&12&12&12&12&12&12& 3
	\end{array}$

	$C_3=C_1^3$ 4 image points? 12,5,21,18$\\
	\begin{array}{crrrrrrrrrrrrrrrrrrrrrrrrrrrrrrr}
    i&      0& 1& 2& 3& 4& 5& 6& 7& 8& 9&10&11&12&13&14&15\\
\gamma(i)& 12& 5&29& 8&16&25&21&10& 3&19&28&30&18& 2& 6&24\\
\gamma'(i)&24&22&21&23& 0&25&18&30& 6& 2&12&10& 4&13&16&27\\
  t(i)&     8& 8& 8& 2& 8& 8& 8& 2& 8& 8& 8& 8& 8& 1& 2& 8\\
    i&     &16&17&18&19&20&21&22&23&24&25&26&27&28&29&30\\
\gamma(i)& & 9&26& 4&23&20&13&14& 0&15&11&17&22& 1&27& 7\\
\gamma'(i)&&14&19&29& 9& 1&26& 8& 3&15&17& 5&11&20&28& 7\\
  t(i)&    & 2& 8& 8& 8& 8& 8& 8& 2& 8& 8& 8& 8& 8& 8& 2
	\end{array}$

	$C_4=C_1^4$ 4 image points? 5,21,18,25$\\
	\begin{array}{crrrrrrrrrrrrrrrrrrrrrrrrrrrrrrr}
    i&      0& 1& 2& 3& 4& 5& 6& 7& 8& 9&10&11&12&13&14&15\\
\gamma(i)&  5&21& 9&17&30&13&18&22&26&10&14&29&25&16&12& 8\\
\gamma'(i)&25&12& 8&16& 5&29&24&23& 0&22&21& 2&26&13& 7&19\\
  t(i)&     6& 6& 6& 3& 6& 6& 6& 3& 6& 6& 6& 6& 6& 1& 3& 6\\
    i&     &16&17&18&19&20&21&22&23&24&25&26&27&28&29&30\\
\gamma(i)& & 7&24&11&28&20& 4& 0& 1& 3& 2&15&23& 6&19&27\\
\gamma'(i)&&30&20&15& 1&10& 6& 4&14&17&28&18& 9&11&27& 3\\
  t(i)&    & 3& 6& 6& 6& 6& 6& 6& 3& 6& 6& 6& 6& 6& 6& 3
	\end{array}$

	$C_5=C_1^6$ 4 image points? 18,25,13,4$\\
	\begin{array}{crrrrrrrrrrrrrrrrrrrrrrrrrrrrrrr}
    i&      0& 1& 2& 3& 4& 5& 6& 7& 8& 9&10&11&12&13&14&15\\
\gamma(i)& 18&25&27& 3& 9&11&13&28& 8&23& 1& 7& 4&29&21&15\\
\gamma'(i)&15& 8&26& 3&24&17&29& 7&18&21& 4&12& 0&13&14&11\\
  t(i)&     4& 4& 4& 1& 4& 4& 4& 1& 4& 4& 4& 4& 4& 1& 1& 4\\
    i&     &16&17&18&19&20&21&22&23&24&25&26&27&28&29&30\\
\gamma(i)& &19&17&16& 0&20& 2& 6&12&24&30&26&14& 5&22&10\\
\gamma'(i)&&16& 9&28& 2&22& 5& 6&23&27&19&25&10& 1&20&30\\
  t(i)&    & 1& 4& 4& 4& 4& 4& 4& 1& 4& 4& 4& 4& 4& 4& 1
	\end{array}$

	$C_6=C_1^8$ 4 image points? 13,4,11,2$\\
	\begin{array}{crrrrrrrrrrrrrrrrrrrrrrrrrrrrrrr}
    i&      0& 1& 2& 3& 4& 5& 6& 7& 8& 9&10&11&12&13&14&15\\
\gamma(i)& 13& 4&10&24&27&16&11& 0&15&14&12&19& 2& 7&25&26\\
\gamma'(i)&28&26& 0&30&29&27&17&14&25& 4& 6& 8&18&13&23& 1\\
  t(i)&     3& 3& 3& 3& 3& 3& 3& 3& 3& 3& 3& 3& 3& 1& 3& 3\\
    i&     &16&17&18&19&20&21&22&23&24&25&26&27&28&29&30\\
\gamma(i)& &22& 3&29& 6&20&30& 5&21&17& 9& 8& 1&18&28&23\\
\gamma'(i)&& 3&10&19&12&21&24& 5& 7&20&11&15&22& 2& 9&16\\
  t(i)&    & 3& 3& 3& 3& 3& 3& 3& 3& 3& 3& 3& 3& 3& 3& 3
	\end{array}$

	$C_7=C_1^{12}$ 4 image points? 16,30,29,9$\\
	\begin{array}{crrrrrrrrrrrrrrrrrrrrrrrrrrrrrrr}
    i&      0& 1& 2& 3& 4& 5& 6& 7& 8& 9&10&11&12&13&14&15\\
\gamma(i)& 16&30&14& 3&23& 7&29& 5& 8&12&25&28& 9&22& 2&15\\
\gamma'(i)&11&18&25& 3&27& 9&20& 7&28& 5&24& 0&15&13&14&12\\
  t(i)&     2& 2& 2& 1& 2& 2& 2& 1& 2& 2& 2& 2& 2& 1& 1& 2\\
    i&     &16&17&18&19&20&21&22&23&24&25&26&27&28&29&30\\
\gamma(i)& & 0&17&19&18&20&27&13& 4&24&10&26&21&11& 6& 1\\
\gamma'(i)&&16&21& 1&26& 6&17&29&23&10& 2&19& 4& 8&22&30\\
  t(i)&    & 1& 2& 2& 2& 2& 2& 2& 1& 2& 2& 2& 2& 2& 2& 1
	\end{array}$

	$C_8$ 4 image points? 0,1,12,19$\\
	\begin{array}{crrrrrrrrrrrrrrrrrrrrrrrrrrrrrrr}
    i&      0& 1& 2& 3& 4& 5& 6& 7& 8& 9&10&11&12&13&14&15\\
\gamma(i)&  0& 1& 3& 5& 2& 4&12&14&11&13&15&17&19&16&18&20\\
\gamma'(i)&10&26&14&18& 5&22& 6& 7& 0& 8& 9& 1&25&30&15&20\\
  t(i)&     4& 5&20&20&20&20& 1& 1& 4& 4& 4& 5&20&20&20&20\\
    i&     &16&17&18&19&20&21&22&23&24&25&26&27&28&29&30\\
\gamma(i)& &22&24&21&23&25&27&29&26&28&30& 7& 9& 6& 8&10\\
\gamma'(i)&&21&19&28& 3&12&11& 2&29&17&23&16&13&27&24& 4\\
  t(i)&    & 5&20&20&20&20& 5&20&20&20&20& 5&20&20&20&20
	\end{array}$

	$C_9=C_8^2$ 4 image points? 0,1,19,23$\\
	\begin{array}{crrrrrrrrrrrrrrrrrrrrrrrrrrrrrrr}
    i&      0& 1& 2& 3& 4& 5& 6& 7& 8& 9&10&11&12&13&14&15\\
\gamma(i)&  0& 1& 5& 4& 3& 2&19&18&17&16&20&24&23&22&21&25\\
\gamma'(i)& 9&16&15&28&22& 2& 6& 7&10& 0& 8&26&23& 4&20&12\\
  t(i)&     2& 5&10&10&10&10& 1& 1& 2& 2& 2& 5&10&10&10&10\\
    i&     &16&17&18&19&20&21&22&23&24&25&26&27&28&29&30\\
\gamma(i)& &29&28&27&26&30& 9& 8& 7& 6&10&14&13&12&11&15\\
\gamma'(i)&&11& 3&27&18&25& 1&14&24&19&29&21&30&13&17& 5\\
  t(i)&    & 5&10&10&10&10& 5&10&10&10&10& 5&10&10&10&10
	\end{array}$

	$C_{10}=C_8^4$ 4 image points? 0,1,26,7$\\
	\begin{array}{crrrrrrrrrrrrrrrrrrrrrrrrrrrrrrr}
    i&      0& 1& 2& 3& 4& 5& 6& 7& 8& 9&10&11&12&13&14&15\\
\gamma(i)&  0& 1& 2& 3& 4& 5&26&27&28&29&30& 6& 7& 8& 9&10\\
\gamma'(i)& 0&11&12&13&14&15& 6& 7& 8& 9&10&21&24&22&25&23\\
  t(i)&     1& 5& 5& 5& 5& 5& 1& 1& 1& 1& 1& 5& 5& 5& 5& 5\\
    i&     &16&17&18&19&20&21&22&23&24&25&26&27&28&29&30\\
\gamma(i)& &11&12&13&14&15&16&17&18&19&20&21&22&23&24&25\\
\gamma'(i)&&26&28&30&27&29&16&20&19&18&17& 1& 5& 4& 3& 2\\
  t(i)&    & 5& 5& 5& 5& 5& 5& 5& 5& 5& 5& 5& 5& 5& 5& 5
	\end{array}$

	$C_{11}$ 4 image points? 0,1,26,14$\\
	\begin{array}{crrrrrrrrrrrrrrrrrrrrrrrrrrrrrrr}
    i&      0& 1& 2& 3& 4& 5& 6& 7& 8& 9&10&11&12&13&14&15\\
\gamma(i)&  0& 1& 5& 4& 3& 2&26&29&27&30&28&11&14&12&15&13\\
\gamma'(i)& 0&11&15&14&13&12& 6& 8&10& 7& 9&16&19&17&20&18\\
  t(i)&     1& 4& 4& 4& 4& 4& 1& 4& 4& 4& 4& 4& 4& 4& 4& 4\\
    i&     &16&17&18&19&20&21&22&23&24&25&26&27&28&29&30\\
\gamma(i)& &21&24&22&25&23& 6& 9& 7&10& 8&16&19&17&20&18\\
\gamma'(i)&&21&22&23&24&25& 1& 4& 2& 5& 3&26&30&29&28&27\\
  t(i)&    & 4& 4& 4& 4& 4& 4& 4& 4& 4& 4& 1& 2& 2& 2& 2
	\end{array}$

	$C_{12}$ 4 image points? 0,6,12,23$\\
	\begin{array}{crrrrrrrrrrrrrrrrrrrrrrrrrrrrrrr}
   i&       0& 1& 2& 3& 4& 5& 6& 7& 8& 9&10&11&12&13&14&15\\
\gamma(i)&  0& 6& 9& 7&10& 8&12&15&13&11&14&24&23&22&21&25\\
\gamma'(i)& 5&26&18&10&22&14& 1& 4& 3& 0& 2& 6&23&15&27&19\\
  t(i)&     5& 5& 5& 5& 5& 5& 5& 5& 5& 5& 5& 5& 5& 5& 5& 5\\
   i&      &16&17&18&19&20&21&22&23&24&25&26&27&28&29&30\\
\gamma(i)& &18&19&20&16&17&30&27&29&26&28& 2& 5& 3& 1& 4\\
\gamma'(i)&&16&12& 8&29&25&11&20&24&28& 7&21& 9&17&30&13\\
  t(i)&    & 1& 5& 5& 5& 5& 5& 5& 5& 5& 5& 5& 5& 5& 5& 5
	\end{array}$

	$C_{13}$ 4 image points? 11,7,2,12$\\
	\begin{array}{crrrrrrrrrrrrrrrrrrrrrrrrrrrrrrr}
    i&      0& 1& 2& 3& 4& 5& 6& 7& 8& 9&10&11&12&13&14&15\\
\gamma(i)& 11& 7&24&28&20& 5& 2&17&29&23&10&13&12&15&14& 0\\
\gamma'(i)&15&27& 6&19&23& 5&30& 1&25&20&10& 0&12&11&14&13\\
  t(i)&     4& 4& 4& 4& 4& 1& 4& 4& 4& 4& 1& 4& 1& 4& 1& 4\\
    i&     &16&17&18&19&20&21&22&23&24&25&26&27&28&29&30\\
\gamma(i)& &25&27&18& 3& 9&19&22& 4&30& 8&26& 1&21&16& 6\\
\gamma'(i)&&29& 7&18&21& 4&28&22& 9& 2&16&26&17& 3& 8&24\\
  t(i)&    & 4& 4& 1& 4& 4& 4& 1& 4& 4& 4& 1& 4& 4& 4& 4
	\end{array}$

	$C_{14}$ 4 image points? 10,25,11,0$\\
	\begin{array}{crrrrrrrrrrrrrrrrrrrrrrrrrrrrrrr}
    i&      0& 1& 2& 3& 4& 5& 6& 7& 8& 9&10&11&12&13&14&15\\
\gamma(i)& 10&25&15& 1&30&20&11&23&17&29& 2& 9& 0& 6& 8& 7\\
\gamma'(i)&29&27&30& 0&28&26& 7&15& 4&18&21&17&23&10&11& 2\\
  t(i)&    31&31&31&31&31&31&31&31&31&31&31&31&31&31&31&31\\
    i&     &16&17&18&19&20&21&22&23&24&25&26&27&28&29&30\\
\gamma(i)& &18&22&14& 5&26& 3&21&28&12&19&27&24& 4&16&13\\
\gamma'(i)&&12& 5&16&25& 8&22& 6&13& 3&20& 1&19&24& 9&14\\
  t(i)&    &31&31&31&31&31&31&31&31&31&31&31&31&31&31&31
	\end{array}$

	$C_{15}$ 4 image points? 13,17,24,12$\\
	\begin{array}{crrrrrrrrrrrrrrrrrrrrrrrrrrrrrrr}
    i&      0& 1& 2& 3& 4& 5& 6& 7& 8& 9&10&11&12&13&14&15\\
\gamma(i)& 13&17&30&21& 9& 5&24&27&16& 4&10&15&12& 0&14&11\\
\gamma'(i)&13&17&30&21& 9& 5&24&27&16& 4&10&15&12& 0&14&11\\
  t(i)&     0& 1& 2& 3& 4& 5& 6& 7& 8& 9&10&11&12&13&14&15\\
    i&     &16&17&18&19&20&21&22&23&24&25&26&27&28&29&30\\
\gamma(i)& & 8& 1&18&28&23& 3&22&20& 6&29&26& 7&19&25& 2\\
\gamma'(i)&& 8& 1&18&28&23& 3&22&20& 6&29&26& 7&19&25& 2\\
  t(i)&    &16&17&18&19&20&21&22&23&24&25&26&27&28&29&30
	\end{array}$

\setcounter{subsection}{89}
	\ssec{Answer to exercises.}

	\sssec{Answer to}
	\vspace{-18pt}\hspace{94pt}{\bf \ref{sec-tinc}.}\\
	First, prove that there are exactly $p+1$ lines through $P,$ more
	generally through any point not on $l.$
	Then, prove that on any line $m$ distinct from $l$ and not incident to
	both $P$ and $Q,$ there are exactly $p+1$ points.  If $Q$ is not on the
	line join $Q$ to all the points on $l$ and determine the intersection
	with $m.$ Then, for $P \times Q$ determine a point on an other line
	through $P$ which is not on $l$ and repeat the argument just given.
	To count the points, observe that any point different from $P$ is on
	a line through $P.$  There are exactly $p+1$ such lines and on each
	there are $p$ points distinct from $P$ hence altogether 
	$(p+1)p + 1$ points.

	\sssec{Answer to}
	\vspace{-18pt}\hspace{94pt}{\bf \ref{sec-epph}.}\\
	Given the hexagon $\{A0, A1, A2, A3, A4, A5,\}$ such that
	the alternate vertices $A0,$ $A2$ and $A4$ are collinear.  The necessary
	and sufficient condition for $A1,$ $A3$ and $A5$ to be collinear is
	that the points $P0,$ $P1$ and $P2$ be collinear.  The necessary
	condition
	follows using 1.5.1. on the hexagon $\{A0, P0, A2, P1, A1, P2\}$.

	\sssec{Answer to}
	\vspace{-18pt}\hspace{94pt}{\bf \ref{sec-eqqc}.}\\
	The construction is\\
	\hth$r_0 := P \times Q_0,$ $r_1 := P \times Q_1,$
	$r_2 := P \times Q_2,$\\
	\hth$p_0 := Q_1 \times Q_2,$ $p_1 := Q_2 \times Q_0,$
	$p_2 := Q_0 \times Q_1,$\\
	\hth$A_0 := p_0 \times r_0,$ $A_1 := p_1 \times r_1,$
	$A_2 := p_2 \times r_2,$\\
	\hth$a_0 := A_1 \times A_2,$ $a_1 := A_2 \times A_0,$
	$a_2 := A_0 \times A_1,$\\
	\hth$P_0 := a_0 \times r_0,$ $P_1 := a_1 \times r_1,$
	$P_2 := a_2 \times r_2,$\\
	\hth$q_0 := P_1 \times P_2,$ $q_1 := P_2 \times P_0,$
	$q_2 := P_0 \times P_1,$\\
	\hth$R_0 := a_0 \times q_0,$ $R_1 := a_1 \times q_1,$
	$R_2 := a_2 \times q_2,$\\
	\hth$p := R_1 \times R_2.$\\
	We have to prove\\
	\hth$R_0$ is on $p$ and $R_i$ is on $p_i.$\\
	 $\ldots$ ..\\
	This gives the configuration\\
	\hth$p$ on $R_i;$ $a_i$ on $P_i, R_i, A_{i+1}, A_{i-1};$
	$p_i$ on $A_i, R_i, Q_{i+1}, Q_{i-1},$\\
	\hth$q_i$ on $R_i, P_{i+1}, P_{i-1};$ $r_i$ on $P, A_i, P_i, Q_i.$\\
	Similarly for lower case and upper case exchanged.\\
	If $p = 3,$ $q_i$ and $p$ must contain a fourth point which is one of
	the 13 known point.  By necessity $P$ is on $p$ and $Q_i$ is on $q_i.$
	See \ref{sec-tproj3}

	\sssec{Answer to}
	\vspace{-18pt}\hspace{94pt}{\bf \ref{sec-treid} and \ref{sec-ereid}.}\\
	I will not repeat the computations of Chapter III.\\
	$b_{01} = b_{03}$ if $\frac{m_0}{m_2} = \frac{-m_1}{-m_0},$ therefore if
	$\overline{M}$ is on the conic $X_0^2 - X_1X_2 = 0,$ which is
	represented by the matrix
	$\matb{|}{ccc}-2&0&0\\0&0&-1\\0&-1&0\mate{|}$.\\
	If $B_{22} = (a,m_1,m_0),$ then $b_{12} = [m_0,0,-a],$
	$b_{22} = [m,-a,0],$ $B_{23} = (a,m_0,m_0),$ $B_{32} = (a,m_1,m_2),$
	$b_{23} = [m_0,-a,0],$ $b_{13} = [m2,0,-a],$ $B_{33} = (a,m_0,m_2),$
	hence $B_{33} \cdot b_{03},$ hence \ref{sec-treid}.\\
	$B_{22} \times B_{33} = [m1m2-m_0^2,-a(m_2-m_0),a(M-0-m_1)],$ line
	incident to $(0,m_0-m_1,m_2-m_0) = AEul_0,$ hence \ref{sec-ereid}.0.
	1, and 2, follow from Chapter III.


	\ssec{Relation between Synthetic and Algebraic Finite
	Projective Geometry.}
\setcounter{subsubsection}{-1}
	\sssec{Introduction.}%
	I start with affine geometry by choosing a particular line $m$ as the
	ideal
	line and 2 ordinary lines $x$ and $y$ which intersect at $O$, as well as
	an ordinary point $M$ on neither $x$ nor $y$.\\
	I will first associate to ordinary points on the line integers
	from 0 to $p-1,$ by defining the successor.  I will then define
	addition of points on the line and prove commutativity
	using the axiom of Pappus. I will then define
	multiplication of points on the line and prove commutativity
	using the axiom of Pappus.  It remains to prove the distributivity law.

	\sssec{Definition.}
	Let $Y := x \times m$, $M_1 := (Y \times M) \times y.$\\
	$A_0 := O,$\\
	$A_{i+1} := (((A_i \times M_1) \times m) \times M) \times x,$ for
	$i = 0,$ 1, \ldots until $A_n = A_0.$\\
	$A_{i+1}$ is called the {\em successor} of $A_0$.

	\sssec{Theorem.}
	\hth$n = p$.

	Proof: The parabolic projectivity associates to $A_0,$ $A_1$, $\sigma$
	with fixed point $Y$ which
	associates to $A_0,$ $A_1$, associates to $A_i,$ $A_{i+1}$.
	By definition $\sigma^n = \epsilon$, the identity mapping.
	If $n < p$, any other ordinary point on the line distinct from $A_i$
	has therefore the same period $n$, after exausting all points
	in the line it follws that $n$ must divide $p$, therefore $n = p$.
	We could also give a group theory proof of this Theorem and use the
	Theorem of Lagrange.

	\sssec{Definition.}
	Given 2 points $A$ and $B$ on $x$, the {\em addition} of the 2
	points, $C := A + B$ is defined as follows,\\
	\ldots..

	\sssec{Theorem.}
	{\em The addition is commutative, in other words, $A + B = B + A$.}

	Let $u := Y \times M$, $A_0 := (M_2 \times B) \times m,$
	$A_1 := (M_1 \times A) \times m,$\\
	$C_0 := (A_1 \times M_1) \times (X \times B_1),$
	$C_1 :=(A_0 \times M_1) \times (X \times B_0).$\\
	The axiom of Pappus applied to the points $A_0,$ $A_1,$ $X$ on $m$
	and $B_0,$ $B_1,$ $M_1$ on $u$ implies that $A_0 \times B_1,$
	$B_0 \times A_1$ intersect on the line $A \times B,$ therefore $C = D$
	or $A + B = B + A.$

	\sssec{Definition.}
	Given 2 points $A$ and $B$ on $x$, the {\em multiplication} of the 2
	points, $C := A \cdot B$ is defined as follows,\\
	Let $X := y \times m,$ $M_0 := (M \times X) \times x,$
	$Z := (M_0 \times M_1) \times m,$\\
	$B' := (B \times Z) \times y,$ $A'' := (A \times M_1) \times m,$
	$C := (A'' \times B') \times x.$

	\sssec{Theorem.}
	{\em The multiplication is commutative, in other words,
	$A \cdot B = B \cdot A$.}

	Proof: We have by definition\\
	$A' := (A \times Z) \times y,$ $B'' := (B \times M_1) \times m,$
	$D := (B'' \times A') \times x,$
	The axiom of Pappus applied to the points $A',$ $B',$ $M'$ on $y$
	and $A'',$ $B'',$ $Z$ on $m$ implies that $A'' \times B',$
	$B'' \times A'$ intersect on the line $A \times B,$ therefore $C = D$ or
	$A \cdot B = B \cdot A.$

	\sssec{Theorem.}
	{\em The distributive law applies, in other words\\
	\hth$ A \cdot (B + C) = (A \cdot B) + (A \cdot C).$}

	Proof: Let $Z$ be a point on $m$ distinct from $X$ and $Y$.\\
	Let $B' := (B \times Z) \times y,$ $C' := (C \times Z) \times y,$
	$BpC' := (BpC \times Z) \times y,$\\
	Let $A''$ be the direction of $M_1 \times A$ or $(M_1 \times A) \times
	m$.\\
	$B \times AtB$, $B \times AtB$ and $B \times AtB$ have the same
	direction therefore, if $U := (B \times Y) \times (ZtC \times X,$ then
	$O B = AtC U$, therefore $AtC AtB+C = O AtB$, therefore\\
	\ldots..

{\tiny\footnotetext[1]{G21.TEX [MPAP], \today}}

	\section{Algebraic Model of Finite Projective Geometry.}

	\setcounter{subsection}{-1}
	\ssec{Introduction.}

	In the general descriptions, I will from time to time give, between
	braces, information to the reader with advanced knowledge.
	This information is not required for the reader without prior
	knowledge, and may be explained in later sections or not. In the next
	paragraph, there are several examples of such use of braces.
	To construct finite Euclidean geometries I will use a model which
	depends on the \{field of\} integers modulo $p.$  The properties of the
	integers ${\Bbb Z}$ are assumed.  The model will be constructed in 4
	steps.\\
	In the first step, described in this section, I will not distinguish
	between points and directions, and use the well known algebraization of
	the finite projective plane \{associated with a Galois field,
	corresponding to the prime $p$\}, see also I.3.\\
	In the second step, (Section 8) I will introduce an ideal line \{which
	plays the role of line at infinity in the Euclidean plane\}, the notion
	of parallelism and of mid-points.\\
	In the third step (III.1), I will introduce the notion of
	perpendicularity \{associated to an involution on the ideal line\}.
	\{All these steps are valid in any field.\}\\
	In the fourth step (III.2 and 3), I will introduce measure of angles
	and of distances, together with a finite trigonometry.

	\ssec{Representation of points, lines and incidence.}

	\sssec{Definition.}
	A point is {\em represented by an ordered triple of integers}
	modulo $p,$ {\em placed between parenthesis}.  Not all 3 integers can be
	simultaneously 0.  Two triples are equivalent  iff  one of them can be
	derived from the other using multiplication, modulo $p,$ by an integer
	which is not zero modulo $p.$

	\sssec{Example.}
	If p = 3, there are 13 points:\\
	\hth (0,0,1), (0,1,0), (0,1,1), (0,1,2), (1,0,0), (1,0,1),\\
	\hth (1,0,2), (1,1,0), (1,1,1), (1,1,2), (1,2,0), (1,2,1), (1,2,2).\\
	(2,2,0) is the same as (1,1,0), (2,1,2) is the same as (1,2,1).

	\sssec{Convention.}\label{sec-crep}
	When I will compute numerically, I will always choose the
	representation of triples in such a way that the first non zero
	integer in the triple is 1.  This representation will be called the
	{\em normal representation}.  When I perform algebraic manipulations,
	I will multiply by the most convenient expression, to simplify the
	components or, if appropriate, to make the symmetry evident.

	\sssec{Notation.}\label{sec-npt}
	A more compact notation for the triples is to use a
	single integer, as follows,\\
	\hth$	(0)$  for  $(0,0,1),$\\
	\hth$	(i+1)$  for  $(0,1,i),$ $0 \leq i < p,$\\
	\hth$	((i+1)p + j + 1)$  for  $(1,i,j),$ $0 \leq i,j < p.$\\
	When there is no ambiguity, I will often drop the parenthesis.

	\sssec{Exercise.}
	Justify the Notation \ref{sec-npt} and therefore check that there
	are $p^2 + p +1$ points in the projective geometry associated to $p.$

	\sssec{Definition.}\label{sec-dpt}
	A line is {\em represented by an ordered triple of integers,}
	modulo $p,$ {\em placed between brackets}.  Again, not all 3 integers
	can be
	simultaneously equal to 0, and 2 triples which can be obtained from
	each other by multiplication of each integer by the same non zero
	integer modulo $p$ are considered equal.

	\sssec{Notation.}
	The notation [0] for [0,0,1], \ldots, similar to \ref{sec-npt}
	will be used for lines.  I will, also drop the bracket around the
	single integer, if there is no ambiguity.

	\sssec{Definition.}
	The point $P = (P_0,P_1,P_2)$ {\em and the line} $l = [l_0,l_1,l_2]$
	{\em are incident}, or $P$ {\em is on} $l$ or $l$ {\em goes through}
	$P$, iff\\
	\hth$P \cdot l := P_0.l_0 + P_1.l_1 + P_2.l_2 = 0 \pmod{p}).$\\
	$P$ {\em and} $l$ {\em are not incident}  iff  $P \cdot l\neq  0.$

	\sssec{Example.}
	For $p = 5,$ (1,0,1) is on [1,2,4], (1,2,3) is on [1,4,2].\\
	The points (5) = (0,1,4), (10) = (1,0,4), (14) = (1,1,3),
	(18) = (1,2,2), (22) = (1,3,1), (26) = (1,4,0) are the 6 points on the
	line [12] = [1,1,1].

	\ssec{Line through 2 points and point through 2 lines.}

	\sssec{Definition.}
	I recall the definition of the {\em cross product of 2 three
	dimensional vectors}.\\
	\hth$X * Y := (X_0, X_1, X_2) * (Y_0, Y_1, Y_2) :=\\
	\hti{16}( X_1Y_2 - X_2Y_1, X_2Y_0 - X_0Y_2, X_0Y_1 - X_1Y_0)$

	\sssec{Notation.}\label{sec-ntimes}
	When I use the cross product of 2 vectors and then
	normalize using the convention \ref{sec-crep}, I will use the symbol
	$``\times"$,
	which recalls the symbol ``$*$", instead of that symbol.  The result is
	unique, if I compute numerically.  It is not unique, if I proceed
	algebraically.  In this case, equality implies that an appropriate
	scaling has been used on either side of the equation or on both sides.
	See Chapter V, for some examples.

	\sssec{Theorem.}
	{\em
	\hth$P \times Q$ is the line through the distinct points $P$ and $Q.$\\
	\hth$p \times q$ is the point on the distinct lines $p$ and $q.$}

	This follows from $(P\times Q)\cdot P = (P\times Q)\cdot Q = 0$ or
	$(p\times q)\cdot p = (p\times q)\cdot q = 0$.

	\sssec{Example.}
	For $p = 5,$ $[1,1,3] := (1,2,4) \times (1,3,2)$ and
	$(1,1,3) := [1,2,4] \times [1,3,2].$

	\sssec{Theorem.}
	\enumb
	\item $k_1 A * B - $k$_2 C * D$ {\em is a line incident
	to $A \times B$ and $C \times D.$}
	\item {\em The lines} k$_0 A * B - $k$_1 C * D,$
		k$_1 C * D - $k$_2 E * F$ {\em and}
		k$_2 E * F - $k$_0 A * B$,
		{\em are incident.}
	\enume

	\sssec{Example.}\label{sec-elinecol}
	For any $p$, let $A$ = (0,1,1), $B$ = (1,2,1), $C$ = (2,1,1),
	$D$ = (1,3,1), $E$ = (2,4,1), $F$ = (4,3,1), then $A * B = [-1,1,-1],$
	$C * D = [-2,-1,5],$ $E * F = [1,2,-10].$ If k$_0$ = k$_1$ = k$_2$ = 1,
	we obtain the lines $[1,2,-6],$ $[-3,-3,15],$ $[2,1,-9]$ incident to
	(4,1,1).

	\sssec{Comment.}
	The algebraic method allows the representation of points or lines by a
	single  symbol.  This method which was well used in 19-th Century text,
	see for instance Salmon, 1879, Chapter XIV, has somehow fallen in
	disfavor.

	\ssec{The model satisfies the axioms of the projective Pappus
	plane of order $p$.}

	\setcounter{subsubsection}{-1}
	\sssec{Introduction.}
	After proving that the algebraic model satisfies the
	axioms of finite projective geometry, I give construction of points
	on a line whose coordinates have a simple algebraic relationship.
	These could be used as a tool for the construction of points whose
	coordinates are known in terms of points contructed earlier.  The
	notation $O + k M$ used in Theorem \ref{sec-tbaker} is partially
	justified in section \ref{sec-Sveccalc}.

	\sssec{Theorem.}
	{\em Each line $l$ contains exactly $p+1$ points, each point $P$ is on
	exactly $p+1$ lines, therefore\\
	The model satisfies axiom \ref{sec-ax1}.3 and its dual.}

	Proof:  We want to find the points $(x,y,z)$ on the line $[a,b,c].$
	At least one of the 3 integers,
	$a,$ $b$ or $c$ is different from 0, let it be $c,$ in this case $x$
	and $y$ cannot both be 0.  Given $x$ and $y$ we can solve
	$ax + by + cz = 0$ for the integer
	$z,$ using the algorithm of Euclid-Aryabatha, $z := -(a\:x + b\:y)/c.$
	(See I.\ref{sec-aaryab})\\
	If $x = 1,$ to each value of $y$ from 0 to $p-1$ corresponds a value
	of $z,$ namely $-(a + b\:y)/c.$\\
	If $x = 0$ and $y = 1,$ we obtain one value of $z,$ namely $-b/c.$\\
	Therefore we obtain altogether $p+1$ points.\\
	Exchanging brackets and parenthesis gives the dual property.

	\sssec{Theorem.}
	{\em The model satisfies the axiom} \ref{sec-ax1}.4. {\em of Pappus.}

	Proof:  I will give the proof in the special case in which the lines
	are [0] = [0,0,1] and [1] = [0,1,0].\\
	The general case can be deduced
	from general considerations on projectivity or can be proven directly.
	This direct proof is left as an exercise.\\
	We choose $A_0 = (1,$a$_0,0),$ $A_1 = (1,$a$_1,0),$
	$A_2 = (1,$a$_2,0)$ and $B_0 = (1,0,$b$_0),$ $B_1 = (1,0,$b$_1),$
	$B_2 = (1,0,$b$_2),$
	with a$_0 $a$_1 $a$_2 $b$_0 $b$_1 $b$_2\neq  0.$  Then\\
	\hth$C_0 = ( $a$_2 $b$_2 - $a$_1 $b$_1,$
	a$_1 $a$_2 ($b$_2-$b$_1),$ b$_1 $b$_2 ($a$_2-$a$_1) ),$\\
	\hth$C_1 = ( $a$_0 $b$_0 - $a$_2 $b$_2,$ a$_2 $a$_0 ($b$_0-$b$_2),$
		b$_2 $b$_0 ($a$_0-$a$_2) ),$\\
	\hth$C_2 = ( $a$_1 $b$_1 - $a$_0 $b$_0,$ a$_0 $a$_1 ($b$_1-$b$_0),$
		b$_0 $b$_1 ($a$_1-$a$_0) ).$\\
	It is easy to verify that a$_0 $b$_0 C_0 + $a$_1 $b$_1 C_1
	+ $a$_2 $b$_2 C_2 = 0,$
	therefore the points $C_0,$ $C_1$ and $C_2$ are collinear as will be
	seen shortly, in \ref{sec-tvecc}.\\
	The special cases, where $A_2 = (0,1,0)$ or $B_2 = (0,0,1)$ or a$_0$ or
	b$_0 = 0,$ can be verified easily.

	\sssec{Theorem.}
	{\em The algebraic model satisfies the axioms} \ref{sec-ax1} {\em of
	finite projective geometry and therefore it can be used to prove all the
	theorems of finite projective geometry.}

	\sssec{Definition.}\label{sec-dpapl}
	Given a triangle $\{a_0,$ $a_1,$ $a_2\},$ and 2 arbitrary
	lines, $x$ and $y,$ the {\em Pappus line of} $x$ {\em and} $y,$ is the
	line $z$ associated
	to the application of the axiom of Pappus to the intersection with $x$
	and $y$ of the lines $a_0,$ $a_1$ and $a_2.$ 

	\sssec{Theorem.}
	{\em If $a_0 = [1,0,0],$ $a_1 = [0,1,0],$ $a_2 = [0,0,1],$
	if $x = [$x$_0,$x$_1,$x$_2]$ and $y = [$y$_0,$y$_1,$y$_2]$ then the
	Pappus line of $x$ and $y$ is}\\
	\hth$  z = [$x$_0$y$_0 ($x$_1$y$_2+$x$_2$y$_1),
		$x$_1$y$_1 ($x$_2$y$_0+$x$_0$y$_2),
		$x$_2$y$_2 ($x$_0$y$_1+$x$_1$y$_0)].$

	The proof is left as an exercise.  Hint: One of the 3 points on $z$ is\\
	        $($x$_0^2$y$_1$y$_2-$x$_1$x$_2$y$_0^2,
	$x$_0$y$_0($x$_2$y$_0-$x$_0$y$_2),
	$x$_0$y$_0($x$_1$y$_0-$x$_0$y$_1)).$

	\sssec{Comment.}
	Definition \ref{sec-dpapl} may be new, it was suggested by one of the
	construction in a triangle of the point of Lemoine from the barycenter
	and orthocenter.  See 4.2.12.
	The operation of deriving $z$ from $x$ and $y$ is commutative but is not
	associative.

	\sssec{Exercise.}
	Verify that if $p = 2$ and $p = 4$ the diagonal points of a
	complete quadrangle are collinear.  For $p = 2,$ choose one such
	quadrangle.  For $p = 4,$ the coordinates of the points are 
	$u + v \xi$ ,
	where $u, v \in {\Bbb Z}_2$ and $\xi^2 + \xi + 1 = 0.$
	Choose a quadrangle which is not in the subspace $v = 0.$

	\sssec{Exercise.}\label{sec-edesar}
	Prove algebraically the 2 cases of the Theorem of Desargues.

	\sssec{Definition.}
	Let the coordinates of the distinct points $O$ = (o$_0$,o$_1$,o$_2$)
	and $M$ = (m$_0$,m$_1$,m$_2$) be normalized, $O + $x$ M$ is the point on
	$O\times M$ whose coordinates are
	(o$_0$+x m$_0$,o$_1$+x m$_1$,o$_2$+x m$_2$).

	\sssec{Comment.}
	The following Theorem relates specific constructions in projective
	geometry to algebraic operations, nothing is claimed as to the
	projective properties of the operation ``+", these will require
	the introduction of preferences associated with affine and Euclidean
	geometries.  In this Theorem we have not used the notation ``$*$"
	which appears in section \ref{sec-Sveccalc}, this the reason why
	the notation $O + k M$ used in Theorem \ref{sec-tbaker} is only
	partially justified in section \ref{sec-Sveccalc}.

	\sssec{Theorem. [Baker]}\label{sec-tbaker}
	{\em Let $A,$ $B,$ $C,$ $E$ be points on the line 
	$a := O \times M$ such that}\\
	$A = O + $a$ M,$ $B = O + $b$ M,$ $C = O + $c$ M,$ $E = O + M,$ {\em the
	following constructions gives points} $O + $x$ M,$ $A'$ {\em for}
	x = $-$a, $D$ {\em for} x = a+b, $D'$ {\em for} x = a+b+c,\\
	$L$ {\em for} x = ab, $I$ {\em for} x = a$^{-1}.$

	{\em Let $P$ be a point not on $a$ and let $Q$ be a point on
	$A \times P,$ distinct from $A$ and $P$.}
	\enumb
	\item {\em Let\\
	\hth$q := O \times Q,$ $U := q \times (P \times M),$\\
	\hth$p := O \times P,$ $V := p \times (Q \times M),$\\
	\hth$	A' := (U \times V) \times a,$\\
	then}\\
	\hth$A' = O + $(-a)$ A.$
	\item {\em Let\\
	\hth$pb := P \times B,$ $R := pb \times (O \times Q),$\\
	\hth$	b := R \times M$,\\
	\hth$pa := A \times P,$ $S := pa \times b,$\\
	\hth$	c := Q \times M,$\\
	\hth$	T := pb \times c,$\\
	\hth$	D := (S \times T) \times a,$\\
	then}\\
	\hth$D = O + $(a+b)$ M.$
	\item {\em Let\\
	\hth$	R' := (O \times T) \times b,$\\
	\hth$	T' := (C \times R') \times c,$\\
	\hth$	D' := (S \times T') \times a,$\\
	then}\\
	\hth$D' = O + $(a+b+c)$ M.$
	\item {\em Let\\
	\hth$J := pa \times (E \times R),$\\
	\hth$K := pb \times (J \times M),$\\
	\hth$L := (Q \times K) \times a,$\\
	then}\\
	\hth$L = O + $ab$ M.$
	\item {\em Let\\
	\hth$	G := (P \times E) \times c,$\\
	\hth$	H = (Q \times E) \times p,$\\
	\hth$	I = (G \times H) \times a,$\\
	then}\\
	\hth$I = O + $a$^{-1} M.$
	\enume

	Proof: Choose the coordinate system such that\\
	$O = (1,0,0),$ $M = (0,1,0),$ $P = (0,0,1),$ then, for some
	a, b, c and q $\neq$ 0,\\
	\hth$A$ = (1,a,0), $B$ = (1,b,0), $C$ = (1,c,0), $E$ = (1,1,0),
	$Q$ = (1,a,q).\\
	For 0, we have $a$ = [0,0,1], $q$ = [0,q,$-$a], $U$ = (0,a,q),
	$p$ = [0,1,0], $V$ = (1,0,q), $A'$ = (1,-a,0).\\
	For 1, we have $pb$ = [b,$-1$,0], $R$ = (a,ab,bq), $b$ = [bq,0,$-$a],
	$pa$ = [a,$-$1,0],\\
	\hth$S$ = (a,a$^2$,bq), $c$ = [q,0,$-$1], $T$ = (1,b,q), 
	$D$ = (1,a+b,0).\\
	For 2, we have $R'$ = (a,b$^2$,bq), $T'$ = (b,bc+b$^2-$ac,bq),
		$D' = (1,a+b+c,0).$\\
	For 3, we have $J$ = (a(b$-1$),a$^2$(b$-1$),bq(a$-1$)),
		$K$ = (a(b$-1$),ab(b$-1$),bq(a$-1$)),\\
	\hth$L$ = (1,ab,0).\\
	For 4, we have $G$ = (1,1,q), $H$ = ($1-$a,0,1), $I$ = (a,1,0).

	\sssec{Exercise.}
	Give constructions
	\enumb
	\item associated with the associativity, commutativity and
	distributivity rules
	(a + b) + c = a + (b + c), (a b) c = a (b c), b + a = a + b,\\
	b a = a b and a (b + c) = a b + a c.
	\item for $F = O + $ab$^{-1} M.$
	\enume

	\ssec{Finite vector calculus and simple applications.}
	\label{sec-Sveccalc}
	\setcounter{subsubsection}{-1}
	\sssec{Introduction.}
	The following properties generalize, to the finite case,
	well known properties of vector calculus.  Capital letters will
	represent points, lower case letters will represent lines, the role
	of points and lines can be interchanged because of duality.  I have
	chosen to give at once the relations which directly apply to geometry
	rather than those which correspond to vector calculus.  These would be
	obtained if all lower case letters are replaced by upper case letters.
	1 or 2 of Theorem \ref{sec-tvecc} justify the representation of any
	point on the line $A * B$ by $k A + l B$.

	\sssec{Theorem.}\label{sec-tvecc}
	\enumb
	\item$     A * B = -B * A.$
	\item$     (A * B) * c = (A \cdot c)\: B - (B \cdot c)\: A.$
	\item$     a * (B * C) = (a \cdot C)\: B - (a \cdot B)\: C.$
	\item$     (A * B) \cdot C = (B * C) \cdot A = (C * A) \cdot B$\\
	\hth$  = - (B * A) \cdot C = - (C * B) \cdot A = - (A * C) \cdot B.$
	\item$     (A * B) \cdot (c * d) = (A \cdot c)(B \cdot d)
		- (A \cdot d) (B \cdot c).$
	\item$     (A * B) * (C * D) = (A \cdot (C * D))B
		- (B \cdot (C * D)) A.$
	\item$     (C * A) * (A * B) = ( (A * B) \cdot C )\: A.$
	\item$     ((A * B) \cdot C) \:P = ((B * C) \cdot P)\: A
		+ ((C * A) \cdot P)\: B + ((A * B) \cdot P)\: C.$
	\item$     (A * B) * C + (B * C) * A + (C * A) * B = 0.$
	\enume

	Proof:  The proof of 0 is immediate, the proof of 1. follows from the
	computation of any of the components of the triples on both sides, for
	the 0-th component,\\
	$(A_2 B_0 - A_0 B_2) c_2 - (A_0 B_1 - A_1 B_0) c_1
		= (A_1 c_1 + A_2 c_2) B_0 - (B_1 c_1 + B_2 c_2) A_0,$\\
	adding and subtracting $A_0$ $B_0$ $c_0$ gives the 0-th component of the
	second member of 1.  2, follows from 0 and 1.
	3, give various expressions of the 3 by 3 determinant constructed
	with the 3 triples as the 3 columns, namely\\
	$A_0 B_1 C_2 + A_1 B_2 C_0 + A_2 B_0 C_1
		- A_0 B_2 C_1 - A_1 B_0 C_2 - A_2 B_1 C_0.$\\
	for 4, $(A * B) \cdot (C * D) = ((C * D) * A) \cdot B
		= (A \cdot C) (D \cdot B)  - (A \cdot D) \cdot  (C \cdot B).$\\
	because of 3 and then 1.\\
	6, follows from 1.\\
	7, follows from 0, 1 and 2 applied to $(B * A) * (C * P).$\\
	8, follows from 1.

	\sssec{Theorem.}\label{sec-tcolldet}
	{\em $A,$ $B$ and $C$ are collinear  iff  $(A * B) \cdot C = 0.$}

	Proof:  This an immediate consequence of the fact that
	$(A * B) \cdot C = 0$ iff  $C$ is on the line $A * B.$

	\sssec{Theorem. [Fano]}
	{\em If $p$ is odd and \{$A,B,C,D$\} is a complete quadrangle, the
	intersections $b * b_1,$ $c * c_1,$ $d * d_1$ of the opposite sides
	are not collinear.}

	Proof:  $(b * b_1) * (c * c_1) * (d * d_1)$\\
	\hth$   = 2 ((B * C) \cdot D) ((C * D) \cdot A) ((D * A) \cdot B)
	((A * B) \cdot C) \neq 0,$\\
	by repeated use of \ref{sec-tvecc}.0 to .3.

	\sssec{Comment.}
	$p = 2$ is excluded, because in this case, the preceding
	Theorem is false, in fact every complete quadrangle has its diagonal
	collinear.  I leave as an exercise the determination of where the above
	theory breaks down.

	\sssec{Comment.}
	In algebraic manipulations, although $A_i = a_{i+1} \times a_{i-1},$
	we can not use this expression when the sum of two or more terms is
	involved, because the scaling to go from ``$*$" to $``\times"$ is
	different for each index $i$.  It is therefore essential to do these
	algebraic manipulations using ``$*$".
	For the various proof, $a_i$ will always denote $A_{i+1} * A_{i-1}$
	and use will often been made of the following Theorem:

	\sssec{Theorem.}
	{\em If $a_i := A_{i+1} * A_{i-1}$ and $t := (A_0 * A_1) \cdot A_2,$
	then\\
	\hth$        a_{i+1} * a_{i-1} = t A_i.$\\
	Again, $i = 0,1,2$ and subscript addition is made modulo 3.}

	Proof:  The conclusion follows from \ref{sec-tvecc}.7 and from
	\ref{sec-tvecc}.3.

	\sssec{Comment.}
	The identity \ref{sec-tvecc}.8 is the fundamental identity in Lie
	algebras.  The set of points or the set of lines form a Lie algebra, if
	we use as multiplication ``$*$".  See, for instance, Cohn, Lie groups.

	\sssec{Notation.}
	$det(A,B,C)$ will denote $(A * B) \cdot C = (B * C) \cdot A =
	 (C * A) \cdot B =$.

	\sssec{Theorem.}\label{sec-tvecpa}
	$det(A,C,E)\:det(B,D,E)\:det(A,B,F)\:det(C,D,F)\\
	\hth	= det(A,C,F)\:det(B,D,F)\:det(A,B,E)\:det(C,D,E)$ iff\\
	$(((A\times E)\times (D\times C))\times ((E\times B)\times (C\times F)))
	\cdot ((B\times D)\times (A\times F)) = 0$.

	Proof:
	By \ref{sec-tvecc}.1 and .2, the second equation is equivalent to\\
	$((A\times E)\times (D\times C))\cdot (C\times F)\:
	((B\times D)\times (A\times F))\cdot (E\times B) =\\
	\hth((A\times E)\times (D\times C))\cdot (E\times B)\:
	((B\times D)\times (A\times F))\cdot (C\times F),$\\
	by \ref{sec-tvecc}.3 this is equivalent to\\
	$((D\times C)\times (C\times F))\cdot (A\times E)\:
	((E\times B)\times (B\times D))\cdot (A\times F) =\\
	\hth((E\times B)\times (A\times E))\cdot (D\times C)\:
	((A\times F)\times (C\times F))\cdot (B\times D),$\\
	by \ref{sec-tvecc}.6 this is equivalent to\\
	$(det(C,F,D)\: C)\cdot (A\times E)\:
	(det(B,D,E)\: B)\cdot (A\times F) =\\
	\hth(det(A,B,E)\: E)\cdot (D\times C)\:
	(det(A,F,C)\: F)\cdot (B\times D).$

	This can be considered as an algebraic form of Pascal's Theorem.
	for the order $A,E,B,D,C,F.$

	\ssec{Anharmonic ratio, harmonic quatern, equiharmonic quatern.}

	\sssec{Convention.}
	In this section, I will use the convention that if point
	is on the line [0,1,0],\\
	\hth$(\infty ,0,1)$ denotes the point $(1,0,0).$

	\sssec{Definition.}
	Given 4 points on the line [0,1,0],\\
	$A_0 = ($m$_0,0,1),$ $A_1 = ($m$_1,0,1),$ $A_2 = ($m$_2,0,1),$
	$A_3 = ($m$_3,0,1)$.\\
	The {\em anharmonic ratio} is defined by\\
	\hth anhr$(A_0,A_1,A_2,A_3)$
	:= anhr$({\rm m}_0,{\rm m}_1,{\rm m}_2,{\rm m}_3)
	:= \frac{({\rm m}_2-{\rm m}_0)
	({\rm m}_3-{\rm m}_1)}{({\rm m}_2-{\rm m}_1)({\rm m}_3-{\rm m}_0)}.$\\
	If m$_i = \infty$  then the 2 factors containing m$_i$ are dropped, e.g.
\\
	if m$_0 = \infty$  then anhr$(A_0,A_1,A_2,A_3)$ :=
	anhr$(\infty,{\rm m}_1,{\rm m}_2,{\rm m}_3) :=
	\frac{{\rm m}_3-{\rm m}_1}{{\rm m}_2-{\rm m}_1}.$ 

	\sssec{Lemma.}
	{\em Let $a$, $b$, $c$ and $d$ be such that $ad-bc \neq 0$, if
	$t({\rm m}) := \frac{a{\rm m}+b}{c{\rm m}+d}$ then\\
	anhr$({\rm m}_0,{\rm m}_1,{\rm m}_2,{\rm m}_3)$ = 
	anhr$(t({\rm m}_0),t({\rm m}_1),t({\rm m}_2),t({\rm m}_3))$.}

	If we project 4 points $A_i$ on $a$ onto 4 points $B_i$ on $b$
	from the point $B$, it is easy to see that each coordinate of $B$
	is a linear functions of {\rm m}, therefore the ratio of 2 of some
	specificate coordinates of $B$ are functions of the form $t({\rm m})$.
	This justifies the following 2 Theorems.

	\sssec{Theorem.}
	{\em Given 4 points $B_i$ on a line $b$, the anharmonic ratio
	of the 4 points is the anharmonic ratio of the 4 ratios obtained
	by dividing the $j$-$th$ coordinate of $B_i$ by the $k$-$th$ coordinate
	for appropriate $j \neq k$.}

	\sssec{Theorem.}
	{\em If 4 points $B_i$ are obtained by successive projections from
	4 points $A_i$, then anhr$(B_0,B_1,B_2,B_3)$ = anhr$(A_0,A_1,A_2,A_3)$.}

	\sssec{Theorem.}
	{\em If $r := $anhr$(A_0,A_1,A_2,A_3),$ then if we permute the points in
	all possible way we obtain, in general 6 different values of the
	anharmonic ratio:}\\
	For \hfill the anharmonic ratio is\\
	$\begin{array}{llllc}
	    A_0,A_1,A_2,A_3 &  A_1,A_0,A_3,A_2 &  A_2,A_3,A_0,A_1 &
	  A_3,A_2,A_1,A_0&   r\\
	    A_0,A_1,A_3,A_2 &  A_1,A_0,A_2,A_3 &  A_2,A_3,A_1,A_0 &  
		A_3,A_2,A_0,A_1&   \frac{1}{r}\\
	    A_0,A_2,A_1,A_3 &  A_1,A_3,A_0,A_2 &  A_2,A_0,A_3,A_1 &
	  A_3,A_1,A_2,A_0&   1-r\\
	    A_0,A_2,A_3,A_1 &  A_1,A_3,A_2,A_0 &  A_2,A_0,A_1,A_3 &
		  A_3,A_1,A_0,A_2&   \frac{1}{1-r}\\ 
	    A_0,A_3,A_1,A_2 &  A_1,A_2,A_0,A_3 &  A_2,A_1,A_3,A_0 &  
		A_3,A_0,A_2,A_1&   \frac{r-1}{r}\\
	    A_0,A_3,A_2,A_1 &  A_1,A_2,A_3,A_0 &  A_2,A_1,A_0,A_3 &
	  A_3,A_0,A_1,A_2&   \frac{r}{r-1} 
	\end{array}$

	\sssec{Theorem.}
	{\em There are 3 cases for which the 6 values are not distinct:\\
	0.\hti{6}$0,\infty,1,1,\infty ,0,$ when 2 points are identical.\\
	1.\hti{6}$-1,-1,2,\frac{1}{2},2,\frac{1}{2}.$\\
	2.\hti{6}$v,\frac{1}{v},\frac{1}{v},v,v,\frac{1}{v},$ with $v^2-v+1 = 0$
	or $v = \frac{1+\sqrt{-3}}{2}$\\
	\hth$v$ is real, if $p \equiv 1 \pmod{6}$.}

	For instance, if $p = 7,$ then $v = -2$ or 3, if $p = 19,$
	then $v = -7$ or 8.

	\sssec{Theorem.}
	{\em Given a complete quadrangle ${A,B,C,D}$, the intersection of 2 of
	lines through opposite vertices and the line through 2 diagonal 
	points make a harmonic quatern with these diagonal points. More
	precisely, let 
	$E := (A \times B)\times(C\times D)$ and
	$F := (A \times D)\times(B\times C)$ be 2 of the diagonal points,
	let $a := E\times F$, $G := a \times (B\times D)$ and
	$H := a \times (A\times C)$, then $ r :=$ anhr$(E,F,G,H) = -1$.}

	Let $I$ be the third diagonal point, projecting the 4 points from $B$
	on $A\times C$ and these from $D$ on $a$ gives $r =$ anhr$(A,C,I,H)$
	= anhr$(F,E,G,H) = \frac{1}{r}$, therefore $r^2 = 1$ but we do not
	have case 0, $r = 1$, therefore $r = -1$ which is case 1.

	\sssec{Definition.}\label{sec-dharmq}
	In the special case of the preceding Theorem:\\
	Case 1, we say that $A_2,$ $A_3$ are {\em harmonic conjugate} of
	$A_0,$ $A_1,$ or that $A_0,A_1,A_2,A_3$ form a {\em harmonic quatern}.\\
	Case 2, we say that $A_0,$ $A_1,$ $A_2,$ $A_3$ form a
	{\em equiharmonic quatern}.

	\sssec{Definition.}\label{sec-dequih}
	The {\em pre-equiharmonic non confined configuration} is defined as
	follows:\\
	Given a complete quadrangle $A,B,F,K,$ determine\\
	\hth$	q := A \times B,$ $p := A \times F,$ $b := B \times F,$
	$r := A \times K,$ $f := B \times K,$\\
	\hth$	H := p \times f,$ $J := b \times r,$\\
	choose $C$ on $q,$ distinct from $A$ and $B,$ determine\\
	\hth$c := C \times K,$ $P := b \times c,$ $R := p \times c,$
		$g := C \times F,$ $Q := f \times g,$ $L := g \times r,$\\
	\hth	$d := P \times L,$ $h := Q \times R,$ $D := d \times h.$

	\sssec{Theorem.}
	{\em Given} \ref{sec-dequih}
	\begin{enumerate}
	\item$D \cdot q = 0.$
	\item$J \cdot h = 0  \Rightarrow   H \cdot d = 0.$
	\item {\em The geometric condition $J \cdot h = 0$ is equivalent to\\
	\hth$	A,B,C,D$ is an equiharmonic quatern.}
	\end{enumerate}

	Proof:  Let $A = (0,0,1),$ $B = (0,1,1),$ $F = (1,0,1),$ $K = (1,1,1).$
	We have $q = [1,0,0],$ $p = [0,1,0],$ $b = [-1,-1,1],$ $r = [-1,1,0],$
	$f = [0,-1,1],$ $H = (1,0,0),$ $J = (1,1,2).$
	Let $C$ = (0,c,1), then $c$ = [1$-$c,$-1$,c], $g$ = [c,1,$-$c],
	$P$ = (1$-$c,1,2$-$c),
	$R$ = (c,0,c$-$1), $Q$ = (c$-$1,c,c), $L$ = (c,c,c+1),
	$d$ = [c$^2-$c+1,2c$-1,-$c$^2$], 
	$h$ = [c$^2-$c,2c$-1,-$c$^2$], $D$ = (0,c$^2$,2c$-$1).\\
	$J\cdot h = 0$ or $H\cdot d = 0$ are equivalent to c$^2 -$c+1 = 0.

	\sssec{Definition.}
	Given \ref{sec-dequih} and $J \cdot h = 0,$ the pre-equiharmonic
	configuration is then called an {\em equiharmonic configuration}.

	\sssec{Theorem.}
	{\em A pre-equiharmonic configuration is of type\\
	\hth$	10 \times 3 + 2 \times 2  \:\&\:  7 \times 4 + 2 \times 3,$\\
	unless it is equiharmonic, in which case, it is of type\\
	\hth$	12 \times 3  \:\&\:  9 \times 4.$\\
	The sets of 4 points on each of the 9 lines is an equiharmonic quatern.}

	Proof:  The configuration, with projections as given below is as
	follows.
\newpage
	Points:\hti{18}     lines:\hti{27}        from\\
	$\begin{array}{cccccccccccccc}
	A :& q,&p,&r,&\hti{4}&       q :& A,&B,&C,&D,\\
	B :& q,&b,&f,&&  p :& A,&F,&R,&H,&\hti{4}&	P&   (Q,H)\\
	C :& q,&c,&g,&&  r :& A,&J,&K,&L,&&	P&   (Q,K)\\
	D :& q,&h,&d,&&  b :& F,&B,&P,&J,&&	R&   (L,J)\\
	P :& b,&c,&d,&&  f :& H,&B,&K,&Q,&&	R&   (L,K)\\
	R :& p,&c,&h,&&  c :& R,&K,&C,&P,&&	H&   (J,K)\\
	H :& p,&f,&d,&&  g :& F,&Q,&C,&L,&&	H&   (J,L)\\
	K :& r,&f,&c,&&  d :& L,&H,&P,&D,&&	K&   (F,H)\\
	Q :& f,&g,&h,&&  h :& J,&Q,&R,&D,&&	K&   (F,R)\\
	L :& r,&g,&d,\\
	J :& r,&b,&h,\\
	F :& p,&b,&g,\\
	\end{array}$\\
	If we project $A,B,C,D$ from $P$ on $p$ we get $A,F,R,H;$
	if we project $A,B,C,D$ from $Q$ on $p$ we get $A,H,F,R;$
	therefore $r = \frac{1}{1-r}.$
	To establish the results for the other sets, it is sufficient to
	project from a point, those of the line $q.$  The points on each line
	have been arranged correspondingly.  For instance, for $f,$ the point of
	projection is $R$ and the lines are $p,$ $c$ and $h$; if the point of
	projection is $L,$ the order is $K,$ $B,$ $Q,$ $H,$
	(the second point corresponds
	to $A$ or $B$ for $p$ and $r$ the others are obtained circularly).

	\sssec{Exercise.}
	Prove that the configuration of \ref{sec-einfcub} is equiharmonic.

	\sssec{Exercise.}\label{sec-eequi8}
	Study the configuration which starts with $P_0,$ $P_1,$ $P_3,$ $P_5$
	and $P_7$ on $P_1 \times P_5$.
	Constructs $l_0 := P_0 \times P_1,$ $l_1 := P_0 \times P_3,$
	$l_3 := P_0 \times P_5,$ $l_5 := P_3 \times P_5,$
	$l_4 := P_1 \times P_5,$ $l_8 := P_3 \times P_7,$
	$P_6 := l_2 \times l_8,$ $l_3 := P_1 \times P_6,$
	$P_4 := l_1 \times l_3,$ $l_7 := P_4 \times P_7,$
	$P_2 := l_0 \times l_5$.\\
	Using a coordinate system such that $P_0 = (1,0,0)$, $P_1 = (0,1,0)$,
	$P_3 = (0,0,1)$ and $P_5 = (1,1,1)$, determine an algebraic condition
	involving the coordinates of $P_7$ for $P_2$ to be on $l_7$.\\
	Prove that if $P_2$ is on $l_7$, the configuration is of type
	$8 \times 3 \:\&\: 8 \times 3$.

	\sssec{Exercise.}
	Study the configuration which starts with $A_0,$ $A_1,$ $B_0,$ $B_1$
	and $A_2$ on $A_0 \times A_1,$ constructs $d_0 := A_0 \times B_0,$
	$d_1 := A_1 \times B_1,$
	$d_3 := A_0 \times B_1,$ $d_4 := A_1 \times B_0,$
	$d_6 := A_2 \times B_0,$
	$d_8 := A_2 \times B_1,$
	$a := A_0 \times A_1,$ $b := B_0 \times B_1,$ $P := a \times b,$
	$C_0 := d_3 \times d_4,$ $C_1 := d_1 \times d_6,$
	$C_2 := d_0 \times d_8,$
	$a_0 := C_0 \times C_1,$ $d_5 := A_0 \times C_1,$
	$d_7 := A_1 \times C_2,$
	$B_2 := d_5 \times d_7.$\\
	Determine a geometric condition on $A_2$ for $P, A_0, A_1, A_2$ to be
	an equiharmonic quatern, prove that in this case $P$, $B_0,$ $B_1,$
	$B_2$ is also a equiharmonic quatern.


	\ssec{Projectivity of lines and involution on a line.}
	\label{sec-Sprojinv}
	\setcounter{subsubsection}{-1}
	\sssec{Introduction.}
	In the next section we will study algebraically the
	isomorphisms of the plane into itself.  The special case of the
	mapping of a line of the plane into a line will be defined here.  The
	justification will follow from the general definition.  Such a mapping
	is called a projectivity.  Special cases will be studied and
	appropriate constructions will be given.  The notion of amicable
	projectivities, which are at the basis of the definition of equality of
	angles is also introduced.  The concept of harmonic conjugates is due
	to LaHire\footnote{Coxeter, p.16}.  The term projectivity will be used
	here
	only for correspondances between points on lines not for correspondance
	of a plane with itself as done by some authors.  Theorem
	\ref{sec-tconproj} gives a construction, when 2 points $A$ and $B$ are
	fixed, and $D$ corresponds to $C.$

	\sssec{Convention.}
	For simplicity, I will assume that the line is [0,0,1],
	the last component of all the points is 0, I will therefore only write
	the first 2 components.

	\sssec{Definition.}
	The mapping which associates to the point $(x_0,x_1)$ the
	point $(y_0,y_1)$ given by\\
	\hth$	(y_0) = (a  b) (x_0),$\\
	\hth$	(y_1) = (c  d) (x_1)$\\
	with $ad-bc\neq$  0, is called a {\em projectivity}.

	\sssec{Theorem.}
	{\em If $C$ is the intersection of $c$ and $A \times B,$ the point $D$
	such that $A,$ $B,$ $C$ and $D$ form a harmonic quatern is given by\\
	\hth$	D := B \cdot c A + A \cdot c B.$}

	\sssec{Definition.}
	$D$ is called the {\em harmonic conjugate of $C$ with respect to}
	$A$ and $B.$

	\sssec{Theorem.}
	\enumb
	\item {\em  If $K = (1,k,0),$ $L = (1,l,0)$ and $M = (1,m,0)$ then
		$N,$ the conjugate of $M$ with respect to $K$ and $L,$ is given
		by}\\
	\hth$	N = (2m-l-k,k\:m+l\:m-2k\:l,0).$
	\item {\em If $N$ is the harmonic conjugate of $M$ with respect to $K$
	 and $L$,}
	\enume
		{\em then\\
	\hth$M$ is the harmonic conjugate of $N$ with respect to $K$ and $L,$\\
	\hth$K$ is the harmonic conjugate of $L$ with respect to $M$ and $N$,
		and\\
	\hth$L$ is the harmonic conjugate of $L$ with respect to $M$ and $N.$}

	\sssec{Theorem.}\label{sec-tconproj}
	{\em If $A := (1,0,0),$ $B := (0,0,1),$ $C := (1,k,0),$ $D := (1,l,0),$
	the projectivity on $c := A \times B$ which associates $A$ to $A,$
	$B$ to $B$ and $C$ to $D,$ can be constructed by choosing $a$ line $b$
	through $a,$ dictinct from $c,$ a line a through $b$ distinct from
	$c,$ a point $P$ on $a$ not on $b$ or $c,$ then\\
	\hth$S := (P \times C) \times b,$ $T := (P \times D) \times b.$\\
	The mapping $N$ of $M := (1,m,0)$ is obtained by the construction\\
	\hth$Q := (M \times S) \times a,$ $N := (Q \times T) \times c$
	 and $N = (k,lm,0).$}

	The proof is left as an exercise.

	\sssec{Theorem.}
	{\em Let $u := \frac{l}{k},$ $\phi(M) = (1,u m,0)$ and 
	${\phi^j}(M) = (1,u^j m,0).$
	If $u$ is a primitive root of $p,$ the projectivity has order $p-1.$}

	The proof is left as an exercise.

	\sssec{Theorem.}
	{\em If $K,$ $L$ and $M$ are distinct, $N$ is distinct from $M.$}

	Proof.  The last theorem would imply that
	$m(2m-l-k) = k\:m+l\:m-2k\:l or (m-l)(k-m) = 0.$

	\sssec{Theorem.}
	{\em If $K$ and $M$ are exchanged and $L$ is replaced by $N,$ then $N$
	is replaced by $L.$}

	Indeed, $n(2m-k-l) = k\:m+l\:m-2k\:l$ can be written
	$l(2k-m-n) = m\:k+n\:k-2m\:n.$

	The following theorem gives a construction of a projectivity on a line
	in which $B_i$ corresponds to $A_i,$ $i = 0,1,2.$

	\sssec{Theorem.}\label{sec-tconproj1}
	{\em Given $A_i,$ $B_i,$ $i = 0,1,2,$ 6 points on a line $u,$ such that
	the $A_i$ are distinct and the $B_i$ are distinct.
	Choose the line $s\neq$  $u$ and the point $S,$ with
	$S \cdot u\neq$  0, $S \cdot s\neq$  0.\\
	Construct\\
	$C_i := (S \times A_i) \times s,$ $D_j := (B_0 \times C_j) \times 
	(C_0 \times B_j),$ $j = 1,2,$
	$d := D_1 \times D_2,$ then for any $A_l$ on $u$, construct\\
	$C_l := (S \times A_l) \times  s,$ $D_l := (B_0 \times C_l) \times d,$
	$B_l := (C_0 \times D_l) \times u.$\\
	The mapping which associates $B_l$ to $A_l$ is a projectivity.}

	\sssec{Theorem.}
	{\em Given $A_0 = (1,0,0),$ $A_1 = (1,a_1,0),$ $A_2 = (1,a_2,0)$ and
	$B_0 = (0,1,0),$ $B_1 = (1,b_1,0),$ $B_2 = (1,b_2,0),$
	then the projectivity which is defined in the preceding theorem and
	associates to $A_i,$ $B_i,$ $i = 0,1,2,$ associates to\\
	\hth$	A_j = (1,a_j,0),$ the point\\
	\hth$B_j = ( (a_2-a_1)a_j, (a_2b_2-a_1b_1)a_j - a_1a_2(b_2-b_1), 0),$
	$j > 2.$}

	The proof is left as an exercise.

	\sssec{Theorem.}
	{\em Let $g = \frac{a_2\:b_2-a_1\:b_1}{2(a_2-a_1})$ and
		$h = a_1\:a_2\frac{b_2-b_1}{a_2-a_1},$\\
	The projectivity, which associates to $(1,u,0),$}
	$(1,2g-\frac{h}{u}),0)$
	\enumb
	\item {\em is an involution  iff  $g = 0,$}
	\item {\em is an hyperbolic projectivity if $g^2-h$ is a quadratic
	residue modulo $p,$}\\
	\item {\em is an elliptic projectivity if $g^2-h$ is a non residue and}
	\item {\em is a parabolic projectivity if $h = g^{2,}$ the fixed point
		being $g.$}
	\enume

	Proof.  If we eliminate $u_1$ from $u_1 = 2g-\frac{h}{u_0}$
	 and $u_0 = 2g-\frac{h}{u_1},$
	we get $2g(u_0^2-g\:u_0+h).$  If this relation is to be satisfied for
	all $u_0,$
	it is necessary that $g = 0.$  The condition is sufficient because if
	$\phi(u) := \frac{h}{u},$ then $\phi \circ \phi$  is the identity.

	\sssec{Theorem.}\label{sec-tconproj2}
	{\em Given 3 distinct points $A_0,$ $A_1$ and $A_2$ on the line $a$ and
	3 distinct points $B_0,$ $B_1$ and $B_2$ on the line $b,$\\
	let $A_2 = r_0 A_0 + r_1 A_1$ and 
	$B_2 = s_0 B_0 + s_1 B_1,$ then\\
	if $A_j = t_0 A_0 + t_1 A_1,$ 
	$B_j = \phi(A_j) := \frac{s_0 t_0}{r_0} B_0 
	+ \frac{s_1 t_1}{r_1 B_1}$\\
	is a projectivity which associates $A_j$ to $B_j$ for all $j.$}

	Proof.  The correspondance clearly associates $A_j$ to $B_j$ for
	$j = 0,$
	using $t_1 = 0,$ for $j = 1$ using $t_0 = 0$ and for $j = 2$ using
	$t_0 = r_0$ and $t_1 = r_1.$ 
	The proof that it is a projectivity is left as an exercise.

	\sssec{Theorem.}\label{sec-tconproj3}
	\enumb
	\item {\em If the lines $a$ and $b$ of the preceeding Theorem coincide,
	    there exists constants $f_0$, $f_1$, $f_2$ and $f_3$ such that\\
	\hth$    B_j = (f_0 t_0 + f_1 t_1) A_0 + (f_2 t_0 + f_3 t_1) A_1.$\\
	\hth    If $B_0 = b_{00} A_0 + b_{01} A_1$ and
	$B_1 = b_{10} A_0 + b_{11} A_1,$
	then}\\
	\hth$f_0 = \frac{s_0 b_{00}}{r_0},$ $f_1 = \frac{s_1 b_{10}}{r_1,}$ 
	    $f_2 = \frac{s_0 b_{01}}{r_0},$ $f_3 = \frac{s_1 b_{11}}{r_1}.$ 
	\item {\em The values $t_0$ and $t_1$ for which $A_j$ is a fixed point,
	in other words, for which $A_j$ = $B_j$ satisfy}\\
	\hth$	f_1 t_1^2 - (f_3-f_0) t_0 t_1 - f_2 t_0^2 = 0.$
	\item {\em The projectivity is hyperbolic, parabolic or elliptic if\\
	\hth$	(f_3-f_0)^2 + 4 f_1f_2$ is positive, zero or negative.}
	\item {\em The projectivity is an involution  iff  $f_0 + f_3 = 0.$}
	\enume

	The proof is left as an exercise.

	\sssec{Definition.}\label{sec-dconproj}
	Using the notation of \ref{sec-tconproj2} and of \ref{sec-tconproj3}
	with primes used
	for an other projectivity, we say that {\em 2 projectivities on the same
	line are amicable}  iff  there exists a constant $k$ different from 0
	such that\\
	\hth$f_1' = k f_1,$ $f_2' = k f_2,$ $f_3'-f_0' = k (f_3-f_0).$

	\sssec{Theorem.}
	{\em Two amicable projectivities are either both hyperbolic, or
	both parabolic or both elliptic.  If they are both hyperbolic,
	they have the same fixed points.}

	\sssec{Example.}\label{sec-eprojp5}
	For $p = 5,$ a projectivity $\phi$  associates to\\
	$A_0 = (5),$ $A_1 = (10),$ $A_2 = (14), (18), (22), (26),$
	$B_0 = (26),$ $B_1 = (5),$ $B_2 = (18), (22), (10), (14).$
	$r_0 = r_1 = 1,$ $s_0 = 1,$ $s_2 = -2,$ $b_{00} = -1,$ $b_{01} = 1,$
	$b_{10} = 1,$
	$b_{11} = 0,$ $f_0 = -1,$ $f_1 = -2,$ $f_2 = 1,$ $f_3 = 0.$
	A second projectivity $\phi'$ associates to\\
	$A'_0 = (5),$ $A'_1 = (10),$ $A'_2 = (14), (18), (22), (26),$
	$B'_0 = (18),$ $B'_1 = (14),$ $B'_2 = (10), (5), (26), (22).$
	$r'_0 = r'_1 = 1,$ $s'_0 = -1,$ $s'_2 = 2,$ $b00' = 2,$ $b01' = 1,$
	$b_{10}' = 1,$
	$b_{11}' = 1,$ $f_0' = -2,$ $f_1' = 2,$ $f_2' = -1,$ $f_3' = 2.$
	$\phi'$ is an involution an $f_0' + f_3' = 0.$
	$\phi$	and $\phi'$ are sympathic, with $k = -1.$
	The fixed points are complex and correspond to $t_0$ = 1 and
	$t_1 = 1 + \sqrt{-2}$ or $t_1 = 1 - \sqrt{-2}.$

	\sssec{Comment.}
	The definition \ref{sec-dconproj} will be used in III.1.3. to define
	equality of angles.

	\sssec{Theorem.}
	{\em If $x$ is one of the coordinates, the projectivity takes the
	form\\
	\hth$	F(x) = \frac{a+bx}{c+dx}$\\
	and the fixed points are the roots of\\
	\hth$	d x^2 + (c-b) x -a = 0.$}

	\sssec{Exercise.}

	Prove that the following construction defines a projectivity on $u$
	in which $A_{i+1}$ corresponds to $A_i,$ the points $A_0$ to $A_3$ being
	given.  Let
	$l_f$ is a line through $A_2$ distinct from $u,$
	$E$ is a point on $l_f$ distinct from $A_2,$ 
	$F$ is a point on $l_f$ distinct from $A_2$ and $E,$
	$l_d$ is a line through $A_0$ distinct from $u,$
	$D_1 := l_f \times l_d,$
	$D$ is a point on $A_1 \times D_1$ distinct from $A_1$ and $D_1,$
	$E_0 := (A_0 \times E) \times (A_1 \times F),$
	$E_2 := (A_3 \times F) \times (E \times (l_d \times (A_2 \times D))),$
	$A_{i+1} := (((((A_i \times D) \times l_d) \times E) \times l_e)
	 \times F) \times i,$ $i = 4,$ \ldots.
	The preceding construction is less efficient than that in
	\ref{sec-tconproj1}.

{\tiny\footnotetext[1]{G22.TEX [MPAP], \today}}
	\ssec{Collineation, central collineation, homology and elation.}
	\label{sec-Scollinalg}
	\setcounter{subsubsection}{-1}
	\sssec{Introduction.}
	Collineations, \{which are isomorphisms of the plane onto
	itself\} have been defined in \ref{sec-dcollin}.  They will now be
	studied algebraically.  The {\em point mapping} which associates points
	to points is represented by a non singular matrix, and so is the
	{\em line mapping} which associates lines to lines. Two matrices which
	can be obtained from each other by multiplication, modulo $p,$ by an
	integer different from 0 correspond to the same collineation.

	\sssec{Theorem.}\label{sec-tcoll}
	{\em Given 2 complete quadrangles $A_j$ and $B_j,$ j = 0,1,2,3,\\
	Let $a_i := A_{i+1}*A_{i-1}$ and $b_i := B_{i+1}*B_{i-1},$\\
	let $A_3 = r_0 A_0 + r_1 A_1 + r_2 A_2,$
	    $B_3 = s_0 B_0 + s_1 B_1 + s_2 B_2$\\
	\hti{4}$q_i := \frac{s_i}{r_i},$ $u_i := q_{i+1} q_{i-1},$ \\
	then, up to a proportionality constant,
	$q_i = \frac{b_i \cdot B_3}{a_i \cdot A_3}.$\\
	Moreover,}
	\enumb
	\item {\em the mapping $\gamma$ defined by\\
	    $B_l := \gamma(A_l) :=
	q_0(a_0 \cdot A_l)B_0 + q_1(a_1 \cdot A_l)B_1 + q_2(a_2 \cdot A_l)B_2$\\
	is the point mapping of a collineation which associates to $A_j$, $B_j$
	for} j = 0 to 3.
	\item {\em the mapping $\gamma'$ defined by\\
	$\gamma'(a_l) := u_0(A_0 \cdot a_l)b_0 + u_1(A_1 \cdot a_l)b_1
		+ u_2(A_2 \cdot a_l)b_2$
	is the corresponding line mapping.}
	\enume

	Proof:  By hypothesis, $r_0 \neq 0$, because $A_3$ is not on
	$A_1 \times A_2,$ similarly, $r_1,$ $r_2,$ as well as $s_0,$ $s_1$ and
	 $s_2$ are $\neq 0$, therefore, $q_0,$ $q_1$ and $q_2$ are well defined
	and $\neq 0$.\\
	$a_0\cdot A_3 = (A_1*A_2)\cdot(r_0A_0+r_1A_1+r_2A_2) = r_0\:det(A_0,A_1,
		A_2),$\\
	similarly, $a_i\cdot A_3 = r_i\:det(A_0,A_1,A_2),$
		$b_i\cdot B_3 = r_i\:det(B_0,B_1,B_2),$\\
	hence the alternate expression for $q_i$.\\
	0. follows from \ref{sec-tvecc} by observing that
	$(A_i * A_3) * a_i = r_{i+1} A_{i+1} - r_{i-1} A_{i-1}.$\\
	The details are left as an exercise.\\
	If $M$ and $N$ are any 2 points on $a_l$ and $a_l = M * N,$\\
	$\gamma'(a_l) = \gamma'(M * N) = \gamma(M) * \gamma(N)$\\
	\hth$         = q_1 q_2 (a_1 \cdot M\: a_2 \cdot N 
		- a_2 \cdot M\: a_1 \cdot N) b_0 +  \ldots $\\
	\hth$         = q_1 q_2( (a_1 * a_2).(M * N)) b_0 +  \ldots ,$\\
	\hth$         = u_0 t (A_0 \cdot a_l) b_0 +  \ldots  ,$\\
	because of \ref{sec-tvecc}.  Dividing by $t,$ we get 1.

	\sssec{Theorem.}\label{sec-tidcoll}
	{\em If a collineation transforms each of the points of a complete
	quadrangle into itself, every point is transformed into itself.}

	\sssec{Definition.}
	The collineation of \ref{sec-tidcoll} is called the {\em identity
	collineation} $\epsilon$.

	\sssec{Comment.}
	Theorem \ref{sec-tconproj2} for 1 dimensional sets and Theorem
	\ref{sec-tcoll} for
	2 dimensional sets generalize by induction to $n$ dimensions.

	\sssec{Example.}
	For p = 5, let  $A_0 = (0) = (0,0,1),$ $A_1 = (1) = (0,1,0),$
	$A_2 = (6) = (1,0,0)$ and $A_3 = (12) = (1,1,1),$\\
	let $B_0 = (1),$ $B_1 = (6),$ $B_2 = (12),$ $B_3 = (3) = (0,1,2)$,
	to obtain the point mapping $\gamma$ which associates to $A_j$, $B_j$,
	$a_0 = [0,0,-1],$ $a_1 = [0,-1,0],$ $a_2 = [-1,0,0]$,
	$b_0 = [0,-1,1],$ $b_1 = [-1,0,1],$ $b_2 = [0,0,-1].$\\
	$q_0 = -1,$ $q_1 = -2,$ $q_2 = 2,$
	$u_0 = -4,$ $u_1 = -2,$ $u_2 = 2,$ therefore\\
	$\gamma\vet{X}{Y}{Z} = \mat{-2}{2}{0}{-2}{0}{1}{-2}{0}{0}\vet{X}{Y}{Z} =
	\matb{[}{c}-2X+2Y\\-2X+Z\\-2X\mate{]},$\\
	$\gamma'\matb{[}{c}x\\y\\z\mate{]} = \mat{0}{2}{0}{0}{0}{}{-1}{-2}{-2}
		\matb{[}{c}x\\y\\z\mate{]} = \vet{2y}{-z}{-2x-2y+z}.$

	\sssec{Comment.}
	The following mapping can be used in certain cases but not
	in all cases:\\
	$\phi(M) := q_0( M \cdot A) \phi (A) + q_1( M \cdot B) \phi (B)
	+ q_2 (M \cdot C) \phi (C),$ with\\
	$q_0 = \frac{1}{A \cdot P}(\phi (B) * \phi (C)) \phi (P),$
	$q_1 = \frac{1}{B \cdot P}(\phi (C) * \phi (A)) \phi (P),$\\
	\hth$q_2 = \frac{1}{C \cdot P}(\phi (A) * \phi (B)) \phi (P).$\\
	Indeed, one of the scalar product $A \cdot P$ or $B \cdot P$ or
	$C \cdot P$ can be 0, and cases exists for which whatever permutation
	of the 4 points $A,$ $B,$ $C$ and
	$P$ is used, the same difficulty occurs.

	\sssec{Theorem.}\label{sec-tcollmat}
	\ref{sec-tcoll} {\em can be rewritten using matrix notation. Let {\bf a}
	be a matrix whose rows are the components of the sides of the
	triangle $A_0,$ $A_1,$ $A_2,$ ${\bf a}_{i,j} := {\bf a}_{j,i},$ etc.\\
	Let {\bf Q} be the matrix ${\bf Q}_{i,i} := q_i,$ ${\bf Q}_{i,j} := 0$
	for $i \neq j,$\\
	let ${\bf U}_{i,i} := q_{i+1} q_{i-1},$ ${\bf U}_{i,j} := 0,$
	$i\neq j.$\\
	Let ${\bf A}_l$ and ${\bf B}_l$ be column vectors,\\
	then ${\bf M} := {\bf B\: Q\: a}^T$ defines the
	collineation ${\bf B}_l = {\bf M\: A}_l$ and ${\bf M}' :=
	{\bf b\: U\: A}^T$ gives ${\bf b}_l = {\bf M}' {\bf a}_l.$\\
	Moreover ${\bf M}' = {\bf M}^{-1}$ is the adjoint matrix.}

	\sssec{Example.}
	Let $A_0 = (7),$ $A_1 = (15),$ $A_2 = (19),$ $A_3 = (28),$\\
	$B_0 = (27),$ $B_1 = (3),$ $B_2 = (10),$ $B_3 = (14).$\\
	${\bf B} =  \mat{1}{0}{1}{-1}{1}{0}{1}{2}{-1},$
	${\bf Q} = \mat{-2}{0}{0}{0}{-1}{0}{0}{0}{2},$
	${\bf a}^T = \mat{0}{1}{1}{2}{2}{-2}{-1}{2}{1},$\\
	${\bf M} = {\bf B\: Q\: a}^T = \mat{-2}{2}{0}{-2}{0}{-1}{-2}{0}{0},$\\
	${\bf b} =  \mat{-1}{-1}{2}{2}{-2}{-2}{-1}{-1}{1},$
	${\bf U} = \mat{-2}{0}{0}{0}{1}{0}{0}{0}{2},$
	${\bf A}^T = \mat{1}{0}{1}{1}{1}{-1}{1}{2}{-2},$\\
	${\bf m = b\: U\: A}^T = \mat{0}{2}{0}{0}{0}{1}{-2}{-2}{-1},$

	\sssec{Definition.}\label{sec-dcentrcoll}
	A collineation is called a {\em central collineation} if
	the collineation transforms every point of a given line into itself,
	and it is not the identity.\\
	The line is called the {\em axis of the central collineation}.

	\sssec{Theorem.}
	{\em Let a collineation be given by 2 complete quadrangles with 2
	fixed points $A_0$ and $B_0$ and 2 other pairs $A_2$, $B_2$ and
	$A_3$, $B_3,$ the necessary and sufficient condition for this
	collineation to be a central collineation is that $A_2 \times A_3,$
	$B_2 \times B_3$ and $A_0 \times A_1$ have a point in common.}

	\sssec{Theorem.}
	{\em In a central collineation, if $B_l$ corresponds to $A_l$ and is
	distinct from $A_l,$ then $A_l \times B_l$ passes through a fixed point
	$F.$}

	\sssec{Definition.}
	$F$ is called the {\em center of the central collineation}.

	\sssec{Definition.}\label{sec-dhomolelat}
	A central collineation is a {\em homology}  iff  its
	center is not on its axis.\\
	A central collineation is an {\em elation}  iff  its center is on
	its axis.

	\sssec{Comment.}
	Theorem \ref{sec-tcoll} or \ref{sec-tcollmat} could serve as an
	alternate definition of collineation.oreo

	\sssec{Exercise.}
	Characterize the matrix of a central collineation, and of an elation.

	\sssec{ Notation.}
	When matrices are used to represent collineations correlations it is
	convenient to have a notation for the inverse matrix scaled by a
	convenient non zero factor, meaning that each entry is multiplied by
	that factor, $N^I$ will be used.
	\ssec{Correlations, polarity.}\label{sec-Scorrel}
	\setcounter{subsubsection}{-1}
	\sssec{Introduction.}
	Correlations have been defined in \ref{sec-dcorrel}.  Their
	algebraic study follows directly from that of collineations.  Their
	importance is due to their intimate relation with conics as will be
	seen in \ref{sec-Sconics}.

	\sssec{Definition.}
	The mapping which associates to the point $(m),$ the line [$m$], for all
	$m = 0 to p^{2k}+p^k$ is called a {\em basic duality}.  It will be
	denoted by $\delta$.

	\sssec{Theorem.}
	{\em The mapping $\delta$  is a correlation.}

	\sssec{Theorem.}\label{sec-tcorr0}
	{\em Given a point collineation $\gamma$  and the corresponding line
	collineation $\gamma'$, then the mapping\\
	\hth$\rho  := \delta  \circ \gamma$\\
	is a point correlation, and the corresponding line correlation is\\
	\hth$\rho' := \delta \circ \gamma'.$\\
	In particular,\\
	if $Q = \gamma(P)$ and $q = \gamma'(p),$ then $\rho(P) = \overline{Q}$
	and $\rho'(p) = \overline{q}$.}

	\sssec{Theorem.}\label{sec-tcorr}
	{\em Given a complete quadrangle $A_j$ and a complete
	quadrilateral $b_j,$ $j = 0,1,2,3,$\\
	Let $a_i := A_{i+1}*A_{i-1}$ and $B_i = b_{i+1}*b_{i-1},$\\
	$A_3 = r_0 A_0 + r_1 A_1 + r_2 A_2,$
	    $b_3 = s_0 b_0 + s_1 b_1 + s_2 b_2$\\
	$q_i := \frac{s_i}{r_i},$ $u_i := q_{i+1} q_{i-1},$\\
	then $q_i := \frac{B_i \cdot b_3}{a_i \cdot A_3}.$\\
	Moreover, the correlation which associates to $A_j,$ $a_j$, $j$ = 0 to
	3, is given by}
	\enumb
	\item {\em the point to line mapping}\\
	$b_l := \rho (A_l) := q_0(a_0.A_l)b_0 + q_1(a_1.A_l)b_1 +
		q_2(a_2.A_l)b_2$,
	\item {\em and the line  to point mapping}
	$\rho '(a_l) := u_0(A_0.a_l)B_0 + u_1(A_1.a_l)B_1 + u_2(A_2.a_l)B_2.$
	\enume

	The proof follows from \ref{sec-tconproj2} and from \ref{sec-tcorr0}.

	\sssec{Example.}
	For p = 5, the correlation $\rho$  defined by\\
	$\rho(0) = [0],$ $\rho(1) = [1],$ $\rho(6) = [12],$ $\rho(12) = [19],$
	= (1,2,3), implies\\
	$a_0 = [0,0,1],$ $a_1 = [0,1,0],$ $a_2 = [1,0,0]$ and\\
	$\rho'(a_0) = (1,0,-1),$ $\rho'(a_1) = (1,-1,0),$ $\rho'(a_2) =
	(-1,0,0).$\\
	$q_0 = 2,$ $q_1 = 1,$ $q_2 = 1,$ 
	$u_0 = 1,$ $u_1 = 2,$ $u_2 = 2,$ therefore\\
	$\rho(X,Y,Z) = [-X, -X-Y, -X-2Z],$
	$\rho'[x,y,z] = (2x-2y-z,2y,z).$

	\sssec{Theorem.}\label{sec-tcorrmat}
	{\em Using the notation of} \ref{sec-tconproj2}, \ref{sec-tconproj3},
	and \ref{sec-tcorr} {\em can be written in matrix notation.\\
	Let ${\bf Q}_{i,i} := Q_i$ and ${\bf Q}_{i,j} := 0$ for $i\neq j,$\\
	let ${\bf R}_{i,i} := q_{i+1} q_{i-1}$ and ${\bf R}_{i,j} := 0$ for
	$i\neq j,$ then\\
	${\bf N} := {\bf b\: R\: a}^T$ defines a correlation ${\bf b}_l =
	{\bf N\: A}_l,$ and ${\bf N}' := {\bf B\: V\: A}^T$ determines
	${\bf B}_l = {\bf N}' {\bf a}_l.$\\
	Moreover ${\bf N}' = {\bf N}^{-1^T}$ is the adjoint matrix.}

%
%
	\sssec{Definition.}
	A {\em polarity} is a correlation which satisfies\\
	\hth$\rho' \circ \rho  = \epsilon$.\\
	In this case $\rho(P)$ is called the {\em polar} of P and
	$\rho'(p)$ is called the {\em pole} of p.

	\sssec{Example.}
	The correlation which associates to $A = (0),$ $B = (1),$
	$C = (6)$ and $P = (13)$ the lines [11], [7], [2] and [15], is
	a polarity and\\
	$\rho(X,Y,Z) = [Y+Z,X+Z,X+Y],$
	$\rho'[x,y,z] = (-x+y+z,x-y+z,x+y-z).$

	\sssec{Theorem.}
	{\em If {\bf M} is a matrix associated to a correlation, then this
	correlation is a polarity  iff  the matrix is symmetric, in
	other words  iff  ${\bf M} = {\bf M}^T.$}

	\sssec{Comment.}
	Theorem \ref{sec-tcorr} or \ref{sec-tcorrmat} could serve as an
	alternate definition of correlations.

	\sssec{Definition.}
	A {\em degenerate line correlation} $\rho_d$ corresponds to a
	function which associates to the set of points in the plane, lines
	which are obtained by multiplying the vector associated to the point to
	the left by the matrix\\
	\hth${\bf D} = \mat{0}{-U_2}{U_1}{U_2}{0}{-U_0}{-U_1}{U_0}{0}$.

	\sssec{Theorem.}
	{\em If $U$ is the point $(U_0, U_1, U_2),$ then {\bf D} associates to
	the point $V = (V_0, V_1, V_2),$ the line $U \times V$.\\
	In the correlation, the image of all points are lines through the point
	$U$ and therefore all lines have $U$ has their image.
	The matrix corresponding to $\rho'_d$ is therefore,\\
	\hth$\mat{U_0}{U_0}{U_0}{U_1}{U_1}{U_1}{U_2}{U_2}{U_2}$.}

	\sssec{Exercise.}
	Prove (Seidenberg, p.193-196)
	\enumb
	\item that a linear transformation is the product of 2 polarities.
	\item that the set of fixed point and fixed lines of a linear
	transformation form a self dual configuration.
	\enume

{\tiny\footnotetext[1]{G23.TEX [MPAP], \today}}

	\ssec{Conics.}\label{sec-Sconics}
	\setcounter{subsubsection}{-1}
	\sssec{Introduction.}
	The following definition was first given by von Staudt.
	The connection between polarity and conics was anticipated already by
	Apollonius and clearly understood by La Hire.

	\sssec{Definition.}
	Given a polarity $\rho$ with inverse $\rho',$ a conic is the set of
	points $P$ such that\\
	\hth$P \cdot \rho(P) = 0.$\\
	and the set of lines $p$ such that\\
	\hth$p \cdot \rho'(p) = 0.$\\
	In other words it is the set of points which are on their polar and
	the set of lines which are on their pole.

	If the polarity corresponds to a symmetric matrix\\
	\hth$\mat{a_0}{b_2}{b_1}{b_2}{a_1}{b_0}{b_1}{b_0}{a_2},$\\
	the equation of the corresponding point conic is\\
	\hth$a_0 X_0^2 + a_1 X_1^2 + a_2 X_2^2
		+ 2 ( b_0 X_1 X_2 + b_1 X_2 X_0 + b_2 X_0 X_1 ) = 0.$

	\sssec{Theorem.}
	\enumb
	\item {\em 5 points no 3 of which are collinear determine a conic.}
	\item the conic through A, B, C, D and E is given by\\
	\hth$k_1 [A\times B]\bigx[C\times D] = k_2 [A\times D]\bigx[B\times C],$
\\
	with\\
	\hth$ k_1 = [A\times D]\cdot E$ . $[B\times C]\cdot E,$
		$k_2 = [A\times B]\cdot E$ . $[C\times D]\cdot E.$
	\enume

	\sssec{Example.}\label{sec-econic5}
	Given the data of \ref{sec-elinecol}, the conic through A, B, C, D and E
	is\\
	\hth$2 X_0^2 - X_1^2 - 4 X_2^2 + 5 X_1 X_2 - 4 X_2 X_0 = 0.$

	\sssec{Exercise.}
	Prove that a conic has $p+1$ points in a finite projective plane
	associate with $p$.

	If we join one point $P$ to the $p$ others we obtain $p$ lines through
	$P$ therefore the left over line is the tangent at $P$.

	\sssec{Comment.}
	For $p = 3,$ the conic has 4 points, hence it cannot be constructed
	by giving 5 points, but it can be constructed if we give
	4 points and a tangent at one of these points or 3 points and
	the tangents at 2 of these points.  See \ref{sec-eqqc3}.\\
	For $p = 2$, a conic can be constructed using 3 non collinear points
	and the tangents at 2 of these points.

	\sssec{Theorem.}
	{\em The pole of} [1,1,1] {\em with respect to the conic\\
	$b_0 X_1 X_2 + b_1 X_2 X_0 + b_2 X_0 X_1 \rho$
	    + $(X_0 + X_1 + X_2) (u_0 X_0 + u_1 X_1 + u_2 X_2)$ = 0,\\
	is\\
	$(b_0(-b_0+b_1+b_2) + 2 u_0b_0 - u_1(b_0+b_1-b_2)
	- u_2(b_0-b_1+b_2),  \ldots  ,  \ldots ).$}

	\sssec{Theorem.}
	{\em The pole of} [1,1,1] {\em with respect to the conic\\
	\hth$	c_0 X_1 X_2 + c_1 X_2 X_0 + c_2 X_0 X_1
		        + u_0 X_0^2 + u_1 X_1^2 + u_2 X_2^2 = 0,$\\
	is\\
	\hth$	(c_0(-c_0+c_1+c_2) - 2 u_1c_1 - 2 u_2c_2 + 4 u_1u_2,  \ldots  ,
	  \ldots  ).$}

	\sssec{Example.}\label{sec-econic13}
	For $p = 13,$
	if $b_0 = 1,$ $b_1 = 6,$ $b_2 = 2,$ $u_0 = -5,$ $u_1 = 4,$ $u_2 = 2,$\\
	then $c_0 = -4,$ $c_1 = 2,$ $c_2 = 5,$ and the pole of [1,1,1] is
	(1,6,3) = (95).

	\sssec{Theorem.}
	{\em Given the conic\\
	\hth$	a_0 X_0^2 + a_1 X_1^2 + a_2 X_2^{}2 + b_0 X_1 X_2 + b_1 X_2 X_0
	+ b_2 X_0 X_1 = 0.$\\
	and a point $(P_0,P_1,P_2),$ with $P_2 \neq 0,$ on the conic, all
	the points are given by\\
	\hth$	X_0 = a_1P_0 u^2 - (2 a_1P_1 + b_0P_2) uv
	- (a_0P_0 + b_2P_1 + b_1P_2) v^2,$\\
	\hth$	X_1 = - (b_2P_0 + a_1P_1 + b_0P_2) u^2
	- (2 a_0P_0 + b_1P_2) uv + a_1P_0 v^2,$\\
	\hth$	X_2 = a_1P_2 u^2 + b_2P_2 uv + a_0 P_2 v^2,$\\
	using the $p+1$ values of the homogeneous pair $(u,v).$}

	Proof:  The points $(v,u,0)$ on [0,0,1] joined to $P$ is the line\\
	\hth$l = [-P_2u,P_2v,P_0u-P_1v].$\\
	$X$ is on $l$  iff
	$P_2X_1v = P_2X_0 u - P_0 X_2 u + P_1X_2 v,$ substituting in the
	equation of the
	conic, if $A$ is the coefficient of $X_0^2$ and $B$ that of $X_2^2,$ 
	using the property of the products of the roots of the equations gives
	$P_2X_2 = A,$ $P_0X_0 = B,$ this gives $X_0$ and $X_2,$ substituting
	in $l$ gives $X_1.$ 

	\sssec{Theorem. [Chasles]}
	{\em Given the configuration of Desargues} \ref{sec-tdes} {\em there
	exists a conic such that $A_i$ is the pole of
	$b_i := B_{i+1} \times B_{i-1}$ and
	vice-versa.\\
	Clearly $B_i$ is also the pole of $a_i := A_{i+1} \times A_{i-1}.$}

	Proof:  Let $A_0 = (1,0,0),$ $A_1  = (0,1,0),$ $A_2 = (0,0,1),$
	$C = (1,1,1),$
	$c = (c_0,c_1,c_2)$ and $B_0 = (b,1,1).$  We have\\
	$C_0 = (0,c_2,-c_1),$ $C_1 = (-c_2,0,c_0),$ $C_2 = (c_1,-c_0,0),$ 
	$b_1 = [c_0,c_1,-bc_0-c_1],$ $b_2 = [c_0,-bc_0-c_2,c_2],$ 
	$B_1 = (c_1,(b-1)c_0+c_1,c_1),$ $B_2 = (c_2,c_2,(b-1)c_0+c_2).$ 
	The transformation which associates to $a_i,$ $k_i B_i,$ with\\
	$k_0$ := $c_1c_2,$ $k_1$ := $c_2,$ $k_2$ := $c_1$ is\\
	\hth$\matb{(}{ccc}bc_1c_2&c_1c_2&c_1c_2\\
		c_1c_2&((b-1)c_0+c_1)c_2&c_1c_2\\
		c_1c_2&c_1c_2&((b-1)c_0+c_2)c_1\mate{)},$\\
	is a line to point polarity because the representative matrix is
	symmetric.  Its inverse can easily be obtained by determining $b_0.$ 
	This is left as an exercise.

	\sssec{Notation.}
	If $u = [u_0,u_1,u_2]$ and $v = [v_0,v_1,v_2]$ are 2 lines, then\\
	$u \bigx v = (u_0 X_0 + u_1 X_1 + u_2 X_2)
	(v_0 X_0 + v_1 X_1 + v_2 X_2).$ 

	\sssec{Definition.}
	Given 2 conics $\alpha$  and $\beta$ if there exist integers $k$ and $l$
	and lines $u$ and $v$ such that\\
	\hth$	k \alpha  + l \beta  = u \bigx v,$\\
	then $v$ is called the {\em radical axis with respect to} $u$
	{\em of $\alpha$  and $\beta$ }.

	\sssec{Lemma.}\label{sec-lsymmmat}
	{\em If {\bf N} is a symmetric matrix and {\bf A} and {\bf B} are 2
	vectors, then ${\bf A} \cdot ({\bf N}\:{\bf B}) = {\bf B} \cdot
	({\bf N}\:{\bf A})$.}

	\sssec{Theorem.}
	{\em A conic or the corresponding polarity determines an involution
	on every line, by associating to each point its conjugate on that line.
	Moreover if $A_0$ and $B_0$ are conjugates as well as $A_1$ and $B_1,$
	if $A_l = t_0 A_0 + t_1 A_1$ its conjugate $B_l$ is given by
		$B_l = ( (A_1 . B_0) t_0 + (A_1 . B_1) t_1 ) A_0 
		     - ( (A_0 . B_0) t_0 + (A_0 . B_1) t_1 ) A_0.$}\\
	This property follows from the notion of conjugates
	and from $B_l = (A_1 * A_0) * (t_0 N A_0 + t_1 N A_1),$
	with $B_0 = N A_0$ and $B_1 = N A_1.$  The Lemma confirms the
	involutive property.

	\sssec{Example.}\label{sec-equaqua}
	For $p = 5,$ starting with $A(0) = (6),$ $A(1) = 1,$ $A(2) = 0,$
	$A(3) = 12,$ the quadrangle-quadrilateral configuration is\\
	$a_{1,2} = [6],$ $a_{2,0} = [1],$ $a_{0,1} = [0],$
	$a_{0,3} = [5],$ $a_{1,3} = [10],$ $a_{2,3} = [26],$
	$D_0 = (2),$ $D_1 = (7),$ $D_2 = (11),$
	$d_0 = [30],$ $d_1 = [27],$ $d_2 = [15],$
	$A_{0,3} = (5),$ $A_{1,3} = (10),$ $A_{2,3} = (26),$
	$A_{1,2} = (24),$ $A_{2,0} = (17),$ $A_{0,1} = (13),$
	$a_0 = [24],$ $a_1 = [17],$ $a_2 = [13],$ $a_3 = [12].$

	\sssec{Example.}
	The points and lines of the extended quadrangle-\nolinebreak
	quadrilateral configuration are those of \ref{sec-equaqua} and
	$B_{1,0} = (9),$ $B_{2,1} = (16),$ $B_{0,2} = (3),$
	$B_{2,0} = (21),$ $B_{0,1} = (4),$ $B_{1,2} = (8),$
	$B_{0,3} = (15),$ $B_{1,3} = (27),$ $B_{2,3} = (30),$
	$B_{3,0} = (18),$ $B_{3,1} = (22),$ $B_{3,2} = (14),$
	$b_{1,0} = [4],$ $b_{2,1} = [8],$ $b_{0,2} = [21],$
	$b_{2,0} = [3],$ $b_{0,1} = [9],$ $b_{1,2} = [16],$
	$b_{0,3} = [2],$ $b_{1,3} = [7],$ $b_{2,3} = [11],$
	$b_{3,0} = [18],$ $b_{3,1} = [22],$ $b_{3,2} = [14].$

	\sssec{Example.}
	The conical points and lines of the extended quadrangle-
	quadrilateral configuration are the 6 points and 6 lines
	$C_{1,0} = (20),$ $C_{2,1} = (29),$ $C_{0,2} = (19),$
	$C_{2,0} = (28),$ $C_{0,1} = (23),$ $C_{1,2} = (25),$
	$c_{1,0} = [20],$ $c_{2,1} = [29],$ $c_{0,2} = [19],$
	$c_{2,0} = [28],$ $c_{0,1} = [23],$ $c_{1,2} = [25].$

	\sssec{Definition.}
	A {\em degenerate conic} is a set of points and lines represented by
	an equation corresponding to a singular 3 by 3 symmetric matrix.

	\sssec{Exercise.}
	Describe all the types of degenerate conics.

	\sssec{Exercise.}
	The number of conics, degenerate or not is $(q^2+q+1)(q^3+1),$\\
	The number of degenerate conics are\\
	$\begin{array}{ll}
	line \bigx line& q^2+q+1,\\
	line_1 \bigx line_2& \frac{1}{2}(q^2+q+1)q(q+1)\\
	non \:real\: line  \bigx  its\: conjugate&
	\frac{1}{2}(q^2+q+1)q(q-1),\\
	\end{array}$\\
	The number of non degenerate conics is $q^5-q^2.$ 

	\sssec{Table.}
	$\begin{array}{lrrrrrr}
	q&                                2&   3&    4 &   5&     7&      11\\
	\hline\\
	line \bigx line&                  7&  13&   21&   31&    57&     133\\
	non\:real\:line\bigx its\:conjugate&7&   39&  126&  310& 1197&    7315\\
	line_1 \bigx line_2&             21&   78&  210&  465&  1596&   8778\\
	non\:degenerate\:conics&         28& 234&  1008& 3100& 16758& 160930\\
	all\:conics&                     63& 364&  1365& 3906& 19608& 177156
	\end{array}$

	\ssec{The general conic.}

	\setcounter{subsubsection}{-1}
	\sssec{Introduction.}
	There is a more general connection between correlations
	and conics, which leads to the concept of a general conic, which is one
	of 4 types, the points of a conic of von Staudt and the lines of an
	other conic of von Staudt.  It has $p+1$ points and $p+1$ lines; a
	degenerate conic consisting of $2p+1$ points on 2 distinct lines and
	$2p+1$
	lines through 2 distinct points; a degenerate conic consisting of $p+1$
	points on 1 line and of $p+1$ lines through 1 point; and finally the
	degenerate conic consisting of 1 point and 1 line.
	In the last case, in complex projective geometry, all the complex
	points are on a pair of complex conjugate lines and all the complex
	lines are through a pair of complex conjugate points.
	To every correlation is associated a general point conic and a general
	line conic.

	\sssec{Definition.}
	A {\em general conic} consists of a point conic which is the set
	of points in a correlation which are on their image and of a line conic
	which is the set of lines in a correlation which are on their image.

	\sssec{Theorem.}
	{\em If {\bf N} is the matrix associated to a correlation, the equation
	of the point conic is\\
	\hth${\bf X}^T {\bf N\:X} = 0,$ where {\bf X} is the vector
	$(X_0,X_1,X_2).$\\
	The equation of the line conic is\\
	\hth${\bf x}^T {\bf N}^{-1} {\bf x} = 0,$ where {\bf x} is the vector
	$(x_0,x_1,x_2).$}

	\sssec{Theorem.}
	{\em Let {\bf A} be the most general antisymmetric matrix,\\
	\hth${\bf A} = \mat{0}{-w}{v}{w}{0}{-u}{-v}{u}{0},$\\
	all the correlations associated to ${\bf N} + {\bf A}$ define the same
	point conic.}

	\sssec{Definition.}
	Given a matrix {\bf N}, its symmetric part ${\bf N}^S$ is defined by\\
	\hth${\bf N}^S := \frac{{\bf N} + {\bf N}^T}{2},$\\
	and its {\em antisymmetric part} ${\bf N}^A$ by\\
	\hth${\bf N}^A := \frac{{\bf N} - {\bf N}^T}{2}.$

	\sssec{Theorem.}
	{\em Given a correlation $\rho,$ $\rho'$.\\
	If $T$ is on the point conic, then $\rho(T)$ is on the line conic.\\
	If $t$ is on the line conic, then $\rho'(t)$ is on the point conic.\\
	The general conic degenerates if $det({\bf N}) = 0.$\\
	The center corresponds to the vector which is the homogeneous
	solution of ${\bf N}^S {\bf C} = 0,$ and the central line, to that of}
	$({\bf N}^S)^{-1} {\bf c} = 0.$

	\sssec{Definition.}
	Given a general conic, the tangent  $t$  at the point  $T$  of
	the point conic is defined by  $t := {\bf N}^S T.$
	The {\em contact}  $T$  of a line  $t$   which belongs to a line conic
	is defined by  $T := ({\bf N}^S)^{-1} t.$

	\sssec{Theorem.}
	{\em If the correlation is a polarity then the tangent at a point
	$T$  of a point conic is on the corresponding line conic.  Similarly,
	the contact of a line $t$ of line conic is on the corresponding point
	conic.}

	\sssec{Theorem.}
	{\em If a conic is non degenerate, the necessary and sufficient
	condition for the set of tangents to a point conic to coincide with the
	set of lines on the line conic is that the correlation be a polarity.}

	Proof:  Let {\bf N} be the matrix associated to the correlation.  The
	line conic is ${\bf x}^T {\bf N}^{-1} {\bf x} = 0.$  The tangents
	${\bf t} = {\bf N}^S {\bf X}$ to the point conic are on\\
	${\bf t}^T ({\bf N}^S)^{-1} {\bf N} ({\bf N}^S)^{-1}\:{\bf t} = 0.$\\
	If $\rho$ is a polarity, ${\bf N} = {\bf N}^S$ and ${\bf N}^{-1} =
	{\bf N}^{S^-1} {\bf N} ({\bf N}^S)^{-1}.$  Vice-versa,\\
	if ${\bf N}^{-1} = ({\bf N}^S)^{-1} {\bf N} ({\bf N}^S)^{-1}$ then
	${\bf N}^S = {\bf N} ({\bf N}^S)^{-1} {\bf N}.$\\
	Therefore, using ${\bf N} = {\bf N}^S + {\bf N}^A,$ 
	$2 {\bf N}^A = - {\bf N}^A ({\bf N}^S)^{-1} {\bf N}^A,$ \\
	transposing, $- 2 {\bf N}^A = - {\bf N}^A ({\bf N}^S)^{-1} {\bf N}^{A,}$
	because ${\bf N}^{A^T} = - {\bf N}^A,$ therefore\\
	${\bf N}^A = 0,$ ${\bf N}$ is symmetric and therefore $\rho$ is a
	polarity.

	\ssec{The Theorem of Pascal and Brianchon.}
	\setcounter{subsubsection}{-1}
	\sssec{Introduction.}
	A fundamental theorem associated to conics was discovered
	by Blaise Pascal.  It allows construction of any point on a conic given
	by 5 points and in particular the other intersection of a line through
	one point of a conic.  See I, \ldots.

	There is a general principle of linear construction that if a point or
	line is uniquely define, that point or line can be obtained by a linear
	construction.  The points of intersection of a conic with a general
	line are not uniquely defined and therefore do not admit a linear
	construction on the other hand if the line passes to a known point
	of the conic, the other intersection of the line and the conic is 
	uniquelly defined.  The Pascal construction of \ref{sec-cpascal} is a
	solution to this problem which follows from the following Theorem.

	\sssec{Theorem [Pascal].}\label{sec-tpascal}
	{\em If 6 points $A_0,$ $A_1,$ $A_2,$ $A_3,$ $A_4,$ $A_5$ are on
	a conic and the Pascal points are defined as\\
	\hth$P_0 := (A_0 \times A_1) \times (A_3 \times A_4),$\\
	\hth$P_1 := (A_1 \times A_2) \times (A_4 \times A_5),$\\
	\hth$P_2 := (A_2 \times A_3) \times (A_5 \times A_0),$\\
	then the points $P_0,$ $P_1,$ $P_2$ are collinear}
	(Pascal, 1639, Lemma 1 and 3).

	There are "degenerate" forms of this theorem in which 2 consecutive
	points coincide and the cord is replaced by the tangent at these
	points for instance {\em if the tangent at $A_0$ is $t_0,$ the Pascal
	points are\\
	\hth$P_0 := t_0 \times (A_3 \times A_4),$\\
	\hth$P_1 := (A_0 \times A_2) \times (A_4 \times A_5),$\\
	\hth$P_2 := (A_2 \times A_3) \times (A_5 \times A_0),$\\
	and the points $P_0,$ $P_1,$ $P_2$ are collinear}.

	Proof:  The Theorem of Pascal will now be proven in the 4 cases,
	6 points, 5 points and the tangents at one of them, 4 points and the
	tangents at 2 of them and finally 3 points and their tangents.  In each
	case, the coordinates will be chosen to simplify the algebra.
	See also \ref{sec-tvecpa}.
	\enumb
	\item Let the 6 points of the conic be $A_0,$ $C_0,$ $A_1,$ $B_1,$
		$A_2,$ $B_2.$ Choose the coordinates such that $A_0 = (1,0,0),$
		$A_1 = (0,1,0)$ and $A_2 = (0,0,1),$ choose the barycenter
		$M = (1,1,1)$ at the intersection of $A_1 \times B_1$ and
		$A_2 \times B_2,$ let the line $A_0 \times B_0$ be $[0,r,-s].$\\
		Because the conic passes through $A_i,$ it has an equation of
		the form\\
		0. $   u X_1 X_2 + v X_2 X_0 + w X_0 X_1 = 0.$\\
		$A_1 \times B_1 = A_1 \times M = [1,0,-1],$
		$A_2 \times B_2 = A_2 \times M = [1,-1,0],$ therefore,\\
	\hth$B_1 = (u+w,-v,u+w), B_2 = (u+v,u+v,-w),$\\
		$C_0 = (-urs,s(vr+ws),r(vr+ws)),$\\
		hence the Pascal points are\\
	\hth$P_0 = (A_0 \times C_0) \times (B_1 \times A_2)
		= (s(u+w),-vs,-vr),$\\
	\hth$P_1 = (C_0 \times A_1) \times (A_2 \times B_2)
		= (urs,urs,-r(vr+ws)),$\\
	\hth$P_2 = (A_1 \times B_1) \times (B_2 \times A_0)
		= (w,-(u+v),w).$\\
		which are all on $[v(ru+rv+sw),w(su+rv+sw),su(u+v+w)].$
	\item Let the 5 points of the conic be $A_0,$ $A_1,$ $B_1,$ $A_2,$
		$B_2,$ and let the tangent $t$ be chosen at $A_0.$\\
		With the coordinates chosen as above and the conic again of the
		form 0.0, the tangent is $[0,w,v]$ and the Pascal points are\\
	\hth$P_0 = t \times (B_1 \times A_2) = (u+w,-v,w),$\\
	\hth$P_1 = (A_0 \times A_1) \times (A_2 \times B_2) = (1,1,0),$\\
	\hth$P_2 = (A_1 \times B_1) \times (B_2 \times A_0) =  (w,-(u+v),w).$\\
		which all are on $[-w,w,u+v+w].$
	\item Let the 4 points be $A_0,$ $A_1,$ $A_2$ and $B_0$ and the
		tangents be $t_1$ at $A_1$ and $t_2$ at $A_2.$ Choose the
		coordinates as above, except for $M$ on $A_0 \times B_0
		= [0,1,-1],$ then $r = s = 1.$
	\hth$B_0 = (-u,v+w,v+w),$\\
		the tangents are $t_1 = [w,0,u]$ and $t_2 = [v,u,0].$\\
		The Pascal points are\\
	\hth$P_0 = (A_0 \times A_1) \times t_2 = (-u,v,0),$\\
	\hth$P_1 = (A_1 \times A_2) \times (A_0 \times B_0) = (0,1,1),$\\
	\hth$P_2 = t_1 \times (A_2 \times B_0) = (-u,v+w,w),$\\
		which all are on $[v,u,-u].$
	\item Let the points be $A_0,$ $A_1$ and $A_2$ and the tangents be
		those at these points, using again 0.0. as the equation of
		the conic, the tangents are\\
	\hth$        t_0 = [0,w,v],$ $t_1 = [w,0,u],$ $t_2 = [v,u,0].$\\
		the Pascal points are\\
	\hth$        P_0 = t_0 \times (A_1 \times A_2) = (0,v,-w),$\\
	\hth$        P_1 = t_1 \times (A_2 \times A_0) = (-u,0,w),$\\
	\hth$        P_2 = t_2 \times (A_0 \times A_1) = (u,-v,0),$\\
		which are all on $[vw,wu,uv].$
	\item The other cases, 4 tangents and 2 points of contact, 5 tangents
	and 1 point  of contact, 6 tangents, can be proven by duality.
	\enume

	\sssec{Theorem [Pascal].}
	The reciprocal of the preceding Theorem is true.  In other words,\\
	{\em if the Pascal points $P_i$ are collinear, the 6 points $A_k$ are
	on a conic.}

	The proof is left as an exercise.

	\sssec{Notation.}\label{sec-nconics}
	The property that 6 points are on a conic $\gamma$ will be noted\\
	\hth$incidenceconic(A,B,C,D,E,F[,\gamma])$ or
		$incidenceconic(A_k[,\gamma]).$\\
	Similar notation will be used for degenerate or for dual forms, for
	instance\\
	\hth$incidenceconic(A,t,B,C,D,E)$, where $t$ is the tangent at A.\\
	\hth$incidenceconic(a,b,c,d,e,f)$ where $a,b,c,d,e,f$ are 6 tangents
	to the conic.\\
	Theorem \ref{sec-tpascal} will be denoted as follows.\\
	No.\hti{6}$Pascal(A_k[,a_k];\langle P_i[,p]\rangle)$\\
	Hy.\hti{6}$incidenceconic(A_k).$\\
	De.\hti{6}$P_i := (A_i \times A_{i+1}) \times (A_{i+3} \times A_{i+4})
$\\
	Co.\hti{6}$\langle P_i,p\rangle.$

	If $t_k$ is the tangent at $A_k$, $A_k$ is followed by $t_k$.

	The Pascal line associated to the points $A_k$ will be denoted by
	p := Pascal$(A_k).$

	\sssec{Definition.}
	The dual of the Theorem of Pascal is called the {\em Theorem of
	Brianchon}.  Brianchon discovered the Theorem before Gergonne discovered
	the important principle of duality.\\
	In the degenerate case of a triangle inscribed in a conic and of
	the triangle outscribed to the conic at these points
	(\ref{sec-tpascal}.3), the line is called the {\em Pascal line of the
	triangle} and the point of the dual Theorem, the {\em Brianchon point
	of the triangle}. von Staudt (1863) calls them, pole and polar of the 
	triangle.

	\sssec{Theorem. [Generalization of von Staudt]}
	{\em If $p$ is the Pascal line of the hexagon $A_0,A_1,A_2,A_3,A_4,A_5$
	inscribed in a conic $\gamma$ and $P$ is the Brianchon point of the
	outscribed hexagon formed by the tangents at $A_j$, then $P$ is the pole
	of $p$.}

	The proof follows at once from the properties of poles and polars.

	\sssec{Corollary. [von Staudt]}
	{\em If $\{A_0,A_1,A_2\}$ is a triangle inscribed in a conic, then its
	Pascal line is the polar of its Brianchon point.}

	\sssec{Theorem. [von Staudt]}
	If 2 triangles $\{A_0,A_1,A_2\}$ and $\{B_0,B_1,B_2\}$ are inscribed in
	a conic $\gamma$ and are perspective with center $C$ and axis $c$,
	and $P$, $Q$ are their Brianchon pointa and $p$ and $q$ are their Pascal
	lines, tehn $\langle P,Q,C;pq\rangle$ and $\langle p,q,c;PQ\rangle$,
	moreover, quatern($P,Q,C,pq\times c$) and quatern($p,q,c,PQ\times C$).

	\sssec{Notation.}
	\hth$(A_{i,j,k}) := det(A_i,A_j,A_k).$

	\sssec{Theorem.}
	{\em If 6 points $A_k,$ $k = 0$ to 5, are on a conic then\\
	$  (A_{k+2,k+3,k+1}) (A_{k+3,k+4,k}) (A_{k+4,k+5,k+1}) (A_{k+5,k,k+2})\\
	\hth -(A_{k+2,k+3,k})(A_{k+3,k+4,k+1})(A_{k+4,k+5,k+1})(A_{k+5,k,k+2})\\
	\hth +(A_{k+3,k+4,k+1})(A_{k+4,k+5,k+2})(A_{k+5,k,k+2})(A_{k,k+1,k+3})\\
	\hth -(A_{k+3,k+4,k+1})(A_{k+4,k+5,k+2})(A_{k+5,k,k+3})(A_{k,k+1,k+2})
	= 0.$\\
	The addition for the subscript is done modulo 6.}

	Proof:  The Theorem will be proven for $k = 0.$
	Let $a_k := A_k \times A_{k+1}.$
	The Pascal points are $P_k := a_k \times a_{k+3},$ we have\\
	\hth$(P_{0,1,2}) := det(P_0,P_1,P_2) = (P_0 * P_1) \cdot P_2 = 0,$\\
	but\\
	$P_0 = (A_0 * A_1) * (A_3 * A_4)
	     = det(A_0,A_3,A_4) A_1 - det(A_1,A_3,A_4) A_0,$\\
	\hth$     = (A_{0,3,4}) A_1 - (A_{1,3,4}) A_0,$\\
	similarly,\\
	\hth$P_1 = (A_{1,4,5}) A_2 - (A_{2,4,5}) A_1,$ \\
	\hth$P_2  = (A_{2,5,0}) A_3 - (A_{3,5,0}) A_2
	     = (A_{0,2,5}) A_3 - (A_{0,3,5}) A_2,$\\
	therefore\\
	$det(P_0,P_1,P_2) =
	 (A_{0,3,4}) (A_{1,4,5}) (A_{0,2,5}) (A_{1,2,3})
	        - (A_{1,3,4}) (A_{1,4,5}) (A_{0,2,5}) (A_{0,2,3})$\\
	\hth$        + (A_{1,3,4}) (A_{2,4,5}) (A_{0,2,5}) (A_{0,1,3})
	        - (A_{1,3,4}) (A_{2,4,5}) (A_{0,3,5}) (A_{0,1,2}) = 0.
	\Box$

	\sssec{Construction.}\label{sec-cpascal}
	\hth point Pascal($A_0,A_1,A_2,A_3,A_4,A'_5;[P_0,P_1,P_2,]A_5)$\\
	is used as an abbreviation for the Pascal construction\\
	$P_0 := (A_0\times A_1)\times(A_3\times A_4),$\\
	$P_1 := (A_1\times A_2)\times(A_4\times A'_5),$\\
	$P_2 := (A_2\times A_3)\times(P_0\times P_1),$\\
	$A_5 := (A_4\times A'_5)\times(P_2\times A_0).$\\
	It gives the point $A_5$ on the conic through $A_0$ to $A_5$ on the line
	$A_4\times A'_5$.
	\hth line Pascal($a_0,a_1,a_2,a_3,a_4,a'_5;[p_0,p_1,p_2,]a_5)$\\
	is used for the dual construction.

	\sssec{Theorem.}\label{sec-tconic2}
	{\em When $p = 2$, the points and lines of a conic2 configuration}
	\ref{sec-dconic2} {\em are the points and lines of a conic.}

	Proof:
	The Pascal points are the diagonal points of the complete
	quadrangle configuration which are collinear because of
	Theorem \ref{sec-tproj2}.

	\sssec{Theorem.}\label{sec-tqqc3}
	{\em Let $p = 3,$ in a quadrangle quadrilateral configuration,
	$Q_i,P,q_i,p$} (\ref{sec-dsd}), {\em there is a conic whose tangent at
	$P$ is $p$ and at $Q_i$ is $p_i.$\\
	In other words the elements of a conic3 configuration}
	(\ref{sec-dconic3} {\em are the points and lines of a conic.}

	Proof:  From \ref{sec-eqqc3}, $Q_i$ is on $q_i,$ the Pascal-Brianchon
	theorem gives\\
	Pascal$(P,p,Q_{i+1},q_{i+1},Q_{i-1},q_{i-1};\langle R_0,P_{i-1},P_{i+1}
	,p\rangle).$

	\sssec{Theorem.}\label{sec-tconical}
	{\em The conical points of Definition} \ref{sec-dconical} {\em are on
	a conic and the conical lines are on a conic.}

	Proof: Pascal$(AF_0,FA_2,AF_0,FA_2,AF_0,FA_2;\langle R_1,R_2,R_3,p
	\rangle.$

	\sssec{Theorem.}
	{\em The following points of the extended Pappus configuration
	are on a conic, 1 point on each of the lines $d,$ $\overline{d},$
	say $M_0$ and $\overline{M}_0$ and
	the intersection with the lines joining the other points say $a_1$
	and $a_2$ 
	with the lines joining $M_0$ or $\overline{M}_0$ with the other points
	on $d$ or $\overline{d}.$\\
	This gives the 18 conics}
	\enumb
	\item$M_i,\overline{N}_{i+1},N_{i+1},
	\overline{M}_i,\overline{N}_{i-1},N_{i-1},$
	\item$M_i,\overline{P}_{i-1},P_{i-1},
	\overline{M}_i,\overline{P}_{i+1},P_{i+1},$
	\item$M_{i+1},\overline{N}_{i+1},L_{i-1},
	\overline{M}_{i-1},\overline{N}_{i-1},L_{i-1},$
	\item$M_{i+1},Q_{i-1},\overline{P}_{i-1},\overline{M}_{i-1},Q_{i+1},
	\overline{P}_{i-1},$
	\item$M_{i-1},Q_{i+1},P_{i+1},\overline{M}_{i+1},Q_{i-1},P_{i-1},$
	\item$M_{i-1},L_{i+1},\overline{N}_{i+1},\overline{M}_{i+1},L_{i-1},
	\overline{N}_{i-1}.$
	\enume

	Proof:  This follows from Pascal's Theorem applied to the points in the
	given order, order which was chosen in such a way that the Pascal line
	was always $\overline{m}_0,$ containing $P_0,$ $\overline{P}_0,$
	$Q_0$ and $\overline{D}.$
	Exchanging $N_{i+1}$ and $\overline{N}_{i-1}$ for 0. gives an other
	Pascal line $m_0.$ 
	The 9 Pappus lines are therefore the Pascal lines of the 18 conics.

	\sssec{Theorem.}
	{\em The conics 0. and 1. have the same tangent at their common
	point.}

	Proof:  The coefficients of conic 0. are, for $i = 0,$
	$a_0 = m_0m_1m_2,$\\
	$a_1 = m_1^2m_2,$ $a_2 = m_2^2m_1,$ 
	$b_0 = m_1m_2(m_1+m_2),$ $b_1 = m_2(m_1^2+m_2m_0),$\\
	$b_2 = m_1(m_2^2+m_0m_1),$\\
	This follows easily because $M_0,$ $\overline{M}_0$ are on $a_0,$
	giving $a_1,$ $a_2$ and $b_0,$
	$N_1,\overline{N}_1$ are on $a_1,$ this gives $a_0$ and
	$b_1,$ $N_2,$ $\overline{N}_2$ are on $a_2,$ giving $b_2.$ 
	The coefficients of conic 1, for $i = 0$ are the same except for
	$a_0 = m_0(m_1^2-m_1m_2+m_2^2).$  The algebra is simplified by
	noting that the
	equation for $P_1$ and $\overline{P}_1,$ gives after subtraction
	$b_1$ from $a_2$ and $b_0,$ 
	and for $P_2$ and $\overline{P}_2,$ gives after subtraction $b_2$
	from $a_1$ and $b_0.$ 

	The Theorem follows at once.  The tangent at $M_0$ is\\
	$[m_0(m_1+m_2-m_1m_2,m_1m_2,m_1m_2]$ and at $\overline{M}_0$ is
	$[m_1+m_2-m_0,m_1,m_2].$

	\sssec{Exercise.}\label{sec-esalm}
	Study the configuration of all 18 conics associated to the
		extended Pappus configuration.

	\ssec{The Theorems of Steiner, Kirkman, Cayley and Salmon.}
	\label{sec-SSKC}
\setcounter{subsubsection}{-1}
	\sssec{Introduction.}
	The set of Theorems given here originates with the work of Steiner
	(1828, 1832 - Werke I, p.451).
	Proofs have been given using Pascal's Theorem and Desargues Theorem
	in the plane or starting with properties of the configuration\\
	5 * 3 \& 5 * 3 in three (Cremona, 1877) or four
	(Richmond, 1894, 1899, 1900, 1903) dimensions, subjected to a linear
	condition. An alternate approach starts with the work of Sylvester
	1844 (Papers, I, p.92), 1862 (II, p265.) For a good summary, see
	Salmon, 1879, p. 379-383, Baker, II, (2d Ed. 1930), p. 219-236
	and Friedrich Levi, 1929, p.192-199..\\
	The cyclic permutation notation allows the results to be given in a
	simple algebraic way and suggests the related synthetic construction.

	\sssec{Definition.}
	Given 6 points $A_j$, $j = 0$ to 5, on a conic, a {\em conical hexagon}
	abbreviated here by {\em hexagon} is a permutation $h$ of 0 to 5.
	Given $h$ this defines a specific Pascal line\\
	$p(h) := Pascal(A_{h(0)},A_{h(1)},A_{h(2)},A_{h(3)},A_{h(4)},A_{h(5)}),
	$\\
	A {\em map} will denote here a permutation which acts on $h$.

	\sssec{Example.}
	Let $h = [013524] = (\begin{array}{c}012345\\013524\end{array})
	= (2354).$ The ordered set of point associated to $h$ is
	$A_0A_1A_3A_5A_2A_4$. The map $\sigma$ = (135) associates to this set,
	the set $(2354)(135) = (15)(234) = [053421] = h'$ or
	$A_0A_5A_3A_4A_2A_1$, for
	instance, $h'(2) = h\sigma(2) = h(2) = 3,$ $h'(3) = h\sigma(3) = h(5) =
	4.$ The multiplication of permutations is done from right to left.

	\sssec{Definition.}
	\enumb
	\item The {\em Steiner map} is $\sigma$ = (135),
	\item the {\em Steiner conjugate map} is $\gamma$ = (35).
	\item the {\em Kirkman map} is $\kappa$ = (021)(345),
	\item the {\em Cayley-Salmon map} is $\chi$ = (14),
	\item the {\em Salmon map} is $\lambda$ = (2354).
	\item the {\em line-Steiner maps} are $\sigma_0$ = (23)
		and $\sigma_1$ = (45).
	\enume

	\sssec{Theorem.}
	{\em Given $r = (012345)$ and $s = (05)(14)23),$\\
	$h = (\ldots ij\ldots),$ $r^{-1}hr = r^{-1}(r\ldots i-1,j-1\ldots)$ and
	$s^{-1}hs = (\ldots s(i),s(j)\ldots)$ have the same Pascal line.}\\
	The permutations $r^{-k}hr^k$ and $s^{-1}hs$ are called {\em Pascal
	equivalent}.

	\sssec{Theorem.}
	\enumb
	\item (024), (042), (153) {\em are Pascal equivalent to the Steiner
		map} (135).
	\item (02), (04), (13), (15), (24) {\em are Pascal equivalent to the
		Steiner conjugate map} (35).
	\item (012)(354),(015)(243), (045)(132), (051)(234), (054)(123)
		{\em are Pascal equivalent to the Kirkman map} (021)(345).
	\item (03), (25) {\em are Pascal equivalent to the Cayley-Salmon map}
		(14).
	\item (0132), (0215), (0451), (0534), (1243) {\em are Pascal equivalent
		to the Salmon map} (2354).
	\item (024), (042), (153) {\em are Pascal equivalent to the line-Steiner
		maps} (23) {\em and} (45).
	\enume

	\sssec{Theorem. [Steiner (Pascal)]}
	\enumb
	\item $\langle p(h), p(h\sigma), p(h\sigma^2);S(h)\rangle,$
		$S(h)$ is called the {\em Steiner point of} $h$.
	\item $S(h\gamma)$, called the {\em Steiner conjugate point of} $h$,
		{\em is on the polar of $S(h)$ with respect to the conic.}
	\item {\em there are 10 pairs of conjugate Steiner points.}
	\enume
	See \ref{sec-tpap1}.

	Proof: Let h = [012345] = (). I will use here the abbreviations\\
	\hth$ij$ for the line $A_i\times A_j$,\\
	\hth$ijkl$ for the Pascal point $(A_i\times A_j\times (A_k\times A_l).$
\\
	$p(h) = P_0\times P'_0,$ with $P_0 = 0134,$
	$P'_0 = 0523,$\\
	$p(h\kappa) = P_1\times P'_1,$ with $P_1 = 0125,$
	$P'_1 = 1423,$\\
	$p(h\kappa^2) = P_2\times P'_2,$ with $P_2 = 2534,$
	$P'_2 = 0514,$\\
	Let $Q_0 := (P_1\times P_2) \times (P'_1\times P'_2) = 2514,$\\
	$Q_1 := (P_2\times P_0) \times (P'_2\times P'_0) = 3450,$\\
	$Q_2 := (P_0\times P_1) \times (P'_0\times P'_1) = 0123,$\\
	Pascal($A_1A_4A_3A_2A_5A_0;\langle Q_0,Q_1,Q_2\rangle),$ therefore\\
	Desargues$^{-1}(\{P_0,P_1,P_2\},\{P'_0,P'_1,P'_2\};
		\langle Q_0,Q_1,Q_2\}\rangle,S(h))$,\\
	or\\
	Desargues$^{-1}(\{2435,0312,0514\},\{p34125,p25135,p35124\}
		\{p13245\times p12354,1245,0523\},\\
	\hti{12}\{p12345,p13245,p12354\};
		\langle S(e),S(\sigma),S(\sigma^2\}\rangle,s(e))$,

	\sssec{Theorem. [Kirkman 1849, 1850]}
	\enumb
	\item $\langle p(h), p(h\kappa), p(h\kappa^2);K(h)\rangle,$
		$K(h)$ is called the {\em Kirkman point of} $h$.
	\item {\em there are 60 Kirkman points which are 3 by 3 on the 60
		Pascal lines, giving a configuration of type} 60 * 3
			$\&$ 60 * 3\footnote{Levi, p. 194}.
	\enume

	Proof: Let h = [012345] = (). The proof, for i = 0 is as follows.\\
	$p(h) = P_0\times P'_0,$ with $P_0 = 0134,$
	$P'_0 = 0523,$\\
	$p(h\kappa) = P_1\times P'_1,$ with $P_1 = 0134,$
	$P'_1 = 1245,$\\
	$p(h\kappa^2) = P_2\times P'_2,$ with $P_2 = 2534,$
	$P'_2 = 0312,$\\
	Let $Q_0 := (P_1\times P_2) \times (P'_1\times P'_2) = 0325,$\\
	$Q_1 := (P_2\times P_0) \times (P'_2\times P'_0) = 0145,$\\
	$Q_2 := (P_0\times P_1) \times (P'_0\times P'_1) = 1234,$\\
	Pascal($A_2A_5A_4A_3A_0A_1;\langle Q_0,Q_1,Q_2\rangle),$ therefore\\
	Desargues$^{-1}(\{P_0,P_1,P_2\},\{P'_0,P'_1,P'_2\};
		\langle Q_0,Q_1,Q_2\}\rangle,S(h))$,\\
	or\\
	Pascal(014235) $\implies \langle 1435,0524,0123;p14235\rangle,$\\
	Desargues$^{-1}(p14235,\{0523,1423,0514\},\{14,05,23\},
		\{0134,0135,1235\},\\
	\hti{12}\{35,p12435,01\};\langle p12345,p14523,p21435\}\rangle,K(e))$,

	\sssec{Theorem [Salmon]}
	{\em If 2 triangles have their vertices on a conic, their sides are
	tangent to a conic}\footnote{Salmon, p. 381}.

	Proof:\\
	Desargues$^{-1}(\{0135,0145,0245\},\{45,p14523,01\},
		\{0234,1234,1235\},\\
	\hti{12}\{12,p21435,34\};\langle 1245,K(e),0134;P\rangle)$.

	\sssec{Exercise.}
	Prove $\langle p12345,p125423,p34215\rangle$.

	\sssec{Theorem. [Steiner]}
	\enumb
	\item $\langle S(h),S(h\sigma_0),S(h\sigma_1);s(h)\rangle$,
		$s(h)$ is called the {\em Steiner line of} $h$,
	\item $S(h\sigma_0\sigma_1) \incid s(h),$
	\item {\em there are 15 Steiner lines} $s(h).$
	\enume

	The proofs follows from \ref{sec-esalm}.0 and from the fact that the
	Brianchon lines of the conic inscribed in the 2 triangles are Pascal
	lines of the original conic.

	\sssec{Theorem. [Cayley and Salmon]}
	\enumb
	\item $\langle K(h\chi),K(h\chi),K(h\chi);cs(h)
	\rangle$, $cs(h)$ is called the {\em Cayley-Salmon line of} $h$,
	\item $S(h) \incid cs(h),$
	\item {\em there are 20 Cayley-Salmon lines.}
	\item {\em The 60 Kirkman points, the 20 Steiner points, the 60 Pascal
		lines and the 20 Cayley-Salmon lines form a $80 * 4 \& 80 * 4$
		configuration} (See Levi, p. 199).
	\enume

	\sssec{Theorem.}
	{\em 18 Pascal points and 12 Pascal lines are used in the preceding
	Theorem and these are vertices and sides of 3 complete quadrilaterals.}

	\sssec{Theorem. [Salmon]}
	\enumb
	\item $\langle cs(h),cs(h\lambda),cs(h\lambda^2);Sa(h)\rangle$,
		$Sa(h)$ is called the {\em Salmon point of} $h$,
	\item $cs(h\lambda^3) \incid Sa(h),$
	\item {\em there are 15 Salmon points} $Sa(h).$
	\enume

	\sssec{Theorem.}
	In the preceding Theorem:
	\enumb
	\item {\em Each of the 24 Pascal lines occurs exactly twice.}
	\item {\em The Pascal points of $h$ occur 4 times, the other 30
	Pascal points occur twice.}
	\item {\em The 3 Pascal points of $h,$ the 8 points $S(h\lambda^i),$
	$K_i(h\chi),$ $i = 0,1,2,3$ and the 12 associated Pascal lines
	form a pseudo configuration of type}\\
	\hth$ 3 * 4 + 8 * 3\: \&\: 12 * 3,\: (11).$
	\enume

	\sssec{Example.}
	In all cases $h = e = [012345] = ().$
	\enumb
	\item The Theorem of Steiner.\\
	$\begin{array}{cccc}
	S(e)\incid&S(\sigma_0)\incid&S(\sigma_1)\incid&
	S(\sigma_0\sigma_1)\incid\\
	() = [012345]&(23) = [013245]&(45) = [012354]&(23)(45) = [013254]\\
	(135) = [032541]&(1235) = [023541]&(1345) = [032451]&(12345) =
		[023451]\\
	(153) = [052143]&(1523) = [053142]&(1453) = [042153]&(14523) = [043152]
	\end{array}\\
	\langle\langle p12345,p14523,p34125\rangle,\langle p13245,p14523,p24135
	\rangle,\langle p12354,p15423,p35124\rangle,\\
	\langle p13254,p15432,15234\rangle\rangle.$ (Fig. 200b)
	\item The Theorem of Cayley-Salmon.\\
	$\begin{array}{cccc}
	S(e)\incid&K(\chi)\incid&K(\sigma\chi)\incid&K(\sigma^2\chi)\incid\\
	() = [012345]&(14) = [042315]&(1435) = [042531]&(1453) = [042153]\\
	(135) = [032541]&(024531) = [204153]&(0241) = [204315]&(024351) =
		[204531]\\
	(153) = [052143]&(043512) = [420531]&(045312) = [420153]&(0412) =
		[420315]
	\end{array}\\
	\langle\langle p12345,p14523,p34125\rangle,\langle p42315,p23514,p24135
	\rangle,\langle p13524,p25134,p15342\rangle,\\
	\langle p35124,p21354,p24513\rangle\rangle.$ (Fig. 200b')
	\item The Theorem of Salmon.\\
	Add to Example 1, (Fig. 200b'' and b4)\\
	$\begin{array}{cccc}
	S(\lambda)\incid&K(\lambda)\incid&K(\lambda\sigma\chi)\incid&
		K(\lambda\sigma^2\chi)\incid\\
	(2354) = [013524]&(12354) = [023514]&(12345) = [023451]&(123) = 
		[023145]\\
	(15)(234) = [053421]&(031) = [302145]&(03541) = [302514]&(03451) =
		[302451]\\
	(1423) = [043125]&(02)(1345) = [230451]&(02)(13) = [230145]&(02)(1354) =
		[230514]
	\end{array}\\
	\langle\langle p13524,p12435,p43125\rangle,\langle p23514,p21453,p32154
	\rangle,\langle p15432,p25143,p14523\rangle,\\
	\langle p23145,p24513,p32415\rangle\rangle. (Fig. 200b1)\\
	\begin{array}{cccc}
	S(\lambda)\incid&K(\lambda)\incid&K(\lambda\sigma\chi)\incid&
		K(\lambda\sigma^2\chi)\incid\\
	(25)(34) = [015432]&(134)(25) = [035412]&(1325) = [035241]&(13)(254) = 
		[035124]\\
	(14325) = [045231]&(054231) = [503124]&(052341) = [503412]&(051)(23) =
		[503241]\\
	(12543) = [025134]&(0322)(15) = [350241]&(031542) = [350124]&(034152) =
		[350412]
	\end{array}\\
	\langle\langle p15432,p13254,p25134\rangle,\langle p21453,p31245,p24135
	\rangle,\langle p14253,p34125,p12435\rangle,\\
	\langle p35124,p32415,p41235\rangle\rangle. (Fig. 200b2)\\
	\begin{array}{cccc}
	S(\lambda)\incid&K(\lambda)\incid&K(\lambda\sigma\chi)\incid&
		K(\lambda\sigma^2\chi)\incid\\
	(2453) = [04253]&(15324) = [054213]&(15)(24) = [054321]&(15243) = 
		[054132]\\
	(1245) = [024351]&(0431)(25) = [405132]&(041)(253) = [405213]&(04251) =
		[405321]\\
	(13)(245) = [034152]&(05142) = [540321]&(052)(143) = [540132]&(0532)(14)
		= [540213]
	\end{array}\\
	\langle\langle p14253,p15342,p25143\rangle,\langle p31245,p42315,p32154
	\rangle,\langle p12345,p43125,p13254\rangle,\\
	\langle p23145,p41235,p21354\rangle\rangle.$ (Fig. 200b3)
	\enume

	\sssec{Exercise.}
	\enumb
	\item Give the geometric interpretation of the Theorems in this section.
	\item Determine the pseudo configuration associated to the Theorem of
		Salmon.
	\enume

	\ssec{B\'{e}zier Curves for drawing Conics, Cubics, \ldots.}
\setcounter{subsubsection}{-1}
	\sssec{Introduction.}
	The drawing of curves is facilitated by the notion of B\'{e}zier curves.
	These originate with the work of de Casteljau at Citro\"{e}n in 1959
	and were popularized and generalized by B\'{e}zier. To describe easily
	complicated curves in 2, 3, \ldots dimensions, we start with a
	B\'{e}zier polygon \ref{sec-dcastel} to construct a parametric
	representation of points on the curve iteratively.  The associated
	theory is briefly given here. The curve can be expressed in terms of the
	B\'{e}zier polygon by means of Bernstein polynomials
	(\ref{sec-tbezier}),
	the derivatives and differences of the curve can be similarly
	expressed and related to each other. The example for a curve
	whose $i$-th coordinates can be approximated by cubic polynomials is
	given in \ref{sec-ebezcubic}.

	\sssec{Theorem.}
	{\em Let $P := (w_0(1-I)^2,2w_1I(1-I),w_2I^2),$ $w_0w_1w_2\neq 0$ then}
	\enumb
	\item $P$ {\em is the parametric equation of a conic, which passes
	through the points} $P(0) = (1,0,0),$ $P(1) = (0,0,1),$
	$P(\infty) = (w_0,-2w_1,w_2),$
	\item {\em the tangent $t$ at $P$ is}\\
	\hth$[w_1w_2I^2,-w_2w_0I(1-I),w_0w_1(1-I)^2],$\\
	{\em in particular, the tangent at $P(0)$ is $[0,0,1],$ at $P(1)$ is
	$[1,0,0],$ (which meet at $U = (0,1,0)$) and at $P(\infty)$ is}\\
	\hth$[w_1w_2,w_2w_0,w_0w_1],$
	\item $t$ {\em meets $t(0)$ at\\
	\hth$ T = [w_0(1-I),w_1I,0],$\\
	in particular,} $T(\infty) = [-w_0,w_1,0],$
	\item {\em the anharmonic ratio}\\
	\hth anhr$(U,P(0),T(\infty),T) = \frac{w_0-w_1}{w_0} I.$
	\enume

	The proof starts with the observation that the coordinates $P_0,$ $P_1,$
	$P_2$, satisfy the equation\\
	\hth$ 4w_1^2 P_0P_2 = w_2w_0 P_1^2,$\\
	which is indeed the equation of a conic with the prescribed properties. 
	The corresponding polarity matrix is\\
	\hth$\mat{0}{0}{2w_1^2}{0}{-w_2w_0}{0}{2w_1^2}{0}{0}.$

	Notice that $T(0) = (1,0,0)$ and is not undefined.

	The tangent can either be obtained from the polarity or using $P\times
	DP$, where its direction $DP = (-w_0(1-I),w_1(1-2I),w_2I)$.

	The last statement of the Theorem is associated with the four tangents
	Theorem of J. Steiner.  It can be used as a method to draw conics.
	In the excellent language Postcript (see Reference Manual), a general
	method is given to draw curves based on rhe work of de Cateljau and
	B\'{e}zier as well as a method to draw ellipses using the {\em Euclidean
	concepts} of rotation and scaling differently in the direction of its
	axis. This method does not allow to draw hyperbolas or parabolas and
	ignores the fact that a conic is a projective concept.  The following
	gives a method which allows to draw conics using 3 points $A,\:B,\:C,$
	and the tangents $t_A,\:t_B$ at two of the 2 points.\\
	It is then generalized to other curves.

	\sssec{Algorithm.}
	If the barycentric coordinates are chosen in such a way that
	$A = (1,0,0),$ $B = (0,0,1),$ $t_A \times t_B = (0,1,0)$ and
	$C = (1,1,1),$ then
	the points on the conic are given by $P$ of the preceding Theorem,
	with, for instance, $w_0 = 2,$ $w_1 = -1,$ $w_2 = 2.$
	In the case of a finite field, we compute $P$ for each
	element of the field or for an appropriate subset of it. In the case
	of the field of reals, we can compute $P$ for $tan(\pi t),$ $t = 0$
	to 1, avoiding $1/2$,  a section of the conic can be obtained by
	appropriately limiting the set \{$t$\}, joining the successive points
	by segments will automatically give the asymptotes for an hyperbola,
	which is appropriate because their directions are indeed points in the
	Euclidean plane, as we prefer to consider it (in its extended form).
	An other approach is to limit the domain of $P$ to $[0,1]$ to obtain
	one section of the conic and to replace $w_1$ by $-w_1$, which is
	equivalent to compose {\bf P} with $\frac{I}{2I-1},$ to obtain the
	complement, see Farin, p.185.

	For some of the Theorems, see Farin.

	In what follows, the superscript of $B,$ $P$ and {\bf P} are indices and
	not exponents.

	\sssec{Definition.}
	The {\em Bernstein polynomials} are\\
	\hth$B_i^n := \ve{n}{i}\:I^i(1-I)^{n-i},$ $0 \leq i \leq n$.\\
	By convention $B_{-1}^n = B_{n+1}^n := 0.$

	In particular,\\
	\hth$B_0^2 = (1-I)^2,\:B_1^2 = 2I(1-I),\:B_2^2 = I^2.$

	\sssec{Theorem.}
	\enumb
	\item$B_0^0 = 1,$ $B_i^n = (1-I)B_i^{n-1} + I B_{i-1}^{n-1},$
	\item$D B_i^n = n(B_{i-1}^{n-1}-B_i^{n-1}).$
	\item$\sum_{j=0}^n B_j^n = 1.$
	\enume

	\sssec{Definition.}
	A {\em weighted point} {\bf P} is a set of 3 non homogeneous coordinates
	which are not all 0.  We can add weighted points and multiply by
	scalars, but two weighted points which differ by a multiplicative
	constant are not equivalent. I will use the notation $P$ for the
	equivalent point.

	\sssec{Definition. [de Casteljau]}\label{sec-dcastel}
	Given $n+1$ weighted points ${\bf P}_0,$ ${\bf P}_1,$ \ldots,
	${\bf P}_n,$ called the {\em B\'{e}zier polygon}, define\\
	\hth${\bf P}_i^0 := {\bf P}_i,$ $0 \leq i \leq n,$\\
	\hth${\bf P}_i^j := (1-I)\:{\bf P}_i^{j-1} + I\: {\bf P}_{i+1}^{j-1},$
	$1 \leq j \leq n,$ $0 \leq i \leq n-j.$\\
	\hth${\bf P}^n := {\bf P}_0^n,$\\
	The curve $P^n$ is called the {\em de Casteljau curve of order $n$}.

	The same curve is also called the {\em B\'{e}zier curve}.

	\sssec{Theorem.}\label{sec-tbezier}
	\enumb
	\item $P^n(0) = {\bf P}_0,$ $P^n(1) = {\bf P}_n$.
	\item  ${\bf P}_i^j = \sum_{k=0}^j {\bf P}_{i+k} B_k^j,$
	$0 \leq j \leq n,$ $0 \leq i \leq n-j,$

	{\em in particular,}
	\item${\bf P}^n = \sum_{k=0}^n {\bf P}_k B_k^n,$
	\item$D {\bf P}^n = n \sum_{k=0}^{n-1}
		({\bf P}_{k+1}-{\bf P}_k)B_k^{n-1}.$
	\enume

	\sssec{Definition.}
	\hth$\Delta {\bf Q}_k = {\bf Q}_{k+1} - {\bf Q}_k,$\\
	\hth$\Delta^{r+1} {\bf Q}_k = \Delta^r{\bf Q}_{k+1}-\Delta^r{\bf Q}_k,$

	\sssec{Theorem.}
	\hth$\Delta {\bf Q}_0 = \sum_{k=0}^r (-1)^{r-j}\ve{r}{j}{\bf Q}_{i+k}.$

	\sssec{Theorem.}
	\hth$D^r{\bf P}^n = \frac{n!}{(n-r)!}\sum_{k=0}^{n-r}
		\Delta^r{\bf P}_k B_k^{n-r}.$\\
	\hth$D^r{\bf P}^n = \frac{n!}{(n-r)!} \Delta^r{\bf P}_0^{n-r}.$\\
	{\em In particular},\\
	\hth$D{\bf P}_n = n({\bf P}_1^{n-1}-{\bf P}_0^{n-1}).$

	\sssec{Curves with cubic parametrization.}\label{sec-ebezcubic}
	For $n = 3,$\\
	\hth${\bf P}^3 = {\bf P}_0(1-I)^3 + 3{\bf P}_1(1-I)^2I
	+ 3{\bf P}_2(1-I)I^2 + {\bf P}_3I^3.$\\
	\hth$D{\bf P}^3(0) = {\bf P}_1 - {\bf P}_0,$\\
	\hth$D{\bf P}^3(1) = {\bf P}_3-{\bf P}_2.$\\
	In other words the direction of the tangents at the end points
	is that of the line joining the end points to the nearest point.

	If the cubic associated with the $i$-th coordinate of the curve $P^3$
	is\\
	\hth$f = c_0 + c_1I+c_2I^2+c_3I^3$,\\
	then the $i$-th coordinate $a_j$ of the B\'{e}zier polygon ${\bf P}_j$
	is given by\\
	\hth$a_0 = c_0,$ $a_1 = c_0+\frac{1}{3}b_1,$
	$a_2 = a_1+\frac{1}{3}(c_1+c_2),$ $a_3 = c_0+c_1+c_2+c_3.$\\
	Indeed, $a_0(1-I)^3 + 3a_1(1-I)^2I + 3a_2(1-I)I^2+a_3I^3 = f.$

	These last formulas allows for the determination of the weighted points
	${\bf P}_i$ of the cubic (approximation) given the 3 non homogeneous
	coordinates of the parametrized curve.  If the cubic associated with the
	$i$-th coordinate reduces to a linear function then
	${\bf P}_i = {\bf P}^3(i/3),$ $i = 0,1,2,3$.

	It is often convenient to choose $-1$ and 1 for the end points instead
	of 0 and 1 by means of a change of variable.  If $g = d_0 + d_1I + d_2 
	I^2 + d_3 I^3$ is the new polynomial, $f = c_0 + c_1I + c_2 I^2 + c_3I^3
	= g \circ \phi,$ with $\phi = 2I-1$. In this case, we obtain the
	symmetric formulas,\\
	$a_0 = d_0 - d_1 + d_2 - d_3,$
	$a_1 = d_0 -\frac{1}{3}d_1  -\frac{1}{3}d_2 + d_3,$
	$a_2 = d_0 +\frac{1}{3}d_1  -\frac{1}{3}d_2 - d_3,$
	$a_2 = d_0 + d_1 + d_2 + d_3,$

	\sssec{Example.}
	For the curve $(I-3I^3,1-I^2,1)$, for the first coordinate,
	$d_0 = d_2 = 0,$ $d_1 = 1,$, $d_3 = -3$, therefore
	$a_0 = 2,$ $a_1 = -\frac{10}{3},$ $a_2 = \frac{10}{3}$, $a_3 = -2.$
	for the second coordinate,
	$d_1 = d_3 = 0,$ $d_0 = 1,$, $d_2 = -1$, therefore
	$a_0 = 0,$ $a_1 = \frac{4}{3},$ $a_2 = \frac{4}{3},$ $a_3 = 0.$
	Therefore the B\'{e}zier polygon is\\
	${\bf P}_0 = (2,0,1),$ ${\bf P}_1 = (-\frac{10}{3},\frac{4}{3},1),$
	${\bf P}_2 = (\frac{10}{3},\frac{4}{3},1),$ ${\bf P}_3 = (-2,0,1).$
 
	This gives the Cartesian coordinates of the following points on the
	curve associated with $i/20,$ $i$ = 0 to 20:\\
	  2.000,0.00;  1.287,0.19;  0.736,0.36;  0.329,0.51;  0.048,0.64;
	 -0.125,0.75, -0.208,0.84, -0.219,0.91, -0.176,0.96, -0.097,0.99; 
	  0.000,1.00;  0.097,0.99;  0.176,0.96;  0.219,0.91;  0.208,0.84; 
	  0.125,0.75; -0.048,0.64; -0.329,0.51; -0.736,0.36; -1.287,0.19;
	 -2.000,0.00.\\
	The complement of the curve using $\frac{i/20}{2i/20-1}$ is\\
	 2.0000,0.0000; 3.0041,-0.2346; 4.6094,-0.5625; 7.3178,-1.0408; \ldots,
	-7.3178,-1.0408;-4.6094,-0.5625;-3.0041,-0.2346;-2.0000, 0.0000.

	\sssec{Problem.}
	Given a curve in the plane, what are the condition for a
	representation of the 3 non-homogeneous coordinates by polynomials
	of degree $n$.  For conics, we have seen that $n = 2.$
	\ssec{Projectivity determined by a conic.}

	\sssec{ Definition.}
	Joining 2 distinct points of a conic, is to determine the
	line through the 2 points.  Joining a point of a conic to itself is to
	determine the {\em tangent to the conic} at that point.

	\sssec{ Example.}
	For $p = 3,$
	The conic $X^2 + 2 Y Z = 0$ has the points
	(0), (1), (13), (17), (25) and (29).\\
	The tangents at $(X_0,Y_0,Z_0)$ is $[X_0,Z_0,Y_0].$\\
	The tangent at (0) is [1] and the tangent at (1) is [0].\\
	These points joined to (0) give the lines
	[6], [26], [16], [21], [11].  These points joined to (1) give
	[0], [8], [10], [7], [9].\\
	These lines determine on the ideal line [12], the projectivity which
	associates to\\
	(26), (5), (14), (18), (22), (26), the points
	(5), (26), (18), (22), (10), (14).\\
	This is precisely the projectivity $\phi$  of \ref{sec-eprojp5}.

	\sssec{Theorem.}
	{\em Let $N$ be a symmetric matrix associated to a conic.}
	\enumb
	\item {\em $P$ is on the conic if $P \cdot NP = 0.$}
	\item {\em If $P$ is on the conic and $C$ is not, the other point on the
		conic, if any, is\\
	\hth $P + y C$, with $y = - 2 \frac{C \cdot NP}{C \cdot NC}.$}
	\enume

	Proof:
	\hth$	(P + yC) \cdot N(P + yC) = 0,$\\
	or\\
	\hth$P \cdot NP + yC \cdot NP + y P \cdot NC + y^2 C \cdot NC = 0,$\\
	but $P \cdot NP = 0$ and $C \cdot NC\neq 0$ and $N$ is symmetric,
	therefore $C \cdot NP = P \cdot NC,$
	hence $y = - 2 \frac{C \cdot NP}{C \cdot NC}.$

	\sssec{ Theorem.}
	{\em Let $l$ be a line and $A,$ $B$ be 2 points on the line but not on
	the conic associated with the symmetric matrix $N$,\\
	let $a := A \cdot NA,$ $b := B \cdot NB,$ $c := A \cdot NB 
	= B \cdot NA,$\\
	let $C$ be an arbitrary point on the line, $C = A + k B.$\\
	Let $P_1$ and $P_2$ be 2 distinct points on the conic,
	let $a_1 = (P_1 \cdot A*P_2),$ $b_1 := (P_1 \cdot B*P_2),$
	$c_1 := (A \cdot B*P_1),$
	$c_2 := (A \cdot B*P_2),$
	$d_1 := A \cdot NP_1,$ $d_2 := B \cdot NP_2.$
	If $C \times P_1$ meets the conic at $Q$ and $P_2 \times Q$ meets $l$
	at $D,$ then}
	\enumb
	\item $D = ( (a a_1) + 2(c a_1 - d_1 c_1) k + (b a_1 +2 d_2 c_1) k^2 )
		B$\\
	\hth$	- ( (a b_1 -2 d_1c_2) + 2 (c b_1 - d_2 c_2) k + b b_1 k^2 A$.
	\item {\em The correspondance berween $C$ and $D$ is a projectivity.}
	\enume

	Proof: $Q = (C \cdot NC) P_1 - 2(C \cdot NP_1) C,$\\
	$D = (A*B) * (P_2*Q) = (A \cdot P_2*Q) B - (B \cdot P_2*Q) A\\
	\hti{4}= (Q \cdot A*P_2) B - (Q \cdot B*P_2) A\\
	\hti{4}= ( (C \cdot NC) (P_1 \cdot A*P_2)
		- 2(C \cdot NP_1) (C \cdot A*P_2) ) B\\
	\hth- ( (C \cdot NC) (P_1 \cdot B*P_2) - 2(C \cdot NP_1)
	(C \cdot B*P_2) ) A,$\\
	but $C \cdot NC = a + 2kc + k^2b$\\
	therefore\\
	$D = ( (a+2ck+bk^{2)} a_1 - 2( d_1 + kd_2)(-c_1k) ) B$\\
	\hti{4}$- ( (a + 2kc + bk^{2)} b_1 -2(d_1+d_2k)(c_2) ) A$\\

%
	\sssec{Theorem.}
	{\em If the line $l$ is $[1,1,1]$ then the conic $X^2 + Y^2 + kZ^2 = 0$
	determines on $l$ the involution $\eta$\\
	\hth$	\eta(1,Y,-1-Y) = (1,f(Y),-1-f(Y)),$ with\\
	\hth$	f(Y) = -(\frac{1+k)+kY}{k+(1+k)Y}.$}

	Proof.  The point $(1,Y,-1-Y)$ on $l$ has the polar $[1,Y,-k(1+Y)],$
	which meets $l$ at $(Y+k(1+Y),-1-k(1+Y),1-Y) = (1,f(Y),-1-f(Y)).$
	\ssec{Cubics.}

	\sssec{Notation.}
	In this section, the cubic is denoted by $\gamma$ , $(I,i)$ will
	denote an inflection point and
	the corresponding tangent, $(A,a)$ a point on the cubic and its tangent,

	\sssec{Theorem.}
	{\em Given $I,$ there exists 3 $(A_i,a_i)$ such that
	$a_i \cdot I = 0$ and $A_i$ are collinear.}

	\sssec{Theorem.}
	{\em Let $B_j,$ $j = 0$ to 5 be on $\gamma$  and a conic $\theta$ , if
	$C_j$ is the third point on $B_j \times B_{j+3},$ then $C_j$ are
	collinear.}

	\sssec{Corollary.}
	{\em Given $(A_i,a_i),$ $i = 0$ to 2, let $B_i$ be the other point on
	$a_i$ then $B_i$ are collinear.}

	\sssec{Corollary.}
	{\em If $B_k,$ $k = 0$ to 3 are on $\gamma$.  Let a conic
	$\theta_l$ meet $\gamma$  also
	at $C_{l,0}$ and $C_{l,1},$ $C_{l,0} \times C_{l,1}$ passes through a
	fixed point $D$ of $\gamma$.}

	\sssec{Theorem.}
	{\em The third point on $I_1 \times I_2$ is an inflection point.}

	\sssec{Theorem.}
	{\em Given $(A_l,a_l),$ $l = 0,$ 1, and $B_l$ is the other point on
	$a_l,$ $(A_0 \times  A_1) \times (B_0 \times B_1)$ is on the cubic.}

	\sssec{Theorem.}
	{\em The anharmonic ratio of the 4 tangents through $A$ distinct from
	$a$ is constant.}

{\tiny\footnotetext[1]{G24.TEX [MPAP], \today}}
	\ssec{Other models for projective geometry.}

	\setcounter{subsubsection}{-1}
	\sssec{Introduction.}
	Many models can be derived from the model given in
	section 1.  This is most easily accomplished by starting with a
	correspondence between points in the plane and adjusting for special
	cases.  One such correspondence is $(x_0,x_1,x_2)$ to 
	$(\frac{1}{x_0},\frac{1}{x_1},\frac{1}{x_2}),$ and
	will be studied in some detail.  It assumes some given triangle
	$\{A_0,A_1,A_2\},$ whose vertices have coordinates $(1,0,0), (0,1,0),
	(0,0,1)$.

	\sssec{Definition.}\label{sec-dmodel}
	In inversive geometry, the {\em "points"} are the points $(x_0,x_1,x_2),$ with
	$x_0 x_1 x_2\neq 0$ together with the lines
	$[0,x_1,x_2],$ $[x_0,0,x_2],$ $[x_0,x_1,0],$\\
	the {\em "lines"} are the point conics\\
	\hth$	a_0 x_1 x_2 + a_1 x_2 x_0 + a_2 x_0 x_1 = 0.$\\
	degenerate or not.
	A {\em "point" is on a "line"}, which is a non degenerate point conic,
	iff it belongs to it or is tangent to it.
	If the point conic degenerates in 2 lines, one which is a side of the
	triangle and the other passes through the opposite vertex, then the
	"points" who belong to it are the two lines and the points on the line
	through the opposite vertex but not on the sides of triangle.
	If the point conic degenerates in 2 lines, which are 2 sides of the
	triangle, the "points" which belong to it are the lines through the
	common vertex.

	\sssec{Example.}
	The "line" $x_1 x_2 + 2 x_2 x_0 + 3 x_0 x_1 = 0$ belongs to the
	"points"\\
	$(-x_1x_2, (2 x_1 + 3 x_2)x_2, (2 x_1 + 3 x_2)x_1),$
	$(2 x_1 + 3 x_2) x_1 x_2\neq 0$
	and to the "points"
	$[0,3,2],$ $[3,0,1],$ $[2,1,0]$ tangent respectively at $A_0,$ $A_1$
	and $A_2.$\\
	The "line" $2 x_2 x_0 + 3 x_0 x_1 = 0$ belongs to the "points"
	$(x_0,2,-3),$ $x_0\neq 0$ and to the "points" [1,0,0] and [0,3,2].\\
	The "line" $x_1 x_2 = 0$ belongs to the "points" $[0,x_1,x_2].$

	\sssec{Theorem.}
	{\em The model} \ref{sec-dmodel} {\em satisfies the axioms}
	\ref{sec-ax1} {\em of projective geometry.}

	This is most easily seen if we associate
	to the point $P = (x_0,x_1,x_2),$ $x_0 x_1 x_2\neq 0$ the "point"
		$P' = (\frac{1}{x_0},\frac{1}{x_1},\frac{1}{x_2})$
	or $(x_1 x_2, x_2 x_0, x_0 x_1)$,
	to the point $Q_0 = (0,x_1,x_2),$ the "point" $Q_0' = [0,x_1,-x_2],$
	to the point $Q_1 = (x_0,0,x_2),$ the "point" $Q_1' = [x_0,0,-x_2],$
	to the point $Q_2 = (x_0,x_1,0),$ the "point" $Q_2' = [x_0,x_1,0],$ and
	to the line $l = [a_0,a_1,a_2],$
		the line $l',$ $a_0 x_1 x_2 + a_1 x_2 x_0 + a_2 x_0 x_1 = 0.$

	Indeed if $P \cdot l = 0,$ $P'$ is on $l'$ and if $Q_0 \cdot l = 0,$
	$a_1 x_1 + a_2 x_2 = 0,$
	while the tangent to $l'$ at $A_0$ is $[0,a_2,a_1] = [0,x_1,-x_2].$

	\sssec{Theorem.}
	{\em The "lines" are the conics through the vertices $A_0,$ $A_1,$
	$A_2.$}

	\sssec{Theorem.}
	{\em The "conics" are the quartics\\
	0.\hti{6}$b_0 x_1^2 x_2^2 + b_1 x_2^2 x_0^2 + b_2 x_0^2 x_1^2 +$\\
	\hth$(c_0 x_0 + c_1 x_1 + c_2 x_2) x_0 x_1 x_2 = 0.$\\
	The quartic has double points (or nodes) at the vertices $A_0,$ $A_1,$
	$A_2.$\\
	The branches through $A_0$ are real if and only if $c_0^2 > 4 b_1 b_2,
	$\\
	the branches through $A_1$ are real if and only if $c_1^2 > 4 b_2 b_0,
	$\\
	the branches through $A_2$ are real if and only if $c_2^2 > 4 b_0 b_1.
	$\\
	Vice versa if a quartic as double points at $A_0,$ $A_1$ and $A_2$ it
	is of the form 0.}

	\sssec{Theorem.}
	{\em If the quartic has double points with real branches at $A_0,$ $A_1$
	and $A_2,$ the tangents $P_0' P_1'$ at $A_0,$ $P_2' P_3'$ at $A_1$ and
	$P_4' P_5'$ at $A_2$ are
	such that if $K_0'$ is the tangent to the conic $(A_0,(A_1,P_3'),
	(A_2,P_4')),$
	if $K_1'$ is the tangent to the conic $(A_1,(A_2,P_5'),(A_0,P_0')),$
	and if $K_2'$ is the tangent to the conic $(A_2,(A_0,P_1'),(A_1,P_2')),$
	then there is a conic through $A_0,$ $A_1,$ $A_2$ with tangents $K_0',$
	$K_1',$ $K_2'.$}

	This is a direct consequence of the Theorem of Pascal associated to
	the model.

	\sssec{Theorem.}
	{\em If a quartic has double points with real branches at
	$A_0$, $A_1$ and $A_2,$ then the 6 tangents at these points belong to
	the same line conic.}

	Proof:  Let the tangents at $A_0,$ $A_1,$ $A_2$ be
		$[0,1,z],$ $[0,1,z'],$ $[x,0,1],$ $[x',0,1],$ $[1,y,0],$
	$[1,y',0],$\\
	the tangents $[0,1,z]$ at $A0$ satisfy $b1 z^2 + c_0 z + b_2,$ therefore
	$z z' = \frac{b_2}{b_1},$ similarly, $y y' = \frac{b_0}{b_2}$ and
	$x x' = \frac{b_1}{b_0}.$\\
	On the other hand, applying Brianchon's theorem to these tangents
	gives the Brianchon lines $[0,x'y,-1],$ $[-1,0,y'z],$ $[z'x,-1,0]$
	and these belong to the same point if $x'y y'z z'x = 1.$\\
	This conjecture was most strongly confirmed by a computer program
	and proven within an hour.

	\sssec{Theorem.}
	{\em If $b_0 = 0,$ then the quartic degenerates in the side
	$A_1 \times A_2$
	and a cubic with double point at $A_0$ passing through $A_1$ and $A_2.$}

	The conic $c_1 c_2 y z + b_1 c_1 z x + b_2 c_2 x y = 0$ plays, in the
	invertible
	geometry, the role of the "line" tangent at the "point" [1,0,0].

	\ssec{Notes.}
	\sssec{Theorem. [Jones]}
	{\em Let $n$ be even. If an $n$-gon is inscribed in a conic and $n$-1 
	sides meet a line at fixed points, then the $n$-th side also meets the 
	line at a fixed point and dually.}

	\sssec{Theorem. [Jones]}
	{\em The preceding Theorem, when $n = 4$ is equivalent to Pascal's
	Theorem.}

{\tiny\footnotetext[1]{G25.TEX [MPAP], \today}}
	\section{Geometric Models on Regular Pythagorean Polyhedra.}

	\setcounter{subsection}{-1}
	\ssec{Introduction.}
	Completely independently, one of my first student at the
	"Universit\'{e} Laval", Quebec City, made the important discovery
	that the regular
	polyhedra can be used as models for finite geometries associated with
	2, 3 and 5.  Then, he introduced the nomenclature of selector
	(s\'{e}lecteur) for the notion of cyclic difference sets, introduced by
	J. Singer, in 1938, to label points and hyperplanes in $N$ dimensional
	projective geometry of order $p^k$ (See Baumert, 1971) and to construct
	an appropriate numbering of the points and lines on the polyhedra.\\
	Except for the fundamental contribution of Singer, the introduction of
	selector polarity (prepared by the use of $f(a+b)$ instead of  $f(a-b)$
	in the definition of incidence), the introduction of auto-polars and
	those on the conics for the dodecahedron, all the results in this
	section are due to Fernand Lemay.

	Clearly we have only to study the tetrahedron, the cube and the
	dodecahedron, because the octahedron is dual to the cube and the
	icosahedron is dual to the dodecahedron.

	\ssec{The selector.}\label{sec-Sselector}

	\setcounter{subsubsection}{-1}
	\sssec{Introduction.}
	The important concept of the cyclic difference sets
	allows for an arithmetization of projective geometry which is as close
	to the synthetic point of view as is possible.  With it, it is not only
	trivial to determine all the points on a line, and lines incident to a
	point, but also the lines through 2 points and points on 2 lines.  This
	concept makes duality explicit through the correlation, which is a
	polarity when $p \geq  5.$  The definitions of selector function and
	selector correlation is implicit in Lemay's work.

	\sssec{Definition.}
	A {\em difference set} associated to $q = p^k$ is a set of $q+1$
	integers $\{s_0,$ $s_1,$  $\ldots$  , $s_q\}$ such that the
	$q^2+q$ diferences $s_i - s_j,$ 
	$i\neq j$ modulo $n :=  q^2+q+1$ are distinct and different from 0.
	When applied to Geometry, I will prefer the terminology of Lemay and
	use the synonym {\em selector}.  The elements of the selector are called
	{\em selector numbers}.

	\sssec{Theorem. [Singer]}
	{\em For any $q = p^k$ there exists difference sets.}

	\sssec{Theorem.}\label{sec-tseltr}
	{\em If $\{s_i\},$ $i = 0$ to $q,$ is a difference set and $k$ is
	relatively prime to $n,$ then}
	\enumb
	\item$\{s'_i = a + k s_{i+1}\},$ {\em is also a difference set.\\
	The indices are computed modulo $q+1$ and the selector numbers,
	modulo $n.$}
	\enume

	Using 0, we can always find a selector for which 0 and 1 are
	selector numbers.

	\sssec{Example. [Singer]}\label{sec-esel}
	The following are difference sets associated with $q = p^k:$\\
	For $p = 2:$  \{0, 1, 3\} modulo 7.\\
	For $p = 3:$  \{0, 1, 3, 9\} modulo 13.\\
	For $q = 2^2:$ \{0, 1, 4, 14, 16\} modulo 21.\\
	For $p = 5:$  \{0, 1, 3, 8, 12, 18\} modulo 31.\\
	For $p = 7:$  \{0, 1, 3, 13, 32, 36, 43, 52\} modulo 57.\\
	For $q = 2^3:$ \{0, 1, 3, 7, 15, 31, 36, 54, 63\} modulo 73.\\
	For $q = 3^2:$ \{0, 1, 3, 9, 27, 49, 56, 61, 77, 81\} modulo 91.\\
	For $q = 11:$ \{0, 1, 3, 12, 20, 34, 38, 81, 88, 94, 104, 109\}
	modulo 133.

	\sssec{Definition.}
	If $a = 1$ and $k = -1,$ the selector $s'_i := 1 - s_i$ is called
	the {\em complementary selector} or {\em co-selector} of $s_i.$ \\
	The selectors obtained using $k = 2,$ $\frac{1}{2},$ $-2,$
	$-\frac{1}{2}$ are called respectively {\em bi-selector},
	{\em semi-selector}, {\em co-bi-selector}, {\em co-semi-selector}.

	\sssec{Example.}
	For $q = 4,$ other selectors are\\
	\hth$	\{10, 12, 17, 18, 21\}, \{0, 5, 20, 7, 17\}$ and
	$\{0, 1, 6, 8, 18\}.$

	\noindent For $p = 7$, if\\
	the selector is		 $\{0, 1, 7,24,36,38,49,54\}$, then\\
	the co-selector is	 $\{0, 1, 4, 9,20,22,34,51\}$,\\
	the bi-selector is	 $\{0, 1, 5,27,34,37,43,45\}$,\\
	the co-bi-selector is	 $\{0, 1,13,15,21,24,31,53\}$,\\
	the semi-selector is	 $\{0, 1, 9,11,14,35,39,51\}$,\\
	the co-semi-selector is  $\{0, 1, 7,19,23,44,47,49\}$.

	\sssec{Program.}
	All selectors derived by multiplication from one of them are
	given in [113]MODP30.

	\sssec{Definition.}
	The {\em selector function $f$ associated to the selector} $\{s_i\}$ 
	is the function from ${\Bbb Z}_n$ to ${\Bbb Z}_n$ \\
	\hth$	f(0) = 0,$ $f(s_j-s_i) = s_i,$ $i\neq  j.$

	\sssec{Example.}
	For $p = 2,$ the selector function associated with
	$\{0, 1, 3\} \umod {7}$ is\\
	$\begin{array}{lrrrrrrr}
	i& 0& 1& 2& 3& 4& 5& 6\\
	f(i)&0& 0& 1& 0&3&3&1
	\end{array}$\\
	For $p = 3,$ the selector function associated with
	$\{0, 1, 3, 9\} \umod {13}$ is\\
	$\begin{array}{lrrrrrrrrrrrrr}
	i& 0& 1& 2& 3& 4& 5& 6& 7& 8& 9&10&11&12\\
	f(i)&0& 0& 1& 0&-4&-4& 3&-4& 1& 0& 3& 3& 1
	\end{array}$\\
	For $p = 5,$ the selector function associated with
	$\{0, 1, 3, 8, 12, 18\} \umod {31}$ is\\
	$\begin{array}{lrrrrrrrrrrrrrrrr}
	i&   0& 1& 2& 3& 4& 5& 6& 7& 8& 9&10&11&12&13&14&15\\
	f(i)&0& 0& 1& 0& 8& 3&12& 1& 0& 3& 8& 1& 0&18&18& 3\\
	\cline{1-17}
	i&   &16&17&18&19&20&21&22&23&24&25&26&27&28&29&30\\
	f(i)&&18& 1& 0&12&12&18&12& 8& 8&18& 8&12& 3& 3& 1
	\end{array}$

	\sssec{Theorem.}
	\enumb
	\item$f(j-i) -i = f(i-j) - j \pmod{n}.$
	\enume

	Points, lines and incidence in the 2 dimensional geometry associated
	with $q = p^k$ and $n := q^2 + q + 1$ are defined as follows.

	\sssec{Definition.}
	The {\em points} are elements of the set \{0, 1,  $\ldots$  , $n-1$\},\\
	The {\em lines} are elements of the set \{0, 1,  $\ldots$  , $n-1$\}.\\
	A point  $a$  is {\em incident} to a line  $b$  iff  $f(a+b) = 0.$

	\sssec{Notation.}
	The points are denoted by a lower case letter or by an
	integer in ${\Bbb Z}_n.$  The lines are denoted by a lower case letter
	or by an integer in ${\Bbb Z}_n$ followed by an asterix.  The line
	incident to the points $a$ and $b$ is denoted $a \times b,$ the point
	incident to the lines $a^*$ and $b^*$ is denoted $a^* \times b^*.$

	\sssec{Theorem.}\label{sec-tselinc}
	{\em Given a selector $\{s_j\}$ associated with $q = p^k$ and the
	corresponding selector function $f:$}
	\enumb
	\item {\em The $q+1$ points incident to or on the line $i^*$ are}
	$s_j-i \umod{n}.$
	\item {\em The $q+1$ lines incident to or on the point $i$ are}
	$(s_j-i \umod{n})^*.$
	\item $a\neq b \implies a \times b = (f(b-a)-a)^*.$
	\item $a\neq b \implies a^* \times b^* = f(b-a) - a.$
	\item {\em $a$ on $b^*$  iff  $b$ on $a^*.$}
	\enume

	The statements in the preceding Theorem reflect the duality in
	projective geometry.

	\sssec{Definition.}
	The {\em selector polarity} is the correlation which
	associates to the point  $i$  the line  $i^*.$  The points  $x$  which
	are on $x^*$ are called {\em auto-polars}.

	The name "polarity" is appropriate because of \ref{sec-tselinc}.4.\\
	The selector polarity and the auto-polars play an important role in a
	natural way of labeling the elements of the Pythagorean solids.

	\sssec{Theorem.}
	{\em The auto-polars are given by\\
	\hth$	a_i = \frac{s_i}{2},$ modulo $n.$}

	Indeed we should have for an auto-polar $x,$ $x = s_i - x.$

	\sssec{Definition.}
	A {\em primitive polynomial of degree} 3 {\em over} $GF(q),$ is an
	irreducible polynomial $P$ of degree 3 such that\\
	\hth$	I^k\neq  \underline{1}$ for $k = 1$ to $q-2,$\\
	where $I$ is the identity function and $\underline{1}$ the constant
	polynomial 1.\\
	The multiplication is done modulo $P$ and polynomials which differ by a
	multiplicative constant$\neq$  0 modulo $q$ are equivalent.

	\sssec{Theorem. [Singer]}
	{\em For each value of $q = p^k$ a selector can be obtained by choosing
	a primitive polynomial of degree 3 over $GF(q).$ The selector is the set
	of exponents of $I$ between 0 and $q-2$ which are of degree less than
	2.}

	\sssec{Example.}
	For $p = 3,$ $P = I^3 - I + 1,$\\
	$I^0 = 1,$ $I^1 = I,$ $I^2 = I^2,$ $I^3 = I - 1,$ $I^4 = I^2 - I,$
	$I^5 = I^2 - I + 1,$
	$I^6 = I^2 + I + 1,$ $I^7 = I^2 - I - 1,$ $I^8 = I^2 + 1,$
	$I^9 = I + 1,$
	$I^{10} = I^2 + I,$ $I^{11} = I^2 + I - 1,$ $I^{12} = I^2 - 1$ and we
	have $I^{13} = 1.$\\
	Therefore the selector is $\{0, 1, 3, 9\}$.

	\ssec{The tetrahedron.}\label{sec-Stetra}

	\setcounter{subsubsection}{-1}
	\sssec{Introduction.}
	I have found useful to introduce the adjectives vertex, edge and in
	later sections, face, to distinguish points and lines which have
	different representation in the Pythagorean solids.

	\sssec{Definition.}
	The {\em points} in the tetrahedron model consist of
	\enumb
	\item The 4 {\em vertex-points}, which are the 4 vertices (or the
		opposite planes or the line through the center $C$ of the
		tetrahedron perpendicular to one of the 4 planes).
	\item The 3 {\em edge-points}, which are the pairs of orthogonal
		edges, (or the mid-points of 3 non orthogonal edges or the
		line through these points and the center $C$).
	\enume
	The {\em lines} in the tetrahedron model consist of
	\begin{enumerate}
	\setcounter{enumi}{1}
	\item The 6 {\em edge-lines}, which are incident to the 2
		vertex-points and to the edge-point on them.
	\item The {\em tetrahedron-line}, which is incident to the 3
		edge-points.
	\enume

	\sssec{Theorem.}
	{\em The model satisfies the axioms of projective geometry for
	$p = 2.$}

	\sssec{Theorem.}
	{\em With the selector \{0,1,3\} $\umod{7},$ the 3 points, 0, 4 and 5
	are auto-polars.  It is therefore natural to associate them to the 3
	edge-points.  These points are on the line $3^*,$ it is natural to
	associate it to the tetrahedral line.  Any of the vertex-points can be
	chosen as the polar 3 of $3^*.$  We will choose the 3 adjacent
	edge-lines as $0^*,$ $4^*$ and $5^*$ such that $0 \cdot 0^*
	= 4 \cdot 4^* = 5 \cdot 5^* = 0.$\\
	The other vertex-points are the third point on $0^*,$ $4^*$ and $5^*,$
	therefore
	2 is on the line and $2^*$ is the line orthogonal to the line associated
	to 5, similarly for 1 and $1^*,$ to 0 and 6 and $6^*,$ to 4.}

	\sssec{Figure.}
	\begin{picture}(200,100)(-140,-10)
	\put(0,0){\line(80,0){80}}
	\put(0,0){\line(3,5){40}}
	\put(80,0){\line(-3,5){40}}
	\put(0,0){\circle*{2}}		\put(-10,4){$A_1$}
	\put(40,0){\circle{2}}		\put(36,4){$P_0$}
	\put(80,0){\circle*{2}}		\put(78,4){$A_2$}
	\put(20,33.3){\circle*{2}}	\put(8,35.2){$P_2$}
	\put(60,33.3){\circle*{2}}	\put(62,35.2){$P_1$}
	\put(40,66.7){\circle{2}}	\put(36,70.4){$A_0$}
	\put(40,22.2){\circle{2}}	\put(36,24.3){$P$}
	\end{picture}

	\sssec{Theorem.}
	{\em A complete quadrangle configuration consists of the 4
	vertex-points $A_0,$ $A_1,$ $A_2,$ P and the 6 edge-lines
	$a_0 = A_1 \times A_2,$ 
	$a_1 = A_2 \times A_0,$ $a_2 = A_0 \times A_1,$ $p_0 = P \times A_0,$
	$p_1 = P \times A_1,$ $p_2 = P \times A_2.$  It
	has the 3 edge-points $P_i = p_i \times a_i$ as its diagonal points, and
	these are on the tetrahedron-line $p.$}

	\sssec{Exercise.}\label{sec-eseltetra}
	\enumb
	\item For $q = 2,$ determine the primitive polynomial giving the
		selector
	$\{0, 1, 3\}$.
	\item Determine the correspondence between the selector notation
		and the homogeneous coordinates for points and lines.  Note
		that these are not the same.
	\item The correspondence $i$ to $i^*$ is a polarity whose fixed
		points are on a line.  Determine the matrix representation and
		the equation satisfied by the fixed points.
	\item Determine the degenerate conic through 0, 1, 2 and 5 with
		tangent $5^*$ at $5$, its matrix representation and its
		equation in homogeneous coordinates.
	\item Determine all the non degenerate conics.
	\enume

	\ssec{The cube.}\label{sec-Scube}

	\sssec{Convention.}
	In what follows we identify elements of the cube,
	which are symmetric with respect to its center $C,$ for instance,
	the parallel faces.  There are therefore 3 independent faces,
	4 independent vertices and 6 independent edges.

	\sssec{Definition.}
	The {\em points} in the cube model consist of
	\enumb
	\item The 3 {\em face-points}, which are the square faces or their
		centers or the lines joining $C$ to these points.
	\item The 4 {\em vertex-points}, which are the vertices or the
		lines joining $C$ to these vertices.
	\item The 6 {\em edge-points}, which are the edges, or the mid-points
		of the edges or the lines joining $C$ to these points.
	\enume
	The {\em lines} in the cube model consist of
	\begin{enumerate}
	\setcounter{enumi}{2}
	\item The 3 {\em face-lines}, corresponding to a face $f,$ which are
		incident to the 2 face-points and to the 2 edge-points in the
		plane through $C$ parallel to $f.$
	\item The 4 {\em vertex-lines}, corresponding to a vertex $V,$ which are
		ncident to the vertex-points $V$ and to the 3 edge-points not
		adjacent to $V.$
	\item The 6 {\em edge-lines}, corresponding to an edge $e,$ which are
		incident to the face-point perpendicular to $e,$ to the 2
		vertex-points and the edge-point on $e.$
	\enume

	\sssec{Theorem.}
	{\em The cube model satisfies the axioms of projective geometry
	for $p = 3.$}

	\sssec{Theorem.}\label{sec-tcube}
	{\em With the selector \ref{sec-esel} for p = 3, the auto-polars are
	0, 7, 8 and 11.  If we examine the quadrangle-quadrilateral
	configuration, we observe that $p$ and $q_i$ are the lines which
	require a 4-th point, it is easy to verify that, with $p = 3,$ $P$ is
	on $p$ and $Q_i$ is on $q_i.$  Moreover\\
	$0^* \cdot 0$ = 0, this suggest to take $P = 0,$ $Q_i = 7, 8, 11.$
	Hence\\
	$r_i := P \times Q_i = 9, 1, 3;$ $p_i := Q_{i+1} \times Q_{i-1} =
	5, 2, 6;$\\
	$A_i := r_i \times p_i =
	4, 12, 10;$ $a_i := A_{i+1} \times A_{i-1} = A_i;$ $P_i
	:= a_i \times r_i = p_i;$ \\
	$q_i := P_{i+1} \times P_{i-1} = Q_i;$ $R_i := a_i \times q_i = r_i;$
	$p := R_1 \times R_2 = p.$}

	\sssec{Theorem.}
	{\em Because of \ref{sec-tcube}, the vertex-points are auto-polars,
	we can choose them as 0, 7, 8 and 11, the other elements of the cube
	follow from \ref{sec-tcube}.  The edge-points are $R_i$ and $P_i,$ the
	face-points are $A_i.$	}

	\sssec{Figure.}
	\begin{picture}(200,160)(-130,-10)
	\put(0,0){\line(1,0){100}}
	\put(0,100){\line(1,0){100}}
	\put(30,130){\line(1,0){100}}
	\put(0,0){\line(0,1){100}}
	\put(100,0){\line(0,1){100}}
	\put(130,30){\line(0,1){100}}
	\put(0,100){\line(1,1){30}}
	\put(100,100){\line(1,1){30}}
	\put(100,0){\line(1,1){30}}
	\put(0,0){\circle*{2}}		\put(0,-8){8}\put(4,4){$Q_1$}
	\put(50,0){\circle{2}}		\put(50,-8){9}\put(46,4){$R_0$}
	\put(100,0){\circle*{2}}	\put(96,-8){11}\put(88,4){$Q_2$}
	\put(0,100){\circle*{2}}	\put(4,92){0}\put(-4,104){$P$}
	\put(50,100){\circle{2}}	\put(54,92){5}\put(46,104){$P_0$}
	\put(100,100){\circle*{2}}	\put(104,92){7}\put(90,104){$Q_0$}
	\put(0,50){\circle{2}}		\put(4,42){2}\put(4,54){$P_1$}
	\put(50,50){\circle{3}}		\put(54,42){10}\put(46,54){$A_2$}
	\put(100,50){\circle{2}}	\put(104,42){1}\put(88,54){$R_1$}
	\put(30,130){\circle*{2}}	\put(34,134){11}  
	\put(80,130){\circle{2}}	\put(84,134){9}  
	\put(130,130){\circle*{2}}	\put(134,134){8}  
	\put(15,115){\circle{2}}	\put(19,111){6}  
	\put(65,115){\circle{3}}	\put(69,111){12}\put(58,119){$A_1$}
	\put(115,115){\circle{2}}	\put(119,111){3}\put(111,119){$R_2$}
	\put(115,15){\circle{2}}	\put(119,7){6}\put(107,19){$P_2$}
	\put(115,65){\circle{3}}	\put(119,57){4}\put(111,69){$A_0$}
	\put(130,30){\circle*{2}}	\put(134,22){0}  
	\put(130,80){\circle{2}}	\put(134,72){2}  
	\end{picture}

	\sssec{Exercise.}\label{sec-eselcube}
	\enumb
	\item For $p = 3,$ determine the primitive polynomial giving the
		selector $\{0, 1, 3, 9\}$.
	\item Determine the correspondence between the selector notation and the
		homogeneous coordinates for points and lines.  Note that these
		are not the same.
	\item The correspondence $i$ to $i^*$ is a polarity whose fixed points
		are on a line.  Determine the matrix representation and the
		equation satisfied by the fixed points.
	\item Determine the degenerate conic through 0, 1, 2 and 5 with tangent
		$4^*$ at $5$, its matrix representation and its the equation in
		homogeneous coordinates.  Hint: use \ref{sec-tqqc3}.
	\item Determine all the conics.
	\enume

	\ssec{The dodecahedron.}\label{sec-Sdodeca}

	\sssec{Convention.}
	In what follows we identify elements of the dodecahedron
	which are symmetric with respect to its center $C,$ for instance,
	the parallel faces.  There are therefore 6 independent faces,
	$\frac{5}{3}6 = 10$ independent vertices and $\frac{5}{2}6 = 15$
	independent edges.

	\sssec{Definition.}\label{sec-ddodeca}
	The {\em points} in the dodecahedron model consist of
	\enumb
	\item The 6 {\em face-points}, which are the pentagonal faces or their
		 center or the lines joining $C$ to these points.
	\item The 10 {\em vertex-points}, which are the vertices or the lines
		joining $C$ to these vertices.
	\item The 15 {\em edge-points}, which are the edges, or the mid-points
		of the edges or the lines joining $C$ to these points.
	\enume
	The {\em lines}  in the dodecahedron model consist of
	\begin{enumerate}
	\setcounter{enumi}{2}
	\item The 6 {\em face-lines}, which are incident to the corresponding
		face-point $F$ and to the 5 edge-points in the plane through $C$
		perpendicular to $CF.$
	\item The 10 {\em vertex-lines}, corresponding to a vertex $V$, which
		are incident to the 3 edge-points in the plane through $C$
		perpendicular to $CV$ and to the 3 vertex-points which joined
		to $V$ form an edge.
	\item The 15 {\em edge-lines}, corresponding to an edge $E,$ which are
		incident to the 2 face-points, the 2 vertex-points and the 2
		edge-points in the plane through $C$ and $E.$
	\enume

	\sssec{Theorem.}
	{\em The dodecahedron model satisfies the axioms of projective
	geometry for $p = 5.$}

	\sssec{Example.}
	For $p = 5,$ the selector function associated with the selector
	$\{0, 1, 3, 8, 12, 18\}$ is\\
	$\begin{array}{lrrrrrrrrrrrrrrrr}
	i&  0& 1& 2& 3& 4& 5& 6& 7& 8& 9&10&11&12&13&14&15\\
	f(i)&0& 0& 1& 0& 8& 3&12& 1& 0& 3& 8& 1& 0&\!-\!13&\!-\!13& 3\\
	type&f& e& e& e& f& v& f& v& e& f& v& v& e& v& e& s\\[10pt]

	i& &16&17&18&19&20&21&22&23&24&25&26&27&28&29&30\\
	f(i)&& \!-\!13& 1& 0&12&12&\!-\!13&12& 8& 8&\!-\!13& 8&12& 3& 3& 1\\
	type&& f& f& e& v& v& v& e& e& v& e& e& e& e& v& e
	\end{array}$

	The auto-polars are $\{0, 4, 6, 9, 16, 17\}$.

	The "type" is explained in the following Theorem.

	\sssec{Theorem.}\label{sec-tdodeca}
	{\em A natural labeling of the points of the dodecahedron and of
	the dodecahedral configuration, associated with the selector
	\ref{sec-esel}, for p = 5, can be obtained as follows.  If we examine
	the dodecahedron configuration,\\
	\hth$	FA_i \cdot fa_i = AF_i \cdot af_i = 0,$\\
	it is therefore natural to choose $FA_i$ and $AF_i$ as the
	auto-polars, but this cannot be done arbitrarily.  Let us choose any
	3 of them as $FA_i,$ 0, 16 and 17.  To obtain $P$ and $A_i,$ we can
	proceed as follows.\\
	$pq_i = FA_i \times FA_{i+1} = 18, 15, 1;$ $PQ_i = pq_i;$\\
	$PR_i = pq_i \times FA_{i-1} = 14, 3, 2;$\\
	$qp_i = PR_i \times PQ_i = 25, 28, 30;$\\
	$AF_i = QP_i \times QP_{i+1} = 6, 4, 9;$\\
	$a_i = FA_{i+1} \times AF_{i-1} = 23, 26, 8;$ $A_i = a_i;$\\
	$p = PR_1 \times PR_2 = 29.$\\
	We therefore choose $P = 29$ and $A_i = 23,26,8.$  We obtain,
	according to \ref{sec-dsd}, \ref{sec-desd} and \ref{sec-dconical}:\\
	$a_i = A_i,$ $r_i = 20, 5,10,$\\
	$P_i = 11,13,24,$ $q_i = 19, 7,21,$ $R_i = r^i,$ $p_i = P_i,$ $p = 29,
	$\\
	$PQ_i = pq_i = 18,15, 1,$ $QP_i = qp_i = 25,28,30,$\\
	$QR_i = qr_i =12,27,22,$ $PR_i = pr_i = 14, 3, 2,$\\
	$AF_i = af_i =  6, 4, 9,$ $FA_i = fa_i =  0,16,17.$}

	Therefore the face-points are $FA_i,$ $AF_i;$ the vertex-points are
	$p_i,$ $q_i,$ 
	$r_i;$	the edge-points are $A_i,$ $PQ_i,$ $QP_i,$ $QR_i,$ $PR_i.$ 

	\sssec{Figure.}
        \begin{verbatim}
                                     Q2 
                                     21
                              QR1 o       o  PQ1 
                              27     FA0     15
                      R2  o           0             P0 
                      10                            11
                    o      QP2                 A1 o    .
                A0         30                  26         PR2 
                23              o P    QP0  R0 o             2
             o                  29 o 25 o 20                   .
         P2      FA2        QP1               QR2       AF1       Q1 
         24      17         28       FA1      22         4         7
                        R1             16          Q0 
           o             5                       19               .
          PQ0                 A2           PQ2 
          18          QR0          8            1         PR1        18
                  12                 P1                3
            .            AF2         13         AF0              .
                          9           o          6
            7                             PR0                       24
                2                     14                    23
                      11               o             10
                             15             27
                                     21
	\end{verbatim}

	\sssec{Comment.}
	We observe that in the dodecahedron, $FA_i$ are adjacent and if $AF_i$
	are constructed as in \ref{sec-tdodeca}, that these are not.  Moreover,
	$FA_i,$	$FA_{i+1}$ and $AF_{i+1}$ are adjacent; $AF_i,$ $AF_{i-1},$
	$FA_{i-1}$ are not.
	Therefore $AF_{i+1}$ and $AF_{i-1}$ are not adjacent to $FA_{i+1}$
	and $AF_i$ and therefore are adjacent to $FA_i$ and $FA_{i-1}$.  There
	is therefore a consistent way to define adjacency of conical points, if,
	given 3 of them named $FA_i$, the 3 others are labelled according to the
	construction \ref{sec-tdodeca}.

	\sssec{Definition.}
	If 3 conical points are labelled $FA_i$ and the 3 others, $AF_i,$ are
	labelled according to the construction \ref{sec-tdodeca}.  The triples
	$FA_i;$	$FA_i,$ $FA_{i+1}$, $AF_{i+1};$ $AF_{i+1},$ $AF_{i-1},$
	$FA_i$ and $AF_{i+1},$ $AF_{i-1},$ $FA_{i-1}$
	are {\em adjacent} and the other triples are {\em not adjacent}.

	The notion of "adjacent" and "not adjacent" can be interchanged.

	\ssec{Difference Sets with a Difference.}

	\setcounter{subsubsection}{-1}
	\sssec{Introduction.}
	After introducing distances in $n$ dimensional affine
	geometry and the associated selector, it occured to me that we could
	consider other difference sets for sets associated to $p^k$
	by choosing polynomials which are not irreducible. I discuss here
	briefly the extension to difference sets appropriate to the study
	of geometries in 2 and higher dimensions.

	\sssec{Definition.}
	A {\em difference set} associated to $q = p^k$ and to a
	polynomial of degree 3 with one root, is a set of $q$
	integers $\{s_0,$ $s_1,$  $\ldots$  , $s_{q-1}\}$ such that the $q^2-q$
	differences $s_i - s_j,$ 
	$i\neq	 j$ modulo $n =  q^2-1$ are distinct and different from 0
	modulo $q+1$.\\
	A {\em difference set} associated to $q = p^k$ and to a
	polynomial of degree 3 with two roots, is a set of $q-1$
	integers $\{s_0,$ $s_1,$  $\ldots$  , $s_{q-2}\}$ such that the
	$q^2-3q+2$\\
	$= (q-1)(q-2)$ differences $s_i - s_j,$ $i\neq j$ modulo $n = q^2-q$ are
	distinct and different from 0 modulo $q$ and modulo $q-1$.

	When applied to Geometry, I will prefer the terminology of Lemay and
	use the synonym {\em selector}.  The elements of the selector are called
	{\em selector numbers}.

	\sssec{Theorem.}
	{\em There exists always a polynomial $P$ of degree 3 with one root or 2
	roots such that $I$ is a generator of the multiplicative group of
	polynomials, of degree at most 2, with coefficients in $Z_p$, normalized
	to have the
	coefficient of the highest power 1, which are relatively prime to $P$.
	The selector numbers are the powers of $I$ modulo $P$ which are
	polynomials of degree at most 1.}

	The proof can be adapted easily from that of the irreducible case and
	is left as an exercise.

	\sssec{Example.}
	For $p = 3,$ $P = I^3 + I + 1,$\\
	$I^0 = 1,$ $I^1 = I,$ $I^2 = I^2,$ $I^3 = I + 1,$\\
	$I^4 = I^2 + I,$ $I^5 = I^2 - I - 1,$
	$I^6 = I^2 - I + 1,$ $I^7 = I^2 + 1,$\\
	and we have $I^8 = 1.$ Therefore the selector is 0, 1, 3.

	\sssec{Example.}
	The following are difference sets associated
	with $q = p^k:$ \\
	For $p = 3:$ $I^3+I+1,$ root 1, selector 0,1,3 $\pmod{8}.$\\
	\hth$I^3+I^2-I-1,$ roots 2,2,1, selector 0,1 $\pmod{6}.$\\
	For $p = 5:$ $I^3 - I -1,$ root 2, selector 0,1,3,11,20 $\pmod{24}.$\\
	\hth$I^3 - 2 I -1,$ roots 3,3,4, selector 0,1,3,14 $\pmod{20}$.\\
	For $p = 7:$ $I^3-I^2-2,$ root 5, selector 0,1,7,11,29,34,46 $\pmod{48}
	.$\\
	\hth$ I^3-3I-2,$ roots 2,6,6, selector 0,1,3,11,16,20 $\pmod{42}.$\\
	For $q = 11:$ $I^3 - I -1,$ root 6,\\
	\hth selector 0,1,3,28,38,46,67,90,101,107,116 $\pmod{120}.$\\
	\hth$I^3 - I^2 - I -1,$ roots 7,7,9,\\
	\hth selector 0,1,9,15,36,38,43,62,94,107 $\pmod{110.}$

%
	\sssec{Definition.}
	The {\em selector function} $f$ associated to the selector $\{s_i\}$ 
	is the function from ${\Bbb Z}_n$ to ${\Bbb Z}_n$ 
		$f(s_j-s_i) = s_i,$ $i\neq j,$
		for all other values $f(l) = -1.$

	\sssec{Example.}
	For $p = 5,$ the selector function associated with
	\{0, 1, 3, 11, 20\} is\\
	$\begin{array}{lcrrrrrrrrrrrr}
	i&\vline&0& 1& 2& 3& 4& 5& 6& 7& 8& 9&10&11\\
	\hline
	f(i) &\vline&-1& 0& 1& 0&20&20&-1&20& 3&11& 1& 0\\[10pt]
	i&\vline& 12&13&14&15&16&17&18&19&20&21&22&23\\
	\hline
	f(i) &\vline&-1&11&11&20&11& 3&-1& 1& 0& 3& 3& 1
	\end{array}$

	\sssec{Theorem.}
	\enumb
	\item {\em If the defining polynomial has 1 root then the selector has
	$n := p$ elements, the selector function has $p^2-1$ elements and is
	$-1$ for the $p-1$ multiples of $p+1.$}
	\item {\em If the defining polynomial has 2 distinct roots then the
	selector has $n := p-1$ elements, the selector function has $p(p-1)$
	elements and is $-1$ for the $2p-1$ multiples of $p$ and $p-1.$}
	\enume

	\sssec{Theorem.}
	\enumb
	\item {\em If $f(i-j)\neq  -1$ then $f(j-i) -i = f(i-j) - j.$}
	\enume

	Points, lines and incidence in the 2 dimensional geometry associated
	with $q = p^k$ and $n := q^2 + q + 1$ are defined as follows.

	\sssec{Definition.}
	The {\em points} are elements of the set \{0, 1,  \ldots , n-1\},
	the {\em lines} are elements of the set \{0, 1,  \ldots  , $n-1$\},
	a point  $a$  is {\em incident} to a line  $b$  iff  $f(a+b) = 0.$

	\sssec{Notation.}
	The points are denoted by a lower case letter or by an
	integer in ${\Bbb Z}_n.$  The lines are denoted by a lower case letter
	or by an
	integer in ${\Bbb Z}_n$ followed by an asterix.  The line incident to
	the points $a$ and $b$ is denoted $a \times b,$ the point incident to
	the lines $a*$ and $b*$ is denoted $a* \times b*.$

%
	We leave as an exercise to state and prove Theorems analogous to
	those in Section \ref{sec-Sselector}.

	\sssec{Definition.}
	The {\em dual affine plane}, is a Pappian plane in which we prefer
	the ``special" points which are those on a line $l$
	and a point $P$ not on $l$ and the ``special" lines which are
	those through $P$ and the line $l$.

	Th dual affine geometry can be studied by associating with it
	a polynomial which has 1 root. I give here some examples of Desargues, 
	Pappus and Pascal configurations.

	I illustrate Pappus and Desargues configurations using the notation of
	\ref{sec-npappus} and of \ref{sec-tdes} and give the points on
	a conic obtained using Pascal's construction.
	
	\sssec{Example.}
	For $p = 11,$ with the polynomial $I^3-I-1$, we have\\
	Pappus$(\langle 89,51,79\rangle,69*,\langle 33,88,110,\rangle,13*;
	\langle 92,71,6\rangle,95*)$,\\
	with $A_{i+1}\times B_{i-1} = (56*,87*,32*)$ and
	$A_{i-1}\times B_{i+1} = (28*,77*,115*).$\\
	Desargues$(98,\{7,1,70\},\{37*,31*,100*\},\{60,98,73\},\{50*,47*,60*\};
\\
	\langle 70,76,7\rangle,\langle 60*,89*,50*\rangle,31*),$\\
	Desargues$(111,\{115,69,13\},\{54*,33*,51*\},\{41,119,10\},
\{91*,80*,117*\};\\
	\langle 67,68,70\rangle,\langle 5*,47*,110*\rangle,53*)$,\\
	From Pascal's construction we obtain the following points are on a
	conic: 9,10,33,51,58,60,74,77,79,87,96,98.

	\sssec{Exercise.}
	Define a geometry corresponding to a polynomial which has 2 roots.
	\ssec{Generalization of the Selector Function for higher dimension.}
	\setcounter{subsubsection}{-1}
	\sssec{Introduction.}
	I will briefly stae one result for dimensions 3 and 4 concerning
	defining polynomials associated to the non irreducible case
	and illustrate for dimnesions 3, 4 and 5.

%
	\sssec{Theorem.}
	{\em If the $P_i$ denotes a primitive polynomial of degree $i.$}
	\enumb
	\item {\em For $k = 3,$ the defining polynomials $P$ can have the
	following form,\\
	\hth$	P_4, P_1P_3, P_1^2P_2,$\\
	there are $p^4+p^3+p^2+p+1,$ $p^4-1,$ $p^4-p$ polynomials relatively
	prime to $P,$ in these respective cases.}
	\item {\em For $k = 4,$ the defining polynomials $P$ can have the
	following form,\\
	\hth$	P_5, P_1P_4, P_1^2P_3, P_2P_3.$\\
	there are $p^5+p^4+p^3+p^2+p+1,$ $(p^3-1)(p+1),$ $p^5-1,$ $p^5-p$
	polynomials relatively prime to $P,$ in these respective cases.}
	\enume

	Proof: The polynomials in the sets are those which are relatively prime
	to the defining polynomial.  There are $p^k$ homogeneous polynomials of
	degree $k.$  If, for instance, $k = 4$ and the defining polynomial $P$
	is $P_2P_3,$ there are $p^2+p+1$ polynomials which are multiple of
	$P^2$ and $p+1,$ which are multiples of $P_3,$ hence
	$p^4+p^3+p^2+p+1-(p^2+p+1)-(p-1)$ 
	polynomials relatively prime to $P.$

	\sssec{Example.}
	$a_0,$ $a_1,$\ldots$a_k$  represent $I^{k+1}-a_0I^k-a_1I^{k-1} -a_k$.\\
	$\begin{array}{ccccrl}
	k	&p	&period	&def.\:pol.  & | sel.| 	&roots\:of\:def.\:pol.\\
	\cline{1-6}
	3	&3	&40	&2,1,1,1	&    13		&--\\
		&	&26	&1,1,1,1	&     9		&1\\
		&	&24	&0,1,1,1	&     8		&2,2\\
		&5	&156	&1,2,0,2	&    31		&--\\
		&	&124	&1,0,0,2	&    25		&4\\
		&	&120	&0,0,1,2	&    24		&4,4\\
		&7	&400	&0,1,1,4	&    57		&--\\
		&	&342	&0,0,1,1	&    49		&3\\
		&	&336	&0,0,3,1	&    48		&5,5\\
		&11	&1464	&0,0,2,5	&    133	&--\\
		&	&1330	&0,0,1,1	&    121	&3\\
		&	&1320	&1,5,2,4	&    120	&1,1\\[5pt]

	4	&3	&121	&2,0,0,0,1	&40	&--\\
		&	&104	&0,1,0,0,1	&35	&
(I^2+I-1)(I^3-I^2+I+1)	\\
		&	& 80	&0,2,0,0,1	&27	&2\\
		&	& 78	&1,0,0,0,1	&26	&2,2\\
		&5	&781	&4,0,0,0,1	&156	&--\\
		&	&744	&2,2,0,0,1	&149	&
(I^2+I+2)(I^3+2I^2-I+2)	\\
		&	&624	&2,0,0,0,1	&125	&3\\
		&	&620	&3,0,1,0,1	&124	&3,3\\
		&7	&2801	&3,0,0,0,1	&400	&--\\
		&	&2736	&6,0,0,0,1	&391	&
(I^2+2I-2)(I^3-I^2-3I-3)\\
		&	&2400	&3,1,0,0,1	&343	&3\\
		&	&2394	&0,3,3,0,1	&342	&5,5\\
		&11\\[5pt]
		&11	&16105	&0,0,10,0,9	&1464	&--\\
		&	&15960	&0,0,0,10,8	&1451	&
(I^2+I-)(I^3-I^2-I-)\\
		&	&14640	&0,0,0,10,9	&1331	&10\\
		&	&14630	&0,0,0,9,7	&1330	&3,3\\
		&13	&30941	&8,0,0,0,1	&2380	&--\\
		&	&30744	&5,0,0,0,1	&2365	&
(I^2-3I+6)(I^3-2I^2+I+2)\\
		&	&28560	&2,0,0,0,1	&2197	&11\\
		&	&28548\\[5pt]

	5	&3	&364	&1,0,0,0,0,1	&121	&--\\
		&	&242	&1,1,0,0,0,1	& 81	&2\\
		&	&240	&1,2,1,0,0,1	& 80	&2,2\\
	\end{array}$

	\sssec{Definition.}
	Given a selector $s,$ the {\em selector function} associates to
	the integers in the set ${\Bbb Z}_n$ a set of $p+1$ integers or $p$
	integers obtained as follows,\\
	\hth$	s(j) \in f_i$  iff  $sel(l) - sel(j) = i$ for some $l.$

	\sssec{Theorem.}
	\enumb
	\item{\em $f(i)$ is the set of points on the line $i^* \times 0^*.$}
	\item{\em $f(i)-j,$ where we subtract $j$ from each element in the
		set, is the set of points in $(i+j)^*\times j^*,$ equivalently}
	\item{\em $f(i-j)-j,$ is the set of points in  $i^* \times j^*.$}
	\item$	a^*\times b^*\times c^* = ((a-i)^*\times (b-i)^*\times (c-i)^*)
		 - i.$
	\enume

%
%
	\sssec{Theorem.}
	\enumb
	\item[0.~~]{\em If the defining polynomial is primitive, then}
	\item[  0.]$| s|  = \frac{p^k-1}{p-1},$
	\item[  1.] if $i \notequiv \frac{p^k-1}{p-1},$ $| f(i)|  = p+1.$
	\item[1.~~]{\em If the defining polynomial has one root, then}
	\item[  0.]$| s|  = p^k,$
	\item[  1.]	if $i \neq  0,$ $| f(i)|  = p,$
	\item[2.~~]{\em If the defining polynomial has a double root, then}
	\item[  0.] $| s|  = p^k-1,$
	\item[  1.] if $i \notequiv p,p^2-1,$ $| f(i)|  = p,$
	\item[  2.]	if $i \equiv p$ and $i\neq 0,$ $| f(i)|  = p-1,$
	\enume

	\sssec{Example.}\label{sec-eric3select}
	\enumb
	\item$ k = 3,$ $p = 3,$ defining polynomial $I^4 - 2 I^3 - I^2 - I - 1
	=(I-1)(I^3-I+1),$\\
	selector: $\{0,1,2,9,10,13,15,16,18,20,24,30,37\}$\\
	selector function:\\
	$\begin{array}{rrrrrcrrrrrcrrrrr}
	 0&-1&-1&-1&-1&\vline&14& 1& 2&10&16&\vline&28& 2& 9&13&30\\
	 1& 0& 1& 9&15&\vline&15& 0& 1& 9&15&\vline&29& 1&13&20&24\\
	 2& 0&13&16&18&\vline&16& 0& 2&24&37&\vline&30& 0&10&20&30\\
	 3&10&13&15&37&\vline&17& 1&13&20&24&\vline&31& 9&10&18&24\\
	 4& 9&16&20&37&\vline&18& 0& 2&24&37&\vline&32& 9&10&18&24\\
	 5&10&13&15&37&\vline&19& 1&18&30&37&\vline&33& 9&16&20&37\\
	 6& 9&10&18&24&\vline&20& 0&10&20&30&\vline&34&15&16&24&30\\
	 7& 2& 9&13&30&\vline&21& 9&16&20&37&\vline&35& 2&15&18&20\\
	 8& 1& 2&10&16&\vline&22& 2&15&18&20&\vline&36& 1&13&20&24\\
	 9& 0& 1& 9&15&\vline&23& 1&18&30&37&\vline&37& 0&13&16&18\\
	10& 0&10&20&30&\vline&24& 0&13&16&18&\vline&38& 2&15&18&20\\
	11& 2& 9&13&30&\vline&25&15&16&24&30&\vline&39& 1& 2&10&16\\
	12& 1&18&30&37&\vline&26&15&16&24&30&\vline\\
	13& 0& 2&24&37&\vline&27&10&13&15&37&\vline
	\end{array}$

	\item$k = 3,$ $p = 3,$ defining polynomial $I^4 - I^3 - I^2 - I - 1
	=(I-1)^2(I^2+I-1),$\\
	selector: $\{0,1,2,8,11,18,20,22,23\}$\\
	selector function:\\
	$\begin{array}{rrrrcrrrrcrrrrcrrrr}
	 0&-1&-1&-1&\vline& 7& 1&11&20&\vline&14& 8&20&23&\vline&21& 1& 2&23\\
	 1& 0& 1&22&\vline& 8& 0&18&20&\vline&15& 8&11&22&\vline&22& 0& 1&22\\
	 2& 0&18&20&\vline& 9& 2&11&18&\vline&16& 2&11&18&\vline&23& 0&11&23\\
	 3& 8&20&23&\vline&10& 1& 8&18&\vline&17& 1&11&20&\vline&24& 2&20&22\\
	 4&18&22&23&\vline&11& 0&11&23&\vline&18& 0& 2& 8&\vline&25& 1& 2&23\\
	 5&18&22&23&\vline&12& 8&11&22&\vline&19& 1& 8&18&\vline\\
	 6& 2&20&22&\vline&13&-1&-1&-1&\vline&20& 0& 2& 8&\vline
	\end{array}$

	\item$  k = 3,$ $p = 3,$ defining polynomial $I^4 - I^2 - I - 1.$\\
	selector: $\{0,1,2,4,14,15,19,21\}$\\
	selector function:\\
	$\begin{array}{rrrrcrrrrcrrrrcrrrr}
	 0&-1&-1&-1&\vline& 6&15&19&-1&\vline&12& 2&14&-1&\vline&18& 1&21&-1\\
	 1& 0& 1&14&\vline& 7&14&19&21&\vline&13& 1& 2&15&\vline&19& 0& 2&19\\
	 2& 0& 2&19&\vline& 8&-1&-1&-1&\vline&14& 0& 1&14&\vline&20& 1& 4&19\\
	 3& 1&21&-1&\vline& 9&15&19&-1&\vline&15& 0& 4&-1&\vline&21& 0& 4&-1\\
	 4& 0&15&21&\vline&10& 4&14&15&\vline&16&-1&-1&-1&\vline&22& 2& 4&21\\
	 5&14&19&21&\vline&11& 4&14&15&\vline&17& 2& 4&21&\vline&23& 1& 2&15
	\end{array}$
	\enume

	\sssec{Example.}
	In the case of Example \ref{sec-eric3select}.0.
	If we denote by $i^\dagger,$ the lines $0^* \times i^*,$ these lines,
	which are sets of 4 points can all be obtained from\\
	$1^\dagger	= \{0,1,9,15\},$ $2^\dagger = \{0,13,16,18\},$
	$4^\dagger = \{9,16,20,37\}$ and\\
	$10^\dagger	= \{0,10,20,30\}$ by adding an integer modulo $n$.\\
	$1^\dagger	+ 0 = 1^\dagger,9^\dagger,15^\dagger,$ $1^\dagger + 1
	= 39^\dagger,8^\dagger,14^\dagger,$\\
	\hth$1^\dagger + 9 = 6^\dagger,31^\dagger,32^\dagger,$ 
		$1^\dagger + 15 = 34^\dagger,25^\dagger,26^\dagger,$ \\
	$2^\dagger	+ 0 = 2^\dagger,24^\dagger,37^\dagger,$ $2^\dagger
		+ 2 = 22^\dagger,35^\dagger,38^\dagger,$
	 $2^\dagger + 37 = 3^\dagger,5^\dagger,27^\dagger,$ \\
	\hth$	2^\dagger + 24 = 13^\dagger,16^\dagger,18^\dagger,$\\
	$4^\dagger	+ 0 = 4^\dagger,21^\dagger,33^\dagger,$ $4^\dagger
		+ 4 = 17^\dagger,29^\dagger,36^\dagger,$
		$4^\dagger + 21 = 12^\dagger,19^\dagger,23^\dagger,$ \\
	\hth$4^\dagger + 33 = 7^\dagger,11^\dagger,28^\dagger,$\\
		$10^\dagger + 0 = 10^\dagger,20^\dagger,30^\dagger.$ 


	\ssec{The conics on the dodecahedron.}

	\setcounter{subsubsection}{-1}
	\sssec{Introduction.}
	The reader may want to skip this section until he has
	become familiar with conics.  In it, we summarize the various types and
	sub-types of conics as they relate to the representation of the
	finite projective plane, for $p = 5,$ on the dodecahedron.  We will see
	later, IV.1.12.  that the dodecahedron can also be used to represent
	the finite polar and the finite non-Euclidean geometry, for $p = 5.$

	\sssec{Definition.}
	The {\em conics are all of the same type} if the
	classification into face-points, vertex-points and edge-points is the
	same.  The {\em conics  are of the same sub-type} if they can be
	derived from each other using any of the 60 collineations which exchange
	face-points.

	\sssec{Notation.}
	In the following Theorem we use the notation\\
	``60 $fffvve,$ 30 $C1$ (6,9,17;11,29;22), 30 $C2$ (9,16,17;7,13;15)."\\
	to indicate that we have 60 conics with 3 face-points, 2 vertex-lines
	and 1 edge-line.  These are of the sub-type $C1$ and $C2.$  An example
	of a
	conic of a given sub-type is provided in parenthesis, ``;" separates
	points of a different classification, these points are given in the
	order, face-point, vertex-point, edge-point, and in the same
	classification in increasing order.\\
	A pictorial representation of the sub-types is given in Figure
	\ref{sec-fdodeca}.

	\sssec{Theorem.}
	{\em The $31.30.25.16.\frac{6}{6!} = 3100$ conics are of the following
	type and sub-type.}\\
	$\begin{array}{rrrllrll}
	  1&ffffff,& 1& A&(0,4,6,9,16,17).\\
	 30&ffffee,&30& B&(6,9,16,17;2,22).\\
	 60&fffvve,&30& C1&(6,9,17;11,29;22),&  30& C2&(9,16,17;7,13;15).\\
	120&fffvee,&30& D1&(9,16,17;7;1,25),&   30& D2&(9,16,17;5;1,18),\\
	&&	    30& D3&(0,4,6,;5;1,4),&     30& D4&(0,4,6;13;18,30).\\
	 30&ffvvvv,&15& E1&(0,4;5,11,13,20),&   15& E2&(0,4;7,19,21,29).\\
	 60&ffvvve,&60& F &(0,4;11,13,29;8).\\
	360&ffvvee,&30& G1&(0,4;11,20;1,14),&   30& G2&(0,4;5,13;3,30),\\
	&&	    30& G3&(0,4;7,21;18,27),&   30& G4&(0,4;19,21;12,28),\\
	&&	    60& G5&(0,4;11,29;1,22),&   60& G6&(0,4;5,11;14,18),\\
	&&	    60& G7&(0,4;21,29;1,27),&   60& G8&(0,4;5,21;2,27).\\
	180&ffveee,&30& H1&(0,4;11;8,22,25),&     30& H2&(0,4;13;8,15,22),\\
	&&	    60& H3&(0,4;21;2,18,28),&   60& H4&(0,4;21;1,22,30).\\
	135&ffeeee,&15& I1&(0,4;2,15,22,25),&   30& I2&(0,4;2,12,14,15),\\
	&&	    30& I3&(0,4;15,18,22,30),&  60 &I4&(0,4;1,15,25,30).\\
	 12&fvvvvv,& 6& J1&(16;7,10,11,21,24),&  6& J2&(16;5,13,19,20,29).\\
	120&fvvvve,&60& K1&(16;5,10,13,19;27),& 60& K2&(6;5,7,10,11;26).\\
	300&fvvvee,&30& L1&(4;10,11,20;1,12),&  30& L2&(0;5,10,13;27,30),\\
	&&	    60& L3&(16;11,21,24;12,18),&60& L4&(0;10,11,20;1,15),\\
	&&	    60& L5&(0;7,10,21,3,22),&   60& L6&(0;5,11,13;12,26).\\
	480&fvveee,&30& M1&(9;11,21;2,14,28),&  30& M2&(0;10,20;25,28,30),\\
	&&	    60& M3&(17;19,21;2,23,25),& 60& M4&(0;19,24;2,14,15),\\
	&&	    60& M5&(9;10,20;1,15,18),&  60& M6&(16;20,29;12,27,28),\\
	&&	    60& M7&(9;7,29;8,12,30),&   60& M8&(16;24,21;1,12,26),\\
	&&	    60& M9&(17;7,21;3,8,14).\\
	480&fveeee,&60& N1&(0;10;3,15,27,30),&  60& N2&(0;10;14,15,18,30),\\
	&&	    60& N3&(0;10;12,14,22,30),& 60& N4&(0;13;1,23,27,30),\\
	&&	    60& N5&(0;10;1,2,27,28),&   60& N6&(0;13;1,26,28,30),\\
	&&	    60& N7&(0;13;2,3,15,22),&   60& N8&(0;13,2,12,22,23).\\
	 12&feeeee,& 6& O1&(16;1,8,22,25,28),&   6& O2&(16;3,12,14,26,30).\\
	 10&vvvvvv,&10& P &(5,7,10,11,21,29).\\
	 60& vvvvee,&30& Q1&(10,11,20,29;15,27),&30& Q2&(29,5,13;21;1,12).\\
	240& vvveee,&30& R1&(10,21,29;3,18,23),& 30& R2&(10,11,20;23,27,30),\\
	&&	    30& R3&(10,21,29;14,26,28),&30& R4&(10,20,24;3,18,23),\\
	&&	    60& R5&(10,21,29;3,25,28),& 60& R6&(11,13,29;1,8,15).\\
	270& vveeee,&15& S1&(11,20;3,18,27,30),& 15& S2&(11,20;2,15,22,25),\\
	&&	    30& S3&(11,20;1,2,12,22),&  30& S4&(11,20;1,14,18,30),\\
	&&	    30& S5&(7,21;18,22,25,27),& 30& S6&(7,21;2,12,14,15),\\
	&&	    30& S7&(7,21;3,12,14,30),&  30& S8&(7,21;8,12,14,26),\\
	&&	    60& S9&(7,21;2,3,18,25).\\
	120& veeeee,&30& T1&(21;3,12,14,22,28),& 30& T2&(5;3,8,22,26,28).\\
	&&	    60& T3&(24;3,12,14,22,27).\\
	 20& eeeeee,& 10& U1&(1,2,14,18,25,30),&  10& U2&(3,12,14,15,25,27).
	\end{array}$

	The proof of the decomposition into types was done using a computer
	program 
	which took 22 minutes to run an an IBM PC.

	\sssec{Figure.}\label{sec-fdodeca}
	The pictorial representation of a conic of a given sub-type
	on the dodecahedron is as follows.
{\footnotesize \begin{verbatim}
A: the 6 faces.


B:     f        .        f

                .
    .        .     .        .
          .           .

           o         o
  f                          f
            .   .   .

        .              .


C1:           .  o  .               C2:    f        .        f

            .         .                             .
                            f           .        .     .        .
           o           .                      o           o
              .     .                               f
      .          .          .                  .         .

  .       f      .      f        .              .   o   .

      .          .          .               .              .
              .     .
           .           o


\end{verbatim}

\begin{verbatim}
D1:             o                D2:           .  .  .

      f         .        f                   .         .
                                                  f
                .                           .           .
    .       .       .        .                 .     .
         .             .               o          o          o
                f
          o           o            .       f      .      f        .

            .   .   .                  .          .          .
                                               .     .
        .              .                    .           .


D3:           .  .  .               D4:           .  .  .

            .         .                         o         o
     f           f          f            f                      f
           .           .                       .           .
              o     o                             .     .
      .          .          .             .          .          .

  .              .               .    .              .               .

      .          o          .             .          o          .
              .     .                             .     .
           .           .                       .           .


\end{verbatim}

\begin{verbatim}
E1:           o  .  o                E2:
                                                o           o
            .         .                            .     .
                                           .          .          .
           .           .
              .     .                  .       f      .      f        .
      .          o          .
                                           .          .          .
  .       f      .      f        .                 .     .
                                                o           o
      .          o          .
              .     .
           .           .


F:            .  o  o

            .         .

           .           .
              .     .
      .          o          .

  .       f      .      f        .

      .          .          .
              .     .
           .           o

\end{verbatim}

\begin{verbatim}

            .         o
G1:                                 G2:            o  .  o
           .           .
              .     .                            .         .
      .          o          .
                                                .           .
  .       f      .      f        .                 .     .
                                           o          .          o
      .          o          .
              .     .               .          f      .      f        .
           .           .
                                           .          .          .
            o         .                            .     .
                                                .           .


                                                .         o
G3:                                 G4:
           o           o                       o           .
              .     .                             .     .
      o          .          o             .          .          .

  .       f      .      f        .    .       f      .      f        .

      .          .          .             .          .          .
              .     .                             .     .
           .           .                       .           o

                                                o         .


\end{verbatim}

\begin{verbatim}
            o         .
G5:                                 G6:            o  .  .
           .           o
              o     .                            .         o
      .          .          .
                                                .           .
  .       f      .      f        .                 .     .
                                           o          o          .
      .          o          .
              .     .               .          f      .      f        .
           .           .
                                           .          .          .
            .         .                            .     .
                                                .           .


            o         .
G7:                                 G8:            o  .  .
           .           o
              .     .                            .         .
      .          .          .
                                                .           o
  .       f      .      f        .                 o     .
                                           .          .          o
      .          .          o
              .     .               .          f      .      f        .
           .           o
                                           .          .          .
            .         .                            .     .
                                                .           .

\end{verbatim}

\begin{verbatim}

H1:            .  o  .              H2:            .  o  o

             .         .                         .         .

            .           .                       .           .
               .     .                             .     o
       .          o          .             .          .          .

.          f      .      f        . .          f      .      f        .

       .          .          .             .          .          .
               o     o                             o     .
            .           .                       .           .


             o         .                         .         .

\end{verbatim}

\begin{verbatim}

H3:            .  .  .              H4:            .  .  .

             .         .                         o         .

            .           o                       .           .
               o     .                             o     .
       o          .          .             .          .          o

.          f      .      f        . .          f      .      f        .

       .          .          .             .          .          .
               .     .                             .     .
            .           .                       .           o


I1:                                 I2:           .  .  .
           .           .
              o     o                           o         o
      .          .          .
                                               .           .
  .       f      .      f        .                o     o
                                          .          .          .
      .          .          .
              o     o                 .       f      .      f        .
           .           .
                                          .          .          .
                                                  .     .
                                               .           .

\end{verbatim}

\begin{verbatim}

                                                 o
I3:                                 I4:
           .           .                       .           .
              .     o                             .     o
      o          .          .             .          .          o

  .       f      .      f        .    .       f      .      f        .

      .          .          o             .          .          .
              o     .                             .     o
           .           .                       .           .

\end{verbatim}

\begin{verbatim}

J1:           .  .  .              J2:              .

            .         .                             o
                                        .       .       .        .
           o           o                     o             o
              .     .                               f
      .          o          .                 .           .

  o              .               o              o   .   o

      .          .          .               .              .
              .     .
           .           .
                 f


K1:              f                  K2:            .  .  .
           o           o
              .     .                            .         .
      .          o          .                         f
                                                .           .
  .              .               .                 .     .
                                           .          o          .
      .          .          .
              o     .                  o              .               o
           o           .
                                           .          .          o
                                                   .     .
                                                o           .

\end{verbatim}

\begin{verbatim}

L1:            o  .  o              L2:            .  .  .

             .         .                         .         .
                  f                                   f
            .           .                       .           .
               .     .                             o     o
       o          .          o             .          o          .

  .               .               . o                 .               o

       .          o          .             .          .          .
               .     .                             .     .
            .           .                       .           .


L3:                                 L4:            o  .  o
           o           o
              .     o                            o         .
      .          .          .                         f
                                                .           .
  o              o               .                 .     .
                                           .          o          .
      .          .          .
              .     .               .                 .               .
           .           .
                 f                         o          .          .
                                                   .     .
                                                .           .

\end{verbatim}

\begin{verbatim}

L5:            o  .  o              L6:            .  .  .

             .         .                         .         .
                  f                                   f
            .           .                       .           .
               .     .                             .     o
       .          .          .             .          o          .

.                 o               o o                 .               .

       .          .          .             .          .          .
               .     o                             p     .
            .           .                       o           .


M1:             o                   M2:             o

                .                                   .
    .       .       .        .          .       o       o        .
         .             .                     o             o
                f                                   f
          o           o                       .           .

            o   .   o                           .   .   .


\end{verbatim}

\begin{verbatim}

             .         .                         o         .

M3:            .  o  .              M4:            .  o  .

             .         .                         .         .
                  f                                   f
            .           .                       .           .
               .     .                             .     .
       o          .          o             .          .          .

.                 .               . o                 .               o

       .          .          .             .          .          .
               .     .                             .     .
            o           o                       .           .

             .         .                         o         .


M5:            o  .  o              M6:            o  .  o

             .         o                         .         o
                  f                                   f
            .           .                       .           .
               .     .                             .     .
       .          o          .             .          .          o

.                 .               . .                 .               .

       .          .          o             .          .          .
               .     .                             o     .
            .           .                       .           .

\end{verbatim}

\begin{verbatim}

M7:            .  .  .              M8:            .  .  .

             o         .                         .         .
                  f                                   f
            .           .                       .           .
               o     .                             .     o
       .          o          .             .          o          o

o                 .               . .                 .               .

       o          .          .             .          .          o
               .     .                             .     .
            .           .                       o           .


M9:            .  .  .               M10:    .  o  .

             .         .                   .         .
                  f                             f
            .           .                 .           .
               .     .                       .     .
       o          o          o       o          .          .

.                 .               .             .               o    o

       .          .          o                  .          .
               .     .                       .     .
            .           o                 .           o

\end{verbatim}

\begin{verbatim}

                                                 o
N1:           .  .  .               N2:
                                               .           .
            o         .                           .     .
                 f                        o          o          o
           .           .
              o     o                 .       f      o               .
      .          o          .
                                          .          .          .
  .              .               .                .     .
                                               .           .
      .          .          .
              o     .
           .           .


N3:            .  .  .              N4:            o  o  .

             .         .                         .         .
                  f
            .           o                       .           .
               .     o                             o     o
       o          .          .             .          .          .

  .               .               .   .       f       o               .

       .          .          .             .          .          .
               .     o                             .     .
            .           .                       .           .


\end{verbatim}

\begin{verbatim}

             o         .                         .         .

N5:            .  .  .              N6:            .  .  .

             .         .                         o         o
                  f                                   f
            .           .                       .           .
               o     .                             .     .
       .          o          o             .          .          .

.                 .               . .                 .               .

       o          .          .             .          o          .
               .     .                             .     o
            .           .                       .           .

             .         .                         o         .


             o         .                         o         .

N7:            .  .  .              N8:            .  .  .

             .         .                         .         .
                  f                                   f
            .           .                       .           .
               o     .                             .     .
       o          .          .             o          .          o

.                 .               . .                 .               .

       .          o          o             o          o          .
               .     .                             .     .
            .           .                       .           .



             o

\end{verbatim}

\begin{verbatim}


       O1:              .                 O2:              o

                .                                  .
    .       o       o        .         o       .       .        o
         .             .                    .             .
                                                   f
          o           o                      .           .

            .   o   .                          .   .   .

        .              .                     o            o
                f


P:
           o           .
              .     .
      .          o          .

  o              .               o

      .          o          .
              .     .
           .           o

\end{verbatim}

\begin{verbatim}


       Q1:              .                  Q2:
                                              .           .
                .                                o     o
    .       o       o        .           .          o          .
         o             o
                                     .              .               .
          .           .
                                         .          o          .
            o   .   o                            .     o
                                              o           .
        .              .


R1:           .  .  .               R2:              o

            .         .                              o
                                         .       o       o        .
           o           o                      .             .
              .     .
      .          o          .                  .           .

  .              o                               o   .   o

      .          .          .                  .           .
              o     o
           .           .


\end{verbatim}

\begin{verbatim}


R3:            .  o  .              R4:            o  .  o

             .         .                         .         .

            o           o                       .           .
               .     .                             .     .
       o          o          o             .          .          .

.                 .               . .                 o               .

       .          .          .             .          o          .
               .     .                             o     o
            .           .                       .           .


R5:           .  .  .               R6:            .  .  .

            .         o                          o         .

           o           o                        o           o
              .     .                              .     o
      .          o          o              o          .          .

  .              .                  .                 .               .

      .          .          .              .          .          .
              o     .                              .     .
           .           .                        .           o

\end{verbatim}

\begin{verbatim}

            .         .
S1:                                 S2:
           .           .                       .           .
              .     .                             o     o
      o          o          o             .          o          .

  .              .               .    .              .               .

      o          o          o             .          o          .
              .     .                             o     o
           .           .                       .           .


\end{verbatim}

\begin{verbatim}

            o         .                         .         o
S3:                                 S4:
           .           .                       .           .
              o     .                             .     .
      .          o          .             o          o          .

  .              .               .    .              .               .

      .          o          .             .          o          o
              o     .                             .     .
           .           .                       .           .

            o         .                         o         .


             o

S5:            .  .  .              S6:              .

             .         .                             .
                                         .       o       o        .
            o           o                     o             o
               .     .
       o          .          o                 o           o

.                 .               .              .   .   .

       .          .          .               .              .
               o     o
            .           .

\end{verbatim}

\begin{verbatim}

S7:         .  .  .                 S8:            .  o  .

          o         o                            o         o

         o           o                          o           o
            .     .                                .     .
    .          .          .                .          .          .

.              .               .    .                 o               .

    o          .          o
            .     .
         .           .

\end{verbatim}

\begin{verbatim}
            .         .
S9:
           o           o
              o     .
      o          .          .

  .              .               .

      o          .          .
              .     o
           .           .

            .         .


T1:           .  .  .               T2:            o  o  .

            o         o                          o         .

           .           .                        .            .
              .     .                              .      o
      o          .          o              .          .           o

  .              o               .  . .               o               .

      .          o          .              .          .          .
              .     .                              .     .
           .           .                        .           .
\end{verbatim}

\begin{verbatim}

T3:           .  .  .

             o         .

            .           .
               .     .
       .          o          o

   .              o

       .          .          o
               o     .
            .           .

U1:           .  .  .               U2:            .  o  .

            o         o                          .         .

           .           .                        .           .
              .     .                              o     o
      o          .          o              .          .          .

  .              .               .  .                 o               .

      .          .          .              o          .          o
              o     o                              .     .
           .           .                        .           .
\end{verbatim}}

	The family of types of conics was determined interactively using
	a computer program.

	\ssec{The truncated dodecahedron.}
	\setcounter{subsubsection}{-1}
	\sssec{Introduction.}
	After defining convex uniform polyhedra, whose notion
	may go back to Archimedes and were fully studied by Kepler, we will
	show that one of them, the truncated dodecahedron can be used as a
	model for the finite projective plane of order $3^2.$ 

	\sssec{Definition.}
	A polyhedron with regular faces, in Euclidean 3-space is
	{\em uniform} if it has symmetry operations taking a given vertex into
	any other vertex, otherwize it is {\em non-uniform}.  If, in addition,
	all faces are congruent, the  polyhedra is {\em regular}.

	\sssec{Theorem. [Euclid]}
	{\em There are 5 convex regular polyhedra.}

	\sssec{Notation. [See Johnson.]}
	In the following Theorem, we use the following notation, developped by
	several Mathematicians. $\{n\}$ denotes a regular polygon with $n$
	sides, $(n.q.n.q)$ denotes a vertex with adjoining faces successively
	with $n,$ $q,$ $n,$ $q$ sides, $<n.q>$ denotes
	an edge ajoining a face with $n$ sides and one with $q$ sides.

	\sssec{Theorem. [Kepler]}
	{\em Besides regular prisms and antiprisms, there are
	13 convex uniform, non-regular polyhedra.

	These are}
\footnotesize{
	Name\hti{20}		Faces\hti{12}	Vertices\hti{10}	Edges\\
	$\begin{array}{llll}
	Cuboctahedron          &8\{3\},6\{4\}	&12(3.4.3.4)& 24<3.4>\\
	Icosidodecahedron      &20\{3\},12\{5\}&30(3.5.3.5)& 60<3.5>\\
	Truncated\:tetrahedron  &4\{3\},4\{6\}	&12(3.6^2)&   12<3.6>,6<6.6>\\
	Truncated\:octahedron   &6\{4\},8\{6\}	&24(4.6^2)&    24<4.6>,12<6.6>\\
	Truncated\:cube         &8\{3\},6\{8\}	&24(3.8^2)&    24<3.8>,12<8.8>\\
	Truncated\:icosahedron  &12\{5\},20\{6\}&60(5.6^2)&    60<5.6>,30<6.6>\\
	Truncated\:dodecahedron&20\{3\},12\{10\}&60(3.10^2)&60<3.10>,30<10.10>\\
	Rhombicuboctahedron    &8\{3\},18\{4\}	&24(3.4^3)&    24<3.4>,24<4.4>\\
	Rhombicosidodecahedron	&20\{3\},30\{4\},&60(3.4.5.4)& 60<3.4>,60<4.5>\\
				& 12\{5\}\\
	Truncated\:cuboctahedron	&12\{4\},8\{6\},	&48(4.6.8)&
	24<4.6>,24<4.8>,\\
				& 6\{8\}	&&		24<6.8>\\
	Truncated		&30\{4\},20\{6\},&120(4.6.10)& 60<4.6>,
	60<4.10>,\\
	\hti{4}icosidodecahedron&12\{10\}		&&	60<6.10>\\
	Snub cuboctahedron	&32\{3\},6\{4\}&24(3^4.4)&    36<3.3>,24<3.4>\\
	Snub icosidodecahedron	&80\{3\},12\{5\}&60(3^4.5)&   90<3.3>,60<3.5>\\
	n-gonal\:prism		&n\{4\},2\{n\}&2n(4^2.n)&    n<4.4>,2n<4.n>\\
	n-gonal\:antiprism	&2n\{3\},2\{n\}&2n(3^3.n)&    2n<3.3>,2n<3.n>
	\end{array}$}

\normalsize
	\sssec{Theorem. [N. W. Johnson]}
	{\em There are 92 convex non-uniform regular-faced polyhedra.}

	The fact that all vertices are of the same type does not insure
	uniformity, as the example of the elongated square gyrobicupola of
	J. C. P. Miller shows.  This non-uniform polyhedra has the same
	characteristics as the rhombicuboctahedron, but has the part below
	the 8 squares turned 45 degrees.

	Before discussing the truncated dodecahedron as a model for the
	Pappian plane associated with $3^2$, I will discuss the pentagonal
	antiprism as a model for the Pappian plane associated with $2^2$.

	\sssec{Notation.}
	I identify elements which are 
	symmetrical with respect to the center of the antiprism.
	For the pentagonal antiprism, with $i = 0,1,2,3,4,$ I will denote by
	$t_i,$ the 5 triangular faces, by $v_i,$ the 5 vertices, by $e_i,$ the
	5 pentagonal-triangular edges, by $f_i,$ the 5 triangular-triangular
	edges and by $p,$ pentagonal face. We have altogether 21 elements
	to represents the 21 points in the plane associated with $2^2$.

	\sssec{Theorem.}
	For $q = 2^2,$
	\enumb
	\item {\em The selector is} $\{0,1,4,14,16\}.$
	\item {\em The corresponding selector function $f$ is, and the 
	representation of the points on the antiprism are}\\
	$\begin{array}{ccrrrrrrrrrrr}
	i&\vline&0&1&2&3&4&5&6&7&8&9&10\\
	\hline
	f(i)&\vline&0&0&14&1&0&16&16&14&14&16&4\\
	repr.\:of\:points&\vline&v_2&e_2&v_4&f_3&f_2&t_3&t_1&v_0&v_1&e_0&e_3\\
	repr.\:of\:lines&\vline&v_2&t_2&v_4&f_3&f_2&e_3&e_1&v_0&v_1&t_0&t_3
\\[10pt]
	i&\vline&&11&12&13&14&15&16&17&18&19&20\\
	\hline
	f(i)&\vline&&14&4&1&0&1&0&4&4&16&1\\
	repr.\:of\:points&\vline&&v_3&f_4&f_1&p&t_0&t_2&e_1&f_0&e_4&t_4\\
	repr.\:of\:lines&\vline&&v_3&f_4&f_1&p&e_0&e_2&t_1&f_0&t_4&e_4\\
	\end{array}$
	\item {\em The incidence properties are}\\
	$e_i^* \incid v_i,e_{i\pm 2},t_{i\pm 1}$,\\
	$t_i^* \incid v_i,t_{i\pm 2},f_{i\pm 1}$,\\
	$f_i^* \incid v_i,e_{i\pm 1},f_{i\pm 1}$,\\
	$v_i^* \incid p,v_i,e_i,t_i,f_i$,\\
	$p \incid v_i.$
	\enume

	Proof: I leave as an exercise the determination of the fundamental
	polynomial and the corresponding selector.\\
	The selector function follows easily from its definition.\\
	The selector polarity which associates $i$ to $i^*$ has the fixed points
	7,8,11,0 and 2 on the line $14^*$. I will associates to 14 and to
	$14^*$ the pentagonal face and its incident points or lines to $v_i$.
	Starting from that, one of the possible solution is given in 1.
	Notice that I use the same correspondance between $e_i$, $v_i$ and $f_i$
	for the points and the lines but exchange $t_i$ and $e_i$ to get
	the corresponding points and lines.

	\sssec{Figure.}
	The corresponding drawing for the Projective plane over $2^2$ is given
	page \ldots.

	\sssec{Notation.}
	For the truncated dodecahedron, I will denote by $t,$ a
	triangular face, by $d,$ a decagonal face, by $v,$ a vertex, by $e,$ a
	$<10,10>$ edge and by $u,$ a $<3.10>$ edge.  The lower case notation is
	used indifferently  for points and lines, the upper case notation for
	points.

	\sssec{Lemma.}
	{\em If $x_Y$ denotes the number of points of type $Y$ incident to a
	line of type $x,$ then\\
	\hth$2 | f_Y,$ $5 | d_Y,$ $3 | t_Y$ for $Y\neq  T,$ $3 | t_T-1.$} 

	Proof:  For instance, there are 10 $T$-points, each is adjacent to 10
	lines; on the other hand, the 30 $e$-lines are adjacent to 30 $e_T$ 
	triangles, the 15 $f$-lines are adjacent to 15 $f_T$ triangular-points,
	 $\ldots$  .  This implies\\
	\hth$30 e_T + 15 f_T + 30 v_T + 10 t_T + 6 d_T = 100,$
	which gives, modulo $2,$ $f_T = 0,$ modulo $5,$ $d_T = 0,$
	modulo $3,$ $t_T = 1.$

	\sssec{Theorem.}\label{sec-ttruncdodeca0}
	{\em For $q = 3^2,$ a primitive polynomial is}
	\enumb
	\item$I^3 - I - \epsilon ,$\\
	{\em with $\epsilon = 1 + \alpha,$ an 8-th root of unity and
	$\alpha^2 = -1.$\\
	The powers of $\epsilon$ are 1, $1 + \alpha,$ $-\alpha,$
	$1 - \alpha,$ $-1,$ $-1 - \alpha,$ $\alpha,$ $-1 + \alpha.$

	The corresponding selector is}
	\item$\{0, 1, 3, 9, 27, 49, 56, 61, 77, 81\}.$

	{\em The corresponding selector function is}
	\item
\footnotesize{
	\begin{verbatim}
	0  1  2  3  4  5  6  7  8  9 10 11 12 13 14 15 16 17 18 19 20 21 22 23
	   0  1  0 77 56  3 49  1  0 81 81 49 81 77 77 61 77  9 81 61 56 27 77
	t  e  e  e  d  e  e  v  v  v  v  e  u  u  e  v  e  u  v  u  u  v  e  u
	0  1  2 56 57 24 21 54 48  9 10 66 65 19 30 16 15 39 18 13 23  6 32 20

	  24 25 26 27 28 29 30 31 32 33 34 35 36 37 38 39 40 41 42 43 44 45 46
	   3 56  1  0 49 27 61 61 49 61 27 56 56 81 56 61  9 77 49 49 56 49  3
	   v  u  u  v  t  u  v  v  v  v  v  v  u  u  u  u  u  v  v  u  v  e  t
	   5 89 80 49 50 69 14 63 22 41 45 88 58 87 74 17 85 33 64 73 53 34 84

	  47 48 49 50 51 52 53 54 55 56 57 58 59 60 61 62 63 64 65 66 67 68 69
	   9  1  0 27 49  9  3 27  1  0 61  3 81  1  0 56 77 27 27 81 27  9 49
	   t  v  v  t  u  u  v  v  d  v  u  u  t  u  v  u  e  v  u  v  u  v  u
	  47  8 27 28 67 83 44  7 60  3  4 36 59 55 77 78 31 42 12 11 51 72 29

	  70 71 72 73 74 75 76 77 78 79 80 81 82 83 84 85 86 87 88 89 90
	  77 81  9 27  3 77  1  0  3 61  1  0  9  9 56  9 61 81  3  3  1
	   t  v  e  u  u  d  t  e  d  u  d  v  v  u  t  u  t  u  e  d  e
	  86 71 68 43 38 79 76 61 62 75 26 81 82 52 46 40 70 37 35 25 90
	\end{verbatim}}
\normalsize
	\item$\{0, 46, 47, 50, 59, 70, 28, 76, 84, 86\}.$
	\enume

	{\em The letters refer to the type.  The last row gives the conjugate,
	for instance, 61 is the conjugate of 77.}

	Proof:  To retrieve the primitive polynomial associated with $S,$ the
	selector 1, because $3 \in S,$ $I^3 = \beta I + \gamma,$
	$\beta$  and $\gamma$  are chosen in such a way that $I^{56}$ has no
	term of second degree.  The computations
	are facilitated by preparing first a table giving $g(i) \ni$\\
	\hth$	1 + \epsilon^i = \epsilon^{g(i)},$ $0 \leq  i \leq  t,$\\
	and by use of the convention $\epsilon^{-1} = 0.$

	The conjugates are obtained when $\alpha$  is replaced by $-\alpha$.

	\sssec{Heuristics.}
	The truncated dodecahedron has 182 faces, vertices and
	edges.  using symmetry with respect to the center we expect that a
	model can be found for the projective geometry of order $3^2,$ with 91
	points and with 10 points on each line .  We will solve simultaneously
	the following problems, discover appropriate incidence properties,
	associate to the vertices, integers from 0 to 90 to take advantage of
	the selector and determine a fundamental projectivity on a line to
	prepare for a representation of finite Euclidean geometry.  I will
	describe here some of the steps which have led me to the solution given
	in \ref{sec-dtruncdodeca} to \ref{sec-ttruncdodeca2}.

	The auto-correlates should be the points of a conic $\gamma$.  I will
	choose
	this conic as a circle in the corresponding Euclidean plane.  The
	intersection of the lines $0 \times 70 = 77^*$ and of
	$46 \times 28 = 72^*,$ which is 75, is chosen as the center of the
	circle.  The points on the polar $75^*$ 
	are 2, 6, 16, 17, 19, 25, 43, 65, 72, 77.

	To obtain a fundamental projectivity, we want to choose 2 points,
	$A,$ $B,$ on the circle and project from them any point $X$ on the
	circle onto $75^*,$ 
	giving $X_A$ and $X_B,$ $X_A$ corresponds to $X_B,$  we want to choose
	$A$ and $B$ such
	that the projectivity is of order 10.  A trial gave a projectivity of
	order 5, it was then easy to obtain one of order 10 using $A = 0$ and
	$B = 50.$  The computations start as follows:
	$0 \times 0 = 0^* \times 75^* = 77 \times 0 = 0^*$ with 0 as the other
	point on $\gamma$.\\
	$50 \times 0 = 27^* \times 75^* = 65 \times 0 = 27^*$ with 50 as the
	other point on $\gamma$.\\
	$50 \times 50 = 50^* \times 75^* = 6 \times 0 = 3^*$ with 46 as the
	other point on $\gamma$.\\
	$50 \times 46 = 31^* \times 75^* = 25 \times 0 = 56^*$ with 84 as the
	other point on $\gamma$.\\
	 $\ldots$  .  Hence the projectivity \ref{sec-ttruncdodeca2} and the
	equidistant points
	0, 50, 46, 84, 47, 70, 86, 28, 76, 59 on $\gamma$ .

	With one $d$-face chosen as $75^*,$ the 10 $t$-faces are subdivided
	into 2 sets, those adjacent to the d-face and those which are not.  The
	vertices of the pentagonal points 0, 46, 47, 86, 76 are chosen for the
	successive triangles adjacent to the $d$-face.  The diametrically
	opposite point, e. g. 70 of 0 is chosen for the triangle not adjacent
	to the $d$-face but adjacent to the triangle 0.

	Because $0 \times 46 = 3,$ $46 \times 47 = 45,$  $\ldots$ ,
	$0 \times 70 = 77,$ $50 \times 86 = 6$, I chose the pentagonal side 3
	for the $e$-point between the $t$-points 0 and 46,  $\ldots$, the
	diameter 77 for the $e$-point between the $t$-points 0 and 70, 
	$\ldots$.  Because $5 | 75_Y,$ we chose these 10 $e$-points as
	incident to $75^*.$	

	These consideration suggest Definition \ref{sec-dtruncdodeca} and
	Theorems \ref{sec-ttruncdodeca1} and \ref{sec-ttruncdodeca2}.

	\sssec{Definition.}\label{sec-dtruncdodeca}
	The points in the truncated dodecahedron model consist of
	\enumb
	\item The 10 {\em triangular face-points} $T.$
	\item The 6 {\em decagonal face-points} $D.$
	\item The 30 {\em vertex-points} $V.$
	\item The 30 {\em triangular-decagonal edge-points} $U.$
	\item The 15 {\em decagonal-decagonal edge-points} $E.$
	\enume

	The {\em lines} in the truncated dodecahedron model consist of
	\enumb
	\item The $10$ {\em triangular face-lines} $t$.  Each is incident to
	itself as a point, to the 3 adjacent $<12.12>$ edge-points $E,$ and to
	the 6 vertex-points $V$ which are the vertices of the triangle adjacent
	to the 3 edge-points which are not themselves adjacent to these
	points.\\
	For instance, for $t = 0,$ the incident points are\\
	$0(T), 1(E),77(E), 3(E),56(V),61(V), 27(V),81(V), 9(V),49(V):$
\newpage
\footnotesize{
	\begin{verbatim}
	                  o 61                       27 o

	                       .   1           77  .
	              56 o         o  .       .  o        o 81
	                                 o _0
	             .                   .                    .

	        .                        o 3                       .     .

	             .                   .                   .

	                 .           o       o           .
	                       .     49      9     .
	\end{verbatim}}
\normalsize
	\item The $6$ {\em decagonal face-lines} $d.$  Each is incident to
	its 5 $<3.12>$ edges $U$, and the 5 $<12.12>$ edges $E$ adjacent to its
	5 adjacent triangles.\\
	For instance 75($d$) is incident to\\
	$25(U), 43(U), 65(U), 17(U), 19(U) and 2(E), 6(E), 16(E), 77(E), 72(E):$
\footnotesize{
	\begin{verbatim}

	                               .
	                               o72
	                               .

	                            .  o  .
	             .                 19                .
	               o .     .               .     . o
	             77       o 17           25 o       2
	                     .         _75       .

	                       .               .
	                         o 65     43 o
	                            .     .
	                      .                 .
	                  16 o                   o 6
	                    .                     .
	\end{verbatim}}
\normalsize
	\item The $30$ {\em vertex-lines} $v.$  Each is incident to
	the $<12.12>$ edge $E_0$ adjacent to it and to the vertex at the other
	end of it, to the 2 triangular points $T_1$ and $T_2$ adjacent to the
	other edges $E_1$
	and $E_2$, to the 2 $<12.12>$ edges $E_3$ and $E_4$ opposite to $E_1$
	or $E_2$ belonging to the same decagon as $v$, to the vertices adjacent
	to $E_3$ or $E_4$ closest to $v$, to the $<3.12>$ edges $U$ belonging
	to the same decagon as $T_1$ or $T_2$ and the triangle opposite $E_0$.\\
	For instance, 9($v$) is adjacent to\\
	$72(E), 68(V), 0(T),47(T), 85(E),40(E), 82(V),40(V), 83(U), 52(U):$
\newpage
\footnotesize{
	\begin{verbatim}
	                  .                             .

	                       .                   .
	                 .  o 47      .       .     0 o   .

	             .                   ._9                  .
	        52 o                                            o 83
	        .                        o 72                      .     .

	             .                   o 68                .
	                    40                        85
	                 .  o        .       .        o  .
	                       o 18             82 o

	                  .                             .
	\end{verbatim}}
\normalsize

	\item The $30$ {\em triangular-decagonal edge-lines} $u.$  Each is
	incident to the adjacent decagonal point $D_0$, to the $<12.12>$ edge
	$E_0$ adjacent to the triangle adjacent to $u$ and to the $<3.12>$ edges
	$U$ adjacent to $E_0$, to the vertices in the same decagons $D_1$ and 
	$D_2$ as $E_0$ opposite the vertex adjacent to $E_0$ and the same
	triangle as $u$, and to the $<3.12>$ edges $U_1$ and $U_2$ adjacent to
	the triangle adjacent to $D_0$ and $D_1$ or $D_2$ not adjacent to these
	decagons and to the vetices adjacent to $D_0$ and the $<12.12>$ edges
	of $D_0$ adjacent to $U_1$ or $U_2$.\\
	For instance, 17($t$) is adjacent to\\
	$75(D), 77(E), 83(T),39(T), 64(V),10(V), 60(U),74(U), 44(V),32(V):$
\footnotesize{
	\begin{verbatim}
	                             .       .
	                    32 o                   o 44

	                  .                             .

	                                 o 75
	                  .                             .

	               74 o    .         _17       .    o 60
	                 .           .   .   .           .

	             .                   .                    .

	        .                        o 77                      .     .

	          10 o                   .                   o 64
	                            39 o   o 83
	                 .           .       .           .
	                       .                   .
	\end{verbatim}}
\normalsize
	\item The $15$ {\em decagonal-decagonal edge-lines} $e.$  Each is
	incident to the 2 decagonal points, the 2 triangular points,
	the 2 $<3.12>$ edges, the 2 vertices, the $<12.12>$ edge, whose
	center in the the equatorial plane through e and the $<12.12>$ edge
	perpendicular to that plane.\\
	For instance, 1($e$) is incident to\\
	$55(D),80(D), 0(T),76(T), 26(U),60(U), 48(V),8(V),90(E), 2(E):$
\newpage
\footnotesize{
	\begin{verbatim}

	                             .   o   .
	                       .         90        .

	                  .                             .

	                                 o 55
	                  .                             .

	                       .         26        .
	                 .           .   o   .           .
	                                 o 0
	             .                   o 48                .       2
	        .                     _1 .                        .  o  .
	             .                   o 8                 .
	                                 o 76
	                 .           .   o   .           .
	                       .         60        .

	                  .                             .

	                                 o 80
	                  .                             .

	                       .                   .

	                             .  (o)  .
	                                (90)
        \end{verbatim}}
        \enume

	\sssec{Theorem.}\label{sec-ttruncdodeca1}
	{\em The truncated docecahedron model satisfies the axioms}
	\ref{sec-ax1} {\em for q = $3^2.$}

	\sssec{Figure.}
	The corresponding drawing for the Projective plane over $3^2$ is given
	page \ldots.

	\sssec{Theorem.}\label{sec-ttruncdodeca2}
	{\em A fundamental projectivity on line $75^*$ is\\
        \hth$(77,65,6,25,72,17,16,43,2,19).$\\
        The elements are alternately of type $v$ and $u$.}

	\sssec{Exercise.}\label{sec-eseltrdodeca0}
	Given the selector function $f$ of \ref{sec-ttruncdodeca0} and the 6
	dodecagonal faces, 4, 55, 75, 78, 80, 89, reconstruct the preceding
	figure using the following rules, which are first examplified,
	\enumb
	\item $4(D) \times 89(D) = 5*(e),$ 5 is the decagonal-decagonal
	edge line which is in the equatorial plane through the center of
	the decagons 4 and 89.
	\item $1(E) \times 80(D) = 60*(u),$ 1 must be adjacent to a triangular
	face 76(T) adjacent to the decagon 80, 60 is then the triangular
	decagonal edge line adjacent to 76 and 80. 
	\item $1(E) \times 3(E) = 0*(t),$ 1 and 3 must be adjacent to the same
	triangular face 0(T), 0 is that face.
	\item $0(T) \times 5(E) = 56*(v),$ 0 must be adjacent to a
	decagonal-decagonal edge line 1(D) which must be adjacent to a
	tringular face 76(T) adjacent to the edge 5, 56 is the vertex
	adjacent to the latter 2.
	\enume
	The integers of the second member follow from the selector function
	for instance $5 = f(89-4)-4$.\\
	The above rules are clearly redundant.\\
	Determine alternate rules, for instance the rule corresponding to
	2 triangular faces or 2 vertices adjacent to the same
	decagonal-decagonal edge.
	Slightly more ambitious is to dermine all the possible rules.

	\sssec{Theorem.}
	{\em There are several configurations which represent a projective
	plane of order 3.  The quadrangle consists of 4 triangular face-points,
	the diagonal points, of 3 decagonal-decagonal edge-points, the
	quadrilateral, of 6 vertex-points.  All the other points on the
	truncated dodecahedron represent complex points, 6 on each of the 13
	lines.}

	The first example is associated with the primitive polynomial 
	\ref{sec-ttruncdodeca0}.0.
\footnotesize{
	\begin{verbatim}
	                        .                .
	                         o              o
	                         47  .   o   . 59
	                       .         20        .

	                  .                             .


	                  .                             .

	                        o 9              18 o
	                 .           .       .           .

	             .                   o 10                .       1
	        .  o 0                   o 90               76 o  .  o  .
	             .                   o 82                .

	                 .           .       .           .
	                       o 81             71 o

	                  .                             .
	\end{verbatim}}

\normalsize
	The conjugates are given in table 
	\ref{sec-ttruncdodeca0}.2.

	A second example is as follows
\footnotesize{
	\begin{verbatim}

	                            82       10
	                             o   o   o
	                       .         90        .

	                  .                             .


	                  .                             .

	                86 o    .                   .  o 84
	                 .           .       .           .

	             .                   o 8                 .       2
	        .                        o 1                      o  o  o
	             .                   o 48                .    7     54

	                 .           .       .           .
	                       .                   .
	                   o 46                        o 70
	                  .                             .

	                                 o 80
	                  .                             .

	                       .                   .

	                             .  (o)  .
	                                (90)
	\end{verbatim}}

\normalsize
	on the $t$-line $70^8,$ the conjugates are $21(V)$ and $11(E),$ $22(E)$
	and $77(E),$ $24 (V)$ and $30 (V)$.\\
	on the $e$-line $1^*,$ the conjugates are $0(T)$ and $76(T),$ $26(U)$
	and $60(U),$ $55(D)$ and $80(D).$\\
	on the $v$-line $8^*,$ the conjugates are $41(V)$ and $83(U),$ $53(V)$
	and $73(U),$ $19(U)$ and $69(U).$

	Proof:
	For the conjugates we use the Pascal construction to determine the
	6-th point on the line on a conic through 4 real points and 1 complex
	point.

	\sssec{Exercise.}\label{sec-etruncdodeca1}
	For $q = 2^2,$
	\enumb
	\item determine the primitive polynomial giving the
	    selector 0, 1, 4, 14, 16.
	\item Determine the correspondance between the selector notation
		and the homogeneous coordinates for points and lines. 
		Note that these are not the same.
	\item The correspondance $i$ to $i^*$ is a polarity whose fixed
	points are on a line.  Determine the matrix representation, the polar
	of $(X,Y,Z)$ and the equation satisfied by the fixed points.
	\item Determine the fundamental projectivity on the line $14^*$
	 using a point conic which has no points on $14^*.$
	\item Illustrate Pascal's Theorem.
	\enume

	\sssec{Exercise.}\label{sec-etruncdodeca2}
	\enumb
	\item Explore the usefulness of the truncated cuboctahedron less the
	    hexagonal faces and the $<4.8>$ edges as a model for
	    the projective geometry of order 7.
	\item Show that the 14-gonal antiprism can be used as a model for
	    the projective geometry of order 7.  More generally,
	\item Show that the $n$-gonal antiprism  can be used as a model for
	    the projective geometry of order $q = p^k$ when $p \equiv -1
	\pmod{4}$, with $n = \frac{q^2+q}{4}.$
	\item Show that the $n$-gonal antiprism  can be used as a model for
	    the projective geometry of order $q = p^k$ when $q \equiv 1
	\pmod{12}$, with $n = \frac{q^2+q}{2}.$  Finally,
	\item Show that the $n$-gonal prism  can be used as a model for
	    the projective geometry of order $q = p^k$ when $q \equiv -1
	\pmod{3}$, with $n = \frac{q^2+q}{3}.$
	\item is there a general theory when using prisms or antiprisms?
	\enume

	\sssec{Exercise.}\label{sec-etruncdodeca3}
	For $q = 2^3.$
	\enumb
	\item to 4. Answer question similar to those of \ref{sec-etruncdodeca1}
	\setcounter{enumi}{4}
	\item Show that the $18$-gonal antiprism can be used as a model for
	    the projective geometry of order $2^3.$ More generally,
	\item Show that the $n$-gonal antiprism  can be used as a model for
	    the projective geometry of order $q = 2^k$, with
	    $n = \frac{q^2+q}{4}.$
	\enume


\setcounter{section}{3}
\setcounter{subsection}{89}
	\sssec{Answer to}
	\vspace{-18pt}\hspace{94pt}{\bf \ref{sec-eseltetra}.}
	\enumb
	\item  For $q = 2,$ the primitive polynomial giving the selector
	0, 1, 3, is\\
	\hth$    I^3 + I + 1.$\\
	    The auto-correlates are 0  11  2  7  8.\\
	    The selector function is\\
	$\begin{array}{cclllllllllllllllllllll}
	i&\vline&0&1&2&3&4&5&6&7&8&9&10&11\\
	f(i)&\vline&&0&14&1&0&16&16&14&14&16\\
	type&\vline&F_0&V_0&F_4&V_2&T_0&T_2&V_1&F_3&F_1&T_3&E_2&F_2\\
	\cline{2-14}
	i&\vline&12&13&14&15&16&17&18&19&20\\
	f(i)&\vline&4&1&0&1&0&4&4&16&1\\
	type&\vline&T_4&E_1&P&V_3&E_0&T_1&E_3&V_4&E_4
	\end{array}$
	\item  The correspondence between the selector notation and the
	    homogeneous coordinates for points and lines is\\
	$\begin{array}{ccll}
	\hth&	i&	I^i&	     i^* \\
	&	0&	1   &    6^*: 1,2,4,\\
	&	1&	I&	1^*:	0,2,6,\\
	&	2&	I^2&	     0^*: 0,1,3,\\
	&	3&	I+1&	5^*:	2,3,5,\\
	&	4&	I^2+I&	   3^*: 0,4,5,\\
	&	5&	I^2+I+1&	 4^*: 3,4,6,\\
	&	6&	I^2+1&	   2^*: 1,5,6.
	\end{array}$
	\item  The matrix representation is\\
	\hth$    M = \mat{1}{0}{1}{0}{1}{0}{1}{0}{0},
		 M^{-1} = \mat{0}{0}{1}{0}{1}{0}{1}{0}{1}.$
	    and the equation satisfied by the fixed points is $(X_0+X_1)^2$ = 0.
	\item  The degenerate conic through 0, 1, 2 and 5 with tangent $5^*$
	at 5, is represented by the matrix\\
	\hth$N = \mat{0}{1}{1}{1}{0}{0}{1}{0}{0}.$\\
	    The polar of 0 is $0^*,$ of 1 is $0^*,$ of 2 is $5^*,$ of 4 is
		$4^*,$ of 5 is $5^*$ of 6 is $6^*$ and of 3 is undefined.
	    The equation in homogeneous coordinates is $X_0(X_1+X_2) = 0.$
	\item  A circle with center 14 can be constructed as follows.  I first
	    observe that a direction must be orthogonal to itself.  Indeed,
	    if 0 is a direction, the others form an angle 1,2,3,4 $\umod{5},$
	    we cannot play favorites and must choose 0.  If $A_0 = 1,$
	    $C \times A_0$ and therefore the tangent has direction 0,
	    $A_0 \times A_{i+1}$ has
	    direction $i \umod{5}$ or are the points 0, 7, 8, 2, 11.

	    It is natural to choose the pentagonal face-point as 14, and the
	    edge-points on the pentagon as 0, 8, 11, 7, 2.  The points on the
	    circle 1, 6, 3, 15, 19 are chosen as the vertex-points opposite
	    the corresponding edge-point, 1 opposite 0,  6 opposite 8, \ldots.
	    This gives the types, with subscripts indicated in 0. and the
	    definition:

	    The points are represented on the 5-anti-prism as follows.
	    The pentagonal face-point, P,
	    the 5 triangular face-points, $T_i,$ 
	    the 5 vertex-points, $V_i,$ 
	    the 5 triangular-triangular edge-points, $E_i,$ 
	    the 5 pentagonal-triangular edge-points $F_i.$ 

	    The lines are represented on the 5-anti-prism as follows.
	    The pentagonal face-line, f, which is incident to $F_i,$ 
	    the 5 triangular face-lines, $t_i,$ which are incident to $F_i,$ 
		$F_i,$ $T_{i+1},$ $T_{i-1},$ $E_{i+2},$ $E_{i-2}$.  If $f$ is
	the pentagonal edge of
		$t_i$	and $V,$ $V'$ are on $f$, $F_i$ is on it, $T_{i+1}$
	$(T_{i-1})$ share $V$ $(V'),$
		$E_{i+2}$ $(E_{i-2})$ are on an edge through $V$ $(V')$ not on
		$t_i$ \\
	    the 5 vertex-lines, $v_i,$ which are incident to\\
		$F_i,$ $V_{i+2},$ $V_{i-2},$ $E_{i+1},$ $E_{i-1}.$  If $t$ is
	the face with $v_i$ on its
		pentagonal edge these are all the vertices, and edge-points on
		it distinct from $v_i.$ \\
	    the 5 triangular-triangular edge-lines, $e_i,$ which are incident to
		$F_i,$ $T_{i+2},$ $T_{i-2},$ $V_{i+1},$ $V_{i-1}.$  $V_{i+1}$
	and $V_{i-1}$ are on the same edge as $e_i,$ the line which joins the
	center $C$ of the antiprism
		to $E_{i}$ is parallel to the edge containing $F_i,$ $T_{i+2}$
	and $T_{i-2}$ are
		the triangular faces which are not adjacent to $E_i$ or $F_i.$ 

	    the 5 pentagonal-triangular edge-lines. $f_i,$ which are incident to
		$P,$ $T_i,$ $V_i,$ $E_i,$ $F_i.$  $T_i$ is adjacent to $f_i,$
	$V_i$ is opposite $f_i,$ 
	$E_i$	joined to the center of the antiprism is parallel to $T_i.$ 
	\enume

	\sssec{Answer to}
	\vspace{-18pt}\hspace{94pt}{\bf \ref{sec-eselcube}.}\\
	For $p = 3,$
	\enumb
	\item  The primitive polynomial giving the selector 0, 1, 3, 9 is
	    $I^3 - I - 1.$
	\item  The correspondence between the selector notation and the
	    homogeneous coordinates for points and lines is\\
	$\begin{array}{rllll}
		&i	&I^i	&     i^*\\
		&0	&1      &12^*:& 1,2,4,10,\\
		&1	&I	&1^*:&	0,2,8,12,\\
		&2	&I^2	&     0^*:& 0,1,3,9,\\
		&3	&I+1	&7^*:&	2,6,7,9,\\
		&4	&I^2+I	&   3^*:& 0,6,10,11,\\
		&5	&I^2+I+1	& 4^*:& 5,9,10,12,\\
		&6	&I^2+2I+1	&1^{0*}:& 3,4,6,12,\\
		&7	&I_2+I+2	& 6^*:& 3,7,8,10,\\
		&8	&I^2+1	&   2^*:& 1,7,11,12,\\
		&9	&I+2    &11^*:& 2,3,5,11\\
		&10	&I_2+2I	&  9^*:& 0,4,5^,7, \\
		&11	&I^2+2I+2	&5^*:& 4,8,9,11,\\
		&12	&I^2+2	&  8^*:& 1,5,6,8.
	\end{array}$

	\item  The matrix representation of the polarity $i$ to $i^*$ is\\
	\hth$M = \mat{1}{0}{1}{0}{1}{0}{1}{0}{0},$
		$M^{-1} = \mat{0}{0}{1}{0}{1}{0}{1}{0}{2}.$\\
	    The equation satisfied by the fixed points is
	$X_0^2 + X_1^2 + 2 X_2 X_0  = 0.$
	\item  The degenerate conic through 0, 1, 2 and 5 with tangent $4^*$ at
	$5$, is obtained by constructing the quadrangle-quadrilateral
	configuration starting with $P = 5$ and $Q_i = \{0,1,2\}.$  We obtain
	$q_i$ = $\{3^*,2^*,7^*\},$  which are the tangents at $Q_i.$  The
	matrix representation is\\
	\hth$N = \mat{0}{1}{1}{1}{0}{1}{1}{1}{0}$ with equation
	$X_1X_2 + X_2X_0 + X_0X_1 = 0.$\\
	    We can check that the polar of $10 = 3^* \times 4^*$ is
		$9^* = 0 \times 5.$
	\enume

	\sssec{Answer to}
	\vspace{-18pt}\hspace{94pt}{\bf \ref{sec-eseltruncdodeca1}.}
	\enumb
	\item  For $q = 2^2,$ the primitive polynomial giving the selector
	    0, 1, 4, 14, 16 is $I^3-I^2-I-\epsilon$ , with\\
	\hth$\epsilon^2+\epsilon +1 = 0.$
	\item  The correspondence between the selector notation and the
	    homogeneous coordinates are as follows, $i^*$ has the homogeneous
	    coordinates associated with $I^i.$\\
	$\begin{array}{clrr}
	\hth&     i&	I^i&		i*\\
	\cline{2-4}
	&     0& 1&	      20^*\\
	&     1& I&	      14^*\\
	&     2& I^2&	      0^*\\
	&     3& I^2+I+\epsilon&		 10^*\\
	&     4& I+\epsilon&	^2    19^*\\
	&     5& I^2+\epsilon&	^2I   4^*\\
	&     6& I^2+\epsilon&	^2I+1       18^*\\
	&     7& I^2+1&	   15^*\\
	&     8& I^2+\epsilon&		    3^*\\
	&     9& I^2+\epsilon^2I+\epsilon& 	5^*	\\
	&    10& I^2+\epsilon I+1&	 9^*\\
	&    11& I^2+\epsilon	^2&   13^*\\
	&    12& I^2+\epsilon I+\epsilon& 	11^*\\
	&    13& I^2+I+\epsilon^2&  6^*\\
	&    14& I+1&	     2^*\\
	&    15& I^2+I&	    1^*\\
	&    16& I+\epsilon&		     12^*\\
	&    17& I^2+\epsilon I&	   16^*\\
	&    18& I^2+\epsilon I+\epsilon^2& 17^*\\
	&    19& I^2+\epsilon^2I+\epsilon^2& 8^*\\
	&   20& I^2+I+1&	  7^* 
	\end{array}$

	    To obtain the last column, for row 9, $[1,\epsilon^2,\epsilon ] =
	(1,1,1) \times (1,\epsilon,0) = 20 \times 17 = 5*.$
	\item  The correspondence $i$ to $i^*$ is a polarity whose fixed points
	are on a line.  The matrix representation is obtained by using the image
	of 4 points.\\
	\hth	 0 = (0,0,1), $M(0) = 0^* = [1,0,0],$\\
	\hth	 1 = (0,1,0), $M(1) = 1^* = [1,1,0],$\\
	\hth	 2 = (1,0,0), $M(2) = 2^* = [0,1,1],$\\
	\hth	18 = $(1,\epsilon,\epsilon^2),$ $M(18) = 18^*
		= [1,\epsilon^2,1].$\\ 
	    The first 3 conditions give the polarity matrix as\\
	    The last condition gives $\beta \epsilon + \alpha\epsilon^2 = 1,$
	$\gamma + \beta \epsilon = \epsilon^2,$ $\gamma = 1.$  Hence 
	$\gamma  = 1,$ $\beta  = 1,$ $\alpha  = 1.$  Therefore\\
	\hth$M = \mat{0}{1}{1}{1}{1}{0}{1}{0}{0},$
		$M^{-1} = \mat{0}{0}{1}{0}{1}{1}{1}{1}{1}.$\\
	    Note that $M$ is real and could have been obtained from the reality
	    and non singularity conditions, giving directly
	$\alpha = \beta = \gamma  = 1.$\\
	    The polar of $(X_0,X_1,X_2)$ is $[X_1+X_2,X_0+X_1,X_0].$\\
	    The fixed points $(X_0,X_1,X_2)$ satisfy $X_1^2 = 0$ corresponding
	to $14^*.$ 
	\item  A point conic with no points on 14 is  1, 3, 4, 5,13,\\
	    the corresponding line conic is       15,19,10,16, 8.\\
	    Projecting from 1 and 3,     	   1, 3, 5,13, 4,\\ we get the
	    fundamental projectivity,		   8, 2,11, 0, 7 on $14^*.$ 

	\item To illustrate Pascal's Theorem, because there are only 5 points on
	    a conic, we need to use the degenerate case.  The conic through
	    0, 1, 2 and the conjugate points 9 and 18 is
	    The last condition gives $\beta \epsilon  + \alpha\epsilon^2 = 1,$
	$\gamma + \beta \epsilon = \epsilon^2,$ $\gamma = 1.$\\
	\hth    Hence $\gamma  = 1,$ $\beta  = 1,$ $\alpha  = 1.$  Therefore\\
	\hth$M = \mat{0}{1}{1}{1}{1}{0}{1}{0}{0},$
		$M^{-1} = \mat{0}{0}{1}{0}{1}{1}{1}{1}{1}.$\\
	    Note that $M$ is real and could have been obtained from the reality
	    and non singularity conditions, giving directly
	$\alpha = \beta = \gamma  = 1.$\\
	    The polar of $(X_0,X_1,X_2)$ is $[X_1+X_2,X_0+X_1,X_0].$\\
	    The fixed points $(X,X_1,X_2)$ satisfy $X_1^2 = 0$ corresponding to
	    $14^*.$ 
	\item A point conic with no points on 14 is  1, 3, 4, 5,13,\\
	\hth$\mat{0}{1}{1}{1}{0}{1}{1}{1}{0}$\\
	    The tangents at (0,0,1), (0,1,0), (1,0,0),
	$(1,\epsilon^2,\epsilon),$ $(1,\epsilon,\epsilon^2)$ 
	    are		    [1,1,0], [1,0,1], [0,1,1],
		$[1,\epsilon^2,\epsilon),$ $(1,\epsilon,\epsilon^2],$ 
	    or		     $1^*,$      $15^*,$     $2^*,$   $5^*,$    $17^*.$ 
	    On the other hand, using Pascal's Theorem, the tangent at 0 is
	    given by\\
	$  ( ( ( (0 \times 1) \times (9 \times 18) ) \times ( (18 \times 0)
	\times (1 \times 2) ) ) \times (2 \times 9) ) \times 0\\
	\hth= ( ( (0^* \times 7^*) \times (4^* \times 20^*) ) \times 12^*)
	 \times 0\\
	\hth= ( ( (14 \times 17) = 8^*) \times 12^* {\rm or} 13) \times 0 = 1^*.$ 
	\enume

	\sssec{Answer to}
	\vspace{-18pt}\hspace{94pt}{\bf \ref{sec-eseltruncdodeca2}.}\\
	For $q = 57,$
	choose the auto-correlates as point on a circle although 0 is on the
	circle draw as it is the center.  With the succession of points $X_i,$\\
	$\begin{array}{llrrrrrr}
	x_i = 0 \times X_i &36,& 1,&52,&43,& 3,&32,&13,\\
	X_i		   &   16,&35,&18,&50,&29,&26,&30,\\
	y_{i+1} = X_{i-1} \times X_{i+1}&      22,&42,& 8,&14,&10,&28,&44,\\
	y_{i+2} = X_{i-2} \times X_{i+2}&      34,& 2,&41,&17,&40,&20,&23,\\
	y_{i+3} = X_{i-3} \times X_{i+3}&       7,&31,& 6,&27,&54,&25,&39,\\
	y_{i+1} \times x_i&  21,&51,& 5,&46,&33,& 4,&45,\\
	y_{i+2} \times x_i&  24,&56,&48,&15,&49,&38,&47,\\
	y_{i+3} \times x_i&  53,&12,&37,& 9,&55,&11,&19.\\
	\end{array}$
	
	This gives all the points in the projective plane of order 7.
	We observe\\
	$\begin{array}{lllllll}
	16^*&21^*&24^*&53^*&22^*&34^*& 7^* \\
	36&36&36&36&36&36&36\\
	16&&&&35,30&18,26&50,29\\
	42,44&22& 8,28&&14,10\\
	41,20&&34&17,40&& 2,23\\
	27,54&31,39&& 7&&& 6,25\\
	&&46,33& 5, 4&21&&51,45\\
	&15,49&&56,47&48,38&24\\
	&37,11&12,19&&& 9,55&53\\
	\cline{1-7}
	35^*&51^*&56^*&12^*&42^*& 2^* 31^*\\
	 1& 1& 1& 1& 1& 1& 1\\
	35&&&&16,18&50,30&29,26\\
	22, 8&42&14,44&&10,28\\
	17,23&& 2&40,20&&34,41\\
	54,25& 7, 6&&31&&&27,39\\
	&&33, 4&46,45&51&&21, 5\\
	&49,38&&24,48&15,47&56\\
	& 9,19&53,37&&&55,11&12\\
	\cline{1-7}
	18^*& 5^* 48^*&37^*& 8^* 14^*& 6^*\\
	52&52&52&52&52&52&52\\
	18&&&&35,50&16,29&26,30\\
	42,14& 8&22,10&&28,44\\
	34,40&&41&20,23&& 2,17\\
	25,39&31,27&& 6&&& 7,54\\
	&& 4,45&21,33& 5&&51,46\\
	&38,47&&56,15&24,49&48\\
	&53,55&12, 9&&&11,19&37
	\end{array}$

	\sssec{Answer to}
	\vspace{-18pt}\hspace{94pt}{\bf \ref{sec-eseltruncdodeca3}.}\\
	For q = $2^3,$ \\
	$\begin{array}{lrrrrrrrrrcrrrrrrrrr}
	36:&		 0&37&38&40&44&52&18&27&68&\hti{2}&1^*& 3^*& 7^*&
	2^*& 4^*& 5^*\\
	36 \times  0 =  0^*:&   0& 1& 3& 7&15&31&36&54&63& & 0 & 0 & 0 & 
	1 & 3 &31\\
	36 \times 37 = 37^*:&  17&26&36&37&39&43&51&67&72& &72 &51 &67 & 
	72 &72 &26\\
	36 \times 38 = 38^*:&  16&25&35&36&38&42&50&66&71& &35 &71 &66 &
	71 &50 &71\\
	36 \times 40 = 40^*:&  14&23&33&34&36&40&48&64&69& &14 &33 &69 &
	34 &69 &69\\
	36 \times 44 = 44^*:&  10&19&29&30&32&36&44&60&65& &30 &60 &29 &
	29 &32 &10\\
	36 \times 52 = 52^*:&   2&11&21&22&24&28&36&52&57& & 2 &28 &24 &
	52 &11 & 2\\
	36 \times 18 = 18^*:&  13&18&36&45&55&56&58&62&70& &62 &70 &56 &
	13 &70 &58\\
	36 \times 27 = 27^*:&   4& 9&27&36&46&47&49&53&61& &53 & 4 &47 &
	61 &27 &49\\
	36 \times 68 = 68^*:&   5& 6& 8&12&20&36&41&59&68& & 6 &12 & 8 &
	 5 &59 &68\\
	\end{array}$

	Conic with no point on 36:  2, 4, 5, 6,13,28,31,46,63\\
	\hth	line conic:	   29,59,31, 9,18,43,28,35,64.\\
	Fundamental projectivity: from 2 and 5 on the conic, the points\\
	\hth	 2, 5, 6,31,13,28, 4,46,63 give the points on $36^*:$\\
	\hth	38, 0,68,27,52,37,40,18,44.

 empty 


\chapter{FINITE PRE INVOLUTIVE GEOMETRY}
%
	\section{An Overview of the Geometry of the Hexal Complete 5-Angles.}

	\setcounter{subsection}{-1}
	\ssec{Introduction.}

	In the geometry of Euclid, not every pair of lines have a point in 
	common, namely the parallel ones.  I call Euclidean Geometry, that
	geometry which consists in completing the plane of Euclid by the ideal 
	points and the lines of Euclid by the ideal line. To each set of
	parallel lines correspond
	its direction or point at infinity or ideal point. The line at
	infinity or ideal line is incident to all ideal points. Figures Pl and 
	St may help the reader to visualize. In Fig. Pl, projecting the line
	$b$ on the line  $c$ from the point $P$ establishes a one to one
	correspondance between the points on these lines, if we include the 
	ideal point $C_i$, on $c$, corresponding to $B_i$ and the ideal point
	$B_\infty$, on $b$, corresponding to $C_\infty$. Replacing lines 
	$b$ and $c$ by planes, perpendicular to the plane ${\Bbb P}$ of figure 
	establishes a one to one correspondance between a line through 
	$C_\infty$ perpendicular to ${\Bbb P}$ and the ideal line through
	$B_\infty$.\\
	This led to the concept of perspectivity, which I have schematized in
	Fig. St. In it, the shading, corresponds to the method used by Chinese
	artists to represent distances in paintings. The tiling
	corresponds to the method used by Western painters. Johannes Vermeer's
	use of perspective in his paintings was so accurate as to allow P. T. A.
	Swillens to reconstruct, from the size of a chair, in the painting, not 
	only the size of the rooms, but also to estimate the height of the
	artist.\\
	Affine geometry is obtained from Euclidean geometry by discarding
	the notions associated with congruences of figures, projective
	geometry is obtained by discarding the notion of parallelism,
	thereby making the properties of any point or line in the plane 
	indistinguishable from that of any other.\\
	I will describe at a later time, how I was lead to the discovery of
	finite Euclidean geometry and to the extension of many of the
	properties of Euclidean geometry.  While working out a proof for these
	results, it occured to me, that the results can be placed in the
	framework of finite projective geometry.  I will, as I proceed, make
	the connection with the results in classical Euclidean geometry.
	The results can be considered as proceeding from an, apparently new,
	configuration consisting of 14 points and 13 lines.
	This configuration is defined starting from an ordered complete
	5-angle, $A_0,$ $A_1,$ $A_2,$ $M$ and $\ov{M},$ in which the
	first 3 points can be rotated and the last 2 points interchanged.
	In other words the configuration is the same if we replace $A_0,$ $A_1$ 
	and $A_2$ by $A_1,$ $A_2$ and $A_0$ and independently $M$ by
	$\ov{M}$ and $\ov{M}$ by $M.$
	In involutive geometry, (the Euclidean geometry without measure of
	angles and distances), we define altitudes and their intersection, the
	orthocenter, we define medians and their intersection, the barycenter.
	In the generalization to projective geometry, the orthocenter and the
	barycenter become two arbitrary points, whose role is interchangeable.
	The proofs are constructive, and the only construction required are
	those of lines through 2 given points and of points at the intersection
	of two given lines, but these constructions must be valid for all $p$.
	They do not involve the construction of an
	arbitrary line through a given point, as required to obtain, for
	instance, an arbitrary point on a conic, by the construction of Pascal
	or of MacLaurin.\\
	No special relations will be assumed here between the points obtained
	during the construction.  The special relation $\ov{M}$ on the
	polar of $M$ with respect to the triangle $A_0,$ $A_1,$ $A_2$ will be
	studied in Chapter IV,\\
	in the section on Cartesian geometry and the special case where $M$ and
	$\ov{M}$ are respectively on the polar of $\ov{M}$ and $M$ with respect
	to the triangle will also be discussed elsewhere.\\
	The beginning of a synthetic proof is given in section 4.3.  Synthetic
	proofs are highly desirable and are from my point of view more elegant,
	but require much more time to develop.\\
	The constructions and statements are given in a compact found using a
	notation which will now be explained.

	\ssec{Notation and application to the special configuration of
	Desargues and to the pole and polar of with respect to a triangle.}

	\setcounter{subsubsection}{-1}
	\sssec{Introduction.}
	In the preceding Chapter, I have introduced a notation for points,
	lines,
	incidence and statements.  Additional notation is given here for
	conics, for points on conics and tangents to conics and a notation
	which allows to describe at once 3 points or 6 points associated to a
	triangle.

	\sssec{Notation.}\label{sec-nalgebraic}
	The identifier for a point conic will be a lower case
	Greek letter or an identifier starting with a lower case preceded by
	a backward quote `` ` ''.  The identifier for a line conic will be an
	upper case Greek letter or an identifier starting with an upper case
	preceded by a backward quote `` ` ''.

	The subscript $i$, will have the values $0,1$ and $2$.
	Hence $A_i$ denotes 3 points $A_0,$ $A_1$ and $A_2.$ 

	If subscripts involve the letter $i$ and addition, the addition is done
	modulo 3, for instance,\\
	\hth$		a_i := A_{i+1} \times A_{i+2}$\\
	is equivalent to\\
	\hth$a_0 := A_1 \times A_2,$ $a_1 := A_2 \times A_0,$
	$a_2 := A_0 \times A_1.$\\
	It represents the construction of the sides $a_0,$ $a_1$ and $a_2$ 
	of a triangle with vertices $A_0,$ $A_1$ and $A_2.$ 

	To indicate that a conic $\gamma$  is constructed as that conic which
	passes through the 5 distinct points $P_0,$ $P_1,$ $P_2,$ $P_3,$ and
	$P_4,$ we write\\
	\hth$	\gamma  := conic(P_0, P_1, P_2, P_3, P_4).$\\
	To indicate that a conic $\gamma_1$ is constructed as that conic whose
	tangent at $P_0$ is $a_0$ and at $P_2$ is $a_2$ and passes also by
	$P_4,$ we write\\
	\hth$	\gamma_1 := conic((P_0,a_0),(P_2,a_2),P_4).$ or
		$\gamma_1 := conic(P_0,a_0,P_2,a_2,P_4).$

	When 3 lines $x_i$ are concurrent, the intersection $X$ can be obtained
	using any of the three pairs.  I have chosen, arbitrarily,\\
	\hth$	X := x_1 \times x_2 (*),$\\
	as a reminder that 2 other definitions of $X$ could have been chosen.
	In the special case, $x_1$ = $x_2,$ the other choice\\
	\hth$	X := x_0 \times x_1,$\\
	will be used. `` (*) '' denotes therefore not only a Definition but
	also a Theorem or Conclusion.  A similar notation will be used for 
	conics.\\
	\hth$X \cdot \gamma  = 0$ and $X_i \cdot \gamma = 0,$ are the notations
	corresponding to the point $X$ is on the conic $\gamma$  and the
	triple $X_0,$ $X_1,$ $X_2$ is on the conic $\gamma$.\\
	\hth $P = Pole(p,\alpha$ ), is the notation for $P$ is the pole of $p$
	with respect to the conic $\alpha$ .\\
	\hth$\gamma\:is\:a\:circle = 0$ or $\gamma\:is\:a\:cocircle$ is either
	an hypothesis, to indicate a prefered conic from which all other
	circles are defined or a Conclusion,\\
	\hth$X = Center(\gamma )$ and $\ov{X} = Cocenter(\gamma )$ is an
	abbreviation for
	$X$ is the center of the conic $\gamma$  (not necessarily a circle)
	and $\ov{X}$ is
	the cocenter of the conic $\gamma$, in other words $X,$
	$(\ov{X})$ is the polar of
	$m$ $(\ov{m})$ with respect to $\gamma$.  See section 4.3.3.

	\sssec{Example.}\label{sec-epoletr}
	With this notation, the special configuration of Desargues
	of 0.4.6. can be defined by\\
	\hth$	a_i := A_{i+1} \times A_{i-1},$	$r_i := P \times A_i,$\\
	\hth$	P_i := a_i \times r_i,$	    $q_i := P_{i+1} \times P_{i-1},$\\
	\hth$	R_i := a_i \times q_i,$	    $p_i := A_i \times R_i,$\\
	\hth$	Q_i := p_{i+1} \times p_{i-1},$	$p := R_1 \times R_2 (*),$\\
	and the conclusions of the special Desargues Theorem are implied by
	the last Definition-Conclusion and by the Conclusion,\\
	\hth$	Q_i \cdot r_i = 0.$\\
	Let $P = (p_0,p_1,p_2),$ and $A_0 = (1,0,0),$ $A_1 = (0,1,0),$
	$A_2 = (0,0,1),$
	then\\
	\hth$	a_0 = [1,0,0],$	 $r_0 = [0,p_2,-p_1],$\\
	\hth$	P_0 = (0,p_1,p_2),$       $q_0 = [-p_1p_2,p_2p_0,p_0p_1],$\\
	\hth$	R_0 = (0,p_1,-p_2),$      $p_0 = [0,p_2,p_1],$\\
	\hth$	Q_0 = (-p_0,p_1,p_2),$     $p = [p_1p_2,p_2p_0,p_0p_1].$

	\sssec{Example.}\label{sec-eqqc3c}
	For $p = 3,$ prove that if $A_0 = (1,0,0),$ $A_1 = (0,1,0),$
	$A_2 = (0,0,1)$ and $P = (1,1,1)$ then the other elements of the
	quadrangle quadrilateral configuration II.\ref{sec-dsd} are\\
	\hth$P_0 = (0,1,1),$ $Q_0 = (-1,1,1),$ $R_0 = (0,1,-1),$ \ldots, and\\
	\hth$a_0 = [1,0,0],$ $p = [1,1,1],$\\
	\hth$p_0 = [0,1,1],$ $q_0 = [-1,1,1],$
	$r_0 = [0,1,-1],$ \ldots and that\\
	the conic of II.\ref{sec-tqqc3} is\\
	\hth$X_0^2 + X_1^2 + X_2^2 = 0.$

	\sssec{Theorem.}
	{\em With the above notation, the polar $p$ can be obtained
	algebraically from the pole $P$ or the pole $P$ from the polar $p$ using
	the first or the second formula:}
	\enumb
	\item$      pA_i = P * A_i,$	       $Pa_i = p * a_i,$
	\item$      pA_i = (P \cdot a_i) A_i - (A_i \cdot a_i) P,$
		$Pa_i = (p \cdot A_i) a_i - (a_i \cdot A_i) p,$
	\item$ p_i = (P \cdot a_{i+1})(P \cdot a_{i-1}) A_{i+1} * A_{i-1}
		+ (P \cdot a_{i+1})(A_{i-1} \cdot a_{i-1}) P * A_{i+1}\\
	\hth	- (P \cdot a_{i-1})(A_{i+1} \cdot a_{i+1}) P * A_{i-1},\\
	\hti{4}P_i = (p \cdot A_{i+1})(p \cdot A_{i-1}) a_{i+1} * a_{i-1}
		+ (p \cdot A_{i+1})(a_{i-1} \cdot A_{i-1}) p * a_{i+1}\\
	\hth	- (p \cdot A_{i-1})(a_{i+1} \cdot A_{i+1}) p * a_{i-1}.$
	\item$  Pa_i = (P \cdot a_{i+1}) A_{i+1} - (P \cdot a_{i-1}) A_{i-1},\\
	\hti{4}PA_i = (p \cdot A_{i+1}) a_{i+1} - (p \cdot A_{i-1}) a_{i-1}.$
	\item$      p = \frac{1}{P \cdot a_0}a_0 + \frac{1}{P \cdot a_1}a_1
			+ \frac{1}{P \cdot a_2}a_2,$
		$P = \frac{1}{p \cdot A_0}A_0 + \frac{1}{p \cdot A_1}A_1
			+ \frac{1}{p \cdot A_2}A_2.$
	\enume

	Proof:  Only the first part of 2 to 4 needs to be proven, because of
	duality.
	To obtain 2, we use $P * P = 0$ and $A_{i+1} * P = - P * A_{i+1}.$
	To obtain 3, we recall that $a_i = A_{i+1} * A_{i-1},$
	we use $A_i * a_j = 0$ when $i\neq j$ and $A_i * a_i = (A_0 * A_1)
	\cdot A_2 = t,$
	then divide by $t$ and by $P \cdot a_i \neq 0$.
	To obtain 4, we use $p = R_{i+1} * R_{i-1}.$  We divide by
	$(P \cdot a_i)(P \cdot a_{i+1})(P \cdot a_{i-1}\neq 0,$ and obtain\\
	$p = A_{i+1} * \frac{1}{P \cdot a_i}A_{i-1} + A_{i-1} *
	\frac{1}{P \cdot a_{i+1}}A_i
		+ A_i * \frac{1}{P \cdot a_{i-1}}A_{i+1}$, or\\
	$p = \frac{1}{P \cdot a_i}a_i + \frac{1}{P \cdot a_{i+1}}a_{i+1}
	+ \frac{1}{P \cdot a_{i-1}}a_{i-1}$.

	\sssec{Example.}
	For $p = 13,$ $A_i = ( 36(1,1,9),27(1,1,0),151(1,10,7) ),$
	$P = (68(1,4,2)),$\\
	$a_i = [ 175(1,12,5), 150(1,10,6), 170(1,12,0) ],$
	$r_i = [77,138,31],$\\
	$S_i = (143,63,33),$ $p_i = [108,46,37],$ $R_i = (48,16,32),$\\
	$p = \frac{1}{7}a_0 + \frac{1}{1}a_1 + \frac{1}{10}a_2
	= 2a_0 + a_1 + 4a_2 = [124],$ $p_i = [140,176,106],$\\
	$P_i = (51,132,84),$ $P = \frac{1}{11}A_0 + \frac{1}{9}A_1 +
	\frac{1}{6}A_2 = 6A_0 + 3A_1 + 11A_2 = (68).$

	\sssec{Definition.}
	An {\em hexal complete 5-angle configuration}, is a configuration
	which starts with an ordered set of 5 points $A_0,$ $A_1,$ $A_2,$ $M$
	and $\ov{M}.$

	In the configuration obtained from it, if a point $X_0$ is constructed,
	5 other points are obtained.  $X_1$ is obtained by replacing $A_0,$ 
	$A_1,$	$A_2$ by $A_1,$ $A_2,$ $A_0;$ $X_2$ is obtained by replacing
	the same points by $A_2,$ $A_0,$ $A_1$ and the points $\ov{X}_i$
	are obtained by exchanging in the construction of $X_i,$ $M$  and
	$\ov{M}.$  The same holds for lines.  The first letter has a
	macron placed above it in the naming of the
	construction which exchanges $M$ and $\ov{M}.$
	\{ In the group of permutation on the 5 points of the complete 5-angle,
	the figure is invariant under the cyclic group generated by the
	permutation\\
	\hth$	( A_0 A_1 A_2  M  \ov{M} )$\\
	\hth$	( A_1 A_2 A_0  \ov{M} M  ) \}.$\\
	In special cases, several of these elements or all of the elements may
	coincide.

	\sssec{Comment.}
	We know from II.1.5.6. that a complete 5-angle requires $p \geq  5,$
	therefore, the definition and results that follow are non vacuous only
	if p $\geq$  5.
	We introduce here a terminology inspired from corresponding terms
	in Euclidean geometry.  In some instances, the correspondence will be
	made explicitly.  For instance,
	the line $m$ which will be constructed corresponds to the ideal line or
	line at infinity in Euclidean geometry, we will therefore call  $m$
	the ideal line.  In the symmetry which exchanges $M$ and
	$\ov{M}$, to $m$ corresponds $\ov{m}$, which will be called the coideal
	line. \{$\ov{m}$ corresponds to the orthic axis. \}
	The conic $\theta$  which will be constructed corresponds to the
	circumcircle and the conic $\gamma$  to the circle of Brianchon-Poncelet
	also called the nine-point circle.

	\sssec{Definition.}
	$\theta$ and any conic $\delta , (\ov{\delta} )$ such that there
	exists a radical axis $u$, $(\ov{u})$ with respect to $m$
	$(\ov{m})$ is called a {\em circle}
	({\em cocircle}) and $u$ is called the {\em radical} ({\em coradical})
	 {\em axis} of $\theta$  and
	$\delta, (\ov{\delta} )$.

	Algebraically, we have, for some integers $k_1,$ $k_2$ and $k_3,$\\
	\hth$	k_1 \delta  + k_2 \theta  = k_3 (m) \bigx (u),$\\
	where $\theta$, $m,$ $\delta$  and $u$ are expressed exactly as in the
	corresponding expressions P0.7, P1.19, 1.20, \ldots, below.

	A {\em triangle} consists of its vertices and its sides.  When we want
	to be
	specific we will use either or both, for instance the given triangle
	can be written as $\{A_i\}$ or $\{a_i\}$ or $\{A_i,a_i\}$.\\
	To each of the section of this Theorem corresponds a sequence of
	theorems in Euclidean geometry which will be given in the corresponding
	sections of Chapter IV.

	We will give separately the construction of the various points and
	lines of the hexal configuration ({\em zetetic part}) and the proof
	that the construction satisfies the given properties ({\em poristic
	part}).

	\ssec{An overview of theorems associated with equality of
	distances and angles.  The ideal line, the
	orthic line, the line of Euler, the circle of Brianchon-Poncelet, the
	circumcircle, the point of Lemoine.}

	\setcounter{subsubsection}{-1}
	\sssec{Introduction.}
	As the generalization proceeds, the 4 points on the line
	of Euler, become 10 points on its generalization.  The 9 (or 12) points
	on the circle of Brianchon-Poncelet (also called circle of Euler)
	become 20 points on the corresponding conic.  New results, which will
	be given in part III, are further consequences.\\
	The definitions are numbered starting with D, the conclusions are
	numbered stating with C, the proofs, which consist of the algebraic
	expressions of the various points, lines and conics, which can easily
	be checked and from which the conclusions can easily be verified, have a
	number corresponding to the definition, starting with P.\\
	The numbering in this overview is the same as the number in the
	complete theory, given in Chapter 5 and 6.

	\sssec{Theorem.}
	{\em If we derive a point $X$ and a line $x$ by a given construction
	from $A_i$, $M$ and $\ov{M}$, with the coordinates as given in
	G0.0 and G0.1, below, and the point $\ov{X}$ and line $\ov{x}$
	are obtain by the same construction interchange $M$ and $\ov{M}$,\\
	\hth$X = (f_0(m_0,m_1,m_2),f_1(m_0,m_1,m_2),f_2(m_0,m_1,m_2)),\\
	\hth x = [g_0(m_0,m_1,m_2),g_1(m_0,m_1,m_2),g_2(m_0,m_1,m_2)],\\
	\implies\\
	\hth \ov{X} = (m_0f_0(m_0^{-1},m_1^{-1},m_2^{-1}),
		m_1f_1(m_0^{-1},m_1^{-1},m_2^{-1}),
		m_2f_2(m_0^{-1},m_1^{-1},m_2^{-1})),\\
	\hth \ov{x} = [m^{-1}_0g_0(m_0^{-1},m_1^{-1},m_2^{-1}),
		m^{-1}_1g_1(m_0^{-1},m_1^{-1},m_2^{-1}),
		m^{-1}_2g_2(m_0^{-1},m_1^{-1},m_2^{-1})].$}

	Proof: The point collineation
	$\bf{C} = \matb{(}{ccc}q_0&0&0\\0&q_1&0\\0&0&q_2\mate{)},$
	associates to (1,1,1), $(q_0,q_1,q_2),$ and to $(m_0,m_1,m_2),$
	$(r_0,r_1,r_2),$ if $r_i = q_i m_i.$\\
	In the new system of coordinates, \\
	$X = (q_0f_0(q_0^{-1}r_0,q_1^{-1}r_1,q_2^{-1}r_2),
		q_1f_1(q_0^{-1}r_0,q_1^{-1}r_1,q_2^{-1}r_2),
		q_2f_2(q_0^{-1}r_0,q_1^{-1}r_1,q_2^{-1}r_2)).$\\
	Exchanging $q_i$ and $r_i$ and then replacing $q_i$ by 1 and $r_i$
	by $m_i$ is equivalent to substituting $m_i$ for $q_i$ and 1 for $r_i,$ 
	which gives $\ov{X}$. $\ov{x}$ is obtained similarly.

	The line collineation is\\
	\hth$\matb{(}{ccc}q^{-1}_0&0&0\\0&q^{-1}_1&0\\0&0&q^{-1}_2\mate{)}.$

	\sssec{Theorem.}
	{\em Given a complete 5-angle, 5 distinct points, no 3 of
	which are on the same line, $A_0,$ $A_1,$ $A_2,$ $M$ and $\ov{M},$
	$A_i$ are called the vertices,
	$M$ is called the barycenter and $\ov{M},$ the orthocenter.}
	\begin{enumerate}
	\item  The ideal line and the orthic line. See Fig. 1,

	H0.0.\hti{3}$A_i,$\\
	H0.1.\hti{3}$M,$ $\ov{M}$,\\
	D0.0.\hti{3}$   a_i := A_{i+1} \times A_{i-1},$\\
	D0.1.\hti{3}$  ma_i := M \times A_i,
	      \ov{m}a_i := \ov{M} \times A_i,$\\
	D0.2.\hti{3}$  M_i := ma_i \times a_i,
	   \ov{M}_i := \ov{m}a_i \times a_i,$\\
	D0.3.\hti{3}$  mm_i := M_{i+1} \times M_{i-1},
       \ov{m}m_i := \ov{M}_{i+1} \times \ov{M}_{i-1},$\\
	D0.4.\hti{3}$  MA_i := a_i \times mm_i,
	  \ov{M}A_i := a_i \times \ov{m}m_i,$\\
	D0.7.\hti{3}$  m := MA_1 \times MA_2 (*),$
	     $\ov{m} := \ov{M}A_1 \times \ov{M}A_2 (*).$

	The nomenclature:\\
	N0.0.\hti{3}$  a_i$ are the {\em sides}.\\
	N0.3.\hti{3}$  ma_i$ are the {\em medians}, $\ov{m}a_i$ are the
	{\em comedians} or\\
	\hth$	\ov{m}a_i$ are the {\em altitudes}, $ma_i$ are the
	{\em coaltitudes},\\
	N0.4.\hti{3}$  M_i$ are the {\em mid-points} of the sides.\\
	\hth$	\ov{M}_i$ are the {\em feet} or the feet of the
	altitudes,\\
	N0.5.\hti{3}$  (M_i, mm_i)$ is the {\em complementary triangle},\\
	\hth$	(\ov{M}_i,$ $\ov{m}m_i)$ is the {\em orthic
	triangle},\\
	N0.6.\hti{3}$  MA_i$ are the {\em directions of the sides},\\
	N0.8.\hti{3}$  m$ is the {\em ideal line} corresponding to the line at
	 infinity,\\
	\hth$	\ov{m}$ is the {\em coideal line} or the {\em orthic
	line}, which is the polar\\
	\hth	    of $\ov{M}$ with respect to the triangle.

	Proof:\\
	P0.0.\hti{3}$   a_0 = [1,0,0],$\\
	P0.1.\hti{3}$  ma_0 = [0,1,-1],$     $\ov{m}a_0 = [0,-m_2,m_1],$\\
	P0.2.\hti{3}$  M_0 = (0,1,1),$	 $\ov{M}_0 = (0,m_1,m_2),$\\
	P0.3.\hti{3}$  mm_0 = [-1,1,1],$
	$\ov{m}m_0 = [-m_1m_2,m_2m_0,m_0m_1],$\\
	P0.4.\hti{3}$  MA_0 = (0,1,-1),$     $\ov{M}A_0 = (0,m_1,-m_2),$\\
	P0.7.\hti{3}$  m = [1,1,1],$
		    $\ov{m} = [m_1m_2,m_2m_0,m_0m_1],$

	\item  The line of Euler and the circle of Brianchon-Poncelet.
	See Fig. 2, 2b.

	Let\\
	D1.0.\hti{3}$   eul := M \times \ov{M}$\\
	D1.20.\hti{3}$  \gamma  := conic(M_0, M_1, M_2, \ov{M}_1,
		\ov{M}_2) (*),$\\
	then\\
	C1.1\hti{4}$\gamma\:is\:a\:circle,$ $\gamma\:is\:a\:cocircle = 0.$

	The nomenclature:\\
	N1.0.\hti{3}$  eul$ is the {\em line of Euler}.\\
	N1.11.\hti{3}$  \gamma$  is the {\em circle of Brianchon-Poncelet}.
	 In Euclidean
		geometry, the circle of Brianchon-Poncelet, is also called the
		circle of 9 points or circle of Feuerbach or, improperly, the
		circle of Euler.  It passes through the midpoints of the sides,
		the feet of the altitudes and the midpoints of the segment
		joining the vertices to the orthocenter.
	    The Definition-Conclusion D1.20. corresponds to the first part of
		the Theorem of Brianchon-Poncelet.

	Proof:\\
	P1.0.\hti{3}$   eul = [m_1-m_2,m_2-m_0,m_0-m_1],$\\
	P1.20.\hti{3}$  \gamma : m_1m_2 X_0^2+ m_2m_0 X_1^2 + m_0m_1 X_2^2$\\
	\hti{10}$- m_0(m_1+m_2) X_1 X_2 - m_1(m_2+m_0) X_2 X_0
			- m_2(m_0+m_1) X_0 X_1 = 0,$\\
	\hth$	\gamma ^{-1}: m_0^2(m_1-m_2)^2 x_0^2 + m_1^2(m_2-m_0)^2 x_1^2
			+ m_2^2(m_0-m_1)^2 x_2^2\\
	\hti{12}		- 2m_1m_2(m_0(3s1-2m_0)+m_1m_2) x_1 x_2\\
	\hti{12}		- 2m_2m_0(m_1(3s1-2m_1)+m_2m_0) x_2 x_0\\
	\hti{12}	- 2m_0m_1(m_2(3s1-2m_2)+m_0m_1) x_0 x_1.$

	\item  The circumcircle. See Fig. 4, 4b.

	Let\\
	D1.6.\hti{2}$ Imm_i := m \times \ov{m}m_i,$
	$\ov{I}mm_i := \ov{m} \times mm_i,$\\
	D1.7.\hti{2}$  ta_i := A_i \times Imm_i,$\\
	D1.19.\hti{2}$ \theta  := conic(A_1,ta_1,A_2,ta_2,A_0),$\\
	then\\
	C1.2.\hti{2}$ \ov{I}mm_i \cdot ta_i = 0.$\\
	H1.1.\hti{2}$  \theta\:is\:a\:circle = \theta\:is\:a\:cocircle = 0.$

	The nomenclature:\\
	N1.4.\hti{2}$ Imm_i$ are the {\em directions of the
	antiparallels} $a_i$ {\em with respect\\
	\hth	to the sides} $a_{i+1}$ and $a_{i-1}.$\\
	N1.5.\hti{2}$ ta_i$ are the {\em tangents at} $A_i$ {\em to the
	circumcircle},\\
	N1.10.\hti{2}$\theta$  is the {\em circumcircle}.

	Proof:\\
	P1.6.\hti{2}$ I\ov{m}m1_0 = (m_0(m_1-m_2),-m_1(m_2+m_0),
	m_2(m_0+m_1)),$\\
	\hth$	\ov{I}\ov{m}m1_0 = (m_0(m_2-m_1),-m_1(m_2+m_0),
	m_2(m_0+m_1)),$\\
	P1.7.\hti{2}$  ta_0 = [0,m_2(m_0+m_1),m_1(m_2+m_0)],$\\
	P1.19.\hti{2}$ \theta :\:m_0(m_1+m_2)X_1X_2 + m_1(m_2+m_0)X_2X_0\\
	\hti{12} + m_2(m_0+m_1)X_0X_1 = 0,$\\
	\hth$	2 \theta  + \gamma  = (m) \bigx (\ov{m}).\\
	\hth	\theta ^{-1}: m_0^2(m_1+m_2)^2 x_0^2 + m_1^2(m_2+m_0)^2 x_1^2
			+ m_2^2(m_0+m_1)^2 x_2^2\\
	\hti{12}- 2m_1m_2(m_2+m_0)(m_0+m_1) x_1 x_2\\
	\hti{12}- 2m_2m_0(m_0+m_1)(m_1+m_2) x_2 x_0\\
	\hti{12}- 2m_0m_1(m_1+m_2)(m_2+m_0) x_0 x_1 = 0,$

	\item  The point of Lemoine. See Fig. 3.

	Let\\
	D1.2.\hti{2}$Maa_i := ma_{i+1}\times \ov{m}a_{i-1},$
		$\ov{M}aa_i := ma_{i-1}\times \ov{m}a_{i+1},$\\
	D1.3.\hti{2}$mMa_i := Maa_i\times \ov{M}aa_i,$\\
	D1.4.\hti{2}$K := mMa_1 \times mMa_2 (*),$\\
	D1.8.\hti{2}$  T_i := ta_{i+1} \times ta_{i-1},$\\
	D12.1.\hti{2}$  at_i:= A_i \times T_i,$

	The nomenclature:\\
	N1.5.\hti{2}$ (T_i,ta_i)$ is the {\em tangential triangle},\\
	N12.1.\hti{2}$ at_i$ are the {\em symmedians}, of d'Ocagne,\\
	N1.2.\hti{2}$ K$ is the {\em point of Lemoine}, also called point of
	Grebe or of Lhuillier.
	\enume

	Proof:\\
	\hti{4}P1.2.\hti{2}$Maa_0 = (m_0,m_1,m_0)$,
		$\ov{M}aa_0 = (m_0,m_0,m_2),$\\
	\hti{4}P1.3.\hti{2}$mMa_0 = [q_0,m_0(m_2-m_0),-m_0(m_0-m_1)],$\\
	\hti{4}P1.4.\hti{2}$K = (m_0(m_1+m_2),m_1(m_2+m_0),m_2(m_0+m_1)),$\\
	\hti{4}P1.8.\hti{2}$T_0 = (-m_0(m_1+m_2),m_1(m_2+m_0),m_2(m_0+m_1)),$\\
	\hti{4}P12.1.\hti{2}$at_0 = [0,m_2(m_0+m_1),-m_1(m_2+m_0)],$

	\ssec{The fundamental $3*4+11*3\:\&\:3*5+10*3$ configuration.}

	\setcounter{subsubsection}{-1}
	\sssec{Introduction.}
	It would be desirable to have a synthetic proof of
	the sequence of Theorems given in this and in the following Chapters.
	In many instances, it is not difficult to obtain it, using the standard
	Theorems of projective geometry, mainly those of Pappus, Desargues and
	Pascal.  In other cases, the proof is less obvious.  Theorem 4.3.1.,
	which can be considered as the starting point, has a first part
	which required additional constructions.  The proof implies the validity
	of the extension of all the Theorems to finite
	projective geometries associated to Galois fields of order $p^{k},$
	$p > 3$ and to the projective geometries associated to the field of
	rationals, the field of reals, the field of complex numbers, the real
	p-adic field, the complex p-adic field,  $\ldots$.  For the second part,
	the proof is synthetic.

	\sssec{Theorem.}
	{\em Let $A_0,$ $A_1,$ $A_2,$ $M$ and $\ov{M}$ be a complete
	5-angle, see Fig. 0,\\
	\hth$		a_i := A_{i+1} \times A_{i-1}$\\
	\hth$    ma_i := M \times A_i,$
	 $ \ov{m}a_i := \ov{M} \times A_i$\\
	\hth$    M_i := ma_i \times a_i,$
	$\ov{M}_i := \ov{m}a_i \times a_i,$\\
	\hth$Maa_i := ma_{i+1}\times\ov{m}a_{i-1},$\\
	\hth$cc_i := A_i\times Maa_i,$\\
	\hth$P := cc_1\times cc_2 (*),$\\
	\hth$CA_i := cc_i\times a_i,$\\
	\hth$caa_i := CA_{i+1}\times CA_{i-1},$\\
	\hth$    c_i := M_{i+1} \times \ov{M}_{i-1},$\\
	\hth$CC_i := a_i \times c_i$\\
	\hth$		p := CC_1 \times CC_2 (*),$\\
	\hth$\gamma$ := conic$(M_i,\ov{M}_1,\ov{M}_2) (*)$.
	then\\
	\hth$    CC_i \cdot caa_i = 0.$

	The configuration involves the 14 points $A_i,$ $M_i,$
	$\ov{M}_i,$ $C_i,$ $M$
	and $\ov{M}$ and the 13 lines $a_i,$ $ma_i,$
	$\ov{m}a_i,$ $c_i$ and $p.$}

	Proof: For the first part, (see Fig. 0')\\
	dual-Pappus$(\langle ma_2,ma_0,ma_1\rangle,\langle \ov{m}a_1,\ov{m}a_2,
	\ov{m}a_0\rangle; \langle cc_0,cc_1,cc_2\rangle,P)$,\\
	therefore $cc_i$ are incident to $P$,\\
	Desargues$(P,\{A_i\},\{CA_i\};\langle CC_i\rangle,p)$,\\
	therefore $caa_i\times a_i$ are incident to $p$,\\
	Desargues$^{-1}(cc_0,\{ma_1,cc_1,a_1\},\{\ov{m}a_2,cc_2,a_2\};
	\langle caa_0,c_0,a_0\rangle,CC_0)$\\
	therefore $caa_0,$ $c_0$ and $a_0$ are incident to $CC_0$.

	For the second part, the Theorem of Pascal implies that the points
	$M_0,$	$\ov{M}_0,$ $M_1,$ $\ov{M}_1,$ $M_2,$
	$\ov{M}_2,$ are on a conic, because the points $C_0,$ $C_1$ 
	and $C_2$ are collinear.\\
	The conic may degenerate in two lines.  This will occur if, for
	instance, $\ov{M}_0$ is on $M_1 \times M_2$ and in this case
	$M_0$ is on
	$\ov{M}_1 \times \ov{M}_2.$  Indeed,

	\sssec{Theorem.}
	{\em Let $A_0,$ $A_1,$ $A_2,$ $M$ and $\ov{M}$ be a complete
	5-angle
	such that $\ov{M}_0,$ $M_1$ and $M_2$ are collinear then $M_0,$
	$\ov{M}_1$ and $\ov{M}_2$ 
	are collinear.}

	Proof:  A synthetic proof is as follows.
	Let $E := M_0 \times (M \times \ov{M}),$ the Theorem of Pappus
	applied to $A_j M M_j$ and
	$A_0 \ov{M}_0 \ov{M}$ for $j = 1$ and $2$ implies that
	$M_0,$ $E,$ $\ov{M}_j$ are
	collinear, therefore $M_0$ $\ov{M}_1$ and $\ov{M}_2$ are
	collinear.

	\sssec{Theorem.}
	{\em If $m = [1,1,1]$ and $\ov{m} = [m_0,m_1,m_2],$ then with
	respect to the line conic\\
	\hth$a_{00} x_0^2 + a_{11} x_1^2 + a_{22} x_2^2
		+ 2 (a_{12} x_1 x_2 + a_{20} x_2 x_0 + a_{01} x_0 x_1) = 0$\\
	the pole of $m$ is\\
	\hth$    (a_{00} + a_{01} + a_{20}, a_{01} + a_{11} + a_{12},
	a_{20} + a_{12} + a_{22})$\\
	and the pole of $\ov{m}$ is\\
	\hth$    (a_{00}m_1m_2 + a_{01}m_2m_0 + a_{02}m_0m_1,
	 a_{01}m_1m_2 + a_{11}m_2m_0 + a_{12}m_0m_1,\\
	\hti{12}	a_{20}m_1m_2 + a_{12}m_2m_0 + a_{22}m_0m_1).$}


	\ssec{An overview of theorems associated with bisected angles.
	The inscribed circle, the point of Gergonne, the point of Nagel.}

	\setcounter{subsubsection}{-1}
	\sssec{Introduction.}
	I will now give a construction associated to a conic inscribed in a
	triangle.  The degenerate case of the Theorem of Brianchon implies
	that if $JJ_i$ are the points of contact on $A_{i+1} \times  A_{i-1},$
	then the lines $A_i \times JJ_i$ pass through a point $J.$  We can
	choose arbitrarily a point $I$ or its polar $i.$  The construction in
	Theorem 3.12 determines a pair of points $M$ and $\ov{M}$ which
	in the case of Euclidean geometry will correspond to the barycenter and
	to the orthocenter.  As will be seen later, the function which
	associates $M,$ $\ov{M}$ to $J,$ $I$ is not one to one.  It is
	therefore
	necessary to start with this construction if we want to extend to
	projective geometry that part of the geometry of the triangle which
	is related to the inscribed circles.  In this case, Part 0. should
	precede Part 1.

	\sssec{Theorem.}
	{\em Given a complete 5-angle, 5 distinct points, no 3 of
	which are on the same line, $A_0,$ $A_1,$ $A_2,$ $J$ and $I,$
	$A_i$ are called the vertices,\\
	$J,$ is the {\em point of Gergonne} and\\
	$I,$ is the center of the {\em inscribed circle}.}

	\enumb
	\item The barycenter and orthocenter derived from the point of 
	Gergonne and the center of the inscribed circle.
	\enume

	Let\\
	H0.0.\hti{3}$	A_i,$ (See Fig. 20b)\\
	H0.2.\hti{3}$	J, I,$\\
	D0.0.\hti{3}$  a_i := A_{i+1} \times A_{i-1},$\\
	D0.8.\hti{3}$  ja_i := J \times A_i,$\\
	D0.9.\hti{3}$  JJ_i := ja_i \times a_i,$\\
	D0.10.\hti{2}$  j_i := JJ_{i+1} \times JJ_{i-1},$\\
	D0.11.\hti{2}$  Ja_i := j_i \times a_i,$\\
	D0.23'.\hti{2}$  ji_i := JJ_i \times I,$\\
	D0.26.\hti{2}$  Jia_i := ji_{i+1} \times a_{i-1},$
		$Ji\ov{a}_i :=  ji_{i-1}\times a_{i+1},$\\
	D0.27.\hti{2}$  jia_i := Jia_{i+1} \times Ja_{i-1},$
		$ji\ov{a}_i := Jia_{i-1} \times Ja_{i+1},$\\
	D0.28.\hti{2}$  Jai_i := jia_{i+1} \times j_{i-1},$
		$Ja\ov{i}_i := jia_{i-1} \times j_{i+1},$\\
	D20.0.\hti{2}$ Ji_i := jai_{i+1} \times ja\ov{i}_{i-1},$\\
	D20.22.\hti{2}$ \iota  := conic(JJ_0,JJ_1,JJ_2,Ji_1,Ji_2) (*),$\\
	D0.5'.\hti{2}$ m_i := Jai_i \times J\ov{a}i_i,$\\
	D0.6.\hti{3}$ MM_i := m_{i+1} \times m_{i-1},$\\
	D0.1'.\hti{2}$ ma_i := A_i \times MM_i,$\\
	D0.4'.\hti{2}$ MA_i := m_i \times a_i,$\\
	D0.7.\hti{3}$ m := MA_1 \times MA_2 (*),$\\
	D0.H.\hti{2}$ M := ma_1 \times ma_2 (*),$\\
	D0.25'.\hti{2}$ IMa_i := m \times ji_i,$\\
	D0.1'.\hti{2}$ \ov{m}a_i := A_i \times I\ov{m}a_i,$\\
	D0.H.\hti{2}$ \ov{M} := \ov{m}a_1 \times \ov{m}a_2 (*),$\\
	then\\
	C0.2.\hti{3}$  a_i \cdot \iota = 0.$\\
	C0.5.\hti{3}$  m_i \cdot A_i = 0.$\\
	C20.3.\hti{2}$ \iota\:is\:a\:circle = 0.$\\
	C20.4.\hti{2}$ I = Center(\iota ).$

	The nomenclature:\\
	N20.3.\hti{3}$\iota$  is the {\em inscribed circle},\\
	N0.12.\hti{3}$JJ_i$ are the {\em Gergonnian points}, these are the
	points of contact of the inscribed circle with the
	sides of the triangle.\\
	$Ja_i$ is the {\em pole of} $ja_i$ {\em with respect to
	the inscribed circle}.\\
	$Ji_i$ is the {\em point of the inscribed circle
		diametrically opposite to} $JJ_i.$\\
	Again, $M$ is the barycenter and $\ov{M}$ is the orthocenter.

	Proof:\\
	For a synthetic proof see section $\ldots$  G272.tex.\\
	Let $A_0 = (1,0,0),$ $A_1 = (0,1,0),$ $A_2 = (0,0,1),$
	$J = (j_0,j_1,j_2),$ $I = (i_0,i_1,i_2),$\\
	$m$ is constructed in such a way that $I$ is the pole of $m$ with
	respect to
	    $\iota$ , therefore if the line $m$ is chosen to be $[1,1,1],$ then
	\enumb
	\item$     I = (j_0(j_1+j_2),j_1(j_2+j_0),j_2(j_0+j_1)),$ therefore
	there is no loss
	    of generality if we set
	\item$     i_0 := j_0(j_1+j_2),$ $i_1 := j_1(j_2+j_0),$
	$i_2 := j_2(j_0+j_1).$
	    We will use the abbreviations for symmetric functions of $j_0,$
	$j_1,$ $j_2$
	    using ``p" instead of ``s" as used for the symmetric functions of
	    $m_0,m_1,m_2.$  For instance,
	\item$     p_{11} = j_1j_2+j_2j_0+j_0j_1.$
	\enume

	\noindent P0.0.\hti{3}$  a_0 = [1,0,0].$\\
	P0.8.\hti{3}$  ja_0 = [0,j_2,-j_1].$\\
	P0.9.\hti{3}$  JJ_0 = (0,j_1,j_2).$\\
	P0.10.\hti{3}$  j_0 = [-j_1j_2,j_2j_0,j_0j_1].$\\
	P0.11.\hti{3}$  Ja_0 = (0,j_1,-j_2).$\\
	P0.23.\hti{3}$  ji_0 = [j_1j_2(j_1-j_2),j_2j_0(j_1+j_2),
		-j_0j_1(j_1+j_2)].$\\
	P0.26.\hti{3}$  Jia_0 = (j_0(j_2-j_0),j_1(j_2+j_0),0),$
		$Ji\ov{a}_0 = (j_0(j_1-j_0),0,j_2(j_0+j_1)).$\\
	P0.27.\hti{3}$  jia_0 = [j_1j_2(j_0+j_1),j_2j_0(j_0+j_1),
		j_0j_1(j_1-j_0)],$\\
	\hth$ ji\ov{a}_0 = [j_1j_2(j_2+j_0),j_2j_0(j_2-j_0),
		j_0j_1(j_2+j_0)].$\\
	P0.28.\hti{3}$  Jai_0 = (j_0(j_1+j_2),-j_1j_2,j_1j_2),$
		$Ja\ov{i}_0 = (j_0(j_1+j_2),j_1j_2,-j_1j_2).$\\
	P20.0.\hti{2}$ Ji_0 = (j_0(j_1+j_2)^2,j_1j_2^2,j_1^2j_2).$\\
	P20.22.\hti{2}$ \iota : j_1^2j_2^2 X_0^2 + j_2^2j_0^2 X_1^2
		+ j_0^2j_1^2 X_0 X_1\\
	\hti{12}- 2j_0j_1j_2( j_0 X_1 X_2 + j_1 X_2 X_0 + j_2 X_0 X_1) = 0.$\\
	\hth$	j^{-1} \theta  + \iota  = - i \bigx
		 [j_1^2j_2^2,j_2^2j_0^2,j_0^2j_1^2].$\\
	\hth$	\iota ^{-1}: j_1j_2 x_1 x_2 + j_2j_0 x_2 x_0 + j_0j_1 x_0 x_1
		 = 0.$

	\noindent P0.25.\hti{2}$ IMa_0 = (j_0(j_1+j_2)^2,
		j_1(j_2^2-p_{11}),j_2(j_1^2-p_{11})).$\\
	P0.12.\hti{2}$ \ov{m}a_0 = [0,j_2(j_1^2-p_{11}),-j_1(j_2^2-p_{11})].$\\
	P0.16.\hti{2}$ \ov{M} = (j_0(j_1^2-p_{11})(j_2^2-p_{11}),
	    j_1(j_2^2-p_{11})(j_0^2-p_{11}), j_2(j_0^2-p_{11})(j_1^2-p_{11})).$

	For $m_i,$ $MM_i,$ $ma_i,$ $MA_i,$ $m$ and $M,$ see 3.7.\\
	The following relations are useful in the derivation of some
	of the formulas either given above or given below:
	\enumb
	\item$	m_0 = j_0(j_1^2-p_{11})(j_2^2-p_{11})$
	\item$	m_1+m_2 = -j_0(j_1+j_2)^2(j_0^2-p_{11}),$
	\item$	m_0(m_1+m_2) = (j_0(j_1+j_2))^2 jp,$ with
	\item$	jp = -(j_0^2-p_{11})(j_1^2-p_{11})(j_2^2-p_{11})$
	\item$	m_1-m_2 = -(j_1-j_2)(j_0^2-p_{11})(p_{11}+j_1j_2)$
	\item$	m_0(m_1-m_2) = -j_0(j_1-j_2)(p_{11}+j_1j_2)(j_0^2-p_{11})$
	\item$	s_1 = 4j_0j_1j_2p_{11},$
	\item$	s_1+m_0 = j_0(...),$
	\item$	j_1m_2-j_2m_1 = (j_2^2-j_1^2)(j_0^2-p_{11})p_{11},$
	\item$	m_1m_2 = -j_1j_2(j_0^2-p_{11})jp.$
	\enume


{\tiny\footnotetext[1]{G27.TEX [MPAP], \today}}

\setcounter{section}{1}
	\section{The Geometry of the Hexal Complete 5-Angles.}

	\setcounter{subsection}{-1}
	\ssec{Introduction.}

	Section \ldots contains a synthesis of a very large number of Theorems
	in Euclidean Geometry, using the presentation introduced in section i.
	This is followed by a proof given also as in section 1.\\
	The set of Theorems includes some which are always valid, some
	which are valid when the given triangle has a tangent circle and
	some which are valid when the point of Steven exits.
	\\
	In the second case I indicate that the definitions and Theorem are
	meaningful by labeling the section with (J). In the third case I label
	the section with (Mu), if neither case apply I label the section
	with (M). Definitions and conclusions contained in sections without
	(M), (J) or (Mu) are always meaningful.\\
	I start with a triangle $\{A_i\}$.\\
	In case (M), I choose the barycenter $M$ and the orthocenter $\ov{M}$.\\
	In case (Mu), I choose the barycenter $M$ = (m$_0$, m$_1$, m$_2$) and
	the point of Steven,\\
	\hth$Mu = (\sqrt{{\mb m}_0},\sqrt{{\mb m}_1},\sqrt{{\mb m}_2}).$
	This assumes that
	$k$m$_0$,$k$m$_1$ and $k$m$_2$ are quadratic residues for some $k$.
	I then determine the orthocenter from $M$ and $Mu$.\\
	In case (J), I assume that the triangle has a tangent circle and
	derive the orthocenter from the barycenter and the point of Gergonne
	$J$. The point $\ov{J}$, if it exist is such that the constructions
	obtained by use $\ov{J}$ instead of $J$ give eventually $\ov{M}$
	instead of $M$ and vice-versa. If the $\ov{J}$ does not exist
	these constructions are meaningless. The construction in the rightmost
	column of the sections marked with (J) should therefore be ignored.\\
	At the end of section 0, whatever the variant, the ideal line and
	orthic line have been constructed as well as the medians and altitudes,
	mid-points, the feet, and the complemantary and anticomplematary 
	triangles.\\
	In section 1, we construct the line $eul$ of Euler, the point $K$
	of Lemoine, the circumcircle $\theta$ and the circle $\gamma$ of
	Brianchon-Poncelet. Hypothesis \ldots

	Because in finite geometry
	If the section starts with (M), (J) or (Mu), it is only to 
	In this Chapter, I will give systematically most of the results which
	generalize the known results of the geometry of the triangle in
	classical Euclidean geometry.  In 5.1. and in 5.4. the corresponding
	constructions can be done with the ruler alone.  In 5.5. the
	corresponding constructions in classical Euclidian geometry would
	require also the compass.  (\#) is used to indicate Theorems obtained
	starting June 10, 1982, by systematically obtaining incidence relations
	on 2 examples and verifying that the conjecture so obtained is indeed a
	Theorem.

	What corresponds to isotropic points and foci of conics in the
	hexal complete 5-angles configuration is given in 5.2. and what
	corresponds to perpendicular directions,in 5.3.
	A summary of all incidence properties obtained in this Chapter is
	given in 5.7. to allow an easier access to the results.

	\ssec{The points of Euler, the center of the circle of
	Brianchon-Poncelet,
	and of the circumcircle, the points of Schr\"{o}ter, the point of
	Gergonne of the orthic triangle, the orthocentroidal circle.}

	\sssec{Theorem.}
	{\em Given a complete 5-angle, 5 distinct points, not 3 of
	which are on the same line, $A_0,$ $A_1,$ $A_2,$ $M$ and $\overline{M},$
	The vertices $A_i$ are those of a triangle, $M$ is the
	{\em barycenter} and $\overline{M}$ {\em is the orthocenter}.}


	\sssec{Proof of Theorem 5.1.1.}
	{\em 
	The algebraic proof will be summarized by giving the coordinates of the
	points and lines constructed in 5.1.1.  The incidence properties follow
	from straightforward computation of scalar products or substitution in
	the equation of the conics.\\
	For triples, the coordinates of the 0-th subscript will be given.  The
	coordinates of subscript 1 and 2 are obtained by applying the mapping
	$\rho$ and $\rho$ $^2$ to it.  $\rho$  is defined as follows,
	we substitute $m_1,m_2,m_0$ for $m_0,m_1,m_2$ in each of the
	components, and rotate these, the 0-th
	coordinate becoming the first, the first coordinate becoming the second
	and the second coordinate becoming the 0-th.
	For instance, from $em_0$ of 5.0. we get\\
	\hth    $em_1$ = $[s_1+m_1,m_2-m_0,-(s_1+m_1)]$ and
	$    em_2 = [-(s_1+m_2),s_1+m_2,m_0-m_1].$}

	The hypothesis imply, $m\neq 0,$ $m\neq 0,$ $m\neq 0,$ $m\neq m_2,$
	 $m\neq m_0$ and $m\neq m_1.$\\
	I will use the usual abbreviations for the symmetric functions,\\
	\hth$	s_1 := m_0 + m_1 + m_2, s_{11} := m_1m_2 + m_2m_0 + m_0m_1,$\\
	\hth$	s_2 := m_0^2 + m_1^2 + m_2^2, etc.$\\
	and\\
	\hth$	q_0 := m_0^2-m_1m_2, q_1 := m_1^2-m_2m_0, q_2 := m_2^2-m_0m_1,$\\
	and the identity\\
	\hth$	(m_1+m_2)(m_2+m_0)(m_0+m_1) = s_{211} + 2s_{111}.$

	The dual of the reciprocal of all elements have also been included.

	\sssec{Comment.}
	To determine in succession the homogeneous coordinates,
	we have used the definition.  To check the results, if for instance
	$x := P \times Q,$ we can simply verify $x \cdot P = x \cdot Q = 0.$
	The construction asserts implicitly that for $x := P \times Q,$ in
	general, $P$ and $Q$ are distinct, in other words for some value of
	$p$ and some $\overline{M},$
	$P$ and $Q$ are distinct.  It may of course happen that for a particular
	example $P = Q,$ 2 cases are possible, the coordinates of $x,$ for this
	example, are not all 0, this means that some alternate construction,
	using for instance one of the conclusions, will determine $x,$ in the
	other case, $x$ can not be constructed.
	For instance, for $p = 37$ and $\overline{M} = (202) = (1,4,16),$
	$Ste = AA_1 = (880),$
	but $stAA_1 = (472) = (1,11,27)$ is well defined but cannot be obtained
	using $Ste \times AA_1.$ 
	On the hand, for any $p,$ if $\overline{M} = (1,p-1,p-1),$
	$\overline{F}_0 = M_0 = (0,1,1)$
	and $\overline{f}m_0 = (0,0,0)$ and is therefore undefined.

	\sssec{Comment.}
	The determination of a conic, with known intersections $X1,$ $X2$
	with $a_0,$ Y1, Y2 with $b_0$ and $Z1,$ $Z2$ with $a_2,$ can be
	obtained easily.

	\ssec{Isotropic points and foci of conics.}

	\setcounter{subsubsection}{-1}
	\sssec{Introduction.}
	The following pairs of point can not be obtained by the
	construction involving only intersection of known lines or lines
	through known points, they are sufficiently important to be defined.

	\sssec{Definition.}
	$I,I' = m \times \gamma,$ $\overline{I},\overline{I}'
	= \overline{m} \times \gamma$.\\
	The first pair corresponds to the {\em isotropic points},
	the second pair to the {\em co-isotropic points}.

	\sssec{Theorem.}
	{\em With the definitions of Theorem 5.1, we have\\
	\hth$\overline{I},\overline{I}' = (m_0(m_1+m_2),-m_1(m_2+j\sigma ),
	-m_2(m_1-j\sigma )),$\\
	\hth$	j = +1$ or $-1,$ $\sigma  := \sqrt{-s_{11}},$\\
	\hth$I,I' = (m_0(m_1+m_2),-m_0m_1-j\tau ,-m_2m_0+j\tau ),$ where\\
	\hth$	j := +1$ or $-1,$ $\tau  := \sqrt{-m_0m_1m_2s_1}.$}

	\sssec{Definition.}
	$F$ is a {\em focus} of a non degenerate conic iff both $F \times I$ and
	$F \times \overline{I}$ are tangent to the conic.

	\sssec{Theorem.}
	{\em If the conic is not a parabola, there are 4 foci, real or
	complex.}

	\ssec{Perpendicular directions.}

	\sssec{Definition.}
	Two directions $IA$ and $IB$ {\em are perpendicular} iff one
	direction is on the polar of the other with respect to any circle.
	We will write  $IA \perp  IB.$

	\sssec{Theorem.}
	{\em 
	Let $(X_0,X_1,X_2)$ be an ideal point, the perpendicular direction is\\
	\hth$(m_0(m_1X_2-m_2X_1),m_1(m_2X_0-m_0X_2),m_2(m_0X_1-m_1X_0)),$}

	\sssec{Theorem.}
	{\em $(X_0,X_1,X_2)$ and $(Y_0,Y_1,Y_2)$ are perpendicular directions
		if}
	\enumb
	\item$    m_1m_2X_0Y_0 + m_2m_0X_1Y_1 + m_0m_1X_2Y_2 = 0.$
	\enume

	\sssec{Theorem.}
	{\em The following are perpendicular directions.  Let
	D0.\hti{5}$   Imeul = (s_1-3m_0,s_1-3m_1,s_1-3m_2),$\\
	then\\
	C0.\hti{5}$MA_i \perp  I\overline{m}a_i.$\\
	C1.\hti{5}$   Im_i \perp  I\overline{m}_i.$\\
	C2.\hti{5}$   I \perp  I, I' \perp  I'.$\\
	C3.\hti{5}$   EUL \perp  Imeul.$\\
	See also C12.4, C12.5, C16.7,}

	\sssec{Exercise.}
	Construct $Imeul$ of the preceding Theorem.

	\ssec{The circle of Taylor, the associated circles, the circle of
	 Brocard
	the points of Tarry and Steiner, the conics of Simson and of Kiepert,
	the associated circumcircles, the circles of Lemoine.}

	\setcounter{subsubsection}{-1}
	\sssec{Introduction.}
	Besides the properties given in Theorem 5.1.1., many
	other properties of Euclidean geometry generalize to projective
	geometry.  These will now be stated.  The numeration started in Theorem
	5.1.1. is continued.

	\sssec{Theorem.}
	{\em Given the hypothesis of Theorem 5.1.1. and the points and
	lines defined (or constructed) in that Theorem.}


	\sssec{Notation.}
	To make some of the algebraic expression less cumbersome,
	we have often used the symmetric functions\\
	\hth$s_1 := m_0+m_1+m_2,$\\
	\hth$s_{11} := m_1m_2+m_2m_0+m_0m_1,$\\
	\hth$s_2 := m_0^2+m_1^2+m_2^2,$\\
	\hth$s_{21} := m_0^2(m_1+m_2)+m_1^2(m_2+m_0)+m_2^2(m_0+m_1).$\\
	and similarly in the equations for conics other symmetric functions.
	We have also used, at times,\\
	\hth$    q_0 := m_0^2-m_1m_2, q_1 := m_1^2-m_2m_0, q_2 := m_2^2-m_0m_1,$\\
	and the following identities in the calculations:\\
	\hth$    m_1q_0 + m_2q_1 + m_0q_2 = 0 and m_2q_0 + m_0q_1 + m_1q_2 = 0.$

	\sssec{Proof of Theorem 5.4.1..}
	{\em The proof is given in the same way as the
	proof of 5.1.1.}


	\sssec{Proof of Theorem 5.4.3., P.15.16..}
	The details of the proof to obtain the equation of the circle of
	Brocard will now be given.
	The equation of a conic which has the radical axis m with the circle
	$\gamma$	 is
	\enumb
	\item$m_0(m_1+m_2) X_1 X_2 + m_1(m_2+m_0) X_2 X_0
	+ m_2(m_0+m_1) X_0 X_1$\\
	\hth$	+ (X_0 + X_1 + X_2) (u_0 x_0 + u_1 X_1 + u_2 X_2) = 0.$\\
	$u_0,$ $u_1$ and $u_2$ are determined in such a way that the conic passes
	through $Br3_i.$ 
	For $Br3_0$ we have\\
	\hth$	X_0 + X_1 + X_2 = 4m_1m_2+m_2m_0+m_0m_1,$\\
	and for the first line of 0.\\
	\hth$	(4m_1m_2+m_2m_0+m_0m_1) (m_1m_2(m_2+m_0)(m_0+m_1)).$\\
	Hence we have to solve\\
	\hth$2 m_1m_2    u_0 + m_1(m_2+m_0) u_1 + m_2(m_0+m_1) u_2
	+ m_1m_2(m_2+m_0)(m_0+m_1) = 0,$\\
	\hth$m_0(m_1+m_2) u_0 + 2 m_2m_0    u_1 + m_2(m_0+m_1) u_2
	+ m_2m_0(m_0+m_1)(m_1+m_2) = 0,$\\
	\hth$m_0(m_1+m_2) u_0 + m_1(m_2+m_0) u_1 + 2 m_0m_1    u_2
	+ m_0m_1(m_1+m_2)(m_2+m_0) = 0.$\\
	Replacing $u_0,$ $u_1$ and $u_2$ in terms of $v0,$ $v1$ and $v2,$ given by
\\
	\hth$v0 := \frac{u_0}{m_1m_2}, v1 := \frac{u_1}{m_2m_0},
	v2 := \frac{u_2}{m_0m_1}, we get$\\
	\hth$  2 m_1m_2    v0 + m_1(m_2+m_0) v1 + m_2(m_0+m_1) v2
	+ (m_2+m_0)(m_0+m_1) = 0,$\\
	\hth$  m_0(m_1+m_2) v0 + 2 m_2m_0    v1 + m_2(m_0+m_1) v2
	+ (m_0+m_1)(m_1+m_2) = 0,$\\
	\hth$  m_0(m_1+m_2) v0 + m_1(m_2+m_0) v1 + 2 m_0m_1    v2
	+ (m_1+m_2)(m_2+m_0) = 0.$\\
	The determinant is $D = -6 s_{222} + 2 s_{33}.$
	The numerator for $v0$ is $E(m_0+m_1)(m_2+m_0)$ with
	$E = (s_{211}-s_{22}).$\\
	Hence the solution for $u_0.$\\
	To obtain P15.16., we have to determine\\
	\hth$u_1+u_2+m_0(m_1+m_2) = \frac{m_0(m_1+m_2)(D + Es_{11} + Em_1m_2)}{D}
$\\
	\hth$	= \frac{m_0(m_1+m_1)(-3s_{222} + s_{33} + Em_1m_2)}{D},$\\
	because $Es_{11} = s_{211}s_{11}-s_{22}s_{11}.$
	But $(s_{211} - s_{22})s_{11} = 3s_{222} - s_{33},$
	hence the equation for PUb.
	\enume

	To obtain the relation between $\theta$  and $\beta$ , knowing that\\
	\hth$	A \theta  + \beta  = B (m) \bigx (lem),$\\
	it is easy to obtain $A$ and $B,$ for instance, $K.\beta$  gives\\
	\hth$	A (3m_0m_1m_2(m_1+m_2)(m_2+m_0)(m_0+m_1))$\\
	\hth$	= B (2s_{11} 3m_0m_1m_2(m_1+m_2)(m_2+m_0)(m_0+m_1))$\\
	therefore with $B = 1,$ $A = 2s_{11}.$

	       ADD PERPENDICULARITY, e.g. $IiI_i$ $\perp$ $Iai_i.$
	 Cross refer. at end of G2705

	\ssec{Theorems associated with bisected angles.  The outscribed
	circles, the circles of Spieker, the point of Feuerbach, the
	barycenter of the excribed triangle.}

	\setcounter{subsubsection}{-1}
	\sssec{Introduction.}
	I will now give a construction associated to a conic
	inscribed in a 	triangle.  The degenerate case of the Theorem of
	Brianchon implies that if $JJ_i$ are the points of contact on
	$A_{i+1} A_{i-1},$
	then the lines $A_i \times JJ_i$ pass through a point $J.$ We can choose
	arbitrarily a point $I$ or its polar $i.$  The construction in
	Theorem 5.5.1. determines a pair of points $M$ and $\overline{M}$
	which in the case
	of Euclidean geometry will correspond to the barycenter and to the
	orthocenter.  As will be seen later, the function which
	associates $M,$ $\overline{M}$ to $J,$ $I$ is not one to one.  It is
	therefore
	necessary to start with this construction if we want to extend to
	projective geometry that part of the geometry of the triangle which
	is related to the inscribed circles.  Part 0. should therefore
	precede Part 1. of Theorem 5.1.1.   Part 20. given next follows
	Part 19. of Theorem 5.1.1.

IN THE NEXT SECTION REVERSE THE ORDER. FIRST SHOW THAT THE point diametrically
	opposed to $JJ_0$ is on the line $m_0 \times j_2\times JJ_1,$
	then that $m_0 \times j_2$ is on the line $j_0 \times a_0$ and 
	$ij_0 \times a_1$\\
STUDY	from which it follows that $m_0 \times j_2$ defines $m_0$ with $A_0$ and
	can be obtained from $JJ_i.$ 

	\sssec{Heuristics.}
	Before giving the construction I will look back at
	Euclidean geometry and determine properties which have guided me in
	the construction given below.  Let I be the center of the circle
	$\iota$	 inscribed in the triangle $(A_0,$ $A_1,$ $A_2),$ let $JJ_i$
	be the point of contact with $a_i,$ let $m_i$ be the parallel to
	$a_i$ through $A_i.$	\\
	First, if $Ja_0 := j_0 \times a_0,$ $Jia_2 := a_1 \times ji_0,$ 
	and $Jai_0 := j_2 \times m_0,$ then $Ja_0,$ $Jia_2$ and $Jai_0$ are
	collinear
	because they are the Pascal points of the hexagon with cords or
	tangents $j_0,$ $a_1,$ $j_2,$ $a_0,$ $ji_0,$ $jai_1.$ 
	If we start from $A_i,$ $J$ and $I$ we can therefore construct $JJ_i,$ 
	$Ja_0,$	$Jia_2,$ $Jai_0$ and $m_0,$ hence $MA_0 := a_0 \times m_0,$ 
	similarly we can construct $MA_1$ and therefore the ideal line
	$m := MA_0 \times MA_1.$  (The construction below is a variant which
	uses a "symmetric" point $J\overline{a}i_0$ also on $m_0.)$ 
	From $m$ we can derive the barycenter $M$ as the polar of $m$ with
	respect to the triangle $\{A_i\}.$ 
	Next, the conic through $JJ_i$ with tangent $a_1,$ $a_2$ can be defined
	as a circle, the altitude $\overline{m}a_0$ can be obtained as parallel
	to $I \times JJ_0$ and therefore the orthocenter $\overline{M}$ can be
	constructed.\\
	Finally, let $jai_1 := Jai_0 \times JJ_2$ and
	$Ji_0 := jai_1 \times ji_0,$ 
	I claim that $Ji_0$ is on the inscribed circle.  Indeed, first the
	triangles $(JJ_{i-1},A_i,JJ_{i+1})$ are isosceles triangles, then,
	for $i = 0,$ the triangle $(JJ_1,A_0,Jai_0)$ which is similar to the
	triangle $(JJ_1,$ $A_2,$ $JJ_0)$ is therefore an isosceles triangle and
	$|A_0,Jai_0| = |A_0,JJ_1| = |A_0,JJ_2|.$ 
	Therefore $angle(A_0,JJ_2,Jai_0) = \frac{1}{2}(\pi - 
	angle(JJ_2,A_0,Jai_0)) = \frac{1}{2}angle(A_0,A_1,A_2)
	= angle(A_0,A_1,I),$  therefore $j_0$ is parallel to $A_1 \times I$
	and therefore perpendicular to $j_1,$ it follows
	that $(Ji_0,JJ_0)$ is a diameter.\\
	We can therefore construct $Ji_0$ on $\iota$ .


	\sssec{Proof of Theorem 5.5.2..}
	\enumb
	\item$  A synthetic proof of \ldots  is as follows.  Pascal's Theorem
	 gives$\\
	\hth$    ct(ji_0,a_0,j_1,a_2,j_0,j\overline{a}i_2?) = 
	(J\overline{i}a_1,Ja_0)$\\
	\hth$	\Rightarrow   Ja\overline{i}_0  \Rightarrow   
	ja\overline{i}_2  \Rightarrow   Ji_0.$\\
	\hth$    ct(ji_0,a_0,j_2,a_1,j_0,jai_1?) = (Jia_2,Ja_0)$\\
	\hth$	\Rightarrow   Jai_0  \Rightarrow   jai_1  \Rightarrow   Ji_0,$
 hence 0.0 and 0.2.\\
	\hth$    ct(j_2,a_1,j\overline{a}i_2,jai_1,a_2,j_1) =
	(Jai_0,A_0,J\overline{a}i_0),$\\
	\hth	which are therefore collinear, hence 0.3.\\
	\hth$    ct(jai_1,j_1,a_0,j_2,j\overline{a}i_2,tangent(Ji_0)) =$\\
	\hth$	(Jai_0,J\overline{a}i_0,MA_0),$ hence 0.4.\\
	\hth$    Im_i$ is the pole of $ji_i,$ therefore $i$ is the polar of $I,$
		 hence 0.5.\\

	    The coordinates of the various points are easy to derive.
	Let $A_0 = (1,0,0),$ $A_1 = (0,1,0),$ $A_2 = (0,0,1),$
	$J = (j_0,j_1,j_2),$ $I = (i_0,i_1,i_2),$\\
	$m$ is constructed in such a way that $I$ is the pole of $m$ with
	respect to
	    $`\iota$ , therefore if the line $m$ is chosen to be $[1,1,1],$ then
	\item$     I = (j_0(j_1+j_2),j_1(j_2+j_0),j_2(j_0+j_1)),$ therefore
	there is no loss
	    of generality if we set
	\item$     i_0 := j_0(j_1+j_2),$ $i_1 := j_1(j_2+j_0),$
	$i_2 := j_2(j_0+j_1).$\\
	    We will use the abbreviations for symmetric functions of $j_0,$
	 $j_1,$ $j_2$
	    using "p" instead of "s" as used for the symmetric functions of
	    $m_0,m_1,m_2.$  For instance,
	\item$     p_{11} = j_1j_2+j_2j_0+j_0j_1.$
	We have also expressed the coordinates in terms of $i_0,$ $i_1$ and
	$i_2.$
	The symmetric functions of $i_0,$ $i_1$ and $i_2$ use "o" instead of
	"s".
	The expression of $j_0,$ $j_1$ and $j_2$ in terms of $i_0,$ $i_1$ and
	 $i_2$ is given by
	\item$     j_0 = \frac{(o_1-2i_1)(o_1-2i_2)}{ip},  \ldots , where$
	\item$     ip^2 = 2(o_1-2i_0)(o_1-2i_1)(o_1-2i_2).$
	This alternate notation has the advantage that the information on the
	associate construction for the excribed circles is obtained by
	replacing either $i_0$ by $-i_0,$ or $i_1$ by $-i_1,$ or $i_2$ by
	$-i_2.$
	\enume

	\sssec{Proof of Theorem 5.5.2., P21.8.}
	The proof or the preceding theorem is straightforward, I will only
	give details for the determination of $\pi$ :
	Let {\bf C} be the symmetric matrix associated to the polarity of $pi$,
	let {\bf M} be the matrix whose i-th column are the coordinates of
	$m_i,$ 
	let {\bf J} be the matrix whose i-th column are the coordinates of
	$Mna_i,$ 
	let {\bf K} be a diagonal matrix of unknown scaling factors $k_0,$
	$k_1,$ $k_2.$\\
	\hth$	{\bf C J} = {\bf M K} or {\bf C} = {\bf M K J}^-1$\\
	expresses the fact that $m_i$ is the polar of $Mna_i.$ \\
	\hth${\bf J}^{-1} = \mat{p_{11}}{2j_0j_1-p_{11}}{2j_2j_0-p_{11}}
		{2j_0j_1-p_{11}}{p_{11}}{2j_1j_2-p_{11}}
		{2j_2j_0-p_{11}}{2j_1j_2-p_{11}}{p_{11}}.$\\
	The problem is now reduced to a set ot 3 homogeneous equations in the
	unknowns $k_0,$ $k_1,$ $k_2,$ which express the symmetry of {\bf C},
	namely, after simplification,\\
	\hth$   - j_0j_1 k_0      + j_0j_1 k_1 + j_2(j_1-j_0) k_2 = 0,$\\
	\hth$j_0(j_2-j_1) k_0      - j_1j_2 k_1      + j_1j_2 k_2 = 0,$\\
	\hth$     j_2j_0 k_0 + j_1(j_0-j_2) k_1      - j_2j_0 k_2 = 0,$\\
	giving $k_0 = j_1j_2,$ $k_1 = j_2j_0,$ $k_2 = j_0j_1.$

	\sssec{Comment.}
	The following alternate definition for Nagel's point which
	is clearly more clumsy:\\
	\hth$    oi := O \times I,$\\
	\hth$    Ioi := i \times oi,$\\
	\hth$    n\overline{m} := \overline{M} \times Ioi,$\\
	\hth$    N := n\overline{m} \times mi,$

	We have\\
	\hth$    oi = [ j_1j_2(j_1-j_2)(j_2+j_0)(j_0+j_1),
j_2j_0(j_2-j_1)(j_0+j_1)(j_1+j_2),$\\
	\hth$		j_0j_1(j_0-j_2)(j_1+j_2)(j_2+j_0) ],$\\
	\hth$    Ioi = (j_0(j_1+j_2)(p21-2j_0p_{11}), \ldots ),$\\
	\hth$    n\overline{m} = [j_0(j_1^2-j_2^2)(j_0^2-p_{11}),$ \ldots ],
 done backward from $N.$\\
	We also have\\
	\hth$	    mj = [j_1-j_2,j_2-j_0,j_0-j_1],$

	\sssec{Comment.}
	An alternate method to obtain quickly the relation between
	the barycentric coordinates of the point of Gergonne and of the
	orthocenter is as follows.\\
	Let $n_0 = m_0(m_1+m_2),$  $\ldots$  , we known that are circles are\\
	\hth$n_0 x_1 x_2 + \ldots - (X_0 + X_1 + X_2)(u_0 X_0 + \ldots ) = 0.$\\
	the line equation of the inscribed circle is\\
	\hth$	j_1j_2 x_1x_2 +  \ldots  = 0$\\
	to express that it is a circle we can use\\
	\hth$	{\bf A} = 2 adjoint({\bf B}),$\\
	where {\bf A} is the polarity matrix associated to the general circle
	and {\bf B} the matrix associated to (2).  The constant is arbitrary
	and reflect the chosen scaling.\\
	\hth${\bf A} = \mat{-2u_0}{n_2-u_0-u_1}{n_1-u_2-u_0}
		{n_2-u_0-u_1}{-2u_1}{n_0-u_1-u_2}
		{n_1-u_2-u_0}{n_0-u_1-u_2}{-2u_2},
	{\bf B} = \mat{0}{j_0j_1}{j_2j_0}
		{j_0j_1}{0}{j_1j_2}
		{j_2j_0}{j_1j_2}{0}.$\\
	This gives at once\\
	\hth$u_0 = (j_1j_2)^2$ and\\
	\hth$n_0 = 2j_0^2j_1j_2+u_1+u_2 = (j_0(j_1+j_2))^2.$

	\sssec{Comment.}
	To obtain the points of contact of the outscribed circle $\iota$ $_0,$ 
	Let $J_0$ be the corresponding point of Gergonne $[g_0,g_1,G_2].$
	We have\\
	\hth$	g_0(g_1+G_2) = -j_0(j_1+j_2),$\\
	\hth$	g_1(G_2+g_0) = j_1(j_2+j_0),$\\
	\hth$	G_2(g_0+g_1) = j_2(j_0+j_1),$\\
	adding 2 equations and subtracting the third gives\\
	\hth$    g_1G_2 = p_{11},$ $G_2g_0 = -j_0j_1,$ $g_0g_1 = -j_2j_0,$ with
	$p_{11} = j_1j_2+j_2j_0+j_0j_1.$\\
	Hence with an appropriate constant of proportionality,\\
	\hth$	(g_0,g_1,G_2) = (-j_0j_1j_2,j_2p_{11},j_1p_{11}).$\\
	Therefore the points of contact with $a_0,$ $a_1$ and $a_2$ are\\
	\hth$    Na_0 = (0,j_2,j_1),$ $Nai_0 = (j_2j_0,0,-p_{11}),$
	$Na\overline{i}_0 = (j_0j_1,-p_{11},0).$

	\sssec{Theorem.}
	{\em The isogonal transformation of J is\\
	\hth$	isog(J) = ( j_0(j_1+j_2)^2,j_1(j_2+j_0)^2,j_2(j_0+j_1)^2),$\\
	\hth$			-j_0j_1(j_2+j_0)(p1+j_2) ),$}

	Proof:\\
	$x_0 = [-j_1j_2,j_2j_0,j_0j_1],$\\
	\hth$X_0 = (j_0(j_2+j_1),j_1(j_2+j_0),-j_2(j_0+j_1)),$\\
	\hth$Y_0 = (j_0(j_2+j_0),j_1(j_2+j_0),2j_2j_0),$\\
	\hth$y_0 = [j_1j_2(j_2+j_0)(j_0+j_1),
	j_2j_0(j_0j_1-j_1j_2+3j_2j_0+j_0^2),$\\
	\hth$Z_0 = (j_0(j_1^2+j_2^2-j_1j_2-j_0j_1),j_1(j_2+j_0)^2,
	j_1(j_0+j_1)^2),$\\
	\hth$	z_0 = [0,j_2(j_0+j_1^2),-j_1(j_2+j_0)^2],$ hence the Theorem.

	\sssec{Definition.}
	Many other constructions can be easily derived from the
	following operation called the {\em dual construction}.
	Instead of the quintuple $A_i$, $M,$ $\overline{M},$ consider instead
	the quintuple $A_i$, $M,$ $\overline{M}' := T\overline{m}m.$\\
	The construction associated to every point $X =(X_0,X_1,X_2),$ a point
	$X'$ whose coordinates are the reciprocal $X' = (X_1X_2,X_2X_0,X_0X_1)$
	and to every line $x = (x_0,x_1,x_2)$ the reciprocal
	$x' = (x_1x_2,x_2x_0,x_0x_1).$\\
	A few of the dual points and lines are not new but most are and
	lead easily to the construction of important points and lines.
	See for instance the exercise on the line of Longchamps.
	We have
	$ma'_i = ma_i,$ $M'_i = M_i,$ $m'_i = m_i,$ $Im'_i = Im_i,$
	$AC'_i = AC_i,$
	$i' = i,$ $O' = K,$ $oa'_i = at_i,$ $\overline{S}' = \overline{S},$
	$Ima'_i = Ima_i,$ $K' = O,$
	$ok' = ok.$\\
	An example is given in ex3.3.

	\sssec{Corollary.}
	{\em We can now summarize incidence properties associated
	with the historically important line of Euler and circle of Brianchon-
	Poncelet.}

	\enumb
	\item The following 14 points are on the line of Euler:
	the barycenter $M,$ the orthocenter $\overline{M},$ the point $PP$
	 of D3.3, the center $EE$ and cocenter $\overline{E}E$ of the circle
	 of Brianchon-Poncelet, the center $O$ and cocenter $\overline{O}$ of
	the circumcircle, the points $Am,$ $\overline{A}m$ of D7.9,
	the points $D_i$ of D8.4, the center $G$ and the cocenter
	 $\overline{G}$ of the orthocentroidal circle.

	\item  The following 24 points are on the conic of Brianchon-Poncelet:
	    the midpoints $M_i$, the feet $\overline{M}_i$, the
	Euler points $E_i$, $\overline{E}_i$, the
	    points $F_i$ and $\overline{F}_i$ of D6.2, points of Schr\"{o}ter
	 $S$ and  $\overline{S},$ the
	    points of Feuerbach, $Fe$ and $Fe_i$.
	\enume

	The complete set of incidence properties are given in detail in
	section 5.7.

	\sssec{Comment.}
	Given the algebraic coordinates of a point it is sometimes
	difficult to obtain a construction starting from $M$ and
	$\overline{M}.$  One
	additional tool is provided by using homologies.  We will give here
	an example, which allows the easy construction of other points on the
	line of Euler.

	\sssec{Definition.}
	A {\em barycentric homology} is a homology with center $M$ and
	axis $m.$

	\sssec{Example.}
	One such homology and its inverse is\\
	\hth$	{\bf D} = \mat{0}{1}{1}{1}{0}{1}{1}{1}{0},
	 {\bf D}^{-1} = \mat{-1}{1}{1}{1}{-1}{1}{1}{1}{-1}.$

	\sssec{Theorem.}
	{\em The transforms of the 14 points on the line of Euler, given
	in 5.10.0 are as follows:\\
	{\bf D}$(M) = M,$\\
	{\bf D}$(\overline{M}) = O,$\\
	{\bf D}$(PP) = (q_1+q_2,q_2+q_0,q_0+q_1),$\\
	{\bf D}$(EE) = (3s_1-m_0,3s_1-m_1,3s_1-m_2),$\\
	{\bf D}$(\overline{E}E) = (3s_{11}-m_1m_2,3s_{11}-m_2m_0,
	3s_{11}-m_0m_1),$\\
	{\bf D}$(O) = EE,$\\
	{\bf D}$(\overline{O}) = (s_{21}-m_0^2(m_1+m_2),s_{21}-m_1^2(m_2+m_0),
	s_{21}-m_2^2(m_0+m_1)),$\\
	{\bf D}$(Am) = \overline{M},$\\
	{\bf D}$(\overline{A}m) = (s_{21}-s_{111}-m_0s_{11},
	s_{21}-s_{111}-m_1s_{11},s_{21}-s_{111}-m_2s_{11}),$\\
	{\bf D}$(D_0) = (m_0(m_1+m_2)-2m_1m_2,m_0(m_2+m_0)-2m_1m_2,
	m_0(m_0+m_1)-2m_1m_2),$\\
	{\bf D}$(G) = (5s_1-3m_0,5s_1-3m_1,5s_1-3m_2),$\\
	{\bf D}$(\overline{G}) = (s_{21}-9s_{111}-m_0s_{11},
	s_{21}-9s_{111}-m_1s_{11},s_{21}-9s_{111}-m_2s_{11}),$}

	\sssec{Exercise.}
	Complete the table for the inverse transform,
	{\bf D}T-1T$(M) = M,$\\
	{\bf D}T-1T$(\overline{M}) = (s_1-3m_0,s_1-3m_1,s_1-3m_2).$\\
	Observe that {\bf D}T-1T$(\overline{M}).m = 0.$

	\ssec{Duality and symmetry for the inscribed circle.}

	\setcounter{subsubsection}{-1}
	\sssec{Introduction.}
	$\ldots$	

	\sssec{Theorem.}
	{\em If $m_1+m_2,$ $m_2+m_0$ and $m_0+m_1$ are all quadratic residues or
	or non quadratic residue, then both the dual of the inscribed circle
	and the symmetric of the inscribed circle are real.  Moreover,
	if $i$ and $j$ are the dual of $I$ and $J$ and if $\overline{I}$ and
	$J\overline{}$ are the symmetric of $I$ and $J$ then}
	\enumb
	\item$    i = [ \sqrt{m_1+m_2}, \sqrt{m_2+m_0}, \sqrt{m_0+m_1} ],$
	\item$j = [ (-i_0+i_1+i_2)^{-1}, (i_0-i_1+i_2)^{-1},
	 (i_0+i_1-i_2)^{-1} ]$\\
	and
	\item$    \overline{I} = (m_0 i_0, m_1 i_1, m_2 i_2),$
	\item$    \overline{J} = (m_0 j_0, m_1 j_1, m_2 j_2),$
	\enume

	Proof:\\
	For the symmetric case, $\overline{I} \times (a_0 \times (A_0 \times 
	\overline{J})),$ $A_0 \times MA_0$ and $\overline{m}$
	are concurrent, moreover $\overline{I}$ is the pole of $\overline{m}$
	with respect to $\iota$ .\\
	Therefore,
	\enumb
	\setcounter{enumi}{3}
	\item$\overline{J}_1(\overline{I}_2 m_0(m_1+m_2)
	+ \overline{I}_0 m_1m_2) =$\\
	\hth$		\overline{J}_2(\overline{I}_0 m_1m_2
	+ \overline{I}_1 m_0(m_1+m_2)),$\\
	and in view of P0.15,
	\item$    \overline{I}_0 = \overline{J}_0\overline{J}_1 m_2m_0
	+ \overline{J}_2\overline{J}_0 m_0m_1.$
	This relation and the 2 others obtained by circularity give\\
	\hth$-\overline{I}_0 m_1m_2 + \overline{I}_1 m_2m_0 
	+ \overline{I}_2 m_0m_1 = 2 \overline{J}_2^{-1} m_2$
	\enume
	Using 4, we get\\
	\hth$	\overline{I}_0^2 \frac{(m_2-m_1)}{m_0^2}
	- \overline{I}_1^2 \frac{(m_1+m_2)}{m_1^2} 
	+ \overline{I}_2^2 \frac{(m_1+m_2)}{m_2^2} = 0,$\\
	as well as 2 other similar equations.  These equations are compatible
	and give using the minors\\
	\hth$(\frac{\overline{I}_0}{m_0})^2 = m_1+m_2,$
	$(\frac{\overline{I}_1}{m_1})^2 = m_2+m_0,$
	$(\frac{\overline{I}_2/m_2})^2 = m_0+m_1.$\\
	For the dual case, it follows from .1 and .6 (G2722), that the
	coordinates of I are proportional to $\sqrt{m_0(m_1+m_2)},$  $\ldots$,
	 those of the
	dual are obtained by replacing $m_0$ by $m_1m_2,$ \ldots.

	\sssec{Theorem.}
	{\em If $m_0,$ $m_1$ and $m_2$ are all quadratic residues or all non
	quadratic residues, then the dual of the symmetric of the inscribed
	circle is real.  Morover if $\overline{j}$ and $\overline{i}$ are the
	dual of the symmetric of $J$ and $I,$ then}
	\enumb
	\item 0.$    \overline{i} = [ \sqrt{\frac{m_1+m_2}{m_0}},
	 \sqrt{\frac{m_2+m_0}{m_1}}, \sqrt{\frac{m_0+m_1}{m_2}} ],$
	{\em and}\\
		1.$    \overline{j} = [ (m_0(-m_0 \overline{i}_0
	+ m_1 \overline{i}_1 + m_2 \overline{i}_2))^{-1},$\\
	\hth$		 (m_1( m_0 \overline{i}_0 - m_1 \overline{i}_1
	+ m_2 \overline{i}_2))^{-1},$\\
	\hth$		 (m_2( m_0 \overline{i}_0 + m_1 \overline{i}_1
	- m_2 \overline{i}_2))^{-1}].$
	\enume

	\sssec{Example.}
	For $p = 29,$ if $M = (60)$ and $\overline{M} = (258) = (1,7,25),$
	$(1,\frac{m_2+m_0}{m_1+m_2},\frac{m_0+m_1}{m_2+m_0}) = (1,-1,-7),$ with
	 a choice of the
	square roots, $i_0$ = 1, $i_1$ = -12, $i_2$ = -14, hence $i = [538]$ and
	$j = [1,-2,-9] = [833].$
	Moreover $\overline{I} = (1,3,-2) = (144),$ and
	$\overline{J} = (1,-14,7) = (472).$
	$(m_0,m_1,m_2) = (1,7,25),$ with a choice of the square roots,
	$(\sqrt{m_0},\sqrt{m_1},\sqrt{m_2}) = (1,6,5),$ hence
	$\overline{i} = [1,-2,3] = [816]$ and
	$\overline{j} = [1,7,-14] = [248].$
	then $J = (164),$ $I = (448),$ $\overline{I} = (144),$
	$ \overline{J} = (472),$
	and $i = [538],$ $j = [833],$ $\overline{i} = [816],$
	$ \overline{j} = [248].$

	\ssec{Summary of the incidence properties obtained so far }

	\setcounter{subsubsection}{-1}
	\sssec{Introduction.}
	The incidence properties of points, lines and conics
	will now be summarized.  There are several reasons for doing this.
	First, having so many elements, it is difficult to keep in ones mind
	at any one time all of the properties given above.  Second, it is
	important to insure that the elements obtained are in general distinct.
	Third, it is important to obtain from the elements defined any
	incidence properties not already discovered.  For this purpose, I
	create a program, which, for given examples, determine all incidence
	properties, by comparison, it was possible to eliminate a few
	incidence properties which were peculiar to a given example, for the
	others attempting an algebraic proof determined if the incidence
	property was indeed general.  Quite a few new Theorems were obtained
	in this way.  They have been indicated by (\#).

	I have ordered them in the order of the definitions.  The notation is
	self explanatory.

	\sssec{Theorem.}
	{\em The incidence properties are as follows:}

	\sssec{Proof of Theorem 6.1.1.}
	{\em 
	\sssec{Exercise.}
	Construct the vertical tangent of the parabola of Kiepert
	and prove that it is\\
	\hth$[ m_0(m_1-m_2)(s_1-3m_1)(s_1-3m_2),
	m_1(m_2-m_0)(s_1-3m_2)(s_1-3m_0),$\\
	\hth$		m_2(m_0-m_1)(s_1-3m_0)(s_1-3m_1) ].$}

	\sssec{Exercise.}
	Construct the conic of Jerabek (Vigarie, N.99),\\
	\hth$	m_0(m_1^2-m_2^2) X_1 X_2 + m_1(m_2^2-m_0^2) X_2 X_0$\\
	\hth$		+ m_2(m_0^2-m_1^2) X_0 X_1 = 0.$


	\sssec{Comment.}
	There exist a large number of conditional theorems.  For
	instance, if $s_{21} + 12s_{11} = 0$ then $G\cdot \ov{i} = \ov{G}\cdot i
	= 0.$\\
	An example is provided by $p = 29,$ $m_0 = 1,$ $m_1 = 6,$ $m_2 = 11,$
	corresponding to $J = 94,$ $I = 315.$

	\sssec{Exercise.}
	The line of Simson.\\
	Let\\
	D.0.	$Y_i := X \times Im_i,$\\
	D.1.	$y := Y_1 \times Y_2,$\\
	H.0.	$X \cdot \theta = 0,$\\
	then\\
	C.0.	$Y_0 \cdot y = 0.$\\[10pt]
	P.0.	$Y_0 = (),$\\
	P.1.	$y = [((m_1+m_2)(X_1+X_2)-m_0X_0)X_1X_2,
			((m_2+m_0)(X_2+X_0)-m_1X_1)X_2X_0,
			((m_0+m_1)(X_0+X_1)-m_2X_2)X_0X_1].$

	\sssec{Exercise.}
	The excribed circles.
\begin{verbatim}
	Let iii[i] := radical axis(iota[i+1],iota[i-1]),
	then
	    iii[i]\cdot En = 0.
	Ex g277, iii[] = [647,435,847],
	Ex4.0, iii[] = [873,651,964],
	Ex5.0, iii[] = [723,837,965].

	The conic of Neuberg.  (Mathesis, Ser.2, Vol.6, p. 95).
	P40.2.	`Neuberg: m_0m_1m_2(m_0 X_0^2 + m_1 X_1^2 + m_2 X_2^2)
			+ s11( m_0(m_1+m_2)X_1X_2 + m_1(m_2+m_0)X_2X_0
				+ m_2(m_0+m_1)X_0X_1) = 0.
		(Bastin, Mathesis, p.97)
\end{verbatim}

	\sssec{Exercise.}
	(Neuberg, see Casey, no 80,81,82)\\
	The barycenter of the triangle $\{Ste,BRa,Abr\}$ (see D16.5,14,16) is
	$M$.\
	$BRa \times Abr = [m_1m_2(m_2-m_0)(m_0-m_1),m_2m_0(m_0-m_1)(m_1-m_2),
		m_0m_1(m_1-m_2)(m_2-m_0)].$

	\sssec{Exercise.}
	Some points on the circumcircle.\\
	Construct the points $\ov{M}iqm_i$ on $\theta$ and $ac_i$ distinct from
	$A_i$ and the points $Miq\ov{m}_i$ on $\theta$ and on $\ov{a}c_i$
	distinct from $A_i.$

	\sssec{Answer to}
	(partial).\\
	$Miq\ov{m}[0] = (m_0(m_1+m_2),m_1(m_1-m_2),m_2(m_2-m_1)),$\\
	$\ov{M}iqm[0] = (m_1+m_2,m_2-m_1,m_1-m_2),$

	\sssec{Example.}
	$p = 29,$ $A_i = (30,1,0),$ $M = (60),$ $\ov{M} = (215),$\\
	$Miq\ov{m}_i = (545,512,699),$ $\ov{M}iqm_i = (115,261,855).$

	\sssec{Exercise.}
	The point of Miquel.\\
	Given an arbitrary line which does not pass through the vertices and is
	neither the ideal or coideal line, $q = [q_0,q_1,q_2],$
	Let $Q_i := q \times a_i.$  Determine\\
	    $\mu iq_i := conic(I',I'',A_i,Q_{i+1},Q_{i-1}),\\
		\ov{\mu}iq_i := conic(\ov{I}',\ov{I}'',A_i,Q_{i+1},Q_{i-1}),$\\
	Construct the point $Miq$ which is on $\mu iq_i$ and $\theta,$\\
	the point $\ov{M}iq$ which is on $\ov{\mu}iq_i$ and $\theta,$\\
	the circle $\mu q$ of $Miquel$ which passes through the center of
	$\mu iq_i$\\
	and the cocircle $\ov{\mu}q$ which passes through the center of
	$\ov{\mu}iq_i.$\\
	The following special cases are of interest:\\
	$q = \ov{m},$ in which case $Miq = \ov{M}iq,$ which we denote
	$Miq\ov{i},$\\
	$q = e,$ we then denote the point and copoint of Miquel by $Miqe$ and
	$\ov{M}iqe,$\\
	$q = \ov{m}_i,$ giving $Miq\ov{m}_i$ of Exercise <...> above,\\
	$q = m_i,$ giving $\ov{M}iqm_i$ of Exercise <...> above.

	\sssec{Answer to}
	(Partial)\\
	$miq_0 = (q0-q1)(q0-q2)(m_0(m_1+m_2) X_1 X_2 + ... )\\
		+ (X_0 + X_1 + X_2) ( m_2(m_0+m_1)q1(q0-q2) X_1
					+ m_1(m_2+m_0)q2(q0-q1) ) = 0.$\\
	$\ov{m}iq_0 = (m_0q0-m_1q1)(m_0q0-m_2q2)(m_0(m_1+m_2)X_0X_2+...)
		+ (m_1m_2 X_0 + m_2m_0 X_1 + m_0m_1 X_2)\\
	\hth ((m_0+m_1)q1(m_0q0-m_2q2) X_1 + (m_2+m_0)q1(m_0q0-m_1q1) X_2)
		= 0.$\\
	$Miq = (m_0(m_1+m_2)q1q2(q0-q2)(q1-q0),m_1(m_2+m_0)q2q0(q1-q0)(q2-q1),
		m_2(m_0+m_1)q0q1(q2-q1)(q0-q2)),$\\
	$\ov{M}iq = ((m_1+m_2)q1q2(m_0q0-m_1q1)(m_0q0-m_2q2),
			(m_2+m_0)q2q0(m_1q1-m_2q2)(m_1q1-m_0q0)),\\
			(m_0+m_1)q0q1(m_2q2-m_0q0)(m_2q2-m_1q1)),$\\
	$Miq\ov{i} = ((m_1+m_2)(m_2-m_0)(m_0-m_1),(m_2+m_0)(m_0-m_1)(m_1-m_2),
		(m_0+m_1)(m_1-m_2)(m_2-m_0)),$\\
	$Miqe = (m_0(m_1+m_2)(m_2-m_0)(m_0-m_1)(m_0-2m_1+m_2)(m_0+m_1-2m_2),
		m_1(m_2+m_0)(m_0-m_1)(m_1-m_2)(m_1-0m_2+m_0)(m_1+m_2-0m_0),\\
	\hth	m_2(m_0+m_1)(m_1-m_2)(m_2-m_0)(m_2-1m_0+m_1)(m_2+m_0-1m_1))$,\\
	$\ov{M}iqe = ((m_1+m_2)(m_2-m_0)(m_0-m_1)(2m_2m_0-m_1(m_2+m_0))
		(2m_0m_1-m_2(m_0+m_1)),\\
	\hth	     (m_2+m_0)(m_0-m_1)(m_1-m_2)(2m_0m_1-m_2(m_0+m_1))
			(2m_1m_2-m_0(m_1+m_2)),\\
	\hth	     (m_0+m_1)(m_1-m_2)(m_2-m_0)(2m_1m_2-m_0(m_1+m_2))
			(2m_2m_0-m_1(m_2+m_0))).$

	\sssec{Exercise.}
	(Sondat, See Mathesis, Ser. 2, Vol.6, pp. 81-83)\\
	Let $B_0 \cdot \mu iq_0 = 0,$ let\\
	D.0.	$b_1 := B_0 \times Q_1,$ $b_2 := Q_2 \times B_0,$\\
	D.1.	$B_1 := muiq_1 \times b_2 - Q_2,$
		$B_2 := muiq_2 \times b_1 - Q_1,$\\
	D.2.	$b_0 := B_1 \times B_2,$\\
	D.3.	$ab_i := A_i \times B_i,$\\
	D.4.	$S := ab_1 \times ab_2,$\\
	Let\\
	$S_1\cdot muiq_1 = 0,$ $S_2\cdot \mu iq_2 = 0,$\\
	D.5.	$sa_1 := S_1 \times A_1,$ $sa_2 := S_2 \times A_2,$\\
	D.6.	$T := sa_1 \times sa_2,$\\
	D.7.	$sa_0 := T \times A_0,$\\
	D.8.	$S_0 := muiq_0 \times sa_0 - A_0,$\\
	D.9.	$\sigma := conic(S_i,S,T),$\\
	C.0.	$Q[0] \cdot b_0 = 0.$\\
	C.1.	$S \cdot ab[0] = 0.$\\
	C.2.	$S \cdot \theta = 0.$\\
	C.3.	$\sigma is a circle,$\\
	C.4.	$T = \ov{M}$ ==> center($\sigma)\cdot q = 0.$\\
	<...> double check the above.

	\sssec{Exercise.}
	Construct the point common to the circumcircle and the
	circle through  $A_i,$ $\ov{M}_i+1$ and $\ov{M}_i-1.$  (See also the
	transformation of Hamilton.

	\sssec{Partial Answer to}
	a) The center is $E_i,$\\
	$x_i := E_i \times I\ov{mm}_i,$\\
	$y_i := O \times Im_i,$\\
	$Z_i := x_i \times y_i,$\\
	$z_i := A_i \times Z_i,$\\
	$z_i$ corresponds to the perpendicular to $O \times E_i,$ hence
	contains the desired intersection $HH_i.$\\
	$x_0 = [m_1m_2(s_1+m_0),-m_2(2m_0s_1+m_1(m_1+m_2)),
		-m_1(2m_0^2+3m_0m_1+m_2s_1)],$\\
	$y_0 = [-(s_1+m_0),m_1+m_2,m_1+m_2],$\\
	$Z_0 = ((m_1^2-m_2^2)(m_0+m_2),m_1(s1+m_0)(m_0+2m_2),-m_2(s1+m_0)s1),$\
	$z_0 = [0,m_2s1,m_1(m_0+2m_2)],$\\
	The rest of the construction is that of Pascal:\\
	$aa_i := A_{i+1} \times AA_{i-1},$\\
	$ZZ_i := aa_i \times z_i,$\\
	$zz_i := Im_{i+1} \times ZZ_{i-1},$\\
	$Y_i := a_{i+1} \times zz_{i-1},$\\
	$yy_i := Y_{i+1} \times AA_{i-1},$\\
	$Miq\ov{m}_i := yy_{i+1} \times z_i.$

	\sssec{Exercise.}
	\enumb
	\item Complete a section on the conic of Nagel, with\\
	$\nu := conic(JJ_0,JJ_1,JJ_2,Na_1,Na_2),$\\
	$\nu_i := conic(\ldots),$
	\item Give a construction for the other intersection of the conic with
	$Ja_i$ and $Na_i.$
	\item Give a construction for the center of the conic.
	\item Are there some other points on this conic which have already been
	constructed or that you can construct?
	\enume

	\sssec{Partial Answer to}
	$\mu = j_0j_1j_2(X_0^2+X_1^2+X_2^2) - j_0(j_1^2+j_2^2)X_1 X_2
		 -\ldots = 0$, 	<...> not checked\\
	other intersection with $Ja_0 = (p_2+j_0^2),j_0j_1,j_2j_0),$\\
	other intersection with $Na_0 = (p_{22}+j_1^2j_2^2,j_0j_1j_2^2,
		j_0j_1^2j_2),$\\
	center (\ldots?)$ (j_0(j_1^4(j_2-j_0)+j_2^4(j_1-j_0)+j_1^2j_2^2
		(2j_0+3j_1+3j_2)+j_0^2j_1^2(j_1+3j_2)+j_2^2j_0^2(3j_1+j_2)),
		\ldots)$

	\sssec{Exercise.}
	The circles of Lemoine-Tucker.\\
	D.0.	$X_i := K - x A_i,$ $x$ is some integer,\\
	D.1.	$x_i := X_{i+1} \times X_{i-1},$\\
	D.2.	$XX_i := x_{i+1} \times a_{i-1},$\\
	D.3.	$X\ov{X}_i := x_{i-1} \times a_{i+1},$\\
	D.4.	$\xi := conic(XX_0,XX_1,XX_2,X\ov{X}_1,X\ov{X}_2),$\\
	then\\
	C.0.	$x_i\cdot Im_i = 0.$\\
	C.1.	$\xi \cdot X\ov{X}_0 = 0.$\\
	C.2.	$\xi$ is a circle.\\
		D.0., can be replaced by a construction which start with
		a point
		$X_0$ on the symmedian $at_0,$ the parallel through $X_0$ to
		the side $a_2$ or $a_1$ intersect the symmedians $at_1$ or
		$at_2$ at $X_1$ or $X_2.$\\
	Proof.\\
	P.0.	$X_0 = (m_0(m_1+m_2)+x,m_1(m_2+m_0),m_2(m_0+m_1)),$\\
	P.1.	$x_0 = [m_0(m_1+m_2),m_0(m_1+m_2),-s_{11}-m_0m_1+x],$\\
	P.2.	$XX_0 = (s_{11}+m_2m_0-x,m_1(m_2+m_0),0),$\\
	P.3.	$X\ov{X}_0 = (s_{11}+m_0m_1-x,0,m_2(m_0+m_1)),$\\
	P.4.	$(2s_{11}-x)^2 \theta - (m) \bigx (u) = 0,$ with\\
		$u_0 = m_1m_2(m_2+m_0)(m_0+m_1)(s_{11}+m_1m_2-x),$\ldots.

	\sssec{Comment.}
	The following are special cases:\\
	$x = 0$ gives the first circle of Lemoine lambda1,\\
	$x = s_{11}$ gives the second circle of Lemoine lambda2,\\
	$x = \ldots$ gives the circle of Taylor,\\
	$x = 2s_{11}$ gives the degenerate circle $(i) \bigx (i),$\\
	$x = \frac{1}{0}$ gives $\theta.$

	\sssec{Exercise.}
	\ssec{The harmonic polygons. [Casey]}

	\sssec{Definition.}
	Given a conic $\theta$ and a point $K$ not on the conic, an
	inscribed polygon $A_i$, $i = 0, ... n-1$ is a {\em harmonic polygon} if
	$(A_{i-1},A_i,A_{i+1},A'_i)$ is harmonic for all $i$, where\\
	$ka_i := K \times A_i,$\\
	$A'_i := \theta \times ka_i - A_i,$\\
	$k := polar(K),$\\
	$B_i := polar(ka_i),$\\
	$K$ is called the point of Lemoine of the polygon,\\
	$k$ is called the line of Lemoine.

	\sssec{Theorem.}
	{\em If $A_i,$ $i = 0$ to $n-1$ is a harmonic polygon then
	$A'_i,$ $i = 0$ to $n-1$ is a harmonic polygon.}

	\sssec{Construction.}
	Given $K,$ $A_0,$ $A_1,$\\
	construct $k,$ $ka_0,$ $B_0,$ $ka_1,$ $B_1,$\\
	\un{for} $i = 1$ \un{to} $n-1$\\
	\un{begin}\\
	\hth$ka_i := polar A_i,$
	\hth$B_i := ka_i \times k,$
	\hth$c_i := A_{i-1} \times B_i,$
	\hth$A_{i+1} := \theta \times c_i - A_{i-1},$
	\un{end}

	\sssec{Construction.}
	Details.\\
	H.0.	$x(A_1-B_1) = y(A_0-B_0) + B_0A_1-A_0B_1$\\
	H.1.	$x^2 + y^2 = 1,$\\
	then\\
	C.0.	$y^2( (A_0-B_0)^2 + (A_1-B_1)^2) + 2y (A_0-B_0)(B_0A_1-A_0B_1)
			+ (B_0A_1-A_0B_1)^2 - (A_1-B_1)^2 = 0,$\\
	C.1.	$y = -2(A_0B_0)(B_0A_1-A_0B_1)/((A_0-B_0)^2+(A_1-B_1)^2)-A_1.$\\
	CHECK THE ABOVE <...>

	\sssec{Example.}
	For $p = 31,$\\
	let $K = (0,-8,1),$ $A_0 = (1,0,1),$ $A_1 = (),$
	then $A_i = (1,0,1),$

	\sssec{Exercise.}
	Complete a section on polars of the vertices with respect
	to the conic of Brianchon-Poncelet.
	\enumb
	\item Give an explicit construction for the tangents to gamma at
		the mid-points and at the feet.
	\item Give an explicit construction for the polar $pp_i$ of $A_i$ with
		respect to $\gamma.$
	\item Verify that the intersections $PP_i := pp_i\times A_i$
		are collinear on $pp.$
	\enume
	This result can be used as the starting point for special results in
	the geometry of the tetrahedron.  (An other approach is suggested by the
	theorem <...>).  The lines $a_i$ and $pp$ correspond to the ideal lines
	in the four faces of a tetrahedron whose opposite vertices are
	perpendicular.  The tetrahedron so obtained have the additional
	properties that $A_i \times M_i$ are concurrent as well as
	$A_i \times \ov{M}_i.$

	\sssec{Comment.}
	An other model of projective geometry within projective
	geometry is suggested by the following.\\
	Associate to the point $(X_0,X_1,X_2),$ the point
	$(X_1X_2,X_2X_0,X_0X_1),$\\
	associate to the line $[a_0,a_1,a_2],$ the conic\\
	\hth$a_0 X_1 X_2 + a_1 X_2 X_0 + a_2 X_0 X_1 = _0.$\\
	The ideal is the conic\\
	\hth$	X_1X_2 + X_2X_0 + X_0X_1 = _0,$\\
	and the coideal is\\
	\hth$	m_0X_1X_2 + m_1X_2X_0 + m_2X_0X_1 = 0.$\\
	Some care has to be exercised because if, for instance, two of the
	coordinates $X_0,$ $X_1,$ $X_2$ are 0, the image is not defined.\\
	In the following definition "Point'" and "Line" is used for the new
	objects which have the properties of "point" and "line" defined above.

	\sssec{Definition.}
	Given a triangle ${(A_0,A_1,A_2),(a_0,a_1,a_2)},$\\
	the Points are
	\begin{itemize}
	\item the points not on the sides of the triangle,
	\item the line through the vertices, (including $a_0,a_1,a_2$),\\[10pt]
	the Lines are
	\item the conics through A0, A1 and A2, including the degenerate conics
	  which consist of one side and a line through the opposite vertex.
\\[10pt]
	A Point is on a Line if
	\item \ldots
	\end{itemize}
	If two of the points are the isotropic points, the lines become the
	circles passing through a given point.
	A large number of properties of circles as well as properties of
	projective geometry can be obtained by pursuing this approach.
	In particular a study of the quartics which are associated to the
	circles is of interest.\\
	An early reference on circular triangles is by Miquel, J. de
	Liouville, Vol. 9, 1844, p. 24.

	Special cases.  2 Special cases are of interest.

	\sssec{Notation.}
	$P :== p_1 \times p_2$ does not denote an actual construction, but
	    a construction in which $p_1$ or $p_2$ are assumed to be known.\\
	"==" was suggested by the mode of drawing using dashed lines
	    rather than continuous ones.\\
	In the example below, $D$ is not known, hence we can not construct
	$A \times D.$

	The following problem is of interest.

	\sssec{Exercise.}
	Given 2 conics with 3 points in common, determine by a linear
	construction the fourth point on both conics.\\
	One solution is the following.\\
	Let $A,B,C$ be the known points and $D$ be the unknown point.
	Let $E$ and $F$ be on the first conic $\gamma,$ $U$ and $V$ on the
	second conic $\gamma'.$  Determine first by the Pascal's construction\\
	point Pascal$(U,V,C,B,A,E;E'),$ and point Pascal$(A,B,C,U,V,E;E'),$
	$E'$ on $\gamma'$ and $A\times E$, $F'$ on $\gamma'$ and $B\times F$,\\
	let $K :== (D \times A) \times (C \times B),$
	$L := (A \times E) \times (B \times F),$
	$M :== (E \times C) \times (F \times D),$
	$M' :== (E' \times C) \times (F' \times D),$\\
	then Pascal($D,A,E,C,B,F;K,L,M),$ and
	Pascal($D,A,E',C,B,F';K,L,M').$ This implies
	incidence$(L,M,M').$\\
	Using Desargues$^{-1}(\langle L,M',M\rangle,\{C,E,E'\},\{D,F,F'\};G),$
	it follows that $D$ is incident to $c\times G$, with\\
	$G := (E \times F) \times (E' \times F'),$\\
	the triangles $\{C,E,E'\}$ and $\{D,F,F'\}$ being perspective.\\
	$D$ follows from point Pascal$(B,A,E,F,C,G;D).$

	\sssec{Construction.}
	The complete construction is the following:\\
	$P_0 := (U \times V) \times (B \times A),$
	$P_1 := (V \times C) \times (A \times E),$
	$P_2 := (C \times B) \times (P_0 \times P1),$
	$E' := (A \times E) \times (P_2 \times U),$\\
	$P_1' := (V \times C) \times (B \times F),$
	$P_2' := (C \times A) \times (P_0 \times P_1'),$
	$F' := (C \times A') \times (P_2' \times U),$\\
	$G := (E \times F) \times (E' \times F'),$\\
	$Q_0 := (B \times A) \times (F \times C),$
	$Q_1 := (A \times E) \times (C \times G),$
	$Q_2 := (E \times F) \times (Q_0 \times Q_1),$
	$D := (C \times G) \times (Q_2 \times B).$

	An other solution is the following

	\sssec{Theorem.}
	{\em Let $t_A$ and $t_B$ be the tangents to the first conic at $A$
	and at $B,$\\
	let $t'_A$ and $t'_B$ the tangents to the second conic at $A$ and $B,$\\
	$O_1 := t_A \times t_B,$
	$O'_1 := t'_A \times t'_B,$
	$oo' := O_1 \times O'_1,$
	$ab := A \times B,$
	$BA' := t_B \times t'_A,$
	$AB' := t_A \times t'_B,$
	$ab' := BA' \times AB',$
	$E := ab \times ab',$
	$cd := C \times E,$
	$bc := B \times C,$
	$F := bc \times oo',$
	$ad := A \times F,$
	$D := ad \times cd,$\\
	then\\
	$D$ is on conic$(A,t_A,B,t_B,C)$ and conic$(A,t'_A,B,t'_B,C)$}.

	Proof:
	Assume that $A$ and $B$ are the isotropic points then the 2 conics
	are circles.  $O$ and $O'$ are their centers.  $cd \perp oo'$.
	Therefore $(A,B,D,ab \times oo')$ is a harmonic quatern.
	$bc$ and $ad$ meet on $oo'.$  This can be checked using
	$A = (1,i,0),$ $B = (1,-i,0),$ $C = (0,1,1),$ $D = (0,-1,1),$
	$oo' = [1,0,0],$ $bc = [-i,-1,1],$ $ad = [i,-1,-1],$ $F = (0,1,1).$



	\ssec{Cubics.}
\setcounter{subsubsection}{-1}
	\sssec{Introduction.}
	Cubics have extensively studied by Newton, MacLaurin, Gergonne, Plucker,
	Salmon, \ldots.\\
	I will give here a few properties, many of which generalize to higher 
	degree curves, most of them  taken from Salmon, 1979, sections 29 to
	31 and 148 to 159?:

	\sssec{Theorem.}
	{\em All cubics which pass trough 8 fixed points pass also through a
	ninth.}

	\sssec{Definition.}
	If 9 points are on a one parameter family of non degenerate cubics,
	we say that they form a {\em cubic configuration.} This configuration is
	not confined.

	\sssec{Theorem. [MacLaurin]}
	{\em Let $A_0,$ to $A_7$ be 8 points of a cubic, such that
	$A_0,$ $A_1,$ $A_2,$ $A_3,$ $A_4,$ $A_5,$ are on a conic $\alpha$ and
	$A_0,$ $A_1,$ $A_2,$ $A_3,$ $A_6,$ $A_7,$ are on a conic $\beta$ then
	$(A_4\times A_5)\times(A_6\times A_7)$ is on the cubic and the 9
	points form a cubic configuration.}

	Proof: This follows when the preceding Theorem is applied to the
	degenerate cubics consisting of the conic $\alpha$ and the line
	$A_6\times A_7$ and the conic $\beta$ and the line $A_4\times A_5$.

	\sssec{Corollary.}\label{sec-c9}
	{\em If 2 lines meet a cubic at points $B_{0,i}$ and $B_{1,i}$
	then the 3 points $B_{2,i}$ on the cubic and on $B_{0,i} \times
	B_{1,i}$ are collinear,}\\
	or equivalently\\
	{\em If 6 points $B_j$, $j = 0$ to 5, are on a cubic and 2 of the points
	$C_i := (B_i\times B_{i+1})\times (B_{i+3}\times B_{i+4})$ are on the
	cubic then the third point is on the cubic and the 9 points form a
	cubic configuration.}

	Proof: This follows when the preceding Theorem is applied to the
	degenerate conics $\alpha$ through $B_{0,i}$ and $B_{1,i}$ and
	$beta$ through $B_{i,0}$ and $B_{i,1}$.

	The alternate form corollary gives Pappus' Theorem when the cubic
	degenerates into 3 lines.

	\sssec{Notation.}
	I will write $C9(B_0,B_1,B_2,B_3,B_4,B_5;C_0,C_1,C_2).$

	\sssec{Theorem. [Salmon]}\label{sec-tsal}
	{\em The 3 parameter family of cubics through the 6 points $A_i$
	and $B_i$, which are not on a conic is}\\
	\hth$\Sigma_{i=0,1,2}s_i(A_i\times B_i)\bigx(A_{i+1}\times A_{i-1})
		\bigx(B_{i+1}\times B_{i-1})\\
	= (A_{i+1}\times B_{i-1})\bigx(A_{i-1}\times B_i)
		\bigx(A_i\times A_{i+1}).$

	Proof: It is easy to verify that each of the points is on each
	of the 4 degenerate cubics and that these are independent.\\
	There are many alternate forms possible, I have chosen the above one
	which displays a useful symmetry property.

	\sssec{Definition.}
	The {\em tangential point} of a point $C$ on a cubic is the third
	intersection of the tangent at $C$ with the cubic.

	\sssec{Corollary.}
	{\em If 3 points of a cubic are on a line $a$, their tangential points
	are on a line $s$.}

	This follows from the degenerate case $B_{0,i} = B_{1,i}$.

	\sssec{Definition.}
	The line $s$ is called the {\em satellite} of the line $a$.

	\sssec{Notation.} Given 2 points $A$ and $B$ on a cubic, the third point
	on the cubic and the line $A\times B$ is denoted $A\star B$.

	\sssec{Theorem.}
	{\em Given 6 lines $a_i$ and $b_i$ and their 9 intersections}\\
	\hth$i\:j := a_i\times b_j$,
	\enumb
	\item {\em These 9 points form a cubic configuration.}
	\item {\em If $C_0$ is a point on $11\times 22,$ the cubic of the family
	through the points $i\:j$ and $C_0$ are such that if we define
	the following points,}\\
	$C_i := i+1,i+1\star i-1,i-1,$ $D_i := i+1,i-1\star i-1,i+1,$\\
	$E_i := i,i+1\star i-1,i,$ $\ov{E}_i := i+1,i\star i,i-1,$\\
	$F_i := i,i\star i+1,i-1,$ $\ov{F}_i := i,i\star i-1,i+1,$\\
	$Cc_i := C_{i+1}\star C_{i-1},$ $Cd_i := C_{i+1}\star D_{i-1},$
	$Dc_i := D_{i+1}\star C_{i-1},$\\
	$CF_i := C_i\star F_i,$ $Cf_i := C_{i+1}\star F_{i-1},$
	$Fc_i := F_{i+1}\star C_{i-1},$\\
	$CF_i := C_i\star F_i,$ $Cf_i := C_{i+1}\star F_{i-1},$
	$Fc_i := F_{i+1}\star C_{i-1},$\\
	$C\ov{F}_i := C_i\star \ov{F}_i,$
	$C\ov{f}_i := C_{i+1}\star \ov{F}_{i-1},$
	$\ov{F}c_i := \ov{F}_{i+1}\star C_{i-1},$\\
	$DE_i := D_i\star E_i,$ $De_i := D_{i+1}\star E_{i-1},$
	$Ed_i := E_{i+1}\star D_{i-1},$\\
	$D\ov{E}_i := D_i\star \ov{E}_i,$
	$D\ov{e}_i := D_{i+1}\star \ov{E}_{i-1},$
	$\ov{E}d_i := \ov{E}_{i+1}\star D_{i-1},$\\
	$DF_i := D_i\star F_i,$ $D\ov{F}_i := D_i\star \ov{F}_i,$\\
	$Ee_i := E_{i+1}\star E_{i-1},$ $EF_i := E_i\star F_i,$
	$Ef_i := E_{i+1}\star F_{i-1},$\\
	$Fe_i := F_{i+1}\star E_{i-1},$ $E\ov{f}_i := E_{i+1}\star\ov{F}_{i-1},$
	$\ov{F}e_i := \ov{F}_{i+1}\star E_{i-1},$\\
	$\ov{Ee}_i := \ov{E}_{i+1}\star \ov{E}_{i-1},$
	$\ov{EF}_i := \ov{E}_i\star \ov{F}_i,$
	$\ov{E}f_i := \ov{E}_{i+1}\star F_{i-1},$\\
	$F\ov{e}_i := F_{i+1}\star \ov{E}_{i-1},$
	$\ov{Ef}_i := \ov{E}_{i+1}\star\ov{F}_{i-1},$
	$\ov{Fe}_i := \ov{F}_{i+1}\star \ov{E}_{i-1},$\\
	$F\ov{F}_i := F_i\star \ov{F}_i,$\\
	$C'_i := C_i\star C_i$, $D'_i := D_i\star D_i$,\\
	$E'_i := E_i\star E_i$, $\ov{E}'_i := \ov{E}_i\star \ov{E}_i$,
	$F'_i := F_i\star F_i$,	$\ov{F}'_i := \ov{F}_i\star \ov{F}_i$,\\
	$K_i := C_i\star i+1,i-1$,
	$\ov{K}_i := C_i\star i-1,i+1$,\\
	$L_i := D_{i+1}\star i-1,i-1$, $\ov{L}_i := D_{i-1}\star i+1,i+1$,\\
	$M_i := E_i\star i,i$, $\ov{M}_i := \ov{E}_i\star i,i$,
	$N_i := E_i\star i-1,i+1$, $\ov{N}_i := \ov{E}_i\star i+1,i-1$,\\
	$P_i := F_{i+1}\star i+1,i$, $\ov{P}_i := \ov{F}_{i+1}\star i,i+1$,\\
	$Q_i := F_{i-1}\star i,i-1$, $\ov{Q}_i := \ov{F}_{i-1}\star i-1,i$,
	\enume
\newpage
{\footnotesize
	{\em We have the following table for the operation $\star$ between
	points on the cubic:}\\
	\noindent$\begin{array}{lclllclllclllclllclll}
	\star&\smv&00&11&22&\smv&12&20&01&\smv&21&02&10&\smv&
	C_0&C_1&C_2&\smv&D_0&D_1&D_2\\
	\hline
	C_0&\smv&D_0&22&11&\smv&K_0&C\ov{F}_2&C\ov{F}_1
	&\smv&\ov{K}_0&CF_2&CF_1&\smv&C'_0&Cc_2&Cc_1&\smv&00&Cd_2&Dc_1\\
	C_1&\smv&22&D_1&00&\smv&C\ov{F}_2&K_1&C\ov{F}_0
	&\smv&CF_2&\ov{K}_1&CF_0&\smv&Cc_2&C'_1&Cc_0&\smv&Dc_2&11&Cd_0\\
	C_2&\smv&11&00&D_2&\smv&C\ov{F}_1&C\ov{F}_0&K_2
	&\smv&CF_1&CF_0&\ov{K}_2&\smv&Cc_1&Cc_0&C'_2&\smv&Cd_1&Dc_0&22\\
	\hline
	D_0&\smv&C_0&L_2&\ov{L}_1&\smv&21&F\ov{e}_2&\ov{E}f_1
	&\smv&12&\ov{F}e_2&E\ov{f}_1&\smv&00&Dc_2&Cd_1&\smv&D'_0&22'&11'\\
	D_1&\smv&\ov{L}_2&C_1&L_0&\smv&\ov{E}f_2&02&F\ov{e}_0
	&\smv&E\ov{f}_2&20&\ov{F}e_0&\smv&Cd_2&11&Dc_1&\smv&22'&D'_1&00'\\
	D_2&\smv&L_1&\ov{L}_0&C_2&\smv&F\ov{e}_1&\ov{E}f_0&10
	&\smv&\ov{F}e_1&E\ov{f}_0&01&\smv&Dc_1&Cd_0&22&\smv&11'&00'&D'_2\\
	\hline
	E_0&\smv&M_0&DE_2&DE_1&\smv&\ov{F}_0&01&20&\smv&N_0&EF_2&EF_1
	&\smv&21'&\ov{E}_2&\ov{E}_1&\smv&DE_0&Ed_2&De_1\\
	E_1&\smv&DE_2&M_1&DE_0&\smv&01&\ov{F}_1&12&\smv&EF_2&N_1&EF_0
	&\smv&\ov{E}_2&02'&\ov{E}_0&\smv&De_2&DE_1&Ed_0\\
	E_2&\smv&DE_1&DE_0&M_2&\smv&20&12&\ov{F}_2&\smv&EF_1&EF_0&N_2
	&\smv&\ov{E}_1&\ov{E}_0&10'&\smv&Ed_1&De_0&DE_2\\
	\hline
	\ov{E}_0&\smv&\ov{M}_0&D\ov{E}_2&D\ov{E}_1
	&\smv&\ov{N}_0&\ov{EF}_2&\ov{EF}_1&\smv&F_0&10&02
	&\smv&12'&E_2&E_1&\smv&D\ov{E}_0&\ov{E}d_2&D\ov{e}_1\\
	\ov{E}_1&\smv&D\ov{E}_2&\ov{M}_1&D\ov{E}_0
	&\smv&\ov{EF}_2&\ov{N}_1&\ov{EF}_0&\smv&10&F_1&21
	&\smv&E_2&20'&E_0&\smv&D\ov{e}_2&D\ov{E}_1&\ov{E}d_0\\
	\ov{E}_2&\smv&D\ov{E}_1&D\ov{E}_0&\ov{M}_2
	&\smv&\ov{EF}_1&\ov{EF}_0&\ov{N}_2&\smv&02&21&F_2
	&\smv&E_1&E_0&01'&\smv&\ov{E}d_1&D\ov{e}_0&D\ov{E}_2\\
	\hline
	F_0&\smv&12&\ov{F}e_2&E\ov{f}_1&\smv&00&Dc_2&Cd_1
	&\smv&\ov{E}_0&P_2&Q_1
	&\smv&CF_0&Fc_2&Cf_1&\smv&DF_0&\ov{F}_2&\ov{F}_1\\
	F_1&\smv&E\ov{f}_2&20&\ov{F}e_0&\smv&Cd_2&11&Dc_0
	&\smv&Q_2&\ov{E}_1&P_0
	&\smv&Cf_2&CF_1&Fc_0&\smv&\ov{F}_2&DF_1&\ov{F}_0\\
	F_2&\smv&\ov{F}e_1&E\ov{f}_0&01&\smv&Dc_1&Cd_0&22
	&\smv&P_1&Q_0&\ov{E}_2
	&\smv&Fc_1&Cf_0&CF_2&\smv&\ov{F}_1&\ov{F}_0&DF_2\\
	\hline
	\ov{F}_0&\smv&21&F\ov{e}_2&\ov{E}f_1&\smv&E_0&\ov{P}_2&\ov{Q}_1
	&\smv&00&Dc_2&Cd_1
	&\smv&C\ov{F}_0&\ov{F}c_2&C\ov{f}_1&\smv&D\ov{F}_0&F_2&F_1\\
	\ov{F}_1&\smv&\ov{E}f_2&02&F\ov{e}_0&\smv&\ov{Q}_2&E_1&\ov{P}_0
	&\smv&Cd_2&11&Dc_0
	&\smv&C\ov{f}_2&C\ov{F}_1&\ov{F}c_0&\smv&F_2&D\ov{F}_1&F_0\\
	\ov{F}_2&\smv&F\ov{e}_1&\ov{E}f_0&10&\smv&\ov{P}_1&\ov{Q}_0&E_2
	&\smv&Dc_1&Cd_0&22
	&\smv&\ov{F}c_1&C\ov{f}_0&C\ov{F}_2&\smv&F_1&F_0&D\ov{F}_2\\
	\end{array}$
	\noindent$\begin{array}{lclllclllclllclll}
	\star&\smv&E_0&E_1&E_2 &\smv&\ov{E}_0&\ov{E}_1&\ov{E}_2
	&\smv&F_0&F_1&F_2 &\smv&\ov{F}_0&\ov{F}_1&\ov{F}_2\\
	\hline
	E_0&\smv&E'_0&Ee_2&Ee_1&\smv&00'&C_2&C_1
	&\smv&EF_0&Ef_2&Fe_1&\smv&12&E\ov{f}_2&\ov{F}e_1\\
	E_1&\smv&Ee_2&E'_1&Ee_0&\smv&C_2&11'&C_0
	&\smv&Fe_2&EF_1&Ef_0&\smv&\ov{F}e_2&20&E\ov{f}_0\\
	E_2&\smv&Ee_1&Ee_0&E'_2&\smv&C_1&C_0&22'
	&\smv&Ef_1&Fe_0&\ov{E}F_2&\smv&E\ov{f}_1&\ov{F}e_0&01\\
	\hline
	\ov{E}_0&\smv&00'&C_2&C_1&\smv&\ov{E}'_0&\ov{Ee}_2&\ov{Ee}_1
	&\smv&21&\ov{E}f_2&F\ov{e}_1&\smv&\ov{EF}_0&\ov{Ef}_2&\ov{Fe}_1\\
	\ov{E}_1&\smv&C_2&11'&C_0&\smv&\ov{Ee}_2&\ov{E}'_1&\ov{Ee}_0
	&\smv&F\ov{e}_2&02&\ov{E}f_0&\smv&\ov{Fe}_2&\ov{EF}_1&\ov{Ef}_0\\
	\ov{E}_2&\smv&C_1&C_0&22'&\smv&\ov{Ee}_1&\ov{Ee}_0&\ov{E}'_2
	&\smv&\ov{E}f_1&F\ov{e}_0&10&\smv&\ov{Ef}_1&\ov{Fe}_0&\ov{EF}_2\\
	\hline
	F_0&\smv&EF_0&Fe_2&Ef_1&\smv&21&F\ov{e}_2&\ov{E}f_1
	&\smv&F'_0&10'&02'&\smv&FF_0&D_2&D_1\\
	F_1&\smv&Fe_2&EF_1&Ef_0&\smv&\ov{E}f_2&02&F\ov{e}_0
	&\smv&10'&F'_1&21'&\smv&D_2&FF_1&D_0\\
	F_2&\smv&Fe_1&Ef_0&\ov{E}F_2&\smv&\ov{F}e_1&E\ov{f}_0&10
	&\smv&02'&21'&F'_2&\smv&D_1&D_0&FF_2\\
	\hline
	\ov{F}_0&\smv&12&\ov{F}e_2&E\ov{f}_1
	&\smv&\ov{EF}_0&\ov{Fe}_2&\ov{Ef}_1
	&\smv&FF_0&D_2&D_1&\smv&\\ov{F}'_0&01'&20'\\
	\ov{F}_1&\smv&E\ov{f}_2&20&\ov{F}e_0
	&\smv&\ov{Ef}_2&\ov{EF}_1&\ov{Fe}_0
	&\smv&D_2&FF_1&D_0&\smv&01'&\ov{F}'_1&12'\\
	\ov{F}_2&\smv&\ov{F}e_1&E\ov{f}_0&01
	&\smv&\ov{Fe}_1&\ov{Ef}_0&\ov{EF}_2
	&\smv&D_1&D_0&FF_2&\smv&20'&12'&\ov{F}'_2\\
	\end{array}$
}
\newpage

	Proof:\\
{\footnotesize
	\noindent$\begin{array}{lclllclllclllc}
	\alpha_0&D9(C_0&D_0&21&20&10&11&;&00&12&22)\\
	\rho\alpha_0&D9(\ov{E}_2&F_2&01&00&20&21&;&10&22&02)\\
	\rho^2\alpha_0&C9(E_1&\ov{F}_1&11&10&00&01&;&20&02&12)\\
	\alpha_1&C9(E_1&\ov{E}_2&21&11&22&12&;&C_0&02&01)\\
	\sigma\alpha_1&C9(\ov{E}_1&E_2&12&11&22&21&;&C_0&20&10)\\
	\beta\alpha_1&C9(F_2&\ov{F}_1&11&21&12&22&;&D_0&02&01)\\
	\sigma\beta\alpha_1&C9(\ov{F}_2&F_1&11&12&21&22&;&D_0&20&10)\\
	\alpha_2&C9(C_1&D_2&01&F_2&20&00&;&Cd_0&10&22)\\
	\sigma\alpha_2&C9(C_1&D_2&10&\ov{F}_2&02&00&;&Cd_0&01&22)\\
	021021\alpha_2&C9(C_2&D_1&02&\ov{F}_1&10&00&;&Dc_0&20&11)\\
	\sigma 021021\alpha_2&C9(C_2&D_1&20&F_1&01&00&;&Dc_0&02&11)\\
	\alpha_3&C9(C_0&F_0&00&C_1&10&11&;&CF_0&12&22)\\
	\sigma\alpha_3&C9(C_0&\ov{F}_0&00&C_1&01&11&;&C\ov{F}_0&21&22)\\
	021021\alpha_3&C9(C_0&\ov{F}_0&00&C_2&20&22&;&C\ov{F}_0&21&11)\\
	\sigma 021021\alpha_3&C9(C_0&F_0&00&C_2&02&22&;&CF_0&12&11)\\
	\alpha_4&C9(D_0&E_0&20&E_2&11&21&;&DE_0&01&12)\\
	\sigma\alpha_4&C9(D_0&\ov{E}_0&02&\ov{E}_2&11&12&;&D\ov{E}_0&10&21)\\
	021021\alpha_4&C9(D_0&\ov{E}_0&10&\ov{E}_1&22&12&;&D\ov{E}_0&02&21)\\
	\sigma 021021\alpha_4&C9(D_0&E_0&01&E_1&22&21&;&DE_0&20&12)\\
	\alpha_5&C9(E_0&F_0&12&E_2&02&01&;&EF_0&00&20)\\
	\sigma\alpha_5&C9(\ov{E}_0&\ov{F}_0&21&\ov{E}_2&20&10&;
	&\ov{EF}_0&00&02)\\
	021021\alpha_5&C9(\ov{E}_0&\ov{F}_0&21&\ov{E}_1&01&02&;
	&\ov{EF}_0&00&10)\\
	\sigma 021021\alpha_5&C9(E_0&F_0&12&E_1&10&20&;&EF_0&00&01)\\
	\alpha_6&C9(E_1&\ov{F}_2&22&F_2&11&12&;&E\ov{f}_0&10&01)\\
	\sigma\alpha_6&C9(E_1&F_2&22&\ov{F}_2&11&21&;&\ov{E}f_0&01&10)\\
	\alpha_7&C9(F_1&\ov{E}_2&02&\ov{F}_1&22&20&;&F\ov{e}_0&21&11)\\
	\sigma \alpha_7&C9(\ov{F}_1&E_2&20&F_1&22&02&;&\ov{F}e_0&12&11)\\
	\alpha_8&C9(C_0&E_0&20&21&21&11&;&21'&01&22)\\
	\sigma \alpha_8&C9(C_0&\ov{E}_0&02&12&12&11&;&12'&10&22)\\
	012210\alpha_8&C9(F_1&F_2&22&21&21&11&;&21'&01&20)\\
	\sigma \alpha_8&C9(\ov{F}_1&\ov{F}_2&22&12&12&11&;&12'&10&02)\\
	210102\alpha_8&C9(\ov{E}_0&E_0&01&00&00&10&;&00'&20&02)\\
	120102\alpha_8&C9(D_2&D_1&01&00&00&20&;&00'&10&02)\\

	\alpha_9&C9(C_1&10&11&C_2&02&22&;&CF_0&12&00)\\
	\sigma \alpha_9&C9(C_1&01&11&C_2&20&22&;&C\ov{F}_0&21&00)\\
	210012\alpha_9&C9(D_1&10&11&F_1&22&02&;&\ov{F}e_0&12&20)\\
	\sigma 210012\alpha_9&C9(D_1&01&11&\ov{F}_1&22&20&;&F\ov{e}_0&21&02)\\
	012102\alpha_9&C9(F_2&11&10&D_2&02&22&;&E\ov{f}_0&12&01)\\
	\sigma 012102\alpha_9&C9(\ov{F}_2&11&01&D_2&20&22&;&\ov{E}f_0&21&10)\\
	102210\alpha_9&C9(E_2&02&01&E_1&10&20&;&EF_0&00&12)\\
	\sigma 102210\alpha_9&C9(\ov{E}_2&20&10&\ov{E}_1&01&02
	&;&\ov{EF}_0&00&21)\\
	120201\alpha_9&C9(E_1&22&20&E_2&11&01&;&DE_0&21&12)\\
	\sigma 120201\alpha_9&C9(\ov{E}_1&22&02&\ov{E}_2&11&10
	&;&D\ov{E}_0&12&21)\\
	102120\alpha_9&C9(F_1&01&02&\ov{F}_1&10&20&;&Dc_0&00&11)\\
	201210\alpha_9&C9(\ov{F}_2&02&01&F_2&20&10&;&Cd_0&00&22)\\
	\end{array}$
}

	\ssec{The cubics of Grassmann.}
	\sssec{Definition.}
	Given 6 lines $a_i$ and $b_i$, among the 15 intersections we choose
	the following 9,\\
	D1.0.	$A_i := a_{i+1}\times a_{i-1},$\\
	D1.1.	$B_i := b_{i+1}\times b_{i-1},$\\
	D1.2.	$E_i := a_i\times b_i,$\\
	the non confined configuration consisting of these 9 points and 6 lines
	each containing 3 of the points is called a {\em Grassmann
	configuration.} It is noted $(\{A_i\},\{B_i\},\{E_i\}).$

	\sssec{Theorem. [Grassmann]}\label{sec-grass}
	{\em Given 2 triangles $\{A_i,a_i\}$ and $\{B_i,b_i\}$, the locus of
	the points $X$ is a cubic, if X is such that the points obtained by
	finding the intersections of the lines joining $X$ to the vertices of
	one of the triangles and the corresponding sides of the second triangle,
	namely $(X \times A_i)\times b_i$, are collinear.}

	\sssec{Theorem.}
	{\em Let}\\
	D2.0.	$aB_i := A_{i+1} \times B_{i-1},$
		$a\ov{B}_i := A_{i-1} \times B_{i+1},$\\
	D2.1.	$AB_i := aB_i \times a\ov{B}_i,$\\
	D2.2.	$abE_i := AB_{i+1} \times E_{i-1},$
		$ab\ov{E}_i := AB_{i-1} \times E_{i+1},$\\
	D2.3.	$DE_i := abE_i \times ab\ov{E}_i,$\\
	D2.4.	$de_i := DE_i \times E_i,$\\
	D2.5.	$ab_i := A_i \times B_i,$\\
	D2.6.	$D_i := de_i \times ab_i,$\\
	D2.6.	$ba_i := AB_{i+1} \times AB_{i-1},$\\
	D2.8.	$abd_i := D_i \times AB_i,$\\
	D2.9.	$C_i := ba_i \times abd_i,$\\
	D4.0.	$ae_i := A_i \times E_i,$ $be_i := B_i \times E_i,$\\
	D4.1.	$ce_i := C_i\times E_i,$\\
	D4.2.	$aab_i := A_i\times AB_i,$ $bab_i := B_i\times AB_i,$\\
	D4.3.	$F_i := be_i\times aab_i,$ $\ov{F}_i := ae_i\times bab_i,$\\
	D4.4.	$f_i := F_{i+1}\times F_{i-1},$
		$\ov{f}_i := \ov{F}_{i+1}\times \ov{F}_{i-1},$\\
	D4.5.	$A'_i := ce_i\times \ov{f}_i,$	$B'_i := ce_i\times f_i,$\\
	D4.6.	$aC_i := A_{i+1}\times C_{i-1},$
		$a\ov{C}_i := A_{i-1}\times C_{i+1},$\\
	D4.6.	$bC_i := B_{i+1}\times C_{i-1},$
		$b\ov{C}_i := B_{i-1}\times C_{i+1},$\\
	D4.7.	$CF_i := bC_i\times b\ov{C}_i,$
		$C\ov{F}_i := aC_i\times a\ov{C}_i,$\\
	D4.8.	$cD_i := C_{i+1}\times D_{i-1},$
		$c\ov{D}_i := C_{i-1}\times D_{i+1},$\\
	D4.9.	$aF_i := A_{i+1}\times F_{i-1},$
		$a\ov{F}_i := A_{i-1}\times F_{i+1},$\\
	D4.10.	$Cd_i := cD_i\times aF_i,$
		$Dc_i := c\ov{D}_i\times a\ov{F}_i,$\\
	D4.11.	$aD_i := A_{i+1}\times D_{i-1},$
		$a\ov{D}_i := A_{i-1}\times D_{i+1},$\\
	D4.12.	$eF_i := E_{i+1}\times F_{i-1},$
		$e\ov{F}_i := E_{i-1}\times F_{i+1},$\\
	D4.13.	$Ef_i := aD_i\times eF_i,$
		$Fe_i := a\ov{D}_i\times e\ov{F}_i,$\\
	D4.14.	$bD_i := B_{i+1}\times D_{i-1},$
		$b\ov{D}_i := B_{i-1}\times D_{i+1},$\\
	D4.15.	$\ov{fE}_i := E_{i+1}\times \ov{F}_{i-1},$
		$\ov{f}E_i := E_{i-1}\times \ov{F}_{i+1},$\\
	D4.16.	$E\ov{f}_i := bD_i\times \ov{fE}_i,$
		$\ov{F}e_i := b\ov{D}_i\times \ov{f}E_i,$\\
	D4.17.	$ef_i := E_i\times F_i,$ $e\ov{f}_i := E_i\times \ov{F}_i,$\\
	D5.0.	$a'D_i := A'_{i+1} \times D_{i-1},$
		$a'\ov{D}_i := A'_{i-1} \times D_{i+1},$\\
	D5.1.	$abe_i := AB_i \times E_i,$\\
	D5.2.	$M_i := a'D_i \times abe_i,$\\
	D5.3.	$dAB_i := D_{i+1} \times AB_{i-1},$
		$dA\ov{B}_i := D_{i-1} \times AB_{i+1},$\\
	D5.4.	$a'E_i := A'_{i+1} \times E_{i-1},$
		$a'\ov{E}_i := A'_{i-1} \times E_{i+1},$\\
	D5.5.	$L_i := dAB_i \times a'E_i,$
		$\ov{L}_i := dA\ov{B}_i \times a'\ov{E}_i,$\\
	D5.6.	$fB_i := F_{i+1} \times B_{i-1},$
		$f\ov{B}_i := F_{i-1} \times B_{i+1},$\\
	D5.7.	$\ov{f}A_i := \ov{F}_{i+1} \times A_{i-1},$
		$\ov{fA}_i := \ov{F}_{i-1} \times A_{i+1},$\\
	D5.8.	$P_i := fB_i \times \ov{f}A_i,$
		$Q_i := f\ov{B}_i \times \ov{fA}_i,$\\
	D5.9.	$ac_i := A_i \times C_i,$
		$bc_i := B_i \times C_i,$\\
	D5.10.	$bde_i := B_i \times DE_i,$
		$ade_i := A_i \times DE_i,$\\
	D5.11.	$K_i := ac_i \times bde_i,$
		$\ov{K}_i := bc_i \times ade_i,$\\
	D5.12.	$dE_i := D_{i+1} \times E_{i-1},$
		$eD_i := E_{i+1} \times D_{i-1},$\\
	D5.13.	$bK_i := B_{i+1} \times K_{i-1},$
		$kB_i := B_{i-1} \times K_{i+1},$\\
	D5.14.	$Ed_i := eD_i \times bK_i,$
		$De_i := dE_i \times kB_i,$\\
	D6.0.	$cL_i := C_{i+1} \times L_{i-1},$
		 $c\ov{L}_i := C_{i-1} \times \ov{L}_{i+1},$\\
	D6.1.	$C'_i := cL_i \times c\ov{L}_i,$\\
	D6.2.	$k\ov{f}_i := K_i \times \ov{F}_i,$
		$\ov{k}f_i := \ov{K}_i \times F_i,$\\
	D6.3.	$D'_i := k\ov{f}_i \times \ov{k}f_i,$\\
	D6.4.	$a'f_i := A'_i \times F_i,$
		$a'\ov{f}_i := A'_i \times \ov{F}_i,$\\
	D6.5.	$df_i := D_i \times F_i,$ $d\ov{f}_i := D_i \times \ov{F}_i,$\\
	D6.6.	$DF_i := a'\ov{f}_i \times df_i,$
		$D\ov{F}_i := a'f_i \times d\ov{f}_i,$\\
	D6.7.	$fC_i := F_{i+1} \times C_{i-1},$
		$f\ov{C}_i := F_{i-1} \times C_{i+1},$\\
	D6.8.	$\ov{f}C_i := \ov{F}_{i+1} \times C_{i-1},$
		$\ov{fC}_i := \ov{F}_{i-1} \times C_{i+1},$\\
	D6.9.	$fDE_i := F_{i+1} \times DE_{i-1},$
		$fD\ov{E}_i := \ov{F}_{i+1} \times DE_{i-1},$\\
	D6.10.	$aM_i := A_{i+1} \times M_{i-1},$
		$bM_i := B_{i+1} \times M_{i-1},$\\
	D6.11.	$Fc_i := fC_i \times fD\ov{E}_i,$
		$\ov{F}c_i := \ov{f}C_i \times fDE_i,$\\
	D6.12.	$Cf_i := f\ov{C}_i \times aM_i,$
		$C\ov{f}_i := \ov{fC}_i \times bM_i,$\\
	D6.13.	$c_i := C_{i+1} \times C_{i-1},$
		$ll_i := L_i \times \ov{L}_i,$\\
	D6.14.	$Cc_i := c_i \times ll_i,$\\
	D6.15.	$ccL_i := Cc_{i+1} \times L_{i-1},$
		$cc\ov{L}_i := Cc_{i-1} \times \ov{L}_{i+1},$\\
	D6.16.	$AB'_i := ccL_i \times cc\ov{L}_i,$\\
	D6.17.	$lQ_i := L_{i+1} \times Q_{i-1},$
		$p\ov{L}_i := P_{i+1} \times \ov{L}_{i-1},$\\
	D6.18.	$F'_i := lQ_i \times p\ov{L}_i,$\\
	D6.19.	$f\ov{f}_i := F_i \times \ov{F}_i,$
		$dab'_i := D_i \times AB'_i,$\\
	D6.20.	$FF_i := f\ov{f}_i \times dab'_i,$\\
	D7.0.	$a'_i := A_i \times A'_i,$\\
	D7.1.	$b'_i := B_i \times A'_i,$\\
	D7.2.	$c'_i := C_i \times C'_i,$\\
	D7.3.	$d'_i := D_i \times D'_i,$\\
	D7.4.	$ab'_i := AB_i \times AB'_i,$\\
	D7.5.	$e'_i := E_i \times AB'_i,$\\
	D7.6.	$f'_i := F_i \times F'_i,$\\
	D7.7.	$\ov{f}'_i := \ov{F}_i \times F'_i,$\\
	{\em then}\\
	C2.0.	$C_i \incid e_i,$\\
	C2.1.	$X_i = Y_i$,\\
	C2.2.	$A'_i = B'_i,$\\
	C2.3.	$A_i \incid e\ov{f}_i,$ $B_i \incid ef_i,$\\
	C2.4.	$A'_{i+1}\times D_{i-1} = M_i,$\\
	C2.5.	$A'_{i+1}\times E_{i-1} = L_i,$\\
	C2.6.	$E_{i+1}\times A'_{i-1} = M_i,$\\
	C2.7.	$F'_i = \ov{F}'_i$,

	\sssec{Theorem.}
	{\em Given a Grassmann configuration,} $(\{A_i\},\{B_i\},\{E_i\}),$
	\enumb
	\item {\em the points
	$C_i,$ $AB_i$, $D_i,$ $F_i,$ $\ov{F}_i,$ $Cd_i,$ $Dc_i,$ $CF_i,$
	$C\ov{F}_i,$ $DE_i$, $Ef_i$, $Fe_i$, $E\ov{f}_i$, $\ov{F}e_i$, 
	are on the cubic $\gamma$ through $A_i$, $B_i$, $E_i$.}
	\item $(A_i\,a'_i)$, $(B_i\,b'_i),$ $(C_i\,c'_i),$ $(D_i\,d'_i),$
	$(AB_i\,ab'_i),$ $(E_i\,e'_i),$ $(F_i\,f'_i),$ $(\ov{F}_i\,\ov{f}'_i)
	\incid\gamma.$
	\item $E_i\star E_i = AB_i\star AB_i$.
	\item {\em The points $A_i$, $B_i$, $AB_i$, are on a cubic
	configuration.}
	\item $\langle A'_{i+1},A'_{i-1},AB'_i,\rangle$,
	\item $\langle B'_{i+1},B'_{i-1},AB'_i,\rangle$,
	\item $\langle AB'_{i+1},AB'_{i-1},C'_i,\rangle$,
	\item $\langle AB'_i,C'_i,D'_i,\rangle$,
	\item $\langle A'_i,AB'_i,F'_i,\rangle$,
	\enume

	Proof:
	To prove that $AB_0$ is on the cubic, we have to prove, because of
	\ref{sec-grass}, $\langle (AB_0\times A_0)\times b_0,
	(AB_0\times A_1)\times b_1,(AB_0\times A_2)\times b_2,\rangle$,
	but the second point is $B_2$ and the third is $B_1$ and both
	are on $b_0$.  The rest of the proof follows from \ref{sec-tgras1}
	given below.

	\sssec{Comment.}
	The preceding Theorem was conjectured in the process of construction
	the third point on a Grassmann cubic and the line through 2 points
	of the cubic, using intersection of conics and lines, (Fig. hd.c)
	using the Theorem of Grassmann, (Henry White, 1925, p. 109, Fig. 27.)
	For instance, the conic through $B = B_{01},$ $C = B_{02},$
	$A' = B_{10},$ $D = (B\times C' = B_{12})\times (C\times B' = B_{11}$)
	and $X,$ intersects $A \times X$ at the third point $A\star X$.

\newpage
	\sssec{Theorem.}\label{sec-tgras1}
	{\em We have the following table for the operation $\star$ between
	points on a Grassmann cubic:}\\
{\footnotesize
	$\begin{array}{lclllclllclllclll}
	\star&\smv&AB_0&AB_1&AB_2&\smv&A_0&A_1&A_2 &\smv&B_0&B_1&B_2\\
	\hline
	AB_0&\smv&AB'_0&C_2&C_1&\smv&F_0&B_2&B_1&\smv&\ov{F}_0&A_2&A_1\\
	AB_1&\smv&C_2&AB'_1&C_0&\smv&B_2&F_1&B_0&\smv&A_2&\ov{F}_1&A_0\\
	AB_2&\smv&C_1&C_0&AB'_2	&\smv&B_1&B_0&F_2&\smv&A_1&A_0&\ov{F}_2\\
	\hline
	A_0&\smv&F_0&B_2&B_1&\smv&A'_0&E_2&E_1&\smv&D_0&AB_2&AB_1\\
	A_1&\smv&B_2&F_1&B_0&\smv&E_2&A'_1&E_0&\smv&AB_2&D_1&AB_0\\
	A_2&\smv&B_1&B_0&F_2&\smv&E_1&E_0&A'_2&\smv&AB_1&AB_0&D_2\\
	\hline
	B_0&\smv&\ov{F}_0&A_2&A_1&\smv&D_0&AB_2&AB_1&\smv&A'_0&E_2&E_1\\
	B_1&\smv&A_2&\ov{F}_1&A_0&\smv&AB_2&D_1&AB_0&\smv&E_2&A'_1&E_0\\
	B_2&\smv&A_1&A_0&\ov{F}_2&\smv&AB_1&AB_0&D_2&\smv&E_1&E_0&A'_2\\
	\end{array}$

	\noindent$\begin{array}{lclllclllclllclllclll}
	\star&\smv&AB_0&AB_1&AB_2&\smv&A_0&A_1&A_2&\smv&B_0&B_1&B_2&\smv&
	C_0&C_1&C_2&\smv&D_0&D_1&D_2\\
	\hline
	C_0&\smv&D_0&AB_2&AB_1&\smv&K_0&C\ov{F}_2&C\ov{F}_1
	&\smv&\ov{K}_0&CF_2&CF_1&\smv&C'_0&Cc_2&Cc_1&\smv&AB_0&Cd_2&Dc_1\\
	C_1&\smv&AB_2&D_1&AB_0&\smv&C\ov{F}_2&K_1&C\ov{F}_0
	&\smv&CF_2&\ov{K}_1&CF_0&\smv&Cc_2&C'_1&Cc_0&\smv&Dc_2&AB_1&Cd_0\\
	C_2&\smv&AB_1&AB_0&D_2&\smv&C\ov{F}_1&C\ov{F}_0&K_2
	&\smv&CF_1&CF_0&\ov{K}_2&\smv&Cc_1&Cc_0&C'_2&\smv&Cd_1&Dc_0&AB_2\\
	\hline
	D_0&\smv&C_0&L_2&\ov{L}_1&\smv&B_0&Fe_2&Ef_1
	&\smv&A_0&\ov{F}e_2&E\ov{f}_1&\smv&AB_0&Dc_2&Cd_1
	&\smv&D'_0&AB_2'&AB_1'\\
	D_1&\smv&\ov{L}_2&C_1&L_0&\smv&Ef_2&B_1&Fe_0
	&\smv&E\ov{f}_2&A_1&\ov{F}e_0&\smv&Cd_2&AB_1&Dc_1
	&\smv&AB_2'&D'_1&AB_0'\\
	D_2&\smv&L_1&\ov{L}_0&C_2&\smv&Fe_1&Ef_0&B_2
	&\smv&\ov{F}e_1&E\ov{f}_0&A_2&\smv&Dc_1&Cd_0&AB_2
	&\smv&AB_1'&AB_0'&D'_2\\
	\hline
	E_0&\smv&M_0&DE_2&DE_1&\smv&\ov{F}_0&A_2&A_1&\smv&F_0&B_2&B_1
	&\smv&A'_0&E_2&E_1&\smv&DE_0&Ed_2&De_1\\
	E_1&\smv&DE_2&M_1&DE_0&\smv&A_2&\ov{F}_1&A_0&\smv&B_2&F_1&B_0
	&\smv&E_2&A'_1&E_0&\smv&De_2&DE_1&Ed_0\\
	E_2&\smv&DE_1&DE_0&M_2&\smv&A_1&A_0&\ov{F}_2&\smv&B_1&B_0&F_2
	&\smv&E_1&E_0&A'_2&\smv&Ed_1&De_0&DE_2\\
	\hline
	F_0&\smv&A_0&\ov{F}e_2&E\ov{f}_1&\smv&AB_0&Dc_2&Cd_1
	&\smv&E_0&P_2&Q_1
	&\smv&CF_0&Fc_2&Cf_1&\smv&DF_0&\ov{F}_2&\ov{F}_1\\
	F_1&\smv&E\ov{f}_2&A_1&\ov{F}e_0&\smv&Cd_2&AB_1&Dc_0
	&\smv&Q_2&E_1&P_0
	&\smv&Cf_2&CF_1&Fc_0&\smv&\ov{F}_2&DF_1&\ov{F}_0\\
	F_2&\smv&\ov{F}e_1&E\ov{f}_0&A_2&\smv&Dc_1&Cd_0&AB_2
	&\smv&P_1&Q_0&E_2
	&\smv&Fc_1&Cf_0&CF_2&\smv&\ov{F}_1&\ov{F}_0&DF_2\\
	\hline
	\ov{F}_0&\smv&B_0&Fe_2&Ef_1&\smv&E_0&P_2&Q_1
	&\smv&AB_0&Dc_2&Cd_1
	&\smv&C\ov{F}_0&\ov{F}c_2&C\ov{f}_1&\smv&D\ov{F}_0&F_2&F_1\\
	\ov{F}_1&\smv&Ef_2&B_1&Fe_0&\smv&Q_2&E_1&P_0
	&\smv&Cd_2&AB_1&Dc_0
	&\smv&C\ov{f}_2&C\ov{F}_1&\ov{F}c_0&\smv&F_2&D\ov{F}_1&F_0\\
	\ov{F}_2&\smv&Fe_1&Ef_0&B_2&\smv&P_1&Q_0&E_2
	&\smv&Dc_1&Cd_0&AB_2
	&\smv&\ov{F}c_1&C\ov{f}_0&C\ov{F}_2&\smv&F_1&F_0&D\ov{F}_2\\
	\end{array}$

	\noindent$\begin{array}{lclllclllclllclll}
	\star&\smv&E_0&E_1&E_2
	&\smv&F_0&F_1&F_2 &\smv&\ov{F}_0&\ov{F}_1&\ov{F}_2\\
	\hline
	E_0&\smv&AB'_0&C_2&C_1
	&\smv&B_0&Ef_2&Fe_1&\smv&A_0&E\ov{f}_2&\ov{F}e_1\\
	E_1&\smv&C_2&AB'_1&C_0
	&\smv&Fe_2&B_1&Ef_0&\smv&\ov{F}e_2&A_1&E\ov{f}_0\\
	E_2&\smv&C_1&C_0&AB'_2
	&\smv&Ef_1&Fe_0&B_2&\smv&E\ov{f}_1&\ov{F}e_0&A_2\\
	\hline
	F_0&\smv&B_0&Fe_2&Ef_1&\smv&F'_0&A'_2&A'_1&\smv&FF_0&D_2&D_1\\
	F_1&\smv&Ef_2&B_1&Fe_0&\smv&A'_2&F'_1&A'_0&\smv&D_2&FF_1&D_0\\
	F_2&\smv&Fe_1&Ef_0&B_2&\smv&A'_1&A'_0&F'_2&\smv&D_1&D_0&FF_2\\
	\hline
	\ov{F}_0&\smv&A_0&\ov{F}e_2&E\ov{f}_1
	&\smv&FF_0&D_2&D_1&\smv&F'_0&A'_2&A'_1\\
	\ov{F}_1&\smv&E\ov{f}_2&A_1&\ov{F}e_0
	&\smv&D_2&FF_1&D_0&\smv&A'_2&F'_1&A'_0\\
	\ov{F}_2&\smv&\ov{F}e_1&E\ov{f}_0&A_2
	&\smv&D_1&D_0&FF_2&\smv&A'_1&A'_0&F'_2\\
	\end{array}$
}	
\newpage
	Proof:\\
{\footnotesize
	$\begin{array}{lclllclllclllc}
	\alpha_0&D9(C_0&D_0&B_0&A_1&B_2&AB_1&;&AB_0&A_0&AB_2)\\
	\rho\alpha_0&D9(E_2&F_2&A_2&AB_0&A_1&B_0&;&B_2&AB_2&B_1)\\
	\rho^2\alpha_0&C9(E_1&\ov{F}_1&AB_1&B_2&AB_0&A_2&;&A_1&B_1&A_0)\\
	\alpha_1&C9(E_1&E_2&B_0&AB_1&AB_2&A_0&;&C_0&B_1&A_2)\\
	\beta\alpha_1&C9(F_2&\ov{F}_1&AB_1&B_0&A_0&AB_2&;&D_0&B_1&A_2)\\
	\sigma\beta\alpha_1&C9(\ov{F}_2&F_1&AB_1&A_0&B_0&AB_2&;&D_0&A_1&B_2)\\
	\alpha_2&C9(C_1&D_2&A_2&F_2&A_1&AB_0&;&Cd_0&B_2&AB_2)\\
	\sigma\alpha_2&C9(C_1&D_2&B_2&\ov{F}_2&B_1&AB_0&;&Cd_0&A_2&AB_2)\\
	021021\alpha_2&C9(C_2&D_1&B_1&\ov{F}_1&B_2&AB_0&;&Dc_0&A_1&AB_1)\\
	\sigma 021021\alpha_2&C9(C_2&D_1&A_1&F_1&A_2&AB_0&;&Dc_0&B_1&AB_1)\\
	\alpha_3&C9(C_0&F_0&AB_0&C_1&B_2&AB_1&;&CF_0&A_0&AB_2)\\
	\sigma\alpha_3&C9(C_0&\ov{F}_0&AB_0&C_1&A_2&AB_1&;&C\ov{F}_0&B_0&AB_2)\\
	021021\alpha_3&C9(C_0&\ov{F}_0&AB_0&C_2&A_1&AB_2&;&C\ov{F}_0&B_0&AB_1)\\
	\sigma 021021\alpha_3&C9(C_0&F_0&AB_0&C_2&B_1&AB_2&;&CF_0&A_0&AB_1)\\
	\alpha_4&C9(D_0&E_0&A_1&E_2&AB_1&B_0&;&DE_0&A_2&A_0)\\
	021021\alpha_4&C9(D_0&E_0&B_2&E_1&AB_2&A_0&;&DE_0&B_1&B_0)\\
	\alpha_5&C9(E_0&F_0&A_0&E_2&B_1&A_2&;&B_0&AB_0&A_1)\\
	\sigma\alpha_5&C9(E_0&\ov{F}_0&B_0&E_2&A_1&B_2&;&A_0&AB_0&B_1)\\
	021021\alpha_5&C9(E_0&\ov{F}_0&B_0&E_1&A_2&B_1&;&A_0&AB_0&B_2)\\
	\sigma 021021\alpha_5&C9(E_0&F_0&A_0&E_1&B_2&A_1&;&B_0&AB_0&A_2)\\
	\alpha_6&C9(E_1&\ov{F}_2&AB_2&F_2&AB_1&A_0&;&E\ov{f}_0&B_2&A_2)\\
	\sigma\alpha_6&C9(E_1&F_2&AB_2&\ov{F}_2&AB_1&B_0&;&Ef_0&A_2&B_2)\\
	\alpha_7&C9(F_1&E_2&B_1&\ov{F}_1&AB_2&A_1&;&Fe_0&B_0&AB_1)\\
	\sigma \alpha_7&C9(\ov{F}_1&E_2&A_1&F_1&AB_2&B_1&;&\ov{F}e_0&A_0&AB_1)\\
	\sigma \alpha_8&C9(C_0&E_0&B_1&A_0&A_0&AB_1&;&A'_0&B_2&AB_2)\\
	\alpha'_8&C9(A_0&A_0&B_1&B_0&B_0&A_1&;&A'_0&AB_2&E_1)\\
	012210\alpha_8&C9(F_1&F_2&AB_2&B_0&B_0&AB_1&;&A'_0&A_2&A_1)\\
	\sigma \alpha_8&C9(\ov{F}_1&\ov{F}_2&AB_2&A_0&A_0&AB_1&;&A'_0&B_2&B_1)\\
	210102\alpha_8&C9(E_0&E_0&A_2&AB_0&AB_0&B_2&;&AB'_0&A_1&B_1)\\
	120102\alpha_8&C9(D_2&D_1&A_2&AB_0&AB_0&A_1&;&AB'_0&B_2&B_1)\\

	210012\alpha_9&C9(D_1&B_2&AB_1&F_1&AB_2&B_1&;&\ov{F}e_0&A_0&A_1)\\
	\sigma 210012\alpha_9&C9(D_1&A_2&AB_1&\ov{F}_1&AB_2&A_1&;
	&Fe_0&B_0&B_1)\\
	012102\alpha_9&C9(F_2&AB_1&B_2&D_2&B_1&AB_2&;&E\ov{f}_0&A_0&A_2)\\
	\sigma 012102\alpha_9&C9(\ov{F}_2&AB_1&A_2&D_2&A_1&AB_2&;
	&Ef_0&B_0&B_2)\\
	120201\alpha_9&C9(E_1&AB_2&A_1&E_2&AB_1&A_2&;&DE_0&B_0&A_0)\\

	\alpha_{10}&C9(AB_0&E_0&A_2&D_2&A'_1&A_1&;&M_0&A_1&B_2)\\
	\alpha_{11}&C9(AB_2&D_1&B_1&A'_1&E_2&A_0&;&L_0&A_1&B_1)\\
	?\alpha_{11}&C9(AB_1&D_2&B_2&A'_2&E_1&A_0&;&\ov{L}_0&A_2&B_2)\\
	\alpha_{12}&C9(F_1&B_2&B_0&A_2&\ov{F}_1&A_1&;&P_0&E_1&AB_1)\\
	?\alpha_{12}&C9(F_2&B_1&B_0&A_1&\ov{F}_2&A_2&;&Q_0&E_2&AB_2)\\
	\alpha_{13}&C9(C_0&C_0&AB_0&C_1&L_2&AB_1&;&C'_0&D_0&AB_2)\\
	?\alpha_{13}&C9(C_0&C_0&AB_0&C_2&\ov{L}_1&AB_2&;&C'_0&D_0&AB_1)\\
	\alpha_{14}&C9(F_0&F_0&B_0&\ov{F}_0&\ov{F}_0&A_0&;&F'_0&E_0&AB_0)\\
	\alpha_{15}&C9(A_0&C_0&E_0&DE_0&B_0&B_0&;&K_0&A'_0&D_0)\\
	?\alpha_{15}&C9(B_0&C_0&E_0&DE_0&A_0&A_0&;&\ov{K}_0&A'_0&D_0)\\
	\alpha_{16}&C9(AB_0&AB_0&AB_1&L_2&Cc_1&C_0&;&AB'_0&C_2&D_0)\\
	?\alpha_{16}&C9(AB_0&AB_0&AB_2&\ov{L}_1&Cc_2&C_0&;&AB'_0&C_1&D_0)\\
	\alpha_{17}&C9(AB_0&AB_0&C_1&Cc_0&C_0&AB_1&;&AB'_0&AB_2&C_2)\\
	\alpha_{18}&C9(C_1&C_2&AB_2&L_0&\ov{L}_0&AB_1&;&Cc_0&D_2&D_1)\\
	\alpha_{19}&C9(F_0&F_0&D_1&\ov{L}_2&P_1&A_0&;&F'_0&\ov{F}_2&AB_0)\\
	\alpha_{20}&C9(D_0&F_0&E_0&\ov{F}_0&A'_0&B_0&;&DF_0&B_0&A_0)\\
	?\alpha_{20}&C9(D_0&\ov{F}_0&E_0&F_0&A'_0&A_0&;&D\ov{F}_0&A_0&B_0)\\
	\end{array}$
\newpage
	\noindent$\begin{array}{lclllclllclllc}
	\alpha_{21}&C9(F_1&C_2&AB_1&\ov{F}_1&DE_2&E_1&;&Fc_0&AB_0&B_1)\\
	?\alpha_{21}&C9(\ov{F}_1&C_2&AB_1&F_1&DE_2&E_1&;&\ov{F}c_0&AB_0&A_1)\\
	\alpha_{22}&C9(C_1&F_2&E_2&M_2&A_1&AB_0&;&Cf_0&B_2&AB_2)\\
	?\alpha_{22}&C9(C_1&\ov{F}_2&E_2&M_2&B_1&AB_0&;&C\ov{f}_0&A_2&AB_2)\\
	\alpha_{23}&C9(F_0&F_0&D_2&L_1&Q_2&A_0&;&F'_0&\ov{F}_1&AB_0)\\
	?\alpha_{23}&C9(F_0&F_0&D_1&\ov{L}_2&P_1&A_0&;&F'_0&\ov{F}_2&AB_0)\\
	\alpha_{24}&C9(D_0&D_0&C_0&K_0&\ov{F}_0&B_0&;&D'_0&AB_0&A_0)\\
	?\alpha_{24}&C9(D_0&D_0&C_0&\ov{K}_0&F_0&A_0&;&D'_0&AB_0&A_0)\\
	\alpha_{25}&C9(E_1&D_2&B_2&B_1&K_2&C_2&;&Ed_0&A_2&E_0)\\
	?\alpha_{25}&C9(E_2&D_1&B_1&B_2&K_1&C_1&;&De_0&A_1&E_0)\\
	\alpha_{26}&C9(C_0&F_0&A_0&B_0&AB'_0&AB_0&;&CF_0&AB_0&D_0)\\
	?\alpha_{26}&C9(C_0&\ov{F}_0&B_0&A_0&AB'_0&AB_0&;&C\ov{F}_0&AB_0&D_0)\\
	\alpha_{27	}&C9(F_0&\ov{F}_0&B_0&D_0&AB'_0&AB_0&;&FF_0&AB_0&A_0)\\
	\end{array}$
}
	\sssec{Theorem.}
	{\em Given a Grassmann configuration $(\{A_i\},\{B_i\},\{E_i\}),$
	THe tangential points at $A_i$ and $B_i$ are the same,
	in the above case $A'_i$, it will be added to the configuration
	after a semi colon.}

	\sssec{Lemma.}\label{sec-lgras}
	{\em If $(\{A_i\},\{B_i\},\{E_i\},\{A'_i\}),$ is a Grassmann
	configuration so is}
	\enumb
	\item $(\{A_i\},\{B_i\},\{E_i\}),$,\\
	{\em where\\
	\hth$AB_i := A_{i+1}\star B_{i-1},$ and
	$C_i := AB_{i+1}\star AB_{i-1}.$}
	\item $(\{F_i\},\{\ov{F}_i\},\{A'_i\}),$,\\
	{\em where}\\
	\hth$F_i := AB_{i}\star A_{i},$ and $\ov{F}_i := AB_{i}\star B_{i},$
	\enume

	This follows at once from grom \ref{sec-tgras1}.

	\sssec{Theorem.}
	{\em Given a Grassmann configuration $(\{A_i\},\{B_i\},\{E_i\};
	\{AB_i\}),$
	the following are also Grassmann configurations:}
	\enumb
	\item $(\{AB_i\},\{E_i\},\{C_i\};\{AB'_i\}),$
	\item $(\{F_i\},\{\ov{F}_i\},\{A'_i\};\{F'_i\}),$
	\item $(\{DE_i\},\{C_i\},\{C_i\star AB'_i\}),$
	\item $(\{D_i\},\{A'_i\},\{AB'_i\};\{d'_i\}).$
	\enume

	Proof: This follows by repeated applications of the Lemma
	\ref{sec-lgras}.
\newpage
	\ssec{Grassmannian cubics in Involutive Geometry.}
	\sssec{Definition.}
	In involutive geometry I will give the name of {\em Grassmannian cubic}
	to the special case where the 6 lines are $m_i$ and $\ov{m}_i$.

	\sssec{Theorem.}
	{\em The correspondence between the elements as given above and those of
	involutive geometry is as follows}\\
	$\begin{array}{cccccccccccccc}
	A_i=CF_i&B_i=C\ov{F}_i&E_i=DE_i&AB_i=C_i=Cc_i&F_i=K_i
	&\ov{F}_i=\ov{K}_i\\
	MM_i&\ov{M}M_i&A_i&AT_i&F_i&\ov{F}_i\\
	\hline
	Ef_i=C\ov{f}_i&E\ov{f}_i=Cf_i&Fe_i=\ov{F}c_i&\ov{F}e_i=Fc_i&
	A'_i=M_i&L_i=Cd_i&\ov{L}_i=Dc_i&P_i&Q_i\\
	Ef_i&E\ov{f}_i&Fe_i&\ov{F}e_i&MM'_i&L_i&\ov{L}_i&P_i&Q_i\\
	\hline
	a_i&b_i&aB_i&a\ov{B}_i&ab_i&e_i&ae_i&be_i&ba_i\\
	mm_i&\ov{m}m_i&c_i&\ov{c}_i&a_i&aeUL_i&nm_i&\ov{n}m_i&eul\\
	\end{array}$

	\sssec{Theorem.}
	\enumb
	\item {\em The Grassmann cubic passes through the points $MM_i,$
	$\ov{M}M_i,$ $A_i,$ $D_i.$}
	\item {\em Its equation is}\\
	\hth$m_0(m_1+m_2)X_1X_2((m_2+m0)X_1+(m_0+m_1)X_2)
		m_1(m_2+m_0)X_2X_0((m_0+m1)X_2+(m_1+m_2)X_0)
		m_2(m_0+m_1)X_0X_1((m_1+m2)X_0+(m_2+m_0)X_1)
		+(s_{21}+2s_{111})X_0X_1X_2 = 0$\\
	\enume

	Proof: Using \ref{sec-tsal} on the points $A_i$ and $AT_i$,
	we obtain the given form, to determine the coefficients $g_i$ of
	$X_1X_2((m_2+m0)X_1+(m_0+m_1)X_2)$, \ldots and $g$ of $X_0X_1X_2$
	we impose the condition that the cubic passes through $MM_i$,
	this gives the system of equations\\
	\hth$\ldots.$

	\sssec{Theorem.}
	{\em Let}\\
	D2.0.	$aB_i := A_{i+1} \times B_{i-1},$
		$a\ov{B}_i := A_{i-1} \times B_{i+1},$\\
	D2.1.	$AB_i := aB_i \times a\ov{B}_i,$\\
	D2.2.	$abE_i := AB_{i+1} \times E_{i-1},$
		$ab\ov{E}_i := AB_{i-1} \times E_{i+1},$\\
	D2.3.	$DE_i := abE_i \times ab\ov{E}_i,$\\
	D2.4.	$de_i := DE_i \times E_i,$\\
	D2.5.	$ab_i := A_i \times B_i,$\\
	D2.6.	$D_i := de_i \times ab_i,$\\
	D2.6.	$ba_i := AB_{i+1} \times AB_{i-1},$\\
	D2.8.	$abd_i := D_i \times AB_i,$\\
	D2.9.	$C_i := ba_i \times abd_i,$\\
	D4.0.	$ae_i := A_i \times E_i,$ $be_i := B_i \times E_i,$\\
	D4.1.	$ce_i := C_i\times E_i,$\\
	D4.2.	$aab_i := A_i\times AB_i,$ $bab_i := B_i\times AB_i,$\\
	D4.3.	$F_i := be_i\times aab_i,$ $\ov{F}_i := ae_i\times bab_i,$\\
	D4.4.	$f_i := F_{i+1}\times F_{i-1},$
		$\ov{f}_i := \ov{F}_{i+1}\times \ov{F}_{i-1},$\\
	D4.5.	$A'_i := ce_i\times \ov{f}_i,$	$B'_i := ce_i\times f_i,$\\
	D4.6.	$aC_i := A_{i+1}\times C_{i-1},$
		$a\ov{C}_i := A_{i-1}\times C_{i+1},$\\
	D4.6.	$bC_i := B_{i+1}\times C_{i-1},$
		$b\ov{C}_i := B_{i-1}\times C_{i+1},$\\
	D4.7.	$CF_i := bC_i\times b\ov{C}_i,$
		$C\ov{F}_i := aC_i\times a\ov{C}_i,$\\
	D4.8.	$cD_i := C_{i+1}\times D_{i-1},$
		$c\ov{D}_i := C_{i-1}\times D_{i+1},$\\
	D4.9.	$aF_i := A_{i+1}\times F_{i-1},$
		$a\ov{F}_i := A_{i-1}\times F_{i+1},$\\
	D4.10.	$Cd_i := cD_i\times aF_i,$
		$Dc_i := c\ov{D}_i\times a\ov{F}_i,$\\
	D4.11.	$aD_i := A_{i+1}\times D_{i-1},$
		$a\ov{D}_i := A_{i-1}\times D_{i+1},$\\
	D4.12.	$eF_i := E_{i+1}\times F_{i-1},$
		$e\ov{F}_i := E_{i-1}\times F_{i+1},$\\
	D4.13.	$Ef_i := aD_i\times eF_i,$
		$Fe_i := a\ov{D}_i\times e\ov{F}_i,$\\
	D4.14.	$bD_i := B_{i+1}\times D_{i-1},$
		$b\ov{D}_i := B_{i-1}\times D_{i+1},$\\
	D4.15.	$\ov{fE}_i := E_{i+1}\times \ov{F}_{i-1},$
		$\ov{f}E_i := E_{i-1}\times \ov{F}_{i+1},$\\
	D4.16.	$E\ov{f}_i := bD_i\times \ov{fE}_i,$
		$\ov{F}e_i := b\ov{D}_i\times \ov{f}E_i,$\\
	D4.17.	$ef_i := E_i\times F_i,$ $e\ov{f}_i := E_i\times \ov{F}_i,$\\
	D5.0.	$a'D_i := A'_{i+1} \times D_{i-1},$
		$a'\ov{D}_i := A'_{i-1} \times D_{i+1},$\\
	D5.1.	$abe_i := AB_i \times E_i,$\\
	D5.2.	$M_i := a'D_i \times abe_i,$\\
	D5.3.	$dAB_i := D_{i+1} \times AB_{i-1},$
		$dA\ov{B}_i := D_{i-1} \times AB_{i+1},$\\
	D5.4.	$a'E_i := A'_{i+1} \times E_{i-1},$
		$a'\ov{E}_i := A'_{i-1} \times E_{i+1},$\\
	D5.5.	$L_i := dAB_i \times a'E_i,$
		$\ov{L}_i := dA\ov{B}_i \times a'\ov{E}_i,$\\
	D5.6.	$fB_i := F_{i+1} \times B_{i-1},$
		$f\ov{B}_i := F_{i-1} \times B_{i+1},$\\
	D5.7.	$\ov{f}A_i := \ov{F}_{i+1} \times A_{i-1},$
		$\ov{fA}_i := \ov{F}_{i-1} \times A_{i+1},$\\
	D5.8.	$P_i := fB_i \times \ov{f}A_i,$
		$Q_i := f\ov{B}_i \times \ov{fA}_i,$\\
	D5.9.	$ac_i := A_i \times C_i,$
		$bc_i := B_i \times C_i,$\\
	D5.10.	$bde_i := B_i \times DE_i,$
		$ade_i := A_i \times DE_i,$\\
	D5.11.	$K_i := ac_i \times bde_i,$
		$\ov{K}_i := bc_i \times ade_i,$\\
	D5.12.	$dE_i := D_{i+1} \times E_{i-1},$
		$eD_i := E_{i+1} \times E_{i-1},$\\
	D5.13.	$bK_i := B_{i+1} \times K_{i-1},$
		$kB_i := B_{i-1} \times K_{i+1},$\\
	D5.14.	$Ed_i := eD_i \times bK_i,$
		$De_i := dE_i \times kB_i,$\\
	D6.0.	$cL_i := C_{i+1} \times L_{i-1},$
		 $c\ov{L}_i := C_{i-1} \times \ov{L}_{i+1},$\\
	D6.1.	$C'_i := cL_i \times c\ov{L}_i,$\\
	D6.2.	$k\ov{f}_i := K_i \times \ov{F}_i,$
		$\ov{k}f_i := \ov{K}_i \times F_i,$\\
	D6.3.	$D'_i := k\ov{f}_i \times \ov{k}f_i,$\\
	D6.4.	$a'f_i := A'_i \times F_i,$
		$a'\ov{f}_i := A'_i \times \ov{F}_i,$\\
	D6.5.	$df_i := D_i \times F_i,$ $d\ov{f}_i := D_i \times \ov{F}_i,$\\
	D6.6.	$DF_i := a'\ov{f}_i \times df_i,$
		$D\ov{F}_i := a'f_i \times d\ov{f}_i,$\\
	D6.7.	$fC_i := F_{i+1} \times C_{i-1},$
		$f\ov{C}_i := F_{i-1} \times C_{i+1},$\\
	D6.8.	$\ov{f}C_i := \ov{F}_{i+1} \times C_{i-1},$
		$\ov{fC}_i := \ov{F}_{i-1} \times C_{i+1},$\\
	D6.9.	$fDE_i := F_{i+1} \times DE_{i-1},$
		$fD\ov{E}_i := \ov{F}_{i+1} \times DE_{i-1},$\\
	D6.10.	$aM_i := A_{i+1} \times M_{i-1},$
		$bM_i := B_{i+1} \times M_{i-1},$\\
	D6.11.	$Fc_i := fC_i \times fD\ov{E}_i,$
		$\ov{F}c_i := \ov{f}C_i \times fDE_i,$\\
	D6.12.	$Cf_i := f\ov{C}_i \times aM_i,$
		$C\ov{f}_i := \ov{fC}_i \times bM_i,$\\
	D6.13.	$c_i := C_{i+1} \times C_{i-1},$
		$ll_i := L_i \times \ov{L}_i,$\\
	D6.14.	$Cc_i := c_i \times ll_i,$\\
	D6.15.	$ccL_i := Cc_{i+1} \times L_{i-1},$
		$cc\ov{L}_i := Cc_{i-1} \times \ov{L}_{i+1},$\\
	D6.16.	$AB'_i := ccL_i \times cc\ov{L}_i,$\\
	D6.17.	$lQ_i := L_{i+1} \times Q_{i-1},$
		$p\ov{L}_i := P_{i+1} \times \ov{L}_{i-1},$\\
	D6.18.	$F'_i := lQ_i \times p\ov{L}_i,$\\
	D6.19.	$f\ov{f}_i := F_i \times \ov{F}_i,$
		$dab'_i := D_i \times AB'_i,$\\
	D6.20.	$FF_i := f\ov{f}_i \times dab'_i,$\\
	D7.0.	$a'_i := A_i \times A'_i,$\\
	D7.1.	$b'_i := B_i \times A'_i,$\\
	D7.2.	$c'_i := C_i \times C'_i,$\\
	D7.3.	$d'_i := D_i \times D'_i,$\\
	D7.4.	$ab'_i := AB_i \times AB'_i,$\\
	D7.5.	$e'_i := E_i \times AB'_i,$\\
	D7.6.	$f'_i := F_i \times F'_i,$\\
	D7.7.	$\ov{f}'_i := \ov{F}_i \times F'_i,$\\
	{\em then}\\
	C2.0.	$C_i \incid e_i,$\\
	C2.1.	$X_i = Y_i$,\\
	C2.2.	$A'_i = B'_i,$\\
	C2.3.	$A_i \incid e\ov{f}_i,$ $B_i \incid ef_i,$\\
	C2.4.	$A'_{i+1}\times D_{i-1} = M_i,$\\
	C2.5.	$A'_{i+1}\times E_{i-1} = L_i,$\\
	C2.6.	$E_{i+1}\times A'_{i-1} = M_i,$\\
	C2.7.	$F'_i = \ov{F}'_i$,\\

	Proof:
	H0.0.	$m_0 = [0,1,1],$ $\ov{m}_0 = [0,m2,m1],$\\
	D1.0.	$MM_0 = (-1,1,1),$ $\ov{M}M_0 = (m0,-m1,-m2)$\\
	D0.10.	$A_0 = (1,0,0)$\\
	D0.10.	$a_0 = [1,0,0]$\\
	D..	$ma_0 = [0,1,-1],$ $\ov{m}a_0 = [0,m2,-m1],$\\
	D..	$eul_0 = [m1-m2,-(m2-m0),-(m0-m1)],$\\
	D..	$y_0 = [m1-m2,-(m2+m0),-(m0+m1)],$\\
	D..	$y~_0 = [m1-m2,m2+m0,m0+m1],$\\
	D..	$AT_0 = (0,m0+m1,-(m2-m0)),$\\
	D..	$k = [m1+m2,m2+m0,m0+m1],$\\
	D..	$\ov{t}AM_0 = [s1+m0,m2+m0,m0+m1],$\\
	D..	$tAM_0 = [s11+m1m2,m0(m2+m0),m0(m0+m1)],$\\
	D..	$F_0 = (s11+m1m2,-m1(s1+m0),-m2(s1+m0)),$\\
	D..	$\ov{F}_0 = (m0(s1+m0),-(s11+m1m2),-(s11+m1m2)),$\\
	D..	$\ov{f}f_0 = [m1(m1-m2),s11+m2m0,m1(s1+m2)],$\\
	D..	$D_0 = (m1+m2)q0,-m1(m1-m2)(m2+m0),m2(m0+m1)(m1-m2)),$\\
	D..	$eul = [m1-m2,m2-m0,m0-m1],$\\
	D..	$aAT_0 = [0,m2+m0,m0+m1)],$\\
	D..	$\ov{f}_0 = [(m1+m2)q0,-(m2+m0)(s11+m2m0),
		-(m0+m1)(s11+m0m1),$\\
	D..	$MM'_0 = (m0(m1-m2)(m2+m0)(m0+m1),-(m0+m1)(m1+m2)q0,
		(m1+m2)(m2+m0)q0),$\\
	D..	$_0 = [(m1-m2)(s11+m2m0),m1(m2-m0)^2,
		m2q1-m0q2+m0m1(m1-m2)],$\\
	D..	$_0 = [s11+m2m0,m0(s1+m1),0],$\\
	D..	$Fe_0 = (m0(s1+m1)(m2q1-m0q2+m0m1(m1-m2)),
		-(s11+m2m0)(m2q1-m0q2+m0m1(m1-m2)),
		(s11+m2m0)(m1q2-m0q1-m2m0(m1-m2)),$\\
	D..	$_0 = [],$\\
	D..	$_0 = [],$\\
	D..	$Ef_0 = (m0(s1+m2)(m1q2-m0q1-m2m0(m1-m2),
		(s11+m0m1)(m2q1-m0q2+m0m1(m1-m2)),
		-(s11+m0m1)(m1q2-m0q1-m2m0(m1-m2))),$\\
	D..	$_0 = [],$\\
	D..	$_0 = [],$\\
	D..	$\ov{F}e_0 = ((s11+m2m0)(m1q2-m0q1-m2m0(m1-m2),
		(s11+m0m1)(m1q2-m0q1-m2m0(m1-m2),
		-(s11+m0m1)(m2q1-m0q2+m0m1(m1-m2)),$\\
	D..	$_0 = [],$\\
	D..	$_0 = [],$\\
	D..	$E\ov{f}_0 = ((s11+m0m1)(m2q1-m0q2+m0m1(m1-m2),
		m1(s1+m2)(m1q2-m0q1-m2m0(m1-m2),
		-m2(s1+m2)(m2q1-m0q2+m0m1(m1-m2)),$\\
	D..	$_0 = [],$\\
	D..	$_0 = [],$\\
	D..	$Ed_0 = (m0(m1+m2)^2(m2+m0)(m0-m1),
		-m1q2(2m1(m1+m2)+(m2+m0)(m0+m1)),s21+2s111),$\\
	D..	$_0 = [],$\\
	D..	$_0 = [],$\\
	D..	$De_0 = (m0(m1+m2)^2(m2-m0)(m0+m1),
		s21+2s111,m2q1(2m2(m1+m2)+(m2+m0)(m0+m1))),$\\
	D..	$_0 = [],$\\
	D..	$_0 = [],$\\
	D..	$L_0 = ((m2+m0)q1(m0^2+m1m2+3m0(m1+m2),
		m1(m1+m2)(m2-m0)(m0^2+m1m2+3m0(m1+m2),
		-(m1+m2)(m2-m0)(s21+2s111)),$\\
	D..	$_0 = [],$\\
	D..	$_0 = [],$\\
	D..	$\ov{L}_0 = ((m0+m1)q2(m0^2+m1m2+3m0(m1+m2),
		(m0-m1)(m1+m2)(s21+2s111),
		-m2(m0-m1)(m1+m2)(m0^2+m1m2+3m0(m1+m2)),$\\
	D..	$_0 = [],$\\
	D..	$_0 = [],$\\
	D..	$P_0 = (m0(m1+m2)^2(m2-m0)(m0+m1),
		(s21+2s111)q1,-m2q1(m0^2+m1m2+3m0(m1+m2)),$\\
	D..	$_0 = [],$\\
	D..	$_0 = [],$\\
	D..	$Q_0 = (m0(m1+m2)^2(m2+m0)(m0-m1),
		-m1q2(m0^2+m1m2+3m0(m1+m2),-(s21+2s111)q2),$

\newpage
	\sssec{Theorem.}\label{sec-tgras2}
	{\em We have the following table for the operation $\star$ between
	points on the Grassmannian cubic:}\\
{\footnotesize
	$\begin{array}{lclllclllclllclllclll}
	\star&\smv&AT_0&AT_1&AT_2&\smv&MM_0&MM_1&MM_2
	&\smv&\ov{M}M_0&\ov{M}M_1&\ov{M}M_2&\smv&A_0&A_1&A_2&\smv&D_0&D_1&D_2\\
	\hline
	AT_0&\smv&D_0&AT_2&AT_1&\smv&F_0&\ov{M}M_2&\ov{M}M_1
	&\smv&\ov{F}_0&MM_2&MM_1&\smv&MM'_0&A_2&A_1
	&\smv&AT_0&\ov{L}_2&L_1\\
	AT_1&\smv&AT_2&D_1&AT_0&\smv&\ov{M}M_2&F_1&\ov{M}M_0
	&\smv&MM_2&\ov{F}_1&MM_0&\smv&A_2&MM'_1&A_0
	&\smv&L_2&AT_1&\ov{L}_0\\
	AT_2&\smv&AT_1&AT_0&D_2	&\smv&\ov{M}M_1&\ov{M}M_0&F_2
	&\smv&MM_1&MM_0&\ov{F}_2&\smv&A_1&A_0&MM'_2
	&\smv&\ov{L}_1&L_0&AT_2\\
	\hline
	MM_0&\smv&F_0&\ov{M}M_2&\ov{M}M_1&\smv&MM'_0&A_2&A_1
	&\smv&D_0&AT_2&AT_1&\smv&\ov{F}_0&MM_2&MM_1
	&\smv&\ov{M}M_0&Ef_2&Fe_1\\
	MM_1&\smv&\ov{M}M_2&F_1&\ov{M}M_0&\smv&A_2&MM'_1&A_0
	&\smv&AT_2&D_1&AT_0&\smv&MM_2&\ov{F}_1&MM_0
	&\smv&Fe_2&\ov{M}M_1&Ef_0\\
	MM_2&\smv&\ov{M}M_1&\ov{M}M_0&F_2&\smv&A_1&A_0&MM'_2
	&\smv&AT_1&AT_0&D_2&\smv&MM_1&MM_0&\ov{F}_2
	&\smv&Ef_1&Fe_0&\ov{M}M_2\\
	\hline
	\ov{M}M_0&\smv&\ov{F}_0&MM_2&MM_1&\smv&D_0&AT_2&AT_1
	&\smv&MM'_0&A_2&A_1&\smv&F_0&\ov{M}M_2&\ov{M}M_1
	&\smv&MM_0&E\ov{f}_2&\ov{F}e_1\\
	\ov{M}M_1&\smv&MM_2&\ov{F}_1&MM_0&\smv&AT_2&D_1&AT_0
	&\smv&A_2&MM'_1&A_0&\smv&\ov{M}M_2&F_1&\ov{M}M_0
	&\smv&\ov{F}e_2&MM_1&E\ov{f}_0\\
	\ov{M}M_2&\smv&MM_1&MM_0&\ov{F}_2&\smv&AT_1&AT_0&D_2
	&\smv&A_1&A_0&MM'_2&\smv&\ov{M}M_1&\ov{M}M_0&F_2
	&\smv&E\ov{f}_1&\ov{F}e_0&MM_2\\
	\hline
	A_0&\smv&MM'_0&A_2&A_1&\smv&\ov{F}_0&MM_2&MM_1
	&\smv&F_0&\ov{M}M_2&\ov{M}M_1
	&\smv&D_0&AT_2&AT_1&\smv&A_0&De_2&Ed_1\\
	A_1&\smv&A_2&MM'_1&A_0&\smv&MM_2&\ov{F}_1&MM_0
	&\smv&\ov{M}M_2&F_1&\ov{M}M_0
	&\smv&AT_2&D_1&AT_0&\smv&Ed_2&A_1&De_0\\
	A_2&\smv&A_1&A_0&MM'_2&\smv&MM_1&MM_0&\ov{F}_2
	&\smv&\ov{M}M_1&\ov{M}M_0&F_2
	&\smv&AT_1&AT_0&D_2&\smv&De_1&Ed_0&A_2\\
	\hline
	D_0&\smv&AT_0&L_2&\ov{L}_1&\smv&\ov{M}M_0&Fe_2&Ef_1
	&\smv&MM_0&\ov{F}e_2&E\ov{f}_1&\smv&A_0&Ed_2&De_1
	&\smv&D'_0&D_2&D_1\\
	D_1&\smv&\ov{L}_2&AT_1&L_0&\smv&Ef_2&\ov{M}M_1&Fe_0
	&\smv&E\ov{f}_2&MM_1&\ov{F}e_0&\smv&De_2&A_1&Ed_1
	&\smv&D_2&D'_1&D_0\\
	D_2&\smv&L_1&\ov{L}_0&AT_2&\smv&Fe_1&Ef_0&\ov{M}M_2
	&\smv&\ov{F}e_1&E\ov{f}_0&MM_2&\smv&Ed_1&De_0&A_2
	&\smv&D_1&D_0&D'_2\\
	\hline
	F_0&\smv&MM_0&\ov{F}e_2&E\ov{f}_1&\smv&AT_0&\ov{L}_2&L_1
	&\smv&A_0&P_2&Q_1&\smv&\ov{M}M_0&Fe_2&Ef_1
	&\smv&FD_0&\ov{F}_2&\ov{F}_1\\
	F_1&\smv&E\ov{f}_2&MM_1&\ov{F}e_0&\smv&L_2&AT_1&\ov{L}_0
	&\smv&Q_2&A_1&P_0&\smv&Ef_2&\ov{M}M_1&Fe_0
	&\smv&\ov{F}_2&FD_1&\ov{F}_0\\
	F_2&\smv&\ov{F}e_1&E\ov{f}_0&MM_2&\smv&\ov{L}_1&L_0&AT_2
	&\smv&P_1&Q_0&A_2&\smv&Fe_1&Ef_0&\ov{M}M_2
	&\smv&\ov{F}_1&\ov{F}_0&FD_2\\
	\hline
	\ov{F}_0&\smv&\ov{M}M_0&Fe_2&Ef_1&\smv&A_0&P_2&Q_1
	&\smv&AT_0&\ov{L}_2&L_1&\smv&MM_0&\ov{F}e_2&E\ov{f}_1
	&\smv&\ov{F}D_0&F_2&F_1\\
	\ov{F}_1&\smv&Ef_2&\ov{M}M_1&Fe_0&\smv&Q_2&A_1&P_0
	&\smv&L_2&AT_1&\ov{L}_0&\smv&E\ov{f}_2&MM_1&\ov{F}e_0
	&\smv&F_2&\ov{F}D_1&F_0\\
	\ov{F}_2&\smv&Fe_1&Ef_0&\ov{M}M_2&\smv&P_1&Q_0&A_2
	&\smv&\ov{L}_1&L_0&AT_2&\smv&\ov{F}e_1&E\ov{f}_0&MM_2
	&\smv&F_1&F_0&\ov{F}D_2\\
	\end{array}$

	\noindent$\begin{array}{lclllclll}
	\star&\smv&F_0&F_1&F_2 &\smv&\ov{F}_0&\ov{F}_1&\ov{F}_2\\
	\hline
	F_0&\smv&F'_0&MM'_2&MM'_1&\smv&D'_0&D_2&D_1\\
	F_1&\smv&MM'_2&F'_1&MM'_0&\smv&D_2&D'_1&D_0\\
	F_2&\smv&MM'_1&MM'_0&F'_2&\smv&D_1&D_0&D'_2\\
	\hline
	\ov{F}_0&\smv&D'_0&D_2&D_1&\smv&F'_0&MM'_2&MM'_1\\
	\ov{F}_1&\smv&D_2&D'_1&D_0&\smv&MM'_2&F'_1&MM'_0\\
	\ov{F}_2&\smv&D_1&D_0&D'_2&\smv&MM'_1&MM'_0&F'_2\\
	\end{array}$
}	
\newpage
	Proof:\\
{\footnotesize
	$\begin{array}{lclllclllclllc}
	\alpha_0&D9(C_0&D_0&\ov{M}M_0&MM_1&\ov{M}M_2&AT_1&;&AT_0&MM_0&AT_2)\\
	\rho\alpha_0&D9(E_2&F_2&MM_2&AT_0&MM_1&\ov{M}M_0&;&\ov{M}M_2&AT_2&\ov{M}M_1)\\
	\rho^2\alpha_0&C9(E_1&\ov{F}_1&AT_1&\ov{M}M_2&AT_0&MM_2&;&MM_1&\ov{M}M_1&MM_0)\\
	\alpha_1&C9(E_1&A_2&\ov{M}M_0&AT_1&AT_2&MM_0&;&AT_0&\ov{M}M_1&MM_2)\\
	\beta\alpha_1&C9(F_2&\ov{F}_1&AT_1&\ov{M}M_0&MM_0&AT_2&;&D_0&\ov{M}M_1&MM_2)\\
	\sigma\beta\alpha_1&C9(\ov{F}_2&F_1&AT_1&MM_0&\ov{M}M_0&AT_2&;&D_0&MM_1&\ov{M}M_2)\\
	\alpha_2&C9(C_1&D_2&MM_2&F_2&MM_1&AT_0&;&L_0&\ov{M}M_2&AT_2)\\
	\sigma\alpha_2&C9(C_1&D_2&\ov{M}M_2&\ov{F}_2&\ov{M}M_1&AT_0&;&L_0&MM_2&AT_2)\\
	021021\alpha_2&C9(C_2&D_1&\ov{M}M_1&\ov{F}_1&\ov{M}M_2&AT_0&;&\ov{L}_0&MM_1&AT_1)\\
	\sigma 021021\alpha_2&C9(C_2&D_1&MM_1&F_1&MM_2&AT_0&;&\ov{L}_0&\ov{M}M_1&AT_1)\\
	\alpha_3&C9(C_0&F_0&AT_0&AT_1&\ov{M}M_2&AT_1&;&MM_0&MM_0&AT_2)\\
	\sigma\alpha_3&C9(C_0&\ov{F}_0&AT_0&AT_1&MM_2&AT_1&;&\ov{M}M_0&\ov{M}M_0&AT_2)\\
	021021\alpha_3&C9(C_0&\ov{F}_0&AT_0&AT_2&MM_1&AT_2&;&\ov{M}M_0&\ov{M}M_0&AT_1)\\
	\sigma 021021\alpha_3&C9(C_0&F_0&AT_0&AT_2&\ov{M}M_1&AT_2&;&MM_0&MM_0&AT_1)\\
	\alpha_4&C9(D_0&A_0&MM_1&A_2&AT_1&\ov{M}M_0&;&A_0&MM_2&MM_0)\\
	021021\alpha_4&C9(D_0&A_0&\ov{M}M_2&A_1&AT_2&MM_0&;&A_0&\ov{M}M_1&\ov{M}M_0)\\
	\alpha_5&C9(E_0&F_0&MM_0&A_2&\ov{M}M_1&MM_2&;&\ov{M}M_0&AT_0&MM_1)\\
	\sigma\alpha_5&C9(E_0&\ov{F}_0&\ov{M}M_0&A_2&MM_1&\ov{M}M_2&;&MM_0&AT_0&\ov{M}M_1)\\
	021021\alpha_5&C9(E_0&\ov{F}_0&\ov{M}M_0&A_1&MM_2&\ov{M}M_1&;&MM_0&AT_0&\ov{M}M_2)\\
	\sigma 021021\alpha_5&C9(E_0&F_0&MM_0&A_1&\ov{M}M_2&MM_1&;&\ov{M}M_0&AT_0&MM_2)\\
	\alpha_6&C9(E_1&\ov{F}_2&AT_2&F_2&AT_1&MM_0&;&E\ov{f}_0&\ov{M}M_2&MM_2)\\
	\sigma\alpha_6&C9(E_1&F_2&AT_2&\ov{F}_2&AT_1&\ov{M}M_0&;&Ef_0&MM_2&\ov{M}M_2)\\
	\alpha_7&C9(F_1&A_2&\ov{M}M_1&\ov{F}_1&AT_2&MM_1&;&Fe_0&\ov{M}M_0&AT_1)\\
	\sigma \alpha_7&C9(\ov{F}_1&A_2&MM_1&F_1&AT_2&\ov{M}M_1&;&\ov{F}e_0&MM_0&AT_1)\\
	\sigma \alpha_8&C9(C_0&A_0&\ov{M}M_1&MM_0&MM_0&AT_1&;&A'_0&\ov{M}M_2&AT_2)\\
	\alpha'_8&C9(MM_0&MM_0&\ov{M}M_1&\ov{M}M_0&\ov{M}M_0&MM_1&;&A'_0&AT_2&A_1)\\
	012210\alpha_8&C9(F_1&F_2&AT_2&\ov{M}M_0&\ov{M}M_0&AT_1&;&A'_0&MM_2&MM_1)\\
	\sigma \alpha_8&C9(\ov{F}_1&\ov{F}_2&AT_2&MM_0&MM_0&AT_1&;&A'_0&\ov{M}M_2&\ov{M}M_1)\\
	210102\alpha_8&C9(E_0&A_0&MM_2&AT_0&AT_0&\ov{M}M_2&;&AB'_0&MM_1&\ov{M}M_1)\\
	120102\alpha_8&C9(D_2&D_1&MM_2&AT_0&AT_0&MM_1&;&AB'_0&\ov{M}M_2&\ov{M}M_1)\\

	210012\alpha_9&C9(D_1&\ov{M}M_2&AT_1&F_1&AT_2&\ov{M}M_1&;&\ov{F}e_0&MM_0&MM_1)\\
	\sigma 210012\alpha_9&C9(D_1&MM_2&AT_1&\ov{F}_1&AT_2&MM_1&;
	&Fe_0&\ov{M}M_0&\ov{M}M_1)\\
	012102\alpha_9&C9(F_2&AT_1&\ov{M}M_2&D_2&\ov{M}M_1&AT_2&;&E\ov{f}_0&MM_0&MM_2)\\
	\sigma 012102\alpha_9&C9(\ov{F}_2&AT_1&MM_2&D_2&MM_1&AT_2&;
	&Ef_0&\ov{M}M_0&\ov{M}M_2)\\
	120201\alpha_9&C9(E_1&AT_2&MM_1&A_2&AT_1&MM_2&;&A_0&\ov{M}M_0&MM_0)\\

	\alpha_{10}&C9(AT_0&A_0&MM_2&D_2&A'_1&MM_1&;&M_0&MM_1&\ov{M}M_2)\\
	\alpha_{11}&C9(AT_2&D_1&\ov{M}M_1&A'_1&A_2&MM_0&;&L_0&MM_1&\ov{M}M_1)\\
	?\alpha_{11}&C9(AT_1&D_2&\ov{M}M_2&A'_2&A_1&MM_0&;&\ov{L}_0&MM_2&\ov{M}M_2)\\
	\alpha_{12}&C9(F_1&\ov{M}M_2&\ov{M}M_0&MM_2&\ov{F}_1&MM_1&;&P_0&A_1&AT_1)\\
	?\alpha_{12}&C9(F_2&\ov{M}M_1&\ov{M}M_0&MM_1&\ov{F}_2&MM_2&;&Q_0&A_2&AT_2)\\
	\alpha_{13}&C9(C_0&AT_0&AT_0&AT_1&L_2&AT_1&;&C'_0&D_0&AT_2)\\
	?\alpha_{13}&C9(C_0&AT_0&AT_0&AT_2&\ov{L}_1&AT_2&;&C'_0&D_0&AT_1)\\
	\alpha_{14}&C9(F_0&F_0&\ov{M}M_0&\ov{F}_0&\ov{F}_0&MM_0&;&F'_0&A_0&AT_0)\\
	\alpha_{15}&C9(MM_0&AT_0&A_0&A_0&\ov{M}M_0&\ov{M}M_0&;&F_0&A'_0&D_0)\\
	?\alpha_{15}&C9(\ov{M}M_0&AT_0&A_0&A_0&MM_0&MM_0&;&ov{F}_0&A'_0&D_0)\\
	\alpha_{16}&C9(AT_0&AT_0&AT_1&L_2&Cc_1&AT_0&;&AB'_0&AT_2&D_0)\\
	?\alpha_{16}&C9(AT_0&AT_0&AT_2&\ov{L}_1&Cc_2&AT_0&;&AB'_0&AT_1&D_0)\\
	\alpha_{17}&C9(AT_0&AT_0&AT_1&Cc_0&AT_0&AT_1&;&AB'_0&AT_2&AT_2)\\
	\alpha_{18}&C9(C_1&AT_2&AT_2&L_0&\ov{L}_0&AT_1&;&Cc_0&D_2&D_1)\\
	\alpha_{19}&C9(F_0&F_0&D_1&\ov{L}_2&P_1&MM_0&;&F'_0&\ov{F}_2&AT_0)\\
	\alpha_{20}&C9(D_0&F_0&A_0&\ov{F}_0&A'_0&\ov{M}M_0&;&DF_0&\ov{M}M_0&MM_0)\\
	?\alpha_{20}&C9(D_0&\ov{F}_0&A_0&F_0&A'_0&MM_0&;&D\ov{F}_0&MM_0&\ov{M}M_0)\\
	\end{array}$
\newpage
	\noindent$\begin{array}{lclllclllclllc}
	\alpha_{21}&C9(F_1&AT_2&AT_1&\ov{F}_1&A_2&A_1&;&Fc_0&AT_0&\ov{M}M_1)\\
	?\alpha_{21}&C9(\ov{F}_1&AT_2&AT_1&F_1&A_2&A_1&;&\ov{F}c_0&AT_0&MM_1)\\
	\alpha_{22}&C9(C_1&F_2&A_2&M_2&MM_1&AT_0&;&MM_0&\ov{M}M_2&AT_2)\\
	?\alpha_{22}&C9(C_1&\ov{F}_2&A_2&M_2&\ov{M}M_1&AT_0&;&\ov{M}M_0&MM_2&AT_2)\\
	\alpha_{23}&C9(F_0&F_0&D_2&L_1&Q_2&MM_0&;&F'_0&\ov{F}_1&AT_0)\\
	?\alpha_{23}&C9(F_0&F_0&D_1&\ov{L}_2&P_1&MM_0&;&F'_0&\ov{F}_2&AT_0)\\
	\alpha_{24}&C9(D_0&D_0&AT_0&F_0&\ov{F}_0&\ov{M}M_0&;&D'_0&AT_0&MM_0)\\
	?\alpha_{24}&C9(D_0&D_0&AT_0&ov{F}_0&F_0&MM_0&;&D'_0&AT_0&MM_0)\\
	\alpha_{25}&C9(E_1&D_2&\ov{M}M_2&\ov{M}M_1&F_2&AT_2&;&Ed_0&MM_2&A_0)\\
	?\alpha_{25}&C9(E_2&D_1&\ov{M}M_1&\ov{M}M_2&F_1&AT_1&;&A_0&MM_1&A_0)\\
	\alpha_{26}&C9(C_0&F_0&MM_0&\ov{M}M_0&AB'_0&AT_0&;&MM_0&AT_0&D_0)\\
	?\alpha_{26}&C9(C_0&\ov{F}_0&\ov{M}M_0&MM_0&AB'_0&AT_0&;&\ov{M}M_0&AT_0&D_0)\\
	\alpha_{27	}&C9(F_0&\ov{F}_0&\ov{M}M_0&D_0&AB'_0&AT_0&;&FF_0&AT_0&MM_0)\\
	\end{array}$
}
	\sssec{Exercise.}\label{sec-egraseuc}
	Study the Grassmannian cubic when the 6 lines are $mm_i$ and
	$\ov{m}m_i$.
\newpage
	\ssec{Answer to}
	\vspace{-18pt}\hspace{94pt}{\bf \ref{sec-egraseuc}.}\\
	\sssec{Definition.}
	In involutive geometry I will give the name of {\em Grassmannian cubic}
	to the special case where the 6 lines are $mm_i$ and $\ov{m}m_i$.

	\sssec{Theorem.}
	{\em The correspondence between the elements as given above and those of
	involutive geometry is as follows}\\
	$\begin{array}{cccccccccc}
	A_i&B_i&E_i&AB_i&C_i=Cc_i\\
	M_i&\ov{M}_i&EUL_i&D_i&Aeul_i\\
	\hline
	a_i&b_i&aB_i&a\ov{B}_i&ab_i&e_i&ae_i&be_i&ba_i\\
	mm_i&\ov{m}m_i&c_i&\ov{c}_i&a_i&aeUL_i&nm_i&\ov{n}m_i&eul\\
	\end{array}$

	\sssec{Theorem.}
	\enumb
	\item {\em The Grassmann cubic passes through the points $M_i,$
	$\ov{M}_i,$ $EUL_i,$ $D_i.$}
	\item {\em Its equation is ?}\\
	\hth$g_0X_0(-X_0+X_1+X_2)(-m_1m_2X_0+m_2m_0X_1+m_0m_1X_2)
		+g_1X_1(X_0-X_1+X_2)(m_1m_2X_0-m_2m_0X_1+m_0m_1X_2)\\
		+g_2X_2(X_0+X_1-X_2)(m_1m_2X_0+m_2m_0X_1-m_0m_1X_2)
		= 8m_0m_1m_2(m_2X_0+m_0X_1-m_0X_2)(-m_1X_0+m_0X_1+m_1X_2)
		(m_2X_0-m_2X_1+m_1X_2)$\\
	where\\
	\hth$g_0 = (m_1-m_2)(s_{21}-2m_0(m_1^2+m_2^2-m_1m_2))$, \ldots
	\enume

	Proof: Using \ref{sec-tsal} on the points $M_i$ and $\ov{M}_i$,
	we obtain the given form, to determine $g_i$ we impose the
	condition that the cubic passes through $EUL_i$, with
	$EUL_0 = (-m_0(m_1-m_2),m_1(m_2-m_0),m_2(m_0-m_1)),$ this gives the
	system of equations\\
	\hth$4m_1m_2(m_0-m_1)s_1\:+\:4m_1m_2(m_2-m_0)s_2
		= (m_0-m_1)(m_2-m_0)(m_1+m_2)^2,\\
	4m_2m_0(m_0-m_1)s_0\:+\:4m_2m_0(m_1-m_2)s_2
		= (m_1-m_2)(m_0-m_1)(m_2+m_0)^2,\\
	4m_0m_1(m_2-m_0)s_0\:+\:4m_0m_1(m_1-m_2)s_1
		= (m_2-m_0)(m_1-m_2)(m_0+m_1)^2,$\\
	the $s_i$ are proportional to
	$(m_1-m_2)(s_{21}-2m_0(m_1^2+m_2^2-m_1m_2))$
	and the constant of proprtionality is easily determined by
	substitution into one of the equations.\\
	Verify that $(0,m_0-m_1,-(m_2-m_0))$ is on the cubic.
\newpage

	\ssec{The cubics of Tucker.}
	\vspace{-18pt}\hspace{171pt}\footnote{Tucker, Messenger of Mathematics,
	 Ser. 2, Vol. 17, 1887-1888, p. 103}\\[8pt]
	\sssec{Lemma.}
	\enumb
	\item $m_2(m_0-m_1)q_1+m_1(m_2-m_0)q_2 = (m_1^2-m_2^2)q_0.$
	\item $m_0(m_1-m_2)q_1q_2+m_1(m_2-m_0)q_2q_0+m_2(m_0-m_1)q_0q_1\\
		= -s_1s_{11}(m_1-m_2)(m_2-m_0)(m_0-m_1).$
	\item $m_1m_2(m_2-m_0)(m_0-m_1)q_0+m_2m_0(m_0-m_1)(m_1-m_2)q_1\\
		+m_0m_1(m_1-m_2)(m_2-m_0)q_2 = -(q_0q_1q_2)^2.$
	\enume

	\sssec{Definition.}
	Let $A_i$ and $Q$ be a complete quadrilateral, the family of cubics
	associated to $A_i$and $Q$ are the cubics through $A_i$, $Q_i$
	and tangent at $A_i$ to $aq_i$, with\\
	$q$, the polar of $Q$ with respect to $\{A_i\}$,\\
	$Q_i := a_i \times q,$ $aq_i := A_i \times Q_i.$

	\sssec{Theorem.}
	{\em If $Q = (T_0,T_1,T_2),$ the Tucker family of cubics is}\\
	$(T_0 X_1 X_2 + T_1 X_2 X_0 + T_2 X_0 X_1)
		(T_1T_2 X_0 + T_2T_0 X_1 + T_0T_1 X_2)\\
	\hti{12} = k T_0T_1T_2 X_0 X_1 X_2.$\\

	{\em Any point $R$ distinct from $A_i$ and $Q$ is on one and only one of
	these cubics, noted Tucker$(Q)(R).$}

	\sssec{Theorem.}
	{\em If $R = (R_0,R_1,R_2)$ is on Tucker$(M),$ so are}\\
	isobaric$(R)$ or $(R_0,R_1,R_2),$ $(R_2,R_0,R_1),$ $(R_1,R_2,R_0),$\\
	semi reciprocal$(R)$ or $(R_0,R_2,R_1),$ $(R_2,R_1,R_0),$
	$(R_1,R_0,R_2),$\\
	reciprocal$(R)$ or $(R_1R_2,R_2R_0,R_0R_1),$ $(R_0R_1,R_1R_2,R_2R_0),$
	$(R_2R_0,R_0R_1,R_1R_2),$\\
	iso reciprocal$(R)$ or $(R_1R_2,R_0R_1,R_2R_0),$
	$(R_0R_1,R_2R_0,R_1R_2),$ $(R_2R_0,R_1R_2,R_0R_1).$

	\sssec{Theorem.}
	{\em The following are special cases of Tucker cubics:\\
	$k = 1/0,$ for $R \cdot a_i = 0,$
	$Tucker(Q)(R) = a_0 \bigx a_1 \bigx a_2.$\\
	$k = 0,$ for $R \cdot m = 0$ or on
	    $conic(Q) := conic(A_1,aq_1,A_2,aq_2,A_0 ),$\\
	    where $aq_i := A_i \times Q_i,$\\
	    Tucker$(Q)(R) = conic(Q) \bigx m.$
	$k = 1$, for $R \cdot aq_i = 0$,
	Tucker$(Q)(R) = aq_0 \bigx aq_1 \bigx aq_2.$\\
	$k = 9,$ Tucker$(Q)(Q).$
	Finally the constant $k$ is the same for Tucker$(Q)(R)$ and for
	Tucker$(R)(Q).$}

	\sssec{Theorem.}
	\enumb
	\item conic$(K) = \theta.$
	\item Tucker$(M)(\ov{M})$ {\em is incident to} $A_i,$ $MA_i,$ $\ov{M},$
	$PO$, $P\ov{O}$, $MAI,$ $P,$ $\ov{P}$, $Atm_i,$
	\item Tucker$(\ov{M})(M)$ {\em is incident to} $A_i,$ $\ov{M}A_i,$ $M,$
	$\ov{A}tm_i,$
	$\ov{Tm}m,$ $\ov{T}mm,$ $Tmm.$
	\item Tucker$(M)(K)$ {\em is incident to} $A_i,$ $MA_i,$ $K,$ $Br1_i,$ 
	$B\ov{r},$ $B\ov{r}$
	\item Tucker$(\ov{M})(K)$ {\em is incident to} $A_i,$ $\ov{I}m_i,$ $K,$
	 $\ov{B}r1_i.$
	\item Tucker$(\ov{M})(O)$ {\em is incident to} $A_i,$ $MA_i,$ $0,$
	$LEM$,
	\item Tucker$(M)(O)$ {\em is incident to \ldots}
	\enume

	\sssec{Theorem.}
	{\em In the finite case, there are $p+1$ such cubics, each has besides
	the 6 vertices $A_i$ and $Q_i$ of the complete quadrilateral, a number
	of points which is
	a multiple of 6 except when $k = \frac{1}{0}$ and 1, when it is
	$3(p-2),$ $k = 0,$ when it is $2p-5-\js{-3}{p},$ $k = 9$ when it is
	$p-5-\js{-3}{p}.$\\
	$\js{-3}{p}$ is the Jacobi symbol = 1 when $p = 1 \pmod{6}$ and = -1
	when $p = 5 \pmod{6}.$}

	\sssec{Construction of the cubic of Tucker$(M)(\ov{M})$ by the ruler 
	only.}
	H0.0.	$A_i,$ See Fig. t and t'\\
	H0.1.	$M,$ $\ov{M},$\\
	D0.0 to .5, construct $a_i,$ $ma_i,$ $\ov{m}a_i,$ $M_i,$ $\ov{M}_i,$
		$mm_i,$ $MA_i,$ $mm_i,$ $m_i,$\\
	D1.2, D3.0, 3.1, D4.12 and D4.26 construct $Maa_i,$ $\ov{M}aa_i,$
	$cc_i,$ $\ov{c}c_i,$ $MMb_i,$\\
	\hth$MM\ov{b}_i,$ $mn_i,$ $m\ov{n}_i,$ $PO,$ $P\ov{O}.$\\
	We then proceed as follows\\
	D80.0.	$Aa_i := \ov{M},P\ov{O},PO,$\\
	D80.0.	$aaA_i := A_{i+1} \times Aa_{i-1},$
		$aa\ov{A}_i := A_{i-1} \times Aa_{i+1},$\\
	D80.0.	$Ad_i := aaA_i \times aa\ov{A}_i,$\\
	D80.0.	$maaA_i := MA_{i+1} \times Aa_{i-1},$
		$maa\ov{A}_i := MA_{i-1} \times Aa_{i+1},$\\
	D80.0.	$Ab_i := maaA_i \times maa\ov{A}_i,$\\
	D80.0.	$aab_i := A_{i+1} \times Ab_{i-1},$
		$aa\ov{b}_i := A_{i-1} \times Ab_{i+1},$\\
	D80.0.	$Ac_i := aab_i \times aa\ov{b}_i,$\\
\\
	D80.1.	$aaaa_i := Aa_{i+1} \times Aa_{i-1},$
		$acac_i := Ac_{i+1} \times Ac_{i-1},$\\
	D80.1.	$Ba_i := aaaa_i \times acac_i,$\\
	D80.1.	$aaac_i := Aa_{i+1} \times Ac_{i-1},$
		$aaa\ov{c}_i := Aa_{i-1} \times Ac_{i+1},$\\
	D80.1.	$Bc_i := aaac_i \times aaa\ov{c}_i,$\\
	D80.1.	$abA_i := A_{i+1} \times Ba_{i-1},$
		$ab\ov{A}_i := A_{i-1} \times Ba_{i+1},$\\
	D80.1.	$Bd_i := abA_i \times ab\ov{A}_i,$\\
	D80.1.	$mabA_i := MA_{i+1} \times Ba_{i-1},$
		$mab\ov{A}_i := MA_{i-1} \times Ba_{i+1},$\\
	D80.1.	$Bb_i := mabA_i \times mab\ov{A}_i,$\\
\\
	D80.2.	$baba_i := Ba_{i+1} \times Ba_{i-1},$
		$bcbc_i := Bc_{i+1} \times Bc_{i-1},$\\
	D80.2.	$Ca_i := baba_i \times bcbc_i,$\\
	D80.2.	$babc_i := Ba_{i+1} \times Bc_{i-1},$
		$bab\ov{c}_i := Ba_{i-1} \times Bc_{i+1},$\\
	D80.2.	$Cc_i := babc_i \times bab\ov{c}_i,$\\
	D80.2.	$acA_i := A_{i+1} \times Ca_{i-1},$
		$ac\ov{A}_i := A_{i-1} \times Ca_{i+1},$\\
	D80.2.	$Cd_i := acA_i \times ac\ov{A}_i,$\\
	D80.2.	$macA_i := MA_{i+1} \times Ca_{i-1},$
		$mac\ov{A}_i := MA_{i-1} \times Ca_{i+1},$\\
	D80.2.	$Cb_i := macA_i \times mac\ov{A}_i,$\\
	D80.3.	$`Tucker := cubic(A_i,mm_i,MA_1,MA_2,\ov{M}),$\\
\\
	C80.0.	$MA_0 \cdot `Tucker = 0,$\\
	C80.0.	$PO \cdot `Tucker = 0,$ $P\ov{O} \cdot `Tucker = 0,$\\
	C80.0.	$\ov{i}OK \cdot `Tucker = 0,$\\
	C80.0.	$Ba_0 \cdot \ov{i}OK = 0,$\\
	C80.0.	$(Ba_i \times Aa_i) \cdot `Tucker = 0,$ at $Aa_i$?\\
	C80.0.	$(Bb_i \times Ab_i) \cdot `Tucker = 0,$ at $Ab_i$?\\
	C80.0.	$(Bc_{i-1} \times Ad_i) \cdot `Tucker = 0, $at $Ad_i$?\\
	C80.0.	$(Bd_{i+1} \times Ac_i) \cdot `Tucker = 0, $at $Ac_i$?\\
	C80.1.	$Ab_i,$ $Ac_i,$ $Ad_i \cdot `Tucker = 0,$\\
	C80.1.	$Ba_i,$ $Bb_i,$ $Bc_i,$ $Bd_i \cdot `Tucker = 0,$\\
	C80.1.	$Ca_i,$ $Cb_i,$ $Cc_i,$ $Cd_i \cdot `Tucker = 0,$
	C80.2.	$Ad_i \cdot tmm_i = 0,$\\
	C80.2.	$Ac_i \cdot mIA_i = 0,$\\
	C80.2.	$Ab_i \cdot mAM_i = 0,$\\
	C80.2.	$Ac_i \cdot (MA_i\times Ad_i) = 0,$\\
	C80.2.	$Ba_i \cdot pOL_i = 0,$\\
	we can continue indefinitely.

	\sssec{The cubic of Tucker$(M)(\ov{M})$.}
	I have determined all the intersections of the following lines
	with the cubic of Tucker$(M)(\ov{M})$.\\
	$mm_i:\:A_i,A_i,MA_i,\\
	a_i:\:A_{i+1},A_{i-1},MA_i,\\
	ma_i:\:A_i,\ov{M},MNa_{i},\\
	mn_i:\:A_i,P\ov{O},MNa_{i+1},\\
	m\ov{n}_i:\:A_i,P{O},MNa_{i-1},\\
	mIA_i:\:A_i,MAI,Atm_i,\\
	cc_i:\:A_i,P,Atm_{i+1},\\
	\ov{c}c_i:\:A_i,\ov{P},Atm_{i+1},\\
	m:\:MA_i,\\
	maM_i:\:MA_i,\ov{M},Atm_i,\\
	aaM_i:\:MA_{i+1},PO,Atm_{i-1},\\
	aa\ov{M}_i:\:MA_{i-1},P\ov{O},Atm_{i+1},\\
	tmm_i:\:MA_i,\ov{P},MNa_{i-1},\\
	tm\ov{m}_i:\:MA_i,\ov{P},MNa_{i+1},\\
	mpo:\:\ov{M},PO,PAM\\
	mp\ov{o}:\:\ov{M},P\ov{O},,PA\ov{M}\\
	\ov{\i}POK:\:\ov{M},MAI,\\
	:\:\ov{M},P,\\
	:\:\ov{M},\ov{P},\\
	pOL:\:PO,P\ov{O},POl,\\
	:\:MA_i,MNA_i,MAI,\\
	pmai:\:P,MAI,,PAM\\
	pma\ov{i}:\:\ov{P},MAI,,PA\ov{M}\\
	pp:\:P,\ov{P},POl,$

	\sssec{Notation.}
	The correspondance between the notation used here and that used in EUC. 
	is as follows:\\
	$\begin{array}{cccccc}
	aaA_i&aa\ov{A}_i&maaA_i&maa\ov{A}_i\\
	m\ov{n}_1,\ov{m}a_2,mn_0&mn_2,m\ov{n}_0,\ov{m}a_1&
	aaM_0,maM_2,aa\ov{M}_1&aa\ov{M}_0,aaM_2,maM_1
	\end{array}$



	\ssec{NOTES}
	Vigari\'{e} (Mathesis. S\'{e}rie 1, Vol. 9, 1889, Suppl. pp. 1-26
	gives the distances to the sides $\delta_a,$ $\delta_b,$ $\delta_c,$
	the normal coordinates $x,$ $y,$ $z$ which are proportional to these,
	and or the barycentric coordinates $\alpha,$ $\beta$, $\gamma,$ which
	are proportional to $ax,$ $by,$ $cz,$ where $a,$ $b,$ $c$ are the
	lengths of the sides.\\
	These are given in terms of $a,$ $b,$ $c$ and the trigonometric
	functions of the angles of the triangle.  To obtain our barycentric
	coordinates it is sufficient to replace in $\alpha,$ $\beta$ and
	$\gamma,$\\
	\hth$a^2$ by $m_0(m_1+m_2),$\htn or $a$ by $a0 = j\:j_0(j_1+j_2),$\\
	\hth$b^2$ by $m_1(m_2+m_0),$\htn or $b$ by $a1 = j\:j_1(j_2+j_0),$\\
	\hth$c^2$ by $m_2(m_0+m_1),$\htn or $c$ by $a2 = j\:j_2(j_0+j_1),$\\
	and to replace the trigonometric functions as follows:\\
	\hth	$sin A$ by $s a_0,$\htn$cos A$ by $c \frac{m_1+m_2}{a_0},
	\htn tan A$ by $t m_0,$\\
	\hth	$sin B$ by $s a_1,$\htn$cos B$ by $c \frac{m_2+m_0}{a_1},
	\htn tan B$ by $t m_1,$\\
	\hth	$sin C$ by $s a_2,$\htn$cos C$ by $c \frac{m_0+m_1}{a_2},
	\htn tan C$ by $t m_2,$\\
	where\\
	$p_{11} = j_1j_2+j_2j_0+j_0j_1,\\
	j^2 = (p_{11}-j_0^2)(p_{11}-j_1^2)(p_{11}-j_2^2),\\
	s^2 = \frac{m_0+m_1+m_2}{(m_1+m_2)(m_2+m_0)(m_0+m_1)}
	    = p_{11} (\frac{2}{j(j_1+j_2)(j_2+j_0)(j_0+j_1)})^2,\\
	c^2 = \frac{m_0m_1m_2}{(m_1+m_2)(m_2+m_0)(m_0+m_1)},$
	$c = \frac{j}{(j_1+j_2)(j_2+j_0)(j_0+j_1)},\\
	t = \frac{s}{c}.$\\
	(Vigari\'{e}'s notation is here given between quotes.\\
	Twice the area "2S" by $a_0 a_1 a_2 s$,\\
	$m_0+m_1+m_2 = 4j_0j_1j_2p_{11}$,\\
	$m_0m_1m_2 = j_0j_1j_2jp^2$,\\
	$("2S")^2 = m_0m_1m_2(m_0+m_1+m_2) = (2j_0j_1j_2\:jp)^2 p_{11},$\\
	the radius of the inscribed circle "r" by
	$r^2 = (j_0j_1j_2)^2/p_{11}.$\\
	Moreover\\
	$a+b+c = 2j\:p_{11}$,\\
	$b+c-a = j\:j_0(j_1+j_2)$,\\
	$b^2 - c^2 = m_0(m_1-m_2),\\
	b + c = j(p_{11} + j_1j_2),\\
	b^2 + c^2 - a^2 = 2 m_1m_2.$

	The coordinates are given for the following points,
	I give first Vigari\'{e}'s notation and under it my notation.
	$\begin{array}{lllllllllll}
	G& K&        H&O&H_o&\Omega_1&\Omega_2&I&I_a&I_b&I_c\\
	M& K& \ov\{M\}&O&MAI&Br&      Br~&     I&I_0&I_1&I_2\\
	\hline
	O_9&I_c({\bf 26})&\nu&\Gamma&\nu`_1&\nu'_b&\nu'_c&\Gamma'_a&
		\Gamma'_b&\Gamma'_c\\
	EE&En&	   EE&      N& J \\
	\hline
	I_o&{\bf 32}&J_\delta& J_\rho&N  &R  &\rho    &\rho'&V&W&\\
	   &        &        &       &Tar&Ste&BR\ov{a}&Bra\\
	\hline
	V_2&W_2&P  &P_2& D&D_2&Z  &A_1  &B_1  &C_1\\
	   &   &Tbb&   &  &Tnn&Bro&Br1_0&Br1_1&Br1_2\\
	\hline
	A_2  &B_2  &C_2  &A_3  &B_3  &C_3  \\
	Br3_0&Br3_1&Br3_2&     &     &     &Br2_0&Br2_1&Br2_2\\
	\hline
	{\bf 58}&{\bf 59}&{\bf 60}&{\bf 61}&{\bf 62}&{\bf 63}&{\bf 64}&
		\delta_0&\delta&{\bf 67}\\
	mm_i    &m~m_i   &j_i     &mf_i    &\ov{m}  &   aia &     at_i&
		lem      &    o&     ok\\
	\hline
	{\bf 68}&{\bf 69}&\Sigma_1&\Sigma'_1&\Sigma''_1&\Sigma_2&IO&
		KH_o&HH_o\\
	bbr     &eul    &     &    &  &    &
		    &    \\
	\end{array}$

	In our notation we have,\\
	$"D" = (m1m2(m2+m0)(m0+m1),\ldots)$,\\
	$"I_0" = (j_1j_2(j_2+j_0)(j_0+j_1),  $\ldots$  )\\
	"J_\delta" = (j_0j_1(j_1+j_2)(j_2+j_0),j_1j_2(j_2+j_0)(j_0+j_1),
		j_2j_0(j_0+j_1)(j_1+j_2)),\\
	"J_\rho" = (j_2j_0(j_0+j_1)(j_1+j_2),j_0j_1(j_1+j_2)(j_2+j_0),
		j_1j_2(j_2+j_0)(j_0+j_1)),\\
	"P2" = inverse("P")\\
	"P2" = ((s_{11}+m_2m_0)(s_{11}+m_0m_1), \ldots )\\
	"41" = (m_1m_2(m_2+m_0)(m_0+m_1), \ldots ).$\\

	The equations are given for the following lines,
	$m_i,$ $\ov{m}_i,$ $X_{i+1} \times X_{i-1},$ where $X_i = a_i\times ai_i,$
	$\ov{m}_i],$ $\ov{\i},$ $I_{i+1} \times I_{i-1},$ $at_i,$ "65", "66",
	ok, "line of Brocard" $:= Br1 \times Br2,$ e, "70", "71", "72".

	The equations of the circles.\\
	$\theta$, $\iota$, $\iota_i$, $\gamma$,
	$"78",polar circle:
	m1m2X0^2+m2m0X1^2+m0m1X2^2 = 0, could only find the obvious I, and I[i] 
		on it
	 "79",anticomplementary of "78":
	(X0+X1+X2)(m_0(m_1+m_2)X0 + m_1(m_2+m_0) X1 + m_2(m_0+m_1)X2)
	-m_0(m_1+m_2)X1X2 + m_1(m_2+m_0) X2X0 + m_2(m_0+m_1)X0X1) = 0,
	m_0(m_1+m_2)X0^2+ m_1(m_2+m_0)X1^2 + m_2(m_0+m_1)X2^2
	+ 2 (m1m2 X1X2 + m2m0 X2X0 + m0m1 X0X1) = 0$
	$\beta,$,"81" family of Circles of Tucker",
	$\lambda1,$ $lambda2,$ "Circle of Taylor `Tay?",
	"Circles of Neuberg" (D35.4),
	"86" "Circles of M'Cay"
	"87":`alphap[i] "Circles of Apolonius",
	"88": family of "Circles of Schoute",

	The equations are given for other curves,

	The conic of Brocard (D36.19), the conic of Lemoine (D36.7),
	also mentioned by Neuberg, m\'{e}moire sur le t\'{e}tra\`{e}dre, p.5 VI,
\\
	$iM := I \times M,$	$\ov{\i}M := \ov{I}\times \ov{M}$\\ 
	$i'M := I' \times M,$	$\ov{\i}':= \ov{I}'\times \ov{M}$\\ 
	$\ldots ^{-1} \cdot iM = \ldots ^{-1} \cdot i'M
		= \ldots ^{-1} \cdot \ov{i} = \ldots ^{-1}\cdot \ov{i}'M$.\\
	Hence the foci of $\ldots$  are M and K and\\
	the cofoci of $\ov{\ldots}$ are $\ov{M}$ and K.\\
	$MK = (5s_{11}-3m_1m_2,5s_{11}-3m_2m_0,5s_{11}-3m_0m_1),$\\
	$\ov{M}K = (m_0(5s_1-3m_0),m_1(5s_1-3m_1),m_2(5s_1-3m_2)),\\
	\ldots^{-1}: (s_{11}+3m_1m_2)x_1x_2 + (s_{11}+3m_2m_0)x_2x_0
		+ (s_{11}+3m_0m_1)x_0x_1 = 0.\\
	\ov{\ldots}^{-1}: (5m_1m_2-m_0^2+m_1^2+m_2^2)x_1x_2 +  \ldots  = 0.$\\
	points of contact$_0 = ((s_{11}+3m_2m_0)(s_{11}+3m_0m_1), \ldots$).\\
	These are the feet of the symedians of $A_1 MA_2.$\\
	points of contact$_0 = ((5m_2m_0+m_0^2-m_1^2+m_2^2)
		(5m_0m_1+m_0^2+m_1^2-m_2^2),$  \ldots).\\[10pt]
	the conic K (D36.2),\\[10pt]
	the conic of Simmons (92),\\[10pt]
	the conic of Steiner (S36.3),\\[10pt]
	the "hyperbola" of Kiepert (D16.19),\\[10pt]
	the first parabolas of Artzt (D36.8),\\[10pt]
	The second parabolas of Artzt (96):\\[10pt]
	Artzt2:\\[10pt]
	The parabolas of Brocard (97):\\[10pt]
	Brocard1:\\[10pt]
	Brocard2:\\[10pt]
	Focus(Kiepert2).theta = 0.\\[10pt]
	The conic of Jerabek (99):\\[10pt]
	Jerabek = inverse(e)\\
	Jerabek: $m_0(m_1^2-m_2^2) X_1 X_2 + m_1(m_2^2-m_0^2) X_2 X_0
			+ m_2(m_0^2-m_1^2) X_0 X_1 = 0.$\\[10pt]
	The conic centrally associated to a point (99'):\\
	conic(X).  Given X = $(X_0,X_1,X_2),$\\
	    Let $X_i := a_i \times (A_i \times X),$\\
	    conic(X) := conic($X_i, \times pole of i$),\\
	conic(X): $(-X_0+X_1+X_2)X_1^2X_2^2$ $X_0^2$ +  $\ldots$ 
		$-2(X_0^3X_1X_2$ $X_1 X_2 +  \ldots  ) = 0.$\\
	I := conic(I),\\[10pt]
	"I" no point on it
	The conic I (100):\\
	supplementary($\theta$) = I,\\
	I: $(j_1j_2)^3 ((j_0+j_1)(j_2+j_0))^2$ $X_0^2$ +  $\ldots$ 
		-2( $(j_0(j_1+j_2))^3 j_1j_2(j_0+j_1)(j_2+j_0) X_1 X_2
	 +  \ldots$  ) = 0\\[10pt]
	or $(p_{11}^2+p_1p_{111})(X_0^2+X_1^2+X_2^2)-2(j_0^2(j_1+j_2)^2 X_1X_2
		+ j_1^2(j_2+j_0)^2 X_2X_0 + j_2^2(j_0+j_1)^2 X_0X_1) = 0$\\
	the conics of Simson (D16.18),\\[10pt]

	m1m2 X1X2 + m2m0 X2X0 + m0m1 X0X1= 0 no point on it

	\ssec{The cubic of 17 points.}
	\setcounter{subsubsection}{-1}
	\sssec{Introduction.}
	The cubic of 17 points is defined without explicit reference by
	Vigari\'{e}. It can be defined as the cubic through the vertices of
	a triangle, its midpoints and the midpoints Mma$_i$ between the vertices
	and the feet. The other 8 points are the barycenter, orthocenter,
	center of the outscribed circle, point of Lemoine, and the 4 centers
	of the tangent circles.  Other points and tangent on it will also be 
	given. In particular, $\ov{K}LL_i$, Flor, $\ov{A}RTM$, are on the cubic
	and at$_i$ is the tangent at A$_i$, 
	mf$_i$ is the tangent at M$_i$, mk is the tangent at M,
	ok is the tangent at K.

	\sssec{Definition.}
	The {\em cubic of 17 points} is defined by\\
	"cubic17 := cubic(A$_i$,M$_i$,Mma$_i$).

	\sssec{Theorem.}
	O $\cdot$"cubic17 = M $\cdot$"cubic17 = $\ov{M}\cdot$"cubic17 = 
	K$\cdot$"cubic17 \\
	\hth = I $\cdot$cubic17 = I$_i\cdot$"cubic17  = 0.\\
	$\ov{K}$LL$_i \cdot$"cubic17 = 0.\marginpar{19.7.82}\\
	$\ov{A}$RTM$ \cdot$"cubic17 = 0.\marginpar{3.7.91}\\
	at$_i\cdot$"cubic17  = 0,
	
	Proof.\\
	"cubic17: $m_0(m_1+m_2) X_1 X_2 (X_1-X_2)
		+ m_1(m_2+m_0) X_2 X_0 (X_2-X_0)\\
	\hth + m_2(m_0+m_1) X_0 X_1 (X_0-X_1) = 0.$

	\sssec{Theorem.}
	$\begin{array}{lcllllllllll}
	&\vline&M&\ov{M}&O&K&Flor&\ov{A}RTM&C17a&C17b&C17c&C17d\\
	\hline
	M&\vline&K&O&\ov{M}&M&\ov{A}RTM&Flor&C17b\\
	\ov{M}&\vline&&C17a&M&\ov{A}RTM&&K&\ov{M}\\
	O&\vline&&&Flor&K&O&&?\\
	K&\vline&&&&O&&\ov{M}&&\\
	Flor&\vline&\ov{A}RTM&?&O&?&M&\\
	\ov{A}RTM&\vline&Flor&K&?&\ov{M}&M&\\
	C17a&\vline&C17b&\ov{M}&?&&&&?&&\\
	C17b&\vline&C17a&&&&&\\
	C17c&\vline&\ov{A}RTM&C17d&?&?&M&?&?\\
	\end{array}$\\
	$\begin{array}{lcllllll}
	&\vline&A_i&M_i&Mma_i&\ov{K}LL_i\\
	\hline
	M&\vline&M_i&A_i&\ov{K}LL_i&Mma_i&\\
	\ov{M}&\vline&Mma_i&C17d_i&A_i&C17b_i&\\
	O&\vline&\ov{K}LL_i&M_i&C17c_i&A_i&\\
	K&\vline&A_i&Mma_i&M_i&C17d_i&\\
	Flor&\vline&C17d_i&C17b_i&Mma_i&C17e_i&\\
	\ov{A}RTM&\vline&C17c_i&\ov{K}LL_i&C17f_i&M_i&\\
	C17a&\vline&C17f_i&C17c_i&?&\ov{K}LL_i\\
	C17b&\vline&C17b_i&C17e_i&C17d_i&C17g_i\\
	C17c&\vline&C17d_i&C17b_i&?&?\\
	\end{array}$\\
	$\begin{array}{lcllllll}
	&\vline&A_i&M_i&Mma_i&\ov{K}LL_i\\
	\hline
	A_i&\vline&K&M&\ov{M}&O&\\
	M_i&\vline&M&O&K&\ov{A}RTM&\\
	Mma_i&\vline&\ov{M}&K&C17c_i&M&\\
	\ov{K}LL_i&\vline&O&\ov{A}RTM&M&C17a&\\
	\end{array}$\\
	$\begin{array}{lcllllll}
	&\vline&A_{i-1}&M_{i-1}&Mma_{i-1}&\ov{K}LL_{i-1}\\
	\hline
	A_{i+1}&\vline&M_i&A_i&\ov{K}LL_i&Mma_i\\
	M_{i+1}&\vline&A_i&Mma_i&M_i&C17d_i&\\
	Mma_{i+1}&\vline&\ov{K}LL_i&M_i&C17c_i&A_i&\\
	\ov{K}LL_{i+1}&\vline&Mma_i&C17d_i&C17b_i&M_i\\
	\end{array}$

	\sssec{Exercise.}\label{sec-ecubic17}
	Construct the tangents to "cubic17 at Mma$_i$ and at $\ov{M}$.

	\sssec{Exercise.}\label{sec-ecocubic17}
	Give properties of
	$\ov{"}$cubic17 := cubic(A$_i$,$\ov{M}_i$,$\ov{M}$ma$_i$).

	\sssec{Partial answer to}
	\vspace{-18pt}\hspace{120pt}\ref{sec-ecubic17}\\[8pt]
	The tangent at Mma$_0$ is\\
	$[(m1-m2)(s1-2m0),(m1+m2)(s1-2m1),-(m1+m2)(s1-2m2)].$\\
	The tangent at $\ov{M}$ is\\
	$[m1m2(m1-m2)(s1-2m0),m2m0(m2-m0)(s1-2m1),\\
	\hth m0m1(m0-m1)(s1-2m2)].$\\
	A$_0\oplus $A$_0$ = K, A$_1\oplus $A$_2$ = M$_0$, A$_0\oplus $M$_0$ = M,
	A$_0\oplus $Mma$_0 = \ov{M}$, A$_1\oplus $Mma$_2 = \ov{K}$LL$_0$,\\ 
	A$_0\oplus \ov{K}$LL$_0$ = O,  
	M$_0\oplus $M$_0$ = K, M$_1\oplus $M$_2$ = Mma$_0$,\\
	M$_0\oplus \ov{K}$LL$_0 = \ov{A}$RTM,
	M$_1\oplus \ov{K}$LL$_2$ = C17a$_0$,\\
	O$\oplus $O = Flor, where\\
	Flor = $((m1+m2)(s1-2m1)(s1-2m2),(m2+m0)(s1-2m2)(s1-2m0),\\
	\hth (m0+m1)(s1-2m0)(s1-2m1))$,\\
	tangent at $\ov{M}$ is $[m1m2(m1-m2)(s1-2m0),\ldots]$.\\
	$\ov{M}\times \ov{M} = (m0(m1+m2)(s1-2m1)(s1-2m2)(s22-2m0^2),\ldots),$\\
	$M\times Flor = [(m1-m2)(s1-2m0),(m2-m0)(s1-2m1),(m0-m1)(s1-2m2)],$\\
	$\ov{M}\times \ov{A}RTM = [(m1-m2)(s1-2m0)(s2-2m0^2),\ldots],$\\
	$\ov{M}\oplus \ov{A}RTM = K,$\\
	$K\times Flor = [(m1-m2)(m2+m0)(m0+m1)(s1-2m0)^2,\ldots],$\\
	$C17a = (m0(m1+m2)(s1-2m1)(s1-2m2)(s2-2m0^2),\ldots),$
	$C17b = ((s1-2m0)(s2-2m1^2)(s2-2m2^2),\ldots),$\marginpar{4.7.91}\\
	$C17c = ((m1+m2)(s1-2m1)(s1-2m2)),\ldots),$\marginpar{4.7.91}\\
	$C17d = ((s2-2m1^2)(s2-2m2^2)(s3-s21-2s2m0),\ldots),$\marginpar{4.7.91}\\
	$C17b_0 = (2m0(m1+m2)(s1-2m1)(s1-2m2),s1(s1-2m1)(s2-2m2^2),\\
	\hth s1(s1-2m2)(s2-2m1^2)),$\marginpar{4.7.91}\\
	$C17c_0 = ((m1+m2)(s1-2m1)(s1-2m2),m1s1(s1-2m1),m2s1(s1-2m2)),$\\
	$C17d_0 = (m0s1,(m2+m0)(s1-2m2),(m0+m1)(s1-2m1)),$\\
	$C17e_0 = ((s2-2m1^2)(s2-2m2^2)(s3-s21-2s111),\\
	\hth -2m1(m2+m0)(s2-2m1^2)(s3-2m2^3-s21+2m2(m0^2+m1^2)+2s111),
	\hth -2m2(m0+m1)(s2-2m2^2)(s3-2m1^3-s21+2m1(m2^2+m0^2)+2s111)),$\\
	$C17f_0 = (s1(s2-2m2^2)(s2-2m1^2),2m1(m2+m0)(s1-2m2)(s2-2m1^2),\\
	\hth -2m2(m0+m1)(s1-2m1)(s2-2m2^2)),$\\
	$C17g_0 = ((m1+m2)s1(s1-2m0)(s1-2m1)(s1-2m2)
		(s3-s21+2s111+2m2(s2-2m2^2))((s3-s21+2s111+2m1(s2-2m2^2))),
	\hth -m1(s3-s21-2s111)(s3-s21+2s111+2m1(s2-2m1^2))
		(s4+2s22-4(m2^4+m0^2m1^2)),
	\hth -m2(s3-s21-2s111)(s3-s21+2s111+2m2(s2-2m2^2))
		(s4+2s22-4(m1^4+m2^2m0^2))$,\marginpar{4.7.91}

	\sssec{Answer to}
	\vspace{-18pt}\hspace{92pt}\ref{sec-ecocubic17}\\[8pt]
	$\ov{O}\cdot\ov{"}$cubic17 = M$\cdot\ov{"}$cubic17  =
	$\ov{M}\cdot\ov{"}$cubic17 = K $\cdot\ov{"}$cubic17= 0.\\
	KLL$_i\cdot\ov{"}$cubic17 = 0.\marginpar{19.7.82}\\
	at$_i\cdot\ov{"}$cubic17  = 0,\\
	$\ov{"}$cubic17:\\
	$m_0^2(m_1+m_2)$ $X_1 X_2$ $(m_2 X_1 - m_1 X_2)$
		+ $m_1^2(m_2+m_0)$ $X_2 X_0$ $(m_0 X_2 - m_2 X_0)$\\
	\hth + $m_2^2(m_0+m_1)$ $X_0 X_1$ $(m_1 X_0 - m_0 X_1) = 0.$\\[10pt]
	\ssec{The cubic of 21 points.}
	cubic21\\[10pt]
	\ssec{The Barbilian Cubics.}
\setcounter{subsubsection}{-1}
	\sssec{Introduction.}
	In posthumously published works of Dan Barbilian, also known in his
	native Roumanian Country as the poet Eon Barbu, the following Theorem is
	proven.  The loci of the pseudo centers of the isotropic cubics which 
	pass through the vertices of a complete quadrilateral and 2 of its
	diagonal elements is a circle.  I observed that in the case where the
	isotropic points  are the fixed points of the involution determined by
	the 3 pairs of opposite sides of the quadrilateral, the third diagonal
	point is also on the cubics.  It is this family of cubics which will be
	studied now, to which I will give the name of the Poet-Mathematician
	Barbilian.

	\sssec{Definition.}
	An {\em isotropic cubic} is a cubic which passes through the isotropic
	points.\\
	The {\em pseudo center} of an isotropic cubic is the intersection of
	its tangents at the isotropic points.\footnote{The isotropic points are
	also called circular points. Barbilian calls a pseudo center, a pseudo
	focus.}

	\sssec{Theorem. [Barbilian]}
	{\em The family of isotropic cubic through the vertices $B_j$ of a
	complete quadrangle and 2 of its diagonal points
	$A_1 := (B_0\times B_2)\times(B_1\times B_3)$ and
	$A_2 := (B_0\times B_1)\times(B_2\times B_3)$ has a circle as the locus
	of the pseudo centers.}  This circle is the {\em Miquel circle} of
	the complete quadrangle and the 2 diagonal points.

	I remind the reader that this circle passes through the center of the
	circumcircles of the triangles $\{B_0,B_1,A_1\}$, $\{B_2,B_3,A_1\}$,
	$\{B_0,B_2,A_2\}$, $\{B_1,B_3,A_2\}$. See g334

	\sssec{Definition.}
	The isotropic cubics through the vertices of a triangle, the feet
	and the orthocenter will be called {\em Barbilian cubics.}

	\sssec{Corollary.}
	{\em The family of Barbilian cubics has a circle as the locus of its 
	pseudo centers.}

	In this case,  $B_0 = A_0$, $B_1 = \ov{M}_1$, $B_2 = \ov{M}_2$,
	$B_3 = \ov{M}$ and the circles circumscribed to $\{B_2,B_3,A_1\}$,
	$\{B_1,B_3,A_2\}$ pass through the point of Miquel, $\ov{M}_0$.

	\sssec{Theorem.}
	{\em The Miquel circle of $B_j$, $A_1$ and $A_2$ is the circle of
	Brianchon-Poncelet.}

	\sssec{Theorem.}\label{sec-tbarb0}
	The following are degenerate Barbilian cubics.
	\enumb
	\item $`Aam_0 \bigx \ov{m}a_0$, its equation is\\
	\hth$(m_1+m_2)\:m_0 X_1X_2(m_2X_1-m_1X_2)
	- m_1\:m_1 X_2X_0(m_0X_2-m_2X_0)\\
	\hti{12}- m_2\:m_2 X_0X_1(m_1X_0-m_0X_1) = 0.$
	\item $`Mma_0 \bigx a_0$, its equation is\\
	\hth $m_1 X_2X_0(m_0X_2-m_2X_0) - m_2 X_0X_1(m_1X_0-m_0X_1) = 0.$
	\enume

	Proof: 0 and 1, follow from the definition of the circles `Aam$_i$ and
	`Mma$_i$ given in section D11.1 and .2.

	\sssec{Theorem.}\label{sec-tbarb1}
	\enumb
	\item {\em The cubics through $A_i$, $\ov{M}_i$, $\ov{M}$, are}\\
	\hth$a_0m_0X_1X_2(m_2X_1-m_1X_2) + a_1m_1X_2X_0(m_0X_2-m_2X_0)\\
	\hti{12}	+ a_2m_2X_0X_1(m_1X_0-m_0X_1) = 0.$
	\item {\em A necessary and condition for the cubics through $A_i$,
	$\ov{M}_i$, $\ov{M}$, to be Barbilian cubics, is}\\
	\hth$a_0 + a_1 + a_2 = 0$.
	\enume

	Proof:
	It is easy to verify 0. For 1, any Barbilian cubic is a linear 
	combination of the degenerate cubics given in the preceding Theorem
	and this satisfy the given condition.

	\sssec{More details on}
	\vspace{-18pt}\hspace{132pt}\footnote{\ref{sec-tbarb1}}\\[8pt]

	Recall that the isotropic points are\\
	\hth$(m_0(m_1+m_2),-m_0m_1-j\tau,-m_2m_0+j\tau),$ with $j = \pm 1$ and
	$\tau^2 = -m_0m_1m_2s_1$.

	\sssec{Theorem.}\label{sec-tbarb2}
	{\em In homogeneous Cartesian coordinates $(X,Y,Z),$ with\\
	$A_0 = (0,h,1),$ $A_1 = (b,0,1),$ $A_2 = (c,0,1)$ and isotropic points
	$(\pm j,1,0)$, we have the following.
	\enumb
	\item {\em The coordinates of the sides, feet, orthocenter and altitudes
	are}\\
	\hth$a_0 = [0,1,0],$ $a_1 = [h,c,-ch],$ $a_2 = [h,b,-bh],$\\
	\hth$\ov{M}_0 = (0,0,1),$ $\ov{M}_1 = (c(h^2+bc),hc(c-b),h^2+c^2),$
		$\ov{M}_2 = (b(h^2+bc,hb(b-c),h^2+b^2),$\\
	\hth$\ov{M} = (0,bc,-h),$\\
	\hth$\ov{m}_0 = [1,0,0]$, $\ov{m}_1 = [c,-h,-bc]$,
		$\ov{m}_2 = [b,-h,-bc]$.
	\item {\em The circle through $A_1$, $A_2$, $\ov{M}_1$, $\ov{M}_2$ is}\\
	\hth$\alpha m_0:\:X^2 + Y^2 + bc Z^2 - (b+c)ZX = 0.$
 	\item {\em The circle through $A_0$, $\ov{M}$, $\ov{M}_1$, $\ov{M}_2$
	is}\\
	\hth$\mu a_0:\:h(X^2 + Y^2) - hbc Z^2 + (bc-h^2)XY = 0.$
	\item {\em The Barbilian cubics are}\\
	\hth$k\alpha m_0\bigx X + l\mu a_0\bigx Y = 0$.
	\item {\em The pseudo center is, with} $d = k(b+c)+l(bc-h^2)$,\\
	\hth$(kd,-lhd,2(k^2+l^2h^2)).$
	\item {\em This is a parametric equation of the circle of
	Brianchon-Poncelet}:\\
	\hth$2h(X^2+Y^2)-(h^2-bc)XY-h(b+c)ZX = 0.$
	\item {\em The transformation from barycentric to Cartesian coordinates
	is\\
	$\mat{0}{b}{c}{h}{0}{0}{1}{1}{1}$ and its inverse is $\frac{1}{h(c-b)}
	\mat{0}{c-b}{0}{-h}{-c}{ch}{h}{b}{-bh},$\\
	where the barycentric coordinates of the orthocenter are given by\\
	\hth$m_0 = \frac{bc(c-b)}{h^2+bc},$ $m_1 = -c,$ $m_2 = b$,\\
	provided} $h^2 = \frac{m_1m_2s_1}{m_0}.$
	\item {\em The barycentric coordinates of the pseudo center are easily
	derived using the value of $h$ and of}\\
	\hth$m_0d = km_0(m_2-m_1)-l(s_{111}-m_1m_2s_1)$.
	\enume

	The details of the proof is left to the reader.

	\sssec{Answer to}
	\vspace{-18pt}\hspace{92pt}\footnote{\ref{sec-tbarb2}}\\[8pt]
	0, is straigthforward.\\
	For 1, $\alpha m_0$ is $u[h,b,-hb]\bigx[h,c,-hc] = v[c,-h,-bc]\bigx
	[b,-h,-bc].$\\
	With $j^2 = -1,$ $u = (cj-h)(bj-h),$ $v = (hj+b)(hj+c) = - u.$\\
	After dividing by $h^2+bc$ we get the equation 1.\\
	For 2, $\mu a_0$ is $u[h,b,-hb]\bigx[c,-h,-bc] = v[h,c,-hc]\bigx
	[b,-h,-bc].$\\
	$u = (hj+c)(bj-h),$ $v = (hj+b)(cj-h) = u.$\\
	After dividing by $c-b$ we get the equation 2.\\
	For 4, the tangent at $(j,1,0)$, obtained by evaluating the first
	partial derivatives at that point, is\\
	\hth$[-2k+2lhj,2kj+2lh,k(b+c)+l(bc-h^2)].$\\
	The tangent at the other isotropic point is obtained by replacing
	$j$ by $-j$.\\
	Their intersection is 4, after dividing by $4j$.
	It is easy to verify that 5, is the equation of a circle through
	$\ov{M}_i$, and that the pseudo center is on it for all values of $k$
	and $l$.\\
	For 5, the coordinates of vertices $A_i$ give the coefficients of the
	matrix. The transform of $(m_0,m_1,m_2)$ is $(m_1b+m_2c,hm_0,s_1)$.
	Comparing with $(0,bc,-h)$ gives 5.
	6, is straightforward.
	For 7, we need to related the cubics in Cartesian and barycentgric
	coordinates.If we use k' and l' for the barycentric case,
	comparison of the coefficients of $X_2^2X0$ in baricentric coordiantes
	gives\\
	$k'm_1^2m_0 = kc(2bc-bc-c^2) = kc^2(b-c),$ and
	$l'm_0m_1 = lh(hc^2-hbc) = lh^2(c-b).$  Using proportionality we can
	therefore write 7.\\
	The pseudo center will be (0,0,1) if $d = 0,$ this gives\\
	$k = h^2-bc = \frac{m_1m_2(s_1+m_0)}{m_0},$ $l = m_2-m_1$.\\
	Substituting, we get, after division by $m_1m_2m_0^{-1}$,
	$k' = s_1+m_0,$ $l' = (m_1-m_2)s_1,$ hence\\
	$a_0 = -(m_1+m_2)(s_1+m_0),$ $a_1 = (2m_1-m_2)s_1+m_0m_1$,
	$a_1 = (2m_1-m_2)s_1+m_0m_1$.
	To check this independently, we should verify that $M_0\times I_0$
	is tangent to the cubic for these values of $k'$ and $l'$.
	$(0,0,1)\times (m_0(m_1+m_2),-m_0m_1-j\tau,-m_2m_0+j\tau)
	= [m_0m_1+j\tau,m_0(m_1+m_2),0].$

	\sssec{Theorem.}
	{\em Given an Barbilian cubic $\Gamma$, there exists a line $l$
	and a circumscribed conic $\phi$ such that}\\
	\hth$\Gamma = \theta \bigx l + m \bigx \phi$.\\
	More specifically, with $l_0$ arbitrary,\\
	\hth$l = [l_0,l_0-a_2,l_0+a_1]$,\\
	\hth$\phi = b_0m_0X_1X_2 + b_1m_1X_2X_0 + b_2m_2X_0X_1 = 0$, with\\
	\hth$b_0 = -m_2a_1+m_1a_2-(m_1+m_2)l_0,$ $b_1 = -m_2a_1-(m_2+m_0) l_0$,
	\hth$b_2 = m_1a_2-(m_0+m_1)l_0.$

	Proof:
	Identification of the coefficients of $X_1^2X_2$ and $X_2^2X_1$ gives\\
	$a_0m_2 = (m_1+m_2)l_1 + b_0,$
	$-a_0m_1 = (m_1+m_2)l_2 + b_0,$\\
	subtracting gives, $a_0 = l_1-l_2,$ and similarly $a_1 = l_2-l_0,$
	and $a_2 = l_0-l_1.$\\
	By substitution, we obtain $b_0$ and similarly $b_1$ and $b_2$, using
	$a_0+a_1+a_2 = 0$.

	\sssec{Definition.}
	$l$ is called a {\em radical axis} of $\Gamma$, $\phi$ is
	called the {\em corresponding radical conic} of $\Gamma$.

	$\theta$ could be replaced by an other circle.

	\sssec{Theorem.}
	\enumb
	\item {\em The non trivial ideal point is} $(a_0,a_1,a_2)$.
	\item {\em The tangent at the non trivial ideal point, or
	asymptote is}\\
	\hth$[m_0a_1a_2(m_2a_1-m_1a_2),m_1a_2a_0(m_0a_2-m_2a_0),
		m_2a_0a_1(m_1a_0-m_0a_1)].$
	\enume

	Proof: This follows from the fact that the non trivial ideal point is
	$m\times l$. The tangent is obtain by taking the partial derivatives
	respectively with respect to $X_0$, $X_1$ and $X_2$ at $(a_0,a_1,a_2)$.
	The first one is\\
	$m_1a_1a_2(m_0a_2-2m_2a_0) - m_2a_1a_2(m_0a_1-2m_1a_0)
	= -m_0a_1a_2(m_2a_1-m_1a_2).$

	\sssec{Comment.}
	Special Barbilian cubics can be obtained by combining the equations
	of Theorem \ref{sec-tbarb0}  For instance,\\
	$\beta a0$ follows from adding the equations 0, for indices 0,1 and 2
	respectively multiplied by $m_0$, $m_1$ and $m_2$.\\
	$\beta a1$, by using on the equations 0, the multipliers $m_1m_2$,
	$m_2m_0$ and $m_0m_1$.\\
	$\beta a2$, by using on the equations 1, equal multipliers.

	\sssec{Theorem.}
	{\em The Barbilian cubic $\beta a0:$\\
	\hth$ m_0(m_1-m_2)X_1X_2(m_2X_1-m_1X_2)
		+ m_1(m_2-m_0)X_2X_0(m_0X_2-m_2X_0)\\
	\hti{12}	+ m_2(m_0-m_1)X_0X_1(m_1X_0-m_0X_1) = 0$\\
	has the following properties:}
	\enumb
	\item {\em A radical axis is $-mai$, with} $mai = [m_0,m_1,m_2].$
	\item {\em The corresponding radical conic has the equation}\\
	\hth $m_0(m_1^2+m_2^2)X_1X_2 + m_1(m_2^2+m_0^2)X_2X_0
		+ m_2(m_0^2+m_1^2)X_0X_1.$
	\item {\em The non trivial ideal point is}
	$MK = (m_1-m_2,m_2-m_0,m_0-m_1).$
	\item {\em The asymptote is}\\
	\hth$[m_0(m_2-m_0)(m_0-m_1)(s_2-m_0s_1),m_1(m_0-m_1)(m_1-m_2)
		(s_2-m_1s_1),m_2(m_1-m_2)(m_2-m_0)(s_2-m_2s_1)].$
	\item {\em The tangent at $A_i$ is $mka_i$, with}\\
	\hth$mka_0 = [0,m_0-m_1,-(m_2-m_0)].$
	\enume

	\sssec{Theorem.}
	{\em The Barbilian cubic $\beta a1:$\\
	\hth$m_0^2(m_1-m_2)X_1X_2(m_2X_1-m_1X_2)
		+ m_1^2(m_2-m_0)X_2X_0(m_0X_2-m_2X_0)\\
	\hti{12}	+ m_2^2(m_0-m_1)X_0X_1(m_1X_0-m_0X_1) = 0$\\
	has the following properties:}
	\enumb
	\item {\em A radical axis is} $-\ov{m}$.
	\item {\em The corresponding radical conic is $2m_0m_1m_2`\ov{S}teiner,
	$ with}\\
	\hth$	`\ov{S}teiner = m_0 X_1X_2 + m_1 X_2X_0 + m_2 X_0X_1 = 0,$
	\item {\em The non trivial ideal point is $EUL$ with}\\
	\hth$ EUL = (m_0(m_1-m_2),m_1(m_2-m_0),m_2(m_0-m_1)).$
	\item {\em The asymptote is $\ov{m}$}.
	\item {\em The tangent at $A_i$ is $mka_i$, with}\\
	\hth$mka_0 = [0,m_0-m_1,-(m_2-m_0)].$
	\enume

	\sssec{Theorem.}
	{\em The Barbilian cubic $\beta a2:$\\
	\hth$m_0(s_1-3m_0)X_1X_2(m_2X_1-m_1X_2)
		+ m_1(s_1-3m_1)X_2X_0(m_0X_2-m_2X_0)
	\hti{12}	+ m_2(s_1-3m_2)X_0X_1(m_1X_0-m_0X_1) = 0$\\
	has the following properties:}
	\enumb
	\item {\em A radical axis is $eul$ with} $eul = [m_1-m_2,m_2-m_0,m_0-m_1].$
	\item {\em The corresponding radical conic is}\\
	\hth$m_0^2(m_1-m_2)X_1X_2 + m_1^2(m_2-m_0)X_2X_0
		+ m_2^2(m_0-m_1)X_0X_1 = 0$,
	\item {\em The non trivial ideal point is $Ieul$, where}
	$Ieul = (s_1-3m_0,s_1-3m_1,s_1-3m_2).$
	\item {\em The asymptote is} $Mkm \times Ieul,$\\
	$[m_0(m_1-m_2)(s_1-3m_1)(s_1-3m_2),m_1(m_2-m_0)(s_1-3m_2)(s_1-3m_0),
		m_2(m_0-m_1)(s_1-3m_0)(s_1-3m_1)].$
	\item {\em The tangent at $A_i$ is} 
	\enume

	Several mappings are defined and these allow an algebraic definition
	of many of the points, these will be given here as theorems for the
	points already defined and as definition for the others.\\
	reciprocal$(X_0,X_1,X_2) := (X_1X_2,X_2X_0,X_0X_1),$\\
	reciprocal$(x_0,x_1,x_2) := (x_1x_2,x_2x_0,x_0x_1),$\\
	inverse$(X_0,X_1,X_2) := (m_0(m_1+m_2)X_1X_2,m_1(m_2+m_0)X_2X_0,
		m_2(m_0+m_1)X_0X_1),$\\
	complementary$(X_0,X_1,X_2) := (X_1+X_2,X_2+X_0,X_0+X_1),$
	[Nagel, 1885]\\
	anticomplementary$(X_0,X_1,X_2) := (-X_0+X_1+X_2,X_0-X_1+X_2,
		X_0+X_1-X_2),$[inverse transformation of de Longchamps, 1886]\\
	supplementary$(X_0,X_1,X_2) := \ldots$\\
	algebraically associated$(X_0,X_1,X_2) :=
			( (-X_0,X_1,X_2),(X_0,-X_1,X_2),(X_0,X_1,-X_2) ),$\\
	Brocardian$(X_0,X_1,X_2) := ( (X_0X_1,X_1X_2,X_2X_0), 
		(X_2X_0,X_0X_1,X_1X_2) ),$\\
	isobaric$(X_0,X_1,X_2) := ( (X_2,X_0,X_1), (X_1,X_2,X_0) ),$\\
	semi reciprocal$(X_0,X_1,X_2) :=
			( (X_0,X_2,X_1), (X_2,X_1,X_0), (X_1,X_0,X_2) ),$\\
	associated$(X_0,X_1,X_2) := (X_1-X_2,X_2-X_0,X_0-X_1).$

	\sssec{Theorem.}
	\enumb
	\item$	K = inverse(M),$
	\item$	O = inverse(\ov{M}) = complementary(\ov{M})$
	\item$	(Br2,Br2) = brocardian(K),$\marginpar{?}
	\item$	I_i = algebraically\: associated(I),$
	\item$	N = anticomplementary(I),$
	\item$	J = reciprocal(N),$
	\item$	N_i = algegraically\: associated(N),$
	\item$	J_i = reciprocal(N_i).$
	\enume

	\sssec{Definition.}
	The following are the definition of other points.
	\enumb
	\item $H_0 := reciprocal(\ov{M})$
	\item $I_c := complementary(I),$
	\item $I_0 := reciprocal(I),$
	\item "Center of equal parallels" := anticomplememtary(I$_0),$
	\item $("J_\delta","J_\rho")$ := Brocardian(I).
	\enume

	\sssec{Exercise.}
	Define in terms of the above functions as many
	points as you can in Theorem $\ldots$ .

	\sssec{Exercise.}
	Determine, for many of the points of Definition $\ldots$ 
	a linear construction and determine their barycentric coordinates.


	\section{Finite Projective Geometry.}
\setcounter{subsection}{-1}
	\ssec{Introduction.}
	The Theorems given here are deduced from Theorems of Involutive
	Geometry.

	\sssec{Theorem.}\label{sec-tgkimh}
	{\em Given 6 points $A_i$ and $B_i$, forming an hexagon inscribed to a
	conic $\alpha$ and outscribed to an other conic $\beta$. Let $C$ be
	the point common to $A_i\times B_i$. Let $T_i$ be the vertices of the 
	tangents to $\alpha$ at $A_i$.}
	\enumb
	\item {\em The lines $B_i\times T_i$ have a point $D$ in common.}
	\item {\em The line $C\times D$ passes through the pole of
		with respect to the triangle $\{A_i\}$ of the Desargues line
		of the perspective triangles $\{T_i\}$ and $\{B_i\}$.}
	\enume

	The Theorem generalizes a Theorem of Kimberling
	\ref{sec-tkimberling1} and \ref{sec-tkimberling2} using
	\ref{sec-tdevkimb}.

	\sssec{Theorem.}\label{sec-tgkim}
	{\em Given the special Desargues configuration with the points $A_i$ on
	the lines of the triangle $\{MM_0,MM_1,MM_2\}$ with center of
	perspectivity $M$. Let $\ov{m}$ be an arbitrary line and $\ov{M}A_i$ be
	its intersection with the side $a_i$ of the triangle $\{A_0,A_1,A_2\}$,
	if $TMa_i$ is the intersection of the line $MM_{i+1}$ $\ov{M}A_{i-1}$
	and the line $MM_{i-1}$ $\ov{M}A_{i+1}$, then the lines joining the
	points $A_i$ to the $TMa_i$ have a point $ARTM$ in common.}

	The Theorem generalizes a Theorem of Kimberling
	\ref{sec-tekim} assuming that the excenters are replaced by the
	vertices of the anti complimentary triangle and the direction of
	the altitudes are replaced by the intersections of the orthic line with
	the side of the triangle.

	\sssec{Definition.}
	The point $ARTM$ is called the {\em point of Luke.}


	\section{Finite Involutive Geometry.}
\setcounter{subsection}{-1}
	\ssec{Introduction.}
	I will now describe the Theorems of involutive Geometry in the
	traditional way, refering for to proofs of the corresponding
	sections of the hexal configuration.\\
	Starting with affine geometry, we obtain an involutive geometry, if
	we choose among all the possible involutions on the ideal line,
	a particular one, called fundamental involution.  We could also
	start directly from projective geometry and choose among all the 
	possible involutions one involution on one of all the possible lines.\\
	This involution can be given in many ways,
	\enumb
	\item by 2 points, the fixed points of the involution,
	\item by 2 pairs of corresponding points on a line,
	\item by a polarity and a line which does not belong to its line conic,
	\item by an hexal complete 5-angle, \ldots. See II.3.
	\enume
	The definitions will be given in terms of the fundamental
	involution.  Because this involution can be elliptic or hyperbolic,
	there are 2 distinct types of real involutive geometries, elliptic and
	hyperbolic.  I will study them together and give theorems, in the
	hyperbolic case, which in some cases can be used as an alternate
	definition of the concepts.  Such theorems will be noted with
	{\bf (H. D.)}.
	When the additional notions of measure of distance and angles
	will have been introduced, the elliptic involutive geometry
	will become the Euclidean Geometry and the hyperbolic one, will be
	that of Minkowski. A third geometry which corresponds to
	the confluence of the 2 fixed points of the involution will be
	considered later, it is the parabolic (involutive) geometry which
	becomes the Galilean geometry.\\
	Among the many ways of starting I will give one.  It is a good
	Exercise to ask students to try other approaches.\\
	I will choose one line $m$ as the ideal line and a conic, given
	by 5 points (again an other set of 5 elements can be chosen) as
	the defining circle.  Mid points of the side of a triangle
	can be obtainred by the construction of the pole of a line
	with respect to a triangle, see \ref{sec-epoletr}.  Using the mid-points
	of the sides we can derive the barycenter.\\
	For perpandicularity, we choose one of the point $A$ on the conic,
	determine its tangent $t_A$, the parallel tangent $t_B,$ by a
	construction which is the dual of that of finding the second point of
	intersection of a line with the conic, and the point of contact $B$.
	$A \times B$ is a diameter.  Perpendicular directions are obtsined as
	follows. If $I_p$ is an ideal point, we determine the second
	intersection $P$ of $A \times I_p$ with the conic, the perpendicular
	direction is then $(P \times B) \times m.$  We can therefore construct
	the altitudes and therefore the orthocenter.

	\ssec{Fundamental involution, perpendicularity, circles.}

	\sssec{Definition.}\label{sec-dfundinv}
	Starting with an affine geometry associated to $p$, a
	particular involution on the ideal line will be called the 
	{\em fundamental involution}.

	\sssec{Definition.}\label{sec-disoptl}
	If the fundamental involution is hyperbolic, its fixed points are 
	called {\em isotropic points}, the other points on the ideal line
	will be called {\em ideal points} or {\em directions}.  (In a 
	hyperbolic involutive
	geometry, the isotropic points are no more called ideal points).
	The lines through the isotropic points, distinct from the ideal line,
	are called {\em isotropic lines}.\\
	The lines which are neither ideal or isotropic are called {\em ordinary}
	{\em lines}, the points which are not ideal or isotropic points are 
	called {\em ordinary points}, ordinary lines or points will 
	abbreviated from now on by lines or points.  On an ordinary line, 
	there are $p$ ordinary points and one ideal point.

	\sssec{Definition.}\label{sec-dperp}
	Corresponding pairs of points in the fundamental involution are called 
	{\em perpendicular ideal points} or {\em perpendicular directions}.
	2 lines whose ideal points are perpendicular directions are called
	{\em perpendicular lines}.\\
	Some obvious results follow from these definitions and from those of
	the corresponding affine geometry.  For instance:

	\sssec{Theorem.}
	{\em All the lines perpendicular to a given line are parallel.}

	\sssec{Definition.}\label{sec-dcircul}
	If the involution defined by a conic on the ideal line
	is the same as the fundamental involution, the corresponding conic
	is called a {\em circle} and the corresponding polarity is called a
	{\em circularity}.

	\sssec{Theorem. (H. D)}
	{\em In a hyperbolic involutive geometry, a necessary and
	sufficient condition for a conic to be a circle is that it passes
	through the isotropic points.}

	\sssec{Definition.}\label{sec-dcenterco}
	The {\em center of a conic} is the pole of the ideal line in the
	corresponding polarity.  (See II.2.3.0).

	\sssec{Theorem. (H. D)}
	{\em In a hyperbolic geometry, the center of a circle is
	the intersection of its isotropic tangents.}

	\sssec{Definition.}\label{sec-ddiam}
	A {\em diameter of a conic} is a line passing through its
	center  (See II.2.3.1).

	\sssec{Definition.}\label{sec-dmediat}
	A {\em mediatrix of 2 points} $A$ and $B$ on a line $l$, which is not
	an isotropic line, is the line perpendicular to $l$ through the 
	mid-point of $A$ $B$.  (See II.6.2.6)

	\sssec{Example.}
	In the examples of involutive and Euclidean geometry,  I will
	make one of 2 choices for the ideal line and for the defining circle.
	\begin{enumerate}
	\setcounter{enumi}{-1}
	\item In the first choice,\\
		0.$[1,1,1]$ is the ideal line, as in affine geometry.\\
		1.$X0^2$ + $X1^2$ + k $X2^2$ = 0, $k\neq -\frac{1}{2}$, is the 
		defining circle,\\
	\hth	$-(1+2k)$ N $p$ for the elliptic case,
		$-(1+2k)$ R $p$ for the hyperbolic case.\\
		2. Let $\delta^2 := -1-2k$.\\
	\hth If $k\neq -1$, the isotropic points $(1,y,-1-y)$ correspond to the 
		roots $y_1$ and $y_2$ of\\
		3. $(1+k)y^2 + 2k y + (1+k) = 0.$\\
	or with\\
		4. $k' = \frac{2k}{1+k)},$\\
	to the roots of\\
		5. $y^2 + k' y + 1 = 0$.\\
	Therefore\\
	\hth$	y_1 = \frac{-k + \delta }{1+k}$ and 
		$y_2 = \frac{-k - \delta }{1+k}.$\\
	\hth If $k = -1,$ the isotropic points are (0,1,-1) and (1,0,-1).\\
	\hth The polar of $(X0,X1,X2)$ is $[X0,X1,kX2]$,\\
	\hth The direction perpendicular to $(X0,X1,-X0-X1)$ is\\
	\hth 	$(kX0 + (1+k)X1, -(1+k)x0 - kX1, X0 - X1)$.\\
	\hth If $k = -\frac{1}{2}$, the conic 0.1. is tangent to the ideal line.
	\item In the second choice,\\
		0.$[0,0,1]$ is the ideal line, as in Euclidean geometry.\\
		1.$k X0^2 + X1^2 = X2^{2}, k\neq 0$, the defining circle,\\
	\hth$		-k$ N $p$ for the elliptic case,\\
	\hth$		-k$ R $p$ for the hyperbolic case.\\
	\hth$	\delta ^2 := -k.$\\
	\hth If $k$ = 0, the conic 1.1. is tangent to the ideal line.\\
	\hth The isotropic points are $(1, \delta , 0)$ and $(1,-\delta ,0)$.\\
	\hth The polar of $(X0,X1,X2)$ is $[kX0,X1,-X2]$.\\
	\hth The direction perpendicular to $(X0,X1,0)$ is $(X1,-kX0,0)$.
	\end{enumerate}

	\ssec{Altitudes and orthocenter.}

	\sssec{Definition.}\label{sec-daltid}
	In a triangle $\{A_i\}$, the {\em altitude} $\ov{m}a_i$
	from $A_i$ is the line through $A_i$ which is perpendicular to the
	opposite side $a_i := A_{i+1} \times A_{i-1}$ (C0.1,N0.3).

	\sssec{Theorem.}
	{\em The altitudes $\ov{m}a_i$ of a triangle are concurrent at a 
	point $\ov{M}$.} (D0.12)

	\sssec{Definition.}\label{sec-dorthoc}
	The point $\ov{M}$ is called the {\em orthocenter} of the 
	triangle.(N0.2)

	\sssec{Theorem.}\label{sec-tright}
	{\em The necessary and sufficient condition for a triangle to be a
	right triangle at $A_i$ is that its orthocenter $\ov{M}$ 
	coincides with $A_i$.}

	\sssec{Theorem.}\label{sec-tisosc}
	{\em The necessary and sufficient condition for a 
	triangle to be an isosceles triangle is that the orthocenter be on 
	the altitude from $A_i$ and distinct from the center of mass.}

	\sssec{Theorem.}\label{sec-tequil}
	{\em The necessary and sufficient condition for a triangle to be an
	equilateral triangle is that the orthocenter and the barycenter
	coincide.}


	\ssec{The geometry of the triangle, I.}

	\setcounter{subsubsection}{-1}
	\sssec{ Introduction.}
	We are now ready to give a large number of results of
	finite involutive geometry associated to a scalene triangle whose
	vertices  $A_0,$ $A_1,$ $A_2$ and whose sides $a_0,$ $a_1,$ $a_2$ are
	ordinary.\\
	Theorem II.6.2.7. determines a point $M,$ the center of mass, at the
	intersection of the {\em medians}$A_0 M_0,$ $A_1 M_1,$ $A_2 M_2,$ the
	points $M_i$ being the mid-points of pairs of vertices.\\
	Theorem 3.1. determines a point $\ov{M},$ the orthocenter, at the
	intersection
	of the {\em altitudes }$A_0 \ov{M}_0,$ $A_1 \ov{M}_1,$
	$A_2 \ov{M}_2,$ the points
	$\ov{M}_i$ being the feet of the altitudes.\\
	In a scalene triangle, $M$ and $\ov{M}$ are distinct, are
	distinct from the
	vertices and are not collinear with any of the vertices.  A large number
	of results can therefore be obtained as direct consequences of
	rephrasing the results of Theorem 3.6. and 4.0.\\
	Similar results can be obtain for right triangles, for isosceles
	triangle and for equilateral triangles.  These will be left as
	exercises.\\
	These results were in fact the starting point of our study of finite
	Euclidean geometry, as explained in section \ldots.\\
	All references will be to Theorems 3.6. and 4.0. unless explicitely
	indicated.

	\sssec{Definition.}
	The {\em ideal points} $MA_i$ {\em of a triangle} are the ideal points
	on its sides.\\
	The {\em orthic points} $\ov{M}A_i$ {\em of a triangle} are the
	points on the corresponding sides $a_i$ of the triangle and
	$\ov{m}m_i$ of the orthic triangle  (D0.13, N0.6). See Fig. 1.

	\sssec{Theorem.}
	{\em The orthic points $\ov{M}A_i$ are on the orthic line
	$\ov{m}$} (D0.14*).

	\sssec{Definition.}\label{sec-dcomplem}
	The triangle $M_i$ is called the {\em complementary triangle.}
	Its sides are denoted $mm_i.$\\
	The triangle $\ov{M}_i$ is called the {\em orthic triangle}.
	Its sides are denoted $\ov{m}m_i$  (D0.18, N0.5).

	\sssec{Definition.}\label{sec-dorthic}
	The {\em orthic line} $\ov{m}$ {\em of a triangle} is the polar
	of its orthocenter with respect to the triangle  (N0.8).\\
	Its direction $EUL$ is called the {\em orthic direction} (N1.1).

	\sssec{Definition.}\label{sec-dleuler}
	The line $eul := M \times \ov{M}$ is called the {\em line of
	Euler} (D1.0, N1.0).

	\sssec{Theorem.}
	{\em The mid-points $M_i$ at the intersection of the medians $ma_i$
	with the sides $a_i$ and the feet $\ov{M}_i$ of the altitudes
	$\ov{m}a_i$ are on a circle $\gamma$} (D1.20, C1.4). See Fig. 2.

	\sssec{Definition.}
	The circle $\gamma$ is called the {\em circle of Brianchon-Poncelet}
	(N1.11).

	\sssec{Theorem.}
	{\em If $Maa_i$ ($\ov{M}aa_i$) is the intersection of the median
	$ma_{i+1}$ ($ma_{i-1}$) with the altitude $\ov{m}a_{i-1}$
	($\ov{m}a_{i+1})$, then the lines $mMa_i$ joining $Maa_i$ and
	$\ov{M}aa_i$ have a point $K$ in common} (D1.2, D1.3, D1.4*).
	See Fig. 3.

	\sssec{Definition.}
	The point $K$ is called the {\em point of Lemoine}  (N1.2).

	\sssec{Definition.}\label{sec-dcircumci}
	The {\em circumcircle $\theta$ of a triangle} $\{A_0, A_1, A_2\}$ is
	the circle passing through the vertices of the triangle  (D1.19, H1.1,
	N1.10).

	\sssec{Theorem.}
	{\em The line $ta_i$ through the vertex $A_i$ parallel to the side
	$\ov{m}_i$ of the orthic triangle is the tangent at $A_i$ to the
	circumcircle} (D1.7, D1.19*). See Fig. 4.

	\sssec{Definition.}\label{sec-dtrtang}
	The triangle with sides $ta_i$ is called the {\em tangential
	triangle. } Its vertices are denoted by $T_i$  (D1.8, N1.5).

	\sssec{Definition.}\label{sec-dtrmixed}
	The {\em mixed triangles} are the triangles with respective sides\\
	\hth$c_i := M_{i+1} \times \ov{M}_{i-1}$ and
	$\ov{c}_i := \ov{M}_{i+1} \times M_{i-1}$ (D1.13, N1.6).\\
	The {\em mixed feet} are the points $CC_i,$ $(\ov{C}C_i)$ on the
	side of the given triangle and the corresponding side $c_i,$
	$(\ov{c}_i)$ of the mixed triangle (D1.14, N1.7). Dee Fig. 5.

	\sssec{Theorem.}
	{\em The mixed feet $CC_i,$ $(\ov{C}C_i)$ of a mixed triangle are
	collinear on the line $p$ $(\ov{p})$} (D1.15*).

	\sssec{Definition.}\label{sec-dmixedl}
	The line $p$ and $\ov{p}$ are called the {\em mixed lines of a
	triangle} (N1.8).

	\sssec{Theorem.}
	{\em The mixed lines $p$ and $\ov{p}$ of a triangle meet at the
	point $PP$ which is on the line of Euler} (D1.16, C1.0).

	\sssec{Definition.}\label{sec-dmixedc}
	$PP$ is called the {\em mixed center of the triangle}  (N1.9).

	\sssec{Definition.}\label{sec-dpmedo}
	The intersection $\ov{I}Ma_i$ of a median with the orthic line
	is called a {\em medorthic point}  (D0.15, N0.9).

	\sssec{Definition.}\label{sec-dpeul}
	The intersection of the lines $\ov{m}e_i$ $(me_i)$ joining the
	medorthic points $\ov{I}Ma_{i+1}$
	$(\ov{I}Ma_{i-1})$ to the foot $\ov{M}_{i-1}$
	$(\ov{M}_{i+1})$ are called the {\em points of Euler} $EE_i$
	({\em Eulerian points} $\ov{E}_i$) (D5.0, D5.1). Fig. 6.

	\sssec{Theorem.}
	The {\em points of Euler} $EE_i$ are the mid-points of the
	segment joining the orthocenter $\ov{M}$ to the vertex $A_i$
	(D5.1, C5.3).

	\sssec{Theorem.}
	{\em The points of Euler $EE_i$ are on the circle of Brianchon-Poncelet}
	(C5.5).\\
	{\em The Eulerian point $\ov{E}E_i$ is on the median $ma_i$ as
	well as on the circle of Brianchon-Poncelet} (C5.0, C5.5).

	\sssec{Theorem.}
	{\em The lines $em_i$ joining the mid-points $M_i$ to the Eulerian
	points $EE_i$ are concurrent at a point $EE.$\\
	The lines $\ov{e}m_i$ joining the feet $\ov{M}_i$ to
	the Eulerian points $\ov{E}E_i$ are concurrent at a point
	$\ov{E}E$} (D5.2, D5.3*).\\
	{\em $EE$ is on the line of Euler and is the center of the circle of
	Brianchon-Poncelet.\\
	$\ov{E}E$ is on the line of Euler and is
	the pole of the orthic line with respect to the circle of
	Brianchon-Poncelet}  (C5.1, C5.4, N5.1).

	\sssec{Definition.}\label{sec-dmedtrix}
	The {\em mediatrix} $mf_i$ is the line through through the
	mid-point $M_i$ perpendicular to the corresponding side $a_i$
	(D6.0, N6.0).

	\sssec{Theorem.}
	{\em The vertex $T_i$ of the tangential triangle is on the mediatrix
	$mf_i$}  (D6.0, C6.8, N6.0).\\
	{\em The mediatrices $mf_i$ are concurrent at a point $O.$\\
	The diameters $\ov{m}f_i$ of the circle of Brianchon-Poncelet
	which pass through the feet of the altitudes pass through the same
	point $\ov{O}$} (D6.4*, N6.1).

	\sssec{Definition.}\label{sec-dcircumce}
	$O$ is the {\em circumcenter} or center of the circumcircle (N6.1).

	\sssec{Theorem.}
	{\em The circumcenter $O$ is on the line of Euler.\\
	The point $\ov{O}$ is on the line of Euler} (C6.1).

	\sssec{Definition.}
	An {\em equilateral conic} is a conic whose ideal points are harmonic
	conjugates of the isotropic points.\\
	An {\em coequilateral conic} is a conic whose points on the orthic line
	are harmonic conjugates of the coisotropic points.

	We leave as an exercise the pproof of the following Theorem and
	Corollary.

	\sssec{Theorem.}
	{\em If an conic passes through the vertices of the triangle}
	\enumb
	\item {\em it is equilateral if and only if it passes through the
		orthocenter.\\
		Its center is on the circle of Brianchon-Poncelet.}
	\item {\em it is coequilateral if and only if it passes through the
		barycenter.\\
		Its cocenter is on the circle of Brianchon-Poncelet.}
	\enume

	\sssec{Corollary.}
	{\em A conic}\\
	\hth$a_0X_0^2+a_1X_1^2+a_2X_2^2+b_0X_1X_2+b_1X_2X_0+b_2X_0X_1 = 0$,\\
	\enumb
	\item {\em is equilateral if and only if}\\
	\hth$m_1m_2(a_1+a_2-b_0)+m_2m_0(a_2+a_0-b_1)+m_0m_1(a_0+a_1-b_2)= 0$.
	\item {\em it is coequilateral if and only if}\\
	\hth$m_0(m_1^2a_1+m_2^2a_2-m_1m_2b_0)+m_1(m_2^2a_2+m_0^2a_0-m_2m_0b_1)
		+m_2(m_0^2a_0+m_1^2a_1-m_0m_1b_2)=0$.
	\enume

	\sssec{Definition.}
	The {\em conic of Kiepert} is the conic circumscribed to the triangle
	passing through the barycenter and the orthocenter. (D3.8.)\\
	The {\em conic of Jerabek is the conic circumscribed to the triangle
	passing through the orthocenter and the point of Lemoine. (D36.16.)\\
	These are therefore equilateral. (C3.3 and C36.).
	The center of one conic is the cocenter of the other and these
	are on the circle of Brianchon-Poncelet (C8.9 and C36.18.)

	\sssec{Definition.}\label{sec-dcneff}
	The circle through the vertices $T_i$ of the tangential triangle is
	the {\em circle of Neff}.

	\sssec{Theorem.}
	\enumb
	\item{\em The circle of Neff is a cocircle.}
	\item{\em The ortic line is the radical axis of the circle of Neff and 
		both the circumcircle and the circle of Brianchon-Poncelet.}
	\enume

	\sssec{Definition.}\label{sec-dtrneff}
	A {\em triangle of Neff} is a triangle whose orthocenter is on the
	conic.
	
	\sssec{Exercise.}\label{sec-etrneff}
	Prove that in a triangle of Neff, one of the sides of the tangential
	triangle is a diameter of the circle of Neff.
	Determine other conditions for this to happen.
xxx
	\sssec{Definition.}\label{sec-dpcomplo}
	The points $EUL_i$ at the intersection of the corresponding
	sides of the complementary triangle and of the orthic triangle are
	called the {\em complorthic points}  (D8.0, N8.0).\\
	The lines $aeUL_i$ joining the complorthic points are called
	{\em complorthic lines}  (D8.3, N8.1).\\
	The triangle whose vertices are the complorthic points is called the
	{\em complorthic triangle} (N8.2).

	\sssec{Definition.}\label{sec-dpmixed}
	The intersections of corresponding sides of the mixed triangles
	are the {\em mixed points} $D_i$ (D8.4, N8.4).

	\sssec{Theorem.}
	\enumb
	\item{\em The mixed points $D_i$ and $\ov{D_i}$ are on the line of
	Euler} (C8.2).
	\item{\em The vertex $A_i$ and the mixed point $D_i$ are on the
	complorthic line $aeUL_i.$} (C8.1, C8.3).
	\item{\em The lines $nm_i$ joining the mid-points of the sides to
	corresponding complorthic points $EUL_i$ are concurrent in a point $S.$}
	\item{\em The lines $\ov{n}m_i$ joining the feet of the altitudes to
	the corresponding complorthic points $EUL_i$ are concurrent in a point
	$\ov{S}$} (D8.1, D8.2*).
	\enume

	\sssec{Definition.}\label{sec-dsch}
	The points $S$ and $\ov{S}$ are the {\em point} and {\em copoint
	of Schr\"{o}ter} (N8.3).

	\sssec{Theorem.}
	{\em $S$ is the point of Miquel of the quadrangle $a_{i+1},$ $a_{i-1},$
	$\ov{m}a_{i+1},$ $\ov{m}a_{i-1}.$  $S$ is therefore also on
	the circles through $A_i,$ $\ov{M},$ $\ov{M}_{i+1},$
	$\ov{M}_{i-1}$ of center $E_i$ and on the circles through
	$A_{i+1},$ $A_{i-1},$ $\ov{M}_i$ of center $M_i.$} (See .)

	\sssec{Theorem.}
	\enumb
	\item{\em The points of Schr\"{o}ter are on the circle of
	Brianchon-Poncelet} (C8.8).
	\item{\em The first point of Schr\"{o}ter $S,$ the Eulerian point
	$\ov{E}E_i$ and the mixed point $D_i$ are on the same line
	$s_i$}
	\item{\em The second point of Schr\"{o}ter $\ov{S},$ the point of
	Euler point $EE_i$ and the mixed point $D_i$ are on the same line
	$\ov{s}_i.$}  (D8.5, C8.4).
	\enume

	\sssec{Theorem.}
	{\em The conic through the barycenter $M$, the orthocenter
	$\ov{M}$
	and the feet $Gm_i$ of the perpendicular $iMA_i$ from $M$ to the
	corresponding altitude $\ov{m}a_i$ are on a circle $`omicron$}
	(D10.3, D10.4, D10.7, C10.7).\\
	{\em This circle passes also through the perpendiculars
	$\ov{G}m_i$
	which are the feet of the perpendiculars $\ov{g}m_i$ from
	$\ov{M}$ to the corresponding median $ma_i$.
	$\{M,\ov{M}\}$ is a diameter whose mid-point is $G$}
	(C10.1, C10.8). See Fig. 9.

	\sssec{Definition.}\label{sec-dorthoce}
	$`omicron$ is called the {\em orthocentroidal circle} (N10.2).

	\sssec{Theorem.}
	{\em If we join \ldots}.\\ 
	(D6.1, D10.1, D10.2, D10.3*).

	\sssec{Definition.}
	The $G$ is the center of the orthocentroidal circle $`omicron$ (N10.13).

	\sssec{Theorem.}
	{\em The line $be_i$ is parallel to the median $ma_i$}
	(D10.5, N10.0, C10.2).

	\sssec{Theorem.}
	{\em The 3 circles, the circumcircle $\theta$, the circle $\gamma$ of
	Brianchon-Poncelet and the orthocentroidal circle $`omicron$ have the
	same radical axis $\ov{m}.$} (C1.5, C10.9)

	\sssec{Definition.}
	An {\em orthocentric quadrangle} is \ldots.

	An example is provided by the circumcentral orthocentric quadrangle
	(N10.1).

	\ssec{The geometry of the triangle. II.}
	
	\sssec{Theorem.}
	{\em The line $tm_i$ through the mid-point $M_i$ parallel to the side
	$\ov{m}_i$ of the orthic triangle is tangent at $M_i$ to the
	circle of Brianchon-Poncelet.} (D12.0, C12.11)

	\sssec{Definition.}\label{sec-dpsymm}
	The line $at_i$ joining the vertex $A_i$ of the triangle to
	the vertex $T_i$ of the tangential triangle are called the
	{\em symmedians} (D12.1, N12.0).

	\sssec{Theorem.}
	{\em The symmedians $at_i$ are concurrent at a point $K$}  (C12.6).

	\sssec{Theorem.}
	{\em The point $K$ of Lemoine, the first point $S$ of Schr\"{o}ter and
	the point $G$ are collinear on the line $gk.$\\
	The point $K$ of Lemoine, the second point $\ov{S}$ of
	Schr\"{o}ter and the point $\ov{G}$ are collinear on the line
	$\ov{g}k$}  (D12.2, C12.7).

	\sssec{Definition.}\label{sec-dpotan}
	Te {\em tangential point} $AMa_i$ ($AM\ov{a}_i$) is the
	intersection od the parallel $am_{i+1}$ ($a\ov{m}_{i-1}$) through
	$A_{i-1}$ ($A_{i+1}$) to the altitude $\ov{m}a_i$ and the parallel
	$am_{i-1}$ ($a\ov{m}_{i+1}$) through $A_i$ to
	$\ov{m}a_{i+1}$ ($\ov{m}a_{i-1}$) (D6.8, D14.4, N14.0).

	\sssec{Definition.}\label{sec-dcitan}
	The {\em tangential circle} $\chi a_i$ ($\chi\ov{a}_i$) is the
	circle though the vertices $A_{i+1}$ and $A_{i-1}$ tangent at $A_{i-1}$
	($A_{i+1}$) to the side $a_{i+1}$ ($a_{i-1}$) (D14.13, C14.8, C14.5).

	\sssec{Theorem.}
	{\em The tangential circle $\chi a_i$ ($\chi\ov{a}_i$) passes
	through the tangential point $AM\ov{a}_{i+1}$ ($AMa_{i-1}$)}
	(D14.13).

	\sssec{Definition.}\label{sec-dparlem}
	The {\em parallels of Lemoine} $kk_i$ are the lines through the
	point $K$ of Lemoine parallel to the sides of the triangle
	(D15.0, N15.0). See Fig. 13.

	\sssec{Definition.}\label{sec-d1trbr}
	The vertices $Br1_i$ of the {\em first triangle of Brocard} are
	the intersections of the mediatrices $mf_i$ with the parallels of
	Lemoine $kk_i$ (D15.1, N15.1).

	\sssec{Theorem.}
	{\em The lines $br0_i$ joining the vertices of a triangle $A_i$ to
	the corresponding vertex $Br1_i$ of the first triangle of Brocard are
	concurrent at a point $BR0$}  (D15.2, D15.3*).

	\sssec{Theorem.}
	{\em The lines $br_i$ ($b\ov{r}_i)$ joining the vertices $A_{i-1}$
	($A_{i+1}$) to the vertices $Br1_{i+1}$ ($Br1_{i-1}$) of the first
	triangle of Brocard are concurrent at a point $Br$
	($B\ov{r}$} (D15.4, D15.5*).

	\sssec{Definition.}\label{sec-dpobro}
	The point $Br$ ($B\ov{r}$) is called the {\em first (second)
	point of Brocard} (N15.4).

	\sssec{Definition.}d2trbr
	The points $Br2_i$ at the intersection of the parallel $ok1_i$ to the
	side $a_i$ through the center $O$ of the circumcircle and the
	perpendicular $km_i$ to $a_i$ through the point $K$ of Lemoine are
	the vertices of the {\em second triangle of Brocard}
	(D13.4, D13.3, D15.6, N15.2).

	\sssec{Theorem.}
	{\em The lines $br3_i$ joining corresponding vertices $Br1_i$ and
	$Br2_i$ of the first and second triangle of Brocard are concurrent at a
	point $Bro$} (D15.9, D15.10*).

	\sssec{Definition.}\label{sec-dcrossta}
	The {\em cross tangential line} $mff_i$ is the line through the 
	tangential points $AMa_{i+1}$ and $AM\ov{a}_{i-1}$ (D15.7, N15.6).

	\sssec{Definition.}\label{sec-d3trbr}
	The {\em vertices of the third triangle of Brocard} $Br3_i$
	are the intersections of the cross tangential line $mff_i$
	and the corresponding symmedian $at_i$ (D15.8, N15.3).

	\sssec{Theorem.}
	{\em The vertices $Br1_i$, $Br2_i$ and $Br3_i$ of the first, second and 
	third triangle of Brocard, the first and second point of Brocard
	$Br1$ and $Br2,$ the center $O$ of the circumcircle and the point $K$
	of Lemoine are on a circle $\beta$ with center $Bro$, the mid-point of
	$\{K,O\}$} (D15.18, C15.17, C15.18, C15.12, C15.13, C15.7).

	\sssec{Definition.}\label{sec-dcibro}
	The circle $\beta$ is called the {\em circle of Brocard} (N15.6).

	\sssec{Definition.}
	The {\em conics of Tarry} $`Tarry[i]$ are the conics through the
	barycenter $M$ 	and through 2 vertices, tangent there to the side
	through the third vertex, $A_i$. (N19.0.)

	\sssec{Theorem.}
	{\em Let $Apt_0$, $(Ap\ov{t}_0)$  be the intersection of the line
	through $A_0$ parallel to the median $ma_1$ $(ma_2)$ with the line
	through $A_2$ $(A_1)$ parallel to the median $ma_0$ and circularly for
	$Apt_1$, $(Ap\ov{t}_1)$, $Apt_2$, $(Ap\ov{t}_2)$, then the line
	$Apt_0\times (Ap\ov{t}_0)$ is the tangent common to $`Tarry_1$ and
	$`Tarry_2$ with $Apt_0$ and $Ap\ov{t}_0$ as point of contact.}
	(D19.7, C19.0, D33.7, C33.5, C33.6.)

	From the coordinates associated with the symmetric Theorem using
	$\ov{M}$ instead of $M$, it is easy to solve the problem of
	C. Bindschelder, El. Math. 1990, p. 56.


	\ssec{Geometry of the triangle. III.}
	\sssec{Definition.}
	The {\em line of Schr\"{o}ter}, $pap$ is \ldots. (N4.1)
	It is tangent to the conics of Steiner, Lemoine and Simmons,
	(P. de Lepiney, Math. 1922-133)(C36.7)
	!dont have def. of Lemoine and Simmons, these are of the form
	$b_0 x_1x_2 + b_1 x_2x_0 + b_2 x_0x_1 = 0,$ with
	$b_0 m_0(m_1-m_2)+ b_1 m_1(m_2-m_0)+ b_2 m_2(m_0-m_1) = 0.$
	!$MK . Center(`Lemoine) = 0,$36.15 no???
	$MK . Center(`Simmons) = 0,$??\\

	\ssec{Geometry of the triangle. IV.}
	\sssec{Theorem. [Kimberling]}\label{sec-tkimberling1}
	\enumb
	\item {\em The lines joining the vertices $T_i$ of the tangential
	triangle to the second intersection $B_i$ of the medians $ma_i$ with the
	circumcircle $\theta$ are concurent at a point $CK$.} (C47.0.)
	\item {\em The lines joining the vertices $T_i$ of the tangential
	triangle to the second intersection $\ov{B}_i$ of the altitudes
	$\ov{m}a_i$ with the circumcircle are concurent at a point $C~K$.}
	(C47.0.)
	\enume

	\sssec{Definition.}
	The points $CK$ and $\ov{C}K$ just defined is called respectively the
	{\em point and copoint of Kimberling}. (N47.0.)

	\sssec{Theorem. [Kimberling]}\label{sec-tkimberling2}
	{\em The point and copoint of Kimberling are on the line of Euler.}
	(C47.4) See \ref{sec-tgkimh}

	\sssec{Theorem.}\label{sec-tdevkimb}
	\enumb
	\item Desargues($M,\{A_i\},\{B_i\},ee$). (D47.21)
	\item Desargues($\ov{M},\{A_i\},\{\ov{B}_i\},\ov{e}e$). (D47.21)
	\item Desargues($CK,\{T_i\},\{B_i\},\ov{m}$). (C47.6)
	\item Desargues($\ov{M},\{A_i\},\{\ov{B}_i\},m$). (C47.6)
	\enume

	\sssec{Theorem. [Sekigichi]}\label{sec-tsek}
	{\em The set of points on a triangle at which the sum of the distances
	to the sides is equal to the arithmetic mean of the lengths of the
	altitudes is a segment of a line through the barycenter.}
	(Amer. Math. Monthly, 1981, 349 and 1984, 257.)

	\sssec{Definition.}
	The line defined in the preceding Theorem is called the {\em line of
	Sekiguchi}.

	\sssec{Theorem.}\label{sec-tsek1}
	{\em The line of Sekiguchi is perpendicular to the line $ok$
	joining the center $O$ of the circumcircle to the point $K$ of Lemoine.}
	
	The segment $[A_0,Sek_0]$ is equal to the segment $\ov{M}_1,A_1$,
	(D18.27), the segment $[\ov{M}_0,Set_1]$ is equal to the segment
	$A_2,\ov{M}_2$, (D18.28), the line $sek_2$ joining $Sek_0$ and $Set_1$
	has the direction of $O\times K$ (C18.23).

	\ssec{Geometry of the triangle. V.}

	\sssec{Definition.}
	The {\em triangle of Nagel}, $\{Na_{i}\}$ has as its vertices
	the point of contact of the $i$-th exscribed circles with the 
	$i$-th side}.(N21.0.)

	\sssec{Definition.}
	The {\em conic of Feuerbach} is the conic through the vertices
	of the triangle, the point of Gergonne $J$ and the incenter $I$.
	(N20.6.)

	\sssec{Theorem. [Feuerbach]}
	{\em The conic of Feuerbach is an equilateral hyperbola,
	it passes through the orthocenter and the point of Nagel,
	it is tangent at $I$ to the line through$I$ and $O$ (Th\'{e}bault), it
	has the point of Feuerbach as its center.} (D20.23., C20.14, C20.15,
	C20.17,C23.8.) See also Neuberg, Math. 1922-51-90.

	\sssec{Theorem. [Kimberling]}\label{sec-tekim}
	{\em If $Kim_0$ is the intersection of the lines from the center $I_1$ 
	and $I_2$ of the excribed circles on the exterior bissectrix through
	$A_0$ perpendicular respectively to the sides $a_2$ and $a_1$, then
	the line $kimc_0$ joining $Kim_0$ to $A_0$ and the similarly obtained
	lines $kimc_1$ and $kimc_2$ have a point $Kim$ in common.} (D21.30.)
	See \ref{sec-tgkim}

	\sssec{Definition.}
	The point $Kim$ is called the {\em excribed point of Kimberling.}
	(N21.5.)

	\sssec{Theorem.}
	{\em The point of Kimberling is on the conic of Feuerbach.}
	(C21.11.)

	\sssec{Theorem.}
	{\em If $Kid_0$ is the intersection of the lines from the center $I_1$ 
	and $I_2$ of the excribed circles on the exterior bissectrix through
	$A_0$ and respectively the points $M\ov{}A_2$ and $\ov{M}_1$ on the
	orthic line $m~$ and the sides $a_2$ and $a_1$, then the line $kidc_0$
	joining $Kid_0$ to $A_0$ and the similarly obtained lines $kidc_1$ and
	$kidc_2$ have a point $Kid$ in common.} (D21.26.)

	\sssec{Definition.}
	The point $Kid$ is called the {\em excribed orthic point.} (N21.4.)

	\sssec{Theorem.}
	{\em The barycentric coordinates of the incenter $I$ are proportional
	to the lengths of the sides of the triangle.}

	\sssec{Exercise.}
	Prove that the point $En$ is the centroid of a wire of uniform
	density forming the sides of the triangle $A_i$. See C. J. Bradley,
	Math. Gazette, 1989, p. 44. for the latter.

	\sssec{Definition. [Mandart]}
	The {\em conic of Nagel} is the conic tangent at the vertices of the
	triangle of Nagel to the sides of the triangle.(N27.0)

	\sssec{Theorem. [Mandart and Neuberg]}
	{\em The center of the conic of Nagel is on the conic of Feuerbach.}
	C27.1. (Math. 1922-125)

	\sssec{Definition. [Mandart]}
	The {\em cercle of Nagel} is the circle circumscribed to the
	triangle of Nagel.(Math. 1922-125)

	\sssec{Theorem.}
	{\em The complimentary point $En$ of the incenter $I$ is the
	center of gravity of the perimeter of the triangle.}
	(See Math. 1889, Suppl. p. 8, {\bf 26})

	\ssec{Sympathic projectivities.}

	\sssec{Introduction.}
	This section discusses is some detail the notion of equality of angles
	in involutive geometry.

	\sssec{Definition.}
	A {\em sympathic projectivity} is one which is 
	amicable with the fundamental involution.  (II, 1.5.10)

	\sssec{Theorem. (H. D.)}
	{\em If the involutive geometry is hyperbolic, a sympathic
	projectivity has 2 fixed points, the isotropic points.}

	\sssec{Theorem.}
	{\em The sympathic projectivities form an Abelian group under 
	composition.\\
	Moreover, using 1.0.10.0.4., if\\
	\hth$	f_b(y) = \frac{-1+by}{b+k'+y},$\\
	then\\
	\hth$	f_{b1} \circ f_{b2} = f_{b3}$, with 
		$b3 = \frac{-1+b1b2}{k'+b1b2}.$}

	Proof.\\
	$f_{b1}(f_{b2}(y))
	= (-(b1+b2+k')+\frac{(-1+b1b2)y}{(b1+k'}(b2+k')-1+(b1+b2+k')y),$\\
	dividing numerator and denominator by $b1+b2+k'$ gives the conclusion of
	the Theorem.\\
	See also \ldots, Section 7.

	\sssec{Example.}
	The method of obtaining sympathic 
	projectivities will be studied in section $\ldots$.  It will be seen 
	that all are powers of a
	sympathic projectivity S which is of order p-1 in the hyperbolic case
	and of power p+1 in the elliptic case.  This generating projectivity
	is not unique, choosing one of these as fundamental sympathic
	projectivity will constitute the next step towards Euclidean geometry,
	the sympathic geometry.  The fundamental involution is 
	$S^{(\frac{p-1}{2})}$ or $S^{(\frac{p+1}{2})}.\\$	

	With $p$ = 7, (elliptic), we will choose $k$ = 0, $\delta^2$ = 6,\\
	\hth The sympathic projectivities are 
		$S^{i},$ $i$ = 0 to 7, with\\
	\hth$	S(1,j,-1-j) = (2-3j,3+2j,-5+j),$\\
	\hth$	S(0,1,6) = (1,4,2), $ or\\
	$\begin{array}{rl}
	S,=,&(\:7,14,20,26,32,38,44,50)\\
		&(38,44,26,14,,7,20,50,32) = S^7,\\
	S^2 = &(\:7,14,20,26,32,38,44,50)\\
		&(20,50,14,44,38,26,32,,7) = S^6,\\
	S^3 = &(\:7,14,20,26,32,38,44,50)\\
		&(26,32,44,50,20,14,,7,38) = S^5.\\
	\end{array}$\\
	    The fundamental involution is\\
	$\begin{array}{rl}
	S^4,=,&(\:7,14,20,26,32,38,44,50)\\
		&(14,\:7,50,32,26,44,38,20).\\
	\end{array}$\\
	    The isotropic points are 
		$(1,\delta,-1-\delta),$ $(1,-\delta,-1+\delta)$\\
	With $p = 7,$ (hyperbolic), we will choose $k = 1,$ $\delta = 2,$\\
	    The sympathic projectivities are $S^{i},$ $i$ = 0 to 5, with\\
	$\begin{array}{rl}
	S = &(26,38,\:7,14,20,32,44,50)\\
		&(26,38,44,20,\:7,14,50,32) = S^5,\\
	S^2 = &(26,38,\:7,14,20,32,44,50)\\
		&(26,38,50,\:7,44,20,32,14) = S^4,\\
	\end{array}$\\
	    The fundamental involution is\\
	$\begin{array}{rl}
	S^3 = &(26,38,\:7,14,20,32,44,50)\\
		&(26,38,32,44,50,\:7,14,20).
	\end{array}$\\
	    The isotropic points are (26) = (1,2,4) and (38) = (1,4,2).\\
	anti\ldots.

	\ssec{Equiangularity.}

	\sssec{Definition.}
	An {\em angle} is an ordered pair $\{a,b\}$ of ordinary lines $a$ and
	$b$.

	\sssec{Definition.}
	{\em Two angles} $\{a,b\}$ and $\{a_1,b_1\}$ {\em are equal} and we 
	write\\
	\hth$	\{a,b\} = \{a_1,b_1\},$\\
	if the ideal points on these lines, $A,$ $B$, $A_1$, $B_1$ are such 
	that there
	exists a sympathic projectivity which associates $A$ to $B$ and $A_1$ 
	to $B_1$.
	Compare with Coxeter, p. 9 and p.125).

	\sssec{Notation.}
	In view of 2.3., we will also use $\{A,B\} = \{A_1,B_1\}$ instead of
	$\{a,b\} = \{a_1,b_1\},$ where $A$, $B$, $A_1$, $B_1$ are the ideal 
	points on $a$, $b$, $a_1$, $b_1$.

	\sssec{Example.}
	For $p$ = 5, starting with Example II.1.5.12. if $\phi'$ is used to
	to define the fundamental involution, then $\phi$  is a sympathic
	projectivity.  We have the equality of angles \{(10),(5)\} 
	= \{(5),(26)\} = \{(26),(14)\} = \{(14),(18)\} = \{(18),(22)\} 
	= \{(22),(10)\} and of angles
	\{(10),(14)\} = \{(5),(18)\} = \{(26),(22)\} = \{(14),(10)\}.

	\sssec{Theorem.}
	{\em If $a$ and $b$ are perpendicular, then $\{a,b\} = \{b,a\}.$
	If $c$ and $d$ are also perpendicular, then $\{a,b\} = \{(c,d\}.$}

	\sssec{Definition.}
	If $a$ and $b$ are perpendicular, the angle $\{a,b\}$ is called
	a {\em right angle}.

	\sssec{Definition.}
	If $a$ and $b$ are not parallel and $c$, through $a \times b$, is such
	that $\{a,c\} = \{c,a\},$ $c$ is called a {\em bisectrix} of $\{a,b\}$.
	If a bisectrix exist, we say that the {\em angle} $\{a,b\}$ {\em can 
	be bisected}.

	\sssec{Theorem.}
	{\em If the ideal points on $a$ and $b$ are $(1,a_1,-1-a_1)$ and
	$(1,b_1,-1-b_1),$ then the ideal point $(1,z,-1-z)$ on the bisectrix $c$
	of $\{a,b\}$ satisfies the second degree equation}
	\enumb
	\item$(k'+a_1+b_1)z^2 - 2(a_1b_1-1)z - (a_1+b_1+k'a_1b_1) = 0,$ 
		{\em with} $k' = \frac{2k}{1+k)}.$
	\item {\em The discriminant of 0. is}\\
	\hth$t' = (a_1^2 + k'a_1 +1) (b_1^2 + k'b_1 +1)$
	\item {\em Moreover,\\
	 0. if $t' \neq  0$ is a quadratic residue, the bisectrices are real and
		perpendicular to each other,\\
	 1. if $t'$ is a non residue, there are no real bisectrices,\\
	 2. if $a$ or $b$ is an isotropic line, $t' = 0,$ the bisectrices 
		coincide with the isotropic line,\\
	 3. if both $a$ and $b$ are isotropic, the bisectrices are undefined,\\
	 4. if $a$ and $b$ are parallel, the bisectrices do not exist but the
		directions given by 0. are that of $a$ and of the perpendicular 
		to $a$.}
	\enume

	Proof:  Let the sympathic projectivity which associates to the ideal
	points on $a$ and $c$ the ideal points on $c$ and $b$, have the form\\
	\hth$	f(x) = \frac{a'+b'x}{c'+d'x},$\\
	then\\
	\hth$	z (c'+d'a_1) = a'+b'a_1,$\\
	\hth$	b_1 (c'+d'z) = a'+b'z.$\\
	Because of 0.0.10., $d' = -a' = 1+k,$ $c'-b' = 2k.$\\
	Substituting and multiplying the first equation by $(c_1-b_1),$ the 
	second by $(c_1-a_1)$ and adding, we obtain the equation 0.\\
	If $a_1$ corresponds to an isotropic point, $a_1^2 + k'a_1 + 1 = 0,$
	$t' = 0,$ the roots of 0. are\\
	\hth$\frac{a_1b_1-1}{k'+a_1+b_1} 
		= a_1\frac{a_1b_1-1}{a_1^{2+k'a_1+b_1a_1}} = a_1.$\\
	If $a_1 = b_1,$ 0. can be written $(z - a)((k'+2a_1)z + (2+k'a) = 0.$
	The perpendicularity follows from 1.9.6.

	\sssec{Example.}
	For $p$ = 7, hyperbolic case, let the ideal point on ``a" be (32)
	and on ``b" be (20), $a_1 = 3,$ $b_1 = 1$, $k = k' = 1,$ 0. is 
	$5z^2$ - 4z = 0, with roots 0 and 5 giving the points (14) and (44).
	If $a_1 = b_1 = 0,$ one root is 0, the other is 5.

	\sssec{Definition.}
	The {\em angle} between distinct non isotropic lines {\em is even}
	if and only if the angle can be bisected.

	\sssec{Theorem.}
	{\em Under the hypothesis of Theorem 0.2.7., an angle is even
	if $a_1^2 + k'a_1 +1$ and $b_1^2 + k'b_1 +1$ are both quadratic residues
	or both non residues.}

	The proof follows at once from $\ldots$ 

	\sssec{Theorem.}
	{\em The relation "even" is a equivalence relation.}\\
	Again this follows from \ldots.

	\sssec{Definition.}
	The sum of two angles \ldots

	\ldots  circle, angle at the center, rotation.

	\ssec{Equidistance, congruence.}
	congruence (translation composed with rotation)

	\sssec{Theorem.}
	{\em Any congruence can be written as the composition of
	a translation  rotation and a translation.}

	\sssec{Definition.}
	A {\em segment} $[A,B]$ is an unordered pair of ordinary points $A$
	and $B$.\\
	\ldots not on the same isotropic line?
	
	\sssec{Definition.}
	{\em Two segments are equal}  iff
       
	\sssec{Theorem.}
	{\em If $[A,B] = [B,C]$ and $C$ is on
	$A \times B,$  then either $A = C$ or $B$ is the mid-point of $[A,C].$}

	\ldots equality of segments on parallel line  iff  equal in the affine 
	sense or $AB = CD$ in affine sense or $BA = CD$
       
	\sssec{Theorem.}
	{\em $\{A,B\} = \{C,D\}$ implies $[A,B] = [C,D].$}

	equal. on non parallel segment using translation and circle
	may have to use tangent to circle\\
	def. of congruence.

	\ssec{Special triangles.}

	\sssec{Definition.}
	A {\em right triangle} is a triangle with 2  perpendicular sides.
	If $a_1$ and $a_2$ are perpendicular we say that the triangle is a
	{\em right triangle at} $A_0$.

	\sssec{Theorem.}
	{\em A necessary and sufficient condition for a triangle to be
	a right triangle at $A_0$ is that $m_1 = m_2 = 0.$}

	\sssec{Exercise.}\label{sec-erighttr}
	If we start with $A_i$, $M$ and $\ov{M}$ in involutive geometry, we
	cannot derive the properties of the right triangles.  Other elements
	have to be prefered.  Make an appropriate choice and construct enough
	elements to determine $\theta$ and $\gamma.$

	\sssec{Answer to}
	\vspace{-18pt}\hspace{94pt}{\bf \ref{sec-erighttr}.}\\
	To obtain the coordinates, we replace\\
	\hth$m_0,$ $m_1,$ $m_2$ by $1,$ $\epsilon m_1,$ $\epsilon m_2,$\\
	and when the coordinates contain terms of different order of $\epsilon,$
	we neglect the terms of higher order.\\
	For instance,\\
	$q_0 = 1,$ $q_1 = -m_2,$ $q_2 = -m_1,$\\
	$\theta: m_0(m_1+m_2) X_1X_2 + m_1(m_2+m_0) X_2 X_0
		+ m_2(m_0+m_1) X_0 X_1 = 0,$\\
	becomes\\
	$\theta: (m_1+m_2) X_1X_2 + m_1 X_2 X_0 + m_2 X_0 X_1 = 0,$\\
	Of course many points or lines will coincide and some of the
	construction which are invalid must be replaced by other constructions
	but the coordinates do not have to be rewritten.
	For instance,\\
	$A_0 = \ov{M}M_0 = \ov{I}Ma_0,$,
	$\ov{m} = \ov{m}_0 = \ov{m}m_0 = ta_0,$
	$\ov{m}a_0 = \ov{m}k,$
	$\ov{M}A_0 = TAa_0,$
	$eul = ma_0$,
	$EUL = Imm_0,$
	$ta_0 = \ov{m},$
	$Aat_0 = \ov{M}_0.$

	We can start, for instance, with $A_i,$ $M$ and $K = (m_1+m_2,m_1,m_2),$
	on $mm_0,$ we construct, as usual, $ma_i,$ $M_i,$ $mm_i,$ $MA_i,$ $m_i,$
	$MM_i,$ $m.$  Then\\
	$mk := M \times K,$ $mk = [m_1-m_2,-m_1,m_2],$\\
	$at_i := K \times A_i,$ $at_0 = [0,m_2,-m_1],$
	$at_1 = [m_2,0,-(m_1+m_2)],$\\
	$Aat_i := at_i \times a_i,$ $Aat_1 = (m_1+m_2,0,m_2),$\\
	$Iat_i := m \times at_i,$ $Iat_0 = (m_1+m_2,-(2m_2+m_1),m_2),$\\
	$ta_i := A_i \times Iat_i,$ $ta_0 = [0,m_2,m_1],$
		$ta_1 = [m_2,0,m_1+m_2],$\\
	$TAa_0 := ta_0 \times a_0,$ $TAa[0] = (0,m_1,-m_2).$\\
	In general, the construction cannot be done for all 3 elements, but
	can done simultaneously for the elements with index 1 and 2, we can
	use $j$ for 1 and $-j$ for 2.\\
	$\theta = conic(A_0,ta_0,A_1,A_2,MM_0),$\\
	$\theta: (m_1+m_2)X_1X_2 + m_1 X_2X_0 + m_2 X_0X_1 = 0,$\\
	$\gamma := conic(M_i,Aat_0,A_0),$\\
	$\gamma: m_2X_1^2 + m_1 X_2^2 - (m_1+m_2)X_1X_2 -m_1 X_2X_0
		-m_2X_0X_1 = 0.$

	When we start with $J$ and $M$, the triangle is a right triangle if
	$I \times J_1\para a_2$ and $I \times J_2\para a_1$.  Moreover,
	$j_0^2 = p_{11},$ $m_1 = j_1(j_2+j_0)(j_2-j_0),$
	$m_2 = j_2(j_0+j_1)(j_0-j_1).$ The usual construction gives
	$\ov{M}_0 = Aat_0$,
	then $at_0 := A_0 \times Aat_0,$ $K := at_0 \times mm_0.$

	\sssec{Definition.}
	An {\em isosceles triangle} $\{A_i\}$ at $A_0$ is a triangle whose
	angles $A_0\: A_1\: A_2$ and $A_1\: A_2\: A_0$ are equal.\\
	What about right isosceles?

	\sssec{Theorem.}
	{\em A necessary and sufficient condition for a triangle to be
	an isosceles triangle at $A_0$ is that $m_1 = m_2.$}

	\sssec{Theorem.}
	{\em If $\{A B C\}$ is an isosceles triangle at $A$, 
	then the sides $A B$ and $A C$ are equal.}

	\sssec{Theorem.}
	{\em If $\{A B C\}$ is an isosceles triangle, then the angle
	$(A \times B,A \times C)$ is even.}

	\sssec{Definition.}
	A triangle  $\{A B C\}$ is an {\em equilateral triangle}  iff  it is
	isosceles at $B$ and $C$.

	\sssec{Theorem.}
	{\em A necessary and sufficient condition for a triangle to be
	an equilateral triangle at $A_0$ is that $m_0 = m_1 = m_2.$}

	\sssec{Theorem.}
	{\em If a triangle is equilateral, then all its angles 
	$A B C,$ $B C A$ and $C A B$ are equal and all its sides are equal.}

	\sssec{Definition.}
	A triangle which is neither a right triangle nor an
	isosceles triangle, and therefore not an equilateral triangle is 
	called a {\em scalene triangle}.

	\sssec{Theorem.}
	{\em A necessary and sufficient condition for a triangle to be
	a scalene triangle is that $m_0$, $m_1$ and $m_2$ be distinct.}

	\sssec{Definition.}
	A triangle which is an isosceles triangle at $A$ 
	but is not an equilateral triangle is called a {\em proper isosceles }
	triangle.

	\sssec{Theorem.}
	{\em If a triangle is equilateral, then all its angles 
	$A B C,$ $B C A$ and $C A B$ are equal and all its sides are equal.}

	The following Definitions and Theorem are only meaningfull in
	Minkowskian Geometry.\footnote{4.3.89}

	\sssec{Definition.}
	An {\em isotropic triangle} is a triangle with one isotropic side.

	\sssec{Definition.}
	A {\em doubly isotropic triangle} is a triangle with 2 isotropic sides

	\sssec{Theorem.}
	\enumb
	\item {\em A necessary and sufficient condition for a triangle to be
		isotropic is that the barycenter be on the complementary
		triangle.}
	\item {\em A necessary and sufficient condition for a triangle to be
		doubly isotropic is that the barycenter be one of the vertices
		of the complementary triangle.}
	\enume

	\sssec{Theorem.}
	{\em If $a[0]$ is an isotropic line, then}
	\enumb
	\item m$_1$ + m$_2$ = 0,
	\item {\em the circumcircle degenerates into $a_0$ and the line}
	[0,m$_0$+m$_1$,-(m$_0$-m$_1$)],
	\enume

	\sssec{Theorem.}
	{\em If $a_1$ and $a_2$ are isotropic lines, then}
	\enumb
	\item m$_0$ = 1, m$_1$, m$_2$ = -1, 
	\item {\em the circumcircle degenerates into $a_1$ and $a_2$.}
	\enume

	\ssec{Other special triangles.}
\setcounter{subsubsection}{-1}
	\sssec{Introduction.}
	There are many other types of triangles that can be defined.
	I will give here 2 examples which allow the constructions of
	configurations of the type 9 * 3 \& 9 * 3, distinct from that of
	Pappus.  For the first one, if we choose $B_1 = M_1,$ $B_2 = M_2$ and
	$C_1 = \ov{M},$ the construction of Section 19 gives
	$C_0 = Mam_2,$ $C_2 = Tara_0,$ $B_0 = tara_2 \times tarb_2.$
	from P19.7, follows that $B_0 \cdot a_0 = 0$ iff $q_0 = 0.$
	This suggest the definition of a triangle of Tarry and the construction
	of the 1-Pappus configuration.\\
	A similar approach determines the construction of the 2-Pappus
	configuration.

	\sssec{Definition.}
	The {\em 1-Pappus configuration} is the set of points\\
	\hth 1-Pappus($A_i,B_i,C_i$),\\
	such that $\{A_i\}$, $\{B_i\}$, $\{C_i\}$ are 3 triangles and
	incidence($A_i,C_i,C_{i-1}$), incidence($B_{i+1},B_{i-1},C_i$).

	\sssec{Definition.}
	The {\em 2-Pappus configuration} is the set of points\\
	\hth 2-Pappus($A_i,B_i,C_i$),\\
	such that $\{A_i\}$, $\{B_i\}$, $\{C_i\}$ are 3 triangles and
	incidence($A_i,B_i,C_i$), incidence($B_i,C_{i+1},C_{i-1}$).

	\sssec{Definition.}
	A {\em triangle of Tarry} is a triangle which is not a right triangle
	whose point of Tarry is well defined and coincides with one of the
	vertices of the triangle.

	\sssec{Theorem.}
	{\em A necessary and sufficient condition for a triangle to be a
	triangle of Tarry at $A_0$ is that $q_0 = 0.$}

	The proof follows at once from $m_i \neq 0$ and from P16.3.\\
	Moreover, if $q_0 = q_1 = 0,$ then the point of Tarry is undefined.

	\sssec{Corollary.}
	{\em A necessary and sufficient condition for a triangle to be a
	triangle of Tarry at $A_0$ is that the orthocenter be on the conic of
	Tarry.}

	\sssec{Theorem.}
	\hth$\ov{M} \cdot `Tarry[0] \Rightarrow$
	1-Pappus($A_i,Tarb_2\:M_1\:M_2,Mam_2\:\ov{M}\:Tara_0$).

	\sssec{Theorem.}
	{\em In Involutive Geometry, $A_0 = (x,y,1),$ $A_1 = (0,0,1),$
	$A_2 = (1,0,1)$, $M = (1+x,y,3),$ $\ov{M} = (xy,x(1-x),y)$ is a triangle
	of Tarry iff\\
	\hth$u^2 -(1+2y^2)u + y^4 = 0$ and $x^2 -x + u = 0.$}

	Proof:
	Assuming that $m = [0,0,1]$ and $X_0^2 + X_1^2 = X_2^2$ is a circle,
	a trivial computation determines $M$ and $\ov{M}$ as given.
	The conic of Tarry is\\
	\hth$a_0 \bigx a_0 + k a_1 \bigx a_2 = 0,$\\
	where $k$ is determined in such a way that $M \cdot `Tarry = 0,$
	this gives $k = 1$. To insure that $\ov{M}$ is on the conic of Tarry
	gives after division by $u := x(1-x),$\\
	$u + (y^2-u)(u-y^2) = 0,$ a simple discussion determines that
	$0 < y \leq \frac{\sqrt{3}}{2}.$

	\sssec{Definition.}
	An {\em Eulerian triangle} is a triangle for which the line of Euler
	is parallel to one of its sides.

	\sssec{Theorem.}
	{\em A necessary and sufficient condition for a triangle to be an
	Eulerian triangle for side $a_0$ is that $2m_0 = m_1+m_2.$}

	The proof follows at once from P1.17.

	\sssec{Theorem.}
	\hth$eul \para a_0 \Rightarrow$
	2-Pappus($A_i,\ov{M}_0\:M_1\:M_2,\ov{M}\:Tarc_0\:Tar\ov{c}_0)$.

	\ssec{Geometry of the triangle. V.}
	(bissectrices)


	\begin{center}
	{\Huge CHAPTER II\\[40pt]
	FINITE PROJECTIVE\\[20pt] GEOMETRY}\\[40pt]
	\end{center}

\setcounter{section}{89}
	\section{Answers to problems and miscellaneous notes.}
	\sssec{Answer to}
	\vspace{-18pt}\hspace{100pt}{\bf \ref{sec-edesar}.}\\
	Let $A_0 = (1,0,0),$ $A_1 = (0,1,0),$ $A_2 = (0,0,1),$ $C = (1,1,1),$
	$B_0 = (b_0,1,1),$ $B_1 = (1,b_1,1),$ $B_2 = (1,1,b_2),$\\
	by hypothesis, $b_i \neq 1,$ $b_2 \neq 0,$ $b_1b_2 \neq 1$ and
	$2-b_0-b_1-b_2+b_0b_1b_2\neq 0$ (because of $\{B_i\}.$\\
	$a_0 = [1,0,0],$ $b_0 = [1-b_1b_2,-(1-b_2),-(1-b_1)],$\\
	$c_0 = [0,1,-1],$ $C_0 = (0,1-b_1,-(1-b_2)),$\\
	$c = [(1-b_1)(1-b_2),(1-b_2)(1-b_0),(1-b_0)(1-b_1)],$ $d = [0,b_2,-1],$
\\
	$D = (b_2,1,b_2),$ $e = [b_2(1-b_1),b_2(1-b_0),-1+b_0-b_2+b_1b_2],$\\
	$E = (1-b_0+b_2b_0-b_1b_2-b_2^2b_0+b_2^2b_0b_1,
		1-b_0+b_2-2b_1b_2-b_2^2+b_2^2b_1+b_0b_1b_2,\\
	\hth	b_2(2-b_0-b_1-b_2+b_0b_1b_2)),$
	$f = [1-b_1b_2,-b_2(1-b_2),-b_2(1-b_1)],$
	$F = (b_2(1-b_1),0,1-b_1b_2),$ $G = (1-b_2,1-b_1b_2,0),$\\
	$g = [1-b_1b_2,-(1-b_2),-b_2(1-b_1)],$\\
	$X = (0,1,b_2),$ $Y = (1-b_0-b_2+b_0b_1b_2,-(1-b_1),-b_2(1-b_1)),$\\
	$Z = (1-b_2+b_2^2-b_2^1b_1,1-b_1b_2,b_2(1-b_2)).$

	\ssec{Answer to exercises.}

	\sssec{Exercise. [Pappus]}
	Define\\
	\hth$\alpha := (A_0 * A_1) \cdot A_2, \beta := (B_0 * B_1) \cdot B_2,$\\
	\hth$\alpha_{1,2} := (A_1 * A_2) \cdot B_1,
		\beta_{1,2} := (B_1 * B_2) \cdot A_1,$\\
	\hth$\alpha_{2,0} := (A_2 * A_0) \cdot B_2,
		\beta_{2,0} := (B_2 * B_0) \cdot A_2,$\\
	\hth$\alpha_{0,1} := (A_0 * A_1) \cdot B_0,
		\beta_{0,1} := (B_0 * B_1) \cdot A_0,$\\
	Using 2.3.17.0,\\
	$C_0 = (A_1 * B_2) * (A_2 * B_1
	 = ( (A_1 * B_2) \cdot B_1) B_2 - ( (B_2 * A_2) \cdot B_1) A_1$\\
	\hth$   = \alpha_{1,2} B_2 - \beta_{1,2} A_1,$\\
	similarly,\\
	$C_1 = \alpha_{2,0} B_0 - \beta_{2,0} A_2,$\\
	$C_2 = \alpha _{0,1} B_1 - \beta_{0,1} A_0,$ therefore\\
	$(C_0 * C_1) \cdot C_2
		= \beta  \alpha_{1,2} \alpha_{2,0} \alpha_{0,1}
			- \alpha  \beta_{1,2} \beta_{2,0} \beta_{0,1}
		+ \alpha_{2,0} \alpha_{1,2} \beta_{2,0} \beta_{0,1}
			- \beta_{2,0} \alpha_{1,2} \alpha_{2,0} \beta_{0,1}$\\
	\hth$	+ \alpha_{0,1} \beta_{1,2} \alpha_{2,0} \beta_{0,1}
			- \beta_{0,1} \beta_{1,2} \alpha_{2,0} \alpha_{0,1}
		+ \alpha_{1,2} \beta_{1,2} \beta_{2,0} \alpha_{0,1}
			- \beta_{1,2} \alpha_{1,2} \beta_{2,0} \alpha_{0,1}$\\
	\hth$	= \beta  \alpha_{1,2} \alpha_{2,0} \alpha_{0,1}
			- \alpha  \beta_{1,2} \beta_{2,0} \beta_{0,1}$.\\
	Therefore, if the points $A_0,$ $A_1,$ $A_2$ and the points $B_0,$
	$B_1,$ $B_2$ are
	collinear, $\alpha = 0$ and $\beta = 0,$ therefore
	$(C_0 * C_1) \cdot C_2 = 0$
	and the points $C_0,$ $C_1,$ $C_2$ are collinear by 2.3.18.

	\sssec{Exercise.}
	(Harmonic quatern).
	$a = [0,0,1],$ let $A = (0,0,1),$ $A \times K = [k,-1,0],$\\
	let $B = (1,k,1),$ $l\neq 0$ and 1.\\
	$B \times L = [l,-1,k-l],$ $A \times M = [m,-1,0],$ $D = (l-k,ml-mk,l-m),$
	$D \times K = [k(l-m), m-l, (l-k)(m-k)],$ $A \times L = [l,-1,0],$
	$C = (m-k,lm-lk,l-m),$ $B \times C = [2kl-km-lm,2m-k-l,(k-l)(k-m)],$
	$N = (2m-l-k,km+lm-2kl,0).$

	\sssec{Exercise.}
	(Projectivity).
	Choose $b = [0,1,0],$ $a = [1,0,0]$ and $P = (0,1,1).$
	$c = [0,0,1],$ $S = (1,0,-k),$ $T = (1,0,-l),$ $Q = (0,m,k),$
	$N = [k,ml,0].$

	\sssec{Exercise.}
	(Projectivity with 3 pairs).
	$C_0 = (1,0,-1),$ $C_j = (1,a_j,-1-a_j),$ $j > 0,$ with obvious
	notation,\\
	$D_j = (1,-a_jb_j, -1-a_j),$ $j = 1,2,$
	$t = ( a_1a_2(b_2-b_1) + a_2b_2-a_1b_1, a_1-a_2, a_2b_2-a_1b_1 ),$
	$D_j = ( a_2-a_1, a_1b_1-a_2b_2 - a_1a_2(b_2-b_1),
	(a_1-a_2)(1+a_j) ),$ $j > 2,$
	hence $B_j$ given above.

	\sssec{Answer to 2.6.14. }
	In part,
	the coefficients of $A_0$ and $A_1$ in $A_l$ and $B_l$ must be
	proportional, therefore,
	$\frac{f_0 t_0 + f_1 t_1}{ f_2 t_0 + f_3 t_1} = \frac{t_0}{t_1},$
	this gives 1.

	\sssec{Answer to 2.3.7.}
	For $p = 2,$ if $A[] = ((1,0,0), (1,0,1), (1,1,0), (1,1,1)),$ and the
	diagonal points are $B_i,$ 
	$A_0 \times A_1 = [0,1,0],$ $A_2 \times A_3 = [1,1,0],$ $B_0 = (0,0,1).$
	$A_1 \times A_2 = [1,1,1],$ $A_3 \times A_0 = [0,1,1],$ $B_1 = (0,1,1).$
	$A_0 \times A_2 = [0,0,1],$ $A_1 \times A_3 = [1,0,1],$ $B_2 = (0,1,0).$
	The diagonal points are on [1,0,0].
	For $p = 4,$ the coordinates are 0, 1, $x,$ $y = 1+x.$
	The addition and multiplication tables are\\
	$\begin{array}{rllccccrllcccc}
	\hth &	+&\vline&  0 &1 &x &y&\htn& \cdot& \vline& 0 &1 &x &y\\
		\cline{2-7}\cline{9-14}
		&0&\vline& 0 &1 &x &y& &0&\vline& 0 &0 &0 &0\\
		&1&\vline& 1 &0 &y &x& &1&\vline& 0 &1 &x &y\\
		&x&\vline& x &y &0 &1& &x&\vline& 0 &x &y &1\\
		&y&\vline& y &x &1 &0& &y&\vline& 0 &y &1 &x
	\end{array}$\\
	If $A[] = ((1,0,0),$ $(1,0,1),$ $(1,x,0),$ $(1,y,1)),$\\
	$A_0 \times A_1 = [0,1,0],$ $A_2 \times A_3 = [x,1,1],$ $B_0 = (1,0,x).$
	$A_1 \times A_2 = [x,1,x],$ $A_3 \times A_0 = [0,1,y],$ $B_1 = (1,1,x).$
	$A_0 \times A_2 = [0,0,1],$ $A_1 \times A_3 = [y,0,y],$ $B_2 = (0,1,0).$
	The diagonal points are on $[x,0,1].$

	\sssec{Answer to 2.5.10. }
	In part,\\
	$C_2 = (c,1-c,1),$ $a_0 = (1,c,-c-1),$ $b = [1,0,0],$
	$B_2 = (c^2-c+1,-c+1,1).$ 
	A geometric condition is $a_0 \cdot C_2 = 0$ or $b \cdot B_2 = 0.$
	The configuration is then of type\\
	\hth$6 * 4 + 3 * 3 + 1 * 2  \:\&\:  2 * 4 + 9 * 3.$
	\sssec{Notes.}
	On \ref{sec-tbaker}:\\
	Commutativity implies that if\\
	\hth$	J' := (B \times P) \times (E \times Q),$
		$J = (b(a-1),b^2(a-1),a(b-1)),$\\
	\hth$	K' := (A \times P) \times (J' \times M),$
		$K = (b(a-1),ab(a-1),a(b-1)),$\\
	\hth then\\
	\hth$	L \cdot (J' \times K') = 0.$

	The construction\\
	\hth$D'' := (R\times T') \times a,$ with D'' = (1,b+c)\\
	is related to the associative property\\
	\hth$	(a+b) + c = a + (b+c).$

	Before\ref{sec-lsymmmat}\\
	\sssec{Theorem.}
	{\em \ldots describe the degenerate conic, perhaps in 2.10.8
	\ldots determine the collineation which leave a general conic fixed
	also special case when it is a conic.}

	\sssec{Examples.}
	For p = 5,\\
	{\bf C0}
	\begin{verbatim}
	 0 1 2 3 4 5 6 7 8 9101112131415161718192021222324252627282930
	 4 4 2 4 4 4 4 4 2 4 4 4 4 4 1 4 4 4 4 4 2 2 4 4 4 1 4 1 4 4 4
	\end{verbatim}
	\hth$ N = \mat{2}{-2}{0}{2}{0}{-1}{2}{0}{0},$
	$N^I = \mat{0}{2}{0}{0}{0}{-1}{-2}{-2}{1}.$\\
	Point conic and its mapping 0  1, 1  6, 10 27, 14  4, 23  8, 27 29,\\
	Line conic and its mapping 1 10, 4 14, 6  0, 8  1, 27 23, 29 27,\\
	Points on line conic and tangent, 0 16, 1 0, 10 22, 14 4, 23 14, 27 29,
\\
	Lines on line conic and contact, 1 7, 4 14, 6 5, 8 13, 27 29, 29 27.\\
	The equation of the point conic is $X^2$ + 2YZ + ZX = 0.\\
	The equation of the line conic is $-2z^2$ + yz - zx + xy = 0.\\
	{\bf C1}
	\begin{verbatim}
	 0 1 2 3 4 5 6 7 8 9101112131415161718192021222324252627282930
	 6 6 3 6 6 6 6 3 6 6 6 6 6 3 6 6 6 1 6 3 6 6 6 6 6 3 3 6 6 6 6
	\end{verbatim}
	\hth$ N = \mat{-2}{2}{0}{-2}{0}{-1}{-2}{0}{0},$
	$N^I = \mat{0}{-2}{0}{0}{0}{-1}{2}{2}{1}.$\\
	Point conic and its mapping 0  1, 1  6, 10 22, 12  5, 18 25, 29  9,\\
	Line conic and its mapping 1 10, 5 18, 6  0, 9  1, 22 12, 25 29,\\
	{\bf C2}
	\begin{verbatim}
	 0 1 2 3 4 5 6 7 8 9101112131415161718192021222324252627282930
	10 51010 2101010 5 21010 2 510101010 510 21010 11010 110 51010
	\end{verbatim}
	\hth$ N = \mat{0}{1}{-2}{1}{1}{0}{-1}{1}{0},$
	$N^I = \mat{0}{0}{2}{-2}{-2}{-1}{2}{-2}{-1}.$\\
	Point conic and its mapping 0  6, 5 18, 6  5, 15 26, 19 22, 21 25,\\
		22 28, 23 14, 24 17, 25  4, 27 10,\\
	The center is (23), the points are on [21] or [2].\\
	Line conic and its mapping 4  6, 5 24, 6  5, 10 22, 14 23, 17 27,\\
		18 25, 22 21, 25 15, 26  0, 28 19,\\
	The central line is [18], the lines pass through (4) or (12).\\
	The equation of the point conic is $Y^2$ + YZ + 2ZX + 2XY = 0.\\
	The equation of the line conic is $y^2$ - $2z^2$ -yz - 2zx + xy = 0.\\
	{\bf C3}
	\begin{verbatim}
	 0 1 2 3 4 5 6 7 8 9101112131415161718192021222324252627282930
	 4 4 4 4 4 2 4 4 4 1 4 4 4 2 4 1 4 2 4 4 4 2 4 4 4 4 4 4 4 4 1
	\end{verbatim}
	Point conic and its mapping 7 10, 9 19, 11 30, 12 18, 29 27, 30  7,\\
	Line conic and its mapping 7 30, 10 12, 18 29, 19  9, 27 11, 30  7,\\
	{\bf C4}
	\begin{verbatim}
	 0 1 2 3 4 5 6 7 8 9101112131415161718192021222324252627282930
	 5 5 5 5 5 5 5 5 5 5 5 5 5 5 5 5 5 1 5 5 5 5 5 5 5 5 5 5 5 5 5
	\end{verbatim}
	Point conic and its mapping 11 29, 14 12, 17 24, 19 28, 26 15, 27 23,\\
	Line conic and its mapping 12 26, 15 19, 23 14, 24 17, 28 11, 29 27,\\
	{\bf C5}
	\begin{verbatim}
	 0 1 2 3 4 5 6 7 8 9101112131415161718192021222324252627282930
	 5 5 5 1 5 5 5 5 5 5 5 5 5 5 5 5 5 5 5 5 5 5 5 5 5 5 5 5 5 5 5
	\end{verbatim}
	Point conic and its mapping 2 27, 3  3, 8  8, 10  1, 16 19, 18 16,\\
	Line conic and its mapping 1  8, 3  3, 8 18, 16 16, 19  2, 27 10,\\
	{\bf C6}
	\begin{verbatim}
	 0 1 2 3 4 5 6 7 8 9101112131415161718192021222324252627282930
	 5 5 5 5 5 5 5 5 5 5 5 5 5 5 5 1 5 5 5 5 5 5 5 5 5 5 5 5 5 5 5
	\end{verbatim}
	Point conic and its mapping 9 14, 10 12, 15 26, 18 29, 23 21, 24 17,\\
	Line conic and its mapping 12 18, 14 23, 17 10, 21 24, 26 15, 29  9,\\
	{\bf C7}
	\begin{verbatim}
	 0 1 2 3 4 5 6 7 8 9101112131415161718192021222324252627282930
	 5 5 5 1 5 5 5 5 5 5 5 5 5 5 5 5 5 5 5 5 5 5 5 5 5 5 5 5 5 5 5
	\end{verbatim}
	Point conic and its mapping 0 16, 3  3, 8  8, 11 28, 16  0, 28 11,\\
	Line conic and its mapping 0 11, 3  3, 8 28, 11  0, 16 16, 28  8,\\
	{\bf C8}
	\begin{verbatim}
	 0 1 2 3 4 5 6 7 8 9101112131415161718192021222324252627282930
	10101010 5 1 51010 2101010 2 51010 1101010 2101010 51010 510 2
	\end{verbatim}
	\hth$ N = \mat{1}{0}{0}{1}{1}{0}{1}{0}{2},$
	$N^I = \mat{1}{-1}{2}{0}{1}{0}{0}{0}{-2}.$\\
	Point conic and its mapping 17 24, the center is (17).\\
	Line conic and its mapping 24 17, the central line is [24].\\
	The equation of the point conic is, with $\delta$ $^2$ = 2,\\
	(X - $(2+2\delta$ )Y - $(2+\delta$ )Z) (X - $(2-2\delta$ )Y
	- $(2-\delta$ )Z) = 0.\\
	The equation of the line conic is
	(x + $(2+2\delta$ )y + $(1-2\delta$ )z) (x + $(2-2\delta$ )y
	+ $(1+2\delta$ )z) = 0.\\
	{\bf C9}
	\begin{verbatim}
	 0 1 2 3 4 5 6 7 8 9101112131415161718192021222324252627282930
	 5 5 1 5 5 5 5 5 5 5 5 5 5 5 5 5 5 5 5 5 5 5 5 5 5 5 5 5 5 5 5
	\end{verbatim}
	Point conic and its mapping 2  5, 5  2, 15 25, 20 30, 26 14, 29 11,\\
	Line conic and its mapping 2 15, 5  2, 11 26, 14 20, 25 29, 30  5,\\
	{\bf C10}
	\begin{verbatim}
	 0 1 2 3 4 5 6 7 8 9101112131415161718192021222324252627282930
	 1 3 6 6 6 6 3 6 6 6 6 3 6 6 6 6 3 6 6 6 6 3 6 6 6 6 3 6 6 6 6
	\end{verbatim}
	Point conic and its mapping 3  3, 4  4, 8 28, 9 29, 13  8, 14  9,\\
	Line conic and its mapping 3 13, 4 14, 8  8, 9  9, 28  4, 29  3,\\
	{\bf C11}
	\begin{verbatim}
	 0 1 2 3 4 5 6 7 8 9101112131415161718192021222324252627282930
	 1 3 3 3 3 3 3 3 3 3 3 3 3 3 3 3 3 3 3 3 3 3 3 3 3 3 3 3 3 3 3
	\end{verbatim}
	Point conic and its mapping 2  5, 5  2, 12 14, 15 13, 17 24, 20 23,\\
	Line conic and its mapping 2 15, 5 12, 13 17, 14 20, 23  2, 24  5,\\
	{\bf C12}
	\begin{verbatim}
	 0 1 2 3 4 5 6 7 8 9101112131415161718192021222324252627282930
	 5 1 5 5 5 5 5 5 5 5 5 5 5 5 5 5 5 5 5 5 5 5 5 5 5 5 5 5 5 5 5
	\end{verbatim}
	\hth$ N = \mat{-2}{-2}{0}{-2}{0}{0}{-2}{0}{-1},$
	$N^I = \mat{0}{2}{0}{2}{-2}{1}{0}{0}{-1}.$\\
	Point conic and its mapping 1  6, 7 15, 8 13, 15 25, 16 18, 19 16,\\
	Line conic and its mapping 6  1, 13 15, 15 19, 16 16, 18  8, 25  7,\\
	{\bf C13}
	\begin{verbatim}
	  0 1 2 3 4 5 6 7 8 9101112131415161718192021222324252627282930
	  1 1 1 1 1 1 1 1 1 1 1 1 1 1 1 1 1 1 1 1 1 1 1 1 1 1 1 1 1 1 1
	\end{verbatim}
	\hth$ N = \mat{0}{1}{1}{1}{0}{1}{1}{1}{0},$
	$N^I = \mat{2}{-2}{-2}{-2}{2}{-2}{-2}{-2}{2}.$\\
	Point conic and its mapping 0 11, 1  7, 6  2, 13 15, 17 27, 24 30,\\
	Line conic and its mapping 2  6, 7  1, 11  0, 15 13, 27 17, 30 24,\\
	The equation of the point conic is YZ + ZX + XY = 0.\\
	The equation of the line conic is $x^2$ + $y^2$ + $z^2$ - 2yz
	- 2zx - 2xy = 0.\\
	{\bf C14}
	\begin{verbatim}
	  0 1 2 3 4 5 6 7 8 9101112131415161718192021222324252627282930
	  2 2 1 2 2 2 2 2 2 1 2 2 2 2 2 1 1 2 2 2 2 2 1 2 2 1 2 2 1 2 2
	\end{verbatim}
	\hth$ N = \mat{1}{-2}{1}{1}{-1}{0}{0}{2}{-1},$
	$N^I = \mat{1}{1}{2}{0}{-1}{-2}{1}{1}{1}.$\\
	Point conic and its mapping 2 15, 9 29, 15  7, 16 18, 22 21, 28  4,\\
	The points are on [19].\\
	Line conic and its mapping 4 28, 7 15, 15  2, 18 16, 21 22, 29  9,\\
	The lines pass through (25).\\
	{\bf C15}
	\begin{verbatim}
	 0 1 2 3 4 5 6 7 8 9101112131415161718192021222324252627282930
	 2 2 2 2 2 1 2 2 2 2 1 2 1 2 1 2 2 2 1 2 2 2 1 2 2 2 1 2 2 2 2
	\end{verbatim}
	\hth$ N = N^I = \mat{1}{-2}{-2}{-2}{1}{-2}{-2}{-2}{1}$\\
	Point conic and its mapping 19 28, 20 23, 23 20, 25 29, 28 19, 29 25,\\
	Line conic and its mapping 19 28, 20 23, 23 20, 25 29, 28 19, 29 25,\\
	The equation of the point conic is $X^2 + Y^2 + Z^2 + YZ + ZX
	+ XY = 0.$\\
	The equation of the line conic is $x^2 + y^2 + z^2 + yz + zx
	+ xy = 0.$

	\sssec{Answer to }
	\vspace{-18pt}\hspace{100pt}{\bf \ref{sec-econic13}.}\\
	For $p = 13,$ points on the conic are, (0,1,1) = (2),
	(0,1,2) = (3), (1,0,1) = (15), (1,0,4) = (18), (1,1,0) = (27),
	(1,2,0) = (40).\\
	The point conic is\\
	  2,  3, 15, 18, 27, 35, 40, 51,133,135,146,151,158,168,\\
	the line conic is\\
	111, 83,156,121,179, 22, 98,148,112,129, 86, 25,165,166.\\
	The representative matrix is\\
	\hth$\mat{1}{-4}{1}{-4}{-6}{-2}{1}{-2}{-3}$\\
	(0,1,2) = (13) and (1,0,12) = (26) are on [1,1,1], the polars are
	[1,6,5] = [97] and [0,1,11] = [12].  Hence the pole of [1,1,1] is
	(1,6,3) = (95).\\
	$1.-6+3.1-4.5-2.(-3) = -17 = -4$\\
	$6.(-3)+8.6-2.(-6)-(-5).5) = 67 = 2,$\\
	$2.5+4.2-(-5).(-3)-4.(-6) = 27 = 1, (-4,2,1) = (1,6,3).$

	\sssec{Answer to}
	\vspace{-18pt}\hspace{100pt}{\bf \ref{sec-econic5}.}\\
	$A\times B = [1,-1,1],$ $C\times D =[2,1,-5],$ $A\times D = [2,-1,1],$
	$B\times C = [1,1,-3],$ $k_1 = 1.3,$ $k_2 = -1.3,$ therefore the conic
	is, after dividing by 3,\\
	$(X_0-X_1+X_2)(2X_0+X_1-5X_2) + (2X_0-X_1+X_2)(X_0+X_1-3X_2) = 0,$\\
	which gives twice the result of the Example.

\setcounter{section}{3}
\setcounter{subsection}{89}
	\sssec{Answer to}
	\vspace{-18pt}\hspace{100pt}{\bf \ref{sec-eseltetra}.}
	\enumb
	\item  For $q = 2,$ the primitive polynomial giving the selector
	0, 1, 3, is\\
	\hth$    I^3 + I + 1.$\\
	    The auto-correlates are 0  11  2  7  8.\\
	    The selector function is\\
	$\begin{array}{cclllllllllllllllllllll}
	i&\vline&0&1&2&3&4&5&6&7&8&9&10&11\\
	f(i)&\vline&&0&14&1&0&16&16&14&14&16\\
	type&\vline&F_0&V_0&F_4&V_2&T_0&T_2&V_1&F_3&F_1&T_3&E_2&F_2\\
	\cline{2-14}
	i&\vline&12&13&14&15&16&17&18&19&20\\
	f(i)&\vline&4&1&0&1&0&4&4&16&1\\
	type&\vline&T_4&E_1&P&V_3&E_0&T_1&E_3&V_4&E_4
	\end{array}$
	\item  The correspondence between the selector notation and the
	    homogeneous coordinates for points and lines is\\
	$\begin{array}{ccll}
	\hth&	i&	I^i&	     i^* \\
	&	0&	1   &    6^*: 1,2,4,\\
	&	1&	I&	1^*:	0,2,6,\\
	&	2&	I^2&	     0^*: 0,1,3,\\
	&	3&	I+1&	5^*:	2,3,5,\\
	&	4&	I^2+I&	   3^*: 0,4,5,\\
	&	5&	I^2+I+1&	 4^*: 3,4,6,\\
	&	6&	I^2+1&	   2^*: 1,5,6.
	\end{array}$
	\item  The matrix representation is\\
	\hth$    M = \mat{1}{0}{1}{0}{1}{0}{1}{0}{0},
		 M^{-1} = \mat{0}{0}{1}{0}{1}{0}{1}{0}{1}.$
	    and the equation satisfied by the fixed points is $(X_0+X_1)^2$ = 0.
	\item  The degenerate conic through 0, 1, 2 and 5 with tangent $5^*$
	at 5, is represented by the matrix\\
	\hth$N = \mat{0}{1}{1}{1}{0}{0}{1}{0}{0}.$\\
	    The polar of 0 is $0^*,$ of 1 is $0^*,$ of 2 is $5^*,$ of 4 is
		$4^*,$ of 5 is $5^*$ of 6 is $6^*$ and of 3 is undefined.
	    The equation in homogeneous coordinates is $X_0(X_1+X_2) = 0.$
	\item  A circle with center 14 can be constructed as follows.  I first
	    observe that a direction must be orthogonal to itself.  Indeed,
	    if 0 is a direction, the others form an angle 1,2,3,4 $\umod{5},$
	    we cannot play favorites and must choose 0.  If $A_0 = 1,$
	    $C \times A_0$ and therefore the tangent has direction 0,
	    $A_0 \times A_{i+1}$ has
	    direction $i \umod{5}$ or are the points 0, 7, 8, 2, 11.

	    It is natural to choose the pentagonal face-point as 14, and the
	    edge-points on the pentagon as 0, 8, 11, 7, 2.  The points on the
	    circle 1, 6, 3, 15, 19 are chosen as the vertex-points opposite
	    the corresponding edge-point, 1 opposite 0,  6 opposite 8, \ldots.
	    This gives the types, with subscripts indicated in 0. and the
	    definition:

	    The points are represented on the 5-anti-prism as follows.
	    The pentagonal face-point, P,
	    the 5 triangular face-points, $T_i,$ 
	    the 5 vertex-points, $V_i,$ 
	    the 5 triangular-triangular edge-points, $E_i,$ 
	    the 5 pentagonal-triangular edge-points $F_i.$ 

	    The lines are represented on the 5-anti-prism as follows.
	    The pentagonal face-line, f, which is incident to $F_i,$ 
	    the 5 triangular face-lines, $t_i,$ which are incident to $F_i,$ 
		$F_i,$ $T_{i+1},$ $T_{i-1},$ $E_{i+2},$ $E_{i-2}$.  If $f$ is
	the pentagonal edge of
		$t_i$	and $V,$ $V'$ are on $f$, $F_i$ is on it, $T_{i+1}$
	$(T_{i-1})$ share $V$ $(V'),$
		$E_{i+2}$ $(E_{i-2})$ are on an edge through $V$ $(V')$ not on
		$t_i$ \\
	    the 5 vertex-lines, $v_i,$ which are incident to\\
		$F_i,$ $V_{i+2},$ $V_{i-2},$ $E_{i+1},$ $E_{i-1}.$  If $t$ is
	the face with $v_i$ on its
		pentagonal edge these are all the vertices, and edge-points on
		it distinct from $v_i.$ \\
	    the 5 triangular-triangular edge-lines, $e_i,$ which are incident to
		$F_i,$ $T_{i+2},$ $T_{i-2},$ $V_{i+1},$ $V_{i-1}.$  $V_{i+1}$
	and $V_{i-1}$ are on the same edge as $e_i,$ the line which joins the
	center $C$ of the antiprism
		to $E_{i}$ is parallel to the edge containing $F_i,$ $T_{i+2}$
	and $T_{i-2}$ are
		the triangular faces which are not adjacent to $E_i$ or $F_i.$ 

	    the 5 pentagonal-triangular edge-lines. $f_i,$ which are incident to
		$P,$ $T_i,$ $V_i,$ $E_i,$ $F_i.$  $T_i$ is adjacent to $f_i,$
	$V_i$ is opposite $f_i,$ 
	$E_i$	joined to the center of the antiprism is parallel to $T_i.$ 
	\enume

	\sssec{Answer to}
	\vspace{-18pt}\hspace{100pt}{\bf \ref{sec-eselcube}.}\\
	For $p = 3,$
	\enumb
	\item  The primitive polynomial giving the selector 0, 1, 3, 9 is
	    $I^3 - I - 1.$
	\item  The correspondence between the selector notation and the
	    homogeneous coordinates for points and lines is\\
	$\begin{array}{rllll}
		&i	&I^i	&     i^*\\
		&0	&1      &12^*:& 1,2,4,10,\\
		&1	&I	&1^*:&	0,2,8,12,\\
		&2	&I^2	&     0^*:& 0,1,3,9,\\
		&3	&I+1	&7^*:&	2,6,7,9,\\
		&4	&I^2+I	&   3^*:& 0,6,10,11,\\
		&5	&I^2+I+1	& 4^*:& 5,9,10,12,\\
		&6	&I^2+2I+1	&1^{0*}:& 3,4,6,12,\\
		&7	&I_2+I+2	& 6^*:& 3,7,8,10,\\
		&8	&I^2+1	&   2^*:& 1,7,11,12,\\
		&9	&I+2    &11^*:& 2,3,5,11\\
		&10	&I_2+2I	&  9^*:& 0,4,5^,7, \\
		&11	&I^2+2I+2	&5^*:& 4,8,9,11,\\
		&12	&I^2+2	&  8^*:& 1,5,6,8.
	\end{array}$

	\item  The matrix representation of the polarity $i$ to $i^*$ is\\
	\hth$M = \mat{1}{0}{1}{0}{1}{0}{1}{0}{0},$
		$M^{-1} = \mat{0}{0}{1}{0}{1}{0}{1}{0}{2}.$\\
	    The equation satisfied by the fixed points is
	$X_0^2 + X_1^2 + 2 X_2 X_0  = 0.$
	\item  The degenerate conic through 0, 1, 2 and 5 with tangent $4^*$ at
	$5$, is obtained by constructing the quadrangle-quadrilateral
	configuration starting with $P = 5$ and $Q_i = \{0,1,2\}.$  We obtain
	$q_i$ = $\{3^*,2^*,7^*\},$  which are the tangents at $Q_i.$  The
	matrix representation is\\
	\hth$N = \mat{0}{1}{1}{1}{0}{1}{1}{1}{0}$ with equation
	$X_1X_2 + X_2X_0 + X_0X_1 = 0.$\\
	    We can check that the polar of $10 = 3^* \times 4^*$ is
		$9^* = 0 \times 5.$
	\enume

	\sssec{Answer to}
	\vspace{-18pt}\hspace{100pt}{\bf \ref{sec-etruncdodeca1}.}\\
	\enumb
	\item  For $q = 2^2,$ the primitive polynomial giving the selector
	    0, 1, 4, 14, 16 is $I^3-I^2-I-\epsilon$ , with\\
	\hth$\epsilon^2+\epsilon +1 = 0.$
	\item  The correspondence between the selector notation and the
	    homogeneous coordinates are as follows, $i^*$ has the homogeneous
	    coordinates associated with $I^i.$\\
	$\begin{array}{clrr}
	\hth&     i&	I^i&		i*\\
	\cline{2-4}
	&     0& 1&	      20^*\\
	&     1& I&	      14^*\\
	&     2& I^2&	      0^*\\
	&     3& I^2+I+\epsilon&		 10^*\\
	&     4& I+\epsilon&	^2    19^*\\
	&     5& I^2+\epsilon&	^2I   4^*\\
	&     6& I^2+\epsilon&	^2I+1       18^*\\
	&     7& I^2+1&	   15^*\\
	&     8& I^2+\epsilon&		    3^*\\
	&     9& I^2+\epsilon^2I+\epsilon& 	5^*	\\
	&    10& I^2+\epsilon I+1&	 9^*\\
	&    11& I^2+\epsilon	^2&   13^*\\
	&    12& I^2+\epsilon I+\epsilon& 	11^*\\
	&    13& I^2+I+\epsilon^2&  6^*\\
	&    14& I+1&	     2^*\\
	&    15& I^2+I&	    1^*\\
	&    16& I+\epsilon&		     12^*\\
	&    17& I^2+\epsilon I&	   16^*\\
	&    18& I^2+\epsilon I+\epsilon^2& 17^*\\
	&    19& I^2+\epsilon^2I+\epsilon^2& 8^*\\
	&   20& I^2+I+1&	  7^* 
	\end{array}$

	    To obtain the last column, for row 9, $[1,\epsilon^2,\epsilon ] =
	(1,1,1) \times (1,\epsilon,0) = 20 \times 17 = 5*.$
	\item  The correspondence $i$ to $i^*$ is a polarity whose fixed points
	are on a line.  The matrix representation is obtained by using the image
	of 4 points.\\
	\hth	 0 = (0,0,1), $M(0) = 0^* = [1,0,0],$\\
	\hth	 1 = (0,1,0), $M(1) = 1^* = [1,1,0],$\\
	\hth	 2 = (1,0,0), $M(2) = 2^* = [0,1,1],$\\
	\hth	18 = $(1,\epsilon,\epsilon^2),$ $M(18) = 18^*
		= [1,\epsilon^2,1].$\\ 
	    The first 3 conditions give the polarity matrix as\\
	    The last condition gives $\beta \epsilon + \alpha\epsilon^2 = 1,$
	$\gamma + \beta \epsilon = \epsilon^2,$ $\gamma = 1.$  Hence 
	$\gamma  = 1,$ $\beta  = 1,$ $\alpha  = 1.$  Therefore\\
	\hth$M = \mat{0}{1}{1}{1}{1}{0}{1}{0}{0},$
		$M^{-1} = \mat{0}{0}{1}{0}{1}{1}{1}{1}{1}.$\\
	    Note that $M$ is real and could have been obtained from the reality
	    and non singularity conditions, giving directly
	$\alpha = \beta = \gamma  = 1.$\\
	    The polar of $(X_0,X_1,X_2)$ is $[X_1+X_2,X_0+X_1,X_0].$\\
	    The fixed points $(X_0,X_1,X_2)$ satisfy $X_1^2 = 0$ corresponding
	to $14^*.$ 
	\item  A point conic with no points on 14 is  1, 3, 4, 5,13,\\
	    the corresponding line conic is       15,19,10,16, 8.\\
	    Projecting from 1 and 3,     	   1, 3, 5,13, 4,\\ we get the
	    fundamental projectivity,		   8, 2,11, 0, 7 on $14^*.$ 

	\item To illustrate Pascal's Theorem, because there are only 5 points on
	    a conic, we need to use the degenerate case.  The conic through
	    0, 1, 2 and the conjugate points 9 and 18 is
	    The last condition gives $\beta \epsilon  + \alpha\epsilon^2 = 1,$
	$\gamma + \beta \epsilon = \epsilon^2,$ $\gamma = 1.$\\
	\hth    Hence $\gamma  = 1,$ $\beta  = 1,$ $\alpha  = 1.$  Therefore\\
	\hth$M = \mat{0}{1}{1}{1}{1}{0}{1}{0}{0},$
		$M^{-1} = \mat{0}{0}{1}{0}{1}{1}{1}{1}{1}.$\\
	    Note that $M$ is real and could have been obtained from the reality
	    and non singularity conditions, giving directly
	$\alpha = \beta = \gamma  = 1.$\\
	    The polar of $(X_0,X_1,X_2)$ is $[X_1+X_2,X_0+X_1,X_0].$\\
	    The fixed points $(X,X_1,X_2)$ satisfy $X_1^2 = 0$ corresponding to
	    $14^*.$ 
	\item A point conic with no points on 14 is  1, 3, 4, 5,13,\\
	\hth$\mat{0}{1}{1}{1}{0}{1}{1}{1}{0}$\\
	    The tangents at (0,0,1), (0,1,0), (1,0,0),
	$(1,\epsilon^2,\epsilon),$ $(1,\epsilon,\epsilon^2)$ 
	    are		    [1,1,0], [1,0,1], [0,1,1],
		$[1,\epsilon^2,\epsilon),$ $(1,\epsilon,\epsilon^2],$ 
	    or		     $1^*,$      $15^*,$     $2^*,$   $5^*,$    $17^*.$ 
	    On the other hand, using Pascal's Theorem, the tangent at 0 is
	    given by\\
	$  ( ( ( (0 \times 1) \times (9 \times 18) ) \times ( (18 \times 0)
	\times (1 \times 2) ) ) \times (2 \times 9) ) \times 0\\
	\hth= ( ( (0^* \times 7^*) \times (4^* \times 20^*) ) \times 12^*)
	 \times 0\\
	\hth= ( ( (14 \times 17) = 8^*) \times 12^* {\rm or} 13) \times 0 = 1^*.$ 
	\enume

	\sssec{Answer to}
	\vspace{-18pt}\hspace{100pt}{\bf \ref{sec-etruncdodeca2}.}\\
	For $q = 57,$
	choose the auto-correlates as point on a circle although 0 is on the
	circle draw as it is the center.  With the succession of points $X_i,$\\
	$\begin{array}{llrrrrrr}
	x_i = 0 \times X_i &36,& 1,&52,&43,& 3,&32,&13,\\
	X_i		   &   16,&35,&18,&50,&29,&26,&30,\\
	y_{i+1} = X_{i-1} \times X_{i+1}&      22,&42,& 8,&14,&10,&28,&44,\\
	y_{i+2} = X_{i-2} \times X_{i+2}&      34,& 2,&41,&17,&40,&20,&23,\\
	y_{i+3} = X_{i-3} \times X_{i+3}&       7,&31,& 6,&27,&54,&25,&39,\\
	y_{i+1} \times x_i&  21,&51,& 5,&46,&33,& 4,&45,\\
	y_{i+2} \times x_i&  24,&56,&48,&15,&49,&38,&47,\\
	y_{i+3} \times x_i&  53,&12,&37,& 9,&55,&11,&19.\\
	\end{array}$
	
	This gives all the points in the projective plane of order 7.
	We observe\\
	$\begin{array}{lllllll}
	16^*&21^*&24^*&53^*&22^*&34^*& 7^* \\
	36&36&36&36&36&36&36\\
	16&&&&35,30&18,26&50,29\\
	42,44&22& 8,28&&14,10\\
	41,20&&34&17,40&& 2,23\\
	27,54&31,39&& 7&&& 6,25\\
	&&46,33& 5, 4&21&&51,45\\
	&15,49&&56,47&48,38&24\\
	&37,11&12,19&&& 9,55&53\\
	\cline{1-7}
	35^*&51^*&56^*&12^*&42^*& 2^* 31^*\\
	 1& 1& 1& 1& 1& 1& 1\\
	35&&&&16,18&50,30&29,26\\
	22, 8&42&14,44&&10,28\\
	17,23&& 2&40,20&&34,41\\
	54,25& 7, 6&&31&&&27,39\\
	&&33, 4&46,45&51&&21, 5\\
	&49,38&&24,48&15,47&56\\
	& 9,19&53,37&&&55,11&12\\
	\cline{1-7}
	18^*& 5^* 48^*&37^*& 8^* 14^*& 6^*\\
	52&52&52&52&52&52&52\\
	18&&&&35,50&16,29&26,30\\
	42,14& 8&22,10&&28,44\\
	34,40&&41&20,23&& 2,17\\
	25,39&31,27&& 6&&& 7,54\\
	&& 4,45&21,33& 5&&51,46\\
	&38,47&&56,15&24,49&48\\
	&53,55&12, 9&&&11,19&37
	\end{array}$

	\sssec{Answer to}
	\vspace{-18pt}\hspace{100pt}{\bf \ref{sec-etruncdodeca3}.}\\
	For q = $2^3,$ \\
	$\begin{array}{lrrrrrrrrrcrrrrrrrrr}
	36:&		 0&37&38&40&44&52&18&27&68&\hti{2}&1^*& 3^*& 7^*&
	2^*& 4^*& 5^*\\
	36 \times  0 =  0^*:&   0& 1& 3& 7&15&31&36&54&63& & 0 & 0 & 0 & 
	1 & 3 &31\\
	36 \times 37 = 37^*:&  17&26&36&37&39&43&51&67&72& &72 &51 &67 & 
	72 &72 &26\\
	36 \times 38 = 38^*:&  16&25&35&36&38&42&50&66&71& &35 &71 &66 &
	71 &50 &71\\
	36 \times 40 = 40^*:&  14&23&33&34&36&40&48&64&69& &14 &33 &69 &
	34 &69 &69\\
	36 \times 44 = 44^*:&  10&19&29&30&32&36&44&60&65& &30 &60 &29 &
	29 &32 &10\\
	36 \times 52 = 52^*:&   2&11&21&22&24&28&36&52&57& & 2 &28 &24 &
	52 &11 & 2\\
	36 \times 18 = 18^*:&  13&18&36&45&55&56&58&62&70& &62 &70 &56 &
	13 &70 &58\\
	36 \times 27 = 27^*:&   4& 9&27&36&46&47&49&53&61& &53 & 4 &47 &
	61 &27 &49\\
	36 \times 68 = 68^*:&   5& 6& 8&12&20&36&41&59&68& & 6 &12 & 8 &
	 5 &59 &68\\
	\end{array}$

	Conic with no point on 36:  2, 4, 5, 6,13,28,31,46,63\\
	\hth	line conic:	   29,59,31, 9,18,43,28,35,64.\\
	Fundamental projectivity: from 2 and 5 on the conic, the points\\
	\hth	 2, 5, 6,31,13,28, 4,46,63 give the points on $36^*:$\\
	\hth	38, 0,68,27,52,37,40,18,44.

\chapter{FINITE INVOLUTIVE SYMPATHIC AND GALILEAN GEOMETRY}
{\tiny\footnotetext[1]{G30.TEX [MPAP], \today}}

\setcounter{section}{-1}
	\section{Introduction.}

	In part II, I have given a construction of a finite projective
	geometry associated to a prime $p.$  In it, there is no notion of
	parallelism, equality of segments or of angles, perpendicularity, etc .
	I have then obtained the well known finite affine geometry.  In it,
	we have the notion of parallel lines, equality of segments on
	a given line or on parallel lines, but we have no circles, no
	notion of equality on non parallel lines, no perpendicularity, etc .
	It is the purpose of Part III to construct a finite Euclidean geometry
	in which these notions as well as measure of angles and distances
	can be obtained.

	In the first step, which I will call {\em involutive geometry}, I choose
	an involution on the ideal line.  This involution either is elliptic,
	in which case it has no real fixed points or is hyperbolic, in which
	case it has 2 real fixed points.  The elliptic case resembles more
	the standard Euclidean geometry, while the hyperbolic case is easier to
	deal with, but the properties of both geometries go hand in hand.
	In it we define circles and perpendicularity.  A principle of
	compensation, which is not evident in the classical case, makes its
	appearance.  For instance, if we consider the lines through the center
	of a circle, half of them do not intersect the circle, but the other
	half do and then at two points.  As an other example, not all triangles
	have an inscribed circle, only roughly one in 4 has, but these have
	4 inscribed circles.  In the involutive geometry, I also define the
	equality of angles and the equality of segments.

	In the second step, I will introduce the {\em sympathic geometry}, in
	which we
	have the notion of measure of angle.  The algebraic development suggests
	a finite trigonometry.  In fact 2 such trigonometries are required for
	each prime, corresponding to the elliptic and to the hyperbolic case.
	The trigonometry for the elliptic case is obtained easily from the
	notion of primitive roots associated to $p.$  The trigonometry for
	the hyperbolic case, requires a generalization.

	In the last step, I introduce the notion of measure of distances and
	obtain the {\em finite Euclidean geometry}.


\setcounter{section}{0}
	\section{Finite involutive geometry.}
    \setcounter{subsection}{8}
	\ssec{Theorems in finite involutive Geometry, which do not 
		correspond to known theorems in Euclidean Geometry.}

	The Theorems in finite Euclidean Geometry fall also in several
	categories.  The first one, \ldots

	The theorems are a direct consequence of \ldots.\\
	The proof follows by assuming like in section \ldots that $m$ 
	corresponds to the line at infinity and \ldots.  The reference in 
	parenthesis is to the section <?> in Theorem \ldots.

	\sssec{Theorem.}
	\enumb
	\item Let $M1 \times H2$ meet $A1 \times A2$ in $C0$,  \ldots, then 
		the points $C0$, $C1$ and $C2$
		are on the same line $p.$
	\item Let $H1 M2$ meet $A1 A2$ in $D0$,  \ldots, then the points 
		$D0$, $D1$ and $D2$\\
		are on the same line $q.$
	\item The intersection $P$ of $p$ and $q$ is on the line $eul$ of Euler.
	\enume

	Proof:  Use AA1, 3.0, 3.1, with $H0 = \overline{M}_0,$ $M0 = M_0,$ 
	$C0 = C_0,$ 
		$D0$ = $C~_0,$ $p = p$ and $q = \overline{p}.$

	\ssec{The geometry of the triangle of degree 2.}

	\ldots   Involves problems of the second degree, bisectrices,
	inscribed circles for even triangles.

	\ssec{Some theorems involving circles.}
\setcounter{subsubsection}{-1}
	\sssec{Introduction.}
	It is not my intention to devlop here the extensive theory on
	circles for involutive geometry over arbitrary fields.
	I will simply give an example which illustrates how the problem
	can be approached effectively.

	\sssec{Definition.}
	Let $\theta$ be a defining circle and $m$, the ideal line,
	any circle $\gamma$ can be written as\\
	\hth$\gamma = \theta + (m) \bigx (r),$\\
	where $(r) = [r_0,r_1,r_2]$ is a given constant times the radical
	axis with $\theta$.
	The {\em 3 dimensional representation of the circle $\gamma$ } is
	defined by the point, with coordinates
	$r_0,r_1$ and $r_2$. I will write $({\bf r})_3 := (r_0,r_1,r_2)_3$
	for that representation. $\theta$ is represented by the origin.
	A degenerate circle $(m) \bigx (r)$ is represented by the direction
	of $r$.

	\sssec{Exercise.} What is the representation of tangent circles.

	\sssec{Lemma.}
	{\em If $\gamma_0$ and $\gamma_1$ are circles, represented by
	$({\bf r}_0)_3$ and $({\bf r}_1)_3$, then the family of circles
	through their intersections is represented by\\
	\hth$({\bf r}_0)_3 + k({\bf r}_1)_3$,\\
	with $k$ an arbitrary element in the field together with $\infty$,
	where $\infty$ represents $\gamma_1$.
	The addition is that of vectors in 3 dimensions and the multiplication
	by $k$ the scalar multiplication.}

	I will also denote the family by\\
	\hth$\gamma_0 + k \gamma_1$.\\
	One can also use the homogeneous representation,\\
	\hth$k_0\gamma_0 + k_1\gamma_1$.\\
	Strictly speaking, this is the representation used in the proofs,
	although I have used the non homogeneous representation to
	simplify the writing.

	\sssec{Theorem. [Bundle]}
	{\em Let $\gamma_j$, $j = 0$ to $3,$ be 4 circles, if there is a circle
	$\alpha$ which passes through the intersection of the circles
	$\gamma_0$ and $\gamma_1$, as well as the intersection of
	$\gamma_2$ and $\gamma_3$, then there is a circle $\beta$ passing
	through the intersections of $\gamma_0$ and $\gamma_2$, as well as
	those of $\gamma_1$ and $\gamma_3$.}

	This is the so called {\em bundle Theorem.}
	\footnote{see Dembosky, p. 256}

	Proof: If ${\bf r}_j)_3$ is the representation of $\gamma_j$. The family
	through the first 2 circles is represented by
	$({\bf r}_0)_3 + k({\bf r}_1)_3$ and that through the last 2 circles
	by\\
	\hth$({\bf r}_2)_3 + l({\bf r}_3)_3$, the hypothesis
	concerning the circle $\alpha$ implies\\
	\hth$({\bf r}_0)_3 + k({\bf r}_1)_3 = u(({\bf r}_2)_3 + l({\bf r}_3)_3$,
	which can be rewritten\\
	\hth$({\bf r}_0)_3 + u({\bf r}_2)_3 = -k(({\bf r}_2)_3 + 
	ul({\bf r}_3)_3)$,
	which gives the conclusion concerning the circle $\beta.$

	\ssec{The parabola, ellipse and hyperbola.}

	\sssec{Introduction.}
	The parabola, ellipse and hyperbola have already be
	defined in affine geometry.  Here we study their properties in
	involutive geometry.

	The parabola.

	The ellipse and hyperbola.

	If we assume that the isotropic points are $(\delta,1,0)$ and 
	$(-\delta,1,0),$
	where $\delta^2$ = d = N p, we will see that by an appropriate \ldots
	transformation, these can be reduced to\\
	\hth$	\frac{X0^2}{A} + \frac{X1^2}{B} = X2^2.$\\
	Recall also that $i^2 = -1$.

	\sssec{Definition.}
	The isotropic tangents are isotropic lines tangent to
	the conic.  The foci are the intersection of 2 isotropic tangents
	through 2 different isotropic points.

	\sssec{Theorem.}
	Given the conic\\
	\hth$	\frac{X0^2}{A} + \frac{X1^2}{B} = X2^2.$
	\enumb
	\item[D0.~~]$    C = A + Bd,$\\
	then
	\item[C0.~~]  The point polarity is\\
	\hth$\mat{B}{0}{0}{0}{A}{0}{0}{0}{-AB}$
	\item[C1.~~]The line polarity is\\
		$\mat{A}{0}{0}{0}{B}{0}{0}{0}{-1}$
	\item[C2.~~] The isotropic tangents through $(\delta,1,0)$ are\\
		$(1,-\delta ,\sqrt{C})$ and $(1,-\delta ,-sqrt{C})$
	\item[C3.~~] The foci are
	\item[C3.0.]$(-\sqrt{C},0,1),$ $(\sqrt{C},0,1),$
	\item[C3.1.]$(0,-\sqrt{C},\delta ),$ $(0,\sqrt{C},\delta ),$
	\item[C4.0.]$C = $R$ p \Rightarrow$ the foci C.3.0. are real, 
		the foci C.3.1. are not.
	\item[C4.1.]$C = $N$ p \Rightarrow$ the foci C.3.1. are real, 
		the foci C.3.0. are not.
	\enume

	\sssec{Theorem.}
	{\em Given the conic}
	\enumb
	\item[D0.~~]$\frac{X0^2}{A} + \frac{X1^2}{B} = X2^2.$
	\item[H1.0.]$A = $R$ p,$ $B = $R$ p,$
	\item[D1.0.]$   a = \sqrt{A},$ $b = \sqrt{B},$
	\item[H1.1.]$  A = $N$ p,$ $B = $N$ p,$
	\item[D1.1.]$  a = \sqrt{\frac{A}{d}},$ $b = \sqrt{\frac{B}{d}},$
	\item[H1.2.]$  A = $R$ p,$ $B = $N$ p,$
	\item[D1.2.]$  a = \sqrt{A},$ $b = \sqrt{\frac{B}{d}},$
	\item[H1.3.]$  A = $N$ p,$ $B = $R$ p,$
	\item[D1.3.]$  a = \sqrt{\frac{A}{d}},$ $b = \sqrt{B},$
		{\em then the conic takes the form}
	\item[C1.0.]$  \frac{x^2}{a^2} + \frac{y^2}{b^2} = 1,$
	\item[C1.1.]$  d\frac{x^2}{a^2} + d\frac{y^2}{b^2} = 1,$
	\item[C1.2.]$  \frac{x^2}{a^2} + d\frac{y^2}{b^2} = 1,$
	\item[C1.3.]$  d\frac{x^2}{a^2} + \frac{y^2}{b^2} = 1,$
	\enume

	\sssec{Theorem.}
	\enumb
	\item[H0.0.]$  p = -1 mod 4 and AB = R p,$\\
	or
	\item[H0.1.]$  p = 1 mod 4 and AB = N p,$\\
	{\em then}
	\item[C0.~~] {\em the conic is an ellipse,}
	\enume

	\sssec{Theorem.}
	\enumb
	\item[H0.0.]{\em$p \equiv 1 \pmod{4}$ and $AB$ R $p,$\\
	or}
	\item[H0.1.]{\em$  p \equiv -1 \pmod{4}$ and $AB$ N $p,$\\
	then}
	\item[C0.~~]{\em the conic is a hyperbola.\\
	The ideal points on it are}
	\item[C1.0.]$(\frac{a}{b} i,1,0), (-\frac{a}{b} i,1,0),$
	\item[C1.1.]$(\frac{a}{b}, 1, 0), (-\frac{a}{b}, 1, 0).$
	\enume

	\ssec{Cartesian coordinates in involutive Geometry.}
	\vspace{-18pt}\hspace{320pt}\footnote{5.7.83}\\[8pt]
	\sssec{Introduction.}

	\sssec{Notation.}
	A pair of reals between parenthesis will denote the
	Cartesian coordinates of a point.  We cannot choose a pair of reals
	between brackets to denote the $x$ and $y$ intercept of a line in the
	Cartesian plane, because we have then no way to represents lines through
	the origin.  We will therefore use triplets, with the last non zero
	coordinate normalized to 1.

	\sssec{Theorem.}
	If we choose as $x$ axis the line [0,1,0] and as $y$ axis [1,0,0],
	then we have the correspondence:\\
	\hth$	C (i,j,1) = (i,j),$\\
	\hth$	C [i,j,k] = [\frac{i}{k},\frac{j}{k},1],$ $k\neq 0,$\\
	\hth$	C [i,j,0] = [\frac{i}{j},1,0],$ $j\neq 0,$\\
	\hth$	C [i,0,0] = [1,0,0].$

	\sssec{Theorem.}
	Given a triangle whose vertices have the Cartesian coordinates\\
	\hth$	(0,a),$ $(b,0),$ $(c,0),$ $a \neq  0,$ $b \neq  c.$
	\enumb
	\item The point whose barycentric coordinates are
	\hth$	(q0,q1,q2),$ with $q0+q1+q2 \neq  0,$\\
	\hth    corresponds to the point whose Cartesian coordinates are\\
	\hth$	( \frac{bq1+cq2}{q0+q1+q2},$ $\frac{aq0}{q0+q1+q2} ).$\\
REDO 1. IN view of the preceding theorem
	\item The line, distinct from the ideal line, whose barycentric\\
	\hth    coordinates are\\
	\hth$	[l0,l1,l2]$\\
	\hth    corresponds to the line whose intercepts are\\
	\hth$	\frac{bl2-cl1}{l2-l1},$ $\frac{bl2-cl1}{(c-b)l0-cl1+bl2) },$\\
	if $bl2-cl1 = 0$\\
	\hth    and $((c-b)l0-cl1+bl2) \neq  0,$ it corresponds to\\
	\hth$	0,\frac{l2-l1}{(c-b)l0-cl1+bl2)},$\\
	\hth    and $((c-b)l0-cl1+bl2) = 0,$ it corresponds to\\
	\hth$	1,0,$
	\item The values of the coordinates of the orthocenter are\\
	\hth$	m0 = bc(b-c),$ $m1 = c(a^2+bc),$ $m2 = -b(a^2+bc).$
	\enume

	\sssec{Definition.}
	\vspace{-18pt}\hspace{94pt}\footnote{13.11.83}\\[8pt]
	The following mapping associates to the non ideal points
	in the finite Euclidean plane associated to $p$, points in the classical
	Euclidean plane.\\
	\hth$	T(i,j) = (i+kp,j+lp),$ where $k$ and $l$ are any integers.\\

	\sssec{Theorem.}
	Let $d = i1j2-i2j1,$ then
	\hth$	(i1,j1) x (i2,j2) = [\frac{j1-j2}{d},$
		$\frac{i2-i1}{d},1],$ $d\neq 0,$\\
	\hth$	(i1,j1) x (i2,j2) = [\frac{j1-j2}{i2-i1},1,0],$
		$d = 0,$ $i2\neq i1,$\\
	\hth$	(i1,j1) x (i2,j2) = [1,0,0],$ $d = 0,$ $i2=i1,$ $j2\neq j1.$

	For the following see ..[1,135]/cartes

	\sssec{Example.}
	For $p$ = 13, let the circles be\\
	Cr:\hth$     x^2 + y^2 = r^2.$\\
	The points on the circles are\\
	\enumb
	\item[C1:]$ (1,0), (-1,0), (0,1), (0,-1), (6,2), (-6,2), (6,-2),
		 (-6,-2), (2,6), (-2,6), (2,-6), (-2,-6),$
	\item[C2:]$ (2,0), (-2,0), (0,2), (0,-2), (4,1), (-4,1), (4,-1),
		 (-4,-1), (1,4), (-1,4), (1,-4), (-1,-4),$
	\item[C3:]$ (3,0), (-3,0), (0,3), (0,-3), (6,5), (-6,5), (6,-5),
		 (-6,-5), (5,6), (-5,6), (5,-6), (-5,-6),$
	\item[C4:]$ (4,0), (-4,0), (0,4), (0,-4), (5,2), (-5,2), (5,-2),
		 (-5,-2), (2,5), (-2,5), (2,-5), (-2,-5),$
	\item[C5:]$ (5,0), (-5,0), (0,5), (0,-5), (4,3), (-4,3), (4,-3),
		 (-4,-3), (3,4), (-3,4), (3,-4), (-3,-4),$
	\item[C6:]$ (6,0), (-6,0), (0,6), (0,-6), (3,1), (-3,1), (3,-1),
		 (-3,-1), (1,3), (-1,3), (1,-3), (-1,-3),$
	\enume

	The isotropic lines through the origin contain the points:\\
	i0:\hti{5}$ (0,0), (1,-5), (2,3), (3,-2), (4,6), (5,1), (6,-4),$\\
	\hth$	   (-1,5), (-2,-3), (-3,2), (-4,-6), (-5,-1), (-6,4),$\\
	i1:\hti{5}$ (0,0), (1,5), (2,-3), (3,2), (4,-6), (5,-1), (6,4),$\\
	\hth$	   (-1,-5), (-2,3), (-3,-2), (-4,6), (-5,1), (-6,-4),$\\
	If we join the origin to the points (1,k) and (1,l) we obtain
	perpendicular directions, with k,l = 0,oo; 1,-1; 2,6; 3,4; -2,-6;
		 -3,-4.

	For $p$ = 13, let the circles be\\
	Cr:\hti{5}$     x^2 - 6 xy + y^2 = r^2.$\\
	The points on the circles are
	\enumb
	\item[C1:]$ (0,1), (1,0), (1,6), (4,4), (4,-6), (6,1), (6,-4),$
	\item[C2:]$ (0,2), (1,-2), (1,-5), (2,0), (2,-1), (5,5), (5,-1),$
	\item[C3:]$ (0,3), (1,1), (1,5), (3,0), (3,5), (5,1), (5,3),$
	\item[C4:]$ (0,4), (2,3), (2,-4), (3,2), (3,3), (4,0), (4,-2),$
	\item[C5:]$ (0,5), (4,5), (4,6), (5,0), (5,4), (6,4), (6,6),$
	\item[C6:]$ (0,6), (2,2), (2,-3), (3,-2), (3,-6), (6,0), (6,-3),$
	\enume
	as well as the points symmetric with respect to the origin.\\
	If we join the origin to the points (1,k) and (1,l) we obtain
	perpendicular directions, with k,l = 0,-4; 1,-1; 2,-5; 3,oo; 4,-2,
	5,-6; 6,-3.

	For $p$ = 11, let the circles be\\
	Cr:\hti{5ex}$     x^2 - 4xy + y^2 = r^2.$\\
	The points on the circles are
	\enumb
	\item[C1:]$ (0,1), (1,4), (4,4),$
	\item[C2:]$ (0,2), (2,-3), (3,3),$
	\item[C3:]$ (0,3), (1,1), (1,3),$
	\item[C4:]$ (0,4), (4,5), (5,5),$
	\item[C5:]$ (0,5), (2,2), (2,-5),$
	\enume

	as well as the points symmetric with respect to the diagonals,
	$(i,j)$ here means $(i,j),$ $(j,i),$ $(-i,-j),$ $(-j,-i).$\\
	The isotropic points are\\
	l0:\hti{5ex}$ (0,0),$ $(1,3),$ $(1,4),$ $(2,-3),$ $(2,-5),$ $(3,-2),$ 
		$(3,1),$\\
	\hth$	$(4,1),$ $(4,5),$ $(5,-2),$ $(5,4),$\\
	l1:\hti{5ex}$ (0,0),$ $(1,-3),$ $(1,-4),$ $(2,3),$ $(2,5),$ $(3,2),$ 
		$(3,-1),$\\
	\hth$	$(4,-1),$ $(4,-5),$ $(5,2),$ $(5,-4),$\\
	If we join the origin to the points $(1,k)$ and $(1,l)$ we obtain
	perpendicular directions, with $k,l$ = 0,-5; 1,-1; 2,$\infty$; 3,5; 
		4, -2.

	For $p$ = 11, let the circles be\\
	Cr:\hti{5ex}$     x^2 + y^2 = r^2.$\\
	The points on the circles are
	\enumb
	\item[C1:] (0,1), (3,5),
	\item[C2:] (0,2), (1,5),
	\item[C3:] (0,3), (2,4),
	\item[C4:] (0,4), (1,2),
	\item[C5:] (0,5), (3,4),
	\enume

	as well as the points symmetric with respect to the 2 axis and the
	diagonals.\\
	$(i,j)$ here means $(i,j),$ $(i,-j),$ $(j,i),$ $(j,-i),$ $(-i,-j),$ 
		$(-i,j),$ $(-j,-i),$ $(-j,i).$\\
	If we join the origin to the points $(1,k)$ and $(1,l)$ we obtain
	perpendicular directions, with $k,l = 0,\infty$; 1,-1; 2,5; 3,-4; 
		4,-3; -2,-5.

	\ssec{Correspondence between circles in finite and classical 
		Euclidean geometry.}

	\sssec{Introduction.}

	\sssec{Theorem.}
	To the point $(x,y),$ in classical geometry, on a circle
	centered at the origin and of radius $r,$ corresponds, if $r$ is not
	congruent to 0 modulo $p$, the point $(x/r mod p, y/r mod p)$ on a 
	circle of
	radius 1 in the finite geometry associated to $p$.\\
	Vice-versa, given a point $P = (x,y)$ on a circle of radius 1 in the
	finite geometry associated to $p$, we can always find a point on a
	circle in the classical geometry which is one of the representatives of
	$P$, in the mapping given in \ldots.

	The proof is left to the reader.  The first part is trivial, the second
	part is not trivial.  See also [135]FINPYT.BAS

	\sssec{Example.}
	For $p$ = 13,\\
	\hth	(2,6) for $r$ = 1 is associated to (15,20) for $r$ = 25.\\
	For $p$ = 29,\\
	\hth	(5,11) for $r$ = 1 is associated to (24,18) for $r$ = 30.\\
	\hth	(8,13) for $r$ = 1 is associated to (108,45) for $r$ = 117.\\
	\hth	(6,9) for $r$ = 1 is associated to (180,96) for $r$ = 204.

	\sssec{Theorem.}
	There exist a circle of radius $u$ in $R^2$ which contains all the
	representatives of a circle in $Z_p^2.$\\
	Indeed, for the radius 1, for instance, if one of the representatives
	is on the circle $x^2 + y^2 = r_1^2$ and if $s_1 = 1/r_1,$ then\\
	\hth$	u = r_i(s_i+k_ip),$ for all i,\\
	by finite induction, if\\
	\hth$	u = r_1(s_1 + p k_1) = r_2(s_2 + p k_2),$\\
	then\\
	\hth$	r_1 k_1 - r_2 k_2 = (r_2s_2 - r_1s_1)/p,$\\
	this gives\\
	\hth$	k_1 = a_1 + r_2 k'_2$ and with $s'_2 = (s_1 + a_1 p)/r_2,$\\
	\hth$	u = r_1 r_2 (s'_2 + p k'_2),  \ldots  .$

	\sssec{Example.}
	For $p$ = 29, for $r$ = 1, we start with\\
	\hth	point in $Z_{29}^2$        $r_i$    $s_i$\\
	\hth	5,11                    5       6\\
	\hth	8,13                    13      9\\
	\hth	6,9                     25      7\\
	then\\
	\hth$	5 k_1 - 13 k_2 = 3,$ $a_1 = -2,$ $s'_2 = -4,$\\
	\hth$	u = 5.13 (-4 + 29 k'_2),$\\
	\hth$	13 k'_2 - 5 k_3 = 3,$ $a_2 = 1,$ $s'_3 = 5,$\\
	\hth$	u = 5.13.5 (5 + 29 k'_3),$\\
	hence the suitable circle in $R^2$ with smallest radii has radius 
		$u$ = 1625
	and contains the points\\
	\hth	in $R^2$          in $Z_{29}^2$\\
	\hth	-1300,-975      5,11;\\
	\hth	-1500,-625      8,13;\\
	\hth	-1560,-455      6,9.\\

	\ssec{Answers to problems.}

	\sssec{Answer to 1.13.1.}
	For the second part.\\
	The problem can be restated successively as follows, given a
	solution of\\
	0.      $x^2 + y^2 = z^2,$ 
	there exist $i,j,k$ such that $(x+ip)^2 + (y+jp)^2 = (z+kp)^2,$ 
	or there exist $u$ and $v$ such that\\
	\hth$	u^2 - v^2 = x,$ $2uv = y,$ $u^2 + v^2 = z,$\\
	eliminating v from the first 2 equations and using 0., gives\\
	\hth$	u^2 = \frac{r+x}{2},$ $v^2 = \frac{r-x}{2},$\\
	$\frac{r+x}{2}$ need not be a quadratic residue, therefore we use 
		instead\\
	\hth$	u^2 = b\frac{r+x}{2},$ $v^2 = b \frac{r-x}{2},$ $c = 1/b,$\\
	this gives $u$ and $v$,\\
	\hth$	x+ip = (u^2-v^2)c,$ $y+ip = 2uv c,$ $z+ip = (u^2+v^2) c.$\\
	A more careful discussion will show that signs may have to be
	changed and the role of $x$ and $y$ interchanged.\\
	For instance, for $p$ = 13 and $2^2 + 6^2 = 1^2,$ 
	$b$ = -2, $c$ = 6, $u^2 = (-5)(-2) = 6^2,$ $v^2 = (6)(-2) = 1^2,$ 
	hence $x+ip = 35,$ $-(y+jp) = 72,$ $z+kp = 222.$\\
	For $p$ = 17, and $x$ = 4, $y$ = 6, $z$ = 1,
	$b$ = 3, $c$ = 6, $u^2 = (6)(3) = 1^2,$ $v^2$ = (12)(3) = $2^2,$ 
	hence after interchange of $x$ and $y$,
	$x+ip = 72,$ $y+jp = 210,$ $z+kp = 222.$\\
	For $p$ = 19 and $3^{2+7^2} = 1^2,$ 
	$b$ = 2, $c$ = 10, $u^2 = (2)(2) = 2^2,$ $v^2 = (-1)(2) = 6^2,$ 
	hence $-(x+ip) = 320,$ $-(y+jp) = 240,$ $z+kp = 400.$



	(AFTER INVOLUTIVE GEOMETRY)

	\sssec{Comment.}
	For the following theorem, I will not give a linear
	construction, although one could be given.  The theorem is a
	generalization of the Theorem of Miquel and can be further
	generalized in the context of Gaussian geometry.

	\sssec{Notation.}
	If $u$ and $v$ are 2 lines and $\xi$ is a conic,\\
	\hth$\xi - u\bigx v = 0$,\\
	is equivalent to\\
	\hth$\xi(X) - (u\cdot X)\: (v\cdot X) = 0$.\\
	This should be moved before the definition of circles.

	\sssec{Theorem.}
	{\em The radical axis of the 2 circles\\
	\hth$\mu_j := \theta - m\bigx u_j$, $j = 0,1$\\
	is $u_1 - u_0$.}

	Indeed, $\mu_1(X) - \mu_0(X) = - (m\cdot X)\: ((u_1-u_0)\cdot X) = 0$
	therefore\\
	\hth$\mu_1 = \mu_0 - m\bigx (u_1-u_0).$

	\sssec{Theorem.}
	{\em The radical axis of each pair of 3 circles are concurrent}.

	Proof: Let the 3 circles be\\
	\hth$\mu_i := \theta - m\bigx u_i$,\\
	$u_i$ are the radical axis of these circles with $\theta$.\\
	the 3 radical axis are $u_2-u_1$, $u_0-u_2$, $u_1-u_0$, but
	$u_2-u_1 = (u_0-u_2) + (u_1-u_0)$, therefore one of the axis passes
	through the intersections of the other 2.

	\sssec{Theorem. [Miquel]}
	H0.\hti{6}$N_i \cdot a_i = 0,$ \\
	H1.\hti{6}$N_i \cdot m \neq 0,$ \\
	D0.\hti{6}$\mu_i := circle(A_i,N_{i+1},N_{i-1}),$\\
	D1.\hti{6}$Miquel := (\mu_1 \times \mu_2) - N_0,$\\
	{\em then}\\
	C0.\hti{6}$Miquel \times \mu_0 = 0.$

	The nomenclature.\\
	N0.\hti{6}Miquel is called the {\em point of Miquel associated to}
	$N_i.$ 

	Proof.  Let\\
	\hth$	N_0 = (0,1,q_0),$ $N_1 = (q_1,0,1),$ $N_2 = (1,q_2,0).$\\
	H1. implies that $1+q_i \neq 0.$\\
	If the equation of $\mu_0$ a circle is\\
	\hth$	\theta  - m \bigx u_0 = 0,$ with
	$u_0 = [u_{0,0},u_{0,1},u_{0,2}],$\\
	$m = [1,1,1]$\\
	and\\
	$\theta = m'_0X_1X_2 + m'_1X_2X_0 + m'_2X_0X_1$\\$m = [1,1,1]$.
	with $m'_0 = m0(m1+m2),$  \ldots).\\
	D0. implies, for $i = 0,$\\
	\hth$	u_0 = [0,\frac{m'_2}{1+q_2},\frac{q_1m'_1}{1+q_1}].$\\
	Let $Miquel = (X_0,X_1,X_2)$.  If the 3 circles have a point in common,
	it is on the intersection of the 3 radical axis $u_i$, it is therefore
	necessary that\\
	\hth$	(u_i-u_{i+1})\cdot Miquel = 0,$ $i = 0, 1, 2,$\\
	therefore\\
	\hth$	Miquel = (u_1-u_2) \times (u_2-u_0),$\\
	this gives after simplification,\\
	\hth$	Miquel = (\frac{m'_0}{1+q_0}(\frac{-q_0m'_0}{1+q_0}
		+\frac{q_0q_1m'_1}{1+q_1}+\frac{m'_2}{1+q_2},  \ldots  ).$\\
	It remains to verify that $Miquel$ belongs to $\mu_i.$ \\
	First, $u_i\cdot Miquel = m'_0m'_1m'_2\frac{1+q_0q_1q_2}
	{1+q_0}(1+q_1)(1+q_2)),$\\
	second, $m\cdot Miquel = -\frac{q_0m'_0m'_0}{(1+q_0)^{2}} + \ldots$ 
	\hth$		+ m'_1m'_2\frac{1+q_1q_2}{(1+q_1)(1+q_2)} +  \ldots $.\\
	It is straigthforward to verify that the product of these two
	expressions is precisely $m'_0X_1X_2 + m'_1X_2X_0 + m'_2X_0X_1.$

	\sssec{Theorem. [Miquel]}
	H0.\hti{6}$n := N_1 \times N_2,\: n\cdot N_0 = 0.$\\
	D0.\hti{6}$n = [n_0,n_1,n_2],$\\
	{\em then}\\
	C0.\hti{6}$Miquel\cdot \theta  = 0.$\\
	C1.\hti{6}$q_0 = -\frac{n_1}{n_2},$ $q_1 = -\frac{n_2}{n_0},$
	$q_2 = -\frac{n_0}{n_1}.$\\
	C2.\hti{6}$Miqnel = (\frac{n_1n_2m'_0}{n_1-n_2},
		\frac{n_2n_0m'_1}{n_2-n_0},
		\frac{n_0n_1m'_2}{n_0-n_1}).$

	The condition that the points $N_i$ be collinear is precisely
	$1 + q_0q_1q_2 = 0,$ but in this case $\theta  = 0$ as follows from the
	expression $u_i\cdot Miquel.$  It is straightforward to verify C1 and
	C2.

	\sssec{Corollary.}
	0.\hti{7}{\em The circles $\mu f_i$ circumscribed to 
	$A_i,\overline{M}A_{i+1},\overline{M}A_{i-1}$
	have a point $\overline{F}ock$ in common.\\
	1.\hti{7}$\overline{F}ock$ is in the circumcircle $\theta$.}\\
	2.\hti{7}$\mu f_i = \theta
		+ m \bigx [0,\frac{m_2m_0(m_0+m_1)}{m_0-m_1},
		-\frac{m_2m_0(m_0+m_1)}{m_0-m_1}]$.\\
	3.\hti{7}$\overline{F}ock =
		(m_0(m_1+m_2)(m_2-m_0)(m_0-m_1),
		 m_1(m_2+m_0)(m_0-m_1)(m_1-m_2),
		 m_2(m_0+m_1)(m_1-m_2)(m_2-m_0))$.
	\footnote{22.12.88}

	This is the special case when $n$ is the orthic line $\overline{m}$ =
	[m$_0$,m$_1$,m$_2$].\\
	The point $\overline{F}ock$ had been constructed before (D38.9) and
	proven to be on $\theta$ (C38.4).

	\sssec{Theorem. [Miquel]}
	D0.\hti{6}$N_{i,j} := midpoint(A_i,N_j),$\\
	D1.\hti{6}$n_{i,j} := mediatrix(A_i,N_j),$\\
	D2.\hti{6}$C_i := n_{i,i+1} \times n_{i,i-1},$\\
	D3.\hti{6}$\phi  := circle(C_0,C_1,C_2),$\\
	{\em then}\\
	C0.\hti{6}$O\cdot \phi  = 0.$

	Proof:\\
	P0.\hti{6}$N_{0,1} = (1+2q_1,0,1),   N_{0,2} = (2+q_2,q_2,0).$\\
	P1.\hti{6}$n_{0,1} = [m2+m0,m0-(1+2q_1)m2,-(1+2q_1)(m2+m0)],$\\
	\hth$	n_{0,2} = [-q_2(m0+m1),(2+q_2)(m0+m1),(2+q_2)m1-q_2m0].$\\
	P2.\hti{6}$C_0 = (m0(1+2q_1+q_1q_2)m0+(2+q_2)(1+q_1)m1
			+(1+2q_1)(1+q_2)m2),$\\
	\hth$		(m2+m0)((1+q_1)m0+(q_1q_2-1)m1),$\\
	\hth$		(m0+m1)((1-q_1q_2)m2+(1+q_2)m0).$\\
	P3.\hti{6}$\phi :  \ldots $

	\sssec{Problem.}
	The following question suggests itself.  Let\\
	\hth$	\nu := circle(N_0,N_1,N_2).$\\
	What relation exists between all circles $\nu$  having the same point of
	Miquel?\\
	Same question in the case for which the point of Miquel is on
	$\theta$.\\
	Theorem \ldots states that all the circles are lines and \ldots that
	one of these
	lines is that of Simson and Wallace.  Again what is the relation between
	these lines?

	\sssec{Theorem. [Simson and Wallace]}
	H0.\hti{6}$X\cdot \theta  = 0,$\\
	D0.\hti{6}$n_i := X \times Im~_i,$\\
	D1.\hti{6}$N_i := n_i \times a_i,$\\
	D2.\hti{6}$n := N_1 \times N_2,$\\
	{\em then}\\
	C0.\hti{6}$N_0\cdot n = 0 (*).$\\
	C1.\hti{6}$(W \times N_i) \perp a_i.?$

	Proof:\\
	P0.\hti{6}$n_0 = [m1X_2-m2X_1,(m1+m2)X_2+m2X_0,-m1X_0-(m1+m2)X_1].$\\
	P1.\hti{6}$N_0 = (0,m1X_0+(m1+m2)X_1,m2X_0+(m1+m2)X_2).$\\
	P2.\hti{6}$n = [X_1X_2(-m0X_0+(m1+m2)(X_1+X_2)),\\
		X_2X_0(-m1X_1+(m2+m0)(X_2+X_0)),$\\
	\hth$		X_0X_1(-m2X_2+(m0+m1)(X_0+X_1))].$\\
	To obtain the last expression we use in each coordinate the relation H0,
	\\$m0(m1+m2)X_1X_2 + m1(m2+m0)X_2X_0 + m2(m0+m1)X_0X_1 = 0.$

	MAY WANT TO REFER HERE TO THE FOLLOWING BUT MOVE IT AS APPLICATION
	OF PARABOLAS.

	\sssec{Theorem.}
	{\em The set of lines having the same point $X$ of Miquel are on
	a line parabola}\footnote{30.12.82}:\\
	C0.\hti{6}$    mup^{-1}(X):$\\
	\hth$X_0 u0(u1-u2)/m'_0 = X_1 u1(u2-u0)/m'_1 = X_2 u2(u0-u1)/m'_2.$\\
	\hth$\mu p(X): (X_1X_2m'_0 U_0)^2 + (X_2X_0m'_1 U_1)^2 
		+ (X_0X_1m'_2 U_2)^2$\\
	\hth$- 2X_0X_1X_2(X_0m'_1m'_2 U_1U_2 + X_1m'_2m'_0 U_2U_0 
		+ X_2m'_0m'_1 U_0U_1).$\\
	C1.\hti{6}$a_i \cdot \mu p(X) = 0.$\\
	C2.\hti{6}{\em The line of Simson and Wallace is the tangent at the
	vertex.}\\
	C3.\hti{6}{\em The point of Miquel is its focus.}

	Proof.  C2 of Theorem \ldots gives C0.



\setcounter{subsection}{8}
	\ssec{The conic of Kiepert.}
\setcounter{subsubsection}{-1}
	\sssec{Introduction.}
	The conic of Kiepert has been constructed in
	5.4.1.D3.8.\footnote{13.1.83}.\\
	Kiepert showed that, in the classical case, if $V_i$ is a
	point on the mediatrix $mf_i$ such that\\
	\hth$angle(A_1,A_2,V_0) = angle(A_2,A_0,V_1 = angle( A_0,A_1,V_2),$\\
	then $v_i := A_i\times V_i$ have a point $V$ in common which is on a
	hyperbola, now known as the hyperbola of Kiepert.  After proving this
	Theorem in the finite case, I will consider several special cases of
	interest, which can be obtained either by a linear or by a second degree
	construction.  In the latter case, if the angle is $\frac{\pi}{4},$
	the point is
	called the point of Vectem, to which is associated a special chapter of
	the classical theory of the geometry of the triangle.  The cases when
	the angle is $\frac{\pi}{3}$ and $\frac{\pi}{6}$ are also discussed
	and a new property is obtained.

	\sssec{Theorem.}
	{\em Let}\footnote{15.9.86}
	\begin{enumerate}
	\setcounter{enumi}{-1}
	\item[H0.0.]$X \cdot \theta  = 0.$
	\item[G0.0.]$X = (X_0,X_1,X_2).$
	\item[D1.0.]$x1 := A_1\times X,$
	\item[P1.0.]$x1 = [X_2,0,-X_0],$
	\item[D1.1.]$V_0 := x1\times mf_0,$
	\item[P1.1.]$V_0 = ((m1+m2)X_0,(m1+m2)X_2-(m1-m2)X_0,(m1+m2)X_0).$
	\item[D1.2.]$v_0 := A_0\times V_0,$
	\item[P1.2.]$v_0 = [0,(m1+m2)X_2,(m1-m2)X_0-(m1+m2)X_2].$
	\item[D1.3.]$x2 := A_2\times V_0,$
	\item[P1.3.]$x2 = [(m1+m2)X_2-(m1-m2)X_0,-(m1+m2)X_0,0].$
	\item[D1.4.]$x3 := Ma_0\times X,$
	\item[P1.3.]$x3 = [-X_1-X_2,X_0,X_0].$
	\item[D1.5.]$Y = x2\times x3.$
	\item[P1.4.]$Y = ((m1+m2)X_0,(m1+m2)X_2-(m1-m2)X_0,\\
		\hti{12} (m1-m2)X_0+(m1+m2)X_1)$

	\noindent\item[D2.0.]$x4 := A_0\times X,$
	\item[P2.0.]$x4 = [0,X_2,-X_1].$
	\item[D2.1.]$X1 := x4\times m,$
	\item[P2.1.]$X1 = (X_1+X_2,-X_1,-X_2).$
	\item[D2.2.]$x5 := A_2\times X1,$
	\item[P2.2.]$x5 = [X_1,X_1+X_2,0].$
	\item[D2.3.]$V_1 := mf_1\times x5,$
	\item[P2.3.]$V_1 = ((m2+m0)(X_1+X_2),(m2+m0)X_1,2m2X_1+(m2+m0)X_2).$
	\item[D2.4.]$v_1 := A_1\times v_1,$
	\item[P2.4.]$v_1 = [2m2X_1+(m2+m0)X_2),0,-(m2+m0)(X_1+X_2)].$

	\noindent\item[D3.0.]$y4 := A_0\times Y,$
	\item[P3.0.]$y4 = [0,Y_2,-Y_1].$
	\item[D3.1.]$Y1 := y4\times m,$
	\item[P3.1.]$Y1 = (Y_1+Y_2,-Y_1,-Y_2).$
	\item[D3.2.]$y5 := A_1\times Y1,$
	\item[P3.2.]$y5 = [Y_2,0,Y_1+Y_2].$
	\item[D3.3.]$V_2 := mf_1\times y5,$
	\item[P3.3.]$V_2 = ((m0+m1)(Y_1+Y_2),(m0+m1)Y_1+2m1Y_2,-(m0+m1)Y_2).$
	\item[D3.4.]$v_2 := A_2\times V_2,$
	\item[P3.5.]$v2 := [(m0+m1)Y_1+2m1Y_2,-(m0+m1)(Y_1+Y_2),0].$

	\noindent\item[D4.0.]$V = v_0\times v_1,$
	\item[P4.0.]$V = (u(X_1+X_2),\\
	\hth	(m2+m0)((m0-m1)(m1-m2)X_0-2m1(m1+m2)X_1+uX_2,
	\hth	2m2(m0+m1)(m1+m2)X_1+uX_2),$\\
	where\\
		$u := (m1+m2)(m2+m0)(m0+m1).$\\
	{\em then}
	\item[C0.0.]$V \cdot v_2 = 0.$
	\item[C0.1.]$V \cdot \kappa iepert = 0.$
	\end{enumerate}

	The construction is based on\\
	$angle(X,A_1,A_2) = angle(A_1,A_2,Y) = angle(X,A_0,A_2)
	= angle(A_0,A_2,V_1),$ \\
	implying the parallelism of $A_0\times X$ and $A_2\times V_1$ and
	 symmetrically for $V_2.$ \\
	For P4.0., after replacing $Y_0,$ $Y_1$ and $Y_2$ by their values from
	P1.4.,
	the equation for $\theta$  is used to express $X_0X_1$ in terms of
	$X_2X_0$ and $X_1X_2.$

	\sssec{Exercise.}
	To complete the proof of x.x.1., the 2 special case $X = A_0$
	and $X = A_1$ should be considered.  This is left as an exercise.\\
	In the first case  $x4$ should be replaced by the tangent $ta_0$ at
	$A,$ in
	the second case $x1$ should be replaced by the tangent $ta_1$ at $A_1.$

	\sssec{Exercise.}
	Proceed in the inverse order and construct $X$ from
	$V.$  Prove
	that if $V$ is on $\kappa iepert$ then $X$ is on $\theta$ .

	\sssec{Exercise.}
	Study the projectivity which associates to $(X_0,X_1,X_2),$
	the point $(V_0,V_1,V_2),$ as given by P4.0. without assuming that
	$(X_0,X_1,X_2)$	is on $\theta$.  Determine 4 points and their images
	and construct
	any of these points if they have not been constructed in this book.

	The following are special cases.\\
	\hth$X = A_2,$ $\alpha  = 0$ gives $V = M.$\\
	\hth$\sigma  = \frac{\pi}{2}$ gives $V = \overline{M}.$\\
	\hth$\sigma  = \frac{\pi}{4}$ gives the point of Vectem (see below).\\
	$\sigma = \frac{\pi}{3}$ gives the equilateral point (see below)
	$\sigma = \frac{\pi}{6}$ gives the hexagonal point (see below)
	$\sigma = angle(A_{i-1},A_i,A_{i+1})$ gives $V = A_i.$

	Other angles give $V = Tar.$  (5.4.1.D16.3.), $V = Br0.$  (5.4.1.D15.3.)
	$V = \overline{B}r0.$  (5.4.1.D15.3.) and $V = En.$  (5.4.1.D21.10)

	\begin{enumerate}
	\setcounter{enumi}{-1}
	\item[D5.0.]$Ma\overline{m}_i := ma_{i+1}\times \overline{m}a_{i-1},$
	$\overline{M}a\overline{m}_i := \overline{m}a_{i+1}\times ma_{i-1},$
	\item[D5.1.]$Ae_i := a_i\times e,$
	\item[D5.2.]$mae_i := Ae_{i+1}\times Ma\overline{m}_{i-1},
	\overline{m}ae_i := Ae_{i+1}\times \overline{M}a\overline{m}_{i-1},$
	\item[D5.3.]$MMa_i := mae_i\times a_i,
	\overline{M}Ma_i := \overline{m}ae_i\times a_i,$
	\item[D5.4.]$  mm := MMa_1\times MMa_2,
	\overline{m}m := \overline{M}Ma_1\times \overline{M}Ma_2,$\\
	{\em then}
	\item[C5.0.]$n_i \cdot Kiepert1 = 0.$
	\item[C5.1.]$mm \cdot Kiepert1
		= \overline{m}m \cdot Kiepert1 = 0.$
	\item[C5.2.]$mm \cdot K = \overline{m}m \cdot K = 0.$
	\item[C5.3.]]$S$ is the center of Kiepert1,
	$ \overline{S}$ is the cocenter\footnote{15.1.83}.\\
	The nomenclature:

	Proof.
	\item[P5.0.]$Ma\overline{m}_0 = (m0,m1,m0),
	\overline{M}a\overline{m}_0 = (m0,m0,m2),$
	\item[P5.1.]$Ae_0 = (0,m0-m1,m0-m2),$
	\item[P5.2.]$mae_0 = [m2(m1-m2),m1(m2-m0),m2(m0-m1)],$\\
	\hth$	\overline{m}ae_0 = [m1(m2-m1),m2(m0-m2),m1(m1-m0)],$
	\item[P5.3.]$MMa_0 = (0,m2(m0-m1),m1(m0-m2)),$\\
	\hth$	\overline{M}Ma_0 = (0,m1(m1-m0),m2(m2-m0)),$
	\item[P5.4.]$mm = [m0(m1-m2),m1(m2-m0),m2(m0-m1)],$\\
	\hth$	\overline{m}m = [m1m2(m1-m2),m2m0(m2-m0),m0m1(m0-m1)],$
	\end{enumerate}

	The tangent at $Tar$ is\\
	\hth$    [m1m2q0^2(m2-m1),m2m0q1^2(m0-m2),m0m1q2^2(m1-m0)].$\\
	The tangent at $Br0$ is\\
	\hth$[m0^3(m1+m2)^2(m1-m2),m1^3(m2+m0)^2(m2-m0),
	\hth	m2^3(m0+m1)^2(m0-m1)].$\\
	The tangent at $\overline{B}r0$ is\\
	\hth$    [m1m2(m1+m2)^2(m1-m2),m2m0(m2+m0)^2(m2-m0),
	\hth	m0m1(m0+m1)^2(m0-m1)].$

	\sssec{Example.}
	With $p = 13,$ $A[] = (14,1,0),$ $M = (28),$
	$ \overline{M} = (44),$
	$Ma\overline{m}_0	= ( 41, 29, 70),$
	$\overline{M}a\overline{m}_0 = ( 31, 42,106),$
			$Ae_0	= (  4, 25, 79),$
	$mae_0	= [138, 81,145],$ $\overline{m}ae_0 = [151, 84,141],$
	$MMa_0	= (  9, 20,131),$ $\overline{M}Ma_0 = (  7, 19,144),$
	$mm = [146],$             $\overline{m}m = [136],$\\
	$V = (\frac{a}{sin(A-\alpha)},\frac{b}{sin(B-\alpha)},
	\frac{c}{sin(C-\alpha)}).$

	\sssec{Theorem.}
	If $V^{\sigma } = (V_0,V_1,V_2)$ is associated with the angle
	$\sigma$, then\\
	$V^{-\sigma} = ((m1-m2)(m2+m0)(m0+m1)V_0+2m0(m1+m2)(m2-m0)V_1\\
	\hti{12}	+2m0(m0-m1)(m1+m2)V_2,\\
	\hth       2m1(m1-m2)(m2+m0)V_0+(m2-m0)(m0+m1)(m1+m2)V_1\\
	\hti{12}	+2m1(m2+m0)(m0-m1)V_2)
			+2m2(m0+m1)(m1-m2)V_0,\\
	\hth	2m2(m2-m0)(m0+m1)V_1+(m0-m1)(m1+m2)(m2+m0)V_2).$

	Proof:\\
	$v_i := V\times A_i,$ 	$v_0 = [0,V_2,-V_1].$ \\
	$V_i := v_i\times mf_i,$ \\
	$V_0 = ((m1+m2)(V_2-V_1),(m1-m2)V_1,(m1-m2)V_2).$ \\
	$va_i = A_{i+1}\times V_{i-1},$	\\
	$va_0 = [(m0+m1)(V_0-V_1),0,(m0-m1)V_0].$ \\
	$v\overline{a}_i = A_{i-1}\times V_{i+1},$\\
	$v\overline{a}_0 = [(m2+m0)(V_2-V_0),(m2-m0)V_0,0].$ \\
	$Va_i	= m\times va_i,$ \\
	$Va_0 = ((m0-m1)V_0,2m1V_0-(m0+m1)V_1,-(m0+m1)(V_0-V_1)).$ \\
	$V\overline{a}_i = m\times v\overline{a}_i,$ \\
	$V\overline{a}_0 = ((m2-m0)V_0,(m2+m0)(V_2-V_0),-2m2V_0+(m2+m0)V_2)).$ \\
	$vb_i := Va_i\times A_i,$\\
	$Vb_0 = [0,(m0+m1)(V_0-V_1),2m1V_0-(m0+m1)V_1].$ \\
	$v\overline{b}_i := V\overline{a}_i\times A_i,$ \\
	$v\overline{b}_0 = [0,2m2V_0-(m2+m0)V_2,(m2+m0)(V_2-V_0)].$ \\
	$V^{-\sigma}_i := vb_{i+1}\times v\overline{b}_{i-1},$\\
	$V^{-\sigma}_0 = ((m1+m2)(V_2-V_1),2m1V_2-(m1+m2)V_1,-2m2V_1
		+(m1+m2)V_2).$ \\
	$v^{-\sigma}_i := V^{-\sigma}_i\times A_i,$ \\
	$v^{-\sigma}_0 = [0,2m2V_1-(m1+m2)V_2,2m1V_2-(m1+m2)V_1].$ \\
	$V^{-\sigma} = v^{-\sigma}_1\times v^{-\sigma}_2,$ \\
	$V^{-\sigma} \cdot v^{-\sigma}_0 = 0 (*).$

	For the determination of $V^{-\sigma}$, I have multiplied the
	components by $m1-m2$\\
	and used the property that $V^{\sigma }$ is on $\kappa iepert$ to
	eliminate $V_1V_2.$  Every
	component is then divisible by $V_0.$

	\sssec{Example.}
	$p = 11,$ $\overline{M} = (1,2,4),$\\
	$V^{\sigma } = (1,-5,-4),$ $V^{-\sigma} = (1,4,3),$\\
	$V_i	= \{(1,-2,5), (1,1,-4), (1,-5,-4),\}$\\
	$V^-_i = \{(1,5,1),(1,4,3), (1,4,2),\}$\\
	the sides of these triangles are\\
	$vv_i = \{[1,0,3], [0,1,7], [1,3,1],\}$
	 $vv^-_i = \{[1,-3,0], [1,-2,-2], [0,1,-5]\}$.
	$vv_i \times vv^-_i = \{(1,4,-4), (1,-4,-1), (1,-1,2),\}$
	 which have [1,4,7] in common.

	\sssec{Exercise}.
	x.x.x. defines a projectivity which fixes the conic of
	$\kappa iepert.$  Determine other properties of this projectivity.

	\sssec{Exercise.}
	Construct $V^{\frac{\pi}{2}-\sigma}$
	and $V^{\frac{\pi}{2}+\sigma}$ and obtain properties involving these
	points and $V^{\sigma},$ $V^{-\sigma}$ and lines derived from these.

	\ssec{The Theorem of Vectem and related results.}
\setcounter{subsubsection}{-1}
	\sssec{Introduction.}
	In classical Euclidean Geometry, the construction
	of the point of Vectem starts with that of squares on the sides of the
	triangle, outside of it.  In the finite case, there is ambiguity and
	the squares need not exist.  It is easy to determine the intersections
	of the circle $\kappa _1$ with the perpendicular through $A_1$ to
	$a_0.$
	This leads to the expression for $X_{1,0}$ given below.  To insure the
	consistency associated to the outside condition of the classical case
	I have started with that definition for $X_{1,0},$ chosen $X_{2,0}$ on
	$\overline{\kappa}_2$ and
	the perpendicular at $A_2$ to $a_0$ in such a way that
	$X_{1,0}\times X_{2,0}$ is
	parallel to $a_0.$  $X_{2,1},$ $X_{0,2}$ and $X_{0,1},$ $X_{1,2},$
	are defined using the
	symmetry operation $\rho$ , defined in section ?.?.?.
	Because $\alpha$  is obtained in section ?.?.?.  using a square root
	operation
	the definitions can be repeated using $-\alpha$  instead of $\alpha$,
	the
	corresponding elements are denoted with a superscript -.  These have
	been given explicitely.  The conclusions have not been written
	explicitely.
	To each conclusion (and proof) coreesponds an other one by exchanging
	no superscript with the superscript $-$ and $\alpha$  by $-\alpha$ .

	Explicit expression for distances and area, if needed, are given uding
	the same notation as in the conclusions, replacing C by F.

	One could also proceed by first choosing one of the intersections of
	$\kappa_1$ with $A_1\times I\overline{m}a_0$ as $X_{1.0}$ and
	 constructing all the
	other points.  For instance, $X_{2,0}$ by
	$(A_2\times I\overline{m}a_0)\times X_{1,0}\times MA_0,$
	$X_{0,1}$ by $(A_0\times I\overline{m}a_1)\times
	(\overline{m}a_2\times (A_0\times X_{1,0})),$ etc.

	\sssec{Theorem.}
	\begin{enumerate}
	\setcounter{enumi}{-1}
	\item[H0.0.]$X_{1,0} := (m0(m1+m2),\alpha -m0m1,-m2m0),$
		 $X_{1+i,i} := \rho^i X_{1,0},$
		$X_{2,0} := (m0(m1+m2),-m0m1,\alpha -m2m0),$
		$X_{2+i,i} := \rho^i X_{2,0},$
	\item[H0.1.]$X^-_{1,0}:= (m0(m1+m2),-\alpha -m0m1,-m2m0),$
		 $X^-_{1+i,i} := \rho^i X^-_{1,0},$
		$X^-_{2,0}	:= (m0(m1+m2),-m0m1,-\alpha -m2m0),$
		 $X^-_{2+i,i} := \rho^i X^-_{2,0},$
	\item[D0.1.]$x_{i,j,k} := A_i\times X_{j,k}, i \neq k.$
		$x^-_{i,j,k} := A_i\times X^-_{j,k}, i \neq k.$
	\item[D0.2.]$V_i := x_{i+1,i-1,i}\times x_{i-1,i+1,i},$
		$V^-_i := x^-_{i+1,i-1,i}\times x^-_{i-1,i+1,i},$
	\item[D0.3.]$W_i := x_{i+1,i,i+1}\times x_{i-1,i,i-1},$
		$W^-_i := x^-_{i+1,i,i+1}\times x^-_{i-1,i,i-1},$
	\item[D0.4.]$U_i := x_{i+1,i-1,i+1}\times x_{i-1,i+1,i-1},$
		$U^-_i := x^-_{i+1,i-1,i+1}\times x^-_{i-1,i+1,i-1},$
	\item[D0.5.]$v_i := A_i v V_i, v^-_i := A_i\times V^-_i,$
	\item[D0.6.]$w_i := X_{i+1,i-1}\times X_{i-1,i+1},$
		$w^-_i := X^-_{i+1,i-1}\times X^-_{i-1,i+1},$
	\item[D0.7.]$V := v_1\times v_2, V^- := v^-_1\times v^-_2,$
	\item[D0.8.]$v : = V\times V^-.$
	\item[D1.0.]$Ix_{i,j,k} := m\times x_{i,j,k}, j \neq k.$
	\item[D1.1.]$Iv_i = m\times v_i,$
	\item[D1.2.]$Iw_i = m\times w_i,$\\
	{\em then}
	\item[C0.0.~]$(X_{i+1,i}\times X_{i-1,i}) \cdot MA_i = 0.$
	\item[C0.1.~]$V \cdot v_0 = 0 (*).$
	\item[C0.2.~]$U_i \cdot \overline{m}a_i = 0.$
	\item[C0.3.~]$W_i \cdot w_i = 0.$
	\item[C0.4.~]$W_i \cdot v_i = 0.$
	\item[C0.5.~]$V^-_i \cdot w_i = 0.$
	\item[C0.6.~]$x_{i+1,i,i+1} ??? x_{i-1,i,i-1}.$
	\item[C0.7.~]$v_i \cdot w_i.$
	\item[C0.8.~]$dist^2(A_{i+1},X_{i+1,i})
		= dist^2(A_{i-1},X_{i-1,i}) = dist^2(A_{1+1},A_{1-1}).$
	\item[C0.9.~]$dist^2(A_{i+1},X_{i,i+1}) = dist^2(A_{i-1},X_{i,i-1}).$
	\item[C0.10.]$dist^2(A_i,V_i) = dist^2(V_{i+1},V_{i-1}).$
	\item[C0.11.]$Area(A_i,X_{i,i+1},X_{i,i-1}) = Area(A_0,A_1,A_2).$

	\noindent\item[P0.1.~]$x_{0,0,1} = [0,m2,m2+m0].\\
		x_{0,0,2} = [0,m0+m1,m1].\\
		x_{0,2,0} = [0,\alpha -m2m0,m0m1].\\
		x_{0,1,0} = [0,m2m0,\alpha -m0m1].\\
		x_{0,2,1} = [0,m1m2-\alpha ,m1(m2+m0)].\\
		x_{0,1,2} = [0,m2(m0+m1),m1m2-\alpha ].$
	\item[P0.2.~]$V_0 = (m0(m1+m2),\alpha -m0m1,\alpha -m2m0).$
	\item[P0.3.~]$W_0 = (m0(m1+m2)\alpha -m0m1m2(s_1+m0),m1m2(\alpha -m0m1),
		m1m2(\alpha -m2m0)).$
	\item[P0.4.~]$U_0 = (-m0m1m2,m1(\alpha -m1m2),m2(\alpha -m1m2)).$
	\item[P0.5.~]$v_0 = [0,m2m0-\alpha ,\alpha -m0m1].$
	\item[P0.6.~]$w_0 = [2m1m2\alpha ,m2m0(m1s_1-\alpha ),
		m0m1(m2s_1-\alpha )].$
	\item[P0.7.~]$V = ((\alpha -m2m0)(\alpha -m0m1),
		(\alpha -m0m1)(\alpha -m1m2),(\alpha -m1m2)(\alpha -m2m0)).$
	\item[P0.8.~]$v =[(m1-m2)(m0s_1-m1m2),(m2-m0)(m1s_1-m2m0),
		(m0-m1)(m0s_1-m1m2)].$

	\noindent
	\item[F0.6.~]$dist^2(A_0,X_{2,0}) = 2 \alpha  - 2m0m1 - m2(m0+m1).$\\
		$dist^2(A_0,X_{1,0}) = 2 \alpha  - 2m2m0 - m1(m2+m0).$

	\item[P1.0.~]$Ix_{0,0,1} = (m0,-(m2+m0),m2).\\
		Ix_{0,0,2} = (m0,m1,-(m0+m1)).\\
		Ix_{0,2,0} = (\alpha -m0(m1+m2),m0m1,m2m0-\alpha ).\\
		Ix_{0,1,0} = (m0(m1+m2)-\alpha ,\alpha -m0m1,-m2m0).\\
		Ix_{0,2,1} = (\alpha -m1(2m2+m0),m1(m2+m0),m1m2-\alpha ).\\
		Ix_{0,1,2} = (m2m0+\alpha ,m1m2-\alpha ,-m2(m0+m1)).$
	\item[P1.1.~]$Iv_0 = (m0(m1+m2)-2\alpha ,\alpha -m0m1,4-m2m0).$
	\item[P1.2.~]$Iw_0 = (m0(m1-m2)\alpha ,m1(m2m0s_1-(2m2+m0)\alpha ,
		m2((2m1+m0)\alpha -m0m1s_1)).$

	The nomenclature:
	\end{enumerate}

	\sssec{Theorem.}
	\begin{enumerate}
	\setcounter{enumi}{-1}
	\item$	2 Area(A_1,A_2,X_{2,1}) = -m0m1.$
	\item$	2 Area(A_1,A_2,X_{0,1}) = \alpha  - m0m1.$
	\item$	dist^2(V_1,V_2) =$
	\end{enumerate}

	\sssec{Comment.}
	The isotropic points are real if $-\alpha$ is a quadratic residue
	(5.5.2.) if $p \equiv  1 \pmod{4}$, there are $p-1$ ordinary points on
	any circle
	and $\pi  = p-1$ is divisible by 4, therefore a square can be
	constructed,
	the diagonals forming the angle $\frac{\pi}{4}$ with the sides, this
	 is consistent
	with the fact that $X_{i,j}$ are integers.  If $\alpha$  is imaginary
	and $p \equiv -1 \pmod{4},$ then $\pi = p+1$ is divisible by 4 and the
	same situation exist.

	\sssec{The equilateral and hexagonal points.}

	Let $\beta = \frac{\alpha}{\sqrt{3}}$,
	\begin{enumerate}
	\setcounter{enumi}{-1}
	\item[H0.0.]$V^e_i = \kappa_{i+1}\times \kappa_{i-1},$\\
		$V^e_0 = (m0(m1+m2),\beta -m0m1,\beta -m2m0).$
	\item[D0.0.]$vea_i := V_{i+1}\times A_{i-1},$\\
		$vea_0 = [m1(m2+m0),m0m1-\beta ,0].$\\
		$ve\overline{a}_i := V_{i-1}\times A_{i+1},$\\
		$ve\overline{a}_0 = [m2(m0+m1),0,m2m0-\beta ].$
	\item[D0.1.]$Vea_i := vea_{i+1}\times MA_{i-1},$\\
		$Vea_0 = (\beta -m0m1,m1(m2+m0),m1m2+\beta )).$\\
		$Ve\overline{a}_i := ve\overline{a}_{i-1}\times MA_{i+1},$\\
		$Ve\overline{a}_0 = (\beta -m2m0),-m1m2-\beta ,m2(m0+m1)).$
	\item[D0.2.]$veb_i := Vea_{i+1}\times M_{i-1},$\\
		$veb_0 = [m2(m0+m1),-m2(m0+m1),2\beta +m2(m0-m1)].$\\
		$ve\overline{b}_i := Ve\overline{a}_{i-1}\times M_{i+1},$\\
		$ve\overline{b}_0 = [m1(m2+m0),2\beta -m1(m2-m0)],-m1(m2+m0).$
	\item[D0.3.]$Veb_i := veb_i\times ve\overline{a}_i,$\\
		$Veb_0 = (\beta -m2m0,3\beta -m1m2,m2(m0+m1)).$\\
		$Ve\overline{b}_i := veb~_i\times vea_i,$\\
		$Ve\overline{b}_0 = (\beta -m0m1,m1(m2+m0),3\beta -m1m2).$
	\item[D0.4.]$vec_i := vVeb_i\times A_i,$\\
		$vec_0 = [0,-m2(m0+m1),3\beta -m1m2].$\\
		$ve\overline{c}_i := vVe\overline{b}_i\times A_i,$\\
		$ve\overline{c}_0 = [0,3\beta -m1m2,-m1(m2+m0)].$
	\item[D0.5.]$V^h_i := vec_{i+1}\times ve\overline{c}_{i-},$\\
		$V^h_0 = (m0(m1+m2),3\beta -m0m1,3\beta -m2m0).$

	\noindent\item[D1.0]$vve_i := V^e_{i+1}\times V^e_{i-1},$\\
		$vve_0 = [\beta +m1m2,2\beta -m2(m0+m1),2\beta -m1(m2+m0)].$
	\item[D1.1.]$vvh_i := V^h_{i+1}\times V^h_{i-1},$\\
		$vvh_0 = [m1m2-\beta ,2\beta -m2(m0+m1),2\beta -m1(m2+m0)].$
	\item[D1.2.]$Veh_i := vve_i\times vvh_i,$\\
		$Veh_0 = (0,-(2\beta -m1(m2+m0)),2\beta -m2*(m0+m1)).$
	\item[D1.3.]$veh := Veh_1\times Veh_2,$\\
		$veh = [(2\beta -m1(m2+m0))(2\beta -m2(m0+m1)),\\
	\hth	(2\beta -m2(m0+m1))(2\beta -m0(m1+m2)),\\
	\hth	(2\beta -m0(m1+m2))(2\beta -m1(m2+m0))].$\\
	{\em then}
	\item[C0.0.]$V^e_i \cdot mf_i = 0.$
	\item[C0.1.]$V^e \cdot v^e_0 = 0 (*).$
	\item[C0.2.]$V^h \cdot \kappa iepert = 0$
	\item[C0.3.]$V^h \cdot v^h_0 = 0 (*).$
	\item[C0.4.]$V^e \cdot \kappa iepert = 0$
	\item[C0.5.]$Veh_i \cdot a_i = 0.$
	\item[C0.6.]$Veh_0 \cdot veh = 0 (*).$

	The nomenclature:
	\item[N0.0.]$V^e_i$ are the {\em equilateral points}, such that\\
	\hth$	angle(V^e_i,A_{i+1},A_{i-1}) = \frac{\pi}{3}.$\\
	$V^e_0$ is therefore on $\kappa_{i+1}$ and $\kappa_{i-1}.$
	\item[N0.1.]$V^h_i$ are the {\em hexagonal points}, such that\\
	\hth$	angle(V^h_i,A_{i+1},A_{i-1}) = \frac{\pi}{6}.$\\
	$V^h_0$ is therefore the barycenter of the equilateral triangle
		$(V^e_i,A_{i+1},A_{i-1}).$
	\end{enumerate}

	\sssec{Answer to x.x.4.}
	\hth$X = A_0$ gives $V = A_1,$\\
	\hth$X = A_1$ gives\\
	\hti{12}$V = ((m2+m0)(m0+m1),2m1(m2+m0),2m2(m0+m1)),$\\
	\hth$X = A_2$ gives $V = M.$\\
	\hth$X = (m1+m2,-(m1-m2),m1-m2)$ gives $V = A_2.$




	\setcounter{section}{1}
	\setcounter{subsection}{10}
	\ssec{Representation of involutive geometry on the dodecahedron.}
	\setcounter{subsubsection}{-1}
	\sssec{Introduction.}
	When $p = 5,$ it is natural to try to represent involutive
	geometry on the dodecahedron.  The most natural choice,  for the ideal
	line, in the hyperbolic case is an edge-line.  We can choose two
	face-points as the isotropic points.  In the elliptic case, the
	simplest choice for the ideal line appears to be a vertex-line.  The
	fundamental involution associates to a vertex-point an edge-point.

	\sssec{Definition.}
	In the case of hyperbolic involutive geometry, the
	{\em isotropic points} are chosen as 2 face-points.

	\sssec{Theorem.}
	{\em With the chosen fundamental involution, the perpendicular
	direction of a vertex-point is a vertex point and to an edge-point is
	an edge-point.}

	\sssec{Example.}
	If the isotropic points are 0 and 4, the perpendicular
	directions are 10 and 24 as well as 23 and 26.

	\sssec{Theorem.}
	There are 100 circles in a hyperbolic involutive geometry.

	$\begin{array}{ccl}
	{\rm Number\: of}&{\rm 	center}&{\rm sub-types}\\
	2&		f&		B,G3;D2,H1.\\
	2&		f&		B,G4;D3,H2.\\
	2&		v&		C1,D4;G1,I3.\\
	2&		v&		C2,D1;G2,I2.\\\
	4&		v&		B,H4;G6,G7.\\
	1&		s&		A,E1;E2,I1.\\
	4&		s&		B,I4;F,H3.\\
	4&		s&		C1,D2;D1,G8.\\
	4&		s&		C2,D3;D4,G5.
	\end{array}$

	Proof:  For the type $ffffss,$ out of 15 quadruples only 6 contain
	2 given ones, therefore the number of conics must be divided by
	$\frac{15}{6}.$\\
	For the type $fffxxx,$ out of 20 triples only 4 contain 2 given ones,
	hence the number of conics must be divided by $\frac{20}{4} = 5.$\\
	For the type $ffxxxx,$ out of 15 pairs, only one is the given one,
	therefore, the number of conics is to be divided by 15.\\
	As a check there are $25*16*\frac{6}{24}$ = 100 conics through 2
	given points.\\
	More precisely the conics are\\
	1 of type $ffffff,$ sub-type $A.$\\
	12 of type $ffffss,$ sub-type $B.$\\
	12 of type $fffvvs,$ 6 of sub-type $C1$, 6 of sub-type $C2.$\\
	24 of type $fffvss,$ 6 each of sub-type $D1,$ $D2,$ $D3$ and $D4.$\\
	2 of type $ffvvvv,$ 1 each of sub-type $E1,$ $E2.$\\
	4 of type $ffvvvs,$ sub-type $F.$\\
	24 of type $ffvvss,$ 2 each of sub-type $G1$ to $G4$ and
		   4 each of type $G5$ to $G8.$\\
	12 of type $ffvsss,$ 2 each of sub-type $H1,$ $H2$ and
			   4 each of sub-type $H3$, $H4.$\\
	9 of type $ffssss,$  1 of sub-type $I1,$ 2 each of sub-type $I2,$ $I3,$
			 and 1 of sub-type\\
	\hth$I4.$\\
	The centers and their relationship to the conic have been determined
	using the program [130]DODECA.

	\sssec{Theorem.}\label{sec-eproj}
	{\em In the case of elliptic involutive geometry, if a vertex-line
	is chosen as the ideal line, there exists an elliptic projectivity
	which associates, alternately, to a vertex-point, an edge-point and to
	an edge-point, a vertex-point.}
 
	\sssec{Definition.}
	The projectivity of Theorem \ref{sec-eproj} is chosen
	as the {\em fundamental projectivity}.

	\sssec{Example.}
	We can choose as ideal line $5^*$ and as fundamental
	projectivity (7,13,23,27,29,26).

	\sssec{Theorem.}
	{\em Given a center, there are 4 circles with 6 ordinary points on
	them.  2 have a diameter in the direction of a ideal face-point and 2
	have a diameter in the direction of the other ideal face-point.}

	\sssec{Theorem.}
	{\em There are 100 circles in an elliptic involutive geometry.}\\
	$\begin{array}{ccl}
	{\rm Number\: of}&{\rm center}&{\rm sub}-{\rm types}\\
	3&		f&		H1, M1; S2,G4.\\
	3&		f&		H2, M4; S1,G3.\\
	1&		v&		A,P;U1,U2.\\
	6&		v&		I4,S5;H3,F.\\
	3&		s&		I2,G1;L1,M3.\\
	3&		s&		I3,G2;L2,M2.\\
	6&		s&		H4,S8;G6,G7.
	\end{array}$

	Proof:  The proof was done using the program [130] EUCLID5
	and the interactive program [130] DODECA.  The semi colon
	separates the circles whose diameter have a different ideal
	face-points.

	The details are on G331.PRN.


{\tiny\footnotetext[1]{G34.TEX [MPAP], \today}}
	\section{Finite Sympathic Geometry.}

	\setcounter{subsection}{-1}
	\ssec{Introduction.}

	See Example 1.10
	\ldots Measure of angles, separate from measure of distances,\\
	2 triangles having the same angles are similar, their sides are
	not equal.\\
	For measure of distances we can do it starting from a unit
	(2 ordinary distinct points) on all lines which have the same
	parity (even or odd), the other parity requires an other unit,
	the two become connected as a subset of sympathic projectivity
	which is Euclidean geometry.\\
	Although we could have subordinated measure of angles to measure of
	distance we prefer to do the reverse.

	\ssec{Trigonometry in a Finite Field for $p.$
		The Hyperbolic Case.}
	\setcounter{subsubsection}{-1}
	\sssec{Introduction.}
	The trigonometry associated to the finite Euclidean
	plane with real isotropic points will first be defined and studied in
	this section for the finite field $Z_p.$  Theorems 1.4. and 1.6.
	determine $sin(1)$ and $cos(1)$ from which all other values can be
	obtained using the addition formulas.
	In section 2, definitions and
	results will be extended, for the finite field associated to $p^e,$ with
	proofs left as an exercise.

	\sssec{Definition.}
	Given the sets
	$Z$ of the integers,
	$Z_p$ of the integers modulo $p,$
	$Z_{p-1}$ of the integers modulo $p-1,$
	let $\delta$  be a square root of a non quadratic residue of $p,$
	I define $\pi$ as follows\\
	\hth$\pi  := p-1.$\\
	Therefore $\frac{\pi}{2}$ is an integer.

	The problem addressed here is to construct 2 functions ${\mathit sine}$ or
	$sin$ and
	${\mathit cosine}$ or $cos$ with domain $Z$ and range $\{Z_p, \delta  Z_p\}$
	which satisfy:

	The Theorem of Pythagoras,
	\begin{enumerate}
	\setcounter{enumi}{-1}
	\item 0.$sin^2(x) + cos^2(x) = 1,$

	The addition formulas,\\
		1.$sin(x+y) = sin(x) cos(y) + cos(x) sin(y),\\
		2. cos(x+y) = cos(x) cos(y) - sin(x) sin(y).$

	The periodicity property
	\item 0.$sin(2\pi +x) = sin(x),$ $cos(2\pi +x) = cos(x),$

	The symmetry properties
	\item 0.$sin(\pi +x) = -sin(x),$ $cos(\pi +x) = -cos(x),\\
		1.sin(-x) = -sin(x),$ $cos(-x) = cos(x),\\
		2.sin(\pi -x) = sin(x),$ $cos(\pi -x) = -cos(x),\\
		3.sin(\frac{\pi}{2}-x) = cos(x),$
			$cos(\frac{\pi}{2}-x) = sin(x),\\
		4.sin(\frac{\pi}{2}+x) = cos(x),$
			$cos(\frac{\pi}{2}+x) = -sin(x),$

	and such that
	\item 0.$sin(0) = 0,$ $cos(0) = 1,\\
		1.sin(\frac{\pi}{2}) = 1, cos(\frac{\pi}{2}) =0,\\
		2.cos(x) \neq 0$ for $0 < x < \frac{\pi}{2}.$
	\end{enumerate}

	\sssec{Theorem.}
	{\em Let}
	\begin{enumerate}
	\setcounter{enumi}{-1}
	\item{\em $g$ be a primitive root of $p,$}
	\item $\gamma  := \sqrt{g},$
	\item $i :=  \gamma ^{\frac{p-1}{2}},$
	\item $e(j) := \gamma ^j,$
	\item $sin(j) = \frac{1}{2i}(e(j) - e(-j)),$
		$cos(j) = \frac{1}{2}(e(j) + e(-j)),$\\
	{\em then}
	\item{\em $i^2 = -1$ and\\
	satisfy 1.1.0.0. to 1.1.3.2..}
	\end{enumerate}

	Proof.  Because $g$ is a primitive root,\\
	\hth	$i^2 = g^{\frac{p-1}{2}} = -1.$\\
	From the definition of $sin(j)$ and $cos(j)$ follows\\
	\hth$	cos(j) + i sin(j) = e(j),$
		$cos(j) - i sin(j) = \frac{1}{e(j)},$\\
	\hth	therefore $cos(j)^2 + sin(j)^2 = 1,$ hence 1.1.0.0.\\
	From the exponentiation properties follows\\
	\hth$	e(j+k) = \gamma ^{j+k} = (cos(j) + i sin(j) )
		(cos(k) + i sin(k))$\\
	\hth$		 = (cos(j) cos(k) - sin(j) sin(k))
			+ i (cos(j) sin(k) + sin(j) cos(k)),$\\
	hence 1.1.0.1.  Because of 1. and 5., $e(\frac{\pi}{2})
	= \gamma^{\frac{p-1}{2}} = i,$ $\frac{1}{e(\frac{\pi}{2})} = -i,$
	hence 1.1.3.1.\\
	0. implies that $\frac{\pi}{2}$ is the smallest exponent of $g$
		which gives $-1$,\\
	hence $\frac{\pi}{2}$ is the smallest exponent of $\gamma$  which
	gives +i or -i,\\
	therefore 1.1.3.2.  The proof of all other properties is left as an
	exercise.

	\sssec{Theorem.}
	{\em Assume $p \equiv  1 \pmod{4}.$  Let}
	\begin{enumerate}
	\setcounter{enumi}{-1}
	\item{\em $g$ be a primitive root of $p,$}
	\item$i := g^{\frac{p-1}{4}},$ $\delta  := \gamma  = sqr(g),$
		$g' := - g^{\frac{p-3}{2}},$\\
	{\em then}
	\item$sin(1) = i \frac{g'-1}{2} \delta,$ 
		$cos(1) = \frac{g'+1}{2} \delta.$
	\end{enumerate}

	Proof:  $g g' = - g^{\frac{p-1}{2}} = 1,$ $
		i^2 = g^{\frac{p-1}{2}} = -1.\\
	\hth	\delta ^{-1} = \delta /g = g' \delta ,$\\
	hence\\
	\hth$	sin(1) = \frac{\delta  - \delta ^{-1}}{2 i} 
		= -i\frac{1-g'}{2} \delta ,$\\
	\hth$	cos(1) = \frac{\delta  + \delta ^{-1}}{2} 
		= \frac{1+g'}{2} \delta .$

	\sssec{Theorem.}
	{\em Assume $p \equiv  -1 \pmod{4}.$  Let}
	\begin{enumerate}
	\setcounter{enumi}{-1}
	\item{\em $g$ be a primitive root of $p,$}
	\item{\em $\delta  := i$ or $\delta ^2 := -1,$
		$g' := -g^{\frac{p-3}{4}},$\\
		then}
	\item$sin(1) = (g-1) g'\frac{1}{2},$
		    $cos(1) = (g+1) g'\frac{1}{2} \delta.$
	\end{enumerate}

	Proof:  $g g'^2 = g^{\frac{p-1}{2}} = -1 = 1/\delta ^2,$
		therefore $\gamma  g' = 1/\delta ,$ \\
		$\gamma ^{-1} = g'\delta$ and $\gamma  = gg'\delta , $\\
	hence\\
	\hth$	sin(1) = \frac{\gamma -\gamma ^{-1}}{2 i} 
		= \frac{1}{2}(g-1)g',$\\
	\hth$	cos(1) = \frac{\gamma +\gamma ^{-1}}{2} 
		= \frac{1}{2}(g+1) g' \delta .$

	\sssec{Example.}
	
	For $p = 13,$ $g = \delta ^2 = 2,$ $i = -5,$ $g' = -6,$ then\\
	$\begin{array}{cccr}
		i&       sin(i)&          cos(i)&          tan(i)\\
		0&        0& 		 1&		 0\\
		1&       -2 \delta& 		 4 \delta&   		 6\\
		2&       -6&              -2&		 3\\
		3&        6 \delta&  		 6 \delta&          	 1\\
		4&       -2&              -6&		-4\\
		5&        4 \delta&  		-2 \delta& 		-2\\
		6&        1&		 0& 		 \infty 
	\end{array}$

	For $p = 11,$ $g = 2,$ $g' = -4,$ $\delta ^2 = -1,$ then\\
	$\begin{array}{cccc}
		i&       sin(i)&          cos(i)&          tan(i)\\
		0&        0&		 1&		 0\\
		1&       -2&		 5 \delta &		-4 \delta \\
		2&        2 \delta& 		 4&		-5 \delta \\
		3&        4&		 2 \delta &		-2 \delta \\
		4&        5 \delta& 		-2&		 3 \delta \\
		5&        1&		 0&		 \infty 
	\end{array}$

	\sssec{Theorem.}
	{\em Given a trigonometric table of $sin$ and $cos,$ all other
	$\phi(p-1)$ tables can be obtained by using}
	\begin{enumerate}
	\setcounter{enumi}{-1}
	\item$sin^{(e)}(j) = sin(je),$ $cos^{(e)}(j) = cos(je),$
		$(e,p-1) = 1,$\\
	\hth{\em with $0 < e < p-1.$}
	\end{enumerate}

	Proof:  We know that there are $\phi(p-1)$ primitive roots.
	If $g^e$ is an other primitive root, then\\
	\hth$	g^{(e)} = g^e,$ $(e,p-1) = 1,$\\
	\hth$	\delta ^{(e)} = g^{\frac{e-1}{2}} \delta ,$\\
	\hth	for $p \equiv  1 \pmod{4},$ $i^{(e)} = g^{e\frac{p-1}{4}},$
		$g'^{(e)} = -g^{e\frac{p-3}{2}}$\\
	\hth	for $p \equiv  -1 \pmod{4},$ $g'^{(e)} = -g^{e\frac{p-3}{4}}.$\\
	Substituting in 2.1.3. and 2.1.4. gives the theorem.\\
	Replacing $\delta$  by $-\delta$  gives tables for which
		$sin(\frac{\pi}{2}) = -1.$

	\sssec{Example.}
	For $p$ = 13, $g$ = 2,\\
	\hth$	e = 5, 7, 11,$\\
	\hth$	g^{(e)} = g^e = 6,$ $-2,$ $-6,$\\
	\hth$	\delta ^{(e)} = 4 \delta ,$ $-5 \delta ,$ $6 \delta ,$\\
	\hth$	i^{(e)} = g^{3e} = -5,$ $5,$ $5,$\\
	\hth$	g'^{(e)} = -g^{5e} = -2,$ $6,$ $2$\\
	\hth$	sin^{(e)}(1) = 4 \delta ,$ $-4 \delta ,$ $2 \delta ,$\\
	\hth$	cos^{(e)}(1) = -2 \delta ,$ $2 \delta ,$ $-4 \delta .$\\
	For $e = 5,$ $sin^{(5)}(1) = 5.\frac{3}{2} \delta ^{(e)} 
		= 1 \delta ^{(e)} = 4 \delta ,$
		$cos^{(5)}(1) = -\frac{1}{2} \delta ^{(e)} 
		= 6 \delta ^{(e)} = -2 \delta .$\\
	The tables are:\\
	$\begin{array}{cccr}
		i&       sin(i)&          cos(i)&          tan(i)\\
		0&	 0&		 1&		 0\\
		1&	 4 \delta& 		-2 \delta &		-2\\
		2&	-6&		 2&		-3\\
		3&	-6 \delta &		-6 \delta &		 1\\
		4&	 2&		-6&		 4\\
		5&	-2 \delta &		 4 \delta& 		 6\\
		6&	 1&		 0&		 \infty 
	\end{array}$

	\ssec{Trigonometry in a Finite Field for $q = p^e$.
		The Hyperbolic Case.}

	\setcounter{subsubsection}{-1}
	\sssec{Introduction.}
	After recalling the definition of Galois fields, for
	$p^2,$ I will state the Theorems which generalize 2.1.2., 2.1.3 and
		2.1.4.

	\sssec{Definition.}
	Let $n$ be a non quadratic residue,  the {\em set of elements
	in the Galois field} $p^2,$ $GF(p^2),$ are the polynomials of degree
	0 or 1,
	for which addition is performed modulo $p$ and multiplication is
	performed modulo $I^2-n.$  More specifically
	\hth$	(u I + v) + (u' I + v') 
		= (u+u' \umod{p}) I + (v+v' \umod{p},$\\
	\hth$	(u I + v) . (u' I + v') 
		= (uv'+u'v \umod{p}) I + (vv'+nuu' \umod{p}).$\\
	\hth Moreover $(u I + v)^{-1} = \frac{-u I + v}{v^2-nu^2}.$\\
	More generally, if $P$ is a primitive polynomial of degree $n,$ i.e. a
	polynomial which has no factors with coefficients in $Z_p,$ the set of
	elements in the Galois field $p^e,$ $GF(p^e),$ are the polynomials of
	degree
	less than $e,$ for which addition is performed modulo $p$ and
	multiplication is performed modulo $P.$

	\sssec{Notation.}
	\hth$u I + v$ will be written $u.v$ or $up+v.\\
	\hth t I^2 + u I + v$
	will be written $t.u.v$ or $t p^2 + u p + v,  \ldots  .$

	\sssec{Example.}
	Let $q = 5^2,$ $n = 3,$ $g = I + 1 = 1. 1 = 6,$ then\footnote{26.10.82}\\
	$g^{-1} = -2. 2 = 3. 2 = 17,$ $g^2 = 2. -1 = 2. 4 = 14,$
	$g^4 = 1. -2 = 1. 3 = 8,$ $g^6 = 0. -2 = 0. 3 = 3,$
	$g^{12} = 0. -1 = 0. 4 = 4,$ hence $-g^{11} = g^{-1}.$

	\sssec{Theorem.}
	{\em 2.1.1. generalizes for $p^e.$}

	The proof as well as the proof
	of the other theorems in this section are left as exercises.

	\sssec{Theorem.}
	{\em Assume $q = p^e \equiv  1 \pmod{4}.$  Let}
	\begin{enumerate}
	\setcounter{enumi}{-1}
	\item{\em $g$ be a primitive root of $p^e,$}
	\item$i := g^{\frac{q-1}{4}},$ $\delta  := sqr(g),$
		$g' := - g^{\frac{q-3}{2}},$\\
	{\em then}
	\item$sin(1) = i (g'-1) \delta \frac{1}{2},
		cos(1) = (g'+1) \delta \frac{1}{2}.$
	\end{enumerate}

	\sssec{Theorem.}
	{\em Assume $q = p^e \equiv  -1 \pmod{4}.$  Let}
	\begin{enumerate}
	\setcounter{enumi}{-1}
	\item{\em $g$ be a primitive root of $p^e,$}
	\item$g' = -g^{\frac{q-3}{4}},$ $\delta ^2 = -1,$
		$g^{-1} = -g^{\frac{q-3}{2}},$\\
	{\em then}
	\item$sin(1) = (g-1) g'\frac{1}{2},$
		$cos(1) = (g+1) g' \delta \frac{1}{2}.$
	\end{enumerate}

	\sssec{Example.}
	
	For $q = 5^2,$ $n = 3,$ $\delta ^2 = g = 6,$ $i = g^6 = 3,$
		$g' = -g^{11} = 17,\\
	sin(1) = 2. 4 \delta  = 14 \delta,$
		$cos(1) = 4. 4 \delta  = 24 \delta .\\
	sin(2) = (-1. 0) . (2. 2) - -2. -1 = 3. 4 = 19,$
		$cos^2(1) = (2. -1) . (1. 1) = 1. 0,\\
	cos(2) = 2 cos^2(1) - 1 = 2. 0 - 0. 1 = 2. -1 = 2. 4 = 14.$\\
	This gives the Table:\\
	$\begin{array}{ccc}
	 k&    sin(k)&          cos(k)\\
	 0&	 0&		 1\\
	 1&      14 \delta &		24 \delta \\
	 2 &     19&		14\\
	 3  &    20 \delta &		21 \delta \\
	 4   &    3&		10\\
	 5    &   4 \delta &		12 \delta \\
	 6     & 20&		20\\
	 7&      12 \delta &		 4 \delta \\
	 8 &     10&		 3\\
	 9  &    21 \delta &		20 \delta \\
	10   &   14&		19\\
	11    &  24 \delta& 		14 \delta \\
	12     &  1&		 0
	\end{array}$

	\sssec{Exercise.}
	
	Verify the following and construct the corresponding trigonometric
	table.
	\begin{enumerate}
	\setcounter{enumi}{-1}
	\item For $q = 13^2,$ $n = -2,$ $\delta ^2 = g = 15,$
	$i = g^{42} = 8,$ $g' = -g^{83} = -147,$ $g^{167} = 35,\\
	sin(1) = 110 \delta ,$ $cos(1) = 18 \delta ,$
	\item For $q = 7^2,$ $n = 3,$ $\delta ^2 = -1,$ $g = 8,$\\
	\hth$    sin(1) = 3. 4. \delta  = 25 \delta ,$
		$cos(1) = 2. 2 \delta  = 18 \delta ,$
	\item For $q = 11^2,$ $\delta ^2 = 13,$\\
	\hth$    sin(1) =  0. 2 \delta  = 2 \delta ,$ $
		cos(1) =  8. 1 \delta  = 89 \delta , $
	\item For $q = 13^2,$ $\delta ^2 = 15,$\\
	\hth$    sin(1) =  11. 0 \delta  = 143 \delta ,$
		$ cos(1) =  3. 1 \delta  = 40 \delta , $
	\item For $q = 17^2,$ $\delta ^2 = 20,$\\
	\hth$    sin(1) = 11. 16 \delta  = 203 \delta ,$
		$cos(1) =  7. 5 \delta  = 124 \delta , $
	\item For $q = 5^3,$ $\delta ^2 = 9,$\\
	\hth$    sin(1) = 3. 3. 0 \delta  = 90 \delta ,$
		$cos(1) = 4. 4. 1 \delta  = 121 \delta ,$\\
	\hth$    sin(2) = 87, cos(2) = 110.$
	\end{enumerate}


{\tiny\footnotetext[1]{G35.TEX [MPAP], \today}}

\setcounter{section}{2}
\setcounter{subsection}{0}
	\ssec{Trigonometry in a Finite Field for p.  The Hyperbolic Case.}

	\setcounter{subsubsection}{-1}
	\sssec{Introduction.}
	The trigonometry associated with the finite Euclidean
	plane with real isotropic points will first be defined and studied in
	this section for the finite field ${\Bbb Z}_p.$  Theorems 
	\ref{sec-ttrig1} and \ref{sec-ttrigm1}
	determine $sin(1)$ and $cos(1)$ from which all other values can be
	obtained using the addition formulas of \ref{sec-dtrigh}.\\
	In section \ref{sec-Strighptoe}, definitions and
	results will be extended, for the finite field associated with
	$p^e,$ with proofs left as an exercise.\\
	The trigonometry associated with the finite Euclidean plane with no real
	isotropic points will obtained in section \ref{sec-Strigeptoe},
	$sin(1)$ and $cos(1)$ will be determined, for the general case $p^e$ in 
	\ref{sec-ttrigeptoe} and \ref{sec-ttrigemptoe}.

	\sssec{Definition.}\label{sec-dtrigh}
	Given the sets\\
	${\Bbb Z}$ of the integers,\\
	${\Bbb Z}_p$	of the integers modulo $p,$\\
	${\Bbb Z}_{p-1}$ of the integers modulo $p-1,$\\
	let $\delta$  be a square root of a non quadratic residue of $p.$\\
	\hth$\pi := p-1.$

	The problem addressed here is to construct 2 functions {\em sine} or
	{\em sin} and {\em cosine} or {\em cos} with domain ${\Bbb Z}$ and
	range $\{{\Bbb Z}_p\cup \delta{\Bbb Z}_p\}$ 
	which satisfy:\\[10ex]
	The Theorem of Pythagoras,
	\enumb
	\item$sin^2(x) + cos^2(x) = 1.$

	The addition formulas,
	\item$0. \:sin(x+y) = sin(x) cos(y) + cos(x) sin(y),$\\
	1.	$cos(x+y) = cos(x) cos(y) - sin(x) sin(y).$

	The periodicity property
	\item $   sin(2\pi +x) = sin(x),$ $cos(2\pi +x) = cos(x),$

	The symmetry properties
	\item$0.\:    sin(\pi +x) = -sin(x),$ $cos(\pi +x) = -cos(x),\\
	1.\:    sin(-x) = -sin(x),$ $cos(-x) = cos(x),\\
	2.\:    sin(\pi -x) = sin(x),$ $cos(\pi -x) = -cos(x),\\
	3.\:    sin(\frac{\pi}{2}-x) = cos(x),$ 
				$cos(\frac{\pi}{2}-x) = sin(x),\\
	4.\:    sin(\frac{\pi}{2}+x) = cos(x),$
		   $cos(\frac{\pi}{2}+x) = -sin(x),$

	and such that
	\item$0.\:    sin(0) = 0, cos(0) = 1,$\\
	$1.\:    sin(\frac{\pi}{2}) = 1,$ $cos(\frac{\pi}{2}) =0,$\\
	$2.\:    sin(x) \neq  \pm 1$ for $0 < x < \frac{\pi}{2}.$
	\enume

	\sssec{Theorem.}\label{sec-ttrig}
	{\em Let}
	\enumb
	\item  {\em $g$ be a primitive root of $p,$}
	\item$  \gamma  := \sqrt{g},$
	\item$  \iota :=  \gamma ^{\frac{p-1}{2}},$
	\item$  e(j) := \gamma ^j,$\\
	{\em then}
	\item$  \iota^2 = -1$ and
	\item$  sin(j) = \frac{e(j) - e(-j)}{2\iota},$
		$cos(j) = \frac{e(j) + e(-j)}{2},$\\
	{\em satisfy} \ref{sec-dtrigh}.0 to .4.
	\enume

	Proof.  Because $g$ is a primitive root,\\
	\hth$\iota^2 = g^{\frac{p-1}{2}} = -1.$\\
	From the definition of sin(j) and cos(j) follows\\
	\hth$	cos(j) + \iota sin(j) = e(j),$
		$cos(j) - \iota sin(j) = e(j)^{-1},$\\
		therefore $cos(j)^2 + sin(j)^2 = 1,$ hence \ref{sec-dtrigh}.0.\\
	From the exponentiation properties follows\\
	\hth$e(j+k) = \gamma^{(j+k)} 
		= (cos(j) + \iota sin(j) )(cos(k) + \iota sin(k))\\
	\hti{12}		 = (cos(j) cos(k) - sin(j) sin(k))
				+ \iota (cos(j) sin(k) - sin(j) cos(k)),$
	hence \ref{sec-dtrigh}.1.\\
	From 5 and \ref{sec-dtrigh}.1 follows \ref{sec-dtrigh}.3.\\
	0, implies that $\frac{\pi}{2}$ is the smallest exponent of $g$
	 which gives -1,
	hence $\frac{\pi}{2}$ is the smallest exponent of $\gamma$  
		which gives $+\iota$ or $-\iota,$
	therefore \ref{sec-dtrigh}.4.2.  The proof of all other properties is
	left as an exercise. The next 2 Theorems give $sin(1)$ and $cos(1)$
	first when $\iota\in Z_p$, then when this is not the case.

	\sssec{Theorem.[Hyperbolic case]}\label{sec-ttrig1}
	{\em Assume $p \equiv 1 \pmod{4}.$  Let}
	\enumb
	\item  {\em $g$ be a primitive root of $p$},
	\item$  i := \iota := g^{\frac{p-1}{4}},$
		$\delta  := \gamma  = \sqrt{g},$
		$ g' := - g^{\frac{p-3}{2}},$\\
	{\em then}
	\item$  sin(1) = i \frac{g'-1}{2} \delta ,$ {\em and}
		$ cos(1) = \frac{g'+1}{2} \delta .$
	\enume

	Proof:  $g g' = - g^{\frac{p-1}{2}} = 1,$
		$i^2 = g^{\frac{p-1}{2}} = -1.$
		$\delta^{-1} = \delta g^{-1} = g' \delta$ ,\\
	hence\\
		$sin(1) = \frac{\delta - \delta^{-1}}{2 i} 
			= -i\frac{1-g'}{2} \delta,$ and
		$cos(1) = \frac{\delta + \delta^{-1}}{2} 
			= \frac{1+g'}{2} \delta$.

	\sssec{Theorem.[Hyperbolic case]}\label{sec-ttrigm1}
	{\em Assume $p \equiv  -1 \pmod{4}.$  Let}
	\enumb
	\item  {\em $g$ be a primitive root of $p,$}
	\item$\delta  := \iota$ or $\delta ^2 := g^{\frac{p-1}{2}} = -1,$
		$g' := -g^{\frac{p-3}{4}},$\\
	{\em then}
	\item$  sin(1) = \frac{(g-1) g'}{2},$ 
		$cos(1) = \frac{(g+1) g'}{2} \delta .$
	\enume

	Proof:  $g g'^2 = g^{\frac{p-1}{2}} = -1 = \delta^{-2}.$ Because
		$\gamma := \sqrt{g},$ by taking square roots,
		$\gamma  g' = \delta^{-1},$
		 $\gamma ^{-1} = g'\delta$  and $\gamma = gg'\delta$ ,\\
	hence\\
		$sin(1) = \frac{\gamma-\gamma^{-1}}{2 \iota} 
			= \frac{(g-1)g'}{2},$
		$cos(1) = \frac{\gamma+\gamma^{-1}}{2} 
			= \frac{(g+1)g'}{2}\delta$.

	\sssec{Example.}
	For $p = 13,$ $g = \delta^2 = 2,$ $i = -5,$ $g' = -6,$ then\\
	$\begin{array}{rcrrr}
	j&\vline&       sin(j)&	  cos(j)&	  tan(j)\\
	\hline
	0&\vline&	0\: 	&	 1\:	&	 0\\
	1&\vline&   -2 \delta& 		 4 \delta &  		 6\\
	2&\vline&       -6\:	    &  -2\:&		 3\\
	3&\vline&6 \delta  &		6 \delta &	 	 1\\
	4&\vline&       -2\:	    &  -6\:&	-4\\
	5&\vline&4 \delta  &		-2 \delta &		-2\\
	6&\vline&	1\:	&	 0\: &		 \infty 
	\end{array}$

	For $p = 11$, $g = 2,$ $g' = -4,$ $\delta^2 = -1.$ then\\
	$\begin{array}{rcrrr}
	j&\vline& sin(j)&	  cos(j)	&  tan(j)\\
	\hline
	0&\vline&0\:	&	 1\:		&0\:\\
	1&\vline&-2\:	&	 5 \delta 	&-4 \delta \\
	2&\vline&2 \delta& 	 4\:		&-5 \delta \\
	3&\vline&4\:	&	 2 \delta	&-2 \delta \\
	4&\vline&5 \delta& 	-2\:		& 3 \delta \\
	5&\vline&1\:	&	 0\:		& \infty \:
	\end{array}$

	\sssec{Theorem.}
	{\em Given a trigonometric table of $sin$ and $cos,$ all other
	$\phi(p-1)$ tables can be obtained by using}
	\enumb
	\item$	sin^{(e)}(j) = sin(je),$ $cos^{(e)}(j) = cos(je),$
		$(e,p-1) = 1,$
		with $0 < e < p-1.$
	\enume

	Proof:  We know that there are $\phi(p-1)$ primitive roots.\\
	If $g^{(e)}$ is an other primitive root, then\\
		$g^{(e)}e = g^e,$ $(e,p-1) = 1,$
		$\delta^{(e)} = g^{\frac{e-1}{2}}\delta$ ,\\
		for $p \equiv  1 \pmod{4},$ $i^e = g^{e\frac{p-1}{4}},$ 
			$g'^e = -g^{e\frac{p-3}{2}}$, \\
		for $p \equiv -1 \pmod{4},$ $g'^e = -g^{e\frac{p-3}{4}}.$ \\
	Substituting in 2.1.3. and 2.1.4. gives the theorem.\\
	Replacing $\delta$  by $-\delta$  gives tables for which 
		$sin(\frac{\pi}{2}) = -1.$

	\sssec{Example.}
	For $p = 13,$ $g = 2,$\\
	$\begin{array}{lcrrr}
	e &\vline& 5\:& 7\:& 11\:\\
	\hline
	g^{(e)} = g^e &\vline& 6\:& -2\:& -6\:\\
	\delta^e &\vline& 4 \delta& -5 \delta& 6 \delta\\
	i^e = g^{3e} &\vline& -5\:& 5\:& 5\:\\
	g'^e = -g^{5e} &\vline& -2\:& 6\:& 2\:\\
	sin^{(e)}(1) &\vline& 4 \delta& -4 \delta& 2 \delta\\
	cos^{(e)}(1) &\vline& -2 \delta& 2 \delta& -4 \delta
	\end{array}$\\
	For $e = 5,$ $sin^{(5)}(1) = 5.\frac{3}{2}\delta^e = 1 \delta^e 
		= 4 \delta,$
		$cos^{(5)}(1) = -\frac{1}{2}\delta^e = 6 \delta^e
		 = -2 \delta$.\\
	The tables are:\\
	$\begin{array}{rcrrr}
	j&\vline& sin(j)	  &cos(j)	&  tan(j)\\
	\hline
	0&\vline& 0\:	&	 1\:	&	 0\\
	1&\vline& 4 \delta& 	-2 \delta &	-2\\
	2&\vline&-6\:	&	 2\:	&	-3\\
	3&\vline&-6 \delta& 	-6 \delta& 	 1\\
	4&\vline& 2\:	&	-6\:	&	 4\\
	5&\vline&-2 \delta& 	 4 \delta &	 6\\
	6&\vline& 1\:	&	 0\:	&	 \infty 
	\end{array}$

	\ssec{Trigonometry in a Finite Field for $q = p^e.$
		  The Hyperbolic Case.}\label{sec-Strighptoe}

	\setcounter{subsubsection}{-1}
	\sssec{Introduction.}
	After recalling the definition of Galois fields, I will generalize
	Theorems \ref{sec-ttrig}, \ref{sec-ttrig1} and \ref{sec-ttrigm1}.

	\sssec{Definition.}
	Let $n$ be a non quadratic residue,  the set of elements in the
	{\em Galois field} $GF(p^2),$ associated with $p^2,$ are the polynomials
	of degree 0 or 1,
	for which addition is performed modulo $p$ and multiplication is
	performed modulo $P := I^2-n.$  More specifically\\
	\hth$	(u I + v) + (u' I + v') = (u+u' \umod{p}) I 
			+ (v+v' \umod{p}),\\
	\hth	(u I + v) \cdot (u' I + v') = (uv'+u'v \umod{p}) I
			 + (vv'+nuu' \umod{p}).$\\
	Moreover $(u I + v)^{-1} = \frac{-u I + v}{v^2-nu^2}.$ \\
	More generally, if $P$ is a primitive polynomial of degree $n,$ i.e. a
	polynomial which has no factors with coefficients in ${\Bbb Z}_p,$
	the set of elements in the Galois field $GF(p^e),$ are the 
	polynomials of degree less than $e$, for which addition and
	multiplication is performed modulo $P.$

	\sssec{Notation.}
	$u I + v$ will be written $u.v$ or $up+v,$  $t I^2 + u I + v$
	will be written $t.u.v$ or $t p^2 + u p + v$.

	\sssec{Example.}
	Let $q = 5^2,$ $n = 3,$ $g = I + 1 = 1.1 = 6,$ then\\
	$g^{-1} = -2.2 = 3.2 = 17,$ $g^2 = 2.-1 = 2.4 = 14,$ 
		$g^4 = 1.-2 = 1.3 = 8,$
	$g^6 = 0.-2 = 0.3 = 3,$ $g^{12} = 0.-1 = 0.4 = 4,$
		 hence $-g^{11} = g^{-1}.$

	\sssec{Theorem.}
	{\em Theorems} \ref{sec-ttrig}, \ref{sec-ttrig1} {\em and}
	\ref{sec-ttrigm1} generalize, with $p$ replaced by $q := p^e.$

	\sssec{Example.}
	For $q = 5^2,$ $n = 3,$ $\delta^2 = g = 6,$ $i = g^6 = 3,$
		$ g' = -g^{11} = 17,$\\
	$sin(1) = 2.4 \delta  = 14 \delta,$
		$cos(1) = 4.4 \delta = 24 \delta$.\\
	$sin(2) = (-1.0) \cdot (2.2) - -2.-1 = 3.4 = 19,$
		 $cos^2(1) = (2.-1) \cdot (1.1)
	    = 1.0,$ $cos(2) = 2 cos^2(1) - 1 = 2.0 - 0.1 = 2.-1 = 2.4 = 14.$\\
	This gives the Table:\\
	$\begin{array}{rcrr}
	 k&\vline&    sin(k)	  &cos(k)\\
	\hline
	 0&\vline&	 0\:	&	 1\:\\
	 1&\vline&     14 \delta& 	24 \delta \\
	 2&\vline&      19\:	&	14\:\\
	 3&\vline&     20 \delta& 	21 \delta \\
	 4&\vline&       3\:	&	10\:\\
	 5&\vline&      4 \delta& 	12 \delta \\
	 6&\vline&      20\:	&	20\:\\
	 7&\vline&     12 \delta& 	 4 \delta \\
	 8&\vline&      10\:	&	 3\:\\
	 9&\vline&     21 \delta& 	20 \delta \\
	10&\vline&      14\:	&	19\:\\
	11&\vline&     24 \delta& 	14 \delta \\
	12&\vline&     1\:	&	 0\:
	\end{array}$

	\sssec{Exercise.}
	Verify the following and construct the corresponding trigonometric
	table.
	\enumb
	\item  For $q = 13^2,$ $n = -2,$ $\delta^2 = g = 15,$
	    $i = g^{42} = 8,$ $g' = -g^{83} = -147,$ $g^{167} = 35,$
	    $sin(1) = 110 \delta$, $cos(1) = 18 \delta$,
	\item  For $q = 7^2,$ $n = 3,$ $\delta ^2 = -1,$ $g = 8,$
	    $sin(1) = 3.4. \delta = 25 \delta,$\\
		$cos(1) = 2.2 \delta = 18 \delta,$
	\item  For $q = 11^2,$ $\delta ^2 = 13,$
	    $sin(1) =  0.2 \delta = 2 \delta$,
		$cos(1) =  8.1 \delta  = 89 \delta$,
	\item  For $q = 13^2,$ $\delta ^2 = 15,$
	    $sin(1) =  11.0 \delta = 143 \delta$, 
		$cos(1) =  3.1 \delta = 40 \delta$,
	\item  For $q = 17^2,$ $\delta^2 = 20,$
	    $sin(1) = 11.16 \delta = 203 \delta$, 
		$cos(1) =  7.5 \delta = 124 \delta$,
	\item  For $q = 5^3,$ $\delta ^2 = 9,$
	    $sin(1) = 3.3.0 \delta = 90 \delta$, 
		$cos(1) = 4.4.1 \delta = 121 \delta,$
	    $sin(2) = 87,$ $cos(2) = 110.$
	\enume
	\ssec{Trigonometry in a Finite Field for $q = p^e.$
	  The Elliptic Case.}\label{sec-Strigeptoe}

	\sssec{Notation.}
	$(GF(q),+,.)$ is a finite field with $q = p^e$ elements,\\
	$(GF(q)_b,+,.)$ for the corresponding extension field
	$GF(q)(\beta),$ with $\beta^2 = b,$
	where $b$ is a non quadratic residue modulo $p.$

	\sssec{Convention.}
	I will heretofore assume that $p$ is a given odd prime.
	The sets $G_b$ and $\overline{G}_b$ depend on $q,$ we could indicate
	that dependence by writing $G_{b,q}$ for $G_b$ and
	$\overline{G}_{b,q}$ for $\overline{G}_b.$ 

	\sssec{Definition.}\label{sec-dtrigcomp}
	Let $G_b = GF(q) \cup \{\infty \}.$\\
	The operation $\circ$ is defined by
	\enumb
	\item$ \infty  \circ a = a,$ $a \in G_b,$
	\item$ -a \circ a = a \circ -a = \infty,$ $a \in GF(q),$
	\item$ a \circ a' = \frac{a.a' + b}{a+a'},$ $a$ and
	$a' \in GF(q),$ $a + a'\neq  0.$\\
	To avoid confusion with the power notation in $GF(q),$ the $k$-th power
	in $G_b$ precedes the exponent with $``o"$.  For instance,
	\item$ a^{o0} = \infty ,$ $a^{o1} = a,$
	$a^{ok} = a \circ a^{o(k-1)}.$
	\enume

	\sssec{Theorem.}
	{\em If $\js{b}{p} = -1,$ in other words, if $b$ is non quadratic
	residue modulo $p,$ then}
	\enumb
	\item{\em  $\{G_b, o\}$ is an Abelian group,}
	\item{\em  $\infty$  is the neutral element,}
	\item{\em the inverse of $a \in GF(q)$ is $-a,$}
	\item$ r \circ s = t \Rightarrow  r \circ (-t) = -s.$
	\enume

	Proof: The associativity property follows from\\
	\hth$	(a \circ a') \circ a''
		= \frac{a.a'.a''+b(a+a'+a'')}{a'.a''+a''.a+a.a'+ b}
			= a \circ (a' \circ a''),$\\
	if $a'\neq -a$ and $a \circ a'\neq  -a'',$
	and from the special cases,\\
	\hth$	(a \circ -a) \circ a' = a' = a \circ (-a \circ a'),$\\
	\hth$	(\infty  \circ a) \circ a' = a \circ a'
	= \infty  \circ (a \circ a').$

	\sssec{Example.}
	With $p = 13,$\\
{\footnotesize
	$\begin{array}{rrcrrrrrrrrrrrrr}
	b&\vline&g&g^{o2}&g^{o3}&g^{o4}&g^{o5}&g^{o6}&g^{o7}&g^{o8}&g^{o9}&
		g^{o10}&g^{o11}&g^{o12}&g^{o13}&g^{o14}\\
	\hline
	2&\vline&2&-5&-6&-4&3&-1&0&1&-3& 4& 6&5&-2&\infty\\
	2&\vline&1&-5&4&-4&5&-1&\infty\\
	\hline
	6&\vline&1&-3& 5& 4&2&-6&0&6&-2&-4&-5&3&-1&\infty\\
	\end{array}$}

	\sssec{Comment.}
	If $p = 2,$ $\js{b}{p} = -1$ is never satisfied, hence the restriction
	$p$ odd.

	\sssec{Definition.}
	If $\js{b}{p} = -1$ and $\beta  = \sqrt{b}$, then\\
	\hth$\overline{G}_b = \{1\} \cup \{ \frac{r + \beta }{r - \beta },$
	$ r \in GF(q)\}.$

	The elements in $\overline{G}_b$ are distinct and $\overline{G}_b$
	is a subset of $GF(q)_b$.
	The operation of multiplication in $GF(q)_b$ induces one in the set
	$\overline{G}_b.$ 

	\sssec{Theorem.}
	{\em $(\overline{G}_b,.)$ and $(G_b,o)$ are isomorphic, with the
	correspondance}
	\enumb
	\item{\em $ 1 \in \overline{G}_b$, corresponds to $\infty \in G_{b}$,}
	\item{\em $ \frac{r + \beta }{r - \beta } \in \overline{G}_b$
	corresponds to $r \in G_b.$}
	\enume

	Proof:	$\frac{r + \beta }{r - \beta } . 
		\frac{s + \beta }{s - \beta }
		= \frac{(rs + b) + (r+s) \beta }{(rs + b) - (r+s) \beta }
			= ( \frac{r \circ s + \beta }{r \circ s - \beta} ),$
	 if $r+s\neq  0.$\\
		If $s = -r,$
	     $\frac{r + \beta }{r - \beta } . \frac{s + \beta }{s - \beta }
			= \frac{(r + \beta ) (-r + \beta )}
		{(-r + \beta ) (-r - \beta )} = 1,$

	\sssec{Theorem.}
	\enumb
	\item{\em $ \overline{G}_{b,p}$ is an Abelian group, of order $p+1.$}
	\item{\em $ ( \frac{r + \beta }{r - \beta } )^{p+1} \equiv  1 \pmod{ p}$
		for any $r \in G_{b,p}.$}
	\enume

	\sssec{Lemma.}
	{\em If $A$ is a cyclic group of order $q+1,$ the number of elements
	of order $d,$ where $d |  q+1,$ is $\phi(d)$ and}\\
	\hth$	q+1 = \Sigma_{d| q+1}  \phi(d).$

	\sssec{Lemma. [Gauss]}
	\vspace{-18pt}\hspace{130pt}\footnote{Herstein, p. 76, 39}\\[8pt]
	{\em If $A$ is an abelian group of order $q$ which is not cyclic then
	there exists a divisor $d$ of $q$ such that the number of solutions of
	$x^d = e$ is larger than $d$.}

	\sssec{Lemma.}
	{\em The polynomial\\
	\hth$	R_d := (r + \beta )^d - (r - \beta )^d$\\
	has at most $d$ roots in $G_d.$}

	Proof:  Dividing by $\beta,$ we obtain a polynomial in $Z_p$ of degree
	$d-1,$
	which has therefore at most $d-1$ roots for $z \in Z_p$ or $d$ roots
	in $G_d$ 
	$(\infty$ being a root).

	\sssec{Example.}
	With $p = 13,$ $b = 2,$ if $S_d$ is the set of roots of $R_d$ \\
	$d = 1,$ $S_1 = \{\infty \}.$\\
	$d = 2,$ $S_2 = \{0\} \cup S_1,$\\ 
	$d = 7,$ $S_7 = \{\pm 1,\pm 4,\pm 5\} \cup S_1,$\\
	$d = 14,$ $S_{14} = \{\pm 2,\pm 3,\pm 6\} \cup S_7 \cup S_2.$ 

	\sssec{Theorem.}
	\hth$ (\overline{G}_{b,p},.)$ {\em is a cyclic group of order} $p+1.$\\
	\hth$(G_{b,p},o )$ {\em is a cyclic group of order} $p+1.$

	\sssec{Example.}
	2 is a generator of $G_{13,2},$  1 is a generator of $G_{13,6}.$ 

	\sssec{Theorem.[Elliptic case]}\label{sec-ttrigeptoe}
	{\em  Given $q = p^e \equiv 1 \pmod{ 4}.$ Let}
	\enumb
	\item{\em$ b$ be a non quadratic residue,}
	\item{\em$ r_b$ be a generator of $G_b$,}
	\item$ i^2 := -1,$
	\item$ \beta ^2 := b,$
	\item$  r := r_b^{o\frac{p-1}{4}},$\\
	{\em then}
	\item$ sin(1) = \frac{r^2+b}{r^2-b},$
	    $cos(1) = \frac{-2r i}{r^2-b} \beta .$
	\enume

	Proof:
	Let $\sigma = \frac{r_b+\beta}{r_b-\beta},$ then $\sigma^{p+1} = 1$
	and $0 < i < p+1 \Rightarrow  \sigma^i\neq  1.$\\
	$\rho^2 = \sigma \Rightarrow \rho^{2(p+1)} = 1$ and 
	$0 < i < 2(p+1) \Rightarrow \rho^i\neq 1.$\\
	If we take square roots twice, $\rho^{\frac{p+1}{2}} = \pm i,$ we want
	$\rho^{\frac{p+1}{2}} = \sigma^{\frac{p+1}{4}} = i,$ then
	$\rho = cos(1) + i\: sin(1),$
	and $\rho^{2(p+1)} = cos(2(p+1)) + i sin(2(p+1)) = 1.$\\
	If $r = r_b^{o\frac{p-1}{4}},$ then $\rho^{\frac{p-1}{2}}
		= \sigma^{\frac{p-1}{4}} = \frac{r+\beta}{r-\beta},$ or
	$i = \rho^{\frac{p+1}{2}} = \rho\rho^{\frac{p-1}{2}} =
	\rho\frac{r+\beta}{r-\beta},$ or
	$cos(1) + i\: sin(1) = \rho  = i \frac{r-\beta }{r+\beta },$
	$cos(1) - i\: sin(1) = \rho^{-1}  = -i \frac{r+\beta}{r-\beta},$\\
	therefore\\
	$cos(1) = -2 i\: r\: \beta\frac{1}{r^2-\beta^2},$ 
	$sin(1) = \frac{r^2+\beta^2}{r^2-\beta^2}.$ 

	\sssec{Example.}
	$q = 13,$ $i = 5$, $b = \beta^2 = 6,$ $r_b = 1$, $r = 5,$ $r^2 = -1,$
		$sin(1) = 3,$ $cos(1) = -4 \beta.$\\
	$\begin{array}{rcrrrr}
	  k&\vline& sin(k)&    cos(k)&   tan(k)& atan(k\beta )\\
	\hline
	  0&\vline&   0\:&	      1\:&  	 0 \beta& 	 0\\
	  1&\vline&   3\:&        -4 \beta&  	-5 \beta& 	-3\\
	  2&\vline&   2 \beta&       -4\:&	 6 \beta& 	 4\\
	  3&\vline&   5\:&        -3 \beta&  	-1 \beta& 	 6\\
	  4&\vline&  -3 \beta& 	5\:&      2 \beta& 	 5\\
	  5&\vline&  -4\:&         2 \beta&   	 4 \beta& 	-1\\
	  6&\vline&  -4 \beta& 	3\:&	 3 \beta& 	 2\\
	  7&\vline&   1\:&         0 \beta&  	 \infty\:& 	-2
	\end{array}\\
	\begin{array}{rcrrrrr}
	 q&\vline& b& i&  r&sin(1)&cos(1)\\
	\hline
	 5&\vline& 2& 2&  2&    -2&  -\beta\\
	17&\vline& 4& 3&  6&    -5& -3\beta\\
	29&\vline&12& 2&  7&    -2&-10\beta\\
	37&\vline& 6& 2&  5&     6&  -\beta\\
	41&\vline& 9& 2&-17&    -5&-19\beta
	\end{array}$

	\sssec{Theorem.[Elliptic case]}\label{sec-ttrigemptoe}
	{\em Given $q = p^e \equiv  -1 \pmod{ 4}.$
	Let}
	\enumb
	\item{\em$ b$ be a non quadratic residue,}
	\item{\em$ r_b$ be a generator of $G_b,$}
	\item$ \iota^2 := -1,$ $\delta  := \iota,$
	\item$ \beta ^2 := b,$
	\item$ r := r_b^{o\frac{p+1}{4}},$\\
	{\em then}
	\item$ cos(1) = \frac{r_b}{\sqrt{b-r_b^2}} \delta ,$
		$sin(1) = \frac{r cos(1)}{r_b}.$
	\enume

	Proof:  The proof proceeds at first as in \ref{sec-ttrigeptoe}.
	$\frac{r+\beta}{r-\beta} = i,$ therefore $r = -\beta i,$
	this establishes the relationship
	between the sign for the square root of -1 and $b.$\\
	$cos(2) = \frac{1}{2}(\sigma + \sigma^{-1} ) 
		= \frac{r_b^2+b}{r_b^2-b}$ and
	$sin(2) = \frac{1}{2i}(\sigma - \sigma^{-1} ) = \frac{2r_b}{r_b^2-b}.$\\
	$2 cos^2(1) = 1+cos(2)  \Rightarrow   
	cos(1) = \frac{r_b}{\sqrt{r_b^2-b}},$\\
	moreover\\
	$\frac{r_b}{b-r_b} = -cos^2(1),$ $sin(1)$ follows from
	$2 sin(1) cos(1) = sin(2)
	= \frac{2 r cos^2(1)}{r_b},$ insuring the consistency between the
	signs of
	$sin(1)$ and $cos(1)$ to insure that $sin(\frac{\pi}{2}) = 1.$

	\sssec{Example.}
	$q = 11,$ $\delta^2 = -1,$ $b = 2,$ $r_b = 1,$ $r = -3,$
	$cos(1) = \delta,$ $sin(1) = -3\delta$.\\
	$\begin{array}{rcrrrr}
	  k&\vline& sin(k)&    cos(k)&   tan(k)& atan(k\beta )\\
	\hline
	  1&\vline&  -3 \delta& 	1 \delta& 	-3&	 3\\
	  2&\vline&  -5\:&        -3\:&       -2&	-2\\
	  3&\vline&   4 \delta& 	4 \delta& 	 1&	-1\\
	  4&\vline&  -3\:&        -5\:&        5&	-5\\
	  5&\vline&   1 \delta&       -3 \delta& 	-4&	 4\\
	  6&\vline&   1\:&         0\:&        \infty& 	-4
	\end{array}\\
	\begin{array}{rcrrrrr}
	 q&\vline& b&r_b& r&   sin(1)&   cos(1)\\
	\hline
	 3&\vline& 2&1&  1&    \delta&   \delta\\
	 7&\vline& 3&1&  2&   3\delta& -2\delta\\
	19&\vline& 2&1&  6&   6\delta&   \delta\\
	23&\vline& 5&1&  8&   4\delta&-11\delta\\
	31&\vline& 3&1&-11& -13\delta&  4\delta\\
	43&\vline& 3&5&-13&  -7\delta&  6\delta\\
	47&\vline& 5&4& 18&   3\delta&-15\delta.
	\end{array}$

	\sssec{Theorem.}
	{\em If as usual, $\pi := p+1$ in the elliptic case or 
	$\pi := p-1$ in the hyperbolic case, then}
	\enumb
	\item$3\:|\:\pi \Rightarrow sin(\frac{\pi}{6}) = \frac{1}{2},$
	$ cos (\frac{\pi}{6}) = \frac{\sqrt{3}}{2}.$
	\item$4\:|\:\pi \Rightarrow cos(\frac{\pi}{4}) = \frac{\sqrt{2}}{2}.$
	\item$5\:|\:\pi \Rightarrow
		cos(\frac{\pi}{5}) = \frac{\sqrt{5}+1}{4},$
		$cos(\frac{2\pi}{5}) = \frac{\sqrt{5}-1}{4},$\\
	\hth	$sin(\frac{\pi}{5}) = \frac{\sqrt{10-2\sqrt{5}}}{4},$
		$sin(\frac{2\pi}{5}) = \frac{\sqrt{10+2\sqrt{5}}}{4}.$
	\enume

	In the classical case, there is no ambiguity of sign, because
	$0<x<\frac{\pi}{2} \implies sin(x), cos(x) >0.$ This is not the case
	in a finite field, the formulas can only give the trigonometric 
	functions up to the sign, or alternately one of the values of
	$\sqrt{3}$, $\sqrt{2}$ $\sqrt{5}$, $\sqrt{10\pm 2\sqrt{5}}$ can be
	derived from $cos(\frac{\pi}{6})$, $cos(\frac{\pi}{4})$,
	$cos(\frac{\pi}{5})$, $sin(\frac{\pi}{5})$ and $sin(\frac{2\pi}{5})$.

	\sssec{Example.}
	\enumb
	\item[p = 11,] elliptic case, $sin(2) = -5,$ $cos(2) = -3\implies
		\sqrt{3} = 5,$\\
		with $\delta^2 = -1,$
		$cos(3) = 4\delta\implies \sqrt{2} = -3 \delta.$
	\item[p = 11,] hyperbolic case, $cos(2) = 4,$ $cos(4) = -2\implies
		\sqrt{5} = 4,$\\
		with $\gamma^2 = -3,$
		$sin(2) = -4 \gamma\implies\sqrt{2} = -5 \gamma,$
		$sin(4) = \gamma\implies\sqrt{-4} = 4 \gamma$.
	\item[p = 19,] elliptic case,
		$cos(5) = -9 \delta\implies\sqrt{2} = \delta,$\\
		with $\delta^2 = 2,$ $cos(4) = -2\implies
	 cos(8) = 7\implies \sqrt{5} = -9,$\\
		$sin(4) = 4\implies\sqrt{9} = -3,$
		 $sin(8) = 3\implies \sqrt{-8} = -7.$
	\enume
\newpage
	\sssec{Definition.}

	\sssec{Lemma.}
	{\em If $r_b$ is a generator of $G_b$ and $b' := \frac{b}{r_b^2},$ then}
	\enumb
	\item{\em $ kr_b$ is a generator of $G_{bk^2},$}
	\item{\em 1 is a generator of $G_{b'}.$}
	\enume

	\sssec{Theorem.}
	\enumb
	\item{\em There exists always fundamental roots.}
	\item{\em There exist $\frac{1}{2}\phi(p+1)$ fundamental roots
	associated with $p.$}
	\enume

	\sssec{Example.}
	6, 7 and 8 are the fundamental roots for $p = 13,$\\
	6, 7 and 12 are the fundamental roots for $p = 17,$\\
	3, 11, 18 and 27 are the fundamental roots for $p = 29,$\\
	6, 14, 15, 18, 19, 20, 23, 24 and 32 are the fundamental roots
	    for $p = 37,$\\
	12, 13, 28, 29, 30 and 35 are the fundamental roots for $p = 41.$

	\sssec{Theorem.}
	{\em Given an involution, $I(x) = \frac{a x + b}{c x - a},$
	$aa + bc \neq  0,$
	an amicable projectivity, in other words a projectivity with the same
	fixed points, real or complex, is given by\\
	\hth$	T(x) = ( \frac{a+f) x + b }{ c x - a + f},$\\
	where $f = \pm  \sqrt{(aa + bc)/d)}.$}

	Proof: If $\js{d}{p} = -1$ then $F_d$ is a fundamental projectivity:\\
	\hth$	F_d = \frac{y + d}{y + 1}.$\\
	In view of \ref{sec-dtrigcomp}, we have $F_d = y \circ 1$.

	\sssec{Comment.}
	\hth$( \frac{r + \beta }{r - \beta } )^e \equiv  1,$ for $e|p+1
	\implies (r + \beta )^e - (r - \beta )^e \equiv  0,$\\
	dividing by $\beta$ , we obtain a polynomial in $Z_p$ of degree $e-1,$
	 which
	has therefore at most $e-1$ roots for $z$ in $Z_p$ or $e$ roots in $G$
	($\infty$  being
	a root).  We want to show that $(G~,.)$ and therefore $(G,\infty )$
	 are cyclic groups.

	\sssec{Theorem.}
	{\em Let}
	\enumb
	\item $r_0 = 1,$ $r_1 = r,$ $r_{i+1} = r_i \circ r,$
	\item $x_i \equiv  1/r_i (\in G~),$
	\item $u_{i+1} \equiv  r u_i + s u_{i-1}, u_0 = 0,$ $u_1 = 1,$\\
		$v_{i+1} \equiv  r v_i + s v_{i-1}, v_0 = 2, v_1 = r,$
	\item $4s \equiv  d -r^2,$
	\item $\alpha  = \frac{r+\sqrt{d}}{2}$ {\em and}
	$\beta  = \frac{r-\sqrt{d}}{2},$\\
	{\em then}
	\item $x_{i+1} \equiv  \frac{r x_i + 1}{d x_i + r} \pmod{p},$ 
		$x_0 = 0,$
	\item $r = \alpha  + \beta , \sqrt{d} = \alpha -\beta ,$
		$ s = -\alpha  \beta .$
	\item $u_i = (\alpha ^i - \beta ^i)/\sqrt{d},$
		$ v_i = \alpha ^i + \beta ^i,$
	\item $2 u_{i+j} = u_i v_j + v_i u_j,$
		$2 v_{i+j} = v_i v_j + b u_i u_j,$
	\item $u_{i+1} (b u_i + r v_i) - v_{i+1} (r u_i + v_i) = 0,$
	\item $x_i v_i = u_i.$
	\enume
	Proof:
	$(r u_i + s u_{i-1}) = \frac{ (\alpha +\beta )(\alpha ^i-\beta ^i) 
	- \alpha  \beta (\alpha ^{i-1}-\beta ^{i-1}) }{\sqrt{d}}
	= \frac{\alpha ^{i+1}-\beta ^{i+1}}{\sqrt{d}} = u_{i+1}.$\\
	Substituting in 2. with $j = 1$ after multiplication by 2 gives\\
	$((u_i v_1 + v_i u_1) (b u_i + r v_i) -
			(v_i v_1 + b u_i u_1) (r u_i + v_i)$\\
	\hth$	= u_i^2 (b v_1 - r b u_1)
			+ u_i v_i (b u_1 + r v_1 - r v_1 - b u_1)
			+ v_i^2 (r u_1 - v_1) = 0.$

	Morover,\\
	$\frac{r u_i + v_i}{d u_i + r v_i}
		= \frac{ (\alpha +\beta )(\alpha ^i-\beta ^i) 
		+ (\alpha -\beta )(\alpha ^i + \beta ^i) }
			 { ((\alpha -\beta )(\alpha ^i-\beta ^i)
		+ (\alpha +\beta )(\alpha ^i+\beta ^i) ) \sqrt{d}}
		= \frac{u_{i+1}}{v_{i+1}} = x_{i+1}.$

	\sssec{Example.}
	For $p = 7,$ elliptic case, $\iota^2 = -1,$\\
	\hth$	A_k = (cos(k), sin(k), 1)$\\
	are points on the conic\\
	\hth$	x^2 + y^2 = z^2.$\\
	If we define it as a circle and $z = 0$ as the ideal line, the
	isotropic points are not real and we have a Euclidean geometry.\\
	The center of the circle, which is the pole of $z = 0$ is (0,0,1).
	There are 8 real points on the circle,\\
	\hth$  A_0 = (1,0,1),$ $A_2 = (-2,-2,1),$ $A_4 = (0,1,1),$
	$ A_6 = (2,-2,1),$\\
	\hth$  A_8 = (-1,0,1),$ $A_{10} = (2,2,1),$ $A_{12} = (0,-1,1),$
	$ A_{14} = (-2,2,1).$\\
	The distances on the lines $a_{2k} = O\times A_{2k},$
	 $k = 0,$ 1, 2, 3, are real.\\
	 $ a_0 = [0,1,0],$ $a_2 = [2,-2,0],$ $a_4 = [1,0,0],$ $a_6 = [1,1,0].$\\
	The other lines through $O$ intersect the circle at complex points:\\
	\hth$    A_1 = (2\iota,4\iota,1),$ $A_3 = (4\iota,2\iota,1),$
		$A_5 = (-4\iota,2\iota,1),$\\
	\hth$	A_7 = (-2\iota,4\iota,1),$ $A_9 = (-2\iota,-4\iota,1),$
		$A_{11} = (-4\iota,-2\iota,1),$\\
	\hth$	A_{13} = (4\iota,-2\iota,1),$ $A_{15} = (2\iota,-4\iota,1).$\\
	These are on the lines $a_{2k+1} = O\times A_{2k+1},$ $k = 0,$ 1, 2 ,
	3,\\
	\hth$    a_1 = [-2,1,0],$ $a_3 = [1,-2,0],$ $a_5 = [1,2,0],$
		$a_7 = [2,1,0].$

	If $B_1 = (1,2,1)$ is a real point on $a_1,$ $B_1$ is on a circle\\
	\hth$	x^2 + y^2 = 5 z^2.$\\
	this circle intersects $a_{2k+1}$ at real points and $a_{2k}$ at complex
	points.  The distances between points on $a_1$ are multiples of $\iota,$
	because $\sqrt{5} = 3\iota.$  The same is true on the lines $a_3,$
	$a_5,$ $ a_7.$

	\sssec{Definition.}
	The smallest $j$ such that $u_j \equiv  0 \pmod{p}$ is called the
	{\em rank of apparition of} $p.$  Hence

	\sssec{Theorem.}
	{\em For a fixed $r$ there are $\frac{1}{2}\phi(p+1)$ values of $s$ in
	$ [1,p-1]$
	for which the rank of apparition is $p+1.$  More generally, there are
	$\phi(\frac{e}{2})$ values of $s$ in $[1,p-1]$ for which the rank of
	apparition is $e,$ $e$ divides $p+1,$ $e > 2.$}

	\sssec{Theorem.}
	{\em If $b$ is a fundamental root modulo $p$, then\\
	\hth$b = {\bf N}p, 1-b = {\bf N}p.$}

	\sssec{Comment.}
	For $p = 17$, 11{\bf N}p, 7{\bf N}p, but 7 is not a fundamental root.

	\sssec{Theorem.}
	{\em For a given $p,$ the sets\\
	$S_h = \{ cos^2(j),$ $j = 1, \ldots , \frac{p-3}{2}\},$ in the
	 hyperbolic case and\\
	$S_e = \{ cos^2(j),$ $j = 1, \ldots , \frac{p-1}{2}\},$ in the
	 elliptic case are\\
	independent of the choice of the primitive root or of the fundamental
	root.}

	\sssec{Example.}
	$p = 11,$ $S_h = \{4,5,7,8\},$ $S_e = \{2,3,6,9,10\}$\\
	$p = 29,$ $S_h = \{3,5,6,7,11,12,15,18,19,23,24,25,27\},$\\
	\hth$S_e = \{2,4,8,9,10,13,14,16,17,20,21,22,26,28\}.$

	\sssec{Theorem.}
	\enumb
	\item{\em Let $sin(j) \in \{{\Bbb Z}_p - \{0\} - \{1\}\}$,\\
	if $p \equiv  -1 \pmod{4},$ then $sin(j) sin(k)\neq  1,$ for all $k,$\\
	if $p \equiv  1 \pmod{4},$ then $sin(j) sin(k) = 1$ for some $k.$}
	\item{\em $\frac{sin(j)}{\gamma} \in \{{\Bbb Z}_p - \{0\}\}$\\
	if $p \equiv  1 \pmod{4},$ then $sin(j) sin(k)\neq  1,$ for all $k,$\\
	if $p \equiv  -1 \pmod{4},$ then $sin(j) sin(k) = 1,$ for some $k.$}
	\enume

\newpage
	\sssec{Example.}
	$\ldots$ Give examples of associated fundamental sympathic
	projectivities, see 1.10.\\
	If $T(r) = \frac{r+d}{r+1}.$  For $p = 7,$ $d = 3,$ (5 is the other
	choice)\\
	$\begin{array}{lcccccccc}
	r&     \infty&    0&  1&  2&  3&  4&  5&  6\\
	T(r)&   1&  3&  2&  4&  5&  0&  6& \infty
	\end{array}$

	\ldots talk about transformations such as $r = 2s$
	leading to the form used in 1.10\\
	$S(s) = \frac{2s+b}{2-3s}.$

	\sssec{Example.}
	For $p = 13,$ (see g35.bas .5.)\\
	$A = (0,1,0),$ $A_j = (1,j,0),$ $j = 0,$  \ldots, $q-1,$
	$ A \times A_j = [0,0,1],$
	$d_j = dist(A,A_j):$ $cos(d_j) = \frac{j}{\sqrt{1+j^2}},$
	$sin(d_j) =  \frac{1}{\sqrt{1+j^2}},$\\
	$\begin{array}{ccccc} 
	j&       \sqrt{1+j^2}&      sin(d_j)& cos(d_j)& d_j\\
	0&      1&               1&       0&       \frac{1}{2}\\
	1&                \ldots.
	\end{array}$

	$tan(dist(A_j,A_k)) = \frac{k-j}{1+jk}$
	$tan(dist(A,A_k) = -\frac{1}{k}$\\
	$tan(dist(A_j,A_l)) = tan(dist(A_j,A_k) + dist(A_k,A_l)).$\\
	Indeed, the second member is\\
	$\frac{ (k-j)(1+lk) + (l-k)(1+jk) }{ (1+jk)(1+lk) - (k-j)(l-k) }
		= \frac{ (l-j) (1+k^2)}{(1+lj)(1+k^2)} = tan(dist(A_j,A_l))$.\\
	and\\
	$tan(dist(A_j,A_l)) = tan(dist(A_j,A) + dist(A,A_l))$\\
	Indeed, the second member is\\
	$\frac{\frac{1}{j} - \frac{1}{l}}{1+ \frac{1}{kl}}
		 = \frac{l-j}{1+kl} = tan(dist(A_j,A_l))$.\\
	$\begin{array}{lcccccccc}
	j\setminus k&   0&     1 &      2&       3&       4&       5&       6\\
	\hline
	0& &\frac{3}{12}& \frac{11}{12}&\frac{10}{12}&\frac{4}{12}&
		&\frac{5}{12}\\
	1& \frac{9}{12}& &\frac{8}{12}& \frac{7}{12}& \frac{1}{12}&
		&\frac{2}{12}\\
	2& \frac{1}{12}& \frac{4}{12}& &\frac{11}{12}&\frac{5}{12}&
		&\frac{6}{12}\\
	3& \frac{2}{12}& \frac{5}{12}& \frac{1}{12}& &\frac{6}{12}&
		&\frac{7}{12}\\
	4& \frac{8}{12}&\frac{11}{12}&\frac{7}{12}& \frac{6}{12}&
		&\frac{1}{12}\\
	5\\
	6& \frac{7}{12}& \frac{10}{12}&\frac{6}{12}&\frac{5}{12}&
		\frac{11}{12}\\
	A& \frac{6}{12}&\frac{9}{12}&\frac{5}{12}& \frac{4}{12}&
		\frac{10}{12}& &\frac{11}{12}
	\end{array}$

	\sssec{Example.}
	$p = 13,$ $h,$ the correspondence between the point
	$A_j = (1,j,0)$ and $d(j)$ is\\
	$\begin{array}{lcccccccccccc}
	j&       0& -2& -3&  1&  4&  6&  \infty & -6& -4& -1&  3&  2\\
	12d(j)&  0&  1&  2&  3&  4&  5&  6&  7&  8&  9& 10& 11
	\end{array}$

	$i = \pm 5$ corresponds to the ideal point, $d(5) = d(-5) = \infty$.\\
	elliptic case:\\
	$\begin{array}{lcccccccccccccc}
	{j}{\delta}&	0&  3& -1& -2&  4& -5&  6&  \infty& -6&5&-4& 2& 1& -3\\
	d(j)&	0&  1&  2&  3&  4&  5&  6&  7&  8&  9& 10& 11& 12& 13
	\end{array}$

	$tan(dist(A_j,A_k)) = \frac{d(k) - d(j)}{1 + d(j)d(k)}$.



	\ssec{Periodicity.}

	\sssec{Definition.}
	Let $f$ be a function, $g(0)$ is arbitrary,\\
	\hth$	g(i+1) = ai + g(i) + f(i) + f(i+1),$\\
	where $a$ is such that\\
	\hth$	g(T) = g(0),$\\
	and we write\\
	\hth$	g = T f.$

	\sssec{Theorem.}
	{\em If $f$ is a periodic function with period $T,$ then\\
	0.\hti{7}$g$ is periodic.\\
	1.\hti{7}$f$ odd $ \Rightarrow   g even.$}

	\sssec{Example.}
	$\begin{array}{lrrrrrrrrrrr}
		i&	0& 1& 2& 3& 4& 5& 6& 7& 8& 9&10\\
		f(i)&	1&-3&-1& 1& 3&-1& 3& 1&-1&-3& 1\\
		g(i)&	0&-2&-6&-6&-2& 0& 2& 6& 6& 2& 0
	\end{array}$\\
	In this example $a = 0.$

	\sssec{Example.}
	$\begin{array}{lrrrrrr}
		i	&0& 1& 2& 3& 4& 5\\
		f(i)	&0& 4& 7& 7& 4& 0\\
		g(i+\frac{1}{2})&0& 9&-8& 9& 0\\
		h(i)-ai	&0& 0& 9& 1&-9&-9\\
		h(i)	&0&-2& 5&-5& 2& 0
	\end{array}$\\
	In this example $a = -\frac{-9}{5}.$

	\sssec{Definition.}
	Let $f$ be a function.
	Let $g$ be defined by\\
	\hth$	g(i+\frac{1}{2}) = u(f_i,f_{i+1}),$\\
	where $u$ is symmetric in its arguments,\\
	we write\\
	\hth$	g = U f.$

	Let $h$ be defined by\\
	\hth$	h(0)$ is arbitrary,\\
	\hth$	h(i+1) = ai + h(i) + g(i+\frac{1}{2}),$\\
	\hth$	h(T) = h(0).$\\
	we write\\
	\hth$	h = M U f.$

	\sssec{Theorem.}
	{\em If be a periodic function with period $T,$ then\\
	0.\hti{7}$     h$ is a periodic function with period $T.$\\
	1.}

	\ssec{Orthogonality.}

	\sssec{Theorem.}
	{\em If $p \equiv  -1 \pmod{4},$ choose the elliptic case and
	$q = \frac{p+1}{2},$ \\
		If $p \equiv 1 \pmod{4},$ choose the hyperbolic case and
	$q = \frac{p-1}{2},$ \\
	0.\hti{7}The trigonometric functions $sin$ and $cos$ are orthogonal.\\
	1.\hti{7}( \ldots ..          ..)\\
	\hth$	( \ldots  sin (i j)  \ldots ),$ $i,j = 1$ to $q-1$\\
	\hth	( \ldots ..          ..)\\
	\hth	is orthogonal, symmetric and $S S = \frac{q}{2} I.$\\
	2.\hti{7}$(\frac{1}{2} \ldots ..  s   \ldots .. \frac{1}{2})$ \\
	\hth$	C =     ( s  cos (i j)  (-1)^i s),$ $i,j = 0$ to $q+1,$\\
	\hth$	(\frac{1}{2} \ldots (-1)^j s .. \frac{1}{2})$\\
		with $s^2 = \frac{1}{2},$ \\
		is orthogonal, symmetric and $C C = \frac{q}{2} I.$}

	\sssec{Example.}
	$p = 7,$ Elliptic case, $q = 4,$\\
	\hti{32}$			(-3 -2 -2 -2 -3)\\
	\hti{16}	(-2  1 -2),\hti{5}(-2 -2  0  2  2)\\
	\hti{12}    S = ( 1  0 -1), C = (-2  0 -1  0 -2)\\
	\hti{16}	(-2 -1 -2),\hti{5}(-2  2  0 -2  2)\\
	\hti{32}			(-3  2 -2  2 -3)\\
		p = 13, hyperbolic case, q = 6,\\
	\hti{40}				(-6  s  s  s  s  s -6)\\
	\hti{16}	(-6  2  1  2 -6)\htn    ( s  2 -6  0  6 -2 -s)\\
	\hti{16}	( 2  2  0 -2 -2)\htn    ( s -6  6 -1  6 -6  s)\\
	\hti{12}    S = ( 1  0 -1  0  1),\: C = ( s  0 -1  0  1  0 -s)\\
	\hti{16}	( 2 -2  0  2 -2)\htn    ( s  6  6  1  6  6  s)\\
	\hti{16}	(-6 -2  1 -2 -6)\htn    ( s -2 -6  0  6  2 -s)\\
	\hti{40}				(-6 -s  s -s  s -s -6)$\\
	with $s = 2 \delta, \delta^2 = 5.$


	\ssec{Conics in sympathic geometry.}

	\sssec{Theorem.}
	{\em 
	$X0_i = a cos(2i),$ $X1_i = b sin(2i),$ $X2_i = 1.$\\
	$X0_i = \frac{a}{\delta} cos(2i+1),$
	$X1_i = \frac{b}{\delta} sin(2i+1),$ $X2_i = 1$.\\
	$X0_i = a cos(2i),$ $X1_i = \frac{b}{\delta} sin(2i),$ $X2_i = 1.$\\
	$X0_i = \frac{a}{\delta} cos(2i+1),$ $X1_i = b sin(2i+1),$
	$X2_i = 1.$}

	\sssec{Example.}
	Let $p = 11,$ $a = 1,$ $b = 2,$ $\delta  = i,$ $i^2 = -1$.\\
	In the elliptic case, \ldots gives
	$X0^2 + X1^2{\frac{2}{4}} = 1:\\
	 (1,0), (-3,1), (-5,5), (0,2), (5,5), (3,1),\\
			   (-1,0), (3,-1), (5,-5), (0,-2), (-5,-5), (-3,-1).$\\
	\ldots gives\\
	$-X0^2 - X1^2{\frac{2}{4}} = 1:\\
	 (1,5), (4,-3), (-3,2), (3,2), (-4,-3), (-1,5),\\
			    (-1,-5), (-4,3), (3,-2), (-3,-2), (4,3), (1,-5).$\\
	In the hyperbolic case, \ldots gives\\
	$X0^2 - \frac{1}{4}X1^2 = 1:\\
	 (1,0), (4,4), (-2,-1), (2,-1), (-4,4),\\
			   (-1,0), (-4,-4), (2,1), (-2,1), (4,-4).$\\
	The asymptotic directions are (5,1,0) and $(-5,1,0)$.\\
	\ldots gives\\
	$- X0^2 + \frac{1}{4}X1^2 = 1:\\
	 (5,-4), (2,-3), (0,2), (-2,-3), (-5,-4),\\
			     (-5,4), (-2,3), (0,-2), (2,3), (5,4).$\\
	The asymptotic directions are (5,1,0) and $(-5,1,0).$

	\sssec{Example.}
	Let $p = 13,$ $a = 1,$ $b = 2,$ $\delta  = 2, i = 5.$\\
	In the elliptic case, \ldots gives\\
	$X0^2 + \frac{1}{2}X1^2 = 1:\\
		(1,0), (-4,-3), (5,-2), (3,6), (-3,6), (-5,-2), (4,-3),\\
		(-1,-0), (4,3), (-5,2), (-3,-6), (3,-6), (5,2), (-4,3).$\\
	\ldots gives\\
	$2 X0^2 + \frac{1}{4}X1^2 = 1:\\
		(3,6), (-1,-3), (5,5), (0,2), (-5,5), (1,-3), (-3,6),\\
		(-3,-6), (1,3), (-5,-5), (0,-2), (5,-5), (-1,3), (3,-6).$\\
	In the hyperbolic case, \ldots gives\\
	$X0^2 + \frac{1}{4}X1^2 = 1:\\
	 (1,0), (-2,1), (-6,-4), (0,2), (6,-4), (2,1),\\
			   (-1,0), (2,-1), (6,4), (0,-2), (-6,4), (-2,-1),$\\
	the asymptotic directions are (4,1,0) and $(-4,1,0).$\\
	\ldots gives\\
	$2 X0^2 + \frac{1}{2}X1^2 = 1:\\
	 (4,-4), (6,-1), (-2,-5), (2,-5), (-6,-1), (-4,-4),\\
			     (-4,4), (-6,1), (2,5), (-2,5), (6,1), (4,4),$\\
	the asymptotic directions are (4,1,0) and $(-4,1,0).$

	\ssec{Regular polygons and Constructibility.}

	\sssec{Definition.}
	A {\em regular polygon} \ldots\\
	because the angles are multiples of $\frac{2r}{p-1}$ or 
	$\frac{2r}{p+1}$.\\
	The only regular polygons are those whose number of sides is a
	divisor of $p-1$ or $p+1.$  If we give the unit angle then we can define
	"convex polygons" and "star polygons", 
	\underline{find appropriate names.}\\
	The constructibility by rule and compass in Euclidean
	geometry corresponds here to those which demand the solution of
	equations of the first and second degree.\\
	The work of Gauss on cyclotomic polynomials extend immediately to the
	finite case because of Theorem \ldots on trigonometric functions.

	\sssec{Theorem.}{\em
	0.\hti{7}For a regular polygon of $n$ sides to exist, $n$ must divide
	 \ldots\\
	1.\hti{7}For a regular polygon to be constructible using only equations
	 of the second degree, $n$ must have the form $2^{i0} p_1^{i1} 
	p_2^{i2} \ldots p_k^{ik},$ where $i0$ is a non negative integer, $i1,$
	$i2,$  $\ldots$  $ik,$ are 0 or 1 and $p_j$ are primes of the form
	$2^k + 1,$ namely, 3, 5, 17, 257, 65537,  \ldots .\\
	All square roots are integers except perhaps the last one.}

	\sssec{Theorem.}
	{\em For triangles.\\
	$cos(\frac{2r}{3}) = \frac{1}{2},$
	$sin(\frac{2r}{3}) = \sqrt{\frac{3}{2}}$.}

	\sssec{Example.}
	$p = 5$, elliptic case,\\
	$cos(\frac{2r}{3}) = 3,$ $sin(\frac{2r}{3}) = \sqrt{2}$.\\
	$p = 7,$ hyperbolic case,\\
	$cos(\frac{2r}{3}) = 4,$
	$sin(\frac{2r}{3}) = \sqrt{6},$\\
	$p = 23,$ elliptic case,\\
	$cos(\frac{2r}{3}) = 12,$ $sin(\frac{2r}{3}) = 8.$

	\sssec{Theorem.}
	{\em For hexagons, we first obtain the triangle and then use\\
	\hth$cos(\frac{2r}{6}) = sin(\frac{2r}{3}),$
	$sin(\frac{2r}{6}) = cos(\frac{2r}{3}).$}

	\sssec{Example.}
	$p = 23,$ elliptic case,\\
	$cos(\frac{2r}{6}) = 8,$ $sin(\frac{2r}{6}) = 12.$

	\sssec{Theorem.}
	{\em For pentagons. The polynomial to solve is\\
	\hth$x^2 - x + 1 = 0,$\\
	$cos(\frac{2r}{5}) = \frac{x}{2},$
	$sin(\frac{2r}{5}) = \sqrt{1 - \frac{x^2}{4}}.$}

	\sssec{Example.}
	$p = 11,$ hyperbolic case, $\gamma^2 = 8.$\\
	$x1 = (1+\sqrt{\frac{5}{2}} = \frac{1+7}{2} = 4,$
	$x2 = \frac{1-7}{2} = -3,$\\
	$cos1(\frac{2r}{5}) = 2,$ $cos2(\frac{2r}{5}) = 4,$\\
	$sin1(\frac{2r}{5}) = \sqrt{-3} = \gamma,$
	$sin2(\frac{2r}{5}) = \sqrt{-4} = -4 \gamma$.\\
	The choice of 1 or 2 is arbitrary as is the choice of the sign of the
	coefficient of $\gamma.$  The second case corresponds to 
	$sin(\frac{2r}{5})$ of the trigonometric table, trig.tab.\\
	The first, corresponds to $sin(\frac{6r}{5})$ of the same table.\\
	$p = 19,$ elliptic case, $\delta^2 = 10.$
	$x = \frac{1+9}{2} = 5,$\\
	$cos1(\frac{2r}{5}) = \frac{5}{2} = -7 = cos(\frac{6r}{5})$
	of trig.tab.\\
	$sin1(\frac{2r}{5}) = \sqrt{9} = 3 = sin(\frac{6r}{5}).$

	\sssec{Theorem.}
	{\em For decagons, we first obtain the pentagon and then use
	$cos(\frac{2r}{10}) = \sqrt{\frac{(1 + cos\frac{2r}{5}}{2}},$
	$sin(\frac{2r}{10}) = sin\frac{\frac{2r}{5}}{2 cos(\frac{2r}{10}}$.}

	\sssec{Example.}
	$p = 19,$ elliptic case,\\
	$cos1(\frac{2r}{10}) = \sqrt{-3} = 4 = cos(\frac{6r}{10}),$
	$sin1(\frac{2r}{10}) = \frac{3}{8} = -2.$

	\sssec{Theorem.}
	{\em For 17 sided polygons.
	The polynomials to solve are in succession:\\
	\hth$u^2 + u + 4 = 0,$ of which we choose 1 root,\\
	\hth$v^2 - u v - 1 = 0,$\\
	\hth$	v' = \frac{-3 + 6 v - v^3}{2},$\\
	\hth$	w^2 - v w + v' = 0,$\\
	\hth$cos(\frac{2r}{p+j}) = \frac{w}{2},$
	$sin(\frac{2r}{p+j}) = \sqrt{1-(\frac{w}{2})^2}.$}

	\sssec{Example.}
	$p = 67,$ elliptic case.\\
	\hth$u = \frac{-1+\sqrt{17}}{2} = 16,$\\
	\hth$v = \frac{u+\sqrt{u^2+4}}{2} = 
	(\frac{16+\sqrt{59}}{2} = \frac{16+40}{2} = 28,$\\
	\hth$v' = 61 = -6,$\\
	\hth$w = \frac{v+\sqrt{v^2-4v'}}{2} = \frac{(28+\sqrt{4}}{2} = 15,
	$\\
	\hth$cos1(\frac{2r}{17}) = \frac{15}{2} = -26 = cos(\frac{16r}{17})$\\
	 of the table obtained using the program trig.bas.\\
	\hth$sin1(\frac{2r}{17}) = \sqrt{1-6} = -14.$\\
	The other choices for the roots of the above equations lead, with
	the right choice of sign, to $cos(\frac{2kr}{17}),$ $k$ = 1,2,3,4,5,6,7.
	From these all the other angles can ge obtained using the trigonometry
	identities \ldots.\\
	$p = 137,$ hyperbolic case,\\
	\hth$u = \frac{-1+47}{2} = 23,$ $v = \frac{23+81}{2} = 52,$\\
	\hth$v' = \frac{-3+312-46}{2} = 63,$ $w = \frac{52+64}{2} = 58,$\\
	$cos1(\frac{2r}{17}) = 29 = cos(\frac{12r}{17})$ of the table obtained
	using the program trig.bas.

	\ssec{Constructibility of the second degree.}

	\setcounter{subsubsection}{-1}
	\sssec{Introduction.}
	In this section we examine the problems which
	correspond or require the intersection of a conic or of a circle
	with a line when one of the intersections is not known.



\setcounter{section}{3}
	\section{Contrast with classical Euclidean Geometry.}
\setcounter{subsection}{-1}
	\ssec{Introduction.}
	To contrast the notions within finite Euclidean geometry
	with those of Euclidean geometry, with have the following summary:

	In finite Euclidean geometry (of the elliptic type),

	The following properties are different in finite Euclidean geometry:
	\begin{enumerate}
	\setcounter{enumi}{-1}
	\item There are $p$ points on each line.
	\item There are $p+1$ lines through every point.
	\item There are $p+1$ points and $p+1$ tangents for each circle.
	\item There are even and odd angles, the even ones can be bisected,
		the odd ones cannot.
	\item There are even triangles, for which there are 4 inscribed
		circles, the others have no inscribed circles.
	\item Angles can be expressed as integers, addition of angles is done
		modulo $p+1.$
	\item line through the center of a circle does not necessarily have
		an intersection with the circle.  The angle between any two
		lines, through the center, which have an intersection is even.
	\item Regular polygons exist only if the number of vertices is a
		divisor of $p+1.$
	\item Distances can be expressed either as integers or as integers
		times an irrational, the addition of distances on the same line
		is done by adding the integers modulo p.  The square of the
		irrational is an integer which is not a square.  For instance,
		for $p = 7,$ the irrational can be chosen to be $\sqrt{3}.$
	\item Trigonometric functions $sin$ and $cos$ can expressed like the
		distances.  The cosine of an even angle is always an integer.
		The cosine of an odd angle is an integer times an irrational.
		If $p \equiv 1 \pmod{4},$ the sine of an even angle is
		an integer that of an odd angles is an integer times an
		irrational, the reverse is true if $p \equiv -1 \pmod{4}.$
	\item Ordering cannot be introduced.  This is replaced by partial
		ordering.
	\end{enumerate}

	Among the properties which are similar, we have the following:
	Incidence, parallelism, equiangularity , equidistance,
	perpendicularity, congruence (of figures), similarity,
	the barycenter, the orthocenter,  the circumcircle,
	the theorem of Pythagoras.

	The constructibility of regular polygons (if they exist in the finite
	case), for instance
	if we replace the field ${\Bbb Z}_p$ by the field $({\Bbb Z}_p,$
	$\sqrt{2},$ we always have regular octogons, if we replace by the
	field $({\Bbb Z}_p,$ $\sqrt{p_i},$ 
	$p_i$ being all the primes, the constructible polygons always exist.
	If we replace the field $Z_p$ by the field of algebraic numbers, it is
	those which are roots of some polynomial, then all regular polygons
	exist.\\
	A similar discussion can be made if we start from the field ${\Bbb Q}$
	of the rationals.  With ${\Bbb Q}$ we can only construct squares,
	extended using \ldots.\\
	If ${\Bbb A}$ is the field of algebraic numbers, \ldots.\\
	The length of the circle as a limit of the length of polygons only make
	sense if we start with $Q.$  The implication of the transcendence of
	\ldots in $A$ is not a number in $A.$\\
	\ldots The field of algebraic numbers, the rational case and the
	existence or non existence of regular polygons.

{\tiny\footnotetext[1]{G39.TEX [MPAP], \today}}

\setcounter{section}{2}
	\section{Parabolic-Euclidean or Cartesian Geometry.}
	\setcounter{subsection}{-1}
	\ssec{Introduction.}
	The Euclidean geometry can be obtained from the projective geometry by
	choosing an appropriate set of elements, namely the ideal line and 2
	complex conjugate points on the line, the isotropic points.
	Similarly, for the hyperbolic geometry, one choose one real conic as the
	ideal.  In the Cartesian geometry, we choose a line, the ideal line and
	a point on that line, the isotropic point.  Definitions and properties
	in this geometry will be stated.  A construction, which allows an
	elementary algebraic proof of the properties, will be given.  The
	transformations leaving invariant the equality of angles and distances
	will be studied in a model of the geometry in the Euclidean plane,
	giving a justification for the name of the geometry.

	We start with a projective plane associated to an arbitrary field.
	A specific line, $i,$ in that plane is chosen, called the ideal line.
	A specific point, $I,$ on that line is also chosen, called the isotropic
	point.  Because a line is chosen, we can use all the concepts of affine
	geometry.  In particular, 2 lines are called parallel if they have the
	same ideal point.  The mid-point of 2 points $A,$ $B$ is the harmonic
	conjugate with respect to $A,$ $B$ of the ideal point on $A \times B.$
	A parabola is a conic tangent to the ideal line.

	I now will define new concepts in the Cartesian geometry.
	To focus the attention on a specific set of properties, I have chosen
	properties which have been inspired by those associated to the geometry
	of the triangle.  Because we want properties which are true in any
	field, it is not appropriate to derive them by a limiting procedure.
	I have therefore stated and proven them independently from the
	corresponding properties in Euclidean geometry and have indicated the
	correspondence by giving the same name as that of the corresponding
	element in Euclidean geometry, but without giving the justification.

	The configuration should give theorems in 2 ways \footnote{2.3.83},
	using \ldots.\\
	The equality of angles is associated with a parabolic
	projectivity (with 2 coincident fixed points).\\
	Recall the construction of a parabolic projectivity on $i,$
	let $I1,$ $I2$ be a pair, the point $I4,$ corresponding to $I3,$
	is obtained as follows,\\
	given $A,$ choose $B$ on $A \times I,$\\
	$C := (A \times I2) \times (B \times I1),$ $c := I \times C,$\\
	$D := (B \times I3) \times c,$ $a := A \times D,$ then\\
	$I4 := a \times i.$

	\sssec{Measure of distances and of angles.}
	The measure of angles and distances play a fundamental role in the
	geometry of Euclid and in the study of non Euclidean geometry by
	de Tilly.  On the other hand, starting from projective geometry, these
	notions are derived notions.  The appropriate definitions for the
	measure of distances and of angles will be given first in the case of a
	real field using a model on the Euclidean plane with given
	perpendicular axis $x$ and $y$ through a point $O$ and the line $l$ with
	equation $y = 1.$  This will justify the name of Cartesian
	geometry.\footnote{6.3.83}

	\ssec{Fundamental Definitions.}

	\sssec{Definition.}
	In in affine geometry let us choose one point on the
	ideal line as a double isotropic point.  This point will be called
	the {\em isotropic point} or {\em sun}.  The associated geometry will
	be called
	{\em parabolic-euclidean} or {\em Cartesian}.

	\sssec{Definition.}
	Any line through the sun is called an {\em isotropic line} or
	{\em solar axis}.  A parabola with the sun as ideal point is called a
	{\em parcircle}.

	\sssec{Comment.}
	There is a configuration which is a special case of the
	hexal configuration which allows the study of the geometry of the
	triangle in the Cartesian geometry.  Indeed it is sufficient to
	choose $\overline{M}$ to be the isotropic point.

	\sssec{Example.}
	$x_{i+1} := x_i + 1 \umod{p}$ is such a projectivity.\\
	With $p = 5,$\\
	$X = \{0, 1, 2, 3, 4, 0,  $\ldots$  \}$\\.
	Hence both angles and distances can be represented by an integer
	modulo $p.$\\
	The circle is replaced by a parcircle,

	\sssec{Notation.}
	Let $\pi$  be the parabolic projectivity with the ideal point
	as fixed point:\\
	0.\hti{6}$pi  = \{(x,x+1)\}.$

	\sssec{Theorem.}
	{\em
	\hth$pi ^i = \{(x,x+i)\},$ $i \in {\Bbb R},$}\\
	therefore, if the coordinates of the point $P$ are $x_P$ and $y_P$ and
	if the parallels to the lines $a$ and $b$ through $O$ meet $l$ at $A$
	and $B,$ we are justified to give the following

	\sssec{Definition.}
	\hth$	dist(P,Q) = y_Q - y_P.$\\
	\hth$	angle(a,b) := x_B - x_A.$

	We have
	\sssec{Theorem.}
	{\em Given any 3 points $A,$ $B$ and $C$ and any 3 lines $a,$ $b$ and
	 $c,$\\
	\hth$	dist(A,C) + dist(C,B) + dist(B,A) = 0,$\\
	\hth$	angle(a,c) + angle(c,b) + angle(b,a) = 0.$}

	\sssec{Comment.}
	Because, in both instances, the notion of measure are
	associated with the coordinates of a 1 dimensional set of points,
	both measure of distances and of angles can be given a sign.  Of course
	what we obtain is not a metric but a semi metric because\\
	\hth$	 |dist(A,B)| \leq |dist(A,C)| + |dist(C,B)|,$\\
	but $dist(P,Q) = 0$ does not imply $P = Q$ but only $y_P = y_Q.$

	In the case of a Gaussian field, if [0,0,1] is the ideal line and
	(0,1,0) is the isotropic point, we can give

	\sssec{Definition.}
	\hth$	dist( (A0,A1,1),(B0,B1,1) ) := B1 - A1 \umod{p^k},$\\
	\hth$	angle( [a0,-1,a2],[b0,-1,b2] ) := b0 - a0 \umod{p^k},$

	\sssec{Theorem.}
	{\em The set of points $Q$ such that angle $AQB$ is constant is a
	parcircle\footnote{3.3.83}.}
	
	\sssec{Theorem.}
	{\em The set of points $Q$ such that angle $QA0$ = angle $AQO$
	isosceles triangle) is a set of 2 lines (which can coincide),
	$A \times S$ and  $B \times S,$
	such that $O$ is equidistant in the Euclidean sense from these 2 lines.}

	Proof:\\
	\hth$A' := (B \times O) \times (A \times S),$
	$ B' := (A \times O) \times (B \times S),$\\
	\hth$I1 := i \times (B \times O),$
	$ I2 := i \times (A \times B) = i \times (A' \times B'),$\\
	\hth$I3 := i \times (A \times O),$\\
	then $(I1,I2) = (I2,I3)$ and $ABO$ is an isosceles triangle.
	$(Ii,Ij)$ denotes the angle of any pair of lines through $Ii$ and $Ij$
	 on $i.$\\
	Any other point $D$ on $B \times S$ is such that $ADO$ is an isosceles
	 triangle.

	\sssec{Definition.}
	$AO = BO$ either if $A,$ $B$ and $I$ are collinear, or if
	$A'B'$ is parallel to $AB$ where\\
	$A' = (B \times O) \times (A \times S) and B' = (A \times O) \times
	(B \times S).$

	\sssec{Definition.}  Two lines are {\em antiparallel} if
	\ldots.

	\ssec{The Geometry of the Triangle in Galilean Geometry.}

	\sssec{Definition.}
	A line in a triangle is a {\em symmedian} if \ldots.

	\sssec{Comment.}
	We could also define measure of distance and angles
	dualy\footnote{4.3.83}.\\
	$(A,B) = (C,D)$ if $( (I \times A),(I \times B) ) and ( (I \times C),
	(I \times D) )$
	are corresponding pairs in an parabolic projectivity with fixed line
	$i.$\\
	If we use as model in the Euclidean plane the line at infinity as
	the ideal line and the point in the direction of the $x$ axis as the sun
	all points on a line through the sun are equidistant from points on
	another line through the sun the measure of distances between the points
	can be chosen as the measure of the distances between the lines.
	Therefore, the distance between $A := (a0,a1,1)$ and $B := (b0,b1,1)$ is
	$b1-a1.$\\
	The angle between $a := [1,-a1,0]$ and $b := [1,-a2,0]$ is $a2-a1.$
	a corresponds to $X = a1 Y,$ if $\alpha$  is the angle with the $y$
	 axis in
	the clockwise direction, $tan(\alpha) = a1,$ the "sun" angle is
	doubled if the tangent is doubled.

	\sssec{Definition.}
	The {\em line of Euler} is \ldots.

	\sssec{Definition.}
	The {\em circumparcircle}

	\sssec{Definition.}
	The {\em first circle of Lemoine.}

	\sssec{Definition.}
	The {\em second circle of Lemoine.}

	\sssec{Definition.}
	The {\em circle of Brocard.}

	The Brianchon-Poncelet-Feurbach theorem becomes\footnote{3.3.83}

	\sssec{Theorem.}
	{\em Given a triangle $\{A_i,a_i\}$ and the parcircle $\iota$  tangent
	to $a_i.$ Let $M$ be any point not on the side of the triangle or on
	$i,$ Let $M_i := (M \times A_i) \times a_i,$ the parcircle
	$\gamma$  through $M_i$ is
	tangent to the parcircle $\iota$ .\\
	By duality, let $m$ be a line not through $A_i$ or $I,$
	let $m_i := (m \times a_i) \times A_i,$ the parcircle tangent to
	$m_i$ is tangent to the parcircle $\iota$.}

	\ssec{The symmetric functions.}
	\sssec{Theorem.}
	{\em The symmetric functions can be expressed in terms of $s_{11}$ and
	$s_{111},$.  More precisely}
	\begin{enumerate}
	\setcounter{enumi}{-1}
	\item[H0.]$s_1 = 0,$ $b := s_{11},$ $c := s_{111},$\\
	{\em then}
	\item[C2.]$s_2 = - 2 b,$
	\item[C3]$s_{21} = - 3 c,$ $s3 = 3 c,$
	\item[C4]$s_{22} = b^2,$ $s_{31} = -2 b^2,$ $ s_4 = 2 b^2,$ 
		$s_{211} = 0,$
	\item[C5]$s_{221} = bc,$ $s_{32} = - bc,$ $ s_{311} = -2 bc,$
		$s_{41} = 5 bc,$ $s_5 = -5 bc,$
	\item[C6]$s_{222} = c^2,$ $s_{33} = b^3 + 3 c^2,$ 
		$s_{321} = -3 c^2,$ $s_{411} = 3 c^2,$
		$s_{42} = -2 b^3 - 3 c^2,$ $s_{51} = 2 b^3 - 3 c^2,$
		$s_6 = -2 b^3 + 3 c^2.$
	\item[C7]$s_{322} = 0,$ $s_{331} = b^2c,$ $s_{421} = - 2 b^2c,$
		$s_{43} = - b^2c,$
		$s_{511} = 2 b^2c,$ $s_{52} = 3 b^2c,$ $s_{61} = - 7 b^2c,$
		$s_7 = 7 b^2c.$
	\item[C8]$s_{332} = bc^2,$ $s_{422} = - 2 bc^2,$
		$s_{431} = - bc^2,$ $s_{44} = b^4 + 4 bc^2,$\\
		$s_{521} = 5 bc^2,$ $s_{53} = - 2 b^4 - 7 bc^2,$
		$s_{611} = - 5 bc^2,$
		$s_{62} = 2 b^4 + 2 bc^2,$ $s_{71} = - 2 b^4 + 8 bc^2,$
		$s_8 = 2 b^4 - 8 bc^2,$
	\item[C9]$s_{333} = c^3,$ $s_{432} = -3 c^3,$
		$s_{441} = b^3c + 3 c^3,$ $s_{522} = 3 c^3,$
		$s_{531} = -2 b^3c - 3 c^3,$ $s_{54} = - b^3c - 3 c^3,$
		$s_{63} = 3 b^3c + 6 c^3,$ $s_{621} = 2 b^3c - 3 c^3,$
		$s_{72} = - 5 b^3c - 3 c^3,$ $s_{711} = - 2 b^3c + 3 c^3,$
		$s_{81} = 9 b^3c - 3 c^3,$ $s_9 = - 9 b^3c + 3 c^3,$
	\item[C10]$s_{433} = 0,$ $s_{442} = b^2c^2,$ 
		$s_{532} = - 2 b^2c^2,$ $s_{541} = - b^2c^2,$\\
		$s_{55} = b^5 + 5 b^2c^2,$ $s_{622} = 2 b^2c^2,$
		$s_{631} = 3 b^2c^2,$
		$s_{64} = - 2 b^5 - 9 b^2c^2,$ $s_{721} = - 7 b^2c^2,$
		$s_{73} = 2 b^5 + 6 b^2c^2,$
		$s_{811} = 7 b^2c^2,$ $s_{82} = - 2 b^5 + b^2c^2,$
		$s_{91} = 2 b^5 - 15 b^2c^2,$
		$s_{10} = - 2 b^5 + 6 b^2c^2,$
	\item[C11]$s_{443} = bc^3,$ $s_{533} = - 2 bc^3,$
		$s_{542} = - bc^3,$ $s_{551} = bc^4 + 4 bc^3,$
		$s_{632} = 5 bc^3,$ $s_{641} = - 2 bc^4 - 7 bc^3,$
		$s_{65} = - bc^4 - 4 bc^3,$
		$s_{722} = - 5 bc^3,$ $s_{731} = 2 bc^4 + 2 bc^3,$ 
		$s_{74} = 3 bc^4 + 11 bc^3,$
		$s_{821} = - 2 bc^4 + 8 bc^3,$ $s_{83} = - 5 bc^4 - 13 bc^3,$
		$s_{911} = 2 bc^4 - 8 bc^3,$ $s_{92} = 7 bc^4 + 5 bc^3,$
		$s_{10,1} = - 11 bc^4 + 11 bc^3,$ $s_{11} = 11 bc^4 - 11 bc^3,$
	\item[C12]$s_{444} = c^4,$ $s_{543} = -3 c^4,$
		$s_{552} = b^3c^2 + 3 c^4,$ $s_{633} = 3 c^4,$
		$s_{642} = - 2 b^3c^2 - 3 c^4,$ $s_{651} = - b^3c^2 - 3 c^4,$
		$s_{66} = b^6 + 6 b^3c^2 + 3 c^4,$ $s_{732} = 2 b^3c^2 - 3 c^4,$
		$s_{741} = 3 b^3c^2 + 6 c^4,$
		$s_{75} = - 2 b^6 - 11 b^3c^2 - 3 c^4,$
		$s_{822} = - 2 b^3c^2 + 3 c^4,$
		$s_{831} = - 5 b^3c^2 - 3 c^4,$
		$s_{84} = 2 b^6 + 8 b^3c^2 - 3 c^4,$
		$s_{921} = 9 b^3c^2 - 3 c^4,$
		$s_{93} = - 2 b^6 - 3 b^3c^2 + 6 c^4,$
		$s_{10,1,1} = - 9 b^3c^2 + 3 c^4,$
		$s_{10,2} = 2 b^6 - 6 b^3c^2 - 3 c^4,$ 
		$s_{11,1} = - 2 b^6 + 24 b^3c^2 - 3 c^4,$
		$s_{12} = 2 b^6 - 24 b^3c^2 + 3 c^4.$
	\end{enumerate}

	\sssec{Theorem.}
	{\em Given a triangle $A_i$ a point $M$ not on the sides of the triangle
	and a point $\overline{M}$ on the polar $m$ of $M$ with respect to the
	triangle.}
	\begin{enumerate}
	\setcounter{enumi}{-1}
	\item{\em $\gamma$ is a parcircle,}
	\item{\em $\theta$ is a parcircle,}
	\item{\em $\chi 1_i$ and $\chi 2_i$ are parcircles.}
	\end{enumerate}

	Proof:
	We will use the abbreviation\\
	\hth$	s11 = m1m2 + m2m0 + m0m1$ and we have\\
	\hth$m0^2 - m1m2 = m1^2 - m2m0 = m2^2 - m0m1$\\
	\hth$	= m1^2 + m2^2 + m1m2 = -s_{11}.$\\
	$s_{11} + m1m2 = -(m1^2+m2^2),$  \ldots.
	\hth$s_{11} + m0^2 = m1m2,$  \ldots.\\
	$2s_{11} - m1m2 = m1m2-2m0^2,$ $(m1-m2)^2 = -(s_{11} + 3 m1m2).$

	COMPARE $Mmm$ and $ji\overline{a},$ $Mm\overline{m}$ and $jia$.

	\sssec{Theorem.}
	{\em The conic\\
	\hth$	m0 X_1 X_2 + m1 X_2 X_0 + m2 X_0 X_1 = 0$\\
	passes through $M,$ $A_i,$ and\\
	\hth$	ZZ_i = (m0,-(m1-m2),m1-m2),$\\
	the tangent at $A_i$ is $\overline{a}c_i,$ \\
	the tangent at $M$ is $mai = [m0,m1,m2],$\\
	the tangent at $ZZ_i$ is $[(m1-m2)^2,m0m1,m2m0]$\footnote{1.3.83}.}

\setcounter{section}{4}

	\section{Transformation associated to the Cartesian geometry.}

\setcounter{subsection}{-1}
	\ssec{Introduction.}
	Such transformation must preserve measure of angles and distances.

	\sssec{Theorem.}
	The transformations associated to the Cartesian geometry
	are represented by unit upper triangular matrices in the Euclidean
	The following are subgroups of these transformation
	The translations\\
	\hth$\matb{(}{ccc}1&0&v\\0&1&w\\0&0&1\mate{)}.$\\
	The shears\\
	\hth$\matb{(}{ccc}1&u&v\\0&1&0\\0&0&1\mate{)},$\\
	and the special shears\\
	\hth$\matb{(}{ccc}1&u&0\\0&1&0\\0&0&1\mate{)}$.

	Indeed,

	\sssec{Definition.}

	\sssec{Theorem.}
	Given a point $P,$ the set of points whose polars with
	respect to a triangle pass through $P$ are on a conic through the
	vertices of the triangle and vice-versa.

	Proof.  The pole of $[q_0,q_1,q_2]$ is $(q_1q_2,q_2q_0,q_0q_1).$
	It is on line
	$(a_0,a_1,a_2)$ if $a_0q_1q_2 + a_1q_2q_0 + a_2q_0q_1 = 0.$

	\sssec{Definition.}
	Given a triangle, the point of Lemoine of a conic
	through the vertices of the triangle is the point $P$ of the preceding
	Theorem.

	\sssec{Definition.}
	Given a triangle, the line of Lemoine of a conic tangent
	to the triangle is the set of points whose polars with respect to the
	triangle are tangent to the conic.

	\sssec{Corollary.}
	The point of Lemoine of the circumcircle is the classical
	point of Lemoine.  The point of Lemoine of the conic of Simsom
	$m_0m_1 X_1 X_2 +  \ldots  = 0$ is $T^{mm} = (m_0m_1,m_1m_2,m_2m_0),$\\
	The line of Lemoine of the inscribed conic is $[j_1j_2,j_2j_0,j_0j_1],$
	it is the line through $Ja_i.$ 


	\ssec{The geometry of the triangle, the standard form.}

	\setcounter{subsubsection}{-1}
	\sssec{Introduction.}
	In this section, we do give only a representative set of
	Theorems using a form similar to that found in works on Geometry
	since Euclid.  Many more Theorems can be deduced from the compact
	form given in section 9.5.  The vertices of the triangle are denoted by
	$A_0,$	$A_1,$ $A_2,$ its sides by $a_0,$ $a_1,$ $a_2.$ 

	\sssec{Definition.}
	A Fano line $p$ of a point $P$ is the line through the
	diagonal points of the quadrangle $A_0,$ $A_1,$ $A_2,$ P.

	\sssec{Definition.}
	The cocenter M of a triangle is the Fano point
	of the ideal line $m.$  (D0.1., .2., .12.)

	\sssec{Construction.}
	Given a triangle $A_i,$ the ideal line m and the center M~,
	we can obtain a conic as follows.  The tangents to the conic $\theta$
	at the vertices of the triangle are be constructed using $m~_i$ = M~ x $A_i.$ 
	Any point on the conic and on a given line through one of its points,
	are be obtained using the construction of Pascal.

	\sssec{Definition.}
	The conic $$\theta$$  constructed in 9.6.3. is by definition a
	{\em circle} the {\em circumcircle} of the triangle.  (D1.12., H1.0.)

	\sssec{Definition.}
	The Euler line of a triangle is the line $eul$ through the
	cocenter and the center of the triangle.  (D1.0.)

	\sssec{Definition.}
	The central parallels $kk_i$ are the lines parallel to
	the sides of a triangle passing through the center $M~.$  (D1.1.)

	\sssec{Theorem.}
	{\em The central parallels $kk_i$ intersect the sides $a_{i+1}$ and

	$a_{i-1}$ at points $K~A~_{i-1}$ and $K~A~_{i+1}$ which are on a circle $\lambda$ .  (D1.2.,
	D2.11., C2.1., C2.2.)

	\sssec{Definition.}
	The circle $\lambda$  is called central circle.

	\sssec{Theorem.}
	{\em The circumcircle and the central circle are tangent at
	a point $LO.$  (C23.0.)

	\sssec{Definition.}
	The central points $M~_i$ of a triangle are the intersection
	of a tangent at a vertex with the opposite sides.  (D0.11.)

	\sssec{Definition.}
	The central line $m~$ of a triangle is the Fano line of its
	center $M~.$  (D0.12.)

	\sssec{Definition.}
	The associated circles $\alpha$ $_i$ are the circles through the
	center $M~$ of the triangle and its vertices $A_{i+1}$ and $A_{i-1}$.
	(D3.6., C3.1.)

	\sssec{Definition.}
	The center-vertex circles $\kappa$ $c_i$ and $\kappa$ $~c_i$ are the circles
	centered at one vertex of a triangle passing through an other vertex.
	(D4.12., C4.0.)

	\sssec{Definition.}
	The bissectrices $i_i$ of a triangle are the lines through a
	vertex $A_i$ such that the lines forms equal angles with the sides of the
	triangle passing through $A_i.$ 

	\sssec{Theorem.}
	{\em The bissectrices have a point I in common.  (9.5.5.,
	D0.3.)

	\sssec{Definition.}
	The bissector is the point I common to the bissectrices
	$i_i.$	 (D.0.3.)  The {\em bissector} line $i$ is the Fano line of the bissector.
	(D20.1)

	\sssec{Comment.}
	The sides of a triangle do not have a point in common,
	therefore, there is no circle tangent to its sides.

	\sssec{Definition.}
	The circles of Apollonius $\alpha$ $p_i$ are the circles centered at
	a central point $M~_i$ through a opposite vertex $A_i.$  (D5.12., C5.0.,
	C5.1.)  They have a common tangent with the circumcircle (C5.3.).  The
	point of contact with the side $a_i$ is on the bissectrix through $A_i.$ 
	(C22.1.)

	\sssec{Theorem.}
	{\em The circles of Apollonius have the same radical axis,
	l~mm, which is the common tangent of the circumcircle and the central
	circle.  (C6.2., C2.6.)

	\sssec{Definition.}
	The sun $MI$ is the direction of the bissector line.
	(D24.0.)

	\sssec{Definition.}
	Any parabola, i.e. a conic tangent to the ideal line,
	whose ideal point is the sun $MI$ is called a {\em sun-parabola}.

	\sssec{Comment.}
	In the isotropic geometry, the center of a parabola is an
	ideal point which is not necessarily its ideal point.

	\sssec{Definition.}
	The center-cocenter conic $\gamma$  is the conic through
	the vertices of the triangle, its center and its cocenter.  (D7.10.)

	\sssec{Theorem.}
	{\em The center-cocenter conic is a sun-parabola.  (C24.1.)

	\sssec{Definition.}
	The tangential circles $\chi$ $t_i$ $(\chi$ $t~_i)$ are the circles
	tangent to $a_{i+1} (a_{i-1})$ at $A_{i-1} (A_{i+1})$ passing through $A_{i+1} (A_{i-1})$.
	(D7.8., C7.0.)

	\sssec{Theorem.}
	{\em The other intersections $K~L_i$ and $K~L~_i$ with the tangential
	circles and the sides of the triangles are on a conic $\xi$  which is
	a sun-parabola.  (D3.1., D7.9., C7.2., D24.2.)

	\sssec{Theorem.}
	{\em The cocircumcircle $\theta$ ~ is the conic through the vertices of
	the triangle for which the tangents are parallel to the opposite side.
	(D1.12., D0.1., C1.0.)

	\sssec{Definition.}
	Let $Eul$ and $E~ul$ be the points of contact of the
	circumcircle and of the co-circumcircle with the line of Euler, the
	conic $\iota$  through these points and circumscribed to the triangle is
	called the {\em bissector conic}.  (D20.19., C20.2., C20.3.)

	\sssec{Theorem.}
	{\em The center of the bissector conic is the bissector point.
	(C20.7.)

	\sssec{Comment.}
	9.6.28. is an alternate definition fron that given in D20.19.

	\ssec{The cubic $\gamma$ a of Gabrielle.}

	\setcounter{subsubsection}{-1}
	\sssec{Introduction.}
	This section and the related section 11. was
	conceived after my daugther asked when I would name a Theorem for her.
	It concerns a general construction which starts from a parabola and
	constructs points on a cubic of which several are assoiated to the
	geometry of the triangle.

	\sssec{Definition.}
	Let x = (x0,x1,x2) be any line of the dual of the
	sun-parabola $\Gamma$ ,
	0.\hti{7}$(m1+m2) x1x2 + (m2+m0)x2x0 + (m0+m1) x0x1 = 0.$\\
	Let $k~k_i$ = [m1+m2,m2,m1].
	The following constructs points X = (X0,X1,X2) of the curve $\gamma$ a called
	the cubic of Gabrielle:
	D1.\hti{6}$X_i := x x k~k_i,$\\
	D2.\hti{6}$x_i := X_i x A_i,$\\
	D3.\hti{6}$X := x_1 x x_2.$\\
	D4.\hti{6}$\gamma a := \{X\}.$\\

	\sssec{Definition.}
	A parametric representation of a curve, with constraint
	arbitrary point are given in terms of 3 homogeneous parameters
	subjected to an homogeneous relation R between these 3 parameters.

	\sssec{Theorem.}
	{\em The curve $\gamma$ a is a point cubic, with axis
	0.\hti{7}$df= [m0^3(m1+m2)^2,m1^3(m2+m0)^2,m2^3(m0+m1)^2].$\\
	It contains the points $A_i,$ $M~_i,$ M, M~, LM.
	A parametric representation, with constraint 9.7.1.0., is
	1.\hti{7}$(x1x2(m0(m1+m2)(x1+x2) + m1m2 x0),$\\
	\hth$	 x2x0(m1(m2+m0)(x2+x0) + m2m0 x1),$\\
	\hth$	 x0x1(m2(m0+m1)(x0+x1) + m0m1 x2)),$\\
	Its equation in homogeneous coordinates is
	P4.\hti{6}$\gamma a: m0 X0(X1^2 + X2^2) + m1 X1(X2^2 + X0^2) + m2 X2(X0^2 + X1^2) = 0.$\\
	Proof:
	Definition 9.7.1. gives
	P1.\hti{6}$X_0 = (m1x1 + m2x2,m1x0 + (m1+m2)x2,m2x0 + (m1+m2)x1),$\\
	P2.\hti{6}$x_0 = [0,m2x0 + (m1+m2)x1,m1x0 + (m1+m2)x2],$\\
	P3.\hti{6}$X = ((m2x0+(m1+m2)x1)(m0x1+(m2+m0)x2),$\\
	\hth$		(m1x0+(m1+m2)x2)(m2x1+(m2+m0)x0),$\\
	\hth$		(m2x0+(m1+m2)x1)(m2x1+(m2+m0)x0)),$\\
	if we multiply all coordinates by x2 and use 9.7.1.0. we get 9.7.3.1.
	By a long algebraic verification it can be shown that the equation $P4.$
	is satisfied by 1.

	\sssec{Theorem.}
	{\em A parametric representation, with constraint
	0.\hti{7}$x0 + x1 + x2 = 0$\\
	is
	1.\hti{7}$(x1x2(m1x1+m2x2),x2x0(m2x2+m0x0),x0x1(m0x0+m1x1)).$\\
	2.\hti{7}$The point 1. is the point on the cubic \gamma a and the line$\\
	\hth$	[x0,x1,x2] through M distinct from M.$\\
	3.\hti{7}$There is one line mi where M is a triple point.$\\

	Proof:  Let (X0,X1,X2) be any point on the line [x0,x1,x2] passing
	through $M,$ eliminating $X0$ between the equation of the cubic $\gamma$ a and
	\hth$X0x0 + X1x1 + X2x2 = 0 gives$\\
	\hth$	(X1+X2)^2 (X1(m0x1(x1+x2)+m1x1^2) + X2(m0x2(x1+x2)+m2x2^2)) = 0$\\
	or
	\hth$	X1 = x0 x2(m0 x0 + m2x2), X2 = x0 x1(m0 x0 + m1x1),$\\
	because of 0., by symmetry we get 1.
	X1+X2 = 0, gives the point M and because (X1+X2) is a double factor
	this point has to be counted twice, hence 2., the point M is a {\em node}.
	The point 1. coincides with M  iff  the 3 coordinates are equal,
	the first 2 give after elimination of $x0,$ $(m0+m1)x1^2$ = $(m2+m0)x2^2,$ 
	hence $x_1$ = i2+i0, x2 = i0+i1, by symmetry x0 = i1+i2, hence 3.

	\sssec{Theorem.}
	{\em If the point of contact of a line
		[x0,x1,x2]
	through $LM$ with the cubic $\gamma$ a is (X0,X1,X2), then
	0.\hti{7}$X0^2 = m0x1x2(m1x1+m2x2), X1^2 = m1x2x0(m2x2+m0x0),$\\
	\hth$	X2^2 = m2x0x1(m0x0+m1x1).$\\
	We have
	1.\hti{7}$m0m1 x1 + m2m0 x1 + m0m1 x2 = 0$\\
	and
	2.\hti{7}$x0 X0 + x1 X1 + x2 X2 = 0.$\\
	Eliminating $X0$ between 2. and the equation of the cubic gives
	\hth$(m2X2+m1X1)(m2x1(m0x0+m1x1)X1^2 + m1x2(m0x0+m2x2)X2^2) = 0, because of 1.$\\
	The first factor corresponds to the point $LM,$ the other factor has a
	double root which gives 0.

	\sssec{Definition.}
	The cubic $\chi$  {\em of Charles} is the cubic through the points
	$M_i,$	$M~_i,$ $LM_i.$ 

	\sssec{Theorem.}
	{\em Let
	0.\hti{7}$a := m0m1m2,$\\
	1.\hti{7}$a0 := m0(m1^2+m2^2+m1m2), a1 := m1(m2^2+m0^2+m2m0),$\\
	\hth$	a2 := m2(m0^2+m1^2+m0m1),$\\
	then
	2.\hti{7}$\chi : a (X0^3 + X1^3 + X2^3)$\\
	\hth$		+ a0 X1X2(X1+X2) + a1 X2X0(X2+X0) + a2 X0X1(X0+X1) = 0.$\\
	3.\hti{7}$The tangent at (X0,X1,X2) is$\\
	\hth$	[a X0^2 + a2 X1^2 + a1 X2^2, a X1^2 + a0 X2^2 + a2 X0^2,$\\
	\hth$		a X2^2 + a1 X0^2 + a0 X1^2],$\\
	4.\hti{7}$The other point (Y0,Y1,Y2) on the tangent [x0,x1,x2] at$\\
	\hth$	(X0,X1,X2), is obtained by eliminating Z0 from 2., where$\\
	\hth$	(X0, X1, X2) is replaced by (Z0, Z1, Z2) and$\\
	\hth$	x0 Z0 + x1 Z1 + x2 Z2 = 0.  The coefficient of Z2^3 is$\\
	\hth$	Y1 X1^2 and that of Z1^3 is Y2 X2^2.$\\

	Proof:  For 4. we observe that the elimination should lead to the
	equation
	\hth$	(X2 Z1 + X1 Z2)^2 (Y2 Z1 + Y1 Z2) = 0.$\\
	An illustration of 4. is given by 12.4.

	\sssec{Conjecture.}
	Given 9 points $A_i,$ $B_i,$ $C_i,$ on a cubic such that $A_i,$ $B_i$ 
	$C_i$	and $(A_0,$ $B_0,$ $C_0),$ $(A_1,$ $B_1,$ $C_1)$ are collinear, then
	$(A_2,$	$B_2,$ $C_2)$ are collinear.

	\sssec{Corollary.}
	If 3 points $A_i$ are on a cubic, the third point $C_i$ on the
	tangent to the cubic at $A_i$ are also collinear.

	\sssec{Example.}
	For q = 16, i0 = 1, i1 = $\epsilon$ $^8,$ i2 = $\epsilon$ $^3,$ 
	\hth$	m0 = 1, m1 = \epsilon , m2 = \epsilon ^6,$\\
	H0.0.	M = 253,		E0.10.	M~ = 184,
	H0.1.\hti{4}$A_i = 2,1,0,$\\
	H0.2.\hti{4}$I = 130,$\\

	E0.0.\hti{4}$a_i = 0^*,1^*,272^*,$\\
	E0.1.\hti{4}$m_i = 253^*,2^*,136^*,      E0.9.   m~_i = 179^*,89^*,90^*,$\\
	E0.2.\hti{4}$M_i = 136,272,137,       E0.11.  M~_i = 15,115,91,$\\
	E0.3.\hti{4}$i_i = 125^*,6^*,233^*,$\\
	E0.4.\hti{4}$I_i = 238,232,234,$\\
	E0.5.\hti{4}$im_i = 44^*,102^*,339^*,$\\
	\hth$	im~_i = 4^*,151^*,168^*,$\\
	E0.6.\hti{4}$Tm_i = 194,14,16,$\\
	\hth$	Tm~_i = 3,30,195,$\\
	E0.7.\hti{4}$tm_i = 61^*,180^*,15^*,$\\
	\hth$	tm~_i = 271^*,203^*,194^*,$\\
	E0.8.\hti{4}$IA_i = 94,240,133,$\\
	\hth$	IA~_i = 101,215,183,$\\
	E0.12.\hti{3}$m = 137^*,			m~ = 189^*,$\\

	E1.0.\hti{4}$eul = 20,$\\
	E1.1.\hti{4}$kk_i = 152^*,205^*,96^*,            kk~_i = 258^*,21^*,147^*,$\\
	E1.2.\hti{4}$KA_i = 234,127,30,               K~A_i = 195,31,237,$\\
	\hth$	KA~_i = 126,16,255,              K~A~_i = 180,128,204,$\\
	E1.3.\hti{4}$kl_i = 203^*,233^*,236^*,           k~l_i = 126^*,194^*,202^*,$\\
	\hth$	kl~_i = 15^*,134^*,237^*,           k~l~_i = 127^*,29^*,14^*,$\\
	E1.4.\hti{4}$B_i = 147,61,102,                B~_i = 101,47,37,$\\
	E1.5.\hti{4}$bb_i =$\\
	E1.6.\hti{4}$Eul_i = 116,254,256,$\\
	E1.7.\hti{4}$Ba_i =$\\
	\hth$	Ba~_i =$\\
	E1.8.\hti{4}$tB_i =$\\
	E1.9.\hti{4}$KK_i = 147,5,200,                K~K_i = 101,199,143,$\\
	E1.10.\hti{3}$eul_i = 88^*,135^*,255^*,$\\
	E1.11.\hti{3}$TT =$\\

	E1.12.\hti{3}$\theta  =   0 ,  1 ,  2 , 40 , 61 ,102 ,123 ,145 ,147 ,$\\
	\hth$	     90^*, 89^*,179^*, 96^*,120^*, 92^*, 71^*, 49^*,216^*,$\\
	\hth$			172 ,197 ,211 ,223 ,229 ,235 ,249 ,250 ,$\\
	\hth$			104^*,270^*,152^*, 10^*,225^*, 20^*,205^*, 54^*,$\\

	E2.0.\hti{4}$tim_i = 118^*,45^*,16^*,            t~im_i = 182^*,101^*,140^*,$\\
	E2.1.\hti{4}$LI_i = 63,193,239,               L~I_i = 181,203,64,$\\
	E2.2.\hti{4}$li_i = 192^*,62^*,238^*,            l~i_i = 13^*,30^*,63^*,$\\
	E2.3.\hti{4}$Atm_i = 211,50,208,              A~tm_i = 151,217,242,$\\
	E2.4.\hti{4}$lt_i = 125^*,254^*,181^*,           l~t_i = 125^*,237^*,31^*,$\\
	E2.5.\hti{4}$LM = 155,			L~M = 163,$\\
	E2.6.\hti{4}$LT_i = 238,135,182,              L~T_i = 238,232,32,$\\
	E2.7.\hti{4}$lm = 98^*,                       l~m = 162^*,$\\
	E2.8.\hti{4}$LMM = 96,			L~MM = 174,$\\
	\hth$	LM~M = 118,			L~M~M = 265,$\\
	E2.9.\hti{4}$tKKL_i =$\\
	\hth$	tKKL~_i =$\\
	E2.10.\hti{3}$lmm = 83^*,                      l~mm = 92^*,$\\
	\hth$	lm~m = 49^*,                     l~m~m = 110^*,$\\

	E3.0.\hti{4}$ka_i = 255^*,61^*,203^*,            k~a_i = 233^*,88^*,126^*,$\\
	\hth$	ka~_i = 193^*,15^*,5^*,             ka~_i = 254^*,127^*,271^*,$\\

	E6.13.\hti{3}$\Gamma  =   0 ,  1 , 33 , 41 , 51 , 93 ,105 ,111 ,129 ,$\\
	\hth$	    116^*,254^*,161^*,235^*,253^*, 43^*, 11^*, 70^*,260^*,$\\
	\hth$			137 ,169 ,171 ,186 ,189 ,241 ,270 ,272 ,$\\
	\hth$			 96^*,107^*,218^*,268^*,174^*,213^*,184^*,256^*,$\\

	E7.1.\hti{4}$TMi = 171,			T~mi = 225,$\\

	E8.0.\hti{4}$dt_i = 139^*,227^*,239^*,$\\
	\hth$	dt~_i = 91^*,151^*,138^*,$\\
	E8.1.\hti{4}$Du_i = 116,6,16,$\\
	\hth$	Du~_i = 90,30,117,$\\

	E9.0.\hti{4}$Eb_i = 133,251,212,              E~b_i = 150,55,167,$\\
	\hth$	Eb~_i = 189,111,261,             E~b~_i = 9,257,146,$\\
	E9.2.\hti{4}$ed_i = 226^*,68^*,256^*,            e~d_i = 183^*,121^*,123^*,$\\
	\hth$	ed~_i = 128^*,44^*,144^*,           e~d~_i = 158^*,185^*,16^*,$\\

	E11.19.\hti{2}$Dh_i = 177,103,75,               D~h_i = 133,253,183,$\\
	\hth$	Dh~_i = 83,159,122,              D~h~_i = 253,34,22,$\\
	E11.20.\hti{2}$di_i = 180,127,13,               d~i_i = 180,3,253,$\\
	\hth$	di~_i = 204,5,193,               d~i~_i = 204,114,2,$\\
	E11.21.\hti{2}$Dj_i = 100,111,261,              D~j_i = 100,124,253,$\\

	E11.22.\hti{2}$dk_i = 27^*,83^*,35^*,              d~k_i = 92^*,71^*,216^*,$\\
	\hth$	dk~_i = 23^*,35^*,224^*,$\\
	E11.23.\hti{2}$du_i = 88^*,232^*,15^*,$\\
	\hth$	du~_i = 114^*,203^*,116^*,$\\
	E11.24.\hti{2}$Dl_i = 101,185,104,$\\
	\hth$	Dl~_i = 101,247,190,$\\
	E11.25.\hti{2}$Dm_i = 204,134,203,$\\
	\hth$	Dm~_i = 14,15,29,$\\
	E11.26.\hti{2}$Dn_i = 30,205,255,$\\
	\hth$	Dn~_i = 16,204,180,$\\
	E11.27.\hti{2}$dn = 33^*,$\\
	\hth$	dn~ = 100^*,$\\
	E11.28.\hti{2}$Do = 200,$\\
	\hth$	Do~ = 174,$\\
	E11.29.\hti{2}$dp = 102^*,$\\
	\hth$	dp~ = 55^*,$\\
	E11.30.\hti{2}$Dq = 178,$\\
	E11.31.\hti{2}$dr = 26^*,$\\

	E11.32.\hti{2}$\gamma a =   0 ,  1 ,  2 , 15 , 40 , 80 , 91 ,100 ,103 ,$\\
	\hth$	     181^*,254^*,125^*,121^*,234^*,208^*,245^*,133^*, 78^*,$\\
	\hth$			115 ,124 ,155 ,169 ,178 ,184 ,253 ,263 ,$\\
	\hth$			118^*,149^*, 39^*,119^*, 26^*, 49^*, 83^*,100^*,$\\

	\hth$	\gamma ~a =   0 ,  1 ,  2 , 53 ,100 ,111 ,136 ,137 ,153 ,$\\
	\hth$	       31^*,237^*,125^*,141^*,173^*, 70^*,200^*,226^*, 41^*,$\\
	\hth$			163 ,184 ,225 ,250 ,253 ,261 ,265 ,272 ,$\\
	\hth$			 18^*, 92^*,111^*,113^*,110^*, 75^*,246^*,117^*,$\\

	E12.0.\hti{3}$Na_i = 89,8,194,                 N~a_i = 272,205,204,$\\
	E12.1.\hti{3}$na_i = 7^*,61^*,115^*,              n~a_i = 204^*,29^*,2^*,$\\
	E12.2.\hti{3}$Nb_i = 248,129,55,               N~b_i = 29,70,251,$\\
	E12.3.\hti{3}$nc_i = 8^*,263^*,144^*,             n~c_i = 158^*,245^*,218^*,$\\
	E12.4.\hti{3}$lMM_i = 192^*,180^*,31^*,           l~MM_i = 114^*,30^*,7^*,$\\
	E12.5.\hti{3}$nd_i = 263^*,58^*,239^*,            n~d_i = 103^*,18^*,66^*,$\\
	E12.6.\hti{3}$Ne_i = 124,243,197,              N~e_i = 208,24,249,$\\
	E12.7.\hti{3}$nf_i = 80^*,157^*,257^*,$\\
	E12.8.\hti{3}$ng_i = 182^*,101^*,140^*,           n~g_i = 118^*,45^*,16^*,$\\
	E12.9.\hti{3}$Nh_i = 106,103,134,              N~h_i = 186,228,222,$\\
	E12.10.\hti{2}$lI = 170^*,                      l~I = 52^*,$\\
	E12.11.\hti{2}$Ni_i = 98,262,54,                N~i_i = 170,18,258,$\\
	\hth$	Ni~_i = 86,210,260,              N~i~_i = 218,35,41,$\\
	E12.12.\hti{2}$nj_i = 179^*,140^*,147^*,            n~j_i = 253^*,239^*,118^*,$\\
	E12.13.\hti{2}$nk_i = 32^*,139^*,200^*,            n~k_i = 221^*,164^*,223^*,$\\
	\hth$	nk~_i = 252^*,144^*,117^*,          n~k~_i = 213^*,48^*,158^*,$\\
	E12.14.\hti{2}$Nl_i = 95,134,189,               N~l_i = 83,124,186,$\\
	\hth$	Nl~_i = 157,136,157,             N~l~_i = 20,228,115,$\\
	E12.15.\hti{2}$nl = 99^*,                       n~l = 150^*,$\\
	\hth$	nl~ = 119^*,                     n~l~ = 161^*,$\\
	E12.16.\hti{2}$nm_i = 268^*,60^*,66^*,             n~m_i = 98^*,152^*,47^*,$\\
	E12.17.\hti{2}$np_i = 47^*,44^*,109^*,             n~p_i = 265^*,268^*,121^*,$\\
	\hth$	np~_i = 76^*,123^*,253^*,           n~p~_i = 140^*,116^*,76^*,$\\
	E12.18.\hti{2}$nq_i = 24^*,139^*,257^*,            n~q_i = 48^*,11^*,223^*,$\\
	\hth$	nq~_i = 146^*,198^*,78^*,           n~q~_i = 82^*,18^*,172^*,$\\
	E12.19.\hti{2}$Nr_i = 157,137,254,              N~r_i = 15,20,15,$\\
	\hth$	Nr~_i = 134,78,198,              N~r~_i = 198,13,83,$\\
	E12.20.\hti{2}$nr = 252^*,                      n~r = 218^*,$\\
	\hth$	nr~ = 202^*,                     n~r~ = 191^*,$\\

	E12.\hti{5}$.	ND_i = 256,186,13,$\\
	E12.\hti{5}$.$\\
	\hth$\chi  =  13 , 15 , 20 , 78 , 83 , 91 , 95 ,103 ,106 ,115 ,116 ,124 ,131 ,$\\
	\hth$    261^*,218^*, 96^*,195^*,197^*,245^*,193^*,  1^*, 10^*,158^*,157^*,131^*,269^*,$\\
	\hth$     15 , 20 ,137 ,272 , 91 ,222 , 83 ,103 ,228 ,116 ,131 ,157 ,198 ,$\\

	\hth$	134 ,136 ,137 ,157 ,186 ,189 ,198 ,222 ,228 ,254 ,256 ,272 ,$\\
	\hth$	146^*,144^*,263^*,119^*, 87^*, 49^*,165^*,178^*, 60^*,257^*, 80^*,  8^*,$\\
	\hth$	254 ,256 , 13 ,136 ,189 , 78 ,115 , 95 ,134 ,106 ,124 ,186 ,$\\
	\hth$	(these are the tangents and the other point on the tangent)$\\

	E19.4.\hti{3}$lI = 170^*,                      l~I = 52^*,$\\
	E19.5.\hti{3}$LJ = 11,$\\

	E22.0.\hti{3}$iA_i = 179^*,237^*,3^*,$\\
	E22.1.\hti{3}$IA = 230,$\\
	E22.2.\hti{3}$ab = 224^*,$\\
	E22.6.\hti{3}$Ia_i = 15,126,4,$\\
	E22.7.\hti{3}$ia = 112^*,$\\

	E25.0.\hti{3}$Dk = 176,$\\
	E25.1.\hti{3}$dl = 100^*,$\\

	Given the center C of a circle and one of its points A, the
	point $X$ on any line $x$ through $A$ (or $y$ through $C)$ can be obtained by
	construction $y$ through $C$ (or $x$ through $A)$ such that the angle XAC =
	angle BCX.
	The above as to be reviewed.  A construction of a point on a given
	tangent (or radius) follows.  There must be a simpler way.
	Let the given points be $A_0,$ $A_1,$ $A_2,$ let the center be C, let
	the radius-tangent be $t,$ Mt := t x m, find $A_4$ on the circle and
	$A_0$	x Mt, find $A_5$ on the circle and $A_1$ x Mt,
	let Y := $(A_0$ x $A_1)$ x $(A_4$ x $A_5),$ (C x Y) x t is the point of contact
	with the circle.
	To find the bissectrix of an angle $A_{-1}$ $A_0$ $A_1$ we use the above
	construction with the tangent-radius $C x (m x (A_{-1} x A_1)$ the
??	point X on the circle is also on the bissectrix.
	\sssec{Notation.}
	Angles and directions will be denoted by an upper case
	letter and a lower case letter underlined.

	\sssec{Theorem.}
	{\em If the angle of the direction of the sides $b_i$ is $nb_i,$ then
 	the angle of the direction of the tangent is $d_{i+1} + d_{i+2} - d_i {\mathit mod} q+1$.

	\sssec{Theorem.}
	{\em If the direction of $a_i$ is $a_i,$ the angles at $A_i$ are
		$A_i	= {\mathit a}_{i+1} - {\mathit a}_{i-1} {\mathit mod} q+1$.

	\sssec{Theorem.}
	{\em 
	0.\hti{7} The angle of the direction of the center of $\chi 1_i$ is\\
	\hth$	{\mathit c1}_i = {\mathit a}_i + {\mathit A}_{i-1}$,\\
	\hth$    that of the center of \chi 2_i is$\\
	\hth$	{\mathit c2}_i = {\mathit a}_i - {\mathit A}_{i+1}.$\\
	1.\hti{7} The center of $\chi 1_i$ is $(A_{i+1} x M   ) x a_{i+1},$\\
	\hth$				  {\mathit c1}_i$\\

	\hth$    that of \chi 2_i is  (A_{i-1} x M   ) x a_{i-1},$\\
	\hth$			    {\mathit c2}_i$\\

%

\section{Problems}.
	\ssec{Problems for Affine Geometry.}
	\sssec{Theorem.}
	{\em If $m$ is the ideal line, and
	$A = (A_0,A_1,A_2)$, $B = (B_0,B_1,B_2)$ then}
	\enumb
	\item {\em the mid-point $A+B$ of $A$ and $B$ is\\
	\hth$A + B = (m\cdot B)A + (m\cdot A)B$}.
	\item {\em the symmetric $2B-A$ of $A$ with respect to $B$ is\\
	\hth$2B - A = 2(m\cdot A)B-(m\cdot B)A$}.
	\enume

	\sssec{Theorem.}\label{sec-t020}
	{\em The mid-points of the diagonals of a complete quadrilateral
	are collinear.} D37.5, C37.15. (020, Chou and Schelter 1986, p. 18)

	\sssec{Definition.}\label{sec-d020}
	The line of the preceding Theorem is called the {\em mid-line} of the
	complete quadrilateral.

	\sssec{Theorem.}\label{sec-t025}
	{\em Given a complete 5-lines, the mid-lines of the 5 complete
	quadrilaterals obtained by suppressiong any of the 5 lines
	have a point in common.} (025, Chou and Schelter 1986, p. 19)

	\sssec{Theorem.}\label{sec-t049}
	{\em Given a triangle $A_i$, and a point $M_0$, let $M_j$ be the
	symmetric of $M_{j-1}$ with respect to $A_{j-1\pmod{3}}$, then}
	\enumb
	\item $M_{i+3}$ {\em is the symmetric of $M_i$ with respect to
		$MM_{i+1}$, vertex of the anticomplementary triangle
		of} $A_i$.
	\item $M_6 = M_0$.
	\item $M_i,M_{i+3},M_{i+1},M_{i+4}$ {\em are parallelograms.}
	\enume

	\ssec{Problems for Involutive Geometry.}
	\sssec{Theorem.}
	{\em The perpendicular direction to $(IX_0,IX_1,IX_2)$ is}
	$(m_0(m_1-m_2)I_0+m_0(m_1+m_2)(I_1-I_2),
		m_1(m_2-m_0)I_1+m_1(m_2+m_0)(I_2-I_0),
		m_2(m_0-m_1)I_2+m_2(m_0+m_1)(I_0-I_1))$.

	\sssec{Theorem. [Buterfly Theorem]}\label{sec-t041}
	{\em If a quadrangle is inscribed in a circle with cent $O$, then a
	diagonal point, $D$, is the midpoint of the intersection with the other
	sides of a perpendicular through $D$ to $O\times D$.(041, Chou
	(1984), p.269.)}

\setcounter{section}{89}
	\section{Answers to problems and miscellaneous notes.}
	\sssec{Answer to}
	\vspace{-18pt}\hspace{100pt}{\bf \ref{sec-t025}.}\\
	Let the lines be $a_i,\ov{m}$ and $\ov{m}' = [m'_0,m'_1,m'_2]$.\\
	The midlines are\\ 
	$l_4 = [s_1-2m_0,s_1-2m_1,s_1-2m_2],$\\
	$l_3 = [s'_1-2m'_0,s'_1-2m'_1,s'_1-2m'_2],$ with
	$s'_1 = m'_0+m'_1+m'_2$,\\
	$\ov{M}A_1+\ov{M}'A_2 = (m'_0(m_2-m_0)-m_0(m'_0-m'_1),
		m'_1(m_0-m_2),m_2(m'_0-m'_1)).$\\
	$\ov{M}'A_1+\ov{M}A_2 = (m_0(m'_2-m'_0)-m'_0(m_0-m_1),
		m_1(m'_0-m'_2),m'_2(m_0-m_1)).$\\
	$l_0 = [m_1m_2m'_0(s'_1-2m'_0)-m'_1m'_2m_0(s_2-2m_0),
		m'_2m_0(2m'_0(m_0-m_1)-m'_1(s_1-2m_1))
		-m_2m'_0(2m_0(m'_0-m'_1)-m_1(s'_1-2m'_1)),
		m'_1m_0(2m'_0(m_2-m_0)+m'_2(s_1-2m_2))
		-m_1m'_0(2m_0(m'_2-m'_0)+m_1(s'_1-2m'_2))],$\\
	The common point is\\
	$P = (m'_0(m_1-m_2)-m_0(m'_1-m'_2),m'_1(m_2-m_0)-m_1(m'_2-m'_0),
		m'_2(m_0-m_1)-m_2(m'_0-m'_1)).$

	\sssec{Answer to}
	\vspace{-18pt}\hspace{100pt}{\bf \ref{sec-t049}.}\\
	Let $M_0 = (m_0,m_1,m_2)$, with $m_0+m_1+m_2 = 1$.\\
	$M_1 = 2A_0-M_0 = (2-m_0,-m_1,-m_2)$,\\
	$M_2 = 2A_1-M_1 = (-2+m_0,2+m_1,m_2)$,\\
	$M_3 = 2A_2-M_2 = (2-m_0,-2-m_1,2-m_2)$,\\
	$M_4 = 2A_0-M_3 = (m_0,2+m_1,-2+m_2)$,\\
	$M_5 = 2A_1-M_4 = (-m_0,-m_1,2-m_2)$,\\
	$M_6 = 2A_2-M_5 = (m_0,m_1,m_2)$.

	\sssec{Answer to}
	\vspace{-18pt}\hspace{100pt}{\bf \ref{sec-t041}.}\\
	Let 3 of the points be $A_i$, let $D := (0,1,x)$, be on $a_0$,
	then the 4-th point is $(y,1,x)$, with
	$y = -\frac{m_0(m_1+m_2)x}{m_2(m_0+m_1)+m_1(m_2+m_0)x}$.
	$O = (m_1+m_2,m_2+m_0,m_0+m_1),$\\
	$D\times O = [(m_2+m_0)x-(m_0+m_1),-(m_1+m_2)x,(m_1+m_2)],$\\
	its direction is $((m_1+m_2)(m_2x+m_1),m_2(m_2+m_0)x-s_{11}-m_2m_0),m_1(m_0+m_1)+x(s_{11}+m_0m_1))$.\\
	The direction perpendicular to $D\times O$ is\\
	$(m_0(m_1-m_2)(m_1+m_2)(m_2x+m_1)
			+m_0(m_1+m_2)(m_2(m_2+m_0)x-s_{11}-m_2m_0)-
			m_1(m_0+m_1)+x(s_{11}+m_0m_1)),
		m_1(m_2-m_0)m_2(m_2+m_0)x-s_{11}-m_2m_0)
			+m_1(m_2+m_0)(m_1(m_0+m_1)+x(s_{11}+m_0m_1)
				-(m_1+m_2)(m_2x+m_1)),
		m_2(m_0-m_1)m_1(m_0+m_1)+x(s_{11}+m_0m_1)
			+m_2(m_0+m_1)((m_1+m_2)(m_2x+m_1)
				-m_2(m_2+m_0)x-s_{11}-m_2m_0)))$.\\
	\ldots.

\begin{figure}
{
\centering
\includegraphics[width=\linewidth]{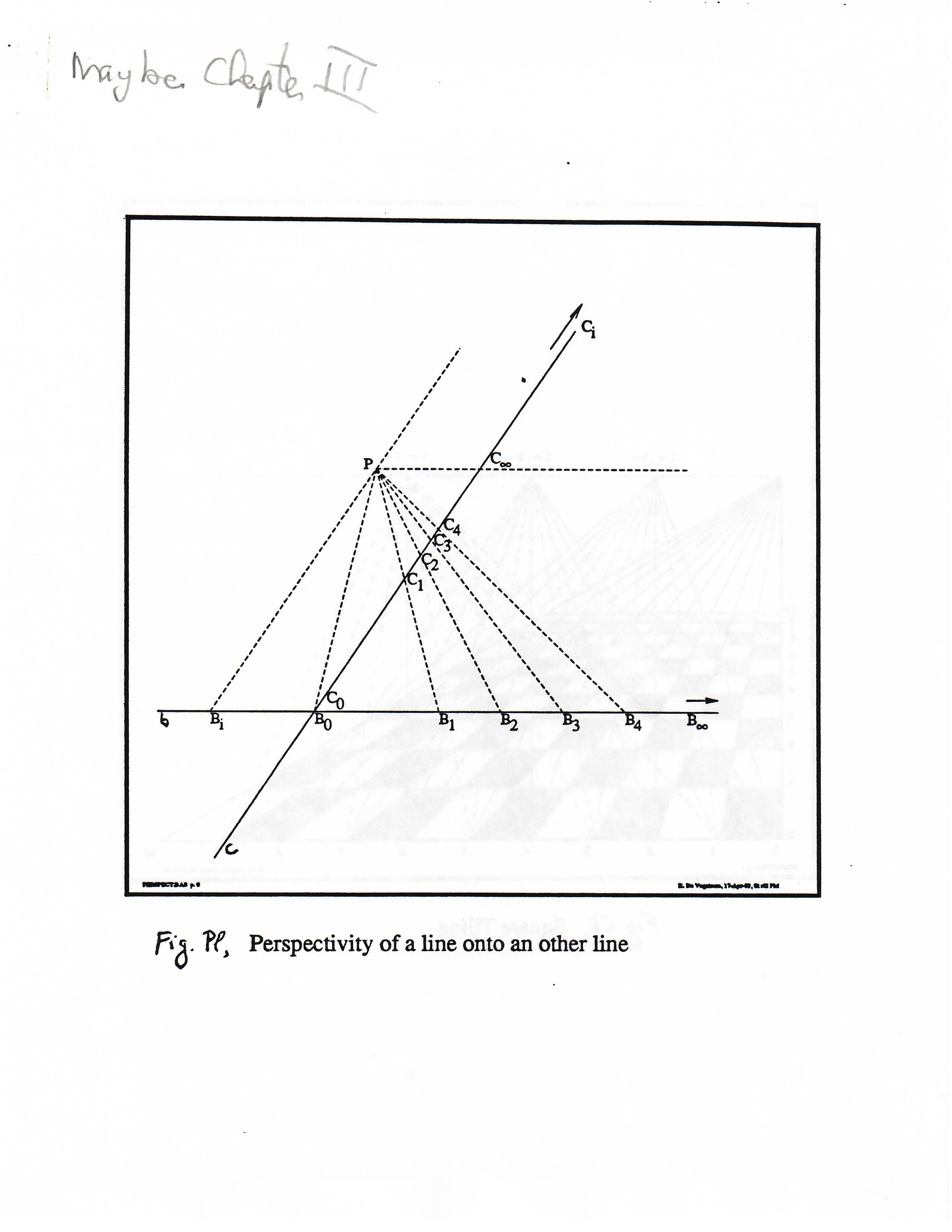}
}
\end{figure}

\begin{figure}
{
\centering
\includegraphics[width=\linewidth]{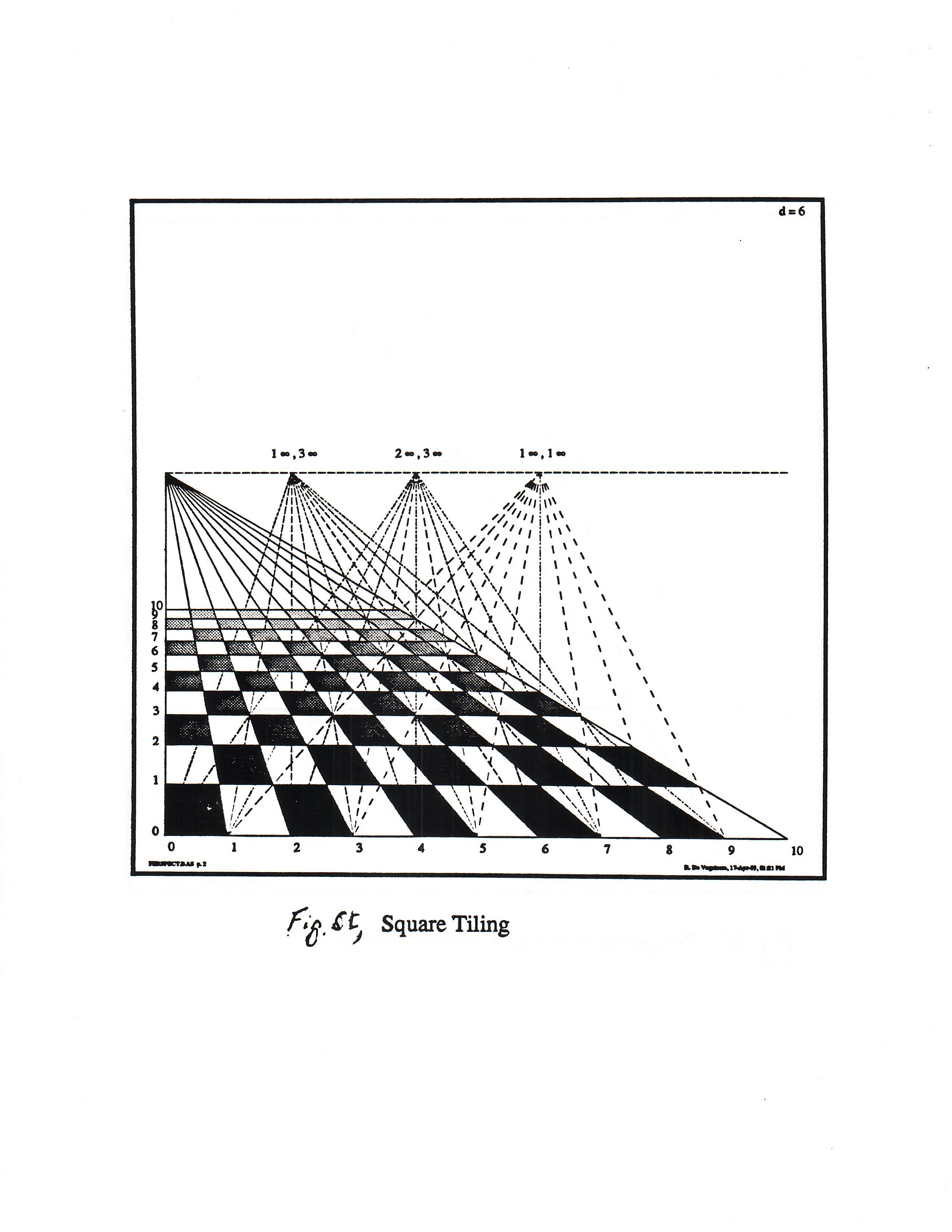}
}
\end{figure}

\begin{figure}
{
\centering
\includegraphics[width=\linewidth]{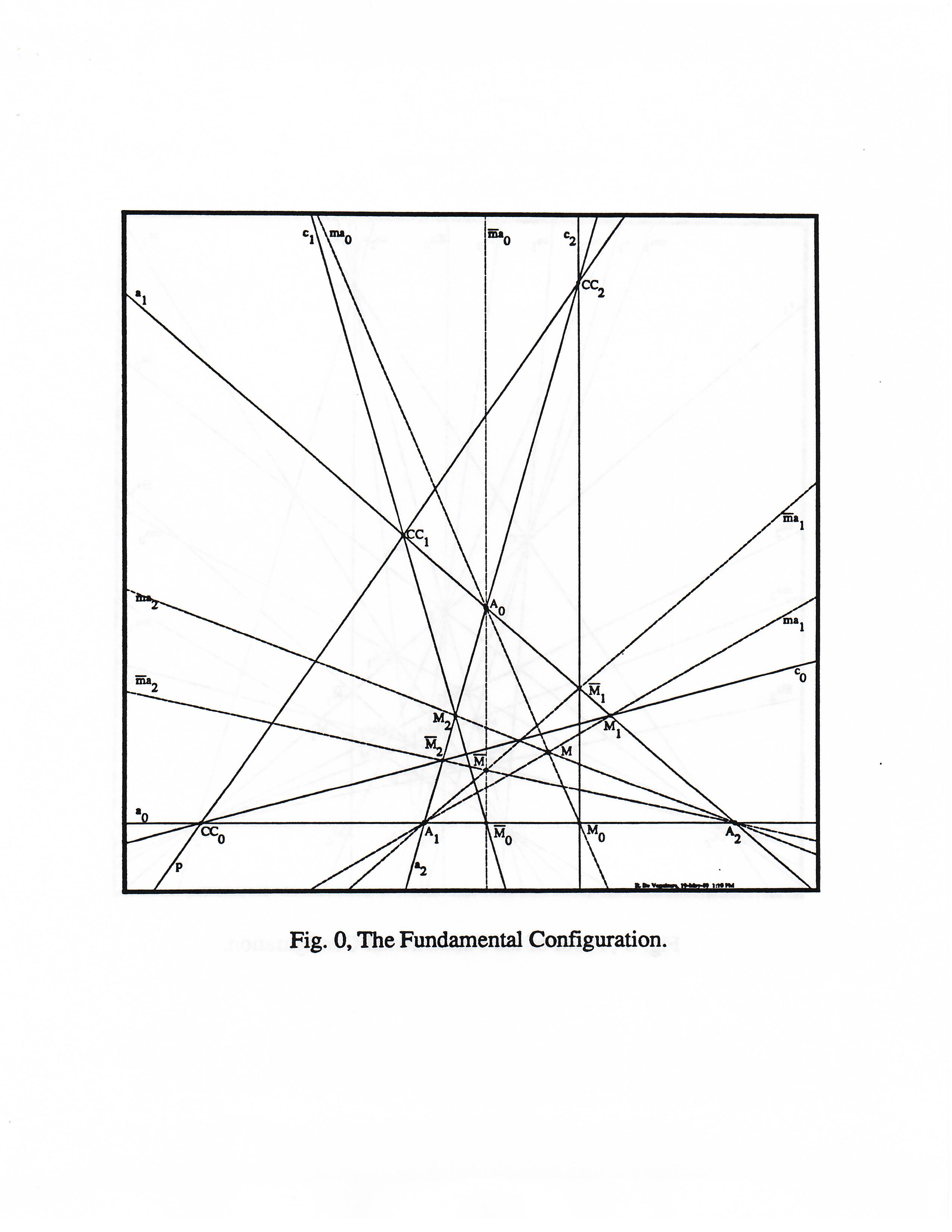}
}
\end{figure}

\begin{figure}
{
\centering
\includegraphics[width=\linewidth]{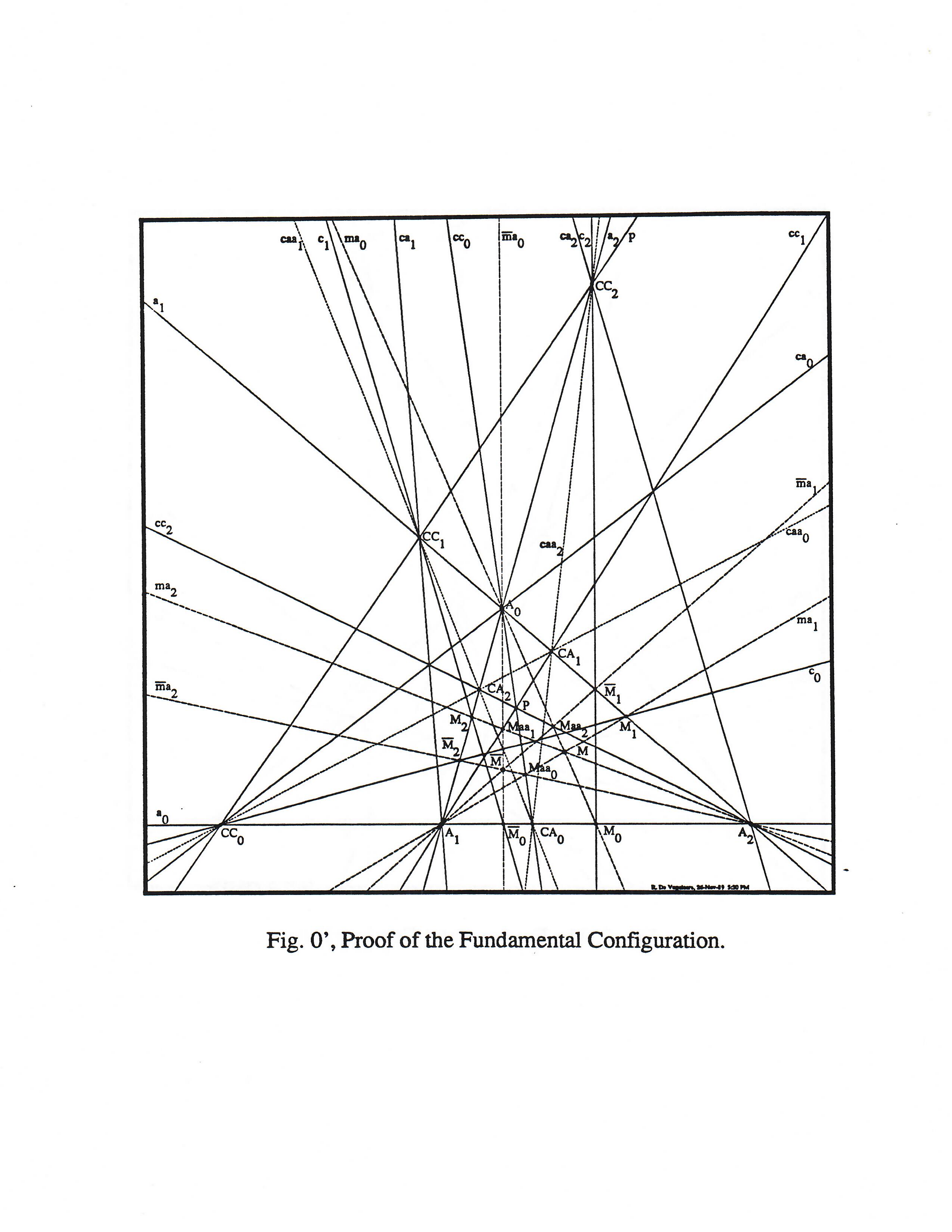}
}
\end{figure}

\begin{figure}
{
\centering
\includegraphics[width=\linewidth]{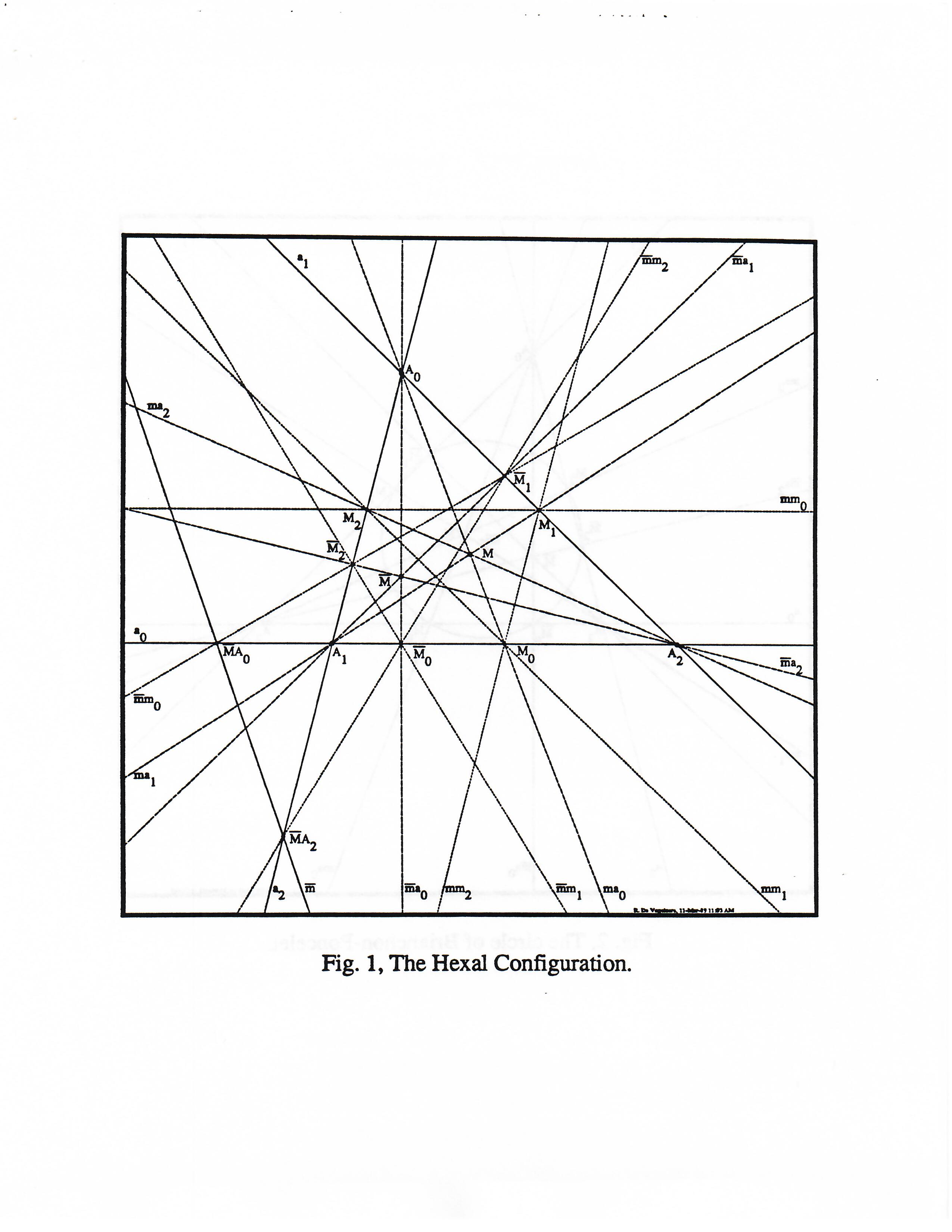}
}
\end{figure}

\begin{figure}
{\centering\includegraphics[width=\linewidth]{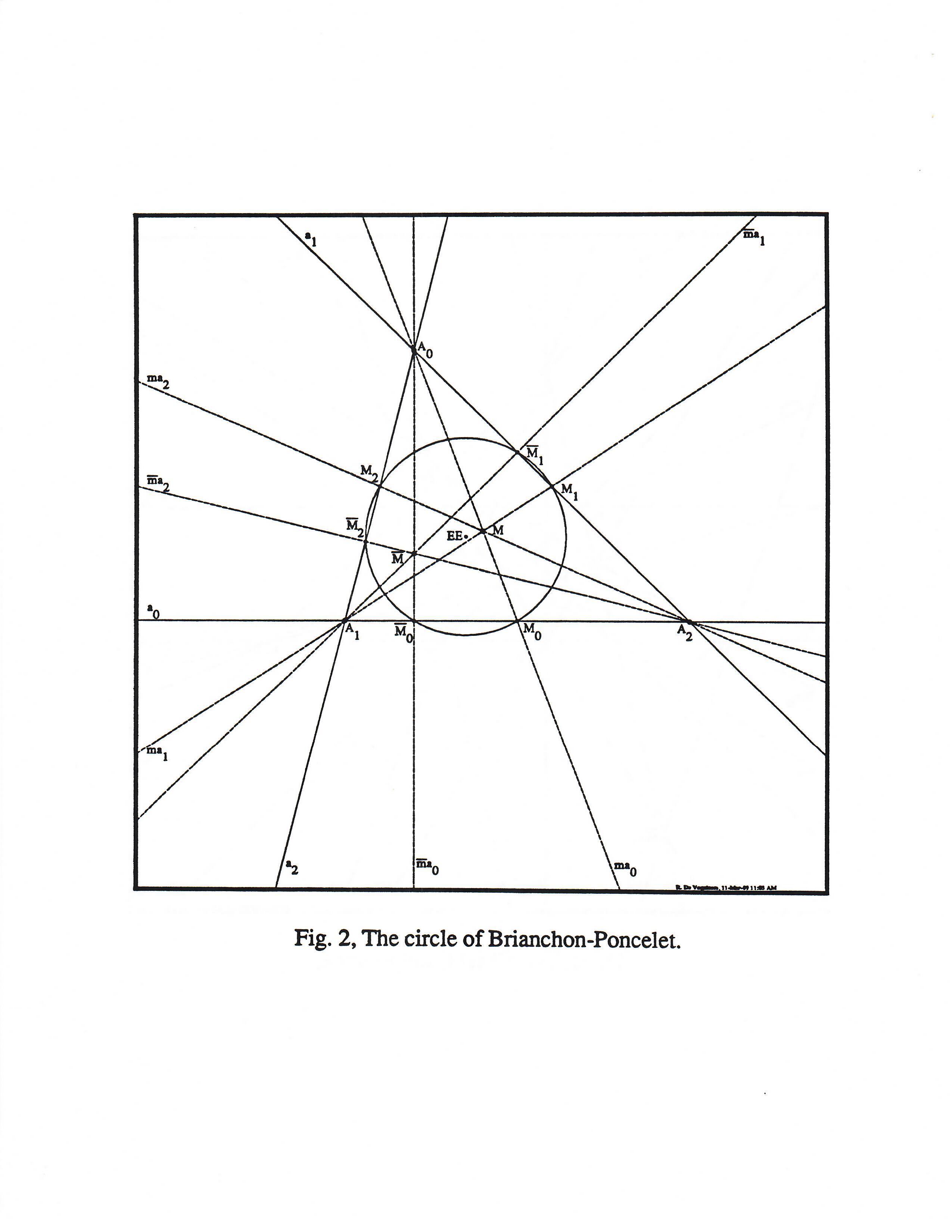}}
\end{figure}

\begin{figure}
{\centering\includegraphics[width=\linewidth]{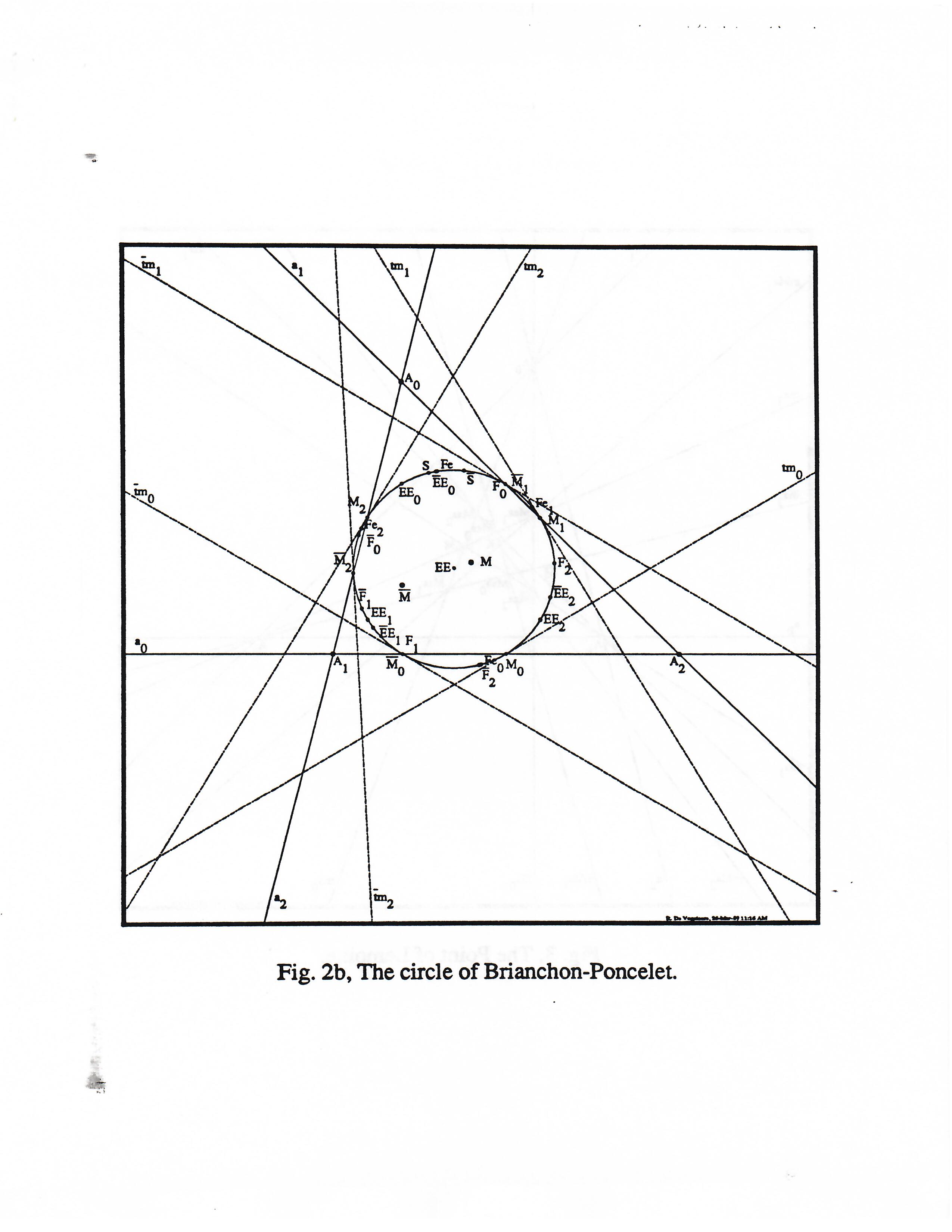}}
\end{figure}

\begin{figure}
{\centering\includegraphics[width=\linewidth]{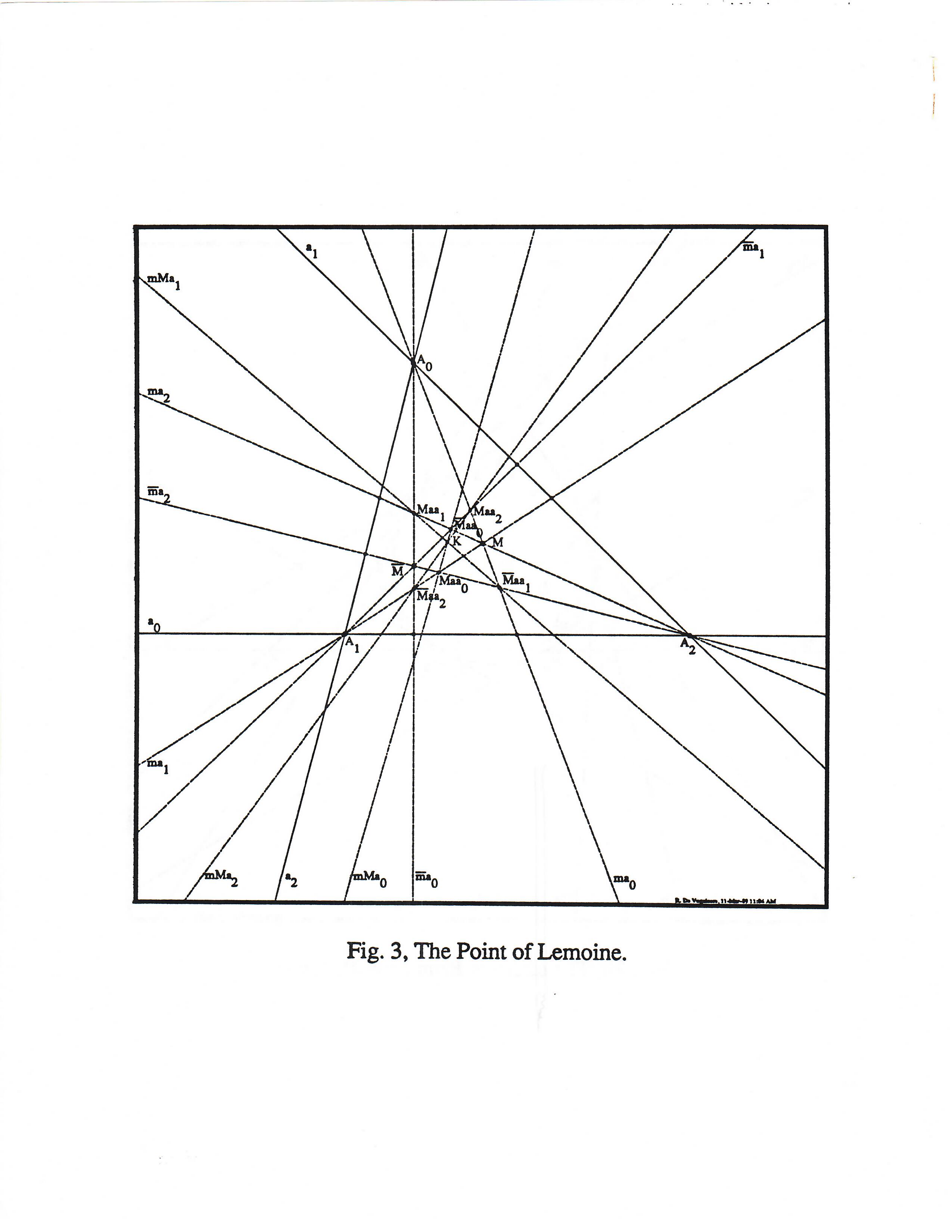}}
\end{figure}

\begin{figure}
{\centering\includegraphics[width=\linewidth]{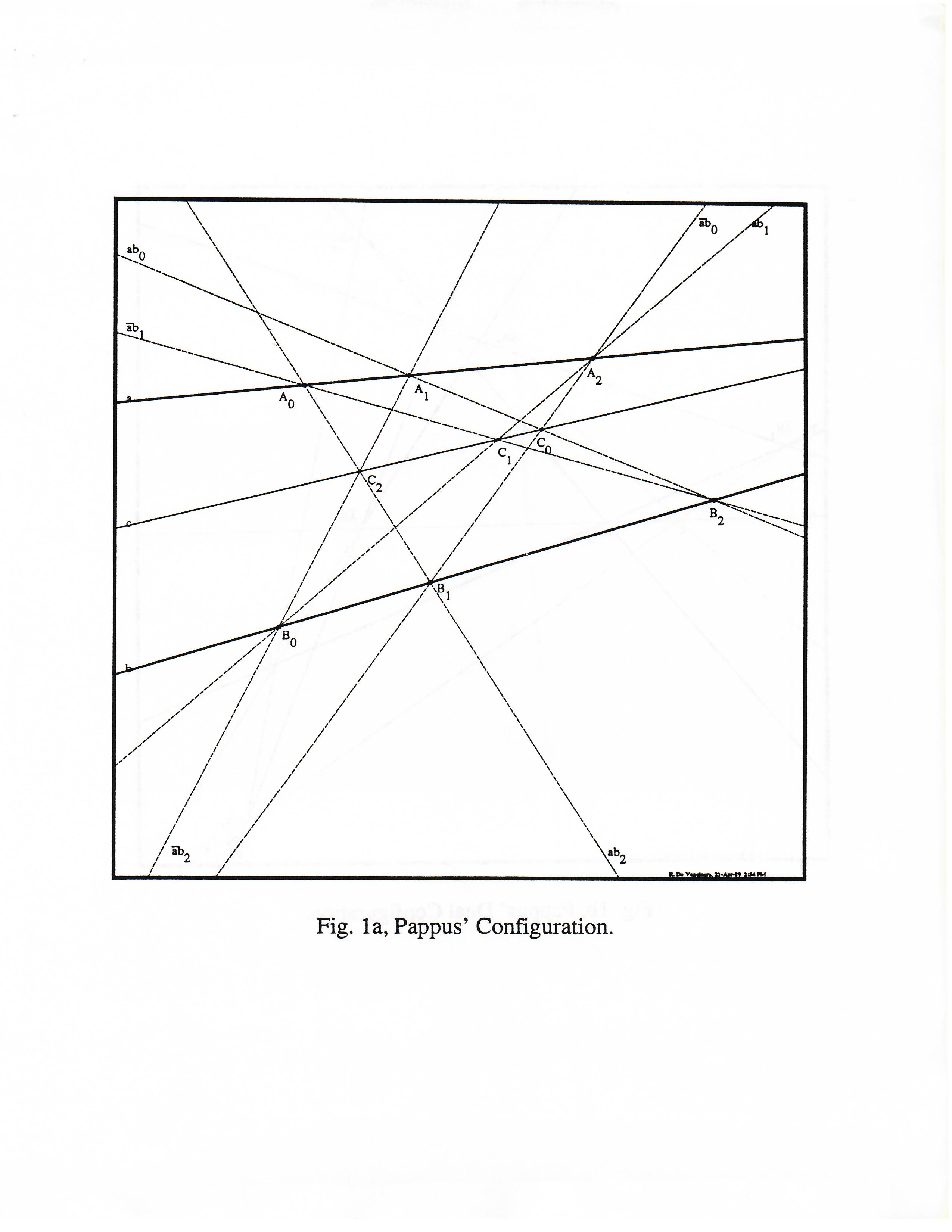}}
\end{figure}

\begin{figure}
{\centering\includegraphics[width=\linewidth]{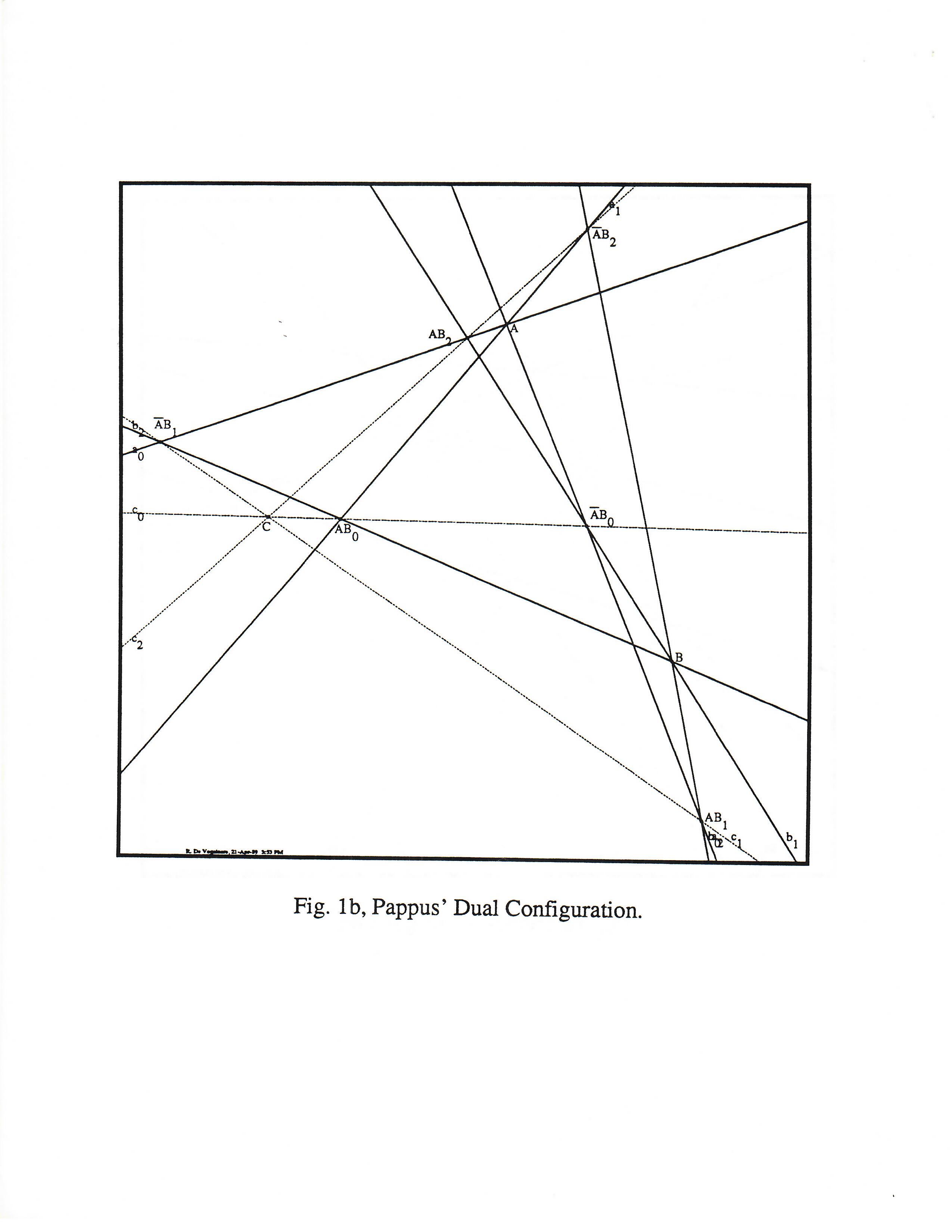}}
\end{figure}

\begin{figure}
{\centering\includegraphics[width=\linewidth]{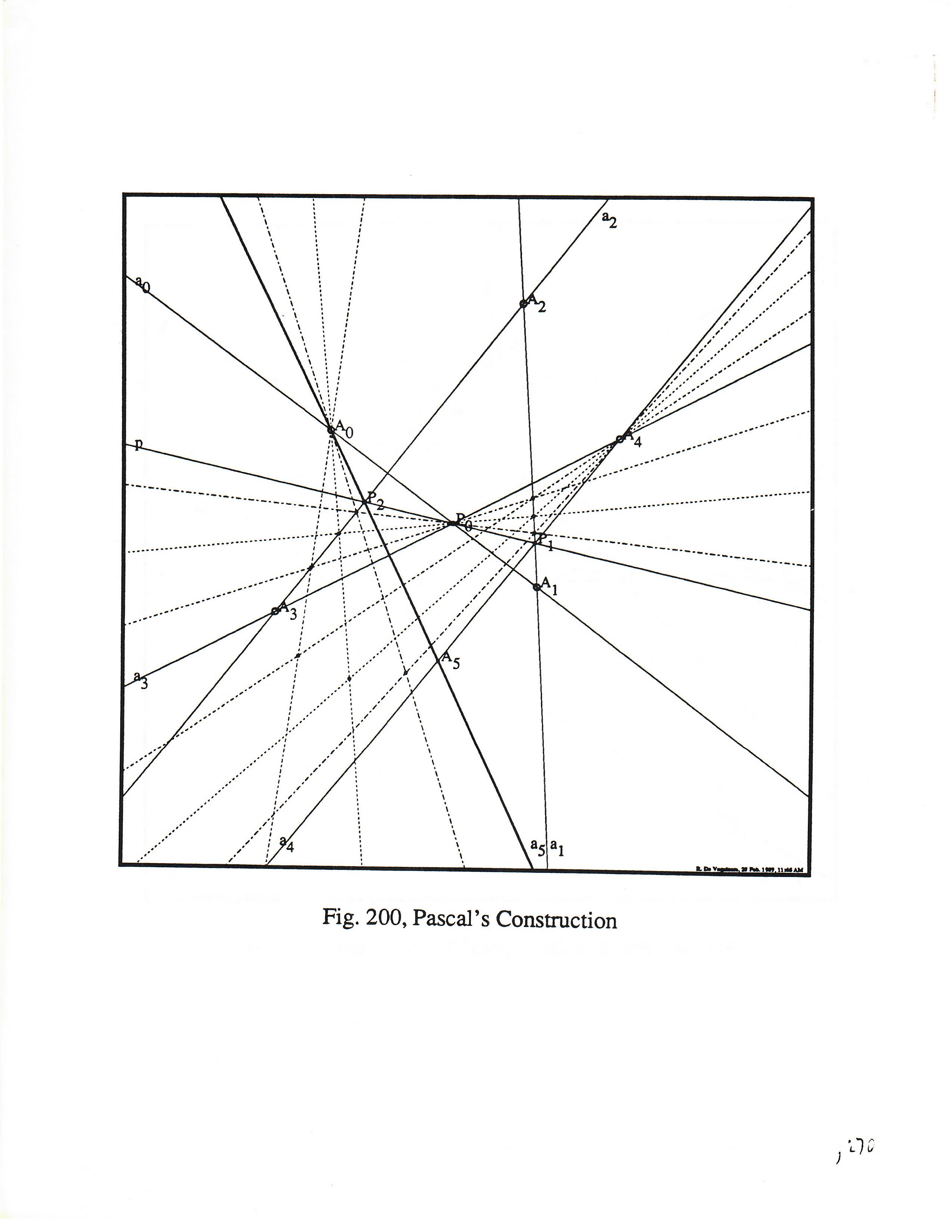}}
\end{figure}

\begin{figure}
{\centering\includegraphics[width=\linewidth]{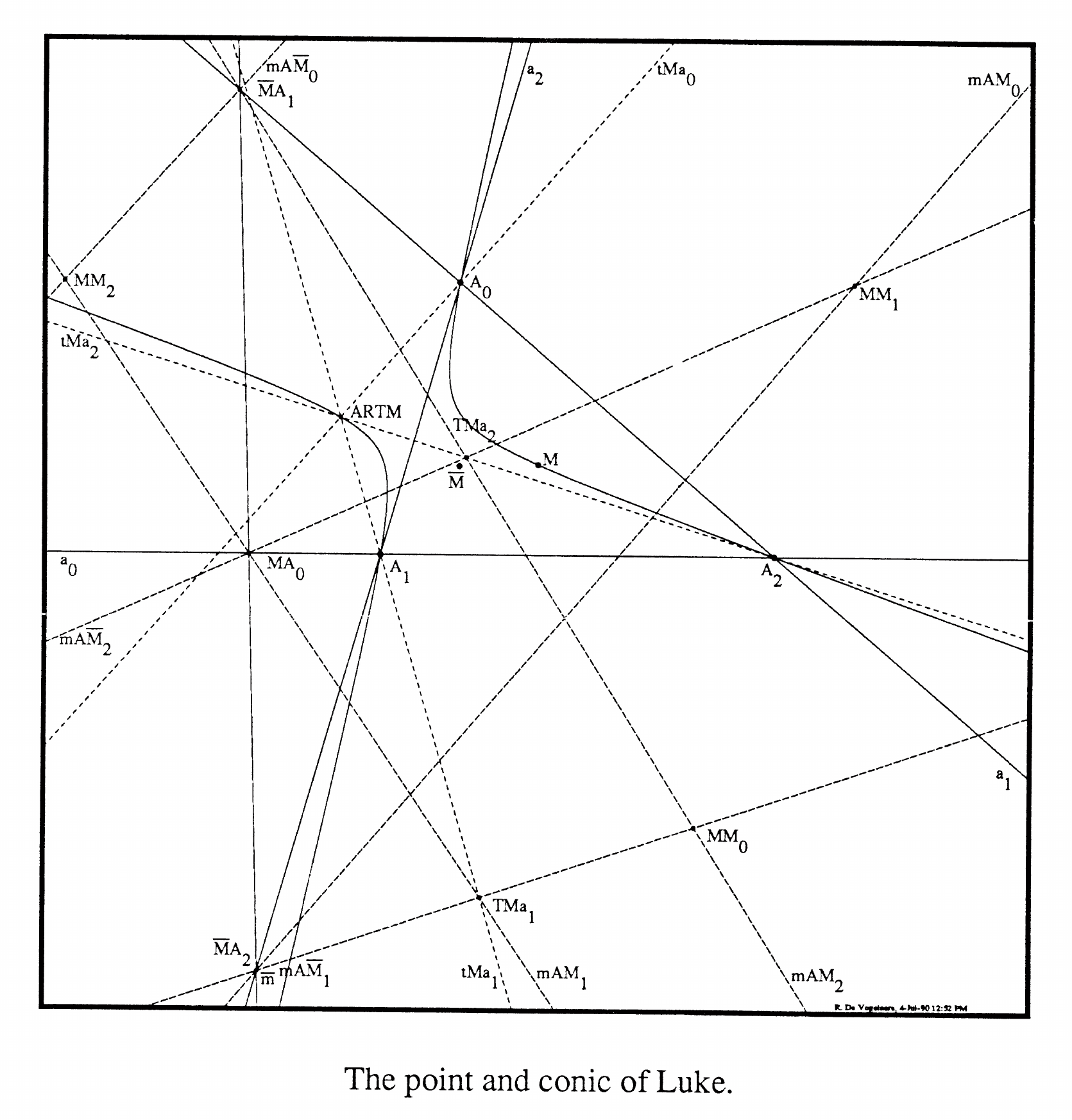}}
\end{figure}

\chapter{FINITE NON-EUCLIDEAN GEOMETRY}


\setcounter{section}{-1}
	\setcounter{section}{-1}
	\section{Introduction.}

	In Chapter IV, Finite Euclidean geometry was constructed.
	In it, we have seen that the angles can be given as integers.
	In the finite hyperbolic Euclidean geometry, the angles can be
	represented by elements in ${\Bbb Z}_{p-1}$ and in finite elliptic
	Euclidean geometry by elements in ${\Bbb Z}_{p+1}$.  The distances can,
	in either case, be represented elements in ${\Bbb Z}_p$ or by $\delta$
	times an element in ${\Bbb Z}_p,$ where $\delta$ is such that $\delta^2$
	is a non quadratic residue in ${\Bbb Z}_p.$\\
	I made many attempts to define angles and distances for a geometry
	which can be considered as the finite form of non-Euclidean geometry.
	The clue was finely provided by the work of Laguerre.
	I will show that using this definition, both angles and distances
	can be treated symmetrically, or to use a mathematical terminology,
	that we have duality between the notions of angle of 2 lines and
	distance of 2 points.

	For those familiar with non-Euclidean geometry, in the classical case,
	there is a distinction between the hyperbolic non-Euclidean geometry
	of Lobatchevski and the elliptic non-Euclidean geometry of Bolyai.
	The axioms, in a form already familiar to Saccheri, are:
	there exists a triangle whose sum of interior angles is equal to
	(Euclidean), smaller than (Lobatchevski) or greater than (Bolyai)
	180 degrees.

	In the hyperbolic case, the set of lines through a point $P$ not on $a$
	line $l$ is subdivided into 2 sets, those which intersect $l$ and those
	who do not.  If we assume continuity, there are 2 lines which form the
	boundary of either set and are called parallels.
	The simplest model is obtained by starting with the 2 dimensional
	projective plane and choosing a given conic as ideal.  We define as
	points those inside the conic and as lines the portion of
	the lines of the projective plane inside the conic.  The parallels to
	$l$ from a point $P$ not on $l$ are those which pass through the
	intersection of $l$ with the conic.

	In the elliptic case, there are no parallels, the lines always
	intersect.  The simplest model is obtained by choosing a sphere in 3
	dimensional Euclidean geometry.  We define as lines the great circles
	of a sphere and as points the points of the sphere, identifying each
	point with its antipode.

	In the finite case, there is no distinction between the elliptic and
	the hyperbolic case.  Indeed in finite projective geometry, the inside
	or the outside of a conic cannot be defined.  Instead, for some lines
	there are no parallels and for others the situation is analogous to
	that described in the classical hyperbolic case.  For those who like to
	refer to some geometric picture, the image of the geometry on the
	sphere will be useful although imperfect.  I will refer to it from time
	to time.  Again, although I would find it more satisfactory to proceed
	synthetically, I will proceed algebraically to reach the goal more
	quickly.

	In finite Euclidean geometry, I proceeded from projective geometry in
	3 steps, affine geometry, involutive geometry and Euclidean geometry.
	Here I will proceed in 2 steps, polar geometry and non Euclidean
	geometry.

	In the involutive geometry, an involution on the ideal line is chosen,
	from which the notions of perpendicularity and circles are derived.
	In finite projective geometry, no conic can be distinguished from any
	other.  To define finite non-Euclidean geometry, I proceed in 2 steps.
	In the first step, I define the finite polar geometry by chosing, or
	better still, prefering a specific polarity, or equivalently a specific
	conic.  From it, the notions of parallelism, circles,
	equality of segments,  $\ldots$  , are derived.  In the second step,
	I introduce the notions of measure of distances and measure of angles,
	in this case also, the ideal conic plays again an essential role.

	\section{Finite Polar geometry.}

	\setcounter{subsection}{-1}
	\ssec{Introduction.}

	After defining the geometry starting from a finite projective geometry
	in which a given polarity is preferred, I define
	elliptic, parabolic and hyperbolic points and lines.
	I then define circles without using the notion of distance,
	equidistance and the dual notion of equiangularity are, as in finite
	Euclidean geometry derived notion.  After defining perpendicularity,
	I define special triangles using equidistance and right angles.  I then
	proceed to define mid-points, medians and mediatrices and finally the
	circumcircles of a triangle.  A new point, which I call the center of
	a triangle is defined using 2 independent methods.  This point also
	exists in classical non-Euclidean geometry, but I have not found any
	reference in the literature.  The intersection of the circumcircles
	of a triangle are obtained and constructed.  Various results obtained
	while studying the center of a triangle are derived.  The circumcircle
	for the special case of a triangle with an ideal vertex is studied and
	finally the properties of the parabola are given in detail.

	\ssec{The ideal conic, elliptic, parabolic and hyperbolic points
	and lines.}

	\sssec{Definition.}
	Among all the conics in the plane, the chosen one is called the
	{\em ideal conic} or the {\em ideal}.  The points on the conic are
	called {\em ideal points} or {\em parabolic points}.  The lines
	tangent to the
	conic are called {\em ideal lines} or {\em parabolic lines}.
	They could also be called isotropic, by analogy with the Euclidean
	case, but I will not use this terminology.\\
	A line which intersects the ideal in 2 real points is called a
	{\em hyperbolic line}, a point which is incident to 2 ideal lines is
	called a {\em hyperbolic point}.\\
	A line which does not contain ideal points or a point which is not on
	an ideal line is called {\em elliptic}.\\
	A point or a line is said to be an {\em ordinary point} if the point
	or the line is either elliptic or hyperbolic.\\
	Two points or two lines are said to be {\em of the same type} if they
	are both either elliptic or hyperbolic.  Points of the same type are
	necessarily ordinary.

	\sssec{Convention.}\label{sec-cidealconic}
	By convention, the conic chosen for the algebraic
	derivation is\\
	\hth$X \cdot X = 0$ or $X_0^2 + X_1^2 + X_2^2 = 0.$

	\sssec{Example.}
	For $p = 13,$ the ideal points are 6, 9, 19, 22, 57, 62, 69,
	76, 79, 118, 134, 141, 148, 153.

	\sssec{Theorem.}
	{\em The polar of  $A = (A_0,A_1,A_2)$ with respect to the ideal conic
	is $\ov{a} = [A_0,A_1,A_2].$}

	\sssec{Notation.}\label{sec-npolar}
	The polar of $A$ will be denoted $\ov{a},$ the pole of $a,$ $\ov{A}$.\\
	This notation should not be confused with the notation in section
	$\ldots$ on finite projective geometry.

	\sssec{Theorem.}
	{\em With $j = +1$ or $-1$, the ideal points on the line
	$a = [a_0,a_1,a_2]$ are}
	\enumb
	\item {\em if $a_1^2 + a_2^2\neq 0,$\\
	\hth$( a_1^2 + a_2^2,$ $-a_0 a_1 + j a_2 \sqrt{d},$
		$-a_0 a_2 - j a_1 \sqrt{d})$,
	where\\
	\hth$d = - (a_0^2 + a_1^2 + a_2^2),$}
	\item {\em if $a_1^2 + a_2^2 = 0$ and $a_1 . a_2\neq  0$,\\
	\hth$( 0, a_1 + j a_2 \sqrt{-1},$ $a_2 - j a_1 \sqrt{-1} ),$}
	\item {\em if $a_1 = a_2 = 0$,\\
	\hth$(0, 1, j \sqrt{-1} ).$}
	\enume

	\sssec{Example.}
	For $p = 13,$
	let $a = [124] = [1,8,6],$ then $d = 3,$ $\sqrt{d} = 4,$ the ideal
	points on $a$
	are $(-4, -8-2, -6-6) = (1,9,3) = (134)$
	and $(-4, -8+2, -6+6) = (1,8,0) = (118).$

	\sssec{Theorem.}
	{\em The point $A = (A_0,A_1,A_2)$ and the line
		$\ov{a} = [A_0,A_1,A_2]$ are}
	\enumb
	\item{\em  parabolic, iff  $A \cdot A = 0,$}
	\item{\em  elliptic, iff  $-A \cdot A$ is a non quadratic residue
		modulo $p$,\\
	    in other words,if there is no integer $x$ such that}\\
	\hth$x^2 = -A \cdot A,$
	\item{\em  hyperbolic, iff  $- A \cdot A$  is a quadratic residue
		modulo $p$.}
	\enume

	\sssec{Example.}
	For $p = 13,$
	(6) = (0,1,5) is parabolic, (172) = (1,12,2) is elliptic and (124) =
	(1,8,6) is hyperbolic.

	\sssec{Theorem.}
	\enumb
	\item{\em There are $p+1$ parabolic or ideal points,}
	\item{\em There are $\frac{p(p-1)}{2}$ elliptic points,}
	\item{\em There are $\frac{p(p+1)}{2}$ hyperbolic points.}
	\enume

	Proof:  Each of the $p+1$ parabolic line meets the other $p$ parabolic
	lines in a hyperbolic point.

	\sssec{Definition.}
	Two {\em lines are parallel} if they have an ideal point in common.\\
	Two {\em points are parallel} if they have an ideal line in common.

	\sssec{Example.}
	For $p = 13,$
	(61) = (1,3,8) and (71) = (1,4,5) are parallel, they are on [134] =
	[1,9,3].

	\sssec{Theorem.}
	{\em The intersections of the sides of a triangle with the polars
	of the opposite vertex with respect to any conic are collinear.}

	This follows at once from II.\ref{sec-tvecc}.7 if we choose the
	coordinates in such a way that the conic is the ideal conic.

	\ssec{Circles in finite polar geometry.}

	\setcounter{subsubsection}{-1}
	\sssec{Introduction.}
	There are 3 kinds of circles in polar geometry.\\
	A hyperbolic circle is a conic tangent to the ideal conic at 2 distinct
	points.  Its center is the intersection of these tangents.\\
	An elliptic circle is a conic tangent to the ideal conic at 2 distinct
	complex conjugate points.\\
	A parabolic circle is one for which the two points of tangency
	coincide.\\
	I will give now the corresponding algebraic definition, when convention
	\ref{sec-cidealconic} is used.\\
	Having introduced the notion of circles, it is natural to define the
	notion of equidistance between points and equiangularity between lines.
	When measure of angles and distances will be introduced, the
	compatibility of the 2 concepts equivalence and measure will be made
	clear.

	\sssec{Definition.}
	The {\em circles of center} $C = (C_0,C_1,C_2)$ are the conics with
	equation\\
	\hth$X \cdot X + k (X \cdot C)^2 = 0.$

	\sssec{Definition.}
	The line $c = [C_0,C_1,C_2]$ is called the {\em central line} of the
	circle.

	\sssec{Theorem.}
	{\em The central line is the polar of the center in the polarity
	associated to the circle as well as in the polarity associated to the
	ideal conic.}

	\sssec{Definition.}
	A circle is called {\em hyperbolic} if its center is hyperbolic,
	{\em elliptic,} if its center is elliptic and {\em parabolic} if its
	center is a parabolic or ideal point.

	\sssec{Theorem.}
	{\em The ordinary points on a circle are all either hyperbolic or
	elliptic.}

	Proof:  If $k$ is a quadratic residue modulo $p,$ then $-X \cdot X$ is
	a quadratic residue and $X$ is necessarily hyperbolic.  If $k$ is a
	non residue, then
	$-X \cdot X$ is a non residue and $X$ is necessarily elliptic.

	\sssec{Theorem.}
	{\em If a circle is hyperbolic, the lines through the center and
	the ideal points on the central line are tangent to both the ideal
	conic and the circle.
	If the circle is parabolic, the center $C$ is an ideal point and its
	polar is the tangent at the ideal point to both the ideal conic and the
	circle.}

	All hyperbolic circles can be constructed using the degenerate form of
	Pascal's construction.  The following Theorem allows the construction
	of parabolic circles and of many elliptic circles.

	\sssec{Theorem.}
	{\em For any circle of center $C$ and central line $c$ through a point
	$X_1$	not on $c,$ if $I_1$ is an ideal point on $C \times X_1$ and
	$I_2$ is a distinct
	ideal point not on $c$ and $I_1 \times I_2$ meets $c$ in $X_0$ then 
	$X_2 := (X_0 \times X_1) \times (C \times I_2)$ is also on the circle.}

	\sssec{Theorem.}
	{\em If a circle of center $C$ is not parabolic, let $A$ and $B$ be
	arbitrary points on the circle, let $M$ and $N$ be the other ideal
	points on $C \times A$ and $C \times B,$ then the central line,
	$A \times B$ and $M \times N$ pass through the same point.}

	\sssec{Example.}
	For $p = 13,$ (see g13.tab)
	\enumb
	\item One of the hyperbolic circles of center (124) has the equation\\
	\hth$6 (x^2 + y^2 + z^{2)} + (x + 8y + 6z)^2 = 0,$\\
	\hth or $7 x^2 + 5 y^2 + 3 z^2 + 5yz - zx + 3 xy = 0.$\\
	    It contains the ideal points 118 and 134 and the elliptic points
	    2, 7, 44, 46, 54, 56, 105, 111, 135, 151, 158, 164.
	\item One of the elliptic circles with center (172) has equation\\
	\hth$2(x^2 + y^2 + z^2 + (x - y + 2z)^2 = 0,$\\
	\hth or $3 x^2 + 3 y^2 + 6 z^2 - 4 yz + 4 zx - 2 xy = 0.$\\
	    It contains the elliptic points 7, 13, 15, 21, 41, 44, 70, 77, 98,
		111, 116, 151, 156, 169.
	\item One of the parabolic circles with center (6) has the equation\\
	\hth$(x^2 + y^2 + z^{2)} - (y + 5z)^2 = 0,$\\
	\hth or $x^2 + 2 z^2 + 3 yz = 0.$\\
	    It contains the ideal point 6 and the hyperbolic points 1, 33, 39,
		81, 89, 100, 101, 109, 110, 121, 129, 171, 177.
	\enume

	\sssec{Definition.}
	Two {\em circles are parallel} if they have one ideal point in
	common.  Two {\em circles are concentric} if they have the same center.

	\sssec{Theorem.}
	{\em Two concentric circles have all their ideal points in common.
	One for the parabolic circles, 2 for the hyperbolic circles.}

	\sssec{Definition.}
	The points $A$ and $B$ {\em are equidistant from the point} $C$  iff
	there exists a circle of center $C$ passing through both $A$ and $B.$

	\sssec{Theorem.}
	{\em $A$ and $B$ are equidistant from $C$  iff\\
	\hth$(A \cdot C)^2\:(B \cdot B) = (B \cdot C)^2\:(A \cdot A).$}

	This suggest the more general definition:

	\sssec{Definition.}\label{sec-ddistneuc}
	The {\em distance between the points} $A$ and $B$ {\em is the same as
	the distance between the points} $C$ and $D$  iff\\
	\hth$(A \cdot B)^2\:(C \cdot C)(D \cdot D)
		= (C \cdot D)^2\:(A \cdot A)(B \cdot B).$\\
	The {\em angle between the lines} $a$ and $b$ {\em is the same as the
	angle between the lines} $c$ and $d$  iff\\
	\hth$(a \cdot b)^2\:(c \cdot c)(d \cdot d)
		= (c \cdot d)^2\:(a \cdot a)(b \cdot b).$

	\sssec{Definition.}
	The angle between $a$ and $b$ is a {\em right angle} iff
	$a \cdot b = 0$ and the distance between $A$ and $B$ is a {\em right
	distance} iff $A \cdot B = 0.$

	\sssec{Comment.}
	Although the distance between 2 points $A$ and $B$ has not yet
	been defined, I will by convention use the notation $d(A,B)$.  This will
	be acceptable, in polar geometry, as long as the notation appears in
	both sides of an equality.  I  will later define the distance between
	2 points and show that it is consistent with \ref{sec-ddistneuc}.

	\sssec{Theorem.}
	{\em The notion of equidistance between pairs of points and the
	notion of equiangularity between pairs of lines is an equivalence
	relation, in other words, the relation is\\
	reflexive: $d(A,B) = d(B,A),$\\
	symmetric: $d(A,B) = d(C,D) \implies d(C,D) = d(A,B),$\\
	transitive: $d(A,B) = d(C,D)$ and $d(C,D) = d(E,F) \implies
	d(A,B) = d(E,F).$}

	\ssec{Perpendicularity.}

	\sssec{Definition.}
	The line $b$ is {\em perpendicular to the line} $a$ iff $b$ passes
	through the pole $A$ of $a$ with respect to the ideal conic, in
	other words, when $a$ and $b$ are conjugates with respect to the ideal
	conic.

	\sssec{Theorem.}
	{\em If the line $b$ is perpendicular to $a,$ then\\
	\hth$b \cdot a = 0.$\\
	In other words, the angle between $a$ and $b$ is a right angle.}

	This follows at once from \ref{sec-ddistneuc}.

	\sssec{Theorem.}\label{sec-tneucortho}
	{\em The perpendicular $h_0,$ $h_1,$ $h_2$ from the vertices
	$A_0,$ $A_1,$ $A_2$ of a triangle to the opposite sides have a point
	$H$ in common.  Moreover,\\
	\hth$h_0 = (A_2 \cdot A_0)\: A_1 - (A_0 \cdot A_1)\: A_2,$\\
	\hth$h_1 = (A_0 \cdot A_1)\: A_2 - (A_1 \cdot A_2)\: A_0,$\\
	\hth$h_2 = (A_1 \cdot A_2)\: A_0 - (A_2 \cdot A_0)\: A_1.$\\
	\hth$H = (A_2 \cdot A_0)(A_0 \cdot A_1)\: A_1*A_2
		+ (A_0 \cdot A_1)(A_1 \cdot A_2)\: A_2*A_0$\\
	\hti{12}$+ (A_1 \cdot A_2)(A_2 \cdot A_0)\: A_0*A_1.$}

	Proof:  $h_0 := A_0 * (A_1*A_2),$ is indeed a line through $A_0$
	perpendicular to $A_1 * _A2$. The results follow easily from
	II.\ref{sec-tvecc}. Related results are obtained in \ref{sec-taltdual}.

	\sssec{Definition.}
	$h_i$ are called the {\em altitudes} of the triangle.  The point $H$ is
	called the {\em orthocenter}.

	\sssec{Example.}\label{sec-eneuch}
	For $p = 13,$
	if $A_0 = (0) = (0,0,1),$ $A_1 = (18) = (1,0,4),$
	$A_2 = (67) = (1,4,1),$
	then $h_0 = (27) = [1,1,0],$ $h_1 = [1,4,3],$ $h_2 = [1,0,12]$
	and $H = (171) = (1,12,1).$

	\sssec{Comment.}
	Section 7 could be placed here, but then the motivation would be absent.

	\ssec{Special triangles.}

	\sssec{Definition.}
	A {\em right}, {\em double right}, {\em polar} {\em triangle} is a
	triangle which has one, two or three right angles.\\
	A {\em right sided} or {\em double right} {\em sided} {\em triangle} is
	a triangle for which the distance between one pair
	or two pairs of vertices is a right distance.

	\sssec{Examples.}
	For $p = 13:$
	The triangle $A = (8) = (0,1,7),$ $B = (17) = (1,0,3),$
	$C = (36) = (1,1,9),$
	with sides [44] = [1,2,4], [150] = [1,10,6], [161] = [1,11,4]
	is a right triangle at $B$.\\
	The triangle $A = (44),$ $B = (17),$ $C = (36),$ with sides [44], [130]
	= [1,8,12], [161] is a double right triangle at $B$ and $C$.
	The triangle $A = (44),$ $B = (17),$ $C = (161),$ with sides [44], [17]
	[161] is a polar triangle.\\
	Exchanging vertices and sides, we obtain, by duality, examples of
	right sided and double right sided triangles.

	\sssec{Definition.}
	A {\em triangle is isosceles} if 2 pairs of vertices are equidistant.\\
	A {\em triangle is equilateral} if all 3 pairs of vertices are
	equidistant.

	\sssec{Theorem.}
	\enumb
	\item{\em If a triangle $\{ABC\}$ is such that $d(A,B) = d(A,C),$ then\\
	\hth$d(a,b) = d(a,c).$}
	\item{\em  If a triangle $\{ABC\}$ is such that $d(A,B) = d(B,C)
		 = d(A,C),$ then}\\
	\hth$d(a,b) = d(a,c) = d(b,c).$
	\enume

	Proof:  The second part follows, by transitivity, from the first part.
	For the first part, let us set $p = A \cdot A,$ $q = B \cdot B,$
	$r = C \cdot C,$ $t = B \cdot C,$
	$u = C \cdot A,$ $v = A \cdot B.$   The hypothesis implies\\
	\hth$v^2 r = u^2 q = s.$\\
	We want to prove that\\
	\hth$w = (a \cdot b)^2 c \cdot c$\\
	does not change when we exchange $b$ and $c$ or $q$ and $r$ as well as
	$u$ and $v$.
	Using II.\ref{sec-tvecc}.2 and .3, $a  \cdot b = u t - v r$ and 
	$c  \cdot  c = p q - v^2,$ \\
	therefore\\
	\hth$w = (u^2 t^2 + v^2 r^2 -2 t u v r) (p q - v^2)$\\
	\hti{12}$= p t^2 s - t^2 u^2 v^2 + p s q r - s^2 - 2 p t q r u v 
			- 2 t s u v.$

	\sssec{Example.}
	For $p = 13:$
	The triangle $A = (172) = (1,12,2),$ $B = (7) = (0,1,6),$ $C = (13)
	= (0,1,12),$ with sides $a = [14] = [1,0,0],$ $b = [182] = [1,12,12],$
	$c = [74] = [1,4,8]$ is an isosceles triangle.
	The triangle $A = (172),$ $B = (7),$ $C = (15) = (1,0,1),$ with sides
	$a = [104] = [1,6,12],$ $b = [182],$ $c = [74]$ is an equilateral
	triangle.

	\sssec{Theorem.}
	{\em In a polar triangle, each vertex is the pole of the opposite
	side and the distance between the vertices is a right distance.}

	\sssec{Definition.}
	Two {\em triangles are dual of each other}  iff the sides of one
	are the polar of the vertices of the other.

	\sssec{Example.}
	The dual of the triangle of example \ref{sec-eneuch} is, with $p = 13,
$\\
	$A_0 = (173) = (1,12,3),$ $A_1 =(53) = (1,3,0),$ $A_2 = (1) = (0,1,0).$

	\sssec{Theorem.}
	{\em A polar triangle is its own dual.}

	\sssec{Theorem.}\label{sec-taltdual}
	{\em The altitudes of a triangle and of its dual coincide.  The
	orthocenter of a triangle and of its dual coincide.}

	The proof is left as an exercise.


\setcounter{section}{1}
\setcounter{subsection}{4}
	\ssec{Mid-points, medians, mediatrices, circumcircles.}

	\setcounter{subsubsection}{-1}
	\sssec{Introduction.}
	For this section, the analogy with the model of the
	non-Euclidean geometry on the sphere is useful.  We recall that each
	point has 2 representations on the sphere, which are antipodes of each
	other.  If we take 2 points $A$ and $B,$ let $A'$ and $B'$ be their
	antipodes, there are 2 points on the great circle, (in the plane
	through the center of the sphere) which are equidistant from $A$ and
	$B,$ namely a point on the arc $AB$ and a point on the arc $A'B,$
	which is the antipode
	of the mid-point on the arc $AB'.$  But the analogy is not complete, in
	the finite case, it is only when the points are of the same type that
	mid-points exist.  There is about 1 chance in 2 that the points are not
	of the same type, there are then no mid-points, there is about 1 chance
	in 2 that they are of the same type, there are then 2 mid-points, this
	is an other example of what I call the law of compensation.\\
	To simplify the algebra, I will introduce a scaling in
	\ref{sec-nneucscal}.
	The scaling contains an arbitrary sign, which may be thought as
	corresponding to the 2 representations on the sphere.  The systematic
	way which is chosen could be replaced by some other one.  The choice
	is influenced by the choice of a primitive root of $p$ and the sign of
	the square root and depends on the rule $\ldots$ .
	Having the concept of mid-points, we can consider those of the vertices
	of a triangle, if the vertices are scaled, we can define interior and
	exterior mid-points.  Again, the choice is arbitrary and depends on the
	rule $\ldots$ .\\
	To each side correspond 2 medians, these meet 3 by 3 in 4 points
	corresponding to the barycenter.  Again the analogy with the geometry
	on the sphere is useful, the 4 barycenters can be considered as
	corresponding to the triangles $\{ABC\},$ $\{A'BC,\}$ $\{AB'C\},$
	$\{ABC'\}$.\\
	A similar treatment can be made for the mediatrices which meet 3 by 3
	in 4 points, each is the center of a circumcircle of the triangle
	$\{ABC\}.$\\
	But, again, the analogy with the geometry on the sphere is not complete.
	Given a triangle, there are about 3 chances in 4 that the 3 vertices
	are not all of the same type, in this case there is no barycenter and
	no circumcircle.  In about 1 chance out of 4, the 3 points are of the
	same type, and there are 4 barycenters and 4 circumcircles.  Again this
	is the compensation.\\
	If the vertices of the triangle are of the same type, the 4 lines
	joining a barycenter to the corresponding center of a circumcircle can
	be considered as generalizations of the line of Euler.  It is natural
	to conjecture that these four lines are concurrent.  This is indeed the
	case.  The surprise is that this point $V$ is not the orthocenter.  The
	coordinates of $V$ are real even if the vertices of the triangle are not
	of the same type.  $V$ must therefore be obtainable in an independent
	way.  One such method is described in section 7.

	I first recall the convention of I.\ref{sec-csqroot}.

	\sssec{Convention.}
	Given $\delta$  a specific square root of a specific non quadratic
	residue of $p$, we choose the square root $a$ of a quadratic residue or
	the square root $a\delta$, of a non residue in such a way that
	$0\leq a<\frac{p-1}{2}.$

	\sssec{Notation.}\label{sec-nneucscal}
	Using the preceding convention, a square root is uniquely defined. It is
	convenient to introduce an other scaling for points and lines different 
	from that given in II.\ref{sec-crep}.\\
	If $A = (A0,A1,A2)$ and $A$ is not an ideal point,\\
	\hth$	       A' = \frac{A}{\sqrt{-A \cdot A}}.$\\
	\hth$| A|  = \sqrt{-A \cdot A}$ is called the length of $A.$
	Either each component is an integer, or each component is an integer
	divided by $\delta$ , in this last case we say that $A'$ is pure
	imaginary.

	\sssec{Theorem.}
	{\em If $A$ is hyperbolic, $A'$ is real, if $A$ is elliptic, $A'$ is
	pure imaginary.  Moreover $A'\cdot A' = -1.$}

	\sssec{Definition.}
	Given 2 points $A$ and $B$ of the same type, $M$ on $A \times B$ is
	called a {\em mid-point} of [A,B]  iff  the distances $M A$ and $M B$
	are equal.

	\sssec{Theorem.}
	{\em The mid-points of $[A,B]$ are $M = A' + B'$ and $M^- = A' - B'.$}

	Proof:  Because of \ref{sec-ddistneuc}, $d(M,A) = d(M,B)$ if\\
	\hth$(M \cdot A')^2 = (M \cdot B')^2,$\\
	or if $(A' \cdot A')^2$ + $(A' \cdot B')^2
		+ 2(A' \cdot A')(A' \cdot B')$\\
	\hth$	       = (B' \cdot B')^2 + (B' \cdot A')^2
	+ 2(B' \cdot B')(B' \cdot A'),$\\
	which is satisfied because of $A' \cdot A' = B' \cdot B' = -1$ and\\
	\hth$	A' \cdot B' = B' \cdot A'.$

	The proof is similar for $M^-$.

	\sssec{Definition.}
	$M$ is called the interior mid-point, $M^-$ is called the
	{\em exterior mid-point}.

	\sssec{Example.}
	For $p = 13,$
	the mid-points of (44) = (1,2,4) and (164) = (1,11,7) are
	(115) = (1,7,10) and (124) = (1,8,6).  Indeed $| 44|  = \sqrt{5},$
	$| 164|  = \sqrt{11},$ hence, $\frac{| 44|}{| 164|}
	= \sqrt{\frac{-8}{-2}} = \sqrt{4} = 2,$
	therefore the mid-points are (1,2,4) + 2(1,11,7) = (3,24,18) = (1,8,6)
	and $(1,2,4) - 2(1,11,7) = (-1,-20,-10) = (1,7,10)$.\\
	With $\delta^2 = 8,$
	$A' = \frac{A}{\delta},$ $B' = \frac{B}{(6 \delta},$ therefore the
	interior mid-point is
	$A + \frac{1}{6}B = A - 2B = (1,7,10)$ and the exterior mid-point is
	$A - \frac{1}{6}B = A + 2B = (1,8,6).$

	\sssec{Definition.}
	$m$ is called a {\em mediatrix} of $[A,B]$  iff  $m$ is perpendicular
	to $A \times B$ and passes through a mid-point of $[A,B].$

	\sssec{Theorem.}
	{\em
	\hth$	m  := A' - B'$ passes through $M = A' + B'$ and\\
	\hth$	       m^- := A' + B'$ passes through $M^- = A' - B'.$}

	Proof:  $m$ is perpendicular to $A \times B$ because\\
	\hth$	       m \cdot (A * B) = m \cdot (A' * B')
		= A' \cdot (A' * B') - B' \cdot (A' * B') = 0.$\\
	$m$ passes through $M,$ because
	$m \cdot M = (A'-B') \cdot (A'+B') = A' \cdot A' - B' \cdot B'
	= -1 - (-1) = 0$.

	\sssec{Theorem.}
	{\em The set of points equidistant from $A$ and $B$ are on $m$ or $m'.$}

	\sssec{Definition.}
	In a triangle, the line joining a vertex to the interior (exterior)
	mid-points of the opposite side is called an {\em interior (exterior)
	median}.

	\sssec{Theorem.}
	{\em If a triangle is isosceles, with $d(A_0,A_1) = d(A_0,A_2),$
	then a median through $A_0$ is also a mediatrix.}

	\ssec{The center  $V$  of a triangle.}

	\sssec{Theorem.}
	{\em Let $M_i$ and $M^-_i$ be the interior and exterior mid-points of
	$A_{i-1}$ and $A_{i+1}$, let $n_i$ and $n^-_i$ be the interior and
	exterior medians associated to $A_i.$}
	\enumb
	\item $G_3 := n_0 \times n_1 \Rightarrow G_3 \cdot n_2 = 0.$ (*)
	\item $G_i := n_i \times n^-_{i+1} \Rightarrow G_i \cdot n^-_{i-1}
	 = 0 $(*).\\
	  {\em Let $m_i$ and $m^-_i$ be the interior and exterior mediatrices of
	    $A_{i+1} A_{i-1}$.}
	\item $O_3 := m_0 \times m_1 \Rightarrow O_3 \cdot m_2 = 0. $(*)
	\item $O_i := m_i \times m^-_{i+1} \Rightarrow O_i \cdot m^-_{i-1}
		= 0. $(*)
	\item $e_j := O_j \times G_j$ {\em and} $V := e_0 \times e_1
		\Rightarrow V \cdot e_2 = V \cdot e_3 = 0. $(*)
	\enume

	The proof follows from Theorem 4.4.12.  As in finite Euclidean geometry
	"*" indicates that there are equivalent definitions, for instance 0
	could be written $G_3 := n_1 \times n_2$ and $G_3 \cdot n_0 = 0.$

	\sssec{Definition.}
	By analogy with Euclidean geometry, the points $G_j$ are
	called the {\em barycenters} of the triangle.\\
	The points $O_j$ are called the {\em centers of the circumcircles} of
	the triangle.\\
	The lines $e_j$ are called the {\em lines of Euler of the triangle}.

	\sssec{Definition.}
	$V$ is called the {\em center of the triangle}.

	\sssec{Theorem.}
	{\em Let $A'_i$ be the normalized coordinates of the vertices of
	a triangle}\footnote{21.9.81}.
	\enumb
	\item{\em The mid points are $A'_{i-1} + j_i A'_{i+1},$ $j_i = +1$ or
		$-1.$}
	\item{\em The mediatrices are} $A'_{i-1} - j_i A'_{i+1}.$
	\item{\em The medians $n_i$ are} $A'_i \times (A'_{i-1} + j_i
		A'_{i+1}).$
	\item{\em Choosing $j_0 j_1 j_2 = 1,$ the medians meet 3 by 3 at the 4
		barycenters which are} $A'_0 + j_2 A'_1 + j_1 A'_2.$
	\item{\em The mediatrices meet 3 by 3 at the 4 centers of circumcircles,
		which are}\\
		$A'_1 * A'_2 + j_2 A'_2*A'_0 + j_1 A'_0*A'_1.$
	\item{\em The Euler lines, joining $G_j$ to $O_j$ are,
		with} $d'_i = A'_{i-1} \cdot A'_{i+1},$
		$(j_0 d'_1 - d'_2) A'_0 + (j_2 d'_2 - j_0 d'_0) A'_1
		+ (d'_0 - j'_2 d'_1) A'_2.$
	\item{\em The Euler lines intersect at $V$ and}\\
	    $V = d'_0 A'_1 * A'_2 + d'_1 A'_2 * A'_0
		+ d'_2 A'_0 * A'_1.$
	\item{\em  Moreover,}
	    $V = A_0 \cdot A_0 A_1 \cdot A_2 A_1*A_2
		+ A_1 \cdot A_1 A_2 \cdot A_0 A_2*A_0
		+ A_2 \cdot A_2 A_0 \cdot A_1 A_0*A_1. $
	\item{\em $V$ exists when the triangle is not a polar triangle.}
	\enume

	Proof:  0. and 1. follow from $\ldots$ .
	For 2., the intersection of $n_1$ and $n_2$ is\\
	\hth$	   n_1*n_2 = - A'_0 (A'_2 \cdot (A'_1*A'_0))$\\
	\hth$		       + j_1 A'_2  (A'_0 \cdot (A'_1*A'_2))$\\
	\hth$		       + j_0 j_1 A'_1  (A'_0 \cdot (A'_1*A'_2))$\\
	\hth$		     = A'_0 + j_0 j_1 A'_1 + j_1 A'_2.$\\
	The same point is obtained if $j_0 j_1 = j_2$ or if
	$j_0 j_1 j_2 = 1.$  There are 4 points corresponding to
	$j_0 = + or -1$ and $j_1 = + or -1.$\\
	The proof of 3. is similar.  The proof of 4. to 8. is left as
	exercises.

	\sssec{Comment.}
	Reality requires, because $A'_i = \frac{A_i}{|A_i|}$  that the
	lengths $|A_i|$  be either all real or all imaginary, hence:

	\sssec{Theorem.}
	{\em The mid-points, mediatrices, medians, barycenter
	and center of circumcircles are real if and only if the vertices are
	either all elliptic or all hyperbolic.  $V$ is always real.}

	\sssec{Example.}\label{sec-neuctr}
	For $p = 13,$ the triangle
	$A_0 = (58) = (1,3,5),$ $A_1 = (51) = (1,2,11),$
	$A_2 = (159) = (1,11,2),$
	has all its vertices hyperbolic.\\
	Let $A'_0 = (6,5,4),$ $A'_1 = (6,-1,1),$ $A'_2 = (6,1,-1).$
	The mid-points of $A_1$ and $A_2$ are (14) = (1,0,0) and
	(13) = (0,1,12).
	All the mid-points are (14), (13); (115), (12); (139), (8).\\
	The mediatrices are [13], [14]; [12], [115]; [8], [139].\\
	The medians are [3], [126]; [91], [176]; [76], [161].\\
	The interior mediatrices [13], [12], [8] meet at $O_3 = (14).$\\
	The centers of the circumcircles are (56), (8), (12) and (14).\\
	The interior medians [3], [91], [76] meet at $G_3 = (33).$\\
	The barycenters are (152), (179), (106), and (33).\\
	The center of the triangle is $V = (152).$

	\ssec{An alternate definition of the center $V$ of a triangle.}

	\sssec{Notation.}\label{sec-nneucv}
	From here on, the following notation will be used
	systematically:
	\enumb
	\item $a_i := A_{i+1} \times A_{i-1}$
	\item $n_i := \frac{A_{i+1}*A_{i-1} }{ a_i,}$
	which means that\\
	\hth$	A_{i+1}*A_{i-1} = n_i a_i,$ defines $n_i,$\\
	$n_i$ is the normalization factor, see 2.3.2. and 2.3.11.
	\item $l_i := A_i \cdot A_i,$
	\item $d_i := A_{i+1} \cdot A_{i-1},$
	\item $t := (A_0 * A_1) \cdot A_2.$\\
	Similarly,
	\item $N_i := \frac{a_{i+1}*a_{i-1} }{A_i,}$
	which means that\\
	\hth$	   a_{i+1}*a_{i-1} = N_i A_i,$ defines $N_i,$
	\item $L_i := a_i \cdot a_i,$
	\item $D_i := a_{i+1} \cdot a_{i-1},$
	\item $T := (a_0 * a_1) \cdot a_2.$
	\enume

	\sssec{Theorem.}\label{sec-tneucv0}
	\enumb
	\item $t = n_1 n_2 N_0 = n_2 n_0 N_1 = n_0 n_1 N_2.$
	\item $n_i^2 L_i = n_i^2 a_i \cdot a_i$
			    $= l_{i+1} l_{i-1} - d_i^2.$
	\item $n_{i+1} n_{i-1} D_i = n_{i+1} n_{i-1} a_{i+1} \cdot a_{i-1}$\\
	\hth$				  = d_{i+1} d_{i-1} - d_i l_i.$
	\item $n_0 n_1 n_2 T = t^2.$

	{\em and the dual relations}
	\item $T = N_1 N_2 n_0 = N_2 N_0 n_1 = N_0 N_1 n_2. $
	\item $N_i^2 l_i = N_i^2 A_i \cdot A_i$\\
	\hth$			   = L_{i+1} L_{i-1} - D_i^2.$
	\item $N_{i+1} N_{i-1} d_i = N_{i+1} N_{i-1} A_{i+1} \cdot A_{i-1}$\\
	\hth$				  = D_{i+1} D_{i-1} - D_i L_i.$
	\item $N_0 N_1 N_2 t = T^2.$
	\enume

	\sssec{Theorem.}\label{sec-tneucv1}
	\enumb
	\item $a_{i+1} * a_{i-1} = t A_i,$
	\item $n_i a_i * A_i = d_{i-1} A_{i-1} - d_{i+1} A_{i+1}.$
	\item $n_i a_i * A_{i+1} = l_{i+1} A_{i-1} - d_i A_{i+1}.$
	\item $n_i a_i * A_{i-1} = d_i A_{i-1} - l_{i-1} A_{i+1}.$
	\enume

	The proof follows easily from 2.3.17. and from 4.6.0.

	\sssec{Example.}\label{sec-eneucv}
	For $p = 13,$ with\\
	$A[] = \{ (0) = (0,0,1),$ $(18) = (1,0,4),$ $(67) = (1,4,1)$ \}, then\\
	$a[] = \{ [173] = [1,12,3],$ $[53] = [1,3,0],$ $[1] = [0,1,0]$ \}.\\
	$l[] = [1,4,5],$ $d[] = [5,1,4],$ $n[] = [10,4,1],$\\
	$L[] = (11,10,1),$ $D[] = (3,12,11),$ $N[] = (1,3,4),$\\
	$t = [4],$ $T = (3).$

	\sssec{Theorem.}
	{\em Let $h$ be the polar of $H$ with respect to the triangle.
	Let $K_i$ be the intersection of $h$ and $a_i,$
	\footnote{13.10.81}\\
	let $v_i$ be the perpendicular at $A_i$ to $A[i] \times K_i$.\\
	Then $v_i$ have a point in common $V$\footnote{19.10.81}.\\
	Moreover, if we define}
	\enumb
	\item $u_i := d_{i+1} d_{i-1} - d_i l_i,$\\
	{\em we have}
	\item $h = u_0 a_0 + u_1 a_1 + u_2 a_2.$
	\item $K_i = u_{i-1}\: A_{i+1} - u_{i+1}\: A_{i-1}.$
	\item $v_i = (d_{i+1}^2 l_{i+1} - d_{i-1}^2 l_{i-1})\: A_i
			- (d_i d_{i+1} l_i - d_{i-1} l_{i-1} l_i)\: A_{i+1}
			+ (d_i d_{i-1} l_i - d_{i+1} l_{i+1} l_i)\: A_{i-1}.$
	\item$V = d_0 l_0\: A_1*A_2 + d_1 l_1\: A_2*A_0 + d_2 l_2\: A_0*A_1.$
	\enume

	Proof:  Because of \ref{sec-tneucortho},\\
	$H = d_1 d_2 a_0 + d_2 d_0 a_1 + d_0 d_1 a_2, $\\
	after simplification,\\
	$H \cdot a_0 = (d_0 d_1 - d_2 l_2) (d_2 d_0 - d_1 l_1), $\\
	using the definition 0, $H \cdot a_0 = u_2 u_1,$
	and because of 2.3.20, after multiplication by $u_0 u_1 u_2,$
	we obtain 1.\\
	$K_i := H * a_i,$ gives 2, after division by $t$.\\
	$A_i * K_i = u_{i+1} a_{i+1} + u_{i-1} a_{i-1},$\\
	therefore,\\
	$V_i = A_i * (A_i * K_i),$\\
	substituting and using 2.3.17.0., we get\\
	\hth$V_i = (u_{i-1} d_{i-1} - u_{i+1} d_{i+1})\: A_i$
		$- u_{i-1} l_i\: A_{i+1} + u_{i+1} l_i\: A_{i-1},$\\
	replacing $u_i$ by its value, we get, 3, from which we obtain\\
	\hth$v_1 * v_2 = d_0 l_0\: A_1*A_2 + d_1 l_1\: A_2*A_0 + d_2 l_2
		\: A_0*A_1,$\\
	after dividing each term by\\
	$(d_2^2 d_1 l_2 + d_1^2 d_0 l_1 + d_0 l_0 l_1 l_2
		       - d_0^3 l_0 - 2 d_1 d_2 l_1 l_2).$

	\ssec{Intersections of the 4 circumcircles.}

	\setcounter{subsubsection}{-1}
	\sssec{Introduction.}
	In this section we study the 4-th point of intersection
	of the 4 circumcircles of a triangle.  The expression for these 6
	points is given in \ref{sec-tneucirc}.  A construction in the case
	where the centers of the 2 circles are given is described in
	\ref{sec-aneuc2circ}.

	\sssec{Notation.}\label{sec-nneucirc}
		${\cal C}_i$ denotes the circumcircle with center $C_i.$\\
		$X_{j,k}$ denotes the intersection of ${\cal C}_j$ and
		${\cal C}_k$
			distinct from the vertices of the triangle $A_i,$
			normalized to $A'_i.$\\
	\hth$	       a_i := A'_i \cdot A'_2,$\\
	\hth$	       d := A'_0 * (A'_1 * A'_2).$

	\sssec{Theorem.}\label{sec-tneucirc}
	\vspace{-18pt}\hspace{90pt}\footnote{Salzbourg-Innsbruck 29-30.9.83}
\\[8pt]
	{\em The intersections of the circumcircles of a triangle $A_i$ are
	given, using}\ref{sec-nneucirc}.\\
	$X_{0,3} = (1 + a_0)\: A'_0 + (a_2 - a_1)\: (A'_1 - A'_2),\\
	X_{1,3} = (1 + a_1)\: A'_1 + (a_0 - a_2)\: (A'_2 - A'_0),\\
	X_{2,3} = (1 + a_2)\: A'_2 + (a_1 - a_0)\: (A'_0 - A'_1),\\
	X_{1,2} = (1 - a_0)\: A'_0 + (a_1 + a_2)\: (A'_1 + A'_2),\\
	X_{2,0} = (1 - a_1)\: A'_1 + (a_2 + a_0)\: (A'_2 + A'_0),\\
	X_{0,1} = (1 - a_2)\: A'_2 + (a_0 + a_1)\: (A'_0 + A'_1).$

	Proof:  Let $B_i := A'_{i+1} * A'_{i-1},$ then\\
	\hth$	       C_0 = - B_0 + B_1 + B_2,$\\
	\hth$	       C_3 = B_0 + B_1 + B_2,$\\
	the circles with these centers are\\
	\hth$	       k X^2 = (X \cdot C_3)^2$ and $k X^2 = (X \cdot C_0)^2,$\\
	they pass through the vertex $A_0$ of the triangle if\\
	\hth$	       k -1 = (A_0 \cdot B_0)^2 = d^2$\\
	and therefore also through the vertices $A_1$ and $A_2$ if $k = - d^2.$
	$X$ is common to the 2 circles if $X \cdot C_0 = j X \cdot C_3$
		($j = +1$ or $-1$),\\
	$j = +1$ leads to the vertices $A'_1$ or $A'_2,$ $j = -1$ gives\\
	\hth$	       X \cdot (C_0 + C_3) = 2 X \cdot (B_1 + B_2) =$\\
	\hth$		       2 X \cdot (A'_0*(A'_1-A'_2)) = 0.$\\
	therefore, for some $l,$ and with $M^-_0 = A'_1 - A'_2,$\\
	\hth$X = l A'_0 + M^-_0,$\\
	hence\\
	\hth$	       X^2 = -l^2 + (M^-_0)^2 + 2l A'_0 \cdot M^-_0$, and\\
	\hth$	       X \cdot C_3 = (lA'_0 + M^-_0) \cdot (B_0+A_0*M^-_0) =
			       lA'_0 \cdot B_0 = l d.$\\
	$X$ is on the circles ${\cal C}_3$ if\\
	\hth$	       -d^{2(-l^2} + (M^-_0)^2 + 2lA'_0 \cdot M^-_0) = l^2 d^2,
$\\
	therefore\\
	\hth$	       l = - \frac{(M^-_0)^2 }{ 2 A'_0 \cdot M^-_0}
			 = -\frac{-2(1+a_0)}{2(a_2-a_1}$\\
	hence the expression for $X_{0,3}$.  The other points are derived
	similarly.

	\sssec{Corollary.}
	{\em The intersections $X_{1,3}$ and $X_{2,3}$ coincide if\\
	\hth$	       1 - a_0 + a_1 + a_2 = 0$\\
	and\\
	\hth$	       X_{1,3} = X_{2,3} = -A'_0 + A'_1 + A'_2,$\\
	with similar expressions for other pairs.}

	The proof is straightforward.

	\sssec{Example.}
	With the triangle of \ref{sec-neuctr},
	$a_i$ = 5,,\\
	$X_{0,3} =\\
	X_{1,3} =\\
	X_{2,3} =\\
	X_{1,2} =\\
	X_{2,0} =\\
	X_{0,1} =$

	\sssec{Theorem.}
	{\em
	Let $M_{0,0},$ $M_{0,1}$ be the mid-points of $a_0,$  $\ldots$  ,
	in the algebraic order defined above,
	then}
	\enumb
	\item{\em  the dual lines are the mediatrices,\\
	   the dual of $M_{0,1}$ passes through $M_{0,1}$ ,  $\ldots$}
	\item{\em  the points $M_{0,1},$ $M_{1,0},$ $M_{2,0}$ are on $o_0,$\\
	   the points $M_{0,0},$ $M_{1,1},$ $M_{2,0}$ are on $o_1,$\\
	   the points $M_{0,0},$ $M_{1,0},$ $M_{2,1}$ are on $o_2,$\\
	   the points $M_{0,1},$ $M_{1,1},$ $M_{2,1}$ are on $o_3.$}
	\item{\em  the dual of $o_i$ is the center $O_i$ of one of the 4
		circumcircles of the triangle $A_i.$}
	\item{\em  By duality, the mediatrices $m_{0,0},$ $m_{1,1},$ $m_{2,1}$
		are on $O_0,$ \ldots.

	    Let $m'_{0,0},$ $m'_{0,1},$ \ldots be the medians
		$A_0 \cdot M_{0,0},$ $A_0 \cdot M_{0,1},$ \ldots then}
	\item{\em the medians $m'_{0,0},$ $m'_{1,1},$ $m'_{2,1}$ are on $G_0,$\\
	   the medians $m'_{0,1},$ $m'_{1,0},$ $m'_{2,1}$ are on $G_1,$\\
	   the medians $m'_{0,1},$ $m'_{1,1},$ $m'_{2,0}$ are on $G_2,$\\
	   the medians $m'_{0,0},$ $m'_{1,0},$ $m'_{2,0}$ are on $G_3.$}
	\enume

	\sssec{Notation.}
	$u * v$ is the vector $(u_1 v_2 - u_2 v_1, u_2 v_0 - u_0 v_2,
				u_0 v_1 - u_1 v_0)$\\
	$u \bigx v$ is the vector $(u_1 v_2 + u_2 v_1, u_2 v_0 + u_0 v_2,$\\
	\hth$			       u_0 v_1 + u_1 v_0)$\\
	$u\: O\: v$ is the vector $(u_0 v_0, u_1 v_1, u_2 v_2))$

	\sssec{Algorithm.}\label{sec-aneuc2circ}
	Given two circles through the points $A_i,$ with centers $C_0$ and
		$C_1,$
	With $j$ in the set $\{0,1\}$ and $i$ in the set $\{0,1,2\},$
	and addition within the indices done modulo 2,\\
	\hth$	       d_j := C_j \: O \: C_j,$\\
	\hth$	       B_i := a_{i+1} X  a_{i+2},$\\
	\hth$	       L := d_0  *  d_1,$\\
	\hth$	       f_i := B_i   \cdot   L,$\\
	\hth$	       G := f \: O\:  f,$\\
	\hth$	       s_i := \frac{a_{i+1} a_{i+2}}{A_i},$\\
	\hth$	       O := \sum_{i=0}^3(s_i G_i A_i).$

	\sssec{Theorem.}
	{\em $O$ is the 4-th point common to the two circles.}



\setcounter{section}{1}
\setcounter{subsection}{8}
	\ssec{Other results in the geometry of the triangle.}

	\setcounter{subsubsection}{-1}
	\sssec{Introduction.}
	The following results were obtained while searching for a construction
	of $V,$ independent from the centers of mass and center of
	circumcircles.

	\sssec{Theorem.}
	{\em Let $I$ be an ideal point on the line $I \times B$.  Let\\
	\hth$	J := (B \cdot B)\: I - 2 (I \cdot B)\: B,$\\
	then $J$ is the other ideal point on $I \times B.$}

	\sssec{Example.}
	For $p = 13,$\\
	Let $I = (22) = (1,0,8)$ and $B = (4) = (0,1,3),$ then\\
	$J = -3(1,0,8) + 4(0,1,3) = (1,3,4) = (57).$\\
	The line $I \times J$ is $[48] = [1,2,8].$

	\sssec{Theorem.}
	{\em Let $a$ be an ordinary line and $B$ an ordinary point not on $a.$\\
	Let $I$ and $K$ be the ideal points on $a$ and $J$ and $L$ be the other
	ideal points on $I \times B$ and $K \times B,$ let $c := J * JJ,$ then\\
	\hth$c = (B \cdot B)\: a - 2 (a \cdot B)\: \ov{B}.$}

	Proof:  With $a = I * K,$\\
	\hth$J = (B \cdot B)\: I - 2(I \cdot B)\: B,\: L = (B \cdot B)\: K
		- 2(K \cdot B)\: B,$ hence\\
	\hth$L * J = (B \cdot B)^2\: K * I + 2(B \cdot B)( (K \cdot B)\: I
		- (I \cdot B)\: K ) * B,$\\
	\hth$= (B \cdot B) (-(B \cdot B)\: a + 2 ((K * I) * B) * B)$\\
	because of 2.3.17.0.,\\
	\hth$= (B \cdot B) (-(B \cdot B)\: a - 2 (a * B) * B),$\\
	\hth$= (B \cdot B) (-(B \cdot B)\: a + 2 (B \cdot B)\: a
		- 2 (a \cdot B)\: B)$\\
	because of 2.3.17.0., hence the Theorem.

	\sssec{Example.}
	For $p = 13,$ let $a = [139] = [1,9,8]$ and $B = (4),$ then
	$c = -3[1,9,8] - 1[0,1,3] = [1,5,9] = [88].$\\
	The ideal points are $I = 22,$ $J = 57,$ $K = 76,$ $L = 79.$

	\sssec{Definition.}
	$c$ as defined in the preceding theorem is called the
	{\em conjugate of} $a$ {\em with respect to} $B.$

	\sssec{Theorem.}
	{\em The lines\\
	\hth$x_{i+1} A_{i+1} - x_{i-1} A_{i-1} + y_{i-1} n_{i-1} a_{i-1}
		- y_{i+1} a_{i+1}$\\
	are concurrent at the point\\
	$\sum_i (y_{i+1} y_{i-1} t + (x_i d_{i+1} - x_{i-1} l_{i-1}) y_{i+1}
			+ (x_i d_{i-1} - x_{i+1} l_{i+1}) y_{i-1} )\: A_i\\
	\hth			+ \sum_i  (x_{i+1} x_{i-1} n_i)\: a_i.$}

	\sssec{Theorem.}
	{\em Given a triangle $A_i$ with sides $a_i,$ let $b_i$ be the
	conjugate of $a_i$ with respect to $A_i.$  Assume that the $b_i$ are not
	collinear.  Let $B_i$ be the vertices of the triangle $b_i,$ then}
	\enumb
	\item{\em$A_i \times B_i$ are concurrent at $W_0.$}
	\item{\em$A_i \times \ov{b}_i$ are concurrent at $H.$}
	\item{\em  $B_i \times \ov{b}_i$ are concurrent at $W_1.$\\
	Moreover,}
	\item $b_i = l_i A_{i+1} * A_{i-1} - 2 t A_i.$ 
	\item $B_i = - 3 l_{i-1} l_{i+1} A_i + 2 l_{i-1} d_{i-1} A_{i+1}$\\
	\hth$		       + 2 l_{i+1} d_{i+1} A_{i-1} + 4 t A_{i+1}
		* A_{i-1}.$
	\item $A_i * B_i = l_{i+1} d_{i+1} A_i * A_{i-1}
		- l_{i-1} d_{i-1} A_{i+1} * A_i$\\
	\hth$		       + 2 t d_{i+1} A_{i+1} - 2 t d_{i-1} A_{i-1}.$
	\item $W_0 = \sum_i  ( (-3 l_{i-1} l_{i+1} d_{i-1} d_{i+1}) A_i
			+ 2 d_i(l_{i-1} d_{i-1}^2
		+ l_{i+1} d_{i+1}^{2)} A_i$ \\
	\hth$		       + 4 t d_{i+1} d_{i-1} A_{i+1} * A_{i-1} ).$
	\item $A_i \times \ov{b}_i = d_{i+1} A_{i+1} - d_{i-1} A_{i-1}.$
	\item $H = d_1 d_2 A_1 * A_2 + d_2 d_0 A_2 * A_0 + d_0 d_1 A_0 * A_1.$ 
	\item {\em With\\
	\hth$x_i = 8 t^2 d_i + 2 l_{i+1} l_{i-1} d_{i+1} d_{i-1}
		- l_0 l_1 l_2 d_i$ and}\\
	\hth$y_i = 4 t l_i d_i,$\\
	\hth$	B_i \times \ov{b}_i = x_{i+1} A_{i+1} - x_{i-1} A_{i-1}
		+ y_{i-1} A_i * A_{i+1} - y_{i+1} A_{i-1} * A_i.$ 
	\item 	$W_1 =$
	\enume

	Proof:  3., is immediate,\\
	Using 2.3.17.0 and 1.7.2., we obtain .4 after division by $t,$\\
	5. and 6. after division by $t$, 7. after division by $l_i.$ \\
	8. is immediate and is indeed the same as 1.3.3.

	\sssec{Example.}\label{sec-eneucw}
	For $p = 13,$ using the same triangle as Example 1.7.4.\\
	$B_i = \{ 9 = (0,1,8), 137 = (1,9,6), 61 = (1,3,8) \},$\\
	$b_i = \{ 180 = [1,12,10], 101 = [1,6,9], 80 = [1,5,1] \}$,\\
	$W_0 = 7 = (0,1,6),$ $H = 171 = (1,12,1),$ $W_1 = 77 = (1,4,11)$.

	\sssec{Exercise.}
	Using the dual triangle, $A_i = \{175, 53, 1\},$
	determine $B_i,$ $b_i$ and $W_0,$ $H$ and $W_1.$

	\ssec{Circumcircle of a triangle with at least one ideal vertex.}

	\setcounter{subsubsection}{-1}
	\sssec{Introduction.}
	In the preceding section I have dealt with circumcircles
	through 3 ordinary points.  I will now discuss the case when 1 or more
	points are ideal points.

	\sssec{Theorem.}
	{\em The only circle through 3 distinct ideal points is the
	ideal conic, $X^2$ = 0.}

	\sssec{Theorem.}
	{\em The only circle through 2 distinct ideal points $A$ and $B$ and
	through an ordinary point C is\\
	$(c \cdot C)^2 X^2 - (c \cdot X)^2 C^2 = 0,$ where\\
	\hth$	c := A \times B.$}

	\sssec{Example.}
	$p = 11,$ $A = (31) = (1,3,2),$ $B = (26) = (1,2,4),$\\
	$C = 54 = (1,6,4).$  The circumcircle through $A,$ $B$ and $C$ is\\
	\hth$	       X_0^2 + X_1^2 + X_2^2 = 4 (X_0 - 2 X_1 - X_2)^2.$

	\sssec{Theorem.}
	{\em There are 2 circles through 1 ideal point $A$ and through
	two distinct ordinary points $B$ and $C$ not collinear with $A.$}

	\sssec{Theorem.}\label{sec-tneucceva}
	{\em Let $X$ be a point not on the sides of a triangle $A_i.$  Let
	$X_i$ be the intersection with $a_i$ of the line $A_i \times X$, or\\
	\hth$	X_i := (A_i \times X) \times a_i,$\\
	Let $\xi$  be a circle through $X_i.$ \\
	Let $Y_i$ be the other intersection of $a_i$ with $\xi$ ,\\
	let $y_i := A_i \times Y_i,$ \\
	then the lines $y_i$ have a point $Y$ in common.}

	\sssec{Comment.}
	The theorem \ref{sec-tneucceva} is analogous to theorem $\ldots$  in
	Euclidean geometry.  It follows from its generalization $\ldots$ 
	to projective geometry.\\
	Because of the clear connection with the Theorem of Ceva, I have
	the following:

	\sssec{Definition.}
	The correspondence $X$ to the various points $Y$ associated to
	the several circumcircles through $X_i,$ is called the {\em Ceva}
	{\em correspondence}.  I will ignore those points $Y$ which happen to
	coincide with a vertex of the triangle.

	\sssec{Comment.}
	Clearly if $Y$ is associated to $X,$ $X$ is associated to $Y$ in a
	Ceva correspondence.  But we cannot call this an involution because the
	correspondence is not one to one or bijective.

	\sssec{Program.}
	The program NETR1.BAS determines the Ceva correspondence.
	It is illustrated in NETR1.HOM.



\setcounter{section}{1}
\setcounter{subsection}{10}
	\ssec{The parabola in polar geometry.}

	\setcounter{subsubsection}{-1}
	\sssec{Introduction.}
	In this section I have defined a parabola for non
	Euclidean geometry and many of the related elements of the parabola,
	by analogy with the definitions of Euclidean geometry.  By duality
	we essentially double the number of these elements, for instance to
	the focus in Euclidean geometry corresponds the focal point and the
	focal line.  The basic equation is given by \ref{sec-tneucparab}.0.

	\sssec{Definition.}
	A {\em parabola} is a conic which is tangent at one point to the ideal
	conic and at one point only.  This point is called the {\em isotropic
	point of the parabola}, the tangent is called the {\em isotropic
	line of the parabola}.

	\sssec{Definition.}
	A {\em focal tangent} $t1$ of a parabola is an ideal tangent to the
	parabola which is not isotropic.\\
	A {\em focal point} $F1$ is an ideal point which is not isotropic.
	There are 2 focal points $F1$ and $F2$ which are either real or
	complex conjugate.

	\sssec{Definition.}
	The {\em focus} $F$ of a parabola is the intersection of the focal
	tangents.  The {\em focal line} $f$ of a parabola is the line through
	the focal points.

	\sssec{Theorem.}\label{sec-tneucfocus}
	{\em The focus is not on the isotropic line.  The focus is not an
	ideal point.}

	Proof:  In the first case, through the focus we could draw 3 tangents to
	the parabola.  In the second case, the parabola would be tangent at a
	second point to the ideal conic and would therefore be a circle.

	\sssec{Definition.}
	The {\em director} $D$ of a parabola is the pole of its focal line with
	respect to the parabola.

	\sssec{Definition.}
	The {\em axis} $a$ of a parabola is the line through its focus and
	its isotropic point.\\
	The {\em axial point} $A$ of a parabola is the point on the focal line
	and on the isotropic line.

	\sssec{Definition.}
	The {\em vertex} $V$ of a parabola is the ordinary point on the
	parabola and its axis.\\
	The {\em vertical line} $v$ of a parabola is the ordinary tangent
	through the axial point.

	\sssec{Theorem.}\label{sec-tneucparab}
	{\em A parabola with isotropic point $I$ and focal tangent $f$ is}
	\enumb
	\item $2 (I \cdot X) (f \cdot X) = t (X \cdot X),$
		$f \cdot I\neq  0,$ $f \cdot f\neq 0,$ $t\neq  0,$
		 $t\neq f \cdot I$.

	{\em The polar of $X0$ is}
	\item $(I \cdot X0) f + (f \cdot X0) I - t X0.$

	{\em The pole of $a$ is}
	\item $(f*I) \cdot a I*f + t (a \cdot f) I + t (a \cdot I) f 
		+ (t^2 - 2t I \cdot f) a.$
	\enume

	Proof:  0,  follows from the general equation of a conic through
	the intersections of the ideal and the lines $\overline{I}$ and $f,$
	see $\ldots$  and
	from \ref{sec-tneucfocus}.11.4.\\
	0, represents a degenerate conic corresponding to the
	lines $\overline{I}$ and $f$ if $t = 0$ and to the lines $I \times F1$
	and $I \times F2$ if $t = f \cdot I$.\\
	The proof of the last fact is left as an exercise.
	1 follows from 0.  See $\ldots$ .
	2 is obtained by chosing 2 points on $a,$ $a * I$ and $f * a,$ and
	determining the intersection of the polars of these lines.\\
	$(f * I) \cdot a$ is a factor of each term.

	\sssec{Example.}
	For $p = 13,$ $I = (1,5,0),$ $f = [1,1,4],$ $t = 4,$\\
	The parabola is\\
	\hth$	       2 (x + 5y)(x + y + 4z) = 4 (x^2 + y^2 + z^{2)}.$\\
	The polar of $X0 =(x0,y0,z0)$ is\\
	\hth$	       [- x0 + 3y0 + 2z0, 3x0 + 3y0 - 3z0, 2x0 - 3y0 - 2z0].$\\
	The pole of $a = [a0,b0,c0]$ is\\
	\hth$	       [a0 + c0, b0 + 5c0, a0 + 5b0 + 6c0].$\\
	The director $D$ is (1,-1,6).\\
	The axial point $A$ is $I * f = (1,5,5).$\\
	With $v^2 = -2,$ the ideal points on the parabola are $I$ and\\
	\hth$	       F1 = (1, 3 + 3v, -1 - 4v), F2 = (1, 3 - 3v, -1 + 4v).$\\
	The tangent at $F1$ is\\
	\hth$f1 = (v - 6) [6 + v, 2 - 5v, -5 -v] = [1, -2 + 6v, 6 + v]$\\
	The ideal lines $t1 = [1,3,4]$ and $t2 = [1,-5,0]$ are tangent to the
	parabola at $T1 = (1,2,-5)$ and $T2 = (2,-5,2),$ they meet at the
	focus $F = (1,-5,-3).$\\
	The directrix is [1,4,4].\\
	The axis $a$ is $I * F = [1,5,5].$\\
	The vertex $V$ is (1,-2,-6).\\
	The tangent at the vertex is [1,4,1].\\
	The vertical line $v$ is [1,4,1].

	\sssec{Theorem.}
	{\em The polar of $X$ is $\sum_j  A_{i,j} X_j,$ with\\
	\hth$       A_{j,j} = 2 I_j f_j - t,$\\
	\hth$       A_{j,k} = I_j f_k + I_k f_j, j\neq  k,$\\
	The pole of $a$ is $\sum_k B_{j,k} a_k,$ with\\
	\hth$       B_{j,j} = t^2 - 2t(I \cdot f - I_j f_j) - (I * f)_j^2,$\\
	\hth$       B_{j,k} = t (I_j f_k + I_k f_j) - (I * f)_j (I * f)_k.$\\
	The dual equation of the conic is\\
	\hth$	       2(I \cdot x) (F \cdot x) = u(x \cdot x),$\\
	The pole of $x$ is $(I \cdot x) F + (F \cdot x) I - u x$ and $F$ and
	$u$ follow from \ldots\\
	if $I_j\neq frac{1}{0},$ $u = (I_0B_{0,0} + I_1 B_{1,1},$
		$- 2 I_0 I_1 B_{0,1}),$\\
	\hth$	       F_j = \frac{u+B_{j,j}}{2 I_j},$\\
	if $I_2$ = 0, and $I_1\neq 0,$ $u = - B_{2,2},$\\
	\hth$       F_0 = \frac{B_{0,0}-B_{2,2}}{2 I_0},$\\
	\hth$       F_1 = \frac{B_{1,1}-B_{2,2}}{2 I_1},$\\
	\hth$       F_2 = \frac{B_{0 \cdot 2} }{ I_0},$\\
	if $I = (1,0,0),$ $u = -B_{1,1}$,\\
	\hth$F_0 = \frac{B_{0,0}-B_{1,1}}{2},$ $F_1 = B_{0,1},$
	$F_2 = B_{0,2}$.}

	Proof:  The matrix $A$ follow from (1), the matrix $B$ is its adjoint
	divided by 2.

	\sssec{Theorem.}
	{\em The director is\\
	\hth$	       D = (f \cdot f) I + (t - f \cdot I) f.$}

	Proof:  Replace $a$ by $f$ in \ref{sec-tneucparab}.2.

	\sssec{Exercise.}
	Complete the following sentences, for the parabola
	\ref{sec-tneucparab}.0:
	\enumb
	\setcounter{enumi}{-1}
	\item The axis $a$ is  .
	\item The axial point $A$ is  .
	\item The vertex $V$ is  .
	\item The vertical line $v$ is  .
	\item The ideal point $J$ is  .
	\enume

	\sssec{Theorem.}
	\enumb
	\item{\em  The vertical line $v$ passes trough the vertex $V.$}
	\item {\em The director $D$, the pole $\overline{f}$ of the focal line
		$f,$ the pole $\overline{d}$ of the directrix $d$ are all on the
		axis $a.$}
	\item{\em  The directrix $d,$ the polar $\overline{F}$ of the focus $F,$
	the polar $\overline{D}$ of
		the director $D$ all pass through the axial point $A.$}
	\item{\em  The other ideal point $J$ on the axis and the other ideal
		line $j$ through the axial point are incident.}
	\enume

	\sssec{Answer }
	\vspace{-18pt}\hspace{94pt}{\bf \ref{sec-eneucw}.}\\
	$B_i = \{144, 64, 88\},$ $ b_i = \{182, 136, 132\},$\\
	$W_0 = (62),$ $H = (171),$ $W_1 = (88).$



\setcounter{section}{1}
\setcounter{subsection}{11}
	\ssec{Representation of polar geometry on the dodecahedron.}

	\setcounter{subsubsection}{-1}
	\sssec{Introduction.}
	When $p = 5,$ the representation of polar geometry
	on the dodecahedron is suggested by the fact that the 6 faces form
	a conic which can be chosen as the ideal.

	\sssec{Definition.}
	Using the dodecahedral representation, the conic which
	consists of the 6 face-points is the {\em ideal conic}.

	\sssec{Theorem.}
	{\em The 15 side-points are hyperbolic and the 10 vertex-points
	are elliptic.}

	This follows at once from the incidence definitions,
	II.\ref{sec-ddodeca}.

	\sssec{Theorem.}
	{\em With the ideal conic of type $A$,\\
	the 3 . 15 conics of type $I3,$ $E1$ and $E2$ are hyperbolic circles,\\
	the 4 . 6 conics of type $J1,$ $J2,$ $O1,$ $O2$ are parabolic circles
	and\\
	the 3 . 10 conics of type $P,$ $U1$ and $U2$ are elliptic circles.}

	Although this should be placed in the Chapter on non-Euclidean
	geometry, we have.

	\sssec{Theorem.}
	{\em With a particular choice of unit, the radii of the various
	sub-types are as follows,\\
	$U1$ and $E2$ are $\frac{\pi}{6},$ $U2$ and $E1$ are $\frac{\pi}{3},$
	$P$ and $I3$ are $\frac{\pi}{3}.$ }

	\sssec{Example.}
	Computations relating to g434.PRN:\\
	If we use as primitive polynomial $I^3-I-2,$ we obtain the
	correspondence:\\
	$\begin{array}{rlrr}
	x&(y_0,y_1,y_2)&(y)&[z]\\
	0&(0,0,1)&(0)	&[6]\\
	1&(0,1,0)&(1)	&[1]\\
	2&(1,0,0)&(6)     &[11]\\
	3&(0,1,2)&(3)     &[21]\\
	4&(1,2,0)&(16)    &[16]\\
	5&(1,3,1)&(22)    &[26]\\
	6&(1,4,4)&(30)    &[7]\\
	7&(1,0,3)&(9)     &[9]\\
	8&(0,1,3)&(4)     &[8]\\
	9&(1,3,0)&(21)    &[10]\\
	10&(1,2,4)&(20)    &[0]\\
	11&(1,0,1)&(7)     &[12]\\
	12&(0,1,1)&(2)     &[24]\\
	13&(1,1,0)&(11)    &[18]\\
	14&(1,1,2)&(13)    &[30]\\
	15&(1,3,2)&(23)    &[2]\\
	16&(1,1,4)&(15)    &[17]\\
	17&(1,0,2)&(8)     &[14]\\
	18&(0,1,4)&(5)     &[28]\\
	19&(1,4,0)&(26)    &[25]\\
	20&(1,4,3)&(29)    &[4]\\
	21&(1,1,3)&(14)    &[22]\\
	22&(1,4,2)&(28)    &[29]\\
	23&(1,2,3)&(19)    &[13]\\
	24&(1,2,1)&(17)    &[20]\\
	25&(1,1,1)&(12)    &[3]\\
	26&(1,2,2)&(18)    &[27]\\
	27&(1,4,1)&(27)    &[19]\\
	28&(1,3,3)&(24)    &[23]\\
	29&(1,3,4)&(25)    &[15]\\
	30&(1,0,4)&(10)    &[5]\\
	\end{array}$

	The second column is $I^x \umod{P},$ the third column is the
	representation\\
	of Chapter II, the fourth column is obtained as follows.

	The ideal conic passes through 0,4,6,9,16 and 17 and is therefore
	represented by the matrix\\
	\hth$Q = \matb{(}{ccc}1&0&1\\0&1&0\\1&0&0\mate{)}$\\
	or the quadratic form $Q(x,y) = x_0 y_0 + x_1 y_1 + x_0 y_2 + x_2 y_0.$ 

	This determines the polar of $(y_0.y_1,y_2)$ as $[y_0+y_2,y_1,y_0],$ 
	but the polar of $x$ is $x^*$, therefore\\
	if $x = (y_0.y_1,y_2)$ and $[z] = [y_0+y_2,y_1,y_0]$ then $x^* = [z].$\\
	For instance, if $x = 7,$ $(y) = (1,0,3),$ $[z] = [-1,0,1] = [10].$

	\sssec{Example.}
	The hyperbolic circles have as center a edge-point.\\
	Those with center $8 = 0^* \times 8^*$ and through the point $u$ can be
	constructed as follows, let $up = 0 \times (4 \times u),$ given any
	line $v$
	through 0, such that $u\cdot v\neq 0$, the Pascal construction gives the
	other point on the conic and $v$ using\\
	\hth$	((((v \times 4) \times up) \times 8 ) \times u ) \times v.$\\
	This gives, with the radii determined below:\\
	0,4;;2,15,22,25 of type $ffssss$ and sub-type $I3,$ radius
		$\frac{\pi}{3}.$\\
	0,4;;5,11,13,20 of type $ffvvvv$ and sub-type $E1,$ radius
		$\frac{\pi}{4}.$\\
	0,4;;7,19,21,29 of type $ffvvvv$ and sub-type $E2.$ radius
		$\frac{\pi}{6}.$

	Before having obtained a synthetic construction of the parabolic and
	elliptic circles we have used the algebraic definition.\\
	The algebraic definition is\\
	\hth$	k x^T Q x - (x^T Q C)^2 = 0,$\\
	with $C$ on the ideal conic, for parabolic circles and $C$ a
	vertex-point
	for elliptic circles.\\
	For $k = 0,$ the circle degenerates in  (a double)
	line, consisting of the points at distance $\frac{\pi}{2}$ from $C.$

	Let $C = 0,$ the polar $0^* = [1,0,0],$ hence the parabolic circles
	are\\
	\hth$	k (x_0^2 - x_1^2 + 2 x_2 x_0) + x_0^2 = 0.$\\
	With $k' = \frac{1}{k},$ the points are (0,0,1) = 0, and
	$(1,x_1,2(1-k'+x_1)^2.$  This
	gives for,\\
	$k =-1,$ (1,0,4) = 30, (1,1,1) = 25, (1,4,1) = 27, (1,2,2) = 26,
		(1,3,2) = 15, in view of the table above.  Hence\\
 	$k = 4,$ \{0,30,25,27,26,15\} of type $fsssss$ and sub-type $O1.$\\
	$k = 2,$ \{0,11,21,20,10,29\} of type $fvvvvv$ and sub-type $J2.$\\
	$k = 3,$ \{0,7,13,19,24,5\} of type $fvvvvv$ and sub-type $J1.$\\
	$k = 1,$ \{0,2,14,22,23,28\} of type $fsssss$ and sub-type $O2.$

	Let $C = 5,$ the polar $5^* = [1,-1,-2],$ hence the circles are\\
	\hth$	k (x_0^2 + x_1^2 + 2 x_2 x_0) - (x_0 - x_1 -2 x_2)^2 = 0.$\\
	The points are,\\
	$(0,1,2 \pm \sqrt{-k}$ and $(1,x_1,
		-k + 2(x_1-1) \pm \sqrt{k^2-k(x_1^2-x_1+2)}.$ \\
	This gives, for $k = -1,$\\
	(0,1,3) = 8, (0,1,1) = 12, (1,4,2) = 22, (1,2,3) = 23, (1,3,2) = 15,\\
	\hth$	(1,3,3) = 22,$ in view of the table above.\\
	Hence, with the radii determined below:\\
	$k = 4,$ \{8,12,22,23,15,28) of type $ssssss$ and sub-type $U2,$ radius 
	$\frac{\pi}{4}.$\\
	$k = 2,$ \{11,21,24,10,19,20\} of type $vvvvvv$ and sub-type $P,$ radius
	$\frac{\pi}{3}$\\
	$k = 3,$ \{5\}, radius 0.\\
	$k = 1,$ \{1,18,2,30,14,25\} of type $ssssss$ and sub-type $U1,$ radius
	$\frac{\pi}{6}.$

	To summarize, we see that, with the ideal conic of type $A,$\\
	the 3 . 15 conics of type $I3,$ $E1$ and $E2$ are hyperbolic circles,\\
	the 4 . 6 conics of type $J1,$ $J2,$ $O1,$ $O2$ are parabolic circles
	and\\
	the 3 . 10 conics of type $P,$ $U1$ and $U2$ are elliptic circles.

	\sssec{Exercise.}
	For a synthetic construction of the parabolic circles and some elliptic
	ones, we can use IV.1.2.7.  This is a good exercise.

	\sssec{Example of Distances.}
	We recall the trigonometric tables for p = 5:\\
	With $\delta$  = 2,\\
	$\begin{array}{rcccrcc}
		x&	sin(x)&	cos(x)&\vline&	x&	sin(x)&	cos(x)\\
		0&	0&	1&\vline&		0&	0&	1\\
		1&	2\delta&2\delta&\vline&	1&	-2&	\delta	\\
		2&	1&	0&	\vline&	2&	\delta&	-2\\
		&	&	&	\vline&	3&	1&	0
	\end{array}\\
	\begin{array}{ccccccc}
		cos^2(d)&\vline		&0&1&2&3&4\\
		\frac{d}{\pi}&\vline	&\frac{1}{2}&0&\frac{1}{6}&\frac{1}{4}
			&\frac{1}{3}
	\end{array}$\\
	where
	\hth$cos(d(C,X)) = \frac{Q(X,C)^2}{Q(X,X)} = \frac{k}{Q(C,C)}.$

	Hence the distances, recorded above.  For instance, in the case of
	elliptic circles, for $C = 5 = (1,3,1)$ and $X = 8 = (0,1,3),$
	$Q(C,C) = 2,$
	$Q(X,X) = 1,$ $Q(X,C) = 1,$ $cos^2(d(C,X)) = 3,$
	$d(C,X) = \frac{\pi}{4}.$ Hence\\ 
	for $k = -1,$ the radius is $\frac{\pi}{4},$ \\
	for $k = 2,$ $cos^2(d(C,X)) = -1,$ $d(C,X) = \frac{\pi}{3},$\\
	for $k = 3,$ $d(C,X) = 0,$\\
	for $k = 1,$ $d(C,X) = \frac{\pi}{6}.$ \\
	In the case of hyperbolic circles with $C = 8$ and $X = 5,$ we have
	$d(8,5) = \frac{\pi}{4},$ the other radii are obtained directly,\\
	$cos^2(d(8,2))	= \frac{(-2)^2}{1 . 1} = -1,$ hence
	$d(8,2) = \frac{\pi}{4},$\\
	$cos^2(d(8,7))	= \frac{(-2)^2}{1 . 2} = 2,$ hence
	$d(8,7) = \frac{\pi}{6}.$ 



\setcounter{section}{1}
\setcounter{subsection}{11}
\setcounter{subsubsection}{14}
	\subsubsection{Example.}
	Computations relating to g434.PRN:
	If we use as primitive polynomial $I^3-I-2,$ we obtain the
	correspondence:\\
	$\begin{array}{rlrr}
	x&(y_0,y_1,y_2)&(y)&[z]\\
	0&(0,0,1)&(0)	&[6]\\
	1&(0,1,0)&(1)	&[1]\\
	2&(1,0,0)&(6)     &[11]\\
	3&(0,1,2)&(3)     &[21]\\
	4&(1,2,0)&(16)    &[16]\\
	5&(1,3,1)&(22)    &[26]\\
	6&(1,4,4)&(30)    &[7]\\
	7&(1,0,3)&(9)     &[9]\\
	8&(0,1,3)&(4)     &[8]\\
	9&(1,3,0)&(21)    &[10]\\
	10&(1,2,4)&(20)    &[0]\\
	11&(1,0,1)&(7)     &[12]\\
	12&(0,1,1)&(2)     &[24]\\
	13&(1,1,0)&(11)    &[18]\\
	14&(1,1,2)&(13)    &[30]\\
	15&(1,3,2)&(23)    &[2]\\
	16&(1,1,4)&(15)    &[17]\\
	17&(1,0,2)&(8)     &[14]\\
	18&(0,1,4)&(5)     &[28]\\
	19&(1,4,0)&(26)    &[25]\\
	20&(1,4,3)&(29)    &[4]\\
	21&(1,1,3)&(14)    &[22]\\
	22&(1,4,2)&(28)    &[29]\\
	23&(1,2,3)&(19)    &[13]\\
	24&(1,2,1)&(17)    &[20]\\
	25&(1,1,1)&(12)    &[3]\\
	26&(1,2,2)&(18)    &[27]\\
	27&(1,4,1)&(27)    &[19]\\
	28&(1,3,3)&(24)    &[23]\\
	29&(1,3,4)&(25)    &[15]\\
	30&(1,0,4)&(10)    &[5]\\
	\end{array}$

	The second column is $I^x \umod{P},$ the third column is the
	representation\\
	of Chapter II, the fourth column is obtained as follows.

	The ideal conic passes through 0,4,6,9,16 and 17 and is therefore
	represented by the matrix\\
	\hth$Q = \matb{(}{ccc}1&0&1\\0&1&0\\1&0&0\mate{)}$\\
	or the quadratic form $Q(x,y) = x_0 y_0 + x_1 y_1 + x_0 y_2 + x_2 y_0.$ 

	This determines the polar of $(y_0.y_1,y_2)$ as $[y_0+y_2,y_1,y_0],$ 
	but the polar of $x$ is $x^*$, therefore\\
	if $x = (y_0.y_1,y_2)$ and $[z] = [y_0+y_2,y_1,y_0]$ then $x^* = [z].$\\
	For instance, if $x = 7,$ $(y) = (1,0,3),$ $[z] = [-1,0,1] = [10].$

	\subsubsection{Example.}
	The hyperbolic circles have as center a edge-point.\\
	Those with center $8 = 0^* \times 8^*$ and through the point $u$ can be
	constructed as follows, let $up = 0 \times (4 \times u),$ given any
	line $v$
	through 0, such that $u\cdot v\neq 0$, the Pascal construction gives the
	other point on the conic and $v$ using\\
	\hth$	((((v \times 4) \times up) \times 8 ) \times u ) \times v.$\\
	This gives, with the radii determined below:\\
	0,4;;2,15,22,25 of type $ffssss$ and sub-type $I3,$ radius
		$\frac{\pi}{3}.$\\
	0,4;;5,11,13,20 of type $ffvvvv$ and sub-type $E1,$ radius
		$\frac{\pi}{4}.$\\
	0,4;;7,19,21,29 of type $ffvvvv$ and sub-type $E2.$ radius
		$\frac{\pi}{6}.$

	Before having obtained a synthetic construction of the parabolic and
	elliptic circles we have used the algebraic definition.\\
	The algebraic definition is\\
	\hth$	k x^T Q x - (x^T Q C)^2 = 0,$\\
	with $C$ on the ideal conic, for parabolic circles and $C$ a
	vertex-point
	for elliptic circles.\\
	For $k = 0,$ the circle degenerates in  (a double)
	line, consisting of the points at distance $\frac{\pi}{2}$ from $C.$

	Let $C = 0,$ the polar $0^* = [1,0,0],$ hence the parabolic circles
	are\\
	\hth$	k (x_0^2 - x_1^2 + 2 x_2 x_0) + x_0^2 = 0.$\\
	With $k' = \frac{1}{k},$ the points are (0,0,1) = 0, and
	$(1,x_1,2(1-k'+x_1)^2.$  This
	gives for,\\
	$k =-1,$ (1,0,4) = 30, (1,1,1) = 25, (1,4,1) = 27, (1,2,2) = 26,
		(1,3,2) = 15, in view of the table above.  Hence\\
 	$k = 4,$ \{0,30,25,27,26,15\} of type $fsssss$ and sub-type $O1.$\\
	$k = 2,$ \{0,11,21,20,10,29\} of type $fvvvvv$ and sub-type $J2.$\\
	$k = 3,$ \{0,7,13,19,24,5\} of type $fvvvvv$ and sub-type $J1.$\\
	$k = 1,$ \{0,2,14,22,23,28\} of type $fsssss$ and sub-type $O2.$

	Let $C = 5,$ the polar $5^* = [1,-1,-2],$ hence the circles are\\
	\hth$	k (x_0^2 + x_1^2 + 2 x_2 x_0) - (x_0 - x_1 -2 x_2)^2 = 0.$\\
	The points are,\\
	$(0,1,2 \pm \sqrt{-k}$ and $(1,x_1,
		-k + 2(x_1-1) \pm \sqrt{k^2-k(x_1^2-x_1+2)}.$ \\
	This gives, for $k = -1,$\\
	(0,1,3) = 8, (0,1,1) = 12, (1,4,2) = 22, (1,2,3) = 23, (1,3,2) = 15,\\
	\hth$	(1,3,3) = 22,$ in view of the table above.\\
	Hence, with the radii determined below:\\
	$k = 4,$ \{8,12,22,23,15,28) of type $ssssss$ and sub-type $U2,$ radius 
	$\frac{\pi}{4}.$\\
	$k = 2,$ \{11,21,24,10,19,20\} of type $vvvvvv$ and sub-type $P,$ radius
	$\frac{\pi}{3}$\\
	$k = 3,$ \{5\}, radius 0.\\
	$k = 1,$ \{1,18,2,30,14,25\} of type $ssssss$ and sub-type $U1,$ radius
	$\frac{\pi}{6}.$

	To summarize, we see that, with the ideal conic of type $A,$\\
	the 3 . 15 conics of type $I3,$ $E1$ and $E2$ are hyperbolic circles,\\
	the 4 . 6 conics of type $J1,$ $J2,$ $O1,$ $O2$ are parabolic circles
	and\\
	the 3 . 10 conics of type $P,$ $U1$ and $U2$ are elliptic circles.

	\subsubsection{Exercise.}
	For a synthetic construction of the parabolic circles and some elliptic
	ones, we can use IV.1.2.7.  This is a good exercise.

	\subsubsection{Example of Distances.}
	We recall the trigonometric tables for p = 5:\\
	With $\delta$  = 2,\\
	$\begin{array}{rcccrcc}
		x&	sin(x)&	cos(x)&\vline&	x&	sin(x)&	cos(x)\\
		0&	0&	1&\vline&		0&	0&	1\\
		1&	2\delta&2\delta&\vline&	1&	-2&	\delta	\\
		2&	1&	0&	\vline&	2&	\delta&	-2\\
		&	&	&	\vline&	3&	1&	0
	\end{array}\\
	\begin{array}{ccccccc}
		cos^2(d)&\vline		&0&1&2&3&4\\
		\frac{d}{\pi}&\vline	&\frac{1}{2}&0&\frac{1}{6}&\frac{1}{4}
			&\frac{1}{3}
	\end{array}$\\
	where
	\hth$cos(d(C,X)) = \frac{Q(X,C)^2}{Q(X,X)} = \frac{k}{Q(C,C)}.$

	Hence the distances, recorded above.  For instance, in the case of
	elliptic circles, for $C = 5 = (1,3,1)$ and $X = 8 = (0,1,3),$
	$Q(C,C) = 2,$
	$Q(X,X) = 1,$ $Q(X,C) = 1,$ $cos^2(d(C,X)) = 3,$
	$d(C,X) = \frac{\pi}{4}.$ Hence\\ 
	for $k = -1,$ the radius is $\frac{\pi}{4},$ \\
	for $k = 2,$ $cos^2(d(C,X)) = -1,$ $d(C,X) = \frac{\pi}{3},$\\
	for $k = 3,$ $d(C,X) = 0,$\\
	for $k = 1,$ $d(C,X) = \frac{\pi}{6}.$ \\
	In the case of hyperbolic circles with $C = 8$ and $X = 5,$ we have
	$d(8,5) = \frac{\pi}{4},$ the other radii are obtained directly,\\
	$cos^2(d(8,2))	= \frac{(-2)^2}{1 . 1}$ = -1, hence
	$d(8,2) = \frac{\pi}{4},$\\
	$cos^2(d(8,7))	= \frac{(-2)^2}{1 . 2}$ = 2, hence
	$d(8,7) = \frac{\pi}{6}.$ 



	\section{Finite Non-Euclidean Geometry.}

	\setcounter{subsection}{-1}
	\ssec{Introduction.}

	\ssec{Trigonometry for the general triangle.}

	\setcounter{subsubsection}{-1}
	\sssec{Introduction.}
	Spherical trigonometry refers to the relation between the
	measure of angles and arcs of a triangle on a sphere.\\
	The formulas of al-Battani (Albategnius, about 920 A.D.) and of
	Jabir ibn Aflah (Geber, about 1130 A.D.) have to be adapted to the
	finite case in which the sine of an angle in the first 2 quadrants
	cannot be considered as positive.  There are several possible
	solutions.  One of these will be given in Theorem \ref{sec-tneuctrig}.\\
	Let $A,$ $B,$ $C$ be the vertices of a triangle, $a,$ $b,$ $c$ be its
	sides.\\
	The measure of the angle between $b$ and $c$ will be denoted $A,$ 
	$\ldots$  .\\
	The distance between the points $B$ and $C$ will be denoted $a,$
	$\ldots$  .

	\sssec{Definition.}
	Let $A,$ $B,$ $C$ be 3 points on the sphere\\
	\hth$	x^2 + y^2 + z^2 = 1,$\\
	of center $O = (0,0,0,1).$\\
	The {\em direction} $DA$ of $A$ is the ideal point on $OA.$\\
	The {\em side} $a = \{B,C\}$ is a section of the circle  ${\cal C}_a$
	which is the intersection of the sphere
	and the plane $O \times B \times C.$
	The {\em spherical distance} of the side $a$, also denoted  $a$  is
	the angle between the directions $DB$ and $DC.$\\
	The {\em angle} $B A C,$ also denoted  $A$
	is the angle of the directions of the tangents at $A$ to the circles
	${\cal C}_b$ and ${\cal C}_c$.

	\sssec{Theorem.}\label{sec-tneuctrig}
	{\em Between the trigonometric functions of the angles and sides
	of a general triangle we have the relations:}

	\enumb
	\item ~~~$ \frac{| sin a| }{sin A} = \frac{| sin b|}{sin B}
		 = \frac{| sin c|}{sin C} = r.$

	\item 0. $cos A = cos B cos C + sin B sin C cos a,$\\
	      1. $cos B = cos C cos A + sin C sin A cos b,$\\
	      2. $cos C = cos A cos B + sin A sin B cos c.$

	\item 0. $cos a = cos b cos c + | sin b|  | sin c|  cos A,$\\
	      1. $cos b = cos c cos a + | sin c|  | sin a|  cos B,$\\
	      2. $cos c = cos a cos b + | sin a|  | sin b|  cos C.$

	\item 0. $sin A = \frac{cos^2 B - cos^2 C}
			{sin B cos C cos c - sin C cos B cos b},$\\
	      1. $sin B = \frac{(cos^2 C - cos^2 A}
			{sin C cos A cos a - sin A cos C cos c},$\\
	      2. $sin C = \frac{(cos^2 A - cos^2 B}
			{sin A cos B cos b - sin B cos A cos a}.$

	\item 0. $| sin a|  = \frac{cos^2 b - cos^2 c}
			{| sin b|  cos c cos C - | sin c|  cos b cos B},$\\
	      1. $| sin b|  = \frac{cos^2 c - cos^2 a}
			{| sin c|  cos a cos A - | sin a|  cos c cos C},$\\
	      2. $| sin c|  = \frac{cos^2 a - cos^2 b}
			{| sin a|  cos b cos B - | sin b|  cos a cos A}.$

	\item 0. $cos A = \frac{sin B cos B cos c - sin C cos C cos b}
			{sin B cos C cos c - sin C cos B cos b},$\\
	      1. $cos B = \frac{sin C cos C cos a - sin A cos A cos c}
			{sin C cos A cos a - sin A cos C cos c},$\\
	      2. $cos C = \frac{sin A cos A cos b - sin B cos B cos a}
			{sin A cos B cos b - sin B cos A cos a}.$

	\item 0.$cos a = \frac{| sin b|  cos b cos C - | sin c|  cos c cos B}
			{| sin b|  cos c cos C - | sin c|  cos b cos B},$\\
	      1.$cos b = \frac{| sin c|  cos c cos A - | sin a|  cos a cos C}
			{| sin c|  cos a cos A - | sin a|  cos c cos C},$\\
	      2.$cos c = \frac{| sin a|  cos a cos B - | sin b|  cos b cos A}
			{| sin a|  cos b cos B - | sin b|  cos a cos A}.$
	\enume

	Proof\footnote{Echo Lake 22.7.84}:\\
	Let the coordinates of the points $A,$ $B,$ $C$ be
	$(A_0,A_1,A_2,1),$ $(B_0,B_1,B_2,1),$ $(C_0,C_1,C_2,1).$\\
	Those of $DA,$ $DB$ and $DC$ are $(A_0,A_1,A_2,0),$ $(B_0,B_1,B_2,0),$
	$(C_0,C_1,C_2,0).$\\
	if $A \cdot B := A_0 B_0 + A_1 B_1 + A_2 B_2$ and
	$A \cdot A := A_0 A_0 + A_1 A_1 + A_2 A_2,$
	by definition (see $\ldots$ )\\
	\hth$	cos a = B \cdot C$\\
	because $B \cdot B = C \cdot C = 1.$\\
	The plane $A \times B \times O$ is
	\hth$\{A_1 B_2 - A_2 B_1, A_2 B_0 - A_0 B_2, A_0 B_1 - A_1 B_0, 0\}$\\
	the tangent to the sphere at $A$ is\\
	\hth$	\{A_0,	    A_1,	     A_2	  ,-1\}$\\
	and the ideal plane is\\
	\hth$	\{0,	     0,	      0,	  , 1\}$\\
	therefore the direction of $A \times B$ is\\
	\hth$	DAB = A \cdot A B_0 - A \cdot B A_0,$
		$A \cdot A B_1 - A \cdot B A_1,$
		$A \cdot A B_2 - A \cdot B A_2,0).$\\
	Similarly the direction of $A \times C$ is\\
	\hth$	DAC = A \cdot A C_0 - A \cdot C A_0,$
		$A \cdot A C_1 - A \cdot C A_1,$
		$A \cdot A C_2 - A \cdot C A_2,0).$\\
		$DAB \cdot DAB = 1 - (A \cdot B)^2 = 1 - cos^2 c
			= sin^2 c,$ and
		$DAC \cdot DAC = sin^2 b.$\\
	Therefore\\
	\hth$	cos A = \frac{B \cdot C + A \cdot B A \cdot C 
		- A \cdot B A \cdot C 
		- A \cdot C A \cdot B}{| sin b|  | sin c| }
		= \frac{cos a - cos c cos b}{| sin b|  | sin c| },$\\
	hence 2.0.\\
	\hth$sin^2 A sin^2 b sin ^2 c = ( 1 - cos^2 A) sin^2 b sin^2 c$\\
	\hti{16}$ =sin^2 b sin^2 c - cos^2 a - cos^2 b cos^2 c
		+ 2 cosa cos b cos c$\\
	\hti{16}$= 1 - cos^2 a - cos^2 b - cos^2 c + 2 cos a cos b cos c$\\
	Therefore, if\\
	\hth$	r := sqrt{\frac{1 - cos^2 a - cos^2 b - cos^2 c 
		+ 2 cos a cos b cos c}
			{sin^2 a sin^2 b sin^2 c}}$\\
	then\\
	\hth$	\frac{sin A }{ | sin a|} = \frac{sin B }{ | sin b| }
		= \frac{sin C }{ | sin c|} = r.$\\
	Simple algebraic manipulations give 3 to 6.\\
	If we eliminate $cos B$ and $cos C$ from 1.1 and 1.2,
\\
	\hth$	cos B = -\frac{sin C cos b + cos A sin B cos c}{sin A},$\\
	\hth$	cos C = -\frac{sin B cos c + cos A sin C cos b}{sin A},$\\
	substituting in 1.0. gives\\
	\hth$	cos A = (sin C cos b + cos A sin B cos c)
			(sin B cos C + cos A sin C cos b)/sin^2 A$\\
	\hti{16}$		- sin B sin C cos a$\\
	\hth$	cos A = (sin c cos b + cos A sin b cos c)
			(sin b cos c + cos A sin c cos b)/sin^2 a$\\
	\hti{16}$		- sin b sin c cos a$\\
	\hth$cos B cos C - cos A = \frac{(cos b - cos c cos a)
		(cos c - cos b cos a)
		- sin^2 a (cos a - cos b cos c) }{sin^2 a | sin b| | sin c| }$\\
	\hti{16}$	= \frac{cos a(-cos^2 b - cos^2 c + 2 cos a cos b cos c + 
		1 -cos^2 a}
			{sin^2 a | sin b|  | sin c| }$\\
	\hti{16}$	= cos a sin B sin C.$\\
	Hence 1.0.

	\sssec{Example.}
	For $p = 13,$ with $\delta^2$ = 2,\\
	let $A = (0,0,1,1),$ $B = (1,2,3,1),$ $C = (6,1,4,1).$\\
	$cos a = 7,$ $cos b = 4,$ $cos c = 3,$
	$| sin a|  = 2,$ $| sin b|  = 5 \delta,$ $| sin c|  = 3 \delta.$\\
	$cos A = 2,$ $cos B = 4 \delta,$ $cos C = 2 \delta,$
	$sin A = 6,$ $sin B = 2 \delta,$ $sin C = -4 \delta.$

	\ssec{Trigonometry for the right triangle.}

	\sssec{Theorem.}
	{\em For a triangle with a right angle at $A,$ let $sin A = 1,$
	$cos A = 0,$ then we have the relations:}

	\enumb
	\item[0.1.] $| sin b|  = | sin a|  sin B,$
	\item[  2.]$| sin c|  = | sin a|  sin C,$

	\item[1.0.] $cos B cos C = sin B sin C cos a,$
	\item[  1.] $cos B = sin C cos b,$
	\item[  2.] $cos C = sin B cos c.$

	\item[2.~~] $cos a = cos b cos c,$
	\enume

	Proof:
	1 and 0.2 follow from \ref{sec-tneuctrig}.0.\\
	0 follows from \ref{sec-tneuctrig}.2.0.\\
	1 and 1.2 follow from \ref{sec-tneuctrig}.1 which gives 1.0, using 2.0.

	\ssec{Trigonometry for other triangles .}

	\sssec{Definition.}
	An {\em auto-dual triangle} is a triangle such that\\
	\hth$	A = a,$ $B = b,$ $C = c.$

	\sssec{Theorem.}\label{sec-tneuctrigadual}
	{\em If a triangle is auto dual, then}
	\enumb
	\item  $cos A = \frac{cos B cos C}{1+ sin B sin C},$
	\item  $sin A = - \frac{sin B + sin C}{1 + sin B sin C}.$
	\enume

	Proof:  0 follows from Theorem \ref{sec-tneuctrig}.
	If we substitute $cos A$ using  0  in $sin^2 A + cos^2A = 1,$ we get
	$sin A = +j \frac{sin B + sin C}{1 + sin B sin C},$ $j = +1$ or $-1$.
	replacing $sin A$ and $cos A$ by their expression in 1.1, gives after
	multiplication by $1 + sin B sin C,$\\
	\hth$ 	1 + sin B sin C = cos^2 C - i (sin B sin C + sin^2 C),$\\
	therefore $j = -1.$

	\sssec{Notation.}
	$A = (s,c),$ is an abbreviation for $sin A = s,$ $cos A = c.$

	\sssec{Example.}
	For $p = 13,$ let
	$a = A = (-4,-4 \delta ),$ $b = B = (-6,-2),$ $c = C = (3, 5 \delta),$
	we easily verify \ref{sec-tneuctrigadual}.0 and .1:\\
	$cos A = - 4 \delta = - \frac{2 . 5 \delta}{1+3 . -6)},$ 
	$sin A = -4 = -\frac{3-6}{1+3 . -6}.$


\section{Tri-Geometry}
	\ssec{The primitive case.}\label{sec-Striprim}

	\setcounter{subsubsection}{-1}
	\sssec{Introduction.}
	To a given polynomial $P_3$ of the third degree, we can
	associate a selector.  The first case I will consider is that when
	the polynomial has no integer roots or is primitive.  To a given
	such polynomial corresponds a selector called the fundamental
	selector and a tri-geometry with non-integer isotropic points and
	lines.  To this fundamental selector we can associate others, see
	g25.prn,  The semi-selector gives conics associated to the auto-polars,
	the co-selector and the bi-selector are associated to the point-conics
	\footnote{17.3.86}
	and line-conics through the isotropic points, the bi-selector and the
	co-selector to the point-conics and line-conics tangent to the
	isotropic lines.  Examples indicated that the other selectors do not
	give lines or conics or, in general, cubics.  It is an open question if
	they have any geometrical significance.

	\sssec{Definition.}
	If $s$ is the selector, the {\em selector function} is a function
	from ${\Bbb Z}_{p^2+p+1}$ to ${\Bbb Z}_{p^2+p+1}$ given by
	\enumb
	\item 	$f(s_j-s_i) := s_i,$ $i\neq j,$ $f(0) = -1.$
	\enume

	\sssec{Theorem.}
	{\em The selector for the $c$-lines is the co-selector of the
	lines.  More precisely,}
	\enumb
	\item 	$s^c(i) = 1-s(i).$\\
	{\em The selector function for $c$-lines is given by}
	\item 	$f^c(i) = 1 - f(-i).$
	\enume

	\sssec{Theorem.}
	\enumb
	\item 	$a \times b = (f(b-a) - a)^*.$ 
	\item 	$a^* \times b^* = f(b-a) - a.$
	\item{\em $a$ is on $b^*$  iff  $f(a+b) = 0$ or $f(a+b) = -1.$}
	\item{\em the points on $a^*$ are $s(i) - a,$ $i = 0$ to $p.$}
	\item 	$a c b = (1 - f(b-a) - b)^c.$ 
	\item 	$a^c c b^c = 1 - f(b-a) - b.$
	\item{\em $a$ is on $b^c$  iff  $f(-a-b) = 1.$}
	\item{\em the points on $a^c$ are $1 - a - s(i),$ $i = 0$ to $p.$}
	\enume

	\sssec{Definition.}
	Let $c^* := a \times b$.  Let $a = s(i) - c,$ let $b = s(j) - c,$
	the {\em gap} of $a$ and $b,$ written\\
	\hth$gap(b,a) := j - i \umod{p+1}.$\\
	Let $c := a^* \times b^*.$  Let $a^* = s(i) - c,$
	let $b^* = s(j) - c,$
	the {\em gap} of $a^*$ and $b^*,$ written\\
	\hth$gap(b^*,a^*) := j - i \umod{p+1}.$

	\sssec{Theorem.}
	{\em Let $a_0$ be a point on $b_0,$ let
	let $a_i$ be on $b_i,$ such that\\
	\hth$gap(a_i,a_0)$ + $gap(b_i,b_0)$ = 0,\\
	the points $a_i$ are on a $c$-line through $a_0$ tangent to $b_0.$ }

	\sssec{Table.}
	The selector for some values of $p$ and equivalent ones
	which are not complementary (obtained by reversing the order are

	$p = 3,$  0: 0,1,3,9.  1: 0,1,4,6.

	$p = 5,$  0: 0,1,3,8,12,18.  1: 0,1,3,10,14,26.  2: 0,1,4,6,13,21.\\
	\hth	3: 0,1,4,10,12,17. 4: 0,1,8,11,13,17.\\
	$p = 7,$  0: 0,1,3,13,32,36,43,52.  1: 0,1,4,9,20,22,34,51.\\
	\hth	2: 0,1,4,12,14,30,37,52.  3: 0,1,5,7,17,35,38,49.\\
	\hth	4: 0,1,5,27,34,37,43,45.  5: 0,1,7,19,23,44,47,49.

	$p = 11,$ 0: 0,1,3,12,20,34,38,81,88,94,104,109.\\
	\hth	1: 0,1,3,15,46,71,75,84,94,101,112,128.\\
	\hth	2: 0,1,3,17,21,58,65,73,100,105,111.\\
	\hth	3: 0,1,3,17,29,61,80,86,91,95,113,126.\\
	\hth	4: 0,1,4,12,21,26,45,+68,84,97,99,127.\\
	\hth	5: 0,1,4,16,50,71,73,81,90,95,101,108.\\
	\hth	6: 0,1,4,27,51,57,79,89,100,118,120,125.\\
	\hth	7: 0,1,5,12,15,31,33,39,56,76,85,98.\\
	\hth	8: 0,1,5,21,24,39,49,61,75,92,125,127.\\
	\hth	9: 0,1,5,24,44,71,74,80,105,112,120,122.\\
	\hth       10: 0,1,5,25,28,68,78,87,89,104,120,126.\\
	\hth       11: 0,1,6,18,39,68,79,82,98,102,124,126.\\
	\hth       12: 0,1,8,21,33,36,47,52,70,74,76,124.\\
	\hth       13: 0,1,9,19,24,31,52,56,58,69,72,98.\\
	\hth       14: 0,1,15,18,20,24,31,52,60,85,95,107.\\
	\hth       15: 0,1,15,25,45,52,58,61,63,80,84,92.\\
	\hth       16: 0,1,16,21,24,49,51,58,62,68,80,94.\\
	\hth       17: 0,1,23,37,57,62,75,83,86,90,92,102.

	\sssec{Example.}
	Let $p = 3.$ If we use the selector 0,1,3,9 and use the
	representation on the cube (g25.prn), the complementary selector
	0,1,5,11 gives the $c$-lines which can be classified as follows:
	3 of type $VVSS,$ 2 vertex-points and 1 side-point through each.\\
	\hth More precisely, two of 2 adjacent vertex-points and 1
	 side-point through each, such that no 2 are in the same face,
	 one of 2 opposite vertex-points and 2 adjacent side-points one
	 through each.\\
	3 of type $FSSS,$ 1 face-point, 1 side-point in it and 2 opposite
		side-points in an other face.\\
	3 of type $FVSS,$ 1 face-point, two adjacent side-points in it and a
		vertex point on one of the side-points.\\
	3 of type $FSVV,$ 1 face-point, two adjacent vertex-points in it and a
		side-point through one of the vertex-points.\\
	1 of type $VFFF,$ 1 vertex-point and the 3 face-points.\\
	Clearly the converse is not true.  For instance, only one of the 4
	vertex-points can serve for the last case given.

	\sssec{Example.}
	Let $p = 7,$ $P_3 = I^3 + 2,$\\
	The powers of $I+3$ are:\\
	$\begin{array}{crcrrrrrrrrrrrrrrr}
	\hth&0&\vline&0&0&1,&\hti{2}0&1&3,&\hti{2}1&-1&2,&\hti{2}1&3&2,&
	\hti{2}1&3&3,\\
	\hth&5&\vline&1&2&0,&1&-3&1,&0&1&-1,&1&2&-3,&1&2&2,\\
	\hth&10&\vline&1&3&-2,&1&0&1,&1&-2&-2,&1&-1&-1,&1&-2&1,\\
	\hth&15&\vline&1&2&1,&1&0&3,&1&1&0,&1&-1&3,&1&0&0,\\
	\hth&20&\vline&1&0&-3,&1&-1&1,&1&-1&-3,&1&-3&-2,&0&1&2,\\
	\hth&25&\vline&1&-2&-1,&1&0&2,&1&3&-1,&1&-1&-2,&1&1&3,\\
	\hth&30&\vline&1&-2&0,&1&1&-2,&1&2&-2,&1&-2&-3,&1&-2&3,\\
	\hth&35&\vline&1&-3&0,&0&1&1,&1&-3&3,&0&1&0,&1&3&0,\\
	\hth&40&\vline&1&-2&2,&1&3&-3,&1&1&-3,&1&0&-1,&1&2&3,\\
	\hth&45&\vline&1&-1&0,&1&2&-1,&1&1&-1,&1&-3&-3,&0&1&-2,\\
	\hth&50&\vline&1&1&1,&1&1&2,&1&3&1,&1&-3&-1,&0&1&-3,\\
	\hth&55&\vline&1&0&-2,&1&-3&2,
	\end{array}$

	\sssec{Example,}
	$p = 7,$ $P_3 = I^3 + 2,$\\
	 ~0$^*$:~~0 ~1 ~7 24 36 38 49 54\\
	 ~1$^*$:~~0 ~6 23 35 37 48 53 56\\
	 ~7$^*$:~~0 17 29 31 42 47 50 51\\
	 24$^*$:~~0 12 14 25 30 33 34 40\\
	 36$^*$:~~0 ~2 13 18 21 22 28 45\\
	 38$^*$:~~0 11 16 19 20 26 43 55\\
	 49$^*$:~~0 ~5 ~8 ~9 15 32 44 46\\
	 54$^*$:~~0 ~3 ~4 10 27 39 41 52\\
	The points 0,3,8,19,21,33,50,56 are on a $c$-line through 0.\\
	The points 0,5,16,18,30,47,53,54 are on a $c$-line through 0.

	The part proving that co-, bi- and semi-selectors are conics was
	proven before this date	\footnote{31.3.86}.
	The equation of the conics through 2
	coordinate points was also obtained earlier.  It remains to prove that
	the 2 are identical.

	\sssec{Lemma.}
	\enumb
	\setcounter{enumi}{-1}
	\item{\em If $i$ is an element of the co-selector, the tangent is
	$(1-2i)^*.$ }
	\item{\em If $i$ is an element of the bi-selector, the tangent is
	$(a-\frac{i}{2})^*,$ for some $a.$}
	\enume

	Proof:\\
	For 0, $(1-2i)^*$ is on $i$ because $f(1-i) = 0$ if $i$ is an element
	of the co-selector.  It remains to prove that it is the only point on
	$(1-2i)^*.$ For 2, by duality?

	\sssec{Theorem.}\label{sec-tgeomprimit}
	Let $S$ be a selector\footnote{17.3.86}.
	\enumb
	\item{\em The points associated to the co-selector are on a conic which
	    passes through the isotropic points.}
	\item{\em The points associated to the bi-selector are on a conic which
	    is tangent to the isotropic lines.}
	\item{\em The points associated to the semi-selector are on a conic
	    for which the isotropic triangle is a polar triangle.}
	\item{\em The conics of the same family are such that 2 distinct points
	    determine a conic and 2 distinct conics have exactly one point in
	    common.}
	\enume

	Proof:  Le $P_3 = I^3 + bI - c$\footnote{2.4.86}.\\
	For 0.  Consider the selector associated to the line through
	$0 = G^0$ and
	$1 = G = I + g.$  Let $G^i = I + h.$  The corresponding point on the
	co-selector is $G^{-i+1}.$  We obtain\\
	\hth$	G^{-i+1} = (g-h) I^2 - h(g-h) I + (g(b+h^2)+c)$\\
	It is easy to check that that point is on the conic\\
	\hth$	(bg+c) X_0^2 + g X_1^2 - g X_0X_1 - X_1X_2$\\
	and that the isotropic points are on this conic.

	Part 1, follows by duality in view of Lemma 3.2.10.1.

	For 2, because the line $i^*$ is on the point $i,$ the correspondance
	which associates $i^*$ to $i$ is a polarity and the points on their
	polars is a conic, the auto-polar conic.  These points are such that
	$f(2i) = 0,$ where $f$ is the selector function and therefore the
	solutions
	$i$ are points corresponding to the semi-selector\footnote{31.3.86}.
	In view of g142.prn, the symmetric matrix $M_2$ which represents the
	auto-conic satisfies for some values of $u,$ $v,$ $w$\\
	\hth$\matb{(}{ccc}u&gv&(g^2-b)w\\0&v&2gw\\0&0&w\mate{)}
	 = \matb{(}{ccc}a_0&b_2&b_1\\b_2&a_1&b_0\\b_1&b_0&a_2\mate{)}\:
	\matb{(}{ccc}0&0&1\\0&1&2g\\1&g&g^2\mate{)}.$\\
	The inverse of the last matrix is\\
	\hth$\matb{(}{ccc}g^2&-g&1\\-2g&1&0\\1&0&0\mate{)}.$\\
	Multiplying the first matrix by this last matrix gives, because
	of the symmetry,\\
	$u = 1,$ $v = 1,$ $w = 1$ and with $b = s_{11},$\\
	\hth$M_2 = \matb{(}{ccc}-s_{11}&0&1\\0&1&0\\1&0&0\mate{)}.$\\
	This matrix clearly associates to the pole $(1,-(\pi_1+\rho_2),
	\rho_1\rho_2)$ the polar, $[\rho_0^2, \rho_0, 1],$ because $s_1 = 0.$

	\sssec{Answer to}
	\vspace{-18pt}\hspace{90pt}{\bf \ref{sec-tgeomprimit}.0.}\\
	I fill in here some of the details:\\
	If $(I + g) * (I + h)^{-1} = uI^2 + vI + w,$\\
	then $((v+uh)I^2 +(w+vh-ub) + (wh+uc) = k (I + g),$\\
	therefore\\
	\hth$	v = -uh,$ $wh + uc = g(w-uh^2-ub)$ or\\
	\hth$	u = g-h,$ $v = -h(g-h),$ $w = g(b+h^2)+c.$\\
	The conic through the isotropic points and through the points 0 and 1
	is of the form $a_0X_0^2 + gX_1^2 - X_1X_2 + b_1X_2X_0 + b_2X_0X_1.$ \\
	To insure that it passes through the isotropic points gives 3 linear
	equations for $a_0,$ $b_0$ and $b_2.$  It is easiest to check a
	posteriori that\\
	\hth$	(bg+c)X_0^2 + g X_1^2 - g X_2X_0 - X_1X_2$\\
	passes through the isotropic
	points, for instance through $(1, \rho_0,$ $r_1\rho_2):$\\
	\hth$	g(\rho_1\rho_2+\rho_2\rho_0+\rho_0\rho_1 + \rho_0^2 
	- \rho_1\rho_2) + (c - \rho_0\rho_1\rho_2) = 0.$

	The point $(u,v,w)$ is on the conic because\\
	\hth$	(bg+c)(g-h)^2 + gh^2(g-h)^2 - g(g-h)w + h(g-h)w$\\
	\hth$	= (g-h)^2(bg+c+gh^2-g(b+h^2)-c) = 0.$

	\sssec{Definition.}
	The mapping which associates to a point $P$ corresponding to $G^k,$
	the point $Q$ corresponding to $G^{-k}$ is called the {\em inversion
	mapping}.

	\sssec{Theorem.}
	{\em If $P_3 = I^3 + bI - c,$}
	\enumb
	\item{\em The inversion mapping $T$ associates to $(x,y,z),$ $(X,Y,Z)$
	with}\\
	\hth$	X = bx^2 + y^2 - xz,$\\
	\hth$	Y = cx^2 - yz,$\\
	\hth$	Z = (bx-z)^2 + by^2 - cxy.$\\
	\hth$    T \circ T (x,y,z)$\\
	\hth$    = (c^2x^3 + bcx^2y + b^2x^2z - 3cxyz - 2bxz^2 + cy^3 
		+ by^2z + z^3).(x,y,z) \:?$\\
	\enume

	\sssec{Example.}
	$p = 7,$ $P_3 = I^3 + 2,$
	\begin{verbatim}
	selector: 0  1  7  24  36  38  49  54
	selector function:
	 0  1  2  3  4  5  6  7  8  9 10 11 12 13 14 15 16 17 18 19 20 21 22 23
	-1  0 36 54 54 49  1  0 49 49 54 38 24 36 24 49 38  7 36 38 38 36 36  1
	24 25 26 27 28 29 30 31 32 33 34 35 36 37 38 39 40 41 42 43 44 45 46 47
	 0 24 38 54 36  7 24  7 49 24 24  1  0  1  0 54 24 54  7 38 49 36 49  7
	48 49 50 51 52 53 54 55 56
	 1  0  7  7 54  1  0 38  1
	c-selector: 0  1  4  9  20  22  34  51
	c-selector fuction:
	 0  1  2  3  4  5  6  7  8  9 10 11 12 13 14 15 16 17 18 19 20 21 22 23
	-1  0 20  1  0  4 51 51  1  0 51  9 22  9 20 51  4 34  4  1  0  1  0 34
	24 25 26 27 28 29 30 31 32 33 34 35 36 37 38 39 40 41 42 43 44 45 46 47
	34  9 51 34 51 22  4 20 34  1  0 22 22 20 20 22 51 20  9 34 22 34 20  4
	48 49 50 51 52 53 54 55 56
	 9  9  1  0  9  4  4 22  1
	\end{verbatim}

	$12 \times 22 = 42^*,$ $12 c 22 = 39^c,$ 
	$52^*$ is tangent at 12 to $39^c,$ $32^*$ is tangent at 22 to $39^c,$ 
	$49^c$	is tangent at 12 to $42^*,$ $29^c$ is tangent at 22 to $42^*,$
	$32^* \times 52^* = 6,$ $29^c c 49^c = 28.$
	$ 6 \times 28 = 30^*,$ $ 6 c 28 = 51^c,$ 
	$52^*$ is tangent at  6 to $51^c,$  $8^*$ is tangent at 28 to $51^c,$ 
	$16^c$	is tangent at  6 to $30^*,$ $29^c$ is tangent at 28 to $30^*,$
	$ 8^* \times 52^* = 14,$ $29^c c 16^c = 50.$  There appears to be no
	connection.

	the $c$-lines (conics through the isotropic points) are\\
	\hth$- 2m X_0^2 + l X_1^2 + 3k X_2^2 - m X_1X_2 - l X_2X_0 - k X_0X_1
	= 0.$

	the $c$-line with $l = m = 0$ is\\
	\hth$	0,1,0; 1,0,0; 1,3,1; 1,3,6; 1,5,2; 1,5,5; 1,6,3; 1,6,4$
	or 12, 18, 19, 22, 27, 38, 40, 52,\\
	the 57 others are all obtained by adding a constant, for instance,
	if we add 38 we get
	50, 56,  0,  3,  8, 19, 21, 33 or
	1,1,1; 0,0,1; 1,4,5; 1,3,2; 1,2,4; 1,0,0; 1,6,1; 1,5,4\\
	which corresponds to $k = m = 0,$\\
	if we add 19 we get 31, 37, 38, 41,46, 0, 2, 14, $(k = l = 0)$\\
	if we add 33 we get 45, 51, 52, 55,  3, 14, 16, 28, $(k = l = m = 1)$\\
	Is there any significance to the fact that 19, 38 are $\frac{57}{3},$
	and $\frac{2.57}{3}$?

	the $c$-points (conics tangent to isotropic lines) are\\
	\hth$3k x_0^2 + l x_1^2 - 2k x_2^2 - k x_1x_2 - l x_2x_0 - k x_0x_1
	= 0.$

	\sssec{Example.}
	$p = 7,$ $P_3 = I^3 + 2,$ selector, 0,1,7,24,36,38,49,54,\\
	$3x_0^2-x_1x_2:$ 18,19,22,27,38,40,52,12, (0,1,4,9,20,22,34,51)\\
	\hti{12}$	    40,38,32,22, 0,53,29,52, (0,1,5,27,34,37,43,45)$\\
	$x_1^2-x_2x_0:$	  56, 0, 3, 8,19,21,33,50,
	$-2x_2^2-x_0x_1:$	37,38,41,46, 0, 2,14,31,

	If we replace $I$ by $I+1$ we obtain $P_3 = I^3+3I^2+3I+3,$ \\
	this gives the same selectors,\\
	$P_3 = I^3+1,$ and $P_3 = I^3+I-1$ give the selectors\\
	(0,1,6,15,22,26,45,55),
	(0,1,3,13,32,36,43,52) and (0,1,5,7,17,35,38,49).

	\sssec{Example.}
	For $p = 3$ and $P_3 = I^3+2,$ the auto-polar conics through two
	of the points (0,0,1), (0,1,0) and (1,0,0) are\\
	$X_0^2 - X_1X_2 = 0$ or $x_0^2 + 3x_1x^2 = 0,$\\
	$X_1^2 - X_2X_0 = 0$ or $x_1^2 + 2x_2x^0 = 0,$\\
	$X_2^2 - X_0X_1 = 0$ or $x_2^2 - x_0x^1 = 0.$\\
	All 3 do not have 3 points or 3 lines in common.

	\sssec{Comment.}
	If 2 conics are in the same family and we known the
	tangents corresponding to the points of one we can obtain those
	of the other.  The sum of the corresponding representation of points
	and lines is a constant.

	\sssec{Program.}
	Examples can be studied using 130$\setminus$ TWODIM.BAS.

	\ssec{The case of 1 root.  Inverse geometry.}
	\label{sec-Stri1root}

	\setcounter{subsubsection}{-1}
	\sssec{Introduction.}
	Let $P_3 = (I^2 + a I + b) (I + c),$ $a^2 - 4 b$ N $p.$\\
	There is one isotropic point $(1,a,b)$ and one isotropic line
	$[c^2,-c,1].$\\
	The isotropic point is not on the line otherwize, $-c$ whould be a root
	of $I^2 + a I + b.$  The $p+1$ ideal points are $(0,1,c)$ and
	$(1,x,c(x-c)).$\\
	The $p+1$ ideal lines are $[0,b,-a]$ and $[1,x,-\frac{1+ax}{b}].$\\
	In the case of the complex field, if $P_3 = (I^2+1) I,$ the
	$c$-lines are the circles through the origin, it is therefore natural
	to call this geometry {\em inverse geometry}.

	\sssec{Definition.}
	The {\em pseudo-bi-selector} is the set $\{2 s_i\},$ \\
	The {\em pseudo-semi-selector} is the set $\{\frac{1}{2} s_i\},$

	\sssec{Example.}
	$p = 5.$ With $P_3 = I^3-I^2-2I-3,$ a generator is $I+2,$
	its powers are\\
	$\begin{array}{cccccccccc}
	0,0,1&0,1,2&1,4,4&1,2,3&0,1,1&1,3,2&1,0,2&1,3,4&1,2,1&0,1,0\\
	1,2,0&0,1,3&1,0,1&1,1,0&1,1,2&1,4,3&1,4,2&1,1,1&1,0,0&1,4,1\\
	1,3,0&1,3,3&1,1,4&1,2,4
	\end{array}$

	The lines are\\
	0 [1,0,0]: \{0,1,4,9,11\}  and (0,1,4)\\
	1 [1,2,4]: \{1,2,5,10,12\} and (1,1,3)
	 $\ldots$ ..\\
	The selector function is\\
	$\begin{array}{ccrrrrrrrrrrrrrrrrrrrrr}
	i&\vline&	1&2&3&4&5&7&8&9&10&11&13&14&15&16&17&19&20&21&22&23\\
	f(i)&\vline&	0&9&1&0&4&4&1&0&1&0&11&11&9&9&11&9&4&4&11&1
	\end{array}$

	The isotropic line [1,1,1], the isotropic point is (1,0,3).
	The ideal lines are $[j] = \{j,j+6,j+12,j+18\},$ for $j = 0$ to 5.

	I will now examine the case when ${\Bbb Z}_p$ is replace by an
	infinite field ${\Bbb R},$ for instance.

	$p = 7,$ $P_3 = I^3 + I,$ $G = I + 3,$\\
	$\begin{array}{cccccccccc}
	0,0,1&0,1,0&0,1,1&0,1,2&0,1,3&0,1,4&0,1,5&0,1,6&1,0,0&1,0,1\\
	 0, i&0^I,0^i&46,12^*&29,28^*&1,20^*&7,44^*& 11,4^*&
	34,36^*&4^I,0^*&I,16^*\\
	1,0,2&1,0,3&1,0,4&1,0,5&1,0,6&1,1,0&1,1,1&1,1,2&1,1,3&1,1,4\\
	16, 8^*&40,32^*&24,40^*&8,24^*& 32, 4^*& 6^I,46^*&12,14^*&
	 14, 6^*& 21,30^*&25,38^*\\
	1,1,5&1,1,6&1,2,0&1,2,1&1,2,2&1,2,3&1,2,4&1,2,5&1,2,6&1,3,0\\
	15,22^*&43,2^i&5^I, 7^*&44,23^*&23,15^*&38,39^*&
	35,47^*&37,31^*&18,3^i&1^I,11^*\\
	1,3,1&1,3,2&1,3,3&1,3,4&1,3,5&1,3,6&1,4,0&1,4,1&1,4,2&1,4,3\\
	 4,27^*&27,19^*&41,43^*& 22, 3^*& 42,35^*&39,7^i&7^I,29^*&
	28,45^*&45,37^*&47,13^*\\
	1,4,4&1,4,5&1,4,6&1,5,0&1,5,1&1,5,2&1,5,3&1,5,4&1,5,5&1,5,6\\
	10,21^*&6, 5^*&33,1^i&3^I, 1^*& 20,17^*&17, 9^*&
	 26,33^*&5,41^*&19,25^*&30,5^i \\
	1,6,0&1,6,1&1,6,2&1,6,3&1,6,4&1,6,5&1,6,6\\
	2^I,34^*&36, 2^*&2,42^*&3,18^*& 31,26^*&9,10^*&	13,6^i
	\end{array}$

	The {\em real isotropic point} is denoted by $I,$ the {\em real
	isotropic line} by $i.$
	$1 \times 9 = (-1) -1 = 7^i,$ $(0,1,3) \times (1,6,5) = [1,3,6].$
	$9^* \times 17^* = (-1) - 1 = 7^I,$ $[1,5,2] \times [1,5,1] = (1,4,0)$
	\footnote{7.4.86}.
	Observe that $k^I$ corresponds to $I G'^I,$ with
	$G' = I+3 \pmod{I_2+1}.$

	The selector is 0,1,7,11,29,34,46.\\
	The co-selector is 0,1,3,15,20,38,42.\\
	The {\em pseudo-bi-selector} is	$\{0, 2 ,14 ,22 ,10 ,20,44 ,4^I\}$.\\
	The corresponding tangents to the bi-conic are
				$\{i,47^*,41^*,37^*,19^*,2^*,14^*,0^*\}$\\
	which is a member of the dual of the co-selector family
	\footnote{9.4.86}.\\
	The pseudo-semi-selector is $\{0 ,17 ,23 ,24 ,41 ,47 ,2^I,6^I\}$.\\
	The corresponding tangents to the semi-conic are
				$\{0^*,17^*,23^*,24^*,41^*,47^*,2^i,6^i\}$.\\	
	The points $2^I,6^I,$ are obtained from the ideal tangents $2^i,6^i.$\\
	The same can be checked if we add 1, $\ldots$ , to the values above,
	we get in this way 24 hyperbolas, e.g.\\
	\hth$\{ 1 ,18 ,24 ,25 ,42 , 0 ,3^I,7^I\}$\\
	\hth$\{47^*,16^*,22^*,23^*,40^*,46^*,1^i,5^i\},$\\
	and 24 ellipses, e.g.\\
	\hth$\{ 1 , 4 , 6 ,15 ,25 ,28 ,30 ,39 \},$\\
	\hth$\{ 0^*, 3^*, 5^*,14^*,24^*,27^*,29^*,38^*\}.$

	Hence, do we also have therefore the Theorem that a selector has
	$\frac{p-1}{2}$ even values and $\frac{p+1}{2}$ odd values?

	\sssec{Comment.}
	In the case of the field ${\Bbb R},$ every polynomial of degree 3
	has necessarily one root.  There is no restriction in assuming that it
	is $P_3 := I^3 + I.$  In this case the isotropic points are
	(1,0,1), and
	the Euclidean isotropic points (1,i,0), (1,-i,0).

	\sssec{Theorem.}
	{\em If the field is ${\Bbb R}$ and $P_3 := I^3 + I,$ the transformation
	associated}
	\enumb
	\item{\em to $k = -1,$ transforms the lines into circles through the
		point (1,0,1).}
	\item{\em to $k = 2,$ transforms the lines into parabolas with focus
	(1,0,1).}
	\item{\em to $k = \frac{1}{2},$ transforms the lines into equilateral
	hyperbolas with center (1,0,1).}
	\enume

	Proof:\\
	For 0, the conics which pass trough $(1,i,0)$ and $(1,-i,0)$ are
	circles. (1,0,1) is the third isotropic point.\\
	For 1, the conics are tangent to the isotropic line through $(1,i,0)$
	and $(1,-i,0)$ which is the ideal line.  Because the focus of a parabola
	is at the intersection of the tangent through the Euclidean isotropic
	point, we have 1.\\
	For 2, because (1,0,1) is the pole of the opposite isotropic line which
	is the ideal line, (1,0,1) is the center of the conic.  Because the
	points on the conic and the ideal line form a harmonic quatern with
	the pole $(1,i,0)$ and the intersection $(1,-i,0)$ with its polar, the
	corresponding directions, which are those of the asymptotes to the
	hyperbola and therefore perpendicular.

	More explicitely:

	\sssec{Theorem.}
	{\em The transformation associated to the case $k = -1$
	associates to to the point $(x,y,1)$ or $(x,y),$ the point $(X,Y,1)$ or
	$(X,Y),$ with}
	\enumb
	\item{\em$(xI^2	+ yI + 1) (XI^2 + YI + 1) = 1 \pmod{P_3}.$
	this gives}
	\item $X - 1 = \frac{x-1}{(x-1)^2 + y^2}, Y = -\frac{y}{(x-1)^2 + y^2}.$
	\item{\em The point $Q = (X,Y)$ and the point $P = (x,y)$ are on a line
	through (1,0), the product of the distances to that point is 1,
	and the points $P$ and $Q$ are separated by (1,0).}
	\enume

	Proof:\\
	0, gives with $I^3$ replaced by $-I$ and $I^4$ replaced by $-I^2,$ \\
	\hth$	-xX + x + X + yY = 0,$\\
	\hth$	-xY - yX + y + Y = 0.$\\
	solving for $X$ and $Y$ gives easily 1.\\
	Moreover,\\
	\hth$	\frac{X-1}{Y} = - \frac{x-1}{y},$

	\sssec{Theorem.}
	{\em The transformation associated to the case $k = 2$
	associates to to the point $(x,y,1)$ or $(x,y),$ the point $(X,Y,1)$ or
	$(X,Y),$ with}
	\enumb
	\item{\em$(XI^2 + YI + 1) = (xI^2 + yI + 1)^2 \pmod{P_3}.$ \\
	this gives}
	\item 	$X - 1 = y^2 - (x-1)^2,$ \\
		$Y = -2y(x-1).$
	\item{\em The line associated to $a(x-1) + by + c$ is the parabola}\\
		$(-2ab^2(X-1)+b(a^2-b^2)Y-2ac^2)^2
		= 4(a^2+b^2)c^2(c^2-b^2(X-1)-abY).$
	\enume

	Proof:
	0, gives with $I^3$ replaced by $-I$ and $I^4$ replaced by $-I^2,$
	we obtain at once 1.
	For 2, to simplify let us write the line as $y = mx'+d,$ with
	$x' := x-1,$
	$m = -\frac{a}{b}$ and $d = -\frac{c}{b}.$  Expressing $y$ in terms of
	$x'$ gives, with $X' := X-1,$\\
	\hth$	2mX' + (m^2-1)Y = 2(m^2+1)dx' + 2md^2$ and\\
	\hth$	X' + mY = -(m^2+1)x'^2 + d^2$\\
	eliminating $x'$ gives 2.

	\sssec{Theorem.}
	{\em The transformation associated to the case $k = \frac{1}{2}$
	associates to the point $(x,y,1)$ or $(x,y),$ the point $(X,Y,1)$ or
	$(X,Y),$ with}
	\enumb
	\item 	$(XI^2	+ YI + 1)^2 = (xI^2 + yI + 1)  \pmod{P_3}.$

	{\em this gives}
	\item 	$x - 1 = Y^2 - (X-1)^2,$ $y = -2Y(X-1).$
	\item{\em The line associated to $a(x-1) + by + c$ is the hyperbola}\\
	\hth$	a(Y^2-(X-1)^2) - 2bY(X-1) + c = 0.$
	\enume

	Proof:
	0, gives with $I^3$ replaced by $-I$ and $I^4$ replaced by $-I^2,$
	we obtain at once 1.

	\sssec{Theorem.}
	\enumb
	\item{\em The lines joining the points associated to the selector
	and their inverse are tangent to a conic.}\\
	\hth$	2y^2 + bz^2 + xz = 0$ or $2bX^2 - Y^2 - 8XZ = 0.$
	\item{\em The lines joining the points associated to the co-selector
	and their inverse are tangent to a conic.}
	\enume

	Proof:  The points of the selector are $I-h,$ their inverse is
	$I^2+hI+h^2+b.$	 The line through these points is $[-2h^2-b,h,1].$

	\sssec{Problem.}
	Complete a set of axioms of inverse geometry using an
	appropriate form of the axiom of Pappus:
	\enumb
	\item   Given 2 distinct points, there exist one and only one line
	    incident to, or passing through, the 2 points, or the {\em points
		 are parallel}.
	\item   Given 2 distinct lines, there exists one and only one point
	    incident to, or on, the 2 lines, or the {\em lines are parallel}.
	\item   There exists at least one line $l$ and two distinct points $P$
		 and $Q$ not incident to $l.$
	\item   On the line $l$ there are exactly $p$ points, $p$ an {\em odd
		prime}.
	\item   Given a line $l$ and a point $P$ not on the line, there exists
	 one and only one line parallel to $l$ through $P.$
	\item   Given a point $P$ and a line $l$ not through the point, there
	exists one and only one point parallel to $P$ on $l.$
	\enume

\setcounter{section}{3}
\setcounter{subsection}{3}
	\ssec{The case of a double root and a single root.}
	\label{sec-Stri2roots}
	\vspace{-18pt}\hspace{450pt}\footnote{16.5.85}\\[8pt]
	\footnote{28.2.86}
	\setcounter{subsubsection}{-1}
	\sssec{Introduction.}
	There is no ambiguity to call this also the case of 2 {\em roots}.

	\sssec{Definition.}
	The {\em selector function} is a function from the set
	${\Bbb Z}_{p(p-1)} - \{0 \pmod{p}\} - \{0 \pmod{p-1}\}$ into
	${\Bbb Z}_{p(p-1)},$
	with\\
	\hth$	f(i-j) = i-j \pmod{p(p-1)},$ for all $i$ and $j$ on $[1,0,0].$

	\sssec{Example.}
	$p = 5$
	The cyclic group is\\
	0,0,1  0,1,1  1,2,1  1,2,4  1,4,1  1,0,1  1,3,3  0,1,3  1,4,3  1,2,3\\
	1,0,2  1,1,1  1,4,2  1,1,2  1,1,4  1,0,3  1,4,4  1,3,2  0,1,2  1,3,2.\\
	The lines are\\
	~0 [1,0,0]: \{0,1,7,18  and (0,1,0), (0,1,4)\}\\
	19 [1,4,1]: \{1,2,8,19  and (1,1,0), (1,0,4)\}\\
	18 [1,2,0]: \{2,3,9,0   and (1,2,0), (1,2,2)\}\\
	17 [1,3,2]: \{3,4,10,1  and (1,3,0), (1,1,3)\}\\
	\ldots.

	The selector is\\
	$\begin{array}{ccrrrrrrrrrrrr}
	i&\vline&	1&2&3&6&7&9&11&13&14&17&18&19\\
	f(i)&\vline&	0&18&18&1&0&18&7&7&7&1&0&1\\
	\end{array}$

	If $f(j-i)$ does not exist, then if $j-i \equiv \pmod{4},$ the points
	are on a line through (1,4,0).\\
	If $f(j-i)$ does not exist, then if $j-i \equiv \pmod{5},$ the points
	are on a line through (1,0,0).\\
	Otherwize, the line if $f(j-i)-i \pmod{20}.$

	There is no restriction in assuming that $P_3 := I^3 - I^2.$ 

	\sssec{Definition.}
	The {\em bi-isotropic point} is $I_0 := (1,-1,0),$\\
	the {\em isotropic} {\em point is} $I_1 := (1,0,0).$\\
	The {\em bi-isotropic line} is $i_0 := [0,0,1],$\\
	the {\em isotropic line} is $i_1 := [1,1,1].$

	\sssec{Theorem.}
	\enumb
	\item{\em The points associated to the co-selector are on a conic, the
	{\rm co-conic}, which passes through the isotropic point $I_1$ and is
	tangent to the isotropic line $i_1$ at the co-isotropic point $I_0.$ }
	\item{\em The points associated to the bi-selector are on a conic, the
	{\rm bi-conic}, which is tangent to the co-isotropic line $i_1$ and is
	tangent to the isotropic line $i_0$ at the isotropic point $I_1.$}
	\item{\em The points associated to the semi-selector are on a conic, the
	{\rm semi-conic}, which is tangent to the isotropic line $i_0$ at the
	co-isotropic point $I_0$ and is such that the polar of the isotropic
	point $I_1$ is the co-isotropic line $i_1.$ }
	\enume

	Proof:\\
	For 0,  The conic of 3.1.8.0. reduces to\\
	\hth$	(k-l) Y^2 + m Z^2 + (m-l) YZ +(m-k) ZX + (k-l) XY = 0,$\\
	which passes through $I_1$ and for which [1,1,1] is the polar of
	$(1,-1,0).$\\
	For 1,\ldots.\\
	For 2,\ldots.

	\ssec{The case of a triple root.  Solar geometry.}
	\label{sec-Stritriple}
	\vspace{-18pt}\hspace{290pt}\footnote{18.2.86}\\[8pt]
	\setcounter{subsubsection}{-1}
	\sssec{Introduction.}
	In this case the  $\ldots$ 
	There is one special point and a line belonging to each other.
	The special point and the special lines are called respectively
	the {\em isotropic point} and the {\em isotropic line}.
	The other points on the isotropic line are called {\em ideal points}.
	The other lines through the isotropic point are called {\em ideal line}.
	The other points and lines are called {\em ordinary}.

	\sssec{Theorem.}
	\enumb
	\item{\em There are $p$ ideal points, $p$ ideal lines.}
	\item{\em There are $p^2$ ordinary points, $p^2$ ordinary lines.}
	\item{\em The isotropic line belongs to 1 isotropic point and $p$ ideal
		points.}
	\item{\em The isotropic point belongs to 1 isotropic line and $p$ ideal
		lines.}
	\item{\em The ideal lines belong to \ldots.}
	\item{\em The ordinary lines belong to \ldots.}
	\enume

	\sssec{Comment.}
	In the parabolic-Euclidean or sun-geometry, among all points
	and lines of projective geometry, one point and a line through it are
	preferred, in this case this is also true, but if we represent this
	geometry in the Cartesian plane and choose the isotropic line as the
	line at infinity and the isotropic point as the direction of the $x$
	axis, the $c$-lines are parabolas which have  $\ldots$ .
	It is therefore natural, by analogy to choose the names
	{\em solar geometry}
	and {\em bi-solar geometry}, for the geometry in question.

	\sssec{Lemma.}
	\enumb
	\item 	$X_1^2	- 2X_2 (X_0 + a X_1) = 0$

	{\em and}
	\item 	$Y_1^2	- 2Y_2 (Y_0 + a Y_1) = 0$

	{\em implies}
	\item 	$(X_1Y_2+X_2Y_1)^2 - 2X_2Y_2 (X_0Y_2+X_2Y_0+X_1Y_1
		+ a(X_1Y_2+X_2Y_1)) = 0.$
	\enume

	Proof:  The first member of 2. is the sum of the first member of 0 and
	1 multiplied respectively by $Y_2^2$ and $X_2^2.$ 

	\sssec{Theorem.}
	{\em Let $a = 0,$ generators of $T$ are $I + 1$ and $I + 2.$
	The cyclic group of order p generated by I+b corresponds to points on
	the conic\\
	\hth$	X_1^2 - 2X_2 (X_0 + a X_1) = 0$\\
	where $a := \frac{1}{2b}.$
	The cyclic group of order $p$ generated by $I^2+I+\frac{1}{2}$
	corresponds to points on the conic\\
	\hth$	X_1^2 - 2X_2 X_0 = 0.$}

	Proof:\\
	The point $(X_0,X_1,X_2)$ corresponds to the polynomial
	$X_0I^2 + X_1I + X_2.$\\ 
	The product of $(X_0I^2 + X_1I + X_2)$ and $(Y_0I^2 + Y_1I + Y_2)$ is
	$(X_0Y_2+X_2Y_0+X_1Y_1,X_1Y_2+X_2Y_1,X_2Y_2).$	 The Theorem folllows
	at once from Lemma  $\ldots$ .

	\sssec{Definition.}
	The {\em selector function} is a function from ${\Bbb Z}_p \bigx
	{\Bbb Z}_p$ to ${\Bbb Z}_p \bigx {\Bbb Z}_p$  $\ldots$ .

	\sssec{Definition.}
	The {\em ideal lines} can be represented by $[i],$ $i \in {\Bbb Z}_p.$\\
	The points on $[i]$ are $(j,i+j \umod{p}).$

	\sssec{Theorem.}
	\enumb
	\item   $(x,y) \times (x',y') = (f(x'-x,y'-y) - (x,y)
	\umod{{\Bbb Z}_p \bigx {\Bbb Z}_p)}.$ 
	\item   $[x,y] \times [x',y'] = [f(x'-x,y'-y) - (x,y)
	\umod{{\Bbb Z}_p \bigx	{\Bbb Z}_p]}.$ 
	\item   $(x,y) \cdot [x',y']$  iff  $f(x'-x,y'-y) = (0,0).$
	\enume

	\sssec{Example.}
	For $p = 3,$ the selector function is\\
	$\begin{array}{cccccccc}
	x&\vline&	0,1&0,2&1,0&1,2&2,0&2,1\\
	f(x)&\vline&	0,0&0,1&0,0&0,1&1,0&1,0
	\end{array}$

	With the exponents in the order exponent of $b$ then exponent of $a,$\\
	\hti{11} points\hti{8} line\\
	$\begin{array}{cccccc}
	\cline{1-6}
	\hth&0,0&0,1&1,0&\vline&0,0\\
	&0,1&0,2&1,1&\vline&0,2\\
	&0,2&0,0&1,2&\vline&0,1\\
	&1,0&1,1&2,0&\vline&2,0\\
	&1,1&1,2&2,1&\vline&2,2\\
	&1,2&1,0&2,2&\vline&2,1\\
	&2,0&2,1&0,0&\vline&1,0\\
	&2,1&2,2&0,1&\vline&1,2\\
	&2,2&2,0&0,2&\vline&1,1
	\end{array}$\\
	The points on the 3 ideal lines are\\
	\hth[0] = \{0,0  1,1  2,2\}\\
	\hth[1] = \{0,1  1,2  2,0\}\\
	\hth[2] = \{0,2  1,0  2,1\}

	\sssec{Example.}
	If $p = 7$ and $P^3 = I^3$ then the group $T$ is\\
	$\begin{array}{ccccccc}
	0,0,1 & 0,1,2 & 1,4,4 & 1,2,6 & 1,6,3 & 1,1,6 & 1,5,4\\
	0,1,1 & 1,3,2 & 1,3,5 & 1,5,2 & 0,1,5 & 1,0,3 & 1,5,3\\
	1,2,1 & 1,3,4 & 1,2,3 & 1,0,5 & 1,6,5 & 1,3,3 & 1,6,4\\
	1,1,5 & 1,0,1 & 1,4,1 & 1,5,5 & 0,1,3 & 1,5,6 & 0,1,6\\
	1,3,6 & 1,1,1 & 1,1,3 & 1,4,2 & 1,4,3 & 1,3,1 & 1,0,6\\
	1,4,5 & 1,1,4 & 1,2,5 & 1,4,6 & 1,0,2 & 1,1,2 & 1,6,6\\
	1,6,1 & 1,6,2 & 1,0,4 & 1,2,4 & 1,2,2 & 1,5,1 & 0,1,4
	\end{array}$\\
	A selector is\\
	$e,$ $a,$ $b,$ $ab^4,$ $a^3b^4,$ $a^3b^6,$ $a^6b^6.$ points on [1,0,0].

	The conics are\\
	$\begin{array}{ccccccc}
	0,0,1 & 0,1,1 & 1,2,1 & 1,1,5 & 1,3,6 & 1,4,5 & 1,6,1\\
	0,0,1 & 0,1,2 & 1,4,4 & 1,2,6 & 1,6,3 & 1,1,6 & 1,5,5\\
	0,0,1 & 0,1,3 & 1,6,2 & 1,3,3 & 1,2,5 & 1,5,3 & 1,4,2\\
	0,0,1 & 0,1,4 & 1,1,2 & 1,4,3 & 1,5,5 & 1,2,3 & 1,3,2\\
	0,0,1 & 0,1,5 & 1,3,4 & 1,5,6 & 1,1,3 & 1,6,6 & 1,2,4\\
	0,0,1 & 0,1,6 & 1,5,1 & 1,6,5 & 1,4,6 & 1,3,5 & 1,1,1\\
	0,0,1 & 1,1,4 & 1,4,1 & 1,5,2 & 1,2,2 & 1,3,1 & 1,6,4
	\end{array}$\\
	The line is [0,1,0] with points
	$( 1,0,1) \cdot (1,0,c) = (1,0,\frac{c}{1+c)} )$\\
	$\begin{array}{ccccccc}
	0,0,1& 1,0,1& 1,0,4& 1,0,5& 1,0,2& 1,0,3& 1,0,6.\\
	\end{array}$\\
	Other points are on [0,0,1].

	They all have a contact of order 2 at (0,0,1) with tangent [1,0,0].

	\ssec{The case of 3 distinct roots.}\label{sec-Stri3roots}
	\vspace{-18pt}\hspace{200pt}\footnote{21.2.86}\\[8pt]
	\sssec{Definition.}
	If the roots are $a,$ $b,$ $c,$
	\enumb
	\item The polynomial which has 2 of the roots corresponds to a point
		called {\em isotropic point}.
	\item The 3 lines through 2 of the 3 isotropic points are called
		{\em isotropic lines}.
	\item Any non isotropic line through an isotropic point is called an
		{\em ideal line}.
	\item Any non isotropic point on an isotropic line is called an
		{\em ideal point}.
	\enume

	\sssec{Example.}
	$p = 7,$ $a = 1,$ $b = 2,$ $c = 4$\footnote{earlier version}.
	\enumb
	\item The isotropic points are $A_0 = (1,1,1),$ $A_1 = (1,2,4),$
	$A_2 = (1,4,2).$
	\item The isotropic lines are $a_0 = [1,1,1],$ $a_1 = [1,4,2],$
	$a_2 = [1,2,4].$
	\item The generators of the group are $\alpha$  = (0,1,2) and
	$\beta = (0,1,1).$
	\item The ideal lines through $A_0$ are\\
	\hth$    [1,3,0] = S_0 = \{e,b^2,b^4, a^3,a^3b^2,a^3b^4\},$\\
	\hth$    [1,4,3] = b S_0,$\\
	\hth$    [0,1,3] = a S_0,$\\
	\hth$    [1,5,6] = ab S_0,$\\
	\hth$    [1,0,5] = a^2 S_0,$\\
	\hth$    [1,6,2] = a^2b S_0.$
	\item The ideal lines through $A_1$ are\\
	\hth$    [1,5,0] = S_1 = \{e,a^2,a^4, ab^3,a^3b^3,a^5b^3\},$\\
	\hth$    [1,2,6] = a S_1,$\\
	\hth$    [1,3,4] = b S_1,$\\
	\hth$    [1,0,3] = ab S_1,$\\
	\hth$    [0,1,5] = b^2 S_1,$\\
	\hth$    [1,6,5] = b^2a S_1.$
	\item The ideal lines through $A_2$ are\\
	\hth$    [1,6,0] = S_2 = \{e,a^2,a^2b^2, a^4b^4,a^4b,b^3,a^2b^5\},$\\
	\hth$    [1,1,5] = ab S_2,$\\
	\hth$    [1,5,1] = a S_2,$\\
	\hth$    [1,0,6] = a^2b S_2,$\\
	\hth$    [0,1,6] = b S_2,$\\
	\hth$    [1,3,3] = ab^2 S_2.$
	\item A selector is $(e,b^5,ab^5,a^2b^3,a^5b^3)$ giving the points\\
	\hth$    (0,0,1), (1,6,1), (1,6,3), (1,6,6), (1,6,4)$ on $[110],$\\
	    the ideal points on this line are (1,6,0), (1,6,5), (1,6,2).\\
	    The 36 other lines are obtained by multiplication by any of the
	    elements in the group, e.g. if we multiply to the left by $a^5b^5,$ 
	    the points are (0,1,4), (1,6,4), (1,2,2), (1,3,6), (1,0,1), and
	    the ideal points are (1,1,5), (1,5,0), (1,4,3).
	\enume

	\sssec{Notes.}
	We have therefore the following operations:\\
	\hth$	l := P \times Q,$ $L := p \times q$\\
	\hth$	R := P \bullet Q,$ $r := p \bullet q?$\\
	where $\bullet$  is done modulo a polynomial of degree 3.
	If the polynomial is primitive the properties are well known, what is
	probably new is what happens when the polynomial is not primitive.
	If it has 3 roots it makes sense to normalize to have the isotropic
	points at (1,0,0), (0,1,0) and (0,0,1), but I do not see how this can
	be done in view of the fact that an isotropic point corresponds to
	$(I-a)(I-b).$

	\sssec{Definition.}
	Given a polynomial of the third degree with 3 distinct roots, a
	{\em line generator} is a generator of a cyclic group of order $p-1$
	whose elements correspond to $p-1$ points of a line through one of the
	isotropic points, the last point is an ideal point on the isotropic
	line which does not belong to the isotropic point.

	\sssec{Definition.}
	Two line generators are said to be {\em independent} if they are
	associated to lines through distinct isotropic points.

		Why did I not worry about this when constructing an example and
		use simply distinct lines?

	\ssec{Conjecture.}
	\vspace{-18pt}\hspace{100pt}\footnote{24.2.86}\\[8pt]
	{\em Given a polynomial of the third degree with 3 distinct
	roots, there exists 2 independent line generators.}

	\sssec{Comment.}  I will choose the roots to be 0, 1 and -1.
	$P_3 = I^3 - I.$\\
	The isotropic points are $A_0 = (1,0,-1),$ $A_1 = (1,1,0),$
	$A_2 = (1,-1,0).$\\
	The isotropic lines are $a_0 = [0,0,1],$ $a_1 = [1,1,1],$
	$a_2 = [1,-1,1].$

	\sssec{Conjecture.}  With the choice just given, there exist
	an $x$ such that\\
	if $y = -2(x+1),$ $(1,0,x),$ $(1,1,1y),$ $(x+y+1,x+1,xy)$ are line
	generators corresponding to lines $[0,1,0],$ $[1,-1,0],$
	$[1,1,-\frac{2x+y+2}{xy}]$ through $A_0,$ $A_1$ and $A_2.$

	\sssec{Example.}
	$\begin{array}{lccccccc}
	p = 3,&\vline&0,0,1 &1,1,2\\
	\hth&\vline& 1,0,1 &1,2,2\\
	\cline{1-4}
	p = 5,&\vline&0,0,1 &1,1,1 &1,1,4 &1,1,2\\
	\hth&\vline&1,0,1 &1,4,2 &1,2,4 &1,3,3\\
	\hth&\vline&1,0,2 &1,2,3 &1,4,4 &0,1,3\\
	\hth&\vline&1,0,3 &0,1,2 &1,3,4 &1,4,1\\
	\cline{1-6}
	p = 7,&\vline &roots &1,2,4,\\
	\hth&\vline&	0,0,1 &1,2,0 &1,2,1 &1,2,5 &1,2,2 &1,2,3\\
	\hth&\vline&	1,0,1 &1,3,2 &1,5,5 &1,4,0 &1,1,6 &1,6,3\\
	\hth&\vline&	1,4,4 &0,6,4 &1,1,4 &1,5,4 &0,1,0 &1,0,4\\
	\hth&\vline&	1,1,3 &1,0,2 &1,6,1 &1,4,1 &0,1,1 &1,3,5\\
	\hth&\vline&	1,4,1 &1,5,3 &1,0,0 &0,1,2 &1,1,2 &1,3,6\\
	\hth&\vline&	1,6,6 &1,5,2 &1,1,0 &1,4,5 &1,3,1 &0,1,4
	\end{array}$ $I^3 -1 = 0$

	The lines are obtained form\\
	\hth$	(0,0,1),(0,1,2),(0,1,0),(0,1,1),(0,1,4)$\\
	or	0,0	4,3	2,4	3,4	5,5\\
	for instance, adding 2,3 modulo 6,6\\
	\hth$	2,3	0,0,	4,1	5,1	1,2$\\
	or	(1,5,4),(0,0,1),(1,5,3),(1,5,2),(1,5,5) on [1,4,0].

	The $c$-lines are all obtained from\\
	\hth$	(0,0,1),(0,1,4),(1,0,2),(1,5,5),(1,6,4)$\\
	or\hti{6}0,0\hti{5}5,5\hti{5}3,1\hti{5}1,2\hti{5}2,1\\
	for instance, adding 3,2 modulo 6,6 gives\\
	\hth$	3,2\hti{5}2,1\hti{5}0,3\hti{5}4,4\hti{5}5,3$\\
	or	(1,6,1),(1,6,4),(1,2,5),(1,1,2),(1,4,5).

	$p = 11,$ line generators: (1,0,1), (1,1,7), (1,10,2).

	$p = 13,$ line generators: (1,0,1), (1,1,9), (1,12,2).

	$p = 17,$ line generators: (1,0,2), (1,1,11), (1,16,4).

	\sssec{Definition.}
	A {\em selector} is a set of $p-2$ elements $P^i_kQ^j_k$ which are
	on an ordinary line.

	\sssec{Theorem.}
	{\em Given 2 independent line generators $P$ and $Q,$
	the isotropic lines are obtained as cosets of the cyclic groups
	generated by $P,$ $Q$ and $P \bullet Q$.\\
	The ordinary lines are obtained by multiplication modulo $P_3$,
	$P^lQ^m$ by the elements of a selector.}

	We may want to put this in a section on {\em triangular geometry}.

	In this geometry we have ordinary, ideal and isotropic points,
	ordinary, ideal and isotropic lines, and $c$-lines.
	These are represented by conics through the isotropic points,
	the ideal $c$-lines are the degenerate conics consisting of an isotropic
	line and an ideal line to the oposite isotropic point.
	The isotropic $c$-lines are the degenerate conics consisting of two
	isotropic lines.
	The lines and the $c$-lines can be interchanged.  If the $c$-lines are
	considered as lines, then the lines are $c$-lines, in other words if
	we start with a geometry where we define the conics through 3 given
	points as lines, the conics are represented by lines.  Pascal's
	Theorem gives the following,  consider a line $l,$ with points $P_{2i}$
	on $a_i,$ and 3 other points $P_1,$ $P_3,$ $P_5,$ this line can be
	considered a $c$-conic,
	indeed, the $c$-lines through successive points are the degenerate
	$c$-lines or ideal $c$-lines,
	$a_0 = (A_0 \times P_1) + a_0,$ $a_1 = (A_2 \times P_1) + a_2,$
	$a_2 = (A_2 \times P_3) + a_2,$
	$a_3 = (A_1 \times P_3) + a_1,$ $a_4 = (A_1 \times P_5) + a_1,$
	$a_5 = (A_0 \times P_5) + a_0.$\\
	The $c$-Pascal points are $Q_0 = (A_1 \times P_3)
	\times (A_0 \times P_1),$
	$Q_1 = (A_2 \times P_1) \times (A_1 \times P_5),$
	$Q_2 = (A_0 \times P_5) \times (A_2 \times P_3). $\\
	These points are on a conic, with $A_0,$ $A_1,$ $A_2$ because the
	Pascal line for the sequence $A_0,$ $Q_0,$ $A_1,$ $Q_1,$ $A_2,$ $Q_2,$
	gives the Pascal points
	$P_1,$	$P_3,$ $P_5.$  This can be used to study what could be called a
	{\em bi-triangular geometry}.

	\sssec{Comment.}  The $c$-lines can be deduced from the line by
	the transformation which associates, in general, to
	$(X_0,$ $X_1,$ $X_2),$ $(X_1 X_2,$ $X_2 X_0,$ $X_0 X_1)$, 
	a conic becomes then a quadric with double points, isolated or not in
	the case or a real field.  A conic through $A_1,$ $A_2,$ but not
	through $A_0,$ becomes a quadric which degenerates in $a_1,$ $a_2$ and
	a conic through $A_1,$ $A_2$ but not $A_0.$ 

	\sssec{Comment.}
	We could choose as isotropic points, in a model of this
	geometry in the Euclidean plane, with Cartesian coordinates, by
	choosing one of them at the origin, and the 2 others at the direction
	of the axis.  The $c$-lines are then hyperbolas passing through the
	origin, with asymptotes in the direction of the axis.

	\sssec{Problem.}
	Study the axiomatic of the triangular geometry and obtain Theorems in
	it.  Circles could be conics through 2 of the isotropic points.

	\sssec{Problem.}
	Study the axiomatic of the triangular bi-geometry and obtain Theorems
	in it.

	\sssec{Comment.}
	The analysis can be repeated in the form of Euclidean
	geometry by considering the non-homogeneous points $(x,y)$ and the
	homogeneous lines $[a,b,c].$  This can be done directly or infered
	from the cases 1, 2 and 4 above, with one of the isotropic lines
	playing the role of the line at infinity in Euclidean geometry.

	\ssec{Notes.}
	In G45, I give a special case of the following Theorem valid when
	$s_1 = a = 0,$ this generalizes the Theorem, with $b = s_{11}$ and
	$c = s_{111}.$\\
	It was obtained earlier.

	\sssec{Theorem.}
	{\em The symmetric functions of the roots are}\\
	\hth$s_1 := \rho_0 + \rho_1 + \rho_2 = a,$\\
	\hth$s_{11} := \rho_1\rho_2 + \rho_2\rho_0 + \rho_0\rho_1 = b,$\\
	\hth$s_{111} : = \rho_0\rho_1\rho_2 = c,$\\
	\hth$s_2 := \rho_0^2 + \rho_1^2 + \rho_2^2 = a^2 - 2b,$\\
	\hth$s_{21} := \rho_0^2(\rho_1+\rho_2) + \rho_1^2(\rho_2+\rho_0)
		+ \rho_2^2(\rho_0+\rho_1) = ab-3c,$\\
	\hth$s_3 := \rho_0^3 + \rho_1^3 + \rho_2^3 = a(a^2-3b) + 3c,$\\
	\hth$s_{211} := ac,$\\
	\hth$s_{22} := b^2-2ac,$\\
	\hth$s_{31} := a(ab-c) - 2 b^2,$\\
	\hth$s_4 := a(a^3-4ab+4c) + 2 b^2.$

	\sssec{Theorem.}
	\enumb
	\item{\em The conic which pass through the isotropic points is}\\
	\hth$	  k ( (b^2-2ac) X_0^2 + b X_1^2 + 3 X_2^2$\\
	\hth$		+ 2a X_1X_2 + 2(a^2-2b) X_2X_0 - (3c-ab) X_0X_1 )$\\
	\hth$	+ l ( bc X_0^2 + 3c X_1^2 + a X_2^2$\\
	\hth$		+ 2b X_1X_2 - (3c-ab) X_2X_0 + 2ac X_0X_1 )$\\
	\hth$	+ m ( 3 c^2 X_0^2 + ac X_1^2 + (a^2-2b) X_2^2$\\
	\hth$	- (3c-ab) X_1X_2 + 2(b^2-2ac) X_2X_0 + 2bc X_0X_1 ) = 0,$
	\item{\em which is tangent to the isotropic lines is}\\
	\hth$	  k ( 3 x_0^2 + (a^2+b) x_1^2 + ac x_2^2$\\
	\hth$		- (3c+ab) x_1x_2 + 2b^{} x_2x_0 - 4a x_0x_1 )$\\
	\hth$	+ l ( a x_0^2 + a(a^2-2b)+3c) x_1^2 + (a^2-2b)c x_2^2$\\
	\hth$	- (a(ab+c)-2b^2) x_1x_2 - (3c-ab) x_2x_0 - 2(a^2-b) x_0x_1 )$\\
	\hth$	+ m ( (a^2-2b) x_0^2 + a(a^3-3ab+4c) x_1^2
		+ (a(a^2-3b)+3c)c x_2^2$\\
	\hth$	+ (a(-a^2b+3b^2-ac)-bc) x_1x_2 + (a(ab-c)-2b^2) x_2x_0$\\
	\hth$			- (a(2a^2-5b)+3c) x_0x_1 ) = 0.$
	\enume

	Proof:\\
	The degenerate conics through the isotropic points are\\
	\hth$	  \alpha _0 ( \rho_1^2 X_0 + \rho_1 X_1 + X_2)
	( \rho_2^2 X_0 + \rho_2 X_1 + X_2) )$\\
	\hti{12}$+ \alpha _1 ( \rho_2^2 X_0 + \rho_2 X_1 + X_2)
	( \rho_0^2 X_0 + \rho_0 X_1 + X_2) )$\\
	\hti{12}$	+ \alpha _2 ( \rho_0^2 X_0 + \rho_0 X_1 + X_2)
	( \rho_1^2 X_0 + \rho_1 X_1 + X_2) ) = 0.$\\
	If we choose, in succession, $\alpha_0$ = $\alpha_1$ = $\alpha_2$ = 1,
	$\alpha_0$ = $\rho_0,$ $\alpha_1$ = $\rho_1,$ $\alpha_2$ = $\rho_2,$
	and $\alpha_0 = \rho_0^2,$ $\alpha_1 = \rho_1^2,$ $\alpha_2 = \rho_2^2,
$\\ 
	we obtain respectively the expressions whose coefficients are $k,$ $l$
	and $m.$  Similarly for the $c$-points we start with the degenerate
	conics tangent to the isotropic lines, which are\\
	\hth$ \alpha_0 (x_0 -(\rho_2 + \rho_0)x_1 + \rho_2 \rho_0 x_2)
	(x_0 -(\rho_0 + \rho_1)x_1 + \rho_0 \rho_1 x_2))$\\
	\hti{12}$+ \alpha_1 (x_0 -(\rho_0 + \rho_1)x_1 + \rho_0 \rho_1 x_2)
	(x_0 -(\rho_1 + \rho_2)x_1 + \rho_1 \rho_2 x_2))$\\
	\hti{12}$+ \alpha_2 (x_0 -(\rho_1 + \rho_2)x_1 + \rho_1 \rho_2 x_2)
	(x_0 -(\rho_2 + \rho_0)x_1 + \rho_2 \rho_0 x_2)) = 0.$

	The following Theorem was develppoed to prove the relation between
	the conics associated to the co, bi and semi-selectors but were
	found not to be needed.  It is now an answer to an exercise.

	\sssec{Answer to exercise   .}
	Let $s_1$ = 0, the conic through (0,1,0), (0,0,1)
	\enumb
	\item   which passes through the isotropic points is

	\hth$	s_{111} X_0^2 - X_1X_2 = 0.$
	\item   which is tangent to the isotropic lines is

	\hth$	x_0^2 - s_{111} x_1x_2 + s_{11} x_2x_0 = 0.$\\
	or\\
	\hth$	(s_{111}X_0 + s_{11} X_1)^2 - 4 s_{111} X_1 X_2 = 0,$
	\item   which has the isotropic triangle as polar triangle is

	\hth$	s_{111} X_0^2 - 2 s_{11} X_0X_1 + 2 X_1X_2 = 0.$
	\enume

	Proof:   Using
	\enumb
	\setcounter{enumi}{2}
	\item 	$\rho_1 +\rho_2 = - \rho_0,$ we can check
	\enume
	for 0, $\rho_0\rho_1\rho_2 - (\rho_1+\rho_2)\rho_1\rho_2 = 0,$\\
	for 1, $\rho_0^4 - \rho_0\rho_1\rho_2\rho_0 +
	(\rho_1\rho_2+\rho_2\rho_0+\rho_0\rho_1)\rho_0^2 = 0,$\\
	for 2, $\vet{\rho_0^3}{\rho_0^2}{\rho_0} =
	\mat{s_{111}}{-s_{11}}{0}{-s_{11}}{0}{1}{0}{1}{0}
	\vet{1}{-\rho_1-\rho_2}{\rho_1\rho_2}.$

	\sssec{Exercise.}
	Let $s_1 = 0.$ Determine the conic through (0,1,0), (0,0,1),
	\enumb
	\item   which passes through the isotropic points,
	\item   which is tangent to the isotropic lines,
	\item   which has the isotropic triangle as polar triangle.
	\enume

	\sssec{Comment.}
	To obtain the statements of the preceding Theorem, I will
	illustrate for the case 2.  If the conic is represented by the
	symmetric matrix\\
	\hth$\matb{(}{ccc}1&\gamma&\beta\\ \gamma&0&\alpha\\
		\beta&\alpha&0\mate{)}$.\\
	The condition that $I_0$ is the pole of $i_0$ gives\\
	\hth$\mu\:\rho_0^2 = 1 - \gamma(\rho_1+\rho_2) + \beta\rho_1\rho_2,$\\
	\hth$\mu\: \rho_0  = \gamma  + \alpha  \rho_1\rho_2,$\\
	\hth$\mu   = \beta  - \alpha (\rho_1+\rho_2).$\\
	Eliminating $\mu$  from the first 2 and the last 2 equations gives\\
	\hth$1 - \gamma(\rho_1+\rho_2) + \beta\rho_1\rho_2 - \gamma\rho_0
	- \alpha  \rho_0\rho_1\rho_2 = 0,$\\
	\hth$\gamma  + \alpha  \rho_1\rho_2 - \beta  \rho_0
	+ \alpha (\rho_0\rho_1+\rho_0\rho_2) = 0,$\\
	or using 3,\\
	\hth$1 + \beta \rho_1\rho_2 - \alpha  s_{111} = 0,$\\
	\hth$\gamma  + \alpha  s_{11} - \beta  \rho_0 = 0.$\\
	Because the conic cannot depend on individual values of $\rho_0,$
	$\rho_1,$ $\rho_2,$ 
	$\beta$	 = 0 and then  $\alpha = \frac{1}{s_{111}}$ and
	$\gamma = - \alpha s_{11}.$

	\ssec{On the tetrahedron.}
	\sssec{Example.}
	Let the roots be 0,1,2,3 and p = 5,\\
	the isotropic points are (1,-1,1,-1), (1,0,1,0), (1,1,-2,0), (1,2,2,0).
\\
	$P_4 = I^4 - I^3 + I^2 - I,$ therefore, $I^5 = I \pmod{P_4}$ and
	$I^6 = I^2 {\mathit mod} P_3.$\\ 
	The cubic surface is given by\\
	$\matb{|}{cccc}
		k	& l	& m	&n\\ 
		Z+T	&Y+Z	&X+Y	&X\\
		X-Z	&T-Y	&Z-x	&Y\\
		Y+Z	&X+Y	&X+T	&Z\mate{|} = 0.$\\
	For instance, a point on $I_0 \times I_1$ is $(u+v,-u,u+v,-u)$ and\\
	$\matb{|}{cccc}
		k	&l	&m	& n \\
		v	&v	&v	&u+v\\
		0	&0	&0	&-u \\
		v	&v	&v	&u+v\mate{|} = 0,$\\
	because 2 rows are equal.\\
	Similarly for a point on $I_2 \times I_3,$ $(u+v,2u+v,2u-2v,0),$\\
	$\matb{|}{cccc}
		  k	&  l	&  m	&  n  \\
		2u-2v	& 4u- v	&3u+2v	& u+ v\\
		-u+3v	&-2u- v	& u-3v	&2u+ v\\
		4u- v	& 3u+2v	& u+ v	&2u-2v\mate{|} = 0,$
	because the sum of the last 3 rows is equal to $0 \pmod{5}.$\\
	No other conditions are needed to obtain a family of cubics with 3
	parameters because if the isotropic points are chosen as (1,0,0,0),
	(0,1,0,0), (0,0,1,0), (0,0,0,1) the cubic is\\
	\hth$d_0 X_1X_2X_3 + d_1 X_2X_3X_0 + d_2 X_3X_0X_1 + d_3 X_0X_1X_2 = 0.$

	\sssec{Example.}
	\footnote{14.4.86}
	For the semi-transformation, let $p = 5$ and
	$P_4 = I^4 - I^3 + I^2 - I,$
	$(X,Y,Z,T)^2 = (x,y,z,t)$ gives\\
	\hth$	x = 2XT + 2YZ + (2XZ+Y^2)$\\
	\hth$	y = 2YT + Z^2 - (2XZ+Y^2) + X^2,$\\
	\hth$	z = 2 ZT + (2XZ+Y^2) + 2 XY,$\\
	\hth$	t = T^2.$\\
	A plane $\{k,l,m,n\}$ is therefore transforned in the quadric
	represented by the symmetric matrix\\
	$\matb{(}{cccc}
		  l&  m&k-l+m&k\\
		  m&k-l+m&k&  l\\
		k-l+m&k&  l&  m\\
		  k&  l&  m&  n\mate{)}.$
	The isotropic points are (1,-1,1,-1), (1,0,1,0), (1,1,-2,0), (1,2,2,0),
\\
	and the corresponding isotropic planes are\\
	\{0,0,0,1\}, \{1,1,1,1\}, \{-2,-1,2,1\}, \{2,-1,-2,1\}.\\
	It is easy to check the latter are the polar of the former,
	independently from $k,$ $l,$ $m,$ $n.$  It is easy to verify that the
	quadric which have the isotropic tetrahedron as polar tetrahedron form a
	3 parameter family and that this generalizes to n dimensions.

	\sssec{Exercise.}
	Study the relation which exist between the
	correspondance between a pair of points and the pair obtained at the
	intersection of the tangents at the 2 points to the $c$-line through
	these points and the intersection of the $c$-lines tangent to the line
	through the 2 points at the two points.\\
	Hint: Study first how to
	obtain from the point on any other line through a point and the $c$-line
	that line through the point which is tangent to the $c$-line.

	\sssec{Definition.}
	Let $l^c$ be a line through $P,$  $\ldots$ .

\chapter{GENERALIZATION TO 3 DIMENSIONS}



\setcounter{section}{-1}
	\setcounter{section}{-1}
	\section{Introduction.}

	I will sketch here part of the generalization to 3 dimensions of what
	has been presented in the preceding parts.  It will be obvious how to
	generalize further to $n$ dimensions.  After a brief look at the
	history, I will review the application of Grassmann algebra to the
	incidence properties of the fundamental objects in 3 dimensions, the
	points, the lines and the planes.

	The finite polar geometry will be introduced in section
	\ref{sec-Spolar}. It is obtain by prefering a plane, the ideal plane,
	to which correspond the notions of affine geometry, parallelism,
	mid-points, equality of segments on parallel lines, and a quadric, the
	fundamental quadric, which, together with the ideal plane, allow for
	the definition of spheres and therefore equality of distances between
	unordered pairs of points as well as orthogonality or more generally
	equality of angles between ordered pair of lines.\\
	To illustrate properties in 3 dimensions, the geometry
	of the triangle in involutive geometry will be generalized, in section
	\ref{sec-Sgentetra} to the study of the general tetrahedron in finite
	polar geometry. In the classical case, the first work on
	the subject is that of Prouhet, this was followed by important memoirs
	of Intrigila  and Neuberg.\\
	We will see that a special case occurs very naturally, that of the
	orthogonal tetrahedron, studied in section \ref{sec-Sorthotetra}.
	We will see that the success of the
	theory of this special case is explained by the generalization to 3
	dimension of the symmetry which exists in 2 dimensions when we exchange
	the barycenter and the orthocenter.\\
	The isodynamic tetrahedron is studied in section \ref{sec-Sisodtetra}.

	The generalization of many other 3 dimensional and $n$ dimensional
	concepts is left to the reader.

	This part ends with an introduction to the anti-polar geometry
	\ref{sec-Santipol}.

	\ssec{Relevant historical background.}

	\setcounter{subsubsection}{-1}
	\sssec{Introduction.}
	In the classical case, the extension to 3 dimensions is already given
	by Euclid.  The earlier definitions of conics derives from the circular
	cone in 3 dimensions.  Of note is also the fact that the 2-dimensional
	Desargues' theorem derives directly from the incidence properties in 3
	dimensions.  Although the algebraic notation of analytic geometry,
	introduced by Descartes immediately extends and at once suggests to go
	beyond the observable to 4 and to $n$ dimensions, it is not suitable if
	we progress from finite polar geometry - where equality of distance and
	angles are defined and are not primary notions - to finite 3
	dimensional Euclidean geometry.  Instead, I will use the notation of
	exterior algebra introduced by Grassmann.

	\ssec{Grassmann algebra applied to incidence properties of
	points, lines and planes}

	\setcounter{subsubsection}{-1}
	\sssec{Introduction.}
	After introducing in \ref{sec-dgras} the algebraic representation of
	points, lines and planes in ${\Bbb Z}_p^3$, I recall the basic concepts
	and properties of the exterior algebra of Grassmann (\ref{sec-ngras}
	to \ref{sec-ngras2}), I define the incidence relations
	(\ref{sec-dgras3}) and derive the associated properties
	(\ref{sec-tgras3}, \ref{sec-tgras4} to \ref{sec-tgras7}).

	\sssec{Definition.}\label{sec-dgras}
	The {\em points} and {\em planes} in 3 dimensions will be represented
	using 4 homogeneous coordinates.  (Not all coordinates are 0, and
	if all coordinates are multiplied modulo $p$ by the same non zero
	element in ${\Bbb Z}_p,$ we obtain the same point or plane.)\\
	Points will be denoted by a capital letter and the coordinates will be
	placed between parenthesis.  Planes will be denoted by a capital letter
	preceded by the symbol ``$|$" or by a calligraphic letter
	and the coordinates will be placed between braces.\\
	The {\em lines} will be represented by 6 homogeneous coordinates
	$[l_0,l_1,l_2,l_3,l_4,l_5],$
	such that $l_0 l_5 + l_1 l_4 + l_2 l_3 = 0.$\\
	This part of the definition
	will be justified in \ref{sec-tgras3}.4 and \ref{sec-tgraspl}.\\
	The normalization will again be such that the leftmost non zero
	coordinate is 1.

	\sssec{Example.}
	Let $p = 7,$\\
	$P := (2,4,6,1) = (1,2,3,4),$ $l := [3,3,1,4,3,5] = [1,1,5,6,1,4],$\\
	${\cal Q} := \{5,1,1,1\} = \{1,3,3,3\}$ are respectively a point, aline
	and a plane.

	\sssec{Notation.}\label{sec-ngras}
	To define algebraically the incidence properties, I will use
	Grassmann algebra with exterior product multiplication.
	If $e_0,$ $e_1,$ $e_2,$ $e_3$ are ``unit" vectors, we write
	\enumb
	\item 	$P := (P_0,$ $P_1,$ $P_2,$ $P_3)
		  := P_0\: e_0 + P_1\: e_1 + P_2\: e_2 + P_3\: e_3.$
	\item 	$l := [l_0,$ $l_1,$ $l_2,$ $l_3,$ $l_4,$ $l_5]$\\
	\hti{4}$:= l_0\: e_0 \vee e_1 + l_1\: e_0 \vee e_2 + l_2\: e_0 \vee
		 e_3 + l_3\: e_1 \vee e_2 + l_4\: e_3 \vee e_1 + l_5\: e_2 \vee
		e_3,$\\
	with
	\item 	  $l_0$ $l_5$ + $l_1$ $l_4$ + $l_2$ $l_3 = 0.$
	\item${\cal Q} := \{{\cal Q}_0,$ ${\cal Q}_1,$ ${\cal Q}_2,$
	${\cal Q}_3\}\\
	\hti{4}:= {\cal Q}_0\: e_1 \vee e_2 \vee e_3 + {\cal Q}_1\: e_3 \vee e_2
	\vee e_0 + {\cal Q}_2\: e_3 \vee e_0 \vee e_1 + {\cal Q}_3\: e_1 \vee
	e_0\vee e_2.$\\
	In each case not all coefficients are zero.\\
	The specific notation for $l$ and ${\cal Q}$ will is justified in
	\ref{sec-tgras2}.  For ${\cal Q}$ the order of the unit vectors is 
	chosen in such a way that the last ones are consecutive, $e_3$, $e_0,$
	$e_1$, $e_2$.
	If condition 2 is not satisfied, the 2-form will be denoted using
	an identifier starting with a lower case letter and followed by
	"$'$ ".
	If an identity is satisfied for the general 2-form $l'$ as well as for
	the line $l,$ I will use the notation $l'$ (see for instance
	\ref{sec-tgras1}).
	\enume

	I recall:

	\sssec{Definition.}\label{sec-dgrasext}
	The {\em exterior product} is defined by using the usual rules
	of algebra, namely, commutativity, associativity, neutral element
	property and distributivity with the exception\\
	\hth$	e_i \vee e_j = - e_j \vee e_i$ which gives, in particular,
	$e_i \vee e_i = 0.$

	\sssec{Lemma.}
	\hth{\em $(e_i \vee e_j) \vee (e_k \vee e_l) = (e_k \vee e_l) \vee
	(e_i \vee e_j).$}

	\sssec{Theorem.}\label{sec-tgras1}
	\hth{\em $P \vee Q = - Q \vee P,\: l' \vee m' =  m' \vee l'.$}

	\sssec{Corollary.}\label{sec-cogras}
	\hth{\em $P \vee P = 0.$}

	\sssec{Definition.}
	Given any expression involving points, lines or planes
	using the Grassmann representation, the {\em dual} of an expression is
	obtained by replacing the coefficient by itself and\\
	\hth$	e_{i_0} \vee \ldots \vee e_{i_{k-1}}$  by  $j\: e_{i_k}
	\vee \ldots \vee e_{i_3}$\\
	where $i_0,$ \ldots, $i_{k-1},$ $i_k,$ \ldots $i_3$ is a permutation of
	0,1,2,3
	and $j = 1$ if the permutation is even, $-1$ if the permutation is odd.

	\sssec{Theorem.}\label{sec-tgras2}
	\enumb
	\item$dual(P) = P_0\:e_1\vee e_2\vee e_3 + P_1\:e_3 \vee e_2 \vee e_0
		+ P_2\: e_3 \vee e_0 \vee e_1 + P_3\: e_1 \vee e_0 \vee e_2.$
	\item$dual(l') = l_0\:e_2\vee e_3 + l_1\:e_3\vee e_1
		+ l_2\:e_1\vee e_2 + l_3\: e_0 \vee e_3
		+ l_4\: e_0 \vee e_2 + l_5\: e_0 \vee e_1.$
	\enume

	Because of the notation \ref{sec-ngras}.1, duality, for a line, simply
	reverses the order of the components of $l.$

	\sssec{Notation.}\label{sec-ngras2}
	\enumb
	\item$P \cdot {\cal Q} := {\cal Q} \cdot P := dual( P \vee {\cal Q} ),$
	\item$l' \wedge {\cal P} := {\cal P} \wedge l'
		:= dual( dual({\cal P}) \vee dual(l') ).$
	\item${\cal P} \wedge {\cal Q} := {\cal Q} \wedge {\cal P}
		:= dual( dual({\cal P}) \vee dual({\cal Q}) ).$
	\enume

	\sssec{Definition.}\label{sec-dgras3}
	A {\em point} $P$ {\em is incident to a line} $l$  iff\\
	\hth$	P \vee l = 0.$\\
	A {\em point} $P$ {\em is incident to a plane} ${\cal Q}$  iff\\
	\hth$	P \vee {\cal Q} = 0.$\\
	A {\em line} $l$ {\em is incident to a plane} ${\cal Q}$  iff\\
	\hth$	l \wedge {\cal Q} = 0.$

	\sssec{Theorem.}\label{sec-tgras3}
	\enumb
	\item~~$P \vee Q = (P_0 Q_1 - P_1 Q_0)\: e_0 \vee e_1
		+ (P_0 Q_2 - P_2 Q_0) \:e_0 \vee e_2$\\
	\hth$+ (P_0 Q_3 - P_3 Q_0)\: e_0 \vee e_3 + (P_1 Q_2 - P_2 Q_1)\: e_1
	\vee e_2$\\
	\hth$+ (P_3 Q_1 - P_1 Q_3)\: e_3 \vee e_1 + (P_2 Q_3 - P_3 Q_2)\: e_2
	\vee e_3.$
	\item~~${\cal P} \vee {\cal Q}
	 = ({\cal P}_2 {\cal Q}_3 - {\cal P}_3 {\cal Q}_2)\: e_0 \vee e_1
	({\cal P}_3 {\cal Q}_1 - {\cal P}_1 {\cal Q}_3)	+  \:e_0 \vee e_2$\\
	\hth$+ ({\cal P}_1 {\cal Q}_2 - {\cal P}_2 {\cal Q}_1)\: e_0 \vee e_3
	 + ({\cal P}_0 {\cal Q}_3 - {\cal P}_3 {\cal Q}_0)\: e_1\vee e_2$\\
	\hth$+ ({\cal P}_0 {\cal Q}_2 - {\cal P}_2 {\cal Q}_0)\: e_3 \vee e_1
	+ ({\cal P}_0 {\cal Q}_1 - {\cal P}_1 {\cal Q}_0)\: e_2	\vee e_3.$
	\item~~$P \vee l' = (P_1 l_5 + P_2 l_4 + P_3 l_3)\:e_1\vee e_2 \vee e_3$\\
	\hth$	+ (-P_0 l_5 + P_2 l_2 - P_3 l_1)\: e_3 \vee e_2 \vee e_0$\\
	\hth$	+ (-P_0 l_4 - P_1 l_2 + P_3 l_0)\:e_3\vee e_0 \vee e_1$\\
	\hth$	+ (-P_0 l_3 + P_1 l_1 - P_2 l_0)\: e_1 \vee e_0 \vee e_2.$
	\item~~${\cal P} \vee l' = ({\cal P}_1 l_0 + {\cal P}_2 l_1
		+ {\cal P}_3 l_2)\:e_1\vee e_2\vee e_3$\\
	\hth$	+ (-{\cal P}_0 l_0 + {\cal P}_2 l_3 - {\cal P}_3 l_4)
		\: e_3 \vee e_2 \vee e_0$\\
	\hth$	+ (-{\cal P}_0 l_1 - {\cal P}_1 l_3 + {\cal P}_3 l_5)
		\:e_3\vee e_0 \vee e_1$\\
	\hth$	+ (-{\cal P}_0 l_2 + {\cal P}_1 l_4 - {\cal P}_2 l_5)
		\: e_1 \vee e_0 \vee e_2.$
	\item~~$l' \vee m' = (l_0 m_5 + l_1 m_4 + l_2 m_3
		 + l_3 m_2 + l_4 m_1 + l_5 m_0)\:e_0\vee e_1 \vee e_2 \vee e_3.$
	\item~~$P \vee {\cal Q} = (P_0 {\cal Q}_0 + P_1 {\cal Q}_1 + P_2
		{\cal Q}_2 + P_3 {\cal Q}_3)\: e_0 \vee e_1 \vee e_2 \vee e_3.$
	\item~~$P \vee (P \vee l') = 0.$
	\item~~${\cal Q} \vee ({\cal Q} \wedge l') = 0.$
	\item 0. $(P \vee l') \wedge l' = - (l_0 l_5+ l_1 l_4 + l_2 l_3) P.$\\
	      1. $(P \vee l) \wedge l = 0.$
	\item 0. $({\cal Q} \wedge l') \vee l' = - (l_0 l_5 + l_1 l_4 + l_2 l_3)
		{\cal Q}.$\\
	      1. $   ({\cal Q} \wedge l) \vee l = 0.$
	\enume

	The proof is straightforward or follows from duality.\\
	The condition \ref{sec-ngras}.2 that a sextuple be a line is precisely
	chosen to insure 8.1 and 9.1.

	\sssec{Example.}\label{sec-egras}
	For $p = 7,$ given
	$P_0 := (1, 2, 3, 4),$ $P_1 := (1,0,1,1),$ $P_2 := (1,1,0,1),$
	$P_3 := (1,0,0,1),$\\
	$l_0 := [1, 1, 5, 6, 1, 4],$ $l_1 := [1, 6, 0, 6, 1, 1],$
	${\cal Q}_0 := \{1, 3, 3, 3\},$\\
	${\cal Q}_1 := \{1, 5, 0, 6\},$ we can easily verify\\
	$P_0$ and $P_1$ are incident to $l_0,$ $P_1$ and $P_2$ are incident to
	$l_1,$\\
	$P_0,$ $P_1,$ $P_2,$ $l_0$ and $l_1$ are incident to ${\cal Q}_0,$
	$P_3$ and $l_0$ are incident to ${\cal Q}_1.$

	\sssec{Notation.}
	As for 2 dimensional finite projective geometry, we will make use of a
	compact notation, assuming that the elements are ordered as if the 4 or
	6 normalized coordinates were forming an integer in base $p.$
	We have the correspondence\\
	$\begin{array}{llcl}
	\hth&	(0) := (0,0,0,1),&\vline	&[0] := [0,0,0,0,0,1],\\
	\hth&	(1) := (0,0,1,0),&\vline	&[1] := [0,0,0,0,1,0],\\
	\hth&      (p+1) := (0,1,0,0),&\vline &[p+1] := [0,0,0,1,0,0],\\
	\hth&   (p^2+p+1) := (1,0,0,0),&\vline &[p^2+p+1] := [0,0,1,0,0,0],\\
	\hth&		&\vline&[p^3+p^2+p+1] := [0,1,0,0,0,0],\\
	\hth&		&\vline&[p^4+p^3+p^2+p+1] := [1,0,0,0,0,0].
	\end{array}$

	\sssec{Example.}
	Continuing Example \ref{sec-egras},\\
	$P_0 = (180),$ $P_1 = (65),$ $P_2 = (107),$ $P_3 = (58),$
	$l_0 = [7222],$
	$l_1 = [17509],$ ${\cal Q}_0 = \{228\},$ ${\cal Q}_1 = \{308\}.$

	\sssec{Theorem.}\label{sec-tgras4}
	{\em $P$ and $Q$ are distinct  iff  $P \vee Q\neq  0.$}

	Proof: By Corollary \ref{sec-cogras}, if $P$ and $Q$ are not distinct,
	$Q = k P,$ $k\neq  0$
	and $P \vee Q = P \vee k P = 0.$  If $P \vee Q = 0,$ let $P_0$ be a
	coefficient of $P$ different from 0, \ref{sec-tgras3}.0 gives
	$P_0 Q_1 = P_1 Q_0,$ $P_0 Q_2 = P_2 Q_0,$ $P_0 Q_3 = P_3 Q_0, $
	therefore if $Q_0 = 0$ then $Q_1 = Q_2 = Q_3 = 0,$ and $Q$ is not
	a point.\\
	If $Q_0\neq 0,$ I can, by homogeneity choose $Q_0 = P_0$ and then
	$Q_1 = P_1,$ $Q_2 = P_2$ and $Q_3 = P_3$ or $Q = P.$

	\sssec{Theorem.}
	{\em Given 2 distinct points $P$ and $Q,$ there exist one and only one
	line $l = P \vee Q$ incident to $P$ and $Q.$}

	Proof: Because of associativity, $P \vee (P \vee Q) = 0$ and
	$(P \vee Q) \vee Q = 0,$
	therefore $l = P \vee Q$ is incident to both $P$ and Q and $l\neq 0$
	because $P$ and $Q$ are distinct.\\
	The line is unique.  Let $P \vee Q\neq 0,$ $P \vee l = Q \vee l = 0$,
	$l\neq 0.$\\
	Because $P \vee Q\neq 0,$ one of the coordinates is different from 0,
	let it be\\
	$P_0 Q_1 - P_1 Q_0.$  Theorem \ref{sec-tgras3}.1 gives 4
	equations associated to $P \vee l = 0$ and 4 equations associated to
	$Q \vee l = 0,$ the last equations are\\
	\hth$	-P_0 l_3 + P_1 l_1 - P_2 l_0 = 0$\\
	\hth$	-Q_0 l_3 + Q_1 l_1 - Q_2 l_0 = 0.$\\
	Multiplying the first by $-Q_0$ and the second by $P_0$ gives
	\begin{enumerate}
	\setcounter{enumi}{-1}
	\item$       (P_0 Q_1 - P_1 Q_0) l_1 = (P_0 Q_2 - P_2 Q_0) l_0,$

	Similarly multiplying by $-Q_1$ and $P_1$ gives
	\item$       (P_0 Q_1 - P_1 Q_0) l_3 = (P_1 Q_2 - P_2 Q_1) l_0,$

	The third equation of each set gives similarly
	\item$       (P_0 Q_1 - P_1 Q_0) l_2 = (P_0 Q_3 - P_3 Q_0) l_0,$
	\item$       (P_0 Q_1 - P_1 Q_0) l_4 = (P_3 Q_1 - P_1 Q_3) l_0,$

	If we add the first equations for $P \vee l = 0$ multiplied by $-Q_0$
	and $-Q_1,$ we get
	\item$       (P_0 Q_1 - P_1 Q_0) l_5 = -Q_1 P_3 l_1 + Q_1 P_2 l_2
				+ Q_0 P_3 l_3 + Q_0 P_2 l_4 = 0.$\\
	Because $l_0\neq 0,$ the first parenthesis is different from 0,
	otherwise,
	it follows from 0, 1, 2 and 3 then $l_1 = l_3 = l_2 = l_4 = 0,$
	and from 4 that $l_5 = 0.$  We can, because of homogeneity write\\
	\hth$	l_0 = P_0 Q_1 - P_1 Q_0,$\\
	it follows that\\
	\hth$	l_1 = P_0 Q_2 - P_2 Q_0.$\\
	\hth$	l_3 = P_1 Q_2 - P_2 Q_1.$\\
	and from 2 and 3,\\
	\hth$	l_2 = P_0 Q_3 - P_3 Q_0,$\\
	\hth$	l_4 = P_3 Q_1 - P_1 Q_3,$\\
	Replacing in 4, gives\\
	\hth$(P_0 Q_1 - P_1 Q_0) l_5 = Q_1 P_2(P_0 Q_3 - P_3) Q_0)
				- Q_1 P_3(P_0 Q_2 - P_2 Q_0)$\\
	\hti{16}$			+ Q_0 P_2(P_3 Q_1 - P_1 Q_3)
				+ Q_0 P_3(P_1 Q_2 - P_2 Q_1)$\\
	\hti{12}$		= (P_0 Q_1 - P_1 Q_0)(P_2 Q_3 - P_3 Q_2).$\\
	hence\\
	\hth$	l_5 = P_2 Q_3 - P_3 Q_2.$\\
	Therefore $l = P \vee Q.$
	\end{enumerate}

	\sssec{Theorem.}\label{sec-tgraspl}
	{\em Given a point $P$ and a line $l,$ not incident to $P,$ there exists
	one and only one plane ${\cal Q} = P \vee l$ incident to $P$ and $l.$}

	Proof:
	Because of associativity, \ref{sec-cogras} and \ref{sec-tgras3}.4.1,
	$P \vee (P \vee l) = (P \vee l) \wedge l = 0,$ therefore
	${\cal Q} = P \vee l$ is incident to both $P$ and $l$ and
	${\cal Q}\neq 0$ because $P$ and $l$ are not incident.\\
	The plane is unique, if $P \vee {\cal Q} = l \vee {\cal Q} = 0$ and
	${\cal Q}\neq 0,$ let ${\cal Q}_0$ be a coefficient of
	${\cal Q}\neq 0.$ $P \vee {\cal Q} = 0$ and $l \vee {\cal Q} = 0$ give\\
	\hth$P_0 {\cal Q}_0 + P_1 {\cal Q}_1 + P_2 {\cal Q}_2
		+ P_3 {\cal Q}_3 = 0,$\\
	\hth$	{\cal Q}_1 l_0 + {\cal Q}_2 l_1 + {\cal Q}_3 l_2 = 0,$\\
	\hth$	-{\cal Q}_0 l_0 + {\cal Q}_2 l_3 - {\cal Q}_3 l_4 = 0,$\\
	\hth$	-{\cal Q}_0 l_1 - {\cal Q}_1 l_3 + {\cal Q}_3 l_5 = 0,$\\
	\hth$	-{\cal Q}_0 l_2 + {\cal Q}_1 l_4 - {\cal Q}_2 l_5 = 0.$\\
	Multiplying the equations respectively by $l_5,$ 0, 0, $-P_3$ and $P_2$
	and adding gives using homogeneity, and the same argument used in the
	preceding Theorem,\\
	\hth$	{\cal Q}_0 = P_1 l_5 + P_2 l_4 + P_3 l_3,$\\
	\hth$	{\cal Q}_1 = -P_0 l_5 - P_3 l_1 + P_2 l_2.$\\
	Similarly, if we multiply respectively by $l_4$, $0,$ $P_3,$ $0$
	and $-P_1$ and then add,\\
	\hth$	{\cal Q}_2 = - P_0 l_4 + P_3 l_0 - P_1 l_2,$\\
	and if we multiply respectively by $l_3,$ 0, $-P_2,$ $P_1$ and 0 and
	then add,\\
	\hth$	{\cal Q}_3 = - P_0 l_3 - P_2 l_0 + P_1 l_2.$\\
	Therefore ${\cal Q} = P \vee l.$

	Using duality, it is easy to deduce from \ref{sec-tgras5} and
	\ref{sec-tgras6}.

	\sssec{Theorem.}\label{sec-tgras5}
	{\em Given 2 distinct planes ${\cal P}$ and ${\cal Q},$ there exist one
	and only one line $l = {\cal P} \wedge {\cal Q}$ incident to
	${\cal P}$ and ${\cal Q}.$}

	\sssec{Theorem.}\label{sec-tgras6}
	{\em Given a plane ${\cal Q}$ and a line $l,$ not incident to
	${\cal Q},$ there exists one and only one point
	$P = {\cal Q} \wedge l$ incident to ${\cal Q}$ and $l.$}

	\sssec{Lemma.}\label{sec-llcrossm}
	{\em If $l$ and $m$ are lines, then}\\
	\hth$(l_1m_5-l_5m_1)(l_0m_5+l_1m_4+l_2m_3)+(l_1m_2-l_2m_1)
		(l_3m_5-l_5m_3)\\
	\hti{12}= l_1m_5(l_0m_5+l_1m_4+l_2m_3+l_3m_2+l_4m_1+l_5m_0).$

	\sssec{Lemma.}
	{\em If $l = [0,0,0,l_3,l_4,l_5]$ and $m = [m_0,m_1,m_2,0,0,0],$ then
	$l$ and $m$ have a point $P$ in common iff $l\vee m = 0$ and}
	$P = (0,m_0,m_1,m_2).$

	Proof:
	The 4 conditions associated with $P\vee l = 0$ give, because not all
	$l_3,$ $l_4$ and $l_5$ can be 0, $P_0 = 0$  and
	$P_1l_5+P_2l_4+P_3l_3 = 0.$  The 4 conditions associated with
	$P\vee m = 0$ give, because not all $m_0,$ $m_1$ and $m_2$ can be 0,
	$P_1 = m_0,$ $P_2 = m_1,$ $P_3 = m_2$, substituting in the remaining
	equation gives the equivalence with $l\vee m = 0$.

	\sssec{Lemma.}
	{\em If $l_5m_5\neq 0$ and $l_1m_5\neq l_5m_1$, if $P$ is on $l$ and
	$m$, then $l$ and $m$ have a point $P$ in common iff $l\vee m = 0$ and}
\\
	\hth $P = (l_0m_5+l_1m_4+l_2m_3,l_3m_4-l_4m_3,l_5m_3-l_3m_5,
		l_4m_5-l_5m_4).$

	Proof:
	The first component of $P\vee l = 0$ and of $P\vee m = 0$ implies
	$P_1 = k_0(l_4m_3-l_3m_4)$, $P_2 = k_0(l_3m_5-l_5m_3)$ and
	$P_3 = k_0(l_5m_4-l_4m_5)$.  Similarly, the second components implies
	$P_0 = k_1(l_1m_2-l_2m_1)$, $P_2 = k_1(l_1m_5-l_5m_1)$ and
	$P_3 = k_1(l_2m_5-l_5m_2)$.  Consistency implies\\
	\hth$(l_1m_5-l_5m_1)(l_5m_4-l_4m_5) = (l_3m_5-l_5m_3)(l_2m_5-l_5m_2).$\\
	or\\
	\hth$l_5m_5(l_3m_2+l_2m_3+l_1m_4+l_4m_1)+m_5^2(-l_1l_4-l_2l_3)\\
	\hti{16}+l_5^2(-m_1m_4-m_2m_3)\\
	\hti{12}= l_5m_5(l_3m_2+l_2m_3+l_1m_4+l_4m_1+l_0m_5+l_5m_0),$\\
	because $l$ and $m$ are lines.  Chosing
	$k_1 = (l_3m_5-l_5m_3)/(l_1m_5-l_5m_1)$ and $k_0 = -1$
	gives the expression for $P$ using Lemma \ref{sec-llcrossm}.

	\sssec{Theorem.}\label{sec-tgras7}
	{\em If 2 distinct lines $l$ and $m$ are such that $l \vee m = 0,$ they
	are incident to a point noted $l \bigx  m$ and to a plane noted
	$l\: {\cal X}\: m.$
	Vice-versa, if 2 distinct lines are incident to the same point or the
	same plane then $l \vee m = 0.$
	Moreover, if $l = [l_0,l_1,l_2,l_3,l_4,l_5]$ and $m = [m_0,m_1,m_2,m_3,
	m_4,m_5],$ then one of the following will give the point $l \bigx m$}
	\enumb
	\item$(l_0m_5+l_1m_4+l_2m_3,l_3m_4-l_4m_3,l_5m_3-l_3m_5,l_4m_5-l_5m_4).$
	\item$(l_1m_2-l_2m_1,l_4m_1+l_3m_2+l_0m_5,l_1m_5-l_5m_1,l_2m_5-l_5m_2),$
	\item$(l_2m_0-l_0m_2,l_0m_4-l_4m_0,l_5m_0+l_3m_2+l_1m_4,l_2m_4-l_4m_2).$
	\item$(l_0m_1-l_1m_0,l_0m_3-l_3m_0,l_1m_3-l_3m_1,l_5m_0+l_4m_1+l_2m_3).$
	\enume
	{\em and the plane $l\: {\cal X}\: m$}
	\begin{enumerate}
	\setcounter{enumi}{3}
	\item$\{l_5m_0+l_4m_1+l_3m_2,l_2m_1-l_1m_2,l_0m_2-l_2m_0,l_1m_0-l_0m_1
\}.$
	\item$\{l_4m_3-l_3m_4,l_5m_0+l_2m_3+l_1m_4,l_4m_0-l_0m_4,l_3m_0-l_0m_3
\}.$
	\item$\{l_3m_5-l_5m_3,l_5m_1-l_1m_5,l_4m_1+l_2m_3+l_0m_5,l_3m_1-l_1m_3
\}.$
	\item$\{l_5m_4-l_4m_5,l_5m_2-l_2m_5,l_4m_2-l_2m_4,l_3m_2+l_1m_4+l_0m_5
\}.$
	\enume

	The proof follows by a judicious application of the Lemmas.

	\sssec{Exercise.}
	If $l$ and $m$ are lines and $l \vee m = 0$
	then $l \vee (l \bigx m) = (l \bigx m) \vee m = 0.$

	\sssec{Exercise.}
	If $l \bigx m$ and $l\: {\cal\: X} m$ are defined by
	\ref{sec-tgras7}.0 and .4, then $(l \bigx m) \vee (l\:{\cal X}\:m) = 0.$



	\section{Affine Geometry in 3 Dimensions.}

	\setcounter{subsection}{-1}
	\ssec{Introduction.}
	To define a 3 dimensional Euclidean geometry, I will start with a
	preferred plane ${\cal I}$ to which are	associated
	the notions of affine geometry. Just as in the case of the Pappian 
	plane, we can define the notions of parallelism, mid-point,
	equality of segments (ordered pair of points) on parallel lines.
	It is convenient to intoduce a matrix notation to express parallelism
	and in later sections polarity and orthogonality.
	Bold faced letters will be used for matrices.  The coordinates
	of points, lines and planes will have associated with them vectors
	which will be considered as row vectors.

	\ssec{The ideal plane and parallelism.}
	\sssec{Definition.}\label{sec-didealpl}
	The preferred plane is called the {\em ideal plane}.
	There is no restriction in chosing the ideal plane
	${\cal I} = \{1,1,1,1\}$ because the coordinates of ${\cal I}$ simply
	corresponds to those of the unit point (1,1,1,1) and can be considered 
	as the polar of the unit point with respect to the tetrahedron
	of the coordinate system (1,0,0,0), (0,1,0,0), (0,0,1,0), (0,0,0,1).

	\sssec{Definition.}
	The points in the ideal plane are called the {\em ideal points},
	the lines in the ideal plane are called the {\em ideal lines}.

	\sssec{Definition.}
	Two {\em lines are parallel}  iff  they are incident to ${\cal I}$ at
	the same point.\\
	Two {\em planes are parallel} iff they are incident to
	${\cal I}$ on the same line.\\
	A {\em plane} ${\cal Q}$ {\em and a line} $l$ {\em are parallel} iff
	the line ${\cal Q} \wedge {\cal I}$ is incident to the point
	${\cal I} \wedge l.$

	\sssec{Definition.}
	\hth${\bf L} := \matb{(}{rrrrrr}1&1&1&0&0&0\\-1&0&0&1&-1&0\\
		0&-1&0&-1&0&1\\0&0&-1&0&1&-1\mate{]},$
	${\bf P} := \matb{[}{rrrr}0&0&1&-1\\0&-1&0&1\\
		0&1&-1&0\\1&0&0&-1\\1&0&-1&0\\1&-1&0&0\mate{\}}.$
 
	\sssec{Theorem.}
	{\em If $l$ is a line then {\bf L}$l^T$ is the direction of the line l
.\\
	If ${\cal P}$ is a plane then {\bf P}${\cal P}^T$ is the direction of
	the plane ${\cal P}$.}

	The proof follows from \ref{sec-tgras3}.3 and 1.
	Notice that {\bf P} is obtained from the transpose of {\bf L} by
	exchanging row $i$ with row $5-i$, $i$ = 0,1,2.

	\sssec{Theorem.}\label{sec-tpolarpar}
	{\em Let ${\cal Q} := \{ {\cal Q}_0,{\cal Q}_1,{\cal Q}_2,{\cal Q}_3 \}$
	and $l := [l_0,l_1,l_2,l_3,l_4,l_5],$\\
	The plane ${\cal Q}$ is parallel to the line l iff}
	\enumb
	\item${\cal Q}{\bf L}l^T = 0,
	or\\
	\hti{4}- Q_0(l_0+l_1+l_2) + Q_1(l_0-l_3+l_4) + Q_2(l_1+l_3-l_5)
	+ Q_3(l_2-l_4+l_5) = 0.$
	\enume

	The proof follows from the property that the direction ${\bf L}l^T$
	of the line $l$ is incident to the plane ${\cal Q}$.  Alternately,
	we can obtain the Theorem by introducing first, directional
	correspondance.

	\sssec{Definition.}
	Let $P := (P_0,P_1,P_2,-(P_0+P_1+P_2))$ be an ideal point.  Let\\
	$m := [m_0,m_1,m_2,m_3,m_4,m_5]$ be an ideal line.  These point and line
	can also be determined by 3 well chosen coordinates.  The coordinates of
	points are placed between double parenthesis and that of lines, between
	double brackets, while the point $P$ as viewed as a point in the plane
	is denoted $(P)$ and the line as $[m]$.\\
	One of the good choices is $[[m_3,-m_1,m_0]],$ indeed,
	in the ideal plane, $m_0 \:e_0 + m_1 \:e_1 + m_2 \:e_2$ is the dual of
	$m_0 \:e_1 \vee e_2 + m_1 \:e_2 \vee e_0 + m_2 \:e_0 \vee e_1,$ while
	the 3 chosen components give $ m_3 \:e_1 \vee e_2 - m_1 \:e_0 \vee
	e_2 + m_0 \:e_0 \vee e_1.$ The other components are
	$m_2 = -m_0-m_1,$ $m_4 = m_3-m_0,$ $m_5 = m_1+m_3.$\\
	We have $P^T = {\bf U}\:(P)^T$ and $[m]^T = {\bf V}\:m^T,$ with\\
	\hth${\bf U} = \matb{(}{rrr}1&0&0\\ 0&1&0\\ 0&0&1\\ -1&-1&-1\mate{)}$
	and ${\bf V} = \matb{(}{rrrrrr}0&0&0&1&0&0\\0&-1&0&0&0&0\\1&0&0&0&0&0
	\mate{)}.$

	The correspondence which associates to $P = (P_0,P_1,P_2,P_3)$ in the 3
	dimensional space the point $(P) = (P_0,P_1,P_2)$ in the 2 dimensional
	plane ${\cal I},$ and which associates to the line
	$m = [m_0,m_1,m_2,m_3,m_4,m_5],$ in the 3 dimensional space the line
	$[m] = [[m_3,-m_1,m_0]],$ in the 2 dimensional plane ${\cal I},$
	is called the {\em directional correspondence}.

	\sssec{Theorem.}
	{\em The directional correspondence is a homomorphism from the 3
	dimensional space onto the ideal plane.}

	\sssec{Theorem.}\label{sec-tvp}
	\hth${\bf V\:P} = {\bf U}^T.$

	\sssec{Theorem.}
	{\em If $P$ and $m$ are in the ideal plane, $P$ is on $m,$ iff $(P)$ is
	on $[m],$ iff}\\
	\hth$(P) \cdot (m) = ((P_0,P_1,P_2)) \cdot [[m_3,-m_1,m_0]]
			  = P_0 m_3 - P_1 m_1 + P_2 m_0 = 0.$

	\sssec{Alternate Proof of}
	\vspace{-18pt}\hspace{150pt}\ref{sec-tpolarpar}.\\[8pt]
	For instance, to the point  $I_l$  and line  $i_{\cal Q}$  of Theorem
	\ref{sec-tpolarpar},
	correspond, the point $(I_l)$ and the line $[i_{\cal Q}]$ in ${\cal I}$.
\\
	\hth$	(I_l) = ((l_0+l_1+l_2,-l_0+l_3-l_4,-l_1-l_3+l_5)),$\\
	\hth$	[i_{\cal Q}] = [[{\cal Q}_0-{\cal Q}_3,{\cal Q}_1-{\cal Q}_3,
		{\cal Q}_2-{\cal Q}_3]],$\\
	${\cal Q}$ is parallel to $l$ iff\\
	$-(l_0+l_1+l_2)({\cal Q}_0-{\cal Q}_3)
		+ (-l_0+l_3-l_4)({\cal Q}_3-{\cal Q}_1)
		- (-l_1-l_3+l_5)({\cal Q}_2-{\cal Q}_3) = 0$\\
	which is \ref{sec-tpolarpar}.0.

	\sssec{Theorem.}
	{\em The mid-point of two points $A = (a_0,a_1,a_2,a_3)$ and
	$B = (b_0,b_1,b_2,b_3)$ is}\\
	\hth$(b_0+b_1+b_2+b_3 )A + (a_0+a_1+a_2+a_3 )B.$

	\sssec{Notation.}
	The {\em mid-point of $A$ and $B$} is denoted by $A + B$.

	\sssec{Exercise.}
	Generalize the construction of the polar $p$ of a point $P$ with respect
	to a triangle to that of the polar {\cal P} of a point $P$ with respect
	to a tetrahedron and prove that if $P = (p_0,p_1,p_2,p_3)$ then
	${\cal P} = \{{\cal P}_1{\cal P}_2{\cal P}_3,{\cal P}_2{\cal P}_3
	{\cal P}_0,{\cal P}_3{\cal P}_0{\cal P}_1,{\cal P}_0{\cal P}_1{\cal P}_2
	\}.$

	\section{Polar Geometry in 3 Dimensions.}\label{sec-Spolar}
\setcounter{subsection}{-1}
	\ssec{Introduction.}
	To define a polar geometry in 3 dimensions, I will start with an
	affine Geometry in 3 dimansions and a preferred non degenerate quadric
	$\theta$ which is not tangent to the ideal plane ${\cal I}$.
	Using the ideal plane and the prefered quadric we can define
	orthogonality, spheres, centers, equality of pairs of points and lines,
	$\ldots$  .

	The preferred quadric is represented by a symmetric 4 by 4 matrix
	{\bf F}, which associates to a 4-vector representing a point (pole), a
	4-vector representing a plane (polar). Its adjoint {\bf G} gives the
	correspondance from polar to pole. From {\bf F}, we can derive a 6 by 6
	matrix {\bf H} which gives the correspondance between a line and its
	polar.  From {\bf F} we can also derive a polarity in the ideal plane
	represented by a 3 by 3 matrix ${\bf J}_3,$ giving the correspondance
	from pole to polar and its adjoint ${\bf K}_3,$ giving the
	correspondance from polar to pole from which we can derive
	perpendicularity between a line and a plane, the direction of the line
	giving the pole and the direction of the plane giving the polar.
	The 6 by 6 matrix {\bf J}, derived from ${\bf J}_3$, allows for a
	direct check of the orthogonality of 2 lines and the 4 by 4 matrix
	{\bf K}, derived from ${\bf K}_3$, allows for a direct check of the
	orthogonality of 2 planes.

	\ssec{The fundamental quadric, poles and polars.}

\setcounter{subsubsection}{-1}
	\sssec{Introduction.}
	The properties of pole and polar are properties in Pappian Geometry.
	They are easily generalized by using a 4 dimension collinearity which 
	transforms the fundamental quadric into an arbitrary quadric or by
	chosing a coordinate system with the four base points on the quadric.

	\sssec{Definition.}\label{sec-dfundqua}
	The preferred quadric is called the {\em fundamental quadric}.
	There is no restriction in chosing the fundamental quadric $\theta$ as
	follows, because it
	simply assumes that the quadric passes through the base points
	(1,0,0,0), (0,1,0,0), (0,0,1,0), (0,0,0,1).
	Let
	\enumb
	\item$ {\bf F} := \matb{(}{cccc}
		0&n_0&n_1&n_2\\n_0&0&n_3&n_4\\
		n_1&n_3&0&n_5\\n_2&n_4&n_5&0\mate{)},$
	\item$ \Theta  := (X0, X1, X2, X3){\bf F}(X0, X1, X2, X3)^T\\
	\hti{2}= n_0 X0 X1 + n_1 X0 X2 + n_2 X0 X3
			+ n_3 X1 X2 + n_4 X3 X1 + n_5 X2 X3 = 0.$\\[10pt]
	The condition of non degeneracy and non tangency are
	\item$d := det({\bf F}) =\\
	\hth n_0^2 n_5^2 + n_1^2 n_4^2 + n_2^2 n_3^2
	- 2 ( n_0 n_1 n_4 n_5 + n_1 n_2 n_3 n_4 + n_2 n_0 n_3 n_5 )\neq 0,$
	\item$t := n_3n_4n_5+n_1n_2n_5+n_0n_2n_4+n_0n_1n_3
		+2n_0n_5(n_0+n_5)\\
	\hth+2n_1n_4(n_1+n_4)+2n_2n_3(n_2+n_3)-(n_0n_5+n_1n_4+n_2n_3)n\neq 0,$\\
		where\\
	\hth $n := n_0 + n_1 + n_2 + n_3 + n_4 + n_5.$
	\enume

	The condition $t\neq 0$ will be verified in \ref{sec-ttangency}.

	\sssec{Definition.}
	In polar geometry, the points in the ideal plane which are not on the
	quadric are called the {\em ideal points}, all lines in the ideal plane
	which are not tangent to the quadric are the {\em ideal lines}.

	\sssec{Definition.}
	The points in the ideal plane and the quadric are the
	{\em isotropic points}.  The lines in the ideal plane tangent to the
	quadric are the {\em isotropic lines}.\\
	If the isotropic points are real, the polar geometry is said to be
	{\em hyperbolic}, if there are no real isotropic points, it is said to
	be {\em elliptic}, if there is exactly one isotropic
	point, the quadric being tangent to the plane, the geometry is said to
	be {\em parabolic}.

	Again as for the involutive geometry, I will not study the parabolic
	case and will study together the elliptic and hyperbolic case.

	\sssec{Definition.}
	The {\em polar} of the point $P = (P_0,P_1,P_2,P_3),$
	is the plane\\
	\hth$	{\cal Q}^T :={\bf F}P^T =
		\{ n_0 P_1 + n_1 P_2 + n_2 P_3, n_0 P_0 + n_3 P_2 + n_4 P_3,\\
	\hti{12} n_1 P_0 + n_3 P_1 + n_5 P_3,n_2 P_0 + n_4 P_1 + n_5 P_2 \}.$\\
	$P$ is called the {\em pole of the plane} ${\cal Q}.$

	\sssec{Theorem.}
	{\em
	\hth$	P^T ={\bf G}{\cal Q}^T$\\
	where {\bf G} is the adjoint of {\bf F}:\\
	{\bf G} =\\
	$\matb{(}{cccccc}2n_3n_4n_5 & n_5(n_0n_5\!-\!n_2n_3\!-\!n_1n_4) &
	 n_4(n_1n_4\!-\!n_0n_5\!-\!n_2n_3) & n_3(n_2n_3\!-\!n_1n_4\!-\!n_0n_5)\\
	n_5(n_0n_5\!-\!n_2n_3\!-\!n_1n_4) &
	 2n_1n_2n_5 & n_2(n_2n_3\!-\!n_1n_4\!-\!n_0n_5) &
	 n_1(n_1n_4\!-\!n_0n_5\!-\!n_2n_3)\\
	n_4(n_1n_4\!-\!n_0n_5\!-\!n_2n_3) &
	 n_2(n_2n_3\!-\!n_1n_4\!-\!n_0n_5) &
	 2n_0n_2n_4 & n_0(n_0n_5\!-\!n_2n_3\!-\!n_1n_4)\\
	n_3(n_2n_3\!-\!n_1n_4\!-\!n_0n_5) &
	 n_1(n_1n_4\!-\!n_0n_5\!-\!n_2n_3) &
	 n_0(n_0n_5\!-\!n_2n_3\!-\!n_1n_4) & 2n_0n_1n_3\mate{)}.$ }

	\sssec{Theorem.}\label{sec-ttangency}
	{\em ${\cal I}$ is tangent to the fundamental quadric iff\\
	\hth$(1,1,1,1){\bf G}\{1,1,1,1\}^T = 0$\\
	or\\
	\hth$t = 0.$\\
	where $t$ is defined in \ref{sec-dfundqua}}.

	$2t$ is simply the sum of all the elements of {\bf G}.

	\sssec{Theorem.}
	{\em If $Q$ is on the polar ${\cal P}$ of $P$ then its polar
	${\cal Q}$ is incident to $P.$}

	The proof is left as an exercise.  The theorem justifies the following
	definition:

	\sssec{Definition.}
	A line $m$ is {\em a polar of a line} $l$ iff it is the
	line common to all the polars of the points of $l.$

	\sssec{Theorem.}
	{\em Let {\bf H} :=\\
	$\matb{(}{cccccc}
	n_1n_4-n_2n_3 & n_5n_1 & -n_5n_2 & n_5n_3 & n_5n_4 & -n_5^2\\
	-n_4n_0 & n_2n_3-n_0n_5 & n_4n_2 & n_4n_3 & -n_4^2 & n_4n_5\\
	n_3n_0 & -n_3n_1 & n_0n_5-n_1n_4 & -n_3^2 & n_3n_4 & n_3n_5\\
	-n_2n_0 & -n_2n_1 & -n_2^2 & n_0n_5-n_1n_4 & n_2n_4 & -n_2n_5\\
	-n_1n_0 & -n_1^2 & -n_1n_2 & -n_1n_3 & n_2n_3-n_0n_5 & n_1n_5\\
	-n_0^2 & -n_0n_1 & -n_0n_2 & n_0n_3 & -n_0n_4 & n_1n_4-n_2n_3
	\mate{)},$\\
	then the polar $m$ of $l$ is given by\\
	\hth$	m^T = {\bf H} l^T.$\\
	Moreover\\
	\hth$	{\bf H}\: {\bf H} = d {\bf I},$\\
	where ${\bf I}$ is the identity matrix and $d$ is the determinant of
	${\bf F}$ given in \ref{sec-dfundqua}.3.}

	The proof is left as an exercise.  As a hint, consider 2 points
	$P$ and $Q$ on $l = P \vee Q$ and their polar  ${\bf F} P$  and
	${\bf  F} Q$.

	\sssec{Definition.}
	The {\em center of a quadric} is the pole of ${\cal I},$

	\sssec{Example.}
	The pole of \{1,0,0,0\} is$\\
	(  2n_3n_4n_5, n_5(n_0n_5-n_2n_3-n_1n_4), n_4(n_1n_4-n_0n_5-n_2n_3),\\
	\hth n_3(n_2n_3-n_1n_4-n_0n_5)).$\\
	The center of the fundamental quadric is$\\
	(2n_3n_4n_5+n_0n_5(n_5-n_3-n_4)+n_1n_4(n_4-n_5-n_3)
		+n_2n_3(n_3-n_4-n_5),\\	
	 \hti{1}2n_5n_1n_2+n_3n_2(n_2-n_5-n_1)+n_4n_1(n_1-n_2-n_5)
		+n_0n_5(n_5-n_1-n_2),\\	
	 \hti{1}2n_2n_4n_0+n_5n_0(n_0-n_2-n_4)+n_1n_4(n_4-n_0-n_2)
		+n_3n_2(n_2-n_4-n_0),\\	
	 \hti{1}2n_0n_1n_3+n_2n_3(n_3-n_0-n_1)+n_4n_1(n_1-n_3-n_0)
		+n_5n_0(n_0-n_1-n_3).$

	\ssec{Orthogonality in space and the ideal polarity.}
\setcounter{subsubsection}{-1}
	\sssec{Introduction.}
	Prefering both an ideal plane and a fundamental quadric allows us
	to define othogonality of lines and planes with lines and planes.
	After defining the polarity in the ideal plane, induced by the
	fundamental quadric, we use it to derive the 3 conditions which
	express the orthogonality of lines and planes and the condition
	which express the orthogonality of 2 lines or of 2 planes.

	\sssec{Definition.}
	A {\em line is orthogonal to a plane} iff the
	polar of its ideal point is incident to the ideal line of the plane.
	A {\em line is orthogonal to a line} iff the polar of its ideal
	point is incident to the ideal point of the other line.
	A {\em plane is orthogonal to a plane} iff the ideal line of one is
	the polar of the ideal line of the other.

	\sssec{Definition.}
	The {\em ideal polarity} is the polarity induced in the ideal plane
	${\cal I}$ by the polarity defined by the quadric $\theta$.

	\sssec{Notation.}
	It is sometimes convenient to use an other notation for the elements of
	the fundamental quadric.\\
	$n_{ij} = n_{ji}$ is the coefficient of $X_iX_j$ in the equation of the 
	fundamental quadric, more specifically,\\
	$n_{01} := n_0,$ $n_{02} := n_1,$ $n_{03} := n_2,$ $n_{12} := n_3,$
	$n_{31} := n_4$ and $n_{23} := n_5.$\\
	The elements of the matices ${\bf K}_3$ and of {\bf K} are more
	easily expressed in terms of $i_{ii}$ and $i_{ij},$ using
	$ijkl$ for permutation of 0123,\\
	\hth$i_{ii} = -(n_{kl}^2+n_{lj}^2+n_{jk}^2+2(n_{lj}n_{kl}+n_{kl}n_{jk}
		+n_{jk}n_{lj}),\\
	\hth i_{ij} = (n_{ik}-n_{il})(n_{jk}-n_{jl})+n_{kl}(3n_{ij}+2n_{kl}-n),$
\\
	with $n := n_{01}+n_{02}+n_{02}+n_{12}+n_{31}+n_{23}.$\\
	For instance,\\
	\hth$i_{00} = -(n_{12}^2+n_{31}^2+n_{23}^2+2(n_{31}n_{23}+n_{23}n_{12}
		+n_{12}n_{31},$\\
	\hth$i_{01} = (n_{02}-n_{03})(n_{12}-n_{13})+n_{23}(3n_{01}+2n_{23}-n),$

	\sssec{Theorem.}
	{\em The point to line ideal polarity in given by the matrix\\
	\hth${\bf J}_3 = {\bf U}^T\:{\bf F\:U}
	= \matb{(}{ccc}-2n_2&n_0-n_2-n_4&n_1-n_2-n_5\\
		n_0-n_2-n_4&-2n_4&n_3-n_4-n_5\\
		n_1-n_2-n_5&n_3-n_4-n_5&-2n_5\mate{)}.$\\[10pt]
	The adjoint matrix ${\bf K}_3$ gives the line to point ideal polarity.}
\\
	\hth${\bf K}_3 = \mat{i_{00}}{i_{01}}{i_{02}}{i_{01}}{i_{11}}{i_{12}}
	{i_{02}}{i_{12}}{i_{22}}.$


	Indeed given a point $(P)$ in ${\cal I}$, ${\bf U}(P)$ gives the
	coordinates of $P$ is space, multiplication to the left by {\bf F}
	determines the polar plane ${\cal P}$.  ${\bf P}{\cal P}$ gives
	the direction $i_{\cal P}$ of ${\cal P}$, muiltiplication to the
	left by {\bf V} gives the 3 coordinates $(i_{\cal P})$ of the direction
	in the ideal plane, using \ref{sec-tvp} gives ${\bf J}_3$.

	For the adjoint matrix,\\
	$i_{00} = 4n_4n_5-(n_3-n_4-n_5)^2 = -(n_3^2+n_4^2+n_5^2)+
		2(n_4n_5+n_5n_3+n_3n_4\\
	\hth= -(n_{12}^2+n_{31}^2+n_{23}^2+2(n_{31}n_{23}+n_{23}n_{12}
		+n_{12}n_{31}.\\
	i_{01} = (n_3-n_4-n_5)(n_1-n_2-n_5)+2n_5(n_0-n_2-n_4)\\
	\hth= (n_{12}-n_{13}-n_{23})(n_{02}-n_{03}-n_{23})
		+2n_{23}(n_{01}-n_{03}-n_{13})\\
	\hth= (n_{02}-n_{03})(n_{12}-n_{13})+n_{23}(3n_{01}+2n_{23}-n).$

	\sssec{Theorem.}
	\enumb
	\item $det({\bf J}_3) = 2t.$
	\item {\em The ideal polarity is not degenerate.}
	\enume

	\sssec{Theorem.}
	{\em If $P$ is on the conic associated with the idea polarity then
	$P$ is on the fundamental quadric.}

	\sssec{Theorem.}\label{sec-tspaceortho}
	{\em Let\\
	\hth$	l_a := l_0+l_1+l_2,$ $l_b := -l_0-l_4+l_3,$
	$l_c := -l_1-l_3+l_5,$ $l_d := -l_2-l_5+l_4,$\\
	\hth$	i_{\cal P} := {\cal P} \wedge {\cal I},
		i_{\cal Q} := {\cal Q} \wedge {\cal I},
		I_l := l \wedge {\cal I}, I_m := m \wedge {\cal I}.$
	\hth then} $[l_a,l_b,l_c,l_d] = I_l = {\bf L}l.$
	\enumb
	\item {\em A line $l = [l_0,l_1,l_2,l_3,l_4,l_5]$ is orthogonal to
	the plane\\
	\hth$	{\cal P} = \{P_0,P_1,P_2,P_3\}$ iff\\
	      0. $    [i_{\cal P}] = {\bf J}_3 (I_l),$
		or if for some $k \neq 0$,\\
	      1. $   k(P_3-P_0) = l_a n_2     + l_b(n_0-n_4) + l_c(n_1-n_5)
	 + l_d n_2,$\\
	\hti{4}$k(P_3-P_1) = l_a(n_2-n_0) + l_b(-n_4) + l_c(n_3-n_5) + l_d n_4,
$\\
	\hti{4}$k(P_3-P_2) = l_a(n_2-n_1) + l_b(n_3-n_4) + l_c(-n_5) + l_d n_5.
$\\
	Other relations can be derived from these, e.g.\\
	      2. $k(P_1-P_0) = l_a n_0 + l_b n_0 + l_c(n_1-n_3) + l_d(n_2-n_4),
	$\\
	\hti{4}$k(P_2-P_0) = l_a n_1 + l_b(n_0-n_3) + l_c n_1 + l_d(n_2-n_5),$}
	\item {\em A line $l = [l_0,l_1,l_2,l_3,l_4,l_5]$ is orthogonal to the
	line\\
	\hth$	m = [m_0,m_1,m_2,m_3,m_4,m_5]$ iff\\
	      0. $   (I_m) {\bf J}_3 (I_l) = 0,$\\
		or}\\
	      1. $- ( (l_0+l_1+l_2)(m_0+m_4-m_3)+(m_0+m_1+m_2)(l_0+l_4-l_3))
			n_0\\
	\hth - ((l_0+l_1+l_2)(m_1+m_3-m_5)+(m_0+m_1+m_2)(l_1+l_3-l_5))n_1\\
	\hth - ((l_0+l_1+l_2)(m_2+m_5-m_4)+(m_0+m_1+m_2)(l_2+l_5-l_4))n_2\\
	\hth + ((l_0+l_4-l_3)(m_1+m_3-m_5)+(m_0+m_4-m_3)(l_1+l_3-l_5))n_3\\
	\hth + ((l_0+l_4-l_3)(m_2+m_5-m_4)+(m_0+m_4-m_3)(l_2+l_5-l_4))n_4\\
	\hth +((l_1+l_3-l_5)(m_2+m_5-m_4)+(m_1+m_3-m_5)(l_2+l_5-l_4))n_5 = 0.$\\
	or\\
	$m{\bf J}l = 0,$ with ${\bf J} = $\\
\tiny{
	\hti{-5}$\matb{(}{cccccc}
	-2n_0&-n_0-n_1+n_3&-n_0-n_2+n_4&n_0-n_1+n_3&-n_0+n_2-n_4
		&n_1-n_2-n_3+n_4\\
	-n_0-n_1+n_3&-2n_1&-n_1-n_2+n_5&n_0-n_1-n_3&-n_0+n_2+n_3-n_5
		&n_1-n_2+n_5\\
	-n_0-n_2+n_4&-n_1-n_2+n_5&-2n_2&n_0-n_1-n_4+n_5&-n_0+n_2+n_4
		&n_1-n_2-n_5\\
	n_0-n_1+n_3&n_0-n_1-n_3&n_0-n_1-n_4+n_5&-2n_3&n_3+n_4-n_5&n_3-n_4+n_5\\
	-n_0+n_2-n_4&-n_0+n_2+n_3-n_5&-n_0+n_2+n_4&n_3+n_4-n_5&-2n_4
		&-n_3+n_4+n_5\\
	n_1-n_2-n_3+n_4&n_1-n_2+n_5&n_1-n_2-n_5&n_3-n_4+n_5&-n_3+n_4+n_5
		&-2n_5\mate{)}.$
}
	\enume
\begin{enumerate}
\setcounter{enumi}{1}
	\item {\em A plane ${\cal P} = \{P_0,P_1,P_2,P_3\}$ is orthogonal to the
plane\\
	\hth$	{\cal Q} = \{Q_0,Q_1,Q_2,Q_3\}$ iff}\\
	      0. $(i_{\cal Q}) {\bf K}_3 (i_{\cal P})^T = 0,$\\
		{\em or}\\
	      1. ${\cal Q} {\bf K}{\cal P}^T = 0,$\\
		{\em where {\bf K} is the symmetric matrix}\\
		${\bf K} = \matb{(}{cccc}
				i_{00}&i_{01}&i_{02}&i_{03}\\
				i_{10}&i_{11}&i_{12}&i_{13}\\
				i_{20}&i_{21}&i_{22}&i_{23}\\
				i_{30}&i_{31}&i_{32}&i_{33}\mate{)}.$
	\enume

	\sssec{Example.}
	The pole of the line ${\cal I} \wedge A_0 = [[1,0,0]],$ is
	$((i_{00},i_{01},i_{02})),$ which gives\\
	$I_{Q_0} = (i_{00},i_{01},i_{02},-(i_{00}+i_{01}+i_{02})),$
	hence if $foot_0 := (A_0 \vee I_{Q_0}) \wedge {\cal A}_0,$ then\\
	$foot_0 = (0,i_{01},i_{02},i_{03}),$ with\\
	$i_{03} = -(i_{00}+i_{01}+i_{02}) = (n_0-n_1-n_3)(n_4-n_3-n_5)
	+2n_3(n_2-n_1-n_5).$

	\sssec{Definition.}
	The defining quadric and any other which has the the same
	ideal polarity is called a {\em sphere}.

	\sssec{Theorem.}
	{\em All spheres degenerate or not are given by\\
	\hth$\Phi := k_1\Theta + k_2\:{\cal I} \bigx {\cal R}.$\\
	with not both $k_1$ and $k_2$ equal to 0 and ${\cal R}$ a plane,
	distinct from the ideal plane.}

	A sphere can be reduced to a point or be degenerate in the ideal
	plane and an other plane, when $k_1 = 0$.

	\sssec{Definition.}
	The plane ${\cal R}$ of the preceding Theorem is called the
	{\em radical plane of the 2 spheres $\Theta$ and $\Phi$}.

	\sssec{Exercise.}
	Give an example of a sphere which reduces to a single point.

	\sssec{Theorem.}
	{\em Given 2 ordinary points $A$ and $B,$ on a sphere and the polar of
	the ideal point on $A \times B$ is incident to $A \times B$ at $C,$
	then $C$ is independent of the sphere.} $C$ is called {\em mid-point}
	of $(A,B).$

	\sssec{Exercise.}
	Generalize the construction of the polar $p$ of a point $P$ with respect
	to a conic to that of the polar ${\cal P}$ of a point $P$ with respect
	to a quadric.

	\sssec{Exercise.}
	Give a construction of the mid-point of 2 points.
	\ssec{The general tetrahedron.}\label{sec-Sgentetra}
\setcounter{subsubsection}{-1}
	\sssec{Introduction.}
	The study of the geometry of the triangle can be generalized in 3
	dimensions to the study of the general tetrahedron.  A special case
	occurs very naturally, that of the orthogonal tetrahedron, studied in
	section \ref{sec-Sorthotetra}.

	\sssec{Notation.}
	Let $IV$ be the set $\{0,1,2,3\}$ and $VI$ be the set $\{0,1,2,3,4,5\}$,
	let $d$ be a function from the set $IV \bigx IV$ to the set $VI$ defined
	by\\
	\hth$	d_{0,1} = 0,$ $d_{0,2} = 1,$ $d_{0,3} = 2,$\\
	\hth$	d_{1,2} = 3,$ $d_{3,1} = 4,$ $d_{2,3} = 5,$\\
	$d_{i,i}$ is undefined, $d_{i,j} = d_{j,i}.$
	$d^{-1}$ denotes the inverse function.\\
	Similarly, let $e$ be a function from the set $IV \bigx IV \bigx IV$ to
	the set $VI$ defined by\\
	\hth$e_{1,2,3} = 0,$ $e_{2,3,0} = 1,$ $e_{3,0,1} = 2,$ $e_{0,1,2} = 3,
$\\
	$e_{i,j,k}$ unchanged when we permute indices and $e_{i,j,k}$ undefined
	if 2 indices are equal.  $e^{-1}$ denotes its inverse.

	\sssec{Notation.}
	$a \bigx b$ indicates first that the lines $a$ and $b$ have a point
	$P$ in common and second define $P.$  The plane through the lines $a$
	and $b$ is similarly denoted by  $a\:{\cal X}\: b.$  See 3.2.D.12.
	In this section, the indices $i,j,k,l$ are in the set
	$\{0,1,2,3\},$ the indices $u,v$ are in the set $\{0,1,2,3,4,5\}.$
	Unless indicated explicitely the indices $i,j,k,l$ or $u,v$ in a
	given statement are distinct.\\
	$l_{i,j}$ or $l_{d_{i,j}}$ represent the same line, the second forms
	indicates explicitely the mapping used to map the 2 dimensional
	array into a 1 dimensional array.\\
	The set $I_J := \{(0,1),(0,2),(0,3),(1,2),(3,1),(2,3)\}.$\\
	The set $J_I := \{(1,0),(2,0),(3,0),(2,1),(1,3),(3,2)\}.$\\
	The definitions only define the object $u(i,j),$ $(i,j)$ $\in I_J$
	and not that $(i,j)$ is in the set $I_J$ unless indicated explicitely.\\
	If $(j,i) \in J_I,$ then $u(i,j) = u(j,i).$\\
	If a definition is followed by $(*),$ this means that one of several
	definitions can be used, those not used are Theorems, for instance in
	D1.12.  $O_1$ can be defined by $facealtitude_{0,1} \bigx
	facealtitude_{3,1},$ and $O_1 \vee facealtitude_{2,1} = 0.$
	A quadric is denoted by a greek letter, $\theta$  say, the point
	quadric is then denoted  by $\Theta,$ the plane quadric by $|\Theta.$

	\sssec{Comment.}
	In this section, I will only give the expression of one of
	the points in a set, the others are obtained as follows, if a point
	$P_{nu}$ is defined symmetrically from $A_1,$ $A_2$ and $A_3$ the
	point $P_{nv}$ is obtained as follows,\\
	\hth	let $nu = n_{d_{i,j}}$ then $nv = n_{d_{i+1,j+1}}.$\\
	where the addition witin the subscripts is done modulo 4.\\
	In particular,\\
	\hth$	n_0 = n_{0,1}$ becomes $n_{1,2} = n_3,$\\
	\hth$	n_1 = n_{0,2}$ becomes $n_{1,3} = n_4,$\\
	\hth$	n_2 = n_{0,3}$ becomes $n_{1,0} = n_0,$\\
	\hth$	n_3 = n_{1,2}$ becomes $n_{2,3} = n_5,$\\
	\hth$	n_4 = n_{3,1}$ becomes $n_{0,2} = n_1,$\\
	\hth$	n_5 = n_{2,3}$ becomes $n_{3,0} = n_2,$\\
	If a line $l_u$ is defined non symmetrically in terms of
	$A_0,$	$A_1,$ $A_2,$ $A_3$ then $l_v$ is obtained by means of a
	permutation P of \{0,1,2,3\}.

	If $l_0 = f(n_0,n_1,n_2,n_3,n_4,n_5),$ then\\
	\hth$l_1 = f(n_1,n_2,n_0,n_5,n_3,n_4),$
		$l_2 = f(n_2,n_0,n_1,n_4,n_5,n_3).$\\
	\hth$l_5 = f(n_5,n_1,n_3,n_2,n_4,n_0),$\\
	\hth$l_4 = f(n_4,n_2,n_5,n_0,n_3,n_1),$
	$l_3 = f(n_3,n_0,n_4,n_1,n_5,n_2).$

	\sssec{Notation.}
	$\{A_i\}$ will denote the tetrahedron with vertices $A_i.$ 
	If we want to indicate explicitely not only the vertices $A_i$ but
	also the edges $a_u$ and the faces ${\cal A}_j$ we will use the more
	elaborate notation $\{A_i,a_u,{\cal A}_j\}.$ 

	\sssec{Comment.}
	For the tetrahedron with vertices $A_0,$ $A_1,$ $A_2,$ $A_3,$ the
	algebra will be done assuming these have the coordinates to be
	(1,0,0,0), (0,1,0,0), (0,0,1,0), (0,0,0,1),
	and that the barycentric point $M$ has coordinates (1,1,1,1).

	\sssec{Theorem.}
	{\em If the coordinates of a point $P$ are $(m_0,m_1,m_2,m_3),$
	$m_0,m_1,m_2,m_3 \neq 0,$
	those of the plane ${\cal P}$, which is its polar with respect to the
	tetrahedron $\{A_i\}$ are $\{m_0^{-1},m_1^{-1},m_2^{-1},m_3^{-1}\}.$ }

	The Euclidean geometry will be defined starting with the ideal plane
	${\cal I}$ which is the polar of $M$ with respect to the tetrahedron
	and starting from the quadric\\
	\hth$\Theta : n_0 X_0 X_1 + n_1 X_0 X_2 + n_2 X_0 X_3 + n_3 X_1 X_2
		+ n_4 X_3 X_1 + n_5 X_2 X_3 = 0.$\\
	as one of the spheres.  Prefering ${\cal I}$ and $\Theta$  allows us to
	define parallelism and orthogonality.

	\sssec{Theorem.}\label{sec-tgentetra}
	{\em Given}\\
	H0.0.\hti{3}$  A_0, A_1, A_2, A_3,$\\
	H0.1.\hti{3}$  M,$ $\ov{M},$\\
	H0.2.\hti{3}$  \Theta.$

	{\em Let}\\
	\noindent D0.0.\hti{3}$  a_{i,j} := A_i \vee A_j,$\\
	D0.1.\hti{3}$  {\cal A}_l := A_i \vee A_j \vee A_k,$\\
	D0.2.\hti{3}$  {\cal I} := polar(M)$ {\em with respect to the
		tetrahedron,}

	\noindent D1.0.\hti{3}$  C := pole({\cal I}),$\\
	D1.1.\hti{3}$  euler := C \vee M,$\\
	D1.2.\hti{3}$  AP_i := pole({\cal A}_i),$\\
	D1.3.\hti{3}$  med_i := C \vee {\cal A}_i,$\\
	D1.4.\hti{3}$  Imed_i := {\cal I}\wedge  med_i,$\\
	D1.5.\hti{3}$  alt_i := A_i \vee Imed_i,$\\
	D1.6.\hti{3}$  Foot_i := {\cal A}_i \wedge alt_i,$\\
	D1.7.\hti{3}$  ipa_{i,j} := Imed_i \vee Imed_j, (i,j) \in I_J,$\\
	D1.8.\hti{3}$  {\cal P}erp_{i,j} := ipa(i,j) \vee A_i, (i,j) \in I_J,$\\
	\hth$	{\cal P}erp_{j,i} := ipa(i,j) \vee A_j, (j,i) \in J_I,$\\
	D1.9.\hti{3}$  Face foot_{i,j} := {\cal P}erp_{i,j}\wedge a_{k,l},
(i,j) \in I_J \:or\: J_I,$\\
	D1.10.\hti{2}$ facealtitude_{i,j} := Face foot_{i,j} \vee A_i,
(i,j) \in I_J\: or\: J_I,$\\
	D1.11.\hti{2}$ O_i := facealtitude_{j,i} \bigx facealtitude_{k,i},$
		$i,j,k\: distinct (*),$\\
	D1.12.\hti{2}$ Mid_i := Foot_i + O_i,$\\
	D1.13.\hti{2}$ mid_i := Mid_i \vee Imed_i,$\\
	D1.14.\hti{2}$ H := mid_0 \bigx mid_1, (*)$\\
	D1.15.\hti{2}$ \eta := ${\em quadric through} $alt_i (*),$\\
	{\em then}\\
	C1.0.\hti{3}$  M = C + H$\\
	C1.1.\hti{3}$  H \vee euler = 0.$\\
	C1.2.\hti{3}$Foot_{i_j} = Foot_{j_i}.$\\
	C1.3.\hti{3}$  O_i \vee \eta = 0.$

	\noindent The nomenclature or alternate definitions:\\
	\hth$	A_i$ are the {\em vertices},\\
	\hth$	M$ is the {\em barycenter},\\
	N0.0.\hti{3}$  a_u$ are the {\em edges},\\
	N0.1.\hti{3}$  {\cal A}_i$ are the {\em faces},\\
	\hth	The {\em tetrahedron} is $(A_i,a_u,{\cal A}_i),$\\
	N0.2.\hti{3}$  {\cal I}$ is the {\em ideal plane},

	\noindent N1.0.\hti{3}$  C$ is the {\em center of the circumsphere},\\
	N1.1.\hti{3}$  euler$ is the {\em line of Euler},\\
	N1.2.\hti{3}$  AP_i$ is the {\em pole of the face} ${\cal A}_i,$\\
	N1.3.\hti{3}$  med_i$ is the {\em mediatrix of the face} ${\cal A}_i,$\\
	N1.4.\hti{3}$  Imed_i$ is the {\em ideal point on the mediatrix}
		$med_i,$\\
	N1.5.\hti{3}$  alt_i$ is the {\em altitude} corresponding to $A_0,$\\
	N1.6.\hti{3}$  Foot_i$ is the {\em foot of} $alt_i,$ corresponding to
		${\cal A}_i,$\\
	N1.7.\hti{3}$  ipa_{i,j}$ is the {\em direction of the planes
		perpendicular to} $a_{k,l},$\\
	N1.8.\hti{3}$  {\cal P}erp_{i,j}, (i,j) \in I_J,$ is the {\em plane
		perpendicular to} $a_{k,l}$ {\em through} $A_i,$\\
	\hth$	{\cal P}erp_{j,i},$ $(j,i) \in J_I,$ is the {\em plane
		perpendicular to} $a_{k,l}$ {\em through} $A_j,$\\
	N1.9.\hti{3}$  Face foot_{i,j},$ $(i,j) \in I_J$ or $J_I,$ is the{\em
		face-foot in the face} ${\cal A}_j,$ {\em on the edge\\
		\hth opposite $A_i$,
		$A_i \vee Facefoot_{i,j}$ is perpendicular to $a_{k,l},$}\\
	N1.10.\hti{2}$ facealtitude_{i,j},$ $(i,j) \in I_J$ or $J_I,$ is the
		{\em face-altitude in the face} ${\cal A}_j$\\
	\hth	{\em through the vertex $A_i,$ perpendicular to} $a_{k,l},$\\
	N1.11.\hti{2}$ O_i$ is the {\em orthocenter of} ${\cal A}_i,$\\
	N1.12.\hti{2}$ mid_i$ is the {\em perpendicular to} ${\cal A}_i$ through
		$M_i,$\\
	N1.13.\hti{2}$ H$ is the {\em center of the hyperboloid} $\eta.$\\
	N1.14.\hti{2}$\eta$ is the {\em hyperboloid of Neuberg.}

	Proof:\\
	G0.0.\hti{3}$  A_0 = (1,0,0,0).$\\
	G0.1.\hti{3}$  M = (1,1,1,1).$\\
	G0.2.\hti{3}$  \Theta: n_0 X_0 X_1 + n_1 X_0 X_2 + n_2 X_0 X_3
			+ n_3 X_1 X_2 + n_4 X_3 X_1 + n_5 X_2 X_3 = 0.$

	\noindent P0.0.\hti{3}$  a_0 = [1,0,0,0,0,0].$\\
	P0.1.\hti{3}$  {\cal A}_0 = \{1,0,0,0\}.$\\
	P0.2.\hti{3}$  {\cal I} = \{1,1,1,1\}.$

	\noindent P1.0.\hti{3}$C = (\\
	\hti{12}2n_3n_4n_5+n_0n_5(n_5-n_3-n_4)
		+n_1n_4(n_4-n_5-n_3)+n_2n_3(n_3-n_4-n_5),\\	
	\hti{12}2n_5n_1n_2+n_3n_2(n_2-n_5-n_1)+n_4n_1(n_1-n_2-n_5)
		+n_0n_5(n_5-n_1-n_2),\\	
	\hti{12}2n_2n_4n_0+n_5n_0(n_0-n_2-n_4)+n_1n_4(n_4-n_0-n_2)
		+n_3n_2(n_2-n_4-n_0),\\	
	\hti{12}2n_0n_1n_3+n_2n_3(n_3-n_0-n_1)+n_4n_1(n_1-n_3-n_0)
		+n_5n_0(n_0-n_1-n_3).$
	P1.1.\hti{3}$  euler = [
		2n_5(n_3n_4-n_1n_2)+(n_1n_4-n_2n_3)(n_2+n_4-n_1-n_3)
			+n_0n_5(n_1+n_2-n_3-n_4),\\
		2n_4(n_5n_3-n_2n_0)+(n_2n_3-n_0n_5)(n_0+n_3-n_2-n_5)
			+n_1n_4(n_2+n_0-n_5-n_3),\\
		2n_3(n_4n_5-n_0n_1)+(n_0n_5-n_1n_4)(n_1+n_5-n_0-n_4)
			+n_2n_3(n_0+n_1-n_4-n_5),\\
		2n_2(n_1n_5-n_0n_4)+(n_0n_5-n_1n_4)(n_4+n_5-n_0-n_1)
			+n_2n_3(n_0+n_4-n_1-n_5),\\
		2n_1(n_0n_3-n_2n_5)+(n_2n_3-n_0n_5)(n_5+n_3-n_2-n_0)
			+n_1n_4(n_2+n_5-n_0-n_3),\\
		2n_0(n_2n_4-n_1n_3)+(n_1n_4-n_2n_3)(n_3+n_4-n_1-n_2)
			+n_5n_0(n_1+n_3-n_2-n_4),$\\
	\noindent P1.2.\hti{3}$  AP_0 = ( 2n_3n_4n_5,n_5(n_0n_5-n_2n_3-n_1n_4),
		n_4(n_1n_4-n_0n_5-n_2n_3),\\
	\hti{12}n_3(n_2n_3-n_1n_4-n_0n_5) ).$\\
	P1.3.\hti{3}$  med_0 = [n_5(n_5-n_3-n_4),n_4(n_4-n_5-n_3),
		n_3(n_3-n_4-n_5),\\
	\hti{12}n_4(n_1-n_2)-n_5(n_0-n_2),n_5(n_0-n_1)-n_3(n_2-n_1),\\
	\hti{12}n_3(n_2-n_0)-n_4(n_1-n_0)].$\\
	P1.4.\hti{3}$  Imed_0 = (n_3^2+n_4^2+n_5^2-2(n_4n_5+n_5n_3+n_3n_4),$\\
	\hth$	-n_5(3n_0+2n_5-n)-(n_1-n_2)(n_3-n_4),$\\
	\hth$	-n_4(3n_1+2n_4-n)-(n_2-n_0)(n_5-n_3),$\\
	\hth$	-n_3(3n_2+2n_3-n)-(n_0-n_1)(n_4-n_5))).$\\
	P1.5.\hti{3}$  alt_0 = [n_5(3n_0+2n_5-n)+(n_1-n_2)(n_3-n_4),\\
	\hti{12}n_4(3n_1+2n_4-n)+(n_2-n_0)(n_5-n_3),\\
	\hti{12}n_3(3n_2+2n_3-n)+(n_0-n_1)(n_4-n_5), 0, 0, 0].$\\
	P1.6.\hti{3}$  Foot_0 = (0,n_5(3n_0+2n_5-n)+(n_1-n_2)(n_3-n_4),\\
	\hti{12}n_4(3n_1+2n_4-n)+(n_2-n_0)(n_5-n_3),\\
	\hti{12}n_3(3n_2+2n_3-n)+(n_0-n_1)(n_4-n_5)],$\\
	P1.7.\hti{3}$  ipa_0 = [2n_5,n_3-n_4-n_5,-n_3+n_4-n_5,-n_1+n_2+n_5,
		-n_1+n_2-n_5,\\
	\hti{12}-n_1+n_2+n_3-n_4.$\\
	P1.8.\hti{3}$  {\cal P}erp_{0,1} = \{ 0,-n_4+n_3+n_2-n_1,-n_5+n_2-n_1,
		n_5+n_2-n_1 \},$\\
	\hth$	{\cal P}erp_{1,0} = \{ n_4-n_3-n_2+n_1,0,-n_5+n_4-n_3,n_5+n_4
		-n_3 \}.$\\
	P1.9.\hti{3}$  Face foot_{0,1} = (0,0,-n_1+n_2+n_5,n_1-n_2+n_5),$\\
	\hth$	Face foot_{1,0} = (0,0,-n_3+n_4+n_5,n_3-n_4+n_5).$\\
	P1.10.\hti{2}$facealtitude_{0,1} = [0,-n_1+n_2+n_5,n_1-n_2+n_5,0,0,0],
		$\\
	\hth$	facealtitude_{1,0} = [0,0,0,n_3-n_4-n_5,n_3-n_4+n_5,0].$\\
	P1.11.\hti{2}$ O_0 = (0, n_5^2-(n_3-n_4)^2, n_4^2-(n_5-n_3)^2,
		n_3^2-(n_4-n_5)^2).$\\
	P1.12.\hti{2}$ Mid_0 = (0,n_5(3n_0+n_5-n)+(n_3-n_4)(n_1-n_2+n_3-n_4),\\
	\hti{12}n_4(3n_1+n_4-n)+(n_5-n_3)(n_2-n_0+n_5-n_3),\\
	\hti{12}n_3(3n_2+n_3-n)+(n_4-n_5)(n_0-n_1+n_4-n_5)),$\\
	P1.13.\hti{2}$ mid_0 = [n_5(3n_0+n_5-n)+(n_3-n_4)(n_1-n_2+n_3-n_4),\\
	\hti{12}n_4(3n_1+n_4-n)+(n_5-n_3)(n_2-n_0+n_5-n_3),\\
	\hti{12}n_3(3n_2+n_3-n)+(n_4-n_5)(n_0-n_1+n_4-n_5),\\
	\hti{12}(n_0-n_1-n_4+n_5)(n_5+n_4-n_3),\\
	\hti{12}(n_2-n_0-n_5+n_3)(n_3+n_5-n_4),
		(n_1-n_2-n_3+n_4)(n_4+n_3+n_5)]$.\\
	P1.14.\hti{2}$ H = (n_3n_4n_5+n_1n_2n_5+n_2n_0n_4+n_0n_1n_3
	+n_0n_5(n_0-n_1-n_2)\\
	\hti{16}+n_1n_4(n_1-n_2-n_0)+n_2n_3(n_2-n_0-n_1),\\
	\hth{12}n_3n_4n_5+n_1n_2n_5+n_2n_0n_4+n_0n_1n_3
		+n_0n_5(n_0-n_3-n_4)\\
	\hti{16}+n_1n_4(n_4-n_0-n_3)+n_2n_3(n_3-n_4-n_0),\\
	\hth{12}n_3n_4n_5+n_1n_2n_5+n_2n_0n_4+n_0n_1n_3
		+n_0n_5(n_5-n_1-n_3)\\
	\hti{16}+n_1n_4(n_1-n_3-n_5)+n_2n_3(n_3-n_1-n_5),\\
	\hth{12}n_3n_4n_5+n_1n_2n_5+n_2n_0n_4+n_0n_1n_3
	 +n_0n_5(n_5-n_2-n_4)\\
	\hti{16}+n_1n_4(n_4-n_5-n_2)+n_2n_3(n_2-n_5-n_4)).$\\
	P1.15.\hti{2}$ \eta: r_0 X_0 X_1 + r_1 X_0 X_2 + r_2 X_0 X_3
			+ r_3 X_1 X_2 + r_4 X_3 X_1 + r_5 X_2 X_3 = 0.$\\
	\hth$r_0 = (n_1-n_2-n_3+n_4)(n_0(3n_5+2n_0-n)+(n_1-n_3)(n_2-n_4),\\
	\hth r_1 = (n_2-n_0-n_5+n_3)(n_1(3n_4+2n_1-n)+(n_2-n_5)(n_0-n_3),\\
	\hth r_2 = (n_0-n_1-n_4+n_5)(n_2(3n_3+2n_2-n)+(n_0-n_4)(n_1-n_5),\\
	\hth r_3 = (n_0-n_4-n_1+n_5)(n_3(3n_2+2n_3-n)+(n_0-n_1)(n_4-n_5),\\
	\hth r_4 = (n_2-n_5-n_0+n_3)(n_4(3n_1+2n_4-n)+(n_2-n_0)(n_5-n_3),\\
	\hth r_5 = (n_1-n_3-n_2+n_4)(n_5(3n_0+2n_5-n)+(n_1-n_2)(n_3-n_4).$

	Details for the computation of P.15 are given in \ref{sec-thypneu}.

	\sssec{Comment.}\label{sec-cgentetra}
	A simple derivation for some of the points in the faces follows from a
	direct application of the results on the geometry of the triangle.
	Indeed, the circumcircle in ${\cal A}_0$ is on the one hand\\
	\hth$	 n_3\: X_1 X_2 + n_4\: X_3 X_1 + n_5\: X_2 X_3 = 0$\\
	and on the other hand\\
	\hth$m_3(m_1+m_2)\: X_1 X_2 + m_2(m_3+m_1)\: X_3 X_1
		+ m_1(m_2+m_3)\: X_2 X_3 = 0,
$\\
	assuming the coordinates of the orthocenter ${\cal A}_0$ to be
	$(0,m_1,m_2,m_3).$\\
	Comparing we get\\
	\hth$m_1m_2 = n_4+n_5-n_3,$ $m_2m_3 = n_3+n_4-n_5,$
	$m_3m_1 = n_5+n_3-n_4.$\\
	$m_1,$ $m_2,$ $m_3$ are proportional to $(m_1m_2)(m_3m_1),$
	$(m_1m_2)(m_2m_3),$ $(m_2m_3)(m_3m_1),$
	therefore using homogeneity\\
	\hth$	m_1 = n_5^2-(n_3-n_4)^2,$\\
	\hth$	m_2 = n_4^2-(n_5-n_3)^2,$\\
	\hth$	m_3 = n_3^2-(n_4-n_5)^2.$\\
	This can therefore be used to derive all the elements in the plane
	directly from Theorem 2.6. of Chapter 2.
	The following are useful,\\
	\hth$	m_2+m_3 = n_5(n_3+n_4-n_5),$\\
	\hth$	m_3+m_1 = n_4(n_5+n_3-n_4),$\\
	\hth$	m_1+m_2 = n_3(n_4+n_5-n_3)$\\
	(The notation is only valid in ${\cal A}_0,$ to have a notation for all
	faces, $m_5,$ $m_4,$ $m_3$ should be replaced by $m_{01},$ $m_{02},$
$m_{03}.$)\\
	For instance, to obtains $Foot_0,$\\
	\hth$ia_0 := | A_0 | I = [0,0,0,1,1,1] = [[1,0,0]],$\\
	\hth$Pia_0 := pole(ia_0) \in {\cal I} = ((i_{00},i_{01},i_{02}))
= (i_{00},i_{01},i_{02},i_{03}),$\\
	with $i_{03} = -i_{00}-i_{01}-i_{02} = (n_0-n_1-n_3)(n_4-n_3-n_5)
+2n_3(n_2-n_1-n_5).$\\
	hence $Foot_0 := (A_0 \vee Pia_0) {\cal A}_0 = (0,i_{01},i_{02},i_{03})
	.$  Hence\\
	\hth$Foot_{i,i} = 0,$\\
	\hth$f_{I_J} := Foot_{i_j} = Foot_{j_i}
		= (n_{i,k}-n_{i,l})(n_{j,k}-n_{j,l}) + n_{k,l}(3n_{i,j}
+2n_{k,l}-n)$\\
	\hth		for $i\neq j,$ and $i,j,k,l$ a permutation of 0,1,2,3.\\
	Hence the Theorem as well as C1.3.\\
	Similarly the center of the circumcircle $\in {\cal A}_0$ is\\
	\hth$	(0,n_5(n_3+n_4-n_5),n_4(n_5+n_3-n_4),n_3(n_4+n_5-n_3)).$

	\sssec{Theorem.}\label{sec-thypneu}
	{\em Let $r_{i,j}$ be the coordinate of $X_iX_j$ in $\eta$.\\
	Let the coordinates of the feet $Foot_0,$ $Foot_1,$ $Foot_2,$ $Foot_3$
	be $(0,f_{01},f_{02},f_{03}),$ $(f_{10},0,f_{12},f_{13}),$
	$ (f_{20},f_{21},0,f_{23}),$ $(f_{30},f_{31},f_{32},0),$\\
	then\\
	\hth$      r_{i,j} = f_{k,l}(f_{i,k}f_{j,l}-f_{i,l}f_{j,k}),$
	where $i,j,k,l$ is an even permutation of} 0,1,2,3.

	Proof:  Let the inverse $f_{J_I}$ of $f_{I_J}$ modulo $p$ be denoted
	$g_{I_J}.$\\
	If all $f_{ij}\neq 0,$ expressing the fact that the quadric contains the
	altitudes gives the equations
\\
	(0)\hti{5}$  f_{01} r_0 + f_{02} r_1 + f_{03} r_2 = 0,$\htn
	(1)\hti{5}$  g_{01} r_5 + g_{02} r_4 + g_{03} r_3 = 0,$\\
	(2)\hti{5}$  f_{10} r_0 + f_{12} r_3 + f_{13} r_4 = 0,$\htn
	(3)\hti{5}$  g_{10} r_5 + g_{12} r_2 + g_{13} r_1 = 0,$\\
	(4)\hti{5}$  f_{20} r_1 + f_{21} r_3 + f_{23} r_5 = 0,$\htn
	(5)\hti{5}$  g_{20} r_4 + g_{21} r_2 + g_{23} r_0 = 0,$\\
	(6)\hti{5}$  f_{30} r_2 + f_{31} r_4 + f_{32} r_5 = 0,$\htn
	(7)\hti{5}$  g_{30} r_3 + g_{31} r_1 + g_{32} r_0 = 0.$\\
	Equations (0) and (7) are obtained by substituting in the
	equation of the quadric, $X_i$ by $k A_i + Foot_i,$\\
	$g_{30},$ $g_{31},$ $g_{32}$ are proportional to $f_{01}f_{02},$
	$f_{03}f_{01},$ $f_{02}f_{03},$
	and therefore to $f_{03}^{-1},$ $f_{02}^{-1},$ $f_{01}^{-1}$.\\
	We solve (0) with respect to $r_0,$ and (3) with respect to $r_5,$
	in terms of $r_1$ and $r_2,$ (5) gives $r_4,$ substitution in (6) gives
	an homogeneous equation in terms of $r_1$ and $r_2$ only, hence
	after division by $f_{01}f_{03} + f_{02}f_{13}$\\
	\hth$r_{0,2} = r_1 =  f_{13}(f_{03}f_{21} - f_{01}f_{23}),$
	$r_{0,3} = r_2 = -f_{12}(f_{02}f_{31} - f_{01}f_{32}).$\\
	equations (0) give $r_0,$ (5) gives $r_4,$ (7) gives $r_3$ and
	(1) gives $r_5.$
	Hence\\
	\hth$	r_{0,1} = r_0 = -f_{23}(f_{03}f_{12} - f_{02}f_{13}),$\\
	\hth$	r_{1,2} = r_3 =  f_{03}(f_{10}f_{23} - f_{13}f_{20}),$\\
	\hth$	r_{3,1} = r_4 =  f_{02}(f_{30}f_{12} - f_{32}f_{10}),$\\
	\hth$	r_{2,3} = r_5 =  f_{01}(f_{20}f_{31} - f_{21}f_{30}),$\\
	equations (2) and (4) can be used as a check.  Summarizing the results
	gives the Theorem.  Simplifying by a common factor, we obtain P1.15.

	\sssec{Theorem.}
	{\em Let}\\
	D2.0.\hti{3}$  Ia_u := a_u \wedge {\cal I},$\\
	D2.1.\hti{3}$  {\cal P}olara_u := polar(Ia_u),$\\
	then\\
	C2.0.\hti{3}$  ipa_u \wedge {\cal P}olar_{5-u} = 0.$\\
	Nomenclature:\\
	N2.0.\hti{3}$  Ia_u$ are the {\em ideal points on the edges} $a_u,$\\
	N2.1.\hti{3}$polara_u$ is the {\em equatorial plane perpendicular to}
		$a_u.$

	Proof.\\
	P2.0.\hti{3}$  Ia_0 = (1,-1,0,0).$\\
	P2.1.\hti{3}$  {\cal P}olara_0 = \{-n_0,n_0,n_1-n_3,n_2-n_4\}.$
	\ssec{The orthogonal tetrahedron.}\label{sec-Sorthotetra}

	\sssec{Definition.}
	A {\em tetrahedron is orthogonal} iff the 3 pairs of opposite sides are
	perpendicular.

	\sssec{Lemma.}
	{\em
	$a_0 \cdot a_5 = 0$ iff $n_1 + n_4 = n_2 + n_3,$\\
	$a_0 \cdot a_1 = 0$ iff $n_3 = n_0 + n_1,$\\
	$a_1 \cdot a_2 = 0$ iff $n_5 = n_1 + n_2,$\\
	$a_2 \cdot a_0 = 0$ iff $n_4 = n_2 + n_0.$}

	The first condition expresses the orthogonality of opposite sides,
	the other conditions the orthogonality of adjacent sides.

	\sssec{Theorem.}
	{\em The tetrahedron is orthogonal iff
	the parameters of the circumsphere satisfy\\
	\hth$	n_0 + n_5 = n_1 + n_4 = n_2 + n_3.$}

	Proof.  The perpendicularity of $A_0 \vee A_1$ and $A_2 \vee A_3$
	implies, because of \ref{sec-tspaceortho}.1, with $l_j = m_j = 0$,
	except for $l_0 = m_5 = 1$,\\
	\hth$n_1-n_3 = n_2-n_4$ or $n_1 + n_4 = n_2 + n_3,$\\
	Similarly that of $A_0 \vee A_2$ and $A_1 \vee A_3$ implies\\
	\hth$	n_0 + n_5 = n_2 + n_3.$

	\sssec{Theorem.}
	{\em Given an orthogonal tetrahedron whose adjacent sied are not
	orthogonal, let}
	\enumb
	\item$	m_0 = (n_0+n_1-n_3)^{-1},$ $m_1 = (n_3+n_4-n_5)^{-1},$\\
		$m_2 = (n_5+n_1-n_2)^{-1},$ $m_3 = (n_2+n_4-n_0)^{-1},$\\[10pt]
	then
	\item$	n_0 = (m_0+m_1)m_2m_3,$ $n_1 = (m_0+m_2)m_3m_1,$
	$n_2 = (m_0+m_3)m_1m_2,$\\
		$n_3 = (m_1+m_2)m_0m_3,$ $n_4 = (m_3+m_1)m_0m_2,$
	$n_5 = (m_2+m_3)m_0m_1.$
	\enume

	Proof:  The non othogonality of adjacent sides implies that the $m_j$
	are well defined. We obtain, because of the orthogonality of opposite
	sides,\\
	\hth$	m_0^{-1} + m_1^{-1} = 2 n_0,$ $m_2^{-1} + m_0^{-1} = 2 n_1,$
	$m_1^{-1} + m_2^{-1} = 2 n_3,$\\
	we also obtain\\
	\hth$	m_0^{-1} + m_3^{-1} = n_1-n_3+n_2+n_4 = 2 n_2,$\\
	\hth$	m_1^{-1} + m_3^{-1} = n_3-n_1+n_2+n_4 = 2 n_4,$\\
	\hth$	m_2^{-1} + m_3^{-1} = n_1+n_2+n_3+n_4-2n_0 = 2 n_5.$\\
	If we multiply by $\frac{1}{2}m_0m_1m_2m_3$ we get 1.

	\sssec{Comment.}
	The indices obey the following rules.\\
	Let $n_{i,j}$ be the coefficient of $X_i X_j,$\\
	we have $n_{0,1} = n_0,$ $n_{0,2} = n_1,$ $n_{0,3} = n_2,
		n_{1,2} = n_3,$ $n_{1,3} = n_4,$ $n_{2,3} = n_5.$\\
	The orthogonality takes the form,\\
	$n_{0,1} + n_{2,3} = n_{0,2} + n_{1,3} = n_{0,3} + n_{1,2}.$\\
	$m_i$ is the inverse of $n_{i,j} + n_{i,k} - n_{j,k},$ where $i,j,k$
	are distinct.\\
	For instance, if $l$ is distinct from $i,j,k,$
	$m_i$ is also the inverse of $n_{i,j} + n_{i,l} - n_{j,l}.$

	\sssec{Definition.}
	A orthogonal tetrahedron is called a {\em special orthogonal
	tetrahedron} at $A_i$ if 2 adjacent sides through $A_i$ are also
	orthogonal.

	\sssec{Theorem.}
	{\em If $A_0 \vee A_1$ is orthogonal to $A_0 \vee A_2$ and
	the tetrahedron is orthogonal then these lines are orthogonal to
	$A_0 \vee A_3$ and\\
	\hth$	n_0 + n_1 - n_3 = n_2 + n_0 -n_4 = n_1 + n_2 - n_5 = 0.$\\
	Vice versa, if $n_0 + n_1 - n_3 = 0$ and the tetrahedron is orthogonal,
	then it is special at $A_0.$ }

	\sssec{Exercise.}
	Discuss the special cases
	\enumb
	\item $n_0+n_1-n_3 = 0,$  $n_3+n_4-n_5 \neq 0$ \ldots.
	\item $n_0+n_1-n_3 = 0,$  $n_3+n_4-n_5= 0$
	\item $n_0 = 0.$
	\enume

	\sssec{Theorem.}
	{\em The coordinates of the points lines and planes defined in
	\ref{sec-tgentetra} are}
	\enumb
	\item[G0.0.~] $A_0 = (1,0,0,0),$
	\item[G0.1.~] $M = (1,1,1,1),$
	\item[G0.2.~] $\Theta:\:(m_0+m_1)m_2m_3 X_0 X_1 +
	(m_0+m_2)m_1m_3 X_0 X_2 + (m_0+m_3)m_1m_2 X_0 X_3\\
	\hti{2} + (m_1+m_2)m_0m_3 X_1 X_2 + (m_3+m_1)m_0m_2 X_3 X_1 + 
	(m_2+m_3)m_0m_1 X_2 X_3 = 0.$
	\item[P0.0.~] $a_0 = a_{0,1} = [1,0,0,0,0,0],$
	\item[P0.1.~] ${\cal A}_0 = \{1,0,0,0\},$
	\item[P0.2.~] ${\cal I} = \{1,1,1,1\},$
	\item[P1.0.~] $C = (-m_0+m_1+m_2+m_3,m_0-m_1+m_2+m_3,
		m_0+m_1-m_2+m_3,\\
	\hth m_0+m_1+m_2-m_3),$
	\item[P1.1.~]$euler = [m_0-m_1,m_0-m_2,m_0-m_3,m_1-m_2,m_3-m_1,m_2-m_3],
$
	\item[P1.2.~] $AP_0 = (m_0(m_1+m_2)(m_3+m_1)(m_2+m_3),
		-m_1(m_2+m_3)(m_0m_1+m_2m_3),\\
	\hth -m_2(m_3+m_1)(m_0m_2+m_3m_1),-m_3(m_1+m_2)(m_0m_3+m_1m_2)),$
	\item[P1.3.~] $med_0 = [m_2+m_3,m_3+m_1,m_1+m_2,m_2-m_1,m_1-m_3,
		m_3-m_2],$
	\item[P1.4.~] $Imed_0 = (-(m_1+m_2+m_3),m_1,m_2,m_3),$
	\item[P1.5.~] $alt_0 = [m_1,m_2,m_3,0,0,0],$
	\item[P1.6.~] $Foot_0 = (0,m_1,m_2,m_3),$
	\item[P1.7.~] $ipa_{0,1} = [m_2+m_3,-m_2,-m_3,m_2,-m_3,0]$
	\item[P1.8.~] ${\cal P}erp_{0,1} = {\cal P}erp_{1,0}
	 = \{0,0,-m_3,m_2\},$
	\item[P1.9.~] $Face foot_{0,1} = Face foot_{1,0} = (m_0,m_1,0,0),$
	\item[P1.10.] $face altitude_{0,1} = [0,m_2,m_3,0,0,0],$
	$face altitude_{1,0} = [0,0,m_0,0,m_1,0],$
	\item[P1.11.] $O_0 = (0,m_1,m_2,m_3),$
	\item[P1.12.] $Mid_0 = (0,m_1,m_2,m_3),$
	\item[P1.13.] $mid_0	= $ $altitude_0$,
	\item[P1.14.] $H = (m_0,m_1,m_2,m_3),$
	\item[P1.15.] $Hyperboloid:\:\\
	\hth  m_2m_3\: X0 X1 - m_1m_3\: X0 X2
		- m_0m_2\: X1 X3 + m_0m_1\: X2 X3 = 0\\
	or\\
	\hth m_2m_3\: X0 X1 - m_1m_2\: X0 X3 - m_0m_3\: X1 X2
		+ m_0m_1\: X2 X3 = 0$.
	\enume

	\sssec{Exercise.}
	Construct a quadric generalizing the conic of Brianchon-Poncelet, and
	verify that its equation is\\
	\hth$m_1m_2m_3\:X_0^2 + \ldots - m_2m_3(m_0+m_1)\:X_0 X_1 + \ldots = 0.
$\\
	Determine points on this quadric by linear constructions which are in
	none of the faces.

	\sssec{Exercise.}
	Construct a quadric which generalizes the sphere of Prouhet,
	passing through the barycenters and orthocenters of the faces and
	verify that its equation is\\
	\hth$	3 (m_0+m_1)m_2m_3 X_0 X_1 +  \ldots $\\
	\hth$    - 2 (X_0 + X_1 + X_2 + X_3) (m_1m_2m_3 X_0 +  \ldots  ) = 0.$\\
	(Coolidge, Treatise, p. 237)

	\ssec{The isodynamic tetrahedron.}\label{sec-Sisodtetra}

	\sssec{Definition.}
	A {\em symmedian} is a line joining a vertex to the point
	of Lemoine of the opposite face.

	\sssec{Definition.}
	An {\em isodynamic tetrahedron} is a tetrahedron in which
	3 of the symmedians are concurrent.

	\sssec{Theorem.}
	{\em A tetrahedron is isodynamic  iff\\
	\hth$	n_0n_5 = n_1n_4 = n_2n_3.$}

	Proof:\\
	Let $K_i$ be the symmedian in the place ${\cal A}_i$, let $k_i :=
	A_i\times K_i,$\\
	\hth$	K_0 = (0,n_5,n_4,n_3),$ $K_1 = (n_5,0,n_2,n_1),$\\
	\hth$	K_2 = (n_4,n_2,0,n_0),$ $K_3 = (n_3,n_1,n_0,0).$\\
	$k_0 = [n_5,n_4,n_3,0,0,0],$ $k_1= [n_5,0,0,n_2,n_1,0]$ and
	$k_2 = [0,n_4,0,n_2,0,n_0]$.\\
	$k_0$ and $k_1$ are coplanar if $n_1n_4 = n_2n_3,$
	$k_0$ and $k_2$ are coplanar if $n_0 n_5 = n_2 n_3,$ hence the theorem.

	\sssec{Theorem.}
	{\em In an isodynamic tetrahedron all 4 symmedians are concurrent.}




\setcounter{section}{1}
\setcounter{subsection}{2}
	\ssec{The orthogonal tetrahedron.}

	\sssec{Definition.}
	A tetrahedron is {\em orthogonal} iff opposite sides are perpendicular.

	\sssec{Lemma.}
	{\em 
	$a_0 \cdot a_5 = 0$  iff  $n_1 + n_4 = n_2 + n_3,$\\
	$a_0 \cdot a_1 = 0$  iff  $n_3 = n_0 + n_1,$\\
	$a_1 \cdot a_2 = 0$  iff  $n_5 = n_1 + n_2,$\\
	$a_2 \cdot a_0 = 0$  iff  $n_4 = n_2 + n_0.$}

	\sssec{Theorem.}
	{\em The tetrahedron is orthogonal  iff
	the parameters of the circumsphere satisfy\\
	\hth$	n_0 + n_5 = n_1 + n_4 = n_2 + n_3.$}

	Proof:  The perpendicularity of $A_0 \vee A_1$ and $A_2 \vee A_3$ 
	implies $(0,0,1,-1) {\bf F} (1,-1,0,0)^T = 0,$ or\\
	$(n_1-n_3) - (n_2-n_4) = 0$ or\\
	\hth$	n_1 + n_4 = n_2 + n_3,$\\
	Similarly that of $A_0 \vee A_2$ and $A_1 \vee A_3$ implies\\
	\hth$	n_0 + n_5 = n_2 + n_3.$

	\sssec{Theorem.}
	{\em If $n_0+n_1-n_3\neq  0,$ $n_3+n_4-n_5\neq  0,$
	$n_5+n_1-n_2\neq 0,$ $n_2+n_4-n_0\neq 0,$
	and the tetrahedron is orthogonal.  Let}
	\enumb
	\item$	m_0 = (n_0+n_1-n_3)^{-1},$ $m_1 = (n_3+n_4-n_5)^{-1},$\\
	\hth$	m_2 = (n_5+n_1-n_2)^{-1},$ $m_3 = (n_2+n_4-n_0)^{-1},$\\
	{\em then}
	\item$	n_0 = (m_0+m_1)m_2m_3,$ $n_1 = (m_0+m_2)m_3m_1,$
	$n_2 = (m_0+m_3)m_1m_2,$\\
	\hth$	n_3 = (m_1+m_2)m_0m_3,$ $n_4 = (m_3+m_1)m_0m_2,$
	$n_5 = (m_2+m_3)m_0m_1.$
	\enume

	Proof:  We obtain at once,\\
	\hth$	m_0^{-1} + m_1^{-1} = 2 n_0,$ $m_2^{-1} + m_0^{-1} = 2 n_1,$
	$m_1^{-1} + m_2^{-1} = 2 n_3,$\\
	using 1.3.3., we also obtain\\
	\hth$	m_0^{-1} + m_3^{-1} = n_1-n_3+n_2+n_4 = 2 n_2,$\\
	\hth$	m_1^{-1} + m_3^{-1} = n_3-n_1+n_2+n_4 = 2 n_4,$\\
	\hth$	m_2^{-1} + m_3^{-1} = n_1+n_2+n_3+n_4-2n_0 = 2 n_5.$\\
	If we multiply by $\frac{1}{2}m_0m_1m_2m_3$ we get 1.

	\sssec{Comment.}
	?? when hyp. of prec theorem replace $\neq 0$ by $= 0,$
	see Theorem 3.7.

	\sssec{Comment.}
	The indices obey the following rules.\\
	Let $n_{0,1}$ be the coefficient of $X_0 X_1,$  $\ldots$  ,\\
	we have $n_{0,1} = n_0,$ $n_{0,2} = n_1,$ $n_{0,3} = n_2,$\\
	\hth$	n_{1,2} = n_3,$ $n_{1,3} = n_4,$ $n_{2,3} = n_5.$\\
	The orthogonality takes the form,\\
	$n_{0,1} + n_{2,3} = n_{0,2} + n_{1,3} = n_{0,3} + n_{1,2}.$\\
	$mi$ is the inverse of $n_{i,j} + n_{i,k} - n_{j,k},$ where $i,j,k$
	are distinct.\\
	For instance, if $l$ is distinct from $i,j,k,$
	$mi$ is also the inverse of $n_{i,j} + n_{i,l} - n_{j,l}.$

	\sssec{Definition.}
	A orthogonal tetrahedron is called a {\em special orthogonal
	tetrahedron} at $A_i$ if 2 adjacent sides through $A_i$ are also
	orthogonal.

	\sssec{Theorem.}
	{\em If $A_0 \vee A_1$ is orthogonal to $A_0 \vee A_2$ and
	the tetrahedron is orthogonal then these lines are orthogonal to
	$A_0 \vee A_3$ and\\
	\hth$	n_0 + n_1 - n_3 = n_2 + n_0 -n_4 = n_1 + n_2 - n_5 = 0.$\\
	Vice versa, if $n_0 + n_1 - n_3 = 0$ and the tetrahedron is orthogonal,
	then it is special at $A_0.$ }

	\sssec{Exercise.}
	Discuss the special cases
	\enumb
	\item $n_0+n_1-n_3 = 0,$  $n_3+n_4-n\neq 0$ \ldots.
	\item $n_0+n_1-n_3 = 0,$  $n_3+n_4-n_5= 0$
	\item $n_0 = 0.$
	\enume

	\sssec{Theorem.}
	{\em The coordinates of the points lines and planes defined in
	in 1.2.7. are}
	H0.0.\hti{3}$A_0 = (1,0,0,0),$\\
	H0.1.\hti{3}$M = (1,1,1,1),$

	P0.0.\hti{3}$a_0 = a_{0,1} = [1,0,0,0,0,0],$\\
	P0.1.\hti{3}$| 0 = \{1,0,0,0\},$\\
	P0.2.\hti{3}$| I = \{1,1,1,1\},$

	P1.0.\hti{3}$C = (-m_0+m_1+m_2+m_3,m_0-m_1+m_2+m_3,m_0+m_1-m_2+m_3,
		m_0+m_1+m_2-m_3),$
	P1.1.\hti{3}$euler = [m_0-m_1,m_0-m_2,m_0-m_3,m_1-m_2,m_3-m_1,m_2-m_3],
$\\
	P1.2.\hti{3}$AP_0 = (m_0(m_1+m_2)(m_3+m_1)(m_2+m_3),
		-m_1(m_2+m_3)(m_0m_1+m_2m_3),$\\
	\hth$		-m_2(m_3+m_1)(m_0m_2+m_3m_1),
		-m_3(m_1+m_2)(m_0m_3+m_1m_2)),$\\
	P1.3.\hti{3}$med_0 = [m_2+m_3,m_3+m_1,m_1+m_2,m_2-m_1,m_1-m_3,m_3-m_2],
$\\
	P1.4.\hti{3}$Imed_0 = (-(m_1+m_2+m_3),m_1,m_2,m_3),$\\
	P1.5.\hti{3}$alt_0 = [m_1,m_2,m_3,0,0,0],$\\
	P1.6.\hti{3}$Foot_0 = (0,m_1,m_2,m_3),$\\
	P1.7.\hti{3}$ipa_{0,1} = [m_2+m_3,-m_2,-m_3,m_2,-m_3,0]$\\
	P1.8.\hti{3}$| Perp._{0,1} = \{0,0,-m_3,m_2\},$ {\em perp. to}
	$a_{0,1}$ through $A_2 \times A_3$\\
	\hth$	 = same ??$\\
	P1.9.\hti{3}$Face foot_{0,1} = (m_0,m_1,0,0),$ {\em on} $a_{0,1}$\\
	\hth$	 = same,$\\
	P1.10.\hti{2}$face altitude_{0,1} = [0,m_2,m_3,0,0,0],$
		$ Face foot_{2,3} \vee A_0$\\
	\hth$	     [1,0] = [0,0,m_0,0,m_1,0],$
		$ Face foot_{2,3} \vee A_1$\\
	P1.11.\hti{2}$O_0 = (0,m_1,m_2,m_3),$\\
	P1.12.\hti{2}$Mid_0 = (0,m_1,m_2,m_3),$\\
	P1.13.\hti{2}$mid_0 = [m_1,m_2,m_3,0,0,0],$\\
	P1.14.\hti{2}$H = (m_0,m_1,m_2,m_3),$\\
	P1.15.\hti{2}$Eta: m_2m_3 X0 X1 - m_1m_3 X0 X2 - m_0m_2 X1 X3 
		+ m_0m_1 X2 X3 = 0\\
		m_2m_3 X0 X1 - m_1m_2 X0 X3 - m_0m_3 X1 X2 + m_0m_1 X2 X3 = 0.$

	${\cal C}oideal = \{m_1m_2m_3,m_2m_3m_0,m_3m_0m_1,m_0m_1m_2\},$\\
	$Cocenter = Barycenter,$

	\sssec{Theorem.}
	{\em If $(0,p_1,p_2,p_3)$ is on the Euler line $eul^{(0)}$ then}
	\enumb
	\item$      p_1(m_2-m_3) + p_2(m_3-m_1) + p_3(m_1-m_2) = 0,$
	\item$      P \vee IC(0)$ {\em intersects the Euler line $eu;$ at}\\
	\hth$  ( (p_1(m_0-m_2) + p_2(m_1-m_0))(m_1+m_2+m_3),$\\
	\hth$    p_1((m_1-m_2)(m_1+m_2+m_3)+m_2(m_1-m_0)) + p_2m_1(m_0-m_1),$\\
	\hth$    p_2((m_1-m_2)(m_1+m_2+m_3)+m_1(m_0-m_2)) + p_1m_2(m_2-m_0),$\\
	\hth$    p_2((m_1-m_3)(m_1+m_2+m_3)+m_1(m_0-m_3))
	+ p_1((m_3-m_2)(m_1+m_2+m_3)+m_2(m_3-m_0)) ),$\\
	{\em or more symmetrically},\\
	\hth$( p_1(s_1(m_0-m_2)+m_0(m_2-m_0) + p_2(s_1(m_1-m_0)+m_0(m_0-m_1),$\\
	\hth$p_1(s_1(m_1-m_2)+m_1(m_2-m_0) + p_2(s_1(m_1-m_1)+m_1(m_0-m_1),$\\
	\hth$p_1(s_1(m_2-m_2)+m_2(m_2-m_0) + p_2(s_1(m_1-m_2)+m_2(m_0-m_1),$\\
	\hth$p_1(s_1(m_3-m_2)+m_3(m_2-m_0) + p_2(s_1(m_1-m_3)+m_3(m_0-m_1) ),$\\
	{\em Moreover if} $p_1 = km_1+s-m_0,$ $p_2 = km_2+s-m_0,$
	$p_3 = km_3+s-m_0,$
	{\em then the point on e is}\\
	$( (k-1)m_0+s,(k-1)m_1+s,(k-1)m_2+s,(k-1)m_3+s ).$\\
	{\em In particular,}\\
	$M = (m_1+m_2+m_3,m_2+m_3+m_0,m_3+m_0+m_1,m_0+m_1+m_2),$
	$P = ( s_1t_0-m_0t_0,s_1(m_1^{2-m_2m_3)-m_1t_0,s_1(m_2^2-m_3m_1)
	-m_2t_0},$\\
	\hth$		s_1(m_3^{2-m_1m_2)-m_3t_0} ),$\\
	\hth$	{\mathit with} t_0 = m_0m_1+m_0m_2+m_0m_3-m_1m_2-m_3m_1-m_2m_3,$\\
	$\overline{O} = ($\\
	\hth$Am = (s_1-3m_0,s_1-3m_1,s_1-3m_2,s_1-3m_3),$\\
	\hth$G = (s_1+2m_0,s_1+2m_1,s_1+2m_2,s_1+2m_3),$\\
	\hth$\overline{A}m = ()$\\
	\hth$D_0 = ()$\\
	\hth$D_1 = ()$\\
	\hth$D_2 = ()$\\
	\hth$G = ()$\\
	\hth$\overline{G} = ()$
	\enume

	{\sssec{Answer (partial).}
	\ldots ?  The polar $pp_0 = [-2m_1m_2,m_2(m_0+m_1),m_1(m_2+m_0],$\\
	\ldots ?  The intersection $PP_0 = (0,-m_1(m_2+m_0),m_2(m_0+m_1)),$\\
	\ldots ?  $pp = [m_1m_2(m_2+m_0)(m_0+m_1),m_2m_0(m_0+m_1)(m_1+m_2),
	m_0m_1(m_1+m_2)(m_2+m_0)].$

	Point ``O", intersection of perpendicular to faces through their
	barycenter\\
	(special case of $\ldots$  with $k = 0$ hence\\
	$``O" = (s_1-m_0,s_1-m_1,s_1-m_2,s_1-m_3).$\\
	``Conjugate tetrahedron",\\
	$``A'_0" = (-2s_1-2m_0,s_1+2m_1,s_1+2m_2,s_1+2m_3),$\\
	barycenter of faces of $[A'[i]]$ are\\
	$``M'_0" = (0,m_1,m_2,m_3),$
	which are the orthocenters of the faces,\\
	Perpendiculars through $M'_0$ to the faces, (which are parallel to
		those of [A[]] meet at\\
	$``O'" = (3s_1-4m_0,3s_1-4m_1,3s_1-4m_2,3s_1-4m_3).$\\
	Hence his theorem:
	Then he generalizes the circle of Brianchon-Poncelet and gives its
	center as the midpoint of $H$ and $``O"$\\
	I believe his $``O"$ is my $G.$\\
	Orthocenter,barycenter and $``O"$ are collinear

	\ssec{The isodynamic tetrahedron.}

	\sssec{Definition.}
	A {\em symmedian} is a line joining a vertex to the point
	of Lemoine of the opposite face.

	\sssec{Definition.}
	An {\em isodynamic tetrahedron} is a tetrahedron in which
	3 of the symmedians are concurrent.

	\sssec{Theorem.}
	{\em A tetrahedron is isodynamic  iff\\
	\hth$	n_0n_5 = n_1n_4 = n_2n_3.$}

	Proof:\\
	\hth$	K_0 = (0,n_5,n_4,n_3),$ $K_1 = (-n_5,0,n_2,-n_1),$\\
	\hth$	K_2 = (-n_4,-n_2,0,n_0),$ $K_3 = ().$\\
	$k_0$ and $k_1$ are coplanar if $n_1n_2 = n_2n_3,$
	$k_0$ and $k_2$ are coplanar if $n_0 n_5 = n_2 n_3,$ hence the theorem.

	\sssec{Theorem.}
	{\em In an isodynamic tetrahedron all 4 symmedians are concurrent.}

	In 3 dimensions start with $A_0 = (1,0,0,0),$ $\ldots$,
	$A_3 = (0,0,0,1),$
	and with $M = (1,1,1,1)$ and $\overline{M} = (m_0,m_1,m_2,m_3).$
	$M$ corresponds to the barycenter, $\overline{M}$ to the intersection
	of the lines joining the orthocenter of the faces to the opposite vertex
	and $A[]$ to the vertices of an orthogonal tetrahedron.  See $\ldots$ 

.       Theorem.  Prove that the tetrahedron is orthogonal.

.       Theorem.  Prove that the circumscribed quadric is given by\\
	\hth$(m_0+m_1)m_2m_3 X_0 X_1 +  \ldots  = 0.$

.       Theorem.  Construct a quadric generalizing the conic of Brianchon-
	Poncelet, and verify that its equation is\\
	\hth$	m_1m_2m_3 X_0^2 +  \ldots  - m_2m_3(m_0+m_1) X_0 X_1 +  \ldots  = 0.$\\
	Determine points on this quadric by linear constructions which are in
	none of the faces.

.       Theorem.  Construct a quadric which generalizes the sphere of Prouhet,
	passing through the barycenters and orthocenters of the faces and
	verify that its equation is\\
	\hth$	3 (m_0+m_1)m_2m_3 X_0 X_1 +  \ldots $\\
	\hth$    - 2 (X_0 + X_1 + X_2 + X_3) (m_1m_2m_3 X_0 +  \ldots  ) = 0.$\\
	(Coolidge, Treatise, p. 237)\\
	Point ``O", intersection of perpendicular to faces through their
	barycenter
	(special case of $\ldots$  with k = 0 hence\\
	$``O" = (s_1-m_0,s_1-m_1,s_1-m_2,s_1-m_3)$\\
	``Conjugate tetrahedron",\\
	$``A'_0" = (-2s_1-2m_0,s_1+2m_1,s_1+2m_2,s_1+2m_3),$\\
	barycenter of faces of $[A'[i]]$ are\\
	$``M'_0" = (0,m_1,m_2,m_3),$\\
	which are the orthocenters of the faces,\\
	Perpendiculars through $M'_0$ to the faces, (which are parallel to
		those of $[A[]]$ meet at\\
	$``O'" = (3s_1-4m_0,3s_1-4m_1,3s_1-4m_2,3s_1-4m_3),$\\
	Hence his theorem:\\
	Then he generalizes the circle of Brianchon-Poncelet and gives its
	center as the midpoint of $H$ and $``O"$\\
	I believe his $``O"$ is my $G$.\\
	Orthocenter,barycenter and $``O"$ are collinear.



\setcounter{section}{1}
\setcounter{subsection}{4}
	\ssec{The antipolarity.}\label{sec-Santipol}

	\sssec{Definition.}
	Consider the 2-form $[l_0,l_1,l_2,l_3,l_4,l_5]$ with
	\enumb
	\item$      l_0 l_5 + l_1 l_4 + l_2 l_3\neq  0,$\\
	the {\em point to plane antipolarity} associates to a point $P$ a plane
	${\cal P} := l' \vee P,$\\
	the {\em plane to point antipolarity} associates to a plane ${\cal P}$
	a point $P := l' {\cal P}.$
	\enume

	\sssec{Theorem.}\label{sec-tantimat}
	{\em The point to plane antipolarity can be represented by an
	antisymmetric matrix\\
	\hth$ {\bf P} = \matb{(}{cccc}
		 0&l_5&l_4&l_3\\
		-l_5& 0&l_2&-l_1\\
		-l_4&-l_2& 0&l_0\\
		-l_3&l_1&-l_0& 0\mate{)}$\\
	The plane to point antipolarity is represented by the antisymmetric
	matrix\\
	\hth${\bf Q} = \matb{(}{cccc}
		 0&l_0&l_1&l_2\\
		-l_0& 0&l_3&-l_4\\
		-l_1&-l_3& 0&l_5\\
		-l_2&l_4&-l_5& 0\mate{)}$\\
	both have determinant $(l_0 l_5 + l_1 l_4 + l_2 l_3)^2\neq0.$}

	\sssec{Theorem.}
	{\em If ${\cal P}$ is associated to $P$ in an antipolarity then $P$ is
	associated to $| P$ and $| P$ is incident to $P.$}

	The proof is left as an exercise.

	\sssec{Theorem.}\label{sec-tantipol}
	{\em Let}
	\enumb
	\item$      d := l_0 l_5 + l_1 l_4 + l_2 l_3,$

	{\em the planes associated in the antipolarity \ref{sec-tantimat}
	to the points of a line $m$  are all incident to a line  $q$  and if}
	\item$      {\bf L} :=  l \: dual(l)^T - d {\bf I},$\\
	{\em then}
	\item$      {\bf L} =\\
		 \matb{(}{cccccc}
		-l_1l_4-l_2l_3&l_0l_4&l_0l_3&l_0l_2&l_0l_1&l_0^2\\
		 l_1l_5&-l_0l_5-l_2l_3&l_1l_3&l_1l_2&l_1^2&l_0l_1\\
		 l_2l_5&l_2l_4&-l_0l_5-l_1l_4&l_2^2&l_2l_1&l_2l_0\\
		 l_3l_5&l_3l_4&l_3^2&-l_0l_5-l_1l_4&l_3l_1&l_3l_0\\
		 l_4l_5&l_4^2&l_4l_3&l_4l_2&-l_0l_5-l_2l_3&l_4l_0\\
		 l_5^2&l_5l_4&l_5l_3&l_5l_2&l_5l_1&-l_1l_4-l_2l_3\mate{)} $
	\item$      q = {\bf L} m.$
	\item$      det({\bf L}) = - d^6.$
	\enume

	The proof is left as an exercise.

	\sssec{Example.}\label{sec-eantipol}
	Let $p = 29,$ $l' = [1,4,3,11,4,10] = [3644209],$\\
	E0.\hti{5}$    d = 1.$\\
	E2.\hti{5}$    {\bf L} = \matb{(}{rrrrrr}
			 9& 4& 11&3&4&1\\
			11&-14&-14& 12&-13&4\\
			 1&12&3&9& 12&3\\
			-6&-14&5&3&-14& 11\\
			11&-13&-14& 12&-14&4\\
			13&11& -6&1& 11&9\mate{)}$\\
	E3.\hti{6}$    m =   732541, 25620,   871,    30,     1,     0.$\\
	\hth$q =     20154561, 8595176, 4156378, 3635799, 3644412, 3644208.$

	\sssec{Theorem.}\label{sec-tantipol2}
	{\em The antipolarity of vertices, faces and edges of a
	tetrahedron follows from what follows.  Given}\\
	H.0.\hti{4}$   l' = [l_0,l_1,l_2,l_3,l_4,l_5],$\\
	H.1.\hti{4}$   A_i,$\\
	{\em let}\\
	D0.0.\hti{3}$  a_{i,j} := A_i \vee A_j,$\\
	D0.1.\hti{3}$  {\cal A}_i := a_{j,k} \vee A_l,$\\
	D1.0.\hti{3}$  {\cal B}_i := l' \vee A_i,
			B_i := l' {\cal A}_i,$\\
	D1.1.\hti{3}$  ba_i := | B_i | A_i,
			b\overline{a}_i := B_i \vee A_i,$\\
	D1.2.\hti{3}$  N_{i,j} := a_{k,l} {\cal B}_j,
			{\cal N}_{i,j} := a_{i,j} \vee B_j,$\\
	D1.3.\hti{3}$  n_{i,j} := B_i \vee A_j,$\\
	D1.4.\hti{3}$  b_{d(i,j)} := | B_i  | B_j,$\\
	{\em then}\\
	C1.0.\hti{3}$  B_i \vee {\cal B}_j = 0.$\\
	C1.1.\hti{3}$  N_{i,j} \vee {\cal N}_{j,i} = 0.$\\
	C1.2.\hti{3}$  N_{i,j} \vee n_{i,j} = 0,$\\
	C1.3.\hti{3}$  n_{i,j} = | B_j  | A_i.$\\
	C1.4.\hti{3}$  b_{d(i,j)} = B_k \vee B_l.$\\
	C1.5.\hti{3}$  b_u = {\bf L} a_u.$

	Proof.\\
	P1.0.\hti{3}$  {\cal B}_0 = \{0,-l_5,-l_4,-l_3\}.$\\
	\hth$	B_0 = (0,-l_0,-l_1,-l_2).$\\
	P1.1.\hti{3}$  ba_0 = [0,0,0,l_3,l_4,l_5].$\\
	\hth$	b\overline{a}_0 = [l_0,l_1,l_2,0,0,0].$\\
	P1.2.\hti{3}$  N_{1,2} = (0,0,l_1,l_2).$\\
	\hth$	\overline{N}_{1,2} = \{0,0,l_4,l_3\}.$\\
	P1.3.\hti{3}$  n_{1,2} = [0,0,0,l_1,-l_2,0].$\\
	P1.4.\hti{3}$  b_0 = [-l_1l_4-l_2l_3,l_1l_5,l_2l_5,l_3l_5,l_4l_5,
	l_5^2].$

	\sssec{Corollary.}
	{\em If $l$ is a line and the definitions of Theorem \ref{sec-tantipol2}
	hold then all the conclusions of \ref{sec-tantipol2} hold.  Moreover
	$b_u = l$ for all $u$.\\
	The mapping is not one to one.  The image of a point $P$ is the plane
	$P \vee l,$ the image of the plane ${\cal Q}$ is the point
	${\cal Q} \vee l.$}

	\sssec{Example.}\label{sec-eantipol1}
	Let $p = 29,$ $l' = [1,4,3,11,4,10] = [3644209],$
	$A = (871,30,1,0),$\\
	E0.0.\hti{3}$  a = [732541,25260,871,30,1,0].$\\
	E0.1.\hti{3}$  {\cal A} = \{871,30,1,0\}.$\\
	E1.0.\hti{3}$  {\cal B} = \{382,1463,7606,22969\}.$\\
	\hth$	B = (149,1397,9293,16386).$\\
	E1.1.\hti{3}$  ba = [    139, 220389, 805824,3570916],$\\
	\hth$	b\overline{a} = [3634832, 741908,  33687,   1203].$\\
	E1.2.\hti{3}$  N =  \matb{(}{rrrr}
			 -- &    9&   33&  146\\
			   27& -- &  875& 1393\\
			   37&  883& -- & 9281\\
			  378& 1248&16009& -- \mate{)}.$\\
	\hth$	{\cal N} = \matb{\{}{rrrr}
			 -- &   11&   34&  378\\
			   19& -- &  883& 1451\\
			   49&  878& -- & 7599\\
			  233& 1103&22737& -- \mate{\}}.$\\
	E1.3.\hti{3}$  n = \matb{[}{rrrr}
			 -- &    639&     56&     26\\
			 659374& -- &  25285&    889\\
			 903264& 732889& -- &   1422\\
			9219913& 745997&  40398& -- \mate{]}.$\\
	E1.4.\hti{3}$b = [20154561,8595176,4156378,3635799,3644412,3644208].$

	\sssec{Example.}
	Let $p = 29,$ $l = [3623186] = [1,4,2,-14,4,12],$ $A = (871,30,1,0),
$\\
	E1.0\hti{4}$   {\cal B} = \{343,1577,13493,18590\}.$\\
	\hth$	B = (148,1281,10148,23752).$\\
	E1.1.\hti{3}$  ba = [286,391098,781435,3574280],$\\
	\hth$	b\overline{a} = [3610443,745272,34514,935].$\\
	E1.2.\hti{3}$  N = \matb{(}{rrrr}
			 -- &   16&   32&  146\\
			   22& -- &  875& 1277\\
			   35&  897& -- &10122\\
			  784& 1045&23578& -- \mate{)}.$\\
	\hth$	{\cal N} = \matb{\{}{rrrr}
			 -- &   12&   53&  320\\
			   28& -- &  881& 1567\\
			   44&  878& -- &13486\\
			  233&  929&18532& -- \mate{\}}$.\\
	E1.3.\hti{3}$  n = \matb{[}{rrrr}
			  -- &    436&     57&     26\\
			  537429& -- &  25285&    885\\
			  854486& 733295& -- &   1393\\
			19121847& 751884&  47967& -- \mate{]}$.\\
	E1.4.\hti{3}$  b_u = [3623186].$

	\sssec{Theorem.}
	{\em An antipolarity can be determined as follows,\\
	Given 4 points $A_i,$ a line $ba_0 \in {\cal A}_0 = A_1 \vee A_2 \vee
	A_3,$\\
	a line $ba_1 \in A_1 := A_2 \vee A_3 \vee A_0$ and a point $B_0$ on
	$n_{0,1} = A_1 \vee (ba_1 \bigx a_5)$ but not on
	$( (ba_0 \bigx a_4) \vee A_0) \bigx (ba_1 \bigx a_2) \vee A_1) )
	  \vee ( (ba_0 \bigx a_3) \vee A_0) \bigx (ba_1 \bigx a_1) \vee A_1) )
	 (= B_2 \vee B_3).$ }

	Proof.  Let us choose the $A_i$ as basis for the coordinate system.\\
	$ba_0 = [0,0,0,l_3,l_4,l_5]$ determines $l_3,l_4,l_5.$\\
	$ba_1 = [0,l_1,l_2,0,0,l_5]$ determines after scaling $l_1$ and $l_2.$\\
	$B_0 =  (0,l_0,l_1,l_2),$ determines, after scaling $l_0.$ Scaling the
	last component should check with $l_2.$

	\sssec{Example.}
	Let $p = 29,$ $A = (871,30,1,0).$\\
	Let $ba_0 = [139] = [0,0,0,1,3,-7],$ $ba_1 = [220389] = [0,1,8,0,0,-12],
$\\
	$N_{0,1} = ba_1 \bigx a_5 = (9),$ $n_{0,1} := N_{0,1} \vee A_1 = [639].$
\\
	Finally let $B_0 = (149) = (0,1,4,3)$ on $[639] = [0,0,0,1,-8,0]$ but
	not on $[20154561] = [1,-2,13,9,-2,-5].$\\
	(For the details of the computations see Example \ref{sec-eantipol1}.)\\
	$ba_0$ gives $l_3 = 1,$ $l_4 = 3,$ $l_5 = -7,$\\
	$ba_1$ gives $l_1 = t,$ $l_2 = 8t,$ $l_5 = -12t,$ $t$ is the scaling
	factor,\\
	hence $l_1 = t = \frac{7}{12} = 3,$ $l_2 = -5.$\\
	$B_0$ gives $l_0 = u,$ $l_1 = 4u,$ $l_2 = 3u,$ hence $l_0 = u =
	\frac{3}{4} = 8.$\\
	$l_2 = 3.8 = -5$ is a check.\\
	Therefore $l' = [8,3,-5,1,3,-7] = [1,4,3,11,4,10].$

	The associated construction is as follows.

	\sssec{Construction.}\label{sec-cantipol}
	Given $A_0,$ $A_1,$ $A_2,$ $A_3,$ $ba_0 \in {\cal A}_0,$
	$ba_1 \in {\cal A}_1,$\\
	$N_{1,0} := ba_0 \bigx a_5,$ $n_{1,0} := N_{1,0} \vee A_0,$\\
	$N_{2,0} := ba_0 \bigx a_4,$ $n_{2,0} := N_{2,0} \vee A_0,$\\
	$N_{3,0} := ba_0 \bigx a_3,$ $n_{3,0} := N_{3,0} \vee A_0.$\\
	$N_{0,1} := ba_1 \bigx a_5,$ $n_{0,1} := N_{0,1} \vee A_1,$\\
	$N_{2,1} := ba_1 \bigx a_2,$ $n_{2,1} := N_{2,1} \vee A_1,$\\
	$N_{3,1} := ba_1 \bigx a_1,$ $n_{3,1} := N_{3,1} \vee A_1.$\\
	$B_2 := n_{2,0} \bigx n_{2,1},$ $n_{2,3} := B_2 \vee A_3,$\\
	$N_{2,3} := n_{2,3} \bigx a_0,$\\
	$B_3 := n_{3,0} \bigx n_{3,1},$ $n_{3,2} := B_3 \vee A_2,$\\
	$N_{3,2} := n_{3,2} \bigx a_0.$\\
	Given $B_0$ on $n_{0,1},$ $not on B_2 \vee B_3$ (otherwize $l'$
	is a line),\\
	$n_{0,2} := B_0 \vee A_2,$ $N_{0,2} := n_{0,2} \bigx a_4,$\\
	$n_{0,3} := B_0 \vee A_3,$ $N_{0,3} := n_{0,3} \bigx a_3.$\\
	$ba_2 := N_{0,2} \vee N_{3,2},$\\
	$N_{1,2} := ba_2 \bigx a_2,$ $n_{1,2} := N_{1,2} \vee A_2,$\\
	$B_1 := n_{1,2} \bigx n_{1,0}.$\\
	$B_i$ are the antipoles of ${\cal A}_i.$\\
	$ba_3 := N_{0,3} \vee N_{2,3}.$\\
	${\cal B}_i := A_i \vee ba_i.$\\
	${\cal B}_i$ are the antipolars of $A_i.$\\
	To complete the construction,\\
	$N_{1,3} := ba_3 \bigx a_1,$ $n_{1,3} := B_1 \vee A_3,$ we can check
	\enumb
	\item$      N_{1,3} \vee n_{1,3} = 0.$
	\enume

	\sssec{Theorem.}\label{sec-tconf2022}
	{\em In the geometry of the triangle if $A_3$ is $\overline{M},$ $ba_0$
	is $m$,
	and $ba_1$ is an arbitrary line and $B_0$ is an arbitrary point on
	$(ba_1 \times (A_2 \times \overline{M})) \times A_1,$ the configuration
	of \ref{sec-cantipol} consisting of 20 points $A,$ $B$ and $N_,$
	and of 22 lines $a,$ $ba,$ $n_,$  satisfies \ref{sec-cantipol}.0.}

	\sssec{Theorem.}
	D1.0.\hti{3}$  {\cal P}a_r := P \vee a_{5-r},$\\
	D1.1.\hti{3}$  Pab_r := {\cal P}a_r  b_{5-r},$\\
	D1.2.\hti{3}$  {\cal P} := Pab_0 \vee Pab_1 \vee Pab_2,$\\
	{\em then}\\
	C1.0.\hti{3}$  {\cal P} \vee P = 0,$\\
	C1.1.\hti{3}$  {\cal P} = P \vee l'.$

	\sssec{Example.}
	With $p,$ $l'$ and $A$ as in Example \ref{sec-eantipol1}.\\
	Let $P = (1742) = (1,1,1,1),$
	then ${\cal P}a = \{24419,1683,899\},$ $Pab = (2357,3443,25116),$
	${\cal P} = \{9747\}.$\\
	Let $P = (5350) = (1,5,9,13),$\\
	then ${\cal P}a = \{20214,1335,891),$ $Pab = (5356,5363,3726),$
	${\cal P} = \{2611\}.$

	\sssec{Exercise.}
	The antipolarity associates to a point quadric $Alpha$, a plane quadric
	${\cal B}eta,$ the points of one are on the tangent of the
	other.  Study this correspondance in detail.




\setcounter{section}{1}
\setcounter{subsection}{5}
	\ssec{Example.}
	\sssec{Case 0.}
	$p = 13,$ $Barycenter = 366,$ $n = 1,3,6,7,10,12,$ $m = 1,2,4,5,$\\
	The tetrahedron is orthogonal.\\
	$Barycenter = (366)$\\
	${\cal I}deal = \{366\}$\\
	$A = (183,14,1,0)$\\
	$a = [30941,2380,183,14,1,0]$\\
	$Center = (1244)$\\
	$Pole of {\cal A} = (1509,525,271,60)$\\
	$mediatrix = [271483,148770,212411,132416]$\\
	$IC = (387,747,591,579)$\\
	$altitude = [107836,32527,2726,330]$\\
	$Foot = (49,240,526,573)$\\
	$ipa = [85840,136570,110191,374757,58939,29111]$\\
	${\cal P}erp. = \{8,24,92,188,222,1197\}$\\
	\hti{4}$ = \{8,24,92,188,222,1197\}$\\
	$Face foot = (12,23,40,188,235,521)$\\
	\hti{4}$ = (12,23,40,188,235,521)$\\
	$face altitude = [26547,50714,88063,31006,187,2718]$\\
	\hti{4}$ = [40,18,12,2388,32462,326]$\\
	$Orthocenters = (49,240,526,573)$\\
	$Mid = (49,240,526,573)$\\
	$mid = [107836,32527,2726,330]$\\
	$Orthocenter = (578)$\\
	${\cal C}oideal = \{1504\}$\\
	$Cocenter = (366)$\\
	$Barycenter $(check)$ = (366)$\\
	$Hyperboloid:\:  11  1  0  0  3  5,\: 4  0  1  12  0  3$

	The coordinates of $IC$ are (1,1,2,9), (1,3,4,5), (1,2,5,5), (1,2,4,6),
\\
	those of $Foot$ are (0,1,2,9), (1,0,4,5), (1,2,0,5), (1,2,4,0),\\
	those of the $Center of Hyperboloid$ are $(1,2,4,5),$\\
	the hyperboloids are\\
	\hti{4}$(-2r1+4r2) X_0 X_1 + r1 X_0 X_2 + r2 X_0 X_3$\\
	\hti{4}$-r2 X_1 X_2 + 3r1 X_3 X_1 + (5r1+3r2) X_2 X_3 = 0.$

	\sssec{Case 1.}
	$p = 13,$ $Barycenter = 1504,$ $n = 1,3,6,7,10,12,$ $m = 1,2,4,5,$\\
	The tetrahedron is orthogonal.\\
	The $Center$ is an ideal point.\\
	$Barycenter = (1504)$\\
	${\cal I}deal = \{578\}$\\
	$A = (183,14,1,0)$\\
	$a = [30941,2380,183,14,1,0]$\\
	$Center = (2368)$

	\sssec{Case 2.}
	$p = 13,$ $Barycenter = 366,$ $n = 1,5,2,6,4,10,$\\
	The $Center$ is an ideal point.\\
	$Barycenter = (366)$\\
	${\cal I}deal = \{366\}$\\
	$A = (183,14,1,0)$\\
	$a = [30941,2380,183,14,1,0]$\\
	$Center = (2248)$

	\sssec{Case 3.}
	$p = 13,$ $Barycenter = 366,$ $n = 1,4,10,9,2,5,$ $m = 1,8,4,2,$\\
	The tetrahedron is orthogonal.\\
	$Barycenter = (366)$\\
	${\cal I}deal = \{366\}$\\
	$A = (183,14,1,0)$\\
	$a = [30941,2380,183,14,1,0]$\\
	$Center = (94)$\\
	$Pole of {\cal A} = (1156,2225,589,942)$\\
	$mediatrix = [208070,207763,208045,207684]$\\
	$IC = (1156,1251,1563,1587)$\\
	$altitude = [252838,32488,3743,252]$\\
	$Foot = (115,237,1537,1587)$\\
	$ipa = [269114,110205,242016,374208,58939,30209]$\\
	${\cal P}erp. = \{12,23,157,189,222,1535\}$\\
	\hti{4}$ = \{12,23,157,189,222,1535\}$\\
	$Face foot = (8,24,105,185,235,1535)$\\
	\hti{4}$ = (8,24,105,185,235,1535)$\\
	$face altitude = [17759,52911,230868,30967,187,3732]$\\
	\hti{4}$ = [92,17,7,2391,32462,248]$\\
	$Orthocenters = (115,237,1537,1587)$\\
	$Mid = (115,237,1537,1587)$\\
	$mid = [252838,32488,3743,-1]$\\
	$Orthocenter = (1589)$\\
	${\cal C}oideal = \{1165\}$\\
	$Cocenter = (299)$\\
	$Barycenter $(check)$ = (366)$\\
	$Hyperboloid:\:  6  1  0  0  10  6,\:     3  0  1  9  0  3$

	\sssec{Case 4.}
	$p = 13,$ $Barycenter = 366,$ $n = 1,4,4,9,2,5,$\\
	$Barycenter = (366)$\\
	${\cal I}deal = \{366\}$\\
	$A = (183,14,1,0)$\\
	$a = [30941,2380,183,14,1,0]$\\
	$Center = (1833)$\\
	$Pole of {\cal A} = (281,150,2341,527)$\\
	$mediatrix = [207817,330382,309075,201913]$\\
	$IC = (1504,2092,1312,1587)$\\
	$altitude = [237459,32774,3396,252]$\\
	$Foot = (108,233,1208,1587)$\\
	$ipa = [269480,110754,242016,375855,163118,348392]$\\
	${\cal P}erp. = \{63,110,157,190,2133,1540\}$\\
	\hti{4}$ = \{330,2215,157,190,227,1847\}$\\
	$Face foot = (2,19,105,194,326,1535)$\\
	\hti{4}$ = (8,24,105,194,235,1873)$\\
	$face altitude = [4577,41926,230868,31084,194,3732]$\\
	\hti{4}$ = [92,17,7,2382,32462,222]$\\
	$Orthocenters = (115,337,1884,1587)$\\
	$Mid = (112,194,194,1587)$\\
	$mid = [246391,173899,316701,-1]$\\
	$Center of hyperb. = (194)$\\
	${\cal C}oideal = \{-2\}$\\
	$Cocenter = (-2)$\\
	$Barycenter $(check)$ = (366)$\\
	$Hyperboloid:\:  9  8  0  0  11  4$

	The coordinates of $Center$ are $(1,9,9,12),$\\
	those of $IC$ are $(1,7,10,8), (1,11,3,11), (1,6,8,11), (1,8,4,0),$\\
	those of $Foot$ are $(0,1,7,3), (1,0,3,11), (1,6,0,11), (1,8,4,0),$\\
	those of $orthocenters$ are $(0,1,7,10),(1,0,11,11),(1,10,0,11),
	(1,8,4,0),$\\
	those of the $Center of Hyperboloid$ are $(1,0,0,11),$\\
	the hyperboloid is\\
	\hti{4}    $X_0 X_1 - 2 X_0 X_2 - 6 X_3 X_1 - X_2 X_3 = 0.$

	\sssec{Case 5.}
	$p = 17,$ $Barycenter = 614,$ $n = 1,2,5,4,11,10,$\\
	$Barycenter = (614)$\\
	${\cal I}deal = {614}$\\
	$A = (307,18,1,0)$\\
	$a = [88741,5220,307,18,1,0]$\\
	$Center = (2954)$\\
	$Pole of {\cal A} = (3872,3441,1230,2256)$\\
	$mediatrix = [174494,1146461,1279312,787119]$\\
	$IC = (4484,2739,4436,1700)$\\
	$altitude = [904299,91648,9268,541]$\\
	$Foot = (184,427,4368,1684)$\\
	$ipa = [170402,88127,1190767,798046,722183,136312]$\\
	${\cal P}erp. = \{121,198,91,548,3911,1187\}$\\
	\hti{4}$ = \{511,2041,3605,1762,407,4829\}$\\
	$Face foot = (12,23,86,318,341,3486)$\\
	\hti{4}$ = (0,0,239,312,477,2908)$\\
	$face altitude = [59263,113306,422825,88928,309,8399]$\\
	\hti{4}$ = [1,0,5,5232,90764,443]$\\
	$Orthocenters = (0,346,2919,3656)$\\
	$Mid = (173,7,29,2687)$\\
	$mid = [850281,37068,146789,323376]$\\
	$Center of hyperb. = (2687)$\\
	${\cal C}oideal = \{-2\}$\\
	$Cocenter = (-2)$\\
	$Barycenter $(check)$ = (614)$\\
	$Hyperboloid:\:  4  2  14  9  1  13$

	The coordinates of the $Center$ are $(1,9,2,12),$\\
	those of $IC$ are $(1,14,7,12), (1,8,7,1), (1,14,4,15), (1,4,13,16),$\\
	those of $Foot$ are $(0,1,9,13), (1,0,7,1), (1,14,0,15), (1,4,13,0),$\\
	those of $orthocenters$ are $(0,0,0,1), (1,0,2,5), (1,9,0,11), 
	(1,11,10,0),$\\
	those of the $Center of Hyperboloid$ are $(1,8,4,0),$\\
	the hyperboloid is\\
	\hti{4}$2 X_0 X_1 + X_0 X_2 + 7 X_0 X_3 
	-4 X_1 X_2 - 8 X_3 X_1 - 2 X_2 X_3 = 0.$

	\sssec{Case 6.}
	$p = 17,$ $Barycenter = 614,$ $n = 1,5,2,11,4,10,$\\
	$Barycenter = (614)$\\
	${\cal I}deal = {614}$\\
	$A = (307,18,1,0)$\\
	$a = [88741,5220,307,18,1,0]$\\
	$Center = (3114)$\\
	$Pole of {\cal A} = (3984,3313,2240,1262)$\\
	$mediatrix = [95281,167230,631754,493841]$\\
	$IC = (4564,2643,1748,4612)$\\
	$altitude = [1218731,93484,6380,373]$\\
	$Foot = (248,331,1476,4608)$\\
	$ipa = [1427515,407150,88127,875734,799871,561605]$\\
	${\cal P}erp. = \{41,107,198,567,1922,4653\}$\\
	\hti{4}$ = \{319,3493,2041,3783,372,1395\}$\\
	$Face foot = (15,22,103,317,392,2908)$\\
	\hti{4}$ = (1,31,1,309,494,3486)$\\
	$face altitude = [74002,108393,506346,88911,312,7821]$\\
	\hti{4}$ = [18,22,0,5235,90475,409]$\\
	$Orthocenters = (1,394,3496,3095)$\\
	$Mid = (61,4,2623,205)$\\
	$mid = [304633,22454,1185019,18414]$\\
	$Center of hyperb. = (2623)$\\
	${\cal C}oideal = \{-2\}$\\
	$Cocenter = (-2)$\\
	$Barycenter $(check)$ = (614)$\\
	$Hyperboloid\:  3  2  10  5  11  14$

	The coordinates of $Center$ are $(1,9,12,2),$\\
	those of $Center of hyperboloid$ are $(1,8,0,4).$\\
	Observe that the last 2 coordinates are exchanged.

	\sssec{Case 7.}
	$p = 17,$ $Barycenter = 614,$ $n = 1,2,5,4,6,10,$\\
	THE QUADRIC IS DEGENERATE

	\sssec{Case 8.}
	$p = 13,$ $Barycenter = 1165,$ $n = 1,4,10,9,2,5,$ $m = 1,8,4,2,$\\
	The tetrahedron is orthogonal.\\
	The $Center$ is an ideal point.\\
	$Barycenter = (1165)$\\
	${\cal I}deal = \{1589\}$\\
	$A = (183,14,1,0)$\\
	$a = [30941,2380,183,14,1,0]$\\
	$Center = (501)$

	\sssec{Case 9.}
	$p = 19,$ $Barycenter = 762,$ $n = 1,5,3,11,15,10,$\\
	$Barycenter = (762)$\\
	${\cal I}deal = \{762\}$\\
	$A = (381,20,1,0)$\\
	$a = [137561,7240,381,20,1,0]$\\
	$Center = (5279)$\\
	$Pole of {\cal A} = (2484,3802,1813,6007)$\\
	$mediatrix = [2260013,884810,642183,1463793]$\\
	$IC = (597,128,4864,6754)$\\
	$altitude = [82689,115,11573,431]$\\
	$Foot = (12,15,4731,6746)$\\
	$ipa = [881367,14830,1375219,884015,884396,386372]$\\
	${\cal P}erp. = \{63,216,149,500,3383,3644\}$\\
	\hti{4}$ = \{738,6887,2775,2190,450,5872\}$\\
	$Face foot = (12,23,77,396,438,1103)$\\
	\hti{4}$ = (15,32,191,395,495,2186)$\\
	$face altitude = [82689,158138,528524,137846,384,7962]$\\
	\hti{4}$ = [115,27,11,7245,142254,647]$\\
	$Orthocenters = (203,452,2201,1217)$\\
	$Mid = (296,541,2557,1882)$\\
	$mid = [233678,882463,2451610,1218391]$\\
	$Center of hyperb. = (2557)$\\
	${\cal C}oideal = \{-2\}$\\
	$Cocenter = (-2)$\\
	$Barycenter $(check)$ = (762)$\\
	$Hyperboloid:\:  10  8  7  0  11  8$

	\sssec{Case 10.}
	$p = 29,$ $Barycenter = 1742,$ $n = 1,5,3,11,4,10,$\\
	$Barycenter = (1742)$\\
	${\cal I}deal = \{1742\}$\\
	$A = (871,30,1,0)$\\
	$a = [732541,25260,871,30,1,0]$\\
	$Center = (7629)$\\
	$Pole of {\cal A} = (13904,3284,4290,19705)$\\
	$mediatrix = [5332964,18533938,781117,7231524]$\\
	$IC = (11875,23216,15207,13667)$\\
	$altitude = [13487988,743909,39576,1283]$\\
	$Foot = (553,1350,15178,13660)$\\
	$ipa = [8254786,6197397,2802094,15763096,16443405,9610961]$\\
	${\cal P}erp. = \{484,433,109,27,11978,11805\}$\\
	\hti{4}$ = \{1317,18560,11630,1716,976,378\}$\\
	$Face foot = (17,40,436,871,1567,17691)$\\
	\hti{4}$ = (26,51,88,878,1422,871)$\\
	$face altitude = [415484,976431,10634475,732541,895,42080]$\\
	\hti{4}$ = [146,38,28,25282,740951,871]$\\
	$Orthocenters = (109,1574,871,18242)$\\
	$Mid = (512,1462,34,24564)$\\
	$mid = [12490883,8959529,853649,15467561]$\\
	$Center of hyperb. = (19403)$\\
	${\cal C}oideal = \{4265\}$\\
	$Cocenter = (759)$\\
	$Barycenter $(check)$ = (1742)$\\
	$Hyperboloid:\:  5  7  14  19  20  25$

	\sssec{Case 11.}
	$p = 29,$ $Barycenter = 19403,$ $n = 1,5,3,11,4,10,$\\
	$Barycenter = (19403)$\\
	${\cal I}deal = \{4265\}$\\
	$A = (871,30,1,0)$\\
	$a = [732541,25260,871,30,1,0]$\\
	$Center = (759)$\\
	$Pole of {\cal A} = (13904,3284,4290,19705)$\\
	$mediatrix = [18523362,18533938,18522525,18526776]$\\
	$IC = (13904,18759,19848,1960)$\\
	$altitude = [21121745,751739,43780,1691]$\\
	$Foot = (866,1098,19384,1944)$\\
	$ipa = [1759206,15692546,12959255,7495237,10224187,21187915]$\\
	${\cal P}erp. = \{399,566,324,1401,8614,8443\}$\\
	\hti{4}$ = \{1194,10989,10499,20218,1650,14791\}$\\
	$Face foot = (17,32,784,889,1567,14327)$\\
	\hti{4}$ = (20,34,523,878,1161,8440)$\\
	$face altitude = [415484,781319,19121847,733063,895,38716]$\\
	\hti{4}$ = [320,55,13,25282,748520,1451]$\\
	$Orthocenters = (527,1574,8458,14617)$\\
	$Mid = (202,901,7622,19518)$\\
	$mid = [4930377,15804147,19711544,3631884]$\\
	$Center of hyperb. = (6873)$\\
	${\cal C}oideal = \{22214\}$\\
	$Cocenter = (17926)$\\
	$Barycenter $(check)$ = (19403)$\\
	$Hyperboloid:\:  21  18  18  12  21  6$




\setcounter{section}{89}
	\section{Answers to problems and miscellaneous notes.}
	\sssec{Theorem.}
	{\em If $(0,p_1,p_2,p_3)$ is on the Euler line $eul$ then}
	\enumb
	\item$      p_1(m_2-m_3) + p_2(m_3-m_1) + p_3(m_1-m_2) = 0,$
	\item{\em$ P \vee IC(0)$ intersects the Euler line $eul$ at\\
	$( (p_1(m_0-m_2) + p_2(m_1-m_0))(m_1+m_2+m_3),\\
	\hth    p_1((m_1-m_2)(m_1+m_2+m_3)+m_2(m_1-m_0)) + p_2m_1(m_0-m_1),\\
	\hth    p_2((m_1-m_2)(m_1+m_2+m_3)+m_1(m_0-m_2)) + p_1m_2(m_2-m_0),\\
	\hth    p_2((m_1-m_3)(m_1+m_2+m_3)+m_1(m_0-m_3))\\
	\hti{12}	 + p_1((m_3-m_2)(m_1+m_2+m_3)+m_2(m_3-m_0)) )$,\\
	or more symmetrically,\\
	$   ( p_1(s_1(m_0-m_2)+m_0(m_2-m_0) + p_2(s_1(m_1-m_0)+m_0(m_0-m_1),\\
	\hth	p_1(s_1(m_1-m_2)+m_1(m_2-m_0) + p_2(s_1(m_1-m_1)+m_1(m_0-m_1),\\
	\hth	p_1(s_1(m_2-m_2)+m_2(m_2-m_0) + p_2(s_1(m_1-m_2)+m_2(m_0-m_1),\\
	\hth p_1(s_1(m_3-m_2)+m_3(m_2-m_0) + p_2(s_1(m_1-m_3)+m_3(m_0-m_1) ).$\\
	Moreover, if $p_1 = km_1+s-m_0,$ $p_2 = km_2+s-m_0,$ $p_3 = km_3+s-m_0,
$\\
	then the point on $eul$ is\\
	$( (k-1)m_0+s,(k-1)m_1+s,(k-1)m_2+s,(k-1)m_3+s ).$\\
	In particular,\\
	$M = (m_1+m_2+m_3,m_2+m_3+m_0,m_3+m_0+m_1,m_0+m_1+m_2),$\\
	$P = ( s_1t_0-m_0t_0,s_1(m_1^2-m_2m_3)-m_1t_0,s_1(m_2^2-m_3m_1)-m_2t_0,
$\\
	\hth		$s_1(m_3^2-m_1m_2)-m_3t_0),$\\
		with $t_0 = m_0m_1+m_0m_2+m_0m_3-m_1m_2-m_3m_1-m_2m_3,$\\
	$\overline{O} = ($\\
	$Am = (s_1-3m_0,s_1-3m_1,s_1-3m_2,s_1-3m_3),$\\
	$G = (s_1+2m_0,s_1+2m_1,s_1+2m_2,s_1+2m_3),$\\
	$\overline{A}m = ($\\
	$D_0 =()$\\
	$D_1 = ()$\\
	$D_2 = ()$\\
	$G = ()$\\
	$\overline{G} = ()$}
	\enume
	\sssec{Answer to}
	\vspace{-18pt}\hspace{94pt}{\bf \ref{sec-e}.}\\
	\subsubsection{Answer (partial).}
	\ldots ?  The polar $pp_0 = [-2m_1m_2,m_2(m_0+m_1),m_1(m_2+m_0],$\\
	\ldots ?  The intersection $PP_0 = (0,-m_1(m_2+m_0),m_2(m_0+m_1)),$\\
	\ldots ?  $pp = [m_1m_2(m_2+m_0)(m_0+m_1),m_2m_0(m_0+m_1)(m_1+m_2),
	m_0m_1(m_1+m_2)(m_2+m_0)].$

	Point ``O", intersection of perpendicular to faces through their
	barycenter\\
	(special case of $\ldots$  with $k = 0$ hence\\
	$``O" = (s_1-m_0,s_1-m_1,s_1-m_2,s_1-m_3).$\\
	``Conjugate tetrahedron",\\
	$``A'_0" = (-2s_1-2m_0,s_1+2m_1,s_1+2m_2,s_1+2m_3),$\\
	barycenter of faces of $[A'[i]]$ are\\
	$``M'_0" = (0,m_1,m_2,m_3),$
	which are the orthocenters of the faces,\\
	Perpendiculars through $M'_0$ to the faces, (which are parallel to
		those of [A[]] meet at\\
	$``O'" = (3s_1-4m_0,3s_1-4m_1,3s_1-4m_2,3s_1-4m_3).$\\
	Hence his theorem:
	Then he generalizes the circle of Brianchon-Poncelet and gives its
	center as the midpoint of $H$ and $``O"$\\
	I believe his $``O"$ is my $G.$\\
	Orthocenter,barycenter and $``O"$ are collinear



\chapter{QUATERNIONIAN GEOMETRY}

\setcounter{section}{-1}
	\section{Introduction.}
	\vspace{-18pt}\hspace{108pt}\footnote{22.1.87}\\[8pt]
	It is a classical result, (see Artin, Harsthorne)
	that if the coordinates which are used to define a projective geometry
	are elements of a non commutative division ring, then Desargues' Theorem
	is true, but Pappus' Theorem is, in general, not true.  More precisely,
	Pappus' theorem implies that the division ring or skew fielf is
	commutative.  I will prove here detailled geometric properties which
	justify the definitions of medians and circumcircular
	polarity in a quaternionian plane.

	The results were conjectured by taking the coordinates in the ring
	with unity associated with quaternions over the finite field
	${\Bbb Z}_p,$ $p$ prime.
	This is not a division ring because a finite division ring is a field.
	In this geometry, not all points define a line and vice-versa.
	The situation is similar to that described by Kn\"{u}ppel and Salow,
	for the case of a commutative ring with unity.  This generalization
	merits to be explored in detail.

	In involutive Geometry, we started with the triangle $\{A_i,a_i\},$
	the barycenter $M$ and the orthocenter $\ov{M}$.
	We constructed the medians $ma_i$, the altitudes $\ov{m}a_i,$ the
	mid-points $M_i$, the feet $\ov{M}_i,$ the complementary triangle
	$(M_i,mm_i),$ the orthic triangle $\{\ov{M}_i,\ov{m}m_i\},$
	the ideal points $MA_i$, the orthic points $\ov{M}A_i$,
	the ideal line $m,$ the orthic line $\ov{m},$  the orthic
	directions, $Imm_i$, the tangential triangle $\{T_i,ta_i\},$ the
	symmedians $at_i$ and the point of Lemoine
	$K.$ Moreover the same tangential triangle $T_i$ can be obtained if we
	interchange the role of $M$ and $\ov{M}.$\\
	I attempted the same construction for
	the Geometry over the quaternion skew field. In this case, however,
	the lines $at_i$ are not concurrent, in general, but form a
	triangle $K_i$  and there exists a polarity in which $K_i$ is the pole
	of $a_i$ and therefore $A_i$ is the pole of $at_i.$  This polarity
	degenerates into all the lines through $K$ in the involutive Geometry.

	The plane corresponding to the involutive plane is defined by chosing a
	complete 5-angle in the quaternionian plane. 3 points are the vertices
	of the basic triangle, 1 is the barycenter, 1 is the cobarycenter.
	We define the ideal line as the polar of the
	barycenter with respect to the triangle and as comedians,
	the lines joining the vertices of the triangle to the cobarycenter.
	It can be shown that there is a polarity which
	associates to the vertices of the triangle 3 of the lines through them,
	corresponding to the tangential lines in Euclidean geometry,
	but, in general, the involution defined by this polarity on the ideal 
	line, does not have the directions of the sides and the direction of the
	comedians as corresponding points.

	\section{Quaternionian Geometry over the reals.}

	\ssec{Points, Lines and Polarity.}
	\sssec{Notation.}
	Identifiers, starting with a lower case letter, will denote quaternions,
	$\ov{q}$ denotes the conjugate of $q,$ $q'$, the conjugate
	inverse, $q^n := q\:\ov{q} = \ov{q}\:q,$
	$q^{-n} := (q^n)^{-1}.$

	\sssec{Definition.}
	The elements and incidence in {\em Quaternionian geometry} in 2
	dimensions are defined as follows.
	\enumb
	\item	The {\em points} are $(q_0,q_1,q_2)$ with right equivalence,
	\item	The {\em lines} are $[l_0,l_1,l_2]$ with right equivalence,
	\item	A point $P$ is {\em incident} to a line $l$  iff\\
	\hth$P \cdot l := \sum_{i=0}^2  \ov{P}_i l_i = 0.$
	\enume

	Condition 2 is consistent with equivalence and can also be written\\
	\hth$	l \cdot P = 0.$\\
	I prefered it, because the usual form $\sum_{i=0}^2 P_i l_i = 0$
	implies $\sum_{i=0}^2  \ov{l}_i \ov{P}_i = 0.$

	We remind the reader of the following

	\sssec{Theorem.}
	{\em In any skew field, if a matrix {\bf A} has a left inverse and a
	right inverse, these are equal.}

	Proof: Let {\bf C} be the left inverse of {\bf A} and {\bf B} be its
	right inverse, by associativity of matrices,\\
	\hth${\bf C} = {\bf C}({\bf A}{\bf B}) = ({\bf C}{\bf A}){\bf B}
		= {\bf B}.$

	\sssec{Lemma.}
	\enumb
	\item $p_i\neq q_i,$ $i= 1,2\implies (1,p_1,p_2) \times (1,q_1,q_2)\\
	\hth= [(p'_1-q'_1)^{-1}(p'_1\ov{p}_2-q'_1\ov{q}_2)
		(\ov{p}_2-\ov{q}_2)^{-1},
		(\ov{p}_1-\ov{q}_1)^{-1},-(\ov{p}_2-\ov{q}_2)^{-1}],$
	\item $(1,p_1,p_2) \times (1,p_1,q_2) = [\ov{p}_1,-1,0],$
	\item $(0,0,1)\times (1,p_1,0) = [\ov{p}_1,-1,0]$,
	\item $(0,0,1)\times (0,1,0) = [1,0,0]$.
	\enume

	Proof: Let the line be $[x_0,x_1,x_2]$, we must have\\
	$x_0 +\ov{p}_1 x_1+\ov{p}_2 x_2 = 0,$ and
	$x_0 +\ov{q}_1 x_1+\ov{q}_2 x_2 = 0,$
	subtracting gives\\
	$(\ov{p}_1-\ov{q}_1)+ x_1 (\ov{p}_2-\ov{q}_2) x_2 = 0,$\\
	hence $x_1$ and $x_2.$ $x_0$ follows from substitution into the second 
	equation.

	\sssec{Theorem.}
	{\em A quaternionian geometry is a perspective geometry.}

	\sssec{Lemma.}\label{sec-ldesar}
	{\em Let $a_i$ and $b_i$ be different from 0.
	The points $P_0 := (0,q_1,r_2),$ $P_1 := (r_0,0,q_2),$
	$P_2 := (q_0,r_1,0)$ are collinear iff\\
	\hth$(\ov{q}_2r'_2)(\ov{q}_1r'_1)(\ov{q}_0r'_0) = -1$
	and the line is given by any of the triples}
	\hth$[r'_0\ov{q}_2,q'_1\ov{r}_2,-1], [-1,r'_1\ov{q}_0,q'_2\ov{r}_0],
		[r'_2\ov{q}_1,-1,q'_0\ov{r}_1].$

	Proof:
	Let $y := [y_0,y_1,-1] := P_0 \times P_1$, we have\\
	$\ov{q}_1y_1+\ov{r}_2 = 0$ and $\ov{r}_0y_0+\ov{q}_2 = 0$,
	this gives the first form $y$. To verify that $P_2$ is on $y$, we need
	$q'_0r'_0\ov{q}_2+\ov{r}_1q'_1\ov{r}_2 = 0.$

	\sssec{Lemma.}
	{\em In a quaternionian geometry the Theorem of Desargues is satisfied.}

	Proof:
	We can always choose the coordinates of $A_i$ and $C$ as follows\\
	\hth$A_0 = (1,0,0),$ $A_1 = (0,1,0),$ $A_2 = (0,0,1),$ $C = (1,1,1).$\\
	It follows that\\
	\hth$a_0 = [1,0,0],$ $a_1 = [0,1,0],$ $a_2 = [0,0,1],$\\
	\hth$c_0 = [0,1,-1],$ $c_1 = [-1,0,1],$ $c_2 = [1,-1,0],$\\
	It follows that $B_i$ are\\
	\hth$B_0 = (q_0,1,1),$ $B_1 = (1,q_1,1),$ $B_2 = (1,1,q_2),$\\
	therefore\\
	\hth$b_0 = [0,\ov{q}_1-1,-(\ov{q}_2-1)],$
	    $b_1 = [-(\ov{q}_0-1),\ov{q}_2-1,0],$\\
	\hti{12}    $b_2 = [\ov{q}_0-1,-(\ov{q}_1-1),0],$\\
	\hth$B_0 = (q_0,1,1),$ $B_1 = (1,q_1,1),$ $B_2 = (1,1,q_2),$\\
	\hth$C_0 = (0,q_1-1,q_2-1),$
	    $C_1 = (q_0-1,q_2-1,0),$\\
	\hti{12}    $C_2 = (q_0-1,q_1-1,0),$\\
	the Theorem follows from Lemma \ref{sec-ldesar}.

	\sssec{Theorem.}
	{\em A quaternionian geometry is a Desarguesian geometry.}

	\sssec{Notation.}
	To give an explicit way of indicating that a 3 by 3 matrix is a
	point or line collineation or a correlation from points to lines or
	from lines to point a parenthesis is used on the side of the point
	and a bracket on the side of a line. This notation is used when we
	give the table of elements. It could also be used in connection
	with the bold face letter representing collineation or correlation.
	This notation is only useful when we apply algebra to geometry.

	\sssec{Theorem.}
	{\em If {\bf C} is a point collineation, the line collineation is
	${\bf C}'^T.$\\
	In particular, the point collineation which associates to $A_i$, $A_i$
	and to $(1,1,1),$ $(q_0,q_1,q_2)$ is\\
	\hth$\matb{(}{ccc}q_0&0&0\\0&q_1&0\\0&0&q_2\mate{)},$\\
	and the line collineation is}\\
	\hth$\matb{[}{ccc}q'_0&0&0\\0&q'_1&0\\0&0&q'_2\mate{]}.$

	Proof: If $Q$ is the image of $P$, and $m$ is the image of $l$, we want
\\
	$0 = P\cdot l = \Sigma \ov{P}_i l_i = \Sigma\ov{\bf C}^{-1}_{ij}Q_jl_i
	= \Sigma\ov{Q}_j{\bf C}'_{ij}l_i = \Sigma\ov{Q}_j m_j = 0,$\\
	for all points $P$ and incident lines $l$. This  requires\\
	\hth$m_j = \Sigma {\bf C}'^T_{ji}l_i.$

	\sssec{Definition.}
	\vspace{-18pt}\hspace{91pt}\footnote{13.12.86}\\[8pt]
	A {\em Hermitian matrix} {\bf M} is a matrix which is equal to its
	conjugate transpose.

	{\bf M} defines a transformation from points to line,

	{\bf M}$^{-1},$ the inverse,
	defines the transformation from lines to points.

	\sssec{Theorem.}
	If {\bf M} is Hermitian and $p = {\bf M} P,$ $q = {\bf M} Q$
	and $P \cdot q = 0$ then $Q \cdot p = 0.$

	Proof:  $Q \cdot p = \sum_i \ov{Q}_i p_i\\
	\hth	= \sum_i \sum_j\ov{Q}_i M_{i,j} P_j\\
	\hth	= \sum_j \sum_i (\ov{Q}_i \ov{M}_{j,i}) P_j\\
	\hth	= \sum_j \sum_i  \ov{(M_{j,i} Q_i)} P_j\\
	\hth	= \sum_j \ov{q}_j P_j = q \cdot P = 0.$

	\sssec{Theorem.}
	\enumb
	\item {\em The transformation defined by a Hermitian matrix is a
		polarity}.
	\item {\em The columns of a polarity ${\bf M}$ are the polars of the
		points $A_i$, with\\
		$A_0 = (1,0,0)$, $A_1 = (0,1,0)$, $A_2 = (0,0,1)$.}
	\item {\em The columns of a inverse polarity ${\bf M}^{-1}$ are the
		poles of the lines $a_i$, with $a_0 = [1,0,0]$, $a_1 = [0,1,0]$,
		$a_2 = [0,0,1]$.}
	\enume

	\sssec{Definition.}
	{\em Polar points} are points incident to their polar.
	{\em Polar lines} are lines incident to their pole.

	\sssec{Comment.}
	A polar line can contain infinitely many polar points. For instance,
	let the polarity be $\matb{[}{rrr}1&-1&-1\\-1&1&-1\\-1&-1&1\mate{)}.$
	$P := (0,1,1),$ has for polar $p = [1,0,0].$ A point $Q := (0,1,q),$
	on $p$ has for polar $[1+q,q,1]$. $Q$ is a polar point
	if $\Re(q) = 0.$ Therefore all the points $(0,1,a_1{\bf i} +
	a_2{\bf j} +a_3{\bf k})$ are polar points on [1,0,0]. 
	
	\sssec{Lemma.}
	Let\\
	\hth$	c_i := q_{i-1} q_i q_{i+1} +  \ov{q}_{i+1}
	\ov{q}_i \ov{q}_{i-1} = 2 Re(q_{i-1} q_i q_{i+1}),$\\
	if $q_i \ov{q}_i \neq  0$ then\\
	\hth$	c_0 = c_1 = c_2$\\
	and we define\\
	\hth$	c := c_0.$

	Proof:\\
	$\ov{q}_2 q_2 c_0 =\ov{q}_2 c_0 q_2
	= \ov{q}_2 (q_2 q_0 q_1 + \ov{q}_1  \ov{q}_0 \ov{q}_2) q_2
	= \ov{q}_2 q_2 (q_0 q_1 q_2) + (\ov{q}_2 \ov{q}_1
	 \ov{q}_0 )\ov{q}_2 q_2\\
	\hth= \ov{q}_2 q_2 (q_0 q_1 q_2 + \ov{q}_2 \ov{q}_1 \ov{q}_0 )
	= \ov{q}_2 q_2 c_1.$

	\sssec{Theorem.}\label{sec-thermquat}
	{\em Let\\
	\hth$a_i = \ov{a}_i,$\\
	\hth$b_i := a_{i+1} a_{i-1} - q_i \ov{q}_i,$\\
	\hth$r_i := \ov{q}_{i-1} \ov{q}_{i+1} - a_i q_i,$\\
	then\\
	\hth$\matb{(}{ccc}a_0&q_2&\ov{q}_1\\
		\ov{q}_2&a_1& q_0\\
		q_1& \ov{q}_0&a_2\mate{)}
	\matb{(}{ccc}b_0&r_2&\ov{r}_1\\
		\ov{r}_2&b_1& r_0\\
		r_1& \ov{r}_0&b_2\mate{)} = d\: {\bf E},$\\
	where {\bf E} is the identity matrix and\\
	\hth$	d := a_0 a_1 a_2 - a_0 q_0 \ov{q}_0 - a_1 q_1
	\ov{q}_1 - a_2 q_2 \ov{q}_2 + 2 Re(q_{2} q_0 q_{1}).$\\
	Moreover,\\
	\hth$b_i = \ov{b}_i,$\\
	\hth$a_i := b_{i+1} b_{i-1} - r_i \ov{r}_i,$\\
	\hth$q_i := \ov{r}_{i-1} \ov{r}_{i+1} - b_i r_i,$}

	Proof:  For instance,\\
	$a_0 b_0 + q_2 \ov{r}_2 + \ov{q}_1 r_1
	= a_0 a_1 a_2 - a_0 q_0 \ov{q}_0 + q_2 q_0 q_1
	- q_2 a_2 \ov{q}_2 + \ov{q}_1 \ov{q}_0
	\ov{q}_2 - \ov{q}_1 a_1 q_1 = d$\\
	\hth and\\
	$a_0 r_2 + q_2 b_1 + \ov{q}_1 \ov{r}_0
	= a_0 \ov{q}_1 \ov{q}_0 - a_0 a_2 q_2 + q_2 a_2 a_0
	 - q_2 q_1 \ov{q}_1 + \ov{q}_1 q_1 q_2
	- \ov{q}_1 a_0 \ov{q}_0 = 0.$\\
	The second part of the proof is obtained similarly or follows
	from \ref{sec-ldesar}.1.

	The 2 parts of the following Lemma use different approaches to the
	problem of constructing polarities.

	\sssec{Lemma.}
	\enumb
	\item {\em If all the components of the lines $x_i$ are non zero and the
	$i$-th component of $x_i$ is real, necessary and sufficient conditions
	for $x_i$ to be polars  of $A_i$ are}
	\begin{enumerate}
	\item [0.]$x_{i+1,i-1}^{-1}\ov{x}_{i-1,i+1} = k_i,$ $k_i$ {\em real}.
	\item [1.]$k_0k_1k_2 = 1.$
	\enume
	\item {Let 3 points have coordinates\\
	\hth$	P_0 = (a_0,q_2,q_1^{-1}),$ $P_1 = (q_2^{-1},a_1,q_0),$
		$P_2 = (q_1,q_0^{-1},a_2),$\\
	where $a_i = \ov{a}_i,$ then, if the norm of $q_0q_1q_2 = 1$
	and if the matrix {\bf P} is
	obtained by multiplying the column vectors $P_i$ respectively by
	$1,$ $q_2^n,$ $\frac{1}{q_1^n},$ this matrix defines a polarity which
	associates the lines $a_i$ to the points $P_i.$}
	\enume

	Proof: For part 0, the condition that the $i$-th component of $x_i$ is
	real can always be satisfied by multiplying the components of $x_i$ by
	$\ov{x}_{ii}.$\\
	It remains to find real numbers which multiplied by $x_i$ give the
	columns of an Hermitian matrix. Considering the elements $x_{01}$
	and $x_{10}$ requires $\ov{x}_{10} = x_{01}k_2$, $k_2$ real,
	or more generally 0.
	Multiplying the first and second column by $k_1$ and $k_1^{-1}$
	requires condition 1, for the elements in position 12 and 21 to be
	conjugates of each other.

	For the second part, the matrix is then\\
	\hth$	\matb{(}{ccc}a_0&\ov{q}_2&q'_1\\
		q_2&q_2^na_1&q_2^n\ov{q}_0\\
		q_1^{-1}&q_2^nq_0&q_1^{-n}a_2\mate{)}.$

	\sssec{Exercise.}
	State and prove the Theorem which extends the preceding Theorem to the
	case where some of the components of the vectors $x_i$ are 0, or some
	of the $q_i$ are 0.

	\sssec{Lemma.}\label{sec-lqua}
	\vspace{-18pt}\hspace{80pt}\footnote{15.1.87}\\[8pt]
	$u_j^n\neq 0,$ $v_j^n\neq 0,$ $j = 1,2$, and
	$d_0 := u_1u_2^{-1} + v_1v_2^{-1},$
	$e_0 := u_2u_1^{-1} + v_2v_1^{-1},$\\
	\hth$	\Rightarrow   u_1^{-1}d_0v_2 = u_2^{-1}e_0v_1$ and\\
	\hti{12}$d_0^nv_2^nu_2^n=e_0^nv_1^nu_1^n.$

	Proof: $u_1^{-1}d_0v_2 = u_2^{-1}v_2 + u_1^{-1}v_1$ and
		$u_2^{-1}e_0v_1 = u_1^{-1}v_1 + u_2^{-1}v_2.$

	\sssec{Lemma.}\label{sec-lforTLinv}
	\hth$(u_0 u_1 u_2 v_0 v_1 v_2)^n \neq 0,$\\
	\hth$d_i := u_{i+1}u_{i-1}^{-1} + v_{i+1}v_{i-1}^{-1},$
	$d := d_0d_1d_2,$\\
	\hth$ 	e_i := u_{i-1}u_{i+1}^{-1} + v_{i-1}v_{i+1}^{-1},$
	$e := e_0e_1e_2,$\\
	\hth$	\Rightarrow   d^n = e^n.$

	Proof: This follows from Lemma \ref{sec-lqua}, taking norms and
	using the fact
	that the norm of a product is the product of the norms.

	\ssec{Quaternionian Geometry of the Hexal Complete 5-Angles.}
	\sssec{Notation.}
	In what follows, I will use the same notation as in involutive Geometry,
	namely,\\
	$l := P \times Q,$ means that the line $l$ is defined as the line
	incident to $P$ and $Q.$\\
	If subscripts are used these have the values 0, 1 and 2 and the
	computation is done modulo 3,\\
	$P \cdot l = 0$ means that the point $P$ is incident to the line $l.$\\
	When 3 lines intersect, this intersection can be defined in 3 ways, this
	has been indicated by using (*) after the definition and implies a
	Theorem.\\
	\hth$\sigma := polarity((M_i,a_i)).$\\
	implies that $\sigma$ is the polarity which associates $M_0$ to $a_0$,
	$M_1$ to $a_1$ and $M_2$ to $a_2$.\\
	\hth$m = polar(\sigma ,M).$\\
	implies that in the polarity $\sigma$, $m$ is the polar of $M.$\\
	The labeling used is ``H," for Hypothesis, ``D", for definitions, ``C",
	for conclusions,``N", for nomenclature, ``P", for proofs, this
	labelling being consistent with that of the corresponding definitions.
	The example given is associated to the quaternions over ${\Bbb Z}_{19}$,
	the labelling is ``E" and is consistent with the corresponding
	definitions.\\
	Because any 3 pairs of points and lines do not necessarily define a
	polarity, if a polarity is defined it implies a conclusion (or Theorem)
	I have therefore replaced ``D" by ``DC".

	\sssec{The special configuration of Desargues.}
	With this notation, the special configuration of Desargues
	can be defined by\\
	\hth$	a_i := A_{i+1} \times A_{i-1},$ $qa_i = Q \times A_i,$\\
	\hth$	Q_i := a_i \times qa_i,$ $qq_i := Q_{i+1} \times Q_{i-1},$\\
	\hth$	QA_i := a_i \times qq_i,$ $q_i := A_i \times QA_i,$\\
	\hth$	QQ_i := q_{i+1} \times q_{i-1},$ $q := QA_1 \times QA_2 (*),$\\
	and the other conclusion of the special Desargues Theorem can be
	written,\\
	\hth$QQ_i \cdot qa_i = 0.$\\
	Let $Q$ and $A_i$ be\\
	\hth$	Q = (q_0,q_1,q_2),$ and $A_0 = (1,0,0),$ $A_1 = (0,1,0),$
	$A_2 = (0,0,1),$\\
	then we have the following results, not obtained in the given order,\\
	\hth$	A_0 = (1,0,0),$ $a_0 = [1,0,0],$\\
	\hth$	Q = (q_0,q_1,q_2),$ $q = [q'_0,q'_1,q'_2],$\\
	\hth$	QA_0 = (0,q_1,-q_2),$ $qa_0 = [0,q'_1,-q'_2],$\\
	\hth$	Q_0 = (0,q_1,q_2),$ $q_0 = [0,q'_1,q'_2],$\\
	\hth$	QQ_0 = (-q_0,q_1,q_2),$ $qq_0 = [-q'_0,q'_1,q'_2],$

	The self duality of the configuration corresponds to the replacement
	of points by lines where upper case letters are replaced by lower case
	letters and coordinates by their conjugate inverse.

	\sssec{Fundamental Hypothesis, Definitions and Conclusions.}
	\label{sec-fundHDC}
	The ideal line and the coideal line.

		Given\\
	H0.0.	$A_i,$\\	
	H0.1.	$M,$ $\ov{M},$\\
		Let\\
	D1.0.	   $a_i := A_{i+1} \times A_{i-1},$\\
	D1.1.   $ma_i := M \times A_i,$
		$\ov{m}a_i := \ov{M} \times A_i,$\\
	D1.2.   $M_i := ma_i \times a_i,$
		$\ov{M}_i := \ov{m}a_i \times a_i,$\\
	D1.3.	$eul = M \times \ov{M},$\\
	DC1.4.	$\sigma := polarity((M_i,a_i)),$
		$\ov{\sigma} := polarity((\ov{M}_i,a_i)),$

	\noindent D2.0.   $mm_i := M_{i+1} \times M_{i-1},$
		$\ov{m}m_i := \ov{M}_{i+1} \times \ov{M}_{i-1},$\\
	D2.1.   $MA_i := a_i \times mm_i,$
	$\ov{M}A_i := a_i \times \ov{m}m_i,$\\
	D2.2.	$m_i := A_i \times MA_i,$
		$\ov{m}_i := A_i \times \ov{M}A_i,$\\
	D2.3.	$MM_i := m_{i+1} \times m_{i-1},$
		$\ov{M}M_i := \ov{m}_{i+1} \times
			 \ov{m}_{i-1},$\\
	D2.4.   $m := MA_1 \times MA_2\:(*),$
		$\ov{m} := \ov{M}A_1 \times \ov{M}A_2\:(*),$\\
	D2.5.	$Ima_i := m \times ma_i,$
		$\ov{I}ma_i := \ov{m} \times \ov{m}a_i,$\\
	D2.6.	$IMa_i := m \times \ov{m}a_i,$
		$\ov{I}Ma_i := \ov{m} \times ma_i,$\\
	D2.7.	$iMA_i := M \times MA_i,$
		$\ov{\imath}MA_i := \ov{M} \times \ov{M}A_i,$\\
	then\\
	C2.0.	$m = polar(\sigma ,M),$ $\ov{m} = polar(\ov{\sigma} ,\ov{M}).$\\
	C2.1.	$mm_i = polar(\sigma ,A_i),$
		$\ov{m}m_i = polar(\ov{\sigma},A_i).$\\
	C2.2.	$ma_i = polar(\sigma ,MA_i),$
		$\ov{m}a_i = polar(\ov{\sigma},\ov{M}A_i).$\\
	C2.3.	$iMA_i = polar(\sigma,Ima_i),$
		$\ov{\imath}MA_i = polar(\ov{\sigma},\ov{I}ma_i).$\\
	Let\\
	D3.0.	$mf_i := M_i \times IMa_i,$
		$\ov{m}f_i := \ov{M}_i \times \ov{I}Ma_i,$\\
	D3.1.	$O := mf_1 \times mf_2 (*),$
		$\ov{O} := \ov{m}f_1 \times \ov{m}f_2 (*),$\\
	D3.2.	$Mfa_i := a_{i+1} \times mf_{i-1},$
		$\ov{M}fa_i := a_{i+1} \times \ov{m}f_{i-1},$\\
	\htp{28}$Mf\ov{a}_i := a_{i-1} \times mf_{i+1},$
		$\ov{M}f\ov{a}_i := a_{i-1} \times
			\ov{m}f_{i+1},$\\
	D3.3.	$mfa_i := Mfa_{i+1} \times A_{i-1},$
		$\ov{m}fa_i := \ov{M}fa_{i+1} \times A_{i-1},$\\
	\htp{28}$mf\ov{a}_i := Mf\ov{a}_{i-1} \times A_{i+1},$
		$\ov{m}f\ov{a}_i := \ov{M}f\ov{a}_{i-1}
			\times A_{i+1},$\\
	D3.4.	$Mfm_i := mfa_i \times m_i,$
		$\ov{M}fm_i := \ov{m}fa_i \times \ov{m}_i,$\\
		then\\
	C3.0.	$O \cdot eul = \ov{O} \cdot eul = 0.$\\
	C3.1.	$Mfm_i \cdot mf\ov{a}_i = \ov{M}fm_i \cdot
			\ov{m}f\ov{a}_i = 0.$\\
	Let\\
	D4.0.  $Imm_i := m \times \ov{m}m_i,$
		$\ov{I}mm_i := \ov{m} \times mm_i,$\\
	D4.1.	  $ta_i := A_i \times Imm_i,$\\
	D4.2.	  $T_i := ta_{i+1} \times ta_{i-1},$\\
	D4.3.	  $at_i:= A_i \times T_i,$\\
	D4.4.	  $K_i := at_{i+1} \times at_{i-1},$\\
	D4.5.		$TAa_i	:= ta_i \times a_i,$\\
	D4.6.		$poK_i	:= Taa_{i+1} \times Taa_{i-1},$\\
	DC4.7.	$\theta  := polarity((A_i,ta_i)),$\\
	DC4.8.	$\lambda := polarity((A_i,at_i)),$\\
	then\\
	C4.0.  $\ov{I}mm_i \cdot ta_i = 0.$\\
	C4.1.	$T_i \cdot mf_i = 0.$\\
	C4.2.	$a_i := polar(\theta,T_i).$\\
	C4.3.	$a_i := polar(\lambda,K_i).$

	The nomenclature:\\
	N0.0.	$A_i$	are the {\em vertices} of the triangle,\\
	N0.1.	$M$ is the {\em barycenter},
		$\ov{M}$ is the {\em cobarycenter}.\\
	N1.0.   $a_i$ are the {\em sides}.\\
	N1.1.   $ma_i$ are the {\em medians},
		$\ov{m}a_i$ are the {\em comedians}\\
	N1.2.   $M_i$ are the {\em mid-points} of the sides.
		$\ov{M}_i$ are {\em the feet of the comedians}\\
	N1.3.	eul is the {\em line of Euler},\\
	N1.4.	$\sigma$ is the {\em Steiner polarity}.
		$\ov{\sigma}$ is the {\em co-Steiner polarity}.

	\noindent N2.0.   $\{M_i,mm_i\}$ is the {\em complementary triangle},\\
	\htp{28}	$\{\ov{M}_i,\ov{m}m_i\}$ is the
			{\em orthic triangle},\\
	N2.1.   $MA_i$ are the {\em directions of the sides},\\
	N2.2.	$\{MM_i,m_i\}$	is the {\em anticomplementary triangle}.\\
	N2.3.   $m$ is the {\em ideal line} corresponding to the line at
		infinity,\\
	\htp{28}	$\ov{m}$ is the {\em orthic line} which is the polar of
		$\ov{M}$ with respect to the triangle.

\noindent N2.4.	$Ima_i$	are the {\em directions of the medians}.\\
	\htp{28}	$IMa_i$ are the {\em directions of the comedians}.

\noindent N3.0.	$mf_i$	are the {\em mediatrices},\\
	N3.1.	$O$ is the {\em center},\\
	N3.2.	$Mfm_i$	are the {\em trapezoidal points},

\noindent N4.0.	$Imm_i$ are the {\em directions of the antiparallels of}
		$a_i${\em with respect to the\\
	\hth	sides} $a_{i+1}$ and $a_{i-1}.$\\
		N4.1.  $(T_i,ta_i)$ is the {\em tangential triangle},\\
	N4.2.  $at_i$ are the {\em symmedians},\\
	N4.3.  $K_i$ is the {\em triangle of Lemoine}.\\
	N4.4.	$\theta$ is the {\em circumcircular polarity}\\
	N4.5.	$\lambda$ is the {\em Lemoine polarity}.

	\sssec{Theorem.}
	{\em If we derive a point $X$ and a line $x$ by a given construction
	from $A_i$, $M$ and $\ov{M}$, with the coordinates as given in
	G0.0 and G0.1, below, and the point $\ov{X}$ and line $\ov{x}$
	are obtain by the same construction interchange $M$ and $\ov{M}$,\\
	\hth$X = (f_0(m_0,m_1,m_2),f_1(m_0,m_1,m_2),f_2(m_0,m_1,m_2)),\\
	\hth x = [g_0(m_0,m_1,m_2),g_1(m_0,m_1,m_2),g_2(m_0,m_1,m_2)],\\
	\implies\\
	\ov{X} = (m_0f_0(m_0^{-1},m_1^{-1},m_2^{-1}),
		m_1f_1(m_0^{-1},m_1^{-1},m_2^{-1}),
		m_2f_2(m_0^{-1},m_1^{-1},m_2^{-1})),\\
	\ov{x} = [m'_0g_0(m_0^{-1},m_1^{-1},m_2^{-1}),
		m'_1g_1(m_0^{-1},m_1^{-1},m_2^{-1}),
		m'_2g_2(m_0^{-1},m_1^{-1},m_2^{-1})].$}

	Proof: The point collineation
	$\bf{C} = \matb{(}{ccc}q_0&0&0\\0&q_1&0\\0&0&q_2\mate{)},$
	associates to (1,1,1), $(q_0,q_1,q_2),$ and to $(m_0,m_1,m_2),$
	$(r_0,r_1,r_2),$ if $r_i = q_i m_i.$\\
	In the new system of coordinates, \\
{\small
	$X = (q_0f_0(q_0^{-1}r_0,q_1^{-1}r_1,q_2^{-1}r_2),
		q_1f_1(q_0^{-1}r_0,q_1^{-1}r_1,q_2^{-1}r_2),
		q_2f_2(q_0^{-1}r_0,q_1^{-1}r_1,q_2^{-1}r_2)).$\\
}
	Exchanging $q_i$ and $r_i$ and then replacing $q_i$ by 1 and $r_i$
	by $m_i$ is equivalent to substituting $m_i$ for $q_i$ and 1 for $r_i,$ 
	which gives $\ov{X}$. $\ov{x}$ is obtained similarly.

	The line collineation is\\
	\hth$\matb{[}{ccc}q'_0&0&0\\0&q'_1&0\\0&0&q'_2\mate{]}.$

	\sssec{Exercise.}
	Prove that if a point to line polarity [{\bf P}) has its $i,j$-th
	element\\
	\hth${\bf P}_{ij} = f_(m_0,m_1,m_2),$\\
	then the $i,j$-th element of the polarity obtained by the same 
	construction, after exchange of $M$ and $\ov{M}$, is\\
	\hth$\ov{\bf P}_{ij} = m'_i f(m_0^{-1},m_1^{-1},m_2^{-1})m_j^{-1}.$\\
	Similarly, for a line to point polarity $({\bf P}^{-1}]$\\
	\hth$({\bf P^{-1}}]_{ij} = g_(m_0,m_1,m_2), \implies
	(\ov{\bf P}^{-1}]_{ij} = m_i g(m_0^{-1},m_1^{-1},m_2^{-1})\ov{m}_j.$

	\sssec{Lemma.}
	\hth$m_1^{-1}(m_0+m_1)(m_0-m_1)^{-1} = -
		(m_0^{-1}+m_1^{-1})(m_0^{-1}-m_1^{-1})^{-1}m_1^{-1}.$

	This Lemma is useful in checking equivalent representations
	of coordinates of points and lines.

	\sssec{Notation.}\label{sec-nqua1}
	\hth $r_i := (m_{i-1}^{-1}+m_i^{-1})^{-1}(m_{i+1}^{-1}-m_{i-1}^{-1}),\\
	\hth s_i := -(m_{i-1}^{-1}+m_i^{-1})^{-1}(m_i^{-1}+m_{i+1}^{-1}),\\
	\hth t_i := s_{i+1}^{-1} s_{i-1}^{-1},\\
	\hth f_i := s_i - s_{i+1}^{-1}s_{i-1}^{-1},\\
	\hth g_i := t_i - t_{i+1}^{-1}t_{i-1}^{-1}.$

	\sssec{Lemma.}\label{sec-lforT}
	\enumb
	\item$	s_0s_1s_2 = -1.$
	\item$	norm(t_0t_1t_2) = 1.$
	\item$s'_2\ov{f}_2s_0^{-1} = -f_2s_1.$
	\item$t'_2\ov{f}_2t_0^{-1} = -f_2t_1.$
	\enume

	Proof: For 0, we use Lemma \ref{sec-lforTLinv} and obtain 1, from the
	definition of $t_i$. For 2, we substitute $f_2$ by its definition and
	compare the terms of both sides of the equality which have the same
	sign.

	\sssec{Proof of }
	\vspace{-19pt}\hspace{84pt}{\bf\ref{sec-fundHDC}}.\\[8pt]
	Let\\
	G0.0.	$A_0 = (1,0,0),$ $A_1 = (0,1,0),$ $A_2 = (0,0,1),$\\
	G0.1.	$M = (1,1,1),$ $\ov{M} = (m_0,m_1,m_2),$\\
	then\\
	P1.0.	$a_0 = (1,0,0),$ $a_1 = (0,1,0),$ $a_2 = (0,0,1),$\\
	P1.1.	$ma_0 = [0,1,-1],$ $\ov{m}a_0 = [0,m'_1,-m'_2],$\\
	P1.2.	$M_0 = (0,1,1),$ $\ov{M}_0 = (0,m_1,m_2),$\\
	P1.3.	$eul = [1,(\ov{m}_1-\ov{m}_2)^{-1}(\ov{m}_2-\ov{m}_0),
		(\ov{m}_1-\ov{m}_2)^{-1}(\ov{m}_0-\ov{m}_1)],$\\
	P1.4.	${\bf S} = \matb{[}{ccc}1&-1&-1\\-1&1&-1\\-1&-1&1\mate{)},$
	${\bf S}^{-1} = \matb{(}{ccc}0&1&1\\1& 0&1\\1&1& 0\mate{]}.\\
	\htp{28}{\bf \ov{S}} = \matb{[}{ccc}m_0^{-n}&-m'_0m_1^{-1}&
		-m'_0m_2^{-1}\\-m'_1m_0^{-1}&m^{-n}_1&-m'_1m_2^{-1}\\
		-m'_2m_0^{-1}&-m'_2m_1^{-1}&m_2^{-n}\mate{)},\\
	\htp{28}{\bf \ov{S}}^{\,-1} = \matb{(}{ccc}0&m_0\ov{m}_1&m_0\ov{m}_2\\
		m_1\ov{m}_0& 0& m_1\ov{m}_2\\
		m_2\ov{m}_0&m_2\ov{m}_1& 0\mate{]}.$

\noindent P2.0.	$mm_0 = [1,-1,-1],$ $\ov{m}m_0 = [m'_0,-m'_1,-m'_2],$\\
	P2.1.	$MA_0 = (0,1,-1),$ $\ov{M}A_0 = (0,m_1,-m_2),$\\
	P2.2.	$m_0 = [0,1,1],$ $\ov{m}_0 = [0,m'_1,m'_2],$\\
	P2.3.	$MM_0 = (1,-1,-1),$ $\ov{M}M_0 = (m_0,-m_1,-m_2),$\\
	P2.4.	$m = [1,1,1],$ $\ov{m} = [m'_0,m'_1,m'_2],$\\
	P2.5.	$Ima_0 = (2,-1,-1)$, $\ov{I}ma_0 = (2m_0,-m_1,-m_2),$\\
	P2.6.	$IMa_0 = (m_1+m_2,-m_1,-m_2),$
		$\ov{I}Ma_0 = (m_0(m_1^{-1}+m_2^{-1}),-1,-1),$\\
	P2.7.	$iMA_0 = [2,-1,-1],$
		$\ov{\imath}MA_0 = [2m'_0,-m'_1,-m'_2],$

%
\noindent P3.0.	$mf_0 = [(m_1+m_2)'(\ov{m}_1-\ov{m}_2),1,-1],\\
	\htp{28}\ov{m}f_0 = [m'_0(m_1^{-1}+m_2^{-1})'(m'_1-m'_2),m'_1,-m'_2,1],$
\\
	P3.1.	$O = (m_1+m_2,m_2+m_0,m_0,m_1),\\
	\htp{28}\ov{O} = (m_0(m_1^{-1}+m_2^{-1}),
			m_1(m_2^{-1}+m_0^{-1}),m_2(m_0^{-1}+m_1^{-1})),$\\
	P3.2.	$Mfa_0 = (1,0,-(m_0+m_1)(m_0-m_1)^{-1}),\\
	\htp{28}\ov{M}fa_0 = (m_0,0,
			-m_2(m_0^{-1}+m_1^{-1})(m_0^{-1}-m_1^{-1})^{-1}),\\
	\htp{28}Mf\ov{a}_0 = (1,m'_2(m_2+m_0)(m_2-m_0)^{-1},0),\\
	\htp{28}\ov{M}f\ov{a}_0 = (m_0,
			m_1(m_2^{-1}+m_0^{-1})(m_2^{-1}-m_0^{-1})^{-1},0),$\\
	P3.3.	$mfa_0 = [(m_1+m_2)'(\ov{m}_1-\ov{m}_2),1,0],\\
	\htp{28}\ov{m}fa_0 = [m'_0(m_1^{-1}+m_2^{-1})'(m'_1-m'_2),m'_1,0],\\
	\htp{28}mf\ov{a}_0 = [(m_1+m_2)'(\ov{m}_1-\ov{m}_2),0,-1],\\
	\htp{28}\ov{m}f\ov{a}_0
			= [m'_0(m_1^{-1}+m_2^{-1})'(m'_1-m'_2),0,-m'_2],$\\
	P3.4.	$Mfm_0 = ((m_1+m_2)(m_1-m_2)^{-1},-1,1),\\
	\htp{28}\ov{M}fm_0 = (m_0(m_1^{-1}+m_2^{-1})(m_1^{-1}-m_2^{-1})^{-1},
			-m_1,m_2),$

\noindent P4.0.	$Imm_0 = (r_0,1,s_0),$ $\ov{I}mm_0 = (-r_0,1,s_0),$\\
	P4.1.	$ta_0 = [0,1,-s'_0],$\\
	P4.2.	$T_0 = (1,s_2,s_1^{-1}),$\\
	P4.3.	$at_0 = [0,s'_2,-\ov{s}_1] = [0,1,-t'_0],$\\
	P4.4.	$K_0 = (1,t_2,t_1^{-1}),$\\
	P4.5.	$Taa_0 = (0,1,s_0),$\\
	P4.6.	$poK_0 = [-1,s'_2,\ov{s}_1],$

\noindent P4.7.	${\bf T} = \matb{[}{ccc}0&\ov{f}_2&-\ov{f}_2s_0^{-1}\\
			f_2&0&-f_2s_1\\-s'_0f_2&-\ov{s}_1\ov{f}_2&0\mate{)},$
		${\bf T}^{-1} = \matb{(}{ccc}1&\ov{s}_2&s'_1\\
			s_2&s^n_2&s^n_2\ov{s}_0\\
			s^{-1}_1&s^n_2s_0&s^{-n}_1\mate{]}.$\\
	P4.8.	${\bf L} = \matb{[}{ccc}0&\ov{g}_2&-\ov{g}_2t_0^{-1}\\
			g_2&0&-g_2t_1\\
			-t'_0g_2&-\ov{t}_1\ov{g}_2&0\mate{)},$
		${\bf L}^{-1} = \matb{(}{ccc}1&\ov{t}_2&t'_1\\
			t_2&t_2^n&t_2^n\ov{t}_0\\
			t_1&t_1^{-n}t'_0&t_1^{-n}\mate{]}.$

	Details of proof:\\
	For P4.0, if the coordinates of $Imm_0$ are $x_0$, 1 and $x_2$,
	we have to solve\\
	\hth$x_0 + 1+x_2 = 0,\\
	\hth -m_0^{-1}x_0 + m_1^{-1}+m_2^{-1}x_2 = 0.$\\
	Multiplying the equations to the left respectively by $m_2^{-1}$ and -1,
	or by $m_0^{-1}$ and 1 and adding gives $x_0$ and $x_2$ using the
	notation \ref{sec-nqua1}.\\
	For P4.7, it is easier to obtain ${\bf T}^{-1}$ first, the columns are
	$T_0,$ $T_1$, $T_2,$ multiplied to the right by 1, $s_2^n,$ $s_1^{-n}.$
	The matrix ${\bf T}$ is then obtained using Theorem \ref{sec-thermquat},
	multiplying by $-s_1^{-n}.$  The equivalence with the matrix whose
	columns are $ta_i$ can be verified using Lemma \ref{sec-lforT}.2.
	A similar proof gives P4.8. It is trivialize by the notationb used for 
	$t$.

	\sssec{Theorem.}
	{\em The product of the diagonal elements of ${\bf T}^{-1}$ and of
	${\bf L}^{-1}$ is the same.}

	This follows from Lemma \ref{sec-lforTLinv}.

	\sssec{Exercise.}
	Prove that the center of the circumcircular polarity is\\
	\hth$(m'_0(m_1^{-1}+m_2^{-1}),m'_1(m_2^{-1}+m_0^{-1}),
		m'_2(m_0^{-1}+m_1^{-1})).$\\
	Therefore, in general, it is distinct from $O$. From this follows,
	that, in general, $mf_i$ is not the polar of $MA_i$ in the
	circumcircular polarity.

%
\newpage
	\section{Finite Quaternionian Geometry.}
	\ssec{Finite Quaternions.}
	\sssec{Definition.}
	{\em Finite Quaternions} over ${\Bbb Z}_p$ are associative elements of
	the form\\
	\hth$q_0 + q_1{\bf i} + q_2{\bf j} + q_3{\bf k},$\\
	where $q_i$ are elements of ${\Bbb Z}_p$ and {\bf i}, {\bf j}, {\bf k}
	are such that\\
	\hth${\bf i}^2 = {\bf j}^2 = -1$ and
	${\bf k} = {\bf i}{\bf j} = -{\bf j}{\bf i}.$

	\sssec{Theorem.}
	{\bf i}, {\bf j}, {\bf k} {\em satisfy}\\
	\hth${\bf k}^2 = -1,$ ${\bf i} = {\bf j}{\bf k} = -{\bf k}{\bf j},$
	${\bf j} = {\bf k}{\bf i} = -{\bf i}{\bf k}.$

	This follows at once from associativity.

	\sssec{Theorem.}
	{\em Finite quaternions in ${\Bbb Z}_p$ can be represented by
	2 by 2 matrices over ${\Bbb Z}_p$.\\
	In particular, if $j_0^2 + j_1^2 = -1$, then we can represent\\
	$1$ by $\ma{1}{0}{0}{1},$ {\bf i} by $\ma{0}{1}{-1}{0},$
	{\bf j} by $\ma{j_0}{j_1}{j_1}{-j_0},$
	{\bf k} by $\ma{j_1}{-j_0}{-j_0}{-j_1}.$}

	\sssec{Comment.}
	If $p \equiv 1 \pmod{4},$ we can find an interger $j_0$ such that
	$j_0^2 = -1$, and choose $j_1 = 0$.

	\sssec{Example.}
	\enumb
	\item $p = 5$, we can represent\\
	$1$ by $\ma{1}{0}{0}{1},$ {\bf i} by $\ma{0}{1}{-1}{0},$
	{\bf j} by $\ma{2}{0}{0}{-2},$ {\bf k} by $\ma{0}{-2}{-2}{0}.$
	\item $p = 7$, we can represent\\
	$1$ by $\ma{1}{0}{0}{1},$ {\bf i} by $\ma{0}{1}{-1}{0},$
	{\bf j} by $\ma{-3}{2}{2}{3},$ {\bf k} by $\ma{2}{3}{3}{-2}.$
	\enume

	Finite quaternions will be represented by an integer using the following
	notation.

	\sssec{Notation.}
	In the example the quaternion over ${\Bbb Z}_p,$\\
	\hth	$q = q_0 + q_1 i + q_2 j + q_3 k$\\
	is represented by\\
	\hth	$q = q_0 + q_1 p + q_2 p^2 + q_3 p^3,$ $0 \leq q_i < p.$

	For instance, when $p$ = 19, the representation of $18+3i+6j$ is 2222,
	of $11+10i+4j+3k$ is 22222, of $16+8i+18j$ is 6666 and of
	$14+12i+13j+9k$ is 66666.

	\ssec{Example in a finite quaternionian geometry.}
{\small	Let $p = 19,$\\
	G0.0.	$A_0 = (1,0,0),$ $A_1 = (0,1,0),$ $A_2 = (0,0,1),$\\
	G0.1.	$M = (1,2222,22222),$ $\ov{M} = (1,6666,66666),$\\
	then\\
	E1.0.	$a_0 = [1,0,0],$ $a_1 = [0,1,0],$ $a_2 = [0,0,1],$\\
	E1.1.	$ma_i = $[0,1,13827], [22219,0,1], [1,6378,0],\\
	\htp{26}	$\ov{m}a_i = $[0,1,41987], [66657,0,1], [1,3333,0],\\
	E1.2.	$M_i = $(0,1,48176), (22219,0,1), (1,2222,0),\\
	\htp{26}	$\ov{M}_i = $(0,1,21174), (70903,0,1), (1,6666,0),\\
	E1.3.	$eul = $[1,35222,126587],\\
	E1.4	${\bf S} = \matb{[}{rrr}1&868&115341\\6378&13&82528\\
			22222&48176&18\mate{)},$
		${\bf S}^{-1} = \matb{(}{rrr}0&5034&115341\\2222&0&116883\\
			22222&13827&0\mate{]},$\\
	\htp{26} $\ov{\bf S} = \matb{[}{rrr}1&1835&33434\\5407&14&88952\\
			104133&48615&11\mate{)},$
		$\ov{\bf S}^{-1} = \matb{(}{rrr}0&443&69658\\6789&0&92894\\
			67891&44653&0\mate{]},$\\
	E2.0.	$mm_i = $[1,6378,22222], [2205,1,13827], [22219,48173,1],\\
	\htp{28}	$\ov{m}m_i = $[1,3333,70894], [6653,1,41987],
			[66657,21177,1],\\
	E2.1.	$MA_i = $(0,1,82525), (115341,0,1), (1,5017,0),\\
	\htp{28}	$\ov{M}A_i = $(0,1,116386), (66657,0,1), (1,573,0),\\
	E2.2.	$m_i = $[0,1,116874], [115341,0,1], [1,861,0],\\
	\htp{28}	$\ov{m}_i = $[0,1,95573], [70903,0,1], [1,3906,0],\\
	E2.3.	$MM_i = $(1,5017,115338), (868,1,82525), (115341,116883,1),\\
	\htp{28}	$\ov{M}M_i = $(1,573,70894), (3903,1,116386),
			(66657,95586,1),\\
	E2.4.	$m = $[1,861,115338],
		$\ov{m} = $[1,3906,66666],\\
	E2.5.	$Ima_i = $(1,6128,61279), (624,1,41443), (61290,123602,1),\\
	\htp{28}	$\ov{I}ma_i = $(1,3906,35637), (2132,1,58383),
		(101928,116383,1),\\
	E2.6.	$IMa_i = $(1,31398,82872), (70470,1,745), (35569,2751,1),\\
	\htp{28}	$\ov{I}Ma_i = $(1,84862,112419),
			(112219,1,114203), (4535,17280,1),\\
	E2.7.	$iMA_i = $[2,6378,22222], [2205,2,13827], [22219,48173,2],\\
	\htp{28}	$\ov{\imath}MA_i = $[2,3333,70894], [6653,2,41987],
			[66657,21177,2],\\
	E3.0.	$mf_i = $[1,57399,96485], (22698,1,119282), (59539,116028,1),\\
	\htp{28}$\ov{m}f_i = $[1,17052,63592), (87814,1,43860),
		(49067,45624,1),\\
	E3.1.	$O = $(1,39571,2622),
		$\ov{O} = $(1,26376,18393).\\
	E3.2.	$Mfa_i = $(1,0,59534), (57399,1,0), (0,119282,1),\\
	\htp{28} $Mf\ov{a}_i = $(1,22693,0), (0,1,116019), (96498,0,1),\\
	\htp{28} $\ov{M}fa_i = $(1,0,49068), (17053,1,0), (0,43863,1),\\
	\htp{28} $\ov{M}f\ov{a}_i = $(1,87803,0), (0,1,45633), (63575,0,1),\\
	E3.3.	$mfa_i = $(1,57399,0], [0,1,119282], [59539,0,1],\\
	\htp{28} $mf\ov{a}_i = $(1,0,96485], [22698,1,0], [0,116028,1],\\
	\htp{28} $\ov{m}fa_i = $(1,17052,0], [0,1,43860], [49067,0,1],\\
	\htp{28} $\ov{m}f\ov{a}_i = $(1,0,63592], [87814,1,0], [0,45624,1],\\
	E3.4.	$Mfm_i = $(1,57399,57265), (10093,1,49647), (92996,18940,1),\\
	\htp{28}	$\ov{M}fm_i = $(1,39191,23604), (90214,1,112020),
			(72715,69295,1),\\
	E3.5.	$iMA_i = $[1,6628,76281], [1112,1,7094], [76270,89247,1],\\
	\htp{28}	$\ov{\imath}MA_i = $[1,5096,35637], [3336,1,21174],
		[101928,79188,1],\\
	E4.0.	$Imm_i = $(1,101541,76547), (91854,1,115568),(74057,64703,1),\\
	\htp{28}$\ov{I}mm_i = $(1,36019,60652),(45706,1,21992),
			(63503,72857,1),\\
	E4.1.	$ta_i = $[0,1,19660], [64952,0,1], [1,51999,0],\\
	E4.2.	$T_i = $(1,115899,64951), (51988,1,114948), (39743,19651,1),\\
	E4.3.	$at_i = $[0,1,86571], [100052,0,1], [1,66787,0],\\
	E4.4.	$K_i = $(1,52716,100037), (66802,1,11323), (84095,86576,1),\\
	E4.5.	$Taa_i = $(0,1,114948), (39743,0,1), (1,115899,0),\\
	E4.6.	$poK_i	= $[1,51999,39734], [115882,1,19660], [64952,1149331],\\
	E4.7.	${\bf T} = \matb{[}{rrr}0&126353&46604\\
			10833&0&10388\\90969&127181&0\mate{)},$
		${\bf T}^{-1} = \matb{(}{rrr}1&21317&72608\\
			115899&14&32443\\64951&105118& 2\mate{]},$\\
	E4.8.	${\bf L} = \matb{[}{rrr}0&66802&84095\\
			70412&0&66737\\53448&70833&0\mate{)},$
		${\bf L}^{-1} = \matb{(}{rrr}1&84484&37508\\
			52716&14&35592\\100037&101959&2\mate{]}.$
}

	Except for interchanges the computation of $\times$ is done as follows\\
	we normalize $l_2$ to 1\\
	$\ov{l}_0.P_0 + \ov{l}_1.P_1 + P_2 = 0$,\\
	$\ov{l}_0.Q_0 + \ov{l}_1.Q_1 + Q_2 = 0$,\\
	Multiplying the first by $P_0^{-1}.Q_0$ to the right and subtract from
	the second equation gives\\
	$\ov{l}_1(Q_1-P_1P_0^{-1}Q_0) + (Q_2-P_2P_0^{-1}Q_0) = 0$,\\
	therefore if $r_3 := (Q_1-P_1P_0^{-1}Q_0)^-1$
	and $r_4 = -(Q_2-P_2P_0^{-1}Q_0)$, then\\
	$l_1 = $ the conjugate of $r_4.r_3$ and\\
	$l_0 = -$ the conjugate of $ \ov(l)_1p_1+p_2)p_0^{-1}.$

	The interchange is done as follows\\
	if $P_0 = 0$, then we exchange $P_i$ and $Q_i$,\\
	if after exchange, $P_0 = 0$, we consider all permutations sub0,
	sub1, sub2, of the subscripts 0, 1 and 2 .

	Correspondance in $Z_{19}$ between representation and quaternion.\\
	$\begin{array}{rcrrrrc}
	representation&\vline&r&i&j&k&for\\
	\cline{1-7}
	2222&\vline&18&2&6&0&M\\
	22222&\vline&11&10&4&3\\
	6666&\vline&16*8&18&0&\ov{M}\\
	66666&\vline&14&12&13&9\\
	35222&\vline&15&10&2&5&eul\\
	126587&\vline&9&12&8&18\\
	&\vline&&&&\\
	\end{array}$

	\section{Miniquaternionian Plane $\Psi$ of Veblen-Wedderburn.}
\setcounter{subsection}{-1}
	\ssec{Introduction.}
	Starting with the work of L. E. Dickson of 1905, non-Desarguesian planes
	of order 9 were discovered by Veblen
	and Wedderburn in 1907, I will here consider only one of these
	which is self dual, and for which non trivial polarities exists,
	and refer to the work of G. Zappa (1957),
	T. G. Ostrom (1964), D. R. Hughes (1957) and T. G. Room and
	P. B. Kirkpatrick (1971) for further reading.

	The synthetic definition used can be traced to Veblen and Wedderburn,
	who first consider points obtained by applying a transformation (see p.
	383), later generalized by J. Singer.  The notation is
	inspired by Room and Kirkpatrick (see Table 5.5.4) using the same
	method I used for the finite plane reversing the indices for lines.\\
	An alternate definition, (5.6.1), is given by Room and Kirkpatrick.

	\ssec{Miniquaternion near-field.}

	\sssec{Definition.}
	A {\em near-field} (${\Bbb N},+,\circ)$ is a set ${\Bbb N}$
	with binary operations such that
	\enumb
	\item ${\Bbb N}$ is finite,
	\item $({\Bbb N},+)$ {\em is an Abelian group, with neutral element} 0,
	\item $({\Bbb N}-\{0\},\circ)$ {\em is an group, with neutral element}
		1,
	\item $\circ$ is right distributive over $+$, or\\
	\hth$(\xi+\eta)\circ \zeta
			= \xi\circ \zeta+\eta\circ \zeta,$
		 for all $\xi,\eta,\zeta\in{\Bbb N}$
	\item $\xi\circ 0 = 0$, for all $\xi\in{\Bbb N}$.
	\enume

	\sssec{Theorem.}
	{\em In any near-field,}
	\enumb
	\item $0\circ\xi = 0,$ for all $\xi\in{\Bbb N}$.
	\item $\xi\circ\eta = 0\implies \xi = 0$ {\em or} $\eta = 0.$
	\item $1,-1\neq 0.$
	\enume

	\sssec{Theorem.}\label{sec-tnearfield9}
	{\em In any near-field of order 9,}
	\enumb
	\item $\{0,1,-1\}\approx {\Bbb Z}_3.$
	\item $\xi+\xi+\xi = 0,$ for all $\xi\in{\Bbb Q}_9,$
	\item $-1\circ\xi = \xi\circ(-1)=\xi,$ {\em for all} $\xi\in{\Bbb Q}_9,$
	\item $(-\xi)\circ\eta = \xi\circ(-\eta) = -(\xi\circ\eta),$
		{\em for all} $\xi,\eta\in{\Bbb Q}_9,$
	\item $(-\xi)\circ(-\eta) = \xi\circ\eta,$ {\em for all}
		$\xi,\eta\in{\Bbb Q}_9,$
	\item {\em Given $\kappa\in {\Bbb Q}_9^*,$
		$\lambda = s - \kappa r$ determines a one to one correspondance
		between the elements $\lambda \in {\Bbb Q}_9$ and the pairs
		$(r,s)$, $r,s\in {\Bbb Z}_3$.}
	\item {\em ${\Bbb Q}_9$ being an other near-field of order 9, the
		groups $({\Bbb Q}_9,+)$ and $({\Bbb Q}'_9,+)$ are isomorphic.}
	\item {\em Besides GF($3^2$) there is only one near-field of order 9,
		which is the smallest near-field which is not a field,}
		(Zassenhaus, 1936).
	\enume

	\sssec{Exercise.}\label{sec-enearfield9}
	Determine the correspondance of \ref{sec-tnearfield9}.5.

	\sssec{Definition.}
	The {\em miniquaternion set} ${\Bbb Q}_9 := \{0,\pm 1, \pm\alpha,\pm
	\beta, \pm \gamma\}$ with the operations of addition and multiplications
	defined from,\\
	\hth$  \xi +\xi +\xi = 0$ for all $\xi \in {\Bbb Q}_9$,\\
	\hth$  \alpha - 1 = \beta ,$ $\alpha + 1 = \gamma ,$\\
	\hth$  \alpha^2 = \beta^2 = \gamma^2 = \alpha\beta\gamma = -1$.\\
	The set ${\Bbb Q}_9^* := \{\pm\alpha,\pm\beta, \pm \gamma\}$.

	\sssec{Theorem.}
	\enumb
	\item $\alpha-\beta = \beta-\gamma = \gamma -\alpha = 1,$
		$\alpha+\beta+\gamma = 0.$
	\item $\beta\gamma = -\gamma\beta = \alpha,$
		$\gamma\alpha = -\alpha\gamma = \beta,$
		$\alpha\beta = -\beta\alpha = \gamma.$
	\item {\em the multiplication is right distributive,
		$(\rho+\sigma)\tau = \rho\tau+\sigma \tau,$\\
		for all} $\rho,\sigma,\tau\in {\Bbb Q}_9$.
	\item $\{{\Bbb Q}_9,+,.\}$ {\em is a near-field}.
	\item $\{{\Bbb Q}_9,+,.\}$ {\em is not a field}, e. g.\\
		$\alpha(\alpha+\beta) = \alpha(-\gamma) = \beta$,
		$\alpha\alpha+\alpha\beta = -1+\gamma = \alpha$.
	\item $\begin{array}{rcrrrrrrrr}
	+&\vline&1&-1&\alpha&-\alpha&\beta&-\beta&\gamma&-\gamma\\
	\hline
	1&\vline&-1&0&\gamma&-\beta&\alpha&-\gamma&\beta&-\alpha\\
	-1&\vline&0&1&\beta&-\gamma&\gamma&-\alpha&\alpha&-\beta\\
	\alpha&\vline&\gamma&\beta&-\alpha&0&-\gamma&1&-\beta&-1\\
	-\alpha&\vline&-\beta&-\gamma&0&\alpha&-1&\gamma&1&\beta\\
	\beta&\vline&\alpha&\gamma&-\gamma&-1&-\beta&0&-\alpha&1\\
	-\beta&\vline&-\gamma&-\alpha&1&\gamma&0&\beta&-1&\alpha\\
	\gamma&\vline&\beta&\alpha&-\beta&1&-\alpha&-1&-\gamma&0\\
	-\gamma&\vline&-\alpha&-\beta&-1&\beta&1&\alpha&0&\gamma
	\end{array}\\[10pt]
	\begin{array}{rcrrrrrrrr}
	\cdot&\vline&1&-1&\alpha&-\alpha&\beta&-\beta&\gamma&-\gamma\\
	\hline
	1&\vline&1&-1&\alpha&-\alpha&\beta&-\beta&\gamma&-\gamma\\
	-1&\vline&-1&1&-\alpha&\alpha&-\beta&\beta&-\gamma&\gamma\\
	\alpha&\vline&\alpha&-\alpha&-1&1&\gamma&-\gamma&-\beta&\beta\\
	-\alpha&\vline&-\alpha&\alpha&1&-1&-\gamma&\gamma&\beta&-\beta\\
	\beta&\vline&\beta&-\beta&-\gamma&\gamma&-1&1&\alpha&-\alpha\\
	-\beta&\vline&-\beta&\beta&\gamma&-\gamma&1&-1&-\alpha&\alpha\\
	\gamma&\vline&\gamma&-\gamma&\beta&-\beta&-\alpha&\alpha&-1&1\\
	-\gamma&\vline&-\gamma&\gamma&-\beta&\beta&\alpha&-\alpha&1&-1
	\end{array}$
	\enume

	\ssec{The miniquaternionian plane $\Psi$.}

	\sssec{Definition.}\label{sec-dpsiinc1}
	With $i,i'\in\{0,1,2\},$ $j\in\{0,1,\ldots,12\},$ and the
	addition being performed modulo $3$ for the first element of a pair,
	and modulo $13$, for the second element in the pair or for the element,
	if single, then the elements and incidence in the
	{\em miniquaternionian plane $\Psi$} are defined as follows.
	\enumb
	\item	The {\em points} $P$ are $(j),(i,j), (i',j)$,
	\item	The {\em lines} $l$ are $[j],[i,j], [i',j]$,
	\item	The {\em incidence} is defined by\\
	\hti{4}$[j] := \{(-j),(1-j),(3-j),(9-j),(i,-j),(i',-j)\},$\\
	\hti{4}$[i,j] := \{(-j),(i,2-j),(i,5-j),(i,6-j),(i',3-j),(i',11-j),\\
	\hth		(i'+1,7-j)(i'+1,9-j),(i'-1,1-j),(i'-1,8-j)\},$\\
	\hti{4}$[i',j] := \{(-j),(i',2-j),(i',5-j),(i',6-j),(i,3-j),\\
	\hth		(i,11-j),(i+1,7-j)(i+1,9-j),(i-1,1-j),(i-1,8-j)\}.$
	\enume

	\sssec{Exercise.}\label{sec-epsiinc1}
	\ref{sec-dpsiinc1}.2 is similar to the use of ordered cosets to
	determine efficiently operations of finite as well as infinite groups.
	In this case, $[j]$ is a subplane, $[i,j]$ and $[i',j]$ are 
	copseudoplanes.
	\enumb
	\item Perform a similar representation of points, lines and incidence
	starting with a subplane which is a Fano plane.
	\item Determine similar representations for non Desarguesian geometries
	of order $5^2$, using a subplane of order 4, or of order 5
	$(651 = 31\cdot 21).$
	\item Determine other such representation for non Desarguesian
	geometries of higher order.
	\enume

	\sssec{Theorem.}
	{\em The same incidence relations obtain, if we
	interchange points and lines in} \ref{sec-dpsiinc1}.2.

	\sssec{Theorem. [see Room and Kirkpatrick]}\label{sec-tminipola}
	\enumb
	\item\begin{enumerate}
		\item [0] {\em The correspondance $(j)$ to $[j]$ and
		    $(i,j)$ to $[i,j]$ and $(i',j)$ to $[i',j]$ is a polarity
			${\cal P}_0$ (${\cal J}^*$).}
		\item [1] {\em The $16$ auto-poles are} (0), (7), (8), (11),
			(0,8), (0,12), (1,4), (1,7), (2,10), (2,11),
			(0',8), (0',12), (1',4), (1',7), (2',10), (2',11).
		\enume

	\item\begin{enumerate}
		\item [0] {\em The correspondance $(j)$ to $[j]$ and $(i,j)$ to
			$[i',j]$ and $(i',j)$ to $[i,j]$ is a polarity
			${\cal P}_1$ (${\cal J}'^*$).}
		\item [1] {\em The $22$ auto-poles are} (0), (7), (8), (11),
			(0,1), (0,3), (0,9), (1,1), (1,3), (1,9), (2,1), (2,3),
			(2,9), (0',1), (0',3), (0',9), (1',1), (1',3), (1',9),
			(2',1), (2',3), (2',9).
		\item [2] (0),  (7), (8), (11), (0,1), (1,9), (2,3), (0',9),
			(1',3), (2',1),\\
			(0), (7), (8), (11), (1,1), (2,9), (0,3), (2',9),
			(0',3), (1',1),\\
			(0), (7), (8), (11), (2,1), (0,9), (1,3), (1',9),
			(2',3), (0',1) {\em are ovals.}
		\enume
	\enume

	\sssec{Exercise.}
	\enumb
	\item Prove that the correspondance $(j)$ to $[j]$ and $(i,j)$ to
			$[(i+1)',j]$ and $(i',j)$ to $[i-1,j]$ is a polarity
			${\cal P}_2$.
	\item Prove that the correspondance $(j)$ to $[j]$ and $(i,j)$ to
			$[(i-1)',j]$ and $(i',j)$ to $[i+1,j]$ is a polarity
			${\cal P}_3$.
	\enume

	\sssec{Exercise.}\label{sec-eminipola}
	\enumb
	\item Determine a configuration in \ref{sec-tminipola}.0.2, which gives
		an example were the Theorem of Pascal is satisfied and an other,
		in which it is not satisfied.
	\item Determine ovals which are subsets of \ref{sec-tminipola}.1.1.
	\enume

	\sssec{Theorem.}
	{\em The polar $m$ of a point $M$ with respect to a triangle is
	incident to that point.}

	Indeed, we can always assume thast the triangle consists of the
	real points $A_0 =(0)$, $A_1 =(1)$, $A_2 =(2)$, and that
	$M = (5) = (1,1,1).$ It follows that $m$ = [4] = [1,1,1] which is
	incident to $M$.

	\sssec{Exercise.}
	Check that the other points and lines of the polar construction are
	$M_i = (4), (8), (3),$ $MA_i = (10), (12), (9),$ $MM_i = (7), (6),
	(11),$ $a_i = [12], [1], [0],$ $ma_i = [9], [8], [11],$
	$m_i = [3], [2], [7],$ $mm_i = [6], [10], [5].$

	\sssec{Theorem. [see Room and Kirkpatrick]}\label{sec-tpsifano}
	\enumb
	\item {\em The planes obtained by taking the complete quadrangle
		associated with $3$ real points $A_0,$ $A_1,$ $A_2,$ and a point
		$\ov{M}$ which such that none of the lines $\ov{M}\times A_i$
		are real are Fano planes associated with ${\Bbb Z}_2$.}
	\item {\em There are $(\frac{1}{6}13.12.9).24 = 5616$ Fano planes that
		contain $3$ real points.}
	\enume

	\sssec{Example.}
	The following is a Fano plane (0), (1), (2), (0,3), (2,1), (0',12),
	(1,0), [12], [1], [0], [0',0], [0,12], [2',11], [0',7].

	\sssec{Exercise.}\label{sec-epsifano}
	Determine the Fano plane associated with (0), (1), (2), (0,7).

{\tiny	\sssec{Comment.}
	The correspondance between the notation of Veblen-Wedderburn
	and Room-Kirkpatrick is\\
	$\begin{array}{lcccccccc}
	Veblen-Wedderburn&\vline&a_j&b_j&c_j&d_j&e_j&f_j&g_j\\
	Room-Kirkpatrick&\vline&k_j&a_j&b_j&c_j&a'_j&b'_j&c'_j\\
	De\:Vogelaere&\vline&[-j]&[0,-j]&[1,-j]&[2,-j]&[0',-j]&[1',-j]&[2',-j]\\
	\hline
	Veblen-Wedderburn&\vline&A_j&B_j&C_j&D_j&E_j&F_j&G_j\\
	Room-Kirkpatrick&\vline&K_j&A'_j&C'_j&B'_j&A_j&C_j&B_j\\
	De\:Vogelaere&\vline&(j)&(0,j)&(1,j)&(2,j)&(0',j)&(1',j)&(2',j)
	\end{array}$}

	\sssec{Example. [Veblen-Wedderburn]}
	With the notation\\
	$(C\langle c_0,c_1,c_2\rangle,\{A_0,A_1,A_2\}\{a_0,a_1,a_2\},
	\{B_0,B_1,B_2\}\{b_0,b_1,b_2\};\\
	\{C_0,C_1,C_2\}\{d_0,d_1,d_2\})$, with $d_i := C_{i+j}\times C_{i-j},$\\
	the following configuration shows that the Desargues axiom is not
	satisfied\\
	$((0)\langle [0,0],[1',0],[2,0]\rangle,
	\{(0,1),(1,7),(1',2\}\{[1,1],[0',8],[0',9]\},
	\{(2,3),(0'3),(2,1)\}\{[0,11],[2',7],[10]\};\\
		\{(2',5),(0,{10}),(1',3)\}\{[2,1],[1,0],[0,4]\}).$

	\sssec{Definition.}
	The {\em Singer} matrix ${\bf G} :=
	\matb{(}{ccc}0&0&1\\1&0&1\\0&1&0\mate{)}.$
	Its powers ${\bf G}^k$ are the columns\\
	$k,k+1,k+2$ of\\
	\hth$\begin{array}{ccrrrrrrrrrrrrr}
	k=&\vline&0&1&2&3&4&5&6&7&8&9&10&11&12\\
	\hline
	&\vline&1&0&0&1&0&1&1&1&-1&-1&0&1&-1\\
	&\vline&0&1&0&1&1&1&-1&-1&0&1&-1&1&0\\
	&\vline&0&0&1&0&1&1&1&-1&-1&0&1&-1&1
	\end{array}$

	\sssec{Problem.}
	Can we characterize the plane $\Psi$ using Theorem
	\ref{sec-tpsifano}.0.
\newpage
	move to g6a.tex:

	\sssec{Answer to}
	\vspace{-18pt}\hspace{94pt}{\bf \ref{sec-enearfield9}.}\\[8pt]
	$\begin{array}{rrcrrrrrrrrr}
	\kappa\hti{2}\vline&\lambda=&\vline&0&1&-1&\alpha&-\alpha&\beta&-\beta
		&\gamma&-\gamma\\
	\hline
	\alpha\hti{3}&r&\vline&0&0&0&-1&1&-1&1&-1&1\\
		&s&\vline&0&1&-1&0&0&-1&1&1&-1\\
	\hline
	\beta\hti{3}&r&\vline&0&0&0&-1&1&-1&1&-1&1\\
		&s&\vline&0&1&-1&1&-1&0&0&-1&1\\
	\hline
	\gamma\hti{3}&r&\vline&0&0&0&-1&1&-1&1&-1&1\\
		&s&\vline&0&1&-1&-1&1-1&0&0\\
	\end{array}$

	\sssec{Definition.}
	The elements and incidence in the {\em miniquaternionian plane $\Psi$}
	are defined as follows.
	\enumb
	\item	The {\em points} are $(\xi_0,\xi_1,\xi_2)$ with right
		equivalence,
	\item	
	\item	A point $P$ is {\em incident} to a line $l$  iff\\
	\hth
	\enume

	\sssec{Definition. [Veblen-Wedderburn]}
	The {\em points} $P$ are $(x,y,1),$ $(x,1,0),$ $(1,0,0),$ the
	{\em lines} $l$ are $[1,b,c],$ $[0,1,c],$ $[0,0,1],$ and the
	{\em incidence} is $P\cdot l = 0$.

	\sssec{Theorem. [Veblen-Wedderburn]}
	\enumb
	\item $[1,b,c]\times [1,b',c'] = (-(yb+c),y,1),$ with $y(b-b')=-(c-c').$
	\item $[1,b,c] \times [0,1,c'] = (c'b-c,-c',1),$
	\item $[1,b,c] \times [0,0,1] = (-b,1,0),$
	\item $[0,1,c] \times [0,1,c'] = (1,0,0),$
	\item $[0,1,c] \times [0,0,1] = (1,0,0),$
	\enume

	\sssec{Theorem. [Veblen-Wedderburn]}
	{\em Let} (a(b+c) = ab+ac)
	\enumb
	\item ${\bf M} := \matb{(}{rrr}1&0&1\\-1&0&0\\0&-1&-1\mate{)}.$
	\item $A_0 := (-1,0,1),$ $B_0 := (-\gamma,\alpha,1),$
		$C_0 := (\beta,-\alpha,1),$ $D_0 := (-\beta,\gamma,1),$
		$E_0 := (\alpha,-\gamma,1),$ $F_0 := (\gamma,-\beta,1),$
		$G_0 := (-\alpha,\beta,1),$
	\item $A_j := {\bf M}^j A_0$, $B_j := {\bf M}^j B_0$, \ldots,
		{\em for $j = 0$ to $12$},
	\item $a_0 := [1,1,1],$ $b_0 := [1,\alpha,1]$, $c_0 := [1,-\alpha,1]$,
		$d_0 := [1,\gamma,1]$, $e_0 := [1,-\gamma,1]$,
		$f_0 := [1,-\beta,1]$, 	$g_0 := [1,\beta,1]$,

	{\em then}
	\item $a_0 = \{A_0,A_1,A_3,A_9,B_0,C_0,D_0,E_0,F_0,G_0\},\\
		b_0 = \{A_0,B_1,B_8,D_3,D_{11},E_2,E_5,E_6,G_7,G_9\},\\
		c_0 = \{A_0,C_1,C_8,E_7,E_9,F_3,F_{11},G_2,G_5,G_6\},\\
		d_0 = \{A_0,B_7,B_9,D_1,D_8,F_2,F_5,F_6,G_3,G_{11}\},\\
		e_0 = \{A_0,B_2,B_5,B_6,C_3,C_{11},E_1,E_8,F_7,F_9\},\\
		f_0 = \{A_0,C_7,C_9,D_2,D_5,D_6,E_3,E_{11},F_1,F_8\},\\
		g_0 = \{A_0,B_3,B_{11},C_2,C_5,C_6,D_7,D_9,G_1,G_8\},$
	\item $A_0 = \{a_0,a_4,a_{10},a_{12},b_0,c_0,d_0,e_0,f_0,g_0\},\\
		B_0 = \{a_0,b_5,b_{12},d_4,d_6,e_7,e_8,e_{11},g_2,g_{10}\},\\
		C_0 = \{a_0,c_5,c_{12},e_2,e_{10},f_4,f_6,g_7,g_8,g_{11}\},\\
		D_0 = \{a_0,b_2,b_{10},d_5,d_{12},f_7,f_8,f_{11},g_4,g_6\},\\
		E_0 = \{a_0,b_7,b_8,b_{11},c_4,c_6,e_5,e_{12},f_2,f_{10}\},\\
		F_0 = \{a_0,c_1,c_{10},d_7,d_8,d_{11},e_4,e_6,f_5,f_{12}\},\\
		G_0 = \{a_0,b_4,b_6,c_7,c_8,c_{11},d_2,d_{10},g_5,g_{12}\},$
	\item $X_j \incid x_k \implies X_{j+l\umod{13}}\incid x_{k+l\umod{13}}.$
	\item 
{\footnotesize  ${\bf M}^2 = \matb{(}{rrr}1&-1&0\\-1&0&-1\\1&1&1\mate{)},$
		${\bf M}^3 = \matb{(}{rrr}-1&0&0\\-1&1&0\\0&-1&0\mate{)},$
		${\bf M}^4 = \matb{(}{rrr}-1&-1&1\\1&0&-1\\1&0&0\mate{)},$
		${\bf M}^5 = \matb{(}{rrr}0&-1&1\\1&1&-1\\1&0&1\mate{)},$
		${\bf M}^6 = \matb{(}{rrr}1&-1&-1\\0&1&-1\\1&-1&0\mate{)},$
		${\bf M}^7 = \matb{(}{rrr}-1&1&-1\\-1&1&1\\-1&0&1\mate{)},$
		${\bf M}^8 = \matb{(}{rrr}1&1&0\\1&-1&1\\-1&-1&1\mate{)},$
		${\bf M}^9 = \matb{(}{rrr}0&0&1\\-1&-1&0\\0&-1&1\mate{)},$
		${\bf M}^{10}= \matb{(}{rrr}0&-1&-1\\0&0&-1\\1&-1&-1\mate{)},$
		${\bf M}^{11}= \matb{(}{rrr}1&1&1\\0&1&-1\\-1&1&1\mate{)},$
		${\bf M}^{12}= \matb{(}{rrr}0&-1&0\\-1&-1&-1\\1&1&0\mate{)},$
		${\bf M}^{13}= \matb{(}{rrr}1&0&0\\0&1&0\\0&0&1\mate{)}.$}
	\enume

	Proof: $[x,y,z]\incid \beta ?$ iff
	$x({\bf M}^k_{00}+{\bf M}^k_{20}) + y({\bf M}^k_{01}+{\bf M}^k_{21})
		+ z({\bf M}^k_{02}+{\bf M}^k_{22})
		+ (z{\bf M}^k_{10}+y{\bf M}^k_{11}+z{\bf M}^k_{22})\beta = 0,$
		\ldots.

	\sssec{Example. [Veblen-Wedderburn]}
	With the notation\\
	$(C\langle c_0,c_1,c_2\rangle,\{A_0,A_1,A_2\}\{a_0,a_1,a_2\},
	\{B_0,B_1,B_2\}\{b_0,b_1,b_2\};\\
	\{C_0,C_1,C_2\}\{d_0,d_1,d_2\})$, with $d_i := C_{i+j}\times C_{i-j},$\\
	the following configuration shows that the Desargues axiom is not
	satisfied\\
	$(A_0\langle b_0,f_0,d_0\rangle,\{B_1,C_7,F_2\}\{c_{12},e_5,e_4\},
	\{D_3,E_3,D_1\}\{b_2,g_6,a_3\};\\
		\{G_5,B_{10},F_3\}\{d_{12},c_0,b_9\}).$

	\sssec{Partial answer to}
	\vspace{-18pt}\hspace{118pt}{\bf \ref{sec-epsiinc1}.}\\[8pt]
	For $n = 7^2$, 2451 = 57.43, for $n = 9^2$, 6643 = 91.73,
	for $n = 11^2$, 14763 = 57.259, For $n = 13^2$, 28731 = 3.9577.  

	\sssec{Answer to}
	\vspace{-18pt}\hspace{94pt}{\bf \ref{sec-epsifano}.}\\[8pt]
	The other points are (1,1), (0',12), (0',0), the lines are [12], [1],
	[0], [1',0], [0,12], [0,11] and the polar of (0,7) is [0',6].

	\sssec{Answer to}
	\vspace{-18pt}\hspace{94pt}{\bf \ref{sec-eminipola}.}\\[8pt]
	$(7)\times(8)=[6], [6]\times[0] = (3)$,
	$(8)\times(0)=[1], [1]\times[7] = (2)$,
	$(0)\times(7)=[9], [9]\times[8] = (5)$,
	$\langle (3),(2),(5);[11]\rangle$.\\
	$(7)\times(8)=[6], [6]\times[0,1] = (1',7)$,
	$(8)\times(0,1)=[0,5], [0,5]\times[7] = (0',6)$,
	$(0,1)\times(7)=[2',6], [2',6]\times[8] = (2,5)$,
	$(2,5)$ is not incident to $(1',7)\times(0'6) = [2,2].$\\
	This has not been checked.

	From Dembowski, p. 129

	\sssec{Definition.}
	A linear ternary ring $(\Sigma,+,\cdot)$ is called a {\em cartesian
	field} iff $(\Sigma,+)$ is associative and is therefore a group.

	\sssec{Definition.}
	A cartesian field is called a {\em quasifield} iff the right 
	distributivity law holds:\\
	\hth$(x+y)z = xz+yz.$\\
	Artzy adds that $xa=xb+c$ has a unique solution, but this is a property
	(28). This is Veblen-Wedderburn.

	\sssec{Definition.}
	A quasifield is called a {\em semifield} iff the left
	distributivity law holds:\\
	\hth$z(x+y) = zx+zy.$

	\sssec{Definition.}
	A quasifield is called a {\em nearfield} iff
	$(\Sigma,\cdot)$ is associative and is therefore a group.

	\sssec{Definition.}
	A semifield is called a {\em alternative field} iff
	$x^2y = x(xy)$ and $xy^2 = (xy)y.$

	\sssec{Theorem.}
	{\em $P$ is $(p,L)$ transitive iff $P$ is $(p,L)$ Desarguesian.
	$p$ is point, $L$ is a line.} Dembowski p.123, 16

	Let $Q_0 = (79)$, $Q_1 = (80),$ $Q_2 = (90),$ and $U = (81)$
	then $q_2 = [79],$ $q_0 = [90],$ $q_1 = [80],$ $v = [88],$ 
	$i = [78],$ $V = (78),$ $I = (82),$ $j = [89]$, $W = (89)$,\\
	Points on $q_2:$ 86,12,25,38,51,64,77\\
	Points on $q_1:$ 85,11,24,37,50,63,76\\
	$11\times 12 = [ 7]:   	84(78,86), 8(51,25),43(77,64),48(86,38),
				52(12,12),54(25,78),66(64,80),72(38,51),\\
	 11\times 25 = [61]:    82(78,80), 0(64,64),18(12,51),28(51,78),
				33(77,12),58(38,86),61(86,25),62(25,38),\\
	 11\times 38 = [75]:	81(78,78), 4(38,25),14(12,80),19(77,86),
				36(51,64),70(64,38),73(25,51),74(86,12),\\
	 11\times 51 = [ 4]:	87(86,86), 1(38,80), 2(77,78),46(25,64),
				55(78,38),57(51,12),69(64,51),75(12,25),\\
	 11\times 64 = [ 8]:	83(86,78), 7(64,25),10(77,51),42(78,12),
				47(12,64),53(51,80),65(38,38),71(25,86),\\
	 11\times 77 = [68]:	88(86,80), 5(38,64),13(51,51),21(77,38),
				30(25,12),32(64,86),67(12,78),68(78,25),$\\
	Coordinates of points:\\
	$\begin{array}{rccccccccc}
	( 0)&\vline&64,64&38,80&77,78&78,51&38,25&38,64&51,86&64,25\\
	( 8)&\vline&51,25&86,77&77,51&80,77&12&51,51&12,80&64,78\\
	(16)&\vline&78,77&12,38&12,51&77,86&51,38&77,38&86,64&64,77\\
	(24)&\vline&80,64&25&77,77&25,80&51,78&78,64&25,12&25,77\\
	(32)&\vline&64,86&77,12&64,12&86,51&51,64&80,51&38&25,25\\
	(40)&\vline&77,80&38,78&78,12&77,64&77,25&12,86&25,64&12,64\\
	(48)&\vline&86,38&38,12&80,38&51&12,12&51,80&25,78&78,38\\
	(56)&\vline&51,77&51,12&38,86&12,77&38,77&86,25&25,38&80,25\\
	(64)&\vline&64&38,38&64,80&12,78&78,25&64,51&64,38&25,86\\
	(72)&\vline&38,51&25,51&86,12&12,25&80,12&77&78&80,80\\
	(80)&\vline& 0&78,78&78,80&86,78&78,86&80,86&86&86,86\\
	(88)&\vline&86,80&80,78&\infty\\
	\end{array}$

	Coordinates of lines:\\
	$\begin{array}{rccccccccc}
	[ 0]&\vline&78,38&38&77,78&12,12&51,77&86,38&12,86&12,77\\
	{[} 8]&\vline&64,77&77,25&64,12&80,25&51,80&78,12&12&64,78\\
	{[}16]&\vline&25,25&77,64&86,12&25,86&25,64&51,64&64,38&51,25\\
	{[}24]&\vline&80,38&77,80&78,25&25&51,78&38,38&64,51&86,25\\
	{[}32]&\vline&38,86&38,51&77,51&51,12&77,38&80,12&64,80&78,77\\
	{[}40]&\vline&77&38,78&51,51&12,38&86,77&51,86&51,38&25,38\\
	{[}48]&\vline&38,64&25,51&80,64&12,80&78,51&51&25,78&64,64\\
	{[}56]&\vline&38,25&86,51&64,86&64,25&12,25&25,77&12,64&80,77\\
	{[}64]&\vline&38,80&78,64&64&12,78&77,77&25,12&86,64&77,86\\
	{[}72]&\vline&77,12&38,12&12,51&38,77&80,51&25,80&78,80&\infty\\
	{[}80]&\vline&80&78,78&86&86,78&86,86&80,86&86,80&78,86\\
	{[}88]&\vline&78&80,78&80,80\\
	\end{array}$

	$[80]:11(80,77),24(80,64),37(80,51),50(80,38),63(80,25),76(80,12),
		79(80,80),85(80,86),89(80,78),90(/infty)$
	$B = A+\alpha,$ $(A,B)\incid[V,Y]$, $ V = (78)$,
		$(76)=(80,12)=(\un{0},\alpha),$
		$Y\times V = (78)\times(76)=[13].$\\
	$[13]:78(78),15(64,78),18(12,51),19(77,86),68(78,25),76(80,12),
		46(25,64),48(86,38),53(51,80),60(38,77)$\\
	hence $12 = 80+\alpha,$ $51 =\alpha+\alpha=-\alpha,$
	$25 =78+\alpha=1+\alpha=\gamma,$ $64=25+\gamma=\gamma+\gamma=-\beta,$
	$38=86+\alpha=-1+\alpha=\beta,$ $77=38+\alpha=\beta+\alpha=-\gamma.$\\
	$\begin{array}{rrrrrrrrrr}
	\infty&0&1&-1&\alpha&-\alpha&\beta&-\beta&\gamma&-\gamma\\
	90&80&78&86&12&51&38&64&25&77
	\end{array}$

	$[\alpha,0]=(12)\times(80,80) = (12)\times(79)=[51],$
	$(a,b)\incid[51]\implies b = a\cdot \alpha.$\\
	$(42)=(78,12)=(\un{1},\alpha)\implies \alpha = 1\times\alpha,$\\
	$(45)=(12,86)-(\alpha,\un{-1})\implies -1 = \alpha\times\alpha,$\\
	$(23)=(64,77)=(-\beta,-\gamma) \implies -\gamma=-\beta\times\alpha,$\\
	$(28)=(51,78) = (-\alpha,\un{1})\implies 1 =-\alpha\times\alpha,$\\
	$(35)=(86,51)=(\un{-1},-\alpha)\implies -\alpha= \un{-1}\times\alpha,$

	Using
	DATA 6,0, 6,4, 6,10, 6,12, 0,0, 1,0, 2,0, 3,0, 4,0, 5,0
	DATA 6,0, 0,7, 0,8, 0,11, 3,2, 3,10, 4,4, 4,6, 5,5, 5,12
	gives the same multiplication table give left not right
	distibutive law
	with $Q_i = 79,81,87$, $U = (83),$ $\alpha = (12)$,
	$q_0 = [87] = \{79,81,82,86, 4,17,\ldots\},$\\
	$q_1 = [81] = \{79,85,87,88,10,23,\ldots\},$\\
	$q_2 = [79] = \{81,87,89,90,12,25,\ldots\},$\\
	with case 7, data 79,81,87,83,12:\\
	$\infty = 87,$ $\un{0} = 81,$ $\un{1} = 89,$ $\un{-1} = 90,$
	$\alpha = 12,$ $-\alpha = 77,$ $\beta = 38,$ $-\beta = 51,$
	$\gamma = 25,$ $-\gamma = 64.$

	This is a try for a section to be included in g19.tex between
	Moufang and Desargues.

	\section{Axiomatic.}
	\ssec{Veblen-MacLagan planes.}
	\sssec{Introduction.}
	The first example of a Veblen-Wedderburn plane was given in 1907
	by Veblen and MacLagan-Wedderburn.  It is associated to the
	algebraic structure of a nearfield, which is a skew field
	which lacks the left distributive law, hence is an other plane
	between the Veblen-Wedderburn plane and the Desarguesian plane.

	\sssec{Axiom. [Da]}
	\vspace{-19pt}\hspace{90pt}\footnote{Da for Desargues leading
	to associativity of multiplication.}\\[8pt]
	Given a Veblen-Wedderburn plane, 2 points $Q_1$ and $Q_2$ on the
	ideal line and an other point $Q_0$ not on it, any 2 parallelograms
	$A_i$ and $B_i$ with directions $Q_1$ and $Q_2$, with
	no sides in common \ldots,???,
	such that $A_j$ and $B_j$ are perspective from $Q_0$ for j = 0 To 2,
	imply that $A_3$ and $B_3$ are perspective from $Q_0$.

	\sssec{Notation.}
	Da$(\{Q_0,Q_1,Q_2\},\{A_j\},\{B_j\}).$

	\sssec{Definition.}
	A {\em Veblen-MacLagan plane} is a Veblen-Wedderburn plane in which
	the axiom Da is satisfied.

	\sssec{Lemma. [For Associativity]}
	H1.0.\hti{3}$A_0,$ $a_{12},$ $x,$ (See Fig. 2?.)\\
	D1.0.\hti{3}$a_{01} := Q_1 \times A_0,$ $a_{02} := Q_2 \times A_0,$\\
	D1.1.\hti{3}$A_1 := a_{01} \times a_{12},$
		$A_2 := a_{02} \times a_{12},$\\
	D1.2.\hti{3}$a_{13} := Q_2 \times A_1,$ $a_{23} := Q_1 \times A_2,$
		$A_3 := a_{13} \times a_{23},$\\
	D2.0.\hti{3}$a_0 := Q_0 \times A_0,$ $a_1 := Q_0 \times A_1,$
		$a_2 := Q_0 \times A_2,$ $a_3 := Q_0 \times A_3,$
	D2.1.\hti{3}$B_0 := a_0 \times y,$
		$b_{01} := Q_1 \times B_0,$ $b_{02} := Q_2 \times B_0,$\\
	D2.1.\hti{3}$B_1 := b_1 \times b_{01},$ $B_2 := b_2 \times b_{02},$
	D2.2.\hti{3}$b_{13} := Q_2 \times B_1,$ $b_{23} := Q_1 \times B_2,$
		$B_3 := b_{13} \times b_{23},$\\
	C1.0.\hti{3}$B_3 \incid b_3,$\\
	{\em Moreover}\\
	$A_0 = (A,B),$ $A_1 = (A',B),$ $A_2= (A,B'),$ $A_3 = (A',B'),$
	$B_0 =$\\ 
	Proof:
	Da$(\{Q_0,Q_1,Q_2\},\{A_j\},\{B_j\}).$

	\sssec{Theorem.}
	{\em In a Veblen-MacLagan plane, the ternary ring $(\Sigma,*)$
	is a nearfield.}:
	\enumb
	\item $(\Sigma,+)$ {\em is an Abelian group,}
	\item $(\Sigma-\{0\},\cdot)$ {\em is a group,}
	\item ($\Sigma,*) = (\Sigma,+,\cdot)$ {\em is right distributive,}
		$(a+b)\cdot c = a \cdot c + b\cdot c$.
	\enume

	\ssec{Examples of Perspective planes.}
	\sssec{Theorem.}
	\enumb
	\item {\em The Cayleyian plane is not a Veblen-MacLagan plane.}
		\marginpar{replace Desarg.?}
	\enume

	\sssec{Definition.}
	A {\em miniquaternion plane} \ldots.

	\sssec{Theorem.}
	\enumb
	\item {\em A miniquaternion plane is a Veblen-MacLagan plane.}
	\item {\em A miniquaternion plane is not a Moufang plane.}
	\enume

	\sssec{Tables.}
	The following are in an alternate notation the known table for $p$ = 3
	and a new table for $p$ = 5.
	The other incidence are obtained by adding one to the subscripts of
	the lines and subtracting one for the subscript of the points.

	{\footnotesize Selectors for $\Psi$ plane, when p = 3:\\
	$(0^0_0):$ $[0^0_2],$ $[0^0_5],$ $[0^0_6],$
		$[1^0_1],$ $[1^0_8],$ $[1^1_7],$ $[1^1_9],$
		$[1^2_3],$ $[1^3_{11}],$\\
	$(0^1_0):$ $[0^1_2],$ $[0^1_5],$ $[0^1_6],$
		$[1^1_3],$ $[1^1_{11}],$ 
		$[1^2_1],$ $[1^3_{8}],$ $[1^0_7],$ $[1^0_9],$\\
	$(0^2_0):$ $[0^2_2],$ $[0^2_5],$ $[0^2_6],$
		$[1^2_7],$ $[1^3_9],$
		$[1^0_3],$ $[1^0_{11}],$ $[1^1_1],$ $[1^1_8],$\\
	\noindent$(1^0_0):$ $[1^0_2],$ $[1^0_5],$ $[1^0_6],$
		$[0^0_1],$ $[0^0_8],$ $[0^1_7],$ $[0^1_9],$
		$[0^2_3],$ $[0^3_{11}],$\\
	$(1^1_0):$ $[1^1_2],$ $[1^1_5],$ $[1^1_6],$
		$[0^1_3],$ $[0^1_{11}],$ 
		$[0^2_1],$ $[0^3_{8}],$ $[0^0_7],$ $[0^0_9],$ \\
	$(1^2_0):$ $[1^2_2],$ $[1^2_5],$ $[1^2_6],$
		$[0^2_7],$ $[0^3_9],$
		$[0^0_3],$ $[0^0_{11}],$ $[0^1_1],$ $[0^1_8],$}

	{\footnotesize Selectors for $\Psi$ plane, when p = 5:\\
	$(0^0_0):$ $[0^0_6],$ $[0^0_{21}],$ $[0^0_{16}],$
		$[1^2_{18}],$ $[1^2_{25}],$ $[1^3_5],$ $[1^3_{13}],$
		$[0^3_4],$ $[1^0_{22}],$ $[0^2_{28}],\\
	\hti{6} [2^0_{12}],$ $[2^0_{23}],$ $[2^0_{24}],$ $[2^0_{26}],$
		$[2^1_1],$ $[2^1_{10}],$ $[2^4_{14}],$ $[2^4_8],$
		$[3^1_9],$ $[3^1_{27}],$ $[3^4_{15}],$ $[3^4_{19}],$
		$[2^2_{29}],$ $[2^3_{20}],$ $[3^0_3],\\
	 (0^1_0):$ $[0^1_{25}],$ $[0^1_{19}],$ $[0^1_6],$
		$[1^3_{12}],$ $[1^3_{26}],$ $[1^4_{18}],$ $[1^4_{21}],$
		$[0^4_{22}],$ $[1^1_{28}],$ $[0^3_{27}],\\
	\hti{6} [2^1_3],$ $[2^1_5],$ $[2^1_{29}],$ $[2^1_{13}],$
		$[2^2_{20}],$ $[2^2_9],$ $[2^0_1],$ $[2^0_{10}],$
		$[3^2_{16}],$ $[3^2_{15}],$ $[3^0_{4}],$ $[3^0_8],$
		$[2^3_{14}],$ $[2^4_{24}],$ $[3^1_{23}],\\
	 (0^2_0):$ $[0^2_{26}],$ $[0^2_8],$ $[0^2_{25}],$
		$[1^4_3],$ $[1^4_{13}],$ $[1^0_{12}],$ $[1^0_6],$
		$[0^0_{28}],$ $[1^2_{27}],$ $[0^4_{15}],\\
	\hti{6} [2^2_{23}],$ $[2^2_{18}],$ $[2^2_{14}],$ $[2^2_{21}],$
		$[2^3_{24}],$ $[2^3_{16}],$ $[2^1_{20}],$ $[2^1_9],$
		$[3^3_{19}],$ $[3^3_4],$ $[3^1_{22}],$ $[3^1_{10}],$ $
		[2^4_1],$ $[2^0_{29}],$ $[3^2_5],\\
	 (0^3_0):$ $[0^3_{13}],$ $[0^3_{10}],$ $[0^3_{26}],$
		$[1^0_{23}],$ $[1^0_{21}],$ $[1^1_3],$ $[1^1_{25}],$
		$[0^1_{27}],$ $[1^3_{15}],$ $[0^0_4],\\
	\hti{6} [2^3_5],$ $[2^3_{12}],$ $[2^3_1],$ $[2^3_6],$
		$[2^4_{29}],$ $[2^4_{19}],$ $[2^2_{24}],$ $[2^2_{16}],$
		$[3^4_8],$ $[3^4_{22}],$ $[3^2_{28}],$ $[3^2_9],$
		$[2^0_{20}],$ $[2^1_{14}],$ $[3^3_{18}],\\
	 (0^4_0):$ $[0^4_{21}],$ $[0^4_9],$ $[0^4_{13}],$
		$[1^1_5],$ $[1^1_6],$ $[1^2_{23}],$ $[1^2_{26}],$
		$[0^2_{15}],$ $[1^4_4],$ $[0^1_{22}],\\
	\hti{6} [2^4_{18}],$ $[2^4_3],$ $[2^4_{20}],$ $[2^4_{25}],$
		$[2^0_{14}],$ $[2^0_8],$ $[2^3_{29}],$ $[2^3_{19}],$
		$[3^0_{10}],$ $[3^0_{28}],$ $[3^3_{27}],$ $[3^3_{16}],$
		$[2^1_{24}],$ $[2^2_1],$ $[3^4_{12}],$\\[10pt]
	\noindent$(1^0_0):$ $[1^0_9],$ $[1^0_{19}],$ $[1^0_{26}],$
		$[0^2_6],$ $[0^2_{12}],$ $[0^3_{21}],$ $[0^3_{23}],$
		$[1^1_{18}],$ $[1^4_5],$ $[0^0_{22}],\\
	\hti{6} [3^0_{27}],$ $[3^0_{15}],$ $[3^0_{24}],$ $[3^0_{16}],$
		$[2^1_4],$ $[2^1_8],$ $[2^4_{28}],$ $[2^4_{10}],$
		$[3^3_{20}],$ $[3^3_{25}],$ $[3^2_{29}],$ $[3^2_{13}],$
		$[3^1_1],$ $[3^4_{14}],$ $[2^0_3],\\
	 (1^1_0):$ $[1^1_{16}],$ $[1^1_8],$ $[1^1_{13}],$
		$[0^3_{25}],$ $[0^3_3],$ $[0^4_6],$ $[0^4_5],$
		$[1^2_{12}],$ $[1^0_{18}],$ $[0^2_{28}],\\
	\hti{6} [3^1_{15}],$ $[3^1_4],$ $[3^1_{29}],$ $[3^1_{19}],$
		$[2^2_{22}],$ $[2^2_{10}],$ $[2^0_{27}],$ $[2^0_9],$
		$[3^4_{24}],$ $[3^4_{26}],$ $[3^3_{14}],$ $[3^3_{21}],$
		$[3^2_{20}],$ $[3^0_1],$ $[2^1_{23}],\\
	 (1^2_0):$ $[1^2_{19}],$ $[1^2_{10}],$ $[1^2_{21}],$
		$[0^4_{26}],$ $[0^4_{23}],$ $[0^0_{25}],$ $[0^0_{18}],$
		$[1^3_3],$ $[1^1_{12}],$ $[0^2_{27}],\\
	\hti{6} [3^2_4],$ $[3^2_{22}],$ $[3^2_{14}],$ $[3^2_8],$
		$[2^3_{28}],$ $[2^3_9],$ $[2^1_{15}],$ $[2^1_{16}],$
		$[3^0_{29}],$ $[3^0_{13}],$ $[3^4_1],$ $[3^4_6],$
		$[3^3_{4}],$ $[3^1_{20}],$ $[2^2_5],\\
	 (1^3_0):$ $[1^3_8],$ $[1^3_9],$ $[1^3_6],$
		$[0^0_{13}],$ $[0^0_5],$ $[0^1_{26}],$ $[0^1_{12}],$
		$[1^4_{23}],$ $[1^2_3],$ $[0^3_{15}],\\
	\hti{6} [3^3_{22}],$ $[3^3_{28}],$ $[3^3_1],$ $[3^3_{10}],$
		$[2^4_{27}],$ $[2^4_{16}],$ $[2^2_4],$ $[2^2_{19}],$
		$[3^1_{14}],$ $[3^1_{21}],$ $[3^0_{20}],$ $[3^0_{25}],$
		$[3^4_{29}],$ $[3^2_{24}],$ $[2^3_{18}],\\
	 (1^4_0):$ $[1^4_{10}],$ $[1^4_{16}],$ $[1^4_{25}],$
		$[0^1_{21}],$ $[0^1_{18}],$ $[0^2_{13}],$ $[0^2_3],$
		$[1^0_5],$ $[1^3_{23}],$ $[0^4_4],\\
	\hti{6} [3^4_{28}],$ $[3^4_{27}],$ $[3^4_{20}],$ $[3^4_9],$
		$[2^0_{15}],$ $[2^0_{19}],$ $[2^3_{22}],$ $[2^3_8],$
		$[3^2_1],$ $[3^2_6],$ $[3^1_{24}],$ $[3^1_{26}],$
		$[3^0_{14}],$ $[3^3_{29}],$ $[2^4_{12}],$\\[10pt]
	\noindent$(2^0_0):$ $[2^0_6],$ $[2^0_{21}],$ $[2^0_{16}],$
		$[3^2_{18}],$ $[3^2_{25}],$ $[3^3_5],$ $[3^3_{13}],$
		$[2^3_4],$ $[3^0_{22}],$ $[2^2_{28}],\\
	\hti{6} [0^0_{12}],$ $[0^0_{23}],$ $[0^0_{24}],$ $[0^0_{26}],$
		$[0^1_1],$ $[0^1_{10}],$ $[0^4_{14}],$ $[0^4_8],$
		$[1^1_9],$ $[1^1_{27}],$ $[1^4_{15}],$ $[1^4_{19}],$
		$[0^2_{29}],$ $[0^3_{20}],$ $[1^0_3],\\
	 (2^1_0):$ $[2^1_{25}],$ $[2^1_{19}],$ $[2^1_6],$
		$[3^3_{12}],$ $[3^3_{26}],$ $[3^4_{18}],$ $[3^4_{21}],$
		$[2^4_{22}],$ $[3^1_{28}],$ $[2^3_{27}],\\
	\hti{6} [0^1_3],$ $[0^1_5],$ $[0^1_{29}],$ $[0^1_{13}],$
		$[0^2_{20}],$ $[0^2_9],$ $[0^0_1],$ $[0^0_{10}],$
		$[1^2_{16}],$ $[1^2_{15}],$ $[1^0_{4}],$ $[1^0_8],$
		$[0^3_{14}],$ $[0^4_{24}],$ $[1^1_{23}],\\
	 (2^2_0):$ $[2^2_{26}],$ $[2^2_8],$ $[2^2_{25}],$
		$[3^4_3],$ $[3^4_{13}],$ $[3^0_{12}],$ $[3^0_6],$
		$[2^0_{28}],$ $[3^2_{27}],$ $[2^4_{15}],\\
	\hti{6} [0^2_{23}],$ $[0^2_{18}],$ $[0^2_{14}],$ $[0^2_{21}],$
		$[0^3_{24}],$ $[0^3_{16}],$ $[0^1_{20}],$ $[0^1_9],$
		$[1^3_{19}],$ $[1^3_4],$ $[1^1_{22}],$ $[1^1_{10}],$ $
		[0^4_1],$ $[0^0_{29}],$ $[1^2_5],\\
	 (2^3_0):$ $[2^3_{13}],$ $[2^3_{10}],$ $[2^3_{26}],$
		$[3^0_{23}],$ $[3^0_{21}],$ $[3^1_3],$ $[3^1_{25}],$
		$[2^1_{27}],$ $[3^3_{15}],$ $[2^0_4],\\
	\hti{6} [0^3_5],$ $[0^3_{12}],$ $[0^3_1],$ $[0^3_6],$
		$[0^4_{29}],$ $[0^4_{19}],$ $[0^2_{24}],$ $[0^2_{16}],$
		$[1^4_8],$ $[1^4_{22}],$ $[1^2_{28}],$ $[1^2_9],$
		$[0^0_{20}],$ $[0^1_{14}],$ $[1^3_{18}],\\
	 (2^4_0):$ $[2^4_{21}],$ $[2^4_9],$ $[2^4_{13}],$
		$[3^1_5],$ $[3^1_6],$ $[3^2_{23}],$ $[3^2_{26}],$
		$[2^2_{15}],$ $[3^4_4],$ $[2^1_{22}],\\
	\hti{6} [0^4_{18}],$ $[0^4_3],$ $[0^4_{20}],$ $[0^4_{25}],$
		$[0^0_{14}],$ $[0^0_8],$ $[0^3_{29}],$ $[0^3_{19}],$
		$[1^0_{10}],$ $[1^0_{28}],$ $[1^3_{27}],$ $[1^3_{16}],$
		$[0^1_{24}],$ $[0^2_1],$ $[1^4_{12}],$\\[10pt]
	\noindent$(3^0_0):$ $[3^0_9],$ $[3^0_{19}],$ $[3^0_{26}],$
		$[2^2_6],$ $[2^2_{12}],$ $[2^3_{21}],$ $[2^3_{23}],$
		$[3^1_{18}],$ $[3^4_5],$ $[2^0_{22}],\\
	\hti{6} [1^0_{27}],$ $[1^0_{15}],$ $[1^0_{24}],$ $[1^0_{16}],$
		$[0^1_4],$ $[0^1_8],$ $[0^4_{28}],$ $[0^4_{10}],$
		$[1^3_{20}],$ $[1^3_{25}],$ $[1^2_{29}],$ $[1^2_{13}],$
		$[1^1_1],$ $[1^4_{14}],$ $[0^0_3],\\
	 (3^1_0):$ $[3^1_{16}],$ $[3^1_8],$ $[3^1_{13}],$
		$[2^3_{25}],$ $[2^3_3],$ $[2^4_6],$ $[2^4_5],$
		$[3^2_{12}],$ $[3^0_{18}],$ $[2^2_{28}],\\
	\hti{6} [1^1_{15}],$ $[1^1_4],$ $[1^1_{29}],$ $[1^1_{19}],$
		$[0^2_{22}],$ $[0^2_{10}],$ $[0^0_{27}],$ $[0^0_9],$
		$[1^4_{24}],$ $[1^4_{26}],$ $[1^3_{14}],$ $[1^3_{21}],$
		$[1^2_{20}],$ $[1^0_1],$ $[0^1_{23}],\\
	 (3^2_0):$ $[3^2_{19}],$ $[3^2_{10}],$ $[3^2_{21}],$
		$[2^4_{26}],$ $[2^4_{23}],$ $[2^0_{25}],$ $[2^0_{18}],$
		$[3^3_3],$ $[3^1_{12}],$ $[2^2_{27}],\\
	\hti{6} [1^2_4],$ $[1^2_{22}],$ $[1^2_{14}],$ $[1^2_8],$
		$[0^3_{28}],$ $[0^3_9],$ $[0^1_{15}],$ $[0^1_{16}],$
		$[1^0_{29}],$ $[1^0_{13}],$ $[1^4_1],$ $[1^4_6],$
		$[1^3_{4}],$ $[1^1_{20}],$ $[0^2_5],\\
	 (3^3_0):$ $[3^3_8],$ $[3^3_9],$ $[3^3_6],$
		$[2^0_{13}],$ $[2^0_5],$ $[2^1_{26}],$ $[2^1_{12}],$
		$[3^4_{23}],$ $[3^2_3],$ $[2^3_{15}],\\
	\hti{6} [1^3_{22}],$ $[1^3_{28}],$ $[1^3_1],$ $[1^3_{10}],$
		$[0^4_{27}],$ $[0^4_{16}],$ $[0^2_4],$ $[0^2_{19}],$
		$[1^1_{14}],$ $[1^1_{21}],$ $[1^0_{20}],$ $[1^0_{25}],$
		$[1^4_{29}],$ $[1^2_{24}],$ $[0^3_{18}],\\
	 (3^4_0):$ $[3^4_{10}],$ $[3^4_{16}],$ $[3^4_{25}],$
		$[2^1_{21}],$ $[2^1_{18}],$ $[2^2_{13}],$ $[2^2_3],$
		$[3^0_5],$ $[3^3_{23}],$ $[2^4_4],\\
	\hti{6} [1^4_{28}],$ $[1^4_{27}],$ $[1^4_{20}],$ $[1^4_9],$
		$[0^0_{15}],$ $[0^0_{19}],$ $[0^3_{22}],$ $[0^3_8],$
		$[1^2_1],$ $[1^2_6],$ $[1^1_{24}],$ $[1^1_{26}],$
		$[1^0_{14}],$ $[1^3_{29}],$ $[0^4_{12}],$}

	An abbreviated form is as follows:


	The array for the indices:\\
	$\begin{array}{ccrrrrr}
	i&\vline&0&1&2&3&4\\
	\hline
	a_i&\vline& 24&29&14& 1&20\\
	b_i&\vline&  3&23& 5&18&12\\
	c_i&\vline& 22&28&27&15& 4\\
	d_i&\vline& 16&19& 8&10& 9\\
	e_i&\vline& 26&13&21& 6&25\\
	\end{array}$
	
	Selectors for the $\Psi$ plane, when p = 5:\\
{\large
	\noindent$(0^0_0):$ $[0^0_{e_2}],$ $[0^0_{e_3}],$ $[0^0_{d_0}],$
		$[1^3_{e_1}],$ $[1^3_{b_2}],$ $[1^2_{e_4}],$ $[1^2_{b_3}],$
		$[0^2_{c_1}]$ $[0^3_{c_4}],$ $[1^0_{c_1}],$\\[4pt]
	\hti{6} $[2^0_{b_1}],$ $[2^0_{b_4}],$ $[2^0_{a_0}],$ $[2^0_{e_0}],$
		$[2^4_{a_2}],$ $[2^4_{d_2}],$ $[2^1_{a_3}],$ $[2^1_{d_3}],$
\\[4pt]
	\hti{6}	$[3^4_{d_1}],$ $[3^4_{c_3}],$ $[3^1_{d_4}],$ $[3^1_{c_2}],$
		$[2^2_{a_1}],$ $[2^3_{a_4}],$ $[3^0_{b_0}],$\\[10pt]
	\noindent$(1^0_0):$ $[1^0_{d_1}],$ $[1^0_{d_4}],$ $[1^0_{e_0}],$
		$[0^3_{b_1}],$ $[0^3_{e_2}],$ $[0^2_{b_4}],$ $[0^2_{e_3}],$
		$[1^4_{b_2}],$ $[1^1_{b_3}],$ $[0^0_{c_0}],$\\[4pt]
	\hti{6} $[3^0_{c_2}],$ $[3^0_{c_3}],$ $[3^0_{a_0}],$ $[3^0_{d_0}],$
		$[3^2_{a_1}],$ $[3^2_{e_1}],$ $[3^3_{a_4}],$ $[3^3_{e_4}],$
\\[4pt]
	\hti{6}	$[2^4_{c_1}],$ $[2^4_{d_3}],$ $[2^1_{c_4}],$ $[2^1_{d_2}],$
		$[3^4_{a_2}],$ $[3^1_{a_3}],$ $[2^0_{b_0}],$\\[10pt]
	\noindent$(2^0_0):$ $[2^0_{e_2}],$ $[2^0_{e_3}],$ $[2^0_{d_0}],$
		$[3^3_{e_1}],$ $[3^3_{b_2}],$ $[3^2_{e_4}],$ $[3^2_{b_3}],$
		$[2^2_{c_1}]$ $[2^3_{c_4}],$ $[3^0_{c_1}],$\\[4pt]
	\hti{6} $[0^0_{b_1}],$ $[0^0_{b_4}],$ $[0^0_{a_0}],$ $[0^0_{e_0}],$
		$[0^4_{a_2}],$ $[0^4_{d_2}],$ $[0^1_{a_3}],$ $[0^1_{d_3}],$
\\[4pt]
	\hti{6}	$[1^4_{d_1}],$ $[1^4_{c_3}],$ $[1^1_{d_4}],$ $[1^1_{c_2}],$
		$[0^2_{a_1}],$ $[0^3_{a_4}],$ $[1^0_{b_0}],$\\[10pt]
	\noindent$(3^0_0):$ $[3^0_{d_1}],$ $[3^0_{d_4}],$ $[3^0_{e_0}],$
		$[2^3_{b_1}],$ $[2^3_{e_2}],$ $[2^2_{b_4}],$ $[2^2_{e_3}],$
		$[3^4_{b_2}],$ $[3^1_{b_3}],$ $[2^0_{c_0}],$\\[4pt]
	\hti{6} $[1^0_{c_2}],$ $[1^0_{c_3}],$ $[1^0_{a_0}],$ $[1^0_{d_0}],$
		$[1^2_{a_1}],$ $[1^2_{e_1}],$ $[1^3_{a_4}],$ $[1^3_{e_4}],$
\\[4pt]
	\hti{6}	$[0^4_{c_1}],$ $[0^4_{d_3}],$ $[0^1_{c_4}],$ $[0^1_{d_2}],$
		$[1^4_{a_2}],$ $[1^1_{a_3}],$ $[0^0_{b_0}],$
}

	\section{Desarguesian Geometry.}
	\vspace{-18pt}\hspace{108pt}\footnote{22.1.87}\\[8pt]
	I will attempt to generalize the results of quaternionian geometry
	to Desarguesian geometry. It is not clear to me now that polarities
	exist in general. Indded, we have seen that we can construct a system of
	homogeneous coordinates over a skew field, for which the incidence
	property is $\Sigma P_i l_i = 0,$ with right equivalence for the lines
	$l$ anf left equivalence for the points $P$. A line collineation can be
	represented by a matrix, $m = {\bf C_l} l$ while a point
	collineation requires, $P^T = Q^T {\bf C_p},$ to allow for
	right and left equivalence.  For a polarity, these equivalences do not 
	appear to be compatible with a matrix transformation.

	It should be kept in mind that every skew field which is not a filed
	has a non trivial subfield generated  by 1, which can be finite (Ore)
	or not.  This implies that given 4 points forming a complete quadrangle,
	there exist a Pappian subgeometry through these 4 points, the
	elements of which are obtained from the linear constructions
	which start from these 4 points. 

	\sssec{Theorem.}
	{\em In any skew field, if a matrix {\bf A} has a left inverse and a
	right inverse, these are equal.}

	Proof: Let {\bf C} be the left inverse of {\bf A} and {\bf B} be its
	right inverse, by associativity of matrices,\\
	\hth${\bf C} = {\bf C}({\bf A}{\bf B}) = ({\bf C}{\bf A}){\bf B}
		= {\bf B}.$

	\sssec{Theorem.}
	{\em IF ${\bf C_p}$ is a point collineation, the line collineation
	${\bf C_l}$ is ${\bf C_p}^{-1}.$\\
	In particular, the point collineation which associates to $A_i$, $A_i$
	and to $(1,1,1),$ $(q_0,q_1,q_2)$ is\\
	\hth$\matb{(}{ccc}q_0&0&0\\0&q_1&0\\0&0&q_2\mate{)},$\\
	and the line collineation is}\\
	\hth$\matb{(}{ccc}q^{-1}_0&0&0\\0&q^{-1}_1&0\\0&0&q^{-1}_2\mate{)}.$

	Proof: If $Q$ is the image of $P$, and $m$ is the image of $l$, we want
\\
	$0 = P\cdot l = \Sigma P_i l_i = \Sigma Q_i{\bf C_p}_{ij}{\bf C_l}_{jk}
	m_k = \Sigma Q_i m_i = 0,$ for all points $P$ and incident lines $l$
	iff ${\bf C_l} = {\bf C_p}^{-1}.$

	\ssec{Desarguesian Geometry of the Hexal Complete 5-Angles.}
	\sssec{Notation.}
	In what follows, I will use the same notation as in involutive Geometry,
	namely,\\
	$l := P \times Q,$ means that the line $l$ is defined as the line
	incident to $P$ and $Q.$\\
	If subscripts are used these have the values 0, 1 and 2 and the
	computation is done modulo 3,\\
	$P \cdot l = 0$ means that the point $P$ is incident to the line $l.$\\
	When 3 lines intersect, this intersection can be defined in 3 ways, this
	has been indicated by using (*) after the definition and implies a
	Theorem.\\
	The labeling used is "H," for Hypothesis, "D", for definitions, "C",
	for conclusions, "N", for nomenclature, "P", for proofs, this
	labelling being consistent with that of the corresponding definitions.

	\sssec{The special configuration of Desargues.}
	With this notation, the special configuration of Desargues
	can be defined by\\
	\hth$	a_i := A_{i+1} \times A_{i-1},$ $qa_i = Q \times A_i,$\\
	\hth$	Q_i := a_i \times qa_i,$ $qq_i := Q_{i+1} \times Q_{i-1},$\\
	\hth$	QA_i := a_i \times qq_i,$ $q_i := A_i \times QA_i,$\\
	\hth$	QQ_i := q_{i+1} \times q_{i-1},$ $q := QA_1 \times QA_2 (*),$\\
	and the other conclusion of the special Desargues Theorem can be
	written,\\
	\hth$QQ_i \cdot qa_i = 0.$\\
	Let $Q$ and $A_i$ be\\
	\hth$	Q = (q_0,q_1,q_2),$ and $A_0 = (1,0,0),$ $A_1 = (0,1,0),$
	$A_2 = (0,0,1),$\\
	then we have the following results, not obtained in the given order,\\
	\hth$	A_0 = (1,0,0),$ $a_0 = [1,0,0],$\\
	\hth$	Q = (q_0,q_1,q_2),$ $q = [q^{-1}_0,q^{-1}_1,q^{-1}_2],$\\
	\hth$	QA_0 = (0,q_1,-q_2),$ $qa_0 = [0,q^{-1}_1,-q^{-1}_2],$\\
	\hth$	Q_0 = (0,q_1,q_2),$ $q_0 = [0,q^{-1}_1,q^{-1}_2],$\\
	\hth$	QQ_0 = (-q_0,q_1,q_2),$ $qq_0 = [-q^{-1}_0,q^{-1}_1,q^{-1}_2],$

	The self duality of the configuration corresponds to the replacement
	of points by lines where upper case letters are replaced by lower case
	letters and coordinates by their inverse.

	\sssec{Fundamental Hypothesis, Definitions and Conclusions.}
	\label{sec-funddesHDC}
	The ideal line and the coideal line.

		Given\\
	H0.0.	$A_i,$\\	
	H0.1.	$M,$ $\ov{M},$\\
		Let\\
	D1.0.	   $a_i := A_{i+1} \times A_{i-1},$\\
	D1.1.   $ma_i := M \times A_i,$
		$\ov{m}a_i := \ov{M} \times A_i,$\\
	D1.2.   $M_i := ma_i \times a_i,$
		$\ov{M}_i := \ov{m}a_i \times a_i,$\\
	D1.3.	$eul = M \times \ov{M},$

	\noindent D2.0.   $mm_i := M_{i+1} \times M_{i-1},$
		$\ov{m}m_i := \ov{M}_{i+1} \times \ov{M}_{i-1},$\\
	D2.1.   $MA_i := a_i \times mm_i,$
	$\ov{M}A_i := a_i \times \ov{m}m_i,$\\
	D2.2.	$m_i := A_i \times MA_i,$
		$\ov{m}_i := A_i \times \ov{M}A_i,$\\
	D2.3.	$MM_i := m_{i+1} \times m_{i-1},$
		$\ov{M}M_i := \ov{m}_{i+1} \times
			 \ov{m}_{i-1},$\\
	D2.4.   $m := MA_1 \times MA_2\:(*),$
		$\ov{m} := \ov{M}A_1 \times \ov{M}A_2\:(*),$\\
	D2.5.	$Ima_i := m \times ma_i,$
		$\ov{I}ma_i := \ov{m} \times \ov{m}a_i,$\\
	D2.6.	$IMa_i := m \times \ov{m}a_i,$
		$\ov{I}Ma_i := \ov{m} \times ma_i,$\\
	D2.7.	$iMA_i := M \times MA_i,$
		$\ov{\imath}MA_i := \ov{M} \times \ov{M}A_i,$\\
	Let\\
	D3.0.	$mf_i := M_i \times IMa_i,$
		$\ov{m}f_i := \ov{M}_i \times \ov{I}Ma_i,$\\
	D3.1.	$O := mf_1 \times mf_2 (*),$
		$\ov{O} := \ov{m}f_1 \times \ov{m}f_2 (*),$\\
	D3.2.	$Mfa_i := a_{i+1} \times mf_{i-1},$
		$\ov{M}fa_i := a_{i+1} \times \ov{m}f_{i-1},$\\
	\htp{28}$Mf\ov{a}_i := a_{i-1} \times mf_{i+1},$
		$\ov{M}f\ov{a}_i := a_{i-1} \times
			\ov{m}f_{i+1},$\\
	D3.3.	$mfa_i := Mfa_{i+1} \times A_{i-1},$
		$\ov{m}fa_i := \ov{M}fa_{i+1} \times A_{i-1},$\\
	\htp{28}$mf\ov{a}_i := Mf\ov{a}_{i-1} \times A_{i+1},$
		$\ov{m}f\ov{a}_i := \ov{M}f\ov{a}_{i-1}
			\times A_{i+1},$\\
	D3.4.	$Mfm_i := mfa_i \times m_i,$
		$\ov{M}fm_i := \ov{m}fa_i \times \ov{m}_i,$\\
		then\\
	C3.0.	$O \cdot eul = \ov{O} \cdot eul = 0.$\\
	C3.1.	$Mfm_i \cdot mf\ov{a}_i = \ov{M}fm_i \cdot
			\ov{m}f\ov{a}_i = 0.$\\
	Let\\
	D4.0.  $Imm_i := m \times \ov{m}m_i,$
		$\ov{I}mm_i := \ov{m} \times mm_i,$\\
	D4.1.	  $ta_i := A_i \times Imm_i,$\\
	D4.2.	  $T_i := ta_{i+1} \times ta_{i-1},$\\
	D4.3.	  $at_i:= A_i \times T_i,$\\
	D4.4.	  $K_i := at_{i+1} \times at_{i-1},$\\
	D4.5.		$TAa_i	:= ta_i \times a_i,$\\
	D4.6.		$poK_i	:= Taa_{i+1} \times Taa_{i-1},$\\
	then\\
	C4.0.  $\ov{I}mm_i \cdot ta_i = 0.$\\
	C4.1.	$T_i \cdot mf_i = 0.$

	The nomenclature:\\
	N0.0.	$A_i$	are the {\em vertices} of the triangle,\\
	N0.1.	$M$ is the {\em barycenter},
		$\ov{M}$ is the {\em orthocenter}.\\
	N1.0.   $a_i$ are the {\em sides}.\\
	N1.1.   $ma_i$ are the {\em medians},
		$\ov{m}a_i$ are the {\em altitudes}\\
	N1.2.   $M_i$ are the {\em mid-points} of the sides.
		$\ov{M}_i$ are {\em the feet of the altitudes}\\
	N1.3.	eul is the {\em line of Euler},

	\noindent N2.0.   $\{M_i,mm_i\}$ is the {\em complementary triangle},\\
	\htp{28}	$\{\ov{M}_i,\ov{m}m_i\}$ is the
			{\em orthic triangle},\\
	N2.1.   $MA_i$ are the {\em directions of the sides},\\
	N2.2.	$\{MM_i,m_i\}$	is the {\em anticomplementary triangle}.\\
	N2.3.   $m$ is the {\em ideal line} corresponding to the line at
		infinity,\\
	\htp{28}	$\ov{m}$ is the {\em orthic line} which is the polar of
		$\ov{M}$ with respect to the triangle.

\noindent N2.4.	$Ima_i$	are the {\em directions of the medians}.\\
	\htp{28}	$IMa_i$ are the {\em directions of the altitudes}.

\noindent N3.0.	$mf_i$	are the {\em mediatrices},\\
	N3.1.	$O$ is the {\em center},\\
	N3.2.	$Mfm_i$	are the {\em trapezoidal points},

\noindent N4.0.	$Imm_i$ are the {\em directions of the antiparallels of}
		$a_i${\em with respect to the\\
	\hth	sides} $a_{i+1}$ and $a_{i-1}.$\\
		N4.1.  $(T_i,ta_i)$ is the {\em tangential triangle},\\
	N4.2.  $at_i$ are the {\em symmedians},\\
	N4.3.  $K_i$ is the {\em triangle of Lemoine}.

	\sssec{Theorem.}
	If we derive a point $X$ and a line $x$ by a given construction
	from $A_i$, $M$ and $\ov{M}$, with the coordinates as given in
	G0.0 and G0.1, below, and the point $\ov{X}$ and line $\ov{x}$
	are obtain by the same construction interchange $M$ and $\ov{M}$,\\
	\hth$X = (f_0(m_0,m_1,m_2),f_1(m_0,m_1,m_2),f_2(m_0,m_1,m_2)),\\
	\hth x = [g_0(m_0,m_1,m_2),g_1(m_0,m_1,m_2),g_2(m_0,m_1,m_2)],\\
	\implies\\
	\ov{X} = (f_0(m_0^{-1},m_1^{-1},m_2^{-1})m_0,
		f_1(m_0^{-1},m_1^{-1},m_2^{-1})m_1,
		f_2(m_0^{-1},m_1^{-1},m_2^{-1})m_2),\\
	\ov{x} = [m^{-1}_0g_0(m_0^{-1},m_1^{-1},m_2^{-1}),
		m^{-1}_1g_1(m_0^{-1},m_1^{-1},m_2^{-1}),
		m^{-1}_2g_2(m_0^{-1},m_1^{-1},m_2^{-1})].$

	Proof: The point collineation
	$\bf{C_p} = \matb{(}{ccc}q_0&0&0\\0&q_1&0\\0&0&q_2\mate{)},$
	associates to (1,1,1), $(q_0,q_1,q_2),$ and to $(m_0,m_1,m_2),$
	$(r_0,r_1,r_2),$ if $r_i = m_i q_i.$\\
	In the new system of coordinates,\\
\small
{	$X = (f_0(q_0^{-1}r_0,q_1^{-1}r_1,q_2^{-1}r_2)q_0,
		f_1(q_0^{-1}r_0,q_1^{-1}r_1,q_2^{-1}r_2q_1),
		f_2(q_0^{-1}r_0,q_1^{-1}r_1,q_2^{-1}r_2)q_2).$\\
}	Exchanging $q_i$ and $r_i$ and then replacing $q_i$ by 1 and $r_i$
	by $m_i$ is equivalent to substituting $m_i$ for $q_i$ and 1 for $r_i,$ 
	which gives $\ov{X}$. $\ov{x}$ is obtained similarly.

	The line collineation is\\
	\hth$\matb{(}{ccc}q^{-1}_0&0&0\\0&q^{-1}_1&0\\0&0&q^{-1}_2\mate{)}.$

	\sssec{Notation.}\label{sec-ndes1}
	\hth$a_i := (m_{i+1}^{-1}-m_{i-1}^{-1})(m_{i-1}^{-1}+m_i^{-1})^{-1},\\
	\hth s_i := -(m_i^{-1}+m_{i+1}^{-1})(m_i^{-1}+m_{i-1}^{-1})^{-1},\\
	\hth t_i := s_{i+2} s_{i+1},\\
	\hth f_i := s_i - s_{i+1}^{-1}s_{i-1}^{-1},\\
	\hth g_i := t_i^{-1} - t_{i+1}t_{i-1}.$

	\sssec{Proof of }
	\vspace{-19pt}\hspace{84pt}{\bf\ref{sec-funddesHDC}}.\\[8pt]
	Let\\
	G0.0.	$A_0 = (1,0,0),$ $A_1 = (0,1,0),$ $A_2 = (0,0,1),$\\
	G0.1.	$M = (1,1,1),$ $\ov{M} = (m_0,m_1,m_2),$\\
	then\\
	P1.0.	$a_0 = (1,0,0),$ $a_1 = (0,1,0),$ $a_2 = (0,0,1),$\\
	P1.1.	$ma_0 = [0,1,-1],$ $\ov{m}a_0 = [0,m_1^{-1},-m_2^{-1}],$\\
	P1.2.	$M_0 = (0,1,1),$ $\ov{M}_0 = (0,m_1,m_2),$\\
	P1.3.	$eul = [1,(m_1-m_2)^{-1}(m_2-m_0),
		(m_1-m_2)^{-1}(m_0-m_1)],$

\noindent P2.0.	$mm_0 = [1,-1,-1],$ $\ov{m}m_0 = [m_0^{-1},-m^{-1}_1,-m_2^{-1}],
$\\
	P2.1.	$MA_0 = (0,1,-1),$ $\ov{M}A_0 = (0,m_1,-m_2),$\\
	P2.2.	$m_0 = [0,1,1],$ $\ov{m}_0 = [0,m_1^{-1},m_2^{-1}],$\\
	P2.3.	$MM_0 = (1,-1,-1),$ $\ov{M}M_0 = (m_0,-m_1,-m_2),$\\
	P2.4.	$m = [1,1,1],$ $\ov{m} = [m_0^{-1},m_1^{-1},m_2^{-1}],$\\
	P2.5.	$Ima_0 = (2,-1,-1)$, $\ov{I}ma_0 = (2m_0,-m_1,-m_2),$\\
	P2.6.	$IMa_0 = (m_1+m_2,-m_1,-m_2),$
		$\ov{I}Ma_0 = ((m_1^{-1}+m_2^{-1})m_0,-1,-1),$\\
	P2.7.	$iMA_0 = [2,-1,-1],$
		$\ov{\imath}MA_0 = [2m_0^{-1},-m_1^{-1},-m_2^{-1}],$

\noindent P3.0.	$mf_0 = [(m_1+m_2)^{-1}(m_1-m_2),1,-1],\\
	\htp{28}\ov{m}f_0 = [m_0^{-1}(m_1^{-1}+m_2^{-1})^{-1}
			(m_1^{-1}-m_2^{-1}),m^{-1}_1,-m_2^{-1},1],$\\
	P3.1.	$O = (m_1+m_2,m_2+m_0,m_0+m_1),\\
	\htp{28}\ov{O} = ((m_1^{-1}+m_2^{-1})m_0,
			(m_2^{-1}+m_0^{-1})m_1,(m_0^{-1}+m_1^{-1})m_2),$\\
	P3.2.	$Mfa_0 = (1,0,-(m_0-m_1)^{-1}(m_0+m_1)),\\
	\htp{28}\ov{M}fa_0 = (m_0,0,
			(m_0^{-1}-m_1^{-1})(m_0^{-1}+m_1^{-1})^{-1}m_2),\\
	\htp{28}Mf\ov{a}_0 = (1,(m_2-m_0)^{-1}(m_2+m_0),0),\\
	\htp{28}\ov{M}f\ov{a}_0 = (m_0,
			(m_2^{-1}-m_0^{-1})^{-1}(m_2^{-1}+m_0^{-1})m_1,0),$\\
	P3.3.	$mfa_0 = [(m_1+m_2)^{-1}(m_1-m_2),1,0],\\
	\htp{28}\ov{m}fa_0 = [m_0^{-1}(m_1^{-1}+m_2^{-1})^{-1}
			(m_1^{-1}-m_2^{-1}),m^{-1}_1,0],\\
	\htp{28}mf\ov{a}_0 = [(m_1+m_2)^{-1}(m_1-m_2),0,-1],\\
	\htp{28}\ov{m}f\ov{a}_0 = [m^{-1}_0(m_1^{-1}+m_2^{-1})^{-1}
			(m_1^{-1}-m_2^{-1}),0,m_2^{-1}],$\\
	P3.4.	$Mfm_0 = ((m_1+m_2)(m_1-m_2)^{-1},-1,1),\\
	\htp{28}\ov{M}fm_0 = ((m_1^{-1}-m_2^{-1})^{-1}(m_1^{-1}+m_2^{-1})m_0,
			-m_1,m_2),$

\noindent P4.0.	$Imm_0 = (a_0,1,s_0),$ $\ov{I}mm_0 = (-a_0,1,s_0),$\\
	P4.1.	$ta_0 = [0,1,-s_0^{-1}],$\\ \marginpar{?}
	P4.2.	$T_0 = (1,s_2,s_1^{-1}),$\\
	P4.3.	$at_0 = [0,s_2^{-1},-s_1] = [0,1,-t_0],$\\
	P4.4.	$K_0 = (1,t_2^{-1},t_1),$\\
	P4.5.	$Taa_0 = (0,1,s_0),$\\
	P4.6.	$poK_0 = [-1,s_2^{-1},s_1],$

	Details of proof:\\
	For P4.0, if the coordinates of $Imm_0$ are $x_0$, 1 and $x_2$,
	we have to solve\\
	\hth$x_0 + 1+x_2 = 0,$ $-x_0m_0^{-1} + m_1^{-1}+x_2m_2^{-1} = 0.$\\
	Multiplying the equations to the right respectively by $m_2^{-1},$ and
	-1 or by $m_0^{-1}$ and 1 and adding gives $x_0$ and $x_2$ using the
	notation \ref{sec-ndes1}.

	\ssec{Perpendicularity mapping.}
	\sssec{Definition.}
	Given $\ov{M}$ and a direction $I_x$, the perpendicular direction $I_y$
	is defined by the following construction\\
	D5.0.	$b := A_0 \times I_x,$\\
	D5.1.	$B := b \times a_0,$\\
	D5.2.	$c := B \times IMa_2,$\\
	D5.3.	$C := c \times \ov{m}a_0,$\\
	D5.4.	$d := C \times A_1,$\\
	D5.5.	$I_y := d \times m,$

	\sssec{Theorem.}
	{\em If $I_x = (-1-q,q,1)$ and $I_y = (1-r,r,1)$ then}\\
	$r = .$

	Proof:\\
	P5.0.	$b = [0,-q^{-1},1],$\\
	P5.1.	$B = (0,q,1),$\\
	P5.2.	$c = [x,-q^{-1},1],$ with $x = m_0^{-1}(m_0+m_1+m_1q^{-1}),$\\
	P5.3.	$C = (y,m_1,m_2),$ with $y = (m_1q^{-1}-m_2)x^{-1},$\\
	P5.4.	$d = [y^{-1},0,-m_2^{-1}],$\\
	P5.5.	$I_y = (y,z,m_2),$ with $z = -y-m_2$,\\
	Therefore $r = -m_2^{-1}(y+m_2)
		= -m_2^{-1}(m_2+(m_1q^{-1}-m_2)(m_0+m_1+m_1q^{-1})^{-1}m_0).$

	If the skew field we therefore have\\
	\hth$-m_2(r+1)m_0^{-1}((m_0+m_1)q+m_1) = m_1-m_2q$.\\
	\hth$-rm_0^{-1}(m_0+m_1)q -rm_0^{-1}m_1-m_0^{-1}m_1q - m_2^{-1}m_1
		= 0$.\\
	which is, in general, not an involution.

	\section{The Hughes Planes.}
\setcounter{subsection}{-1}
	\ssec{Introduction.}
	There are essentially 2 methods to algebraize a plane. The first one
	which start with the work of Desargues coordinatized the plane using
	2 coordinates, the difficulty of representing the ideal points or
	points at infinity can be dealt with by using 3 homogeneous
	coordinates. This approach has been generalized to perspective
	planes, for which the only axioms are those of incidence, by
	using as coordinates, elements of a ternary ring instead of elements in
	a field.  This generalization was given by Marshall Hall in 1943,
	but its origin can be found, for the case of nearfields, introduced
	by Dickson (1905), in the most remarkable paper of 
	Veblen and MacLagan-Wedderburn in 1907 (p. 380-382).\\
	In this paper they give, independently from Vahlen the first example
	of non Pappian Geometry. The indpendent result consited in
	showing that quaternions could be used as coordinates for such
	a geometry.\\
	The second approach, which can be used in finite planes, is to construct
	a difference set, of $q = p^k$ integers as a subset $\{0,\ldots,q^2+q\}$
	from which the points incident to each line, and the lines
	incident to every point can be completely derived.  This approach
	was fully examined for the finite Pappian planes by J. Singer in
	1938, but it again can be traced in the paper of 1907 (p. 383 and 385).
	Moreover, the generalization to non Desarguesian planes is given
	explicitely for a plane of order 9, called $\Psi$ plane
	by Room and Kirkpatrick.\\
	I do prefer, when applying the notion of difference sets to geometry,
	to use, instead of it, the terminology of selector introduced by Fernand
	Lemay, in 1979. (See his most accessible paper of 1983.)\\
	It is the second approach, that I am exploring in this paper, gives many
	of the results in the form of conjectures.\\
	We will see that to give the incidence properties for planes of the
	$\Psi$ type and order $p^2$ we have to give $p$ selectors of $p$
	elements, and in a particular notation the points which are incident
	to $\frac{p-1}{2}$ lines from which all other incidences can be derived.
	The notation is such that the same incidence tables are valid for
	the points on any line, giving rise to a fundamental polarity.\\
	One of the advantages of the selector approach is to eliminate the need
	of addition and multiplication tables in the particular nearfield
	which greatly simplifies the exploration of new properties with
	a computer.  Many of the planes are special case of Hughes planes,
	hence the tiltle of the section.

	I will assume that $p$ is an odd prime.

	\ssec{Nearfield and coordinatization of the plane.}
	\sssec{Definition. [Dickson]}
	A {\em left nearfield} (${\Bbb N},+,\circ)$ is a set ${\Bbb N}$
	with binary operations such that
	\enumb
	\item ${\Bbb N}$ is finite,
	\item $({\Bbb N},+)$ {\em is an Abelian group, with neutral element} 0,
	\item $({\Bbb N}-\{0\},\circ)$ {\em is an group, with neutral element}
		1,
	\item $\circ$ is left distributive over $+$, or\\
	\hth$\zeta\circ (\xi+\eta)
			= \zeta\circ\xi+\zeta\circ\eta,$
		 for all $\xi,\eta,\zeta\in{\Bbb N}$
	\item $0\circ\xi = 0,$ for all $\xi\in{\Bbb N}$.
	\enume

	For a right nearfield, the left distributive law is replaced by the 
	right one and $0\circ\xi = 0,$ is replaced by $\xi\circ 0 = 0$.

	\sssec{Theorem.}
	{\em In any left nearfield,}
	\enumb
	\item $\xi\circ 0 = 0$ for all $\xi\in{\Bbb N}$.
	\item $\xi\circ\eta = 0\implies \xi = 0$ {\em or} $\eta = 0.$
	\item $1,-1\neq 0.$
	\enume

	\sssec{Definition. [Dickson]}
	Let $n$ be a non residue of $p$.
	A {\em Dickson left nearfield} (${\Bbb N},+,\circ)$ is a set ${\Bbb N}$
	with the operations\\
	\hth$(a_0+b_0\alpha)+(a_1+b_1\alpha) := ((a_0+a_1)+(b_0+b_1)\alpha),\\
	\hth (a_0+b_0\alpha)\cdot(a_1+b_1\alpha)
		 := ((a_0a_1+e\:n\:b_0b_1)+(a_1b_0+e\:a_0b_1)\alpha),$\\
	\hti{12} where $e = +1$ if $a_1^2-n\:b_1^2$ is a quadratic residue
	of $p$ and $e = -1$ otherwize.\\
	A Dickson right nearfield is obtained by the replacement of
	$(a_1b_0+e\:a_0b_1)\alpha),$ by $(a_0b_1+e\:a_1b_0)\alpha)$.

	\sssec{Theorem.}
	{\em A Dickson left nearfield is a left nearfield.}

	\sssec{Definition.}
	Let $\beta,$ $\gamma$, $\xi$ and $\eta$ are elements of a Dickson 
	nearfield. In a {\em Hughes plane}, the {\em points} are the triples\\
	\hth$(1,\beta,\gamma),$ $(0,1,\gamma),$ $(0,0,1),$\\
	and the {\em lines} are the triples\\
	\hth$[\eta,\xi,-1],$ $[\eta,-1,0],$ $[1,0,0]$\\
	A point $(P_0,P_1,P_2)$ is {\em incident} to a line $[l_0,l_1,l_2]$ if
	and only if\\
	\hth$P_0l_0+P_1l_1+P_2l_2 = 0$.\\
	A point or a line is {\em real} if the coefficient of $\alpha$ in its
	coordinates are 0. A point or line is {\em complex}, otherwize.

	The notation is used to indicate the close relationship with the
	corresponding coordinates in a ternary ring, see for instance
	Artzy, p. 203-203,\\
	for the points: $(b,c) = (1,b,c),$ $(c) = (0,1,c),$
		 $(\infty) = (0,0,1),$\\
	for the lines:  $[x,y] = [y,x,-1],$ $[y] = [y,1,0],$
		$[\infty] = [1,0,0]$\\
	indeed $1\cdot y+b\cdot x-c = 0$ corresponds to $c = b\cdot x + y$,
	giving the ternary ring conditions of incidence.

	\sssec{Theorem. [Hughes]}
	{\em A Hughes plane is of Lenz-Barlotti type I.1}

	See Hughes, Rosati and Dembowski, p. 247.  The simplest case p = 3, is
	given by Veblen
	and MacLagan-Wedderburn p. 383, it is called in this case a $\Psi$
	plane by Room and Kirkpatrick.
 
	\sssec{Theorem.}
	{\em A real line has $p+1$ real and $p^2-p$ complex points incident to
	it.\\
	A complex line has 1 real and $p^2$ complex points incident to it.}

	See, for instance, Room and Kirkpatrick.

	\sssec{Theorem.}
	{\em The $p$ selectors and the negative inverses of the fundamental
	selector modulo $p^2+p+1$ form a partition of the set
	$\{1,\ldots , p^2+p\}$.}

%
	\sssec{Theorem.}
	{\em The $p^4+p^2+1$ points are partitioned into $p^2+p+1$ real points
	and $p-1$ phyla of complex points.  Each phylum consists of $p$
	classes. Each class consists of $p^2+p+1$ points, which form by
	definition a coplane.}

	Starting with the work of L. E. Dickson of 1905, non-Desarguesian planes
	of order 9 were discovered by Veblen
	and Wedderburn in 1907, I will here consider only one of these
	which is self dual, and for which non trivial polarities exists,
	and refer to the work of G. Zappa (1957),
	T. G. Ostrom (1964), D. R. Hughes (1957) and T. G. Room and
	P. B. Kirkpatrick (1971) for further reading.

	The synthetic definition used can be traced to Veblen and Wedderburn,
	who first consider points obtained by aplying a transformation (see p.
	383), later generalized by J. Singer.  The notation is
	inspired by Room and Kirkpatrick (see Table 5.5.4) using the same
	method I used for the finite plane reversing the indices for lines.\\
	An alternate definition, (5.6.1), is given by Room and Kirkpatrick.

	\ssec{Miniquaternion nearfield.}

	\sssec{Theorem.}\label{sec-tnearfield9}
	{\em In any left nearfield ${\Bbb Q}_9,$ of order 9,}
	\enumb
	\item $\{0,1,-1\}\approx {\Bbb Z}_3.$
	\item $\xi+\xi+\xi = 0,$ for all $\xi\in{\Bbb Q}_9,$
	\item $-1\circ\xi = \xi\circ(-1)=\xi,$ {\em for all} $\xi\in{\Bbb Q}_9,$
	\item $(-\xi)\circ\eta = \xi\circ(-\eta) = -(\xi\circ\eta),$
		{\em for all} $\xi,\eta\in{\Bbb Q}_9,$
	\item $(-\xi)\circ(-\eta) = \xi\circ\eta,$ {\em for all}
		$\xi,\eta\in{\Bbb Q}_9,$
	\item {\em Given $\kappa\in {\Bbb Q}_9^* := {\Bbb Q}_9 - {0},$
		$\lambda = s - \kappa r$ determines a one to one correspondance
		between the elements $\lambda \in {\Bbb Q}_9$ and the pairs
		$(r,s)$, $r,s\in {\Bbb Z}_3$.}
	\item {\em ${\Bbb Q}_9$ being an other nearfield of order 9, the
		groups $({\Bbb Q}_9,+)$ and $({\Bbb Q}'_9,+)$ are isomorphic.}
	\item {\em Besides GF($3^2$) there is only one nearfield of order 9,
		which is the smallest nearfield which is not a field,}
		(Zassenhaus, 1936).
	\enume

	\sssec{Exercise.}\label{sec-enearfield9}
	Determine the correspondance of \ref{sec-tnearfield9}.5.

	\sssec{Definition.}
	The {\em left miniquaternions} is the set
	 ${\Bbb Q}_9 := \{0,\pm 1, \pm\alpha,\pm \beta, \pm \gamma\}$
	with the operations of addition and multiplications defined from,\\
	\hth$  \xi +\xi +\xi = 0$ for all $\xi \in {\Bbb Q}_9$,\\
	\hth$  \alpha - 1 = \beta ,$ $\alpha + 1 = \gamma ,$\\
	\hth$  \alpha^2 = \beta^2 = \gamma^2 = -\alpha\beta\gamma = -1$.\\
	The set ${\Bbb Q}_9^* := \{\pm\alpha,\pm\beta, \pm \gamma\}$.

	For the right miniquaternions, we replace $\alpha\beta\gamma = 1$
	by $\alpha\beta\gamma = -1$.

	\sssec{Theorem.}
	\enumb
	\item $\alpha-\beta = \beta-\gamma = \gamma -\alpha = 1,$
		$\alpha+\beta+\gamma = 0.$
	\item $-\beta\gamma = \gamma\beta = \alpha,$
		$-\gamma\alpha = \alpha\gamma = \beta,$
		$-\alpha\beta = \beta\alpha = \gamma.$
	\item {\em the multiplication is left distributive,
		$\tau(\rho+\sigma) = \tau\rho+\tau\sigma,$\\
		for all} $\rho,\sigma,\tau\in {\Bbb Q}_9$.
	\item $\{{\Bbb Q}_9,+,.\}$ {\em is a left nearfield}.
	\item $\{{\Bbb Q}_9,+,.\}$ {\em is not a field}, e. g.\\
		$\alpha(\alpha+\beta) = \alpha(-\gamma) = \beta$,
		$\alpha\alpha+\alpha\beta = -1+\gamma = \alpha$.
	\item $\begin{array}{rcrrrrrrrr}
	+&\vline&1&-1&\alpha&-\alpha&\beta&-\beta&\gamma&-\gamma\\
	\hline
	1&\vline&-1&0&\gamma&-\beta&\alpha&-\gamma&\beta&-\alpha\\
	-1&\vline&0&1&\beta&-\gamma&\gamma&-\alpha&\alpha&-\beta\\
	\alpha&\vline&\gamma&\beta&-\alpha&0&-\gamma&1&-\beta&-1\\
	-\alpha&\vline&-\beta&-\gamma&0&\alpha&-1&\gamma&1&\beta\\
	\beta&\vline&\alpha&\gamma&-\gamma&-1&-\beta&0&-\alpha&1\\
	-\beta&\vline&-\gamma&-\alpha&1&\gamma&0&\beta&-1&\alpha\\
	\gamma&\vline&\beta&\alpha&-\beta&1&-\alpha&-1&-\gamma&0\\
	-\gamma&\vline&-\alpha&-\beta&-1&\beta&1&\alpha&0&\gamma
	\end{array}\\[10pt]
	\begin{array}{rcrrrrrrrr}
	\cdot&\vline&1&-1&\alpha&-\alpha&\beta&-\beta&\gamma&-\gamma\\
	\hline
	1&\vline&1&-1&\alpha&-\alpha&\beta&-\beta&\gamma&-\gamma\\
	-1&\vline&-1&1&-\alpha&\alpha&-\beta&\beta&-\gamma&\gamma\\
	\alpha&\vline&\alpha&-\alpha&-1&1&-\gamma&\gamma&\beta&-\beta\\
	-\alpha&\vline&-\alpha&\alpha&1&-1&\gamma&-\gamma&-\beta&\beta\\
	\beta&\vline&\beta&-\beta&\gamma&-\gamma&-1&1&-\alpha&\alpha\\
	-\beta&\vline&-\beta&\beta&-\gamma&\gamma&1&-1&\alpha&-\alpha\\
	\gamma&\vline&\gamma&-\gamma&-\beta&\beta&\alpha&-\alpha&-1&1\\
	-\gamma&\vline&-\gamma&\gamma&\beta&-\beta&-\alpha&\alpha&1&-1
	\end{array}$
	\enume

	For the right miniquaternions, we change the sign of the products in 1.
	and exchange rows and columns for the multiplication table,
	e.g. $\alpha\beta = \gamma.$

	\ssec{The first non-Pappian plane, by Veblen and Wedderburn.}
	
	\sssec{Definition. [Veblen-Wedderburn]}
	The {\em points} $P$ are $(x,y,1),$ $(x,1,0),$ $(1,0,0),$ the
	{\em lines} $l$ are $[1,b,c],$ $[0,1,c],$ $[0,0,1],$ and the
	{\em incidence} is $P\cdot l = 0$, where $x,$ $y,$ $b$ and $c$ are
	elements of a left nearfield.

	\sssec{Theorem. [Veblen-Wedderburn]}
	{\em With $b$, $c$, $b'$, $c'$ in ${\Bbb Q}_9$,}
	\enumb
	\item $[1,b,c]\times [1,b',c'] = (-(yb+c),y,1),$ with $y(b-b')=-(c-c').$
	\item $[1,b,c] \times [0,1,c'] = (c'b-c,-c',1),$
	\item $[1,b,c] \times [0,0,1] = (-b,1,0),$
	\item $[0,1,c] \times [0,1,c'] = (1,0,0),$
	\item $[0,1,c] \times [0,0,1] = (1,0,0),$
	\enume

	\sssec{Theorem. [Veblen-Wedderburn]}
	{\em Let} $(a(b+c) = a\:b+a\:c)$
	\enumb
	\item ${\bf M} := \matb{(}{rrr}1&0&1\\-1&0&0\\0&-1&-1\mate{)},$
	\item $A_0 := (-1,0,1),$ $B_0 := (-\gamma,\alpha,1),$
		$C_0 := (\beta,-\alpha,1),$ $D_0 := (-\beta,\gamma,1),$
		$E_0 := (\alpha,-\gamma,1),$ $F_0 := (\gamma,-\beta,1),$
		$G_0 := (-\alpha,\beta,1),$
	\item $A_j := {\bf M}^j A_0$, $B_j := {\bf M}^j B_0$, \ldots,
		{\em for $j = 1$ to $12$},
	\item $a_j := \{M^jX_i\}$, $X_i\in a_0$, and similarly for $b_j$ to
		$g_j$.
	\item $a_0 := [1,1,1],$ $b_0 := [1,\alpha,1]$, $c_0 := [1,-\alpha,1]$,
		$d_0 := [1,\gamma,1]$, $e_0 := [1,-\gamma,1]$,
		$f_0 := [1,-\beta,1]$, 	$g_0 := [1,\beta,1]$,\\
	{\em ~~~~then}
	\item $M$ is of order 13.
	\item $a_0 = \{A_0,A_1,A_3,A_9,B_0,C_0,D_0,E_0,F_0,G_0\},\\
		b_0 = \{A_0,B_1,B_8,D_3,D_{11},E_2,E_5,E_6,G_7,G_9\},\\
		c_0 = \{A_0,C_1,C_8,E_7,E_9,F_3,F_{11},G_2,G_5,G_6\},\\
		d_0 = \{A_0,B_7,B_9,D_1,D_8,F_2,F_5,F_6,G_3,G_{11}\},\\
		e_0 = \{A_0,B_2,B_5,B_6,C_3,C_{11},E_1,E_8,F_7,F_9\},\\
		f_0 = \{A_0,C_7,C_9,D_2,D_5,D_6,E_3,E_{11},F_1,F_8\},\\
		g_0 = \{A_0,B_3,B_{11},C_2,C_5,C_6,D_7,D_9,G_1,G_8\},$
	\item $A_0 = \{a_0,a_4,a_{10},a_{12},b_0,c_0,d_0,e_0,f_0,g_0\},\\
		B_0 = \{a_0,b_5,b_{12},d_4,d_6,e_7,e_8,e_{11},g_2,g_{10}\},\\
		C_0 = \{a_0,c_5,c_{12},e_2,e_{10},f_4,f_6,g_7,g_8,g_{11}\},\\
		D_0 = \{a_0,b_2,b_{10},d_5,d_{12},f_7,f_8,f_{11},g_4,g_6\},\\
		E_0 = \{a_0,b_7,b_8,b_{11},c_4,c_6,e_5,e_{12},f_2,f_{10}\},\\
		F_0 = \{a_0,c_2,c_{10},d_7,d_8,d_{11},e_4,e_6,f_5,f_{12}\},\\
		G_0 = \{a_0,b_4,b_6,c_7,c_8,c_{11},d_2,d_{10},g_5,g_{12}\},$
	\item $X_j \incid x_k \implies X_{j+l\umod{13}}\incid x_{k+l\umod{13}}.$
	\enume

	Given $a_0$, to $g_0$, it is easy to verify that $A_0$ is on all these
	lines and determine $B_0$ to $G_0$ all on $a_0$ and $B_0$ on $b_0,$
	$C_0$ on $c_0$, \ldots.\\
	Having determined, the other points using 2, it is easy to verify
	which points are on $b_0$, \dots.\\
	The notation helps gretaly in justifyong that 2 points have one and
	only one line in common and 2 liners have only one point in common.
	The notation can be made even more compact.  See \ref{sec-dpsiinc1}.

	The following are the powers of $M$.\\
{\footnotesize  ${\bf M}^2 = \matb{(}{rrr}1&-1&0\\-1&0&-1\\1&1&1\mate{)},$
		${\bf M}^3 = \matb{(}{rrr}-1&0&0\\-1&1&0\\0&-1&0\mate{)},$
		${\bf M}^4 = \matb{(}{rrr}-1&-1&1\\1&0&-1\\1&0&0\mate{)},$
		${\bf M}^5 = \matb{(}{rrr}0&-1&1\\1&1&-1\\1&0&1\mate{)},$
		${\bf M}^6 = \matb{(}{rrr}1&-1&-1\\0&1&-1\\1&-1&0\mate{)},$
		${\bf M}^7 = \matb{(}{rrr}-1&1&-1\\-1&1&1\\-1&0&1\mate{)},$
		${\bf M}^8 = \matb{(}{rrr}1&1&0\\1&-1&1\\-1&-1&1\mate{)},$
		${\bf M}^9 = \matb{(}{rrr}0&0&1\\-1&-1&0\\0&-1&1\mate{)},$
		${\bf M}^{10}= \matb{(}{rrr}0&-1&-1\\0&0&-1\\1&-1&-1\mate{)},$
		${\bf M}^{11}= \matb{(}{rrr}1&1&1\\0&1&-1\\-1&1&1\mate{)},$
		${\bf M}^{12}= \matb{(}{rrr}0&-1&0\\-1&-1&-1\\1&1&0\mate{)},$
		${\bf M}^{13}= \matb{(}{rrr}1&0&0\\0&1&0\\0&0&1\mate{)}.$}

	\sssec{Example. [Veblen-Wedderburn]}\label{sec-eVWDes}
	With the notation\\
	$non-Desargues(C\langle c_0,c_1,c_2\rangle,\{A_0,A_1,A_2\}
	\{a_0,a_1,a_2\},\{B_0,B_1,B_2\}\{b_0,b_1,b_2\};\\
	\{C_0,C_1,C_2\}\{d_0,d_1,d_2\})$, with $d_i := C_{i+j}\times C_{i-j},$\\
	the following configuration shows that the Desargues axiom is not
	satisfied\\
	$non-Desargues(A_0\langle b_0,f_0,d_0\rangle,\{B_1,C_7,F_2\}
	\{c_{12},e_8,e_9\},\{D_3,E_3,D_1\}\{b_{11}g_7,a_3\};\\
		\{G_5,B_{10},F_3\}\{d_1,c_0,b_9\}).$

	\ssec{The miniquaternionian plane $\Psi$.}

	\sssec{Definition.}\label{sec-dpsiinc1}
	With $i\in\{0,1,2\},$ $i'\in\{0',1',2'\},$ $j\in\{0,1,\ldots,12\},$
	and the addition being performed modulo $3$ for the first element of a
	pair, and modulo $13$, for the second element in the pair or for the
	element, if single, then the elements and incidence in the
	{\em miniquaternionian plane $\Psi$} are defined as follows.
	(See \ref{sec-cVWRKDV}
	\enumb
	\item	The 91 {\em points} $P$ are $(j),(i,j), (i',j)$,
	\item	The 91 {\em lines} $l$ are $[j],[i,j], [i',j]$,
	\item	The {\em incidence} is defined by\\
	\hti{4}$[j] := \{(-j),(1-j),(3-j),(9-j),(i,-j),(i',-j)\},$\\
	\hti{4}$[i,j] := \{(-j),(i,2-j),(i,5-j),(i,6-j),(i'+1,3-j),
	(i'+1,11-j),\\
	\hth		(i'-1,7-j)(i'-1,9-j),(i',1-j),(i',8-j)\},$\\
	\hti{4}$[i',j] := \{(-j),(i',2-j),(i',5-j),(i',6-j),(i-1,3-j),
		(i-1,11-j),\\
	\hth		(i+1,7-j)(i+1,9-j),(i,1-j),(i,8-j)\}.$
	\enume
	giving the 10 points on each line. $i$ and $i'$ in the same definiton 
	correspond to the same integer, 0, 03-' or 1, 1' \ldots.

	\sssec{Exercise.}\label{sec-epsiinc1}
	\ref{sec-dpsiinc1}.2 is similar to the use of ordered cosets to
	determine efficiently operations of finite as well as infinite groups.
	In this case, $[j]$ is a subplane, $[i,j]$ and $[i',j]$ are 
	copseudoplanes.
	\enumb
	\item Perform a similar representation of points, lines and incidence
	starting with a subplane which is a Fano plane.
	\item Determine similar representations for non Desarguesian geometries
	of order $5^2$, using a subplane of order 4, or of order 5
	$(651 = 31\cdot 21).$
	\item Determine other such representation for non Desarguesian
	geometries of higher order.
	\enume

	\sssec{Theorem.}
	{\em The same incidence relations obtain, if we
	interchange points and lines in} \ref{sec-dpsiinc1}.2.

	\sssec{Theorem. [see Room and Kirkpatrick]}\label{sec-tminipola}
	\enumb
	\item\begin{enumerate}
		\item [0] {\em The correspondance $(j)$ to $[j]$ and
		    $(i,j)$ to $[i,j]$ and $(i',j)$ to $[i',j]$ is a polarity
			${\cal P}_0$ (${\cal J}^*$).}
		\item [1] {\em The $16$ auto-poles are} (0), (7), (8), (11),
			(0,8), (0,12), (1,4), (1,7), (2,10), (2,11),
			(0',8), (0',12), (1',4), (1',7), (2',10), (2',11).
		\enume
	\item\begin{enumerate}
		\item [0] {\em The correspondance $(j)$ to $[j]$ and $(i,j)$ to
			$[i',j]$ and $(i',j)$ to $[i,j]$ is a polarity
			${\cal P}_1$ (${\cal J}'^*$).}
		\item [1] {\em The $22$ auto-poles are} (0), (7), (8), (11),
			(0,1), (0,3), (0,9), (1,1), (1,3), (1,9), (2,1), (2,3),
			(2,9), (0',1), (0',3), (0',9), (1',1), (1',3), (1',9),
			(2',1), (2',3), (2',9).
		\item [2] (0),  (7), (8), (11), (0,1), (1,9), (2,3), (0',9),
			(1',3), (2',1),\\
			(0), (7), (8), (11), (1,1), (2,9), (0,3), (2',9),
			(0',3), (1',1),\\
			(0), (7), (8), (11), (2,1), (0,9), (1,3), (1',9),
			(2',3), (0',1) {\em are ovals.}
		\enume
	\enume

	\sssec{Exercise.}
	\enumb
	\item Prove that the correspondance $(j)$ to $[j]$ and $(i,j)$ to
			$[(i+1)',j]$ and $(i',j)$ to $[i-1,j]$ is a polarity
			${\cal P}_2$.
	\item Prove that the correspondance $(j)$ to $[j]$ and $(i,j)$ to
			$[(i-1)',j]$ and $(i',j)$ to $[i+1,j]$ is a polarity
			${\cal P}_3$.
	\enume

	\sssec{Exercise.}\label{sec-eminipola}
	\enumb
	\item Determine a configuration in \ref{sec-tminipola}.0.2, which gives
		an example were the Theorem of Pascal is satisfied and an other,
		in which it is not satisfied.
	\item Determine ovals which are subsets of \ref{sec-tminipola}.1.1.
	\enume

	\sssec{Theorem.}
	{\em The polar $m$ of a point $M$ with respect to a triangle is
	incident to that point.}

	Indeed, we can always assume thast the triangle consists of the
	real points $A_0 =(0)$, $A_1 =(1)$, $A_2 =(2)$, and that
	$M = (5) = (1,1,1).$ It follows that $m$ = [4] = [1,1,1] which is
	incident to $M$.

	\sssec{Exercise.}
	Check that the other points and lines of the polar construction are
	$M_i = (4), (8), (3),$ $MA_i = (10), (12), (9),$ $MM_i = (7), (6),
	(11),$ $a_i = [12], [1], [0],$ $ma_i = [9], [8], [11],$
	$m_i = [3], [2], [7],$ $mm_i = [6], [10], [5].$

	\sssec{Theorem. [see Room and Kirkpatrick]}\label{sec-tpsifano}
	\enumb
	\item {\em The planes obtained by taking the complete quadrangle
		associated with $3$ real points $A_0,$ $A_1,$ $A_2,$ and a point
		$\ov{M}$ which such that none of the lines $\ov{M}\times A_i$
		are real are Fano planes associated with ${\Bbb Z}_2$.}
	\item {\em There are $(\frac{1}{6}13.12.9).24 = 5616$ Fano planes that
		contain $3$ real points.}
	\enume

	\sssec{Notation.} For a Fano subplane with 7 elements, I will use
	the notation associated with the selector ${0,1,3}$ and construction:\\
	Given a complete quadrangle $0, 1, 2, 5$,\\
	$0^* := 0 \times 1$, $6^* := 1 \times 2$, $1^* := 2 \times 0$,
	$5^* := 2 \times 5$, $3^* := 0 \times 5$, $2^* := 1 \times 5$,\\
	$3 := 0^* \times 5^*$, $4 := 6^* \times 3^*$, $6 := 1^* \times 2^*$,
	$4^* := 3 x 6$. The Fano plane propery implies $4 \iota 4^*$.\\
	The configuration is denoted by
	$Fano(0,1,2,3,4,5,6,0^*,1^*,2^*,3^*,4^*,5^*,6^*).$
	
	\sssec{Example.}
	The following is a Fano plane configuration:
	$Fano((0), (1), (2), (2,0), (1,1), (0,3), (1',12),$\\
	$[0], [1], [0,12], [1',0], [1',6], [0',11], [12]).$

	\sssec{Exercise.}\label{sec-epsifano}
	Determine the Fano plane associated with (0), (1), (2), (0,7).

{\tiny	\sssec{Comment.}\label{sec-cVWRKDV}
	The correspondance between the notation of Veblen-Wedderburn
	and Room-Kirkpatrick is\\
	$\begin{array}{lcccccccc}
	Veblen-Wedderburn&\vline&a_j&b_j&c_j&d_j&e_j&f_j&g_j\\
	Room-Kirkpatrick&\vline&k_j&a_j&b_j&c_j&a'_j&b'_j&c'_j\\
	De\:Vogelaere&\vline&[-j]&[0,-j]&[2',-j]&[1,-j]&[0',-j]&[1',-j]&[2,-j]\\
	\hline
	Veblen-Wedderburn&\vline&A_j&B_j&C_j&D_j&E_j&F_j&G_j\\
	Room-Kirkpatrick&\vline&K_j&A'_j&C'_j&B'_j&A_j&C_j&B_j\\
	De\:Vogelaere&\vline&(j)&(0',j)&(2,j)&(1',j)&(0,j)&(1,j)&(2',j)
	\end{array}$}

	\sssec{Example. [Veblen-Wedderburn]}
	The example \ref{sec-eVWDes} becomes with the above notation\\
	$((0)\langle [0,0],[1',0],[2,0]\rangle,
	\{(0,1),(1,7),(1',2\}\{[1,1],[0',8],[0',9]\},
	\{(2,3),(0'3),(2,1)\}\\
	\hth\{[0,11],[2',7],[10]\};\\
		\{(2',5),(0,{10}),(1',3)\}\{[2,1],[1,0],[0,4]\}).$

	\sssec{Problem.}
	Can we characterize the plane $\Psi$ using Theorem
	\ref{sec-tpsifano}.0.

\newpage
	\sssec{Definition.}
	The {\em Singer} matrix ${\bf G} :=
	\matb{(}{ccc}0&0&1\\1&0&1\\0&1&0\mate{)}.$
	Its powers ${\bf G}^k$ are the columns\\
	$k,k+1,k+2$ of\\
	\hth$\begin{array}{ccrrrrrrrrrrrrr}
	k=&\vline&0&1&2&3&4&5&6&7&8&9&10&11&12\\
	\hline
	&\vline&1&0&0&1&0&1&1&1&-1&-1&0&1&-1\\
	&\vline&0&1&0&1&1&1&-1&-1&0&1&-1&1&0\\
	&\vline&0&0&1&0&1&1&1&-1&-1&0&1&-1&1
	\end{array}$

	move to g6a.tex:

	\sssec{Answer to}
	\vspace{-18pt}\hspace{94pt}{\bf \ref{sec-enearfield9}.}\\[8pt]
	$\begin{array}{crrcrrrrrrrrr}
	\kappa =&\vline&\lambda=&\vline&0&1&-1&\alpha&-\alpha&\beta&-\beta
		&\gamma&-\gamma\\
	\hline
	\alpha&\vline&r&\vline&0&0&0&-1&1&-1&1&-1&1\\
		&s&\vline&0&1&-1&0&0&-1&1&1&-1\\
	\hline
	\beta&\vline&r&\vline&0&0&0&-1&1&-1&1&-1&1\\
		&s&\vline&0&1&-1&1&-1&0&0&-1&1\\
	\hline
	\gamma&\vline&r&\vline&0&0&0&-1&1&-1&1&-1&1\\
		&s&\vline&0&1&-1&-1&1&1-1&0&0\\
	\end{array}$
	For $-\alpha,$ $-\beta$, $-\gamma$, change the sigh of $r$.

	\sssec{Definition.}
	The elements and incidence in the {\em miniquaternionian plane $\Psi$}
	are defined as follows.
	\enumb
	\item	The {\em points} are $(\xi_0,\xi_1,\xi_2)$ with right
		equivalence,
	\item	
	\item	A point $P$ is {\em incident} to a line $l$  iff\\
	\hth
	\enume

\newpage
	\sssec{Partial answer to}
	\vspace{-18pt}\hspace{128pt}{\bf \ref{sec-epsiinc1}.}\\[8pt]
	For $n = 7^2$, 2451 = 57.43, for $n = 9^2$, 6643 = 91.73,
	for $n = 11^2$, 14763 = 57.259, For $n = 13^2$, 28731 = 3.9577.  

	\sssec{Answer to}
	\vspace{-18pt}\hspace{94pt}{\bf \ref{sec-epsifano}.}\\[8pt]
	We have $Fano((0),(1),(2),(1'0),(2,1),(0,7),(1',12),
	[0],[1],[0,12],[2',0],[1',6],[0,-2],[12]).$

	\sssec{Answer to}
	\vspace{-18pt}\hspace{94pt}{\bf \ref{sec-eminipola}.}\\[8pt]
	$(7)\times(8)=[6], [6]\times[0] = (3)$,
	$(8)\times(0)=[1], [1]\times[7] = (2)$,
	$(0)\times(7)=[9], [9]\times[8] = (5)$,
	$\langle (3),(2),(5);[11]\rangle$.\\
	$(7)\times(8)=[6], [6]\times[0,1] = (1',7)$,
	$(8)\times(0,1)=[0,5], [0,5]\times[7] = (0',6)$,
	$(0,1)\times(7)=[2',6], [2',6]\times[8] = (2,5)$,
	$(2,5)$ is not incident to $(1',7)\times(0'6) = [2,2].$\\
	This has not been checked.

	From Dembowski, p. 129

	\sssec{Definition.}
	A linear ternary ring $(\Sigma,+,\cdot)$ is called a {\em cartesian
	field} iff $(\Sigma,+)$ is associative and is therefore a group.

	\sssec{Definition.}
	A cartesian field is called a {\em quasifield} iff the right 
	distributivity law holds:\\
	\hth$(x+y)z = xz+yz.$\\
	Artzy adds that $xa=xb+c$ has a unique solution, but this is a property
	(28). This is Veblen-Wedderburn.

	\sssec{Definition.}
	A quasifield is called a {\em semifield} iff the left
	distributivity law holds:\\
	\hth$z(x+y) = zx+zy.$

	\sssec{Definition.}
	A quasifield is called a {\em nearfield} iff
	$(\Sigma,\cdot)$ is associative and is therefore a group.

	\sssec{Definition.}
	A semifield is called a {\em alternative field} iff
	$x^2y = x(xy)$ and $xy^2 = (xy)y.$

	\sssec{Theorem.}
	{\em $P$ is $(p,L)$ transitive iff $P$ is $(p,L)$ Desarguesian.
	$p$ is point, $L$ is a line.} Dembowski p.123, 16

	Let $Q_0 = (79)$, $Q_1 = (80),$ $Q_2 = (90),$ and $U = (81)$
	then $q_2 = [79],$ $q_0 = [90],$ $q_1 = [80],$ $v = [88],$ 
	$i = [78],$ $V = (78),$ $I = (82),$ $j = [89]$, $W = (89)$,\\
	Points on $q_2:$ 86,12,25,38,51,64,77\\
	Points on $q_1:$ 85,11,24,37,50,63,76\\
	$11\times 12 = [ 7]:   	84(78,86), 8(51,25),43(77,64),48(86,38),
				52(12,12),54(25,78),66(64,80),72(38,51),\\
	 11\times 25 = [61]:    82(78,80), 0(64,64),18(12,51),28(51,78),
				33(77,12),58(38,86),61(86,25),62(25,38),\\
	 11\times 38 = [75]:	81(78,78), 4(38,25),14(12,80),19(77,86),
				36(51,64),70(64,38),73(25,51),74(86,12),\\
	 11\times 51 = [ 4]:	87(86,86), 1(38,80), 2(77,78),46(25,64),
				55(78,38),57(51,12),69(64,51),75(12,25),\\
	 11\times 64 = [ 8]:	83(86,78), 7(64,25),10(77,51),42(78,12),
				47(12,64),53(51,80),65(38,38),71(25,86),\\
	 11\times 77 = [68]:	88(86,80), 5(38,64),13(51,51),21(77,38),
				30(25,12),32(64,86),67(12,78),68(78,25),$\\
	Coordinates of points:\\
	$\begin{array}{rccccccccc}
	( 0)&\vline&64,64&38,80&77,78&78,51&38,25&38,64&51,86&64,25\\
	( 8)&\vline&51,25&86,77&77,51&80,77&12&51,51&12,80&64,78\\
	(16)&\vline&78,77&12,38&12,51&77,86&51,38&77,38&86,64&64,77\\
	(24)&\vline&80,64&25&77,77&25,80&51,78&78,64&25,12&25,77\\
	(32)&\vline&64,86&77,12&64,12&86,51&51,64&80,51&38&25,25\\
	(40)&\vline&77,80&38,78&78,12&77,64&77,25&12,86&25,64&12,64\\
	(48)&\vline&86,38&38,12&80,38&51&12,12&51,80&25,78&78,38\\
	(56)&\vline&51,77&51,12&38,86&12,77&38,77&86,25&25,38&80,25\\
	(64)&\vline&64&38,38&64,80&12,78&78,25&64,51&64,38&25,86\\
	(72)&\vline&38,51&25,51&86,12&12,25&80,12&77&78&80,80\\
	(80)&\vline& 0&78,78&78,80&86,78&78,86&80,86&86&86,86\\
	(88)&\vline&86,80&80,78&\infty\\
	\end{array}$

	Coordinates of lines:\\
	$\begin{array}{rccccccccc}
	[ 0]&\vline&78,38&38&77,78&12,12&51,77&86,38&12,86&12,77\\
	{[} 8]&\vline&64,77&77,25&64,12&80,25&51,80&78,12&12&64,78\\
	{[}16]&\vline&25,25&77,64&86,12&25,86&25,64&51,64&64,38&51,25\\
	{[}24]&\vline&80,38&77,80&78,25&25&51,78&38,38&64,51&86,25\\
	{[}32]&\vline&38,86&38,51&77,51&51,12&77,38&80,12&64,80&78,77\\
	{[}40]&\vline&77&38,78&51,51&12,38&86,77&51,86&51,38&25,38\\
	{[}48]&\vline&38,64&25,51&80,64&12,80&78,51&51&25,78&64,64\\
	{[}56]&\vline&38,25&86,51&64,86&64,25&12,25&25,77&12,64&80,77\\
	{[}64]&\vline&38,80&78,64&64&12,78&77,77&25,12&86,64&77,86\\
	{[}72]&\vline&77,12&38,12&12,51&38,77&80,51&25,80&78,80&\infty\\
	{[}80]&\vline&80&78,78&86&86,78&86,86&80,86&86,80&78,86\\
	{[}88]&\vline&78&80,78&80,80\\
	\end{array}$

	$[80]:11(80,77),24(80,64),37(80,51),50(80,38),63(80,25),76(80,12),
		79(80,80),85(80,86),89(80,78),90(/infty)$
	$B = A+\alpha,$ $(A,B)\incid[V,Y]$, $ V = (78)$,
		$(76)=(80,12)=(\un{0},\alpha),$
		$Y\times V = (78)\times(76)=[13].$\\
	$[13]:78(78),15(64,78),18(12,51),19(77,86),68(78,25),76(80,12),
		46(25,64),48(86,38),53(51,80),60(38,77)$\\
	hence $12 = 80+\alpha,$ $51 =\alpha+\alpha=-\alpha,$
	$25 =78+\alpha=1+\alpha=\gamma,$ $64=25+\gamma=\gamma+\gamma=-\beta,$
	$38=86+\alpha=-1+\alpha=\beta,$ $77=38+\alpha=\beta+\alpha=-\gamma.$\\
	$\begin{array}{rrrrrrrrrr}
	\infty&0&1&-1&\alpha&-\alpha&\beta&-\beta&\gamma&-\gamma\\
	90&80&78&86&12&51&38&64&25&77
	\end{array}$

	$[\alpha,0]=(12)\times(80,80) = (12)\times(79)=[51],$
	$(a,b)\incid[51]\implies b = a\cdot \alpha.$\\
	$(42)=(78,12)=(\un{1},\alpha)\implies \alpha = 1\times\alpha,$\\
	$(45)=(12,86)-(\alpha,\un{-1})\implies -1 = \alpha\times\alpha,$\\
	$(23)=(64,77)=(-\beta,-\gamma) \implies -\gamma=-\beta\times\alpha,$\\
	$(28)=(51,78) = (-\alpha,\un{1})\implies 1 =-\alpha\times\alpha,$\\
	$(35)=(86,51)=(\un{-1},-\alpha)\implies -\alpha= \un{-1}\times\alpha,$

	Using
	DATA 6,0, 6,4, 6,10, 6,12, 0,0, 1,0, 2,0, 3,0, 4,0, 5,0
	DATA 6,0, 0,7, 0,8, 0,11, 3,2, 3,10, 4,4, 4,6, 5,5, 5,12
	gives the same multiplication table give left not right
	distibutive law
	with $Q_i = 79,81,87$, $U = (83),$ $\alpha = (12)$,
	$q_0 = [87] = \{79,81,82,86, 4,17,\ldots\},$\\
	$q_1 = [81] = \{79,85,87,88,10,23,\ldots\},$\\
	$q_2 = [79] = \{81,87,89,90,12,25,\ldots\},$\\
	with case 7, data 79,81,87,83,12:\\
	$\infty = 87,$ $\un{0} = 81,$ $\un{1} = 89,$ $\un{-1} = 90,$
	$\alpha = 12,$ $-\alpha = 77,$ $\beta = 38,$ $-\beta = 51,$
	$\gamma = 25,$ $-\gamma = 64.$

	This is a try for a section to be included in g19.tex between
	Moufang and Desargues.

	\section{Axiomatic.}
	\ssec{Veblen-MacLagan planes.}
\setcounter{subsubsection}{-1}
	\sssec{Introduction.}
	The first example of a Veblen-Wedderburn plane was given in 1907
	by Veblen and MacLagan-Wedderburn.  It is associated to the
	algebraic structure of a nearfield, which is a skew field
	which lacks the left distributive law, hence is an other plane
	between the Veblen-Wedderburn plane and the Desarguesian plane.

	\sssec{Axiom. [Da]}
	\vspace{-19pt}\hspace{90pt}\footnote{Da for Desargues leading
	to associativity of multiplication.}\\[8pt]
	Given a Veblen-Wedderburn plane, 2 points $Q_1$ and $Q_2$ on the
	ideal line and an other point $Q_0$ not on it, any 2 parallelograms
	$A_i$ and $B_i$ with directions $Q_1$ and $Q_2$, with
	no sides in common \ldots,???,
	such that $A_j$ and $B_j$ are perspective from $Q_0$ for j = 0 To 2,
	imply that $A_3$ and $B_3$ are perspective from $Q_0$.

	\sssec{Notation.}
	Da$(\{Q_0,Q_1,Q_2\},\{A_j\},\{B_j\}).$

	\sssec{Definition.}
	A {\em Veblen-MacLagan plane} is a Veblen-Wedderburn plane in which
	the axiom Da is satisfied.

	\sssec{Lemma. [For Associativity]}
	H1.0.\hti{3}$A_0,$ $a_{12},$ $x,$ (See Fig. 2?.)\\
	D1.0.\hti{3}$a_{01} := Q_1 \times A_0,$ $a_{02} := Q_2 \times A_0,$\\
	D1.1.\hti{3}$A_1 := a_{01} \times a_{12},$
		$A_2 := a_{02} \times a_{12},$\\
	D1.2.\hti{3}$a_{13} := Q_2 \times A_1,$ $a_{23} := Q_1 \times A_2,$
		$A_3 := a_{13} \times a_{23},$\\
	D2.0.\hti{3}$a_0 := Q_0 \times A_0,$ $a_1 := Q_0 \times A_1,$
		$a_2 := Q_0 \times A_2,$ $a_3 := Q_0 \times A_3,$
	D2.1.\hti{3}$B_0 := a_0 \times y,$
		$b_{01} := Q_1 \times B_0,$ $b_{02} := Q_2 \times B_0,$\\
	D2.1.\hti{3}$B_1 := b_1 \times b_{01},$ $B_2 := b_2 \times b_{02},$
	D2.2.\hti{3}$b_{13} := Q_2 \times B_1,$ $b_{23} := Q_1 \times B_2,$
		$B_3 := b_{13} \times b_{23},$\\
	C1.0.\hti{3}$B_3 \incid b_3,$\\
	{\em Moreover}\\
	$A_0 = (A,B),$ $A_1 = (A',B),$ $A_2= (A,B'),$ $A_3 = (A',B'),$
	$B_0 =$\\ 
	Proof:
	Da$(\{Q_0,Q_1,Q_2\},\{A_j\},\{B_j\}).$

	\sssec{Theorem.}
	{\em In a Veblen-MacLagan plane, the ternary ring $(\Sigma,*)$
	is a nearfield.}:
	\enumb
	\item $(\Sigma,+)$ {\em is an Abelian group,}
	\item $(\Sigma-\{0\},\cdot)$ {\em is a group,}
	\item ($\Sigma,*) = (\Sigma,+,\cdot)$ {\em is right distributive,}
		$(a+b)\cdot c = a \cdot c + b\cdot c$.
	\enume

	\ssec{Examples of Perspective planes.}
	\sssec{Theorem.}
	\enumb
	\item {\em The Cayleyian plane is not a Veblen-MacLagan plane.}
		\marginpar{replace Desarg.?}
	\enume

	\sssec{Definition.}
	A {\em miniquaternion plane} \ldots.

	\sssec{Theorem.}
	\enumb
	\item {\em A miniquaternion plane is a Veblen-MacLagan plane.}
	\item {\em A miniquaternion plane is not a Moufang plane.}
	\enume

	\sssec{Tables.}
	The following are in an alternate notation the known table for $p$ = 3
	and a new table for $p$ = 5.
	The other incidence are obtained by adding one to the subscripts of
	the lines and subtracting one for the subscript of the points.

	\section{Bibliography.}
\begin{enumerate}
\item Artzy, Rafael, Linear Geometry, Reading Mass., Addison-Wesley, 1965,
	273 pp.
\item Baumert, Leonard D., Cyclic Difference Sets, N. Y., Springer, 1971.
\item Dembowski, Peter, Finite Geometries, Ergebnisse der Mathematik und ihrer
	Grenzgebiete, Band 44, Springer, New-York, 1968, 375 pp.
\item Dickson, G\"{o}ttingen Nachrichten, 1905, 358-394.
\item Hall, Marshall, Projective Planes , Trans. Amer. Math. Soc., Vol. 54,
	1943, 229-277.
\item Hughes, D. R., A class of non-Desarguesian projective planes,
	Canad. J. of Math., Vol. 9, 1957, 378-388. (I,9.p.1)
\item Lemay, Fernand, Le dod\'{e}ca\`{e}dre et la g\'{e}om\'{e}trie projective
	d'ordre 5, p. 279-306 of Johnson, Norman L., Kallaher, Michael J.,
	 Long Calvin T., Edit.,	Finite Geometries, N. Y., Marcel Dekker Inc.
	1983.
\item Maclagan-Wedderburn, J. H., {\em A theorem on finite algebras,}
	Trans. Amer. Math. Soc., Vol. 6, 1905, 349. 
	(A finite skew-field is a field)
\item Moore, E. H., Mathematical Papers, Chicago Congress, 1893, 210-226.
	(Def. of Galois Fields.)
\item Room, Thomas Gerald and Kirkpatrick P. B., Miniquaternion geometry; an
	introduction to the study of projective planes, Cambridge [Eng.]
	University Press, 1971.
\item Rosati, L. A., I gruppi di collineazioni dei pliani di Hughes. Boll. Un.
	Mat. Ital. Vol. 13, 505-513.
\item Singer, James, A Theorem in Finite Projective Geometry and some 
	applications to number Theory, Trans. Amer. Math. Soc., Vol. 43, 1938,
	377-385 
\item Vahlen, Karl Theodor, {\em Abstrakte Geometrie, Untersuchungen
	\"{u}ber die Grundlagen der Euklidischen und nicht-Euklidischen
	Geometrie}, Mit 92 abbildungen im text, 2., neubearb. aufl.
	Leipzig, S. Hirzel, 1940, Series title: Deutsche mathematik, im
	auftrage der Deutschen forschungsgemeinschaft, 2. beiheft.
\item Veblen, Oswald \& Bussey, W. H., {\em Finite Projective Geometries},
	Trans. Amer. Math. Soc., Vol. 7, 1906, 241-259. (On PG($n,p^k$))
\item Veblen, Oswald, and Wedderburn, Jospeh Henri MacLagan,
	Non-Desarguesian and non-Pascalian geometries, Trans. Amer. Math.
	Soc., Vol. 8, 1907, 379-388.
\item Zappa, G. Sui gruppi di collineazioni dei paini di Hughes, Boll. Un.
	Mat. Ital. (3), Vol. 12, 1957, p. 507-516.
\enume

\newpage
\setcounter{section}{89}
	\section{Answer to problems and Comments.}

	\sssec{Notation.}\label{sec-nqua1a}
	\hth$u_{ij} := q_iq_j^{-1} - r_ir_j^{-1},\\
	\hth v_{ij} := q_iq_j^{-1} + r_ir_j^{-1},\\
	\hth s_i := -v_{i,i-1}^{-1}v_{i,i+1},\\
	\hth a_i := -v_{i-1,i}^{-1}u_{i-1,i+1},\\
	\hth t_i := s_{i+2} s_{i+1},\\
	\hth f_i := s_i - s_{i+1}^{-1}s_{i-1}^{-1},\\
	\hth g_i := t_i^{-1} - t_{i+1}t_{i-1},\\
	\hth k_i := q'_{i-1}\ov{q}_{i+1},\\
	\hth l_i := r'_{i-1}\ov{r}_{i+1}.$

	\sssec{Exercise.}
	Prove $q_1^{-1}u_{12}r_2 + q_2^{-1}u_{20}r_0 + q_0^{-1}u_{01}r_1 = 0,$\\
	associated with $M \cdot eul = 0.$

	The proof follows form substitution of $u_ij$ by their definition.

	\sssec{Exercise.}
	Prove $u_{12}u_{02}^{-1}u_{01} = -u_{10}u_{20}^{-1}u_{21},$
	associated with 2 equivalent forms of $eul$ one for which the first
	coordinate is one and the other obtain by "rotation", the second
	coordinate being one.

	Form the definition of $u_{02}$ it follows by multiplication to the
	right or left by $u_{02}^{-1}$, that\\
	\hth$ q_0q_2^{-1}u_{02}^{-1} - r_0r_2^{-1}u_{02}^{-1} = 1,\\
	\hth u_{02}^{-1}q_0q_2^{-1} - u_{02}^{-1}r_0r_2^{-1} = 1.$\\
	Moreover,\\
	\hth$u_{20} = q_2q_0^{-1} - r_2r_0^{-1} = -q_2q_0^{-1}u_{02}r_2r_0^{-1}$
\\
	or\\
	\hth$u_{20}^{-1} = r_0r_2^{-1} u_{02}^{-1} q_0q_2^{-1}.$

	If we substitute in the identity to prove, with both terms in the second
	member, $u_{12},$ $u_{01},$ $u_{10}$ and $u_{21},$ by their definition,
	we get\\
	$q_1q_2^{-1}u_{02}^{-1}q_0q_1^{-1} - q_1q_2^{-1}u_{02}^{-1}r_0r_1^{-1}\\
	-r_0r_1^{-1}u_{02}^{-1}q_1q_2^{-1} + r_0r_1^{-1}u_{02}^{-1}r_1r_2^{-1}\\
	- q_1q_0^{-1}r_0r_2^{-1}u_{02}^{-1}q_0q_2^{-1}q_2q_1^{-1}
	+ q_1q_0^{-1}r_0r_2^{-1}u_{02}^{-1}q_0q_2^{-1}r_2r_1^{-1}\\
	+ r_1r_0^{-1}r_0r_2^{-1}u_{02}^{-1}q_0q_2^{-1}q_2q_1^{-1}
	- r_1r_0^{-1}r_0r_2^{-1}u_{02}^{-1}q_0q_2^{-1}r_2r_1^{-1} = 0,$\\
	because terms 3 and 7 cancel, terms 1 and 5 as well as 4 and 8 give 1
	and -1, terms 2 and 6 give 0 by application of the identities
	given at the begginning of the proof.

	\sssec{Lemma.}\label{sec-lforTa}
	\enumb
	\item$	norm(s_0s_1s_2) = 1.$
	\item$	norm(t_0t_1t_2) = 1.$
	\item$s'_2\ov{f}_2s_0^{-1} = -f_2s_1.$
	\item$\ov{t}_2\ov{g}_2t_0 = -g_2t_1^{-1}.$
	\enume

	Proof: For 0, we use Lemma \ref{sec-lforTLinv} and obtain 1, from the
	definition of $t_i$. For 2, we substitute $f_2$ by its definition and
	compare the terms of both sides of the equality which have the same
	sign.

	\sssec{Proof of }
	\vspace{-19pt}\hspace{84pt}{\bf\ref{sec-fundHDC}}.\\[8pt]
	Let\\
	G0.0.	$A_0 = (1,0,0),$ $A_1 = (0,1,0),$ $A_2 = (0,0,1),$\\
	G0.1.	$M = (q_0,q_1,q_2),$ $\ov{M} = (r_0,r_1,r_2),$\\
	then\\
	P1.0.	$a_0 = (1,0,0),$ $a_1 = (0,1,0),$ $a_2 = (0,0,1),$\\
	P1.1.	$ma_0 = [0,q'_1,-q'_2],$ $\ov{m}a_0 = [0,r'_1,-r'_2],$\\
	P1.2.	$M_0 = (0,q_1,q_2),$ $\ov{M}_0 = (0,r_1,r_2),$\\
	P1.3.	$eul = [1,-u'_{12}\ov{u}_{02},-u'_{21}\ov{u}_{01}],$\\
	P1.4.	${\bf S} = \matb{(}{ccc}q_0^{-n}&-q'_0q_1^{-1}&-q'_0q_2^{-1}\\
		-q'_1q_0^{-1}&q^{-n}_1&-q'_1q_2^{-1}\\
		-q'_2q_0^{-1}&-q'_2q_1^{-1}&q_2^{-n}\mate{)},$
	${\bf S}^{-1} = \matb{(}{ccc}0&q_0\ov{q}_1&q_0\ov{q}_2\\
		q_1\ov{q}_0& 0& q_1\ov{q}_2\\
		q_2\ov{q}_0&q_2\ov{q}_1& 0\mate{)}.\\
	\htp{28}{\bf \ov{S}} = \matb{(}{ccc}r_0^{-n}&-r'_0r_1^{-1}&
		-r'_0r_2^{-1}\\-r'_1r_0^{-1}&r^{-n}_1&-r'_1r_2^{-1}\\
		-r'_2r_0^{-1}&-r'_2r_1^{-1}&r_2^{-n}\mate{)},$
	${\bf \ov{S}}^{\,-1} = \matb{(}{ccc}0&r_0\ov{r}_1&r_0\ov{r}_2\\
		r_1\ov{r}_0& 0& r_1\ov{r}_2\\
		r_2\ov{r}_0&r_2\ov{r}_1& 0\mate{)}.$

\noindent P2.0.	$mm_0 = [-q'_0,q'_1,q'_2],$ $\ov{m}m_0 = [-r'_0,r'_1,r'_2],$\\
	P2.1.	$MA_0 = (0,q_1,-q_2),$ $\ov{M}A_0 = (0,r_1,-r_2),$\\
	P2.2.	$m_0 = [0,q'_1,q'_2],$ $\ov{m}_0 = [0,r'_1,r'_2],$\\
	P2.3.	$MM_0 = (-q_0,q_1,q_2),$ $\ov{M}M_0 = (-r_0,r_1,r_2),$\\
	P2.4.	$m = [q'_0,q'_1,q'_2],$ $\ov{m} = [r'_0,r'_1,r'_2],$\\
	P2.5.	$Ima_0 = (-2q_0,q_1,q_2)$,
		$\ov{I}ma_0 = (-2r_0,r_1,r_2),$\\
	\htp{28}$I\ov{m}a_0 = (-q_0(q_1^{-1}r_1+q_2^{-1}r_2),r_1,r_2),$
		$\ov{Im}a_0 = (-r_0(r_1^{-1}q_1+r_2^{-1}q_2),q_1,q_2),$\\
	P2.6.	$iMA_0 = [2q'_0,-q'_1,-q'_2],$
		$\ov{\imath}MA_0 = [2r'_0,-r'_1,-r'_2],$

\noindent P3.0.	$mf_0 = [k_1v'_{21}\ov{u}_{21},k_0^{-1},-1],$
		$\ov{m}f_0 = [l_1v'_{21}\ov{u}_{21},-l_0^{-1},1],$\\
	P3.1.	$O = [],$
		$\ov{O} = [],$\\
	P3.2.	$Mfa_0 = (\ov{k}_2,0,-k'_{0}v_{10}u_{10}^{-1}),$
		$\ov{M}fa_0 = (\ov{l}_2,0,l'_{0}v_{10}u_{10}^{-1}),$\\
	\htp{28}	$Mf\ov{a}_0 = (1,k'_2\ov{v}_{02}u'_{02},0),$
		$\ov{M}f\ov{a}_0 = (1,-l'_2\ov{v}_{02}u'_{02},0),$\\
	P3.3.	$mfa_0 = (\ov{k}_1v'_{21}\ov{u}_{21},k_0^{-1},0),$
		$\ov{m}fa_0 = (\ov{l}_1v'_{21}\ov{u}_{21},-l_0^{-1},0),$\\
	\htp{28}	$mf\ov{a}_0 = (\ov{k}_1v_{21}^{-1}\ov{u}_{21},0,-1),$
		$\ov{m}f\ov{a}_0 = (\ov{l}_1v_{21}^{-1}\ov{u}_{21},0,1),$\\
	P3.4.	$Mfm_0 = (k_1^{-1}v_{21}u_{21}^{-1},-\ov{k}_0,1),$
		$\ov{M}fm_0 = (l_1^{-1}v_{21}u_{21}^{-1},\ov{l}_0,-1),$

\noindent P4.0.	$Imm_0 = (a_0,1,s_0),$ $\ov{I}mm_0 = (-a_0,1,s_0),$\\
	P4.1.	$ta_0 = [0,1,-s'_0],$\\
	P4.2.	$T_0 = (1,s_2,s_1^{-1}),$\\
	P4.3.	$at_0 = [0,s'_2,-\ov{s}_1] = [0,1,-\ov{t}_0],$\\
	P4.4.	$K_0 = (1,t_2^{-1},t_1),$\\
	P4.5.	$Taa_0 = (0,1,s_0),$\\
	P4.6.	$poK_0 = [-1,s'_2,\ov{s}_1],$

\noindent P4.7.	${\bf T} = \matb{(}{ccc}0&\ov{f}_2&-\ov{f}_2s_0^{-1}\\
			f_2&0&-f_2s_1\\-s'_0f_2&-\ov{s}_1\ov{f}_2&0\mate{)},$
		${\bf T}^{-1} = \matb{(}{ccc}1&\ov{s}_2&s'_1\\
			s_2&s^n_2&s^{-n}_1s^{-1}_0\\
			s^{-1}_1&s^n_2s_0&s^{-n}_1\mate{)}.$\\
	P4.8.	${\bf L} = \matb{(}{ccc}0&\ov{g}_2&-\ov{g}_2t_0\\
			g_2&0&-g_2t_1^{-1}\\
			-\ov{t}_0g_2&-t'_1\ov{g}_2&0\mate{)},$
		${\bf L}^{-1} = \matb{(}{ccc}1&t'_2&\ov{t}_1\\
			t_2^{-1}&t_2^{-n}&t_1^nt_0\\
			t_1&t_1^n\ov{t}_0&t_1^n\mate{)}.$

	Proof:\\
	For P4.0, if the coordinates of $Imm_0$ are $x_0$, 1 and $x_2$,
	we have to solve\\
	\hth$q_0^{-1}x_0 + q_1^{-1}+q_2^{-1}x_2 = 0,\\
	\hth -r_0^{-1}x_0 + r_1^{-1}+r_2^{-1}x_2 = 0.$\\
	Multiplying the equations to the left respectively by $q_2$ and $-r_2,$
	or by $q_0$ and $r_0$ and adding gives $x_0$ and $x_2$ using the
	notation \ref{sec-nqua1a}.\\
	For P4.7, it is easier to obtain ${\bf T}^{-1}$ first, the columns are
	$T_0,$ $T_1$, $T_2,$ multiplied to the right by 1, $s_2^n,$ $s_1^{-n}.$
	The matrix ${\bf T}$ is then obtained using Theorem \ref{sec-thermquat},
	multiplying by $-s_1^{-n}.$  The equivalence with the matrix whose
	columns ate $ta_i$ can be verified using Lemma \ref{sec-lforTa}.2.
	A similar proof gives P4.8.

	\sssec{Theorem.}
	{\em The product of the diagonal elements of ${\bf T}^{-1}$ and of
	${\bf L}^{-1}$ is the same.}

	This follows from Lemma \ref{sec-lforTLinv}.

	The correspondance between the definitions in EUC and here is as
	follows
	\begin{tabbing}
	DD0.0 \= DD1.0\=\vline\=DD0.1\=DD1.1\=\vline\=DD0.2\=DD1.2\=\vline
	\=DD1.0\=DD1.3\=\vline\=DD6.12\=DD1.4\kill
	D0.0 \> D1.0\>\vline\>D0.1\>D1.1\>\vline\>D0.2\>D1.2\>\vline
	\>D1.0\>D1.3\>\vline\>D36.12\>DC1.4\\
	D0.3 \> D2.0\>\vline\>D0.4\>D2.1\>\vline\>D0.5\>D2.2\>\vline
	\>D0.6\>D2.3\>\vline\>D0.7\>D2.4\\
	D10.3 \> D2.5\>\vline\>D0.25\>D2.6\>\vline\>D10.3\>D2.7\>\vline
	\>D6.0\>D3.0\>\vline\>D6.4\>D3.1\\
	? \> D3.2\>\vline\>?\>D3.3\>\vline\>D14.0?\>D3.4\>\vline
	\>D1.6\>D4.0\>\vline\>D1.7\>D4.1\\
	D1.8 \> D4.2\>\vline\>D12.1\>D4.3\>\vline\>D1.2\>D4.4\>\vline
	\>D1.4\>D4.4\>\vline\>D1.9\>D4.5\\
	D15.12? \> D4.6\>\vline\>D1.19\>DC4.7\>\vline\>\>DC4.8\\
	\end{tabbing}

\chapter{FUNCTIONS OVER FINITE FIELDS}

\setcounter{section}{-1}
	\section{Introduction.}

	\sssec{Notation.}
	The first notation is standard, the second is useful for Theorem
	\ref{sec-tfinortho1}.2.1.\\
	\hth$(a)_i := \prod_0^{i-1} (a+i) = a(a+1)\ldots(a+i-1)$.\\
	\hth$[a]_i := \prod_0^{i-1} (a+2i) = a(a+2)\ldots(a+2i-2)$.

	\sssec{Notation.}
	The following notation, favored on the European
	continent, but seldom used elsewhere, is quite useful:\\
	\hth$	0 !!	  := 1$\\
	\hth$	2n !!	  := 2 . 4 .  \ldots  . 2n.$\\
	\hth$	(2n+1) !! := 1 . 3 .  \ldots  . (2n+1).$

	\section{Polynomials over Finite Fields.}

	\ssec{Definition and basic properties.}

	\setcounter{subsubsection}{-1}
	\sssec{Introduction.}
	In a finite field, we can define polynomials of degree
	up to $p-1.$  These are determined by their values at $i$ in ${\Bbb Z}_p.$  If
	these are defined in the real field with rational coefficients, the
	definition and properties automatically extend to the finite field.

	\sssec{Definition.}\label{sec-dfinpol}
	A {\em polynomial} is a function
	\hth$a_0 I^{p-1} + a_1 I^{p-2} +  \ldots  + a_{p-1}$\\
	which associates to
	$x \in {\Bbb Z}_p$ the integer\\
	\hth$	a_0 x^{p-1} + a_1 x^{p-2} +  \ldots  + a_{p-1}.$\\
	The {\em polynomial is of degree} $k$ iff $a_0 \ldots a_{p-k-1}$
	are congruent to $0$ modulo $p$ and $a_{p-k}$ is not.

	\sssec{Theorem. [Lagrange]}
	{\em Given $k+1$ distinct integers $x_i$ (modulo $p$), $k<p,$ and given
	$k+1$ integers $f_i,$ $\exists$  a polynomial $P$ of degree $k$ $\ni$ \\
	\hth$	P(x_i) = f_i,$ $i = 0$ to $k.$}

	\ssec{Derivatives of polynomials.}

	\sssec{Definition.}
	The {\em derivative of the polynomial} \ref{sec-dfinpol}.0 is
	\hth$(p-1)a_0 I^{p-2} + (p-2)a_1 I^{p-3} +  \ldots  + a_{p-2}.$
	\section{Orthogonal Polynomials over Finite Fields.}

	\setcounter{subsection}{-1}
	\ssec{Introduction.}
	The main purpose of writing this Chapter is connected with interesting
	symmetry properties of the orthogonal polynomials, in ${\Bbb Z}_p.$
	In the classical theory there is a scaling factor which is
	arbitrairely chosen for each of the families of orthogonal polynomials.
	For some time now, the same scaling factor, in each case, is universely
	used.  When determining values of the Chebyshev polynomials for some
	small values of $p,$ I was struck by the symmetry properties given
	in \ref{sec-tfincheb} and \ref{sec-tfincheb1}.  These properties are
	dependent on the
	scaling factors and it turns out that the unanimously accepted ones are
	essentially the only ones giving this property.  The same property has
	been found for the polynomials of Legendre and of Laguerre.  For the
	polynomial of Hermite this is not the case.  I have succeeded in
	obtaining some scaling, given in \ref{sec-dfinherm} for which a
	symmetry can be
	obtained.  This scaling is given by expressions which are different for
	the even and for the odd Hermite polynomials, therefore the recurrence
	relation has a constant whose expression differs for even and odd
	indices.  It is therefore possible to give an a-posteriori
	justification of the scaling factor for the classical polynomial, and
	there is some reason to introduce a different scaling for the Hermite
	polynomials.  The case of the Jacobi polynomials with 2 parameters $a$ 
	and $b$ is left as an exercise. With $a = b$, again a scaling is
	required to obtain symmetry.

	\ssec{Basic Definitions and Theorems.}

	\setcounter{subsubsection}{-1}
	\sssec{Introduction.}
	For orthogonal polynomials, recurrence relations, differential
	equations and values of the coefficients generalize automatically, from
	the classical case.  Therefore, we have the definitions
	\ref{sec-dfinortho} and the theorems \ref{sec-tfinortho1}
	to \ref{sec-tfinortho3}.

	\sssec{Definition.}\label{sec-dfinortho}
	The {\em polynomials of Chebyshev of the first} $(T_n)$ and {\em of
	the second kind} $(U_n),$ {\em of Legendre} $(P_n),$ {\em of Laguerre}
	$(L_n)$ and {\em of Hermite} $(H_n)$ are
	defined by the recurrence relations:
	\enumb
	\item[0.~~]$	T_0 := 1,$ $T_1 := I,$
		$T_{n+1} := 2(2 I - 1) T_n - T_{n-1},$
	\item[1.~~]$	U_0 := 1,$ $U_1 := 2I,$
		$U_{n+1} := 2(2 I - 1) U_n - U_{n-1},$
	\item[2.~~]$P_0^{(a)} := 1,$ $P_1^{(a)} := I,$\\
	\hth$(n+2a+1)\: P_{n+1}^{(a)} := (2n+2a+1) I\: P_n^{(a)}
		- n\: P_{n-1}^{(a)},$ $a \geq 0,$ $n < p-1-2a,$\\
	\item[3.0.]$	L_0 := 1,$ $L_1 := 1 - I,$\\
	\hth$	(n+1) L_{n+1} := (2n+1-I) I L_n - n L_{n-1},$ $n < p-1,$\\
	\item[  1.]$	L_0^{(a)} := 1,$ $L_1^{(a)} := a+1 - I,$\\
	\hth$(n+1) L_{n+1}^{(a)} := (2n+a+1-I) L_n^{(a)} - (n+a) L_{n-1}^{(a)},$
	$n < p-1,$
	\item[4.~~]$H_0 := 1,$ $H_1 := 2I,$ $H_{n+1} := 2I H_n - 2n H_{n-1}.$
	\enume

	The Legendre polynomial is $P_n := P_n^{(0)}$ and
	$P_n^{(a)} := \frac{n!a!}{(n+a)!}\:P_n^{(a,b)},$ where
	$P_n^{(a,b)},$ are the polynomials of Jacobi, scaled so that
	$P_n^{(a)}(1) = 1.$\\
	$L_n^{(a)}$ are the generalized Laguerre polynomials and
	$L_n = L_n^{(0)}$.

	See for instance Handbook of Mathematical functions, p. 782.

	\sssec{Theorem.}\label{sec-tfinortho1}
	{\em If $X_{n,j}$ denotes the coefficient of $I^j$ in the polynomial
	$X_n,$ }
	\enumb
	\item[0.~~]$T_{n,n-2j} = \frac{1}{2}n 2^{(n-2j)} 
		(-1)^j \frac{(n-j-1)!}{j! (n-2j)!},$
	\item[1.~~]$U_{n,n-2j} = \frac{1}{2}n 2^{(n-2j)} 
	(-1)^j \frac{(n-j)!}{j! (n-2j)!},$
	\item[2.0.]$P_{n,n-2j} = 2^{(-n)}
	(-1)^j \frac{(2n-2j)!}{j! (n-j! (n-2j)!},$
	\item[  1.]$P_{2n,2j}^{(a)} = (-1)^(n-j)\ve{n}{j}
		\frac{[2a+2n+1]_j[2j+1]_{n-j}}{[2a+2]_n},\\
		P_{2n,2j}^{(a)} = (-1)^(n-j)\ve{n}{j}
		\frac{[2a+2n+3]_j[2j+3]_{n-j}}{[2a+2]_n},$
	\item[3.0.]$L_{n,j} = (-1)^j \frac{n!}{(n-j)! j!^2},$\\
	\item[	1.]$L_{n,j}^{(a)} = (-1)^j \frac{(n+a)!}{(n-j)! (a+j)! j!},$
	\item[4.~~]$H_{n,n-2j} = n! 2^{(n-2j)} (-1)^j \frac{1}{j! (n-2j)!},$
	\enume

	See for instance Handbook of Mathematical functions, p. 775.

	\sssec{Theorem.}\label{sec-tfinortho2}
	{\em The polynomials of Chebyshev, of the first $(T_n)$ and of
	the second kind $(U_n),$ of Legendre $(P_n),$ of Laguerre
	$(L_n)$ and of Hermite $(H_n)$ satisfy by the differential equations}
	\enumb
	\item[0.~~]$(1-I^2) D^2 T_n - I D T_n + n^2 T_n = 0,$
	\item[1.~~]$(1-I^{2)} D^2 U_n - 3 I D U_n + n(n+2) U_n = 0,$
	\item[2.~~]$(1-I^{2)}\: D^2 P_n^{(a)} - 2(a+1) I) \:D P_n^{(a)}
				   + n(n+2a+1)\: P_n^{(a)} = 0,$
	\item[3.0.]$I D^2 L_n + (1 - I) D L_n + n L_n = 0.$\\
	\item[	1.]$I D^2 L_n^{(a)} + (a+1 - I) D L_n^{(a)} + n L_n^{(a)} = 0.$
	\item[4.~~]$D^2 H_n - 2 I DH_n + 2n H_n = 0.$
	\enume

	See for instance Handbook of Mathematical functions, p. 781.

	\sssec{Theorem.}\label{sec-tfinortho3}
	\enumb
	\item[0.~~]$T_n(1) = 1,$ $DT_n(1) = n^2.$\\
	\item[1.~~]$U_n(1) = n+1,$ $DU_n(1) = \frac{n(n+1)(n+2)}{3}.$
	\item[2.~~]$P_n^{(a)}(1) = 1,$
			$DP_n^{(a)}(1) = - \frac{n(n+2a+1)}{2(a+1)}.$
	\item[3.0.]$L_n(0) = 1,$ $DL_n(0) = -n.$\\
	\item[	1.]$L_n^{(a)} = \ve{n+a}{n},$ $DL_n^{(a)} = -\ve{n+a}{n-1}.$
	\item[4.~~]$H_{2n}(0) = (-1)^n \frac{(2n)!}{n!},$ $DH_{2n}(0) = 0.$
	\hth$	H_{2n+1}(0) = 0,$ $DH_{2n+1}(0) = (-1)^{n+1} \frac{(2n)!}{n!}.$
	\enume

	\sssec{Comment.}
	It is easy to verify that, contrary to the classical case, the roots of
	the orthogonal polynomials are not necessarily in ${\Bbb Z}_p.$ 
	For instance, $T_2$ has a root in ${\Bbb Z}_p$ iff $2$ R $p$ or 
	$p \equiv \pm 1 \pmod{8}.$

	\sssec{Program.}
	[m130]FIN\_ORTHOG.HOM illustates the use of the program
	[m130]FIN\_ORTHOG.BAS, which determines these various orthogonal
	polynomials.  [m130]FIN\_ORTHOG.NOT are notes tracing some of the steps
	leading to the conjectures proven here.

	\ssec{Symmetry properties for the Polynomials of Chebyshev of the
	first and second kind.}

	\sssec{Theorem.}\label{sec-tfincheb}
	\enumb
	\item$	T_{p+i,j} = T_{p-i,j}.$
	\item$	T_{i+2pk,j} = -T_{i+pk,j} = T_{i,j},$ $j<p.$
	\enume

	Proof:\\
	\hth$	T_{p+i,j} = (-1)^{\frac{1}{2} (p+i-j)} 2^j
	\frac{(\frac{1}{2} (p+i+j)-1)! \frac{1}{2} (p+i)}
	{(\frac{1}{2} (p+i-j))! j!}$\\
	\hti{12}$	= (-1)^{\frac{1}{2} (p+i-j)} (-1)^{\frac{1}{2} (p-i-j)}
	(-1)^{\frac{1}{2} (p-i+j+i)}
		\frac{2^j(\frac{1}{2} (p-i+j-2))! \frac{1}{2} (p-i)}
		{(\frac{1}{2} (p-i-j))! j!}$\\
	\hti{12}$	= (-1)^{\frac{1}{2} (p-i-j)} 2^j 
	\frac{(\frac{1}{2} (p-i+j-2))! \frac{1}{2} (p-i)}
	{(\frac{1}{2} (p-i-j))! j!}$\\
	\hti{12}$		= T_{p-i,j}.$

	\sssec{Example.}
	For $p = 5,$\\
	\hth$	       T_0 = - T_{10} = T_{20} =  1.$\\
	\hth$	T_1 = - T_9 = - T_{11} = T_{19} =     +  I.$\\
	\hth$	T_2 = - T_8 = - T_{12} = T_{18} = -1       +2 I^2.$\\
	\hth$	T_3 = - T_7 = - T_{13} = T_{17} =     +2 I       -  I^3.$\\
	\hth$	T_4 = - T_6 = - T_{14} = T_{16}
	=  1       +2 I^2       -2 I^4.$\\
	\hth$		       T_5 = T_{15} =  0.$

	\sssec{Theorem.}\label{sec-tfincheb1}
	\enumb
	\item$	U_{p-1+i,j} = U_{p-1-i,j}.$
	\item$	U_{i+2pk,j} = -U_{i+pk,j} = U_{i,j},$ $j<p.$
	\enume

	Proof:\\
	\hth$	U_{p-1+i,j} = (-1)^{(\frac{1}{2} (p-1+i-j)}
	\frac{(\frac{1}{2} (p-1+i+j))! 2^j}{(\frac{1}{2} (p-1+i-j))! j!}$\\
	\hti{12}$	 = (-1)^{(\frac{1}{2} (p-1+i-j)} (-1)^{(\frac{1}{2}
	 (p-i-j-1)})$\\
	\hth$\frac{(-1)^{(\frac{1}{2}(p-i+j-1))} (\frac{1}{2} (p-i+j-1))! 
	2^j}{(\frac{1}{2}(p-1+i-j))! j!}$\\
	\hti{12}$	 = (-1)^{(\frac{1}{2} (p-1+i-j)} 2^j 
	\frac{(\frac{1}{2} (p-1-i+j))!}{\frac{1}{2} ((p-1-i-j))! j!}$\\
	\hti{12}$			= U_{p-1-i,j}.$

	\sssec{Example.}
	For $p = 5,$\\
	\hth$	U_0 = U_8 = - U_{10} = - U_{18} = 1.$\\
	\hth$	U_1 = U_7 = - U_{11} = - U_{17} =      2 I.$\\
	\hth$	U_2 = U_6 = - U_{12} = - U_{16} = -1       -  I^2.$\\
	\hth$	U_3 = U_5 = - U_{13} = - U_{15} =	I       -2 I^3.$\\
	\hth$	     U_4 = - U_{14} =       =  1       -2 I^2       +  I^4.$\\
	\hth$	     U_9 = - U_{19} =       =  0.$

	\ssec{Symmetry properties for the Polynomials of Legendre.}

	\sssec{Introduction.}

	\sssec{Theorem.}
	\vspace{-18pt}\hti{90}\footnote{24.11.83 and 17.2.89}\\[8pt]
	\hth$P_{p-1-2a-n}^{a)} = P_n^{(a)},$ $n\leq \frac{p-1}{2}-a.$

	Proof:\\
	Let $p' = \frac{1}{2} (p-1).$  The recurrence relations
	\ref{sec-dfinortho}.2 imply\\
	\hth$	(p'+1) P_{p'+1} = - p' P_{p'-1},$\\
	hence $P_{p'+1} = P_{p'-1}.$\\
	They can also be written,\\
	\hth$	(n+1)P_{p-n-2} = -(2n+1)P_{p-n-1} - n P_{p-n}.$\\
	Therefore, starting from $P_{p'}$ and from $P_{p'-1}$ and $P_{p'+1},$
	we obtain by induction $P_{p-1-n} = P_n.$ 

	\sssec{Example.}
	$\begin{array}{l}
	For p = 11,\\
	\hth       P_0 = P_{10} =  1,\\
	\hth	P_1 = P_9 =     I,\\
	\hth	P_2 = P_8 =  5    - 4 I^2,\\
	\hth	P_3 = P_7 =    4I       - 3 I^3,\\
	\hth	P_4 = P_6 = -1    -   I^2	 + 3 I^4,\\
	\hth	P_5      =   -5I       + 5 I^3	  + I^5,
	\end{array}
	\begin{array}{l}
	For p= 13,\\
	\hth	P_0 = P_{12} = 1,\\
	\hth	P_1 = P_{11} =   I,\\
	\hth	P_2 = P_{10} = 6	 - 5 I^2,\\
	\hth	P_3 = P_{ 9} =       5 I	- 4 I^3,\\
	\hth	P_4 = P_{ 8} = 2	 + 6 I^2	+ 6 I^4,\\
	\hth	P_5 = P_{ 7} =     - 3 I	+   I^3       + 3 I^5,\\
	\hth	P_6 =       - 6       - 4 I^2	-   I^4       - I^6,
	\end{array}$

	\sssec{Theorem.}
	\footnote{20.11.87}
	{\em
	$P_{p-1-n} = P_n,$ $n<p.?$}\\
	\hth$P_0^{(a)} = 1,$\\
	\hth$P_1^{(a)} =	I,$\\
	\hth$P_2^{(a)} = \frac{- 1 + (2a+3) I^2}{2(a+1)},$\\
	\hth$P_3^{(a)} = - \frac{-3 I + (2a+5) I^3}{2^2(a+1)(a+2)},$\\
	\hth$P_4^{(a)} = \frac{ 3 - 6(2a+5) I^2+ (2a+5)(2a+7) I^4}{2^2(a+1)(a+2)},$\\
	\hth$P_5^{(a)} = \frac{15 I - 10(2a+7) I^3 + (2a+7)(2a+9) I^5}
		{2^2(a+1)(a+2)},$\\

	\sssec{Example.}
	$\begin{array}{l}
	For p = 11, a = 2\\
	\hth	P_0^{(2)} = 1,\\
	\hth	P_1^{(2)} = I,\\
	\hth	P_2^{(2)} = -2 + 3 I^2,\\
	\hth	P_3^{(2)} = 5 I - I^3,\\
	\hth	P_4^{(2)} = -2 + 3I^2,\\
	\hth	P_5^{(2)} = I,\\
	\hth	P_6^{(2)} = 1,\\
	\end{array}
	\begin{array}{l}
	For p= 13, a = 2\\
	\hth	P_0^{(2)} = 1,\\
	\hth	P_1^{(2)} = I,\\
	\hth	P_2^{(2)} = 2 - I^2,\\
	\hth	P_3^{(2)} = 6 I - 5 I^3,\\
	\hth	P_4^{(2)} = -4 - 6 I^2 - 2 I^4,\\
	\hth	P_5^{(2)} = 6 I	- 5 I^3,\\
	\hth	P_6^{(2)} = 2 - I^2,\\
	\hth	P_7^{(2)} = I,\\
	\hth	P_8^{(2)} = 1,\\
	\end{array}$

	\ssec{Symmetry properties for the Polynomials of Laguerre.}

	\sssec{Theorem (La).}
	\enumb
	\item$      L_{p-1-i,j} = (-1)^j L_{i+j,j},$ $0 \leq  i,j,$ $ i+j < p.$
	\enume

	Proof:
	\begin{enumerate}
	\item$      L_{n,j} = (-1)^j \frac{1}{j!} \ve{n}{j}.$\\
	See for instance, Handbook p.775.
	\end{enumerate}

	\sssec{Example.}
	For $p = 7,$\\
	\hth$	L_0 = 1.$\\
	\hth$	L_1 = 1 -   I.$\\
	\hth$	L_2 = 1 - 2 I - 3 I^2.$\\
	\hth$	L_3 = 1 - 3 I - 2 I^2 +   I^3.$\\
	\hth$	L_4 = 1 + 3 I + 3 I^2 - 3 I^3 - 2 I^4.$\\
	\hth$	L_5 = 1 + 2 I - 2 I^2 + 3 I^3 - 3 I^4 -   I^5.$\\
	\hth$	L_6 = 1 +   I - 3 I^2 +   I^3 - 2 I^4 +   I^5 -   I^6.$

	\sssec{Theorem (La).}
	\enumb
	\item$      L_{p-a-1-i,j}^{(a)} = (-1)^{j+a}\:{ L_{i+j,j}}(a),$
	$0 \leq  i,j,$ $ i+j < p-a.$
	\item$	L_{j,j}^{(a)} = - (a-1)!\:I^{p-a},$ $0 < a,$ $p-a \leq  j < p.$
	\item$	L_{i,j}^{(a)} = 0,$ $a > 0,$ $j < p-a \leq  i < p.$
	\enume

	The proof is left to the reader.

	\sssec{Example.}
	For $p = 13,$ a = 5\\
	\hth$	L_0^{(5)} = 1.$\\
	\hth$	L_1^{(5)} = 6 -   I.$\\
	\hth$	L_2^{(5)} =-5 + 6 I - 6 I^2.$\\
	\hth$	L_3^{(5)} = 4 - 2 I + 4 I^2 + 2 I^3.$\\
	\hth$	L_4^{(5)} =-4 - 6 I + 5 I^2 + 5 I^3 + 6 I^4.$\\
	\hth$	L_5^{(5)} = 5 - 2 I - 5 I^2 -   I^3 - 5 I^4 + 4 I^5.$\\
	\hth$	L_6^{(5)} =-6 + 6 I - 4 I^2 + 5 I^3 + 5 I^4 + 5 I^5 - 5 I^6.$\\
	\hth$	L_7^{(5)} =-1 -   I + 6 I^2 + 2 I^3 - 6 I^4 + 4 I^5 + 5 I^6 
	- 3 I^7.$\\
	\hth$	L_8^{(5)} =	  	  2 I^8.$\\
	\hth$	L_9^{(5)} =	 	  2 I^8 - 6 I^9.$\\
	\hth$	L_{10}^{(5)} =	  2 I^8 +   I^9 - 2 I^{10}.$\\
	\hth$	L_{11}^{(5)} =	  2 I^8 - 5 I^9 - 6 I^{10} -   I^{11}.$\\
	\hth$	L_{12}^{(5)} =	  2 I^8 + 2 I^9 +   I^{10} - 4 I^{11}
	-   I^{12}.$

	\ssec{Symmetry properties for the Polynomials of Hermite.}

	\sssec{Definition.}\label{sec-dfinherm}
	The {\em scaled Hermite polynomials} are defined by
	\enumb
	\item$      H^s_0 = 1,$
	\item$      H^s_1 = I,$
	\item$[\frac{1}{2} n] H^s_n = a_n I H^s_{n-1}
	- \frac{1}{2} (n-1) H^s_{n-2},$\\
	\hth	where $a_n = 1$ if $n$ is even and $a_n = [\frac{1}{2} n],$
	the largest
		integer in $\frac{1}{2} n$ if $n$ is odd.
	\enume

	\sssec{Example.}
	In the fields ${\Bbb Q}$ or ${\Bbb R}$,\\
	\hth$	H^s_2 = -\frac{1}{2} + I^2,$\\
	\hth$	H^s_3 = -\frac{3}{2} I + I^3,$\\
	\hth$	H^s_4 = \frac{3}{8} - \frac{3}{2} I^2 + \frac{1}{2} I^4,$\\
	\hth$	H^s_5 = \frac{15}{8} I - \frac{5}{2} I^3 + \frac{1}{2} I^5,$\\
	\hth$	H^s_6 = -\frac{5}{16} + \frac{15}{8} I^2 - \frac{5}{4} I^4
	+ \frac{1}{6} I^6,$\\
	\hth$	H^s_7 = -\frac{35}{16} I + \frac{35}{8} I^3 - \frac{7}{4} I^5
	+ \frac{1}{6} I^7.$

	\sssec{Theorem.}
	\hth$H^s_{2n}(0) = (-1)^n \frac{(2n-1)!!}{(2n)!!},$ $DH^s_{2n}(0) = 0.$\\
	\hth$H^s_{2n+1}(0) = 0,$ $DH^s_{2n+1}(0) = (-1)^n \frac{(2n+1)!!}{(2n)!!}.
$\\

	\sssec{Lemma.}
	{\em In} ${\Bbb Z}_p,$ $p>2,$
	\enumb
	\item$      (p-1)! = -1.$
	\item$      (p-1-i)! = (-1)^{(i+1)} \frac{1}{i!},$ $0 \leq  i < p.$
	\item$      \ve{p-1-i}{j} = (-1)^j \ve{i+j}{j},$
		$0 \leq  i,j,$ $i+j < p.$
    	\item$      \ve{kp+i}{j} = \ve{i}{j},$ $j < p.$
	\item$(p-2-i)!!\:i!! = (-1)^{\frac{1}{2} k} (p-1-k-i)!! (k+i-1)!!$
				$0 \leq  i < p-1,$ $0 < k+i < p.$
	\enume

	Proof:  0, is the well known Theorem of Wilson.
	1, can be considered as a generalization.\\
	\hth$	(p-1-i)! = (-1)^i  (p-1) \ldots (i+1)$\\
	\hti{12}$		 = (-1)^i \frac{(p-1)!}{i!}$\\
	\hti{12}$		 = (-1)^{(i+1)} \frac{1}{i!}.$\\
	For 2,  $\frac{(p-1-i)!}{(p-1-i-j)! j!}
	= (-1)^{(i+1)} \frac{(i+j)!}{(-1)^{i+j+1} i! j!}
	= (-1)^j \ve{i+j}{j}.$\\

	\sssec{Lemma.}
	{\em Modulo} $p,$ $p>2,$
	\enumb
	\item[0.~~]$( (p-2)!! )^2 = (-1)^{\frac{1}{2} (p-1)}$
	\item[1.~~]$(p-1)!! (p-2)!! = -1.$
	\item[0.0.]$(p-2-i)!! i!! = (-1)^s (p-2)!!,$\\
		{\em where s = $\frac{1}{2}$ i when i is even
		and   s = $\frac{1}{2}$ (p-2-i) when i is odd.}\\
	\item[	1.]{\em or where} $s = [ \frac{1}{2} ([\frac{1}{2} p] + 1 + i)]
	+ [\frac{1}{4} (p+1)].$
	\enume

	0 and 1 are well known and are given for completeness.
	for 2, if $i$ is even,\\
	\hth$	(p-2-i)!! i!! = (p-i)!! (i-2)!! (\frac{i}{p-i}$ or $-1)$\\
	\hti{12}$		= (-1)^{(\frac{1}{2} i)} (p-2)!! 0!!.$\\
	   if $i$ is odd,\\
	\hth$	(p-2-i)!! i!! = (p-4-i)!! (i+2)!! (\frac{p-2-i}{i+2}$ or $-1)$\\
	\hti{12}$		= (-1)^{\frac{1}{2} (p-2-i)} (0)!! p-2!!.$
	2.1, can be verified by choosing $p = 1,3,5,7$ and $i = 0,1,2,3,4.$

	\sssec{Theorem.}
	{\em For scaled Hermite}
	\enumb
	\item$H^s_{i+j,j} = 0,$ $0 \leq  i,j,$ $i$ odd.
	\item$H^s_{p-1-i,j} = H^s_{i+j,j},$ $0 \leq i,j,$ $i${\em even, $j$ even,}
	$i+j < p.$
	\item$H^s_{p-2-i,j} = H^s_{i+j,j},$ $0 \leq i,j,$ $i${\em even, $j$ odd,}
	$i+j < p.$
	\enume

	The proof is left as an exercise.

	\sssec{Example.}
	For $p = 11,$\\
	\hth$	H^s_0 =  1.$\\
	\hth$	H^s_1 =      I.$\\
	\hth$	H^s_2 =  5     +  I^2$\\
	\hth$	H^s_3 =    4 I       +  I^3$\\
	\hth$	H^s_4 = -1     +4 I^2       -5 I^4$\\
	\hth$	H^s_5 =   -5 I       +3 I^3       -5 I^5$\\
	\hth$	H^s_6 = -1     -5 I^2       -4 I^4       +2 I^6$\\
	\hth$	H^s_7 =    4 I       +3 I^3       +  I^5       +2 I^7$\\
	\hth$	H^s_8 =  5      4 I^2       -4 I^4       +4 I^6      -5 I^8$\\
	\hth$	H^s_9 =      I       +  I^3       -5 I^5       +2 I^7     
	-5 I^9.$\\
	\hth$       H^s_{10} =  1      + I^2       -5 I^4       +2 I^6    
	  -5 I^8      -I^{10}.$

	\sssec{Problem.}
	The Jacobi polynomials can be defined by\\
	\hth$Ps_0^{(a,b)} := 1, P_1^{(a,b)} := \frac{a-b + (a+b+2)I}{2(a+1)},$\\
	\hth$	2(n+1)(n+a+b+1)(2n+a+b)\: P_{n+1}^{(a,b)} :=\\
	\hti{12}	((2n+a+b+1)(a^2-b^2)\\
	\hti{16}	 + (2n+a+b)(2n+a+b+1)(2n+a+b+2) I)\: P_n^{(a,b)}\\
	\hti{12}		- 2(n+a)(n+b)(2n+a+b+2)\: P_{n-1}^{(a,b)}.$
	Determine an appropriate scaling for the Jacobi polynomials that
	gives symmetry properties which generalize those of the special case
	where $a = b$.




\setcounter{section}{2}
	\section{Addition Formulas for Functions on a Finite Fields.}

	\setcounter{subsection}{-1}
	\ssec{Introduction.}
	Ungar, gave recently the addition formulas associated
	with a generalization of the trigonometric and hyperbolic functions by
	Ricatti.  This suggested the extension to the finite case.  Section 1,
	is the Theorem of Ungar, the special case for 3 functions is given in
	\ref{sec-Sric3}, with the associated invariant \ref{sec-tric3}.2.
	The invariant defines the distances, addition, which in fact
	corresponds to the multiplication of associated Toeplitz matrices
	gives the angles.  For 3 dimensions we have 2 special cases,
	$p \equiv 1 \pmod{3}$ and $p \equiv -1 \pmod{3}.$ In the latter case all
	non isotropic direction form a cycle.  In the former case, we can
	consider that the set of $(p-1)^2$ non isotropic directions
	corresponds to a direct product of 2 cyclic groups of order $p-1,$
	I conjecture (14) that there are always pairs of generators which are
	closely related called special generators (13).
	This is extended to more than 3 functions in \ref{sec-Sric4}.
	The connection with difference sets is given at the end of that
	chapter.

	\ssec{The Theorem of Ungar.}

	\sssec{Theorem. [Ungar]}
	{\em If $f$ is a solution of}
	\enumb
	\item$	D^{n+1} f + a_n D^n f +  \ldots  + a_0 f = 0,$ $a_{n+1} = 1,$
	\item$	D^n f(0) = 1,$ $D^k f(0) = 0,$ $ 0 \leq  k < n,$

	{\em then}
	\item$f(x+y) = \sum_{m=0}^n a_{m+1} 
	\sum_{k=0}^m  D^k f (x) D^{m-k} f (y).$

	{\em More generally, if}
	\item$	D^k f(0) = d_k,$ $0 \leq  k \leq  n,$

	{\em then}
	\item$	 \sum_{m=0}^n  a_{m+1} \sum_{k=0}^m  D^k f (x) D^{m-k} f (y)\\
	\hth= \sum_{m=0}^n  a_{m+1} \sum_{k=0}^m  d_k D^{m-k} f (x+y).$
	\enume

	Proof: See Abraham Ungar, 1987.

	\sssec{Example.}
	\enumb
	\item[0.0.]$	a_m = 0,$ $0 \leq  m \leq  n,$ $f = \frac{I^n}{n!},$
	\item[	1.]$	(x+y)^n = \sum_{k=0}^n  \ve{n}{k} x^k y ^{n-k}.$
	\item[1.0.]$	n = 0,$ $a_1 = 1,$ $a_{-1} = 0,$ $f = e^I,$
	\item[  1.]$	e^{x+y} = e^x e^y.$
	\item[2.0.]$	n = 1,$ $a_2 = 1,$ $a_1 = 0,$ $a_0 = 1,$ $f = sin,$
	\item[  1.]$	sin(x+y) = sin(x) cos(y) + cos(x) sin(y).$
	\item[3.0.]$	n = 1,$ $a_2 = 1,$ $a_1 = 0,$ $a_0 = -1,$ $f = sinh,$
	\item[  1.]$	sinh(x+y) = sinh(x) cosh(y) + cosh(x) sinh(y).$
	\item[4.~~]$	a_{n+1} = 1,$ $a_k = 0,$ $0 < k \leq  n,$ $a_0 = -j,$
	where $j = \pm 1,$\\
		$D^if = R^{(jn,0)} = \sum_{k=0}^{\infty} 
	\frac{I^{rk-i}}{(rk-i)!},$ where $r = n+1,$\\
		$(j=-1?)$\\
		1. $	R^{(-n,0)}(x+y) = R^{(-n,0)}(x) R^{(-n,n)}(y)$\\
	\hth$			+ R^{(-n,1)}(x) R^{(-n,n-1)}(y) +  \ldots .$\\
	\hth$			+ R^{(-n,0)}(x) R^{(-n,n)}(y).$\\
		$(j=-1?)$\\
		These are, with my notation, the functions of Vincenzo Ricatti.
	In particular, when $n = 2,$? we have the following Theorem.
	\enume

	\sssec{Theorem.}
	{\em If $n$ is odd, then}
	\enumb
	\item $R^{(-n,0)} = R^{(n.0)}(-I).$
	\item$R^{(-n,j)} = (-1)^j R^{(n.j)}(-I).$
	\enume

	\ssec{The case of 3 functions.}\label{sec-Sric3}
	\sssec{Theorem.}\label{sec-tric3}
	{\em Let}
	\enumb
	\item$	f = R^{(2,0)},$ $g = R^{(2,1)} h = R^{(2,2)},$

	{\em then}
	\item$	f(x+y) = f(x) h(y) + g(x) g(y) + h(x) f(y),$\\
	$	g(x+y) = g(x) h(y) + h(x) g(y) + f(x) f(y),$\\
	$	h(x+y) = h(x) h(y) + f(x) g(y) + g(x) f(y),$
	\item$	f^3 + g^3 + h^3 - 3 f g h = 1.$
	\item$	  (f(x) h(y) + g(x) g(y) + h(x) f(y))^3$\\
	\hth$	+ (g(x) h(y) + h(x) g(y) + f(x) f(y))^3$\\
	\hth$	+ (h(x) h(y) + f(x) g(y) + g(x) f(y))^3$\\
	\hth$	- 3 (f(x) h(y) + g(x) g(y) + h(x) f(y))$\\
	\hti{11}$	    (g(x) h(y) + h(x) g(y) + f(x) f(y))$\\
	\hti{11}$	    (h(x) h(y) + f(x) g(y) + g(x) f(y))$\\
	\hti{4}$	= (f(x)^3 + g(x)^3 + h(x)^3 - 3 f(x) g(x) h(x))$\\
	\hth$	  (f(y)^3 + g(y)^3 + h(y)^3 - 3 f(y) g(y) h(y)).$
	\enume

	Proof:\\
	  $g = Df,$ $h = Dg = D^2f,$ $f = Dh = D^2g = D^3f,$\\
	$f,$ $g$ and $h$ satisfy the same differential equation, whose Wronskian
	is constant this gives\\
	\hth$\det\left|\begin{array}{ccc}f&g&h\\g&h&f\\h&f&g
	\end{array}\right| = -1.$

	\sssec{Theorem.}
	{\em The solution of \ref{sec-tric3}, $ f,$ $g,$ $h$ is given by}
	\enumb
	\item	$f = A e^I + B     e^{\beta I} + C     e^{\beta^{-1} I},$
	\item	$g = A e^I + B \beta    e^{\beta I} 
			+ C \beta ^{-1} e^{\beta   I},$
	\item	$h = A e^I + B \beta ^{-1} e^{\beta I}
			+ C \beta^{-1}e^{\beta   I},$\\
	{\em where}
	\item$	\beta ^2 + \beta  + 1 = 0,$
	\item$	A = \frac{1}{3},$ $B = \frac{1}{3} \beta ,$
		$C = \frac{1}{3} \beta^{-1}.$
	\enume

	\sssec{Corollary.}
	\hth$   f = e^{-\frac{1}{2} I} cos(\frac{\sqrt{3}}{2}I),$\\
	\hth$ g = e^{-\frac{1}{2} I} cos((\frac{\sqrt{3}}{2}
		+ \frac{\pi}{3})I),$\\
	\hth$ h = e^{-\frac{1}{2} I} cos((\frac{\sqrt{3}}{2}
	- \frac{\pi}{3})I),$\\
	{\em is a solution of} $D^3 f = f.$

	I examined the more general case\footnote{7.12.87} starting from
	$f_1,g_1,h_1$	and using the addition formulas, it appears that
	the period is always $p-1$ and that if\\
	\hth$f^3 + g^3 + h^3 - 3 f g h = 1$ then we have\\
	$f(\frac{\pi}{3}) = g(\frac{2\pi}{3}) = 0$ when $p \equiv  1 \pmod{6}.$

	Application to the case of 3 dimensional Affine geometry associated
	to $p$\footnote{8.12.87}.

	\sssec{Lemma.}\label{sec-lric3}
	{\em Let}
	\enumb
	\item[0.~~]$T(x,y,z) := x^3+y^3+z^3-3xyz,$ $L(x,y,z) : = x+y+z,$\\
	\hth$	S(x,y,z) := x^2+y^2+z^2-yz-zx-xy,$

	{\em then}\\
	\hth$	T(x,y,z) = L(x,y,z) S(x,y,z).$
	\item[0.0.]{\em If $-3$N$p$ or $p \equiv  2 \pmod{3},$ then the
		only points of $S = 0$ are}\\
	\hth$	(a,a,a).$\\
	\item[	1.]{\em If $-3$R$p$ or $p \equiv  1 \pmod{3},$ then}\\
	\hth$		S(x,y,z) = (x+\tau y+\tau 'z)(x+\tau 'y+\tau z),$ with\\
	\hth$		\tau  := \frac{1}{2} (-1+\sqrt{-3}\tau') 
	:= \frac{1}{2} (-1-\sqrt{-3}).$
	\item[2.~~]{\em The number of lines through the origin on
		$T(x,y,z) = 0$ is\\
	\hth	$3p$ if $p \equiv 1 \pmod{3}$\\
	and\\
	\hth $p+2$
	if $p \equiv -1 \pmod{3}.$}
	\enume

	Proof:  The number of lines through the origin is the same as the
	number of points in a plane not through the origin which are on the
	ideal line or on $S.$\\
	If $p \equiv 1 \pmod{3},$ this gives
	$(p+1) + 1,$\\
	if $p \equiv -1 \pmod{3},$ this gives $3(p+1) - 3.$

	\sssec{Definition.}
	Let ${\cal T}$  be the set of points $(x,y,z) \ni$
	\enumb
	\item$	x^3 + y^3 + z^3 - 3 x y z = 1.$

	Let the addition in {\cal T}  be defined by
	\item$	(x,y,z) + (x',y',z')
		:= (y y' + x z' + z x', x x' + y z' + z y', z z' + x y' + y x')$
	\enume

	\sssec{Theorem.}
	\enumb
	\item$({\cal T} ,+) ${\em is an Abelian group with neutral element}
	$(0,0,1).$
	\item$	(x,y,z) + (x',y',z') + (x'',y'',z'')$\\
	\hti{6}$= (x(x'y''+y'x''+z'z'') + y(x'x''+y'z''+z'y'')
		+ z(x'z''+y'y''+z'x''),$\\
	\hth$x(x'z''+y'y''+z'x'') + y(x'y''+y'x''+z'z'')
		+ z(x'x''+y'z''+z'y''),$\\
	\hth$x(x'x''+y'z''+z'y'') + y(x'z''+y'y''+z'x'')
		+ z(x'y''+y'x''+z'z'')).$
	\item$	(x,y,z) + (y^2-zx,x^2-yz,z^2-xy) = (0,0,1).$
	\enume

	\sssec{Corollary.}
	\enumb
	\item$	2(x,y,z) = (y^2+2zx,x^2+2yz,z^2+2xy).$
	\item$	3(x,y,z) = (3(x^2y+y^2z+z^2x),3(x^2z+y^2x+z^2y,1+9xyz)).$
	\enume

	\sssec{Theorem.}
	{\em If}
	\enumb
	\item$	(x_n,y_n,z_n) := n(x,y,z)$
	{\em then}
	\item$	x_n+y_n+z_n = (x+y+z)^n.$
	\item$	x_{k(p-1)+i} + y_{k(p-1)+i} + z_{k(p-1)+i} = x_i+y_i+z_i.$
	\enume

	\sssec{Theorem.}
	{\em Let}
	\enumb
	\item$u^2 = -3 x^2 + 2 s x - \frac{s^2}{3} + \frac{4}{3s}$
	\item$y = \frac{1}{2} (s - x \pm  u),$
	\item$z = s - x - y,$
	{\em then}$(x,y,z) \in {\cal T}$.
	\enume

	Proof:\\
	Substitute $z$ by $s - x - y \in x^3 + y^3 + z^3 - 3 xyz -1 = 0$ gives\\
	$3 s \:y^2 - 3 s(s-x) y + (3s\:x^2 - 3 s^2\:x + s^3-1) = 0,$\\
	dividing by $3s,$ the discriminant is the second member of 0. $\Box$

	\sssec{Definition.}
	\enumb
	\item The {\em distance} $d$ between 2 points $(x,y,z)$ and
		$(x',y',z')$ is given by\\
	\hth$d^3(x,x') := (x'-x)^3 + (y'-y)^3 + (z'-z)^3 - 3(x'-x)(y'-y)(z'-z).$
	\item If the distance between 2 distinct points is 0, the line incident
	    to the 2 points is called {\em isotropic}.
	\enume

	\sssec{Theorem.}
	{\em The isotropic lines are those on the surface T(x,y,z) = 0.}

	\sssec{Lemma.}
	{\em
	\hth$	d^3(x,x') = d^3(0,x) - d^3(0,x')
		-3(x(x'^2-y'z')+ y(y'^2-z'x')+z(z'^2-x'y') )$\\
	\hti{16}$+ 3 (x'(x^2-yz)+y'(y^2-zx)+z'(z^2-xy) ).$}

	\sssec{Theorem.}
	\enumb
	\item$  d(P,Q) = - d(Q,P).$
	\item{\em If $P = (0,0,0,1),$ then $P \times Q$ is isotropic iff $Q$ is
		on the line $l$ joining
	$P$ to (1,1,1,1) or on a line through $P$ perpendicular to $l.$}
	\enume

	The ideal points on the surface satisfy\footnote{11.12.87}\\
	\hth$	x^3 + y^3 + z^3 - 3 x y z = 0.$

	\sssec{Definition.}
	The {\em normal to the surface} ${\cal T}$  at (a,b,c) is\\
	\hth$	[a^2-bc,b^2-ca,c^2-ab].$

	\sssec{Notation.}
	If $p \equiv 1 \pmod{6},$ $\delta := \frac{p-1}{3}.$ 

	\sssec{Theorem.}
	\hth If $p \equiv  1 \pmod{6},$ $({\cal T} ,+) 
		\sim C_{p-1} \bigx  C_{p-1}.$

	Proof:  The order of the group follows from Lemma \ref{sec-lric3}?.

	\sssec{Lemma.}
	{\em If an Abelian group is isomorphic to  $C_q \bigx  C_q$ 
	if $u$ and $v$ are of order $p$ and $u^i \neq v^j$ for all $i$ and $j$
	between $1$ and $q,$ $u$ and $v$ are generators of the group.}

	\sssec{Lemma.}\label{sec-lric31}
	{\em If $p \equiv 1 \pmod{6},$ $g$ is a primitive root of $p$ and
	$(a,b,c)$ and $(a',b',c')$ are obtained using}
	\enumb
	\item$b,c = \frac{g-a \pm\sqrt{\frac{-g^2+6ag-9a^2+4g^{-1}}{3}}}{2}$\\
	{\em if their $i$-th iterates are distinct, $0<i<p-1,$ then
	(a,b,c) and (a',b',c') are generators.\\
	In particular, if $h^3 = 1$ then}
	\item$	b,c = \frac{h-a \pm  \sqrt{(h+3a)(h-a)}}{2}$
	\enume

	Proof:  the fact that g is primitive insures that the sum of the
	components of the $i$-th iterate of $(a,b,c)$ is $g^i,$ because these
	are distinct for $i = 1$ to $p-2,$ the Lemma follows.

	\sssec{Definition.}
	If the pair $(a,b,c)$ and $(b,a,c)$ are pairs of generators of
	(T,+) then $(a,b,c)$ is called a {\em special generator} of $T.$

	\sssec{Conjecture.}
	Given a primitive root of $p \equiv 1 \pmod{6},$ there exists always
	special generators $(a,b,c) \ni  a+b+c = g \pmod{p}$
	\footnote{21.12.87}.

	\sssec{Theorem.}
	{\em If $(a_1,b_1,c_1)$ is a special generator, the period is\\
	0    1    2  \ldots $\delta$ $\delta$+1 $\delta$+2 \ldots$2\delta$
		 $2\delta$ +1 $2\delta$ +2 \ldots\\ 
	0   $a_1$   $a_2$  \ldots   1   $c_1$   $c_2$  \ldots     
		0   $b_1$   $b_2$  \ldots\\
	0   $b_1$   $b_2$  \ldots   0   $a_1$   $a_2$  \ldots
		     1   $c_1$   $c_2$  \ldots \\
	1   $c_1$   $c_2$  \ldots   0   $b_1$   $b_2$  \ldots     
		0   $a_1$   $a_2$  \ldots.\\
	For $(b_1,a_1,c_1)$ the period, scaled again is\\
	0    1    2  \ldots     $\delta$   $\delta$ +1  $\delta$ +2  \ldots
		$2\delta$  $2\delta$ +1 $2\delta$ +2 \ldots \\
	0   $b_1$   $b_2$  \ldots   0   $a_1$   $a_2$  \ldots
		1   $c_1$   $c_2$  \ldots \\
	0   $a_1$   $a_2$  \ldots   1   $c_1$   $c_2$  \ldots
		0   $b_1$   $b_2$  \ldots \\
	1   $c_1$   $c_2$  \ldots   0   $b_1$   $b_2$  \ldots
		0   $a_1$   $a_2$  $\ldots$  .}

	\sssec{Algorithm.}
	For a given $p,$ we determine the smallest positive
	primitive root $g,$ then for increasing values of $a,$ we determine $b$ and $c$
	using \ref{sec-lric31}.0, if $c(delta) = 1$ we permute $a,b,c,$ in the order\\
	\hth$	b,a,c,$~~$c,b,a,$~~$a,c,b,$~~$c,a,b,$~~$b,c,a,$\\
	unless $p \equiv 1 \pmod{9},$ in which case we try a new $p,$ if c(delta) = 0,
	we save the period and permute, if $c(delta) = 0$ and the values of
	$a(i),$ $b(i),$ $c(i)$ are not distinct from some $i,$ from the corresponding
	saved values, we permute again, if we exaust the permutations, we
	ignore this value of $a.$  When we have obtained $(a,b,c)$ such that
	the first a(delta) = -1 and and the second = 1 we exchange.

	\sssec{Example.}
	\enumb
	\item$  p = 7,$ $g^i = 1,3,2,6,4,5,$	${\cal T}  =$\\
	$\begin{array}{lllllllll}
	0	&1	&2	&3	&4	&5	&6	&7	&8\\
	(0,0,1)&(0,0,2)&(0,0,4)&(0,1,0)&(0,2,0)&(0,4,0)&(1,0,0)&(1,3,6)
		&(1,5,5)\\
	(0,0,1)&(0,0,1)&(0,0,1)&(0,1,0)&(0,1,0)&(0,1,0)&(1,0,0)&(5,1,2)
		&(2,3,3)\\
	\cline{1-9}
	9	&10	&11	&12	&13	&14	&15	&16	&17\\
	(1,6,3)&(2,0,0)&(2,3,3)&(2,5,6)&(2,6,5)&(3,1,6)&(3,2,3)&(3,3,2)
		&(3,4,5)\\
	(5,2,1)&(1,0,0)&(2,3,3)&(5,2,1)&(5,1,2)&(1,5,2)&(3,2,3)&(3,3,2)
		&(2,5,1)\\
	\cline{1-9}
	18	&19	&20	&21	&22	&23	&24	&25	&26\\
	(3,5,4)&(3,6,1)&(4,0,0)&(4,3,5)&(4,5,3)&(4,6,6)&(5,1,5)&(5,2,6)
		&(5,3,4)\\
	(2,1,5)&(1,2,5)&(1,0,0)&(5,2,1)&(5,1,2)	&(2,3,3)&(3,2,3)&(2,5,1)
		&(1,2,5)\\
	\cline{1-9}
	27	&28	&29	&30	&31	&32	&33	&34	&35\\
	(5,4,3)&(5,5,1)&(5,6,2)&(6,1,3)&(6,2,5)&(6,3,1)&(6,4,6)&(6,5,2)
		&(6,6,4)\\
	(1,5,2)	&(3,3,2)&(2,1,5)&(2,5,1)&(1,5,2)&(2,1,5)&(3,2,3)&(1,2,5)
		&(3,3,2)
	\end{array}\\
	\begin{array}{rcrrrrrrr}
	+	&\vline& 	& 0	& 9	& 4	&29   	&20	&27\\
	\cline{1-9}
	 0	&\vline& 	& 0	& 9	& 4	&29	&20	&27\\
	30	&\vline& 	&30	&35	&13	&24	&26	&11\\
	10	&\vline& 	&10	&31	& 2	&21	& 3	&32\\
	34	&\vline& 	&34	& 8	&17	&16	& 7	&33\\
	 5	&\vline& 	& 5	&18	& 6	&14	& 1	&12\\
	22	&\vline& 	&22	&15	&19	&23	&25	&28
	\end{array}$

	When scaled the group is isomorphic to $C_6 \bigx C_2,$ we have the
	equivalences,\\
	0,1,2; 3,4,5; 6,10,20; 7,13,22; 8,11,23; 9,12,21; 14,27,3; 15,24,33;
	16,28,35; 17,25,30; 18,29,32; 19,26,34.\\
	We have the table\\
	$\begin{array}{rcrrrrrr}
	 + &\vline& 0& 9& 3&18& 6&14\\
	\cline{1-8}
	 0 &\vline& 0& 9& 3&18& 6&14\\
	16 &\vline&16& 7&15&19& 8&17
	\end{array}$

	\item$  p = 13,$ $g^i = 1,2,4,8,3,6,12,11,9,5,10,7,$ the scaled period
 is\\
	$\begin{array}{cccc}
	( 0,0,1)&( 7,1,6)&( 7,9,11)&(11,10,6),\\
	( 1,0,0)&( 6,7,1)&(11,7,9) &( 6,11,10),\\
	( 0,1,0)&( 1,6,7)&( 9,11,7)&(10,6,11),
	\end{array}$

	\item$ p = 19,$ $g^i = 1,2,4,8,16,13,7,14,9,18,17,15,11,3,6,12,5,10,$
		the scaled period is\\
	$\begin{array}{cccccc}
	( 0, 0, 1),&( 3, 9, 8),&(15, 1, 4),&( 8,13,18),&( 7, 5, 8),&(11, 0, 9)\\
	( 1, 0, 0),&( 8, 3, 9),&( 4,15, 1),&(18, 8,13),&( 8, 7, 5),&( 9,11, 0)\\
	( 0, 1, 0),&( 9, 8, 3),&( 1, 4,15),&(13,18, 8),&( 5, 8, 7),&( 0, 9,11)
	\end{array}$
	\enume

	\sssec{Example.}
	The following are special generators, for the given primitive
	root, which is the smallest positive one:\\
	$\begin{array}{ccrclcccrclcccrcl}
	  p  &\vline&   g &\vline&  sp. gen.	&\vline\vline&
	  p  &\vline&   g &\vline&  sp. gen.    	&\vline\vline&
	  p  &\vline&   g &\vline&  sp. gen.\\
	\cline{1-17}
	  7 &\vline&   3 &\vline&  (6,1,3)	&\vline\vline&
	283 &\vline&   3 &\vline&  (3,158,125)&\vline\vline&
	631 &\vline&   3 &\vline&  (324,4,306)\\
	 13 &\vline&   2 &\vline&  (1,2,12)	&\vline\vline&
	307 &\vline&   5 &\vline&  (4,192,116)&\vline\vline&
	643 &\vline&  11 &\vline&  (152,0,502)\\
	 19 &\vline&   2 &\vline&  (6,18,16)	&\vline\vline&
	313 &\vline&  10 &\vline&  (5,21,297)	&\vline\vline&
	661 &\vline&   2 &\vline&  (2,134,527)\\
	 31 &\vline&   3 &\vline&  (30,4,0)	&\vline\vline&
	331 &\vline&   3 &\vline&  (0,237,97)	&\vline\vline&
	673 &\vline&   5 &\vline&  (3,52,623)\\
	 37 &\vline&   2 &\vline&  (18,7,14)	&\vline\vline&
	337 &\vline&  10 &\vline&  (180,0,167)	&\vline\vline&
	691 &\vline&   3 &\vline&  (425,0,269)\\
	 43 &\vline&   3 &\vline&  (8,35,3)	&\vline\vline&
	349 &\vline&   2 &\vline&  (50,299,2)	&\vline\vline&
	709 &\vline&   2 &\vline&  (424,3,284)\\
	 61 &\vline&   2 &\vline&  (3,2,58)	&\vline\vline&
	367 &\vline&   6 &\vline&  (22,346,5)	&\vline\vline&
	727 &\vline&   5 &\vline&  (377,352,3)\\
	 67 &\vline&   2 &\vline&  (2,12,55)	&\vline\vline&
	373 &\vline&   2 &\vline&  (53,6,316)	&\vline\vline&
	733 &\vline&   6 &\vline&  (1,541,197)\\
	 73 &\vline&   5 &\vline&  (4,12,62)	&\vline\vline&
	379 &\vline&   2 &\vline&  (5,200,176)	&\vline\vline&
	739 &\vline&   3 &\vline&  (0,400,342)\\
	 79 &\vline&   3 &\vline&  (5,42,35)	&\vline\vline&
	397 &\vline&   5 &\vline&  (8,22,372)	&\vline\vline&
	751 &\vline&   3 &\vline&  (4,426,324)\\
	 97 &\vline&   5 &\vline&  (24,3,75)	&\vline\vline&
	409 &\vline&  21 &\vline&  (390,38,2)	&\vline\vline&
	757 &\vline&   3 &\vline&  (5,122,632)\\
	103 &\vline&   5 &\vline&  (79,25,4)	&\vline\vline&
	421 &\vline&   2 &\vline&  (6,5,412)	&\vline\vline&
	769 &\vline&  11 &\vline&  (1,404,375)\\
	109 &\vline&   6 &\vline&  (13,0,102)	&\vline\vline&
	433 &\vline&   5 &\vline&  (3,273,162)	&\vline\vline&
	787 &\vline&   2 &\vline&  (411,3,375)\\
	127 &\vline&   3 &\vline&  (77,0,53)	&\vline\vline&
	439 &\vline&  15 &\vline&  (0,264,190)	&\vline\vline&
	811 &\vline&   3 &\vline&  (0,188,626)\\
	139 &\vline&   2 &\vline&  (107,31,3)	&\vline\vline&
	457 &\vline&  13 &\vline&  (14,456,0)	&\vline\vline&
	823 &\vline&   3 &\vline&  (15,0,811)\\
	151 &\vline&   6 &\vline&  (106,51,0)	&\vline\vline&
	463 &\vline&   3 &\vline&  (0,335,331)	&\vline\vline&
	829 &\vline&   2 &\vline&  (7,572,252)\\
	157 &\vline&   5 &\vline&  (4,39,119)	&\vline\vline&
	487 &\vline&   3 &\vline&  (39,0,551)	&\vline\vline&
	853 &\vline&   2 &\vline&  (155,698,2)\\
	163 &\vline&   2 &\vline&  (2,29,134)	&\vline\vline&
	499 &\vline&   7 &\vline&  (1,87,418)	&\vline\vline&
	859 &\vline&   2 &\vline&  (228,625,8)\\
	181 &\vline&   2 &\vline&  (2,36,145)	&\vline\vline&
	523 &\vline&   2 &\vline&  (310,6,209)	&\vline\vline&
	877 &\vline&   2 &\vline&  (10,5,864)\\
	193 &\vline&   5 &\vline&  (122,3,73)	&\vline\vline&
	541 &\vline&   2 &\vline&  (93,3,447)	&\vline\vline&
	883 &\vline&   2 &\vline&  (147,6,732)\\
	199 &\vline&   5 &\vline&  (30,5,167)	&\vline\vline&
	547 &\vline&   2 &\vline&  (335,2,212)	&\vline\vline&
	907 &\vline&   2 &\vline&  (553,2,354)\\
	211 &\vline&   2 &\vline&  (33,2,178)	&\vline\vline&
	571 &\vline&   3 &\vline&  (7,8,559)	&\vline\vline&
	919 &\vline&   7 &\vline&  (129,0,727)\\
	223 &\vline&   3 &\vline&  (138,0,88)	&\vline\vline&
	577 &\vline&   5 &\vline&  (4,300,278)	&\vline\vline&
	937 &\vline&   5 &\vline&  (7,493,442)\\
	229 &\vline&   6 &\vline&  (1,168,66)	&\vline\vline&
	601 &\vline&   7 &\vline&  (138,463,7)	&\vline\vline&
	967 &\vline&   5 &\vline&  (3,661,308)\\
	241 &\vline&   7 &\vline&  (20,4,224)	&\vline\vline&
	607 &\vline&   3 &\vline&  (0,441,169)	&\vline\vline&
	991 &\vline&   6 &\vline&  (8,228,761)\\
	271 &\vline&   6 &\vline&  (159,7,111)	&\vline\vline&
	613 &\vline&   2 &\vline&  (3,50,562)	&\vline\vline&
	997 &\vline&   7 &\vline&  (0,625,379)\\
	277 &\vline&   5 &\vline&  (11,5,266)	&\vline\vline&
	619 &\vline&   2 &\vline&  (5,65,551) 	&\vline\vline&	
	\end{array}$

	\sssec{Lemma.}
	{\em If $p \equiv  -1 \pmod{6},$}
	\enumb
	\item$s_i := P_{i0}+P_{i1}+P_{12}  \Rightarrow   s_i = s_1^i.$
	\item$f_i := P_{i0}^2+P_{i1}^2+P_{12}^2-(P_{i1}P_{12}+P_{12}P_{i0}
		+P_{i0}P_{i1})  \Rightarrow   f_i = f_1^i.$
	\enume

	\sssec{Lemma.}
	{\em If $p \equiv -1 \pmod{6},$ if $g$ is a primitive root of $p$,}

	Proof:  Let  $\ldots$\\ 
	To determine what happens for the solutions for $i = 1, 2, \ldots 
	p-1$:\\
	let $g$ be a generator for $p,$ we want
	To determine what happens for the solutions for $i \equiv 0 \pmod{p-1}$
:\\
	\hth$	a + b + c = 1,$ $a^2+b^2+c^2-(bc+ca+ab) = 1,$ $ \Rightarrow $\\
	\hth$	(a+b+c)^2 = 1,$ $bc+ca+ab = 0,$\\
	given $a,$ $b+c = 1-a,$ $bc = a(a-1),$\\
	\hth$	b,c = \frac{1-a \pm  \sqrt{(a-1)(1-3a)}}{2}$\\
	special solutions $(1,0,0),$ $(-\frac{1}{3},\frac{2}{3},\frac{2}{3}),$
	$2(-\frac{1}{3},\frac{2}{3},\frac{2}{3}) = (0,1,0),$
	$2(\frac{2}{3},-\frac{1}{3},\frac{2}{3}) = (1,0,0),$
	there should be $\frac{p+1 - 3 - 3}{2} = \frac{p-5}{2}$ possible values
	of $a$.

	\sssec{Notation.}
	$\epsilon := (0,0,1),$ $\alpha := (1,0,0),$ $\beta := (0,1,0).$

	\sssec{Theorem.}\label{sec-tric3cyclic}
	{\em If $p \equiv -1 \pmod{6},$ and $3\delta$  = $p^2-1,$ then}
	\enumb
	\item$	({\cal T} ,+) \sim C_{3\delta }.$
	\item{\em If $h$ is a generator of this group then $h^{\delta }
		= (0,1,0)$ or
	$(1,0,0),$ in the former case, we will choose $g = h^{-1}$
	otherwize\\
		we will choose $g = h.$}
	\item{\em If $g^i = (a,b,c),$ then $g^{i+\delta }  = (c,a,b),$
	$ g^{i+2\delta }  = (b,c,a),$
		$g^{pi} = (b,a,c),$ $g^{pi+\delta } = (c,b,a),$
	$g^{pi+2\delta } = (a,c,b)$\footnote{13.12.87}}.\\
	\item$d(P_i,P_{p+(p+1)l+(p-1)kj+i}) = d(P_i,P_{1+(p+1)l-(p-1)kj+i}).$ 
	\item$	d(P_i,P_{(p+1)l+(p-1)kj+i}) = d(P_i,P_{(p+1)l-(p-1)kj+i}).$
	\enume

	Proof:\\
	If $p \equiv -1 \pmod{6},$ to any of the $p^2-1$ line through the origin
	which does not pass through (1,1,1,1) and is not in the plane
	perpendicular to this last line, associate a point $a,b,c,1$ say
	let $u := a^3 + b^3 + c^3 - 3abc,$ $u \neq 0$ and $u$ has a unique
	cube root v in ${\Bbb Z}_p,$ therefore the point
	$(\frac{a}{v},\frac{b}{v},\frac{c}{v}) \in  {\cal T}.$\\
	\hth$(a,b,c) + (1,0,0) = (c,a,b),$ $(a,b,c) + (0,1,0) = (b,c,a).$

	\sssec{Lemma.}

	Proof: $\phi$ $(p^2-1)$ = 2 $\phi$ (p-1) $\phi$ (p+1),

	\sssec{Example.}
	$p = 11,$ $g = 7,$

	\sssec{Lemma.}
	\enumb
	\item$  (a+b+c)(a^2+b^2+c^2-bc-ca-ab) = a^3+b^3+c^3-3abc.$
	\item{\em Given $a$ and $g,$ a primitive root of $p,$ then}\\
	$b,c = \frac{g-a \pm \sqrt{\frac{-g^2+6ag-9a^2+4g^{-1}}{3}}}{2}$
	\enume

	Proof:  $a+b+c = g$ and $a^2+b^2+c^2-bc-ca-ab = g^{-1}  \Rightarrow
	(a+b+c)^2 = g^2$ and $bc+ca+ab = 0,$ therefore $b+c = g-a$ and
	$bc = \frac{g^2-g^{-1}}{3} -a(g-a),$ hence
	$b$ and $c$ are roots of a quadratic equations, this gives 1.

	\sssec{Theorem.}\label{sec-tric3gen}
	{\em If $p \equiv  -1 \pmod{6},$}
	\enumb
	\item[0.~~]{\em A necessary condition for $(a,b,c)$ in ${\cal T}$  to be a
	generator
	    is that $a+b+c$ be a primitive root for $p.$}
	\item[1.~~]{\em If $(a,b,c)$ is a generator such that $(a,b,c)^{\delta }
	= (1,0,0),$}
	\item[2.0.]$p \equiv  2 \pmod{9},$ {\em or}
	$p \equiv  11 \pmod{18} \Rightarrow$\\
		$(b,c,a)^{\delta } = \epsilon,$ $(c,a,b)^{\delta } = \beta,$
		$(c,b,a)^{\delta } = \epsilon,$ $(a,c,b)^{\delta } = \alpha,$
		$(b,a,c)^{\delta } = \beta,$\\
		\item$p \equiv  5 \pmod{9},$
	or $p \equiv  5 \pmod{18} \Rightarrow$\\
		$(b,c,a)^{\delta } = \beta$, $(c,a,b)^{\delta } = \epsilon$,
		$(c,b,a)^{\delta } = \alpha$, $(a,c,b)^{\delta } = \epsilon,$
	$(b,a,c)^{\delta } = \beta$,\\
		\item$p \equiv  8 \pmod{9},$
	 {\em or} $p \equiv  17 \pmod{18} \Rightarrow$\\
		$(b,c,a)^{\delta } = \alpha$, $(c,a,b)^{\delta } = \alpha$,
		$(c,b,a)^{\delta } = \beta$, $(a,c,b)^{\delta } = \beta$,
	$(b,a,c)^{\delta } = \beta,$\\
	\item[  1.]$	p(u,v,w) = (v,u,w).$
	\enume

	\sssec{Definition.}
	Given a generator $(a,b,c)$ and a non isotropic scaled
	direction $(u,v,w)$ the corresponding {\em angular direction} is the
	multiplier $i$ such that i(a,b,c) = (u,v,w).

	\sssec{Conjecture.}
	\vspace{-18pt}\hspace{90pt}\footnote{28.12.87}\\[8pt]
	\enumb
	\item$angular\: direction(P_{i+k},P_i) = i + angular \:
		direction(P_k,P_0)$\\
		$\pmod{p^2-1}.$ 
	\item$	angular\: direction(O,M_i) = i + angular\: direction(0,M_0),$
		where $M_i$ is the mid-point of $(P_i, P_{i+1})$.
	\item$	angular\: direction(O,N_i) = i + angular\: direction(0,N_0),$
		where $N_i$ is the mid-point of $(P_{i-1}, P_{i+1})$.
	\item$	angular\:direction(P_i,N_i) = i + angular\:direction(P_0,N_0).$
	\enume

	\sssec{Example.}
	$p = 17,$ generator $(13,4,3),$\\
	$angular\:direction((0,0,1),(13,4,3)) = 164,$\\
	$angular\:direction((0,0,1),(9,6,1)) = 224,$\\
	$M_0 = (15,2,2),$ $angular\:direction(O,M_0) = 60,$\\
	$N_1 = (13,3,6),$ $angular\:direction(O,N_0) = 33,$\\
	$angular\:direction(P_1,N_1) = 40.$

	\sssec{Corollary.}
	{\em The coordinates of the normal to the surface ${\cal T}$  are those
	of $-p(a,b,c)$}\footnote{28.12.87}.

	\sssec{Lemma.}
	{\em If $(a,b,c)$ is a generator and $g^i$ is a primitive root of $p,$
	then}
	\enumb
	\item $i,$ prime, $\equiv  1 \pmod{3}  \Rightarrow
	   (g^i)^{\delta } = \alpha .$?
	\item $i,$ prime, $\equiv  2 \pmod{3}  \Rightarrow
	   (g^i)^{\delta } = \beta .?$
	\item$i$ {\em is not a prime} $\Rightarrow (g^i)^{\delta} = \epsilon.$
	\enume

	\sssec{Example.}
	The following table gives generators $(a,b,c)$ for
	the given values of $p$ and $g,$\\
		$\alpha$,$\epsilon$ $,\beta$ $-\epsilon$ $,\alpha$ $,\beta$
 		$\alpha$	$,\beta$
	 $,\epsilon$ $-\alpha$ $,\epsilon$ $,\beta$ 		$\alpha$	
	$,\alpha$ $,\alpha$ $-\beta$ $,\beta$ $,\beta$ \\
	$\begin{array}{llcllcll}
	p = 11,& 4,9,0	&\vline&p = 5, & 0,4,3	&\vline&p = 17,& 13,4,3\\
	g = 2,	&1,7,5  &\vline&g = 2,&	 	&\vline&g = 3,&	15,16,6\\
	\cline{1-2}\cline{4-5}\cline{7-8}
	p = 29, &0,17,14 &\vline&p = 23,& 11,1,16&\vline&p = 53,& 17,38,0\\
	g = 2,	&4,24,3	&\vline&g = 5,&	5,21,2  &\vline&g = 2,&2,34,19\\
		&7,16,8  &\vline&&	6,13,9  &\vline&&	23,25,7\\
	 	&11,26,23&\vline&& 	18,10,0	&\vline&&	48,52,8\\
	\cline{1-2}\cline{4-5}
	p = 47,	&0,30,22&\vline&p = 41,& 36,1,10&\vline&&	11,24,20\\
	g = 5,	&1,12,39&\vline&g = 6,&	35,3,9	&\vline&&12,49,47\\
		&2,4,46	&\vline&&	30,12,5	&\vline&p = 71,& 69,7,2\\
		&6,25,21&\vline&&	17,16,14&\vline&g = 7,&36,39,3\\
		&18,7,27&\vline&&	40,33,15&\vline&&	4,43,31\\
		&13,45,41&\vline&& 	31,29,28&\vline&&	21,48,9\\
		\cline{4-5}
		&17,15,20&\vline&p = 59,& 53,0,8&\vline&&	13,49,16\\
	 	&19,44,36&\vline&g = 2,&44,15,2	&\vline&&	17,67,65\\
	\cline{1-2}
	p = 83,	&62,23,0&\vline&&	30,28,3	&\vline&&	56,58,35\\
	g = 2,	&2,20,63&\vline&&	48,4,9	&\vline&& 	50,54,45\\
			\cline{7-8}
		&4,11,70&\vline&&	19,6,36	&\vline&p = 89,&39,48,5\\
		&50,29,6&\vline&&	32,17,12&\vline&g = 3,&	8,13,71\\
		&8,16,61&\vline&&	27,21,13&\vline&&	37,45,10\\
		&27,9,49&\vline&& 	18,52,50&\vline&&	29,44,19\\
		\cline{4-5}
		&77,10,81&\vline&&		&\vline&&	20,83,78\\
		&21,13,51&\vline&&		&\vline&&	43,86,52\\
		&25,18,42&\vline&&		&\vline&&	46,81,54\\
		&35,31,19&\vline&&		&\vline&& 	47,72,62\\
			\cline{7-8}
		&45,64,59&\vline&&		&\vline\\
	 	&67,53,48&\vline&&		&\vline\\
	\cline{1-2}
	\end{array}$

	\sssec{Example.}
	\enumb
	\item$  p = 5, {\cal T}  =$\\
	$\begin{array}{rrrrrrrrrrrr}
	0& 1& 2& 3& 4& 5& 6& 7& 8& 9&10&11\\
	\cline{1-12}
	 0& 0& 1&-1&-1&-1& 2&-1& 1&-2&-1& 1\\
	 0&-1&-1&-2&-2& 0& 2& 1& 0& 0& 1&-1\\
	 1&-2&-1& 1&-1&-2& 0&-2& 0&-1&-1&-2\\

	12&13&14&15&16&17&18&19&20&21&22&23\\
	\cline{1-12}
	 -1&-2& 0&-2& 0&-1&-1&-2&-2& 0& 2& 1\\
	 -1&-1& 2&-1& 1&-2&-1& 1&-1&-2& 0&-2\\
	 -2& 0& 2& 1& 0& 0& 1&-1&-1&-1& 2&-1
	\end{array}$

	The ideal points are (last coord. 0),
	A	  B	    C	      D	       E	F	G
	(0,1,-1), (1,0,-1), (1,-1,0), (1,1,1), (1,1,-2),(1,2,2),(1,-2,1).\\
	successive powers, A,G,E,F,B,C,0; B,F,E,G,A,C,0; C,C; D,D;
		E,D,E; F,G,D,F.

	\item$  p = 11, g = 2, {\cal T}  =$\\
	$\begin{array}{rrrrrrrrrrrrrrrrr}
	0& 1& 2& 3& 4& 5& 6& 7& 8& 9&10&11&12&13&14&15&16\\
	\cline{1-17}
	0& 1&-5& 0&-2& 4& 5& 1&-3&-4& 4& 3&-3& 5& 4& 4&-4\\
	0& 3& 3&-4& 3&-3&-3&-4& 0&-1&-5& 1&-3&-3& 5&-2&-2\\
	1& 4& 0&-1& 3&-2& 1& 5&-3& 1& 2& 4& 4& 4&-5&-3&-2\\
	\cline{1-17}
	17&18&19&20&21&22&23&24&25&26&27&28&29&30&31&32&33\\
	\cline{1-17}
	-2& 4& 1&-3&-5& 3&-3&-4&-5& 5&-2&-4& 2& 2&-3&-2&-4\\
	4& 3&-1&-4& 5&-5& 5&-4& 5& 1& 0&-1& 0&-5& 5& 2& 0\\
	0&-2&-4&-3&-3& 0& 4& 1&-1&-3& 4&-1& 5& 4&-5&-2&-1\\
	\cline{1-17}
	34&35&36&37&38&39&40&41&42&43&44&45&46&47&48&49&50\\
	\cline{1-17}
	5& 5& 4& 1& 4& 2& 1& 4& 0&-1& 3&-2& 1& 5&-3& 1& 2 \\
	4&-5& 4& 5&-2& 5& 0& 1&-5& 0&-2& 4& 5& 1&-3&-4& 4 \\
	-5&-1&-5&-4& 3& 0& 0& 3& 3&-4& 3&-3&-3&-4& 0&-1&-5
	\end{array}$\\
	\ldots 

	\item$ p = 17,$ $g = 3,$ ${\cal T}  = \{$\\
	$\begin{array}{cccccc}
	( 0, 0, 1),&( 1, 2, 3),&(10,13,13),&( 1, 7, 4),&( 4,13, 4),
		&( 8, 0,16),\\
	( 6, 6,13),&( 9,16, 6),&(14, 1, 1),&(11, 2,15),&( 1,13, 1),
		&(13, 8, 1),\\
	( 5, 5, 3),&(11, 9, 7),&( 7, 1, 1),&( 7,12, 1),&(12,11,12),
		&( 2, 1, 3),\\
	(11,11,14),&( 1, 4, 7),&( 1,10,10),&(16, 0, 8),&( 5,15, 5),
		&(16, 9, 6),\\
	( 4, 4, 8),&(11,15, 2),&(14, 9, 9),&( 1, 8,13),&(15, 0,15),
		&( 9,11, 7),\\
	( 5, 5,16),&( 7, 1,12),&( 1, 0, 0),&( 3, 1, 2),&(13,10,13),&( 4, 1, 7),
	\end{array}$\\
	\ldots.\}
	\enume

	\sssec{Example.}
	A generator associated to the given primitive root is such
	that $(a,b,c)^{\delta } = \alpha = (1,0,0),$ with
	$\delta := \frac{p^2-1}{3}.$ \\
	$\begin{array}{rcrcl}
	 p &\vline&  g &\vline generator\\
	\cline{1-5}
	  5 &\vline&  2 &\vline&  (0,4,3)\\
	 11 &\vline&  2 &\vline&  (1,3,4)\\
	 17 &\vline&  3 &\vline&  (3,4,13)\\
	 23 &\vline&  5 &\vline&  (0,7,12)\\
	 29 &\vline&  2 &\vline&  (0,4,7)\\
	 41 &\vline&  6 &\vline&  (0,2,17)\\
	 47 &\vline&  5 &\vline&  (0,7,37)\\
	 53 &\vline&  2 &\vline&  (0,6,20)\\
	 59 &\vline&  2 &\vline&  (0,3,11)\\
	 71 &\vline&  7 &\vline&  (0,8,54)\\
	 83 &\vline&  2 &\vline&  (0,6,60)\\
	 89 &\vline&  3 &\vline&  (0,6,77)\\
	101 &\vline&  2 &\vline&  (0,7,60)\\
	107 &\vline&  2 &\vline&  (0,2,29)\\
	113 &\vline&  3 &\vline&  (0,3,-12)???\\ 
	131 &\vline&  2 &\vline&  (0,4,18)\\
	137 &\vline&  3 &\vline&  (0,2,-25)\\
	149 &\vline&  2 &\vline&  (0,3,-53)\\
	167 &\vline&  5 &\vline&  (0,7,-2)\\
	173 &\vline&  2 &\vline&  (0,3,71)\\
	179 &\vline&  2 &\vline&  (0,4,36)\\
	191 &\vline& 19 &\vline&  (0,5,-83)\\
	197 &\vline&  2 &\vline&  (0,19,61)\\
	\end{array}$\\
	See [M130] RICATTI. for more.



\setcounter{section}{3}
\setcounter{subsection}{2}
	\ssec{The case of 4 Functions.}\label{sec-Sric4}

	\sssec{Definition.}
	The set $R_4$ is the set of elements
	\enumb
	\item$  (x,y,z,t)  \ni  x,y,z,t \in  Z_p$ and\\
	$-(x^2-z^2)^2 + (y^2-t^2)^2 + 4((x^2+z^2)yt - (y^2+t^2)xz) = 1,$\\
	with addition
	\item$  (x,y,z,t) + (x',y',z',t') = (xt'+tx'+yz'+zy',xx'+zz'+yt'+ty',$\\
	\hth$				xy'+yx'+zt'+tz',yy'+tt'+xz'+zx').$
	\enume

	\sssec{Theorem.}
	{\em $(R^4,+)$ is an Abelian group.}

	\sssec{Conjecture.}
	\enumb
	\item{\em If $p \equiv  1 \pmod{4}$ then the maximum period is $p-1.$}
	\item{\em If $p \equiv -1 \pmod{4}$ then the maximum period is $p^2-1$}
	\footnote{22.12.87}.
	\enume

	\sssec{Example.}
	\enumb
	\item	If $p = 3,$ $(1,2,0,1),$ is of period 8.
	\item	If $p = 5,$ $(1,2,0,1),$ is of period 4.
	\item	If $p = 7,$ $(3,3,0,4),$ is of period 48.
	\item	If $p = 11,$ $(3,2,6,3),$ is of period 120.
	\item	If $p = 13,$ $(1,2,3,9),$ is of period 12.
	\item	If $p = 17,$ $(15,13,4,8),$ is of period 16.
	\item	If $p = 19,$ $(12,13,14,1),$ is of period 360.
	\item	If $p = 23,$ $(2,3,6,17),$ is of period 528.
	\item	If $p = 29,$ $(15,16,11,18),$ is of period 28.
	\enume

	\sssec{Lemma.}
	$\matb{|}{ccc}z&y&x\\x&z&y\\y&x&z\mate{|}
	= \matb{|}{ccc}x+y+z&y&x\\x+y+z&z&y\\x+y+z&x&z\mate{|}
	= (x+y+z)\matb{|}{ccc}1&y&x\\1&z&y\\1&x&z\mate{|}\\
	\hth= (x+y+z)\matb{|}{ccc}1&y&x\\0&z-y&y-x\\0&x-z&z-y\mate{|}
	= (x+y+z) ((z-y)^2-(x-z)(y-x)).$

	\sssec{Lemma.}
	$\matb{|}{cccc}t&z&y&x\\x&t&z&y\\y&x&t&z\\z&y&x&t\mate{|} =
	(x+y+z+t)\matb{|}{cccc}1&z&y&x\\1&t&z&y\\1&x&t&z\\1&y&x&t\mate{|}\\
	\hth= (x+y+z+t)\matb{|}{ccc}t-z&z-y&y-x\\x-t&t-z&z-y\\y-x&x-t&t-z
	\mate{|}\\
	\hth =
	(x+y+z+t)\matb{|}{ccc}t-z+y-x&z-y&y-x\\x-t+z-y&t-z&z-y\\y-x+t-z&x-t&t-z
	\mate{|}\\
	\hth = (x+y+z+t)(-x+y-z-t)
	\matb{|}{ccc}1&z-y&y-x\\-1&t-z&z-y\\1&x-t&t-z\mate{|}\\
	\hth = (x+y+z+t)(-x+y-z+t)
	\matb{|}{ccc}1&z-y&y-x\\0&t-y&z-x\\x-z&t-y\mate{|}\\
	\hth = (x+y+z+t)(-x+y-z+t) ((t-y)^2+(x-z)^2).$ 

	\ssec{The case of 5 functions.}\label{sec-Sric5}
	\sssec{Definition.}
	\hth$det(x,y,z,t,u) = \matb{|}{ccccc}
	u&t&z&y&x\\x&u&t&z&y\\y&x&u&t&z\\z&y&x&u&t\\t&z&y&x&u\mate{|}$

	\sssec{Definition.}
	The set $R_5$ is the set of elements\\
	\enumb
	\item$  (x,y,z,t,u)  \ni  x,y,z,t,u \in  Z_p$ and\\
	\hth$	det(x,y,z,t,u) = 1$\\
	with addition
	\item$  (x,y,z,t,u) + (x',y',z',t',u')$\\
	$ = (xt'+yu'+zt'+tz'+uy',xy'+yx'+zu'+tt'+uz',$\\
	\hth$	xz'+yy'+zx'+tu'+ut',xt'+yz'+zy'+tx'+uu',xu'+yt'+zz'+ty'+ux').$
	\enume

	\sssec{Lemma.}
	{\em
	\hth$\matb{|}{cccc}u&t&z&y\\x&u&t&z\\y&x&u&t\\z&y&x&u\mate{|}\\
	\hth= (u^2-xt)^2 - (tu-xz)*(xu-yt)\\
	\hti{12} + (zu-xy)*(x^2-yu) + (t^2-zu)*(yu-zt)\\
	\hti{12} - (zt-yu)*(xy-zu) + (z^2-yt)(y^2-xz)$\\
	\hth$= u^4 - x^3y - y^3t - z^3x - t^3z + x^2t^2 + y^2z^2$\\
	\hti{12}$+ 2x^2zu + 2y^2xu + 2z^2tu + 2t^2uy - 3u^2xt - 3u^2yz - xyzt$}

	\sssec{Theorem.}
	{\em $det(x,y,z,t,u) = s
		(2(x^4+y^4+z^4+t^4+u^4) - s(x^3+y^3+z^3+t^3+u^3)$\\
	\hti{12}$	+ x^2(y(2z+2t-3u)-3zt+2tu+2uz) +  \ldots $\\
	\hti{12}$	- yztu - ztux - tuxy - uxyz - xyzt),$\\
	with $s = x+y+z+t+u.$}

	Proof:  We use\\
	\hth$\matb{|}{ccccc}
	u&t&z&y&x\\x&u&t&z&y\\y&x&u&t&z\\z&y&x&u&t\\t&z&y&x&u\mate{|}
	= s\matb{|}{ccccc}
	1&t&z&y&x\\1&u&t&z&y\\1&x&u&t&z\\1&y&x&u&t\\1&z&y&x&u\mate{|}$\\
	and then the Lemma.

	\sssec{Theorem.}
	{\em $(R^5,+)$ is an Abelian group.}

	\sssec{Conjecture.}
	\enumb
	\item{\em If $p \equiv  1 \pmod{10}$ then the maximum period is $p-1.$}
	\item{\em If $p \equiv  9 \pmod{10}$ then the maximum period is $p^2-1$}
	\footnote{24.12.87}.
	\item{\em If $p \equiv  \pm 3 \pmod{10}$ then the maximum period is
	$p^4-1$}\footnote{22.12.87}.
	\enume

	For examples see \ref{sec-ericn}.




\setcounter{section}{3}
	\section{Application to geometry.}
	\setcounter{subsection}{-1}
	\ssec{Introduction.}
	To define distances in a sub geometry of affine
	$k$-dimensional geometry, we have to define a homogeneous function
	$f(P)$ of degree $k.$  We can then either define the distance between
	2 points $P$ and $Q$ by the $k$-th root of $f(Q-P)$ or the hypercube
	between 2 points $P$ and $Q$ by $f(Q,P)$.  I will not discuss here the
	extension of a 2-dimensional distance to $n$-dimension as is done in
	Euclidean geometry.\\
	To define angles, we can associate to a point $P,$ a $k$ by $k$ matrix
	by a bijection, if the set of these matrices, which are of determinant
	1, form a subset of an Abelian group under matrix multiplication, with
	generator $G_0,$ \ldots $G_l,$ we can define then angular direction of
	a point associated to the matrix $G_0^{i_0}$ \ldots $G_l^{i_l}$ by
	$(i_0,$ \ldots $,i_l).$\\
	We can also define $f(P)$ as the determinant of the associated matrix.
	If in the 2 dimensional real affine geometry, we associate to $(x,y)$
	the matrix\\
	\hth$\ma{y}{x}{-x}{y},$\\
	then $f(x,y) = x^2+y^2,$ the matrices of determinant 1 for an Abelian
	group with generator $x = sin(1),$ $y = cos(1),$ and we obtain the
	$2$-dimensional Euclidean distance and angle.\\
	If in the $2$-dimensional real affine geometry, we associate to $(x,y)$
	the matrix\\
	\hth$\ma{y}{x}{x}{y},$\\
	then $f(x,y) = y^2-x^2,$ the matrices of determinant 1 for an Abelian
	group with generator $x = sinh(1),$ $y = cosh(1),$ and we obtain the
	$2$-dimensional Minkowskian distance and angle.\\
	This will now be extended using the generalization of the hyperbolic
	functions by Ricatti.

	\ssec{k-Dimensional Affine Geometry.}
	\sssec{Definition.}
	In $k$-dimensional affine geometry we define the {\em Ricatti
	function} as the function which associates to the point
	$P = (P_0, \ldots ,P_{k-1}),$ the Toeplitz matrix ${\bf T}$, defined by
	\hth$	{\bf T}_{i,j} = P_{k-1-i+j},$ $0 \leq  i,j < k,$\\
	where the subscripts computation is done modulo $k.$

	\sssec{Theorem.}
	{\em The matrix multiplication defines an addition (which is a
	convolution) for the points as follows,
	if ${\bf T}$ is associated to $P$ and ${\bf U}$, to $Q,$
	${\bf T U}$ is associated to $R$ with\\
	\hth$	P \circ Q := R_i = \sum_j  P_j Q_{i-1-j}.$\\
	For instance,}
	\enumb
	\item$  k = 3,$\\
	\hth$	(P \circ Q)_0 = P_0 Q_2 + P_1 Q_1 + P_2 Q_0,$\\
	\hth$	(P \circ Q)_1 = P_0 Q_0 + P_1 Q_2 + P_2 Q_1,$\\
	\hth$	(P \circ Q)_2 = P_0 Q_1 + P_1 Q_0 + P_2 Q_2.$
	\item$  k = 4,$\\
	\hth$	(P \circ Q)_0 = P_0 Q_3 + P_1 Q_2 + P_2 Q_1 + P_3 Q_0,$\\
	\hth$	(P \circ Q)_1 = P_0 Q_0 + P_1 Q_3 + P_2 Q_2 + P_3 Q_1,$\\
	\hth$	(P \circ Q)_2 = P_0 Q_1 + P_1 Q_0 + P_2 Q_3 + P_3 Q_2,$\\
	\hth$	(P \circ Q)_3 = P_0 Q_2 + P_1 Q_1 + P_2 Q_0 + P_3 Q_3.$
	\enume

	\sssec{Corollary.}
	{\em The set of matrices, associated to all the non ideal points of
	$k$-dimensional affine geometry with determinant 1, form an
	abelian group under matrix multiplication.}

	\sssec{Theorem.}
	{\em If $p = k$ then $P_i = \delta_{i,0}$ has period $p$.  Moreover, the
	$j$-th iterate $P^{(j)}$ is such that $P^{(j)}_i = \delta_{j,i}.$}

	\sssec{Theorem.}
	{\em Let $det(\ldots,z,y,x)$ denote the determinant of the Toeplitz
	matrix associated with $P = (x,y,z,\ldots),$ let $s$ be the sum of the
	components of $P,$ then}
	\enumb
	\item$	k = 3,$\\
	\hti{4}$det(z y x) = x^3+y^3+z^3-3xyz = s ((z-y)^2-(x-z)(y-x)).$
	\item$	k = 4,$\\
	\hti{4}$det(t z y x) = s (-x+y-z+t) ((t-y)^2+(x-z)^2).$
	\item$	k = 5,$\\
	\hti{4}$det(x,y,z,t,u) = s (2(x^4+y^4+z^4+t^4+u^4)
		- s(x^3+y^3+z^3+t^3+u^3)$\\
	\hth$		+ x^2(y(2z+2t-3u)-3zt+2tu+2uz) +  \ldots $\\
	\hth$		- yztu - ztux - tuxy - uxyz - xyzt).$
	\enume

	In the following examples we have obtained what a cyclic generator
	of what appears to be the longest period, without examining the details
	of the structure of the solution.

	\sssec{Example.}\label{sec-ericn}
	\enumb
	\item$	k = 4.$\\
	$\begin{array}{rrl}
	p&	period&	cyclic\: generator\\
	\cline{1-3}
	3&	8&	(1,2,0,1),\\
	5&	4&	(1,2,0,1),\\
	7&	48&	(3,3,0,4),\\
	11&	120&	(3,2,6,3),\\
	13&	12&	(1,2,3,9),\\
	17&	16&	(15,13,4,8),\\
	19&	360&	(12,13,14,1),\\
	23&	528&	(2,3,6,17),\\
	29&	28&	(15,16,11,18).
	\end{array}$

	\item$	k = 5.$\\
	$\begin{array}{rrl}
	p&	period&	cyclic\: generator\\
	\cline{1-3}
	3&	80&	(1,1,0,1,2),\\
	5&	5&	(1,0,0,0,0),\\
	7&	2400&	(1,2,4,1,4),\\
	11&	10&	(4,2,1,4,2).\\
	13&	28560&	(3,5,1,11,8)\\
	17&	83520&	(9,7,8,2,11),\\
	19&	18&	(7,16,15,2,0),\\
	23&	279840&	(14,12,4,7,14),\\
	29&	840&	(13,8,25,5,9),\\
	31&	30&	(26,30,11,2,27).
	\end{array}$

	\item$	k = 6,$\\
	$\begin{array}{rrl}
	p&	period&	cyclic\: generator\\
	\cline{1-3}
	3&	6&	(0,1,1,2,1,0),\\
	5&	24&	(0,2,4,1,0,0),\\
	7&	6&	(3,5,0,6,1,2),\\
	11&	120&	(5,9,2,4,2,2),\\
	13&	12&	(3,8,4,2,1,10)\\
	17&	288&	(2,12,5,14,4,3),\\
	19&	18&	(2,12,5,14,4,3),\\
	23&	528&	(3,16,20,4,13,18),\\
	29&	840&	(10,2,8,22,14,4),\\
	31&	30&	(1,11,7,4,29,13),
	\end{array}$

	\item$	k = 7,$\\
	$\begin{array}{rcrl}
	p&	period&	cyclic\: generator\\
	\cline{1-4}
	3&&	728&	(1,2,2,1,0,1,1,1),\\
	5&&	15624&	(0,0,4,3,0,0,0),\\
	7&&	7&	(1,0,0,0,0,0,0),\\
	11&&	1330&	(0,1,8,4,3,4,4),\\
	13&&	168&	(6,3,11,2,4,2,0),\\
	17&&	24137568&(3,10,13,16,3,4,5),\\
	19&&	47045880&(18,10,17,5,14,9,5),\\
	23&&	12166&	(0,17,4,7,3,15,5),\\
	29&&	28&	(26,7,9,10,15,7,15),\\
	31& (6)&887503680&(26,26,15,18,26,22,19),\\
	37& (3)&50652&	(9,9,18,31,0,27,25),\\
	41& (2)&1680&	(7,10,22,32,27,2,27).
	\end{array}$

	\item$	k = 8,$\\
	$\begin{array}{rrl}
	p&	period&	cyclic\: generator\\
	\cline{1-3}
	3&	8&	(1,2,2,1,0,1,1,1),\\
	5&	24&	(2,0,3,3,1,4,4,0),\\
	7&	48&	(1,3,0,0,3,3,4,3),\\
	11&	120&	(5,1,9,9,4,8,2,8),\\
	13&	168&	(7,11,5,4,12,9,5,1),\\
	17&	16&	(9,13,7,10,0,15,4,3),\\
	19&	360&	(18,11,4,13,8,1,7,16),\\
	23&	528&	(9,10,22,4,8,17,16,11),\\
	29&	840&	(28,1,14,21,9,26,14,5),\\
	31&	960&	(28,6,30,20,25,1,30,18),\\
	37&	1368&	(0,21,5,5,28,36,9,9),\\
	41&	40&	(16,30,30,27,14,18,18,17),
	\end{array}$

	\item$	k = 9,$\\
	$\begin{array}{rrll}
	p	&period&	cyclic\: generator\\
	\cline{1-3}
	3	&18&	(0,2,2,1,2,2,0,1,1),\\
	5	&15624&	(1,1,2,0,1,1,4,1,1),& may\: not\: be\: largest\:period\\
	7\\
	11\\
	13\\
	17	&288&	(15,8,16,15,9,13,7,10,12),\\
	19	&18&	(2,4,15,11,6,11,4,0,6),\\
	23\\
	29\\
	31\\
	37	&36&	(5,27,14,9,28,24,20,12,11)
	\end{array}$
	\enume

	\sssec{Example.}
	Here we have written $j$ when the maximum period is $p^j-1,$ 
	unless the number is underlined in which case the period is given.

	$\begin{array}{rcllllllllllll}
k\setminus p&\vline&3&5&7&11&13&17&19&23&29&31&37&41\\
	\cline{1-14}
	3&\vline&\underline{k}&2&1&2&1&2&1&2&2&1&1&2\\
	4&\vline&2&1&2&2&1&1&2&2&1&2&1&1\\
	5&\vline&4&\underline{k}&4&1&4&4&2&4&2&1&4&1\\
	6&\vline&\underline{k}&2&1&2&1&2&1&2&2&1&1&2\\
	7&\vline&6&6&\underline{k}&3&2&6&6&3&1&6&3&2\\
	8&\vline&2&2&2&2&2&1&2&2&2&2&2&1\\
	9&\vline&\underline{k}&6&3&6&3&2&1&6&6&3&1&6\\
	10&\vline&4&\underline{2k}&4&1&4&4&2&4&2&1&4&1\\
	11&\vline&5&5&10&\underline{k}&10&10&10&1&10&5&5&\geq 10\\
	12&\vline&\underline{2k}&2&2&2&1&2&2&2&2&2&1&2\\
	13&\vline&3&4&12&12&\underline{k}&6&\geq 12&6&3&4&\geq 10&\geq 10\\
	14&\vline&6&6&\underline{3k}&3&2&6&6&3&1&6&3&2\\
	15&\vline&\underline{16k}&\underline{8k}&4&2&4&4&2&4&2&1&4&2\\
	16&\vline&4&4&2&4&4&1&4&2&4&2&4&2\\
	17&\vline&16&16&16&\geq 16&4&-&8&\geq 16&\geq 16&\geq 16&?&?\\
	19&\vline&18&5&3&3&\geq 18&9&-&9&\geq 18&6&2&?
	\end{array}$

	$\begin{array}{rcllllllllllll}
k\setminus p&\vline&43&47&53&59&61&67&71&73&79&83&89&97\\
	\cline{1-14}
	3&\vline&1&2&2&2&1&1&2&1&1&2&2&1\\
	4&\vline&2&2&1&2&1&2&2&1&2&2&1&1\\
	5&\vline&4&4&4&2&1&4&1&4&2&4&2&4\\
	6&\vline&1&2&2&2&1&1&2&1&1&2&2&1\\
	7&\vline&1&6&3&6&6&3&1&6&3&2&6&2\\
	8&\vline&2&2&2&2&2&2&2&1&2&2&1&1\\
	9&\vline&3&6&2&6&3&3&2&1&3&6&2&3\\
	10&\vline&4&4&4&2&1&4&1&4&2&4&2&4\\
	11&\vline&2&5&5&5&\geq 10&1&5&\geq 10&\geq 10&\geq 10&1&5\\
	12&\vline&2&2&2&2&1&2&2&1&2&2&2&1\\
	13&\vline&6&4&1&\geq 8&3&\geq 8&\geq 8&4&1&4&\geq 7&\geq 7\\
	14&\vline&1&6&3&6&6&3&1&6&3&2&6&2\\
	15&\vline&4&4&4&2&1&4&2&4&2&4&2&4\\
	16&\vline&4&2&4&4&4&4&2&2&2&4&2&1\\
	17&\vline&8&?&?&8&\geq 16&2&?&?&?&\geq 8&4&?\\
	19&\vline&?&9&\geq 9&\geq 9&?&?&?&\geq 9&101&\geq 9&103&\geq 9
	\end{array}$

	Moreover it appears that
	\begin{verbatim}
	{\em k |      18 |     20 |     21 |     22 |     24 |     25 |
	     26 |     27}
	  |  3	 k |  5	 k | 3 104k |  11 5k |  3	 k |  5	 k |
	  13 2k |  3   k
	  |   	   |   	   |  7  6k |        |   	   |   	   |        | 

	{\em k |      28 |     30 |     33 |     34 |     35 |     36 |
	     38 |     39}
	  |  7  12k |  3	8k |  3 22k |  17 $8kL\hti{-1}/$ |
	  5	.k |  3	.k |  19 9k |  3  .k
	  |   	   |  5  4k | 11 40k |        | 7 560k |   	   |
	        | 13  .k

	{\em k |      40 |     42 |     44 |     45 |     46 |     48 |
	     49 |     50}
	  |  5	3k |  3	.k | 11 15k |  3  .k | 23 11k |  3	5k |
	  7  .k |  5   k
	  |   	   |  7	.k |        |  5  .k |   	   |   	   |        | 

	{\em k |      51 |     52 |     54 |     55 |     56 |     57 |
	     58 |     60}
	  |  $3\geq$ 200k | 13	3k |  3     |  11 5k |  7	8k |  3	.k |
	  29 .k |  3  .k
	  | 17	.k |   	   |        |        |   	   | 19	.k |        |
	  5  .k
	\end{verbatim}

	\sssec{Conjecture.}
	If $k = 4,$ then
	\enumb
	\item	if $p \equiv  1 \pmod{4}$ then the maximum period is $p-1.$
	\item 	if $p \equiv  -1 \pmod{4}$ then the maximum period is $p^2-1$
	\footnote{22.12.87}.
	\enume

	\sssec{Conjecture.}
	If $k = 5,$
	\enumb
	\item	if $p \equiv 1 \pmod{10}$ then the maximum period is $p-1,$
	\item 	if $p \equiv 9 \pmod{10}$ then the maximum period is $p^2-1$
	\footnote{24.12.87},
	\item	if $p \equiv \pm 3 \pmod{10}$ then the maximum period is
	$p^4-1$\footnote{24.12.87}.
	\enume

	The above examples may lead to other conjectures perhaps for all $k.$

	\sssec{Conjecture.}
	Let $p |\hti{-1}/ k.$  The maximum period is $p^e-1,$ where $e$ depends
	on $p$ and $k,$ \footnote{2.2.88}
	\enumb
	\item$	e(p^i,p') = e(p^i,p") if p' \equiv  p" \pmod{p^i}.$
	\item$	(q_1,q_2) = 1 \Rightarrow e(k,q_1q_2) = lcm(e(k,p_1),e(k,p_2).$
	\item$	e(p^i,p') = order(p') \in {\Bbb Z}_{p,.}.$\\
	In view of 0, we can define $e_u := e(p^i,u)$ for $u \in
	{\Bbb Z}_{p^i,.}$ 
	\item$	e(p^i,1) = 1,$ $e(p^i,p^i-1) = 2,$
	$k = 5,$ $e_2 = e_3 = 4,$\\
	$k = 7,$ $e_2,e_4 = 3,$ $e_3,e_5 = 6,$\\
	$k = 2^3,$ $e_3,e_5 = 2,$\\
	$k = 3^2,$ $e_4,e_7 = 3,$ $e_2,e_5 = 6,$\\
	$k = 11,$ $e_3,e_4,e_5,e_9 = 5,$ $e_2,e_6,e_7,e_8 = 10,$\\
	$k = 13,$ $e_3,e_9 = 3,$ $e_5,e_8 = 4,$ $e_4,e_{10} = 6,$
		$e_2,e_6,e_7,e_{11} = 12,$
	$k = 17,$ $e_3,$ $e_4,e_5,e_9 = 5,$ $e_2,e_6,e_7,e_8 = 10,?$
	\enume

	\sssec{Theorem.}
	{\em If, for $k = 3,$ $p \equiv 5 \pmod{6},$ we construct a period
	associated to a generator and determine the coplanar directions to
	the directions associated to 0 and 1, we obtain a difference sets
	For the set ${\Bbb Z}_{p^2,.}$ of the numbers from $0$ to $p^2$
	relatively prime to $p.$\\
	The sets have $p(p-1)$ elements.}

	Proof:  The proof is similar to that of Singer.  In this case, the
	directions are the non isotropic ones and 2 non isotropic directions
	determine exactly one plane, through the origin, which contains $p+1$
	directions.

	This Theorem extends to any dimension.  We should check if these
	difference sets are also obtained by some other method.

	\sssec{Example.}\label{sec-eric3gen}
	$k = 3,$  ([130$\setminus$RIC.BAS] $p,$ then diff. set then generator)\\
	$\begin{array}{rll}
	p&   gen.&	difference\: set \pmod{p^2-1}\: of\: p(p-1)\: elements\\
	\cline{1-3}
	5 & (0,4,3)&	0,1,14,16,21\\
	11& (1,3,4)&	0,1,9,28,30,34,41,44,83,98,103\\
	17& (3,4,13)&	0,1,10,13,34,45,59,86,112,114,129,134,191,195,251,259,\\
	&&		282\\
	23& (0,7,12)&	0,1,60,91,134,142,148,203,249,253,266,269,271,298,305,\\
	&&		333,342,352,363,375,450,488,503\\
	29& (0,4,7)&	0,1,134,147,153,228,246,316,326,328,373,411,432,435,\\
	&&		452,457,484,488,521,549,560,575,589,623,719,774,790,\\
	&&		797,832\\
	41& ( 0,2,17 )&0,1,24,199,208,230,424,470,522,525,533,604,682,684,694,\\
	&&	698,748,775,805,823,872,879,915,941,975,1014,1061,1120,\\
	&&	1133,1161,1178,1248,1263,1283,1316,1527,1548,1567,1592,\\
	&&		1643,1675,\\
	47& ( 0,7,37 )&	0,1,8,115,147,253,373,401,412,447,693,714,716,765,889,\\
	&&		923,964,982,994,1095,1124,1182,1185,1258,1303,1308,\\
	&&	1322,1339,1419,1472,1519,1655,1744,1757,1782,1822,1826,\\
	&&		1842,1848,1910,1925,1934,1967,1977,2004,2099,2153,\\
	53& ( 0,6,20 )&	0,1,28,42,59,133,183,194,218,239,339,385,404,497,499,\\
	&&		548,695,721,773,783,805,820,843,849,922,958,962,1048,\\
	&&	1056,1226,1251,1256,1290,1333,1623,1680,1854,1872,1925,\\
	&&	1941,2022,2102,2191,2194,2203,2266,2314,2321,2334,2417,\\
	&&		2450,2554,2621,\\
	59& ( 0,3,11 )&	0,1,243,331,362,386,448,469,488,598,625,734,814,816,\\
	&&	825,839,912,915,969,1012,1134,1227,1484,1626,1633,1667,\\
	&&	1744,1761,1773,1819,2083,2151,2275,2320,2364,2379,2435,\\
	&&	2527,2543,2549,2596,2717,2737,2798,2802,2840,2850,2868,\\
	&&	2876,3022,3071,3101,3106,3138,3233,3272,3305,3417,3430,\\
	\end{array}$

	\newpage
	$\begin{array}{rll}
	p&   gen.&	difference\:set \pmod{p^2-1}\: of\: p(p-1)\: elements\\
	71& ( 0,8,54 )&	0,1,339,345,406,542,687,821,907,989,1171,1294,1429,\\
	\cline{1-3}
	&&	1443,1502,1522,1553,1617,1628,1650,1691,1737,1792,1828,\\
	&&	1946,2108,2125,2229,2237,2247,2266,2281,2461,2500,2503,\\
	&&	2550,2655,2743,2768,2966,2970,3019,3028,3035,3127,3195,\\
	&&	3280,3328,3360,3405,3426,3431,3617,3912,3996,4019,4031,\\
	&&	4162,4273,4343,4400,4460,4490,4514,4590,4592,4630,4673,\\
	&&		4686,4836,5013,\\
	83& (0,6,60)&	0,1,182,187,214,255,500,503,565,590,596,827,1353,1389,\\
	&&	1406,1456,1501,1555,1577,1629,1690,1720,1900,2039,2067,\\
	&&	2136,2250,2261,2265,2336,2645,2704,2737,2783,2785,2792,\\
	&&	2984,3250,3271,3479,3641,3711,3723,3746,3760,3868,3902,\\
	&&	3953,4053,4063,4071,4194,4296,4309,4353,4459,4568,4592,\\
	&&	4611,4675,4722,4738,4764,4896,4973,5013,5093,5191,5230,\\
	&&	5346,5366,5490,5550,5616,5654,5710,5844,5922,6279,6337,\\
	&&	6611,6683,6712,\\
	89& (0,6,77)&0,1,11,323,584,613,697,739,804,940,1052,1256,1273,1430,\\
	&&	1535,1816,1820,1871,1896,2030,2280,2347,2566,2598,2648,\\
	&&	2743,2781,3096,3352,3496,3624,3790,3831,3868,3887,3922,\\
	&&	3927,3953,3974,4115,4179,4293,4397,4445,4478,4561,4736,\\
	&&	4815,4885,4971,5074,5082,5098,5268,5280,5369,5426,5479,\\
	&&	5556,5679,5830,5858,6067,6135,6138,6184,6259,6303,6683,\\
	&&	6783,6822,6852,7024,7047,7195,7197,7255,7269,7289,7501,\\
	&&		7544,7562,7589,7682,7691,7697,7704,7822,7858,\\
	101& (0,7,60)&	0,1,40,354,640,885,888,1015,1031,1072,1120,1125,1217,\\
	&&	1273,1361,1461,1487,1569,1580,1634,1638,1683,1754,1993,\\
	&&	2069,2128,2223,2321,2656,2773,2837,2872,3052,3180,3383,\\
	&&	3458,3548,3830,3987,4019,4093,4385,4676,4688,4719,4942,\\
	&&	4957,4975,5449,5477,5647,5765,5874,5947,5970,6030,6142,\\
	&&	6194,6264,6349,6489,6621,6790,6800,6901,6923,7064,7315,\\
	&&	7317,7528,7657,7665,7686,7695,7720,7737,7799,7886,7970,\\
	&&	8148,8198,8225,8408,8474,8598,8634,8795,8931,9038,9052,\\
	&&	9099,9177,9190,9214,9258,9389,9408,9475,9856,9876,10194\\
	107& (0,2,29)&0,1,29,224,230,300,471,497,538,789,1049,1190,1193,1276,\\
	&&	1467,1509,1566,1709,1774,1919,2067,2598,2834,2859,3009,\\
	&&	3023,3028,3230,3334,3395,3450,3474,3571,3732,3856,3941,\\
	&&	4166,4292,4329,4369,4381,4449,4561,4595,4615,4713,5053,\\
	&&	5388,5395,5743,5747,6086,6276,6298,6345,6752,6788,6848,\\
	&&	6901,6922,7031,7033,7327,7602,7632,7696,7704,7739,7958,\\
	&&	8096,8211,8238,8249,8366,8688,8704,8779,8823,8872,8956,\\
	&&	9001,9019,9034,9051,9107,9173,9232,9346,9355,9436,9482,\\
	&&	9802,9850,9860,9873,9960,10247,10446,10549,10735,10827,\\
	&&	10866,10928,11033,11084,11115,11289\\
	\end{array}$

	\newpage
	$\begin{array}{rll}
	p&   gen.&	difference\:set \pmod{p^2-1}\: of\: p(p-1)\: elements\\
	\cline{1-3}
	131& (0,4,18)&	0,1,8,49,136,674,699,811,843,954,1044,1198,1217,1338,\\
	&&	1376,1615,1753,2201,2215,2225,2309,2321,2443,2635,2662,\\
	&&	2702,2704,2843,2936,3284,3782,4495,4881,4947,5006,5039,\\
	&&	5042,5304,5386,5433,5513,5623,5629,5794,6032,6133,6183,\\
	&&	6198,6353,6611,6648,6828,6954,7168,7365,7417,7437,7468,\\
	&&	7567,7621,8051,8160,8343,8389,8411,9030,9048,9242,9300,\\
	&&		9323,9339,9885,10100,10173,10330,10642,10924,10959,\\
	&&		11195,11266,11295,11380,11440,11526,11571,11628,11792,\\
	&&		12096,12159,12272,12488,12644,12688,12923,12934,13220,\\
	&&		13425,13446,13588,13649,13934,13938,14393,14511,14704,\\
	&&		14721,14819,14893,14971,15041,15118,15146,15295,15325,\\
	&&		15359,15414,15582,15744,15749,15931,16022,16294,16401,\\
	&&		16427,16480,16489,16802,16845,16858,16962,17085\\
	137& (0,2,112)&	0,1,61,213,288,306,353,531,568,652,686,755,900,1118,\\
	&&	1175,1185,1763,2101,2179,2322,2473,2489,2578,2763,2785,\\
	&&	2920,3102,3142,3155,3339,3468,3509,3538,3776,4101,4157,\\
	&&	4320,4403,4436,4479,4569,4575,4601,4737,4829,5239,5250,\\
	&&	5277,5486,5822,5881,5945,6056,6339,6430,6791,7095,7107,\\
	&&	7278,7366,7535,7636,7996,8116,8182,8226,8262,8491,8591,\\
	&&	9106,9164,9250,9295,9577,9703,10031,10034,10059,10138,\\
	&&	10187,10235,10524,10747,10801,11185,11302,11309,11326,\\
	&&	11570,11781,11790,11820,12163,12461,12512,12567,12586,\\
	&&	12649,12654,13083,13168,13374,13394,13489,13694,14017,\\
	&&	14432,14800,14894,15147,15256,15386,15534,15681,15683,\\
	&&	15923,16392,16695,16875,16890,16898,17071,17092,17106,\\
	&&	17205,17627,17804,17978,18264,18326,18378,18424,18428,\\
	&&	18505,18536,18578,18697
	\end{array}$


\setcounter{section}{4}
\setcounter{subsection}{1}
	\ssec{Ricatti geometry.}

	\setcounter{subsubsection}{-1}
	\sssec{Introduction.}
	It occured to me that just like in 3 dimensional
	Euclidean geometry, the geometry on the sphere can be used as a model
	for the non euclidean geometry of elliptic type, in the same way the
	geometry on the surface {\cal T} : $x^3 + y^3 + z^3 - 3 xyz = 1,$ can
	be used as a model for an other geometry, if $p \equiv -1 \pmod{6}.$
	I will call this geometry, {\em Ricatti geometry}.
	It turns out that this geometry is more akin to an Euclidean geometry.
	It can be considered as starting from a dual affine geometry in which
	we prefer a line (the ideal line) and a point (the ideal point) which
	is {\em not} on the line.  The line corresponds to the intersection
	with $t = 0$ of the plane $x + y + z = 0,$ the point to the direction
	of the line through the origin and the point $(1,1,1)$.

	\sssec{Definition.}
	Given $p \equiv -1 \pmod{6},$ the group $({\cal T} ,+)$ is cyclic
	(\ref{sec-tric3cyclic}. We determine a generator $(a,b,c)$ of the 
	groupt using in part \ref{sec-tric3gen}.  The {\em points}
	on ${\cal T}$  are labelled according to $i(a,b,c)$ from 0 to $p^2-2.$ 
	The {\em lines} are the set of points on ${\cal T}$ and a plane through
	the origin distinct from $[1,1,1,1]$ and which does not contain the line
	from the origin to $(1,1,1,1)$.  The line through $i$ and $i+1$ is
	labelled $-i^*.$	

	\sssec{Notation.}
	Points are denoted by a lower case letter or integer
	$\umod{p^2-1},$ lines by the same followed by a "$^*$". 

	\sssec{Definition.}
	If 2 points do not determine a line they are called {\em parallel}.\\
	If 2 lines are not incident to a point they are called {\em parallel}.

	\sssec{Definition.}
	There is a correspondence between the point $i$ and the line $i^*$,
	called {\em polarity}.

	\sssec{Theorem.}
	\enumb
	\item{\em  There are $p^2-1$ points and lines.}
	\item{\em  A line is incident to $p$ points and a point to $p$ lines.}
	\item{\em  A point is parallel to $p-2$ points and a line is parallel to
		$p-2$ lines.}
	\item{\em  There is duality in this geometry.}
	\enume

	\sssec{Theorem.}
	{\em If $i^*$ is incident to $i_0,i_1,i_2,i_3$ and $i_4,$ then
	$(i+j)^*$ is incident to $i_0-j,i_1-j,i_2-j,i_3-j$ and $i_4-j.$ }

	\sssec{Definition.}
	Let $D$ be a difference set associated to the integers in
	$Z_{p^2,.},$ between 0 and $p^2-1,$ relatively prime to $p,$
	$D = \{d_0,d_1, \ldots ,d_{p-1}\}.$
	\enumb
	\item  The {\em selector function} is defined as follows,\\
	\hth$	f(k(p+1)) = -1,$\\
	\hth$	f(d_j-d_i) = d_i.$
	\item  With points represented by elements in $Z_{p^2-1}$ and lines
		similarly represented but followed by "$^*,$" the incidence
		relation is defined by\\
	\hth$    i$ is on $j^*$  iff  $f(i+j) = 0.$
	\enume

	\sssec{Theorem.}
	\enumb
	\item{\em$i$ is parallel to $j$ or $i^*$ is parallel to $j^*$  iff
		  $f(i-j) = -1,$}
	\item{\em the line $(i \times j)^*$ incident to $i$ and $j$, not
		parallel, is $(f(i-j)-j)^*.$} 
	\item{\em the line $k^*$ incident to $i$ parallel to $j^*$ is  \ldots}
	\enume

	Proof:  For 2, we want $k$ to be $\equiv j \pmod{p+1}$ such that
		$f(k+i) = 0,$ \ldots.

	\sssec{Example.}
	$p = 5,$ $D = \{0,1,14,16,21\}$,\\
	$\begin{array}{ccrrrrrrrrrrrr}
	 i&\vline  &0& 1& 2& 3& 4& 5& 6& 7& 8& 9&10&11\\
	f(i)&\vline&-1& 0&14&21&21&16&-1&14&16&16&14&14\\[5pt]

	i&\vline&   12&13&14&15&16&17&18&19&20&21&22&23\\
	f(i)&\vline&-1& 1& 0& 1& 0&21&-1&21& 1& 0&16& 1
	\end{array}$

	Examples of such differences sets are given in \ref{sec-eric3gen}.

	\sssec{Example.}
	$p = 5,$ the computations are done $\umod{24}.$
	\enumb
	\item The coordinates of the $i$-th point on ${\cal T}$  are
	$a_i,b_i,c_i,$ the distance between $j$ and $j+i$ is
	$d_i = ^3\sqrt{a_i^3+b_i^3+(c_i-1)^3-3a_ib_i(c_i-1)}.$\\
	$\begin{array}{ccrrrrrrrrrrrr}
	i&\vline&0& 1& 2& 3& 4& 5& 6& 7& 8& 9&10&11\\
	\cline{1-14}
	a_i&\vline&0& 0& 1&-1&-1&-1& 2&-1& 1&-2&-1& 1\\
	b_i&\vline&0&-1&-1&-2&-2& 0& 2& 1& 0& 0& 1&-1\\
	c_i&\vline&1&-2&-1& 1&-1&-2& 0&-2& 0&-1&-1&-2\\
	d_i&\vline&0& 3& 1& 1& 0& 3& 3& 4& 0& 4& 1& 4\\[5pt]

	 i&\vline&12&13&14&15&16&17&18&19&20&21&22&23\\
	\cline{1-14}
	a_i&\vline&-1&-2& 0&-2& 0&-1&-1&-2&-2& 0& 2& 1\\
	b_i&\vline&-1&-1& 2&-1& 1&-2&-1& 1&-1&-2& 0&-2\\
	c_i&\vline&-2& 0& 2& 1& 0& 0& 1&-1&-1&-1& 2&-1\\
	d_i&\vline& 0& 2& 4& 1& 0& 1& 2& 2& 0& 4& 4& 2
	\end{array}$
	\item The line $i^*$ is incident to the points $-i,1-i,14-i,16-i,21-i.$
	\item  The point $i$ is parallel to $i+6,$ $i+12$ and $i-6.$
	\item  The angle between lines $j^*$ and $(j+i)^*$ is $d_i.$
	\enume

	In particular,\\
	$0^*:$	$0,1,14,16,21,$ is parallel to $6^*,12^*,18^*,$ \\
	$1^*:$	$23,0,13,15,20,$ is parallel to $7^*,13^*,19^*,$ \\
	$14^*:$	$10,11,0,2,7,$\\
	$16^*:$	$8,9,22,0,5,$\\
	$21^*:$	$3,4,17,19,0,$\\
	$23^*:$	$1,2,15,17,22.$\\
	on $0^*,$ the distances are\\
	$\begin{array}{rcrrrrr}
	   &\vline&0&1&14&16&21\\
	\cline{1-7}
	 0 &\vline&0&3&4&0&4\\
	 1 &\vline&2&0&1&1&0\\
	14 &\vline&1&4&0&1&4\\
	16 &\vline&0&4&4&0&3\\
	21 &\vline&1&0&1&2&0
	\end{array}$

	The "circles" of radius $r$ and center 0 are\\
	$\begin{array}{ll}
	r\hti{4}&	points\:on\\
	1&	2,3,10,13,15,17\\
	4&	7,9,11,14,21,22\\
	2&	18,19,23\\
	3&	1,5,6\\
	0&	0,4,8,12,16,20
	\end{array}$

	\sssec{Example }
	Let $p = 11,$\\
	$\begin{array}{ccrrrrrrrrrrrrrrrr}
	i&\vline& 0&1&2&3&4&5&6&7&8&9&10&11&12&13&14&15\\
	\cline{1-18}
	a_i&\vline& 0&1&6&0&9&4&5&1&8&7&4&3&8&5&4&4\\
	b_i&\vline&0&3&3&7&3&8&8&7&0&10&6&1&8&8&5&9\\
	c_i&\vline& 1&4&0&10&3&9&1&5&8&1&2&4&4&4&6&8\\
	d_i&\vline& 0&8&10&3&2&1&10&3&2&1&0&8&7&6&9&8\\[5pt]

	i&\vline&16&17&18&19&20&21&22&23&24&25&26&27&28&29&30&31\\
	\cline{1-18}
	a_i&\vline&7&9&4&1&8&6&3& 8&7&6&5&9&7&2&2&8\\
	b_i&\vline&9&4&3&10&7&5&6& 5&7&5&1&0&10&0&6&5\\
	c_i&\vline&9&0&9&7&8&8&0& 4&1&10&8&4&10&5&4&6\\
	d_i&\vline&8&4&6&9&0&10&10& 6&5&6&1&2&8&8&0&2\\[5pt]

	i&\vline&32&33&34&35&36&37&38&39&40&41&42&43&44&45&46&47\\
	\cline{1-18}
	a_i&\vline&9&7&5&5&4&1&4&2&1&4&0&10&3&9& 1&5\\
	b_i&\vline&2&0&4&6&4&5&9&5&0&1&6&0&9&4&5&1\\
	c_i&\vline&9&10&6&10&6&7&3&0&0&3&3&7&3&8&8&7\\
	d_i&\vline&9&3&9&6&7&10&7&2&0&3&5&8&2&8&1&10\\[5pt]
	\end{array}$

	$\begin{array}{ccrrrrrrrrrrrrrrrr}
	i&\vline&48&49&50&51&52&53&54&55&56&57&58&59&60&61&62&63\\
	\cline{1-18}
	a_i&\vline&8&1&2&4&4&4&6&8&9&0&9&7&8&8&0&4\\
	b_i&\vline&8&7&4&3&8&5&4&4&7&9&4&1&8&6&3&8\\
	c_i&\vline&0&10&6&1&8&8&5&9&9&4&3&10&7&5&6&5\\
	d_i&\vline&3&5&0&9&3&7&1&1&8&2&7&5&0&6&4&9\\[5pt]

	i&\vline&64&65&66&67&68&69&70&71&72&73&74&75&76&77&78&79\\
	\cline{1-18}
	a_i&\vline&1&10&8&4&10& 5&4&6&9&10&6&10&6&7&3&0\\
	b_i&\vline&7&6&5&9&7& 2&2&8&9&7&5&5&4&1&4&2\\
	c_i&\vline&7&5&1&0&10& 0&6&5&2&0&4&6&4&5&9&5\\
	d_i&\vline&3&10&10&4&8& 2&0&6&8&1&10&3&9&3&6&8\\[5pt]

	i&\vline&80&81&82&83&84&85&86&87&88&89&90&91&92&93&94&95\\
	\cline{1-18}
	a_i&\vline&0&3&3&7&3&8&8&7&0&10&6&1& 8&8&5&9\\
	b_i&\vline&1&4&0&10&3&9&1&5&8&1&2&4& 4&4&6&8\\
	c_i&\vline&0&1&6&0&9&4&5&1&8&7&4&3& 8&5&4&4\\
	d_i&\vline&0&9&4&1&4&5&2&8&2&9&0&3& 3&9&10&5\\[5pt]

	i&\vline&96&97&98&99&100&101&102&103&104&105&106&107&108&109&110&111\\
	\cline{1-18}
	a_i&\vline&9&4&3&10&7&5&6&5&7&5&1&0&10&0&6&5\\
	b_i&\vline&9&0&9&7&8&8&0&4&1&10&8&4&10&5&4&6\\
	c_i&\vline&7&9&4&1&8&6&3&8&7&6&5&9&7&2&2&8\\
	d_i&\vline&6&5&1&1&0&2&5&7&3&3&2&5&4&3&0&10\\[5pt]

	i&\vline&112&113&114& 115&116&117&118&119&120\\
	\cline{1-11}
	a_i&\vline&2&0&4& 6&4&5&9&5&0\\
	b_i&\vline&9&10&6&10&6&7&3&0&0\\
	c_i&\vline&9&7&5& 5&4&1&4&2&1\\
	d_i&\vline&9&8&1&10&9&8&1&3&0
	\end{array}$

	The "circles" of radius $r$ and center 0 are\\
	$\begin{array}{rlc}
	r&	points\:on\\
	1&	5,9,26,46,54,55,73,83,98,99,114,118&	(12)\\
	10&	2,6,21,22,37,47,65,66,74,94,111,115&	(12)\\
	2&	4,8,27,31,39,44,57,69,86,88,101,106&	(12)\\
	9&	14,19,32,34,51,63,76,81,89,93,112,116&	(12)\\
	3&	3,7,33,41,48,52,64,75,77,91,92,104,105,109,119&	(15)\\
	8&	1,11,15,16,28,29,43,45,56,68,72,79,87,113,117&	(15)\\
	4&	17,62,67,82,84,108	&	(6)\\
	7&	12,36,38,53,58,103	&	(6)\\
	5&	24,42,49,59,85,95,97,102,107&	(9)\\
	6&	13,18,23,25,35,61,71,78,96&	(9)\\
	0&	0,10,20,30,40,50,60,70,80,90,100,110&	(12)
	\end{array}$

	p = 17; circles of center 0 with given radius:\\
	$\begin{array}{rlc}
	0:&0,16,32,48,64,80,96,112,128,144,160,176,192,208,224,240,256,272&(18)\\
	1:&29,60,62,73,89,111,118,133,145,147,156,159,161,162,190,195,202,205,\\
	&	216,235,245,251,266,278					&(24)\\
	2:& 24,47,67,120,151,186,223,253,263,269,275,282		&(12)\\
	3:& 1,8,17,38,66,70,136,219,236,258,267,268			&(12)\\
	4:& 18,39,74,78,87,100,106,109,113,125,134,174,193,260,262 	&(15)\\
	5:& 23,59,92,103,110,124,139,142,165,180,213,234		&(12)\\
	6:& 3,33,36,41,42,44,50,51,55,56,71,77,88,101,121,138,157,172,273,274,
		277							&(21)\\
	7:& 2,4,7,9,34,61,63,68,91,107,119,153,173,198,207		&(15)\\
	8:&5,27,31,40,45,46,57,58,76,84,85,97,104,105,122,130,140,171,189,194,\\
	&	206,209,239,276						&(24)\\
	9:& 12,49,79,82,94,99,117,148,158,166,183,184,191,203,204,212,230,231,\\
	&	242,243,248,257,261,283					&(24)\\
	10:& 81,90,115,135,169,181,197,220,225,227,254,279,281,284,286	&(15)\\
	11:& 11,14,15,116,131,150,167,187,200,211,217,232,233,237,238,244,246,\\
	&	247,252,255,285						&(21)\\
	12:& 54,75,108,123,146,149,164,178,185,196,229,265		&(12)\\
	13:& 26,28,95,114,154,163,175,179,182,188,201,210,214,249,270	&(15)\\
	14:& 20,21,30,52,69,152,218,222,250,271,280,287			&(12)\\
	15:& 6,13,19,25,35,65,102,137,168,221,241,264			&(12)\\
	16:&10,22,37,43,53,72,83,86,93,98,126,127,129,132,141,143,155,170,177,\\
	&	199,215,226,228,259					&(24)
	\end{array}$

	The points common to the circles given above and noted $i:$ if the
	radius is $i$ are given below if there are points which are common with
	\newpage
\begin{verbatim}
	center 1, radius 1: 30,61,63,74,90,112,119,134,146,148,157,160,162,
		163,191,196,203,206,217,236,246,252,267,279
	    0: 112,160; 1: 162; 3: 236,267; 4: 74,134; 6: 157;
	    7: 61,63,119; 8: 206; 9: 148,191,203; 10: 90,279;
	   11: 217,246,252; 12: 146,196; 13: 163; 14: 30;
	center 1, radius 2: 25,48,68,121,152,187,224,254,264,270,276,283
	    0: 48,224; 6: 121; 7: 68; 8: 276; 9: 283; 10: 254; 11: 187;
	   13: 270; 14: 152; 15: 25,264;
	center 1, radius 3: 2,9,18,39,67,71,137,220,237,259,268,269
	    2: 67,269; 3: 268; 4: 18,39; 6: 71; 7: 2,9; 10: 220; 11: 237;
	   15: 137; 16: 259;
	center 1, radius 4: 19,40,75,79,88,101,107,110,114,126,135,175,194,261,
		263
	    2: 263; 5: 110; 6: 88,101; 7: 107; 8: 40,194; 9: 79,261; 10: 135;
	   12: 75; 13: 114,175; 15: 19; 16: 126;
	center 1, radius 5: 24,60,93,104,111,125,140,143,166,181,214,235
	    1: 60,111,235; 2: 24; 4: 125; 8: 104,140; 9: 166; 10: 181;
	   13: 214; 16: 93,143;
	center 1, radius 6: 4,34,37,42,43,45,51,52,56,57,72,78,89,102,122,139,
		158,173,274,275,278
	    1: 89,278; 2: 275; 4: 78; 5: 139; 6: 42,51,56,274; 7: 4,34,173;
	    8: 45,57,122; 9: 158; 14: 52; 15: 102; 16: 37,43,72;
	center 1, radius 7: 3,5,8,10,35,62,64,69,92,108,120,154,174,199,208
	    0: 64,208; 1: 62; 2: 120; 3: 8; 4: 174; 5: 92; 6: 3; 8: 5; 12: 108;
	   13: 154; 14: 69; 15: 35; 16: 10,199;
	center 1, radius 8: 6,28,32,41,46,47,58,59,77,85,86,98,105,106,123,131,
		141,172,190,195,207,210,240,277
	    0: 32,240; 1: 190,195; 2: 47; 4: 106; 5: 59; 6: 41,77,172,277;
	    7: 207; 8: 46,58,85,105; 11: 131; 12: 123; 13: 28,210; 15: 6;
	   16: 86,98,141;
\end{verbatim}
\begin{verbatim}
	center 1, radius 9: 13,50,80,83,95,100,118,149,159,167,184,185,192,204,
		205,213,231,232,243,244,249,258,262,284
	    0: 80,192; 1: 118,159,205; 3: 258; 4: 100,262; 5: 213; 6: 50;
	    9: 184,204,231,243; 10: 284; 11: 167,232,244; 12: 149,185;
	   13: 95,249; 15: 13; 16: 83;
	center 1, radius 10: 82,91,116,136,170,182,198,221,226,228,255,280,282,
		285,287
	    2: 282; 3: 136; 7: 91,198; 9: 82; 11: 116,255,285; 13: 182;
	   14: 280,287; 15: 221; 16: 170,226,228;
\end{verbatim}
	\newpage
\begin{verbatim}
	center 1, radius 11: 12,15,16,117,132,151,168,188,201,212,218,233,234,
		238,239,245,247,248,253,256,286
	    0: 16,256; 1: 245; 2: 151,253; 5: 234; 8: 239; 9: 12,117,212,248;
	   10: 286; 11: 15,233,238,247; 13: 188,201; 14: 218; 15: 168; 16: 132;
	center 1, radius 12: 55,76,109,124,147,150,165,179,186,197,230,266
	    1: 147,266; 2: 186; 4: 109; 5: 124,165; 6: 55; 8: 76; 9: 230;
	   10: 197; 11: 150;  13: 179;
	center 1, radius 13: 27,29,96,115,155,164,176,180,183,189,202,211,215,
		250,271
	    0: 96,176; 1: 29,202; 5: 180; 8: 27,189; 9: 183; 10: 115; 11: 211;
	   12: 164; 14: 250,271; 16: 155,215;
	center 1, radius 14: 21,22,31,53,70,153,219,223,251,272,281,0
	    0: 0,272; 1: 251; 2: 223; 3: 70,219; 7: 153; 8: 31; 10: 281;
	   14: 21; 16: 22,53;
	center 1, radius 15: 7,14,20,26,36,66,103,138,169,222,242,265
	    3: 66; 5: 103; 6: 36,138; 7: 7; 9: 242; 10: 169; 11: 14; 12: 265;
	   13: 26; 14: 20,222;
	center 1, radius 16: 11,23,38,44,54,73,84,87,94,99,127,128,130,133,142,
		144,156,171,178,200,216,227,229,260
	    0: 128,144; 1: 73,133,156,216; 3: 38; 4: 87,260; 5: 23,142; 6: 44;
	    8: 84,130,171; 9: 94,99; 10: 227; 11: 11,200; 12: 54,178,229;
	   16: 127;
\end{verbatim}
\begin{verbatim}
	center 2, radius 1: 31,62,64,75,91,113,120,135,147,149,158,161,163,164,
		192,197,204,207,218,237,247,253,268,280
	    0: 64,192; 1: 62,147,161; 2: 120,253; 3: 268; 4: 113; 7: 91,207;
	    8: 31; 9: 158,204; 10: 135,197; 11: 237,247; 12: 75,149,164;
	   13: 163; 14: 218,280;
	center 2, radius 2: 26,49,69,122,153,188,225,255,265,271,277,284
	    6: 277; 7: 153; 8: 122; 9: 49; 10: 225,284; 11: 255; 12: 265;
	   13: 26,188; 14: 69,271;
	center 2, radius 3: 3,10,19,40,68,72,138,221,238,260,269,270
	    2: 269; 4: 260; 6: 3,138; 7: 68; 8: 40; 11: 238; 13: 270;
	   15: 19,221; 16: 10,72;
	center 2, radius 4: 20,41,76,80,89,102,108,111,115,127,136,176,195,262,
		264
	    0: 80,176; 1: 89,111,195; 3: 136; 4: 262; 6: 41; 8: 76; 10: 115;
	   12: 108; 14: 20; 15: 102,264; 16: 127;
	center 2, radius 5: 25,61,94,105,112,126,141,144,167,182,215,236
	    0: 112,144; 3: 236; 7: 61; 8: 105; 9: 94; 11: 167; 13: 182; 15: 25;
	   16: 126,141,215;
\end{verbatim}
	\newpage
\begin{verbatim}
	center 2, radius 6: 5,35,38,43,44,46,52,53,57,58,73,79,90,103,123,140,
		159,174,275,276,279
	    1: 73,159; 2: 275; 3: 38; 4: 174; 5: 103; 6: 44;
	    8: 5,46,57,58,140,276; 9: 79; 10: 90,279; 12: 123; 14: 52; 15: 35;
	   16: 43,53;
	center 2, radius 7: 4,6,9,11,36,63,65,70,93,109,121,155,175,200,209
	    3: 70; 4: 109; 6: 36,121; 7: 4,9,63; 8: 209; 11: 11,200; 13: 175;
	   15: 6,65; 16: 93,155;
	center 2, radius 8: 7,29,33,42,47,48,59,60,78,86,87,99,106,107,124,132,
		142,173,191,196,208,211,241,278
	    0: 48,208; 1: 29,60,278; 2: 47; 4: 78,87,106; 5: 59,124,142;
	    6: 33,42; 7: 7,107,173; 9: 99,191; 11: 211; 12: 196; 15: 241;
	   16: 86,132;
\end{verbatim}
\begin{verbatim}
	center 2, radius 9: 14,51,81,84,96,101,119,150,160,168,185,186,193,205,
		206,214,232,233,244,245,250,259,263,285
	    0: 96,160; 1: 205,245; 2: 186,263; 4: 193; 6: 51,101; 7: 119;
	    8: 84,206; 10: 81; 11: 14,150,232,233,244,285; 12: 185; 13: 214;
	   14: 250; 15: 168; 16: 259;
	center 2, radius 10: 83,92,117,137,171,183,199,222,227,229,256,281,283,
		286,0
	    0: 0,256; 5: 92; 8: 171; 9: 117,183,283; 10: 227,281,286; 12: 229;
	   14: 222; 15: 137; 16: 83,199;
	center 2, radius 11: 13,16,17,118,133,152,169,189,202,213,219,234,235,
		239,240,246,248,249,254,257,287
	    0: 16,240; 1: 118,133,202,235; 3: 17,219; 5: 213,234; 8: 189,239;
	    9: 248,257; 10: 169,254; 11: 246; 13: 249; 14: 152,287; 15: 13;
	center 2, radius 12: 56,77,110,125,148,151,166,180,187,198,231,267
	    2: 151; 3: 267; 4: 125; 5: 110,180; 6: 56,77; 7: 198;
	    9: 148,166,231; 11: 187;
	center 2, radius 13: 28,30,97,116,156,165,177,181,184,190,203,212,216,
		251,272
	    0: 272; 1: 156,190,216,251; 5: 165; 8: 97; 9: 184,203,212; 10: 181;
	   11: 116; 13: 28; 14: 30; 16: 177;
	center 2, radius 14: 22,23,32,54,71,154,220,224,252,273,282,1
	    0: 32,224; 2: 282; 3: 1; 5: 23; 6: 71,273; 10: 220; 11: 252;
	   12: 54; 13: 154; 16: 22;
	center 2, radius 15: 8,15,21,27,37,67,104,139,170,223,243,266
	    1: 266; 2: 67,223; 3: 8; 5: 139; 8: 27,104; 9: 243; 11: 15; 14: 21;
	   16: 37,170;
	center 2, radius 16: 12,24,39,45,55,74,85,88,95,100,128,129,131,134,
		143,145,157,172,179,201,217,228,230,261
	    0: 128; 1: 145; 2: 24; 4: 39,74,100,134; 6: 55,88,157,172;
	    8: 45,85; 9: 12,230,261; 11: 131,217; 13: 95,179,201;
	   16: 129,143,228;
\end{verbatim}



\setcounter{section}{4}
\setcounter{subsection}{2}
	\ssec{3 - Dimensional Equidistance Curves.}

	\setcounter{subsubsection}{-1}
	\sssec{Introduction.}
	On the surface {\cal T} , for $p \equiv -1 \pmod{5}$, we can define
	besides lines (intersection with a plane through the origin), circles
	(pts equidistant using the cubic function from a given point),
	line-circle (set of tangents in space to the circles), podars (set of
	points where tangents in space intersect {\cal T} ), mediatrices (set of
	points equidistant from 2 points).  This section describes those
	curves.

	\sssec{Definition.}
	The {\em circles} are the set of points on {\cal T}  such that the
	cubic distance from a given point on {\cal T} , called the {\em center}
	of the circle, is a given integer $r,$ called the {\em radius} of the
	circle.

	\sssec{Theorem.}\label{sec-t3dimeq}
	{\em The circles of radius $r$ and center $(0,0,1),$ are the points
	$(x,y,z)$ which satisfy 0. and 1. or 0. and 2.}
	\enumb
	\item$	x^3 + y^3 + z^3 - 3xyz = 1.$
	\item$	x^3 + y^3 + z^3 - 3z^2 + 3z - 1 - 3xyz + 3xy = r^3.$
	\item$	3z^2 - 3z - 3xy = - r^3.$
	\enume

	\sssec{Definition.}
	A {\em line-circle} is the set of lines tangent in space to a circle.

	\sssec{Theorem.}
	{\em The line-circle associated to the circles in \ref{sec-t3dimeq}
	have at $(x,y,z)$ the direction $(\Delta x,\Delta y,\Delta z)$ given by}
	\enumb
	\item$\Delta x = xz_2+z'y_2,$ $\Delta y = -yz_2-z'x_2,$
		$\Delta z = -xx_2+yy_2,$\\
	{\em where}
	\item$	x_2 = x^2-yz,$ $y_2 = y^2-zx,$ $z_2 = z^2-xy,$ $z' = 2z-1.$
	\enume

	Proof:  the component of the direction satisfy\\
	\hth$	(x^2-yz)\Delta x + (y^2-zx)\Delta y + (z^2-xy)\Delta z = 0,$\\
	\hth$	    -y \Delta x -      x \Delta y +  (2z-1)\Delta z = 0,$\\
	hence 0.

	\sssec{Definition.}
	A {\em podar} is set of points where tangents in space intersect
	{\cal T}.

	\sssec{Theorem.}
	{\em The coordinates of points on the podar associated to the circle in
	\ref{sec-t3dimeq} are the points $(x+t\Delta x,y+t\Delta y,z+t\Delta z)$
	where $t$ satisfies}
	\enumb
	\item$	t = - 3\frac{x\Delta x_2+y\Delta y_2+z\Delta z_2}
	{(\Delta x+\Delta y+\Delta z}(\Delta x_2+\Delta y_2+\Delta z_2)).$\\
	{\em where}
	\item$\Delta x_2 = \Delta x^2-\Delta y\Delta z, 
	\Delta y_2 = \Delta y^2-\Delta z\Delta x, \Delta z_2 
	= \Delta z^2-\Delta x\Delta y.$
	\enume

	Proof.  A point $(x+t\Delta x,y+t\Delta y,z+t\Delta z)$ on the line
	$(x,y,z)$ with
	direction $(\Delta x,\Delta y,\Delta z)$ is on {\cal T}  if $t$
	satisfies the cubic equation,\\
	$(\Delta x+\Delta y+\Delta z)(\Delta x_2+\Delta y_2+\Delta z_2)t^3 +
	3(x\Delta x_2+y\Delta y_2+z\Delta z_2)t^2\\
	\hth +3(x_2\Delta x+y_2\Delta y+z_2\Delta z)t + (x^3+y^3+z^3-3xyz-1)+1
 = 0,$\\
	the coefficient of $t^0$ is 0 because $(x,y,z)$ is on {\cal T} , that
	of $t$ is 0 because it is $x_2(xz_2+z'y_2) + y_2(-yz_2-z'x_2) +
	z_2(-xx_2+yy_2) = 0.$

	\sssec{Theorem.}
	\enumb
	\item{\em If the tangent $k^*$ at $i$ to the circle, centered at 0 of
	radius $r,$ meets {\cal T}  at $j,$ then the tangent $(k^{+2i)*}$ at
	$-i$ to the circle, centered at 0 of radius $-r \pmod{p},$ meets
	{\cal T}  at $j-2i.$}
	\item{\em If the tangent at $i$ to the circle, centered at 0 of radius
		$r,$ meets {\cal T}  at $j$ parallel to $i,$ then the tangent
		at $-i$ to the circle, centered at 0 of radius $-r \pmod{p},$
		meets {\cal T}  at $j-2i$ parallel to $-i.$}
	\enume

	\sssec{Example.}
	Let $p = 11,$\\
	$\begin{array}{ccrrrrrrrrrrrrrrr}
	r = 1,\\
	circle &\vline&5& 9&26&46&54&55&73&83&98&99& 114& 118\\
	podar&\vline& 92&15& 113&43&51&52& 4&44& 107&45&81&97\\
	line\!-\!circle&\vline& 29^*&19^*&8^*&118^*&110^*&109^*&30^*&0^*&22^*&119^*&
	 40^*&32^*\\[5pt]
	r = 10,\\
	circle&\vline&2&6&21&22&37&47&65&66&74&94&111&115\\
	podar&\vline&101&93&87&31&118&98&62&63&71&61&117&82\\
	line\!-\!circle&\vline& 28^*&28^*&77^*&98^*&46^*&56^*&99^*&98^*&90^*&60^*&
	 37^*&39^*\\[5pt]
	r = 2\\
	circle&\vline&4&8&27&31&39&44&57&69&86&88&101&106\\
	podar&\vline&28&59&36&118& 6&68&36&66&95&49&98&85\\
	line\!-\!circle&\vline& \!--\!&95^*&93^*&3^*&115^*&\!--\!&93^*&95^*&34^*&115^*&63^*&
	44^*&\\[5pt]
	r = 9\\
	circle&\vline&14&19&32&34&51&63&76&81&89&93&112&116\\
	podar&\vline&113&16&113&43&48&42&100&48&56&102&43&20\\
	line\!-\!circle&\vline&16^*&25^*&51^*&86^*&113^*&87&\!--\!&73^*&65^*&27^*&111^*&
	\!--\\[5pt]
	r = 3\\
	circle&\vline&3&7&33&41&48&52&64&75&77&91104& 105& 109& 119\\
	podar&\vline&9&88&99&107&48&58&73&84& 8&97\\
	line\!-\!circle&\vline& 25^*&76^*&65^*&57^*&\!--\!&96^*&56^*&45^*&26^*
	&57^*\\
	r = 3\\
	circle&\vline&&92& 104& 105& 109& 119\\
	podar&\vline&&38&83&84&40&80\\
	line\!-\!circle&&\vline&6^*&46^*&45^*&114^*&84^*\\[5pt]
	r = 8\\
	circle&\vline&1&11&15&16&28&29&43&45&56&68\\
	podar&\vline&82&62&114&115&94&35&94&54&65&74\\
	line\!-\!circle&\vline& 82^*&92^*&15^*&14^*&70^*&119^*&60^*&75^*&64^*
		&80^*\\
	r = 8\\
	circle&\vline&&72&79&87&113&117\\
	podar&\vline&&72&25&33&74& 3\\
	line\!-\!circle&\vline&&\!--&19^*&11^*&90^*&31^*\\[5pt]
	r = 4\\
	circle&\vline&17&62&67&82&84&108\\
	podar&\vline&41&113&91&43&84&108\\
	line\!-\!circle&\vline& \!--\!&41^*&\!--\!& 1^*&\!--\!&\!--\\[5pt]
	r = 7\\
	circle&\vline&12&36&38&53&58&103\\
	podar&\vline&12&36&119&77&109& 7\\
	line\!-\!circle&\vline& \!--\!&\!--\!&45^*&\!--\!&45^*&\!--\\[5pt]
	r = 5\\
	circle&\vline&24&42&49&59&85&95&97&102&107\\
	podar&\vline&24&21&73&83&16&56&103&111&53\\
	line\!-\!circle&\vline& \!--\!&108^*&\!--\!&\!--\!&18^*&108^*&51^*&18^*&111^*&\\[5pt]
	r = 6\\
	circle&\vline&13&18&23&25&35&61&71&78&96\\
	podar&\vline& 79&27&29&106&86&85&95&57&96\\
	line\!-\!circle&\vline& 85^*&102^*&5^*&58^*&68^*&\!--\!&\!--\!&72^*&\!--
	\end{array}$

	\sssec{Definition.}
	A {\em mediatrix} of 2 points is the set of points equidistant
	from them.

	\sssec{Conjecture.}
	The number of points on the mediatrix is $\equiv 0 \pmod{4}$,
	unless the points are $i$ and $i+k(p-1)$ in which case it is  \ldots 
	namely these points are for 0 and $p-1,$ $(p-1)k$ and $(p+1)k-1.$
	When $p = 11,$ the multiples are 8, 12 and 16; when $p = 17,$ 12, 16, 20
	and 24; when $p = 23,$ 0, 16, 20, 24, 28, 32.

	\sssec{Example }
	of mediatrices for $p = 11:$\\
	For 0 and 3, 0 and 6, 0 and 9, no points.\\
	For 0 and 1: 16,22,29,55,66,92,99,105.\\
	For 0 and 2: 25,34,38,45,77,84,88,97.\\
	For 0 and 4: 1,3,6,7,8,9,15,31,52,72,93,109,115,116,117,118.\\
	For 0 and 5: 16,18,19,23,44,58,67,81,102,106,107,109.\\
	For 0 and 7: 2,5,25,48,49,78,79,102.\\
	For 0 and 8: 1,7,39,41,54,74,87,89.\\
	For 0 and 10: $10i$ and $12i-1$.\\
	For 0 and 11: 2,5,6,9,52,56,75,79.\\
	For 0 and 12: 25,28,35,39,63,64,68,69,93,97,104,107.\\
	For 0 and 13: 28,29,32,44,56,57,76,77,89,101,404,105.\\
	For 0 and 14: 3,11,15,29,43,91,105,119.\\
	For 0 and 15: 6,9,16,21,34,37,43,48,53,82,87,92,98,101,114,119.\\
	For 0 and 16: 22,37,45,64,72,91,99,114.\\
	For 0 and 17: 26,28,35,44,45,51,53,59.\\
	For 0 and 18: 3,4,7,11,14,15,29,32,42,57,65,73,81,96,106,109.\\
	For 0 and 19: 3,16,21,27,51,52,66,73,87,88,112,118.

	\sssec{Definition.}
	The {\em horizon} of a point $P$ on {\cal T}  is the set of points on
	{\cal T} and the tangent plane through $P.$

	\sssec{Theorem.}
	{\em The coordinates of points on the tangent at $P = (x,y,z)$ in
	the plane {\cal P}  through $O$ and $Q = (x'.y',z')$ which is also on
	{\cal T}  is\\
	\hth$(x+t\Delta x,y+t\Delta y,z+t\Delta z),$ where $t$ satisfies}
	\enumb
	\item$	t = - 3\frac{x\Delta x_2+y\Delta y_2+z\Delta z_2}
	{(\Delta x+\Delta y+\Delta z}(\Delta x_2+\Delta y_2+\Delta z_2)).$\\
	{\em where}
	\item$\Delta x_2 = \Delta x^2-\Delta y\Delta z,
	\Delta y_2 = \Delta y^2-\Delta z\Delta x,
	\Delta z_2 = \Delta z^2-\Delta x\Delta y.$
	\enume

	Proof.  The direction of the normal to {\cal P}  at $P$ is
	($a := yz'-zy',$ $b := zx'-xz',$ $c := xy'-yx'$).  The direction
	$(\Delta x,\Delta y,\Delta z)$ satisfies $a\Delta x+b\Delta y+c\Delta z
	= 0$ and $x_2\Delta x+y_2\Delta y+z_2\Delta z = 0,$ where\\
	\hth$	x_2 = x^2-yz,$ $y_2 = y^2-zx,$ $z_2 = z^2-xy,$\\
	therefore $\Delta x = y_2c-z_2b,$ $\Delta y = z_2a-x_2c,$
	$\Delta z = x_2b-y_2a.$\\
	A point $(x+t\Delta x,y+t\Delta y,z+t\Delta z)$ on the line $(x,y,z)$
	with direction $(\Delta x,\Delta y,\Delta z)$ is on {\cal T}  if $t$
	satisfies the cubic equation,\\
	$(\Delta x+\Delta y+\Delta z)(\Delta x_2+\Delta y_2+\Delta z_2)t^3 +
	3(x\Delta x_2+y\Delta y_2+z\Delta z_2)t^2 
		+3(x_2\Delta x+y_2\Delta y+z_2\Delta z)t + (x^3+y^3+z^3-3xyz-1)
	+1 = 0,$\\
	the coefficient of $t^0$ is 0 because $(x,y,z)$ is on {\cal T}, that of
	t is 0 because it is $x_2(y_2c-z_2b) + y_2(z_2a-x_2c) + z_2(x_2b-y_2a)
	 = 0.$

	\sssec{Algorithm.}
	To determine the horizon as 0 we determine for each point with $z = 1,$
	on which the line $0 \times sel(i)$ it is located, if $x+y = 0,$ the
	point is the ideal point, if $x = 0$ or $y = 0,$ 0 is a triple contact,
	if it is on no line $0 \times sel(i),$ then it corresponds to points
	parallel to it.  This is implement in [$\setminus 130\setminus$RIC.BAS]
	option 12.

	Proof: The horizon of {\cal P}  of 0 are the points on {\cal T} and
	$z = 1$ or $x^3 + y^3 - 3xy = 0,$ those in the plane $x = kt,$ $y = lt$ 
	satisfy $(k^3+l^3) t - 3kl = 0$ and $t = 0,$ twice.  If $k = 0$ or
	$l = 0$ then $t = 0$ is a triple root, if $k = -l,$ or $x+y = 0,$ then
	the point is an ideal point.
	In all other cases, $t = \frac{3kl}{k^3+l^3)}.$ 

	\sssec{Example.}
	For $p = 11,$ the points $H$ on the horizon of 0 have their tangent $t$
	and the points $Q$ on {\cal T} for which the tangent is $t^*$ given by\\
	$\begin{array}{rcrcl}
	   H   &\vline&    t^*  &\vline& 	     Q\\
	\cline{1-5}
	  0   &\vline&    1^*  &\vline&8,27,29,33,40,43,82,97,102,119  (y = 0)\\
	      &\vline&   41^*  &\vline&3,42,57,62,79,80,88,107,109,113 (x = 0)\\
	  6   &\vline&   28^*  &\vline&   2,6,13,16,55,70,75,92,93,101\\
	  9   &\vline&    0^*  &\vline&   1,9,28,30,34,41,44,83,98,103\\
	 24   &\vline&   --   &\vline&   12,24,36,48,60,72,84,96,108\\
	 51   &\vline&  103^*  &\vline&   17,18,26,45,47,51,58,61,100,115\\
	 66   &\vline&   98^*  &\vline&   5,22,23,31,50,52,56,63,66,105\\
	 81   &\vline&   83^*  &\vline&   15,20,37,38,46,65,67,71,78,81\\
	 87   &\vline&   34^*  &\vline&   7,10,49,64,69,86,87,95,114,116\\
	 99   &\vline&   30^*  &\vline&   4,11,14,53,68,73,90,91,99,118\\
	117   &\vline&   44^*  &\vline&   39,54,59,76,77,85,104,106,110,117\\
	  \infty    &\vline&    9^*  &\vline&   19,21,25,32,35,74,89,94,111,112
	\end{array}$

	\sssec{Conjecture.}
	{\em The points on the horizon of 0 are multiples of 3.}

	\sssec{Example.}
	The horizon of 0 is for
	\enumb
	\item$ p = 5,$\\
	$\begin{array}{ccrrrrrr}
	selector&\vline& 0& 1&14&16&21&--\\
	horizon&\vline&0&15& \infty&0& 3&18
	\end{array}$
	\item$ p = 11,$\\
	$\begin{array}{ccrrrrrrrrrrr}
	selector&\vline& 0& 1& 9&28&30&34&41&44&83&98 103&--\\
	horizon&\vline&9& 0& \infty&6&99&87& 0 117&81&66&51&24
	\end{array}$
	\item$ p = 17,$\\
	$\begin{array}{ccrrrrrrrrrrrrrrr}
	selector&\vline&0&1&10&13&34&45&59&86&112&114&129&134\\
	horizon&\vline&114&111&24&246&225&150&0&165&210&81&159&213\\
	selector&\vline&191&195&251&259&282&--\\
	horizon&\vline&\infty&93&0&141&120&108
	\end{array}$
	\item$ p = 23,$\\
	$\begin{array}{ccrrrrrrrrrrrrrrr}
	selector&\vline&0&1&60&91&134&142&148&203&249&253&266&269\\
	horizon&\vline&0&270&273&180&69&387&471&285&279&366&3&219\\
	selector&\vline&271&298&305&333&342&352&363&375&450&488&503&--\\
	horizon&\vline&81&231&426&444&33&0&498&402&453&\infty&294&408
	\end{array}$
	\enume




	\ssec{Generalization of the Selector Function.}

	\setcounter{subsubsection}{-1}
	\sssec{Introduction.}
	The selector function was introduced by Fernand
	Lemay to determine easily from the selector, points on 2 lines, lines
	incident to 2 points, points on lines or lines incident to points.
	This notion is generalized to 3 and more? dimensions.

	\sssec{Definition.}
	{\em defining polynomial}

	\sssec{Theorem.}
	{\em If the $P_i$ denotes a primitive polynomial of degree $i,$
	for $k = 3,$ the defining polynomials $P$ can have the following form,\\
	\hth$	P_4, P_1P_3, P_1^2P_2,$\\
	there are $p^4+p^3+p^2+p+1,$ $p^4-1,$ $p^4-p$ polynomials relatively
	prime to P, in these respective cases.\\
	For $k = 4,$ the defining polynomials $P$ can have the following form,\\
	\hth$	P_5, P_1P_4, P_1^2P_3, P_2P_3.$\\
	there are $p^5+p^4+p^3+p^2+p+1,$ $(p^3-1)(p+1),$ $p^5-1,$ $p^5-p$
	polynomials relatively prime to $P,$ in these respective cases.

	Proof: The polynomials in the sets are those which are relatively prime
	to the defining polynomial.  There are $p^k$ homogeneous polynomials of
	degree $k.$  If, for instance, $k = 4$ and the defining polynomial $P$
	is $P_2P_3,$ there are $p^2+p+1$ polynomials which are multiple of
	$P^2$ and $p+1,$ which are multiples of $P_3,$ hence
	$p^4+p^3+p^2+p+1-(p^2+p+1)-(p-1)$  polynomials relatively prime to $P.$}

	\sssec{Example.}
	$a_0,a_1,$ \ldots represents $I^{k+1}-a_0I^k-a_1I^{k-1} - \ldots.$\\
	$\begin{array}{rcrlrl}
	k&	p&period&def. pol.&sel.&roots of def. pol. or prim. pol.\\
	3&	3&	40&	2,1,1,1&    13&		--\\
	&&		26&	1,1,1,1&     9&		1\\
	&&		24&	0,1,1,1&     8&		2,2\\
	&	5&	156&	1,2,0,2&    31&		--\\
	&&		124&	1,0,0,2&    25&		4\\
	&&		120&	0,0,1,2&    24&		4,4\\

	4&	3&	121&	2,0,0,0,1&	40&	--\\
	&&		104&	0,1,0,0,1&	35&	(I^2+I-1)(I^3-I^2+I+1)\\
	&&		 80&	0,2,0,0,1&	27&	2\\
	&&		 78&	1,0,0,0,1&	26&	2,2\\
	&	5&	781&	4,0,0,0,1&	156&	--\\
	&&		744&	2,2,0,0,1&	149&(I^2+I+2)(I^3+2I^2-I+2)\\
	&&		624&	2,0,0,0,1&	125&	3\\
	&&		620&	3,0,1,0,1&	124&	3,3\\
	&	7&	2801&	3,0,0,0,1&	400&	--\\
	&&		2736&	6,0,0,0,1&	391& (I^2+2I-2)(I^3-I^2-3I-3)\\
	&&		2400&	3,1,0,0,1&	343&	3\\
	&&		2394&	0,3,3,0,1&	342&	5,5\\
	&	11\\
	&	13&	30941&	8,0,0,0,1&	2380&	--\\
	&&		30744&	5,0,0,0,1&	2365&(I^2-3I+6)(I^3-2I^2+I+2)\\
	&&		28560&	2,0,0,0,1&	2197&	11\\
	&&		28548\\

	5&	3&	364&	1,0,0,0,0,1&	121&	--\\
	&&		242&	1,1,0,0,0,1&	 81&	2\\
	&&		240&	1,2,1,0,0,1&	 80&	2,2\\
	&	5
	\end{array}$

	\sssec{Definition.}
	Given a selector $s,$ the {\em selector function} associates to
	the integers in the set ${\Bbb Z}_n$ a set of $p+1$ integers or $p$ integers
	obtained as follows,\\
	\hth$	s(j) \in f_i$  iff  $sel(l) - sel(j) = i$ for some $l.$

	\sssec{Theorem.}
	\enumb
	\item$	f(i)$ {\em is the set of points on the line} $i^* \times 0^*.$
	\item 0.$f(i)-j,$ {\em where we subtract $j$ to each element is the set,
		is the set of points in $(i+j)^* \times j^*,$ equivalently}
	\item 1.$f(i-j)-j,$ {\em is the set of points in}  $i^* \times j^*.$
	\item$	a^* \times b^* \times c^* = ((a-i)^* \times (b-i)^* \times
		(c-i)^*) - i.$
	\enume

	\sssec{Definition.}

	\sssec{Theorem.}
	{\em }

	\sssec{Theorem.}
	\enumb
	\item	{\em If the defining polynomial is primitive, then}
	\item 0.		$| s|  = \frac{p^k-1}{p-1},$
	\item 1.		{\em if} $i \neq  0,$ $| f(i)|  = p+1$
	\item	{\em If the defining polynomial has one root, then}
	\item 0.		$| s|  = p^k,$
	\item 1.		{\em if} $i \neq  0,$ $| f(i)|  = p,$
	\item	{\em If the defining polynomial has one root, then}
	\item 0.		$| s|  = p^k-1,$
	\item 1.		{\em if} $i ?,$ $| f(i)|  = p,$
	\item 2.		{\em if} $i ?,$ $| f(i)|  = p-1,$
	\item	{\em If the defining polynomial has one quadratic factor, then}
	\item 0.		$| s|  = (p^k-1)(p+1)?,$
	\item 1.		{\em if} $i ?,$ $| f(i)|  = p,$
	\enume

	\sssec{Example.}
	\enumb
	\item$  k = 3,$ $p = 3,$ defining polynomial
		$I^4 - 2 I^3 - I^2 - I - 1.$\\
	selector: 0  1  2  9  10  13  15  16  18  20  24  30  37\\
	selector function:\\
	$\begin{array}{rrrrrcrrrrrcrrrrr}
	 0&-1&-1&-1&-1&\vline& 14& 1& 2&10&16&\vline& 28& 2& 9&13&30\\
	 1& 0& 1& 9&15&\vline& 15& 0& 1& 9&15&\vline& 29& 1&13&20&24\\
	 2& 0&13&16&18&\vline& 16& 0& 2&24&37&\vline& 30& 0&10&20&30\\
	 3&10&13&15&37&\vline& 17& 1&13&20&24&\vline& 31& 9&10&18&24\\
	 4& 9&16&20&37&\vline& 18& 0& 2&24&37&\vline& 32& 9&10&18&24\\
	 5&10&13&15&37&\vline& 19& 1&18&30&37&\vline& 33& 9&16&20&37\\
	 6& 9&10&18&24&\vline& 20& 0&10&20&30&\vline& 34&15&16&24&30\\
	 7& 2& 9&13&30&\vline& 21& 9&16&20&37&\vline& 35& 2&15&18&20\\
	 8& 1& 2&10&16&\vline& 22& 2&15&18&20&\vline& 36& 1&13&20&24\\
	 9& 0& 1& 9&15&\vline& 23& 1&18&30&37&\vline& 37& 0&13&16&18\\
	10& 0&10&20&30&\vline& 24& 0&13&16&18&\vline& 38& 2&15&18&20\\
	11& 2& 9&13&30&\vline& 25&15&16&24&30&\vline& 39& 1& 2&10&16\\
	12& 1&18&30&37&\vline& 26&15&16&24&30&\vline&\\
	13& 0& 2&24&37&\vline& 27&10&13&15&37&\vline&\\
	\end{array}$
	\item$  k = 3, p = 3, defining polynomial I^4 - I^3 - I^2 - I - 1.$\\
	selector: 0  1  2  8  11  18  20  22  23\\
	selector function:\\
	$\begin{array}{rrrrcrrrrcrrrrcrrrr}
	0&-1&-1&-1 &\vline& 7& 1&11&20 &\vline&14& 8&20&23 &\vline&21& 1& 2&23\\
	1& 0& 1&22 &\vline& 8& 0&18&20 &\vline&15& 8&11&22 &\vline&22& 0& 1&22\\
	2& 0&18&20 &\vline& 9& 2&11&18 &\vline&16& 2&11&18 &\vline&23& 0&11&23\\
	3& 8&20&23 &\vline&10& 1& 8&18 &\vline&17& 1&11&20 &\vline&24& 2&20&22\\
	4&18&22&23 &\vline&11& 0&11&23 &\vline&18& 0& 2& 8 &\vline&25& 1& 2&23\\
	5&18&22&23 &\vline&12& 8&11&22 &\vline&19& 1& 8&18 &\vline\\
	6& 2&20&22 &\vline&13&-1&-1&-1 &\vline&20& 0& 2& 8 &\vline\\
	\end{array}$

	\item$  k = 3, p = 3, defining polynomial I^4 - I^2 - I - 1.$\\
	selector: 0  1  2  4  14  15  19  21\\
	selector function:\\
	$\begin{array}{rrrrcrrrrcrrrrcrrrr}
	 0&-1&-1&-1 &\vline& 6&15&19&-1 &\vline& 12& 2&14&-1
		&\vline& 18& 1&21&-1\\
	 1& 0& 1&14 &\vline&7&14&19&21 &\vline& 13& 1& 2&15
		&\vline& 19& 0& 2&19\\
	 2& 0& 2&19 &\vline&8&-1&-1&-1 &\vline& 14& 0& 1&14
		&\vline& 20& 1& 4&19\\
	 3& 1&21&-1 &\vline&9&15&19&-1 &\vline& 15& 0& 4&-1
		&\vline& 21& 0& 4&-1\\
	 4& 0&15&21 &\vline& 10& 4&14&15 &\vline& 16&-1&-1&-1
		&\vline& 22& 2& 4&21\\
	 5&14&19&21 &\vline& 11& 4&14&15 &\vline& 17& 2& 4&21
		&\vline& 23& 1& 2&15
	\end{array}$
	\enume

	\sssec{Example.}
	In the case of Example 3.6.x.0.
	if we denote by $i^\%,$ the lines $0^* \times i^*,$ these lines, which
	are sets of 4 points can all be obtained from\\
	$1^\% = \{0,1,9,15\},$ $2^\% = \{0,13,16,18\},$ $4^\% = \{9,16,20,37\}$ and
\\
	$10^\% = \{0,10,20,30\}$ by adding an integer modulo $n.$\\
	$1^\% + 0 = 1^\%,9^\%,15^\%,$ $1^\% + 1 = 39^\%,8^\%,14^\%,$
		$1^\% + 9 = 6^\%,31^\%,32^\%,$ $1^\% + 15 = 34^\%,25^\%,26^\%,$ \\
	$2^\% + 0 = 2^\%,24^\%,37^\%,$ $2^\% + 2 = 22^\%,35^\%,38^\%,$
		$2^\% + 37 = 3^\%,5^\%,27^\%,$  $2^\% + 24 = 13^\%,16^\%,18^\%,$\\
	$4^\% + 0 = 4^\%,21^\%,33^\%,$ $4^\% + 4 = 17^\%,29^\%,36^\%,$
		$4^\% + 21 = 12^\%,19^\%,23^\%,$ \\
	\hth$	4^\% + 33 = 7^\%,11^\%,28^\%,$\\
	\hth$10^\% + 0 = 10^\%,20^\%,30^\%.$

	\sssec{Definition.}

	\sssec{Conjecture.}




\setcounter{section}{4}
	\section{Generalization of the Spheres in Riccati Geometry.}
	\ssec{Dimension $k$.}
	\setcounter{subsubsection}{-1}
	\sssec{Introduction.}
	If we choose the ``sphere" $x^3+y^3+z^3-3xyz = 1$ in
	3 dimension we do not obtain for a given prime all periods as we do
	with the selector.  We have to generalyze using what is derived from
	differential equations with constant coefficients in which the
	coefficient of the $k-1$-th derivative is zero to obtain a constant
	Wronskian.  But, just as in the case of 2 dimensions, to obtain all
	sets of trigonometric functions, corresponding to the circular and
	hyperbolic functions, for all $p,$ we have to introduce in 3 dimension
	a cubic non residue if there is any,  $\ldots$ 
	I will first recall same well known definitions and Theorems of
	linear differential equations.

	\sssec{Definition.}\label{sec-wronsky}
	Given a linear differential equation
	\hth$	D^k x = C_0 x + C_1 Dx +  \ldots  C_{k-2} D^{k-2}x,$
	and $k$ solutions $y_i$ of these equations, the {\em Wronskian} is the
	matrix of functions whose $j$-th row are the $j$-th derivatives of
	$y_i,$  for $i = 0$ to $k-1.$

	\sssec{Theorem.}
	\enumb
	\item{\em  The functions $y_i$ are independent solution iff the
		determinant of the Wronskian is different from 0 for a
		particular value of the independent variable.}
	\item{\em  The determinant of the Wronskian is a constant function.}
	\item{\em  If $W(0) = E,$ then $W(x+y) = W(x) W(y).$}
	\enume

	\sssec{Theorem.}
	{\em If the linear differential equation \ref{sec-wronsky}.0. is such
	that $C_i$ are constant functions then any linear combination of $x$ 
	and its derivatives is also a solution of \ref{sec-wronsky}.0.}

	\sssec{Comment.}
	If we choose $x$ such that its derivatives are 0 except the $k-1$-th,
	chosen equal to 1, it is easy to obtain independent solutions using
	linear combination of $x$ and its derivatives to insure $W(0) = E.$
	If $det(W(t)) = 1,$ then $det(W(nt)) = 1$ and the surface
	$D^ix(nt),$ $i = 0$ to $k-1,$ can be chosen as a ``sphere" in
	$k$-dimension
	and $n$ as the angle between the directions joining the origin to the
	points $(D^ix((n+a)t))$ and $(D^ix(at)).$ 

	\sssec{Notation.}
	Given 2 solutions $x$ and $x'$ of \ref{sec-wronsky}.0. and a parameter
	$t$ let $x_i$ := $\{D^jx(it)\}$ and $x'_i$ := $\{D^jx'(it)\},$\\
	\hth$	y_{i,j} := x_i x'_j + x_j x'_i,$ $i \neq  j$\\
	\hth$	y_{i,i} := x_ix'_i,$

	\ssec{Dimension 3.}
	\sssec{Theorem.}
	{\em For $k = 3,$ let}
	\enumb
	\item[0.0.]$D^3 x_0 = C_0 x_0 + C_1 Dx_0,$\\
	{\em with}\\
	\item[0.1.]$Dx_0(0) = 0,$ $Dx_0 = x_1(0) = 0,$ $D^2x_0(0) = x^2(0) = 1,
	$\\
	{\em then}
	\item[1.0.] the set of functions $x_2 - C_1 x_0,$ $x_1,$ $x_0$ are
	independent
	\item[  1.] their Wronskian is\\
	\hth$	W = \matb{|}{ccc}
		x_2 - C_1 x_0&	  x_1&		x_0\\
		      C_0 x_0 &  x_2 &	x_1\\
		       C_0 x_1 &      C_0 x_0 + C_1 x_1& x_2\mate{|}$
	\item[  3.] $	W(0) = E.$
	\item[2.~~]	The distance from $(0,0,0)$ to $(x_0,x_1,x_2)$ is\\
	\hth$C_0^2 x_0^3 + C_0 x_1^3 + x_2^3 - 3 x_0x_1x_2 - C_1 x_1^2 x_2
	- C_1 x_2^2x_0 + 2 C_0C_1 x_0^2x_1\\
	\hti{12} + C_1^2 x_0x_1^2.$
	\item[3.~~] The addition formulas are\\
	\hth$	x''_0 = y_{0,2} + y_{1,1} - C_1 y_{0,0},$\\
	\hth$	x''_1 = y_{1,2} + C_0 y_{0,0},$\\
	\hth$	x''_2 = y_{2,2} + C_0 y_{0,1} + C_1 y_{1,1}.$
	\item[4.~~] The tangent plane at $(x_0,x_1,x_2)$ is\\
	\hth$	[ 3 C_0^2 x_0^2 - 3 C_0 x_1x_2 - C_1 x_2^2 + 4 C_0C_1 x_0x_1 
		+ C_1^2 x_1^2,\\
	\hti{9}	  3 C_0 x_1^2 - 3 C_0 x_2x_0 + 2 C_0C_1 x_0^2 - 2 C_1 x_1x_2
		+ 2 C_1^2 x_0x_1, \\
	\hti{9}	  3 x_2^2 - 3 C_0 x_0x_1 - 2 C_1 x_2x_0  - C_1 x_1^2 ].$
	\enume

	\sssec{Definition.}
	If\\[-10pt]
	\enumb
	\item$	p \equiv  1 \pmod{6},$ then $\nu^3 = n$ is a non cubic residue
		and the functions are not necessary real, we therefore denote
		then by $\xi_i$ instead of $x_i$ and express $\xi_i$ in terms
		of a power of $\nu$  and an integer $x_i$ as follows,
	\item$	x_0(3i) = x_0(3i),$ $\xi _1(3i) = x_1(3i) \nu ,$
	$\xi _2(3i) = x_2(3i) \nu ^2,$
		$\xi_0(3i+1) = x_0(3i) \nu^2,$ $\xi_1(3i+1) = x_1(3i+1),$
	$\xi_2(3i) = x_2(3i)\nu$,
		$\xi_0(3i+2) = x_0(3i) \nu$, $\xi_1(3i+2) = x_1(3i+2) \nu^2,$
	 $\xi_2(3i) = x_2(3i).$\\
		Moreover
		$C_1$ is replaced by $C_1 \nu^2.$

		The {\em addition formulas} become, for instance,\\
	\item $x_0(3i) = x_0(1)x_2(3i-1)+x_2(1)x_0(3i-1)n
		+ x_1(1)x_1(3i-1)$\\
	\hth$		 - C_1 x_0(1)x_0(3i-1) n,$\\
	\hth$	x_1(3i) = x_1(1)x_2(3i-1)+x_2(1)x_1(3i-1)
		+ C_0 x_0(1)x_0(3i-1),$\\
	\hth$	x_2(3i) = x_2(1)x_2(3i-1) n + C_0 (x_0(1)x_1(3i-1)
		+x_1(1)x_0(3i-1)n)$\\
	\hth$		 + C_1 x_1(1)x_1(3i-1) n.$
	\item $x_0(3i+1) = x_0(1)x_2(3i)+x_2(1)x_0(3i)n + x_1(1)x_1(3i)n$\\
	\hth$		 - C_1 x_0(1)x_0(3i) n,$\\
	\hth$	x_1(3i+1) = x_1(1)x_2(3i)+x_2(1)x_1(3i)n + C_0 x_0(1)x_0(3i),$\\
	\hth$	x_2(3i+1) = x_2(1)x_2(3i) + C_0 (x_0(1)x_1(3i)+x_1(1)x_0(3i))$\\
	\hth$		 + C_1 x_1(1)x_1(3i) n.$
	\item $x_0(3i+2) = x_0(1)x_2(3i+1)+x_2(1)x_0(3i+1)
		+ x_1(1)x_1(3i+1)$\\
	\hth$		 - C_1 x_0(1)x_0(3i+1),$\\
	\hth$	x_1(3i+2) = (x_1(1)x_2(3i+1)+x_2(1)x_1(3i+1)n)
		+ C_0 x_0(1)x_0(3i+1),$\\
	\hth$	x_2(3i+2) = x_2(1)x_2(3i+1) n + C_0 (x_0(1)x_1(3i+1)
		+x_1(1)x_0(3i+1))$\\
	\hth$		 + C_1 x_1(1)x_1(3i+1) n.$
	\enume

	\sssec{Theorem.}
	(on the period special case for type 1 and 2 and $p \equiv  1
	\pmod{6}$)\\
	{\em The period for type 0, 1 and 2 is respectively $p^2+p+1,$ $p^2-1,$
	$p^2-p.$ }

	\sssec{Notation.}
	The period in $k$-dimension, which depends on the type is
	denoted by $\pi$ $_k.$ 

	\sssec{Theorem.}
	{\em (on the selector)}

	\sssec{Example.}
	For $k = 3,$ (See $\setminus 130$ RIC.BAS)
	\enumb
	\item$  p = 5,$\\
	$\begin{array}{cllr}
		type	&period&	C_0,C_1	&x_0(1),x_1(1),x_2(1)\\
		   0	&    31	& 1,3	  	&0,3,1\\
		   1	&    24	& 1,0		&0,4,3\\
		   2	&    20	& 1,2		&0,2,1
	\end{array}$

	\item$  p = 7,$ $\nu  = 2,$\\
	$\begin{array}{cllr}
		type	&period	&C_0,C_1 &x_0(1),x_1(1),x_2(1)	\\
		   0	&    57	& 1,0		&0,5,6\\
		   1	&    48	& 1,1		&0,5,4\\
		   2	&    42	& 3,3		&0,3,0
	\end{array}$

	\item$  p = 11,$\\
	$\begin{array}{cllr}
		type	&period	&C_0,C_1 &x_0(1),x_1(1),x_2(1)	\\
		   0	&   133	& 1,3	  	&0,1,6\\
		   1	&   120	& 1,0		&0,2,5\\
		   2	&   110	& 1,5		&0,1,4
	\end{array}$

	\item$  p = 13,$ $\nu  = 2,$\\
	$\begin{array}{cllr}
		type	&period	&C_0,C_1 &x_0(1),x_1(1),x_2(1)	\\
		   0	&   183	& 1,0		&0,7,3\\
		   1	&   168	& 1,1		&1,1,10\\
		   2	&   156	& 4,1		&0,7,0
	\end{array}$

	\item$  p = 17,$\\
	$\begin{array}{cllr}
		type	&period	&C_0,C_1 &x_0(1),x_1(1),x_2(1)	\\
		   0	&   307	& 1,4	  	&0,1,2\\
		   1	&   288	& 1,0		&0,2,3\\
		   2	&   272	& 1,12		&0,2,10
	\end{array}$

	\item$  p = 19,$ $\nu  = 2,$\\
	$\begin{array}{cllr}
		type	&period	&C_0,C_1 &x_0(1),x_1(1),x_2(1)	\\
		   0	&   381	& 1,0		&0,3,7\\
		   1	&   360	& 1,1		&0,5,6\\
		   2	&   342	& 4,2		&0,2,11
	\end{array}$

	\item$  p = 23,$\\
	$\begin{array}{cllr}
		type	&period	&C_0,C_1 &x_0(1),x_1(1),x_2(1)	\\
		   0	&   553	& 1,3	  	&0,1,16\\
		   1	&   528	& 1,0		&0,2,9\\
		   2	&   506	& 1,1		&0,1,1
	\end{array}$

	\item$  p = 29,$\\
	$\begin{array}{cllr}
		type	&period	&C_0,C_1 &x_0(1),x_1(1),x_2(1)	\\
		   0	&   871	& 1,1	  	&0,1,1\\
		   1	&   840 & 1,0		&0,2,13\\
		   2	&   812	& 1,10		&0,3,7
	\end{array}$

	\item$  p = 31,$ $\nu  = 3,$\\
	$\begin{array}{cllr}
		type	&period	&C_0,C_1 &x_0(1),x_1(1),x_2(1)	\\
		   0	&   993	& 2,0		&0,3,24\\
		   1	&   960	& 1,2		&0,4,5\\
		   2	&   930 & 6,3		&0,3,29
	\end{array}$
	\enume

	\sssec{Example.}
	For $k = 3,$ (See [m130] WRONSKI.BAS)\\
	The table also includes the coordinates of a line,\\
	    e.g., for $p = 5,$ type 0, $3^* = [2,0,3].$
	\enumb
	\item$  p = 5,$ type 0, $C_0 = 1,$ $C_1 = 3,$\\
	$\begin{array}{ccrrrrrrrrrrrrrrrr}
	i&\vline&0&1&2&3&4&5&6&7&8&9&10&11&12&13&14&15\\
	\cline{1-18}
	x_0&\vline&	0&0&4&2&2&3&4&1&2&4&1&1&1&3&4&0\\
	x_1&\vline&	0&3&1&0&2&2&4&2&4&4&0&0&4&2&2&3\\
	x_2&\vline&	1&1&3&4&0&4&1&4&0&2&0&3&1&0&2&2\\
	l* &\vline& 11&14&28&29&19&20&2&23&24&8&0&10&27&15&4&5\\[5pt]
	i&\vline&16&17&18&19&20&21&22&23&24&25&26&27&28&29&30\\
	\cline{1-18}
	x_0&\vline&4&1&4&0&	2&0&3&1&0&2&2&4&2&4&4\\
	x_1&\vline&4&1&2&4&	1&1&1&3&4&0&4&1&4&0&2\\
	x_2&\vline&4&2&4&4&	0&0&4&2&2&3&4&1&2&4&1\\
	l* &\vline&6&18&30&94&1&21&13&12&22&3&16&17&7&25&26\\
	\end{array}$\\
	selector: $\{0,1,15,19,21,24\}$\\
	selector function:\\
	$\begin{array}{ccrrrrrrrrrrrrrrrr}
	& \vline& -1&0&19&21&15&19&15&24&24&15&21&21&19&19&1&0\\
	& \vline&15&15&1&0&1&0&24&1&0&21&24&19&24&21&1
	\end{array}$

	\item$p = 5,$ type 1, $C_0 = 1,$ $C_1 = 0,$\\
	$\begin{array}{ccrrrrrrrrrrrr}
	i&\vline&0&1&2&3&4&5&6&7&8&9&10&11\\
	\cline{1-14}
	x_0&\vline&	0&0&1&4&4&4&2&4&1&3&4&1\\
	x_1&\vline&	0&4&4&3&3&0&2&1&0&0&1&4\\
	x_2&\vline&	1&3&4&1&4&3&0&3&0&4&4&3\\[5pt]
	i&\vline&12&13&14&15&16&17&18&19& 20&21&22&23\\
	\cline{1-14}
	x_0&\vline&4&3&0&3&0&4&4&3&3&0&2&1\\
	x_1&\vline&4&4&2&4&1&3&4&1&4&3&0&3\\
	x_2&\vline&3&0&2&1&0&0&1&4&4&4&2&4
	\end{array}$\\
	selector: $\{0,1,14,16,21,--\}$\\
	selector function:\\
	$\begin{array}{ccrrrrrrrrrrrr}
	&\vline&-1&0&14&21&21&16& -1&14&16&16&14&14\\
	&\vline&-1&1&0&1&0&21&-1&21& 1&0&16&1
	\end{array}$

	\item$p = 5,$ type 2, $C_0 = 1,$ $C_1 = 2,$\\
	$\begin{array}{ccrrrrrrrrrrrrrrrrrrrr}
	i&\vline&0&1&2&3&4&5&6&7&8&9&10&11&12&13&14&15&16&17&18&19\\
	\cline{1-22}
	x_0&\vline&	0&0&4&2&1&2&3&0&2&4&2&3&3&1&1&1&4&3&0&4\\
	x_1&\vline&	0&2&4&2&3&3&1&1&1&4&3&0&4&0&0&4&2&1&2&3\\
	x_2&\vline&	1&1&4&3&0&4&0&0&4&2&1&2&3&0&2&4&2&3&3&1
	\end{array}$\\
	selector: $\{0,1,7,18,--,--\}$\\
	selector function:\\
	$\begin{array}{ccrrrrrrrrrrrrrrrrrrrr}
	& \vline& -1&0&18&18&-1&-1&1&0&-1&18&-1&7&-1&7&7&-1&-1&1&0&1\\
	\end{array}$

	\item$p = 7,$ type 0, $n = 5,$ $C_0 = 1,$ $C_1 = 0,$\\
	$\begin{array}{ccrrrrrrrrrrrrrrrrrrrr}
	i&\vline&0&1&2&3&4&5&6&7&8&9&10&11&12&13&14&15&16&17&18&19\\
	\cline{1-22}
	x_0&\vline&	0&0&4&3&5&2&1&0&2&3&1&4&3&3&1&4&6&4&3&0\\
	x_1&\vline&	0&5&6&5&0&5&3&6&4&4&1&6&1&5&6&3&2&6&2&2\\
	x_2&\vline&	1&6&5&5&3&3&0&5&3&0&1&0&2&6&6&2&4&3&1&0\\[5pt]
	i&\vline& 20&21&22&23&24&25&26&27&28&29&30&31&32&33&34&35&36&37&38&39\\
	\cline{1-22}
	x_0&\vline&	3&5&1&0&3&2&5&1&5&2&5&3&3&5&2&4&5&6&4&1\\
	x_1&\vline&	4&3&3&2&6&0&1&6&0&3&0&4&4&5&4&5&6&2&0&0\\
	x_2&\vline&	0&5&6&3&6&2&0&6&6&2&5&6&6&3&1&5&5&6&0&2\\[5pt]
	i&\vline& 40&41&42&43&44&45&46&47&48&49&50&51&52&53&54&55&56\\
	\cline{1-19}
	x_0&\vline&	2&6&4&4&0&2&4&0&6&0&1&1&1&1&3&4&4\\
	x_1&\vline&	3&4&6&5&6&0&5&4&4&3&4&5&6&3&3&3&4\\
	x_2&\vline&	3&2&0&6&4&1&2&3&6&3&6&2&3&4&5&3&5
	\end{array}$\\
	selector: $\{0,1,7,19,23,44,47,49\}$\\
	selector function:\\
	$\begin{array}{ccrrrrrrrrrrrrrrrrrrrr}
	&\vline& -1&0&47&44&19&44&1&0&49&49&47&47&7&44&44&49&7&47&1&0\\
	&\vline& 44&23&1&0&23&19&23&49&19&47&19&49&44&47&23&23&44&7&19&19\\
	&\vline& 7&23&7&1&0&19&1&0&1&0&7&7&49&23&47&49&1
	\end{array}$

	\item$p = 7,$ type 1, $n = 2,$ $C_0 = 1,$ $C_1 = 1,$\\
	$\begin{array}{ccrrrrrrrrrrrrrrrr}
	i&\vline& 0&1&2&3&4&5&6&7&8&9&10&11&12&13&14&15\\
	\cline{1-18}
	x_0&\vline&	0&0&4&5&4&4&1&3&1&6&6&0&3&1&0&6\\
	x_1&\vline&	0&5&3&2&6&5&3&2&1&0&5&2&4&2&4&2\\
	x_2&\vline&	1&4&5&5&2&5&4&2&2&1&6&2&1&3&0&5\\[5pt]
	i&\vline&16&17&18&19& 20&21&22&23&24&25&26&27&28&29&30&31\\
	\cline{1-18}
	x_0&\vline&5&1&3&6&	5&1&4&0&6&1&2&1&2&1&3&3\\
	x_1&\vline&6&6&1&6&	2&1&1&4&3&1&4&5&0&6&0&4\\
	x_2&\vline&0&1&1&1&	0&0&1&3&1&1&2&6&2&5&5&0\\[5pt]
	i&\vline&32&33&34&35&36&37&38&39&40&41&42&43&44&45&46&47\\
	\cline{1-18}
	x_0&\vline&4&3&1&4&4&2&5&4&	2&6&0&3&0&2&2&2\\
	x_1&\vline&4&4&0&0&4&5&4&4&	1&3&1&6&6&0&3&1\\
	x_2&\vline&6&2&0&5&3&2&6&5&	3&2&1&0&5&2&4&2
	\end{array}$\\
	selector: $\{0,1,11,14,23,42,44,--\}$\\
	selector function:\\
	$\begin{array}{ccrrrrrrrrrrrrrrrr}
	&\vline& -1&0&42&11&44&44&42&42&-1&14&1&0&11&1&0&44\\
	&\vline&-1&42&44&23& 42&23&1&0&-1&23&23&44&14&42&14&11\\
	&\vline&-1&11&14&14&23&11&11&23& -1&1&0&1&0&14&44&1
	\end{array}$

	\item$p = 7,$ type 2, $n = 5,$ $C_0 = 3,$ $C_1 = 3,$\\
	$\begin{array}{ccrrrrrrrrrrrrrrrr}
	i&\vline & 0&1&2&3&4&5&6&7&8&9&10&11&12&13&14&15\\
	\cline{1-18}
	x_0&\vline&	0&0&2&0&6&5&5&2&5&6&6&0&6&1&5&3\\
	x_1&\vline&	0&3&0&6&4&4&2&4&2&6&0&2&1&4&1&0\\
	x_2&\vline&	1&0&2&6&4&3&6&2&2&0&2&5&6&1&0&4\\[5pt]
	i&\vline&16&17&18&19&20&21&22&23&24&25&26&27&28&29&30&31\\
	\cline{1-18}
	x_0&\vline&0&1&4&3&	0&2&4&4&6&1&1&3&2&5&1&6\\
	x_1&\vline&5&6&3&0&	3&4&6&2&1&5&1&2&4&5&6&6\\
	x_2&\vline&6&1&0&3&	6&2&2&5&4&1&3&6&5&2&2&6\\[5pt]
	i& \vline&32&33&34&35&36&37&38&39&40&41\\
	\cline{1-12}
	x_0&\vline&4&4&6&4&1&6&5&4&	2&3\\
	x_1&\vline&6&6&6&5&6&4&6&2&	1&0\\
	x_2&\vline&2&2&5&2&6&6&3&5&	0&0
	\end{array}$\\
	selector: $\{ 0,1,3,11,16,20,--,--\}$\\
	selector function:\\
	$\begin{array}{ccrrrrrrrrrrrrrrrrr}
	 && \vline&-1&0&1&0&16&11&-1&-1&3&11&1&0& -1&3&-1&1\\
	 &&\vline&0&3&-1&1&-1&20&20&-1&20&16&16& -1&16&-1&11&11\\
	 &&\vline&20&11&-1&-1&16&20&3&3&1
	\end{array}$
	\enume

	\sssec{Definition.}
	\enumb
	\item The {\em direction} $dir(i,j)$ {\em of 2 points} $i$ and $j$ on
	the ``sphere" is the
	    direction of the line associated to the 2 points.
	\item A triangle $(i,j,k)$ is {\em isosceles}  iff  $j-i = k-j.$
	\item The {\em planar direction} $pl(i,j)$ {\em of 2 points} $i$ and $j$
	on the ``sphere" is
	    when $C_0 = 1$ and $C_1 = 0$ that of the normal to the plane passing
	    through the origin, $i$ and $j.$
	\item The {\em $t$ plane} $t^*(i)$ {\em at the point $i$} is the plane
	through the origin parallel to the tangent plane at $i.$
	\enume

	\sssec{Theorem.}
	\vspace{-18pt}\hspace{90pt}\footnote{10.2.88}
	\enumb
	\item[0.~~]$	dir(i,j) = dir(i+k,j+k) - k.$
	\item[1.~~]{\em if $c_0 = 1$ and $c_1 = 0$ then $pl(i,j) = pl(i+k,j+k)
		 + pk.$}
	\item[2.~~]{\em if $c_0 = 1$ and $c_1 = 0$ then $n(i) = -pi,$}
	\item[3.0]{\em  For types 0 and 1, the correspondance $i,$
		$n(i)$ is a bijection.}
	\item[  1.]{\em  For type 2, there are $2p-3$ values of $n$ which are
		undefined because the length of the normal is 0 (ideal point).}
	\item[4.~~]{\em For types 0 and 1, $t^*(i) + i = t^*(j) + j.$}
	\enume

	\sssec{Corollary.}
	{\em If a triangle $(i,j,k)$ is isosceles, then\\
	\hth$	dir(j,k) = dir(i,j) + k-i$}

	\sssec{Example.}
	$p = 5,$ type 0,
	$dir(0,1) = 21,$ $dir(1,2) = 22,$ $dir(0,2) = 5,$ $dir(0,3) = 25,$
	$dir(0,4) = 17.$\\
	$t^*(0) = 3^*,$ $t^*(1) = 2^*$.\\ 
	In the triangle $\{0,2,4\},$ $dir(2,4) = 5,$ $dir (4,0) = 17,$
	$dir(2,4) = 7 = 5+2.$\\
	$p = 5,$ type 1, pl(0,1) = 8, pl(1,2) = 3, pl(2,3) = 22, pl(0,2) = 6,
	$pl(1,2) = 1.$\\
	$t^*(0) = 8^*,$ $t^*(1) = 7^*.$\\
	$p = 11,$ type 1, $pl(0,1) = 40,$ $pl(1,2) = 29,$ $pl(2,3) = 18.$

	\ssec{Dimension 4.}

	\sssec{Theorem.}
	{\em For $k = 4,$ let}
	\enumb
	\item[0.0.]$D^4x_0 = C_0 x_0 + C_1 Dx_0 + C_2 D^2x_0,$\\
	{\em with
	\item[0.1.]$Dx_0(0) = 0,$ $Dx_0 = x_1(0) = 0,$
		$D^2x_0(0) = x_2(0) = 0,$
		$D^3x_0(0) = x_3(0) = 1,$\\
	then}
	\item[1.0.]{\em  the functions $x_3 - C_2x_1 - C_1 x_0,$
		$x_2 - C_2 x_0,$ $x_1$ and $x_0$ are independent}\\
	\item[  1.]{\em  their Wronskian is}\\
	\hth$	W = \matb{|}{cccc}
		x_3 - C_2x_1 - C_1x_0&    x_2 - C_2x_0&     x_1	&x_0\\
		     C_0x_0 &   x_3 - C_2x_1 &    x_2 &	x_1\\
		     C_0x_1 & C_0x_0 + C_1x_1&     x_3& 	x_2\\
		C_0x_2 & C_0x_1 + C_1x_2&C_0x_0 + C_1x_1 + C_2x_2&x_3\mate{|}
	$\\
	\item[  2.] $W(0) = E.$
	\item[2.~~]{\em The addition formulas are}\\
	\hth$	x''_0 = y_{0,3} + y_{1,2} - C_1 y_{0,0} - C_2 y_{0,1}$\\
	\hth$	x''_1 = y_{1,3} + y_{2,2} - C_0 y_{0,0} - C_2 y_{1,1}$\\
	\hth$	x''_2 = y_{2,3} + C_0 y_{0,1} + C_1 y_{1,1}$\\
	\hth$x''_3 = y_{3,3} + C_0(y_{0,2} + y_{1,1}) + C_1 y_{1,2} 
		+ C_2 y_{2,2}$
	\enume

	\sssec{Example.}


 
\chapter{FINITE ELLIPTIC  FUNCTIONS}


\setcounter{section}{-1}
\section{Introduction.}

	The success of the study of the harmonic polygons of Casey (II.6.1),
	suggested the study of the polygons of Poncelet.  After having
	conjectured that the Theorem of Poncelet, as given in I.2.2.
	generalized to the finite case, and because one of the proof of this
	Theorem, in the classical case, is by means of elliptic functions,
	this suggested that these too could be generalized to the finite case.
	Just as the additions properties were used to define the trigonometric
	functions, the same properties were generalized to the finite case.  It
	was soon realized that the poles of the elliptic functions correspond
	to values, which in the finite case are outside of the finite field.
	The basic definitions and properties of section 1 do not give directly
	functions but an abelian group structure on a set $E,$ whose elements
	are, in general, triplets of integers modulo $p.$  In section 2, this
	structure will be described as the direct product of the Klein 4-group
	and an abelian group which can be used as seen in section 3 to define 3
	functions which generalize, in the finite case, the functions $sn,$ $cn$
	and $dn$ of Jacobi.

	In this Chapter, $j$ and $j'$ will denote $+1$ or $-1.$

\section{The Jacobi functions.}

	\ssec{Definitions and basic properties of the Jacobian elliptic
	group.}

	\setcounter{subsubsection}{-1}
	\sssec{Introduction.}
	Given $p$ and $m$ different from 0 and 1, we will define in
	3.1.1, the set $E = E(p,m)$ and, in 3.1.7., an operation ``+'' from
	$E \bigx E$
	into $E.$  The basic result that $(E,+)$ is an abelian group is given in
	3.1.15.

	\sssec{Definition.}
	Given $s,c,d \in {\Bbb Z}_p.$ 
	The {\em elements of} $E$ are\\
	\hth$	(s,c,d)$\\
	such that\\
	D0.\hti{5}$   s^2 + c^2 = 1$ and $d^2 + m\: s^2 = 1.$\\
	as well as, when $-1$ and $-m$ are quadratic residues,\\
	\hth$(\infty ,c\:\infty,d\:\infty ),$ where $c^2 = -1$ and $d^2 = -m.$

	\sssec{Notation.}
	\hth$i := \sqrt{-1},$ $m_1 := 1-m,$ $k := \sqrt{m},$ $k_1 := \sqrt{m_1}.$

	\sssec{Theorem.}
	H0.\hti{5}$   (s,c,d),$ $(s_1,c_1,d_1),$ $(s_2,c_2,d_2) \in E,$\\
	H1.\hti{5}$   j = +1$ or $-1,$\\
	{\em then}\\
	C0.\hti{5}$   d^2 - m\: c^2 = m_1.$\\
	C1.\hti{5}$   c^2 + m_1 s^2 = d^2.$\\
	C2.\hti{5}$   c^2 + s^2 d^2 = d^2 + m\: s^2 c^2 = 1 - m\: s^4.$\\
	C3.\hti{5}$   m(1 - c)(1 + c) = (1 - d)(1 + d).$\\
	C4.\hti{5}$   (c + d)(1 + j d) = (1 + j\: c)(j\: m_1 + d + m\:c).$\\
	C5.\hti{5}$   m(c + d)(1 - j\: c) = (1 - j\: d)(j\: m_1 + d + m\:c)$\\
	C6.\hti{5}$   d^2 - m\: s^2 c^2 = d^2 + c^2(d^2 - 1).$\\
	C7.\hti{5}$   d_1^2 d_2^2 + m\:m_1 s_1^2 s_2^2 = m_1 - m\: c_1^2 c_2^2.$\\
	C8.\hti{5}$(d_1 s_1 c_2 + d_2 s_2 c_1)(d_1 d_2 - m\: s_1 s_2 c_1 c_2)
	 = (s_1 c_2 d_2 + s_2 c_1 d_1)(d_1^2 d_2^2 + m m_1 s_1^2 s_2^2).$

	Proof:  Each of the identities can easily be verified using
	Definition 1.1.  If C3, is written\\
	C3'\hti{6}$   m(1 - j\: c)(1 + j\: c) = (1 - j\: d)(1 + j\: d),$\\
	then C5, follows from C4.

	\sssec{Lemma.}
	H0.\hti{5}$m\: s_0^2 s_1^2 = 1,$\\
	{\em then}\\
	C0.\hti{5}$d_0^2 = -m\: s_0^2 c_1^2,$ $d_1^2 = -m\: s_1^2 c_0^2,$\\
	C1.\hti{5}$(s_0 c_1 d_1)^2 = (s_1 c_0 d_0)^2,$\\
	C2.\hti{5}$(c_0 c_1)^2 = (d_0 s_0 d_1 s_1)^2,$\\
	C3.\hti{5}$(d_0 d_1)^2 = (m\: s_0 c_0 s_1 c_1)^2,$\\
	C4.\hti{5}$s_0 s_1 \neq  0.$

	\sssec{Lemma.}
	H0.\hti{5}$m\: s_0^2 s_1^2 = 1,$\\
	H1.\hti{5}$s_0 c_1 d_1 = -j\: s_1 c_0 d_0,$\\
	{\em then}\\
	C0.\hti{5}$c_0 c_1 = j\: d_0 s_0 d_1 s_1,$\\
	C1.\hti{5}$d_0 d_1 = j\: m\: s_0 c_0 s_1 c_1.$\\
	C2.\hti{5}$c_0 = 0  \Rightarrow   d_1 = 0$ and $c_1 \neq  0.$\\
	\hti{9}$c_1 = 0  \Rightarrow   d_0 = 0$ and $c_0 \neq  0.$

	Proof:  If $c_0 = c_1 = 0$ then $s_0^2 = s_1^2 = 1$ hence $m = 1,$
	which is excluded.

	\sssec{Lemma.}
	H0.\hti{5}$   (s_0 c_1 d_1)^2 = (s_1 c_0 d_0)^2,$\\
	{\em then}\\
	C0.\hti{5}{\em$s_0 = j\: s_1$ or $m\: s_0^2 s_1^2 = 1.$}
\:
	\sssec{Definition.}
	The {\em addition} is defined as follows:\\
	Let     $D = 1 - m\: s_0^2 s_1^2.$ \\
	\hth    If $D \neq  0,$ then\\
	D0.\hti{5}$   (s_0,c_0,d_0) + (s_1,c_1,d_1) =
	(\frac{s_0c_1d_1 + s_1c_0d_0}{D},\frac{c_0c_1 - d_0s_0d_1s_1}{D},
		\frac{d_0d_1 - m\: s_0c_0s_1c_1}{D} ),$\\
	\hth If $D = 0,$ $s_0c_1d_1 = s_1c_0d_0,$ $c_0 \neq  0$ and
	$c_1 \neq  0,$ then\\
	D1.\hti{5}$(s_0,c_0,d_0) + (s_1,c_1,d_1) = (\infty,c\:\infty,d\:\infty),
$\\
	\hth where $c = \frac{c_1}{s_1 d_0}$ and $d = \frac{d_1}{s_1 c_0},$\\
	\hth If $D = 0,$ $s_0c_1d_1 = -s_1c_0d_0, c_0 \neq  0$ and
	$c_1 \neq  0,$ then\\
	D2.0.\hti{3}$ (s_0,c_0,d_0) + (s_1,c_1,d_1)
	= (\frac{s_0^2 - s_1^2}{2 s_0 c_1 d_1},
	 \frac{c_0^2 + c_1^2}{2 c_0 c_1}, \frac{d_0^2 + d_1^2}{2 d_0 d_1}),$\\
	\hth If $D = 0,$ $s_0 c_1 d_1 = j\: s_1 c_0 d_0,$ $c_0 = 0$ and
	$c_1 \neq  0,$ then\\
	D2.1.\hti{3}$ (s_0,c_0,d_0) + (s_1,c_1,d_1)
		= (\infty,c\:\infty,d\:\infty ),$\\
	\hth where $c = \frac{-d_0 s_1}{c_1}$ and
	$d = \frac{d_0^3}{m\: s_0 c_1^3}.$\\
	\hth If $D = 0,$ $s_0c_1d_1 = j\: s_1c_0d_0,$ $c_0 \neq  0$ and
	$c_1 = 0,$ then\\
	D2.2.\hti{3}$ (s_0,c_0,d_0) + (s_1,c_1,d_1)
		= (\infty ,c\:\infty ,d\:\infty ),$\\
	\hth where $c = \frac{-d_1 s_0}{c_0}$ and
		$d = \frac{d_1^3}{m\: s_1 c_0^3}.$\\
	\hth If $s_0 \neq  0,$ then\\
	D3.0.\hti{3}$ (\infty ,c\:\infty ,d\:\infty ) + (s_0,c_0,d_0)$\\
	\hth$	= (s_0,c_0,d_0) + (\infty ,c\:\infty ,d\:\infty )
		 = (\frac{-c d}{m\: s_0}, \frac{d d_0}{m\: s_0},
		 \frac{c c_0}{m\: s_0}).$\\
	\hth If $s_0 = 0,$ then\\
	D3.1.\hti{3}$ (\infty ,c \:\infty ,d \:\infty ) + (0,c_0,d_0)$\\
	\hth$= (0,c_0,d_0) + (\infty ,c \:\infty ,d \:\infty )
	 = (\infty ,c \:d_0 \:\infty ,d\: c_0 \:\infty ).$\\
	D4.\hti{5}$(\infty,c_0 \:\infty,d_0 \:\infty) + (\infty,c_1 \:\infty ,
	d_1 \:\infty ) = (0,\frac{d_0 d_1}{m},c_0 c_1).$

	\sssec{Example.}
	With $p = 11,$ $m = 3,$ $(-\frac{1}{11}) = (-\frac{3}{11}) = -1,$\\
	\hth$E = \{(0,1,1),(0,1,-1),(0,-1,1),(0,-1,-1),\\
	\hti{12}	(1,0,3),(1,0,-3),(-1,0,3),(-1,0,-3),\\
	\hti{12}	(5,3,5),(5,3,-5),(5,-3,5),(5,-3,-5),\\
	\hti{12}	(-5,3,5),(-5,3,-5),(-5,-3,5),(-5,-3,-5)\}.$\\
	If the elements of $E$ in the above order are abbreviated 0,1,2, \ldots
	,15, then the addition table is\\
	$\begin{array}{ccrrrrrrrrrrrrrrrr}
	 + &\vline& 0&1&2&3&4&5&6&7&8&9&10&11&12&13&14&15\\
	\cline{1-18}
	 0 &\vline& 0&1&2&3&4&5&6&7&8&9&10&11&12&13&14&15\\
	 1 &\vline& 1&0&3&2&7&6&5&4&13&12&15&14&9&8&11&10\\
	 2 &\vline& 2&3&0&1&6&7&4&5&14&15&12&13&10&11&8&9\\
	 3 &\vline& 3&2&1&0&5&4&7&6&11&10&9&8&15&14&13&12\\
	 4 &\vline& 4&7&6&5&2&1&0&3&10&13&14&9&8&15&12&11\\\
	 5 &\vline& 5&6&7&4&1&2&3&0&9&14&13&10&11&12&15&8\\
	 6 &\vline& 6&5&4&7&0&3&2&1&12&11&8&15&14&9&10&13\\
	 7 &\vline& 7&4&5&6&3&0&1&2&15&8&11&12&13&10&9&14\\
	 8 &\vline& 8&13&14&11&10&9&12&15&4&1&2&5&0&7&6&3\\
	 9 &\vline& 9&12&15&10&13&14&11&8&1&6&7&2&5&0&3&4\\
	10 &\vline&10&15&12&9&14&13&8&11&2&7&6&1&4&3&0&5\\
	11 &\vline&11&14&13&8&9&10&15&12&5&2&1&4&3&6&7&0\\
	12 &\vline&12&9&10&15&8&11&14&13&0&5&4&3&6&1&2&7\\
	13 &\vline&13&8&11&14&15&12&9&10&7&0&3&6&1&4&5&2\\
	14 &\vline&14&11&8&13&12&15&10&9&6&3&0&7&2&5&4&1\\
	15 &\vline&15&10&9&12&11&8&13&14&3&4&5&0&7&2&1&6
	\end{array}$

	\sssec{Example.}
	With $p = 13,$ $m = 3,$ $(-\frac{1}{13}) = (-\frac{3}{13}) = 1,$\\
	\hth$E = \{(0,1,1),(0,1,-1),(0,-1,1),(0,-1,-1),\\
	\hti{12}(\infty,5 \:\infty,6 \:\infty),(\infty,5 \:\infty,-6 \:\infty),
	(\infty,-5 \:\infty,6\:\infty),(\infty,-5\:\infty,-6\:\infty),\\
	\hti{12}(6,2,6),(6,2,-6),(6,-2,6),(6,-2,-6),\\
	\hti{12}(-6,2,6),(-6,2,-6),(-6,-2,6),(-6,-2,-6)\}.$\\
	If the elements of $E$ in the above order are abbreviated 0,1,2, \ldots,
	15, then the addition table is\\
	$\begin{array}{ccrrrrrrrrrrrrrrrr}
	 + &\vline& 0&1&2&3&4&5&6&7&8&9&10&11&12&13&14&15\\
	\cline{1-18}
	 0 &\vline& 0&1&2&3&4&5&6&7&8&9&10&11&12&13&14&15\\
	 1 &\vline& 1&0&3&2&6&7&4&5&13&12&15&14&9&8&11&10\\
	 2 &\vline& 2&3&0&1&5&4&7&6&14&15&12&13&10&11&8&9\\
	 3 &\vline& 3&2&1&0&7&6&5&4&11&10&9&8&15&14&13&12\\
	 4 &\vline& 4&6&5&7&3&1&2&0&12&14&13&15&11&9&10&8\\
	 5 &\vline& 5&7&4&6&1&3&0&2&10&8&11&9&13&15&12&14\\
	 6 &\vline& 6&4&7&5&2&0&3&1&9&11&8&10&14&12&15&13\\
	 7 &\vline& 7&5&6&4&0&2&1&3&15&13&14&12&8&10&9&11\\
	 8 &\vline& 8&13&14&11&12&10&9&15&7&1&2&4&0&5&6&3\\
	 9 &\vline& 9&12&15&10&14&8&11&13&1&4&7&2&6&0&3&5\\
	10 &\vline& 10&15&12&9&13&11&8&14&2&7&4&1&5&3&0&6\\
	11 &\vline& 11&14&13&8&15&9&10&12&4&2&1&7&3&6&5&0\\
	12 &\vline& 12&9&10&15&11&13&14&8&0&6&5&3&4&1&2&7\\
	13 &\vline& 13&8&11&14&9&15&12&10&5&0&3&6&1&7&4&2\\
	14 &\vline& 14&11&8&13&10&12&15&9&6&3&0&5&2&4&7&1\\
	15 &\vline& 15&10&9&12&8&14&13&11&3&5&6&0&7&2&1&4
	\end{array}$

	\sssec{Theorem.}
	C0.\hti{5}$(s_0,c_0,d_0) + (j's_0,jc_0,j\:j'd_0) = (0,j,j\:j').$\\
	C1.\hti{5}$(\infty ,c \:\infty ,d \:\infty ) + (\infty ,j\:c \:\infty ,
		j'\:d\: \infty ) = (0,-j',-j).$\\
	C2.\hti{5}$-(s_0,c_0,d_0) = (-s_0,c_0,d_0).$\\
	C3.\hti{5}$-(\infty ,c \:\infty ,d \:\infty )
		= (\infty ,-c \:\infty ,-d \:\infty ).$\\
	C4.\hti{5}$(s_0,c_0,d_0) + (0,j,j\:j') = (js_0,j'c_0,j\:j'd_0).$\\
	C5.\hti{5}$(\infty ,c \:\infty ,d \:\infty ) + (0,j',j)
		= (\infty ,jc \:\infty ,j'd \:\infty ).$

	\sssec{Theorem.}
	H0.\hti{5}$   (s_0,c_0,d_0) + (s_1,c_1,d_1) = (s_2,c_2,d_2),$\\
	{\em then}\\
	C0.\hti{5}$   (s_0,c_0,d_0) + (-s_1,-c_1,d_1) = (-s_2,-c_2,d_2).$\\
	C1.\hti{5}$   (s_0,c_0,d_0) + (s_1,-c_1,-d_1) = (s_2,-c_2,-d_2).$\\
	C2.\hti{5}$   (s_0,c_0,d_0) + (-s_1,c_1,-d_1) = (-s_2,c_2,-d_2).$

	\sssec{Notation.}
	We will use the notation, which is customary in abelian
	groups with addition as operation symbol,\\
	\hth$	n(s_0,c_0,d_0) = (n-1)(s_0,c_0,d_0) + (s_0,c_0,d_0), n \in Z$\\
	using induction starting with $n = 0$ or $n = -1.$

	\sssec{Theorem.}
	C0.\hti{5}{\em $n(-s_0,-c_0,-d_0) = j\: n(s_0,c_0,d_0),$
		with $j = (-1)^n.$

	\sssec{Theorem.}
	D0.\hti{5}$   D_2 := 1 - m\: s^4,$\\
	{\em then}\\
	C0.\hti{5}$   2(s,c,d) = (\frac{2s\:c\:d}{D_2},\frac{c^2-s^2d^2}{D_2},
		\frac{d^2-m\: s^2c^2}{D_2}).$

	\sssec{Theorem.}
	{\em $(E,+)$ is an abelian group.  Its order is divisible by 4.}

	The proof, although tedious, is straigthforward.
	The closure follows from the definition .6.
	Associativity follows, non trivially from .6.
	The neutral element is $(0,1,1).$
	The additive inverse element of $(s,c,d)$ is given by 1.10. C2. and C3.

	\sssec{Definition.}
	The group $(E,+)$ is called the {\em Jacobian elliptic group}
	associated to the prime $p$ and the integer $m \in {\Bbb Z}_p.$ 

	\sssec{Corollary.}
	{\em The following constitute special cases.}
	\begin{enumerate}
	\setcounter{enumi}{-1}
	\item	{\em For $m = 0,$ the elements of the group are\\
	\hth$	(sin k ,cos k ,1) and (sin k ,cos k ,-1),$\\
		    and the addition formulas reduce to\\
	\hth$    	(sin k ,cos k ,j_1) + (sin l ,cos l ,j_2) =$\\
	\hth$	(sin(j_1k+j_2l),cos(j_2k+j_2l),j_1j_2),$ $j_1$ and $j_2$ are
		$+1$ or $-1.$}
	\item	{\em For $m = 1,$ the elements of the group are\\
	\hth$	(tanh k,cosech k,cosech k), (tanh k,cosech k,-cosech k)$\\
	\hth	    and if $c^2 = -1,$\\
	\hth$(\infty,c\:\infty,c\:\infty),$ $(\infty,-c\:\infty ,c \:\infty),$
	$(\infty,c\:\infty,-c\:\infty),$ $(\infty,-c\:\infty,-c\:\infty ).$\\
	\hth	    and the addition formulas correspond to\\
	\hth$	tanh k_0 + tanh k_1
		= \frac{tanh k_0 + tanh k_1}{1 + tanh k_0 tanh k_1}.$\\
	\hth$	cosech(k_0+k_1)
		= \frac{cosech k_0 cosech k_1}{1 + tanh k_0 tanh k_1 }.$}
	\end{enumerate}

	\sssec{Comment.}
	To remove some of the mystery associated with some of the
	formulas just given, assume that the finite field is replaced by
	the field of reals.  For instance,
	D4, is obtained by replacing in D0, $c$ by $i\: s,$ $d$ by $i\:k\: s,$
	$c_1$ by $i\:s_1,$ $d_1$ by $i\:k\:s_1$ and letting $s$ and $s_1$
	tend to infinity.

	\ssec{Finite Jacobian elliptic groups for small p.}

	\setcounter{subsubsection}{-1}
	\sssec{Introduction.}
	It can be shown that $(E,+)$ is isomorphic to the direct
	product of the Klein 4-group and the group E associated to the
	finite Weierstrass p function introduced by Professor Tate and that the
	kernel of a homomorphism between the 2 groups is the subgroup of $(E,+)$
	of elements with $s = 0.$  A less precise form of this Theorem is given
	in 2.1. and is illustrated by the examples given in this section and
	prepares for the definition of finite Jacobi elliptic functions.
	In many cases the generator of the larger group allows the inclusion
	of one of the generators of the Klein 4-group.

	\sssec{Theorem.}
	{\em $(E,+)$ is isomorphic to ${\Bbb Z}_2 \bigx {\Bbb Z}_{2n}$
	or to ${\Bbb Z}_4 \bigx {\Bbb Z}_{4n}.$}

	\sssec{Example.}
	In the example the generators of the factor groups will be
	given.  The additional information in the second column will be
	explained in the Chapter on isomorphisms and homomorphisms.\\
	$\begin{array}{cccl}
	p&   m&		   E\:is\:isomorphic&     generator\\
	&&to\:{\Bbb Z}_i \bigx{\Bbb Z}_n&   of\:{\Bbb Z}_i\:and\:{\Bbb Z}_n \\
	\cline{1-4}
	\cline{1-4}
	3&   2&	   {\Bbb Z}_2 \bigx{\Bbb Z}_2  &(0,1,-1), (0,-1,1)\\
	\cline{1-4}
	5&   3&	   {\Bbb Z}_2 \bigx{\Bbb Z}_2  &(0,1,-1), (0,-1,1)\\
	&    2 = m'(2)	   {\Bbb Z}_2 \bigx{\Bbb Z}_4  &(0,1,-1), (1,0,2)\\
	&    4 = m_j(2)&	''	  &(0,1,-1), (\infty ,2 \:\infty ,\infty )\\
	\cline{1-4}
	7&   3&	   {\Bbb Z}_2 \bigx{\Bbb Z}_2  &(0,1,-1), (0,-1,1)\\
	&    2&	   {\Bbb Z}_2 \bigx{\Bbb Z}_4  &(0,-1,1), (2,2,0)\\
	&    4 = m''(2)&	''	  &(0,1,-1), (1,0,2)\\
	&    6 = m'(4)&	''	  &(0,1,-1), (1,0,3)\\
	&    5&	   {\Bbb Z}_2 \bigx{\Bbb Z}_6  &(0,1,-1), (2,-2,3)\\
	\cline{1-4}
	11&  5&	   {\Bbb Z}_2 \bigx{\Bbb Z}_4  &(0,-1,1), (3,5,0)\\
	&    8 = m'(9)&	''	  &(0,1,-1), (1,0,2)\\
	&    9 = m''(5)&	''	  &(0,1,-1), (1,0,5)\\
	&    2&	   {\Bbb Z}_2 \bigx{\Bbb Z}_6  &(0,1,-1), (3,5,4)\\
	&    6&		''	  &(0,1,-1), (5,3,4)\\
	&   10&		''	  &(0,-1,1), (5,3,-2)\\
	&    3&	   {\Bbb Z}_2 \bigx{\Bbb Z}_8  &(0,1,-1), (5,3,5)\\
	&    4 = m''(3)&	''	  &(0,-1,1), (3,5,3)\\
	&    7 = m'(3)&	''	  &(0,1,-1), (3,5,2)\\
	\cline{1-4}
	13&  2 = m'(2)&	   {\Bbb Z}_2 \bigx{\Bbb Z}_4  &(0,1,-1), (1,0,5)\\
	&   12 = m_j(2) = m''(12)&      ''	   &(0,1,-1),
		(\infty ,5 \:\infty ,\infty )\\
	&    4 = m'(10) = m''(10)& {\Bbb Z}_4 \bigx{\Bbb Z}_4 
		 &(\infty ,5 \:\infty ,3 \:\infty ), (1,0,6)\\
	&   10 = m_j(4)&	''	  &(\infty ,5 \:\infty ,4 \:\infty ), (1,0,2)\\
	&    6&	   {\Bbb Z}_2 \bigx{\Bbb Z}_6  &(0,1,-1), (2,6,4)\\
	&    8 = m_j(6)&	''	  &(0,-1,1), (6,2,5)\\
	&    3 = m_j(11)&	  {\Bbb Z}_2 \bigx{\Bbb Z}_8  &(0,1,-1), (6,2,6)\\
	&    5 = m_j(9)&	''	  &(0,1,-1), (6,2,4)\\
	&    9 = m''(3)&	''	  &(0,1,-1), (2,6,2)\\
	&   11 = m'(5)&	''	  &(0,1,-1), (2,6,3)\\
	&    7 = m_j(7)& {\Bbb Z}_2 \bigx{\Bbb Z}_{10}	&(0,-1,1), (2,6,-5)\\
	\cline{1-4}
	17&  2 = m'(2)&{\Bbb Z}_4 \bigx{\Bbb Z}_4&(\infty,4\:\infty,7\:\infty),
		 (1,0,4)\\
	&    9 = m''(2) = m_j(9)&	''	  &(\infty ,4\:\infty ,5 \:\infty ),
		 (1,0,3)\\
	&  16 = m_j(2) = m''(16)&      '' &(\infty,4\:\infty,1\:\infty ), (1,0,6)\\
	&    6&	   {\Bbb Z}_2 \bigx{\Bbb Z}_6  &(0,1,-1), (3,3,7)\\
	&    12 = m_j(6)&	       ''	  &(0,1,-1), (4,6,8)\\
	&    4&	   {\Bbb Z}_2 \bigx{\Bbb Z}_8  &(0,1,-1), (3,3,4)\\
	&    5 = m_j(13)&	       ''	  &(0,1,-1), (6,4,5)\\
	&    13 = m''(4) = m_j(5)&       ''	  &(0,1,-1), (6,4,3)\\
	&    14 = m_j(4) = m'(5)&       ''	  &(0,1,-1), (4,6,7)\\
	&    7&	   {\Bbb Z}_2 \bigx{\Bbb Z}_{10}	&(0,1,-1), (4,-6,5)\\
	&    11 = m_j(7)&	       ''	  &(0,1,-1), (3,3,-2)\\
	&    3&	   {\Bbb Z}_2 \bigx{\Bbb Z}_{12}	&(0,1,-1), (4,6,2)\\
	&    8&		''	  &(0,1,-1), (4,6,3)\\
	&    10 = m'(3) = m_j(8)&       ''	  &(0,1,-1), (3,3,8)\\
	&    15 = m_j(3) = m''(8)&       ''	  &(0,1,-1), (3,3,6)\\
	\end{array}$\\
	$\begin{array}{cccl}
	p&   m&		   E\:is\:isomorphic&     generator\\
	&&to\:{\Bbb Z}_i \bigx{\Bbb Z}_n&   of\:{\Bbb Z}_i\:and\:{\Bbb Z}_n \\
	\cline{1-4}
	19& 12&	   {\Bbb Z}_2 \bigx{\Bbb Z}_6  &(0,1,-1), (3,7,-8)\\
	&   3&	   {\Bbb Z}_2 \bigx{\Bbb Z}_8  &(0,1,-1), (7,3,5)\\
	&   11 = m'(3)&	''	  &(0,1,-1), (3,7,4)\\
	&    7 = m''(11)&	       ''	  &(0,-1,1), (2,4,7)\\
	&    4&		''	  &(0,1,-1), (2,4,2)\\
	&    5 = m''(4)&	''	  &(0,-1,1), (4,2,4)\\
	&   14 = m'(4)&	''	  &(0,1,-1), (4,2,9)\\
	&    2&	   {\Bbb Z}_2 \bigx{\Bbb Z}_{10}	&(0,1,-1), (4,-2,8)\\
	&   10&		''	  &(0,1,-1), (3,-7,5)\\
	&   18&		''	  &(0,1,-1), (2,4,9)\\
	&    6&	   {\Bbb Z}_2 \bigx{\Bbb Z}_{12}	&(0,-1,1), (7,3,7)\\
	&   15 = m'(16)&	       ''	  &(0,1,-1), (2,4,6)\\
	&   16 = m''(6)&	''	  &(0,1,-1), (3,7,3)\\
	&    9&		''	  &(0,1,-1), (4,2,3)\\
	&   13 = m'(9)&	''	  &(0,1,-1), (3,7,6)\\
	&   17 = m''(9)&	''	  &(0,-1,1), (7,3,2)\\
	&    8&	   {\Bbb Z}_2 \bigx{\Bbb Z}_{14}	&(0,1,-1), (2,4,8)\\
	\cline{1-4}
	23&  4&	   {\Bbb Z}_2 \bigx{\Bbb Z}_8  &(0,-1,1), (4,10,11)\\
	&    6 = m''(4)&	''	  &(0,1,-1), (8,11,10)\\
	&   15 = m'(6)&	''	  &(0,1,-1), (11,8,7)\\
	&    5&	   {\Bbb Z}_2 \bigx{\Bbb Z}_{10}	&(0,1,-1), (4,10,6)\\
	&   10&		''	  &(0,1,-1), (4,-10,5)\\
	&   17&		''	  &(0,1,-1), (9,-9,2)\\
	&    2&	   {\Bbb Z}_2 \bigx{\Bbb Z}_{12}	&(0,-1,1), (11,8,9)\\
	&   12 = m''(2)&	''	  &(0,1,-1), (9,9,8)\\
	&   22 = m'(12)&	       ''	  &(0,1,-1), (10,4,3)\\
	&    3&		''	  &(0,-1,1), (8,11,4)\\
	&    8 = m''(3)&	''	  &(0,1,-1), (10,4,11)\\
	&   11 = m'(8)&	''	  &(0,1,-1), (11,8,2)\\
	&   13&		''	  &(0,-1,1), (9,9,11)\\
	&   16 = m''(13)&	       ''	  &(0,1,-1), (8,11,9)\\
	&   21 = m'(16)&	       ''	  &(0,1,-1), (9,9,5)\\
	&   7&	    {\Bbb Z}_2 \bigx{\Bbb Z}_{14}	&(0,1,-1), (4,10,2)\\
	&   14&		''	  &(0,1,-1),  (8,-11,5)\\
	&   19&		''	  &(0,1,-1), (8,-11,2)\\
	&    9&	   {\Bbb Z}_2 \bigx{\Bbb Z}_{16}	&(0,-1,1), (4,10,8)\\
	&   18 = m''(9)&	''	  &(0,1,-1), (10,4,8)\\
	&   20 = m'(18)&	       ''	  &(0,1,-1), (4,10,7)
	\end{array}$




	\ssec{Finite Jacobian Elliptic Function.}

	\sssec{Definition.}
	Given a prime $p$ and an integer $m$ in ${\Bbb Z}_p,$ 3.2.1. defines an
	cyclic group of order $2n$ and $4n.$  If we choose a generator
	g := $(s_1,c_1,d_1)$ of this group, we obtain by successive addition
	n g = $n(s_1,c_1,d_1)$ = $(s_n,c_n,d_n).$  The {\em finite Jacobi
	elliptic functions}
	$sn$, $cn$ and $dn$, $scd$ are defined by\\
	\hth$ sn(n) := s_n,$ $cn(n) := c_n,$ $dn(n) := d_n,$
	$scd(n) := (s_n,c_n,d_n).$\\
	The {\em period} is denoted by $4K.$

	\sssec{Example.}
	For $p = 11,$ $m = 3,$ $K = 2,$	|
	For $p = 13,$ $m = 3,$ $K = 2,$\\
	$\begin{array}{ccccccccc}
	i&sn(i)&cn(i)&dn(i)&\vline&i&sn(i)&cn(i)&dn(i)\\
	0& 0& 1& 1&\vline&0& 0& 1& 1\\
	1&-5& 3&-5&\vline&1& 6& 2& 6\\
	2& 1& 0& 3&\vline&2& \infty &-5 \infty &-6 \infty \\
	3&-5&-3&-5&\vline &3&-6&-2&-6\\
	4& 0&-1& 1&\vline &4& 0&-1&-1\\
	5& 5&-3&-5&\vline &5& 6&-2&-6\\
	6&-1& 0& 3&\vline &6& \infty & 5 \infty & 6 \infty \\
	7& 5& 3&-5&\vline &7&-6& 2& 6\\
	8& 0& 1& 1&\vline &8& 0& 1& 1
	\end{array}$

	\sssec{Definition.}
	\begin{enumerate}
	\setcounter{enumi}{-1}
	\item$	ns := \frac{1}{sn}, nc := \frac{1}{cn},
	nd := \frac{1}{dn},$
	\item$	sc := \frac{sn}{cn}, cd := \frac{cn}{dn},
	ds := \frac{dn}{sn},$
	\item$	cs := \frac{cn}{sn}, dc := \frac{dn}{cn},
	sd := \frac{sn}{dn}.$
	\end{enumerate}

	The notation is due to Glaisher,
	Glaisher, J.W.L., On elliptic functions, Messenger of Mathematics,
	Vol. 11, 1881, 81-95.

	\ssec{Identities and addition formulas for finite elliptic
	functions.}

	\setcounter{subsubsection}{-1}
	\sssec{Introduction.}
	The formulas given in this section are for the most part the same as in
	the real case.  Theorem \ref{sec-taddjac} gives the addition
	formulas.  Theorem \ref{sec-tadduvw}, which may be new is needed to
	prove the
	addition formula for the Jacobi {\em Zeta function}.  Theorem
	\ref{sec-taddcnuvw} is given for sake of completeness.  It is clearly
	less elegant than \ref{sec-tadduvw}.

	\sssec{Lemma.}
	{\em
	\hth$1 - m s_0^2 s_1^2 = c_1^2 + d_0^2 s_1^2 = d_0^2 + m s_0^2 c_1^2.$}

	\sssec{Theorem.}
	{\em 
	\begin{enumerate}
	\setcounter{enumi}{-1}
	\item$      sn^2(u) cn^2(v) dn^2(v) - sn^2(v) cn^2(u) dn^2(u)$\\
	\hth$	  = (1 - m sn^2(u) sn^2 v) (sn^2(u) - sn^2(v).$
	\item$      cn^2(u) cn^2(v) - sn^2(u) sn^2(v) dn^2(u) dn^2(v$\\
	\hth$	  = (1 - m sn^2(u) sn^2 v) (1 - sn^2(u) - sn^2(v).$
	\item$      dn^2(u) dn^2(v) - m^2 sn^2(u) sn^2(v) cn^2(u) cn^2(v$\\
	\hth$  = (1 - m sn^2(u) sn^2 v)(1 - m sn^2(u) - m sn^2(v)
		+ m sn^2(u) sn^2(v).$
	\end{enumerate}}

	\sssec{Theorem.}\label{sec-taddjac}
	\begin{enumerate}
	\setcounter{enumi}{-1}
	\item$      sn(u+v) = \frac{sn^2(u) - sn^2(v)}
				{sn(u) cn(v) dn(v) -sn(v) cn(u) dn(u)}.$
	\item$      cn(u+v) = \frac{1 - sn^2(u) - sn^2(v)}
				{cn(u) cn(v) + sn(u) sn(v) dn(u) dn(v)}.$
	\item$dn(u+v) = \frac{1 - m sn^2(u) - m sn^2(v) + m sn^2(u) sn^2(v)}
			{dn(u) dn(v) + m sn(u) sn(v) cn(u) cn(v)}.$
	\item$      cn(u+v) = \frac{sn(u) cn(u) dn(v) - sn(v) cn(v) dn(u)}
		{sn(u) cn(v) dn(v) - sn(v) cn(u) dn(u)},$ for $u \neq (v) .$
	\item$      dn(u+v) = \frac{sn(u) dn(u) cn(v) - sn(v) dn(v) cn(u)}
		{sn(u) cn(v) dn(v) - sn(v) cn(u) dn(u)},$ for $u \neq (v) .$
	\end{enumerate}

	Formulas 0., 3. and 4. are due to Cayley (1884).

	\sssec{Theorem.}
	{\em 
	\begin{enumerate}
	\setcounter{enumi}{-1}
	\item$      - m cn(u) cn(v) cn(u+v) + dn(u) dn(v) dn(u+v) = 1 - m.$
	\item$      dn(v) dn(u+v) + m cn(u) sn(v) sn(u+v) = dn(u).$
	\item$      sn(v) dn(u) sn(u+v) + cn(v) cn(u+v) = cn(u).$
	\end{enumerate}}

	\sssec{Theorem.}\label{sec-tadduvw}
	{\em
	\begin{enumerate}
	\setcounter{enumi}{-1}
	\item$      sn(u+v+w) (sn(v) sn(u+w)-sn(w) sn(u+v))$\\
	\hth$		= sn(u)( sn(v) sn(u+v) - sn(w) sn(u+w) ).$\\
	\item$	sn(a_0-a_1) sn(a_1-a_2) sn(a_2-a_0)$\\
	\hth$	- sn(a_1-a_2) sn(a_2-a_3) sn(a_3-a_1)$\\
	\hth$	+ sn(a_2-a_3) sn(a_3-a_0) sn(a_0-a_2)$\\
	\hth$	- sn(a_3-a_0) sn(a_0-a_1) sn(a_1-a_3) = 0$
	\footnote{3.11.83}.
	\end{enumerate}}

	Proof:  If we write $u = a_0-a_1,$ $v = a_1-a_2,$ $w = a_3-a_0,$ then
	$u+v = a_0-a_2,$\\
	$u+w = a_3-a_1,$ $u+v+w = a_3-a_2$ and we obtain 0, from 1.
	To prove 1, let us introduce the notation\\
	\hth$s_0 := sn a_0,$ $s_1 := sn a_1,$ $s_2 := sn a_2,$ $s_3 := sn a_3.
$\\
	and similarly for $c_i$ and $d_i.$ 
	Let\\
	\hth$	B_0 := (s_1^2 - s_2^2) (s_2^2 - s_3^2) (s_3^2 - s_1^2)$\\
	\hti{12}$(s_0 c_1 d_1 + s_1 c_0 d_0) (s_0 c_2 d_2 + s_2 c_0 d_0)$\\
	\hti{12}$		(s_0 c_3 d_3 + s_3 c_0 d_0).$\\
	Let $B_1,$ $B_2,$ $B_3$ be obtained by adding 1, 2, 3 modulo 4 to each
	digit, using \ref{sec-taddjac}.0 in 1. and reducing to the same
	denominator, we have to prove that\\
	\hth$	B_0 - B_1 + B_2 - B_3 = 0.$\\
	Using .0.D0.,\\
	$B_0 = (s_1^4 (s_3^2 - s_2^2) + s_2^4 (s_1^2 - s_3^2)
		+ s_3^4 (s_2^2 - s_1^2))$ \\
	\hth$	(  s_1 s_2 s_3 c_0 d_0 (1 - s_0^2) (1 - m s_0^2)$\\
	\hti{12}$	 + s_2 s_3 s_0 c_1 d_1 (1 - s_0^2) (1 - m s_0^2)$\\
	\hti{12}$	 + s_3 s_0 s_1 c_2 d_2 (1 - s_0^2) (1 - m s_0^2)$\\
	\hti{12}$	 + s_0 s_1 s_2 c_3 d_3 (1 - s_0^2) (1 - m s_0^2)$\\
	\hti{12}$	 + s_0^3  c_1 d_1 c_2 d_2 c_3 d_3$\\
	\hti{12}$	 + s_0^2 s_1 c_2 d_2 c_3 d_3 c_0 d_0$\\
	\hti{12}$	 + s_0^2 s_2 c_3 d_3 c_0 d_0 c_1 d_1$\\
	\hti{12}$	 + s_0^2 s_3 c_0 d_0 c_1 d_1 c_2 d_2)$\\
	therefore\\
	$B_0 - B_1 + B_2 - B_3 =$\\
	\hth$(s_1 s_2 s_3 c_0 d_0 (s_0^4s_1^4(s_3^2-s_2^2)(m-m) +  \ldots $\\
	\hti{12}$	+ s_0^4s_1^2s_2^2(1+m-1-m) +  \ldots $\\
	\hti{12}$	+ s_0^4s_1^2(1-1) +  \ldots ) +  \ldots  )$\\
	\hti{12}$+ (s_0 c_1 d_1 c_2 d_2 c_3 d_3
		 (s_0^4 s_1^2 s_2^2 (1-1) + \ldots $\\
	\hti{12}$	+ s_0^2 s_1^4 s_2^2 (-1+1) +  \ldots $\\
	\hti{12}$+ s_1^4 s_2^2 s_3^2 (-1+1) +  \ldots  ) +  \ldots  ) = 0.$\\
	The given terms come from\\
	\hth$s_1^4(s_3^2-s_2^2)s_1s_2s_3c_0d_0ms_0^4$ in $B_0$ and from the
	term in $B_1$\\
	\hth corresponding to the term $s_3^4(s_2^2-s_1^2)s_0s_1s_2c_3d_3
		ms_0^4$ in $B_0,$\\
	the term in $B_2$ corresponding to $-s_2^4s_3^2s_3s_0s_1c_2d_2
	(-1-m)s_0^2$ in $B_0$ and
	the term in $B_1,$ to $-s_3^4s_1^2s_0s_1s_2c_3d_3(-1-m)s_0^2$ in $B_0,$
	the term in $B_2$ corresponding to $-s_2^4s_3^2s_3s_0s_1c_2d_2$ in
	$B_0$ and the term in $B_3,$ to $-s_1^4s_2^2s_2s_3s_0c_1d_1$ in $B_0$,
	the term in $B_2$ corresponding to $-s_2^4s_3^2
		s_0^2s_2c_3d_3c_0d_0c_1d_1$ in $B_0$ and
		the term in $B_1$, to $-s_3^4s_1^2
		s_0^2s_3c_0d_0c_1d_1c_2d_2$ in $B_0,$
	the term $-s_1^4s_2^2$ $s_0^3c_1d_1c_2d_2c_3d_3$ in $B_0$ and
	the term in $B_2,$ to $s_3^4s_2^2s_0^2s_2c_3d_3c_0d_0c_1d_1$ in $B_0$
	and the term in $B_2$ corresponding to $-s_3^4s_1^2s_0^2s_2c_3d_3c_0d_0
		c_1d_1$ in $B_0$ and
	the term in $B_3,$ to $-s_2^4s_3^2s_0^2s_1c_2d_2c_3d_3c_0d_0$ in $B_0,$
	The reduction involves $4(6 . 2 . 2+4 . 3 . 4+4 . 3 . 2) + 4(3+6+3)2$
	terms, which exausts the list of $4(6(4 . 4+4)) = 480$ terms in
	$B_0 - B_1 + B_2 -B_3$.

	\sssec{Comment.}
	Formula \ref{sec-tadduvw}.0, should be compared with the formula of
	Jacobi, (Crelle Vol. 15)\\
	\hth$	sn(u+v+w) sn(u) (1 - m sn(v) sn(w) sn(u+v) sn(u+w))$\\
	\hth$		= sn(u+v) sn(u+w) - sn(v) sn(w).$\\
	Formulas \ref{sec-tadduvw}.0. to 1. should also be compared with the
	formulas of Glaisher (1881) and of Cayley (Crell Vol. 41),

	\sssec{Corollary.}
	{\em
	\begin{enumerate}
	\setcounter{enumi}{-1}
	\item$      sn(u+1) = \frac{sn(1) (sn(1) sn(2) - sn(u-1) sn(u))}
		{sn(1) sn(u) -sn(2) sn(u-1)},$ $u = 3,$  \ldots.\\
	\end{enumerate}}

	Proof:  Use \ref{sec-tadduvw}.0. with $u = v = 1$ and $w = u-1,$

	\sssec{Theorem.}\label{sec-taddcnuvw}
	\begin{enumerate}
	\setcounter{enumi}{-1}
	\item$      cn(u+v+w) = \frac{
		\begin{array}{c}
		sn(u) dn(v) dn(w) (cn(v) cn(u+v) - cn(w)cn(u+w))\\
		- dn(u)  (sn(v) cn(v) dn(w)- sn(w)cn(w)dn(v) )\end{array}}
		{dn(u)  (sn(v) dn(w)cn(u+w) - sn(w)dn(v) cn(u+v))}.$\\
	\item$	sd(a_1-a_2) cn(a_1-a_2)  - sd(a_3-a_0) cn(a_3-a_0)$\\
	\hth$ = sd(a_0-a_1) (cn(a_0-a_2) cn(a_2-a_1)
		- cn(a_0-a_3) cn(a_3-a_1))$\\
	\hti{12}$- cn(a_2-a_3) (cn(a_2-a_0) sd(a_0-a_3)
		- sd(a_2-a_1) cn(a_1-a_3))$
	\end{enumerate}

	Proof:  One proof is to derive first \ref{sec-taddcnuvw}.1. using the
	same method as in \ref{sec-tadduvw}, the other is to set $a_3 = 0$ and
	 derive the corresponding formula using \ref{sec-taddjac}.

	\ssec{Double and  half arguments.}

	\sssec{Theorem.}
	\begin{enumerate}
	\setcounter{enumi}{-1}
	\item$	sn(2u) = \frac{2 sn(u)cn(u)dn(u)}{1 - m sn^4(u)}.$
	\item$	cn(2u) = \frac{cn^2(u) - sn^2(u)dn^2(u)}{1 - m sn^4(u)}$\\
	\hth$       = \frac{cn^2(u)- sn^2(u)dn^2(u)}{cn^2(u)+ sn^2(u)dn^2(u)}.$
	\item$	dn(2u) = \frac{dn^2(u)- m sn^2(u)cn^2(u)}{1 - m sn^4(u)}.$\\
	\hth$  = \frac{dn^2(u)+ cn^2(u)(dn^2(u)-1)1 - m sn^4(u)}
			{dn^2(u)- cn^2(u)(dn^2(u)-1)}.$
	\end{enumerate}

	\sssec{Theorem.}
	D0.\hti{4}$   s_1 := \sqrt{\frac{1 - c}{1 + d}},$\\
	D1.\hti{4}$   c_1 := \sqrt{\frac{1 + d}{1 + d}},$\\
	D2.\hti{4}$   d_1 := \frac{s (c + d)}{(1 + c) (1 + d) s_1 c_1},$\\
	{\em then}\\
	C0.\hti{4}$   d_1^2 = \frac{c+d}{1+c},$\\
	C1.\hti{4}$   2 (s_1,c_1,d_1) = (s,c,d).$

	Proof.  C0. follows directly from D0. to D1.  It is not used to
	define $d_1$ to insure that $2(s_1,c_1,d_1)$ is $(s,c,d)$ not
	$(-s,c,d).$\\
	The formulas can be derived starting from\\
	\hth$	d_1^2 - m s_1^2 c_1^2 = d (1 - m s_1^4).$\\
	Expressing $c_1^2$ and $d_1^2$ in terms of $s_1^2$ gives\\
	\hth$	m (1+d) s_1^4 - 2 m s_1^2 + 1 - d = 0,$ hence\\
	\hth$	s_1^2 = \frac{m + j \sqrt{m^2 - m (1 - d^2)}}{m (1 + d)},$\\
	where $j = +1$ or $-1,$ hence\\
	\hth$	s_1^2 = \frac{1 - j c}{1 + d}.$\\
	therefore\\
	\hth$	c_1^2 = \frac{j c + d}{1 + d}$ and
	$d_1^2 = \frac{m1 + d + j m c}{1 + d} = \frac{j c + d}{1 + j c}.$\\
	It remains to verify, by substitution, for $c$ and $s.$\\
	For $c,$\\
	\hth$	c_1^2 - s_1^2d_1^2 = \frac{2 j c (j c + d)}{(1 + d) (1 + j c)},$
		therefore $j = 1.$\\
	For $s,$\\
	\hth$	1 + m s_1^4 = 1 + \frac{m (1 - c)^2}{(1 + d)^2}
		= 1 + \frac{1 - c}{1 + d}  \frac{1 - d}{1 + c}$\\
	\hth$	= \frac{2 (c + d)}{(1 + d)(1+c)} = \frac{2 s_1 c_1 d_1}{s},$\\
	\hth$	2 s_1 c_1 d_1 = \sqrt{\frac{1 - j\:c}{1 + j\:c}}
		  \frac{j\:c + d}{1 + d} = \frac{2 s (c + d)}{(1 + c)(1 + d)}$\\
	\hth$		= \frac{2 j_1\: s\: (j c + d)}{(1 + j\: c) (1 + d)},$\\
	$j_1 = +1$ or $-1$ has to be determined once the square roots have been
	chosen unambiguously.

	\sssec{Example.}
	$p = 19,$ $m = 2,$ $\delta$ $^2$ = 2,\\
	Let $(s,c,d) = (4,2,8),$\\
	$s_1^2	= 2,$ $c_1^2 = -1,$ $d_1^2 = -3,$ therefore\\
	$s_1 = \delta$, $c_1 = 3 \delta$  or $-3 \delta$, $d_1 = 4$ or $-4.$

	\sssec{Theorem.}
	{\em If $dn(u) \neq -1,$ then at $u$,
	\begin{enumerate}
	\setcounter{enumi}{-1}
	\item$	sn \circ \frac{1}{2} I = \sqrt{\frac{1 - cn}{1 + dn}}.$
	\item$	cn \circ \frac{1}{2} I = \sqrt{\frac{cn + dn}{1 + dn}}.$
	\item$	dn \circ \frac{1}{2} I
		= \sqrt{\frac{1 - m + m cn + dn }{1 + dn }}.$
	\end{enumerate}}

	\sssec{Theorem.}
	{\em If $sn(u) = 0,$ $cn(u) = 1$ and $dn(u) = -1,$ then
	\begin{enumerate}
	\setcounter{enumi}{-1}
	\item$	sn (\frac{u}{2}) = \sqrt{\frac{1}{m}}$
	\item$	cn (\frac{u}{2}) = \sqrt{\frac{m - 1}{m}}.$
	\item$	dn (\frac{u}{2}) = 0.$
	\end{enumerate}}

	\sssec{Theorem.}
	{\em If $sn(u)= 0,$ $cn(u)= -1$ and $dn(u)= -1,$ then
	\begin{enumerate}
	\setcounter{enumi}{-1}
	\item$	sn (\frac{u}{2}) = \infty .$
	\item$	cn (\frac{u}{2}) = \sqrt{-1} \infty .$
	\item$	dn (\frac{u}{2}) = \sqrt{-m} \infty .$
	\end{enumerate}}

	\sssec{Conjecture.}
	\begin{enumerate}
	\setcounter{enumi}{-1}
	\item$	(1-m)$ R $p  \Rightarrow   scd(2K) = (1,0,\sqrt{1-m}).$
	\item$	-1$ R $p$ and $-m$ R $p  \Rightarrow  
	scd(K) = (\sqrt{-1}\infty , \sqrt{-m}\infty , \infty )$ and\\
	\hth$			      scd(2K) = (0, -1, -1).$
	\item$	m$ R $p and (m-1)$ R $p  \Rightarrow
	   scd(K) = (\sqrt{\frac{1}{m}}, \sqrt{1-\frac{1}{m}}, 0)$ and\\
	\hth$			      scd(2K) = (0, 1, -1).$
	\end{enumerate}

	\sssec{Conjecture.}\label{sec-conjsnkmu}
	{\em If $(1-m)$ R $p$ then}
	\hth$      sn(K-u) = cd(u).$

	\ssec{The Jacobi Zeta function.}

	\setcounter{subsubsection}{-1}
	\sssec{Introduction.}
	Definitions \ref{sec-dzsubi} are inspired by the relation which exist,
	in the real case, between the Jacobi {\em Zeta function}, the
	$\theta$  functions and the Weierstrass $\zeta$ function. 
	See Handbook p.578, 16.34 and p.650, 18.10.7.

	\sssec{Definition.}
	The {\em function $u$} is defined by:\\
	\hth$	u(1) := 0,$\\
	\hth$	u(i+1) := u(i) - m sn(1) sn(i) sn(i+1).$

	\sssec{Definition.}
	The {\em Jacobi Zeta function $Z$} is defined by\\
	\hth$	Z(1) := - \frac{u(K)}{K}.$\\
	\hth$	Z(i) := u(i) + Z(1) i,$ $i\neq 1.$

	\sssec{Theorem.}
	{\em 
	\begin{enumerate}
	\setcounter{enumi}{-1}
	\item$      Z(u+v) = Z(u) + Z(v) - m sn(u) sn(v) sn(u+v).$
	\item$      Z(u+v) = Z(u) + Z(v) + m sd(u) (cn(v) cn(u+v) - cn(u) )$
	\item$      Z(u+v) = Z(u) + Z(v) + sc(u) (dn(v) dn(u+v) - dn(u) )$
	\end{enumerate}}

	Proof of 0.  The formula is true, by definition, for $v = 1.$
	It follows by induction on $v$ and from $\ldots$ \ref{sec-tadduvw}
	Indeed,\\
	\hth$	Z(u+v+1) = Z(u+v) + Z(1) - m sn(1) sn(u+v) sn(u+v+1)$\\
	\hth$		= Z(u) + Z(v) + Z(1) - m sn(u) sn(v) sn(u+v)
				- m sn(1) sn(u+v) sn(u+v+1)$\\
	\hth$		= Z(u) + Z(v+1) + m sn(1) sn(v) sn(v+1)
		    - m sn(u) sn(v) sn(u+v) - m sn(1) sn(u+v) sn(u+v+1)$\\
	\hth$		= Z(u) + Z(v+1) - m sn(u) sn(v+1) sn(u+v+1).$

	The proof of 1. and 2. is left as an exercise.
	Hint:  Use \ref{sec-tadduvw}.

	\sssec{Theorem.}
	{\em 
	\begin{enumerate}
	\setcounter{enumi}{-1}
	\item$      Z(K-u) = -Z(u) + m sn(u) cd(u).$
	\item$Z(\frac{1}{2} K) = \frac{m}{2} sn(\frac{K}{2})
		 cd(\frac{K}{2}),$ if $K$ is even.
	\item$      Z(K) = 0.$
	\item$      Z(K+u) = - Z(K-u).$
	\item$      Z(2K-u) = - Z(2K+u).$
	\item$      Z(2K+u) = Z(u).$
	\end{enumerate}}

	Proof:
	0, follows from the additional formula for ${\mathit Zeta}(K-u)$
	and from \ldots \ref{sec-conjsnkmu}
	See Example 3.1.1. and $\setminus$130 elliptic.bas

	\sssec{Definition.}\label{sec-dzsubi}
	\begin{enumerate}
	\setcounter{enumi}{-1}
	\item$ z_1(u) := Z(u) + cn(u) ds(u)$
	\item$ z_2(u) := Z(u) - dn(u) sc(u).$
	\item$ z_3(u) := Z(u) + m sn(u) cd(u).$
	\item$ z_4(u) := Z(u).$
	\end{enumerate}


	\ssec{Example.}
	Several examples of Jacobian elliptic functions follow.\\[10pt]

$\begin{array}{cllcllcll}
p =  5&m =  3\:\:\delta^2 =  2&&p =  5&m =  2\:\:\delta^2 =  2&&p =  5&m =  4
\:\:\delta^2 = 2\\
.5 &( 1,-1, 2\delta)&&	   .5 &( 1\delta, 2, 1\delta)
&&	.5 &( 1\delta, 1\delta, 1)\\
 1 &( 0,-1, 1)&&	       1 &( 1, 0, 2)
&&	       1 &(\infty , 2\infty , 1\infty )\\
&2K = 1&&		      2 &( 0,-1, 1)&&	       2 &( 0,-1,-1)\\
	&&&		    3 &(-1, 0, 2)
&&	       3 &(\infty ,-2\infty ,-1\infty )\\
	&		   &&&2K = 2 		      &&&2K = 2\\

p =  7&m =  3\:\:\delta^2 =  3&&p =  7&m =  2\:\:\delta^2 =  3&&p =  7&m =  4
\:\:\delta^2 = 3\\
.5 &( 1,-1, 2\delta)&&	   .5 &( 3\delta, 3, 1\delta)
&&	.5 &( 2\delta, 1\delta, 3)\\
 1 &( 0,-1, 1)&&	       1 &( 2, 2, 0)&&	       1 &( 1, 0, 2)\\
&2K = 1&&		      2 &( 0, 1,-1)&&	       2 &( 0,-1, 1)\\
	&&&		    3 &(-2, 2, 0)&&	       3 &(-1, 0, 2)\\
	&		   &&&2K = 2 		      &&&2K = 2\\

p =  7&m =  6\:\:\delta^2 =  3&&p =  7&m =  5\:\:\delta^2 =  3\\
.5 &( 3, 3\delta, 1\delta)&&	.5 &( 3\delta, 3, 3\delta)\\
 1 &( 1, 0, 3)&&	       1 &( 2,-2, 3)\\
 2 &( 0,-1, 1)&&	       2 &(-2, 2,-3)\\
 3 &(-1, 0, 3)&&	       3 &( 0,-1,-1)\\
&2K = 2&&		      4 &( 2, 2,-3)\\
&&&			    5 &(-2,-2, 3)\\
&			   &&&2K = 3\\

p = 11&m =  5\:\:\delta^2 =  2&&p = 11&m =  8\:\:\delta^2 =  2&&p = 11&m =  9
\:\:\delta^2 =  2\\
.5 &(-3\delta, 4, 4\delta)&&       .5 &( 2, 2\delta, 1\delta)
&&       .5 &( 1\delta, 4\delta, 4)\\
 1 &( 3, 5, 0)&&	      1 &( 1, 0, 2)&&	      1 &( 1, 0, 5)\\
 2 &( 0, 1,-1)&&	      2 &( 0,-1, 1)&&	      2 &( 0,-1, 1)\\
 3 &(-3, 5, 0)&&	      3 &(-1, 0, 2)&&	      3 &(-1, 0, 5)\\
&2K = 2		    &&&2K = 2		    &&&2K = 2\\

p = 11&m =  2\:\:\delta^2 =  2&&p = 11&m =  6\:\:\delta^2 =  2&&p = 11&m = 10
\:\:\delta^2 =  2\\
.5 &( 2\delta, 2,-3\delta)&&       .5 &( 2, 2\delta, 4\delta)
&&       .5 &( 1\delta, 4\delta, 5)\\
 1 &( 3, 5, 4)&&	      1 &( 5, 3, 4)&&	      1 &( 5, 3,-2)\\
 2 &(-3,-5,-4)&&	      2 &( 5,-3, 4)&&	      2 &( 5, 3, 2)\\
 3 &( 0,-1,-1)&&	      3 &( 0,-1, 1)&&	      3 &( 0, 1,-1)\\
 4 &( 3,-5,-4)&&	      4 &(-5,-3, 4)&&	      4 &(-5, 3, 2)\\
 5 &(-3, 5, 4)&&	      5 &(-5, 3, 4)&&	      5 &(-5, 3,-2)\\
&2K = 3		    &&&2K = 3		    &&&2K = 3\\
\end{array}$

$\begin{array}{cllcllcll}
p = 11&m =  3\:\:\delta^2 =  2&&p = 11&m =  4\:\:\delta^2 =  2&&p = 11&m =  7
\:\:\delta^2 =  2\\
.5 &( 5\delta, 5\delta, 4)&&.5 &( 4\delta, 1\delta, 4)
&&.5 &( 5\delta, 5\delta, 5)\\
 1 &(-5, 3,-5)&&	      1 &( 3, 5, 3)&&	      1 &(-3, 5, 2)\\
 2 &( 1, 0, 3)&&	      2 &( 5, 3, 0)&&	      2 &( 1, 0, 4)\\
 3 &(-5,-3,-5)&&	      3 &( 3, 5,-3)&&	      3 &(-3,-5, 2)\\
 4 &( 0,-1, 1)&&	      4 &( 0, 1,-1)&&	      4 &( 0,-1, 1)\\
 5 &( 5,-3,-5)&&	      5 &(-3, 5,-3)&&	      5 &( 3,-5, 2)\\
 6 &(-1, 0, 3)&&	      6 &(-5, 3, 0)&&	      6 &(-1, 0, 4)\\
 7 &( 5, 3,-5)&&	      7 &(-3, 5, 3)&&	      7 &( 3, 5, 2)\\
&2K = 4		    &&&2K = 4 		    &&&2K = 4\\

p = 13&m =  2\:\:\delta^2 =  2&&p = 13&m = 12\:\:\delta^2 =  2&&p = 13&m =  4
\:\:\delta^2 =  2\\
.5 &(-5\delta, 4, 3\delta)&&       .5 &( 3\delta, 6\delta, 2)
&&       .5 &( 1\delta,-5, 4\delta)\\
 1 &( 1, 0, 5)&&	      1 &(\infty , 5\infty , 1\infty )
&&       1 &( 1, 0, 6)\\
 2 &( 0,-1, 1)&&	      2 &( 0,-1,-1)&&	      2 &( 0,-1, 1)\\
 3 &(-1, 0, 5)&&	      3 &(\infty ,-5\infty ,-1\infty )
&&       3 &(-1, 0, 6)\\
&2K = 2		    &&&2K = 2		    &&&2K = 2\\

p = 13&m = 10\:\:\delta^2 =  2&&p = 13&m =  6\:\:\delta^2 =  2&&p = 13&m =  8
\:\:\delta^2 =  2\\
.5 &( 3, 3\delta, 1\delta)&&       .5 &(-5, 1\delta, 6\delta)
&&       .5 &( 1\delta,-5,-5\delta)\\
 1 &( 1, 0, 2)&&	      1 &( 2, 6, 4)&&	      1 &( 6, 2, 5)\\
 2 &( 0,-1, 1)&&	      2 &( 2,-6, 4)&&	      2 &( 6, 2,-5)\\
 3 &(-1, 0, 2)&&	      3 &( 0,-1, 1)&&	      3 &( 0, 1,-1)\\
&2K = 2&&		     4 &(-2,-6, 4)&&	      4 &(-6, 2,-5)\\
&&&			   5 &(-2, 6, 4)&&	      5 &(-6, 2, 5)\\
&			  &&&2K = 3		    &&&2K = 3\\

p = 13&m =  3\:\:\delta^2 =  2&&p = 13&m =  5\:\:\delta^2 =  2&&p = 13&m =  9
\:\:\delta^2 =  2\\
.5 &(-5\delta, 4, 6\delta)&&       .5 &( 3\delta, 3, 1\delta)
&&       .5 &( 6\delta, 6\delta, 4)\\
 1 &( 6, 2, 6)&&	      1 &(-6, 2, 4)&&	      1 &( 2, 6, 2)\\
 2 &(\infty ,-5\infty ,-6\infty )&&       2 &( 1, 0, 3)
	      2 &(\infty ,-5\infty ,-2\infty )\\
 3 &(-6,-2,-6)&&	      3 &(-6,-2, 4)&&	      3 &(-2,-6,-2)\\
 4 &( 0,-1,-1)&&	      4 &( 0,-1, 1)&&	      4 &( 0,-1,-1)\\
 5 &( 6,-2,-6)&&	      5 &( 6,-2, 4)&&	      5 &( 2,-6,-2)\\
 6 &(\infty , 5\infty , 6\infty )&&       6 &(-1, 0, 3)
	      6 &(\infty , 5\infty , 2\infty )\\
 7 &(-6, 2, 6)&&	      7 &( 6, 2, 4)&&	      7 &(-2, 6, 2)\\
&2K = 4		    &&&2K = 4		    &&&2K = 4\\

p = 13&m = 11\:\:\delta^2 =  2&&p = 13&m =  7\:\:\delta^2 =  2\\
.5 &( 1\delta,-5, 3\delta)&&       .5 &(-5\delta, 4, 1\delta)\\
 1 &( 2, 6, 3)&&	      1 &( 2, 6,-5)\\
 2 &( 1, 0,-4)&&	      2 &( 6,-2, 3)\\
 3 &( 2,-6, 3)&&	      3 &( 6,-2,-3)\\
 4 &( 0,-1, 1)&&	      4 &( 2, 6, 5)\\
 5 &(-2,-6, 3)&&	      5 &( 0, 1,-1)\\
 6 &(-1, 0,-4)&&	      6 &(-2, 6, 5)\\
 7 &(-2, 6, 3)&&	      7 &(-6,-2,-3)\\
&2K = 4&&		     8 &(-6,-2, 3)\\
&&&			   9 &(-2, 6,-5)\\
&			  &&&2K = 5
\end{array}$


	\ssec{Other results.}

	\sssec{Theorem.}
	H0.\hti{5}$   \js{1-m}{p} = 1,$\\
	{\em then}
	C0.\hti{5}$   (1,0,k_1) \in E,$\\
	C1.\hti{5}$   (1,0,k_1) + (1,0,k_1) = (0,-1,1).$

	\sssec{Definition.}
	Let $e = (s,c,d),$ $s \neq \infty,$\\
	\hth$	sin(2e) := 2s\:c,$ $cos(2e) := c^2-s^2.$

	Can this be justified?\\
	This is done better using $sn = sin \circ am,$ $cn = cos \circ am.$

	\sssec{Theorem?.}
	$sin(2e_0+2e_1) = \ldots$.

	\sssec{Theorem. [Landen]}
	{\em Let}
	H0.\hti{5}$   e_0 := (s_0,c_0,d_0) \in E,$\\
	H1.\hti{5}$   s_0\: c_0\: d_0 \neq  0,$ $s_0 \neq  oo,$\\
	H2.\hti{5}$   e_1 := (s_1,c_1,d_1) := (s_0,c_0,d_0) + (1,0,k_1),$\\
	H3.\hti{5}$l(2e_0) := \frac{sin(2e_1)\:cos(2e_0)-cos(2e_1)\:sin(2e_0)}
				{sin(2e_0)-sin(2e_1)},$\\
	{\em then}
	C0.\hti{5}$   l(e_0) = \frac{1-k_1}{1+k_1}.$

	$l$ or $l(p,m)$ is the Landen constant associated to $p$ and $m.$

	Proof.\\
	P0.\hti{5}$   e_1 = (\frac{c}{d},-\frac{k_1\:s}{d},\frac{k_1}{d}).$\\
	P1.\hti{5}$l(e_0) = \frac{k_1(s^2-c^2) + c^2-k_1^2\:s^2}{d^2+k_1}
		      = \frac{1-k_1}{1+k_1}.$

	\sssec{Comment.}
	We can replace in the above Theorem $(1,0,k_1)$ by
	$(-1,0,k_1),$ this gives the same constant $l.$\\
	We can also replace $k_1$ by $-k_1,$ this gives the constant\\
	\hth$	l_1 = \frac{1}{l}.$\\


	\ssec{Isomorphisms and homomorphisms.}

	\sssec{Theorem.}
	{\em If $k_1$ is real, there exists an isomorphism  $\phi'$ between the
	elliptic group associated to $m$ and that associated to\\
	\hth$	m' = \frac{m}{m-1)},$\\
	\hth$\phi'(s,c,d) := (\frac{k_1\:s}{d},\frac{c}{d},\frac{1}{d})$\\
	$\phi'(s,c,0) := (\:\infty,\frac{c}{k_1s)}\:\infty,
	\frac{1}{k_1s}\:\infty),$\\
	$\phi'(\:\infty,c \:\infty,d \:\infty)
	:= (\frac{k_1}{d},\frac{c}{d},0).$}

	\sssec{Corollary.}
	{\em If $k_1$ is real the order of the group associated to $m$ and
	to $\frac{m}{m-1)}$ are the same.}

	\sssec{Theorem. [Jacobi]}
	{\em If $k$ is real, there exists an isomorphism $\phi''$
	between the elliptic group associated to $m$ and that associated to\\
	\hth$	m'' = \frac{1}{m},$\\
	\hth$\phi''(s,c,d) := (k\:s,d,c)$\\
	$\phi''(\:\infty,c \:\infty ,d \:\infty )
	:= (\:\infty ,\frac{d}{k} \:\infty ,\frac{c}{k} \:\infty)$.}

	\sssec{Corollary.}
	{\em If $k$ is real the order of the group associated to $m$ and
	to $\frac{1}{m}$ are the same.}

	\sssec{Theorem. [Jacobi]}
	{\em If $p \equiv 1 \pmod{4},$ there exists an isomorphism $\phi_1$ 
	between the elliptic group associated to $m$ and that associated to\\
	\hth$	m_j = m_1,$\\
	$\phi_1(s,c,d) := (\sqrt{-1}\frac{s}{c},\frac{1}{c},\frac{d}{c}),$\\
	$\phi_1(s,0,d) := (\:\infty,\frac{c}{\sqrt{-1}}s)\:\infty,
	\frac{d}{\sqrt{-1}}s)\:\infty),$\\
	$\phi_1(\:\infty,c \:\infty ,d \:\infty)
	:= (\frac{\sqrt{-1}}{c},0,\frac{d}{c}).$}

	\sssec{Corollary.}
	{\em If $p \equiv 1 \pmod{4}$, the order of the group associated to $m$
	and to $m_1$ are the same.}

	\sssec{Theorem. [Gauss]}
	{\em If $k$ is real, there exist a homomorphism $\phi_G$ from the
	elliptic group associated to $m$ and that associated to\\
	\hth$	m_G = \frac{4k}{(1+k)^2},$\\
	$\phi_G(s,c,d) := (\frac{(1+k)s}{D},\frac{c\:d}{D},\frac{2}{D}-1),$
	where $D = 1+k\:s^2,$\\
	$\phi_G(s,c,d) := (\:\infty,\frac{c\:d}{(1+k)s}\:\infty,
	\frac{2}{(1+k}s)\:\infty),$ if $1+ks^2 = 0.$\\
	$\phi_G(\:\infty,c \:\infty,d \:\infty) := (0,c\:d,-1).$ CHECK cd\\
	The kernel of the homomorphism is \{(0,1,1), (0,-1,-1)\}.\\
	The image of the homomorphism is a subgroup of index 2.}

	\sssec{Corollary.}
	{\em If $k$ is real the order of the group associated to $m$ and
	to $\frac{4k}{(1+k)^2}$ are the same.}

	\sssec{Theorem. [Landen]}
	{\em If $k_1$ is real, there exist a homomorphism $\phi_L$ from
	the elliptic group associated to $m$ and that associated to\\
	\hth$	m_L = (\frac{1-k_1}{1+k_1})^2,$\\
	$\phi_L(s,c,d) := (\frac{(1+k_1)s\:c}{d},
	\frac{1+k_1}{m}\frac{(d^2-k_1)}{d},\frac{1-k_1}{m}\frac{d^2+k_1}{d}),$\\
	$\phi_L(s,c,0) := (\:\infty,-\frac{k}{m\:s\:c}\:\infty,
		\frac{k}{(1+k_1}^2s\:c)\:\infty),$\\
	$\phi_L(\:\infty,c \:\infty,d \:\infty) := (\:\infty,
	\frac{d^2}{m\:c)}\:\infty,\frac{d^2}{(1+k_1)}^2c)\:\infty),$\\
	The kernel of the homomorphism is \{(0,1,1), (0,-1,1)\}.\\
	The image of the homomorphism is a subgroup of index 2.}

	\sssec{Corollary.}
	{\em If $k_1$ is real the order of the group associated to $m$ and
	$+o (\frac{1-k_1}{1+k_1})^2$ are the same.}

	\sssec{Definition.}
	The {\em amplitude function} is defined by\\
	\hth$	sin \circ am = sn,$ $cos \circ am = cn.$\\
	The example below gives $sin(2k),$ $cos(2k)$ under $sin$ and $cos.$

	\sssec{Theorem.}
	{\em 
	C0.\hti{5}$   D\: am = dn,$\\
	C1.\hti{5}$   D\: sin = cos,$ $D\: cos = -sin.$\\
	C2.\hti{5}$   D\: sn = cn\: dn,$ $D\: cn = -sn\: dn,$
	$D\: dn = -m \:sn \:cn.$\\
	C3.\hti{5}$   D^2 (2 am) = -m \:sin \circ (2 am)$}

	Proof:  Using the derivative of composition of functions,\\
	\hth$D\: sn = D (sin \circ am) = cos \circ am D\: am = cn\: dn.$\\
	The other relations in C2, follow from $sn^2 + cn^2 = 1$ and
	$dn^2 + m\: sn^2 = 1.$  C3, follows from\\
	$D^2(2 am) = 2 D\: dn = -2m\: sn\: cn = -m 2 sin \circ am cos \circ am =
		-m\: sin \circ (2 am).$

	\sssec{Comment.}
	The derivatives will have to be defined in a separate
	section.  Somehow the connection with $p$-adic functions will have to be
	involved.\\
	If $| h|  < 1,$ then\\
	$sin(x + h) = sin(x)\: cos(h) + cos(x)\: sin(h),$\\
	$sin(x+h) - sin(x) = sin(x)\: (cos(h)-1) + cos(x)\: sin(h),$\\
	but\\
	$sin(h) = h + o(h)$ and $cos(h)-1 = h^2 + o(h),$\\
	hence\\
	$\frac{sin(x+h)-sin(x)}{h} = cos(x) + o(1)$ and $D\: sin = cos.$

	For the elliptic functions, we have, see for instance Handbook, l. c.,
	p. 575, 16.22.1 to .3 and\\
	$am(h) = h - \frac{h^3}{3!}m + \frac{h^5}{5!}m(4+m) -  \ldots,$\\
	$sn(h) = h + o(h),$\\
	$cn(h) = 1 + o(h),$\\
	$dn(h) = 1 + o(h),$\\
	$am(h) = h + o(h).$\\
	Hence $sn(x+h)-sn(x) = sn(x)(cn(h)\: dn(h) - 1) + cn(x)\: dn(x)\: sn(h)
		= h (cn(x)\: dn(x)) + o(h),$ therefore\\
	$D\: sn = cn\: dn.$\\
	$D\: (sin \circ am) = cos \circ am D\: am = cn \:dn,$ therefore
	$D\: am = dn.$

\section{Applications.}

	\ssec{The polygons of Poncelet.}

	\sssec{Definition.}
	Let us associated to $e = (s,c,d)$ the point $P(e) = (sin(e),cos(e),1).
$\\
	The set of points $P(e_0+j\:e),$ $j = 0,1, \ldots$  are the vertices of
	a polygon called the {\em polygon of Poncelet}.

	\sssec{Theorem. [Poncelet]}
	{\em The sides $P(e_0+j\:e) \times P(e_0+j\:e+e)$ are tangent to a
	circle.}

	Proof.  Let\\
	\hth$	e_1 := e_0 + j\:e,$  $e_2 := e_1 + e,$\\
	\hth$	P_1 := P(e_1),$ $P_2 := P(e_2),$ then\\
	\hth$P_1 = (2s_1c_1,c_1^2-s_1^2,1),$ \\
	\hth$	e_1 = (s_1,c_1,d_1),$\\
	\hth$	e_2 = ( \frac{s\:c_1d_1+s_1c\:d}{D},
		\frac{c\:c_1-d\:s\:d_1s_1}{D}, 
		\frac{d\:d_1-m\:s\:c\:s_1c_1}{D} ),$\\
	with $D = 1 - m\:s^2s_1^2,$ \\
	\hth$P_2 = (2\frac{(s\:c_1d_1+s_1\:c\:d)(c\:c_1-d\:s\:d_1s_1)}{D^2},$\\
	\hti{12}$\frac{(c\:c_1-d\:s\:d_1s_1)^2-(s\:c_1d_1+s_1c\:d)^2}{D^2},1),
$\\
	\hth$	P_1 \times P_2 = [\ldots].$

	\sssec{Corollary. [Landen]}
	{\em The lines $P_{\ldots} \times P_{\ldots}$ pass through a fixed
	point $L := (0,l,1)$ called the {\em point of Landen}.}

	This is a special case of the Theorem of Poncelet, when\\
	\hth$	(s_0,c_0,d_0) = (1,0,k_1).$

	\sssec{Construction.}
	Determination of Poncelet's polygons.\\
	Let the outscribed circle $\theta$  be\\
	\hth$	X_0^2 + (X_1 - d_1)^2 = S^2,$\\
	let the inscribed circle $\gamma$  be\\
	\hth$	R^2(t_0^2+t_1^2) = (c_1t_1 + t_2)^2,$\\
	Given a point $P = (P_0,P_1,P_2)$ on $\theta$  and a tangent
	$t = (t_0,t_1,t_2)$ to $\gamma$  through $P,$ the other tangent
	$u = (u_0,u_1,u_2)$ is given by\\
	\hth$	u_1 = t_2(P_0^2 - R^2P_2^2),
		u_2 = t_1(P_0^2c_1^2-R^2(P_0^2+P_1^2),
		u_0 = -\frac{P_1u_1+P_2u_2}{P_0}.$\\
	Given a tangent $t$ to $\gamma$  and a point on it and $\theta$, the
	other point $Q = (Q_0,Q_1,Q_2)$ common to $t$ and $\theta$ is given by\\
	\hth$Q_2 = 1,$ $Q_1 = 2\frac{t_0^2d_1-t_1t_2}{t_0^2+t_1^2},
		Q_0 = -\frac{Q_1t_1+t_2}{t_0}.$

\section{The Weierestrass functions.}

	\ssec{Complex elliptic functions.}

	\sssec{Definition.}
	Given $g,$ a non residue of $p,$ or $\js{g}{p} = -1,$ then\\
	$C_{p,g} = C_p$ is the set of pairs $(a,b),$ $a,b \in {\Bbb Z}_p$ such
	that\\
	\hth$(a,b) + (c,d) = (a+b,c+d),$\\
	\hth$(a,b) . (c,d) = (a\:c+b\:d\:g,a\:d+b\:c).$\\
	We could also write $(a,b)$ as $a + b \gamma$ , with $\gamma^2 = g.$

	\sssec{Definition.}
	Given $s,c,d \in C_p,$ we can repeat definition $\ldots$ .

	\sssec{Definition.}
	2 of the functions are pure imaginary, the third is real,\\
	3 types $S,$ $C,$ $D.$

	\sssec{Theorem.}
	H0.\hti{5}$   \delta ^2 = d,$\\
	D0.\hti{5}$   e_1 := (s_1 \delta ,c_1 \delta ,d_1),$
		$ e_2 := (s_2,c_2,d_2),$\\
	\hth$	e3 := (s_3 \delta,c_3 \delta, d_3),$\\
	D1.\hti{5}$   D_4 := 1 - m\: d\: s_1^2 s_2^2,$\\
	D2.\hti{5}$   s_4 := \frac{s_1c_2d_2 + s_2c_1d_1}{D_4},$
		$ c_4 := \frac{c_1c_2 - s_1s_2d_1d_2}{D_4},$\\
	\hth$	d_4 := \frac{d_1d_2 - m\: d\: s_1s_2c_1c_2}{D_4},$\\
		Hence replace $m d$ by $m'$\\
	D3.\hti{5}$   D_5 := 1 - m\: d^2 s_1^3 s_3^2,$\\
		$s_5 := \frac{(s_1c_3d_3 + s_3c_1d_1) d}{D_5},$
		$c_5 := \frac{(c_1c_3s_1 + s_3d_1d_3) d}{D_5},$\\
	\hth$	d_5 := \frac{d_1d_3 - m\: d^2 s_1s_3c_1c_3}{D_5},$\\
	H1.\hti{5}$   D_4 \neq  0,$ $D_5 \neq  0,$\\
	{\em then}\\
	C0.\hti{5}$   e_1 + e_2 = (s_4 \delta , c_4 \delta , d_4),$\\
	C1.\hti{5}$   e_1 + e3 = (s_5, c_5, d_5).$\\
	C2.\hti{5}$   2n\: e_1 \in E,$ $(2n+1) e_1 \in S.$

	\sssec{Definition.}
	$(a,b) > 0$ if $b = 0$ and $0 < a < \frac{p}{2}$\\
	\hth		or if $0 < b < \frac{p}{2}.$

	\sssec{Comment.}
	All that has been said above can be repeated.




	\ssec{Weiertrass' elliptic curves and the Weierstrass elliptic
	 functions.}

	\setcounter{subsubsection}{-1}
	\sssec{Introduction.}
	The modern work on elliptic curves starts and ends with
	the elliptic curves of Weierstrass.  I refer the reader to Lang S.

	\sssec{Theorem.}
	{\em Let}
	D0.\hti{5}$   e_3 := -\frac{1+m}{3},$\\
	D1.\hti{5}$   g_2 = 4\frac{m^2-m+1}{3},$
		$g_3 = \frac{4}{27}(m+1)(m-2)(2m-1).$\\
	D2.\hti{5}$  \Delta := g_2^3 - 27 g_3^2,$ $J := \frac{g_2^3}{\Delta},$\\
	D3.\hti{5}$   pn := e_3 + \frac{1}{s^2},$ $Dpn := -2\frac{cd}{s^3},$\\
	D4.\hti{5}$   e_2 := \frac{2-4m}{3},$\\
	D5.\hti{5}$   g'_2 := \frac{4}{3}16m^2-16m+1,$
		$g'_3 := \frac{8}{27}(2m-1)(32m^2-32m-1),$\\
	D6.\hti{5}$\Delta' := g_2^{'3} - 27 g_3^{'2},$
		$J' = \frac{g_2^{'3}}{\Delta'},$\\
	D7.\hti{5}$   qn := e_2 + \frac{1+c}{1-c},$
		$Dqn := -4\frac{d(1+c)^2}{s^3},$\\
	{\em then}\\
	C0.\hti{5}$   Dpn^2 = 4 pn^3 - g_2 pn - g_3.$\\
	C1.\hti{5}$   Dqn^2 = 4 qn^3 - g'_2 qn - g'_3.$\\
	C2.\hti{5}$  \Delta = 3(12m(m-1))^2,$
		$J = (m^2-m+1)^3 (\frac{2}{27m(m-1)})^2.$\\
	C3.\hti{5}$   \Delta' = 256m(m-1),$
		$J' = (16m^2-16m+1)^3\frac{1}{108m(m-1)}.$

	The proof is straithforward.  Substituting D2 in C0, multiplying by
	$s^4$ and expressing $c^2$ and $d^2$ in terms of $s^2$ gives a
	polynomial of the second degree in $s^2.$  The coefficients of 1,
	$s^2$ and $s^4$ give in turn, $e_3,$ $g_2$ and $g_3.$\\ 
	Substituting D5 in C1, multiplying by $(1-c)^2$ gives similarly a
	polynomial of the second degree in x := 1-c.  The coefficients of 1, $x$
	and $x^3$ give in turn $e_2,$ $g'_2$ and $g'_3.$\\
	$pn$ corresponds to the Weierstrass $P$ function and $Dpn$ to its
	derivative.\\
	The formulas correspond to those of real elliptic functions with the
	ratio of the period $\omega$ and the complete elliptic integral $K$ set
	to 1.\\
	(See for instance Handbook for Mathematical functions, p649, 18.9.1,
	2,3,4,5,8,9 and 11).

	\sssec{Example.}
	With $p = 7,$ $m = 2,$ then $e_3 = -1,$ $g_2 = 3,$ $g_3 = -1,$\\
		with $p = 7,$ $m = 6,$ then $e_3 = 0,$ $g_2 = 3,$ $g_3 = 0.$\\
	with $p = 19,$ $k^2 = 2,$ then $e_3 = 18,$ $g_2 = 4,$ $g_3 = 0,$\\
	\hti{18}$			   e_2 = 17$, $g'_2 = 6,$ $g'_3 = -1,$\\
	$\begin{array}{clcccccccccc}
	&\vline& sn&cn&dn&sin&cos&am&pn&Dpn&qn&Dqn\\
	\cline{1-12}
	  0&\vline&  0&1&1&0&1&0&\infty& \infty& \infty& \infty \\
	  1&\vline&  4&2&8&7&3&8&5&9&-5&5\\
	  2&\vline&  7&3&6&4&-2&1&6&2&-6&-8\\
	  3&\vline&  7&3&-6&-2&4&1&6&-2&-6&8\\
	  4&\vline&  4&2&-8&3&7&8&5&-9&-5&-5\\
	  5&\vline&  0&1& -1&1&0&0&\infty& \infty& \infty& \infty
	\end{array}$

	Jacobi $Z$ function = 0, -3,  -2,  2,  3,  0\\
	Weierstrass $\zeta$ function = $\infty$,  1,   6,  -6,  -1, $ \infty$. 

	with $p = 19,$ $k^2 = 3,$ then $e_3 = 5,$ $g_2 = 3,$ $g_3 =-9,$\\
	\hti{18}		   $e_2 = 3,$ $g'_2 = 9,$ $g'_3 = 5,$\\
	$\begin{array}{clcccccccccc}
	&\vline& sn&cn&dn&sin&cos&am&pn&Dpn&qn&Dqn\\
	\cline{1-12}
	  0&\vline&  0&1&1&0&1&0&\infty& \infty& \infty& \infty\\
	  1&\vline&  7&3&5&7&3&1&-7&8&1&3\\
	  2&\vline& -1&0&6&4&-2&-4&6&0&-7&7\\
	  3&\vline&  7&-3&5&-2&4&9&-7&-8&4&8\\
	  4&\vline&  0&-1&1&3&7&-9&-19&-19&3&-19\\
	  5&\vline& -7&-3&5&1&0&-8&-7&8&4&-8\\
	  6&\vline&  1&0&6&3&-7&5&6&0&-7&-7\\
	  7&\vline&  -7&3&5&-2&-4&0&-7&-8&1&-3\\
	  8&\vline&  0&1&1&4&2&0&\infty& \infty& \infty& \infty 
	\end{array}$

	Jacobi $Z$ function = 0,  -7,  0,  7,  0,  -7,  0,  7,  0.\\
	Weierstrass $\zeta$ function = $\infty$,  6, 0, -6, $\infty$,
		 6, 0, -6, $\infty$.. 

	If $2T$ is the period for the Jacobi functions, $T$ is the period of
	$pn,$ which is even and of $Dpn$ which is odd.
	See the next section for the last 2 columns.

	RERUN LAST EXAMPLE using ..[130]/ELLIPT

	\sssec{Theorem.}
	{\em Let $a := e_3.$}\\
	D0.\hti{5}$   (pn3,Dpn3) := (pn1,Dpn1) + (pn2,Dpn2),$\\
	D1.\hti{5}$   Q := (pn1-a)(pn2-1),$ $Q' := (pn1-a-1)(pn2-a-1),$\\
	{\em then}\\
	C0.\hti{5}$   pn3 - a = 4(pn1-a)(pn2-a)$\\
	C1.\hti{5}$   pn3 - a = ( \frac{(pn1-a)(pn2-a)-m}
		{(pn1-a)Dpn2+(pn2-a)Dpn1} )^2$\\
	C2.\hti{5}$   Dpn3 = \frac{Q(Q-m)(Dpn1Dpn2-4QQ')(Dpn1Dpn2-4mQ)}
		{Q(pn0-a)Dpn1+(pn1-a)Dpn0}.$\\
	RECHECK THIS (the line above)\\
	{\em If $pn1 \neq pn2$ then}\\
	C3.0.\hti{3}$ pn3 = ( \frac{Dpn1-Dpn2}{2(pn1-pn2)} )^2 - (pn1+pn2).$\\
	C4.0.\hti{3}$ Dpn3 = \frac{ pn3(Dpn1-Dpn2) + (pn1Dpn2-pn2Dpn1)}
		{pn2-pn1}.$\\
	{\em If $pn1 = pn2$ and $Dpn1 = -Dpn2$ then}\\
	C3.1.\hti{3}$ pn3 = \infty .$\\
	C4.1.\hti{3}$ Dpn3 = \infty .$\\
	{\em If $pn1 = pn2$ and $Dpn1 = Dpn2$ then}\\
	C3.2.\hti{3}$ pn3 = - 2pn1 +(\frac{3pn1^2-g_2\frac{1}{4}}{Dpn1})^2,$\\
	C4.2.\hti{3}$ Dpn3 = \frac{(3pn1^2-\frac{1}{4}g_2)(pn1-pn3)}{Dpn1}
		- Dpn1.$\\
	{\em If $pn1 = pn2$ and $Dpn1 = Dpn2 = 0$ then}\\
	C3.3.\hti{3}$ pn3 = \infty .$\\
	C4.3.\hti{3}$ Dpn3 = \infty .$

	The proof of C3 and C4 follow from the addition formulas for
	the Jacobi functions.  That of C2 and C3 is analoguous of the
	formulas for the real case (Handbook p.635, 18.4.1, 18.4.2)

	\sssec{Theorem.}
	{\em If}
	H0.\hti{5}$   (pn1,Dpn1) + (pn2,Dpn2) = (pn3,Dpn3),$\\
	{\em then}\\
	C0.\hti{5}$   (pn1,Dpn1),$ $(pn2,Dpn2),$ $(pn3,-Dpn3)$ are collinear.

	This follows at once from 3.10.3. and is the geometric interpretation
	of C3.

	\sssec{Theorem.}
	{\em If}
	D0.\hti{5}$   pn(ti,t) := t^{-2} pn(i),$\\
	D1.\hti{5}$   Dpn(ti,t) := t^{-3} Dpn(i),$\\
	D2.\hti{5}$   g_2(t) := t^{-4} g_2,$\\
	D3.\hti{5}$   g_3(t) := t^{-6} g_3,$\\
	{\em then}\\
	C0.\hti{5}$   Dpn(ti,t)^2 = 4 pn(ti,t)^3 - g_2(t) pn(ti,t) - g(3).$

	The proof follows at once from \ldots.\\
	ellinv.tab gives a table of the invariants $g_2(t)$ and $g_3(t)$ for
	$p = 19$ and $m = 2$ to 18.

	\sssec{Theorem.}
	{\em Making explicit the dependence of $g_2$ and $g_3$ on $m$,}
	C0.\hti{5}$   g2(m+1,t) = g2(-m,t) = g2(m+1,-t) = g2(-m,-t).$\\
	C1.\hti{5}$   g3(m+1,t) = -g3(-m,t) = -g3(m+1,-t) = g3(-m,-t).$\\
	C2.\hti{5}$   g2'(m+1,t) = g2'(-m,t) = g2'(m+1,-t) = g2'(-m,-t).$\\
	C3.\hti{5}$   g3'(m+1,t) = -g3'(-m,t) = -g3'(m+1,-t) = g3'(-m,-t).$

	The sections 3.10.7. to 10.12. were inspired from the formulas on
	complex elliptic functions.\\
	(See for instance, Handbook, p.635 18.4.3.,18.4.8)\\
	In the classical case, the Weierstrass $p$ function is defined in
	such a way that the constant term in the Maclaurin expension is 0.
	The Weierstrass $\zeta$ function is defined as its integral and the
	constant term is again chosen as zero, which is natural if we want
	$\zeta$ to be an odd function.  $\zeta$ is not periodic.\\
	In the finite case, we have chosen $pn$ and $\zeta$ using the same
	definition
	in terms of the Jacobi functions as in the real case, but now
	$\zeta$ as an odd periodic function.  This should be contrasted with
	the classical case in which the Weierstrass function is not periodic.
	Theorem 3.10.9. gives an interesting property of $\zeta(1).$

	\sssec{Definition.}
	\vspace{-18pt}\hspace{100pt}\footnote{5.12.83}\\[8pt]
	The {\em Weierstrass $\zeta$ function} is defined by\\
	\hth$	\zeta(u) = $Z$(u) + \frac{cn(u)\:dn(u)}{sn(u)}.$

	\sssec{Definition.}
	\vspace{-18pt}\hspace{100pt}\footnote{15.11.82}\\[8pt]
	The {\em function $u$} is defined as follows:\\
	Given $u(1) = 0,$\\
	D0.0.$  pn(i) \neq pn(j),$ \\
	\hth$u(i+j) := u(i) + u(j) + \frac{Dpn(i)-Dpn(j)}{2(pn(i)-pn(j))}.$\\
	D0.1.\hti{3}$ pn(i) = pn(j) and Dpn(i)\neq Dpn(j),$ $u(i+j) = \infty.$\\
	D1.0.\hti{3}$ Dpn(i) \neq  0,$
	$u(2i) := 2u(i) + \frac{3pn(i)^2 - \frac{1}{4}g_2}{Dpn(i)}.$\\
	D1.1.\hti{3}$ Dpn(i) = 0,$ $12pn(i)^2 = g(2),$ $u(2i) = \infty .$\\
	D2.0.\hti{3}$ u(0) = u(T) = \infty .$\\
	$T$ is the period \ldots.

	\sssec{Theorem.}
	{\em There exist a constant $\zeta(1)$ such that\\
	\hth$	u(j) + j\;\zeta(1) = u(T-j) + (T-j)\zeta(1),$
	$j = 1$ to $\frac{T}{2} -1.$}

	Proof.  $\ldots$ ?

	\sssec{Theorem.}
	{\em The Weierstrass $\zeta$ function is related to the $u$ function by
\\
	\hth$\zeta(j) := u(j) + j\:\zeta(1).$}

	\sssec{Comment.}
	The definitions and theorems can be repeated replacing
	respectively $pn,$ $Dpn,$ $u,$ $\zeta$ by $qn,$ $Dqn,$ $v,$ $\zeta',$
	but we have the additional property of the next Theorem.

	\sssec{Theorem.}
	{\em $\zeta'(\frac{T}{2})$ = 0.}

	Proof.  $\ldots$ ?

	\sssec{Theorem.}
	{\em $\zeta$ and $\zeta'$ are odd functions and their period is either
	$\frac{T}{2}$	or $T.$}

	\sssec{Notation.}
	$card(X)$ denotes the cardinality of the set $X$.

	\sssec{Theorem.}
	{\em If $p \equiv -1 \pmod{4},$ then}
	\begin{enumerate}
	\setcounter{enumi}{-1}
	\item    $card(g2,g3) + card(g2,-g3) = 2(p+1).$\\
	{\em If $p \equiv 1 \pmod{4},$ then}
	\item    $card(g2,g3) = card(g2,-g3).$
	\end{enumerate}

	\sssec{Corollary.}
	{\em 
	If $p \equiv-1 \pmod{4},$ then  $card(g2,0) = p+1.$}

	\sssec{Comment.}
	Examples can be obtained using P.BAS see P.HOM.

	\sssec{Definition.}
	The {\em function} $Ke$ is defined by\\
	\hth$	-p < Ke(m) \equiv K(m) < p,$ $Ke(m)$ is even.

	\sssec{Theorem. [Hasse conjectured by Artin]}
	\hth$	-2 \sqrt{p} < Ke(m) < 2 \sqrt{p}.$

	\sssec{Conjecture.}
	\vspace{-18pt}\hspace{100pt}\footnote{4.1.84}\\[8pt]
	\begin{enumerate}
	\setcounter{enumi}{-1}
	\item  {\em Given an integer  $x$  in the range}\\
	\hth$	-2 \sqrt{p} < x < 2 \sqrt{p}.$\\
	{\em then there exist a pair $(g2,g3)$ such that the corresponding
	Weierstrass elliptic curve $W2$ has $card\:p+1+x$ and the
	corresponding group is abelian.}\\
	This has been verified up to $p = 47$\footnote{4.1.84}.
	\item{\em If the cardinality of $W2$ is divisible by 4 and $W2$ is not
	abelian, then
	there exist a $J2$ or a $J3$ isomorphic to $W2$\footnote{6.1.84}.}
	\item{\em If $e_1+e_2+e_3 = 0,$ $e_i-e_k$ are all non quadratic residue,
	and $j' \neq 0$, then the elliptic group is isomorphic to\\
	\hth$	C_{4l+2} \bigx C_2$ for some $4l.$}\\
	This has been verified up to $p = 97.$  See g7622, Example.
	\end{enumerate}

.       \sssec{Comment.}
	If for a given m we obtain $J3(m)$ and $J2(m)$ and the corresponding
	$W2$ and $W2',$ $card(W2) = card(W2')$ but $W2$ and $W2'$ are not
	necessarily isomorphic.
	E.g. $p = 17,$ $J3(-2) = C_2 \bigx C_{12},$ $J2(-2) = C_{24}.$

	\sssec{Comment}
	\vspace{-18pt}\hspace{95pt}\footnote{6.1.84}\\[8pt]
	Excluding $(g2,g3) = (0,0),$ if $j' \neq 0,1$ there exists
	2 sets of elliptic curves, corresponding to $(g2,g3)$ and $(g2,-g3).$
	What is the connection between the structure if any?\\
	None except that concerning cardinality.  E.g. $p = 31,$
	$W2(1,5) \sim C_{37},$ $W2(1,-5) \sim C_9 \bigx C_3,$ 
	$W2(3,11) \sim C_{28},$ $W2(3,-11) \sim C_{12} \bigx C_3.$ 


%
	\ssec{The isomorphism between the elliptic curves in 3 and 2 dimensions.}

	This should be integrated with 3.1.

	\setcounter{subsubsection}{-1}
	\sssec{Introduction.}
	The usual correspondance in the real field between the functions $sn,$
	 $cn$ and $dn$ of Jacobi and $P,$ $DP$ of Weierstrass should
	be modified to insure an isomorphism between the 3 dimensional
	elliptic curve associated to $(sn,cn,dn)$ and the 2 dimensional
	elliptic curve associated to $(P,DP).$  This requires in fact
	to associate, in the real case to $sn(t),$ $P(2t).$

.       \sssec{Theorem.}
	The curve $(P,P')$ has a singularity when $m = 0$ and $m = 1,$
	When $m = 0,$ the singularity is $(-\frac{1}{3},0),$ because
	$-\frac{1}{3}$ is a double root
	of $4p^3$ - $g_2$ p - $g_3,$ the regular solution when $P' = 0$ is
	$(\frac{2}{3},0).$\\
	When $m = 1,$ the singularity is $(\frac{1}{3},0),$ because
	$\frac{1}{3}$ is a double root
	of $4p^3 - g_2 p - g_3,$ the regular solution when $P' = 0$ is
	$(-\frac{2}{3},0).$ 

	\sssec{Definition.}\label{sec-dell32}
	Let $(s,c,d)$ in $E,$ if $m = 0$ we add the restriction d = 1,
	if $m = 1,$ we add the restriction c = d.
	Let $e_3$ := $-\frac{1+m}{3},$ $e_2$ := $e_3$ + 1, $e_1$ := $e_3$ + m,
	\begin{enumerate}
	\setcounter{enumi}{-1}
	\item$      T(s,c,d) := ( e_3 + \frac{1+d}{1-c},
	2\frac{(c+d)(1+d)}{s(c-1)}),$\\
	\hth$		if s\neq frac{1}{0}, \infty .$
	\item 0.$    T(0,1,1) := ( \infty , \infty  ),$
		1.$    T(0,1,-1) := ( e_1, 0 ),$ $m \neq  0,1,$
		item 2.$    T(0,-1,1) := ( e_2, 0 ),$
		item 3.$    T(0,-1,-1) := ( e_3, 0 ),$\\
	\hth	(when $m = 0,$ $d = 1,$ when $m = 1,$ $d = -1),$
	\item$      T(\infty ,c \infty ,d \infty ) := ( e_3 - \frac{d}{c},
		 2\frac{(c+d)d}{c} ),$\\
	\hth$	(m \neq  0,$ for $m = 1$ and $c = \sqrt{-1},
	 T(\infty ,c \infty ,c \infty ) = (-\frac{5}{3},4c)).$
	\item$ e_1 = \frac{2m-1}{3},$ $e_2 = \frac{2-m}{3}.$\\
	\hth\ldots D0. gives $s_1 = \frac{1-c}{1+d},$
	$c_1 = \sqrt{\frac{c+d}{1+d}},$\\
	\hth$d1 = \frac{s(c+d)}{(1+c)(1+d)\:s_1\:c_1}.$
	\end{enumerate}

	\sssec{Theorem.}\label{sec-tell32}
	{\em Let $a := \frac{p'}{2(p-e_3)(p-e_3-1)},$ then
	\begin{enumerate}
	\setcounter{enumi}{-1}
	\item$ a^2 \neq  0, -1 \Rightarrow T^{-1}(p,p') = (s,c,d),$ where\\
	\hth$	c := \frac{1-a^2}{1+a^2},$ $s := \frac{c-1}{a},$
	 $d := (1-c)(p-e_3) -1.$
	\item$      a^2 = -1 \Rightarrow T^{-1}(p,p')
		= (\infty ,a \infty ,a(e_3-p) \infty ).$
	\item 0.$    T^{-1}(e[3],0) = (0,-1,-1),$ {\em for} $m \neq  0.$
		1.$    T^{-1}(e[3]+1,0) = (0,-1,1),$ {\em for} $m \neq  1.$
		2.$    T^{-1}(e[3]+m,0) = (0,1,-1),$ {\em for} $m \neq  0,1.$
	\item$      T^{-1}(\infty ,\infty ) = (0,1,1).$
	\end{enumerate}

	Proof:  If $a^2 \neq  0, -1,$ solving
	$(1-c)(p-3_3) = 1+d$ and $sp' = -2(c+d)(p-e_3)$ for $c$ gives
	$c-1 = a\:s$.  Hence $(c-1)^2 = a^2(1-c)^2,$ but $c \neq 1,$
	$1-c = a^2(1+c),$ hence
	$c, s = \frac{c-1}{a} = -2\frac{a}{1+a^2}$ and $d.$

	\sssec{Theorem.}\label{sec-tell32a}
	{\em Let $g_2:= 4\frac{m^2-m+1}{3}$ and
	$g_3 := \frac{4}{27}(m-2)(2m-1)(m+1),$ 
	let $T(s,cd) = (p,p'),$ then}
	\begin{enumerate}
	\setcounter{enumi}{-1}
	\item$      p'^2 = 4 p^3 - g_2 p - g_3.$
	\end{enumerate}

	\sssec{Definition.}
	The bijection defined by \ref{sec-dell32} and justified by
	\ref{sec-tell32}
	defines an [\em isomorphism between the 3 dimensional elliptic curve
	of $\ldots$  and the 2 dimensional elliptic curve} \ref{sec-tell32a}.0.

	\sssec{Definition.}\label{sec-dell32b}
	The invariant $j$ of {sec-tell32a}.0. is\\
	\hth$	j := 2^6 3^3 j',$ with\\
	\hth$	j' := \frac{g_2^3 }{g_2^3 - g_3^2}.$

	\sssec{Theorem.}
	\begin{enumerate}
	\setcounter{enumi}{-1}
	\item$      g_2^3 - g_3^2 = (m(m-1))^2.$
	\item$      j = 2^6 3^3 \frac{(m^2-m+1)^3}{(m(m-1))^2}.$
	\item{\em If $j$ is given and $M$ is a solution of 1, the other
		solutions
		are $1-m,$ $\frac{1}{m},$ $1-\frac{1}{m},$ $\frac{1}{1-m},$
		$\frac{m}{1-m}.$}
	\end{enumerate}

	\sssec{Example.}
	For $p = 11,$ let $u := m(m-1),$
	$j' = -3$ corresponds to $u = 1,-2,-5$ or $m = -3,4,3,-2,-4,5$ giving\\
	\hth$	K(j) = 4$ or $-4.$\\
	$j' = 4$ corresponds to $u = 2, -3$ or $m = 2,-1,-5$ giving $K(j) = 0.
$\\
	$j' = \infty$ corresponds to $u = 0$ or $m = 0, 1$ giving
	$K(j) = 1$ or $-1.$\\
	$j' = 0$ corresponds to $u = -1$ or $m = -5 + 2 \delta$, giving
	$ K(j) = 0.$

	\sssec{Theorem.}
	{\em Given $e_1,$ $e_2$ and $e_3$ such that\\
	H.0.\hti{4}$   e_1 + e_2 + e_3 = 0,$\\
	H.1.\hti{4}$   e_i$ are distinct,\\
	Let $m := \frac{e_1-e_3}{e_2-e_3},$ \\
	\hth$	e_2 - e_3 = d,$\\
	\hth$	g2 := 4 ( e_3^2 - e_1 e_2 ),$\\
	\hth$	g3 := 4 e_1 e_2 e_3.$\\
	then\\
	C.0.\hti{4}$   (d_p) = 1 \Rightarrow J3(m) \sim W2(g2,g3).$\\
	C.1.\hti{4}$   (d_p) = -1 \Rightarrow (p+1 - | J3(m)| ) 
	- (-1_p) (p+1 - |W2(g2,g3)|) = 0.$}

	Proof:  Let $c^2 = \frac{1}{d},$ $e'_i = c^2$ $e_i,$ then
	$e'_1 = e'_3 + m$ and $e'_2 = e'_3 + 1$\ldots.

	\ssec{Correspondance between the Jacobi elliptic curve $(cn,sd)$
	 and the  Weierstrass elliptic curve}

	\sssec{Definition.}\label{sec-djacwei}
	Let $m_1 := 1-m,$ $e_2 := 2\frac{1-2m}{3}.$ 
	\begin{enumerate}
	\setcounter{enumi}{-1}
	\item$      T(cn,sd) := ( e_2 + \frac{1+cn}{1-cn},
	 -4sd\frac{m_1+m\:cn^2}{(1-cn}^2),$\\
	\hth		if $cn \neq  1.$
	\item$      T(1,0) := (\infty , \infty ),$
	\item$      T(\infty ,sd) := ( e_2 - 1,-4m sd),$
	\item$      If -\frac{m_1}{m} = b^2 then$\\
	\hth$	T(b,\infty ) = (e_2 + \frac{1+b}{1-b}, 0).$
	\end{enumerate}

	\sssec{Theorem.}\label{sec-tjacwei}
	{\em Let $a := p - e_2 + 1,$ then}
	\begin{enumerate}
	\setcounter{enumi}{-1}
	\item$      a \neq  0,$ $cn := \frac{a-2}{a},$\\
	\hth{\em $m_1 + m cn^2 \neq  0 then T^{-1}(p,p') = (cn, sd),$\\
	\hth		where $sd := P' \frac{(1-cn)^2}{-4(m_1 + m cn^2)},$\\
	\hth$	m_1 + m cn^2 = 0$ then $T^{-1}(p,p') = (cn,\infty ),$}
	\item$      T^{-1}(\infty ,\infty ) = (1,0).$
	\item$      a = 0 \Rightarrow
	T^{-1}(p,p') = (\infty , \frac{p'}{-4m}).$
	\item$      T(-1,0) = (e_2,0).$
	\item{\em If $T^{-1}(x,0) = (y,\infty )$ and $x \neq  e_2$ then\\
	\hth$	m(m-1) = d^2$ and\\
	\hth$	x = e_1$ or $e_3 = -\frac{e_2}{2} \pm 2 d.$}
	\end{enumerate}

	\sssec{Theorem.}\label{sec-tjacweia}
	{\em Let $g_2 := \frac{4}{3}(16m^2-16m+1)$ and\\
	\hth$	g_3 := \frac{8}{27}(2m-1)(32m^2-32m-1),$\\
	let $T(cn,sd) = (P,P'),$ then}
	\begin{enumerate}
	\setcounter{enumi}{-1}
	\item$      P'^2 = 4 P^3 - g_2 P - g_3.$
	\end{enumerate}

	\sssec{Definition.}
	The bijection defined by \ref{sec-djacwei} and justified by
	\ref{sec-tjacwei}
	defines an {\em isomorphism between the 2 dimensional elliptic curve
	of $\ldots$  and the 2 dimensional elliptic curve} {sec-tell32a}.0.

	\sssec{Theorem.}
	{\em Using \ref{sec-dell32b},}
	\begin{enumerate}
	\setcounter{enumi}{-1}
	\item$      g_2^3 - g_3^2 = 256 m(m-1).$
	\item$      j = 16 \frac{(16m^2-16m+1)^3}{m(m-1)}.$
	\end{enumerate}

	\sssec{Corollary.}
	{\em 
	Let $c := \sqrt{1-\frac{1}{m}},$ $p = e_2 + \frac{1+c}{1-c},$ using
	$e_2$ from \ref{sec-djacwei} and
	$g_2,$	$g_3$ from {sec-tjacweia}, then\\
	\hth$	p^3 - g_2p - g_3  = -1.$\\
	$\ldots$	 DOUBLE CHECK THIS.}

	The proof follows from the fact that the denominator $m_1+m\:cn$ in
	\ref{sec-tjacwei}.0.
	cannot be zero.

	\sssec{Theorem.}
	{\em Given $e_1,$ $e_2$ and $e_3$ such that
	H.0.\hti{4}$   e_1 + e_2 + e_3 = 0,$\\
	H.1.\hti{4}$   1 - (\frac{e_1 - e_3}{3 e_2})^2 = f^2,$\\
	H.2.\hti{4}$   e_i$ are distinct,\\
	Let $m := \frac{1}{2}(1 + \frac{1}{d}),$\\
	\hth$	d = 2\frac{1-2m}{3e_2},$\\
	\hth$	g2 := 4 ( e_3^2 - e_1 e_2 ),$\\
	\hth$	g3 := 4 e_1 e_2 e_3.$\\
	then\\
	C.0.\hti{4}$   (d_p) = 1 \Rightarrow J3(m) \sim W2(g2,g3).$\\
	C.1.\hti{4}$   (d_p) = -1 \Rightarrow (p+1 - | J3(m)| ) 
	- (-1-p) (p+1 - |W2(g2,g3)|) = 0.$}

	Proof:  Let $e'_i = c^2$ $e_i.$  We have a $J2(m)$ if
	$e'_2 = 2\frac{1-2m}{3}$ and $e'_1-e'_3 = 4d = 4 \sqrt{m(m_1-1)}$
	because of \ldots.

	\sssec{Example.}
	$\begin{array}{lcrrccc}
	p&\vline&   e_1,e_2,e_3&  g2,g3&   j'&  m'&  structure\\
	\cline{1-7}
	13&\vline&  1,3,9&	   0,4& 0&   -5&  C_6 \bigx C_2\\ \\
	29&\vline&  1,9,19&	  -13,-12& 13&  3&   C_{14} \bigx C_2 \\
	37&\vline&  1,6,30&	  -13,17 & -6&  -14& C_{22} \bigx C_2 \\
	&\vline&	2,15,20&	 0,-5&0&   -5&  C_{14} \bigx C_2 \\
	41&\vline&  1,13,27&	 -6,10&   -6&  -12& C_{18} \bigx C_2 \\
	53&\vline&  1,19,33&	 -13,17&  -17& -22& C_{30} \bigx C_2 \\
	&\vline&	1,20,32&	 -12,16&  -21& 3&   C_{22} \bigx C_2 \\
	61&\vline&  1,8,52&	  -13,17&  15&  -23& C_{30} \bigx C_2 \\
	&\vline&	1,29,31&	 7,-3&26&  11&  C_{30} \bigx C_2 \\
	&\vline&	2,26,33&	 0,-29&   0&   -28& C_{38} \bigx C_2 \\
	73&\vline&  1,8,64&	  0,4& 0	&   C_{42} \bigx C_2 \\
	&\vline&	1,15,57&	 15,-11&  13&  -21& C_{42} \bigx C_2 \\
	&\vline&	1,29,43&	 -20,24&  -7&  -29& C_{34} \bigx C_2 \\
	89&\vline&  1,13,75&	 20,-16&  44&  -26& C_{50} \bigx C_2 \\
	&\vline&	1,25,63&	 23,-19&  -12& -30& C_{42} \bigx C_2 \\
	&\vline&	1,29,59&	 13,-9 &  15&  24&  C_{42} \bigx C_2 \\
	97&\vline&  1,8,88&	  1,3& 2&   -39& C_{50} \bigx C_2 \\
	&\vline&	1,18,78&	 14,-10&  -17& -20& C_{50} \bigx C_2 \\
	&\vline&	1,35,61&	 0,4& 0	&   C_{42} \bigx C_2 \\
	&\vline&	1,38,58&	 15,-11&  7&   -19& C_{42} \bigx C_2 \\
	&\vline& \ldots ..\\
	113&\vline&  \ldots \\
	&\vline&	1,22,90&	 -6,10&	   21&  C_{58} \bigx C_2 \\
	&\vline&	 \ldots 
	\end{array}$



\section{Complete elliptic integrals of the first and second kind.}

	\setcounter{subsubsection}{-1}
	\sssec{Introduction.}
	Several conjectures appear to justify the terminology of
	complete elliptic integrals of the first and second kind for the
	functions $K$ and $E$ defined below.
	The definitions are inspired from the definitions in the real and
	complex fields. Their importance is associated with Conjecture
	\ref{sec-cjellint}.
	By convention in this section I will use\\
	\hth$	q := [\frac{p}{2}] = \frac{p-1}{2}.$\\
	Again $I$ denotes the identity function (I(i) = i), $D$ denotes the
	derivative operator and f $\circ$ g denotes the function which is
	obtained by composition from the functions $f$ and $g.$

	\sssec{Definition.}\label{sec-dellintke}
	\vspace{-18pt}\hspace{90pt}\footnote{8-10.12.83}\\[8pt]
	\enumb
	\item$K_j := ( \frac{(2j-1)!!}{2j}!! )^{2,} \pmod{p},$
		$0 \leq  j \leq  q.$
	\item$E'_j := K_j \frac{2j}{2j-1)},$ $j = 0$ to $q,$
	\item$E_0 := 1,$ $E_j := - \frac{K_j}{2j-1},$ $j = 1$ to $q,$
	\enume

	\sssec{Definition.}
	\enumb
	\item$      K :=  \sum_{j=0}^q ( K_j I^j),$
	\item$      E' := \sum_{j=0}^q ( E'_j I^j),$
	\item$      E :=  \sum_{j=0}^q ( E_j I^j),$
	\enume

	\sssec{Example.}
	For $p = 11,$\\
	$\begin{array}{ccrrrrrrrl}
	 j &\vline& K& E'&E& D& B& C& K''&\hti{8}all [j]\\
	\cline{1-9}
	 0 &\vline& 1& 0& 1&-5&-5&-4&( 1)\\
	 1 &\vline& 3&-5&-3& 5&-2&-3& 1\\
	 2 &\vline& 1& 5&-4&-1& 2&-2&(-5)\\
	 3 &\vline& 1&-1& 2& 5&-4&-2& -4\\
	 4 &\vline& 3& 5&-2&-5&-3&-1&(-4)\\
	 5 &\vline& 1&-5&-5& 0& 1& 0& -1
	\end{array}$

	For $p = 13,$\\
	$\begin{array}{ccrrrrrrrl}
	 j &\vline& K& E'&E& D& B& C& K''&\hti{8}all [j]\\
	\cline{1-9}
	 0 &\vline& 1& 0& 1&-6&-6& 5& 1\\
	 1 &\vline&-3&-6& 3& 1&-4&-6&( 1)\\
	 2 &\vline& 4& 1& 3&-1& 5& 1& -6\\
	 3 &\vline&-3&-1&-2&-1&-2&-2&(-5)\\
	 4 &\vline& 4&-1& 5& 1& 3& 4& -6\\
	 5 &\vline&-3& 1&-4&-6& 3&-1&( 4)\\
	 6 &\vline& 1&-6&-6& 0& 1& 0& 2
	\end{array}$

	See \ref{sec-dellintkpp} and \ref{sec-dellintdbc}.

	\sssec{Lemma.}\label{sec-lellintke}
	\enumb
	\item$      K_j = K_{q-j},$ $0 \leq  j \leq  q.$
	\item$      E'_j = E'_{q+1-j},$ $0 < j \leq  q.$
	\item$      E'_j = \frac{2j-1}{2j} K_{j-1},$ $0 < j \leq  q.$
	\item$      E_{j+1} = \frac{(2j-1)(2j-3)}{4j^2}E_j,$ $0 < j \leq q.$
	\item$      2j K_j - (2j+1)K_{j-1} = 2j E_j,$ $0 < j \leq  q.$
	\item$      2j E_j + E'_j = 0,$ $0 \leq  j \leq  q.$
	\item$      2(j+1) E'_{j+1} - 2j E'_j - E_j = 0,$ $0 \leq  j < q.$
	\enume
	
	The proof\footnote{10.12.83} follows at once from Lemma 3.0.8.4 (g730)
	and from \ref{sec-dellintke}.

	\sssec{Theorem.}
	(see KE.NOT)\\
	\hth$	K(\frac{1}{m}) = (m_p) K(m), 0 < m < p.$

	Proof.\\
	$K(\frac{1}{m}) = m^q K(m),$ because of \ref{sec-lellintke}.0
	and $m^q = (m_p)$ 
	by the Theorem of Euler.  See for instance Adams and Goldstein, p. 107.

	\sssec{Corollary.}
	\hth$	p \equiv  1 \pmod{4}  \Rightarrow   K(-1) = 0.$

	\sssec{Theorem.}\label{sec-tellintde}
	\enumb
	\item$      2 I DK + K - 2DE' = 0.$
	\item$      2(1-I) DK - K - 2DE = 0.$
	\item$      2I DE + E' = 0.$
	\item$      2(1-I) DE' - E = 0.$
	\enume

	Proof:  This follow immediately from the definitions, for instance,
	the coefficient of $I^j$ is\\
	\hth$	2j K_j + K_j - 2(j+1)E'_{j+1} = 0 for 0 \leq  j < q,$\\
	and that of $I^q = 0$ because $2q+1 = p = 0.$

	\sssec{Corollary.}
	\enumb
	\item$      K(0) = 2DE'(0) = E(0) = 1.$
	\item$      DK(0) - DE(0) = \frac{1}{2}.$
	\item$      E'(0) = 0.$
	\item$      K(1) = -2DE(1) = 2DE'(1) - 2DK(1) = -E'(1).$
	\item$      E(1) = 0.$
	\enume

	\sssec{Theorem.}\label{sec-tellintde1}
	\enumb
	\item$      4 D( I(1-I) DK) - K := 0.$
	\item$      4 I D( (1-I) DE') + E' := 0.$
	\item$      4 (1-I) D( I DE) + E := 0.$
	\item$      K = E + E'.$
	\enume

	This derives from \ref{sec-tellintde} by elimination, for instance,
	eliminating E' from 2 and 3 gives\\
	\hth$	4(1-I) D(I DE) = -2(1-I) DE' = -E$\\
	\hth$	4 I D( (1-I) DE') = 2I DE = -E'.$\\
	1, times I gives\\
	\hth$	4I(1-I) DK - 2I K - 4 I DE = 0,$\\
	using 2 and 0 gives\\
	\hth$	4I(1-I)DK - 2I K + 2 E' = 0,$\\
	\hth$	4 D( I(1-I) DK) - 2 D(I K) + 2I DK + K = 0.$\\
	Finally, it follows from \ref{sec-tellintde}.0 and 1 that
	$D(K - E - E') = 0,$ but $K(0) = E(0) + E'(0) = 1,$ hence 3.

	\sssec{Definition.}\label{sec-dellintkpp}
	\enumb
	\item$      K''_j := \frac{( (2j-3)(2j-7) \ldots  )^2}{j!} \pmod{p},$
	\item$      K'' :=  \sum_{j=0}^{\frac{q}{2}} ( K_{q-2j} I^{(q-2j)} ).$
	\enume

	\sssec{Lemma.}\label{sec-lellintde}
	\enumb
	\item[0.0.]$    p \equiv  1 \pmod{4}  \Rightarrow   K'' is even,$
	\item[	1.]$    p \equiv  -1 \pmod{4}  \Rightarrow   K'' is odd.$
	\item[1.~~]$(q-j+2)(q-j+3) K''_{q-j+2} = (2q-4j+1)^2 K''_{q-j},$
		$1 < j \leq  q.$
	\item[2.~~]$      D( (1-4I^2) DK'') - K'' = 0.$
	\item[3.~~]$      c K''$ {\em are, for arbitrary constant $c,$ the only
		solutions of 2}.
	\enume

	To prove 3., we substitute\\
	\hth$	 \sum_{j=0}^{q} ( X_{q-j} I^{(q-j)}),$\\
	in 2, this gives\\
	\hth$	- (2q+1)^2 X_q) I^{(q)} - (2q-1)^2 X_{q-1}) I^{(q-1})$\\
	\hti{12}$+  \sum_{j=2}^{q} ( (q-j+2)(q-j+3) X_{q-j+2}
		- (2q-4j+1)^2 X_{q-j}) I^{q-j} = 0.$\\
	The coefficient of $I^q$ is zero because $2q+1 = 0,$ hence $X_q$ is
	arbitrary.  The coefficient of $I^{(q-1)}$ must be zero, therefore
	$X_{q-1} = X_{q-3} = \ldots = 0.$

	\sssec{Lemma.}
	\enumb
	\item$K \circ (1-I)$ {\em satisfies} \ref{sec-tellintde1}.0.
	\item$E \circ (1-I)$ {\em satisfies} \ref{sec-tellintde1}.1.
	\item$K \circ (I+\frac{1}{2})$ {\em satisfies} \ref{sec-lellintde}.2.
	\item$K'' = s K'' \circ (-I).$
	\item$K = s K \circ (1-I).$
	\item$E = s E' \circ (1-I).$
	\enume

	For instance,\\
	\hth$	K \circ (1-I) = 4 (D( I(1-I) DK) \circ (1-I)$\\
	\hth$		  = -4 D( (1-I)I DK \circ (1-I) )$\\
	\hth$		  = 4 D( (1-I)I D( K \circ (1-I))).$  Hence 0.
	3, from \ref{sec-lellintde}.0 and $s = (-1)^{\frac{p-1}{2}}.$
	$K = c K'' \circ (I-\frac{1}{2}) = s c K'' \circ (\frac{1}{2} -I)
		 = s K \circ (1-I),$ hence 4.
	Finally, $E \circ (1-I)$ and $E'$ satisfy the same differential
	equation, of second order, moreover $E(1) = s E'(0) = 0$ and
	$DE(1) = -s DE'(0)$ because of \ldots and of $K(1) = s K(0)$ hence 5.

	\sssec{Corollary.}
	\enumb
	\item$K(m) = s K(1-m),$ {\em with} $s = (-1)^q$.
	\item$K(\frac{m}{m-1)}) = (\frac{m^2}{1-m)}) K(m),$ $2 < m < p-1.$
	\item$E(m) + s E(1-m) = K(m)$
	\item$\frac{E(m)}{K(m)} + \frac{E(1-m)}{K(1-m)} = 1.$
	\item$m $R$ p,$ $k^2 = m  \Rightarrow K(\frac{4k}{1+k)}^{2)} = K(m),$
		$2 < m < p-1.$
	\item$1-m $R$ p, k^2 = 1-m  \Rightarrow  
		 K(((1-k)(1+k))^{2)} = K(m),$ $2 < m < p-1.$
	\enume

	4 and 5 are still conjectures.

	\sssec{Corollary.}
	\hth$(-1)^k \sum_{k=0}^j (\frac{(2j-1)!!}{(2j)!!})^2 \ve{j}{k}
		= (-1)^q (\frac{(2k-1)!!}{(2k)!!})^2.$ 

	Exchanging $k$ and $j,$ 1. follows from Lemma 3.0.x. of g730.
	Theorem 7.

	\sssec{Corollary.}
	\enumb
	\item$      p \equiv  -1 \pmod{4}  \Rightarrow   K(\frac{p+1}{2}) = 0,$
	\item$      p \equiv  1 \pmod{4}  \Rightarrow   K(\frac{p+1}{2})
		= 2 E(\frac{p+1}{2}).$
	\enume

	\sssec{Conjecture.}\label{sec-cjellint}
	With the exception of $p = 7,$ $K(3) = 3 (= -4),$\\
	\hth$	m \neq  0,1,$ $| K(m)|  < \frac{p}{2}  \Rightarrow
		K(m) \equiv  p+1 \pmod{4}.$

	\sssec{Corollary.}
	\hth$p \equiv 1 \pmod{4},  \Rightarrow   K(m) \neq 0.$

	\sssec{Example.}
	For $p = 11,$\\
	$\begin{array}{ccrrrrrrrrl}
	 m &\vline& K& E'&E& D& B& C& K''&\frac{E}{K}&\hti{8}all (m)\\
	\cline{1-9}
	 0 &\vline& 1& 0& 1& 1& 1& 1& 0& 1\\
	 1 &\vline&-1&-1& 0& 5&-5& 5&-4& 0\\
	 2 &\vline& 0&-1& 1&-4&-5&-4& 4&\infty\\
	 3 &\vline& 4&-4&-3& 2& 3& 2& 4& 2\\
	 4 &\vline& 4&-5&-2& 3& 4& 3& 0& 5\\
	 5 &\vline&-4& 4& 3&-3&-2&-3&-1& 2\\
	 6 &\vline& 0&-5& 5& 0&-1& 0& 1&\infty\\
	 7 &\vline& 4&-3&-4& 2& 1& 2& 0&-1\\
	 8 &\vline&-4& 2& 5& 0&-1& 0&-4&-4\\
	 9 &\vline&-4& 3& 4&-2&-1&-2&-4&-1\\
	10 &\vline& 0&-1& 1&-4&-5&-4& 4&\infty
	\end{array}$

	For $p = 13,$\\
	$\begin{array}{ccrrrrrrrrl}
	 m &\vline& K& E'&E& D& B& C& K''&\frac{E}{K}&\hti{8}all (m)\\
	\cline{1-9}
	 0 &\vline& 1& 0& 1& 1& 1& 1& 1& 1\\
	 1 &\vline& 1& 1& 0& 6&-6& 6& 4& 0\\
	 2 &\vline& 6& 6& 0&-6& 6&-6&-4& 0\\
	 3 &\vline&-2& 3&-5& 1& 2& 1&-4&-4\\
	 4 &\vline&-2& 1&-3& 3& 4& 3&-4&-5\\
	 5 &\vline&-2&-3& 1&-5&-6&-5&-1& 6\\
	 6 &\vline& 2& 4&-2& 5& 4& 5& 2&-1\\
	 7 &\vline&-6&-3&-3& 4& 3& 4& 2&-6\\
	 8 &\vline& 2&-2& 4&-2&-3&-2&-1& 2\\
	 9 &\vline&-2& 1&-3& 3& 4& 3&-4&-5\\
	10 &\vline&-2&-3& 1&-6&-5&-6&-4& 6\\
	11 &\vline&-2&-5& 3&-3&-4&-3&-4& 5\\
	12 &\vline& 6& 0& 6&-1& 0&-1& 4& 1
	\end{array}$

	\sssec{Definition.}\label{sec-dellintdbc}
	By analogy with the case of the real or complex field,
	\enumb
	\item$      D_j := K_{j+1} - E_{j+1},$ $0 \leq  j < q,$ $D_q := 0.$
	\item$      B_j := K_j - D_j,$ $0 \leq  j \leq  q.$
	\item$      C_j := D_{j+1} - B_{j+1},$ $0 \leq  j < q,$ $C_q := 0.$
	\item$      D :=  \sum_{j=0}^{q} ( D_j I^j),$
	\item$      B :=  \sum_{j=0}^{q} ( B_j I^j),$
	\item$      C :=  \sum_{j=0}^{\frac{q}{2}} ( C_j I^j).$
	\enume

	\sssec{Theorem.}
	\enumb
	\item$      D_j = E_{j+1},$ $0 \leq  j < q.$
	\item$      B_j = E_{q-j},$ $0 \leq  j \leq  q.$
	\item$      D = \frac{E'}{I}.$
	\item$      I B = I E + (I-1) E'.$
	\item$      I^2 C = (2-I) E' - I E.$
	\item$      I^2 C = 2 E' - I K.$	?
	\enume




\section{P-adic functions, polynomials, orthogonal polynomials.}

	\sssec{Comment.}
	In a $p$-adic field, we can define polynomials of degree up to $p-1.$
	These are determined by their values at $i$ in ${\Bbb Z}_p.$  If these
	are defined in the real field with rational coefficient, the definition
	and properties are automatically extended to the $p$-adic field.
	For orthogonal polynomials, recurrence relations, differential
	equations and values of the coefficients generalize automatically.
	Therefore, we have the definitions 1. and the theorems 2. and 3.

	\sssec{Definition.}
	The {\em polynomials of Chebyshev of the first} $(T_n)$ {\em and of
	 the second kind} $(U_n),$ {\em of Legendre} $(P_n),$
	{\em of Laguerre} $(L_n)$ and {\em of Hermite} $(H_n)$ are
	defined by the differential equations:
	\begin{enumerate}
	\setcounter{enumi}{-1}
	\item$      (1-I^{2)} D^2 T_n - I D T_n + n^2 T_n \equiv  0,$\\
	\hth$	T_n(0)\equiv  1, DT_n(0) \equiv  .$
	\item$      (1-I^{2)} D^2 U_n - 3 I D U_n + n(n+2) U_n \equiv  0,$\\
	\hth$	U_n(0) \equiv  1, DU_n(0) \equiv  .$
	\item$      (1-I^{2)} D^2 P_n - 2 I D P_n + n(n+1) P_n \equiv  0,$\\
	\hth$	P_n(0) \equiv  1, DP_n(0) \equiv  .$
	\item$      I D^2 L_n + (1-I) D L_n + n \equiv  0.$
	\item$      D^2 H_n - 2 I H_n + 2n \equiv  0.$
	\end{enumerate}

	\sssec{Theorem.}
	{\em If $X_{n,j}$ denotes the coefficient of $I^j$ in the polynomial
	$X_n,$ then}
	\begin{enumerate}
	\setcounter{enumi}{-1}
	\item$      T_{n,n-2j} \equiv  \frac{n}{2} 2^{(n-2j)} (-1)^j
		 \frac{(n-j-1)!}{j!(n-2j)!},$
	\item$      U_{n,n-2j} \equiv  \frac{n}{2} 2^{(n-2j)} (-1)^j
		  \frac{(n-j)!}{j! (n-2j)!}),$
	\item$      P_{n,n-2j} \equiv  2^{(-n)} (-1)^j
		\frac{(2n-2j)!}{j!(n-j)! (n-2j)!},$
	\item$      L_{n,j} \equiv  (-1)^j \frac{n!}{(n-m)! (m!)^2},$
	\item$      H_{n,n-2j} \equiv  n! 2^{(n-2j)} (-1)^j
		 \frac{1}{j!(n-2j)!},$
	\end{enumerate}

	See for instance Handbook of Mathematical functions, p. 775.

	\sssec{Theorem.}
	\begin{enumerate}
	\setcounter{enumi}{-1}
	\item$T_0  1, T_1 \equiv I, T_{n+1} \equiv 2(2I-1) T_n - T_{n-1},$
	\item$U_0 \equiv  1,$ $U_1 \equiv  2I,$
		$U_{n+1} \equiv  2(2I-1) U_n - U_{n-1},$
	\item$P_0 \equiv  1,$ $P_1 \equiv  I,$
		$(n+1)P_{n+1} \equiv  -(2n+1)I P_n - nP_{n-1},$
	\item$P_0 \equiv  1,$ $P_1 \equiv  I,$
		$(n+1)P_{n+1} \equiv  (2n+1-I) P_n - nP_{n-1},$
	\item$H_0 \equiv  1,$ $H_1 \equiv  2I,$
		$H_{n+1} \equiv  2I H_n - 2n H_{n-1}.$
	\end{enumerate}

	See for instance Handbook of Mathematical functions, p. 782.

	\sssec{Theorem (T).}
	\hth$	T_{i+2pk,j} = -T_{i+pk,j} = T_{i,j}, j<p.$

	Proof:\\
	\hth$	T_{p+i,j} \equiv  (-1)^{\frac{p+i-j}{2}} 2^j 
	\frac{\frac{p+i+j}{2}-1)! \frac{p+i}{2}}
			{(\frac{p+i-j}{2})!  j! }$\\
	\hti{12}$		\equiv  (-1)^{\frac{p+i-j}{2}}
		 (-1)^{\frac{p-i-j}{2}} (-1)^{\frac{p-i+j+i}{2}} 2^j
				\frac{(\frac{p-i+j-2}{2})! \frac{p-i}{2}}
				{(\frac{p-i-j}{2})!  j! }$\\
	\hti{12}$		\equiv  (-1)^{\frac{p-i-j}{2}} 2^j
		\frac{\frac{p-i+j-2}{2})! \frac{p-i}{2}}
				{  (\frac{p-i-j}{2})!  j! }$\\
	\hti{12}$		\equiv  T_{p-i,j}.$

	\sssec{Theorem (U).}
	\begin{enumerate}
	\setcounter{enumi}{-1}
	\item$	U_{i+2pk,j} = -U_{i+pk,j} = U_{i,j},$ $j<p.$
	\item$	U_{p-1+i,j} \equiv  (-1)^{\frac{p-1+i-j}{2}} 
		\frac{(\frac{p-1+i+j}{2})! 2^j}
				{ (\frac{p-1+i-j}{2})!  j!}$
	\item$	U_{p-1+i,j} \equiv  (-1)^{\frac{p-1+i-j}{2}} 
		(-1)^{\frac{p-i-j-1}{2}}$
	\item$			(-1)^{\frac{p-i+j-1}{2}} 
		\frac{(\frac{p-i+j-1}{2})! 2^j}
				{ (\frac{p-1+i-j}{2})!  j!}$
	\item$	U_{p-1-i,j} \equiv  (-1)^{\frac{p-1+i-j}{2}} 2^j
		 \frac{(\frac{p-1-i+j}{2})!}
				{(\frac{p-1-i-j}{2})!  j!}
				\equiv  U_{p-1-i,j}.$
	\end{enumerate}

	\sssec{Theorem (Le).}
	\vspace{-18pt}\hspace{100pt}\footnote{24.11.83}\\[8pt]
	\hth$P_{p-1-n} = P_n,$ $n<p.$

	Proof:  The polynomials can be defined by the recurrence relations,\\
	\hth$	P_0 = 1,$ $P_1 = I,$
		$(n+1)P_{n+1} = (2n+1) I P_n - n P_{n-1}, n < p-1.$\\
	The last equation is valid for $n+1 = p$ and therefore $P_p$ can be
	considered as 0 as far as the proof of the theorem is concerned.  They
	satisfy the Rodrigues' formula\\
	\hth$	P_n = \frac{1}{2^n n!}D^n (I^{2-1)^n}$\\
	Therefore\\
	\hth$	P_{p-1} = \frac{1}{ (2^{p-1} ( \frac{p-1}{2}) ! )^2}
		D^{(p-1)} (-I^{2)^(\frac{p-1}{2})}$\\
		
	\hth$		= (-1)^{\frac{p-1}{2}}  \frac{(p-1)!}
		 { 2^{p-1} ( \frac{p-1}{2} ! )^2},$\\
	but $2^{(p-1)} = 1,$ $(p-1)! = -1,$
	and\\
	\hth$	( \frac{p-1}{2} ! )^2 = (-1)^{(\frac{p-1}{2})} (p-1)!$\\
	because $\frac{p-i}{2} = - \frac{p+i}{2},$
	hence\\
	\hth$	P_{p-1} = 1.$\\
	By convention we can write $P_p = P_{-1} = 0,$ if we replace in the
	recurrence relation $n$ by $p-n-1$ we obtain\\
	\hth$	(n+1)P_{p-n-2} = (2n+1)P_{p-n-1} - n P_{p-n},$\\
	and therefore, by induction,\\
	\hth$	P_{p-1-n} = P_n.$

	\sssec{Example.}
	$p = 11,$ see orthog, 120.\\
	\hth$	P_0 = 1,$\\
	\hth$	P_1 =   I,$\\
	\hth$	P_2 = 5 -4 I^2,$\\
	\hth$	P_3 =  4I -3I^3,$\\
	\hth$	P_4 =-1 -  I^2  + 3 I^4,$\\
	\hth$	P_5 = -5I + 5I^3 + I^5,$\\
	\hth$	P_6 =-1 -  I^2  + 3 I^4,$\\
	\hth$	P_7 =  4I -3I^{}3$\\
	\hth$	P_8 = 5 - 4I^2,$\\
	\hth$	P_9 =   I,$\\
	\hth$	P_{10} = 1.$

	For $p= 13,$\\
	\hth$	P_0 = P_{12} = 1,$\\
	\hth$	P_1 = P_{11} =   I,$\\
	\hth$	P_2 = P_{10} = 6 -5 I^2,$\\
	\hth$	P_3 = P_{ 9} =  5I -4I^3,$\\
	\hth$	P_4 = P_{ 8} = 2 + 6 I^2  + 6 I^4,$\\
	\hth$	P_5 = P_{ 7} = -3I + I^3 + 3I^5,$\\
	\hth$	P_6 =	-6 -4 I^2  -  I^4 -I^6,$

	\sssec{Theorem (La).}
	{\em }

	\sssec{Theorem (H).}
	{\em }

	\sssec{Definition.}
	The {\em scaled Hermite polynomials} are defined by
	\begin{enumerate}
	\setcounter{enumi}{-1}
	\item$      H_0 = 1,$
	\item$      H_1 = I,$
	\item$      [\frac{n}{2}] H_n = a_n H_{n-1} - \frac{n-1}{2} H_{n-2},$\\
	\hth where $a_n = 1$ if $n$ is even and $a_n = [n\frac{1}{2}],$
		the largest integer in $\frac{n}{2}$ if $n$ is odd.
	\end{enumerate}

	\sssec{Example.}
	\hth$	H_2 = -\frac{1}{2} + I^2,$\\
	\hth$	H_3 = -\frac{3}{2} I + I^3,$\\
	\hth$	H_4 = \frac{3}{8} - \frac{3}{2} I^2 + \frac{1}{2} I^4,$\\
	\hth$	H_5 = \frac{15}{8} I - \frac{5}{2} I^3 + \frac{1}{2} I^5,$\\
	\hth$H_6 = -\frac{5}{16} + \frac{15}{8} I^2 - \frac{5}{4} I^4
		+ \frac{1}{6} I^6,$\\
	\hth$H_7 = -\frac{35}{16} I + \frac{35}{8} I^3 - \frac{7}{4} I^5
		+ \frac{1}{6} I^7.$

	\sssec{Lemma.}
	{\em Modulo} $p,$ $p>2,$
	\begin{enumerate}
	\setcounter{enumi}{-1}
	\item$      (p-1)! \equiv  -1.$
	\item$      (p-1-i)! \equiv  (-1)^{(i+1)} \frac{1}{i!},$ $0 \leq i < p.$
	\item$\ve{p-1-i}{j} \equiv  (-1)^j \ve{i+j}{j},$ $0 \leq  i,j,$
		$ i+j < p.$
	\item$\ve{kp+i}{j} \equiv \ve{i}{j},$ $j < p.$
	\item$      (p-2-i)!!  i!!  \equiv  (-1)^{k\frac{1}{2}} (p-1-k-i)!!
		  (k+i-1)!! $\\
	\hth$			0 \leq  i < p-1,$ $0 < k+i < p.$
	\end{enumerate}

	Proof:  0. is the well known Theorem of Wilson.
	1, can be considered as a generalization.\\
	\hth$	(p-1-i) \equiv  (-1)^i  (p-1) \ldots (i+1)$\\
	\hth$		\equiv  (-1)^i \frac{(p-1)!}{ i!} $\\
	\hth$		\equiv  (-1)^{i+1} \frac{1}{i!}.$\\
	For 2.  $\frac{(p-1-i)!}{ (p-1-i-j)! j!}$\\
	\hti{12}$\equiv  (-1)^{(i+1)}\frac{ (i+j)!}{ (-1)^{i+j+1} i! j!}
			\equiv  (-1)^j \ve{i+j}{j}.$

	\sssec{Lemma.}
	{\em Modulo $p,$ $p>2,$}
	\begin{enumerate}
	\setcounter{enumi}{-1}
	\item$      ( (p-2)!!  )^2 \equiv  (-1)^{(\frac{p-1}{2}})$
	\item$      (p-1)!!  (p-2)!!  \equiv  -1.$
	\item 0.$    (p-2-i)!!  i!!   (-1)^s (p-2)!! ,$\\
	\hth {\em where $s = i\frac{1}{2} w \equiv  n\: i$ is even\\
	\hth and   $s = \frac{p-2-i}{2}$ when $i$ is odd.
		1.   or where $s = [ \frac{[\frac{p}{2}] + 1 + i}{2}]
		+ [\frac{p+1}{4}].$}
	\end{enumerate}

	0, and 1 are well known and given for completeness.
	2, if $i$ is even,\\
	\hth$	(p-2-i)!!  i!!   (p-i)!!  (i-2)!!  (i\frac{1}{p-i)}$ or $-1)$\\
	\hth$			 (-1)^{(i\frac{1}{2})} (p-2)!!  0!! .$\\
	   if $i$ is odd,\\
	\hth$(p-2-i)!! i!! (p-4-i)!! (i+2)!! (\frac{p-2-i}{i+2}$ or $-1)$\\
	\hti{12}$	\equiv (-1)^{(\frac{p-2-i}{2}))} (0)!!  p-2!! .$\\
	2.1, can be verified by choosing $p = 1,3,5,7$ and $i = 0,1,2,3,4.$

	\sssec{Theorem (La).}
	\begin{enumerate}
	\setcounter{enumi}{-1}
	\item$L_{p-1-i,j} \equiv  (-1)^j L_{i+j,j},$ $0 \leq  i,j,$ $i+j < p.$
	\end{enumerate}

	Proof:\\
	For 0,$ L_{n,j} = (-1)^j \ve{n}{j} \frac{1}{j!}$.
	  See for instance, Handbook p.775.

	\sssec{Lemma.}
	\vspace{-18pt}\hspace{80pt}\footnote{10.12.83}\\[8pt]
	\hth$\sum ((2j-1)!! (2k)!! \frac{1}{ (2j)}!! (2k-1)!! ))^2
		 j!\frac{1}{k!(j-k)}!)$\\
	\hth$	 = (-1)^([p\frac{1}{2}]-k),$ $j = k$ to $[p\frac{1}{2}].$

	This is needed for g761, \ldots. 
	Not yet proven.\\
	The expression which is summed in the first member can be replaced by\\
	\hth$	\frac{ (\frac{(2j)!}{(2k!)} )^2 k! }{  j! (j-k)! 2^{2(j-k)}}.$




	\ssec{Trigonometric Functions.}
	\setcounter{subsubsection}{-1}
	\sssec{Introduction.}
.       Connection with $p$-adic fields.
	In $p$-adic fields, introduced by Kurt Hensel,
	trigonometric functions are defined.  The connection between these
	and those obtained in finite fields has to be explored.
	To that effect, 2 programs have been written, padic.bas and sin.bas.
	The first program obtains the functions sin and cos for arguments
	which are congruent to 0 modulo $p.$
	For instance to $7^4$

	\sssec{Example.}
	$\begin{array}{cclll}
	&	    x&   sin(x)&  cos(x)\\
	\cline{1-4}
	&\frac{1}{2}&     0.4333&  0.4343& 1.0605\\
	&	0.1&     0.1011&  1.0331\\
	&	0.2&     0.2012&  1.0562\\
	&	0.3&     0.3062&  1.0622\\
	&	0.4&     0.4013&  1.0655\\
	&	0.5&     0.5062&  1.0516\\
	&	0.6&     0.6061&  1.0344\\
	&	0.01&    0.0100&  1.0003
	\end{array}$

	\sssec{Example.}
	In base 7 and $7^4,$ we have for the elliptic case\\
	$\begin{array}{rccccccc}
		y&\vline& sin(y)&cos(y)&sin(y)&cos(y)\\
	\cline{1-8}
		0&\vline& 0& 1& 0.000& 1.000& 0.000& 1.000\\
		2&\vline& 2& 2& 2.126& 2.406& 2.653& 2.653\\
		4&\vline& 1& 0& 1.053& 0.514& 1.000& 0.000\\
		6&\vline& 2& 5& 2.054& 5.143& 2.653& 5.013\\
		8&\vline& 0& 6& 0.332& 6.606& 0.000& 6.000\\
		10&\vline&5& 5& 5.444& 5.352&$\ldots$ ..\\
		12&\vline&6& 0& 6.631& 0.614\\
		14&\vline&5& 2& 5.030& 2.601\\
		16&\vline&0& 1& 0.101& 1.033
	\end{array}$

	If we observe that\\
	\hth$	0 .332\hti{2}= - 0 .434$ and $ 6 .606 = -1.060,$\\
	we get the first clue for the relation between the functions in the
	$p$-adic field and in base $7^4.$
	\newpage
	For $13^n,\\$ 
	$x\hti{16}sin(x)\hti{32}   cos(x)$\\
	$\begin{array}{rrrcrrrrrrrrcrrrrrrrr}
	\cline{1-21}
	0.&0&&\vline&\:\:0.&0&0&0&0&0&0&0&\:\:&1.&0&0&0&0&0&0&0\\
	0.&1&&\vline&0.&1&0&2&2&11&9&10&&1.&0&6&61&2&5&4&0\\
	0.&2&&\vline&0.&2&0&3&4&6&3&4&&1.&0&11&12&4&4&12&8\\
	0.&3&&\vline&0.&3&0&2&6&9&3&6&&1.&0&2&6&11&1&0&3\\
	0.&4&&\vline&0.&4&0&11&7&7&7&10&&1.&0&5&12&1&5&9&0\\
	0.&5&&\vline&0.&5&0&3&9&3&12&5&&1.&0&7&5&12&7&7&2\\
	0.&6&&\vline&0.&6&0&3&10&4&2&0&&1.&0&8&11&1&4&8&10\\
	0.&7&&\vline&0.&7&0&10&10&9&1&1&&1.&0&8&4&8&0&8&8\\
	0.&\frac{1}{2}&&\vline&0.&7&6&3&&&&	&&1.&0&8&11	\\
	0.&8&&\vline&0.&8&0&10&10&10&8&5&&1.&0&7&10&5&4&0&7\\
	0.&9&&\vline&0.&9&0&2&10&6&9&8&&1.&0&5&3&8&9&1&10\\
	0.&10&&\vline&0.&10&0&11&8&4&3&5&&1.&0&2&9&4&10&11&8\\
	0.&11&&\vline&0.&11&0&10&6&7&3&10&&1.&0&11&1&11&0&11&12\\
	0.&12&&\vline&0.&12&0&11&3&2&11&4&&1.&0&6&7&5&1&10&3\\
	0.&0&1&\vline&0.&0&1&0&0&0&2&2&&1.&0&0&0&6&6&6&6\\
	&\frac{\pi}{6}&&\vline&7.&6&6&6&6&6&6&6&&2.&4&3&4&6&1&7&5\\
	&\frac{\pi}{3}&&\vline&2.&4&3&4&6&1&7&5&&7.&6&6&6&6&6&6&6\\
	&\frac{\pi}{2}&&\vline&1.&0&0&0&0&0&0&0&&0.&0&0&0&0&0&0&0
	\end{array}$

	\sssec{Example.}
	\begin{verbatim}	
p,first e,s%? 13,2,7
c% = 2
  1    7    2   7.  0    2.  0
  2    2    7   2.  0    7.  0
  3    1    0   1.  0    0.  0
  4    2    6   2.  0    6.  0
p,first e,s%? 169,2,137
c% = 41
  1  137   41 137.  0   41.  0
  2   80  150  80.  0  150.  0
  3    1   91   1.  0   91.  0
  4    2   45   2.  0   45.  0
  5  163   50 163.  0   50.  0
  6   13  168  13.  0  168.  0
  7   58   37  58.  0   37.  0
  8   11  162  11.  0  162.  0
  9  168   65 168.  0   65.  0
 10   76   98  76.  0   98.  0
 11  149   28 149.  0   28.  0
 12  143    1 143.  0    1.  0
 13   85   54  85.  0   54.  0
 14   67   33  67.  0   33.  0
 15    1  117   1.  0  117.  0
 16   15   97  15.  0   97.  0
 17   46   63  46.  0   63.  0
 18   39  168  39.  0  168.  0
 19  110   24 110.  0   24.  0
 20   24  110  24.  0  110.  0
 21  168   39 168.  0   39.  0
 22   63   46  63.  0   46.  0
 23   97   15  97.  0   15.  0
 24  117    1 117.  0    1.  0
 25   33   67  33.  0   67.  0
 26   54   85  54.  0   85.  0
 27    1  143   1.  0  143.  0
 28   28  149  28.  0  149.  0
 29   98   76  98.  0   76.  0
 30   65  168  65.  0  168.  0
 31  162   11 162.  0   11.  0
 32   37   58  37.  0   58.  0
 33  168   13 168.  0   13.  0
 34   50  163  50.  0  163.  0
 35   45    2  45.  0    2.  0
 36   91    1  91.  0    1.  0
 37  150   80 150.  0   80.  0
 38   41  137  41.  0  137.  0
 39    1    0   1.  0    0.  0
 40   41   32  41.  0   32.  0
	\end{verbatim}

	\sssec{Comment.}
	Using $s_1 = x,$ $s_{2i+1} = s_ix^2\frac{(2i-1)^2}{2i(2i+1)},$ \\
	and $arcsin(x) = s_1 + s_3 +  \ldots,$\\
	$arcsin(.7,6,6,6,6,6,6) = .7,6,9,12,9,6,5$\\
	and\\
	$sin(.7,6,9,12,9,6,5) = .7,6,6,6,6,6,6 = (.1,0,0,0,0,0,0)$\\
	modulo 169, at 13 we read  85 and 54 corresponding to\\
				.76   and .24\\
	hence $\frac{\pi}{6}$ correspond to $arcsin(.7,6,6) = .7,6,9$
	and $p$ to $.3,0,5,11,7,1,7$\\
	If we use the program padic.bas we can, given
	$sin (\alpha)$ and $cos(\alpha),$ obtained using 115a n2adic,
	determine $sin(i \alpha)$ and $cos (i \alpha)$.\\
	For $p = 13,$ using $sin(\alpha) \equiv 7 \umod{13},$
	we get $\pmod{p^3}$ all distinct values except\\
	$(0, 2, 5, 6, 7, 8, 11)\:169 \pm 1.\\$
	The non zero numbers in the parenthesis are the non residues.\\
	This indicates that the connection is tenuous.

	\sssec{Theorem.}
	\vspace{-18pt}\hspace{90pt}\footnote{11.10.82}\\[8pt]
	If $s(1) \equiv g \pmod{p^n},$ where $g$ is a primitive root of $p,$
	then \ldots\\ 
	Let $s1 = c2 = \frac{1}{2} = 6. 6 6 6  \ldots,$
	\hti{4}$c1 = s2 = \sqrt{\frac{3}{2}} = 2.12  \ldots,$\\
	given $\xi$, determine $s = sin(\xi),$ $c = cos(\xi),$ in
	the $p$-adic field, $(\xi \equiv 0 \pmod{p},$\\
	determine,\\
	$s(i) = sin(i \xi) cp(i) + cos(i \xi) sp(i),$
	$c(i) = cos(i \xi) cp(i) - sin(i \xi) sp(i).$\\
	($\xi$ was \{ control H $-$)?

	\ssec{Integration.}

	\sssec{Definition.}
	\vspace{-18pt}\hspace{90pt}\footnote{17.12.82}
	\begin{enumerate}
	\item$ Int_{2i}^{2j} \:sin = cos(2j) - cos(2i).$
	\item$	Mid_{2i}^{2j} \:sin
			= sin(2i+1) + sin (2i+3) +  \ldots  + sin(2j-1).$
	\item$	Trap_{2i}^{2j} \:sin
		= \frac{1}{2}sin(2) + sin (2i+2) + \ldots + sin(2j-2)
			+ \frac{1}{2}sin(2j) .$ 
	\item$Simpson_{2i}^{2j} \:sin
			= sin(2i) + 4 sin (2i+2) + 2 sin(2i+4) + 4 sin(2i+6)\\
	\hti{16}			+  \ldots + 4 sin(2j-2) + sin(2j).$
	\end{enumerate}

	\sssec{Theorem.}
	\hth$Int_{2i}^{2j} \:sin
			= -2 sin(1)\: Mid_{2i}^{2j} \:sin\\
	\hti{16}	= -2 tan(1)\: Trap_{2i}^{2j} \:sin\\
	\hti{16}	= -2\frac{sin(2)}{2+cos(2)}\: Simpson_{2i}^{2j} \:sin$.




	\section{P-adic field.}

	\ssec{Generalities.}

	\sssec{Notation.}
	Writing $x = x_0 + x_1 p + x_2 p^2 + x_3 p^3 +  \ldots$  in the form\\
	$x = x_0.x_1 x_2 x_3  \ldots$, we will also write\\
	\hth$x \equiv  x_0 \pmod{p},$ $x \equiv  x_0.x_1 \pmod{p^2},$ \ldots.\\
	We have, for instance,\\
	\hth 0. $\frac{1}{2}$ = 0. 7 6 6 6 6 6 6.

	\sssec{Example.}
	$p = 5,$ up to $p^6,$ \\
	\hth$	\frac{1}{2} = 3.22222,$ indeed, $2\:3 = 6 \equiv  1 \pmod{5},$\\
	\hth$	2\:(3 . 2) = 2 . 13 = 26 \equiv  1 \pmod{5^2},$  \ldots.\\
	\hth$	-1/2 = 2.22222,$ $-.1/2 = .22222,$ $ -.01/2 = .02222.$

	\sssec{Definition.}
	In the $p$-adic field, the {\em exponetial, logarithmic and
	trigonometric functions} are defined by:\\
	$exp(x) := 1 + x + \frac{1}{2!}x^2 + \frac{1}{3!}x^3 +  \ldots,$
		for $| x| < \ldots$.\\ 
	$log(x) := x - \frac{1}{2}x^2 + \frac{1}{3}x^3 +  \ldots,$
	 for $| x| < $.\\
	$sin(x) := x - \frac{1}{3!}x^3 + \frac{1}{5!}x^5 -  \ldots$,
	  for  $| x|  \leq p^{-1}$\\
	$cos(x) := 1 - \frac{1}{2!}x^2 + \frac{1}{4!}x^4 -  \ldots,$
	 for  $| x|  \leq p^{-1}.$

	\sssec{Example.}
	$p = 5,$ $x = 0.1,$
	$\begin{array}{lllclll}
	1	&= 1.00000& 00\\
	x	&= 0.10000& 00&\vline&	 x		&= 0.10000& 00\\
	x^2/2	&= 0.03222& 22&\vline&-x^2\frac{1}{2}	&= 0.02222& 22\\
	x^3/1 . 1	&= 0.00140& 40&\vline& x^3\frac{1}{3} 	&= 0.00231& 31\\
	x^4/4 . 4&= 0.00044& 34&\vline&-x^4\frac{1}{4}		&= 0.00011& 11\\
	x^5/.44	&= 0.00044& 34&\vline& x^5/.1 	&= 0.00010& 00\\
	x^6/.4301&= 0.00004& 03&\vline&-x^6/1 . 1	&= 0.00000& 40\\
	x^7/.31031&= 0.00000& 24&\vline& x^7/2 . 1 	&= 0.00000& 03\\
	\cline{6-7}
	x^8/.422422&= 0.00000& 04\\
	\cline{2-3}
	exp(0 . 1)&= 1.13341& 24&\vline&log(1 .1)	&= 0 .12320& 14\\
	\cline{1-7}
	 x	&= 0.10000& 00&\vline& 1		&= 1.00000& 00\\
	-x^3/1 .1&= 0.00404& 04&\vline&    -x^2\frac{1}{2} &= 0.02222& 22\\
	 x^5/.44&= 0.00044& 34&\vline& x^4/4 	&= 0.00044& 34\\
	-x^7/.31031&= 0.00000& 30&\vline&-x^6/1 .1	&= 0 .00001& 41\\
	\cline{2-3}
		&&		&\vline& x^8/.422422 &= 0.00000& 04\\
	\cline{6-7}
	sin(0 .1)&= 0 .10443& 24&\vline&cos(1 .1)	&= 1 .02213& 03\\
	\end{array}$

	\sssec{Example.}
	For $p = 13,$\\
	$x\hti{16}sin(x)\hti{32}   cos(x)$\\
	$\begin{array}{rrrrcrrrrrrrrcrrrrrrrr}
	\cline{1-22}
	0.&1&&&\vline&\:\:0.&1&0&2&2&11&9&10,&\:\:&1.&0&6&6&12&5&4&0\\
	0.&2&&&\vline&0.&2&0&3&4&6&3&4,&&1.&0&11&12&4&4&12&8\\
	0.&3&&&\vline&0.&3&0&2&6&9&3&6,&&1.&0&2&6&11&1&0&3\\
	0.&4&&&\vline&0.&4&0&11&7&7&7&10,&&1.&0&5&12&1&5&9&0\\
	0.&5&&&\vline&0.&5&0&3&9&3&12&5,&&1.&0&7&5&12&7&7&2\\
	0.&6&&&\vline&0.&6&0&3&10&4&2&0,&&1.&0&8&11&1&4&8&10\\
	0.&7&&&\vline&0.&7&0&10&10&9&1&1,&&1.&0&8&4&8&0&8&8\\
	0.&\frac{1}{2}&&&\vline&0.&7&6&3&0&5&5&7,&&1.&0&8&1&10&0&1&0\\
	0.&8&&&\vline&0.&8&0&10&10&10&8&5,&&1.&0&7&10&5&4&0&7\\
	0.&9&&&\vline&0.&9&0&2&10&6&9&8,&&1.&0&5&3&8&9&1&10\\
	0.&10&&&\vline&0.&10&0&11&8&4&3&5,&&1.&0&2&9&4&10&11&8\\
	0.&11&&&\vline&0.&11&0&10&6&7&3&10,&&1.&0&11&1&11&0&11&12\\
	0.&12&&&\vline&0.&12&0&11&3&2&11&4,&&1.&0&6&7&5&1&10&3\\
	0.&0&1&&\vline&0.&0&1&0&0&0&2&2,&&1.&0&0&0&6&6&6&6\\
	0.&11&1&&\vline&0.&11&1&10&4&10&9&1,&&1.&0&11&3&3&10&11&1\\
	0.&0&2&&\vline&0.&0&2&0&0&0&3&4,&&1.&0&0&0&11&12&12&12\\
	0.&10&2&&\vline&0.&10&2&11&12&2&1&8,&&1.&0&2&2&1&1&3&3\\
	0.&0&3&&\vline&0.&0&3&0&0&0&2&6,&&1.&0&0&0&2&6&6&6\\
	\end{array}$

	\sssec{Definition.}
	The {\em Chebyshev polynomials} are defined by the recurrence
	relation\\
	\hth$T_{i+1}(x) := -T_{i-1}(x) + 2x T_i(x),$ $i = 1, 2,$ \ldots, with\\
	\hth$T_0(x) := 1,$ $T_1(x) := x.$

	\sssec{Definition.}
	For $p \equiv 1 \pmod{4},$ a root $c1$ of $T_i$ is called a
	{\em primitive root modulo $p$}, if all the roots of $T_i,$ \\
	\hth	$c_1,$	$c_3,$  \ldots, $c_{2i-1}$\\
	can be obtained from it using the addition formulas,\\
	\hth$c_1 := c1,$ $c_3 := 4 c1^3 - 3 c1,$
	\hth$c_{2i+1} := - c_{2i-3} + 2 c_{2i-1} (2c1^2-1),$ $i = 2,$ \ldots,
	$i-1.$

	\sssec{Example.}
	For $p = 13,$ the roots of $T_3 = 4x^3 - 3x$ are 0, 2 and $-2.$
	2 and $-2$ are primitive roots.\\
	If $c1 = c_1 = 2$ then $c_3 = 0,$ $c_5 = -2.$\\
	For $p = 17,$ $T_4(x) = 8 x^4 - 8 x^2 + 1,$ which has the primitive
	roots $4, -4, 6, -6.$

	\sssec{Notation.}
	For $p \equiv 1 \pmod{4},$ the roots of $T_{\frac{p-1}{4}}$ will be
	denoted,\\
	\hth$cos(\alpha),$ $cos(3\alpha),$ \ldots, $cos(\frac{p-3}{2}\alpha).$\\
	with $\alpha = \frac{\pi}{\frac{p-1}{2}}.$\\
	We will also define\\
	\hth$	cos(0 \alpha) := 1,$\\
	\hth$	cos(2k\alpha) := -cos((2k-2)\alpha) + 2 c1 cos((2k-1)\alpha),$\\
	\hth$		k = 1,$  \ldots, $\frac{p-1}{4}.$\\
	\hth$	sin(k \alpha) = cos( (\frac{p-1}{4}-k) \alpha).$

	\sssec{Example.}
	For $p = 13,$ $\alpha = \frac{\pi}{6},$ $\delta^2 = 2,$\\
	\hth(the lines $k = \frac{1}{2},$ $\frac{3}{2},$  \ldots  will be
	explained in 1.1, 1.2),\\
	$k\hti{16}sin(k \alpha)\hti{32}   cos(k \alpha)$\\
	$\begin{array}{rcrrrrrrrrrcrrrrrrrrr}
	\cline{1-21}
	0&\vline&\:\:0.&0&0&0&0&0&0&0&\:\:&1.&0&0&0&0&0&0&0\\
	\frac{1}{2}&\vline&9.&7&1&1&0&9&12&0&\delta&&2.&1&8&7&6&2&6&7&\delta\\
	1&\vline&7.&6&6&6&6&6&6&6&&2.&4&3&4&6&1&7&4\\
	\frac{3}{2}&\vline&6.&6&6&6&6&6&6&6&\delta&&6.&6&6&6&6&6&6&6&\delta\\
	2&\vline&2.&4&3&4&6&1&7&4&&7.&6&6&6&6&6&6&6\\
	\frac{5}{2}&\vline&2.&1&8&7&6&2&6&7&\delta&&9.&7&1&1&0&9&12&0&\delta\\
	3&\vline&1.&0&0&0&0&0&0&0&&0.&0&0&0&0&0&0&0\\
	\frac{7}{2}&\vline&2.&1&8&7&6&2&6&7&\delta&&4.&5&11&11&12&3&0&12&
	\delta\\
	4&\vline&2.&4&3&4&6&1&7&4&&6.&6&6&6&6&6&6&6\\
	\frac{9}{2}&\vline&6.&6&6&6&6&6&6&6&\delta&&7.&6&6&6&6&6&6&6&\delta\\
	5&\vline&7.&6&6&6&6&6&6&6\\&\vline&&11.&8&9&8&6&11&5&8\\
	\frac{11}{2}&\vline&9.&7&1&1&0&9&12&0&\delta\\&11.&11&4&5&6&10&6&5&
	\delta\\
	6&\vline&0.&0&0&0&0&0&0&0\\&\vline&&12.&12&12&12&12&12&12&12\\
	\frac{13}{2}&\vline&4.&5&11&11&12&3&0&12&\delta\\&11.&11&4&5&6&10&6&5&
	\delta\\
	7&\vline&6.&6&6&6&6&6&6&6\\&\vline&&11.&8&9&8&6&11&5&8\\
	\frac{15}{2}&\vline&7.&6&6&6&6&6&6&6&\delta&&7.&6&6&6&6&6&6&6&\delta\\
	8&\vline&11.&8&9&8&6&11&5&8&&6.&6&6&6&6&6&6&6\\
	\frac{17}{2}&\vline&11.&11&4&5&6&10&6&5&\delta&&4.&5&11&11&12&3&0&12&
	\delta\\
	9&\vline&12.&12&12&12&12&12&12&12&&0.&0&0&0&0&0&0&0\\
	\frac{19}{2}&11.&11&4&5&6&10&6&5&\delta&&9.&7&1&1&0&9&12&0&\delta\\
	10&\vline&11.&8&9&8&6&11&5&8&&7.&6&6&6&6&6&6&6\\
	\frac{21}{2}&7.&6&6&6&6&6&6&6&\delta&&6.&6&6&6&6&6&6&6&\delta\\
	11&\vline&6.&6&6&6&6&6&6&6&&2.&4&3&4&6&1&7&4\\
	\frac{23}{2}&4.&5&11&11&12&3&0&12&\delta&&2.&1&8&7&6&2&6&7&\delta\\
	12&\vline&0.&0&0&0&0&0&0&0&&1.&0&0&0&0&0&0&0
	\end{array}$

	In this particular case, $sin(\alpha) = \frac{1}{2},$
	$cos(\alpha) = \frac{\sqrt{3}}{2},$ are sufficient to obtain the
	entries 0, 1, 2,  \ldots, 12, in the table.\\
	We have therefore a first method of obtaining tables of trigonometric
	functions in a finite field\footnote{12.10.82}.
	We choose $x,$ such that $| x|  = p^{-1},$
	therefore $|  k x |  \leq p^{-1}.$  ( If $| x|  < p^{-1},$ primitivity
	is not insured. See 1.3.)\\
	We compute $sin(k\: x)$ and $cos(i\: x)$ by Maclaurin series, (see 0.1.)
	and use the addition formulas\\
	$sin(k(\alpha + x)) = sin(kp \alpha) cos(kx) + cos(kp \alpha) sin(kx),$
	$cos(k(\alpha + x)) = cos(kp \alpha) cos(kx) - sin(kp \alpha) sin(kx),$
	where $kp$ is $k \pmod{p}.$

	\sssec{Example.}
	For $p = 13,$ and $x$ = 0. 1,\\
	\hth(The lines $k = \frac{1}{2}, \frac{3}{2},  \ldots  will be
	explained in 1.3).$\\
	$k\hti{16}sin(k (\alpha+x))\hti{32}   cos(k (\alpha+x))$\\
	$\begin{array}{rcrrrrrrrrrcrrrrrrrrrll}
	\cline{1-21}
	\frac{1}{2}&&\vline&\:\:9.&8&2&11&7&1&5&9&\delta&&\:\:2.&3&7&12&
12&12&11&1&\delta\\
	1&&\vline&7.&8&0&4&3&7&4&5&&&2.&10&8&7&2&10&10&6\\
	\frac{3}{2}&&\vline&6.&2&6&10&8&11&2&2&\delta&&6.&10&12&0&3&3&5&7&
\delta\\
	2&&\vline&2.&5&12&3&5&5&1&0&&&7.&2&10&6&12&2&11&4\\
	\frac{5}{2}&&\vline&2.&4&11&4&4&7&6&4&\delta&&9.&2&11&4&4&2&10&1&
\delta\\
	3&&\vline&1.&0&2&6&11&1&0&3&&&0.&10&12&10&6&3&9&6\\
	\frac{7}{2}&&\vline&2.&2&11&4&2&8&8&2&\delta&&4.&11&8&8&0&12&4&5&
\delta\\
	4&&\vline&2.&2&0&11&5&1&11&10&&&6.&11&6&10&6&9&12&0\\
	\frac{9}{2}&&\vline&6.&5&4&1&7&4&9&2&\delta&&7.&5&2&2&7&5&9&0&\delta\\
	5&&\vline&7.&9&8&12&12&9&1&1&&&11.&12&1&7&5&3&7&5\\
	&\frac{11}{2}&&\vline&9.&9&3&0&3&0&5&5&\delta&&11.&7&2&1&8&7&3&3&
\delta\\
	6&&\vline&0.&7&12&9&2&8&10&12&&&12.&12&4&1&11&8&4&2\\
	&0.&1&\frac{1}{2}&4.&5&10&4&8&7&1&10&\delta&&11.&11&2&9&10&12&1&2&
\delta\\
	&7&\vline&6.&5&12&12&1&12&6&5&&&11.&5&0&5&11&0&7&5\\
	&2.&1&\frac{1}{2}&7.&0&0&4&6&1&11&11&\delta&&7.&12&5&8&7&11&2&6&
\delta\\
	&8&\vline&11.&4&8&6&4&5&11&6&&&6.&9&3&2&8&6&8&1\\
	&4.&1&\frac{1}{2}&11.&6&8&1&3&6&3&2&\delta&&4.&9&6&0&7&10&12&2&
\delta\\
	&9&\vline&12.&12&7&9&4&3&11&2&&&0.&9&0&2&10&6&9&8\\
	&6.&1&\frac{1}{2}&11.&12&8&7&3&11&6&6&\delta&&9.&0&11&7&11&0&9&2&
\delta\\
	\end{array}$

	$\begin{array}{rcrrrrrrrrrcrrrrrrrrrll}
	&10&\vline&11.&0&6&7&4&10&9&1&&&7.&0&10&7&2&3&6&0\\
	&8.&1&\frac{1}{2}&7.&4&5&9&12&6&4&3&\delta&&6.&4&0&12&2&12&2&12&
\delta\\
	&11&\vline&6.&2&1&10&3&11&10&1&&&2.&3&6&11&0&10&4&12\\
	&10.&1&\frac{1}{2}&4.&2&7&8&9&4&10&4&\delta&&2.&7&6&7&7&12&9&11\\
	&12&\vline&0.&12&0&11&3&2&11&4&&&1.&0&6&7&5&1&10&3\\
	&11.&1&\frac{1}{2}&9.&6&3&1&9&11&1&4&\delta&&2.&12&5&5&0&2&9&11\\
	&0.&1&\vline&7.&6&8&10&12&0&7&10&&&2.&4&9&10&11&5&4&1\\
	&1.&2&\frac{1}{2}&6.&9&11&0&4&6&7&5&\delta&&6.&3&6&4&12&2&2&1&\delta\\
	&1.&1&\vline&2.&11&2&1&9&11&6&8&&&7.&4&3&4&3&6&1&11\\
	&3.&2&\frac{1}{2}&2.&8&3&4&5&2&3&6&\delta&&9.&4&2&4&12&3&0&9&\delta\\
	&2.&1&\vline&1.&0&11&10&10&7&2&5&&&0.&11&11&9&10&7&2&10\\
	&5.&2&\frac{1}{2}&2.&11&2&0&10&2&0&5&\delta&&4.&0&7&0&0&12&5&9&\delta\\
	&3.&1&\vline&2.&9&6&4&3&12&4&7&&&6.&0&4&6&9&10&2&9\\
	&7.&2&\frac{1}{2}&6.&11&2&0&2&12&9&0&\delta&&7.&11&4&4&1&12&8&3&
\delta\\
	&4.&1&\vline&7.&11&9&3&11&12&4&10&&&11.&6&5&5&1&10&10&2\\
	&9.&2&\frac{1}{2}&9.&11&8&0&12&9&2&12&\delta&&11.&3&4&7&9&10&3&7&
\delta\\
	&5.&1&\vline&0.&8&11&9&9&12&5&11&&&12.&12&5&12&7&1&5&7\\
	&11.&2&\frac{1}{2}&4.&7&7&12&11&6&6&4&\delta&&11.&2&5&1&12&10&10&3&
\delta\\
	&6.&1&\vline&6.&7&1&6&3&1&8&4&&&11.&11&0&6&6&2&8&6\\
	&0.&3&\frac{1}{2}&7.&6&10&9&6&10&5&4&\delta&&7.&6&2&3&0&4&2&10&
\delta
	\end{array}$

	A second method to obtain trigonometric tables in finite fields is
	to start with $sin(\alpha)$ such that $x_0.x_1$ is a primitive root for
	$p^2$ and $cos(\alpha) = \sqrt{1-sin^2(\alpha)}$ and use the addition
	formulas to obtain in succession\\
	$sin(k \alpha)$ and $cos(k \alpha)$ for $k = 2, 3,$ \ldots.  For
	instance, we have the following

	\sssec{Example.}
	For $q = 13^4$ and $sin(\alpha)$ = 7. 8 0 4, then\\
	$cos(\alpha) = 2$.10 8 7, the period is $12 . 13^3$ and\\
	$ k$       $sin(k \alpha)$     $cos(k \alpha)$\\
	$\begin{array}{llcllllcllll}
	1&&\vline&\:\:7.&8&0&4&&\:\:2.&10&8&7\\
	2&&\vline&2.&5&12&3&&7.&2&10&6\\
	3&&\vline&1.&0&2&6&&0.&10&12&10\\
	4&&\vline&2.&2&0&11&&6.&11&6&10\\
	5&&\vline&7.&9&8&12&&11.&12&1&7\\
	6&&\vline&0.&7&12&9&&12.&12&4&1\\
	7&&\vline&6.&5&12&12&&11.&5&0&5\\
	8&&\vline&11.&4&8&6&&6.&9&3&2\\
	9&&\vline&12.&12&7&9&&0.&9&0&2\\
	10&&\vline&11.&0&6&7&&7.&0&10&7\\
	11&&\vline&6.&2&1&10&&2.&3&6&11\\
	12&&\vline&0.&12&0&11&&1.&0&6&7\\
	0.&1&\vline&7.&6&8&10&&2.&4&9&10\\
	11.&1&\vline&0.&11&1&10&&1.&0&11&3\\
	0.&2&\vline&2.&4&4&4&&7.&6&2&11\\
	10.&2&\vline&0.&10&2&11&&1.&0&2&2\\
	0.&3&\vline&1.&0&0&0&&0.&0&10&12\\
	\end{array}$

	\sssec{Theorem.}
	{\em 
	If $sinh(x) = x + \frac{1}{3!}x^3 + \frac{1}{5!}x^5 +  \ldots,$
		$| x| \leq p^{-1},$\\
	and $cosh(x) = 1 + \frac{1}{2!}x^2 + \frac{1}{4!}x^4 +  \ldots,$
		$| x| \leq p^{-1},$\\
	then\\
	\hth$sin(x i) = i sinh(x),$ $cos(x i) = cosh(x),$}

	\sssec{Example.}
	With $p = 11,$\\
	$sin(0.\frac{1}{2} i)$ = 0. 6 5 810 3 8 0 $i,$
	$cos(0.\frac{1}{2} i)$ = 1. 0 7 9 5 7 0 2.\\
	$sin(0 .1 i)  = $0. 1 0 2 9 0 9 6 $i$, $cos(0 .1 i) = $1. 0 6 5 0 5 9 6.

	\sssec{Example.}
	For $p = 11$ and $x$ = . 1,
	Using 6.1.13 and\\
	$sin(\frac{\alpha}{2})$ = 9. 3 5 9 2 2 410  and
	$cos(\frac{\alpha}{2})$ = 5. 7 6 3 0 4 6 9 $i,$\\
	$sin(\alpha)$ = 2.10 5 5 6 5 6 6 $i$ and
	$cos(\alpha)$   = 4. 919 3 8 7 9 4\\
	of 1.6, we obtain\\
	$sin((\frac{\alpha+1}{2}) i)$ = 9. 6 910 9 2 1 0 $i,$
	$cos((\frac{\alpha+1}{2}) i)$ = 5. 8 610 0 0 7 2,
	$sin((\alpha+1) i)$     = 2. 3 5 7 8 5 5 2 $i,$
	$cos((\alpha+1) i)$     = 4. 0 1 210 6 810.\\
	The table can be computed from the first values the second are
	used here as a check:\\
	$k\hti{16}sin(k (\alpha+x)/2)\hti{32}   cos(k( \alpha+x)/2)$\\
	$\begin{array}{rrcrrrrrrrrrrcrrrrrrrrrr}
	\cline{1-23}
	0.&0&\vline&0.&0&0&0&0&0&0&0&&&1.&0&0&0&0&0&0&0\\
	1.&0&\vline&9.&6&9&10&9&2&1&0&&&5.&8&6&10&0&0&7&2&i\\
	2.&0&\vline&2.&3&5&7&8&5&5&2&i&&4.&0&1&2&10&6&8&10\\
	3.&0&\vline&4.&6&5&4&8&0&9&8&&&2.&4&2&1&0&3&3&1&i\\
	4.&0&\vline&5.&3&2&8&0&5&3&7&i&&9.&2&5&0&0&9&10&5\\
	5.&0&\vline&1.&0&10&9&8&0&6&4&&&0.&3&5&4&6&0&0&9&i\\
	6.&0&\vline&5.&2&1&8&5&2&4&9&i&&2.&0&4&10&4&9&10&0\\
	7.&0&\vline&4.&5&10&1&4&2&6&8&&&9.&8&6&7&7&3&4&1&i\\
	8.&0&\vline&2.&5&6&4&3&2&10&4&i&&7.&9&8&2&0&10&10&8\\
	9.&0&\vline&9.&9&4&10&3&10&1&9&&&6.&1&1&0&9&3&4&1&i\\
	10.&0&\vline&0.&6&10&2&7&2&7&10&i&&&10.&10&3&4&6&3&9&9\\
	0.&1&\vline&2.&7&2&5&3&7&0&2&&&6.&3&3&9&5&2&10&4&i\\
	1.&1&\vline&9.&9&0&9&0&0&5&0&i&&7.&0&10&5&5&7&2&7\\
	2.&1&\vline&7.&3&9&9&9&5&4&8&&&9.&4&4&6&2&4&1&7&i\\
	3.&1&\vline&6.&6&3&4&0&5&9&8&i&&2.&5&2&6&2&7&4&1\\
	4.&1&\vline&10.&10&8&9&6&0&8&7&&&0.&2&6&10&2&2&3&3&i\\
	5.&1&\vline&6.&9&6&7&1&2&5&7&i&&9.&7&2&7&3&2&10&8\\
	6.&1&\vline&7.&6&5&9&9&7&0&2&&&2.&0&1&5&8&0&0&1&i\\
	7.&1&\vline&9.&3&8&0&9&9&6&8&i&&4.&2&4&4&1&4&9&10\\
	8.&1&\vline&2.&9&5&8&7&1&7&6&&&5.&10&5&5&8&1&1&10&i\\
	9.&1&\vline&0.&10&0&9&7&9&9&0&i&&1.&0&6&4&6&8&8&6\\
	10.&1&\vline&9.&0&1&8&8&0&2&7&&&5.&6&1&3&4&9&9&5&i\\
	\end{array}$

	$\begin{array}{rrcrrrrrrrrrrcrrrrrrrrrr}
	0.&2&\vline&2.&10&9&3&7&3&9&5&i&&4.&9&1&3&5&1&9&6\\
	1.&2&\vline&4.&8&9&5&5&6&3&4&&&2.&8&5&8&0&0&2&6&i\\
	2.&2&\vline&5.&5&5&6&8&8&7&9&i&&9.&8&5&0&3&10&3&2\\
	3.&2&\vline&1.&0&8&5&3&9&2&5&&&0.&4&4&1&9&9&0&6&i\\
	4.&2&\vline&5.&0&0&2&7&6&10&0&i&&2.&6&7&10&3&4&0&6\\
	5.&2&\vline&4.&3&1&0&4&4&6&5&&&9.&1&7&10&8&2&10&0&i\\
	6.&2&\vline&2.&9&4&8&7&7&9&0&i&&7.&7&3&6&9&3&10&7\\
	7.&2&\vline&9.&4&10&6&3&3&0&7&&&6.&10&4&2&0&9&10&0&i\\
	8.&2&\vline&0.&7&9&3&9&1&3&4&i&&10.&10&2&6&8&4&8&8\\
	9.&2&\vline&2.&2&1&2&3&7&1&10&&&6.&5&0&10&10&2&1&5&i\\
	10.&2&\vline&9.&2&6&3&4&6&7&5&i&&7.&2&7&6&0&1&10&3\\
	0.&3&\vline&7.&1&3&0&7&10&10&1&&&9.&0&0&3&10&10&2&10&i\\
	1.&3&\vline&6.&4&3&0&9&3&9&2&i&&2.&10&2&3&7&6&3&0\\
	2.&3&\vline&10.&10&4&9&7&9&9&0&&&0.&1&7&7&1&6&1&9&i\\
	3.&3&\vline&6.&0&0&10&3&3&8&3&i&&9.&1&9&3&1&1&6&3\\
	4.&3&\vline&7.&8&1&5&1&4&6&8&&&2.&7&1&5&10&2&10&5&i\\
	5.&3&\vline&9.&10&8&5&1&5&10&4&i&&4.&4&0&1&0&4&5&2\\
	6.&3&\vline&2.&3&1&0&6&0&8&0&&&5.&1&10&0&0&4&9&0&i\\
	7.&3&\vline&0.&9&1&6&7&2&5&2&i&&1.&0&2&7&9&4&0&2
	\end{array}$

	\ssec{Extension to the half argument.}
	\setcounter{subsubsection}{-1}
	\sssec{Introduction.}
	The tables of trigonometric functions can be extended to the half
	arguments.  These are required for the angles in finite Euclidean
	geometry.

	\sssec{Theorem.}
	{\em If $g$ is a primitive root of $p,$ and $\delta^2$ = g,
	then $c1' = cos(\alpha\frac{1}{2}) \delta{-1}$ is a primitive root of\\
	\hth$	S'_{\frac{p-1}{2}} = T_{\frac{p-1}{2}} \circ (\delta I),$\\
	where $I$ is the identity function.\\
	Indeed, $T_{2n} = T_n \circ (2 I^2 -1).$  The other roots are
	denoted by $c2',$ $c3',$  \ldots.}

	\sssec{Example.}
	For $p = 13,$ g = 2,\\
	\hth$S'_6 = 256 I^6 -192 I^4 + 36 I^2 -1,$\\
	$c1' = cos(\alpha\frac{1}{2}) \delta{-1}$ = 2. 1 8 7 6 2 6 7, from
	which we derive the
	values in 0.6. for $i = \frac{1}{2},$ $\frac{3}{2},$  \ldots,
	$\frac{11}{2}.$ 

	\sssec{Comment.}
	The method given at the end of section 0.6. enables to
	complete the table of Example 0.7.  Alternately, if $g$ is a primitive
	root for $p^2,$ $p \equiv 1 \pmod{4},$ we know that $g$ is a primitive
	root for $p^e,$ $e = 3, 4, \ldots.$\\
	If $\delta^2 = g,$ $sin(\alpha/2)
	= \delta \sqrt{\frac{1-cos(\alpha)}{2g}}$
	and $cos(\alpha\frac{1}{2}) = \delta \sqrt{\frac{1+cos(\alpha)}{2g}}.$

	\sssec{Example.}
	For $p = 13,$\\
	$sin(\alpha\frac{1}{2}) \delta{-1} = \sqrt{}$(3. 7 7 412 310 4)
	= 9. 8 211 7 1 5 9 $\delta,$
	$cos(\alpha\frac{1}{2}) \delta{-1} = \sqrt{}$($4$.1211 1 7 2 9 1)
	= 2. 3 712121211 1 $\delta.$\\
	One of the signs of the square roots can be chosen arbitrarily,
	the other must be chosen in such a way that\\
	$sin(\alpha) = 2g\: sin(\alpha\frac{1}{2}) cos(\alpha\frac{1}{2}).$

	\sssec{Theorem.}
	{\em For $p \equiv -1 \pmod{4},$ with $\delta^2 = -1,$\\
	$cos(\alpha/2) \delta{-1}$ is a primitive root of\\
	\hth$V_{\frac{p-3}{2}} = (T_{\frac{p-1}{2}} I^{-1}) \circ (\delta I).$\\
	$sin(\alpha/2)$ is a primitive root of\\
	\hth$U_{\frac{p-3}{2}} = (T_{\frac{p-1}{2}} I^{-1}) \circ
	\sqrt{1 - I^2}.$}

	\sssec{Example.}
	For $p = 11,$ $\delta^2 = -1,$ $\alpha = \frac{\pi}{5},$ \\
	\hth$	V_4 = 16 I^4 + 20 I^2 + 5,$\\
	with roots $5$. 7 6 3 0 4 6 9 and $2$.10 5 5 6 5 6 6,\\
	\hth$	U_4 = 16 I^4 - 12 I^2 + 1,$\\
	with roots 9. 3 5 9 2 2 410 and 4. 910 3 8 7 9 4.
	Hence,
	$k\hti{16}sin(k \alpha)\hti{32}   cos(k \alpha)$\\
	$\begin{array}{rcrrrrrrrrrcrrrrrrrrr}
	\cline{1-21}
	0&\vline&\:\:0.&0&0&0&0&0&0&0&&\:\:1.&0&0&0&0&0&0&0\\
	\frac{1}{2}&\vline&9.&3&5&9&2&2&4&10&&5.&7&6&3&0&4&6&9&\delta\\
	1&\vline&2.&10&5&5&6&5&6&6&\delta&4.&9&10&3&8&7&9&4\\
	\frac{3}{2}&\vline&4.&9&10&3&8&7&9&4&&2.&10&5&5&6&5&6&6&\delta\\
	2&\vline&5.&7&6&3&0&4&6&9&\delta&9.&3&5&9&2&2&4&10\\
	\frac{5}{2}&\vline&1.&0&0&0&0&0&0&0&&0.&0&0&0&0&0&0&0\\
	3&\vline&5.&7&6&3&0&4&6&9&\delta&2.&7&5&1&8&8&6&0\\
	\frac{7}{2}&\vline&4.&9&10&3&8&7&9&4&&9.&0&5&5&4&5&4&4&\delta\\
	4&\vline&2.&10&5&5&6&5&6&6&\delta&7.&1&0&7&2&3&1&6\\
	\frac{9}{2}&\vline&9.&3&5&9&2&2&4&10&&6.&3&4&7&10&6&4&1&\delta\\
	5&\vline&0.&0&0&0&0&0&0&0,&\vline&10.&10&10&10&10&10&10&10,
	\end{array}$

	\sssec{Tables.}
	These can be found in the Handbook for Mathematical functions.
	Table 22.3 gives $T$ and $V$ and, by a simple transformation, $S'.$
	Table 22.5 gives $U.$\\
	$U$ and $V$ can be obtained by recurebces:\\
	$U_0 = 1,$ $U_2 = 4 I^2 - 1,$ $U_{2i+2}
		= 2 (2 I^2 - 1) U_{2i} - U_{2i-2}.$\\
	$V_0 = 1,$ $V_2 = 4 I^2 + 3,$ $V_{2i+2}
		= 2 (2 I^2 + 1) V_{2i} - V_{2i-2}.$\\
	We have\\
	$p = 5,$ $g = 2,$ $S'_2 = 4 I^2 - 1,$\\
	\hth$	c2 = s1$ = 2. 2 2 2 2 2 2 2, $c1$ = 2. 2 2 2 2 2 2 2.\\
	$p = 7,$  $U_2 = 4 I^2 - 1,$ $s1$ = 4. 3 3 3 3 3 3 3,\\
	\hth	$V_2 = 4 I^2 + 3,$ $c1$ = 1. 6 3 6 2 1 4 0.\\
	$p = 11,$ $U_4 = 16 I^4 - 12 I^2 + 1,$ $s1$ = 9. 3 5 9 2 2 410.\\
	\hth	$V_4 = 16 I^4 + 20 I^2 + 5,$ $c1$ = 5. 7 6 3 0 4 6 9.\\
	$p = 13,$ $g = 2,$ $S'_6 = 256 I^6 -192 I^4 + 36 I^2 -1,$\\
	\hth	$c3' = s1$ = 9. 7 1 1 0 912 0, $c1' = c1$ = 2. 1 8 7 6 2 6 7.\\
	$p = 17,$ $g = 3,$ $S'_8 = 10368 I^8 - 6912 I^6 + 1440 I^4 - 96 I + 1,
$\\
	\hth	$c5' = s1$ = 10. 8 8 4 3 514 1,
	\hth	$c1' = c1$ = 5.15151513 2 3 9.\\
	$p = 19,$ $U_8 = 256 I^8 - 448 I^6 + 240 I^4 - 40 I^2 + 1,$\\
	\hth	$s1$ = 14.13 0 618 9 118,\\
	\hth	$V_8 = 256 I^8 + 576 I^6 + 432 I^4 + 120 I^2 + 9,$\\
	\hth	$c1$ = 9.101611 61512 1.

	\sssec{Comment.}
	The following values may be useful,\\
	$\sqrt{-1}_5$	 = 2. 1 2 1 3 4 2 3 0 3 2\\
	$\sqrt{-1}_{13}$ = 5. 5 1 0 5 5 1 0 1 8 8,\\
	$\sqrt{-1}_{17}$ = 4. 210 5121612 813 314,\\
	$\sqrt{-1}_{29}$ =12. 112 1181615 3 92425.




	\ssec{The logarithm.}

	\sssec{Definition.}
	The {\em exponentiaonal function} and the {\em logarithmic function}
	are defined by the following $p$-adic expansion.\\
	D.0.\hti{4}$exp(x) = 1 + x +  \ldots  + \frac{1}{n!}x^n + \ldots,$
		$| x| < 1.$\\
	D.1.\hti{4}$log(1+x) = x - \frac{1}{2}x^2 +  \ldots 
		+ (-1)^n \frac{1}{n}x^n +  \ldots,$ $| x|  < 1.$

	The classical theorem is (See for instance Koblitz, 1977.)

	\sssec{Theorem.}
	{\em 
	\ldots	}

	\sssec{Motivation.}
	For p = 5,\\
	\hth$	log(-4) = log(1-5) = 0.41041,$\\
	\hth$	log(1 .1) = log(1+5) = 0.12420,$\\
	\hth$	log(-4.1) = log(1-10) = 0.32314.$\\
	If we want $log(x\:y) = log(x)+log(y)$ to hold, we have 3 equations
	to determine $log(-1),$ $log(2)$ and $log(3).$
	$log(1 .3)$ and $log(-4.4)$ can be used as check.  This gives\\
	\hth$	log(-1) = 0,$\\
	\hth$	log(2) = 0.23240,$\\
	\hth$	log(3) = 0.43134.$\\
	Clearly we now have a function which is not a bijection, for instance,\\
	\hth$	log(1 .20230) = 0 .23420.$\\
	This suggest that we can extend the range of definition of the
	logarithm function.
	The equation $x^{p-1} \equiv 1 \pmod{p}$ has $p-1$ roots,
	 1,2, \ldots,$p-1,$
	therefore the equation $x^{p-1} = 1$ has $p-1$ roots in the $p$-adic
	field, with first digit 1,2, \ldots, $p-1.$\\
	In general the roots are $1,$ $x1,$ $x2,$ \ldots, $-x2,$ $-x1,$ $-1.$

	\sssec{Algorithm.}\label{sec-newtonpadic}
	If $g$ is a primitive root in ${\Bbb Z}_p,$ the corresponding
	primitive root $g' \in \ldots$  can be obtained by Newton's method:\\
	\hth$	y := g,$\\
	\hth$	y = y - \frac{1}{\frac{p-1}{2}}
		(y + \frac{1}{y}^{\frac{p-3}{2}}$ for $i = 1$ to $n$.

	\sssec{Theorem.}
	{\em 
	Algorithm \ref{sec-newtonpadic} determines the first $2n$ digits of
	$g'$.}

	\sssec{Theorem. }
	{\em Given a prime $p,$ a primitive root $g \in {\Bbb Z}_p$ and $x,$}
	H0.\hti{5}$   | x|  = 1,$\\
	D0.\hti{5}$   x0 = ind(int(x)),$\\
	{\em where $ind$ is the index function in ${\Bbb Z}_p$ associated to
	$g,$}\\
	D1.\hti{5}$   y := x/g^{'x0} - 1$\\
	D2.\hti{5}$   z := y - \frac{1}{2}y^2 +  \ldots  + \frac{1}{n}(-y)^n$\\
	{\em then}\\
	C0.\hti{5}$   | z|  < 1.$\\
	C1.\hti{5}$   log_{p,g}(x) = \{x0, z\}$

	\sssec{Example.}
	For $p = 5,$ $x1 = 2 .1213423.$\\
	(All the roots are 1, 2 .1213423, $-2 .1213423 = 3 .3231021,$
	$-1 = 4.4444444$)

	For $p = 7,$ $x1 = 2.4630262,$ $x2 =3.4630262.$

	For $p = 11,$ $x1$ = 2.10 4 9 1 2 3 9, $x2$ = 3. 0 1 2 3 610 8,\\
	\hth$	    x3$ = 4. 7 9 5 2 9 8 0, $x4$ = 5. 2 5 1 7 8 510.

	For $p = 13,$ $x1$ = 2. 6 2 2 4 2 5 8, $x2$ = 3.11 6 9 7 2 4 4,\\
	\hth$x3$ = 4.11 6 9 7 2 4 4, $x5$ = 5. 5 1 0 5 5 1 0\\
	\hth$x6$ = 6. 1 910 3 5 6 4.

	For $p = 17,$ $x1$ = 2. 9 312 914 1 5, $x2$ = 3.13 2 3 011 4 0,\\
	\hth$x3$ = 4. 210 5121612 8, $x4$ = 5. 9 0 516 9 1 5,\\
	\hth$x5$ = 6. 214 4 1 6 2 3, $x6$ = 7. 4 216 11514 2,\\
	\hth$x7$ = 8. 6 1 415 116 2.

	For $p = 19,$ $x1$ = 2. 614 414131014, $x2$ = 3.16 7 8161815 1,\\
	\hth$x3$ = 4. 51717 5 614 0, $x4$ = 5. 3 3131113 716,\\
	\hth$x5$ = 6.12 21714181716, $x6$ = 7.15 7 0 118 0 4,\\
	\hth$x7$ = 8.15 7 0 118 0 4, $x8$ = 9. 118 21712 5 1.

	For $p = 23,$ $x1$ = 2.11211015 2 912, $x2$ = 3. 51717 71821 7,\\
	\hth	    $x3$ = 4.2122 4 7 81622, $x4$ = 5. 1 219 8 2 919,\\
	\hth	    $x5$ = 6.201517 114 720, $x6$ = 7.1519 5 8 81519,\\
	\hth	    $x7$ = 8.171710171911 3, $x8$ = 9. 713 1 7191512,\\
	\hth	    $x9$ =10.11 72112221517, $x10$=11. 81017 3 31922.

	The logarithmic functions as defined is not one to one. $p-1$ arguments
	give the same value.  To make it one to one we give the following

	\sssec{Definition.}
	Given a primitive root $g,$\\
	\hth$    log_{p,g}(x) = \{i0, x0, log_p(x)\},$\\
	\hth$	i0 \in (Z,+),$ $x0 \in ({\Bbb Z}_{p-1},+),$
		$log_p(x) \in \ldots .$\\
	where\\
	$| x|  = p^{-i0},$ $g^{x0} = int(x),$ and $log_p(x)$ is the $p$-adic
	logarithm.\\
	$i0$ is called the {\em characteristic}, $x0,$ the {\em index}
	 and $log_p(x),$ the {\em mantissa.}

	\sssec{Example.}
	For $p = 5$ and $g = 2,$\\
	\hth$log_{5,2}(1 .0000000) = \{0,0.\},$
		$log_{5,2}(2 .1213423) = \{1,0.\},$\\
	\hth$log_{5,2}(3 .3231021) = \{3,0.\},$
		$log_{5,2}(4 .4444444) = \{2,0.\}.$

	\sssec{Theorem.}
	H0.\hti{5}$   | x| ,$ $| y|  = 1.$\\
	C0.\hti{5}$   log_{p,g}(x * y) = log_{p,g}(x) + log_{p,g}(y).$

	\sssec{Problem.}
	Extend the definition to allow $| x|$, $|y|$ to be anything
	using the relation C.0. and the idea of mantissa.




	\ssec{P-adic Geometry and Related Finite Geometries.}

	\setcounter{subsubsection}{-1}
	\sssec{Introduction.}
	Some thought on finite geometry for different powers of $p$ and
	the $p$-adic geometry.

	Given $p,$ if we take the points on a line, with coordinates of the form
	$x0.$ \ldots these are indistinguishable if we use the $p$-adic
	valuation to the precision 1.\\
	If we consider the points $x0.x1,$ these
	are distinguishable to the precision 1 but not to the precision
	 $\frac{1}{p}.$\\
	They can be considered, if $p$ is 5 say, to have a color associated to
	the various digits of $x0$ the colors\\
	"purple, blue, green, yellow, orange, red."
	As we proceed to the next digit all yellows become sub-colours
	proceeding from yellow to orange.\\
	We can with more discrimination distinguish them.
	We can proceed further \ldots.\\
	We observe also that, in this scheme of things, what is a shade of
	yellow for one is a shade of green for some one else.
	This is associated to a change of origin.\\
	When this is transfered to angles this will give different
	trigonometric tables associated to different values of $x$.

	\sssec{Comment.}
	\vspace{-18pt}\hspace{94pt}\footnote{19.10.82}\\[8pt]
	The trigonometric functions work as follows.\\
	In the hyperbolic case,\\
	for $p,$ the period is $2(p-1),$\\
	for $p^2,$ the period is $2p(p-1),$  \ldots. \\
	The factor 2 corresponds to the fact, that just as in Euclidean geometry
	the total angle is $2 \pi,$ lines which form an angle $\pi$ correspond
	to the same direction.\\
	For the hyperbolic case, there are 2 real isotropic points,\\
	There are on the ideal line\\
	$p + 1$ points,\\
	$p - 1$ ideal (non isotropic) points for the $p$-geometry,\\
	$p^2 - 1$ ideal points for the $p^2$-geometry,
	of which $p^2 - p$ are not included in the preceding set.\\
	The hyperbolic trigonometry associated to $p^2$
	is presumably for these $p^2 - p$ directions \ldots.\\
	For the elliptic case, the period is $2(p+1)$ as one would expect
	for the $p$-geometry.\\
	For the $p^2$ case, the period of the trigonometric functions is
	$p^2 + p.$  I do not yet understand how this comes into the picture.

	\sssec{Definition.}
	\vspace{-18pt}\hspace{96pt}\footnote{21.10.82}\\[8pt]
	$p$-ADIC GEOMETRY

	Let $p \equiv -1 \pmod{4}$, let $z = 0$ be the ideal line,
	let $x^2 + y^2 = z^2$ be a circle.\\
	There are no solutions of $x^2 + y^2 = 0,$ therefore the
	isotropic points are not real.\\
	If $z \neq 0,$ then we can normalize using $z = 1,$ consider
	\hth$      x^2 + y^2 = 1.$\\
	Let $| x|  \geq | y|.$  If $| x|  > 1$ then there are no solutions.\\
	(Hint: divide by $p^{| x| }$ and work modulo $p$)

	\sssec{Lemma.}
	{\em 
	If $s = sin(\alpha)$ is a root of \ldots  then\\
	\hth$sin(p \alpha ) = sin(\alpha ),$ $cos(p \alpha ) = cos(\alpha ).$}

	\sssec{Theorem.}
	{\em 
	Let $s = sin(\alpha)$ be a root of \ldots.\\
	If $|sin(\beta)-sin(\alpha)|<1$ then\\
	\hth$	lim_{n\rightarrow\infty} sin(p^n \beta ) = sin(\alpha ).$\\
	If $| x|  = 1,$ there are solutions, $x = x0.x1x2$ \ldots.
	Let $| x0|$  be a primitive solution of the polynomial \ldots.\\ 
	let $x1.$ be \ldots.\\
	To $x$ corresponds the pair $(x,y) = (sin(\alpha),cos(\alpha)),$\\
	Let $X(x) = sin(p \alpha),$
	then the sequence $x,X(x),X^{2(x)},$  \ldots  convergences to $x'.$
	Moreover $x'$ is a root of $\ldots$ }

	\sssec{Example.}
	$p = 11,$ let $x = sin(\alpha) = 3,$
	X(3. 0 0 0 0 0 0 0) = 3. 0 2 8 011 4 6,\\
	\hth$	cos(p \alpha ) = 10.10 312 5 7 4 1,$\\
	X(3. 0 2 0 0 0 0 0) = 3. 0 2 6 0 71212,\\
	X(3. 0 2 6 0 0 0 0) = 3. 0 2 6 7 612 2,\\
	X(3. 0 2 6 7 0 0 0) = 3. 0 2 6 71211 2,\\
	X(3. 0 2 6 712 0 0) = 3. 0 2 6 71212 1,\\
	X(3. 0 2 6 71212 0) = 3. 0 2 6 71212 2.

	\newpage
	\sssec{Example.}
	$p = 11,$ elliptic case, $g = -1$\\
	The following is a table of $sin,$ with
	$sin(k=(2l+1) \alpha )\delta^{-1},$ $sin(k=2l \alpha)$\\
	$\begin{array}{llccccccc}
	1&\vline& 4.000&4.100&4.200&4.300&4.400&4.500&4.600\\
	2&\vline& 5.325&5.112&5.665&5.450&5.245&5.034&5.512\\
	3&\vline& 2.350&2.241&2.115&2.030&2.624&2.506&2.431\\
	4&\vline& 1.063&1.031&1.033&1.064&1.050&1.000&1.054\\
	5&\vline& 2.632&2.322&2.033&2.451&2.130&2.522&2.243\\
	6&\vline& 5.524&5.430&5.303&5.200&5.130&5.026&5.635\\
	0 .1&\vline&  4.566&4.550&4.540&4.536&4.524&4.511&4.504\\
	\cline{1-9}
	1&\vline& 4.500&4.510&4.520&4.530&4.540&4.550&4.560\\
	2&\vline& 5.034&5.013&5.061&5.040&5.026&5.005&5.053\\
	3&\vline& 2.506&2.562&2.556&2.543&2.530&2.524&2.511\\
	4&\vline& 1.000&1.000&1.000&1.000&1.000&1.000&1.000\\
	5&\vline& 2.522&2.562&2.533&2.504&2.544&2.515&2.555\\
	6&\vline& 5.026&5.013&5.000&5.063&5.050&5.044&5.031\\
	0 .1&\vline&  4.511&4.510&4.516&4.515&4.514&4.513&4.512\\
	\cline{1-9}
	1&\vline& 4.510&4.511&4.512&4.513&4.514&4.515&4.516\\
	2&\vline& 5.013&5.011&5.016&5.014&5.012&5.010&5.015\\
	3&\vline& 2.562&2.561&2.560&2.566&2.565&2.564&2.563\\
	4&\vline& 1.000&1.000&1.000&1.000&1.000&1.000&1.000\\
	5&\vline& 2.562&2.566&2.563&2.560&2.564&2.561&2.565\\
	6&\vline& 5.013&5.012&5.011&5.010&5.016&5.015&5.014\\
	0 .1&\vline&  4.510&4.510&4.510&4.510&4.510&4.510&4.510
	\end{array}$



\chapter{DIFFERENTIAL EQUATIONS AND FINITE MECHANICS}


\setcounter{section}{-1}
	\section{Introduction.}
	In the context of Finite Geometry, we should examine the subject of
	Differential Equations, their approximation and the application to
	finite mechanics.\\
	I will describe the first success associated with the harmonic
	polygonal motion, then \ldots

	\section{The first Examples of discrete motions.}
	\ssec{The harmonic polygonal motion.}
	\setcounter{subsubsection}{-1}
	\sssec{Introduction.}

	In classical Euclidean geometry as well as in finite Euclidean geometry
	I define the harmonic polygonal motion as the motion which associates
	to linearly increasing time successive points of the harmonic polygon.
	I will determine, for the classical case, the differential equation of
	the motion, by considering first points which are close to each other
	\ref{sec-tharmpolyg}. I will then prove that this equation is satisfied
	when points are not close to each other \ref{sec-tharmpolyg1}.
	The equation bears ressemblance with the equation
	of Kepler.  Because the method uses derivatives of functions of the
	trigonometric functions only and in view of the method of Hensel for
	$p$-adic functions, the result extend automatically to the finite case.

	\sssec{Definition.}
	Given a conic
	\enumb
	\item$	A(E) = (a cos(E), b sin(E), 1),$\\
	and the point of Lemoine $K = (q,0),$ given the correspond harmonic
	polygon of Casey (g2734, $p.5\ldots$ ), $A_i,$ 
	we define the {\em harmonic polygonal motion} by
	\item$	A(E(t_0+i\:h)) = A_i.$
	\enume

	\sssec{Theorem.}\label{sec-tharmpolyg}
	{\em Let}
	\enumb
	\item$	r = \frac{q}{a},$\\
	{\em if $h$ is small, the motion satisfies the differential equation}
	\item$	C\:DE = 1 - r\:cos(E),$\\
	{\em for some constant of integration $C.$}

	Proof:
	The polar $k$ of $K$ is
	\item$	k = [\frac{q}{a^2},0,-1],$\\[10pt]
	the polar $a(E)$ of $A(E)$ is
	\item$	a(E) = [\frac{cos(E)}{a},\frac{sin(E)}{b},-1],$\\[10pt]
	it meets $k$ at
	\item$	B(E) = (a^2sin(E),b(q-acos(E)),qsin(E)).$\\[10pt]
	the condition that $A_{i-1} \times A_{i+1}$ passes through $B(E(t))$
	gives
	\item$\matb{|}{ccc}a^2 sin(E(t))&   b(q-acos(E(t)))& qsin(E(t))\\
			acos(E(t-h))&   bsin(E(t-h))&    1\\
			acos(E(t+h))&   bsin(E(t+h))&    1\mate{|} = 0.$\\[10pt]
	Let
	\item$	s(t) = \frac{1}{2}(E(t+h)+E(t-h)),
		d(t) = \frac{1}{2}(E(t+h)-E(t-h)),$
	then to the second order in $h,$ with $k = \frac{1}{2}h^2,$\\
	\hth$	s = E + k D^2 E,$ $d = h DE,$
	\item	$cos(s) = cos(E) - k sin(E) D^2 E,$
		$sin(s) = sin(E) + k cos(E) D^2 E,$
		$cos(d) = 1 - k (DE)^2.$\\
	Replacing the determinant by that obtained by using instead of the
	last 2 lines their half sum and their half difference, gives after
	division of the first row and first column by $a$ and the second
	column by $b,$\\
	\hth$\matb{|}{ccc}sin(E)&   r-cos(E)&	rsin(E)\\
		cos(s)cos(d)&   sin(s)cos(d)&    1\\
		-sin(s)sin(d)&  cos(s)sin(d)&    0\mate{|} = 0.$\\
	Dividing the last row by $sin(d)$ and expanding with respect to the
	first row gives, after changing sign,
	\item$	sin(E)cos(s)+(r-cos(E))sin(s)-rsin(E)cos(d) = 0,$\\
	or, after using 7 and dividing by $k,$
	\item$	(1 - r cos(E)) D^2 E = r sin(E) (DE)^2,$\\
	integrating gives 1, with some appropriate constant $C,$
	\item$	(1 - e cos(E)) DE = C.$
	This is tantalizing close to Kepler's equation.
	\enume

	\sssec{Theorem.}\label{sec-tharmpolyg1}
	{\em The harmonic polygonal motion associated to the point of
	Lemoine $(r\:a,0,1)$ is described on the ellipse by the differential
	equation}
	\enumb
	\item$	C\: DE = 1 - r\: cos(E).$

	Proof:  We have to show that if we take the derivative of the relation
	between 3 points equidistant in time, namely \ref{sec-tharmpolyg}.8,
	this derivative is 0 if the differential equation 0. is satisfied.
	We can assume that $C = 1.$\\
	From 0 and from \ref{sec-tharmpolyg}.6. follows
	\item$	Ds = 1 - \frac{1}{2}r(cos(E(t+h))+cos(E(t-h)))$\\
	\hth$	   = 1 - r cos(s) cos(d),$
	\item$	Dd = - \frac{1}{2}r(cos(E(t+h))-cos(E(t-h)))$\\
	\hth$	   = r sin(s) sin(d).$\\
	Taking the derivative of \ref{sec-tharmpolyg}.8, gives\\
	$cos(s) cos(E) (1-rcos(E)) - sin(s) sin(E) (1-rcos(s)cos(d))\\
	\hth+ sin(s)sin(E) (1-rcos(E)) + cos(s) (r - cos(E)) (1-rcos(s)cos(d))\\
	\hth- r cos(d) cos(E) (1-rcos(E)) + r^2 sin(s) sin(E) sin(d)^2.$\\ 
	We would like to prove that this expression is identically zero.
	0, gives
	\item$	r cos(d) = cos(s) + sin(s) \frac{r - cos(E) }{ sin(E)},$
	substituting in the expression gives\\
	$cos(s)cos(E)(1-rcos(E))$\\
	\hth$    - r sin(s)sin(E)cos(E) + cos(s)(r - cos(E))$\\
	\hth$    + r^2 sin(s)sin(E)$\\
	\hth$    + (sin(s)cos(s)sin(E) - r cos^2(s) + cos^2(s)cos(E) - cos(E)$\\
	\hth$	+ r cos^2(E)) (cos(s) + sin(s)\frac{r - cos(E)}{sin(E)})$\\
	\hth$    - sin(s)sin(E)(cos(s)+sin(s)(r - cos(E)))/ sin(E))^2.$
	\enume

	The coefficient of $r^2$ is\\
	$sin(s)sin(E) + \frac{(cos^2(E) -cos^2(s))sin(s)}{sin(E)}
	- \frac{sin^3(s)}{sin(E)}$\\
	\hth$= \frac{sin(s)(sin^2(E) + cos^2(E) - cos^2(s) - sin^2(s)}
	{sin(E)} = 0.$\\
	The coefficient of $r$ is\\
	$-cos(s)cos^2(E) - sin(s)sin(E)cos(E) + cos(s)$\\
	\hth$    + (cos(s) - sin(s)\frac{cos(E)}{sin(E)}(cos^2(E) - cos^2(s))$\\
	\hth$    + sin(s)(sin(s)cos(s)sin(E)+cos^2(s)cos(E)-cos(E))/sin(E)$\\
	\hth$    - 2 sin(s)sin(E)(cos(s)-sin(s)cos(E)/sin(E))sin(s)/sin(E)$\\
	\hti{4}$= sin(s)(-cos^3(E)+cos(E)cos^2(s)-cos(E)+sin(s)cos(s)sin(E)$\\
	\hth$    +cos^2(s)cos(E)-2cos(s)sin(s)sin(E)+2sin^2(s)cos(E))/sin(E)$\\
	\hth$    - cos(s)cos^2(E)-sin(s)sin(E)cos(E)+cos(s)
		+cos(s)cos^2(E)-cos^3(s)$\\
	\hti{4}$= sin(s)(cos(E)(-cos^2(E)+cos^2(s)-1+cos^2(s)+2sin^2(s))/sin(E)
$\\
	\hth$    -sin(s)sin(E) cos(E) - sin(s)sin(E)cos(E)+cos(s)sin^2(s)\\
	\hti{4}= sin(s)cos(E)sin(E)-sin^2(s)cos(s)-sin(s)sin(E)cos(E)
		+cos(s)sin^2(s) = 0.$\\
	The term independent of $r$ is\\
	$cos(s)cos(E) - cos(s)cos(E)$\\
	\hth$    + (sin(s)cos(s)sin(E)+cos^2(s)cos(E)-cos(E))$\\
	\hth$	(cos(s)-sin(s)cos(E)/sin(E)$\\
	\hth$    - sin(s)sin(E)(cos(s)-sin(s)cos(E)/sin(E))^2$\\
	\hti{4}$= (cos(s)-sin(s)cos(E)/sin(E))$\\
	\hth$    (sin(s)cos(s)sin(E)+cos^2(s)cos(E)-cos(E)-sin(s)cos(s)sin(E)$\\
	\hth$	+sin^2(s)cos(E)) = 0.$

	\sssec{Theorem.}
	\vspace{-18pt}\hspace{90pt}\footnote{26.6.83}\\[8pt]
	{\em Let $e'^2 = 1 - e^2,$ then}
	\enumb
	\item$	e'\:tan(\frac{1}{2}e'M) = (1+e)\:tan (\frac{1}{2}E).$
	\enume

	\sssec{Theorem.}
	\enumb
	\item {\em If $t(M) = tan(\frac{1}{2}E)$ and}
	\footnote{29.6.83}
	\item$	k = \frac{e+1}{e-1},$
	\item$ t'_i = \sqrt{-k}t(M_i)$,\\[5pt]
	{\em then}
	\item$	t'_{1+2} = \frac{t'_1 + t'_2}{1 - t'_1t'_2}.$
	\enume

	\sssec{Example.}
{\small	For $p = 13$ and $e = 2,$ $k = 3$ and\\
	$\begin{array}{ccrrrrrrrrrrrrr}
	M &\vline& 0&1&2&3&4&5&6&7&8&9&10&11&12\\
	t(M) &\vline& 0&1&7&11&8&4&\infty&9&5&2&6&12&0\\
	t'_M&\vline&0&6&3&1&-4&-2&\infty&2&4&-1&-3&-6&0\\
	\frac{1}{2}E &\vline& 0&3&7&1&4&6&8&11&5&9
	\end{array}$

	For $p = 13$ and $e = 3,$ $k = 2$ and\\
	$\begin{array}{ccrrrrrrrrrrrrrr}
	M &\vline& 0&1&2&3&4&5&6&7&8&9&10&11&12&13\\
	t(M) &\vline& 0&2&12&4&8&7&3&\infty&10&6&5&9&1&11\\
	\frac{1}{2}E &\vline& 0
	\end{array}$
}
	\sssec{Programs.}
	The programs pl.bas and planet.bas  $\ldots$ 

	\sssec{Exercise.}
	Prove that the acceleration is\\
	$\frac{1-r\:cosE}{C^2}(-(a\:cosE,b\:sinE,0)+r(a\:cos(2E),b\:sin(2E),0))$
	\ssec{The Parabolic Motion.}
\setcounter{subsubsection}{-1}
	\sssec{Introduction.}
	The parabola has been studied in g33.
	Galileo Galilei was the first to show that the motion of a particle
	in a uniform gravitional field is a parabola.  (Love, p.45)
	The result extend to the finite case.

	\sssec{Theorem.}
	{\em In both the infinite and finite cases, the solution of\\
	\hth$	m D^2 x = 0$ and $m D^2 y = -mg$\\
	is\\
	\hth$	x(t) = v_0 t,$ $y(t) = -\frac{1}{2}g t^2 + v_1 t,$ or\\
	\hth$	y = a x^2 + b x,$ with\\
	\hth$	a := -\frac{g}{2v_0^2},$ $b := \frac{v_1}{v_0}.$}

	Proof: Comparing the equation in the form\\
	\hth$	(x + \frac{b}{2a})^2 = \frac{y}{a}+(\frac{b}{2a})^2.$\\
		with the standard equation $y^2 = 4c x$
	shows that the vertex $V$ and the directrix $d$ are\\
	\hth$	V = (\frac{b}{2a}, (\frac{b}{2a})^2),$\\
	\hth$	d: \:y = \frac{v_0^2 + v_1^2}{2g} = \frac{v^2}{2g}$\\
	corresponding to the Torricelli law.

	\sssec{Example.}
{\small	For $p = 7,$ $g = 1$ and $v_0 = v_1 = 4,$ then $a = -2$, $b = 1$,\\
		$y(x) = -2x^2 + x = 1-2(x-2)^2,$\\
	$\begin{array}{ccrrrrrrrr}
	x&\vline&0&1&2&3&-3&-2&-1&0\\
	y&\vline&0&-1&1&-1&0&-3&-3&1\\
	z&\vline&1&1&1&1&1&1&1&0\\
	t&\vline&0&2&-3&-1&1&3&-2
	\end{array}$
}
	\ssec{Attempts to Generalize Kepler's Equation.}

	\setcounter{subsubsection}{-1}
	\sssec{Introduction.}
	I have made many attempts to generalize Kepler's
	equation or the simple planetary motion to the finite case.
	In section $\ldots$ , I examine the use of $p$-adic function to obtain
	a solution in the neighbourhood of a circular motion.

	\ssec{The circular motion.}

	\sssec{Definition.}
	The {\em circular motion} is defined by\\
	\hth$x(t) = cos(t),$ $y(t) = sin(t),$\\
	\hth$Dx(t) = -sin(t),$ $Dy(t) = cos(t).$\\
	This assumes that the unit of distance is chosen as the radius of the
	circle and the unit of time is chosen in such a way that the period
	is $2\pi$.

	\section{Approximation to the Solution of Differential Equations.}

	\setcounter{subsection}{-1}
	\ssec{Introduction.}
	To approximate the solution of differential equations
	it is important to insure that essential properties are preserved.
	In particular, for conservative systems, the same should hold.
	In this connection, I developed in 1956 a method of first order
	and a method of second order which are contact transformations and
	therefore preserve the essential properties of conservative systems.
	These will be applied to the finite case.

	\ssec{Some Algorithms.}
	\sssec{Algorithm.}
	The first order algorithm is defined by

	\sssec{Theorem.}
	{\em }

	\sssec{Algorithm.}\label{sec-adeq}
	The second order algorithm for the solution of the
	differential equation\\
	\hth$D^2 {\bf x} = {\bf f} \circ {\bf x},$ ${\bf x}(0) = {\bf x}_0,$
	${\bf Dx}(0) = {\bf Dx}_0,$ \\
	is defined by\\
	\hth${\bf x}_{i+1} = {\bf x}_i + h {\bf Dx}_i
		+ \frac{1}{2}h^2 {\bf f}_i, 
		{\bf Dx}_{i+1} = {\bf Dx}_i + \frac{1}{2}h ({\bf f}_{i+1}
		+ {\bf f}_i),$ \\
	where\\
	\hth$	{\bf f}_i := {\bf f}({\bf x}_i).$ 

	\sssec{Definition.}
	A {\em mapping is reversible} iff

	\sssec{Theorem.}
	{\em Given the Algorithm \ref{sec-adeq}, the mapping is reversible.}

	Proof:  If we solve for ${\bf x}_i$ and ${\bf Dx}_i,$ we get\\
	\hth$	{\bf Dx}_i = {\bf Dx}_{i+1} - \frac{h}{2} ({\bf f}_{i+1}
		 + {\bf f}_i),$\\
	\hth$	{\bf x}_i = {\bf x}_{i+1} - h {\bf Dx}_i
		- \frac{1}{2}h^2 {\bf f}_i,\\
	\hti{12}		= {\bf x}_{i+1} - h {\bf Dx}_{i+1}
	+ \frac{1}{2}h^2 ({\bf f}_i + {\bf f}_{i+1}).$

	\sssec{Definition.}
	A {\em mapping is symplectic} iff

	\sssec{Theorem.}
	{\em The mapping defined in algorithm \ref{sec-adeq} is symplectic.}

	Proof:

	\sssec{Example.}
{\small	Let $x$ and $f$ be one dimensional, let\\
	\hth$	f(x) = - x - 2 x^3,$ \\
	let $Dx_0 = 0,$ we have the following solutions, for $h = 1$
	and various initial conditions $(x(0),Dx(0)= 0)$.\\
	$p = 11,$\\
	$\begin{array}{ccccccccccc}
	i&\vline&0&1&2&3&4&5&6&7&8\\
	\hline
	(x,Dx)_i&\vline&1,0&5,3&-4,2&-2,0&-4,-2&5,-3&1,0\\
	&\vline&3,0&2,1&5,2&-5,2&-2,1&-3,0&-2,-1&-5,-2&5,-2\\
	&\vline&4,0&4,0\\
	&\vline&5,0&4,-1&3,-2&0,-3&-3,-2&-4,-1&-5,0&-4,1&-3,2\\
	\end{array}$\\
	$p = 13,$\\
	$\begin{array}{ccccccccccc}
	i&\vline&0&1&2&3&4&5&6&7&8\\
	\hline
	(x,Dx)_i&\vline&1,0&6,-6&2,0&6,6&1,0\\
	&\vline&3,0&-6,2&-6,-2&3,0\\
	&\vline&4,0&3,3&-3,3&-4,0&-3,-3&3,-3&4,0\\
	&\vline&5,0&1,1&-6,4&-4,3&0,4&4,3&6,4&-1,1&-5,0\\
	&\vline&6,0&-5,6&5,6&-6,0&5,-6&-5,-6&6,0\\
	\end{array}$
}
	\sssec{Theorem.}
	{\em If we apply the mapping \ref{sec-adeq} to\\
	\hth$D^2x = -x,$ $x(0) = 0,$ $Dx(0) = 1,$\\
	we obtain, up to a scaling factor the trigonometric functions.}

	\sssec{Example.}
{\small	With $p = 13$ and $h = 1,$\\
	$\begin{array}{cccccccccc}
	i&\vline&0&1&2&3&4&5&6&7\\
	&&&        8&9&10&11&12&13&14\\
	(x,Dx)_i&\vline&0,1&1,-5&3,-3&-5,-4&-5,4&3,3&1,5&0,-1
	\end{array}$

	$?(\delta x,Dx)((i) = (sin,cos)(-10 i),$ if $sin(1) = 3$ and
	$cos(1) = 3\delta,$
	with $\delta^2 = 2.$
}
	\sssec{Program.}
	[130] PENDUL(um)



	\section{The Parabolic Motion.}
\setcounter{subsection}{-1}
	\subsection{Introduction.}
	The parabola has been studied in g33.
	Galileo Galilei was the first to show that the motion of a particle
	in a uniform gravitional field is a parabola.  (Love, p.45)
	The result extend to the finite case.

	\subsubsection{Theorem.}
	{\em In both the infinite and finite cases, the solution of\\
	\hth$	m D^2 x = 0$ and $m D^2 y = -mg$\\
	is\\
	\hth$	x(t) = v0 t,$ $y(t) = -\frac{1}{2}g t^2 + v1 t,$ or\\
	\hth$	y(x) = a x^2 + b x,$ with\\
	\hth$	a := \frac{g}{2v0^2},$ $b := \frac{v1}{v0}.$}

	Proof: Comparing the equation in the form\\
	\hth$	(x - \frac{b}{2a)})^2 = -\frac{(y-\frac{b^2}{4a)})}{a}$\\
		with the standard equation $y^2 = 4c x$
	shows that the vertex $V$ and the directrix $d$ are\\
	\hth$	V = (\frac{b}{2a)}, \frac{b^2}{4a)}),$\\
	\hth$	d: y = \frac{v0^2 + v1^2}{2g} = \frac{v^2}{2g)}$\\
	corresponding to the Torricelli law.

	\subsubsection{Example.}
	For $p = 7,$ $g = 1$ and $v0 = v1 = 4,$\\
		y(x) = $-2x^2$ + x,\\
	$\begin{array}{ccrrrrrrrr}
	x&\vline&0&1&2&3&-3&-2&-1&0\\
	y&\vline&0&-1&1&-1&0&-3&-3&1\\
	z&\vline&1&1&1&1&1&1&1&0\\
	t&\vline&0&2&-3&-1&1&3&-2
	\end{array}$

	\section{Attempts to Generalize Kepler's Equation.}

	\setcounter{subsubsection}{-1}
	\subsubsection{Introduction.}
	I have made many attempts to generalize Kepler's
	equation or the simple planetary motion to the finite case.
	In section $\ldots$ , I examine the use of $p$-adic function to obtain
	a solution in the neighbourhood of a circular motion.

	\subsection{The circular motion.}

	\subsubsection{Definition.}
	The {\em circular motion} is defined by\\
	\hth$x(t) = cos(t),$ $y(t) = sin(t),$\\
	\hth$Dx(t) = -sin(t),$ $Dy(t) = cos(t).$\\
	This assumes that the unit of distance is chosen as the radius of the
	circle and the unit of time is chosen in such a way that the period
	is $2\pi$.

	\section{Approximation to the Solution of Differential Equations.}

	\setcounter{subsubsection}{-1}
	\subsubsection{Introduction.}
	To approximate the solution of differential equations
	it is important to insure that essential properties are preserved.
	In particular, for conservative systems, the same should hold.
	In this connection, I developed in 1956 a method of first order
	and a method of second order which are contact transformations and
	therefore preserve the essential properties of conservative systems.
	These will be applied to the finite case.

	\subsubsection{Algorithm.}
	The first order algorithm is defined by

	\subsubsection{Theorem.}
	{\em }

	\subsubsection{Algorithm.}
	The second order algorithm for the solution of the
	differential equation\\
	\hth$D^2 {\bf x} = {\bf f} \circ {\bf x},$ ${\bf x}(0) = {\bf x}_0,$
	${\bf Dx}(0) = {\bf Dx}_0,$ \\
	is defined by\\
	\hth${\bf x}_{i+1} = {\bf x}_i + h {\bf Dx}_i
		+ \frac{1}{2}h^2 {\bf f}_i, 
		{\bf Dx}_{i+1} = {\bf Dx}_i + \frac{1}{2}h ({\bf f}_{i+1}
		+ {\bf f}_i),$ \\
	where\\
	\hth$	{\bf f}_i := {\bf f}({\bf x}_i).$ 

	\subsubsection{Definition.}
	A {\em mapping is reversible} iff

	\subsubsection{Theorem.}
	{\em Given the Algorithm 4.1.3., the mapping is reversible.}

	Proof:  If we solve for ${\bf x}_i$ and ${\bf Dx}_i,$ we get\\
	\hth$	{\bf Dx}_i = {\bf Dx}_{i+1} - \frac{h}{2} ({\bf f}_{i+1}
		 + {\bf f}_i),$\\
	\hth$	{\bf x}_i = {\bf x}_{i+1} - h {\bf Dx}_i
		- \frac{1}{2}h^2 {\bf f}_i,
	\hti{12}		= {\bf x}_{i+1} - h {\bf Dx}_{i+1}
	+ \frac{1}{2}h^2 ({\bf f}_i + {\bf f}_{i+1}).$

	\subsubsection{Definition.}
	A {\em mapping is  } iff

	\subsubsection{Theorem.}
	{\em The mapping defined in algorithm 4.1.3. is}

	Proof:

	\subsubsection{Example.}
	Let $x$ and $f$ be one dimensional, let\\
	\hth$	f(x) = - x - 2 x^3,$ \\
	let $Dx_0 = 0,$ we have the folowing solutions\\
	$p = 11,$\\
	$\begin{array}{cccccccccc}
	i&\vline&0&1&2&3&4&5&6&7\\
	&&8&9&10&11\\
	(x,Dx)_i&\vline&1,0&5,3&-4,2&-2,0&-4,-2&5,-3\\
	&3,0&2,1&5,2&-2,1&-3,0&-2,-1&-5,-2&5,-2\\
	&&2,-1\\
	&4,0\\
	&5,0&4,-1&3,-2&0,-3&-3,-2&-4,-1&-5,0&-4,1\\
	&&-3,2&0,3&3,2&4,1
	\end{array}$

	$p = 13,$\\
	$\begin{array}{cccccccccc}
	i&\vline&0&1&2&3&4&5&6&7\\
	&&&8&9&10&11&12&13&14\\
	&&&15\\
	(x,Dx)_i&\vline&1,0&6,-6&2,0&6,6\\
	&3,0&-6,2&-6,-2\\
	&4,0&3,3&-3,3&-4,0&-3,-3&3,-3\\
	&5,0&1,1&-6,4&-4,3&0,4&4,3&6,4&-1,1\\
	&&-5,0&-1,-1&6,-4&4,-3&0,-4&-4,-3&-6,-4\\
	&&1,-1\\
	&6,0&-5,6&5,6&-6,0&5,-6&-5,-6\\
	\end{array}$

	\subsubsection{Theorem.}
	{\em If we apply the mapping 0.3. to
	\hth$D^2x = -x,$ $x(0) = 0,$ $Dx(0) = 1,$\\
	we obtain, up to a scaling factor the trigonometric functions.}

	\subsubsection{Example.}
	With $p = 13$ and $h = 1,$\\
	$\begin{array}{cccccccccc}
	i&\vline&0&1&2&3&4&5&6&7\\
	&&&        8&9&10&11&12&13&14\\
	(x,Dx)_i&\vline&0,1&1,-5&3,-3&-5,-4&-5,4&3,3&1,5&0,-1
	\end{array}$

	$(\delta x,Dx)((i) = (sin,cos)(-10 i),$ if $sin(1) = 3$ and
	$cos(1) = 3\delta,$ with $\delta^2 = 2.$

	\subsubsection{Program.}
	[130] PENDUL(um)




	\ssec{On the existence of primitive roots.}
	\setcounter{subsubsection}{-1}
	\sssec{Introduction.}
	I will first give a non constructive proof of the existence of primitive
	roots and the give a construction. The first proof insures that the
	construction is always successful.

	\sssec{Theorem.}
	\begin{enumerate}
	\setcounter{enumi}{-1}
	\item$  d = ord_p(x),$ $(d,p) = g,$ $0 < l < d
		\Rightarrow ord_p(x^{l)} = \frac{d}{g},$
	\item$  d = ord_p(x),$ $0 \leq i,j < d,$ $x^i \equiv x^j \pmod{p}
		\Rightarrow i = j.$
	\item{\em If $d | p-1$ then $x^d \equiv 1 \pmod{p}$ has $\phi (d)$
		solutions of order $d.$}  Hint 2.25.
	\item{\em In particular, there are $\phi(p-1)$ primitive roots of $p.$}
	\item$  d = ord_p(z),$ $e = ord_p(y),$ $(d,e) = 1
		\Rightarrow ord_p(z.y) = d.e.$
	\end{enumerate}

	What follows is a Theorem which gives a constructive method of\\
	determining primitive roots or more generally of solutions of\\
	\hth$	d = ord(x),$ where $d | p-1.$\\
	The construction is inspired by Gauss, 1801, section 55.

	\sssec{Theorem.}
	{\em Let $\Pi_{j = 1}^n p_j^{i_j}$ be a prime factorization of $q-1.$\\
	Let $\frac{a_j^{(q-1)}}{p_j} -\neq  1 \pmod{q}$
	and $a_j^{(q-1)} \equiv 1 \pmod{q},$ for $j = 1,2, \ldots n,$ then}
	\begin{enumerate}
	\setcounter{enumi}{-1}
	\item {\em  $p_j^k = ord_q( a_j^{\frac{q-1}{P_j}} ),$
		$0 \leq k_j \leq i_j.$\\
	Let $P_j = p_j^{i_j,}$ let $h_j \equiv a_j^{\frac{q-1}{P_j}} \pmod{q},$
		then,}
	\item {\em  in particular, $p_j = ord_q(h_j).$\\
	Let $h_j^{k_j} = a_j^{\frac{q-1}{i_j^{k_j}}} \pmod{q},$
	$ 0 \leq k_j < i_j,$ then}
	\item$  \Pi_{j = 1}^n h_j^{k_j}$\\
	\hth$	= ord_q ( \Pi_{j = 1}^n h_j^{k_j} ).$\\
	{\em Let $h = \Pi_{j = 1}^n h_j \pmod{q},$ then,}
	\item  {\em in particular, $q-1 = ord_q(h),$}
	\item  {\em $q$ is prime,}
	\item  {\em $h$ is a primitive root of $q.$}
	\end{enumerate}


\chapter{COMPUTER IMPLEMENTATION}



\setcounter{section}{-1}
	\section{Introduction.}
	One of the tradition of Mathematicians is to discover properties
	by working on special cases or examples, this is especially so at the
	beginning of many branches of Mathematics, geometry, number theory,
	algebra, \ldots.  This was certainly the tradition kept up by Euler, see
	\ldots, by Gauss, see \ldots.\\
	In Euclidean geometry, the special cases were obtained by drawing a
	reasonably accurate figure, in number theory by numerical computation,
	and in algebra by algebraic manipulations.  All three can now be done
	accurately and with great speed using computers and these are now
	becoming more and more available to every one.\\
	Depending on our training or, I believe, on the structure of our
	individual brain, such experimentation is almost essentail for many
	to obtain a thourough understanding of basic concepts.\\
	To help in the understanding of the material given above and, I hope,
	to help the reader in the discovery of new properties, it is becoming
	essentail to provide him with the tools to realize quickly computer
	programs.\\
	When the subject matter is well settled and the experimentation is not
	at the basic level, a higher level non interactive language such as
	FORTRAN, ALGOL, PASCAL, PL1, ADA, is an excellent choice.  When this is
	not the case, an interactive language such as BASIC or APL is by far
	preferable.
	BASIC, BASIC+, BASIC+ extended.

	Hardware, operating system, files, interaction, language, compiler,
	interpreter.


\clearpage

\chapter*{REFERENCES}

\begin{enumerate}
\item Adobe Systems, Postscript Language, Tutorial and Cookbook, N. Y., Addison- Wesley, 1985, 243 pp. \\ 
\item	Adobe Systems, Postscript Language, Reference Manual, N. Y., Addison-Wesley, 1985, 319 pp. \\ 
\item	Apollonius, A treatise on Conic Sections, ed. Th. L. Heath, Cambridge, 1896. \\ 
\item	Apollonius, Les Coniques d'Apollonius de Perg\'{e}, trans. P. ver Ecke, Bruges, Belgique, 1923. \\ 
\item	Artin, E., Geometric Algebra, N. Y., Interscience, 1957. \\ 
\item	Artzy, Rafael, Linear Geometry, Reading Mass., Addison-Wesley, 1965, 273 pp. \\ 
\item	Aryabhata I, The Aryabhatiya of Aryabhata, Tranlated with notes by Walter Eugen Clark, Chicago, Ill. Univ. of Chicago Press, 1930. \\ 
\item	Aryabhata I, Aryabhatiya, Ed. by Kripa Shandar Shukla, New Delhi, Indian Nat. Sc. Acad.,1976. 219 pp. \\ 
\item	Baker, Henry, Frederick, Principles of Geometry, Vol. 1 to 4, Cambridge Univ. Press 1, 1929, 195 pp. 2, 1930, 259 pp. \\ 
\item	Barbilian, Dan, (or Barbu, Ion), Pagini Inedite, Vol. 2, Bucarest, Ed. Albatros, 1984, 292 pp. \\ 
\item	Baumert, Leonard D., Cyclic Difference Sets, N. Y., Springer, 1971. \\ 
\item	B\'{e}zier, P., Definition num\'{e}rique des courbes et surfaces, I, II, Automatismes, Vol. 11, 1966, 625-632 and 12, 1967, 17-21. Also Vol. 13, 1968 and Thesis, Univ. of Paris VI, 1977. \\ 
\item	B\'{e}zier, P., The Mathematical Basis of UNISURF CAD System. London, Butterworths, 1986. \\ 
\item	B\'{e}zier, P., Numerical Control; Mathematics and Applications (transl. by R. Forrest). New-York, John Wiley, l972. \\ 
\item	Bolyai, Farkas, Tentamen Juventutem Studiosam ein Elementa Mathiseos Parae  introducendi, Maros-Vasarhely, 1829, see Smith D. E. p. 375. \\ 
\item	Bolyai, Janos, Appendix of Bolyai, Farkas. \\ 
\item	Bolyai, Janos, The science absolute of space independent of the truth and falsity of Euclid's Axiom XI, translated by Dr. George Brus Halstead, Austin,Texas, The Neomon, Vol. 3, 71 pp, 1886. \\ 
\item	Borsuk, Karol and Szmielew Wanda, Foundations of Geometry, Amsterdam, North-Holland, 1960, 444 pp. \\ 
\item	Boubals, J. de Math. Elem. (de Longchamps et Bourget), 1891, p.218. (points of, on the circle of Brianchon-Poncelet) \\ 
\item	Brahmegupta and Bhascara, Algebra with Arithmetic and Mensuration, from the Sanscript of Brahmegupta and Bhascara, translated by Henry Thomas Cole-brooke, London, John Murray, 1817.\\ 
\item	Braikenridge, William, Exercitatio geometrica, London, 1733. \\ 
\item	Brianchon, Charles, Poncelet, Jean, Recherches sur la d\'{e}termination d'une hyperboie \'{e}quilatère au moyen de quatre conditions donn\'{e}es, Ann. de Math., Vol. 11, 1820-1821 , p. 205-220, see Smith, D. E., p. 337. \\ 
\item	Buchheim Arthur, An extension of Pascal's theorem to space of three dimensions. Messenger of Mathematics, Ser. 2, Vol. 14, 1984, 74-75. \\ 
\item	Casey John, Sequel to Euclid, London, 1881, p.101 for III.4.4.0. \\ 
\item	Charles, Michel, Apper\c cu historique sur l'origine et le developement des m\'{e}thodes en g\'{e}om\'{e}trie, 2e Edition, Paris, 1875. \\ 
\item	Ch'in, Chiu-Shao, see Libbrecht, Ulrich. \\ 
\item	Chou, Shang-Ching, Proving Elementary Geometry Theorems using Wu's Algorithm. Contemporary Mathematics, Bledsoe, W. W., and Loveland, D. W. Ed., Amer. Math. Soc. Vol. 29, 1984 243-286. \\ 
\item	Clebsch, Rudolf Frederich Alfred, Vorlesungen uber Geometrie, p.312 \\ 
\item	Coolidge, Julian, The Elements of non-Euclidean Geometry. Oxford Clarendon Press,1909, 291 pp. \\ 
\item	Coolidge, Julian, A treatise on the Circle and the Sphere, Oxford, Clarendon Press, 1916. 603 pp., III.4.4.0. \\ 
\item	Coolidge, Julian, A History of geometrical methods. Oxford Clarendon Press, 1940, 451 pp. \\ 
\item	Coolidge, Julian, A History of the Conic Sections and Quadric Surfaces. Oxford Clarendon Press, 1945, 214 pp. \\ 
\item	Coxeter, H. S. M., The Real projective Plane, New-York, McGraw-Hill, 1949, 196 pp. \\ 
\item	Coxeter H. S. M. and Greitzer, S. L., Geometry Revisited, N. Y. Random House, 1967, 193 pp. \\ 
\item	Coxeter H. S. M. and Moser, W. O. J., Generators and relations for discrete groups. Springer, 1957. \\ 
\item	Dalle A. et De Waele C., G\'{e}om\'{e}trie plane. Namur, Belgique, Wesmael-Charlier, 1936, 408 pp. \\ 
\item	David Antoine, 2000 Th\'{e}orèmes et Problèmes de G\'{e}om\'{e}trie avec Solutions. Namur, Belgique, Wesmael-Charlier, 1956, 1055 pp. \\ 
\item	de Casteljau, P., Outillages m\'{e}thodes calcul. Technical Report, Citroën, Paris 1959, See also 1963. \\ 
\item	de Casteljau, P., Shapes Mathematics and CAD. Kogan Page, London, 1986. \\ 
\item	Dembowski, Peter, Finite Geometries, Ergebnisse der Mathematik und ihrer Grenzgebiete, Band 44, Springer, New-York, 1968, 375 pp. \\ 
\item	Desargues, G\'{e}rard, Brouillon d'un projet d'une atteinte aux \'{e}vènements des rencontres d'un cône avec un plan, Paris, 1639. See Smith D. E., p. 307. \\ 
\item	Descartes, Ren\'{e}, La G\'{e}om\'{e}trie, Nouv. Ed., Paris Hermann, 1886, 91 pp.\\ 
\item	De Vogelaere, R., Finite Euclidean and non-Euclidean Geometry with application to the finite Pendulum and the polygonal harmonic motion. A first step to finite Cosmology. The Big Bang and Georges Lemaitre, Proc. Symp. in honor of 50 years after his initiation of Big-Bang Cosmology, Louvain-la-Neuve, Belgium, October 1983., D. Reidel Publ. Co, Leyden, the Netherlands. 341-355. \\ 
\item	De Vogelaere, R., G\'{e}om\'{e}trie Euclidienne finie. Le cas p premier impair. La Gazette des Sciences Math\'{e}matiques du Qu\'{e}bec, Vol. 10, Mai 1986. \\ 
\item	Dieudonn\'{e}, Jean, La g\'{e}om\'{e}trie des groupes classiques, Berlin, Springer, Ergebnisse Der Math. und ihrer Grenzgebiete, 1963, 125 pp. \\ 
\item	Donath, E. Die merkwurdigen Punkte und linien des Dreiecks, Berlin, VEB Deutscher Verlag der Wissenshaften, 1968. \\ 
\item	Emmerich A., Die Brocarsschen Gebilde, Berlin, Verlag Georg Reimer, 1891, \\ 
\item	Engle and Staeckel, Theorie der Parallellinien von Euklid bis auf Gauss, Leipzig, 1895. (See Mathesis, S\'{e}r. 2, Vol. 6, 1896, Suppl. pp. 1-11 or Rev. des quest. scient. S\'{e}r. 2, Vol. 8, 1895, pp. 603-612) \\ 
\item	Enriques, Frederigo, Lezioni di geometria proiettiva, Bologna, 1904, French Translation, Paris 1930 \\ 
\item	Euclides, Les oeuvres en grec, en latin et en fran\c cais, par Peyrard, Paris, Patris, 1814, 519 pp. \\ 
\item	Evans, Anthony B., On planes of prime order with translations and homologies, J. of Geometry, 34, 1989, 36-41. (Desarguesian planes) \\ 
\item	Eves, Howard, An Introduction to the History of Mathematics, New-York, Holt, Reinehart and Winston, 1953, 588 pp. \\ 
\item	Fano, Sui postulati fondamentali della geometria proiettiva, Giorn. di mat., Vol. 30, 1892, 106-132. (PG(n,p)) \\ 
\item	Farin, Gerald, Curves and Surfaces for Computer Aided and Geometric Design, New-York, Acad. Press, 1988, 334 pp.. \\ 
\item	Feuerbach, Karl, Grundriss zu analytischen Untersuchungen der dreyeckigen Pyramide, Nuremberg, 1827. \\ 
\item	Feuerbach, Karl, Eigenschaften einiger merkwurdigen Punkte des geradlinigen Dreiecks. N\"{u}rnberg, Riegel und Wiessner, 1822, 16+62 pp. \\ 
\item	Fonten\'{e}, G., Extension du Th\'{e}orème de Feuerbach. Nouv. Ann. de Math., S\'{e}r.4, Vol. 5, 1905. \\ 
\item	Forder, Henry, George, The Foundations of Euclidean Geometry, Cambridge Univ. P., 1927, repr. N. Y. Dover P., 1958, 349 pp. \\ 
\item	Forder, Henry, George, Higher Course Geometry, Cambridge Univ. P., 1931, 264 pp. \\ 
\item	Forder, Henry, George, The Calculus of Extension, New-York, Chelsea Pub. Co., 1960. \\ 
\item	Freudenthal, Oktaven, Ausnahmen gruppen und Oktavengeometgrie, Uttrecht, Utrecht Univ. 1960. \\ 
\item	Fritz, Kurt von, The discovery of incommensurability by Hippasus of Metapontum, Annals of Math., Vol. 46, 1945, 242-264. Also Studies in Presocratic Philosophy, Furley David and Allen R. E. Rdit. New-York, Humanities P., 1970, pp.382-412.\\ 
\item	Gauss, Carl, Disquisitiones Arithmeticae, Lipsiae, Fleicher, 1801. Translated by Clarke, Arthur, S.J., New Haven, Yale Univ. P. 1966, 473 pp.\\ 
\item	Gergonne, Joseph, Diaz, circle inscrit Nagel cerele exinscrit.\\ 
\item	Gergonne, Joseph, Diaz, Annales de Math\'{e}matiques, 1827, Vol. 17, p. 220 and 1829, Vol. 19, p. 97 and 129.\\ 
\item	Greenberg, M., Euclidean and Non-Euclidean Geometries, San Francisco, Freeman, 1974.\\ 
\item	Hagge, Der Fuhrrnannsche Kreis und der Brocardsche Kreis, Zeitschrift f\"{u}r rnathematische Urnterricht, vol. 38, 1907.\\ 
\item	Hall, Marshall, Projective Planes , Trans. Amer. Math. Soc., Vol. 54, 1943, 229-277\\ 
\item	Hall, Marshall, Jr, Projective Planes and related Topics, Calif. Inst. of Technology, April 1954, 77 pp.\\ 
\item	Hartshorne, Robin C., Foundation of Projective Geometry, N. Y. Benjamin, 1967, 161 pp.\\ 
\item	Heath, Sir Thomas, The thirteen books of Euclid's elements, Vol. 1, Cambridge University Press, 1908. 424 pp.
Other Edition, "The Classics of the St John's program," Annapolis, The St. John's College Press, 1947.\\ 
\item	Heath, Sir Thomas, Diophantus of Alexandria, 2nd Edition, Cambridge University Press, 1910.\\ 
\item	Heath, Sir Thomas, A Manual of Greek Mathematics, Oxford Univ. P., 1931, 551 pp.\\ 
\item	Heidel, W. A . The Pythagoreans and Greek Mathematics, Amer. J. of Philology, Vol. 61, 1940, 1-33. Also
Studies in Presocratic Philosophy, Furley David and Allen R. E. Edit. New York, Humanities P., 1970, pp. 350-381.\\ 
\item	Hensel, Kurt Theorie der algebraischen Zahlen, Berlin, Teubner, 1908, 349 pp.\\ 
\item	Hensel, Kurt, Zahlentheorie, Berlin, Goeschen'sche Verslaghandlung, 1913, 356 pp.\\ 
\item	Hessenberg G., Math. Ann., Vol. 61, 1905, pp. 161-172.\\ 
\item	Hilbert D., Grurdlagen der Geometrie, 1899, tr. by E. J. Townsend, La Salle, Ill., Open Court Publ. Cp., 1962, 143 pp.\\ 
\item	Hilbert D. und Cohn-Vossen S., Auschauliche Geometrie, Berlin, Springer, 1932, 310 pp.\\ 
\item	Hirschfeld, J. W. P., Projective geometries over finite fields, Oxford, Clarendon Press, 1979. 474 pp.\\ 
\item	Hughes, D. R., A class of non-Desargesian projective planes, Canad. J. of Math., 1957, VoI. 9, 378-388. (I.9.p.1) \\ 
\item	lntrigila, Carmelo, Sul Tretraedro, Rend. della R. Accad. delle Scienze di Napoli, Vol. 22,1883, pp. 69-92. \\ 
\item	Iversen, Birger, An Invitation to Geometry, Math. Inst. Aarhus Univ. Lect. Notes Series, No 59, 1989, 186 pp. \\ 
\item	Jacobi, Karl Gustav, Crelle J. fur Reine und Angewandte Mathernatik, Vol. 15, 199-204, Werke, I, p.336, (4)). \\ 
\item	J\"{a}rnefelt G. and Kustaanheimo Paul, An Observation on Finite Geometries. Den II. Skandinavische matematikerkongress, Trondheim, August 1949, 166-182. \\ 
\item	J\"{a}rnefelt G. , Reflections on a finite Approximation to Euclidean Geometry. Physical and Astronomical Prospects. Suomalaisen Tiedeakatemian Toimituksia, Ser. A, 1951, No 96. \\ 
\item	Johnson, Norman L., Kallaher, Michael J., Long Calvin T., Edit. Finite Geometries, N. Y., Marcel Dekker Inc. 1983. \\ 
\item	Johnson, Roger A., Modern Geometry, Houghton Mifflin Co, 1929, 319 pp. \\ 
\item	Karteszi, F. Introduction to finite geometries. Amsterdam, North Hollard Publ. Co., 1976, 266 pp. \\ 
\item	Kirk, G. S., Popper on Science and the Presocratics. Mind, Vol 69, 1960, 318- 339. 
Also Studies in Presocratic Philosophy, Furley David and Allen P. E. Edit. New-York, Humanities P., 1970, pp.154-177. \\ 
\item	Klein Felix, Geometry, N.Y. MacMillan, 1939. \\ 
\item	Klein, Felix, Famous Problems of Elementary Geometry, tr. by W.W. Beman and D. E. Smith, 1897, N.Y. Ginn, Reprinted in Famous Problems and other Monographs. New York, Chelsea, 1955. \\ 
\item	Kn\"{u}ppel, Frieder and Salow, Edzard, Plane elliptic geometry over rings. Pacific Journal of Mathematics. Vol. 123, (1986), 337-384. \\ 
\item	Koblitz Neal, A Short Course on Some Current Research in p-adic Analysis (Talks at Hanoi Math. Inst., July 1978), 66 pp. (Prof. Ogus) \\ 
\item	Koblitz Neal, p-adic Numbers, p-adic Analysis and Zeta-Functions, Springer- Verlag, N.Y., 1977, 124 pp. \\ 
\item	Lachlan, On Poristic Systems of Circles, Messenger of Mathematics, vol. 16, 1887.\\ 
\item	Laguerre, Edmond Nicolas, Oeuvres, 2 Vol. Paris, Gauthier-VilIars, 1898-1905. (I,9.,p1) \\ 
\item	Lebesgue, Henri, Le\c cons sur les constructions g\'{e}om\'{e}tliques, Paris, Gauthier- Villars, 1950, 304 pp. \\ 
\item	Lehmer, D. H., An elementary course in synthetic projective geometry. Boston, 1917 and Berkeley, California, Univ. of Calif. Pr., 1933, 123 pp. \\ 
\item	Lema\^{i}tre Georges. l'Hygrothèse de l'Atome Primitif, Essai de Cosmogonie, Neucgatel, Ed. du Griffon, 1946, 201 pp.\\ 
\item	Lema\^{i}tre Georges. The Primeval Atom, A Hypothesis of the Origin of the Universe. Trans. by Betty and Serge Korff, van Nostran, N. Y., 1950, 186 pp. \\ 
\item	Lemay, Fernand, Imagination dissidente, Bull. de I'APAME mars 1979. \\ 
\item	Lemay, Fernand, Motivation intrinsèque, Bull. de I'APAME, novembre 1979. \\ 
\item	Lemay, Fernand, Le dod\'{e}caèdre et la g\'{e}om\'{e}trie projective d'ordre 5, see Johnson N. L., 279-306 \\ 
\item	Lemoine, Emile, Propri\'{e}t\'{e}s relatives a deux points du plan d'un triangle qui se d\'{e}duisent d'un point K quelconque du plan comme les points de Brocard se d\'{e}duisent du point de Lemoine. Mathesis, S\'{e}r.1, Vol. 6, 1886, Suppl. 1-27. \\ 
\item	Lemoine, Emile, 1902, G\'{e}om\'{e}trogtaphie, C. Naud, Paris \\ 
\item	Lemoine, Emile, J. de Math. Elem. (de Longchamps et Bourget), 1889, p.93 1890, p.118, (point of, on the circle of Brianchon-Poncelet) \\ 
\item	Libbrecht, Ulrich, Chinese Mathematics in the Thirteenth Century, The Shu- Shu-Chiui-Chang of Ch'in, Chiu-Shao, Cambridge, MIT Press, 1973, 555 pages. \\ 
\item	Lobachevskii, Nikolai, Ivanovich, see Nolden. A., Elementare Einfulrung in die Lobachewskische Geometrie, Berlin, VEB Deuscher Verlag der Wissenschaften, 1958, 259 pp. \\ 
\item	MacLaurin, Colin,Phil. Trans. Roy. Soc. London, 1735. (On Pascal Constr) \\ 
\item	Mansion, Paul, Premiers Principes de !a M\'{e}tag\'{e}om\'{e}trie ou G\'{e}orn\'{e}trie g\'{e}n\'{e}rale, Mathesis, Ser. 2, Vol. 6, 1896, Suppl. 1-46. \\ 
\item	Mascheroni L. G\'{e}om\'{e}trie du Compas, translated in French Carette A.M., Paris 1798\\ 
\item	Maxwell, E. A., The methods of plane projective geometry based on the use of general homogeneous coordinates. Cambridge Univ. Press, 1952. 230 pp. \\ 
\item	Menger, K., Untersunchunten uber allgemeine Metrik, Math. Ann. Vol. 100, 1928, 75-163. \\ 
\item	Michel, Charles, Compl\'{e}ments de g\'{e}om\'{e}trie moderne, Paris, Vuibert, 1926. harmonic polygons, p. 272. \\ 
\item	Michel, Paul-Henri, De Pythagorea Euclide, Paris, Les Belles Lettres, 1950, 699 pp. \\ 
\item	Miquel, Auguste, Th\'{e}orèmes de g\'{e}orn\'{e}trie, J. de Liouvilie, Vol. 3, 1838, p.486. \\ 
\item	Miquel, Auguste, M\'{e}moires de G\'{e}orm\'{e}trie, J. de Math\'{e}matiques Pures et Appliqu\'{e}es, (J. de Liouville), Vol. 9, 1844, p.24. \\ 
\item	Möbius August, Werke, (Calcul Barycentrique) \\ 
\item	Moise, Edwin E., Elementary Geometry from an advanced standpoint, Palo Alto, Addison-Wesley, 1963, 419 pp. \\ 
\item	Moufang, Ruth, Alternatievkörper und der Satz vom vellständigen Vierseit, Abh. Math. Sem. Hamburg, Vol. 9, 1933, 207-222. \\ 
\item	Moulton, F. R., A simple non-Desarguesian plane Geometry, Trans. Amer. Math. Soc. Vol. 3, 1902, 192-195. \\ 
\item	 Nagel, Chretien Henty, Untersuchungen uber die wichtigsten zum Dreicke gehoerige Kreise, 1836. \\ 
\item	 Neuberg, Joeeph, M\'{e}moire sur le t\'{e}traèdre, M\'{e}moires couronn\'{e}s de l'Acad\'{e}mie de Belgique, Vol. .37, 1886, pp. 3-72. \\ 
\item	O'Hara, C. W. and Ward, D. R., An introduction to Projective Geometry, London, Oxford Uni, P., 1937, 298 pp. \\ 
\item	Ostrom, T, G., Finite translation planes, Lecture Notes in Math., Number 158, Berlin, Springer, 1970. \\ 
\item	Ostrom, T. G., Some translation planes that are not well known, Technical Report, N. 13, Department of Math. Washington State Univ., 1968, 49 pp. \\ 
\item	Pascal, Blaise, Pens\'{e}es, Nouv. Ed., Philippe Sellier, 1976, 543 pp. \\ 
\item	Pascal, Blaise, Essay sur les couiques, 1639, see Smith, D. E., p, 326. \\ 
\item	Pascal, Blaise, Oeuvres, Ed. Brunschvig et Boutroux, I, p. 252, (on Pascal, Constr.) \\ 
\item	Pasch, Moritz, Vorlesungen uber neuere Geometrie, Leipzig, Teubner, 1882, 202 pp. \\ 
\item	Pasch, Moritz und Dehn, Max, Vorlesungen uber neuere Geouretrie, Berlin, \\ 
\item	Pickert. G., Projektive Ebenen, Berlin, Springer, 1955, 343 pp. \\ 
\item	Pieri, Un sistema di postulati per la geometria proieitiva, Rev. Math\'{e}m. Torino, Vol 6, 1896. See also Atti Torino, 1904, 1906. \\ 
\item	Pieri, I principii della geornetria di posizione, composti in sistema logico deduttivo, Mem. della Reale Acad. delle Scienze di Torino, serie 2, Vol.48, 1899, pp 1-62. \\ 
\item	Playfair, John, Elements of Geometry, Philadelphia, Lippincoot \& Co, 1864, 318 pp \\ 
\item	Pl\"{u}cker, Julius, Analitische geometrie Entwickelungen, Voi 1 and 2 1828-1831 Crelle, Vol.5, 1830, Vol.12 (1834). Springer, 1976, 275 pp. \\ 
\item	Pl\"{u}cker, Julius, Theorie der algebraischen Curven 1839\\ 
\item	Poncelet, Jean, Victor, Application d'Analyse et de G\'{e}om\'{e}trie, Paris, Mallet- Bachelier, I, 1862, 563 pp., II, 1864, 602 pp. \\ 
\item	Poncelet, Jean. Victor, Trait\'{e} des propri\'{e}t\'{e}e projectives des figures, Paris, Gauthier- Villars, I, 1865, 428+xii pp. II, 1866, 452+vi pp. \\ 
\item	Popper, Sir Karl, Back to the Presocratics, Proc. of the Aristotelian Society, Vol. 59, 1958-9, 1-24, also
Studies in Presocratic Philosophy, Furley David and Allen R. E. Edit. New York, Humanities P., 1970, pp.130-153. \\ 
\item	Prouhet, , Analogies du triangle et du t\'{e}traèdre, Nouv. Ann. de Math., S\'{e}rie 2, Vol. 2, 1863, p.138. \\ 
\item	Reidemeister, K,, Grundlagen der Geometrie, Berlin, Springer, Grundl, der math,Wissens in Einz., Vol. 32, 1968, \\ 
\item	Robert, Alain; Elliptic curves, Lecture Notes in Mathematics, Berlin, Springer, Vol. 326, 1973, 264 pp. \\ 
\item	Roberts, Michael, On the analogues of the Nine-Point Circle in the Space of Three Dimensions, Proc. London Math. Soc., Vol. 19, 1878.  \\ 
\item	Roberts, Samuel, Proc. London Math Soc., Vol 12, 117. (Generalization of Miquel to tetrahedron) \\ 
\item	Robinson, A., Non-standard Analysis, North-Holland, Amsterdam, 1974, 293 pp.\\ 
\item	Robinson, G. de B., The foundations of Geometry, Toronto, 1940. \\ 
\item	Saccheri, Giovanni Girolamo, Euclides ab omni naevo vindicatus Milan, 1732. Tr. George Halstead, London Open Court Pr. 1920, 246 pp. See Engel and St\"{a}ckel.  \\ 
\item	Salmon, George, A treatise on Conic. Sections, 6-th ed. London 1879. \\ 
\item	Salmon, George, A treatise on the higher plane curves, 3d ed. Dublin, Hodges, Foster and Figgis, 1879, 395 pp. \\ 
\item	Schwabhäuser W., Szmielew W., Tarski.A, Metamathdmatische Methoden in der Geometrie, N. Y, Springer, 1980, 482 pp. \\ 
\item	Segre, B, Lectures on modern Geometry. Rome, Cremonese, 1961, 479 pp. \\ 
\item	Segre, C, Un,nuovo campo di ricerche geometriche, Atti R. Acad. Sc. Torino, Vol 25, 1889, 430-457. \\ 
\item	Seidenberg, Lectures in Projective Geometry, Princeton N. J., van Nostrand, 1962, 230 pp. \\ 
\item	Shively, Levi S., An Introduction to Modern Geometry, N. Y., John Wiley, 1939. 167 pp. \\ 
\item	Smith, David, Eugene, History of Mathematics, Vol. I, II. \\ 
\item	Smith, David, Eugene, A Source Book of Mathematics, N. Y. McGraw Hill, 1929,701 pp. \\ 
\item	Smogorzhevskii, A. S., The ruler in Geometrical Constructions, tr. by Halina Moss, New York, Blaisdell, 1961.\\ 
\item	Sommerville, Duncan, Bibliography of non-Euclidean Geometry, London, Harrison, 1911, 403 pp. \\ 
\item	Spieker, Ein merkwurdiger Kreis um der Schwerpunkt des Perimeters des geradlinigen Dreiecks als Analogon des Kreises der neun. Punkte, Grunert's Archiv, Vol. 51,1870. \\ 
\item	Staudt, K. G. C. von, Geometrie der Lage, Nuremberg 1847 \\ 
\item	Staudt, K. G. C. von, Beitrage zur Geometrie der Lage, Nuremberg 1857. \\ 
\item	Steiner, Jacob, Geometrical Constructions with a Ruler, tr. by M. E. Stark, ed. by R C. Archibald, New York Scripta Mathematica, 1950. \\ 
\item	Steiner, Jakob, Collected. Works; Vol. I,. pp. 43 and 135 for III.4.4.0. \\ 
\item	Stevenson, F.W., Projective Planes, W.E. Freeman and Co, 1972, 416 pp. \\ 
\item	Stroeker, R. J., Brocard Points, Circulant Matrices, and Descates' Folium; Math. Magazine, Vol. 61,1988, 172-187 \\ 
\item	Tarski, Alfred, What is Elementary Geometry, The axiomatic method with special reference to Geometry and Physics, Studies in Logic and the Foundation of Mathematics, North-Holland, Amsterdam, 1959, 16-29, Collected. Works, IV, 17-32. \\ 
\item	Taurinus, Theorie der Pallellinien. 1825, 102pp. \\ 
\item	Taurinus, Geometriae primia Elementa, 1826, 76pp. \\ 
\item	Taylor, H. M., On a six point circle connected with a triangle, Messenger of Mathematics, Vol 11, 177-179. (Circle of Taylor). \\ 
\item	Taylor, H. M., The Porism of the ring of circles touching two circles, Messenger of Mathematics, Vol. 7, 1878. III.4.4.0. \\ 
\item	Taylor, W. W., On the ring of circles touching two circles, Messenger of Mathematics, Vol. 7, 1878. III.4.4.0. \\ 
\item	Terquem, Orly, Consideration sur le triangle rectiligne, Nouv. Ann. de Math., Serie 1, Vol. 1, 1842, 196-200 \\ 
\item	Thomas, Ivor, Greek Mathematics, Cambridge, Mass., Harvard Univ. P., Vol. 1. 1939, 505 pp. \\ 
\item	Thureau-Dangin, F. Textes Mathematiques Babyloniens, Leiden 1938. \\ 
\item	Tilly, Joseph Marie de, Essai sur les Principes fondamentaux de G\'{e}om\'{e}trie et de Mecanique, Bruxelles, Mayolez, 1879 192 pp. Also, Mem. Soc. scinc. phys. et natur. de Bordeaux, Vol III Ser. 2, cahier l. \\ 
\item	Tilly, Joseph Marie de, Essai de G\'{e}om\'{e}trie analytique g\'{e}n\'{e}rale, Bruxelles 1892. Blumenthal considers than in this paper Tilly makes a fundamental contribution by introducing n-point relations to characterize a space metrically. \\ 
\item	Tucker, R., The "cosine" orthocenters of a triangle and a cubic through them Messenger of Mathematics, Ser 2, Vol. 17, pp. 97-103. (10 distances between 5 points!) \\ 
\item	Vahlen, Ueber Steinersche Kugelketten, Zeitschrift f\"{u}r Mathematik und Physik, Vol. 41, 1896, III.4.4.0. \\ 
\item	van der Waerden, Mathematics and Astronomy in Mesopotamia, Dict. of Scientific Bibliography, Vol 15.\\ 
\item	Veblen, Oswald, and Bussey, W. H., Trans. Amer. Math. Soc., Vol. 7, 1906, 241-259. (PG(n,pk)) \\ 
\item	Veblen, Oswald, and Young, John, Projective Geometry, Wesley, Boston, I, 1910, II, 1918 \\ 
\item	Ver Eecke, Paul, Proclus de Lycie, Les commentaires sur le premier livre des \'{e}l\'{e}ments d'Euclide, Bruges, Desclee De Brouwer, 1948, 372 pp. \\ 
\item	Verriest, Gustave, El\'{e}ments de G\'{e}om\'{e}trie Projective, Louvain, Feyaerts, 1930, 412 pp.\\ 
\item	Vigari\'{e}, Emile, Premier inventaire de la g\'{e}om\'{e}trie du triangle, Mathesis, S\'{e}r. 1, Vol. 9, 1889, Suppl. pp. 1-26. \\ 
\item	Vigari\'{e}, Emile, La bibliographie de Ia g\'{e}om\'{e}trie du triangle, Mathesis, S\'{e}r. 2, VoI. 6, 1896, Suppl. 1-14. (603 articles) \\ 
\item	Vuibert, Sur la G\'{e}om\'{e}trie classique du Triangle. \\ 
\item	Walker, R., Cartesian and Projective Geometry., London, Edward Arnold and Co, 1953, 320 pp. \\ 
\item	Whitehead, Alfred, The Axioms of Projective Geometry, Cambridge, 1906. \\ 
\item	Wu, Wen-Tsun, On the Decision Problem and the Mechanization of Theorem-Proving in Elementary Geometry, Contemporary Mathematics, Bledsoe, W. W., and Loveland, D. W. Ed., Amer. Math, Soc. VoI. 29, 1984 213-234. \\ 
\item	Wu, Wen-Tsun, Some Recent Advances in Mechanical Theorem-Proving of Geometries. Contemporary Mathematics, Bledsoe, W. W., and Loveland, D.W. Ed., Amer. Math. Soc. Vol. 29, 1984 215-242. \\ 
\item	Young, John, Projective geometry, 4th Carus Monograph, Chicago, 1930. \\
\end{enumerate}



\end{document}